\documentclass[12pt]{amsart}
\usepackage{amssymb}
\usepackage{graphics}  

\def\SENOCNumNarrowMatrices{416}
\def\SENOCNumNarrowMatricesRejectedRankTwo{1}

\def\SEOCNumNarrowMatrices{3030}
\def\SEOCNumNarrowMatricesRejectedRankTwo{8}
\def\SEOCNumNarrowMatricesRejectedPositiveDefinite{5}

\def\SPNOCNumNarrowMatrices{9150}

\def\SPNOCNumNarrowMatricesRejectedImprimitiveRoot{4153}

\def\SPOCNumNarrowMatrices{168075}

\def\SPOCNumNarrowMatricesRejectedImprimitiveRoot{84357}

\def\CPNumNarrowMatrices{137235}

\def\CPNumNarrowMatricesRejectedImprimitiveRoot{101189}

\def\TotalNumEnlargements{512207}
\def\TotalNumNarrowMatricesRemaining{128193}
\def\NumCandidates{214928}
\def\NumCandidatesFailRootSimplicity{10408}
\def\NumRemainingCandidates{204520}
\def\NumCandidatesSSF{2497}
\def\NumCandidatesDistinctSSF{857}

\def\above#1#2#3#4{\genfrac{}{}{0pt}{}{\lower#1\hbox{\smash{$\scriptstyle #3$}}}{\raise#2\hbox{\smash{$\scriptstyle #4$}}}}

\newdimen{\twowidth}
\newdimen{\halftwo}
\def\slashtwo{2\llap{\hbox to 0pt{\hss$|$\hss}\kern\halftwo}}
\newdimen{\threewidth}
\newdimen{\halfthree}
\def\slashthree{3\llap{\hbox to 0pt{\hss$|$\hss}\kern\halfthree}}
\newdimen{\inftywidth}
\newdimen{\halfinfty}
\def\slashinfty{\infty\llap{\hbox to 0pt{\hss$|$\hss}\kern\halfinfty}}

\def\BreakWithGlue#1#2#3{\nolinebreak\hskip #1\hbox{}\penalty #2\hbox{}\nolinebreak\hskip #3}
\def\EasyButWeakLineBreak{\BreakWithGlue{0pt plus 1 fill}{500}{0pt plus 1 fill}}
\def\HardButStrongLineBreak{\BreakWithGlue{0pt plus 1 filll}{600}{0pt plus -1 filll}}
\newcommand\spacer{\ }
\def\instructions#1{$\langle #1\rangle$}

\def\pfillwitharg#1{\mathop{\hbox{$#1$-fill}}\nolimits}
\def\main{\mathop{\rm main}\nolimits}
\def\iso{\cong}

\long\def\GramMatrix#1{#1} 
\long\def\GramMatrix#1{}   
\def\NumLattices{8595}
\def\NumPrintedExplicitly{704}
\def\NumPrintedImplicitly{715}

\def\NumMainLattices{122}

\def\NumWeylGroups{374}
\def\LargestArea{$12\pi$}
\def\SmallestArea{$\pi/12$}
\def\WeylsWithLargestArea{$W_{329}$, $W_{365}$ and $W_{371}$}
\def\WeylsWithSmallestArea{$W_{3}$}
\def\MostLatticesInWeyl{176}
\def\LeastLatticesInWeyl{1}
\def\WeylsWithMostLattices{$W_{68}$}
\def\WeylsWithLeastLattices{$W_{286}$}
\def\MostEdges{28}
\def\WeylsWithMostEdges{$W_{365}$}
\def\LargestCyclic{3}
\def\WeylsWithLargestCyclic{$W_{174}$}
\def\WeylsWithDihedralEight{$W_{182}$, $W_{184}$, $W_{238}$, $W_{285}$, $W_{286}$, $W_{307}$, $W_{324}$, $W_{349}$, $W_{359}$ and $W_{365}$}
\def\WeylsWithDihedralTwelve{$W_{234}$, $W_{320}$ and $W_{350}$}
\def\NumRegularWeyl{2}
\def\RegularWeyls{$W_{234}$ and $W_{286}$}
\def\NumIdealWeyl{3}
\def\IdealWeyls{$W_{247}$, $W_{285}$ and $W_{286}$}
\def\NumRightAngledWeyl{163}

\def\NumCompactChambers{269}

\def\NumCommClasses{39}
\def\NumGenera{8488}
\def\MostReflectivesInGenus{2}
\def\NumGeneraWithMostLattices{107}

\def\GenusWithTwoWeyls{{$1\above{1pt}{1pt}{1}{7}8\above{1pt}{1pt}{1}{1}256\above{1pt}{1pt}{1}{1}$}}
\def\GenusWithTwoWeylsFirstLattice{{L_{341.1}}}
\def\GenusWithTwoWeylsSecondLattice{{L_{342.1}}}

\def\GenusWithTwoWeylsSecondCornerSymbols{{${256}\above{1pt}{1pt}{r}{2}{4}\above{1pt}{1pt}{b}{2}{256}\above{1pt}{1pt}{l}{2}{8}\above{1pt}{1pt}{32,17}{\infty}{32}\above{1pt}{1pt}{32,25}{\infty z}{8}\above{1pt}{1pt}{}{2}$}%
}
\def\GenusWithTwoWeylsFirstWeyl{$W_{341}$}
\def\GenusWithTwoWeylsSecondWeyl{$W_{342}$}
\def\MostNegativeDet{$-29811600$}
\def\MostNegativeDetOfMinimalInDualityClass{$-108000$}
\def\LargestPrimeInDet{97}
\def\LargestRootNorm{14400}
\def\LargestRootNormofMinimalInDualityClass{7200}
\def\NumDualityClasses{1419}

\def\ImplicitGeneraTable{%
\begingroup
\parindent=0pt
\parskip=0pt\small
\raggedbottom
\leavevmode%
$\pfillwitharg{2}({L_{10.1}})$%
\ shares genus with 5-dual\nopagebreak\par%
\nopagebreak\par%
{$1\above{1pt}{1pt}{-2}{{\rm II}}2\above{1pt}{1pt}{1}{1}{\cdot}1\above{1pt}{1pt}{-}{}5\above{1pt}{1pt}{-}{}25\above{1pt}{1pt}{-}{}$}%
\hbox{}\par\smallskip%
\leavevmode%
$\main({L_{248.2}})$%
\ shares genus with {$ {L_{248.1}}$}%
\nopagebreak\par%
\nopagebreak\par%
{$[1\above{1pt}{1pt}{1}{}2\above{1pt}{1pt}{1}{}]\above{1pt}{1pt}{}{0}128\above{1pt}{1pt}{1}{1}$}%
\hbox{}\par\smallskip%
\leavevmode%
$\pfillwitharg{2}({L_{194.1}})$%
\ shares genus with 5-dual\nopagebreak\par%
\nopagebreak\par%
{$1\above{1pt}{1pt}{-2}{{\rm II}}2\above{1pt}{1pt}{1}{1}{\cdot}1\above{1pt}{1pt}{2}{}9\above{1pt}{1pt}{-}{}{\cdot}1\above{1pt}{1pt}{-}{}5\above{1pt}{1pt}{-}{}25\above{1pt}{1pt}{-}{}$}%
\hbox{}\par\smallskip%
\leavevmode%
$\pfillwitharg{2}({L_{195.1}})$%
\ shares genus with 5-dual\nopagebreak\par%
\nopagebreak\par%
{$1\above{1pt}{1pt}{-2}{{\rm II}}2\above{1pt}{1pt}{1}{1}{\cdot}1\above{1pt}{1pt}{-2}{}9\above{1pt}{1pt}{1}{}{\cdot}1\above{1pt}{1pt}{-}{}5\above{1pt}{1pt}{-}{}25\above{1pt}{1pt}{-}{}$}%
\hbox{}\par\smallskip%
\leavevmode%
$\pfillwitharg{2}({L_{73.1}})$%
\ shares genus with 5-dual\nopagebreak\par%
\nopagebreak\par%
{$1\above{1pt}{1pt}{2}{{\rm II}}2\above{1pt}{1pt}{1}{7}{\cdot}1\above{1pt}{1pt}{-}{}5\above{1pt}{1pt}{-}{}25\above{1pt}{1pt}{-}{}{\cdot}1\above{1pt}{1pt}{-2}{}11\above{1pt}{1pt}{-}{}$}%
\hbox{}\par\smallskip%
\leavevmode%
$\pfillwitharg{2}({L_{366.1}})$%
\ shares genus with 13-dual\nopagebreak\par%
\nopagebreak\par%
{$1\above{1pt}{1pt}{-2}{{\rm II}}2\above{1pt}{1pt}{1}{1}{\cdot}1\above{1pt}{1pt}{-}{}13\above{1pt}{1pt}{-}{}169\above{1pt}{1pt}{-}{}$}%
\hbox{}\par\smallskip%
\leavevmode%
$\pfillwitharg{2}({L_{345.1}})$%
\ shares genus with 3-dual; isometric to own %
7-dual\nopagebreak\par%
\nopagebreak\par%
{$1\above{1pt}{1pt}{-2}{{\rm II}}2\above{1pt}{1pt}{1}{1}{\cdot}1\above{1pt}{1pt}{-}{}3\above{1pt}{1pt}{1}{}9\above{1pt}{1pt}{-}{}{\cdot}1\above{1pt}{1pt}{1}{}7\above{1pt}{1pt}{1}{}49\above{1pt}{1pt}{1}{}$}%
\hbox{}\par\smallskip%
\endgroup
}%

\allowdisplaybreaks

\def\rh{{\rm rh}} 
\def\pdual{\mathop{\hbox{\rm $p$-dual}}\nolimits}
\def\pfill{\mathop{\hbox{\rm $p$-fill}}\nolimits}
\def\diag{\mathop{\hbox{diag}}\nolimits}
\def\II{{\rm II}}
\def\I{{\rm I}}

\def\Z{\mathbb{Z}} 
\def\Q{\mathbb{Q}} 
\def\R{\mathbb{R}} 
\def\rp{\R P} 
\def\D{\Delta} 
\def\O{{\rm O}} 

\def\corner#1#2{^{#1}_{#2}}
\def\gluez{z}
\def\gluea{a}
\def\glueboth{b}

\def\a{\alpha}
\def\b{\beta}
\def\c{\gamma}
\def\e{\varepsilon}
\def\ahat{\hat{a}}
\def\bhat{\hat{b}}
\def\rhat{\hat{r}}
\def\shat{\hat{s}}
\def\that{\hat{t}}
\def\Ctilde{\tilde{C}}
\def\DD{\mathcal{D}}
\def\FF{\mathcal{F}}

\def\del{\partial} 
\def\tensor{\otimes}
\let\iso=\cong
\let\cong=\equiv
\def\sset{\subseteq}
\def\sat{{\rm\scriptstyle sat}} 
\def\aut{\mathop{\rm Aut}\nolimits}
\def\notperp{\not\perp}
\def\legendre#1#2{\left(\frac{#1}{#2}\right)}

\def\set#1#2{\{#1\,:\,#2\}}

\def\spanof#1{\langle#1\rangle}

\def\centeroverfull#1{\setbox0=\hbox{#1}\newdimen\foo\foo=\wd0\advance\foo by -\hsize\divide\foo by 2\advance\foo by\hsize\rlap{\kern\foo\llap{\box0}}}

\newcommand\breakok{\discretionary{}{}{}}

\def\Xdefinition#1#2#3{#1\marginpar{\rm\scriptsize #2}\index{#3}}
\def\Xdefinition#1#2#3{#1} 

\def\defn#1{\Xdefinition{#1}{#1}{#1}}
\def\bfdefn#1{\marginpar{\rm\scriptsize #1}{\bf #1}}
\def\bfdefn#1{{\bf #1}}
\def\notation#1{\marginpar{#1}#1}
\def\notation#1{#1} 

\newtheorem{theorem}{Theorem}
\newtheorem{thmdef}[theorem]{Theorem (and Definition)}
\newtheorem{lemma}[theorem]{Lemma}
\theoremstyle{remark}

\newtheorem*{remark}{Remark}
\newtheorem*{remarks}{Remarks}

\begin{document}

\title[Reflective Lorentzian lattices]{The reflective Lorentzian lattices of rank 3}
\author{Daniel Allcock}
\address{Department of Mathematics\\U.T. Austin}
\email{allcock@math.utexas.edu}
\urladdr{http://www.math.utexas.edu/\textasciitilde allcock}
\thanks{Partly supported by NSF grant DMS-0600112.}
\date{November 19, 2010}

\begin{abstract}
We classify all the symmetric integer bilinear forms of signature
$(2,1)$ whose isometry groups are generated up to finite index by
reflections.  There \NumLattices{} of them up to scale, whose
\NumWeylGroups{} distinct Weyl groups fall into \NumCommClasses{}
commensurability classes.  This extends Nikulin's enumeration of the
strongly square-free cases.  Our technique is an analysis of the shape
of the Weyl chamber, followed by computer work using Vinberg's
algorithm and our ``method of bijections''.  We also correct a minor
error in Conway and Sloane's definition of their canonical $2$-adic
symbol.
\end{abstract}

\maketitle

\section{Introduction}
\label{sec-intro}

Lorentzian lattices, that is, integral symmetric bilinear forms of
signature $(n,1)$, play a major role in K3 surface theory and the
structure theory of hyperbolic Kac-Moody algebras.  In both cases, the
lattices which are reflective, meaning that their isometry groups are
generated by reflections up to finite index, play a special role.  In
the KM case they provide candidates for root lattices of KMA's with
hyperbolic Weyl groups that are large enough to be interesting.  For
K3's they are important because a K3 surface has finite automorphism
group if and only if its Picard group is ``2-reflective'', meaning
that its isometry group is generated up to finite index by reflections
in classes with self-intersection $-2$
(see \cite{Piatetski-Shapiro-Shafarevich}).  If the Picard group is
not 2-reflective but is still reflective, then one can often describe
the automorphism group of the surface very
explicitly \cite{Vinberg-most-algebraic}\cite{Borcherds-Coxeter-groups-Lorentzian-lattices-and-K3s}\cite{Nikulin-K3s-with-interesting-symmetry-groups}.
The problem of classifying all reflective lattices of given rank, here
$n+1$, is also of interest in its own right from the perspectives of
Coxeter groups and arithmetic subgroups of $\O(n,1)$.  

In this paper we give the classification of reflective Lorentzian
lattices of rank~$3$.  There are \NumLattices{} of them up to scale,
but this large number makes the result seem more complicated than it
is.  For example,
\MostLatticesInWeyl{} lattices
share a single Weyl group (\WeylsWithMostLattices{} in the table).

Our work is logically independent of Nikulin's monumental paper
\cite{Nikulin-rk3}, but follows it conceptually.  His work had
the same motivations as ours, and solves a problem both more general
and more specific.  He considered Lorentzian lattices that satisfy the
more general condition of being ``almost reflective'', meaning that
the Weyl chamber is small enough that its isometry group is $\Z^l$ by
finite.  This is a very natural class of lattices from both the K3 and
KM perspectives.  He also restricts attention to lattices that are
strongly square-free (SSF), which in rank~$3$ means they have
square-free determinant.  The full classification is desirable from
the KM perspective, because for example the root lattice even of a
finite-dimensional simple Lie algebra need not be SSF.  On the other
hand, every lattice canonically determines a SSF lattice, so that his
classification does yield all lattices whose reflection groups are
maximal under inclusion.  The SSF condition  also enabled Nikulin to
use some tools from the theory of quadratic forms that we avoid.
There are 1097 lattices in Nikulin's 160 SSF lattices' duality classes
(defined in section~\ref{sec-background}), while our
\NumLattices\ lattices fall into \NumDualityClasses\ duality classes.
Most of our techniques run parallel to or are refinements of Nikulin's
methods.

We have strived to keep the table of lattices to a manageable size,
while remaining complete in the sense that the full classification can
be mechanically recovered from the information we give.  So we have
printed only \NumPrintedExplicitly\ of them and given operations to
obtain the others.  Even after this reduction, the table is still
large.  The unabridged list of lattices is available separately
as \cite{Allcock-full-tables-rk3}, including the Gram matrices if
desired.  Anyone wanting to work seriously with the lattices will need
computer-readable data, which is very easy to extract from the TeX
source file of this paper.  Namely, by arcane computer tricks we have
arranged for the TeX source to be simultaneously a Perl script that
prints out all \NumLattices\ lattices (or optionally just some of
them) in computer-readable format.  If the file is saved as say {\tt
file.tex} then simply enter {\tt perl file.tex} at the unix command
line.

For each reflective lattice $L$, we worked out the conjugacy class
of its Weyl group $W(L)$ in $\O(2,1)$, the normalizers of $W(L)$ in
$\aut L$ and $\O(2,1)$, the area of its Weyl chamber, its genus as a
bilinear form over $\Z$, one construction of $L$ for each corner of
its chamber, and the relations among the lattices under $p$-duality
and $p$-filling (defined in section~\ref{sec-background}).  Some
summary information appears in section~\ref{sec-how-to-read}, for example
largest and smallest Weyl chambers, most and fewest lattices with
given Weyl group, most edges (\MostEdges), numbers of chambers with
interesting properties like regularity or compactness, genera
containing multiple reflective lattices, and the most-negative
determinants and largest root norms.

Here is a brief survey of related work.  Analogues of Nikulin's work
have been done for ranks $4$ and~$5$ by R.\ Scharlau \cite{Scharlau}
and C.\ Walhorn \cite{Walhorn} respectively; see also
\cite{Scharlau-Walhorn}.  In rank${}\geq21$, Esselmann
\cite{Esselmann} proved there are only $2$ reflective lattices, both
in rank~$22$---the even sublattice of $\Z^{21,1}$ (recognized as
reflective by Borcherds \cite{Borcherds-I21-1}) and its dual lattice.
Also, Nikulin \cite{Nikulin-finiteness} proved that there are only finitely many
reflective lattices in each rank, so that a classification is
possible.  The lattices of most interest for K3 surfaces are 
the 2-reflective ones, and these have been completely classified by
Nikulin \cite{Nikulin-2-reflective-rank-5+}\cite{Nikulin-2-reflective-rank-3}
and Vinberg \cite{Vinberg-2-reflective-rank-4}.

Beyond these results, the field is mostly a long list of examples found by various authors,
with only glimmers of a general theory, like Borcherds' suggested
almost-bijection between ``good'' reflection groups and ``good''
automorphic forms \cite[\S12]{Borcherds-almost-bijection}.  To choose a few
highlights, we mention Coxeter and Whitrow's analysis \cite{Coxeter-Whitrow} of
$\Z^{3,1}$, Vinberg's famous paper \cite{Vinberg} treating $\Z^{n,1}$ for
$n\leq17$, its sequel with Kaplinskaja \cite{Vinberg-Kaplinskaja} extending this to
$n\leq19$, Conway's analysis \cite{Conway-II25-1} of the even unimodular Lorentzian
lattice of rank~$26$, and Borcherds' method of deducing
reflectivity of a lattice from the existence of a suitable automorphic
form \cite{Borcherds-reflection-groups-of-Lorentzian-lattices}.  Before this, Vinberg's algorithm \cite{Vinberg} was the only
systematic way to prove a given lattice reflective.

Here is an overview of the paper.  After providing background material
in section~\ref{sec-background}, we consider the shape of a hyperbolic
polygon $P$ in section~\ref{sec-shape-of-polygon}.  The main result is
the existence of $3$, $4$ or $5$ consecutive edges such that the lines
in the hyperbolic plane $H^2$ containing them are ``not too far
apart''.  This is our version of Nikulin's method of narrow parts of
polyhedra \cite[\S4]{Nikulin-rk3}\cite{Nikulin-finiteness}, and some of the same
constants appear.  Our technique is based the bisectors of angles and edges of $P$, while his is based on
choosing a suitable point $p$ of $P$'s interior and considering how
the family of lines through $p$ meet $\del P$.  Our method is simpler
and improves certain bounds, but has no extension to
higher-dimensional hyperbolic spaces, while his method works in
arbitrary dimensions.  For a more detailed comparison, please see the
end of section~\ref{sec-shape-of-polygon}.

In section~\ref{sec-shape-of-2-dimensional-Weyl-chamber} we convert
these bounds on the shape of the Weyl chamber into an enumeration of
inner product matrices of $3$, $4$ or $5$ consecutive simple roots of
$W(L)$.  This is mostly a matter of organizing the parameters and
calculations, since the computations themselves are just arithmetic.
In section~\ref{sec-corner-symbols} we introduce ``corner symbols''.
Essentially, given $L$ and two consecutive simple roots $r$ and $s$,
we define a symbol like ${2}\above{1pt}{1pt}{*}{4}{4}$ or
${26}\above{1pt}{1pt}{4,3}{\infty b}{104}$ displaying their norms,
with the super- and subscripts giving enough information to recover
the inclusion $\spanof{r,s}\to L$.  So two lattices, each with a
chamber and a pair of consecutive simple roots, are isometric by an
isometry identifying these data if and only if the corner symbols and
lattice determinants coincide.  Some device of this sort is required,
even for sorting the reflective lattices into isometry classes,
because the genus of a lattice is not a strong enough invariant to
distinguish it.  In one case (see section~\ref{sec-how-to-read}), two
lattices with the same genus even have different Weyl groups.  We also
use corner symbols in the proof of the main theorem
(theorem~\ref{thm-main-thm} in section~\ref{sec-rank-3}) and for
tabulating our results.

Besides the proof of the main theorem, section~\ref{sec-rank-3} also
introduces what we call the ``method of bijections'', which is a
general algorithm for determining whether a finite-index subgroup of a
Coxeter group is generated by reflections up to finite index.
Section~\ref{sec-how-to-read} explains how to read the table at the
end of the paper.  That table displays \NumPrintedExplicitly{}
lattices explicitly, together with operations to apply to obtain the
full list of lattices.  Finally, in section~\ref{sec-Conway-Sloane} we
discuss the effects of various ``moves'' between lattices on the
Conway-Sloane genus symbol \cite[ch.~15]{SPLAG}.  We also resolve a
small error in their definition of the canonical $2$-adic symbol.

Computer calculations were important at every stage of this project,
as explained in detail in section~\ref{sec-rank-3}.  We wrote our programs in
C++, using the PARI software library \cite{PARI}, and the full
classification runs in around an hour on the author's laptop computer.  The source code is available
on request.  Although there is no a posteriori check, like the mass formula used when enumerating
unimodular lattices \cite{Conway-Sloane-unimodular-enumeration},  we can offer several remarks in favor of
the reliability of our results.  First, in the SSF case we recovered
Nikulin's results exactly; second, many steps of the calculation
provided proofs of their results (e.g., that a lattice is reflective
or not); and third, our list of lattices is closed under the
operations of $p$-duality and $p$-filling (see section~\ref{sec-background}) and
``mainification'' (section~\ref{sec-how-to-read}).  
This is a meaningful check because our enumeration methods have
nothing to do with these relations among lattices.

I am grateful to the Max Planck Institute in Bonn, where some of this
work was carried out.

\section{Background}
\label{sec-background}

For experts, the key points in this section are (i) lattices need not
be integral and (ii) the definition of ``reflective hull''.
The notation $\spanof{x,\dots,z}$ indicates the $\Z$-span of some
vectors.

{\it Lattices:\/} A lattice means a free abelian group $L$ equipped
with a $\Q$-valued symmetric bilinear form, called the inner product.
If $L$ has signature $(n,1)$ then we call it Lorentzian.  If the inner
product is $\Z$-valued then $L$ is called integral.  In this case we
call the gcd of $L$'s inner products the scale of $L$.  We call $L$
unscaled if it is integral and its scale is~$1$.  The norm $r^2$ of a
vector means its inner product with itself, and $r$ is called
primitive if $L$ contains none of the elements $r/n$, $n>1$.  $L$ is
called \defn{even} if all its vectors have even norm, and
\defn{odd} otherwise.  If $L$ is nondegenerate then its \defn{dual}, written \notation{$L^*$},
means the set of vectors in $L\tensor\Q$ having integer inner products
with all of $L$.  So $L$ is integral just if $L\sset L^*$.  In this
case $L$'s
\defn{discriminant group} is defined as \notation{$\D(L)$}$:=L^*/L$, a
finite abelian group.   The determinant $\det L$ is the determinant
of any inner product matrix, and when $L$ is integral we have $|\det
L|=|\D(L)|$. 

{\it Roots:\/} A root means a primitive lattice vector $r$ of positive
norm, such that reflection in $r$, given by 
$$
x\mapsto x-2\frac{x\cdot r}{r^2}r
$$ 
preserves $L$.  Under the primitivity hypothesis, this
condition on the reflection is
equivalent to $L\cdot r\sset\frac{1}{2}r^2\Z$.  Now suppose
$r_1,\dots,r_n$ are roots in a lattice $L$, whose span is
nondegenerate.  Then their reflective hull
$\spanof{r_1,\dots,r_n}^{\rh}$ means 
$$
\spanof{r_1,\dots,r_n}^\rh:=\set{v\in
\spanof{r_1,\dots,r_n}\tensor\Q}{\textstyle\hbox{$v\cdot r_i\in\frac{1}{2}r_i^2\Z$ for
    all $i$}}.
$$ The key property of the reflective hull is that it contains the
projection of $L$ to $\spanof{r_1,\dots,r_n}\tensor\Q$.  In the special
case that $r_1,\dots,r_n$ generate $L$ up to finite index, we have
$$
\spanof{r_1,\dots,r_n}\sset L\sset\spanof{r_1,\dots,r_n}^\rh.
$$
Then there are only finitely many possibilities for $L$, giving us
some control over a lattice when all we know is some of its roots.  As
far as I know, the reflective hull was first introduced
in \cite{Scharlau-positive-definite-reflective-lattices}, in terms of
the ``reduced discriminant group''.  For use in the proof of
theorem~\ref{thmdef-finite-corner-symbols}, we mention that if $r$ and
$s$ are simple roots for an $A_1^2$, $A_2$, $B_2$ or $G_2$ root
system, then $\spanof{r,s}^\rh/\spanof{r,s}$ is $(\Z/2)^2$, $\Z/3$,
$\Z/2$ or trivial, respectively.

{\it Weyl Group:\/} We write $W(L)$ for the Weyl group of $L$, meaning
the subgroup of $O(L)$ generated by reflections in roots.  When $L$ is
Lorentzian, $O(L)$ acts on hyperbolic space $H^n$, which is defined as
the image of the negative-norm vectors of $L\tensor\R$ in
$P(L\tensor\R)$.  For $r$ a root, we refer to $r^\perp\sset
L\tensor\R$ and the corresponding subset of $H^n$ as $r$'s mirror.
The mirrors form a locally finite hyperplane arrangement in $H^n$, and
a Weyl chamber (or just chamber) means the closure in $H^n$ of a
component of the complement of the arrangement.  Often one chooses a
chamber $C$ and calls it ``the'' Weyl chamber.  The general theory of
Coxeter groups (see e.g. \cite{Vinberg-general-theory}) implies that
$C$ is a fundamental domain for $W(L)$.  Note that $C$ may have
infinitely many sides, or even be all of $H^n$ (if $L$ has no roots).
A Lorentzian lattice $L$ is called reflective if $C$ has finite
volume.  The purpose of this paper is to classify all reflective
lattices of signature $(2,1)$.  For a discrete group $\Gamma$ of
isometries of $H^n$, not necessarily generated by reflections, its
reflection subgroup means the group generated by its reflections, and
the Weyl chamber of $\Gamma$ means the Weyl chamber of this subgroup.

{\it Simple roots:\/}
The preimage of $C$ in $\R^{2,1}-\{0\}$ has two components; choose one and
call it $\Ctilde$.  For every facet of $\Ctilde$, there is a unique
outward-pointing root $r$ orthogonal to that facet (outward-pointing means
that $r^\perp$ separates $r$ from $\Ctilde$).  These are called the
simple roots, and their pairwise inner products are${}\leq0$.
Changing the choice of component simply negates to simple roots, so
the choice of component will never matter to us.   These
considerations also apply to subgroups of $W(L)$ that are generated by
reflections.

{\it Vinberg's algorithm:\/}
Given a Lorentzian lattice $L$, a negative-norm vector $k$ (called the
controlling vector) and a set of simple roots for the stabilizer of
$k$ in $W(L)$, there is a unique extension of these roots to a set of
simple roots for $W(L)$.  Vinberg's algorithm \cite{Vinberg} finds this
extension, by iteratively adjoining batches of new simple roots.  If
$L$ is reflective then there are only finitely many simple roots, so
after some number of iterations the algorithm will have found them
all.  One can recognize when this happens because the polygon defined
by the so-far-found simple roots has finite area.  Then the algorithm
terminates.  On the other hand, if $C$ has infinite area then the
algorithm will still find all the simple roots, but must run forever
to get them.  In this case one must recognize $L$ as
non-reflective, which one does by looking for automorphisms of $L$
preserving $C$.  See section~\ref{sec-rank-3} for the details of the method we
used.  There is a variation on Vinberg's algorithm using a null
vector as controlling vector \cite{Conway-II25-1}, but the original algorithm was sufficient
for our purposes.

{\it $p$-duality:\/}
For $p$ a prime, the $p$-dual of a nondegenerate lattice $L$ is the
lattice in $L\tensor\Q$ characterized by
$$
\pdual(L)\tensor\Z_q=
\begin{cases}
L^*\tensor\Z_q&\hbox{if $q=p$}\\
L\tensor\Z_q&\hbox{if $q\neq p$}
\end{cases}
$$
for all primes $q$, where $\Z_q$ is the ring of $q$-adic integers.
Applying $p$-duality twice recovers $L$, so $L$ and its $p$-dual have
the same automorphism group.  When $L$ is integral one can be  much
more concrete: $\pdual(L)$ is the sublattice of $L^*$ corresponding to
the $p$-power part of $\D(L)$.  The effect of $p$-duality on the
Conway-Sloane genus symbol is explained in section~\ref{sec-Conway-Sloane}.

In our applications $L$ will be  unscaled, so the
$p$-power part of $\D(L)$ has the form $\oplus_{i=1}^n\Z/p^{a_i}$
where $n=\dim L$ and at least one of the $a_i$'s is $0$.  In this case
we define the ``rescaled $p$-dual'' of $L$ as $\pdual(L)$ with all
inner products multiplied by $p^a$, where $a=\max\{a_1,\dots,a_n\}$.
This lattice is also unscaled, and the $p$-power part of
its discriminant group is $\oplus_{i=1}^n\Z/p^{a-a_i}$.  Repeating the
operation recovers $L$.  We say that unscaled lattices are in
the same duality class if they are related by a chain of
rescaled $p$-dualities.  Our main interest in rescaled $p$-duality is that all
the lattices in a duality class have the same isometry group, so we
may usually replace any one of them by a member of the class with
smallest $|\det|$.

{\it $p$-filling:\/} Suppose $L$ is an integral lattice, $p$ is a
prime, and $\D(L)$ has some elements of order $p^2$.  Then we define
the $p$-filling of $L$ as the sublattice of $L^*$ corresponding to
$p^{n-1}A\sset\D(L)$, where $A$ is the $p$-power part of $\D(L)$ and
$p^n$ is the largest power of $p$ among the orders of elements of
$\D(L)$.  One can check that $\pfill(L)$ is integral, and obviously
its determinant is smaller in absolute value.  The operation is not
reversible, but we do have the inclusion $\aut L\sset\aut(\pfill(L))$.
As for $p$-duality, the effect of $p$-filling on the Conway-Sloane
genus symbol is explained in section~\ref{sec-Conway-Sloane}.

{\it SSF lattices:\/}
Given any integral lattice $L$, one can apply $p$-filling operations
until no more are possible, arriving at a lattice whose discriminant
group is a sum of cyclic groups of prime order.  Then applying
rescaled $p$-duality for some primes $p$ leaves a lattice whose discriminant
group has rank (i.e., cardinality of its smallest generating set) at
most $\frac{1}{2}\dim L$.  Such a lattice is called strongly square-free
(SSF).  For rank~$3$ lattices this is equivalent to $|\D(L)|$ being
square-free.  Nikulin's paper \cite{Nikulin-rk3} classified the SSF reflective
Lorentzian lattices of rank~$3$.
Watson \cite{Watson-I}\cite{Watson-II} introduced
operations leading from any lattice to a SSF one, and our operations
are presumably special cases of his.

\section{The shape of a hyperbolic polygon}
\label{sec-shape-of-polygon}

In this section we develop our version of Nikulin's method of narrow
parts of polyhedra \cite{Nikulin-finiteness}\cite{Nikulin-rk3}.  We will prove the following theorem,
asserting that any finite-sided finite-area polygon $P$ in $H^2$ has
one of three geometric properties, and then give numerical
consequences of them.  These consequences have the general form: $P$
has several consecutive edges, such that the distances between the lines
containing them are less than some a priori bounds.  See the remark at
the end of the section for a detailed comparison of our method with
Nikulin's.  

\begin{thmdef}
\label{thmdef-shape-of-polygons}
Suppose $P$ is a finite-sided finite-area polygon in $H^2$.  Then one
of its edges is \bfdefn{short}, meaning that the angle bisectors based
at its endpoints cross.  
Furthermore, one of the following holds:
\begin{enumerate}
\item
\label{case-there-exists-short-edge}
$P$ has a short edge
orthogonal to at most one of its neighbors, with its neighbors  not
orthogonal to each other.
\item
\label{case-there-exists-short-pair}
$P$ has at least $5$ edges and a \bfdefn{short pair} $(S,T)$, meaning:
$S$ is a short edge orthogonal to its neighbors, $T$ is one of these
neighbors, and the perpendicular bisector of $S$ meets the angle
bisector based at the far vertex of $T$.
\item
\label{case-there-exists-close-pair}
$P$ has at least $6$ edges and a \bfdefn{close pair} of short edges
$\{S,S'\}$, meaning: $S$ and $S'$ are short edges orthogonal to their
neighbors, their perpendicular bisectors cross, and some \bfdefn{long}
(i.e., not short) edge is adjacent to both of them.
\end{enumerate}
\end{thmdef}

For our purposes, the perpendicular bisector of a finite-length edge
of a  polygon $P$  means the ray departing from the midpoint of
that edge, into the interior of $P$.  
The angle bisector at
a finite vertex is a ray defined in the usual way, and the angle
bisector at an ideal vertex means the set of points equidistant from
(the lines containing) the two edges meeting there.  We think of it as a ray based at the ideal vertex.

\begin{lemma}
\label{lem-rays-into-interior-of-P}
Suppose $P$ is a convex polygon in $H^2\cup\partial H^2$ and $x_1,\dots,x_{n\geq2}$
are distinct points of $\del P$, numbered consecutively in the
positive direction around $\del P$, with subscripts read modulo~$n$.
Suppose a ray $b_i$ departs from each $x_i$, into the interior of $P$.
Then for some neighboring pair $b_i,b_{i+1}$ of these rays, either
they meet or one of them meets the open arc $(x_i,x_{i+1})$ of $\del
P$.
\end{lemma}

\begin{remarks}
 When speaking of arcs of $\del P$, we always mean those traversed in
 the positive direction.  Also,  the proof works for the
 closure of any open convex subset of $H^2$;  one should interpret
 $\del P$ as the topological boundary of $P$ in $\rp^2$.
\end{remarks}

\begin{proof}
By convexity of $P$ and the fact that $b_i$ enters the
interior of $P$, 
its closure $\bar b_i$ meets $\del P$ at one other point $y_i$.  (We
refer to $\bar b_i$ rather than $b_i$ because $y_i$ might lie in $\del
H^2$.)  Now, $y_i$ lies in exactly one arc of $\del P$ of the form
$(x_j,x_{j+1}]$, and we define $s_i$ as the number of steps from $x_i$
 to $x_j$.  Formally: $s_i$ is the number in $\{0,\dots,n-1\}$ such
 that $i+s_i\cong j$ mod $n$.

Now choose $i$ with $s_i$ minimal.  If $s_i=0$ then $b_i$ meets
$(x_i,x_{i+1})$ or contains $x_{i+1}\in b_{i+1}$ and we're done, so suppose
$s_i>0$.  Now, either $b_{i+1}$ meets $(x_i,x_{i+1})$, so we're done,
or it meets $b_i$, so we're done, or else $y_{i+1}$ lies in
$(x_{i+1},y_i]\sset(x_{i+1},x_{j+1}]$.  In the last case $s_{i+1}<s_i$, so
   we have a contradiction of minimality.
\end{proof}

\begin{proof}[Proof of theorem~\ref{thmdef-shape-of-polygons}.]
To find a short edge, apply
lemma~\ref{lem-rays-into-interior-of-P} with $x_1,\dots,\breakok x_n$ the
vertices of $P$ and $b_i$ the angle bisector at $x_i$.   No
$b_i$ can meet $[x_{i-1},x_i)$ or $(x_i,x_{i+1}]$, so some neighboring
pair of bisectors must cross, yielding a short edge.

Suppose for the moment that some short edge is non-orthogonal to at
least one of its neighbors.  Then conclusion~\eqref{case-there-exists-short-edge} applies, except if
its neighbors, call them $T$ and $T'$, are orthogonal to each other.
This can only happen if $P$ is a triangle.  In a triangle all
edges are short and at most one angle equals $\pi/2$.  Therefore
conclusion~\eqref{case-there-exists-short-edge} applies with $T$ or $T'$ taken as the short edge.

We suppose for the rest of the proof that every short edge of $P$ is
orthogonal to its neighbors, because this is the only case where
anything remains to prove.  

Next we show that $P$ has${}\geq5$ sides; suppose otherwise.  We
already treated triangles, so suppose $P$ is a quadrilateral.  Because
$H^2$ has negative curvature, $P$ has a vertex with
angle${}\neq\pi/2$; write $E$ and $F$ for the edges incident there.
The angle bisector based there meets one of the other two edges
(perhaps at infinity), say the one next to $E$.  Drawing a picture
shows that $E$ is short.  So $P$ has a short edge involved in an
angle${}\neq\pi/2$, a contradiction.

Next we claim that if $S$ and $S'$ are adjacent short edges, say with $S'$ at
least as long as $S$, then $(S,S')$ is a short pair.  Write $b,b'$ for
the perpendicular bisectors of $S$, $S'$,  write $\c$ for the angle
bisector at $S\cap S'$, and write $\a$ (resp. $\a'$) for the angle bisector
at the other vertex of $S$ (resp. $S'$).  We hope the left half of
figure~\ref{fig-shape-of-polygon} will help the reader.  Now, $\a$
meets $\c$ because $S$ is short, and $b$ passes though $\a\cap\c$
because it is the perpendicular bisector of an edge whose two angles
are equal.  The same conclusions apply with $\a',b',S'$ in place of
$\a,b,S$.  Furthermore, the intersection point $\a'\cap b'\cap\c$ is
at least as far away from $S\cap S'$ as $\a\cap b\cap\c$ is, because
$S'$ is at least as long as $S$.  If $S'$ is strictly longer than $S$,
then $b$ enters the interior of the triangle bounded by $\c$, $\a'$
and $S'$ by passing through $\c$.  Since $b$ cannot meet $S'$, it must
meet $\a'$ as it exits the triangle.  So $(S,S')$ is a short pair.
(In the limiting case where $S$ and $S'$ have the same length, $\a$,
$b$, $\c$, $b'$ and $\a'$ all meet at a single point.)

\begin{figure}
\scalebox{.45}{\includegraphics{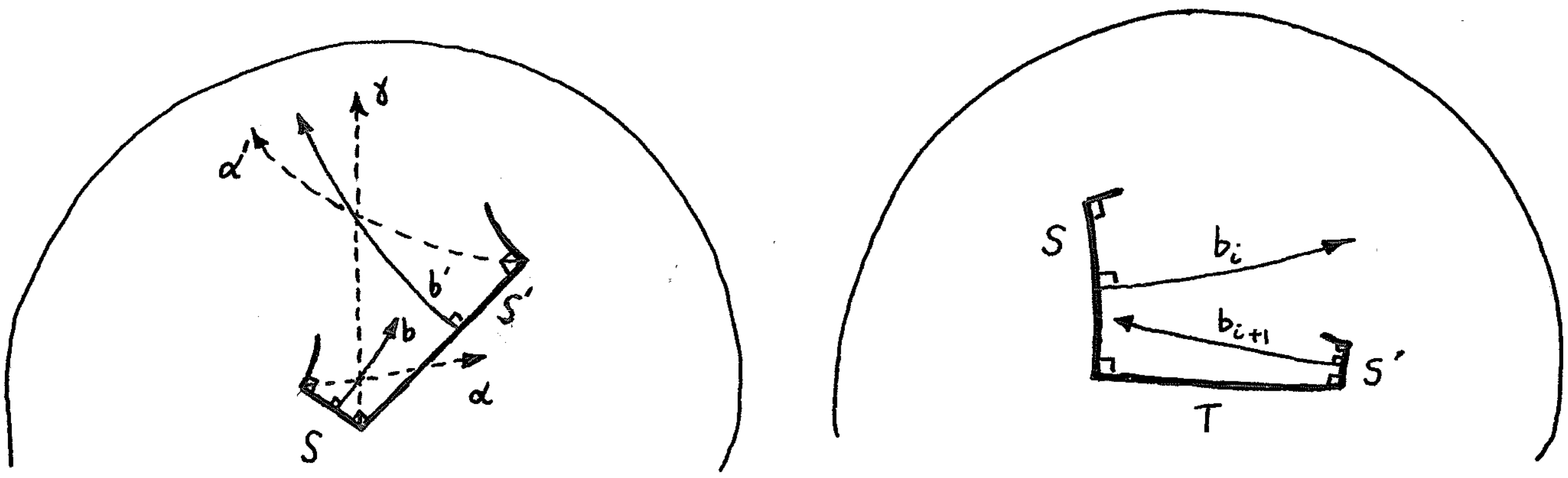}}
\caption{See the proof of theorem~\ref{thmdef-shape-of-polygons}.}
\label{fig-shape-of-polygon}
\end{figure}

Now an argument similar to the one used for quadrilaterals shows that
a pentagon has a
short pair.  Namely, suppose $S$ is a short edge, and write $T,T'$ for
its neighbors and $U,U'$ for their other neighbors.  The perpendicular
bisector of $S$ cannot meet $T$ or $T'$ (since $S\perp T,T'$), so it
meets $U$ or $U'$ (perhaps at infinity), say $U$.
Then the angle bisector at $T\cap U$ either meets $b$, so $(S,T)$ is a
short pair, or it meets $S$.  In the latter case, $T$ is also short,
and the previous paragraph assures us that either $(S,T)$ or $(T,S)$
is a short pair.  This finishes the treatment of a pentagon, so we
take $P$ to have${}\geq6$ edges in the rest of the proof.

At this point it suffices to prove that $P$ has a short pair, or a close pair of
short edges.  
We apply lemma~\ref{lem-rays-into-interior-of-P} again.
Let $x_1,\dots,x_n$ be the midpoints of the short edges and the
vertices not on short edges, in whatever order they occur around $\del
P$.  If $x_i$ is the midpoint of an edge, let $b_i$ be the
perpendicular bisector of that edge.  If $x_i$ is a vertex not on a
short edge, let $b_i$ be the angle bisector there.
Lemma~\ref{lem-rays-into-interior-of-P} now gives us some $i$ such
that either $b_i$ and $b_{i+1}$ meet, or one of them meets
$(x_i,x_{i+1})$.  We now consider the various possibilities:

If both are angle bisectors then $x_i$ and $x_{i+1}$ must be the
endpoints of one edge, or else there would be some other $x_j$ between
them.  The angle bisector based at an endpoint of the edge cannot
cross the interior of the edge, so $b_i$ and $b_{i+1}$ must cross, so
the edge must be short.  But this is impossible because then $b_i$ and
$b_{i+1}$ would not be among $b_1,\dots,b_n$.  So this case cannot
occur.

If one, say $b_i$, is the perpendicular bisector of a short edge $S$,
and the other is an angle bisector, then we may suppose the following holds: the edge $T$ after $S$ is long, and  $b_{i+1}$ is the angle bisector based
at the endpoint of $T$ away from $S$.  Otherwise, some other $x_j$
would come between $x_i$ and $x_{i+1}$.  If $b_i$ and $b_{i+1}$ meet
then $(S,T)$ is a short pair.  Also, $b_i$ cannot meet
$(x_i,x_{i+1})$ since $S\perp T$.  The only other possibility is that
$b_{i+1}$ crosses $S$ between $x_i$ and $S\cap T$.  But in this case
the angle bisectors based at the endpoints of $T$ would cross, so $T$
would also be short, a contradiction.

Finally, suppose $b_i$ and $b_{i+1}$ are the perpendicular bisectors
of short edges $S,S'$.  
If $S,S'$ are adjacent then one of the pairs $(S,S')$, $(S',S)$ is a
short pair; so we suppose they are not adjacent.  There
must be exactly one edge between them, which must be long, for otherwise some other $x_j$
would come between $x_i$ and $x_{i+1}$.  
If $b_i$ and $b_{i+1}$ cross then $S$ and $S'$ forms a close pair of
short edges and we are done.  The only other possibility (up to
exchanging $b_i\leftrightarrow b_{i+1}$, $S\leftrightarrow S'$) is
pictured in the right half of figure~\ref{fig-shape-of-polygon}.
Namely, $b_{i+1}$ meets $S$ between $S\cap T$ and the midpoint of $S$.
In this case, the angle bisector based at $S\cap T$ either meets
$b_{i+1}$, so $(S',T)$ is a short pair, or it meets $S'$, which
implies that $T$ is short, a contradiction.
\end{proof}

In lemmas~\ref{lem-numerical-consequences-of-short-edge}--\ref{lem-numerical-consequences-of-close-short-edges} we give numerical inequalities expressing the
geometric configurations from theorem~\ref{thmdef-shape-of-polygons}.  We
have tried to formulate the results as uniformly as possible.  In each
of the cases (short edge, short pair, close pair of short edges), we
consider three edges $R,T,R'$ of a convex polygon in $H^2$, consecutive
except that there might be a short edge $S$ (resp. $S'$) between $R$
(resp. $R'$) and $T$.  We write $\rhat,\that,\rhat'$ for the
corresponding outward-pointing unit vectors in $\R^{2,1}$ (and $\shat,\shat'$ when
$S,S'$ are present).  We define $\mu=-\rhat\cdot\that$ and
$\mu'=-\rhat'\cdot\that$, and the hypotheses of the lemmas always
give upper bounds on them (e.g., if $R,T$ are adjacent then
$\mu\leq1$).  The main assertion of each lemma is a constraint on
$\lambda:=-\rhat\cdot\rhat'$, namely
\begin{equation}
\label{eq-BASIC-INEQUALITY}
\begin{split}
\lambda<K:={}&1+\mu+\mu'+2\sqrt{1+\mu}\sqrt{1+\mu'}\\
={}&\bigl(\sqrt{1+\mu}+\sqrt{1+\mu'}\,\bigr)^2-1.
\end{split}
\end{equation}
This may be regarded as giving an upper bound on the hyperbolic
distance between (the lines containing) $R$ and $R'$.  The bound
\eqref{eq-BASIC-INEQUALITY} plays a crucial role in section~\ref{sec-shape-of-2-dimensional-Weyl-chamber}.

\begin{lemma}
\label{lem-numerical-consequences-of-short-edge}
Suppose $R,T,R'$ are consecutive edges of a polygon in $H^2$.  Then
$K\leq7$, and $T$ is short if and only if $\lambda<K$.
\end{lemma}

\begin{proof}
Observe that $\rhat-\that$ and $\rhat'-\that$ are vectors orthogonal
to the angle bisectors based at $R\cap T$ and $R'\cap T$; in fact they
point outward from the region bounded by $T$ and these two bisectors.
The  bisectors intersect if and only if
$$
\bigl(\rhat-\that\,\bigr)\cdot\bigl(\rhat'-\that\,\bigr)
<-\sqrt{\bigl(\rhat-\that\,\bigr)^2}\sqrt{\bigl(\rhat'-\that\,\bigr)^2}.
$$ 
(This is the condition that $\rhat-\that$ and $\rhat'-\that$ have
positive-definite inner product matrix.)
Using $\rhat^2=\rhat'^{\,2}=\that^2=1$ and the definitions of $\mu$ and
$\mu'$, then rearranging, yields $\lambda<K$.  And $K$ is largest when
$\mu$ and $\mu'$ are, namely when $T$ is parallel to $R,R'$, yielding
$\mu=\mu'=1$ and 
$K= 7$.
\end{proof}

If a short edge is orthogonal to its neighbors, then its neighbors'
unit vectors' inner product lies in $(-3,-1)$.  This is because $K=3$
when $\mu=\mu'=0$ in lemma~\ref{lem-numerical-consequences-of-short-edge}.  We will use this in the next two lemmas.

\begin{lemma}
\label{lem-numerical-consequences-of-short-pair}
Suppose $R,S,T,R'$ are distinct consecutive edges of a polygon in
$H^2$, with $S$ short and orthogonal to $R$ and $T$.  Then
$\mu\in(1,3)$ and $K<5+4\sqrt2\approx10.657...$, and a necessary
condition for $(S,T)$ to be a short pair is $\lambda<K$.  Furthermore,
the projection $s_0$ of $\rhat'$ to the $\Q$-span of $\shat$ satisfies
\begin{equation}
\label{eq-properties-of-s0}
s_0^2=s_0\cdot\rhat'
=1+\frac{\lambda^2+2\lambda\mu\mu'+\mu'^{\,2}}{\mu^2-1}.
\end{equation}
\end{lemma}

\begin{remark}
In this lemma, $\lambda<K$ is necessary but not sufficient, while in
lemma~\ref{lem-numerical-consequences-of-short-edge} it was both.  The
reason for the difference is that here, if $T$ is short enough then
the angle bisector at $T\cap R'$ will exit $P$ before it can cross the
perpendicular bisector of $S$.  A similar phenomenon occurs in the
next lemma also.
\end{remark}

\begin{proof}
First we note that $\mu>1$ since $R$ and $T$ are ultraparallel.  Also,
$\mu<3$ by the remark before the lemma.  The proof that $\lambda<K$ is exactly
the same as in the previous proof, using the fact that $\rhat-\that$
is orthogonal to the perpendicular bisector of $S$.  (Lemmas \ref{lem-numerical-consequences-of-short-edge}--\ref{lem-numerical-consequences-of-close-short-edges}
acquire a certain unity when the perpendicular bisector of $S$ is
regarded as the ``angle bisector'' of the ultraparallel edges $R$ and
$T$.)  The inequality $K<5+4\sqrt2$ follows from $\mu<3,\mu'\leq1$.

Next, the $\Q$-span of $\shat$ is $\spanof{\rhat,\that\,}^\perp$, so $s_0$ is got
from $\rhat'$ by adding suitable multiples of $\rhat$ and $\that$.
Computation reveals
$$
s_0=\rhat'
+\frac{\lambda+\mu\mu'}{1-\mu^2}\rhat
+\frac{\lambda\mu+\mu'}{1-\mu^2}\that,
$$ 
and then one can work out $s_0\cdot\rhat'$.  Finally,
$s_0^2=s_0\cdot\rhat'$ because $s_0$ was defined as an orthogonal
projection of $\rhat'$. 
\end{proof}

\begin{lemma}
\label{lem-numerical-consequences-of-close-short-edges}
Suppose $R,S,T,S',R'$ are distinct consecutive edges of a polygon in
$H^2$, and assume $S$ and $S'$ are short and orthogonal to their
neighbors.  Then $\mu,\mu'\in(1,3)$, $K<15$, and a necessary condition
for $(S,S')$ to be a close pair of short edges is $\lambda<K$.
Furthermore, the projection $s_0$ of $\rhat'$ to the $\Q$-span of ${\shat}$
satisfies \eqref{eq-properties-of-s0}, and similarly with primed and
unprimed letters exchanged.  Finally,
\begin{equation}
\label{eq-s0-dot-s0prime}
s_0\cdot s_0'=
\lambda+\mu\mu'-\frac{(\lambda+\mu\mu')^3}{(\mu^2-1)(\mu'^{\,2}-1)}.
\end{equation}
\end{lemma}

\begin{proof}
Essentially the same as the previous one.  The computation of
$s_0\cdot s_0'$ is tedious but mechanical.
\end{proof}

\begin{remarks}
The constants $7$, $5+4\sqrt2$ and $15$ appearing in lemmas
\ref{lem-numerical-consequences-of-short-edge}--\ref{lem-numerical-consequences-of-close-short-edges} are essentially the same as the constants $14$,
$10+8\sqrt2$ and $30$ in Nikulin's (4.1.53), (4.1.55) and (4.1.58)
in \cite[Theorem~4.1.8]{Nikulin-rk3}.  The factor of $2$ comes from
the fact that he uses norm~$2$ vectors where we use unit vectors.  And
our short edge, short pair, close pair cases in
theorem~\ref{thmdef-shape-of-polygons} correspond to his type I, II
and III narrow parts of polygons.  In the short edge case, our
condition that the short edge be non-orthogonal to one of its
neighbors corresponds to Nikulin's
condition \cite[Theorem~4.1.4]{Nikulin-rk3} that the Dynkin diagram of
the $3$ roots be connected.

Our results simplify his slightly in that we can suppose $P$
has${}\geq5$ (resp.~$6$) edges in the short pair (resp. close pair)
cases.  More important is our bound on $\rhat\cdot\rhat'$ in terms of
$\rhat\cdot\that$ and $\rhat'\cdot\that$, which we believe is optimal.
The analogous bounds in Nikulin's paper are the second part of
(4.1.4), formula (4.1.10) or (4.1.55), and formula (4.1.13).  Our
bounds are better and simpler in the short pair and close pair cases.
They limit the size of the enumeration in the next section enough that
we can avoid considering class numbers of imaginary quadratic fields,
which Nikulin uses in his section~3.2.
\end{remarks}

\section{The shape of a 2-dimensional Weyl chamber}
\label{sec-shape-of-2-dimensional-Weyl-chamber}

In this section we consider a Lorentzian lattice $L$ of rank~3, and
assume that its Weyl chamber $C$ possesses one of the features
\eqref{case-there-exists-short-edge}--\eqref{case-there-exists-close-pair}
from theorem~\ref{thmdef-shape-of-polygons}.  This always holds if $C$
has finite area.  The goal is to restrict $L$ to one of finitely many
possibilities (up to scale).  This is achieved in
section~\ref{sec-rank-3}, but we do most of the work here, by showing
that the simple roots corresponding to 3, 4 or 5 suitable consecutive
edges have one of finitely many inner product matrices (up to scale).
The point is not really this finiteness, which was already known by Nikulin's
finiteness theorem \cite{Nikulin-finiteness}.  What is important is that all
possibilities are worked out explicitly.  The main tool is
lemma~\ref{lem-integerization-of-inner-products} below, which allows
one to obtain information about two roots of $L$, given information
about the corresponding unit vectors.

In fact we will work with ``quasiroots'' of $L$ rather than roots; the
point is that the methods of this section can be used to gain control
of the 3-dimensional lattice associated to any finite-area
2-dimensional face of any Weyl chamber.  
The author intends to use this in higher-dimensional analogues of this paper.  For
purposes of classifying rank~3 reflective lattices, we will not need
this generalization, so the reader may read ``root'' in place of
``quasiroot'' and assume that $\rho$, $\sigma$, $\tau$, $\sigma'$ and
$\rho'$ are always $2$.  For use in the general case, here is the
definition.  If $L$ is any lattice then we say that $a\in L$ is a
\defn{quasiroot} of $L$ of \defn{weight} $\a\in\{1,2,\dots\}$ if $a$
is primitive, $a^2>0$, and $L\cdot a\sset\frac{1}{\a}a^2\Z$.  
Notice that scaling a lattice up or down doesn't affect the quasiroot
weights.

When a
roman letter represents a quasiroot we will
use the corresponding greek letter for its weight.  A quasiroot of weight~2 is the same thing as a root,
and quasiroot of weight~1 is the same thing as a generator of a
1-dimensional positive-definite summand of $L$.  Of course, every primitive
positive-norm vector of $L$ is a quasiroot of weight equal to its norm.
But sometimes one has an a priori bound on the
weights of quasiroots involved, independent of their norms.

\begin{lemma}
\label{lem-integerization-of-inner-products}
Suppose $a,b$ are quasiroots of weights $\a$, $\b$ in a lattice $L$,
whose associated unit vectors satisfy $\ahat\cdot\bhat=-K$ for some $K>0$.  Then the inner product
matrix of $a$ and $b$ is a rational number times
\begin{equation}
\label{eq-integerization-of-inner-products}
\begin{pmatrix}
\a u&-uv\\
-uv&\b v
\end{pmatrix},
\end{equation}
for some positive integers $u,v$ satisfying $uv=K^2\a\b$.
\end{lemma}

We will also use versions of this lemma with both $=$ signs replaced
by $\leq$ or both by $<$.  These are formal consequences of this
version. 

\begin{proof}
We know $a\cdot b\in\frac{1}{\a}a^2\Z$ and $a\cdot
b\in\frac{1}{\b}b^2\Z$, so
\begin{equation}
\label{eq-foo}
a\cdot b=-va^2/\a
\qquad\hbox{and}\qquad
a\cdot b=-ub^2/\b
\end{equation}
for some positive integers $u,v$.  This gives
$$
-a\cdot b
=
\sqrt{(-a\cdot b)^2}
=
\sqrt{(va^2/\a)(ub^2/\b)}
=
\sqrt{a^2}\sqrt{b^2}\sqrt{uv/\a\b}.
$$
Rewriting $-\ahat\cdot\bhat=K$ as $-a\cdot b=K\sqrt{a^2}\sqrt{b^2}$
yields $uv=K^2\a\b$.  For convenience we scale the inner products so
that $a\cdot b=-uv$, and then \eqref{eq-foo} gives $a^2=\a u$ and
$b^2=\b v$.
\end{proof}

The hypothesis $a\cdot b\neq0$ in the lemma is essential and annoying:
when two quasiroots are orthogonal, their norms need not satisfy any
relation at all.  For example, in the proof of
lemma~\ref{lem-shape-of-2d-SE-integral-polygon}, the hypothesis
$R\notperp T,R'$ will be used to relate the norms of $r$, $t$ and
$r'$, an essential part of the proof.  If we could do away with the
non-orthogonality hypothesis in lemma~\ref{lem-integerization-of-inner-products} then we could carry out
our enumeration of lattices using only the existence of a short edge.
The possibility of orthogonality is what forces us to consider the
more-complicated configurations of theorem~\ref{thmdef-shape-of-polygons}.

The setup for
lemmas~\ref{lem-shape-of-2d-SE-integral-polygon}--\ref{lem-shape-of-2D-CP-integral-polygon}
below is similar to that for
lemmas~\ref{lem-numerical-consequences-of-short-edge}--\ref{lem-numerical-consequences-of-close-short-edges}.
Their hypotheses are exactly the ``outputs'' of the short edge, short
pair and close pair cases of theorem~\ref{thmdef-shape-of-polygons}.
$P$ is a finite-sided finite-area polygon in $H^2$, with angles at
most $\pi/2$, whose edges are the orthogonal complements of quasiroots
of a rank~3 Lorentzian lattice $L$.  $R$, $T$ and $R'$ are distinct
edges of $P$, consecutive except that there might be a short edge $S$
(resp. $S'$) between $R$ (resp. $R'$) and $T$.  We write $r,t,r'$ (and
$s,s'$ when $S,S'$ are present) for the outward-pointing quasiroots of
$L$ orthogonal to $R,\dots,R'$.  We suppose that $r$, $t$ and $r'$
(and $s$ and $s'$ if present) are quasiroots of weights $\rho$, $\tau$
and $\rho'$ (and $\sigma$ and $\sigma'$ if relevant).  As usual, we use hats to
indicate the corresponding unit vectors.  As in section~\ref{sec-shape-of-polygon} we define
$\mu=-\rhat\cdot\that$, $\mu'=-\rhat'\cdot\that$,
$$
K=1+\mu+\mu'+2\sqrt{1+\mu}\sqrt{1+\mu'},
$$
and $\lambda=-\rhat\cdot\rhat'$.  In
lemmas~\ref{lem-shape-of-2d-SE-integral-polygon}--\ref{lem-shape-of-2D-CP-integral-polygon}
$\mu$, $\mu'$ and $\lambda$ will be given by the formulas
$\mu=\sqrt{AB/\rho\tau}$, $\mu'=\sqrt{A'B'/\rho'\tau}$ and
$\lambda=\sqrt{CC'/\rho\rho'}$, where $A,B,A',B',C,C'$ are
non-negative integers described in the lemmas.  In particular, $K$ is
a function of $A,B,A',B'$.  We have tried to make the treatment as
uniform as possible, but the essential difference between orthogonal
and non-orthogonal quasiroots forces separate treatment of the cases
$T\perp R'$, $T\notperp R'$ in each of
lemmas \ref{lem-shape-of-2d-SE-integral-polygon}
and~\ref{lem-shape-of-2d-SP-integral-polygon}.

\begin{lemma}
\label{lem-shape-of-2d-SE-integral-polygon}
Suppose $R,T,R'$ are consecutive edges of $P$, $T$ is short, and
$R\notperp T,R'$.  Then the inner product matrix of $r,t,r'$ is one of
finitely many possibilities, up to scale.

More precisely, let $(A,B)$ vary over all pairs of positive integers
satisfying $AB\leq\rho\tau$ and let $(A',B')$ vary over $(0,0)$ and
all pairs of positive integers satisfying $A'B'\leq\rho'\tau$.  For
fixed such $A,B,A',B'$, let $(C,C')$ vary over all pairs of positive
integers satisfying
\begin{displaymath}
CC' < K^2\rho\rho'
\end{displaymath}
and
\begin{equation}
\label{eq-def-of-beta-short-side-case}
AB'C' = A'BC\hbox{; call the common value $\b$.}
\end{equation}
Then
for some such $A,B,A',B',C,C'$,
the inner product matrix of
$r,t,r'$ is a rational multiple of
\begin{gather}
\label{eq-shape-of-2d-SE-integral-polygon-NONorthogonality-case}
\begin{pmatrix}
\rho AB'&   -ABB'   &    -\b   \\
-ABB'   & \tau BB'  &  -A'B'B  \\
-\b     &  -A'B'B   & \rho'A'B
\end{pmatrix}
\hbox{ if $A',B'\neq0$}\\
\label{eq-shape-of-2d-SE-integral-polygon-orthogonality-case}
\begin{pmatrix}
\rho AC &   -ABC   &   -ACC'   \\
-ABC    & \tau BC  &     0     \\
-ACC'   &     0    & \rho'AC'
\end{pmatrix}
\hbox{ if $A',B'=0$.}
\end{gather}
\end{lemma}

\begin{proof}
We have $\mu\in(0,1]$ by hypothesis, so lemma~\ref{lem-integerization-of-inner-products} says there exist
positive integers $A,B$ with $AB\leq\rho\tau$, such that $r,t$ have
inner product matrix a rational multiple of
\begin{equation}
\label{eq-r-t-inner-product-matrix}
\begin{pmatrix}
\rho A& -AB\\
-AB&\tau B
\end{pmatrix}.
\end{equation}

Now suppose $T\notperp R'$.  Then the argument of the previous
paragraph applies with primed and unprimed symbols exchanged.  And by
shortness of $T$, $\lambda<K$ (lemma~\ref{lem-numerical-consequences-of-short-edge}).  Then another
application of lemma~\ref{lem-integerization-of-inner-products} (using the hypothesis $R\notperp R'$)
shows there exist positive integers $C,C'$ with $CC'<K^2\rho\rho'$, such
that $r,r'$ have inner product matrix a rational multiple of
\begin{equation}
\label{eq-r-rprime-inner-product-matrix}
\begin{pmatrix}
\rho C &-CC'\\
-CC'&\rho'C'
\end{pmatrix}.
\end{equation}
Now we see
$$
\frac{\rho'C'}{\rho C}
=
\frac{r'^{\,2}}{r^2}
=
\frac{r'^{\,2}}{t^2}\frac{t^2}{r^2}
=
\frac{\rho'A'}{\tau B'}\frac{\tau B}{\rho A},
$$
so $AB'C'=A'BC$.  This establishes \eqref{eq-def-of-beta-short-side-case}.  Now we put all three
$2\times2$ matrices together by choosing the scale at which $t^2=\tau
BB'$, leading to \eqref{eq-shape-of-2d-SE-integral-polygon-NONorthogonality-case}.

It remains to consider the case $T\perp R'$.  Then we take $A'=B'=0$
and note that $\sqrt{A'B'/\rho'\tau}$ is indeed equal to
$\mu'=-\rhat'\cdot\that=0$, so this is consistent with the definition
of $\mu'$.  The same argument as before shows that
$r,r'$ have inner product matrix a rational multiple of
\eqref{eq-r-rprime-inner-product-matrix}.  This time
\eqref{eq-def-of-beta-short-side-case} is trivially true because both
sides vanish.  Now we put together \eqref{eq-r-t-inner-product-matrix}
and \eqref{eq-r-rprime-inner-product-matrix} by choosing the scale at
which $r^2=\rho AC$.  Finally, using $r'\cdot t=0$, the resulting
inner product matrix is
\eqref{eq-shape-of-2d-SE-integral-polygon-orthogonality-case}.
\end{proof}

\begin{lemma}
\label{lem-shape-of-2d-SP-integral-polygon}
Suppose $P$ has at least $5$ edges, and $R,S,T,R'$ are consecutive
edges with $(S,T)$ a short pair.  Then the inner product matrix of
$r,s,t,r'$ is one of finitely many possibilities, up to scale. 

More precisely, let $(A,B)$ vary over all pairs of positive integers
with $\rho\tau<AB<9\rho\tau$ and let $(A',B')$ vary over $(0,0)$ and
all pairs of positive integers with $A'B'\leq\rho'\tau$.  For fixed
such $A,B,A',B'$, let $(C,C')$ vary over all pairs of positive
integers satisfying
\begin{gather}
\label{eq-constraint-on-CCprime-SP-case}
\rho\rho'<CC'<K^2\rho\rho'
\\
\label{eq-def-of-beta-short-pair-case}
AB'C' {}= A'BC\hbox{; call the common value $\b$, and}
\\
\label{eq-def-of-N-short-pair-case}
N:=\rho'\sigma+\sigma
\frac{CC'\tau+2\b+A'B'\rho}{AB-\rho\tau}
\hbox{ is an integer}.
\end{gather}
For fixed such $A,B,A',B',C,C'$, let $k$ vary over all positive
integers dividing $N$.  Then for some such $A,B,A',B',C,C',k$, the
inner product matrix of $r,s,t,r'$ is a rational multiple of
\begin{gather}
\label{eq-shape-of-2d-SP-integral-polygon-NONorthogonality-case}
\begin{pmatrix}
\rho AB'&        0      &    -ABB'    &     -\b         \\
0       &\sigma A'BN/k^2&      0      &   -A'BN/k   \\
-ABB'   &        0      & \tau BB'   &    -A'B'B       \\
-\b     &  -A'BN/k      &    -A'B'B   &   \rho'A'B
\end{pmatrix}
\hbox{ if $A',B'\neq0$}\\
\label{eq-shape-of-2d-SP-integral-polygon-orthogonality-case}
\begin{pmatrix}
\rho AC &   0             &  -ABC   &   -ACC' \\
0       & \sigma AC'N/k^2 &    0    & -AC'N/k \\
-ABC    &   0             & \tau BC &     0   \\
-ACC'   & -AC'N/k         &    0    & \rho'AC'
\end{pmatrix}
\hbox{ if $A',B'=0$.}
\end{gather}
\end{lemma}

\begin{proof}
Lemma~\ref{lem-numerical-consequences-of-short-edge} says that $1<\mu<3$.  So apply the previous proof with
this inequality in place of $0<\mu\leq1$.  The result is that $r,t,r'$
have inner products (up to scale) as in \eqref{eq-shape-of-2d-SP-integral-polygon-NONorthogonality-case} or \eqref{eq-shape-of-2d-SP-integral-polygon-orthogonality-case}, where
$A,B,A',B',C,C'$ have the stated properties up to and including
\eqref{eq-def-of-beta-short-pair-case}.  For the treatment of the inner product matrix of $r,r'$
in that argument, note that $R\notperp R'$ by the fact that $P$ has at
least $5$ edges.  Indeed $R$ and $R'$ are not even adjacent.  Since $P$
has all angles${}\leq\pi/2$, this says that the lines containing
$R$ and $R'$ are ultraparallel.  This gives $\lambda>1$, which
proves the $CC'>\rho\rho'$ inequality in \eqref{eq-constraint-on-CCprime-SP-case}.

Now we consider $s$.  As in lemma~\ref{lem-numerical-consequences-of-short-pair}, let $s_0$ be the projection
of $\rhat'$ to $\Q\cdot{s}$, so the projection of $r'$ is
$\sqrt{r'^{\,2}}s_0$, which lies in $\frac{1}{\sigma}s\Z$ since 
$s$ is a quasiroot of weight $\sigma$.
Therefore there is a
positive integer $k$ such that $\sqrt{r'^{\,2}}s_0=-ks/\sigma$, i.e.,
\begin{equation}
\label{eq-formula-for-s}
s=-\sigma\sqrt{r'^{\,2}}s_0/k.
\end{equation}
On the other hand, $r'$ is a quasiroot of weight $\rho'$, so $s\cdot
r'\in\frac{1}{\rho'}r'^{\,2}\Z$.  To work with this condition we define
the integer
\begin{equation}
\begin{split}
\label{eq-formula-for-M}
M:&{}=\rho'\frac{s\cdot r'}{r'^{\,2}}
=\rho'\frac{s\cdot\sqrt{r'^{\,2}}s_0}{r'^{\,2}}
=\rho'\frac{-\sigma\sqrt{r'^{\,2}}\sqrt{r'^{\,2}}}{k r'^{\,2}}s_0^2
=-\frac{\rho'\sigma}{k}s_0^2\\
&{}=-\frac{\rho'\sigma}{k}
\Bigl(1+\frac{\lambda^2+2\lambda\mu\mu'+2\mu'^{\,2}}{\mu^2-1}\Bigr),
\end{split}
\end{equation}
where the last step uses \eqref{lem-numerical-consequences-of-short-pair}.  A calculation shows $M=-N/k$ where
$N$ is defined in \eqref{eq-def-of-N-short-pair-case}.  Since $M$ is
an integer we must have
$N\in\Z$ and $k|N$.

We have shown that $A,B,A',B',C,C',k$ have all the properties claimed,
and all that remains is to work out the remaining entries in the inner
product matrix.  We have $s\cdot r=s\cdot t=0$ by the short pair
hypothesis.  Referring to \eqref{eq-formula-for-M} gives
$$
s\cdot r'=Mr'^{\,2}/\rho'=-Nr'^{\,2}/k\rho'
$$
and referring to \eqref{eq-formula-for-s} and \eqref{eq-formula-for-M} gives
$$
s^2=\frac{\sigma^2}{k^2}r'^{\,2}s_0^2
=\frac{\sigma^2}{k^2}r'^{\,2}\Bigl(-\frac{Mk}{\rho'\sigma}\Bigr)
=\frac{\sigma}{k}r'^{\,2}\Bigl(-\frac{-N/k}{\rho'}\Bigr)
=N\sigma r'^{\,2}/\rho'k^2.
$$
This allows one to fill in the missing entries of \eqref{eq-shape-of-2d-SP-integral-polygon-NONorthogonality-case} and
\eqref{eq-shape-of-2d-SP-integral-polygon-orthogonality-case}, using the values for $r'^{\,2}$ given there.
\end{proof}

\begin{lemma}
\label{lem-shape-of-2D-CP-integral-polygon}
Suppose $P$ has at least $6$ edges and $R,S,T,S',R'$ are consecutive
edges, with $S$ and $S'$ forming a close pair of short edges.  Then
the inner product matrix for $r,s,t,s',r'$ is one of finitely many
possibilities, up to scale.

More precisely,  let $(A,B)$ vary over all pairs of positive integers
satisfying $\rho\tau<AB<9\rho\tau$, and let $(A',B')$ vary over all
pairs of positive integers satisfying $\rho'\tau<A'B'<9\rho'\tau$.
For fixed such $A,B,A',B'$, let $(C,C')$ vary over all pairs of
positive integers satisfying
\begin{gather}
\rho\rho'<CC'<K^2\rho\rho'
\\
\label{eq-def-of-beta-close-pair-case}
AB'C'= A'BC\hbox{; call the common value $\b$, and}
\\
\label{eq-N-and-Nprime-in-Z}
\hbox{$N$ and $N'$ are integers,}
\end{gather}
where $N$ is defined as in \eqref{eq-def-of-N-short-pair-case} and $N'$ similarly, with primed
and unprimed letters exchanged.
For fixed such $A,B,A',B',C,C'$, let $(k,k')$ 
vary over all pairs of positive integers with $k|N$ and $k'|N'$, such that
\begin{equation}
\label{eq-integrality-conditions-for-s-and-sprime}
{\c k^2}/{A'BN}
\hbox{ and }
{\c k'^{\,2}}/{AB'N'}
\hbox{ are integers,}
\end{equation}
where
\begin{equation}
\label{eq-def-of-gamma}
\c=
\frac{\sigma\sigma'\b}{kk'\tau}
\Biggl(\tau+\frac{\b}{CC'}-
\frac{(\tau CC'+\b)^3}{(AB-\rho\tau)(A'B'-\rho'\tau)C^2C'^{\,2}}
\Biggr).
\end{equation}
Then for some such $A,B,A',B',C,C',k,k'$, the inner product matrix of
$r,s,t,s',r'$ is a rational multiple of
\begin{equation}
\label{eq-shape-of-2d-CP-integral-polygon}
\begin{pmatrix}
\rho AB' &       0         &  -ABB'   &  -AB'N'/k'        &     -\b    \\
     0   & \sigma A'BN/k^2 &    0     &       \c          & -A'BN/k \\
   -ABB' &       0         & \tau BB' &        0          &    -A'B'B   \\
-AB'N'/k' &      \c         &    0     & \sigma'AB'N'/k'^{\,2} &       0   \\
   -\b   &  -A'BN/k        & -A'B'B   &        0          & \rho'A'B
\end{pmatrix}.
\end{equation}
\end{lemma}

\begin{proof}
 Mostly this is a repeat of the $\mu'\neq0$ cases of the previous
 proofs.  Follow the proof of
 lemma~\ref{lem-shape-of-2d-SE-integral-polygon}, using
 $\mu,\mu'\in(1,3)$ from
 lemma~\ref{lem-numerical-consequences-of-short-edge} and
 $\lambda\in(1,K)$ from
 lemma~\ref{lem-numerical-consequences-of-close-short-edges}.  The
 result is that after rescaling, $r,t,r'$ have inner products as in
\eqref{eq-shape-of-2d-SE-integral-polygon-NONorthogonality-case}, where $A,B,A',B',C,C'$ satisfy the conditions up to and
 including \eqref{eq-def-of-beta-close-pair-case}.  Then lemma~\ref{lem-shape-of-2d-SP-integral-polygon}'s analysis of $N$ and $k$
 applies, and also applies with primed and unprimed symbols
 exchanged.  This shows that $N,N',k,k'$ have all the properties
 claimed except \eqref{eq-integrality-conditions-for-s-and-sprime}, and justifies all
 entries in \eqref{eq-shape-of-2d-CP-integral-polygon} except the $s\cdot s'$ entry.  To work this out
 we use $s=-\sigma\sqrt{r'^{\,2}}s_0/k$ by \eqref{eq-formula-for-s},
 $s'=-\sigma'\sqrt{r^2}s_0'/k'$ similarly, and formula \eqref{eq-s0-dot-s0prime} for
 $s_0\cdot s_0'$.  The result is
\begin{align*}
s\cdot s'
&{}=\frac{\sigma\sigma'}{kk'}\sqrt{r^2\,r'^{\,2}}\,s_0\cdot s_0'\\
&{}=\frac{\sigma\sigma'}{kk'}\sqrt{r^2\,r'^{\,2}}
\Bigl(\lambda+\mu\mu'-\frac{(\lambda+\mu\mu')^3}{(\mu^2-1)(\mu'^{\,2}-1)}\Bigr).
\end{align*}
A calculation shows that this equals $\c$, defined in \eqref{eq-def-of-gamma}.  This
gives the last entry in the matrix.  Finally, \eqref{eq-integrality-conditions-for-s-and-sprime} is the
condition that $\c=s\cdot s'$ lies in $\frac{1}{\sigma}s^2\Z$ and in
$\frac{1}{\sigma'}s'^{\,2}\Z$, which holds because $s$ and $s'$ are
quasiroots.
\end{proof}

\section{Corner symbols}
\label{sec-corner-symbols}

In this section we introduce a compact notation called a corner symbol for describing
a corner in $H^2$ of the Weyl chamber of a lattice $L$ of signature
$(2,1)$.  It is a formalization of the ideas in
\cite[pp. 59--62]{Nikulin-rk3}, and is 
needed for the proof of our main theorem in section~\ref{sec-rank-3} and presenting our results in
the table.
But the reader should skip it until it is needed.

Throughout this section suppose $L$ is a 3-dimensional Lorentzian
lattice and $e,f$ are consecutive edges of some Weyl chamber.  We will
define the ``corner symbol'' $S(L,e,f)$, whose key property is that
together with $D:=\det L$ it characterizes $(L,e,f)$ up to isometry.
First, there is a unique chamber $C\sset H^2$ of which $e,f$ are
edges.  Choose a component $\Ctilde$ of its preimage in
$L\tensor\R-\{0\}$.  Take $x,y$ to be simple roots corresponding to
$e,f$; implicit in this definition is that they point away from
$\Ctilde$.  Define $z$ as the primitive lattice vector in
$\spanof{x,y}^\perp$ that lies in (the closure of) $\Ctilde$.  We will
define the corner symbol in terms of $x,y,z$, so strictly speaking we
should write $S(L,e,f,\Ctilde)$ or $S(L,x,y,z)$.  But in the end the
choice of $\Ctilde$ is irrelevant because changing it simply negates
$x,y,z$.  So $S(L,e,f)$ will be well-defined.  We will write $X$ and
$Y$ for $x^2$ and $y^2$.  First we will define corner symbols for
corners in $H^2$, and then for ideal vertices.

The definitions involve some standard lattice-theoretic terms, which
we define here for clarity.  If $L$ is a lattice and $M$ a sublattice,
then $M$'s saturation $M_\sat$ is defined as $L\cap(M\tensor\Q)$, and
$M$ is called saturated if $M=M_\sat$.  Now suppose $L$ is integral
and $M$ is saturated.  Then $L$ lies between $M\oplus
M^\perp$ and $(M\oplus M^\perp)^*$, so it corresponds to a subgroup of $\D(M)\oplus\D(M^\perp)$.
This subgroup is called the glue group of the inclusion $M\oplus
M^\perp\to L$.  Because of saturation, it meets the two summands $\D(M)$,
$\D(M^\perp)$ trivially.

\begin{thmdef}
\label{thmdef-finite-corner-symbols}
 Suppose $L$, $x$, $y$ and $z$ are as stated, with $z$ defining a
 finite vertex of $L$'s Weyl chamber.  Then for exactly one row of
 table~\ref{tab-finite-corner-symbols} do they satisfy the conditions listed in the middle
 three columns.  We define their \bfdefn{corner symbol} as the entry in
 the first column.  The corner symbol and $D$ determine the quadruple
 up to isometry, and $z^2$ is given in the last column in terms of
 $D$, $X$, $Y$ and $M:=\min\{X,Y\}$.
\end{thmdef}

\begin{table}
\centeroverfull{%
\begin{tabular}{ccccc}%
{\bf Corner}
&
&$\spanof{x,y}\to\spanof{x,y}_\sat$
&$\spanof{x,y}_\sat\oplus\spanof{z}\to L$\\
{\bf symbol}
&{\bf Angle}  
&{\bf by adjoining}
&{\bf by adjoining}
&$z^2$\\
\hline
\noalign{\smallskip}
$X_2Y$
&$\pi/2$
&nothing
&nothing
&$D/XY$\\
\noalign{\smallskip}
$X_2^sY$
&$\pi/2$
&$(x+y)/2$
&nothing
&$4D/XY$\\
\noalign{\smallskip}
$X_2^bY$
&$\pi/2$
&nothing
&$(x+y+z)/2$
&$4D/XY$\\
\noalign{\smallskip}
$X_2^lY$
&$\pi/2$
&nothing
&$(x+z)/2$
&$4D/XY$\\
\noalign{\smallskip}
$X_2^rY$
&$\pi/2$
&nothing
&$(y+z)/2$
&$4D/XY$\\
\noalign{\smallskip}
$X_2^*Y$
&$\pi/2$
&$(x+y)/2$
&$(x+z)/2$, $(y+z)/2$
&$16D/XY$\\
%
\noalign{\smallskip}
\noalign{\hrule}
\noalign{\smallskip}
$X_3^\pm X$
&$\pi/3$
&nothing
&$(2x+y\pm z)/3$
&$12D/X^2$\\
\noalign{\smallskip}
$X_4Y$
&$\pi/4$
&nothing
&nothing
&$D/M^2$\\
\noalign{\smallskip}
$X_4^*Y$
&$\pi/4$
&nothing
&$(z+\hbox{longer of $x,y$})/2$
&$4D/M^2$\\
\noalign{\smallskip}
$X_6Y$
&$\pi/6$
&nothing
&nothing
&$4D/3M^2$\\
\end{tabular}%
}
\medskip
\caption{Corner symbols for finite vertices of the Weyl chamber of 
 a 3-dimensional Lorentzian lattice $L$;  see theorem~\ref{thmdef-finite-corner-symbols}.
Here $D=\det L$ and
 $M=\min\{X,Y\}$.   When the angle is $\pi/4$ (resp. $\pi/6$), one
 of $X$ and $Y$ is 2 (resp. 3) times the other.}
\label{tab-finite-corner-symbols}
\end{table}

Here are some mnemonics for the notation.  The subscript always gives
the angle at that corner.  The absence of a superscript indicates the
absence of saturation and glue, so that $\spanof{x,y}$ is a summand of
$L$.  An asterisk indicates the most complicated saturation and/or
glue possible.  The remaining superscripts are
$s$ (saturation but no glue), $l$ (glue involving
the left root), $r$ (glue involving the right root),  $b$ (glue
involving both roots) and  $\pm$ (self-explanatory).

\begin{proof}
 At a finite vertex of a Weyl chamber, the angle is $\pi/n$ for
 $n\in\{2,3,4,6\}$, with $x$ and $y$ being simple roots for an
 $A_1^2$, $A_2$, $B_2$ or $G_2$ root system.   The
 projection of $L$ to $\spanof{x,y}\tensor\Q$ lies in the reflective
 hull of $x$ and $y$ (see section~\ref{sec-background}), and $x$ and
 $y$ are primitive in $L$.  This
 constrains the possible enlargements
 $\spanof{x,y}\to\spanof{x,y}_\sat$ and
 $\spanof{x,y}_\sat\oplus\spanof{z}\to L$ so tightly that all
 possibilities are listed.  (In the $A_2$ case there
 must be trivial saturation and nontrivial glue because otherwise
 every element of $W(G_2)=\aut\spanof{x,y}$ would extend to $L$, so
 $x$ and $y$ would not be simple roots.  Similarly, most of the
 $A_1^2$ cases can only occur if $X\neq Y$.)  Also, $z^2$ is easy to work
 out in each case.  The table shows how to reconstruct $(L,x,y,z)$
 from knowledge of $D$ and the corner symbol.
\end{proof}

Corner symbols at ideal vertices are more complicated.  Most of
the information is contained in the structure of $x^\perp$, which we
encode as follows.

\begin{lemma}
 \label{lem-lemma-for-ideal-corner-symbols}
 Suppose $L$ is a $2$-dimensional integral Lorentzian lattice and
 $z\in L$ is primitive isotropic.  Then there is a unique primitive
 null vector $u\in L$ with $z\cdot u<0$, and we define $\lambda$ as
 $z\cdot u$ and $s$ as $[L:\spanof{z,u}]$.  There exists a unique
 $a\in L$ of the form $a=\frac{1}{s}(\sigma z+ u)$ with
 $\sigma\in\{0,\dots,s-1\}$.  The triple $(\lambda,s,\sigma)$
 characterizes $(L,z)$ up to isometry.  The triples that arise by
 this construction are exactly the triples of integers with
 $\lambda<0$, $s>0$, $s^2|2\lambda$, $0\leq\sigma<s$ and
 $(\sigma,s)=1$.
\end{lemma}

\begin{proof}
 The uniqueness of $u$ is clear, and $\lambda<0$ by construction.
 Since $[L:\spanof{z,u}]=s$, $L$ lies between $\frac{1}{s}\spanof{z,u}$ and
 $\spanof{z,u}$, so it corresponds to a subgroup of order $s$ in
 $\frac{1}{s}\spanof{z,u}/\spanof{z,u}\iso(\Z/s)^2$.  By the
 primitivity of $z$ and $u$, this subgroup meets the subgroups
 $\spanof{\frac{1}{s}z}/\spanof{z}$ and
 $\spanof{\frac{1}{s}u}/\spanof{u}$ trivially, so it must be the
 graph of an isomorphism between them.  The
 existence and uniqueness of $a$ follows, as does the condition
 $(\sigma,s)=1$.  The inner product matrix of $L$ with respect to the basis $z$, $a$ is
$$
\begin{pmatrix}0&\lambda/s\\ \lambda/s&2\lambda\sigma/s^2\\\end{pmatrix}
$$
so integrality forces $s^2|2\lambda\sigma$, hence $s^2|2\lambda$.  If
$s$ satisfies this, then the rest of the matrix is automatically
integral, which justifies our claim about which triples can arise.
That $(\lambda,s,\sigma)$ determines $(L,z)$ follows from the
intrinsic nature of our constructions.
\end{proof}

\begin{thmdef}
\label{thmdef-ideal-corner-symbols}
Suppose $(L,x,y,z)$ are as in theorem~\ref{thmdef-finite-corner-symbols}, except that
$z$ defines an ideal vertex of $L$'s Weyl chamber.  Suppose also that
$z$ and $a$ are the basis of $x^\perp$ associated to the pair
$(x^\perp,z)$ by lemma~\ref{lem-lemma-for-ideal-corner-symbols}.  Then for exactly one row of
table~\ref{tab-ideal-corner-symbols} do the conditions in the last two columns hold.  We
define the \bfdefn{corner symbol} of $(L,x,y,z)$ to be the entry in the first
column.  The quadruple can be recovered up to isometry from its symbol
and knowledge of $D=\det L$.
\end{thmdef}

\begin{table}
\begin{tabular}{ccc}%
{\bf Corner}
&{\bf invariants}
&{$\spanof{x}\oplus x^\perp\to L$}\\
{\bf symbol}
&{\bf of $(x^\perp,z)$}
&{\bf by adjoining}\\
\hline
\noalign{\smallskip}
$X\corner{s,\sigma}{\infty}Y$
&$(-s\sqrt{-D/X},s,\sigma)$
&nothing
\\
\noalign{\smallskip}
$X\corner{s,\sigma}{\infty\gluez}Y$
&$(-2s\sqrt{-D/X},s,\sigma)$
&$(x+z)/2$
\\
\noalign{\smallskip}
$X\corner{s,\sigma}{\infty\gluea}Y$
&$(-2s\sqrt{-D/X},s,\sigma)$
&$(x+a)/2$
\\
\noalign{\smallskip}
$X\corner{s,\sigma}{\infty\glueboth}Y$
&$(-2s\sqrt{-D/X},s,\sigma)$
&$(x+z+a)/2$
\\
\end{tabular}
\medskip
\caption{Corner symbols for ideal vertices of the Weyl chamber of 
 a 3-dimensional Lorentzian lattice of determinant $D$; 
 see theorem~\ref{thmdef-ideal-corner-symbols}.
 Note that the values of $s$ and $\sigma$ appear, while the 
 letters $z,a,b$ appear literally when relevant, for example 
 ${26}\above{1pt}{1pt}{4,3}{\infty b}{104}$.}
\label{tab-ideal-corner-symbols}
\end{table}

\begin{proof}
Because $x$ is a root, the glue group of $\spanof{x}\oplus x^\perp\to L$ is trivial
 or $\Z/2$.
 The determinant of $x^\perp$ is $D/X$ or $4D/X$ in these two cases.
  If $(\lambda,s,\sigma)$ are the invariants of
 $(x^\perp,z)$, then $\det x^\perp=\lambda^2/s^2$ gives us an
 equation for $\lambda$ in terms of $D$, $X$ and $s$.  If
 $\lambda=-s\sqrt{-D/X}$ then there is no glue, so the conditions of
 the top row are satisfied.  If $\lambda=-2s\sqrt{-D/X}$ then there is
 $\Z/2$ glue, and the last three rows of the table list all the
 possible gluings of $\spanof{x}$ to $x^\perp$.  This shows that the
 conditions of exactly one row of the table are satisfied, so the
 corner symbol is defined.  

Now we show that it and $D$ determine $(L,x,y,z)$.  By referring to
 the table we see that they determine $(x^\perp,z)$ up to isometry.
 Since the norm of $x$ is part of the symbol, they determine
 $\bigl(\spanof{x}\oplus x^\perp,x,z\bigr)$ up to isometry.  Since the
 symbol also indicates the glue, they determine $(L,x,z)$ up to
 isometry.  These determine the Weyl chamber uniquely, since of all
 the Weyl chambers of $L$ incident at $z$, there is only one for which
 $x$ is a simple root.  This in turn determines $y$, because only one
 simple root of that chamber besides $x$ is orthogonal to $z$.
\end{proof}

\section{The classification theorem}
\label{sec-rank-3}

Here is the main theorem of the paper; the list of lattices appears at
the end of the paper.

\begin{theorem}
\label{thm-main-thm}
Up to scale, there are \NumLattices{} lattices of signature $(2,1)$
whose isometry groups are generated up to finite index by reflections.
\end{theorem}

\begin{proof}
Suppose $L$ is such a lattice, with Weyl group $W$.  Since its Weyl
chamber $C$ has finitely many edges, theorem~\ref{thmdef-shape-of-polygons} applies, so that
$C$ has one of the three features described there: a short edge that is
non-orthogonal to at least one of its neighbors (which are
non-orthogonal to each other), at least $5$ edges and a short pair, or
at least $6$ edges and a
close pair.  Depending on which of those conditions holds, one of
lemmas~\ref{lem-shape-of-2d-SE-integral-polygon}--\ref{lem-shape-of-2D-CP-integral-polygon} applies, with the quasiroot weights
$\rho,\dots,\rho'$ all $2$ because roots are quasiroots of weight~$2$.
The conclusion is that $L$ has $3$, $4$ or $5$ consecutive simple
roots $r,\dots,r'$ whose inner product matrix is one of finitely many
possibilities, up to scale.  These possibilities are listed explicitly
there, and we enumerated them using a computer: the
numbers of matrices of the forms
\eqref{eq-shape-of-2d-SE-integral-polygon-NONorthogonality-case},
\eqref{eq-shape-of-2d-SE-integral-polygon-orthogonality-case},
\eqref{eq-shape-of-2d-SP-integral-polygon-NONorthogonality-case},
\eqref{eq-shape-of-2d-SP-integral-polygon-orthogonality-case} and
\eqref{eq-shape-of-2d-CP-integral-polygon} are
\SENOCNumNarrowMatrices{}, \SEOCNumNarrowMatrices{},
\SPNOCNumNarrowMatrices{}, \SPOCNumNarrowMatrices{} and
\CPNumNarrowMatrices{} respectively.

Some of these matrices are not really possibilities: of those of
form \eqref{eq-shape-of-2d-SE-integral-polygon-NONorthogonality-case},
\SENOCNumNarrowMatricesRejectedRankTwo{} has only rank~$2$.  Of those
of form \eqref{eq-shape-of-2d-SE-integral-polygon-orthogonality-case},
\SEOCNumNarrowMatricesRejectedRankTwo{} have rank~$2$ and
\SEOCNumNarrowMatricesRejectedPositiveDefinite{} are
positive-definite.  Matrices of the form
\eqref{eq-shape-of-2d-SP-integral-polygon-NONorthogonality-case} are
$4\times4$, so the span of $r,\dots,r'$ is obtained by quotienting
$\Z^4$ by the kernel of the inner product.  After
quotienting, $r,\dots,r'$ may no longer be primitive vectors.  When
this happens we reject the matrix because roots are always primitive.
In this way we reject \SPNOCNumNarrowMatricesRejectedImprimitiveRoot{}
matrices, and similarly we reject
\SPOCNumNarrowMatricesRejectedImprimitiveRoot{} and
\CPNumNarrowMatricesRejectedImprimitiveRoot{} of the forms
\eqref{eq-shape-of-2d-SP-integral-polygon-orthogonality-case} and
\eqref{eq-shape-of-2d-CP-integral-polygon}.  This leaves
\TotalNumNarrowMatricesRemaining{} possible inner product matrices
for
$r,\dots,r'$, up to scale.

For $M$ any one of of these matrices, say size $n\times n$, we regard
$M$ as an inner product on $\Z^n$ and define $L_M$ as the quotient by
the kernel of the form, equipped with the induced inner product.  Then
$r,\dots,r'$ are the images in this quotient of the standard basis
vectors of $\Z^n$ and are roots of $L_M$.  Let $E$ vary over all lattices that lie between
the span of these roots and their reflective hull (defined in section~\ref{sec-background}).  There are only
finitely many such enlargements of $L_M$, because $E$ corresponds to a
subgroup of the finite abelian group
$\spanof{r,\dots,r'}^\rh/\spanof{r,\dots,r'}$.  We enumerated the
possible pairs $(M,E)$ with a computer: a total of
\TotalNumEnlargements{}.  Because $r,\dots,r'$ are roots of $L$, $L$
lies between their span and their reflective hull, so $(L,r,\dots,r')$
is isometric to one of our candidates $(E,r,\dots,r')$, up to scale.
Because roots are always primitive, we may discard any
$(E,r,\dots,r')$ for which $r,\dots,r'$ are not primitive, leaving a
total of \NumCandidates{} candidate tuples $(E,r,\dots,r')$ to consider.
We rescaled each candidate to make it unscaled.  

Vinberg's algorithm is the essential tool for determining whether a
given lattice is reflective.  But we are not ready to apply it,
because we need to choose a controlling vector $k$ and find a set of
simple roots for its stabilizer $W_k$ in $W$.  A natural choice is to
take $k$ to represent one of the $2$, $3$ or $4$ vertices of $C$
determined by consecutive pairs among $r,\dots,r'$, and one would
expect $W_k$ to be generated by the reflections in those two roots.
But for some of the candidates $(E,r,\dots,r')$, $W_k$ is actually
larger.  For example, the two roots might be simple roots for an $A_2$
root system, but their span might contain additional roots of $L$,
enlarging the $A_2$ to $G_2$.  We discard the candidate if this occurs
at one of the $2$, $3$ or $4$ vertices.  The reason is that
$r,\dots,r'$ are not simple roots of $E$, because some other mirror of
$E$ passes through the interior of the angle at that vertex.  Since we
started with the assumption that $r,\dots,r'$ are simple roots of $L$,
the candidate $(E,r,\dots,r')$ cannot be isometric to
$(L,r,\dots,r')$.  In this manner we discard
\NumCandidatesFailRootSimplicity{} candidates.  This does not much
reduce the number of cases that must be analyzed, but does allow us to
apply Vinberg's algorithm in the remaining cases.

In principle we could have applied Vinberg's algorithm to the
\NumRemainingCandidates{} remaining candidates $(E,r,\dots,r')$.  In
practice we considered the \NumCandidatesSSF{} SSF cases first.
Because SSF lattices may be distinguished by their genera, one can
immediately sort them into isometry classes and keep only one
candidate from each class.  There turn out to be
\NumCandidatesDistinctSSF{} distinct SSF lattices among our
candidates.  To these we applied Vinberg's algorithm, obtaining
Nikulin's 160 lattices \cite{Nikulin-rk3} and simple roots for them.  We checked
that our results agree with his.

We mentioned in section~\ref{sec-background} that the original form of Vinberg's
algorithm (i.e., with negative-norm controlling vector) was enough for
us.  This is because in every one of the SSF cases, at least one of
the $2$, $3$ or $4$ possible starting corners lies in $H^2$,
rather than being an ideal vertex.  In fact one can show a priori that
this must happen: the Weyl chamber of a $3$-dimensional SSF lattice
cannot have consecutive ideal vertices.  But this isn't really needed since
one can just look at the corners in the 
\NumCandidatesDistinctSSF{} cases.

When applying the algorithm to one of our candidates $(E,r,\dots,r')$,
we checked after each new batch of simple roots whether the polygon
defined by the so-far-found simple roots closed up.  If this happened
then $E$ is reflective and we had a complete set of simple roots.
Otherwise, we looked at the corner symbols of $E$ at the so-far-found
corners.  If some corner symbol appears twice then $E$ admits an
automorphism carrying one corner to the other and preserving the
chamber $C$ being studied.  (Note that each corner can give two corner
symbols since its incident edges can be taking in either order.)
If we find an infinite-order automorphism
this way, or finite-order automorphisms with distinct fixed points,
then $C$ has infinite area and $E$ is non-reflective.  Conversely, if
$E$ is non-reflective then eventually enough corners of $C$ will be
found to show this.  With these stopping conditions, the algorithm is
guaranteed to terminate, either proving non-reflectivity or finding a
set of simple roots.

We remark that a nontrivial automorphism of $E$ preserving $C$ must be
an involution, by the argument for \cite[Prop.~3.4]{Scharlau}.  When we found such an
involution before finishing the search for simple roots, we used it to
speed up the search.

To analyze reflectivity in the non-SSF cases, we sorted them in
increasing order of $|\det E\,|$ and proceeded inductively.  Suppose when
studying a particular case $E$ that we have found all
the reflective lattices of smaller $|\det|$ and obtained simple roots
for them.  Since $E$ is not SSF, its $|\det|$ may be reduced by
$p$-filling or (rescaled) $p$-duality, as explained in
section~\ref{sec-background}.  Write $F$ for the resulting lattice.  If
$F$ is not reflective then $E$ cannot be either, since $W(E)\sset
W(F)$.  And if $F$ is reflective then we can determine whether $E$ is
by a much faster method than Vinberg's algorithm (explained below).  

But it is a nontrivial problem to determine whether $F$ is on the list
of known reflective lattices.  We solved this as follows.  If say $x$
and $y$ are consecutive simple roots of $E$'s Weyl chamber, then they
(or scalar multiples of them) are certainly roots of $F$.  These two
roots might not be extendible to a set of simple roots for $F$,
because of phenomena like the $A_2\sset G_2$ inclusion encountered
above.  But we can look in their rational span to obtain a pair of
simple roots for $F$.  Then we have a corner of $F$'s Weyl chamber and
we compute its corner symbol, say $S$.  We compare this corner symbol
with all the corner symbols for all the reflective lattices of $F$'s
determinant.  If $S$ does not appear then $F$ is not reflective and we
are done.  If $S$ does appear then we obtain an isomorphism of $F$
with a known reflective lattice, extending the two known simple roots
to a full set of simple roots.

We have reduced to the case that applying $p$-filling or rescaled
$p$-duality to $E$ yields a reflective lattice $F$ with a known set of
simple roots.  In the $p$-duality case, the isometry groups of $E$ and
$F$ are equal, so $E$ is reflective and we can obtain simple roots for
it by rescaling $F$'s.  So we assume the $p$-filling case.  Then $E$
is a known sublattice of $F$, lying between $F$ and $pF$, and we must
determine whether $E$ is reflective and obtain simple roots for it
when it is.  We did this using what we call the ``method of
bijections'', which we explain in a more general context below.

In this manner we determined all reflective lattices and found simple
roots for them.  There might be some redundancy in the list, but this
can be easily removed using the lists of corner symbols: two
reflective lattices are isometric if and only if they have the same
determinant and share a corner symbol.
\end{proof}

We now explain the ``method of bijections'', referred to in the proof.
The general setup is that $W$ is a Coxeter group with chamber $C$, $X$
is a finite set on which it acts, $x\in X$ and we want to know whether
the stabilizer $W_x$ is generated up to finite index by reflections.
In the proof of the theorem we used the case that $F$ is a Lorentzian
lattice, $W$ its Weyl group, $X$ the set of subgroups of $F/pF$ for
some prime $p$, $x$ is the subgroup corresponding to the sublattice
$E$ of $F$, and $W_x$ is the subgroup of $W(F)$ preserving $E$.  It
contains $W(E)$ because $\aut E\sset\aut(\pfill(E))$.

The method is a minor variation on traditional coset-enumeration
algorithms.  The idea is to build up a
bijection between a set of chambers and the orbit of $x$ in $X$.
Then the union of the chambers is a fundamental domain for $W_x$.  I imagine a little man
walking from chamber to chamber, being careful never to cross a mirror
of $W_x$.  
For any
chamber $D$ we write $w_D$ for the element of $W$ that sends $D$ to $C$.  
Two chambers $D$, $D'$
are $W_x$-equivalent if and only
if $w_D(x)=w_{D'}(x)$.

During the algorithm we maintain a list $\DD$ of chambers, and also
the list $\FF$ of all pairs $(D,f)$ where $D\in\DD$ and $f$ is an facet of
$D$.  We mark each $(D,f)\in\FF$ by one of the following:
\begin{center}
\begin{tabular}{rl}
``unexplored'':& no assertions are made about $(D,f)$;\\

``mirror'':& the reflection $R_f$ across $f$ lies in $W_x$;\\

``interior'':& $R_f\notin W_x$ and $R_f(D)\in\DD$;\\

``matched'':& $R_f\notin W_x$ and $R_f(D)\notin\DD$ but\\
&$w_D\circ R_f(x)=w_{D'}(x)$ for some $D'\in\DD$.
\end{tabular}
\end{center}
In the initial state, $\DD=\{C\}$ and all pairs $(C,f)$ are marked
``unexplored''.

At each iteration of the algorithm we choose an ``unexplored'' pair
$(D,f)\in\FF$ if one exists.  If none does then the algorithm
terminates.  If $R_f\in W_x$ then we mark the pair ``mirror''.  Otherwise, if $R_f(D)\in\DD$ then
we mark the pair ``interior''.  Otherwise, if $w_D\circ
R_f(x)=w_{D'}(x)$ for some $D'\in\DD$ then we mark the pair
``matched''.  Otherwise we mark $(D,f)$ as ``interior'' and
add $D':=R_f(D)$ to $\DD$ and all pairs
$(D',\hbox{facet of $D'$})$ to $\FF$, and mark them 
``unexplored''.   Then we continue to the next iteration.

In the last case of this chain of otherwise's, note that $w_{D'}=w_D\circ R_f$, and by the
conditions under which this case applies, $w_{D'}$ lies in a coset of
$W_x$ not represented by any previous member of $\DD$.  Since there
are only finitely many cosets of $W_x$ in $W$, the algorithm must
terminate.  Furthermore, at termination $R_f(D)$ is $W_x$-equivalent
to some member of $\DD$, for every $(D,f)\in\FF$, or else $R_f(D)$
would get adjoined to $\DD$ and the algorithm would continue.  So
$\DD$ is a complete set of coset representatives for $W_x$ in $W$ and
the union $U=\cup_{D\in\DD}D$ is a fundamental domain for $W_x$.  
Also, whenever the algorithm
adjoins a chamber $D'=R_f(D)$ to $\DD$, it is already known that $R_f\notin
W_x$.  It follows that both $D$, $D'$ lie in one Weyl chamber of
$W_x$, and by induction that all of $U$ does.

At termination,
each $(D,f)$ is marked ``interior'', ``matched'' or ``mirror''.
Each facet marked ``interior'' really is in $U$'s interior, except perhaps
for its boundary.  
If $(D,f)$ is marked ``matched'', then there is a unique $D'\in\DD$
for which $w_D\circ R_f(x)=w_{D'}(x)$, and we define
$\c_{D,f}=w_{D'}^{-1}\circ w_D\circ R_f$.  This carries $R_f(D)$ to
$D'$, and lies in $W_x$ because it preserves $x\in X$.
Our main claim about the algorithm is the following:

\begin{theorem}
\label{thm-method-of-bijections}
$W_x$ is generated by the reflections across those facets $f$ for
which $(D,f)$ is marked ``mirror'' and the $\c_{D,f}$ for which
$(D,f)$ is marked ``matched''.  The latter generate the stabilizer in
$W_x$ of $W_x$'s Weyl chamber.  In particular, $W_x$ is generated by
reflections up to finite index if and only if the group $\Gamma$
generated by all the $\c_{D,f}$ is finite.
\end{theorem}

\begin{proof}
The algorithm only marks $(D,f)$ as ``mirror'' if $R_f\in W_x$, and
all the $\c_{D,f}$ lie in $W_x$ by construction.  So the supposed
generators do in fact lie in $W_x$.  Next we observe that $\c_{D,f}$
preserves the Weyl chamber of $W_x$ that contains $U$; this is just
the fact that $\c_{D,f}$ sends $R_f(D)$ to $D'$, and $R_f(D)$ lies in
the same Weyl chamber of $W_x$ as $D$ because $R_f\notin W_x$.  It
follows that $\Gamma$ preserves this chamber of $W_x$.  Now consider
the polytope $V$ which is the union of the $\Gamma$-translates of $U$;
it lies in this single chamber of $W_x$.  The reflection across any
one of its facets is a $\Gamma$-conjugate of one of the reflections in
our generating set, so lies in $W_x$.  It follows that $V$ is the
whole of the Weyl chamber for $W_x$ that contains it.  All our claims
follow from this.
\end{proof}

In our situation, checking whether $\Gamma$ is finite amounts to
checking whether the $\c_{D,f}$ have a common fixed point in
hyperbolic space, which is
just linear algebra.  When $\Gamma$ is finite, we take the
$\Gamma$-images of the facets marked ``mirror''.  They bound $W_x$'s
chamber, yielding simple roots for $W_x$, which in our case is $W(E)$.

\section{How to read the table}
\label{sec-how-to-read}

At the end of the paper we exhibit the \NumLattices{} unscaled
reflective Lorentzian lattices of rank~3.  Only
\NumPrintedExplicitly{} lattices are displayed explicitly, and the
rest are obtained by some simple operations.  They are grouped
according to the conjugacy classes of their Weyl groups in
$O(2,1;\R)$. The first \NumMainLattices{} correspond to the main
lattices in Nikulin \cite{Nikulin-rk3} and are listed in the same
order.  Then come the Weyl groups of the lattices of Nikulin's table~2
whose Weyl groups have not been listed yet.  The remainder are listed
in no particular order.  Though we do not give details here, we
classified the \NumWeylGroups{} Weyl groups up to commensurability;
they represent \NumCommClasses{} classes.

\def\GeneralInformation{General Information}
\def\ChamberAngles{Chamber Angles}
\def\ChamberIsometries{Chamber Isometries}
\def\Genus{Genus}
\def\Determinant{Discriminant group and determinant}
\def\ImplicitlyPrintedLattices{Implicitly printed lattices}
\def\CornerSymbols{Corner Symbols and description of the lattices}
\def\SharingGenus{Distinct reflective lattices sharing a genus}
\def\InnerProductMatrix{Inner product matrix}
\def\DiagramAutomorphism{Diagram automorphism}
\def\SimpleRoots{Simple roots}
\def\TableCompleteness{Completeness of the table}
\def\ApplyingOperations{Applying filling, duality and mainification}

{\it\GeneralInformation:\/} For each Weyl group, its Euler
characteristic is $-\chi/24$ and the area of its Weyl chamber is
$\pi\chi/12$, where $\chi$ is given in the table.  The largest area is
\LargestArea{} (\WeylsWithLargestArea) and the smallest is
\SmallestArea{} (\WeylsWithSmallestArea).  The number of lattices
having that Weyl group is listed; the largest number is
\MostLatticesInWeyl{} (\WeylsWithMostLattices) and the smallest is
\LeastLatticesInWeyl{} (\WeylsWithLeastLattices).  For each $W_n$ we
have numbered the printed lattices having that Weyl group as
$L_{n.1}, L_{n.2},\dots$

{\it\ChamberAngles:\/}
At the right side
of the page we give the list
of the chamber's angles $\pi/(\hbox{2, 3, 4, 6 or $\infty$})$ in cyclic order.
The largest number of edges is \MostEdges{} (\WeylsWithMostEdges).  Of
the \NumWeylGroups{} different Weyl chambers, \NumRegularWeyl{} are
regular (\RegularWeyls), \NumIdealWeyl{} have all ideal vertices
(\IdealWeyls), \NumRightAngledWeyl{} have all right angles, and
\NumCompactChambers{} are compact.

{\it\ChamberIsometries:\/}
If the chamber has nontrivial
isometry group then the list of angles is followed by
``$\rtimes C_n$'' or ``$\rtimes D_n$''.  The subscript $n$ gives the number of isometries and
we write $D$ (resp.~$C$) when there are (resp.\ are not) reflections in this
finite group.  
Note the distinction between $C_2$ and $D_2$.
The largest $C_n$ appearing is $C_{\LargestCyclic}$
(\WeylsWithLargestCyclic) and the maximal dihedral groups appearing
are $D_{12}$ (\WeylsWithDihedralTwelve) and $D_8$
(\WeylsWithDihedralEight).  It happens very often that the
chamber admits reflections, but in this case the reflections don't
preserve any of the lattices having that Weyl group.  In the dihedral
case, we have also indicated how the mirrors of the reflections meet
the boundary of the chamber.  A vertical line through a digit (i.e.,
$\slashtwo$, $\slashthree$, $\slashinfty$) means that a mirror passes through that vertex, and a vertical line between two digits
means that a mirror bisects the edge
joining the two corresponding vertices.

{\it\Genus:\/}
For each lattice $L$ listed explicitly we give its genus, in one
version of the notation of Conway and Sloane \cite[ch.~15]{SPLAG}.  That is, for
each prime $p$ dividing twice the determinant, we give the $p$-adic
symbol for $L\tensor\Z_p$.  
We have printed the symbol for the $p^0$ constituent
explicitly to allow easy application of $p$-duality and $p$-filling
(see section~\ref{sec-Conway-Sloane}).
For $p=2$ we give their canonical $2$-adic symbol
in the ``compartment'' notation, but see section~\ref{sec-Conway-Sloane} because there is a
minor error in their paper which requires adjusting the meaning of
``canonical''.  
We have also suppressed the initial $\I$ or $\II$ that they
use because it can be read from the displayed $2$-adic symbol, and we
have suppressed the signature because it is always $(2,1)$.  For
example, they would write $\II_{2,1}(2\above{1pt}{1pt}{1}{1})$ in
place of our $1\above{1pt}{1pt}{2}{\II}2\above{1pt}{1pt}{1}{1}$.

The \NumLattices{} reflective lattices represent \NumGenera{} genera,
the largest number of reflective lattices in a genus is
\MostReflectivesInGenus{}, and \NumGeneraWithMostLattices{} genera
contain this many.  When a listed lattice is not the only reflective
lattice in its genus, this is noted in the table and we give
additional information; see ``\SharingGenus''
below.

{\it\Determinant:\/} The structure of the discriminant group $\D(L)$,
and therefore the determinant of $L$, can be read from the genus
symbol.  To do this one simply reads each symbol $q^{\pm k}$ as the
group $(\Z/q)^k$, where $q$ is a prime power, and the string of
symbols as the direct sum of groups.  One ignores the sign, and when
$q$ is even one also ignores any subscripts.  For example, if $L$ has
genus $1_{\II}^{-2}4^{1}_{1}{\cdot}3^{2}{\cdot}5^{-}125^{1}$ then
discriminant group is $(\Z/4)\oplus(\Z/3)^2\oplus(\Z/5)\oplus(\Z/125)\iso(\Z/1500)\oplus(\Z/15)$.
The most-negative determinant among all reflective lattices is
\MostNegativeDet.  But a better measure of complexity considers only
lattices which are minimal in their duality class,
yielding \MostNegativeDetOfMinimalInDualityClass.  The largest prime
factor among the determinants of our lattices is
\LargestPrimeInDet.

{\it\ImplicitlyPrintedLattices:\/} For some $L$'s in the table,
\instructions{\ldots} appears after the genus;
it indicates how to construct some other lattices, not explicitly
printed, from $L$.  A digit indicates the operation of $p$-filling
(see section~\ref{sec-background}) and $m$ indicates the operation of
``mainification''.  This has meaning when $L$ is odd and the $2^0$
Jordan constituent of the $2$-adic lattice $L\tensor\Z_2$ has
dimension $1$ or~$2$.  It means passage to $L$'s even sublattice,
followed by multiplying inner products by $\frac{1}{4}$ or
$\frac{1}{2}$ to obtain an unscaled lattice.  Every automorphism of
$L$ preserves $L$'s mainification.  For rank $3$ lattices, not
admitting mainification is equivalent to being ``main'' in Nikulin's
terminology \cite{Nikulin-rk3}.  We dislike our name
for the operation, because the mainification of a lattice might not be
main.

A string of these symbols indicates that these operations should be
applied in order, but usually the order is not important because the
operations commute, except for mainification and $2$-filling.  Since
all the primes involved are less than~10, we have printed them without
spaces.  So for example \instructions{3,32,3m,m} indicates that $L$'s
line on the table implicitly represents $4$ additional lattices.  One
is obtained by $3$-filling, one by $3$-filling followed by
$2$-filling, one by $3$-filling followed by mainification, and one by
mainification.  A total of \NumPrintedImplicitly{} lattices are
implicitly printed in this way.  We only use this convention when each
lattice obtained by the stated operations has the same Weyl group as
$L$, so that simple roots for the derived lattices are scalar
multiples of $L$'s.  See ``\ApplyingOperations'' below for how to easily
construct the inner product matrix and simple roots for lattices
obtained in this way.

It turns out that all of Nikulin's lattices are implicitly printed in
this fashion.  We have indicated this by for example
\instructions{\hbox{$23\to N_4$},3,2} for $L_{4.1}$.  This means that the
lattice got by applying $2$-filling and then $3$-filling to $L_{4.1}$
is numbered $N=4$ in Nikulin's \cite[table~1]{Nikulin-rk3}.  We replace $N$
by $N'$ to indicate lattices in his table~2.  We prepared our tables
so that our simple roots come in the same order as his.  Also, some
implicitly printed lattices share their genera with other reflective
lattices; we indicate this with an asterisk ``$*$'' and provide some
additional information in ``\SharingGenus'' below.

{\it\CornerSymbols:\/}
At the right we give the corner symbols at all the corners of the
chamber in the same cyclic order used to list the angles of the Weyl
chamber.  When $\aut L$ contains a nontrivial isometry preserving the
chamber, this isometry is always an involution, and we have written
``$(\times2)$'' after the corner symbols.  The full sequence of corner
symbols is then got by concatenating two copies of the ones listed.
By giving the corner symbols we have given enough information to
reconstruct $L$ in several different ways.  For example, the genus
\GenusWithTwoWeyls{} mentioned below contains
$\GenusWithTwoWeylsSecondLattice={}$\GenusWithTwoWeylsSecondCornerSymbols,
and to construct it we find a consecutive
pair of roots, say ${8}\above{1pt}{1pt}{}{2}{256}$ (the
last root followed by the first one).  So $L$ has roots of norms $8$
and~$256$; the subscript says they are orthogonal, and the absence of
a superscript says that $L$ is generated by them and their orthogonal
complement.  Since the determinant is $-8\cdot256$, a generator for
the complement has norm $-1$.  So $L$ has inner product matrix
$\diag[-1,8,256]$.  Refer to section~\ref{sec-corner-symbols} for full
instructions on constructing $L$ from one of its corner symbols.  The
largest root norm among all reflective lattices (resp. those minimal
in their duality classes) is
\LargestRootNorm{} (resp. \LargestRootNormofMinimalInDualityClass).

{\it\SharingGenus:\/}  When one of our listed lattices $L$ shares its
genus with another reflective lattice, we state what the other lattice
is, either as a lattice in $L$'s duality class or by describing it in
terms of another $L_{n.i}$.  We also give all chains of $p$-dualities
that carry $L$ to a lattice isometric to $L$.  This information is
required to reconstitute the full table from what we have printed.
See ``\TableCompleteness'' below.  A few implicitly-listed lattices
also share their genera with other reflective lattices.  The
corresponding information for these lattices is the following:
\par\smallskip
\ImplicitGeneraTable
\smallskip
\noindent
Reflective lattices having the same genus usually have
the same Weyl group.  But the genus \GenusWithTwoWeyls{} contains the
reflective lattices  $\GenusWithTwoWeylsFirstLattice$ and
$\GenusWithTwoWeylsSecondLattice$, which have Weyl groups
\GenusWithTwoWeylsFirstWeyl{} and \GenusWithTwoWeylsSecondWeyl{}.
The only other genus with
this property is the $2$-dual of this one.  

{\it\InnerProductMatrix:\/} This is the first displayed matrix.  We
did not seek bases in which the matrix entries are small, or bases
related to the corners of the chamber.  Rather, they are in Smith
normal form relative to their duals, which is useful when applying
$p$-filling and $p$-duality.  Namely, let $a,b,c$ be the elementary
divisors of the inclusion $L\to L^*$, where $c|b|a$.  These can be
read from the structure of $\D(L)$, which can be read from the genus
as explained under ``\Determinant''.  With respect to
the basis having the displayed inner product matrix, $L^*$ has basis
$(1/a,0,0)$, $(0,1/b,0)$, $(0,0,1/c)$.  Our lattices are unscaled, so
$c$ is always~$1$.

{\it\DiagramAutomorphism:\/} If $L$ admits a diagram automorphism,
that is, an isometry preserving the Weyl chamber (other than
negation), then this is an involution and we display it right after
the inner product matrix.

{\it\SimpleRoots:\/} Next we list the simple roots of $L$.  If $L$
has no diagram automorphism then all the simple roots are listed, in
the order used for the corner symbols.  If $L$ has a diagram
automorphism then only the first half of the simple roots are
displayed and the rest are got by applying the diagram automorphism.

{\it\TableCompleteness:\/} The table of lattices is complete
in the following sense.  There are \NumDualityClasses{} duality
classes of reflective lattices, and each has exactly one member
appearing in the table,
\NumPrintedExplicitly{} explicitly and \NumPrintedImplicitly{}
implicitly.  
The lattices printed explicitly always have minimal (i.e., least negative)
determinant in their duality classes.
The full list of reflective lattices is obtained from
our explicitly printed lattices by applying the stated operations and
then applying all possible rescaled $p$-dualities.

As stated, this is not quite explicit, because some chains of
$p$-dualities may lead back to $L$.  A fancy way to say this is that
$L$'s duality class is a homogeneous space for $(\Z/2)^\Omega$ where
$\Omega$ is the set of primes dividing $\det L$.  To find the other
lattices in the duality class, without duplicates, one must find a
subgroup of $(\Z/2)^\Omega$ acting simply-transitively, i.e., a
complement to the subgroup $T\sset(\Z/2)^\Omega$ acting trivially.  So
one needs to know $T$ in order to reconstruct the full list of
lattices without duplication.  If $L$ is the unique reflective lattice
in its genus (which is almost always the case) then one can determine
whether a given sequence of rescaled $p$-dualities yields a lattice
isometric to $L$ by computing the genus of the result and comparing
that to $L$'s.  (See section~\ref{sec-Conway-Sloane} for the easy
computational details.)  Otherwise, there is no easy way to determine this, and
it is for this reason that we have listed the nontrivial elements of $T$
explicitly when $L$ shares its genus with another reflective lattice.
We say ``elements'' but it turns out that $T$ is never larger than $\Z/2$.

{\it\ApplyingOperations:\/} Because of our choice of basis for the
lattices (see ``\InnerProductMatrix''), applying $p$-filling is easy.
Let $p^k$ be the $p$-part of the first elementary divisor $a$, and
adjoin $(1/p,0,0)$ to $L$.  This has the effect of dividing the first
row and column of the inner product matrix by $p$, and multiplying the
first coordinate of all roots by $p$.  (If $p^k|b$ then we would treat
the second coordinate similarly, but this never happens because the
explicitly printed lattices are all minimal in their duality classes.)
The result is the inner product matrix for $\pfill(L)$, and the simple
roots of $L$ expressed in a basis for $\pfill(L)$.  It may be that
some of $L$'s roots are divisible by $p$ in $\pdual(L)$, so if
necessary we divide by $p$ to make them primitive.  They will be roots
of $\pfill(L)$ since $\aut L\sset\aut\pfill(L)$.  When the filling
operation is indicated in our tables, these roots are a set of simple roots for
$\pfill(L)$.  This is because we only indicate filling operations when
they don't change the Weyl group.  Similar but slightly more
complicated operations treat mainification and rescaled $p$-duality.

\section{The Conway-Sloane genus symbol}
\label{sec-Conway-Sloane}

Genera of integer bilinear forms present some complicated phenomena,
and Conway and Sloane introduced a notation for working with them that
is as simple as possible \cite[ch.~15]{SPLAG}.  In this section we
state how to apply our operations of filling, duality and
mainification in terms of their notation.  We also correct a minor
error in their formulation of the canonical $2$-adic symbol.

{\it $p$-filling:\/} This affects only the $p$-adic symbol.  One
simply takes the Jordan constituent of largest scale $p^{a+2}$ and
replaces $a+2$ by $a$.  (Recall that $p$-filling is only defined if
$a\geq0$.)  The scale of the resulting summand of
$\pfill(L)\tensor\Z_p$ may be the same as that of another constituent
of $L$.  In this case one takes the direct sum of these sublattices of
$\pfill(L)$ to obtain the $p^a$-constituent.  The direct sum is
computed as follows, where we write $q$ for $p^a$:
\begin{align*}
q^{\e_1n_1}\oplus q^{\e_2n_2}&{}= q^{(\e_1\e_2)(n_1+n_2)}\\
\noalign{\noindent if $q$ is odd, while}
q^{\e_1n_1}_{t_1}\oplus q^{\e_2n_2}_{t_2}&{}= q^{(\e_1\e_2)(n_1+n_2)}_{t_1+t_2}\\
q^{\e_1n_1}_{t_1}\oplus q^{\e_2n_2}_{\II}&{}= q^{(\e_1\e_2)(n_1+n_2)}_{t_1}\\
q^{\e_1n_1}_{\II}\oplus q^{\e_2n_2}_{\II}&{}= q^{(\e_1\e_2)(n_1+n_2)}_{\II}
\end{align*}
if $q$ is even.  
Here the $\e_i$ are signs, $n_i\geq0$ are the ranks and $t_i\in\Z/8$
are the oddities.
In short, signs multiply and  ranks and
subscripts add, subject to the special rules involving $\II$.  When
$p=2$, sometimes there is additional oddity fusion or opportunity for
sign walking.  Here are some
examples of $2$-filling:
\begin{align*}
1^1_18^1_7256^1_1{\cdot}1^29^1
&{}\to
1^1_18^1_764^1_1{\cdot}1^29^1
\\
1^{-2}_28^-_3
&{} \to 
1^{-2}_22^-_3 
= 
[1^{-2}2^-]_5
=
[1^22^1]_1
\\ 
1^1_{-1}[4^18^1]_2
&{}\to
1^1_{-1}[4^12^1]_2
=
[1^12^14^1]_1
\\
1^1_14^1_716^1_1
&{}\to
1^1_14^1_7\oplus4^1_1
=
1^1_14^2_0
\\
1^2_{\II}4^1_1
&{}\to
1^2_{\II}\oplus1^1_1
=
1^3_1.
\end{align*}

{\it Rescaling:\/} Multiplying all inner products by a rational number
$s$ has the following effect on the $l$-adic symbols for all primes
$l$ that divide neither the numerator nor denominator of $s$.  The
sign of a constituent $(\hbox{power of $l$})^{\pm n}$  is flipped just
if $n$ is odd and the Legendre symbol $\legendre{s}{l}$ is~$-1$.  If
$l=2$ we must supplement this rule by describing how to alter the
subscript.  A subscript $\II$ is left alone, while a subscript
$t\in\Z/8$ is multiplied by $s$.

{\it Rescaled $p$-duality:\/} Say $q=p^a$ is the largest scale
appearing in the $p$-adic symbol.  One obtains the rescaled $p$-dual
by an alteration that doesn't affect $L\tensor\Z_l$ for $l\neq p$,
followed by rescaling by $p^a$.  So for $l\neq p$ the effect of
rescaled $p$-duality is the same as that of rescaling.  The $p$-adic
symbol is altered by replacing the scale $p^{a_i}$ of each constituent
by $p^{a-a_i}$, leaving superscripts and subscripts alone.  Since one
usually writes the constituents in increasing order of scale, this
also reverses the order of the constituents.  A single example of
$2$-duality should illustrate the process:
$$
1^1_{-1}2^1_{-1}8^1_1{\cdot}1^-3^19^-{\cdot}1^15^125^-
\to
1^1_14^1_{-1}8^1_{-1}{\cdot}1^13^-9^1{\cdot}1^15^125^-.
$$

{\it Mainification:\/}
This is an alteration of $L$ that doesn't affect $L\tensor\Z_l$ for
odd $l$, followed by rescaling by $\frac{1}{4}$ or $\frac{1}{2}$.  So
for odd $l$ mainification has the same effect as rescaling.  Recall
that the mainification is only defined if $L$ is odd with
$2^0$-constituent of rank~$1$ or~$2$.  Referring to (33)
and (34) of \cite[ch. 15]{SPLAG}, this means that the constituent is
one of $1^1_{\pm1}$, $1^{-1}_{\pm3}$, $1^2_{\pm2}$, 
$1^{-2}_{\pm2}$, $1^2_0$ or $1^{-2}_4$.  To find the $2$-adic symbol of
$\main(L)$, begin by altering the $2^0$-constituent as follows:
$$
1^{\pm1}_t{}\to 4^{\pm1}_t
\qquad
1^{\pm2}_{\pm2}{}\to 2^{\pm2}_{\pm2}
\qquad
1^2_0{}\to2^2_{\II}
\qquad
1^{-2}_4{}\to2^{-2}_{\II}.
$$
The scale of the result may be the same as that of another constituent
present, in which case one takes their direct sum as above.  Finally,
one divides every scale by the smallest scale present ($2$ or $4$),
leaving all superscripts and subscripts alone.  This scaling factor is also
the one to use to work out the effect on the odd $l$-adic symbols.
Here are some examples:
\begin{gather*}
1^2_28^1_7
\to
2^2_28^1_7
\to
1^2_24^1_7
\\
1^2_68^-_5{\cdot}1^-3^-9^1
\to
2^2_68^-_5{\cdot}1^-3^-9^1
\to
1^2_64^-_5{\cdot}1^13^19^-
\\
1^1_18^1_7256^1_1{\cdot}1^29^1
\to
4^1_18^1_7256^1_1{\cdot}1^29^1
=
[4^18^1]_0256^1_1{\cdot}1^29^1
\to
[1^12^1]_064^1_1{\cdot}1^29^1
\\
[1^12^-4^1]_3{\cdot}1^{-2}3^-
\to
[2^-4^2]_3{\cdot}1^{-2}3^-
\to
[1^-2^2]_3{\cdot}1^{-2}3^1.
\end{gather*}

\medskip
{\it The canonical $2$-adic symbol:\/} There is an error in the
definition of the canonical $2$-adic symbol by
Conway-Sloane \cite[ch.~15]{SPLAG}.  Distinct $2$-adic symbols can
represent isometric $2$-adic lattices, and Conway and Sloane give a
calculus for reducing every $2$-adic symbol to a canonical form.
First one uses ``sign walking'' to change some of the signs (and
simultaneously some of the oddities), and then one uses ``oddity
fusion''.  ``Compartments'' and ``trains'' refer to subchains of the
$2$-adic symbol that govern these operations. (The compartments are
delimited by brackets $[\cdots]$.)  Conway and Sloane state that one
can use sign walking to convert any $2$-adic symbol into one which has
plus signs everywhere except perhaps the first positive-dimensional
constituent in each train.  This is not true: sign walking would
supposedly take the $2$-dimensional lattice $1^1_12^-_3$ with total
oddity $1+3=4$ to one of the form $1^-_{\cdots} 2^1_{\cdots}$ with
total oddity $0$.  But the oddity of the first (resp. second)
constituent must be $\pm3$ (resp. $\pm1$) by (33)
of \cite[ch.~15]{SPLAG}, so the total oddity cannot be $0$.

So the canonical form we use is slightly more complicated: minus signs
are allowed either on the first positive-dimensional constituent of a train, or on
the second constituent of a compartment of the form $q^1_{\cdots}
r^-_{\cdots}$ with total oddity~$4$ or $q^-_{\cdots} r^-_{\cdots}$
with total oddity~$0$.  Here $q$ and $r$ are consecutive powers of
$2$.  Oddity fusion remains unchanged.  
Note that the second form $[q^-r^-]_0$ can only occur if the
$q$-constituent is the first positive-dimensional constituent of its
train. 
Using \cite[ch.~15, thm.~10]{SPLAG} we checked
that two $2$-adic lattices are isomorphic if and only if their
canonical symbols (in our sense) are identical.

\vfill\eject

\begin{center}
\bf Table of Rank 3 Reflective Lorentzian Lattices
\end{center}
\medskip
\begingroup
\parindent=0pt
\parskip=0pt\small
\raggedbottom
%
\leavevmode\llap{}%
$W_{1\phantom{0}\phantom{0}}$%
\qquad\llap{11} lattices, $\chi=3$%
\hfill%
$\infty24$%
\nopagebreak\smallskip\hrule\nopagebreak\medskip%
%
%
\leavevmode%
${L_{1.1}}$%
{} : {$1\above{1pt}{1pt}{2}{{\rm II}}4\above{1pt}{1pt}{1}{1}$}\spacer%
\instructions{2\rightarrow N_{1}}%
\EasyButWeakLineBreak%
{${4}\above{1pt}{1pt}{1,0}{\infty b}{4}\above{1pt}{1pt}{r}{2}{2}\above{1pt}{1pt}{*}{4}$}%
\nopagebreak\par%
\nopagebreak\par\leavevmode%
{$\left[\!\llap{\phantom{%
\begingroup \smaller\smaller\smaller\begin{tabular}{@{}c@{}}%
0\\0\\0
\end{tabular}\endgroup%
}}\right.$}%
\begingroup \smaller\smaller\smaller\begin{tabular}{@{}c@{}}%
4\\0\\0
\end{tabular}\endgroup%
\kern3pt%
\begingroup \smaller\smaller\smaller\begin{tabular}{@{}c@{}}%
0\\0\\-1
\end{tabular}\endgroup%
\kern3pt%
\begingroup \smaller\smaller\smaller\begin{tabular}{@{}c@{}}%
0\\-1\\-2
\end{tabular}\endgroup%
{$\left.\llap{\phantom{%
\begingroup \smaller\smaller\smaller\begin{tabular}{@{}c@{}}%
0\\0\\0
\end{tabular}\endgroup%
}}\!\right]$}%
\EasyButWeakLineBreak%
{$\left[\!\llap{\phantom{%
\begingroup \smaller\smaller\smaller\begin{tabular}{@{}c@{}}%
0\\0\\0
\end{tabular}\endgroup%
}}\right.$}%
\begingroup \smaller\smaller\smaller\begin{tabular}{@{}c@{}}%
-1\\-2\\0
\end{tabular}\endgroup%
\HardButStrongLineBreak\kern3pt%
\begingroup \smaller\smaller\smaller\begin{tabular}{@{}c@{}}%
1\\0\\0
\end{tabular}\endgroup%
\HardButStrongLineBreak\kern3pt%
\begingroup \smaller\smaller\smaller\begin{tabular}{@{}c@{}}%
0\\2\\-1
\end{tabular}\endgroup%
{$\left.\llap{\phantom{%
\begingroup \smaller\smaller\smaller\begin{tabular}{@{}c@{}}%
0\\0\\0
\end{tabular}\endgroup%
}}\!\right]$}%
%
%
\hbox{}\par\smallskip%
%
%
\leavevmode%
${L_{1.2}}$%
{} : {$1\above{1pt}{1pt}{-2}{2}8\above{1pt}{1pt}{-}{3}$}\spacer%
\instructions{2\rightarrow N'_{1}}%
\EasyButWeakLineBreak%
{${2}\above{1pt}{1pt}{4,3}{\infty a}{8}\above{1pt}{1pt}{s}{2}{4}\above{1pt}{1pt}{*}{4}$}%
\nopagebreak\par%
\nopagebreak\par\leavevmode%
{$\left[\!\llap{\phantom{%
\begingroup \smaller\smaller\smaller\begin{tabular}{@{}c@{}}%
0\\0\\0
\end{tabular}\endgroup%
}}\right.$}%
\begingroup \smaller\smaller\smaller\begin{tabular}{@{}c@{}}%
-168\\16\\32
\end{tabular}\endgroup%
\kern3pt%
\begingroup \smaller\smaller\smaller\begin{tabular}{@{}c@{}}%
16\\-1\\-3
\end{tabular}\endgroup%
\kern3pt%
\begingroup \smaller\smaller\smaller\begin{tabular}{@{}c@{}}%
32\\-3\\-6
\end{tabular}\endgroup%
{$\left.\llap{\phantom{%
\begingroup \smaller\smaller\smaller\begin{tabular}{@{}c@{}}%
0\\0\\0
\end{tabular}\endgroup%
}}\!\right]$}%
\EasyButWeakLineBreak%
{$\left[\!\llap{\phantom{%
\begingroup \smaller\smaller\smaller\begin{tabular}{@{}c@{}}%
0\\0\\0
\end{tabular}\endgroup%
}}\right.$}%
\begingroup \smaller\smaller\smaller\begin{tabular}{@{}c@{}}%
0\\2\\-1
\end{tabular}\endgroup%
\HardButStrongLineBreak\kern3pt%
\begingroup \smaller\smaller\smaller\begin{tabular}{@{}c@{}}%
-1\\-4\\-4
\end{tabular}\endgroup%
\HardButStrongLineBreak\kern3pt%
\begingroup \smaller\smaller\smaller\begin{tabular}{@{}c@{}}%
1\\-2\\6
\end{tabular}\endgroup%
{$\left.\llap{\phantom{%
\begingroup \smaller\smaller\smaller\begin{tabular}{@{}c@{}}%
0\\0\\0
\end{tabular}\endgroup%
}}\!\right]$}%
%
%
\hbox{}\par\smallskip%
%
%
\leavevmode%
${L_{1.3}}$%
{} : {$1\above{1pt}{1pt}{2}{2}8\above{1pt}{1pt}{1}{7}$}\spacer%
\instructions{m}%
\EasyButWeakLineBreak%
{${2}\above{1pt}{1pt}{4,3}{\infty b}{8}\above{1pt}{1pt}{l}{2}{1}\above{1pt}{1pt}{}{4}$}%
\nopagebreak\par%
\nopagebreak\par\leavevmode%
{$\left[\!\llap{\phantom{%
\begingroup \smaller\smaller\smaller\begin{tabular}{@{}c@{}}%
0\\0\\0
\end{tabular}\endgroup%
}}\right.$}%
\begingroup \smaller\smaller\smaller\begin{tabular}{@{}c@{}}%
-328\\32\\56
\end{tabular}\endgroup%
\kern3pt%
\begingroup \smaller\smaller\smaller\begin{tabular}{@{}c@{}}%
32\\-3\\-6
\end{tabular}\endgroup%
\kern3pt%
\begingroup \smaller\smaller\smaller\begin{tabular}{@{}c@{}}%
56\\-6\\-7
\end{tabular}\endgroup%
{$\left.\llap{\phantom{%
\begingroup \smaller\smaller\smaller\begin{tabular}{@{}c@{}}%
0\\0\\0
\end{tabular}\endgroup%
}}\!\right]$}%
\EasyButWeakLineBreak%
{$\left[\!\llap{\phantom{%
\begingroup \smaller\smaller\smaller\begin{tabular}{@{}c@{}}%
0\\0\\0
\end{tabular}\endgroup%
}}\right.$}%
\begingroup \smaller\smaller\smaller\begin{tabular}{@{}c@{}}%
2\\15\\3
\end{tabular}\endgroup%
\HardButStrongLineBreak\kern3pt%
\begingroup \smaller\smaller\smaller\begin{tabular}{@{}c@{}}%
-3\\-24\\-4
\end{tabular}\endgroup%
\HardButStrongLineBreak\kern3pt%
\begingroup \smaller\smaller\smaller\begin{tabular}{@{}c@{}}%
-1\\-7\\-2
\end{tabular}\endgroup%
{$\left.\llap{\phantom{%
\begingroup \smaller\smaller\smaller\begin{tabular}{@{}c@{}}%
0\\0\\0
\end{tabular}\endgroup%
}}\!\right]$}%

\medskip%
%
\leavevmode\llap{}%
$W_{2\phantom{0}\phantom{0}}$%
\qquad\llap{4} lattices, $\chi=2$%
\hfill%
$\infty23$%
\nopagebreak\smallskip\hrule\nopagebreak\medskip%
%
%
\leavevmode%
${L_{2.1}}$%
{} : {$1\above{1pt}{1pt}{-2}{{\rm II}}8\above{1pt}{1pt}{-}{5}$}\spacer%
\instructions{2\rightarrow N_{2}}%
\EasyButWeakLineBreak%
{${2}\above{1pt}{1pt}{4,1}{\infty b}{8}\above{1pt}{1pt}{b}{2}{2}\above{1pt}{1pt}{+}{3}$}%
\nopagebreak\par%
\nopagebreak\par\leavevmode%
{$\left[\!\llap{\phantom{%
\begingroup \smaller\smaller\smaller\begin{tabular}{@{}c@{}}%
0\\0\\0
\end{tabular}\endgroup%
}}\right.$}%
\begingroup \smaller\smaller\smaller\begin{tabular}{@{}c@{}}%
-24\\8\\16
\end{tabular}\endgroup%
\kern3pt%
\begingroup \smaller\smaller\smaller\begin{tabular}{@{}c@{}}%
8\\-2\\-5
\end{tabular}\endgroup%
\kern3pt%
\begingroup \smaller\smaller\smaller\begin{tabular}{@{}c@{}}%
16\\-5\\-10
\end{tabular}\endgroup%
{$\left.\llap{\phantom{%
\begingroup \smaller\smaller\smaller\begin{tabular}{@{}c@{}}%
0\\0\\0
\end{tabular}\endgroup%
}}\!\right]$}%
\EasyButWeakLineBreak%
{$\left[\!\llap{\phantom{%
\begingroup \smaller\smaller\smaller\begin{tabular}{@{}c@{}}%
0\\0\\0
\end{tabular}\endgroup%
}}\right.$}%
\begingroup \smaller\smaller\smaller\begin{tabular}{@{}c@{}}%
0\\2\\-1
\end{tabular}\endgroup%
\HardButStrongLineBreak\kern3pt%
\begingroup \smaller\smaller\smaller\begin{tabular}{@{}c@{}}%
-1\\-4\\0
\end{tabular}\endgroup%
\HardButStrongLineBreak\kern3pt%
\begingroup \smaller\smaller\smaller\begin{tabular}{@{}c@{}}%
1\\-1\\2
\end{tabular}\endgroup%
{$\left.\llap{\phantom{%
\begingroup \smaller\smaller\smaller\begin{tabular}{@{}c@{}}%
0\\0\\0
\end{tabular}\endgroup%
}}\!\right]$}%
%
%
%
%
%
%
%
%
%
%
%
%
%
%

\medskip%
%
\leavevmode\llap{}%
$W_{3\phantom{0}\phantom{0}}$%
\qquad\llap{6} lattices, $\chi=1$%
\hfill%
$426$%
\nopagebreak\smallskip\hrule\nopagebreak\medskip%
%
%
\leavevmode%
${L_{3.1}}$%
{} : {$1\above{1pt}{1pt}{-2}{{\rm II}}4\above{1pt}{1pt}{1}{7}{\cdot}1\above{1pt}{1pt}{2}{}3\above{1pt}{1pt}{-}{}$}\spacer%
\instructions{2\rightarrow N_{3}}%
\EasyButWeakLineBreak%
{${2}\above{1pt}{1pt}{*}{4}{4}\above{1pt}{1pt}{b}{2}{6}\above{1pt}{1pt}{}{6}$}%
\nopagebreak\par%
\nopagebreak\par\leavevmode%
{$\left[\!\llap{\phantom{%
\begingroup \smaller\smaller\smaller\begin{tabular}{@{}c@{}}%
0\\0\\0
\end{tabular}\endgroup%
}}\right.$}%
\begingroup \smaller\smaller\smaller\begin{tabular}{@{}c@{}}%
-516\\36\\72
\end{tabular}\endgroup%
\kern3pt%
\begingroup \smaller\smaller\smaller\begin{tabular}{@{}c@{}}%
36\\-2\\-5
\end{tabular}\endgroup%
\kern3pt%
\begingroup \smaller\smaller\smaller\begin{tabular}{@{}c@{}}%
72\\-5\\-10
\end{tabular}\endgroup%
{$\left.\llap{\phantom{%
\begingroup \smaller\smaller\smaller\begin{tabular}{@{}c@{}}%
0\\0\\0
\end{tabular}\endgroup%
}}\!\right]$}%
\EasyButWeakLineBreak%
{$\left[\!\llap{\phantom{%
\begingroup \smaller\smaller\smaller\begin{tabular}{@{}c@{}}%
0\\0\\0
\end{tabular}\endgroup%
}}\right.$}%
\begingroup \smaller\smaller\smaller\begin{tabular}{@{}c@{}}%
0\\2\\-1
\end{tabular}\endgroup%
\HardButStrongLineBreak\kern3pt%
\begingroup \smaller\smaller\smaller\begin{tabular}{@{}c@{}}%
1\\-2\\8
\end{tabular}\endgroup%
\HardButStrongLineBreak\kern3pt%
\begingroup \smaller\smaller\smaller\begin{tabular}{@{}c@{}}%
-1\\-3\\-6
\end{tabular}\endgroup%
{$\left.\llap{\phantom{%
\begingroup \smaller\smaller\smaller\begin{tabular}{@{}c@{}}%
0\\0\\0
\end{tabular}\endgroup%
}}\!\right]$}%
%
%
%
%
%
%
%
%
%
%
%
%
%
%
%
%
%
%
%
%
%
%

\medskip%
%
\leavevmode\llap{}%
$W_{4\phantom{0}\phantom{0}}$%
\qquad\llap{44} lattices, $\chi=6$%
\hfill%
$\infty222$%
\nopagebreak\smallskip\hrule\nopagebreak\medskip%
%
%
\leavevmode%
${L_{4.1}}$%
{} : {$1\above{1pt}{1pt}{2}{{\rm II}}4\above{1pt}{1pt}{-}{3}{\cdot}1\above{1pt}{1pt}{1}{}3\above{1pt}{1pt}{1}{}9\above{1pt}{1pt}{-}{}$}\spacer%
\instructions{23\rightarrow N_{4},3,2}%
\EasyButWeakLineBreak%
{${12}\above{1pt}{1pt}{3,2}{\infty b}{12}\above{1pt}{1pt}{r}{2}{18}\above{1pt}{1pt}{b}{2}{4}\above{1pt}{1pt}{*}{2}$}%
\nopagebreak\par%
\nopagebreak\par\leavevmode%
{$\left[\!\llap{\phantom{%
\begingroup \smaller\smaller\smaller\begin{tabular}{@{}c@{}}%
0\\0\\0
\end{tabular}\endgroup%
}}\right.$}%
\begingroup \smaller\smaller\smaller\begin{tabular}{@{}c@{}}%
2124\\0\\-108
\end{tabular}\endgroup%
\kern3pt%
\begingroup \smaller\smaller\smaller\begin{tabular}{@{}c@{}}%
0\\-6\\3
\end{tabular}\endgroup%
\kern3pt%
\begingroup \smaller\smaller\smaller\begin{tabular}{@{}c@{}}%
-108\\3\\4
\end{tabular}\endgroup%
{$\left.\llap{\phantom{%
\begingroup \smaller\smaller\smaller\begin{tabular}{@{}c@{}}%
0\\0\\0
\end{tabular}\endgroup%
}}\!\right]$}%
\EasyButWeakLineBreak%
{$\left[\!\llap{\phantom{%
\begingroup \smaller\smaller\smaller\begin{tabular}{@{}c@{}}%
0\\0\\0
\end{tabular}\endgroup%
}}\right.$}%
\begingroup \smaller\smaller\smaller\begin{tabular}{@{}c@{}}%
-1\\-10\\-18
\end{tabular}\endgroup%
\HardButStrongLineBreak\kern3pt%
\begingroup \smaller\smaller\smaller\begin{tabular}{@{}c@{}}%
3\\28\\60
\end{tabular}\endgroup%
\HardButStrongLineBreak\kern3pt%
\begingroup \smaller\smaller\smaller\begin{tabular}{@{}c@{}}%
1\\9\\18
\end{tabular}\endgroup%
\HardButStrongLineBreak\kern3pt%
\begingroup \smaller\smaller\smaller\begin{tabular}{@{}c@{}}%
-1\\-10\\-20
\end{tabular}\endgroup%
{$\left.\llap{\phantom{%
\begingroup \smaller\smaller\smaller\begin{tabular}{@{}c@{}}%
0\\0\\0
\end{tabular}\endgroup%
}}\!\right]$}%
%
%
\hbox{}\par\smallskip%
%
%
\leavevmode%
${L_{4.2}}$%
{} : {$1\above{1pt}{1pt}{2}{6}8\above{1pt}{1pt}{-}{5}{\cdot}1\above{1pt}{1pt}{-}{}3\above{1pt}{1pt}{-}{}9\above{1pt}{1pt}{1}{}$}\spacer%
\instructions{3m,3,2}%
\EasyButWeakLineBreak%
{${6}\above{1pt}{1pt}{12,5}{\infty a}{24}\above{1pt}{1pt}{s}{2}{36}\above{1pt}{1pt}{*}{2}{8}\above{1pt}{1pt}{b}{2}$}%
\nopagebreak\par%
\nopagebreak\par\leavevmode%
{$\left[\!\llap{\phantom{%
\begingroup \smaller\smaller\smaller\begin{tabular}{@{}c@{}}%
0\\0\\0
\end{tabular}\endgroup%
}}\right.$}%
\begingroup \smaller\smaller\smaller\begin{tabular}{@{}c@{}}%
-3096\\-1440\\432
\end{tabular}\endgroup%
\kern3pt%
\begingroup \smaller\smaller\smaller\begin{tabular}{@{}c@{}}%
-1440\\-669\\198
\end{tabular}\endgroup%
\kern3pt%
\begingroup \smaller\smaller\smaller\begin{tabular}{@{}c@{}}%
432\\198\\-49
\end{tabular}\endgroup%
{$\left.\llap{\phantom{%
\begingroup \smaller\smaller\smaller\begin{tabular}{@{}c@{}}%
0\\0\\0
\end{tabular}\endgroup%
}}\!\right]$}%
\EasyButWeakLineBreak%
{$\left[\!\llap{\phantom{%
\begingroup \smaller\smaller\smaller\begin{tabular}{@{}c@{}}%
0\\0\\0
\end{tabular}\endgroup%
}}\right.$}%
\begingroup \smaller\smaller\smaller\begin{tabular}{@{}c@{}}%
-15\\35\\9
\end{tabular}\endgroup%
\HardButStrongLineBreak\kern3pt%
\begingroup \smaller\smaller\smaller\begin{tabular}{@{}c@{}}%
17\\-40\\-12
\end{tabular}\endgroup%
\HardButStrongLineBreak\kern3pt%
\begingroup \smaller\smaller\smaller\begin{tabular}{@{}c@{}}%
31\\-72\\-18
\end{tabular}\endgroup%
\HardButStrongLineBreak\kern3pt%
\begingroup \smaller\smaller\smaller\begin{tabular}{@{}c@{}}%
-5\\12\\4
\end{tabular}\endgroup%
{$\left.\llap{\phantom{%
\begingroup \smaller\smaller\smaller\begin{tabular}{@{}c@{}}%
0\\0\\0
\end{tabular}\endgroup%
}}\!\right]$}%
%
%
\hbox{}\par\smallskip%
%
%
\leavevmode%
${L_{4.3}}$%
{} : {$1\above{1pt}{1pt}{-2}{6}8\above{1pt}{1pt}{1}{1}{\cdot}1\above{1pt}{1pt}{-}{}3\above{1pt}{1pt}{-}{}9\above{1pt}{1pt}{1}{}$}\spacer%
\instructions{32\rightarrow N'_{2},3,m}%
\EasyButWeakLineBreak%
{${6}\above{1pt}{1pt}{12,5}{\infty b}{24}\above{1pt}{1pt}{l}{2}{9}\above{1pt}{1pt}{}{2}{8}\above{1pt}{1pt}{r}{2}$}%
\nopagebreak\par%
\nopagebreak\par\leavevmode%
{$\left[\!\llap{\phantom{%
\begingroup \smaller\smaller\smaller\begin{tabular}{@{}c@{}}%
0\\0\\0
\end{tabular}\endgroup%
}}\right.$}%
\begingroup \smaller\smaller\smaller\begin{tabular}{@{}c@{}}%
-9144\\-1800\\-1080
\end{tabular}\endgroup%
\kern3pt%
\begingroup \smaller\smaller\smaller\begin{tabular}{@{}c@{}}%
-1800\\-354\\-213
\end{tabular}\endgroup%
\kern3pt%
\begingroup \smaller\smaller\smaller\begin{tabular}{@{}c@{}}%
-1080\\-213\\-127
\end{tabular}\endgroup%
{$\left.\llap{\phantom{%
\begingroup \smaller\smaller\smaller\begin{tabular}{@{}c@{}}%
0\\0\\0
\end{tabular}\endgroup%
}}\!\right]$}%
\EasyButWeakLineBreak%
{$\left[\!\llap{\phantom{%
\begingroup \smaller\smaller\smaller\begin{tabular}{@{}c@{}}%
0\\0\\0
\end{tabular}\endgroup%
}}\right.$}%
\begingroup \smaller\smaller\smaller\begin{tabular}{@{}c@{}}%
1\\-5\\0
\end{tabular}\endgroup%
\HardButStrongLineBreak\kern3pt%
\begingroup \smaller\smaller\smaller\begin{tabular}{@{}c@{}}%
3\\-8\\-12
\end{tabular}\endgroup%
\HardButStrongLineBreak\kern3pt%
\begingroup \smaller\smaller\smaller\begin{tabular}{@{}c@{}}%
-4\\15\\9
\end{tabular}\endgroup%
\HardButStrongLineBreak\kern3pt%
\begingroup \smaller\smaller\smaller\begin{tabular}{@{}c@{}}%
-5\\16\\16
\end{tabular}\endgroup%
{$\left.\llap{\phantom{%
\begingroup \smaller\smaller\smaller\begin{tabular}{@{}c@{}}%
0\\0\\0
\end{tabular}\endgroup%
}}\!\right]$}%

\medskip%
%
\leavevmode\llap{}%
$W_{5\phantom{0}\phantom{0}}$%
\qquad\llap{22} lattices, $\chi=9$%
\hfill%
$\infty242$%
\nopagebreak\smallskip\hrule\nopagebreak\medskip%
%
%
\leavevmode%
${L_{5.1}}$%
{} : {$1\above{1pt}{1pt}{2}{{\rm II}}4\above{1pt}{1pt}{-}{5}{\cdot}1\above{1pt}{1pt}{2}{}5\above{1pt}{1pt}{1}{}$}\spacer%
\instructions{2\rightarrow N_{5}}%
\EasyButWeakLineBreak%
{${20}\above{1pt}{1pt}{1,0}{\infty b}{20}\above{1pt}{1pt}{r}{2}{2}\above{1pt}{1pt}{*}{4}{4}\above{1pt}{1pt}{*}{2}$}%
\nopagebreak\par%
\nopagebreak\par\leavevmode%
{$\left[\!\llap{\phantom{%
\begingroup \smaller\smaller\smaller\begin{tabular}{@{}c@{}}%
0\\0\\0
\end{tabular}\endgroup%
}}\right.$}%
\begingroup \smaller\smaller\smaller\begin{tabular}{@{}c@{}}%
-15020\\360\\680
\end{tabular}\endgroup%
\kern3pt%
\begingroup \smaller\smaller\smaller\begin{tabular}{@{}c@{}}%
360\\-8\\-17
\end{tabular}\endgroup%
\kern3pt%
\begingroup \smaller\smaller\smaller\begin{tabular}{@{}c@{}}%
680\\-17\\-30
\end{tabular}\endgroup%
{$\left.\llap{\phantom{%
\begingroup \smaller\smaller\smaller\begin{tabular}{@{}c@{}}%
0\\0\\0
\end{tabular}\endgroup%
}}\!\right]$}%
\EasyButWeakLineBreak%
{$\left[\!\llap{\phantom{%
\begingroup \smaller\smaller\smaller\begin{tabular}{@{}c@{}}%
0\\0\\0
\end{tabular}\endgroup%
}}\right.$}%
\begingroup \smaller\smaller\smaller\begin{tabular}{@{}c@{}}%
-3\\-50\\-40
\end{tabular}\endgroup%
\HardButStrongLineBreak\kern3pt%
\begingroup \smaller\smaller\smaller\begin{tabular}{@{}c@{}}%
13\\200\\180
\end{tabular}\endgroup%
\HardButStrongLineBreak\kern3pt%
\begingroup \smaller\smaller\smaller\begin{tabular}{@{}c@{}}%
2\\32\\27
\end{tabular}\endgroup%
\HardButStrongLineBreak\kern3pt%
\begingroup \smaller\smaller\smaller\begin{tabular}{@{}c@{}}%
-3\\-46\\-42
\end{tabular}\endgroup%
{$\left.\llap{\phantom{%
\begingroup \smaller\smaller\smaller\begin{tabular}{@{}c@{}}%
0\\0\\0
\end{tabular}\endgroup%
}}\!\right]$}%
%
%
\hbox{}\par\smallskip%
%
%
\leavevmode%
${L_{5.2}}$%
{} : {$1\above{1pt}{1pt}{2}{2}8\above{1pt}{1pt}{-}{3}{\cdot}1\above{1pt}{1pt}{2}{}5\above{1pt}{1pt}{-}{}$}\spacer%
\instructions{2\rightarrow N'_{5}}%
\EasyButWeakLineBreak%
{${10}\above{1pt}{1pt}{4,3}{\infty a}{40}\above{1pt}{1pt}{l}{2}{1}\above{1pt}{1pt}{}{4}{2}\above{1pt}{1pt}{b}{2}$}%
\nopagebreak\par%
\nopagebreak\par\leavevmode%
{$\left[\!\llap{\phantom{%
\begingroup \smaller\smaller\smaller\begin{tabular}{@{}c@{}}%
0\\0\\0
\end{tabular}\endgroup%
}}\right.$}%
\begingroup \smaller\smaller\smaller\begin{tabular}{@{}c@{}}%
-8360\\240\\360
\end{tabular}\endgroup%
\kern3pt%
\begingroup \smaller\smaller\smaller\begin{tabular}{@{}c@{}}%
240\\-6\\-11
\end{tabular}\endgroup%
\kern3pt%
\begingroup \smaller\smaller\smaller\begin{tabular}{@{}c@{}}%
360\\-11\\-15
\end{tabular}\endgroup%
{$\left.\llap{\phantom{%
\begingroup \smaller\smaller\smaller\begin{tabular}{@{}c@{}}%
0\\0\\0
\end{tabular}\endgroup%
}}\!\right]$}%
\EasyButWeakLineBreak%
{$\left[\!\llap{\phantom{%
\begingroup \smaller\smaller\smaller\begin{tabular}{@{}c@{}}%
0\\0\\0
\end{tabular}\endgroup%
}}\right.$}%
\begingroup \smaller\smaller\smaller\begin{tabular}{@{}c@{}}%
-2\\-25\\-30
\end{tabular}\endgroup%
\HardButStrongLineBreak\kern3pt%
\begingroup \smaller\smaller\smaller\begin{tabular}{@{}c@{}}%
9\\100\\140
\end{tabular}\endgroup%
\HardButStrongLineBreak\kern3pt%
\begingroup \smaller\smaller\smaller\begin{tabular}{@{}c@{}}%
1\\12\\15
\end{tabular}\endgroup%
\HardButStrongLineBreak\kern3pt%
\begingroup \smaller\smaller\smaller\begin{tabular}{@{}c@{}}%
-1\\-11\\-16
\end{tabular}\endgroup%
{$\left.\llap{\phantom{%
\begingroup \smaller\smaller\smaller\begin{tabular}{@{}c@{}}%
0\\0\\0
\end{tabular}\endgroup%
}}\!\right]$}%
%
%
\hbox{}\par\smallskip%
%
%
\leavevmode%
${L_{5.3}}$%
{} : {$1\above{1pt}{1pt}{-2}{2}8\above{1pt}{1pt}{1}{7}{\cdot}1\above{1pt}{1pt}{2}{}5\above{1pt}{1pt}{-}{}$}\spacer%
\instructions{m}%
\EasyButWeakLineBreak%
{${10}\above{1pt}{1pt}{4,3}{\infty b}{40}\above{1pt}{1pt}{s}{2}{4}\above{1pt}{1pt}{*}{4}{2}\above{1pt}{1pt}{s}{2}$}%
\nopagebreak\par%
\nopagebreak\par\leavevmode%
{$\left[\!\llap{\phantom{%
\begingroup \smaller\smaller\smaller\begin{tabular}{@{}c@{}}%
0\\0\\0
\end{tabular}\endgroup%
}}\right.$}%
\begingroup \smaller\smaller\smaller\begin{tabular}{@{}c@{}}%
-34760\\200\\1560
\end{tabular}\endgroup%
\kern3pt%
\begingroup \smaller\smaller\smaller\begin{tabular}{@{}c@{}}%
200\\-1\\-9
\end{tabular}\endgroup%
\kern3pt%
\begingroup \smaller\smaller\smaller\begin{tabular}{@{}c@{}}%
1560\\-9\\-70
\end{tabular}\endgroup%
{$\left.\llap{\phantom{%
\begingroup \smaller\smaller\smaller\begin{tabular}{@{}c@{}}%
0\\0\\0
\end{tabular}\endgroup%
}}\!\right]$}%
\EasyButWeakLineBreak%
{$\left[\!\llap{\phantom{%
\begingroup \smaller\smaller\smaller\begin{tabular}{@{}c@{}}%
0\\0\\0
\end{tabular}\endgroup%
}}\right.$}%
\begingroup \smaller\smaller\smaller\begin{tabular}{@{}c@{}}%
-2\\0\\-45
\end{tabular}\endgroup%
\HardButStrongLineBreak\kern3pt%
\begingroup \smaller\smaller\smaller\begin{tabular}{@{}c@{}}%
-1\\-20\\-20
\end{tabular}\endgroup%
\HardButStrongLineBreak\kern3pt%
\begingroup \smaller\smaller\smaller\begin{tabular}{@{}c@{}}%
1\\2\\22
\end{tabular}\endgroup%
\HardButStrongLineBreak\kern3pt%
\begingroup \smaller\smaller\smaller\begin{tabular}{@{}c@{}}%
0\\6\\-1
\end{tabular}\endgroup%
{$\left.\llap{\phantom{%
\begingroup \smaller\smaller\smaller\begin{tabular}{@{}c@{}}%
0\\0\\0
\end{tabular}\endgroup%
}}\!\right]$}%

\medskip%
%
\leavevmode\llap{}%
$W_{6\phantom{0}\phantom{0}}$%
\qquad\llap{6} lattices, $\chi=2$%
\hfill%
$2223$%
\nopagebreak\smallskip\hrule\nopagebreak\medskip%
%
%
\leavevmode%
${L_{6.1}}$%
{} : {$1\above{1pt}{1pt}{-2}{{\rm II}}4\above{1pt}{1pt}{1}{1}{\cdot}1\above{1pt}{1pt}{-2}{}5\above{1pt}{1pt}{-}{}$}\spacer%
\instructions{2\rightarrow N_{6}}%
\EasyButWeakLineBreak%
{${2}\above{1pt}{1pt}{b}{2}{10}\above{1pt}{1pt}{l}{2}{4}\above{1pt}{1pt}{r}{2}{2}\above{1pt}{1pt}{+}{3}$}%
\nopagebreak\par%
\nopagebreak\par\leavevmode%
{$\left[\!\llap{\phantom{%
\begingroup \smaller\smaller\smaller\begin{tabular}{@{}c@{}}%
0\\0\\0
\end{tabular}\endgroup%
}}\right.$}%
\begingroup \smaller\smaller\smaller\begin{tabular}{@{}c@{}}%
-1660\\80\\100
\end{tabular}\endgroup%
\kern3pt%
\begingroup \smaller\smaller\smaller\begin{tabular}{@{}c@{}}%
80\\-2\\-5
\end{tabular}\endgroup%
\kern3pt%
\begingroup \smaller\smaller\smaller\begin{tabular}{@{}c@{}}%
100\\-5\\-6
\end{tabular}\endgroup%
{$\left.\llap{\phantom{%
\begingroup \smaller\smaller\smaller\begin{tabular}{@{}c@{}}%
0\\0\\0
\end{tabular}\endgroup%
}}\!\right]$}%
\EasyButWeakLineBreak%
{$\left[\!\llap{\phantom{%
\begingroup \smaller\smaller\smaller\begin{tabular}{@{}c@{}}%
0\\0\\0
\end{tabular}\endgroup%
}}\right.$}%
\begingroup \smaller\smaller\smaller\begin{tabular}{@{}c@{}}%
-1\\-1\\-16
\end{tabular}\endgroup%
\HardButStrongLineBreak\kern3pt%
\begingroup \smaller\smaller\smaller\begin{tabular}{@{}c@{}}%
-2\\-5\\-30
\end{tabular}\endgroup%
\HardButStrongLineBreak\kern3pt%
\begingroup \smaller\smaller\smaller\begin{tabular}{@{}c@{}}%
1\\0\\16
\end{tabular}\endgroup%
\HardButStrongLineBreak\kern3pt%
\begingroup \smaller\smaller\smaller\begin{tabular}{@{}c@{}}%
1\\2\\15
\end{tabular}\endgroup%
{$\left.\llap{\phantom{%
\begingroup \smaller\smaller\smaller\begin{tabular}{@{}c@{}}%
0\\0\\0
\end{tabular}\endgroup%
}}\!\right]$}%

\medskip%
%
\leavevmode\llap{}%
$W_{7\phantom{0}\phantom{0}}$%
\qquad\llap{16} lattices, $\chi=4$%
\hfill%
$\infty26$%
\nopagebreak\smallskip\hrule\nopagebreak\medskip%
%
%
\leavevmode%
${L_{7.1}}$%
{} : {$1\above{1pt}{1pt}{-2}{{\rm II}}8\above{1pt}{1pt}{1}{7}{\cdot}1\above{1pt}{1pt}{1}{}3\above{1pt}{1pt}{-}{}9\above{1pt}{1pt}{-}{}$}\spacer%
\instructions{23\rightarrow N_{7},3,2}%
\EasyButWeakLineBreak%
{${6}\above{1pt}{1pt}{12,7}{\infty a}{24}\above{1pt}{1pt}{b}{2}{18}\above{1pt}{1pt}{}{6}$}%
\nopagebreak\par%
\nopagebreak\par\leavevmode%
{$\left[\!\llap{\phantom{%
\begingroup \smaller\smaller\smaller\begin{tabular}{@{}c@{}}%
0\\0\\0
\end{tabular}\endgroup%
}}\right.$}%
\begingroup \smaller\smaller\smaller\begin{tabular}{@{}c@{}}%
-1224\\360\\-144
\end{tabular}\endgroup%
\kern3pt%
\begingroup \smaller\smaller\smaller\begin{tabular}{@{}c@{}}%
360\\-102\\39
\end{tabular}\endgroup%
\kern3pt%
\begingroup \smaller\smaller\smaller\begin{tabular}{@{}c@{}}%
-144\\39\\-14
\end{tabular}\endgroup%
{$\left.\llap{\phantom{%
\begingroup \smaller\smaller\smaller\begin{tabular}{@{}c@{}}%
0\\0\\0
\end{tabular}\endgroup%
}}\!\right]$}%
\EasyButWeakLineBreak%
{$\left[\!\llap{\phantom{%
\begingroup \smaller\smaller\smaller\begin{tabular}{@{}c@{}}%
0\\0\\0
\end{tabular}\endgroup%
}}\right.$}%
\begingroup \smaller\smaller\smaller\begin{tabular}{@{}c@{}}%
-1\\-7\\-9
\end{tabular}\endgroup%
\HardButStrongLineBreak\kern3pt%
\begingroup \smaller\smaller\smaller\begin{tabular}{@{}c@{}}%
3\\20\\24
\end{tabular}\endgroup%
\HardButStrongLineBreak\kern3pt%
\begingroup \smaller\smaller\smaller\begin{tabular}{@{}c@{}}%
-1\\-3\\0
\end{tabular}\endgroup%
{$\left.\llap{\phantom{%
\begingroup \smaller\smaller\smaller\begin{tabular}{@{}c@{}}%
0\\0\\0
\end{tabular}\endgroup%
}}\!\right]$}%

\medskip%
%
\leavevmode\llap{}%
$W_{8\phantom{0}\phantom{0}}$%
\qquad\llap{6} lattices, $\chi=6$%
\hfill%
$4242\rtimes C_{2}$%
\nopagebreak\smallskip\hrule\nopagebreak\medskip%
%
%
\leavevmode%
${L_{8.1}}$%
{} : {$1\above{1pt}{1pt}{-2}{{\rm II}}4\above{1pt}{1pt}{-}{3}{\cdot}1\above{1pt}{1pt}{2}{}7\above{1pt}{1pt}{-}{}$}\spacer%
\instructions{2\rightarrow N_{8}}%
\EasyButWeakLineBreak%
{${2}\above{1pt}{1pt}{*}{4}{4}\above{1pt}{1pt}{b}{2}$}\relax$\,(\times2)$%
\nopagebreak\par%
\nopagebreak\par\leavevmode%
{$\left[\!\llap{\phantom{%
\begingroup \smaller\smaller\smaller\begin{tabular}{@{}c@{}}%
0\\0\\0
\end{tabular}\endgroup%
}}\right.$}%
\begingroup \smaller\smaller\smaller\begin{tabular}{@{}c@{}}%
-308\\56\\28
\end{tabular}\endgroup%
\kern3pt%
\begingroup \smaller\smaller\smaller\begin{tabular}{@{}c@{}}%
56\\-10\\-5
\end{tabular}\endgroup%
\kern3pt%
\begingroup \smaller\smaller\smaller\begin{tabular}{@{}c@{}}%
28\\-5\\-2
\end{tabular}\endgroup%
{$\left.\llap{\phantom{%
\begingroup \smaller\smaller\smaller\begin{tabular}{@{}c@{}}%
0\\0\\0
\end{tabular}\endgroup%
}}\!\right]$}%
\hfil\penalty500%
{$\left[\!\llap{\phantom{%
\begingroup \smaller\smaller\smaller\begin{tabular}{@{}c@{}}%
0\\0\\0
\end{tabular}\endgroup%
}}\right.$}%
\begingroup \smaller\smaller\smaller\begin{tabular}{@{}c@{}}%
-1\\0\\-28
\end{tabular}\endgroup%
\kern3pt%
\begingroup \smaller\smaller\smaller\begin{tabular}{@{}c@{}}%
0\\-1\\5
\end{tabular}\endgroup%
\kern3pt%
\begingroup \smaller\smaller\smaller\begin{tabular}{@{}c@{}}%
0\\0\\1
\end{tabular}\endgroup%
{$\left.\llap{\phantom{%
\begingroup \smaller\smaller\smaller\begin{tabular}{@{}c@{}}%
0\\0\\0
\end{tabular}\endgroup%
}}\!\right]$}%
\EasyButWeakLineBreak%
{$\left[\!\llap{\phantom{%
\begingroup \smaller\smaller\smaller\begin{tabular}{@{}c@{}}%
0\\0\\0
\end{tabular}\endgroup%
}}\right.$}%
\begingroup \smaller\smaller\smaller\begin{tabular}{@{}c@{}}%
0\\1\\-3
\end{tabular}\endgroup%
\HardButStrongLineBreak\kern3pt%
\begingroup \smaller\smaller\smaller\begin{tabular}{@{}c@{}}%
-1\\-6\\0
\end{tabular}\endgroup%
{$\left.\llap{\phantom{%
\begingroup \smaller\smaller\smaller\begin{tabular}{@{}c@{}}%
0\\0\\0
\end{tabular}\endgroup%
}}\!\right]$}%

\medskip%
%
\leavevmode\llap{}%
$W_{9\phantom{0}\phantom{0}}$%
\qquad\llap{22} lattices, $\chi=12$%
\hfill%
$\infty2222$%
\nopagebreak\smallskip\hrule\nopagebreak\medskip%
%
%
\leavevmode%
${L_{9.1}}$%
{} : {$1\above{1pt}{1pt}{2}{{\rm II}}4\above{1pt}{1pt}{1}{7}{\cdot}1\above{1pt}{1pt}{-2}{}7\above{1pt}{1pt}{1}{}$}\spacer%
\instructions{2\rightarrow N_{9}}%
\EasyButWeakLineBreak%
{${28}\above{1pt}{1pt}{1,0}{\infty b}{28}\above{1pt}{1pt}{r}{2}{2}\above{1pt}{1pt}{s}{2}{14}\above{1pt}{1pt}{b}{2}{4}\above{1pt}{1pt}{*}{2}$}%
\nopagebreak\par%
\nopagebreak\par\leavevmode%
{$\left[\!\llap{\phantom{%
\begingroup \smaller\smaller\smaller\begin{tabular}{@{}c@{}}%
0\\0\\0
\end{tabular}\endgroup%
}}\right.$}%
\begingroup \smaller\smaller\smaller\begin{tabular}{@{}c@{}}%
-21252\\532\\1204
\end{tabular}\endgroup%
\kern3pt%
\begingroup \smaller\smaller\smaller\begin{tabular}{@{}c@{}}%
532\\-12\\-33
\end{tabular}\endgroup%
\kern3pt%
\begingroup \smaller\smaller\smaller\begin{tabular}{@{}c@{}}%
1204\\-33\\-62
\end{tabular}\endgroup%
{$\left.\llap{\phantom{%
\begingroup \smaller\smaller\smaller\begin{tabular}{@{}c@{}}%
0\\0\\0
\end{tabular}\endgroup%
}}\!\right]$}%
\EasyButWeakLineBreak%
{$\left[\!\llap{\phantom{%
\begingroup \smaller\smaller\smaller\begin{tabular}{@{}c@{}}%
0\\0\\0
\end{tabular}\endgroup%
}}\right.$}%
\begingroup \smaller\smaller\smaller\begin{tabular}{@{}c@{}}%
-31\\-602\\-280
\end{tabular}\endgroup%
\HardButStrongLineBreak\kern3pt%
\begingroup \smaller\smaller\smaller\begin{tabular}{@{}c@{}}%
3\\56\\28
\end{tabular}\endgroup%
\HardButStrongLineBreak\kern3pt%
\begingroup \smaller\smaller\smaller\begin{tabular}{@{}c@{}}%
5\\98\\45
\end{tabular}\endgroup%
\HardButStrongLineBreak\kern3pt%
\begingroup \smaller\smaller\smaller\begin{tabular}{@{}c@{}}%
1\\28\\7
\end{tabular}\endgroup%
\HardButStrongLineBreak\kern3pt%
\begingroup \smaller\smaller\smaller\begin{tabular}{@{}c@{}}%
-9\\-172\\-82
\end{tabular}\endgroup%
{$\left.\llap{\phantom{%
\begingroup \smaller\smaller\smaller\begin{tabular}{@{}c@{}}%
0\\0\\0
\end{tabular}\endgroup%
}}\!\right]$}%
%
%
\hbox{}\par\smallskip%
%
%
\leavevmode%
${L_{9.2}}$%
{} : {$1\above{1pt}{1pt}{2}{6}8\above{1pt}{1pt}{1}{1}{\cdot}1\above{1pt}{1pt}{-2}{}7\above{1pt}{1pt}{1}{}$}\spacer%
\instructions{2\rightarrow N'_{7}}%
\EasyButWeakLineBreak%
{${14}\above{1pt}{1pt}{4,1}{\infty a}{56}\above{1pt}{1pt}{s}{2}{4}\above{1pt}{1pt}{l}{2}{7}\above{1pt}{1pt}{}{2}{8}\above{1pt}{1pt}{r}{2}$}%
\nopagebreak\par%
\nopagebreak\par\leavevmode%
{$\left[\!\llap{\phantom{%
\begingroup \smaller\smaller\smaller\begin{tabular}{@{}c@{}}%
0\\0\\0
\end{tabular}\endgroup%
}}\right.$}%
\begingroup \smaller\smaller\smaller\begin{tabular}{@{}c@{}}%
-39928\\224\\1176
\end{tabular}\endgroup%
\kern3pt%
\begingroup \smaller\smaller\smaller\begin{tabular}{@{}c@{}}%
224\\-1\\-7
\end{tabular}\endgroup%
\kern3pt%
\begingroup \smaller\smaller\smaller\begin{tabular}{@{}c@{}}%
1176\\-7\\-34
\end{tabular}\endgroup%
{$\left.\llap{\phantom{%
\begingroup \smaller\smaller\smaller\begin{tabular}{@{}c@{}}%
0\\0\\0
\end{tabular}\endgroup%
}}\!\right]$}%
\EasyButWeakLineBreak%
{$\left[\!\llap{\phantom{%
\begingroup \smaller\smaller\smaller\begin{tabular}{@{}c@{}}%
0\\0\\0
\end{tabular}\endgroup%
}}\right.$}%
\begingroup \smaller\smaller\smaller\begin{tabular}{@{}c@{}}%
3\\126\\77
\end{tabular}\endgroup%
\HardButStrongLineBreak\kern3pt%
\begingroup \smaller\smaller\smaller\begin{tabular}{@{}c@{}}%
1\\28\\28
\end{tabular}\endgroup%
\HardButStrongLineBreak\kern3pt%
\begingroup \smaller\smaller\smaller\begin{tabular}{@{}c@{}}%
-1\\-42\\-26
\end{tabular}\endgroup%
\HardButStrongLineBreak\kern3pt%
\begingroup \smaller\smaller\smaller\begin{tabular}{@{}c@{}}%
-1\\-35\\-28
\end{tabular}\endgroup%
\HardButStrongLineBreak\kern3pt%
\begingroup \smaller\smaller\smaller\begin{tabular}{@{}c@{}}%
1\\48\\24
\end{tabular}\endgroup%
{$\left.\llap{\phantom{%
\begingroup \smaller\smaller\smaller\begin{tabular}{@{}c@{}}%
0\\0\\0
\end{tabular}\endgroup%
}}\!\right]$}%
%
%
\hbox{}\par\smallskip%
%
%
\leavevmode%
${L_{9.3}}$%
{} : {$1\above{1pt}{1pt}{-2}{6}8\above{1pt}{1pt}{-}{5}{\cdot}1\above{1pt}{1pt}{-2}{}7\above{1pt}{1pt}{1}{}$}\spacer%
\instructions{m}%
\EasyButWeakLineBreak%
{${14}\above{1pt}{1pt}{4,1}{\infty b}{56}\above{1pt}{1pt}{l}{2}{1}\above{1pt}{1pt}{r}{2}{28}\above{1pt}{1pt}{*}{2}{8}\above{1pt}{1pt}{b}{2}$}%
\nopagebreak\par%
\nopagebreak\par\leavevmode%
{$\left[\!\llap{\phantom{%
\begingroup \smaller\smaller\smaller\begin{tabular}{@{}c@{}}%
0\\0\\0
\end{tabular}\endgroup%
}}\right.$}%
\begingroup \smaller\smaller\smaller\begin{tabular}{@{}c@{}}%
-655704\\94472\\11312
\end{tabular}\endgroup%
\kern3pt%
\begingroup \smaller\smaller\smaller\begin{tabular}{@{}c@{}}%
94472\\-13611\\-1630
\end{tabular}\endgroup%
\kern3pt%
\begingroup \smaller\smaller\smaller\begin{tabular}{@{}c@{}}%
11312\\-1630\\-195
\end{tabular}\endgroup%
{$\left.\llap{\phantom{%
\begingroup \smaller\smaller\smaller\begin{tabular}{@{}c@{}}%
0\\0\\0
\end{tabular}\endgroup%
}}\!\right]$}%
\EasyButWeakLineBreak%
{$\left[\!\llap{\phantom{%
\begingroup \smaller\smaller\smaller\begin{tabular}{@{}c@{}}%
0\\0\\0
\end{tabular}\endgroup%
}}\right.$}%
\begingroup \smaller\smaller\smaller\begin{tabular}{@{}c@{}}%
48\\287\\385
\end{tabular}\endgroup%
\HardButStrongLineBreak\kern3pt%
\begingroup \smaller\smaller\smaller\begin{tabular}{@{}c@{}}%
-5\\-28\\-56
\end{tabular}\endgroup%
\HardButStrongLineBreak\kern3pt%
\begingroup \smaller\smaller\smaller\begin{tabular}{@{}c@{}}%
-7\\-42\\-55
\end{tabular}\endgroup%
\HardButStrongLineBreak\kern3pt%
\begingroup \smaller\smaller\smaller\begin{tabular}{@{}c@{}}%
3\\14\\56
\end{tabular}\endgroup%
\HardButStrongLineBreak\kern3pt%
\begingroup \smaller\smaller\smaller\begin{tabular}{@{}c@{}}%
29\\172\\244
\end{tabular}\endgroup%
{$\left.\llap{\phantom{%
\begingroup \smaller\smaller\smaller\begin{tabular}{@{}c@{}}%
0\\0\\0
\end{tabular}\endgroup%
}}\!\right]$}%

\medskip%
%
\leavevmode\llap{}%
$W_{10\phantom{0}}$%
\qquad\llap{16} lattices, $\chi=6$%
\hfill%
$\infty222$%
\nopagebreak\smallskip\hrule\nopagebreak\medskip%
%
%
\leavevmode%
${L_{10.1}}$%
{} : {$1\above{1pt}{1pt}{-2}{{\rm II}}8\above{1pt}{1pt}{1}{1}{\cdot}1\above{1pt}{1pt}{-}{}5\above{1pt}{1pt}{-}{}25\above{1pt}{1pt}{-}{}$}\spacer%
\instructions{25\rightarrow N_{10},5,2*}%
\EasyButWeakLineBreak%
{${10}\above{1pt}{1pt}{20,9}{\infty a}{40}\above{1pt}{1pt}{b}{2}{50}\above{1pt}{1pt}{l}{2}{8}\above{1pt}{1pt}{r}{2}$}%
\nopagebreak\par%
shares genus with 5-dual\nopagebreak\par%
\nopagebreak\par\leavevmode%
{$\left[\!\llap{\phantom{%
\begingroup \smaller\smaller\smaller\begin{tabular}{@{}c@{}}%
0\\0\\0
\end{tabular}\endgroup%
}}\right.$}%
\begingroup \smaller\smaller\smaller\begin{tabular}{@{}c@{}}%
-95800\\1200\\-10800
\end{tabular}\endgroup%
\kern3pt%
\begingroup \smaller\smaller\smaller\begin{tabular}{@{}c@{}}%
1200\\-10\\165
\end{tabular}\endgroup%
\kern3pt%
\begingroup \smaller\smaller\smaller\begin{tabular}{@{}c@{}}%
-10800\\165\\-1042
\end{tabular}\endgroup%
{$\left.\llap{\phantom{%
\begingroup \smaller\smaller\smaller\begin{tabular}{@{}c@{}}%
0\\0\\0
\end{tabular}\endgroup%
}}\!\right]$}%
\EasyButWeakLineBreak%
{$\left[\!\llap{\phantom{%
\begingroup \smaller\smaller\smaller\begin{tabular}{@{}c@{}}%
0\\0\\0
\end{tabular}\endgroup%
}}\right.$}%
\begingroup \smaller\smaller\smaller\begin{tabular}{@{}c@{}}%
-14\\-443\\75
\end{tabular}\endgroup%
\HardButStrongLineBreak\kern3pt%
\begingroup \smaller\smaller\smaller\begin{tabular}{@{}c@{}}%
15\\476\\-80
\end{tabular}\endgroup%
\HardButStrongLineBreak\kern3pt%
\begingroup \smaller\smaller\smaller\begin{tabular}{@{}c@{}}%
28\\885\\-150
\end{tabular}\endgroup%
\HardButStrongLineBreak\kern3pt%
\begingroup \smaller\smaller\smaller\begin{tabular}{@{}c@{}}%
-3\\-96\\16
\end{tabular}\endgroup%
{$\left.\llap{\phantom{%
\begingroup \smaller\smaller\smaller\begin{tabular}{@{}c@{}}%
0\\0\\0
\end{tabular}\endgroup%
}}\!\right]$}%

\medskip%
%
\leavevmode\llap{}%
$W_{11\phantom{0}}$%
\qquad\llap{6} lattices, $\chi=5$%
\hfill%
$4223$%
\nopagebreak\smallskip\hrule\nopagebreak\medskip%
%
%
\leavevmode%
${L_{11.1}}$%
{} : {$1\above{1pt}{1pt}{-2}{{\rm II}}4\above{1pt}{1pt}{1}{7}{\cdot}1\above{1pt}{1pt}{2}{}11\above{1pt}{1pt}{-}{}$}\spacer%
\instructions{2\rightarrow N_{11}}%
\EasyButWeakLineBreak%
{${2}\above{1pt}{1pt}{*}{4}{4}\above{1pt}{1pt}{b}{2}{22}\above{1pt}{1pt}{s}{2}{2}\above{1pt}{1pt}{+}{3}$}%
\nopagebreak\par%
\nopagebreak\par\leavevmode%
{$\left[\!\llap{\phantom{%
\begingroup \smaller\smaller\smaller\begin{tabular}{@{}c@{}}%
0\\0\\0
\end{tabular}\endgroup%
}}\right.$}%
\begingroup \smaller\smaller\smaller\begin{tabular}{@{}c@{}}%
-33924\\748\\1012
\end{tabular}\endgroup%
\kern3pt%
\begingroup \smaller\smaller\smaller\begin{tabular}{@{}c@{}}%
748\\-14\\-23
\end{tabular}\endgroup%
\kern3pt%
\begingroup \smaller\smaller\smaller\begin{tabular}{@{}c@{}}%
1012\\-23\\-30
\end{tabular}\endgroup%
{$\left.\llap{\phantom{%
\begingroup \smaller\smaller\smaller\begin{tabular}{@{}c@{}}%
0\\0\\0
\end{tabular}\endgroup%
}}\!\right]$}%
\EasyButWeakLineBreak%
{$\left[\!\llap{\phantom{%
\begingroup \smaller\smaller\smaller\begin{tabular}{@{}c@{}}%
0\\0\\0
\end{tabular}\endgroup%
}}\right.$}%
\begingroup \smaller\smaller\smaller\begin{tabular}{@{}c@{}}%
2\\16\\55
\end{tabular}\endgroup%
\HardButStrongLineBreak\kern3pt%
\begingroup \smaller\smaller\smaller\begin{tabular}{@{}c@{}}%
3\\22\\84
\end{tabular}\endgroup%
\HardButStrongLineBreak\kern3pt%
\begingroup \smaller\smaller\smaller\begin{tabular}{@{}c@{}}%
-4\\-33\\-110
\end{tabular}\endgroup%
\HardButStrongLineBreak\kern3pt%
\begingroup \smaller\smaller\smaller\begin{tabular}{@{}c@{}}%
-2\\-15\\-56
\end{tabular}\endgroup%
{$\left.\llap{\phantom{%
\begingroup \smaller\smaller\smaller\begin{tabular}{@{}c@{}}%
0\\0\\0
\end{tabular}\endgroup%
}}\!\right]$}%

\medskip%
%
\leavevmode\llap{}%
$W_{12\phantom{0}}$%
\qquad\llap{22} lattices, $\chi=36$%
\hfill%
$\infty222\infty222\rtimes C_{2}$%
\nopagebreak\smallskip\hrule\nopagebreak\medskip%
%
%
\leavevmode%
${L_{12.1}}$%
{} : {$1\above{1pt}{1pt}{2}{{\rm II}}4\above{1pt}{1pt}{-}{3}{\cdot}1\above{1pt}{1pt}{-2}{}11\above{1pt}{1pt}{1}{}$}\spacer%
\instructions{2\rightarrow N_{12}}%
\EasyButWeakLineBreak%
{${44}\above{1pt}{1pt}{1,0}{\infty b}{44}\above{1pt}{1pt}{r}{2}{2}\above{1pt}{1pt}{b}{2}{4}\above{1pt}{1pt}{*}{2}$}\relax$\,(\times2)$%
\nopagebreak\par%
\nopagebreak\par\leavevmode%
{$\left[\!\llap{\phantom{%
\begingroup \smaller\smaller\smaller\begin{tabular}{@{}c@{}}%
0\\0\\0
\end{tabular}\endgroup%
}}\right.$}%
\begingroup \smaller\smaller\smaller\begin{tabular}{@{}c@{}}%
-2348148\\9768\\70268
\end{tabular}\endgroup%
\kern3pt%
\begingroup \smaller\smaller\smaller\begin{tabular}{@{}c@{}}%
9768\\-40\\-293
\end{tabular}\endgroup%
\kern3pt%
\begingroup \smaller\smaller\smaller\begin{tabular}{@{}c@{}}%
70268\\-293\\-2102
\end{tabular}\endgroup%
{$\left.\llap{\phantom{%
\begingroup \smaller\smaller\smaller\begin{tabular}{@{}c@{}}%
0\\0\\0
\end{tabular}\endgroup%
}}\!\right]$}%
\hfil\penalty500%
{$\left[\!\llap{\phantom{%
\begingroup \smaller\smaller\smaller\begin{tabular}{@{}c@{}}%
0\\0\\0
\end{tabular}\endgroup%
}}\right.$}%
\begingroup \smaller\smaller\smaller\begin{tabular}{@{}c@{}}%
42041\\1285284\\1225224
\end{tabular}\endgroup%
\kern3pt%
\begingroup \smaller\smaller\smaller\begin{tabular}{@{}c@{}}%
-154\\-4709\\-4488
\end{tabular}\endgroup%
\kern3pt%
\begingroup \smaller\smaller\smaller\begin{tabular}{@{}c@{}}%
-1281\\-39162\\-37333
\end{tabular}\endgroup%
{$\left.\llap{\phantom{%
\begingroup \smaller\smaller\smaller\begin{tabular}{@{}c@{}}%
0\\0\\0
\end{tabular}\endgroup%
}}\!\right]$}%
\EasyButWeakLineBreak%
{$\left[\!\llap{\phantom{%
\begingroup \smaller\smaller\smaller\begin{tabular}{@{}c@{}}%
0\\0\\0
\end{tabular}\endgroup%
}}\right.$}%
\begingroup \smaller\smaller\smaller\begin{tabular}{@{}c@{}}%
-63\\-1826\\-1848
\end{tabular}\endgroup%
\HardButStrongLineBreak\kern3pt%
\begingroup \smaller\smaller\smaller\begin{tabular}{@{}c@{}}%
-157\\-4796\\-4576
\end{tabular}\endgroup%
\HardButStrongLineBreak\kern3pt%
\begingroup \smaller\smaller\smaller\begin{tabular}{@{}c@{}}%
-16\\-500\\-465
\end{tabular}\endgroup%
\HardButStrongLineBreak\kern3pt%
\begingroup \smaller\smaller\smaller\begin{tabular}{@{}c@{}}%
-33\\-1040\\-958
\end{tabular}\endgroup%
{$\left.\llap{\phantom{%
\begingroup \smaller\smaller\smaller\begin{tabular}{@{}c@{}}%
0\\0\\0
\end{tabular}\endgroup%
}}\!\right]$}%
%
%
\hbox{}\par\smallskip%
%
%
\leavevmode%
${L_{12.2}}$%
{} : {$1\above{1pt}{1pt}{2}{6}8\above{1pt}{1pt}{-}{5}{\cdot}1\above{1pt}{1pt}{-2}{}11\above{1pt}{1pt}{-}{}$}\spacer%
\instructions{2\rightarrow N'_{8}}%
\EasyButWeakLineBreak%
{${22}\above{1pt}{1pt}{4,1}{\infty b}{88}\above{1pt}{1pt}{s}{2}{4}\above{1pt}{1pt}{*}{2}{8}\above{1pt}{1pt}{b}{2}$}\relax$\,(\times2)$%
\nopagebreak\par%
\nopagebreak\par\leavevmode%
{$\left[\!\llap{\phantom{%
\begingroup \smaller\smaller\smaller\begin{tabular}{@{}c@{}}%
0\\0\\0
\end{tabular}\endgroup%
}}\right.$}%
\begingroup \smaller\smaller\smaller\begin{tabular}{@{}c@{}}%
-238040\\528\\4840
\end{tabular}\endgroup%
\kern3pt%
\begingroup \smaller\smaller\smaller\begin{tabular}{@{}c@{}}%
528\\-1\\-11
\end{tabular}\endgroup%
\kern3pt%
\begingroup \smaller\smaller\smaller\begin{tabular}{@{}c@{}}%
4840\\-11\\-98
\end{tabular}\endgroup%
{$\left.\llap{\phantom{%
\begingroup \smaller\smaller\smaller\begin{tabular}{@{}c@{}}%
0\\0\\0
\end{tabular}\endgroup%
}}\!\right]$}%
\hfil\penalty500%
{$\left[\!\llap{\phantom{%
\begingroup \smaller\smaller\smaller\begin{tabular}{@{}c@{}}%
0\\0\\0
\end{tabular}\endgroup%
}}\right.$}%
\begingroup \smaller\smaller\smaller\begin{tabular}{@{}c@{}}%
-1\\10032\\-1672
\end{tabular}\endgroup%
\kern3pt%
\begingroup \smaller\smaller\smaller\begin{tabular}{@{}c@{}}%
0\\-31\\5
\end{tabular}\endgroup%
\kern3pt%
\begingroup \smaller\smaller\smaller\begin{tabular}{@{}c@{}}%
0\\-192\\31
\end{tabular}\endgroup%
{$\left.\llap{\phantom{%
\begingroup \smaller\smaller\smaller\begin{tabular}{@{}c@{}}%
0\\0\\0
\end{tabular}\endgroup%
}}\!\right]$}%
\EasyButWeakLineBreak%
{$\left[\!\llap{\phantom{%
\begingroup \smaller\smaller\smaller\begin{tabular}{@{}c@{}}%
0\\0\\0
\end{tabular}\endgroup%
}}\right.$}%
\begingroup \smaller\smaller\smaller\begin{tabular}{@{}c@{}}%
5\\330\\209
\end{tabular}\endgroup%
\HardButStrongLineBreak\kern3pt%
\begingroup \smaller\smaller\smaller\begin{tabular}{@{}c@{}}%
1\\44\\44
\end{tabular}\endgroup%
\HardButStrongLineBreak\kern3pt%
\begingroup \smaller\smaller\smaller\begin{tabular}{@{}c@{}}%
-1\\-66\\-42
\end{tabular}\endgroup%
\HardButStrongLineBreak\kern3pt%
\begingroup \smaller\smaller\smaller\begin{tabular}{@{}c@{}}%
-3\\-160\\-132
\end{tabular}\endgroup%
{$\left.\llap{\phantom{%
\begingroup \smaller\smaller\smaller\begin{tabular}{@{}c@{}}%
0\\0\\0
\end{tabular}\endgroup%
}}\!\right]$}%
%
%
\hbox{}\par\smallskip%
%
%
\leavevmode%
${L_{12.3}}$%
{} : {$1\above{1pt}{1pt}{-2}{6}8\above{1pt}{1pt}{1}{1}{\cdot}1\above{1pt}{1pt}{-2}{}11\above{1pt}{1pt}{-}{}$}\spacer%
\instructions{m}%
\EasyButWeakLineBreak%
{${22}\above{1pt}{1pt}{4,1}{\infty a}{88}\above{1pt}{1pt}{l}{2}{1}\above{1pt}{1pt}{}{2}{8}\above{1pt}{1pt}{r}{2}$}\relax$\,(\times2)$%
\nopagebreak\par%
\nopagebreak\par\leavevmode%
{$\left[\!\llap{\phantom{%
\begingroup \smaller\smaller\smaller\begin{tabular}{@{}c@{}}%
0\\0\\0
\end{tabular}\endgroup%
}}\right.$}%
\begingroup \smaller\smaller\smaller\begin{tabular}{@{}c@{}}%
-1848\\-1848\\968
\end{tabular}\endgroup%
\kern3pt%
\begingroup \smaller\smaller\smaller\begin{tabular}{@{}c@{}}%
-1848\\-1826\\945
\end{tabular}\endgroup%
\kern3pt%
\begingroup \smaller\smaller\smaller\begin{tabular}{@{}c@{}}%
968\\945\\-483
\end{tabular}\endgroup%
{$\left.\llap{\phantom{%
\begingroup \smaller\smaller\smaller\begin{tabular}{@{}c@{}}%
0\\0\\0
\end{tabular}\endgroup%
}}\!\right]$}%
\hfil\penalty500%
{$\left[\!\llap{\phantom{%
\begingroup \smaller\smaller\smaller\begin{tabular}{@{}c@{}}%
0\\0\\0
\end{tabular}\endgroup%
}}\right.$}%
\begingroup \smaller\smaller\smaller\begin{tabular}{@{}c@{}}%
527\\-1320\\-1584
\end{tabular}\endgroup%
\kern3pt%
\begingroup \smaller\smaller\smaller\begin{tabular}{@{}c@{}}%
472\\-1181\\-1416
\end{tabular}\endgroup%
\kern3pt%
\begingroup \smaller\smaller\smaller\begin{tabular}{@{}c@{}}%
-218\\545\\653
\end{tabular}\endgroup%
{$\left.\llap{\phantom{%
\begingroup \smaller\smaller\smaller\begin{tabular}{@{}c@{}}%
0\\0\\0
\end{tabular}\endgroup%
}}\!\right]$}%
\EasyButWeakLineBreak%
{$\left[\!\llap{\phantom{%
\begingroup \smaller\smaller\smaller\begin{tabular}{@{}c@{}}%
0\\0\\0
\end{tabular}\endgroup%
}}\right.$}%
\begingroup \smaller\smaller\smaller\begin{tabular}{@{}c@{}}%
-105\\209\\198
\end{tabular}\endgroup%
\HardButStrongLineBreak\kern3pt%
\begingroup \smaller\smaller\smaller\begin{tabular}{@{}c@{}}%
-21\\44\\44
\end{tabular}\endgroup%
\HardButStrongLineBreak\kern3pt%
\begingroup \smaller\smaller\smaller\begin{tabular}{@{}c@{}}%
11\\-22\\-21
\end{tabular}\endgroup%
\HardButStrongLineBreak\kern3pt%
\begingroup \smaller\smaller\smaller\begin{tabular}{@{}c@{}}%
77\\-160\\-160
\end{tabular}\endgroup%
{$\left.\llap{\phantom{%
\begingroup \smaller\smaller\smaller\begin{tabular}{@{}c@{}}%
0\\0\\0
\end{tabular}\endgroup%
}}\!\right]$}%

\medskip%
%
\leavevmode\llap{}%
$W_{13\phantom{0}}$%
\qquad\llap{22} lattices, $\chi=42$%
\hfill%
$\infty242\infty242\rtimes C_{2}$%
\nopagebreak\smallskip\hrule\nopagebreak\medskip%
%
%
\leavevmode%
${L_{13.1}}$%
{} : {$1\above{1pt}{1pt}{2}{{\rm II}}4\above{1pt}{1pt}{-}{5}{\cdot}1\above{1pt}{1pt}{2}{}13\above{1pt}{1pt}{1}{}$}\spacer%
\instructions{2\rightarrow N_{13}}%
\EasyButWeakLineBreak%
{${52}\above{1pt}{1pt}{1,0}{\infty b}{52}\above{1pt}{1pt}{r}{2}{2}\above{1pt}{1pt}{*}{4}{4}\above{1pt}{1pt}{*}{2}$}\relax$\,(\times2)$%
\nopagebreak\par%
\nopagebreak\par\leavevmode%
{$\left[\!\llap{\phantom{%
\begingroup \smaller\smaller\smaller\begin{tabular}{@{}c@{}}%
0\\0\\0
\end{tabular}\endgroup%
}}\right.$}%
\begingroup \smaller\smaller\smaller\begin{tabular}{@{}c@{}}%
-14425996\\261508\\523016
\end{tabular}\endgroup%
\kern3pt%
\begingroup \smaller\smaller\smaller\begin{tabular}{@{}c@{}}%
261508\\-4740\\-9481
\end{tabular}\endgroup%
\kern3pt%
\begingroup \smaller\smaller\smaller\begin{tabular}{@{}c@{}}%
523016\\-9481\\-18962
\end{tabular}\endgroup%
{$\left.\llap{\phantom{%
\begingroup \smaller\smaller\smaller\begin{tabular}{@{}c@{}}%
0\\0\\0
\end{tabular}\endgroup%
}}\!\right]$}%
\hfil\penalty500%
{$\left[\!\llap{\phantom{%
\begingroup \smaller\smaller\smaller\begin{tabular}{@{}c@{}}%
0\\0\\0
\end{tabular}\endgroup%
}}\right.$}%
\begingroup \smaller\smaller\smaller\begin{tabular}{@{}c@{}}%
896635\\125112\\24667916
\end{tabular}\endgroup%
\kern3pt%
\begingroup \smaller\smaller\smaller\begin{tabular}{@{}c@{}}%
-16383\\-2287\\-450723
\end{tabular}\endgroup%
\kern3pt%
\begingroup \smaller\smaller\smaller\begin{tabular}{@{}c@{}}%
-32508\\-4536\\-894349
\end{tabular}\endgroup%
{$\left.\llap{\phantom{%
\begingroup \smaller\smaller\smaller\begin{tabular}{@{}c@{}}%
0\\0\\0
\end{tabular}\endgroup%
}}\!\right]$}%
\EasyButWeakLineBreak%
{$\left[\!\llap{\phantom{%
\begingroup \smaller\smaller\smaller\begin{tabular}{@{}c@{}}%
0\\0\\0
\end{tabular}\endgroup%
}}\right.$}%
\begingroup \smaller\smaller\smaller\begin{tabular}{@{}c@{}}%
-1159\\-130\\-31902
\end{tabular}\endgroup%
\HardButStrongLineBreak\kern3pt%
\begingroup \smaller\smaller\smaller\begin{tabular}{@{}c@{}}%
-2351\\-312\\-64688
\end{tabular}\endgroup%
\HardButStrongLineBreak\kern3pt%
\begingroup \smaller\smaller\smaller\begin{tabular}{@{}c@{}}%
-258\\-38\\-7097
\end{tabular}\endgroup%
\HardButStrongLineBreak\kern3pt%
\begingroup \smaller\smaller\smaller\begin{tabular}{@{}c@{}}%
-31\\-10\\-850
\end{tabular}\endgroup%
{$\left.\llap{\phantom{%
\begingroup \smaller\smaller\smaller\begin{tabular}{@{}c@{}}%
0\\0\\0
\end{tabular}\endgroup%
}}\!\right]$}%
%
%
\hbox{}\par\smallskip%
%
%
\leavevmode%
${L_{13.2}}$%
{} : {$1\above{1pt}{1pt}{-2}{2}8\above{1pt}{1pt}{1}{7}{\cdot}1\above{1pt}{1pt}{2}{}13\above{1pt}{1pt}{-}{}$}\spacer%
\instructions{2\rightarrow N'_{10}}%
\EasyButWeakLineBreak%
{${26}\above{1pt}{1pt}{4,3}{\infty b}{104}\above{1pt}{1pt}{s}{2}{4}\above{1pt}{1pt}{*}{4}{2}\above{1pt}{1pt}{s}{2}$}\relax$\,(\times2)$%
\nopagebreak\par%
\nopagebreak\par\leavevmode%
{$\left[\!\llap{\phantom{%
\begingroup \smaller\smaller\smaller\begin{tabular}{@{}c@{}}%
0\\0\\0
\end{tabular}\endgroup%
}}\right.$}%
\begingroup \smaller\smaller\smaller\begin{tabular}{@{}c@{}}%
5304\\104\\2080
\end{tabular}\endgroup%
\kern3pt%
\begingroup \smaller\smaller\smaller\begin{tabular}{@{}c@{}}%
104\\-9\\-65
\end{tabular}\endgroup%
\kern3pt%
\begingroup \smaller\smaller\smaller\begin{tabular}{@{}c@{}}%
2080\\-65\\-198
\end{tabular}\endgroup%
{$\left.\llap{\phantom{%
\begingroup \smaller\smaller\smaller\begin{tabular}{@{}c@{}}%
0\\0\\0
\end{tabular}\endgroup%
}}\!\right]$}%
\hfil\penalty500%
{$\left[\!\llap{\phantom{%
\begingroup \smaller\smaller\smaller\begin{tabular}{@{}c@{}}%
0\\0\\0
\end{tabular}\endgroup%
}}\right.$}%
\begingroup \smaller\smaller\smaller\begin{tabular}{@{}c@{}}%
-11233\\-527280\\54912
\end{tabular}\endgroup%
\kern3pt%
\begingroup \smaller\smaller\smaller\begin{tabular}{@{}c@{}}%
468\\21969\\-2288
\end{tabular}\endgroup%
\kern3pt%
\begingroup \smaller\smaller\smaller\begin{tabular}{@{}c@{}}%
2196\\103090\\-10737
\end{tabular}\endgroup%
{$\left.\llap{\phantom{%
\begingroup \smaller\smaller\smaller\begin{tabular}{@{}c@{}}%
0\\0\\0
\end{tabular}\endgroup%
}}\!\right]$}%
\EasyButWeakLineBreak%
{$\left[\!\llap{\phantom{%
\begingroup \smaller\smaller\smaller\begin{tabular}{@{}c@{}}%
0\\0\\0
\end{tabular}\endgroup%
}}\right.$}%
\begingroup \smaller\smaller\smaller\begin{tabular}{@{}c@{}}%
-412\\-19344\\2015
\end{tabular}\endgroup%
\HardButStrongLineBreak\kern3pt%
\begingroup \smaller\smaller\smaller\begin{tabular}{@{}c@{}}%
-1893\\-88868\\9256
\end{tabular}\endgroup%
\HardButStrongLineBreak\kern3pt%
\begingroup \smaller\smaller\smaller\begin{tabular}{@{}c@{}}%
-225\\-10562\\1100
\end{tabular}\endgroup%
\HardButStrongLineBreak\kern3pt%
\begingroup \smaller\smaller\smaller\begin{tabular}{@{}c@{}}%
-26\\-1220\\127
\end{tabular}\endgroup%
{$\left.\llap{\phantom{%
\begingroup \smaller\smaller\smaller\begin{tabular}{@{}c@{}}%
0\\0\\0
\end{tabular}\endgroup%
}}\!\right]$}%
%
%
\hbox{}\par\smallskip%
%
%
\leavevmode%
${L_{13.3}}$%
{} : {$1\above{1pt}{1pt}{2}{2}8\above{1pt}{1pt}{-}{3}{\cdot}1\above{1pt}{1pt}{2}{}13\above{1pt}{1pt}{-}{}$}\spacer%
\instructions{m}%
\EasyButWeakLineBreak%
{${26}\above{1pt}{1pt}{4,3}{\infty a}{104}\above{1pt}{1pt}{l}{2}{1}\above{1pt}{1pt}{}{4}{2}\above{1pt}{1pt}{b}{2}$}\relax$\,(\times2)$%
\nopagebreak\par%
\nopagebreak\par\leavevmode%
{$\left[\!\llap{\phantom{%
\begingroup \smaller\smaller\smaller\begin{tabular}{@{}c@{}}%
0\\0\\0
\end{tabular}\endgroup%
}}\right.$}%
\begingroup \smaller\smaller\smaller\begin{tabular}{@{}c@{}}%
-2600\\-2600\\1352
\end{tabular}\endgroup%
\kern3pt%
\begingroup \smaller\smaller\smaller\begin{tabular}{@{}c@{}}%
-2600\\-2574\\1325
\end{tabular}\endgroup%
\kern3pt%
\begingroup \smaller\smaller\smaller\begin{tabular}{@{}c@{}}%
1352\\1325\\-675
\end{tabular}\endgroup%
{$\left.\llap{\phantom{%
\begingroup \smaller\smaller\smaller\begin{tabular}{@{}c@{}}%
0\\0\\0
\end{tabular}\endgroup%
}}\!\right]$}%
\hfil\penalty500%
{$\left[\!\llap{\phantom{%
\begingroup \smaller\smaller\smaller\begin{tabular}{@{}c@{}}%
0\\0\\0
\end{tabular}\endgroup%
}}\right.$}%
\begingroup \smaller\smaller\smaller\begin{tabular}{@{}c@{}}%
-17109\\33488\\31304
\end{tabular}\endgroup%
\kern3pt%
\begingroup \smaller\smaller\smaller\begin{tabular}{@{}c@{}}%
-15463\\30267\\28294
\end{tabular}\endgroup%
\kern3pt%
\begingroup \smaller\smaller\smaller\begin{tabular}{@{}c@{}}%
7191\\-14076\\-13159
\end{tabular}\endgroup%
{$\left.\llap{\phantom{%
\begingroup \smaller\smaller\smaller\begin{tabular}{@{}c@{}}%
0\\0\\0
\end{tabular}\endgroup%
}}\!\right]$}%
\EasyButWeakLineBreak%
{$\left[\!\llap{\phantom{%
\begingroup \smaller\smaller\smaller\begin{tabular}{@{}c@{}}%
0\\0\\0
\end{tabular}\endgroup%
}}\right.$}%
\begingroup \smaller\smaller\smaller\begin{tabular}{@{}c@{}}%
-150\\299\\286
\end{tabular}\endgroup%
\HardButStrongLineBreak\kern3pt%
\begingroup \smaller\smaller\smaller\begin{tabular}{@{}c@{}}%
-25\\52\\52
\end{tabular}\endgroup%
\HardButStrongLineBreak\kern3pt%
\begingroup \smaller\smaller\smaller\begin{tabular}{@{}c@{}}%
13\\-26\\-25
\end{tabular}\endgroup%
\HardButStrongLineBreak\kern3pt%
\begingroup \smaller\smaller\smaller\begin{tabular}{@{}c@{}}%
2\\-5\\-6
\end{tabular}\endgroup%
{$\left.\llap{\phantom{%
\begingroup \smaller\smaller\smaller\begin{tabular}{@{}c@{}}%
0\\0\\0
\end{tabular}\endgroup%
}}\!\right]$}%

\medskip%
%
\leavevmode\llap{}%
$W_{14\phantom{0}}$%
\qquad\llap{6} lattices, $\chi=12$%
\hfill%
$222222\rtimes C_{2}$%
\nopagebreak\smallskip\hrule\nopagebreak\medskip%
%
%
\leavevmode%
${L_{14.1}}$%
{} : {$1\above{1pt}{1pt}{-2}{{\rm II}}4\above{1pt}{1pt}{1}{1}{\cdot}1\above{1pt}{1pt}{-2}{}13\above{1pt}{1pt}{-}{}$}\spacer%
\instructions{2\rightarrow N_{14}}%
\EasyButWeakLineBreak%
{${4}\above{1pt}{1pt}{r}{2}{2}\above{1pt}{1pt}{b}{2}{26}\above{1pt}{1pt}{l}{2}$}\relax$\,(\times2)$%
\nopagebreak\par%
\nopagebreak\par\leavevmode%
{$\left[\!\llap{\phantom{%
\begingroup \smaller\smaller\smaller\begin{tabular}{@{}c@{}}%
0\\0\\0
\end{tabular}\endgroup%
}}\right.$}%
\begingroup \smaller\smaller\smaller\begin{tabular}{@{}c@{}}%
-572\\-260\\52
\end{tabular}\endgroup%
\kern3pt%
\begingroup \smaller\smaller\smaller\begin{tabular}{@{}c@{}}%
-260\\-118\\25
\end{tabular}\endgroup%
\kern3pt%
\begingroup \smaller\smaller\smaller\begin{tabular}{@{}c@{}}%
52\\25\\6
\end{tabular}\endgroup%
{$\left.\llap{\phantom{%
\begingroup \smaller\smaller\smaller\begin{tabular}{@{}c@{}}%
0\\0\\0
\end{tabular}\endgroup%
}}\!\right]$}%
\hfil\penalty500%
{$\left[\!\llap{\phantom{%
\begingroup \smaller\smaller\smaller\begin{tabular}{@{}c@{}}%
0\\0\\0
\end{tabular}\endgroup%
}}\right.$}%
\begingroup \smaller\smaller\smaller\begin{tabular}{@{}c@{}}%
103\\-208\\-104
\end{tabular}\endgroup%
\kern3pt%
\begingroup \smaller\smaller\smaller\begin{tabular}{@{}c@{}}%
49\\-99\\-49
\end{tabular}\endgroup%
\kern3pt%
\begingroup \smaller\smaller\smaller\begin{tabular}{@{}c@{}}%
4\\-8\\-5
\end{tabular}\endgroup%
{$\left.\llap{\phantom{%
\begingroup \smaller\smaller\smaller\begin{tabular}{@{}c@{}}%
0\\0\\0
\end{tabular}\endgroup%
}}\!\right]$}%
\EasyButWeakLineBreak%
{$\left[\!\llap{\phantom{%
\begingroup \smaller\smaller\smaller\begin{tabular}{@{}c@{}}%
0\\0\\0
\end{tabular}\endgroup%
}}\right.$}%
\begingroup \smaller\smaller\smaller\begin{tabular}{@{}c@{}}%
17\\-36\\0
\end{tabular}\endgroup%
\HardButStrongLineBreak\kern3pt%
\begingroup \smaller\smaller\smaller\begin{tabular}{@{}c@{}}%
7\\-15\\2
\end{tabular}\endgroup%
\HardButStrongLineBreak\kern3pt%
\begingroup \smaller\smaller\smaller\begin{tabular}{@{}c@{}}%
-6\\13\\0
\end{tabular}\endgroup%
{$\left.\llap{\phantom{%
\begingroup \smaller\smaller\smaller\begin{tabular}{@{}c@{}}%
0\\0\\0
\end{tabular}\endgroup%
}}\!\right]$}%

\medskip%
%
\leavevmode\llap{}%
$W_{15\phantom{0}}$%
\qquad\llap{8} lattices, $\chi=8$%
\hfill%
$\infty232$%
\nopagebreak\smallskip\hrule\nopagebreak\medskip%
%
%
\leavevmode%
${L_{15.1}}$%
{} : {$1\above{1pt}{1pt}{-2}{{\rm II}}8\above{1pt}{1pt}{-}{3}{\cdot}1\above{1pt}{1pt}{-2}{}7\above{1pt}{1pt}{1}{}$}\spacer%
\instructions{2\rightarrow N_{15}}%
\EasyButWeakLineBreak%
{${14}\above{1pt}{1pt}{4,3}{\infty a}{56}\above{1pt}{1pt}{b}{2}{2}\above{1pt}{1pt}{-}{3}{2}\above{1pt}{1pt}{b}{2}$}%
\nopagebreak\par%
\nopagebreak\par\leavevmode%
{$\left[\!\llap{\phantom{%
\begingroup \smaller\smaller\smaller\begin{tabular}{@{}c@{}}%
0\\0\\0
\end{tabular}\endgroup%
}}\right.$}%
\begingroup \smaller\smaller\smaller\begin{tabular}{@{}c@{}}%
-25256\\392\\672
\end{tabular}\endgroup%
\kern3pt%
\begingroup \smaller\smaller\smaller\begin{tabular}{@{}c@{}}%
392\\-6\\-11
\end{tabular}\endgroup%
\kern3pt%
\begingroup \smaller\smaller\smaller\begin{tabular}{@{}c@{}}%
672\\-11\\-14
\end{tabular}\endgroup%
{$\left.\llap{\phantom{%
\begingroup \smaller\smaller\smaller\begin{tabular}{@{}c@{}}%
0\\0\\0
\end{tabular}\endgroup%
}}\!\right]$}%
\EasyButWeakLineBreak%
{$\left[\!\llap{\phantom{%
\begingroup \smaller\smaller\smaller\begin{tabular}{@{}c@{}}%
0\\0\\0
\end{tabular}\endgroup%
}}\right.$}%
\begingroup \smaller\smaller\smaller\begin{tabular}{@{}c@{}}%
-2\\-105\\-14
\end{tabular}\endgroup%
\HardButStrongLineBreak\kern3pt%
\begingroup \smaller\smaller\smaller\begin{tabular}{@{}c@{}}%
11\\560\\84
\end{tabular}\endgroup%
\HardButStrongLineBreak\kern3pt%
\begingroup \smaller\smaller\smaller\begin{tabular}{@{}c@{}}%
1\\52\\7
\end{tabular}\endgroup%
\HardButStrongLineBreak\kern3pt%
\begingroup \smaller\smaller\smaller\begin{tabular}{@{}c@{}}%
-1\\-51\\-8
\end{tabular}\endgroup%
{$\left.\llap{\phantom{%
\begingroup \smaller\smaller\smaller\begin{tabular}{@{}c@{}}%
0\\0\\0
\end{tabular}\endgroup%
}}\!\right]$}%

\medskip%
%
\leavevmode\llap{}%
$W_{16\phantom{0}}$%
\qquad\llap{24} lattices, $\chi=4$%
\hfill%
$2262$%
\nopagebreak\smallskip\hrule\nopagebreak\medskip%
%
%
\leavevmode%
${L_{16.1}}$%
{} : {$1\above{1pt}{1pt}{-2}{{\rm II}}4\above{1pt}{1pt}{-}{3}{\cdot}1\above{1pt}{1pt}{-}{}3\above{1pt}{1pt}{-}{}9\above{1pt}{1pt}{1}{}{\cdot}1\above{1pt}{1pt}{-2}{}5\above{1pt}{1pt}{1}{}$}\spacer%
\instructions{23\rightarrow N_{16},3,2}%
\EasyButWeakLineBreak%
{${36}\above{1pt}{1pt}{*}{2}{20}\above{1pt}{1pt}{b}{2}{6}\above{1pt}{1pt}{}{6}{2}\above{1pt}{1pt}{b}{2}$}%
\nopagebreak\par%
\nopagebreak\par\leavevmode%
{$\left[\!\llap{\phantom{%
\begingroup \smaller\smaller\smaller\begin{tabular}{@{}c@{}}%
0\\0\\0
\end{tabular}\endgroup%
}}\right.$}%
\begingroup \smaller\smaller\smaller\begin{tabular}{@{}c@{}}%
-324180\\-9360\\3600
\end{tabular}\endgroup%
\kern3pt%
\begingroup \smaller\smaller\smaller\begin{tabular}{@{}c@{}}%
-9360\\-66\\69
\end{tabular}\endgroup%
\kern3pt%
\begingroup \smaller\smaller\smaller\begin{tabular}{@{}c@{}}%
3600\\69\\-34
\end{tabular}\endgroup%
{$\left.\llap{\phantom{%
\begingroup \smaller\smaller\smaller\begin{tabular}{@{}c@{}}%
0\\0\\0
\end{tabular}\endgroup%
}}\!\right]$}%
\EasyButWeakLineBreak%
{$\left[\!\llap{\phantom{%
\begingroup \smaller\smaller\smaller\begin{tabular}{@{}c@{}}%
0\\0\\0
\end{tabular}\endgroup%
}}\right.$}%
\begingroup \smaller\smaller\smaller\begin{tabular}{@{}c@{}}%
-5\\-138\\-810
\end{tabular}\endgroup%
\HardButStrongLineBreak\kern3pt%
\begingroup \smaller\smaller\smaller\begin{tabular}{@{}c@{}}%
9\\250\\1460
\end{tabular}\endgroup%
\HardButStrongLineBreak\kern3pt%
\begingroup \smaller\smaller\smaller\begin{tabular}{@{}c@{}}%
3\\83\\486
\end{tabular}\endgroup%
\HardButStrongLineBreak\kern3pt%
\begingroup \smaller\smaller\smaller\begin{tabular}{@{}c@{}}%
-4\\-111\\-649
\end{tabular}\endgroup%
{$\left.\llap{\phantom{%
\begingroup \smaller\smaller\smaller\begin{tabular}{@{}c@{}}%
0\\0\\0
\end{tabular}\endgroup%
}}\!\right]$}%

\medskip%
%
\leavevmode\llap{}%
$W_{17\phantom{0}}$%
\qquad\llap{44} lattices, $\chi=6$%
\hfill%
$22222$%
\nopagebreak\smallskip\hrule\nopagebreak\medskip%
%
%
\leavevmode%
${L_{17.1}}$%
{} : {$1\above{1pt}{1pt}{2}{{\rm II}}4\above{1pt}{1pt}{1}{7}{\cdot}1\above{1pt}{1pt}{2}{}3\above{1pt}{1pt}{1}{}{\cdot}1\above{1pt}{1pt}{-2}{}5\above{1pt}{1pt}{1}{}$}\spacer%
\instructions{2\rightarrow N_{17}}%
\EasyButWeakLineBreak%
{${4}\above{1pt}{1pt}{*}{2}{12}\above{1pt}{1pt}{*}{2}{20}\above{1pt}{1pt}{b}{2}{2}\above{1pt}{1pt}{s}{2}{30}\above{1pt}{1pt}{b}{2}$}%
\nopagebreak\par%
\nopagebreak\par\leavevmode%
{$\left[\!\llap{\phantom{%
\begingroup \smaller\smaller\smaller\begin{tabular}{@{}c@{}}%
0\\0\\0
\end{tabular}\endgroup%
}}\right.$}%
\begingroup \smaller\smaller\smaller\begin{tabular}{@{}c@{}}%
-22020\\360\\240
\end{tabular}\endgroup%
\kern3pt%
\begingroup \smaller\smaller\smaller\begin{tabular}{@{}c@{}}%
360\\-4\\-5
\end{tabular}\endgroup%
\kern3pt%
\begingroup \smaller\smaller\smaller\begin{tabular}{@{}c@{}}%
240\\-5\\-2
\end{tabular}\endgroup%
{$\left.\llap{\phantom{%
\begingroup \smaller\smaller\smaller\begin{tabular}{@{}c@{}}%
0\\0\\0
\end{tabular}\endgroup%
}}\!\right]$}%
\EasyButWeakLineBreak%
{$\left[\!\llap{\phantom{%
\begingroup \smaller\smaller\smaller\begin{tabular}{@{}c@{}}%
0\\0\\0
\end{tabular}\endgroup%
}}\right.$}%
\begingroup \smaller\smaller\smaller\begin{tabular}{@{}c@{}}%
-1\\-28\\-50
\end{tabular}\endgroup%
\HardButStrongLineBreak\kern3pt%
\begingroup \smaller\smaller\smaller\begin{tabular}{@{}c@{}}%
-1\\-30\\-48
\end{tabular}\endgroup%
\HardButStrongLineBreak\kern3pt%
\begingroup \smaller\smaller\smaller\begin{tabular}{@{}c@{}}%
3\\80\\150
\end{tabular}\endgroup%
\HardButStrongLineBreak\kern3pt%
\begingroup \smaller\smaller\smaller\begin{tabular}{@{}c@{}}%
1\\28\\49
\end{tabular}\endgroup%
\HardButStrongLineBreak\kern3pt%
\begingroup \smaller\smaller\smaller\begin{tabular}{@{}c@{}}%
1\\30\\45
\end{tabular}\endgroup%
{$\left.\llap{\phantom{%
\begingroup \smaller\smaller\smaller\begin{tabular}{@{}c@{}}%
0\\0\\0
\end{tabular}\endgroup%
}}\!\right]$}%
%
%
\hbox{}\par\smallskip%
%
%
\leavevmode%
${L_{17.2}}$%
{} : {$1\above{1pt}{1pt}{-2}{6}8\above{1pt}{1pt}{-}{5}{\cdot}1\above{1pt}{1pt}{2}{}3\above{1pt}{1pt}{-}{}{\cdot}1\above{1pt}{1pt}{-2}{}5\above{1pt}{1pt}{-}{}$}\spacer%
\instructions{2\rightarrow N'_{13}}%
\EasyButWeakLineBreak%
{${8}\above{1pt}{1pt}{b}{2}{6}\above{1pt}{1pt}{l}{2}{40}\above{1pt}{1pt}{}{2}{1}\above{1pt}{1pt}{r}{2}{60}\above{1pt}{1pt}{*}{2}$}%
\nopagebreak\par%
\nopagebreak\par\leavevmode%
{$\left[\!\llap{\phantom{%
\begingroup \smaller\smaller\smaller\begin{tabular}{@{}c@{}}%
0\\0\\0
\end{tabular}\endgroup%
}}\right.$}%
\begingroup \smaller\smaller\smaller\begin{tabular}{@{}c@{}}%
9960\\3240\\-120
\end{tabular}\endgroup%
\kern3pt%
\begingroup \smaller\smaller\smaller\begin{tabular}{@{}c@{}}%
3240\\1054\\-39
\end{tabular}\endgroup%
\kern3pt%
\begingroup \smaller\smaller\smaller\begin{tabular}{@{}c@{}}%
-120\\-39\\1
\end{tabular}\endgroup%
{$\left.\llap{\phantom{%
\begingroup \smaller\smaller\smaller\begin{tabular}{@{}c@{}}%
0\\0\\0
\end{tabular}\endgroup%
}}\!\right]$}%
\EasyButWeakLineBreak%
{$\left[\!\llap{\phantom{%
\begingroup \smaller\smaller\smaller\begin{tabular}{@{}c@{}}%
0\\0\\0
\end{tabular}\endgroup%
}}\right.$}%
\begingroup \smaller\smaller\smaller\begin{tabular}{@{}c@{}}%
9\\-28\\-8
\end{tabular}\endgroup%
\HardButStrongLineBreak\kern3pt%
\begingroup \smaller\smaller\smaller\begin{tabular}{@{}c@{}}%
-1\\3\\0
\end{tabular}\endgroup%
\HardButStrongLineBreak\kern3pt%
\begingroup \smaller\smaller\smaller\begin{tabular}{@{}c@{}}%
-13\\40\\0
\end{tabular}\endgroup%
\HardButStrongLineBreak\kern3pt%
\begingroup \smaller\smaller\smaller\begin{tabular}{@{}c@{}}%
0\\0\\-1
\end{tabular}\endgroup%
\HardButStrongLineBreak\kern3pt%
\begingroup \smaller\smaller\smaller\begin{tabular}{@{}c@{}}%
29\\-90\\-30
\end{tabular}\endgroup%
{$\left.\llap{\phantom{%
\begingroup \smaller\smaller\smaller\begin{tabular}{@{}c@{}}%
0\\0\\0
\end{tabular}\endgroup%
}}\!\right]$}%
%
%
\hbox{}\par\smallskip%
%
%
\leavevmode%
${L_{17.3}}$%
{} : {$1\above{1pt}{1pt}{2}{6}8\above{1pt}{1pt}{1}{1}{\cdot}1\above{1pt}{1pt}{2}{}3\above{1pt}{1pt}{-}{}{\cdot}1\above{1pt}{1pt}{-2}{}5\above{1pt}{1pt}{-}{}$}\spacer%
\instructions{m}%
\EasyButWeakLineBreak%
{${8}\above{1pt}{1pt}{r}{2}{6}\above{1pt}{1pt}{b}{2}{40}\above{1pt}{1pt}{*}{2}{4}\above{1pt}{1pt}{l}{2}{15}\above{1pt}{1pt}{}{2}$}%
\nopagebreak\par%
\nopagebreak\par\leavevmode%
{$\left[\!\llap{\phantom{%
\begingroup \smaller\smaller\smaller\begin{tabular}{@{}c@{}}%
0\\0\\0
\end{tabular}\endgroup%
}}\right.$}%
\begingroup \smaller\smaller\smaller\begin{tabular}{@{}c@{}}%
-126840\\-3360\\600
\end{tabular}\endgroup%
\kern3pt%
\begingroup \smaller\smaller\smaller\begin{tabular}{@{}c@{}}%
-3360\\-89\\16
\end{tabular}\endgroup%
\kern3pt%
\begingroup \smaller\smaller\smaller\begin{tabular}{@{}c@{}}%
600\\16\\-1
\end{tabular}\endgroup%
{$\left.\llap{\phantom{%
\begingroup \smaller\smaller\smaller\begin{tabular}{@{}c@{}}%
0\\0\\0
\end{tabular}\endgroup%
}}\!\right]$}%
\EasyButWeakLineBreak%
{$\left[\!\llap{\phantom{%
\begingroup \smaller\smaller\smaller\begin{tabular}{@{}c@{}}%
0\\0\\0
\end{tabular}\endgroup%
}}\right.$}%
\begingroup \smaller\smaller\smaller\begin{tabular}{@{}c@{}}%
-3\\112\\-8
\end{tabular}\endgroup%
\HardButStrongLineBreak\kern3pt%
\begingroup \smaller\smaller\smaller\begin{tabular}{@{}c@{}}%
-2\\75\\-3
\end{tabular}\endgroup%
\HardButStrongLineBreak\kern3pt%
\begingroup \smaller\smaller\smaller\begin{tabular}{@{}c@{}}%
7\\-260\\20
\end{tabular}\endgroup%
\HardButStrongLineBreak\kern3pt%
\begingroup \smaller\smaller\smaller\begin{tabular}{@{}c@{}}%
3\\-112\\6
\end{tabular}\endgroup%
\HardButStrongLineBreak\kern3pt%
\begingroup \smaller\smaller\smaller\begin{tabular}{@{}c@{}}%
2\\-75\\0
\end{tabular}\endgroup%
{$\left.\llap{\phantom{%
\begingroup \smaller\smaller\smaller\begin{tabular}{@{}c@{}}%
0\\0\\0
\end{tabular}\endgroup%
}}\!\right]$}%

\medskip%
%
\leavevmode\llap{}%
$W_{18\phantom{0}}$%
\qquad\llap{88} lattices, $\chi=18$%
\hfill%
$\infty22222$%
\nopagebreak\smallskip\hrule\nopagebreak\medskip%
%
%
\leavevmode%
${L_{18.1}}$%
{} : {$1\above{1pt}{1pt}{2}{{\rm II}}4\above{1pt}{1pt}{1}{7}{\cdot}1\above{1pt}{1pt}{1}{}3\above{1pt}{1pt}{-}{}9\above{1pt}{1pt}{-}{}{\cdot}1\above{1pt}{1pt}{2}{}5\above{1pt}{1pt}{-}{}$}\spacer%
\instructions{23\rightarrow N_{18},3,2}%
\EasyButWeakLineBreak%
{${60}\above{1pt}{1pt}{3,2}{\infty b}{60}\above{1pt}{1pt}{r}{2}{18}\above{1pt}{1pt}{b}{2}{10}\above{1pt}{1pt}{s}{2}{6}\above{1pt}{1pt}{b}{2}{4}\above{1pt}{1pt}{*}{2}$}%
\nopagebreak\par%
\nopagebreak\par\leavevmode%
{$\left[\!\llap{\phantom{%
\begingroup \smaller\smaller\smaller\begin{tabular}{@{}c@{}}%
0\\0\\0
\end{tabular}\endgroup%
}}\right.$}%
\begingroup \smaller\smaller\smaller\begin{tabular}{@{}c@{}}%
-18654660\\443520\\38700
\end{tabular}\endgroup%
\kern3pt%
\begingroup \smaller\smaller\smaller\begin{tabular}{@{}c@{}}%
443520\\-10542\\-921
\end{tabular}\endgroup%
\kern3pt%
\begingroup \smaller\smaller\smaller\begin{tabular}{@{}c@{}}%
38700\\-921\\-80
\end{tabular}\endgroup%
{$\left.\llap{\phantom{%
\begingroup \smaller\smaller\smaller\begin{tabular}{@{}c@{}}%
0\\0\\0
\end{tabular}\endgroup%
}}\!\right]$}%
\EasyButWeakLineBreak%
{$\left[\!\llap{\phantom{%
\begingroup \smaller\smaller\smaller\begin{tabular}{@{}c@{}}%
0\\0\\0
\end{tabular}\endgroup%
}}\right.$}%
\begingroup \smaller\smaller\smaller\begin{tabular}{@{}c@{}}%
-53\\-1750\\-5490
\end{tabular}\endgroup%
\HardButStrongLineBreak\kern3pt%
\begingroup \smaller\smaller\smaller\begin{tabular}{@{}c@{}}%
3\\100\\300
\end{tabular}\endgroup%
\HardButStrongLineBreak\kern3pt%
\begingroup \smaller\smaller\smaller\begin{tabular}{@{}c@{}}%
13\\429\\1350
\end{tabular}\endgroup%
\HardButStrongLineBreak\kern3pt%
\begingroup \smaller\smaller\smaller\begin{tabular}{@{}c@{}}%
2\\65\\220
\end{tabular}\endgroup%
\HardButStrongLineBreak\kern3pt%
\begingroup \smaller\smaller\smaller\begin{tabular}{@{}c@{}}%
-10\\-331\\-1026
\end{tabular}\endgroup%
\HardButStrongLineBreak\kern3pt%
\begingroup \smaller\smaller\smaller\begin{tabular}{@{}c@{}}%
-21\\-694\\-2168
\end{tabular}\endgroup%
{$\left.\llap{\phantom{%
\begingroup \smaller\smaller\smaller\begin{tabular}{@{}c@{}}%
0\\0\\0
\end{tabular}\endgroup%
}}\!\right]$}%
%
%
\hbox{}\par\smallskip%
%
%
\leavevmode%
${L_{18.2}}$%
{} : {$1\above{1pt}{1pt}{2}{6}8\above{1pt}{1pt}{1}{1}{\cdot}1\above{1pt}{1pt}{-}{}3\above{1pt}{1pt}{1}{}9\above{1pt}{1pt}{1}{}{\cdot}1\above{1pt}{1pt}{2}{}5\above{1pt}{1pt}{1}{}$}\spacer%
\instructions{3m,3,2}%
\EasyButWeakLineBreak%
{${30}\above{1pt}{1pt}{12,5}{\infty b}{120}\above{1pt}{1pt}{s}{2}{36}\above{1pt}{1pt}{*}{2}{20}\above{1pt}{1pt}{l}{2}{3}\above{1pt}{1pt}{}{2}{8}\above{1pt}{1pt}{r}{2}$}%
\nopagebreak\par%
\nopagebreak\par\leavevmode%
{$\left[\!\llap{\phantom{%
\begingroup \smaller\smaller\smaller\begin{tabular}{@{}c@{}}%
0\\0\\0
\end{tabular}\endgroup%
}}\right.$}%
\begingroup \smaller\smaller\smaller\begin{tabular}{@{}c@{}}%
-116280\\-4680\\-5040
\end{tabular}\endgroup%
\kern3pt%
\begingroup \smaller\smaller\smaller\begin{tabular}{@{}c@{}}%
-4680\\-186\\-201
\end{tabular}\endgroup%
\kern3pt%
\begingroup \smaller\smaller\smaller\begin{tabular}{@{}c@{}}%
-5040\\-201\\-217
\end{tabular}\endgroup%
{$\left.\llap{\phantom{%
\begingroup \smaller\smaller\smaller\begin{tabular}{@{}c@{}}%
0\\0\\0
\end{tabular}\endgroup%
}}\!\right]$}%
\EasyButWeakLineBreak%
{$\left[\!\llap{\phantom{%
\begingroup \smaller\smaller\smaller\begin{tabular}{@{}c@{}}%
0\\0\\0
\end{tabular}\endgroup%
}}\right.$}%
\begingroup \smaller\smaller\smaller\begin{tabular}{@{}c@{}}%
-2\\-145\\180
\end{tabular}\endgroup%
\HardButStrongLineBreak\kern3pt%
\begingroup \smaller\smaller\smaller\begin{tabular}{@{}c@{}}%
-1\\-40\\60
\end{tabular}\endgroup%
\HardButStrongLineBreak\kern3pt%
\begingroup \smaller\smaller\smaller\begin{tabular}{@{}c@{}}%
1\\72\\-90
\end{tabular}\endgroup%
\HardButStrongLineBreak\kern3pt%
\begingroup \smaller\smaller\smaller\begin{tabular}{@{}c@{}}%
1\\50\\-70
\end{tabular}\endgroup%
\HardButStrongLineBreak\kern3pt%
\begingroup \smaller\smaller\smaller\begin{tabular}{@{}c@{}}%
0\\-10\\9
\end{tabular}\endgroup%
\HardButStrongLineBreak\kern3pt%
\begingroup \smaller\smaller\smaller\begin{tabular}{@{}c@{}}%
-1\\-88\\104
\end{tabular}\endgroup%
{$\left.\llap{\phantom{%
\begingroup \smaller\smaller\smaller\begin{tabular}{@{}c@{}}%
0\\0\\0
\end{tabular}\endgroup%
}}\!\right]$}%
%
%
\hbox{}\par\smallskip%
%
%
\leavevmode%
${L_{18.3}}$%
{} : {$1\above{1pt}{1pt}{-2}{6}8\above{1pt}{1pt}{-}{5}{\cdot}1\above{1pt}{1pt}{-}{}3\above{1pt}{1pt}{1}{}9\above{1pt}{1pt}{1}{}{\cdot}1\above{1pt}{1pt}{2}{}5\above{1pt}{1pt}{1}{}$}\spacer%
\instructions{32\rightarrow N'_{12},3,m}%
\EasyButWeakLineBreak%
{${30}\above{1pt}{1pt}{12,5}{\infty a}{120}\above{1pt}{1pt}{l}{2}{9}\above{1pt}{1pt}{}{2}{5}\above{1pt}{1pt}{r}{2}{12}\above{1pt}{1pt}{*}{2}{8}\above{1pt}{1pt}{b}{2}$}%
\nopagebreak\par%
\nopagebreak\par\leavevmode%
{$\left[\!\llap{\phantom{%
\begingroup \smaller\smaller\smaller\begin{tabular}{@{}c@{}}%
0\\0\\0
\end{tabular}\endgroup%
}}\right.$}%
\begingroup \smaller\smaller\smaller\begin{tabular}{@{}c@{}}%
-1572120\\-310680\\109800
\end{tabular}\endgroup%
\kern3pt%
\begingroup \smaller\smaller\smaller\begin{tabular}{@{}c@{}}%
-310680\\-61395\\21696
\end{tabular}\endgroup%
\kern3pt%
\begingroup \smaller\smaller\smaller\begin{tabular}{@{}c@{}}%
109800\\21696\\-7663
\end{tabular}\endgroup%
{$\left.\llap{\phantom{%
\begingroup \smaller\smaller\smaller\begin{tabular}{@{}c@{}}%
0\\0\\0
\end{tabular}\endgroup%
}}\!\right]$}%
\EasyButWeakLineBreak%
{$\left[\!\llap{\phantom{%
\begingroup \smaller\smaller\smaller\begin{tabular}{@{}c@{}}%
0\\0\\0
\end{tabular}\endgroup%
}}\right.$}%
\begingroup \smaller\smaller\smaller\begin{tabular}{@{}c@{}}%
-219\\1315\\585
\end{tabular}\endgroup%
\HardButStrongLineBreak\kern3pt%
\begingroup \smaller\smaller\smaller\begin{tabular}{@{}c@{}}%
43\\-260\\-120
\end{tabular}\endgroup%
\HardButStrongLineBreak\kern3pt%
\begingroup \smaller\smaller\smaller\begin{tabular}{@{}c@{}}%
44\\-264\\-117
\end{tabular}\endgroup%
\HardButStrongLineBreak\kern3pt%
\begingroup \smaller\smaller\smaller\begin{tabular}{@{}c@{}}%
-14\\85\\40
\end{tabular}\endgroup%
\HardButStrongLineBreak\kern3pt%
\begingroup \smaller\smaller\smaller\begin{tabular}{@{}c@{}}%
-113\\680\\306
\end{tabular}\endgroup%
\HardButStrongLineBreak\kern3pt%
\begingroup \smaller\smaller\smaller\begin{tabular}{@{}c@{}}%
-195\\1172\\524
\end{tabular}\endgroup%
{$\left.\llap{\phantom{%
\begingroup \smaller\smaller\smaller\begin{tabular}{@{}c@{}}%
0\\0\\0
\end{tabular}\endgroup%
}}\!\right]$}%

\medskip%
%
\leavevmode\llap{}%
$W_{19\phantom{0}}$%
\qquad\llap{12} lattices, $\chi=3$%
\hfill%
$4222$%
\nopagebreak\smallskip\hrule\nopagebreak\medskip%
%
%
\leavevmode%
${L_{19.1}}$%
{} : {$1\above{1pt}{1pt}{-2}{{\rm II}}4\above{1pt}{1pt}{-}{3}{\cdot}1\above{1pt}{1pt}{2}{}3\above{1pt}{1pt}{1}{}{\cdot}1\above{1pt}{1pt}{2}{}5\above{1pt}{1pt}{-}{}$}\spacer%
\instructions{2\rightarrow N_{19}}%
\EasyButWeakLineBreak%
{${2}\above{1pt}{1pt}{*}{4}{4}\above{1pt}{1pt}{b}{2}{10}\above{1pt}{1pt}{l}{2}{12}\above{1pt}{1pt}{r}{2}$}%
\nopagebreak\par%
\nopagebreak\par\leavevmode%
{$\left[\!\llap{\phantom{%
\begingroup \smaller\smaller\smaller\begin{tabular}{@{}c@{}}%
0\\0\\0
\end{tabular}\endgroup%
}}\right.$}%
\begingroup \smaller\smaller\smaller\begin{tabular}{@{}c@{}}%
-27060\\300\\540
\end{tabular}\endgroup%
\kern3pt%
\begingroup \smaller\smaller\smaller\begin{tabular}{@{}c@{}}%
300\\-2\\-7
\end{tabular}\endgroup%
\kern3pt%
\begingroup \smaller\smaller\smaller\begin{tabular}{@{}c@{}}%
540\\-7\\-10
\end{tabular}\endgroup%
{$\left.\llap{\phantom{%
\begingroup \smaller\smaller\smaller\begin{tabular}{@{}c@{}}%
0\\0\\0
\end{tabular}\endgroup%
}}\!\right]$}%
\EasyButWeakLineBreak%
{$\left[\!\llap{\phantom{%
\begingroup \smaller\smaller\smaller\begin{tabular}{@{}c@{}}%
0\\0\\0
\end{tabular}\endgroup%
}}\right.$}%
\begingroup \smaller\smaller\smaller\begin{tabular}{@{}c@{}}%
1\\27\\35
\end{tabular}\endgroup%
\HardButStrongLineBreak\kern3pt%
\begingroup \smaller\smaller\smaller\begin{tabular}{@{}c@{}}%
-1\\-26\\-36
\end{tabular}\endgroup%
\HardButStrongLineBreak\kern3pt%
\begingroup \smaller\smaller\smaller\begin{tabular}{@{}c@{}}%
-2\\-55\\-70
\end{tabular}\endgroup%
\HardButStrongLineBreak\kern3pt%
\begingroup \smaller\smaller\smaller\begin{tabular}{@{}c@{}}%
1\\24\\36
\end{tabular}\endgroup%
{$\left.\llap{\phantom{%
\begingroup \smaller\smaller\smaller\begin{tabular}{@{}c@{}}%
0\\0\\0
\end{tabular}\endgroup%
}}\!\right]$}%

\medskip%
%
\leavevmode\llap{}%
$W_{20\phantom{0}}$%
\qquad\llap{22} lattices, $\chi=27$%
\hfill%
$\infty222242$%
\nopagebreak\smallskip\hrule\nopagebreak\medskip%
%
%
\leavevmode%
${L_{20.1}}$%
{} : {$1\above{1pt}{1pt}{2}{{\rm II}}4\above{1pt}{1pt}{1}{1}{\cdot}1\above{1pt}{1pt}{2}{}17\above{1pt}{1pt}{1}{}$}\spacer%
\instructions{2\rightarrow N_{20}}%
\EasyButWeakLineBreak%
{${68}\above{1pt}{1pt}{1,0}{\infty b}{68}\above{1pt}{1pt}{r}{2}{2}\above{1pt}{1pt}{b}{2}{34}\above{1pt}{1pt}{l}{2}{4}\above{1pt}{1pt}{r}{2}{2}\above{1pt}{1pt}{*}{4}{4}\above{1pt}{1pt}{*}{2}$}%
\nopagebreak\par%
\nopagebreak\par\leavevmode%
{$\left[\!\llap{\phantom{%
\begingroup \smaller\smaller\smaller\begin{tabular}{@{}c@{}}%
0\\0\\0
\end{tabular}\endgroup%
}}\right.$}%
\begingroup \smaller\smaller\smaller\begin{tabular}{@{}c@{}}%
-53788\\1156\\2652
\end{tabular}\endgroup%
\kern3pt%
\begingroup \smaller\smaller\smaller\begin{tabular}{@{}c@{}}%
1156\\-24\\-59
\end{tabular}\endgroup%
\kern3pt%
\begingroup \smaller\smaller\smaller\begin{tabular}{@{}c@{}}%
2652\\-59\\-126
\end{tabular}\endgroup%
{$\left.\llap{\phantom{%
\begingroup \smaller\smaller\smaller\begin{tabular}{@{}c@{}}%
0\\0\\0
\end{tabular}\endgroup%
}}\!\right]$}%
\EasyButWeakLineBreak%
{$\left[\!\llap{\phantom{%
\begingroup \smaller\smaller\smaller\begin{tabular}{@{}c@{}}%
0\\0\\0
\end{tabular}\endgroup%
}}\right.$}%
\begingroup \smaller\smaller\smaller\begin{tabular}{@{}c@{}}%
-7\\-170\\-68
\end{tabular}\endgroup%
\HardButStrongLineBreak\kern3pt%
\begingroup \smaller\smaller\smaller\begin{tabular}{@{}c@{}}%
75\\1768\\748
\end{tabular}\endgroup%
\HardButStrongLineBreak\kern3pt%
\begingroup \smaller\smaller\smaller\begin{tabular}{@{}c@{}}%
17\\402\\169
\end{tabular}\endgroup%
\HardButStrongLineBreak\kern3pt%
\begingroup \smaller\smaller\smaller\begin{tabular}{@{}c@{}}%
132\\3128\\1309
\end{tabular}\endgroup%
\HardButStrongLineBreak\kern3pt%
\begingroup \smaller\smaller\smaller\begin{tabular}{@{}c@{}}%
17\\404\\168
\end{tabular}\endgroup%
\HardButStrongLineBreak\kern3pt%
\begingroup \smaller\smaller\smaller\begin{tabular}{@{}c@{}}%
4\\96\\39
\end{tabular}\endgroup%
\HardButStrongLineBreak\kern3pt%
\begingroup \smaller\smaller\smaller\begin{tabular}{@{}c@{}}%
-5\\-118\\-50
\end{tabular}\endgroup%
{$\left.\llap{\phantom{%
\begingroup \smaller\smaller\smaller\begin{tabular}{@{}c@{}}%
0\\0\\0
\end{tabular}\endgroup%
}}\!\right]$}%
%
%
\hbox{}\par\smallskip%
%
%
\leavevmode%
${L_{20.2}}$%
{} : {$1\above{1pt}{1pt}{2}{2}8\above{1pt}{1pt}{1}{7}{\cdot}1\above{1pt}{1pt}{2}{}17\above{1pt}{1pt}{1}{}$}\spacer%
\instructions{2\rightarrow N'_{15}}%
\EasyButWeakLineBreak%
{${34}\above{1pt}{1pt}{4,3}{\infty b}{136}\above{1pt}{1pt}{l}{2}{1}\above{1pt}{1pt}{}{2}{17}\above{1pt}{1pt}{r}{2}{8}\above{1pt}{1pt}{l}{2}{1}\above{1pt}{1pt}{}{4}{2}\above{1pt}{1pt}{s}{2}$}%
\nopagebreak\par%
\nopagebreak\par\leavevmode%
{$\left[\!\llap{\phantom{%
\begingroup \smaller\smaller\smaller\begin{tabular}{@{}c@{}}%
0\\0\\0
\end{tabular}\endgroup%
}}\right.$}%
\begingroup \smaller\smaller\smaller\begin{tabular}{@{}c@{}}%
-2045576\\8160\\11424
\end{tabular}\endgroup%
\kern3pt%
\begingroup \smaller\smaller\smaller\begin{tabular}{@{}c@{}}%
8160\\-30\\-47
\end{tabular}\endgroup%
\kern3pt%
\begingroup \smaller\smaller\smaller\begin{tabular}{@{}c@{}}%
11424\\-47\\-63
\end{tabular}\endgroup%
{$\left.\llap{\phantom{%
\begingroup \smaller\smaller\smaller\begin{tabular}{@{}c@{}}%
0\\0\\0
\end{tabular}\endgroup%
}}\!\right]$}%
\EasyButWeakLineBreak%
{$\left[\!\llap{\phantom{%
\begingroup \smaller\smaller\smaller\begin{tabular}{@{}c@{}}%
0\\0\\0
\end{tabular}\endgroup%
}}\right.$}%
\begingroup \smaller\smaller\smaller\begin{tabular}{@{}c@{}}%
-4\\-289\\-510
\end{tabular}\endgroup%
\HardButStrongLineBreak\kern3pt%
\begingroup \smaller\smaller\smaller\begin{tabular}{@{}c@{}}%
59\\4216\\7548
\end{tabular}\endgroup%
\HardButStrongLineBreak\kern3pt%
\begingroup \smaller\smaller\smaller\begin{tabular}{@{}c@{}}%
7\\501\\895
\end{tabular}\endgroup%
\HardButStrongLineBreak\kern3pt%
\begingroup \smaller\smaller\smaller\begin{tabular}{@{}c@{}}%
56\\4012\\7157
\end{tabular}\endgroup%
\HardButStrongLineBreak\kern3pt%
\begingroup \smaller\smaller\smaller\begin{tabular}{@{}c@{}}%
15\\1076\\1916
\end{tabular}\endgroup%
\HardButStrongLineBreak\kern3pt%
\begingroup \smaller\smaller\smaller\begin{tabular}{@{}c@{}}%
2\\144\\255
\end{tabular}\endgroup%
\HardButStrongLineBreak\kern3pt%
\begingroup \smaller\smaller\smaller\begin{tabular}{@{}c@{}}%
-2\\-143\\-256
\end{tabular}\endgroup%
{$\left.\llap{\phantom{%
\begingroup \smaller\smaller\smaller\begin{tabular}{@{}c@{}}%
0\\0\\0
\end{tabular}\endgroup%
}}\!\right]$}%
%
%
\hbox{}\par\smallskip%
%
%
\leavevmode%
${L_{20.3}}$%
{} : {$1\above{1pt}{1pt}{-2}{2}8\above{1pt}{1pt}{-}{3}{\cdot}1\above{1pt}{1pt}{2}{}17\above{1pt}{1pt}{1}{}$}\spacer%
\instructions{m}%
\EasyButWeakLineBreak%
{${34}\above{1pt}{1pt}{4,3}{\infty a}{136}\above{1pt}{1pt}{s}{2}{4}\above{1pt}{1pt}{*}{2}{68}\above{1pt}{1pt}{s}{2}{8}\above{1pt}{1pt}{s}{2}{4}\above{1pt}{1pt}{*}{4}{2}\above{1pt}{1pt}{b}{2}$}%
\nopagebreak\par%
\nopagebreak\par\leavevmode%
{$\left[\!\llap{\phantom{%
\begingroup \smaller\smaller\smaller\begin{tabular}{@{}c@{}}%
0\\0\\0
\end{tabular}\endgroup%
}}\right.$}%
\begingroup \smaller\smaller\smaller\begin{tabular}{@{}c@{}}%
-62696\\-15640\\952
\end{tabular}\endgroup%
\kern3pt%
\begingroup \smaller\smaller\smaller\begin{tabular}{@{}c@{}}%
-15640\\-3901\\237
\end{tabular}\endgroup%
\kern3pt%
\begingroup \smaller\smaller\smaller\begin{tabular}{@{}c@{}}%
952\\237\\-14
\end{tabular}\endgroup%
{$\left.\llap{\phantom{%
\begingroup \smaller\smaller\smaller\begin{tabular}{@{}c@{}}%
0\\0\\0
\end{tabular}\endgroup%
}}\!\right]$}%
\EasyButWeakLineBreak%
{$\left[\!\llap{\phantom{%
\begingroup \smaller\smaller\smaller\begin{tabular}{@{}c@{}}%
0\\0\\0
\end{tabular}\endgroup%
}}\right.$}%
\begingroup \smaller\smaller\smaller\begin{tabular}{@{}c@{}}%
-40\\170\\153
\end{tabular}\endgroup%
\HardButStrongLineBreak\kern3pt%
\begingroup \smaller\smaller\smaller\begin{tabular}{@{}c@{}}%
-175\\748\\748
\end{tabular}\endgroup%
\HardButStrongLineBreak\kern3pt%
\begingroup \smaller\smaller\smaller\begin{tabular}{@{}c@{}}%
-21\\90\\94
\end{tabular}\endgroup%
\HardButStrongLineBreak\kern3pt%
\begingroup \smaller\smaller\smaller\begin{tabular}{@{}c@{}}%
-87\\374\\408
\end{tabular}\endgroup%
\HardButStrongLineBreak\kern3pt%
\begingroup \smaller\smaller\smaller\begin{tabular}{@{}c@{}}%
1\\-4\\0
\end{tabular}\endgroup%
\HardButStrongLineBreak\kern3pt%
\begingroup \smaller\smaller\smaller\begin{tabular}{@{}c@{}}%
7\\-30\\-32
\end{tabular}\endgroup%
\HardButStrongLineBreak\kern3pt%
\begingroup \smaller\smaller\smaller\begin{tabular}{@{}c@{}}%
-1\\4\\-1
\end{tabular}\endgroup%
{$\left.\llap{\phantom{%
\begingroup \smaller\smaller\smaller\begin{tabular}{@{}c@{}}%
0\\0\\0
\end{tabular}\endgroup%
}}\!\right]$}%

\medskip%
%
\leavevmode\llap{}%
$W_{21\phantom{0}}$%
\qquad\llap{6} lattices, $\chi=18$%
\hfill%
$422422\rtimes C_{2}$%
\nopagebreak\smallskip\hrule\nopagebreak\medskip%
%
%
\leavevmode%
${L_{21.1}}$%
{} : {$1\above{1pt}{1pt}{-2}{{\rm II}}4\above{1pt}{1pt}{1}{7}{\cdot}1\above{1pt}{1pt}{2}{}19\above{1pt}{1pt}{-}{}$}\spacer%
\instructions{2\rightarrow N_{21}}%
\EasyButWeakLineBreak%
{${2}\above{1pt}{1pt}{*}{4}{4}\above{1pt}{1pt}{b}{2}{38}\above{1pt}{1pt}{s}{2}$}\relax$\,(\times2)$%
\nopagebreak\par%
\nopagebreak\par\leavevmode%
{$\left[\!\llap{\phantom{%
\begingroup \smaller\smaller\smaller\begin{tabular}{@{}c@{}}%
0\\0\\0
\end{tabular}\endgroup%
}}\right.$}%
\begingroup \smaller\smaller\smaller\begin{tabular}{@{}c@{}}%
-26372\\-18468\\1672
\end{tabular}\endgroup%
\kern3pt%
\begingroup \smaller\smaller\smaller\begin{tabular}{@{}c@{}}%
-18468\\-12930\\1171
\end{tabular}\endgroup%
\kern3pt%
\begingroup \smaller\smaller\smaller\begin{tabular}{@{}c@{}}%
1672\\1171\\-106
\end{tabular}\endgroup%
{$\left.\llap{\phantom{%
\begingroup \smaller\smaller\smaller\begin{tabular}{@{}c@{}}%
0\\0\\0
\end{tabular}\endgroup%
}}\!\right]$}%
\hfil\penalty500%
{$\left[\!\llap{\phantom{%
\begingroup \smaller\smaller\smaller\begin{tabular}{@{}c@{}}%
0\\0\\0
\end{tabular}\endgroup%
}}\right.$}%
\begingroup \smaller\smaller\smaller\begin{tabular}{@{}c@{}}%
3191\\-1368\\35112
\end{tabular}\endgroup%
\kern3pt%
\begingroup \smaller\smaller\smaller\begin{tabular}{@{}c@{}}%
2233\\-958\\24563
\end{tabular}\endgroup%
\kern3pt%
\begingroup \smaller\smaller\smaller\begin{tabular}{@{}c@{}}%
-203\\87\\-2234
\end{tabular}\endgroup%
{$\left.\llap{\phantom{%
\begingroup \smaller\smaller\smaller\begin{tabular}{@{}c@{}}%
0\\0\\0
\end{tabular}\endgroup%
}}\!\right]$}%
\EasyButWeakLineBreak%
{$\left[\!\llap{\phantom{%
\begingroup \smaller\smaller\smaller\begin{tabular}{@{}c@{}}%
0\\0\\0
\end{tabular}\endgroup%
}}\right.$}%
\begingroup \smaller\smaller\smaller\begin{tabular}{@{}c@{}}%
-4\\2\\-41
\end{tabular}\endgroup%
\HardButStrongLineBreak\kern3pt%
\begingroup \smaller\smaller\smaller\begin{tabular}{@{}c@{}}%
7\\-2\\88
\end{tabular}\endgroup%
\HardButStrongLineBreak\kern3pt%
\begingroup \smaller\smaller\smaller\begin{tabular}{@{}c@{}}%
139\\-57\\1558
\end{tabular}\endgroup%
{$\left.\llap{\phantom{%
\begingroup \smaller\smaller\smaller\begin{tabular}{@{}c@{}}%
0\\0\\0
\end{tabular}\endgroup%
}}\!\right]$}%

\medskip%
%
\leavevmode\llap{}%
$W_{22\phantom{0}}$%
\qquad\llap{12} lattices, $\chi=4$%
\hfill%
$2622$%
\nopagebreak\smallskip\hrule\nopagebreak\medskip%
%
%
\leavevmode%
${L_{22.1}}$%
{} : {$1\above{1pt}{1pt}{-2}{{\rm II}}4\above{1pt}{1pt}{1}{1}{\cdot}1\above{1pt}{1pt}{2}{}3\above{1pt}{1pt}{-}{}{\cdot}1\above{1pt}{1pt}{-2}{}7\above{1pt}{1pt}{-}{}$}\spacer%
\instructions{2\rightarrow N_{22}}%
\EasyButWeakLineBreak%
{${4}\above{1pt}{1pt}{r}{2}{6}\above{1pt}{1pt}{}{6}{2}\above{1pt}{1pt}{b}{2}{42}\above{1pt}{1pt}{l}{2}$}%
\nopagebreak\par%
\nopagebreak\par\leavevmode%
{$\left[\!\llap{\phantom{%
\begingroup \smaller\smaller\smaller\begin{tabular}{@{}c@{}}%
0\\0\\0
\end{tabular}\endgroup%
}}\right.$}%
\begingroup \smaller\smaller\smaller\begin{tabular}{@{}c@{}}%
-25788\\588\\252
\end{tabular}\endgroup%
\kern3pt%
\begingroup \smaller\smaller\smaller\begin{tabular}{@{}c@{}}%
588\\-10\\-7
\end{tabular}\endgroup%
\kern3pt%
\begingroup \smaller\smaller\smaller\begin{tabular}{@{}c@{}}%
252\\-7\\-2
\end{tabular}\endgroup%
{$\left.\llap{\phantom{%
\begingroup \smaller\smaller\smaller\begin{tabular}{@{}c@{}}%
0\\0\\0
\end{tabular}\endgroup%
}}\!\right]$}%
\EasyButWeakLineBreak%
{$\left[\!\llap{\phantom{%
\begingroup \smaller\smaller\smaller\begin{tabular}{@{}c@{}}%
0\\0\\0
\end{tabular}\endgroup%
}}\right.$}%
\begingroup \smaller\smaller\smaller\begin{tabular}{@{}c@{}}%
-1\\-20\\-56
\end{tabular}\endgroup%
\HardButStrongLineBreak\kern3pt%
\begingroup \smaller\smaller\smaller\begin{tabular}{@{}c@{}}%
-1\\-21\\-54
\end{tabular}\endgroup%
\HardButStrongLineBreak\kern3pt%
\begingroup \smaller\smaller\smaller\begin{tabular}{@{}c@{}}%
1\\20\\55
\end{tabular}\endgroup%
\HardButStrongLineBreak\kern3pt%
\begingroup \smaller\smaller\smaller\begin{tabular}{@{}c@{}}%
2\\42\\105
\end{tabular}\endgroup%
{$\left.\llap{\phantom{%
\begingroup \smaller\smaller\smaller\begin{tabular}{@{}c@{}}%
0\\0\\0
\end{tabular}\endgroup%
}}\!\right]$}%

\medskip%
%
\leavevmode\llap{}%
$W_{23\phantom{0}}$%
\qquad\llap{88} lattices, $\chi=48$%
\hfill%
$\infty2222\infty2222\rtimes C_{2}$%
\nopagebreak\smallskip\hrule\nopagebreak\medskip%
%
%
\leavevmode%
${L_{23.1}}$%
{} : {$1\above{1pt}{1pt}{2}{{\rm II}}4\above{1pt}{1pt}{-}{5}{\cdot}1\above{1pt}{1pt}{1}{}3\above{1pt}{1pt}{1}{}9\above{1pt}{1pt}{-}{}{\cdot}1\above{1pt}{1pt}{-2}{}7\above{1pt}{1pt}{-}{}$}\spacer%
\instructions{23\rightarrow N_{23},3,2}%
\EasyButWeakLineBreak%
{${84}\above{1pt}{1pt}{3,2}{\infty b}{84}\above{1pt}{1pt}{r}{2}{18}\above{1pt}{1pt}{b}{2}{12}\above{1pt}{1pt}{*}{2}{4}\above{1pt}{1pt}{*}{2}$}\relax$\,(\times2)$%
\nopagebreak\par%
\nopagebreak\par\leavevmode%
{$\left[\!\llap{\phantom{%
\begingroup \smaller\smaller\smaller
\begin{tabular}{@{}c@{}}%
0\\0\\0
\end{tabular}\endgroup%
}}\right.$}%
\begingroup \smaller\smaller\smaller\begin{tabular}{@{}c@{}}%
-16363116\\948024\\53676
\end{tabular}\endgroup%
\kern3pt%
\begingroup \smaller\smaller\smaller\begin{tabular}{@{}c@{}}%
948024\\-54906\\-3111
\end{tabular}\endgroup%
\kern3pt%
\begingroup \smaller\smaller\smaller\begin{tabular}{@{}c@{}}%
53676\\-3111\\-176
\end{tabular}\endgroup%
{$\left.\llap{\phantom{%
\begingroup \smaller\smaller\smaller\begin{tabular}{@{}c@{}}%
0\\0\\0
\end{tabular}\endgroup%
}}\!\right]$}%
\hfil\penalty500%
{$\left[\!\llap{\phantom{%
\begingroup \smaller\smaller\smaller\begin{tabular}{@{}c@{}}%
0\\0\\0
\end{tabular}\endgroup%
}}\right.$}%
\begingroup \smaller\smaller\smaller\begin{tabular}{@{}c@{}}%
505343\\4563888\\73432800
\end{tabular}\endgroup%
\kern3pt%
\begingroup \smaller\smaller\smaller\begin{tabular}{@{}c@{}}%
-29696\\-268193\\-4315200
\end{tabular}\endgroup%
\kern3pt%
\begingroup \smaller\smaller\smaller\begin{tabular}{@{}c@{}}%
-1632\\-14739\\-237151
\end{tabular}\endgroup%
{$\left.\llap{\phantom{%
\begingroup \smaller\smaller\smaller\begin{tabular}{@{}c@{}}%
0\\0\\0
\end{tabular}\endgroup%
}}\!\right]$}%
\EasyButWeakLineBreak%
{$\left[\!\llap{\phantom{%
\begingroup \smaller\smaller\smaller\begin{tabular}{@{}c@{}}%
0\\0\\0
\end{tabular}\endgroup%
}}\right.$}%
\begingroup \smaller\smaller\smaller\begin{tabular}{@{}c@{}}%
-333\\-3010\\-48342
\end{tabular}\endgroup%
\HardButStrongLineBreak\kern3pt%
\begingroup \smaller\smaller\smaller\begin{tabular}{@{}c@{}}%
-899\\-8120\\-130620
\end{tabular}\endgroup%
\HardButStrongLineBreak\kern3pt%
\begingroup \smaller\smaller\smaller\begin{tabular}{@{}c@{}}%
-403\\-3639\\-58572
\end{tabular}\endgroup%
\HardButStrongLineBreak\kern3pt%
\begingroup \smaller\smaller\smaller\begin{tabular}{@{}c@{}}%
-259\\-2338\\-37656
\end{tabular}\endgroup%
\HardButStrongLineBreak\kern3pt%
\begingroup \smaller\smaller\smaller\begin{tabular}{@{}c@{}}%
-49\\-442\\-7130
\end{tabular}\endgroup%
{$\left.\llap{\phantom{%
\begingroup \smaller\smaller\smaller\begin{tabular}{@{}c@{}}%
0\\0\\0
\end{tabular}\endgroup%
}}\!\right]$}%
%
%
\hbox{}\par\smallskip%
%
%
\leavevmode%
${L_{23.2}}$%
{} : {$1\above{1pt}{1pt}{-2}{2}8\above{1pt}{1pt}{1}{7}{\cdot}1\above{1pt}{1pt}{-}{}3\above{1pt}{1pt}{-}{}9\above{1pt}{1pt}{1}{}{\cdot}1\above{1pt}{1pt}{-2}{}7\above{1pt}{1pt}{-}{}$}\spacer%
\instructions{3m,3,2}%
\EasyButWeakLineBreak%
{${42}\above{1pt}{1pt}{12,11}{\infty b}{168}\above{1pt}{1pt}{s}{2}{36}\above{1pt}{1pt}{*}{2}{24}\above{1pt}{1pt}{b}{2}{2}\above{1pt}{1pt}{s}{2}$}\relax$\,(\times2)$%
\nopagebreak\par%
\nopagebreak\par\leavevmode%
{$\left[\!\llap{\phantom{%
\begingroup \smaller\smaller\smaller
\begin{tabular}{@{}c@{}}%
0\\0\\0
\end{tabular}\endgroup%
}}\right.$}%
\begingroup \smaller\smaller\smaller\begin{tabular}{@{}c@{}}%
-8168328\\-56448\\-59472
\end{tabular}\endgroup%
\kern3pt%
\begingroup \smaller\smaller\smaller\begin{tabular}{@{}c@{}}%
-56448\\-390\\-411
\end{tabular}\endgroup%
\kern3pt%
\begingroup \smaller\smaller\smaller\begin{tabular}{@{}c@{}}%
-59472\\-411\\-433
\end{tabular}\endgroup%
{$\left.\llap{\phantom{%
\begingroup \smaller\smaller\smaller\begin{tabular}{@{}c@{}}%
0\\0\\0
\end{tabular}\endgroup%
}}\!\right]$}%
\hfil\penalty500%
{$\left[\!\llap{\phantom{%
\begingroup \smaller\smaller\smaller\begin{tabular}{@{}c@{}}%
0\\0\\0
\end{tabular}\endgroup%
}}\right.$}%
\begingroup \smaller\smaller\smaller\begin{tabular}{@{}c@{}}%
6299\\88200\\-945000
\end{tabular}\endgroup%
\kern3pt%
\begingroup \smaller\smaller\smaller\begin{tabular}{@{}c@{}}%
43\\601\\-6450
\end{tabular}\endgroup%
\kern3pt%
\begingroup \smaller\smaller\smaller\begin{tabular}{@{}c@{}}%
46\\644\\-6901
\end{tabular}\endgroup%
{$\left.\llap{\phantom{%
\begingroup \smaller\smaller\smaller\begin{tabular}{@{}c@{}}%
0\\0\\0
\end{tabular}\endgroup%
}}\!\right]$}%
\EasyButWeakLineBreak%
{$\left[\!\llap{\phantom{%
\begingroup \smaller\smaller\smaller\begin{tabular}{@{}c@{}}%
0\\0\\0
\end{tabular}\endgroup%
}}\right.$}%
\begingroup \smaller\smaller\smaller\begin{tabular}{@{}c@{}}%
-4\\-133\\672
\end{tabular}\endgroup%
\HardButStrongLineBreak\kern3pt%
\begingroup \smaller\smaller\smaller\begin{tabular}{@{}c@{}}%
-27\\-448\\4116
\end{tabular}\endgroup%
\HardButStrongLineBreak\kern3pt%
\begingroup \smaller\smaller\smaller\begin{tabular}{@{}c@{}}%
-13\\-156\\1926
\end{tabular}\endgroup%
\HardButStrongLineBreak\kern3pt%
\begingroup \smaller\smaller\smaller\begin{tabular}{@{}c@{}}%
-9\\-68\\1296
\end{tabular}\endgroup%
\HardButStrongLineBreak\kern3pt%
\begingroup \smaller\smaller\smaller\begin{tabular}{@{}c@{}}%
-1\\1\\136
\end{tabular}\endgroup%
{$\left.\llap{\phantom{%
\begingroup \smaller\smaller\smaller\begin{tabular}{@{}c@{}}%
0\\0\\0
\end{tabular}\endgroup%
}}\!\right]$}%
%
%
\hbox{}\par\smallskip%
%
%
\leavevmode%
${L_{23.3}}$%
{} : {$1\above{1pt}{1pt}{2}{2}8\above{1pt}{1pt}{-}{3}{\cdot}1\above{1pt}{1pt}{-}{}3\above{1pt}{1pt}{-}{}9\above{1pt}{1pt}{1}{}{\cdot}1\above{1pt}{1pt}{-2}{}7\above{1pt}{1pt}{-}{}$}\spacer%
\instructions{32\rightarrow N'_{16},3,m}%
\EasyButWeakLineBreak%
{${42}\above{1pt}{1pt}{12,11}{\infty a}{168}\above{1pt}{1pt}{l}{2}{9}\above{1pt}{1pt}{}{2}{24}\above{1pt}{1pt}{r}{2}{2}\above{1pt}{1pt}{b}{2}$}\relax$\,(\times2)$%
\nopagebreak\par%
\nopagebreak\par\leavevmode%
{$\left[\!\llap{\phantom{%
\begingroup \smaller\smaller\smaller
\begin{tabular}{@{}c@{}}%
0\\0\\0
\end{tabular}\endgroup%
}}\right.$}%
\begingroup \smaller\smaller\smaller\begin{tabular}{@{}c@{}}%
-22810536\\3272976\\-990864
\end{tabular}\endgroup%
\kern3pt%
\begingroup \smaller\smaller\smaller\begin{tabular}{@{}c@{}}%
3272976\\-469623\\142176
\end{tabular}\endgroup%
\kern3pt%
\begingroup \smaller\smaller\smaller\begin{tabular}{@{}c@{}}%
-990864\\142176\\-43039
\end{tabular}\endgroup%
{$\left.\llap{\phantom{%
\begingroup \smaller\smaller\smaller\begin{tabular}{@{}c@{}}%
0\\0\\0
\end{tabular}\endgroup%
}}\!\right]$}%
\hfil\penalty500%
{$\left[\!\llap{\phantom{%
\begingroup \smaller\smaller\smaller\begin{tabular}{@{}c@{}}%
0\\0\\0
\end{tabular}\endgroup%
}}\right.$}%
\begingroup \smaller\smaller\smaller\begin{tabular}{@{}c@{}}%
-6504625\\-39054960\\20738592
\end{tabular}\endgroup%
\kern3pt%
\begingroup \smaller\smaller\smaller\begin{tabular}{@{}c@{}}%
932578\\5599369\\-2973324
\end{tabular}\endgroup%
\kern3pt%
\begingroup \smaller\smaller\smaller\begin{tabular}{@{}c@{}}%
-283932\\-1704780\\905255
\end{tabular}\endgroup%
{$\left.\llap{\phantom{%
\begingroup \smaller\smaller\smaller\begin{tabular}{@{}c@{}}%
0\\0\\0
\end{tabular}\endgroup%
}}\!\right]$}%
\EasyButWeakLineBreak%
{$\left[\!\llap{\phantom{%
\begingroup \smaller\smaller\smaller\begin{tabular}{@{}c@{}}%
0\\0\\0
\end{tabular}\endgroup%
}}\right.$}%
\begingroup \smaller\smaller\smaller\begin{tabular}{@{}c@{}}%
-2617\\-15715\\8337
\end{tabular}\endgroup%
\HardButStrongLineBreak\kern3pt%
\begingroup \smaller\smaller\smaller\begin{tabular}{@{}c@{}}%
-13491\\-81004\\43008
\end{tabular}\endgroup%
\HardButStrongLineBreak\kern3pt%
\begingroup \smaller\smaller\smaller\begin{tabular}{@{}c@{}}%
-2972\\-17844\\9477
\end{tabular}\endgroup%
\HardButStrongLineBreak\kern3pt%
\begingroup \smaller\smaller\smaller\begin{tabular}{@{}c@{}}%
-3747\\-22496\\11952
\end{tabular}\endgroup%
\HardButStrongLineBreak\kern3pt%
\begingroup \smaller\smaller\smaller\begin{tabular}{@{}c@{}}%
-338\\-2029\\1079
\end{tabular}\endgroup%
{$\left.\llap{\phantom{%
\begingroup \smaller\smaller\smaller\begin{tabular}{@{}c@{}}%
0\\0\\0
\end{tabular}\endgroup%
}}\!\right]$}%

\medskip%
%
\leavevmode\llap{}%
$W_{24\phantom{0}}$%
\qquad\llap{44} lattices, $\chi=9$%
\hfill%
$22224$%
\nopagebreak\smallskip\hrule\nopagebreak\medskip%
%
%
\leavevmode%
${L_{24.1}}$%
{} : {$1\above{1pt}{1pt}{2}{{\rm II}}4\above{1pt}{1pt}{-}{5}{\cdot}1\above{1pt}{1pt}{2}{}3\above{1pt}{1pt}{-}{}{\cdot}1\above{1pt}{1pt}{2}{}7\above{1pt}{1pt}{1}{}$}\spacer%
\instructions{2\rightarrow N_{24}}%
\EasyButWeakLineBreak%
{${4}\above{1pt}{1pt}{*}{2}{28}\above{1pt}{1pt}{b}{2}{6}\above{1pt}{1pt}{b}{2}{14}\above{1pt}{1pt}{s}{2}{2}\above{1pt}{1pt}{*}{4}$}%
\nopagebreak\par%
\nopagebreak\par\leavevmode%
{$\left[\!\llap{\phantom{%
\begingroup \smaller\smaller\smaller\begin{tabular}{@{}c@{}}%
0\\0\\0
\end{tabular}\endgroup%
}}\right.$}%
\begingroup \smaller\smaller\smaller\begin{tabular}{@{}c@{}}%
-278796\\2436\\2856
\end{tabular}\endgroup%
\kern3pt%
\begingroup \smaller\smaller\smaller\begin{tabular}{@{}c@{}}%
2436\\-20\\-27
\end{tabular}\endgroup%
\kern3pt%
\begingroup \smaller\smaller\smaller\begin{tabular}{@{}c@{}}%
2856\\-27\\-26
\end{tabular}\endgroup%
{$\left.\llap{\phantom{%
\begingroup \smaller\smaller\smaller\begin{tabular}{@{}c@{}}%
0\\0\\0
\end{tabular}\endgroup%
}}\!\right]$}%
\EasyButWeakLineBreak%
{$\left[\!\llap{\phantom{%
\begingroup \smaller\smaller\smaller\begin{tabular}{@{}c@{}}%
0\\0\\0
\end{tabular}\endgroup%
}}\right.$}%
\begingroup \smaller\smaller\smaller\begin{tabular}{@{}c@{}}%
-3\\-198\\-124
\end{tabular}\endgroup%
\HardButStrongLineBreak\kern3pt%
\begingroup \smaller\smaller\smaller\begin{tabular}{@{}c@{}}%
-3\\-196\\-126
\end{tabular}\endgroup%
\HardButStrongLineBreak\kern3pt%
\begingroup \smaller\smaller\smaller\begin{tabular}{@{}c@{}}%
4\\264\\165
\end{tabular}\endgroup%
\HardButStrongLineBreak\kern3pt%
\begingroup \smaller\smaller\smaller\begin{tabular}{@{}c@{}}%
13\\854\\539
\end{tabular}\endgroup%
\HardButStrongLineBreak\kern3pt%
\begingroup \smaller\smaller\smaller\begin{tabular}{@{}c@{}}%
3\\196\\125
\end{tabular}\endgroup%
{$\left.\llap{\phantom{%
\begingroup \smaller\smaller\smaller\begin{tabular}{@{}c@{}}%
0\\0\\0
\end{tabular}\endgroup%
}}\!\right]$}%
%
%
\hbox{}\par\smallskip%
%
%
\leavevmode%
${L_{24.2}}$%
{} : {$1\above{1pt}{1pt}{2}{2}8\above{1pt}{1pt}{-}{3}{\cdot}1\above{1pt}{1pt}{2}{}3\above{1pt}{1pt}{1}{}{\cdot}1\above{1pt}{1pt}{2}{}7\above{1pt}{1pt}{1}{}$}\spacer%
\instructions{2\rightarrow N'_{19}}%
\EasyButWeakLineBreak%
{${2}\above{1pt}{1pt}{b}{2}{56}\above{1pt}{1pt}{*}{2}{12}\above{1pt}{1pt}{*}{2}{28}\above{1pt}{1pt}{l}{2}{1}\above{1pt}{1pt}{}{4}$}%
\nopagebreak\par%
\nopagebreak\par\leavevmode%
{$\left[\!\llap{\phantom{%
\begingroup \smaller\smaller\smaller\begin{tabular}{@{}c@{}}%
0\\0\\0
\end{tabular}\endgroup%
}}\right.$}%
\begingroup \smaller\smaller\smaller\begin{tabular}{@{}c@{}}%
-1129128\\-560616\\11088
\end{tabular}\endgroup%
\kern3pt%
\begingroup \smaller\smaller\smaller\begin{tabular}{@{}c@{}}%
-560616\\-278347\\5504
\end{tabular}\endgroup%
\kern3pt%
\begingroup \smaller\smaller\smaller\begin{tabular}{@{}c@{}}%
11088\\5504\\-107
\end{tabular}\endgroup%
{$\left.\llap{\phantom{%
\begingroup \smaller\smaller\smaller\begin{tabular}{@{}c@{}}%
0\\0\\0
\end{tabular}\endgroup%
}}\!\right]$}%
\EasyButWeakLineBreak%
{$\left[\!\llap{\phantom{%
\begingroup \smaller\smaller\smaller\begin{tabular}{@{}c@{}}%
0\\0\\0
\end{tabular}\endgroup%
}}\right.$}%
\begingroup \smaller\smaller\smaller\begin{tabular}{@{}c@{}}%
74\\-151\\-99
\end{tabular}\endgroup%
\HardButStrongLineBreak\kern3pt%
\begingroup \smaller\smaller\smaller\begin{tabular}{@{}c@{}}%
151\\-308\\-196
\end{tabular}\endgroup%
\HardButStrongLineBreak\kern3pt%
\begingroup \smaller\smaller\smaller\begin{tabular}{@{}c@{}}%
-197\\402\\264
\end{tabular}\endgroup%
\HardButStrongLineBreak\kern3pt%
\begingroup \smaller\smaller\smaller\begin{tabular}{@{}c@{}}%
-645\\1316\\854
\end{tabular}\endgroup%
\HardButStrongLineBreak\kern3pt%
\begingroup \smaller\smaller\smaller\begin{tabular}{@{}c@{}}%
-75\\153\\98
\end{tabular}\endgroup%
{$\left.\llap{\phantom{%
\begingroup \smaller\smaller\smaller\begin{tabular}{@{}c@{}}%
0\\0\\0
\end{tabular}\endgroup%
}}\!\right]$}%
%
%
\hbox{}\par\smallskip%
%
%
\leavevmode%
${L_{24.3}}$%
{} : {$1\above{1pt}{1pt}{-2}{2}8\above{1pt}{1pt}{1}{7}{\cdot}1\above{1pt}{1pt}{2}{}3\above{1pt}{1pt}{1}{}{\cdot}1\above{1pt}{1pt}{2}{}7\above{1pt}{1pt}{1}{}$}\spacer%
\instructions{m}%
\EasyButWeakLineBreak%
{${2}\above{1pt}{1pt}{l}{2}{56}\above{1pt}{1pt}{}{2}{3}\above{1pt}{1pt}{}{2}{7}\above{1pt}{1pt}{r}{2}{4}\above{1pt}{1pt}{*}{4}$}%
\nopagebreak\par%
\nopagebreak\par\leavevmode%
{$\left[\!\llap{\phantom{%
\begingroup \smaller\smaller\smaller\begin{tabular}{@{}c@{}}%
0\\0\\0
\end{tabular}\endgroup%
}}\right.$}%
\begingroup \smaller\smaller\smaller\begin{tabular}{@{}c@{}}%
-514248\\3024\\1848
\end{tabular}\endgroup%
\kern3pt%
\begingroup \smaller\smaller\smaller\begin{tabular}{@{}c@{}}%
3024\\-17\\-12
\end{tabular}\endgroup%
\kern3pt%
\begingroup \smaller\smaller\smaller\begin{tabular}{@{}c@{}}%
1848\\-12\\-5
\end{tabular}\endgroup%
{$\left.\llap{\phantom{%
\begingroup \smaller\smaller\smaller\begin{tabular}{@{}c@{}}%
0\\0\\0
\end{tabular}\endgroup%
}}\!\right]$}%
\EasyButWeakLineBreak%
{$\left[\!\llap{\phantom{%
\begingroup \smaller\smaller\smaller\begin{tabular}{@{}c@{}}%
0\\0\\0
\end{tabular}\endgroup%
}}\right.$}%
\begingroup \smaller\smaller\smaller\begin{tabular}{@{}c@{}}%
1\\119\\83
\end{tabular}\endgroup%
\HardButStrongLineBreak\kern3pt%
\begingroup \smaller\smaller\smaller\begin{tabular}{@{}c@{}}%
15\\1792\\1232
\end{tabular}\endgroup%
\HardButStrongLineBreak\kern3pt%
\begingroup \smaller\smaller\smaller\begin{tabular}{@{}c@{}}%
1\\120\\81
\end{tabular}\endgroup%
\HardButStrongLineBreak\kern3pt%
\begingroup \smaller\smaller\smaller\begin{tabular}{@{}c@{}}%
-1\\-119\\-84
\end{tabular}\endgroup%
\HardButStrongLineBreak\kern3pt%
\begingroup \smaller\smaller\smaller\begin{tabular}{@{}c@{}}%
-1\\-120\\-82
\end{tabular}\endgroup%
{$\left.\llap{\phantom{%
\begingroup \smaller\smaller\smaller\begin{tabular}{@{}c@{}}%
0\\0\\0
\end{tabular}\endgroup%
}}\!\right]$}%

\medskip%
%
\leavevmode\llap{}%
$W_{25\phantom{0}}$%
\qquad\llap{24} lattices, $\chi=6$%
\hfill%
$22222$%
\nopagebreak\smallskip\hrule\nopagebreak\medskip%
%
%
\leavevmode%
${L_{25.1}}$%
{} : {$1\above{1pt}{1pt}{-2}{{\rm II}}4\above{1pt}{1pt}{1}{1}{\cdot}1\above{1pt}{1pt}{-}{}3\above{1pt}{1pt}{1}{}9\above{1pt}{1pt}{1}{}{\cdot}1\above{1pt}{1pt}{2}{}7\above{1pt}{1pt}{1}{}$}\spacer%
\instructions{23\rightarrow N_{25},3,2}%
\EasyButWeakLineBreak%
{${12}\above{1pt}{1pt}{*}{2}{252}\above{1pt}{1pt}{b}{2}{2}\above{1pt}{1pt}{l}{2}{36}\above{1pt}{1pt}{r}{2}{14}\above{1pt}{1pt}{b}{2}$}%
\nopagebreak\par%
\nopagebreak\par\leavevmode%
{$\left[\!\llap{\phantom{%
\begingroup \smaller\smaller\smaller\begin{tabular}{@{}c@{}}%
0\\0\\0
\end{tabular}\endgroup%
}}\right.$}%
\begingroup \smaller\smaller\smaller\begin{tabular}{@{}c@{}}%
1925028\\-70560\\2520
\end{tabular}\endgroup%
\kern3pt%
\begingroup \smaller\smaller\smaller\begin{tabular}{@{}c@{}}%
-70560\\2586\\-93
\end{tabular}\endgroup%
\kern3pt%
\begingroup \smaller\smaller\smaller\begin{tabular}{@{}c@{}}%
2520\\-93\\2
\end{tabular}\endgroup%
{$\left.\llap{\phantom{%
\begingroup \smaller\smaller\smaller\begin{tabular}{@{}c@{}}%
0\\0\\0
\end{tabular}\endgroup%
}}\!\right]$}%
\EasyButWeakLineBreak%
{$\left[\!\llap{\phantom{%
\begingroup \smaller\smaller\smaller\begin{tabular}{@{}c@{}}%
0\\0\\0
\end{tabular}\endgroup%
}}\right.$}%
\begingroup \smaller\smaller\smaller\begin{tabular}{@{}c@{}}%
-5\\-134\\72
\end{tabular}\endgroup%
\HardButStrongLineBreak\kern3pt%
\begingroup \smaller\smaller\smaller\begin{tabular}{@{}c@{}}%
-47\\-1260\\630
\end{tabular}\endgroup%
\HardButStrongLineBreak\kern3pt%
\begingroup \smaller\smaller\smaller\begin{tabular}{@{}c@{}}%
0\\0\\-1
\end{tabular}\endgroup%
\HardButStrongLineBreak\kern3pt%
\begingroup \smaller\smaller\smaller\begin{tabular}{@{}c@{}}%
17\\456\\-216
\end{tabular}\endgroup%
\HardButStrongLineBreak\kern3pt%
\begingroup \smaller\smaller\smaller\begin{tabular}{@{}c@{}}%
6\\161\\-70
\end{tabular}\endgroup%
{$\left.\llap{\phantom{%
\begingroup \smaller\smaller\smaller\begin{tabular}{@{}c@{}}%
0\\0\\0
\end{tabular}\endgroup%
}}\!\right]$}%

\medskip%
%
\leavevmode\llap{}%
$W_{26\phantom{0}}$%
\qquad\llap{8} lattices, $\chi=24$%
\hfill%
$\infty22\infty22\rtimes C_{2}$%
\nopagebreak\smallskip\hrule\nopagebreak\medskip%
%
%
\leavevmode%
${L_{26.1}}$%
{} : {$1\above{1pt}{1pt}{-2}{{\rm II}}8\above{1pt}{1pt}{1}{7}{\cdot}1\above{1pt}{1pt}{-2}{}11\above{1pt}{1pt}{-}{}$}\spacer%
\instructions{2\rightarrow N_{26}}%
\EasyButWeakLineBreak%
{${22}\above{1pt}{1pt}{4,3}{\infty b}{88}\above{1pt}{1pt}{b}{2}{2}\above{1pt}{1pt}{s}{2}$}\relax$\,(\times2)$%
\nopagebreak\par%
\nopagebreak\par\leavevmode%
{$\left[\!\llap{\phantom{%
\begingroup \smaller\smaller\smaller\begin{tabular}{@{}c@{}}%
0\\0\\0
\end{tabular}\endgroup%
}}\right.$}%
\begingroup \smaller\smaller\smaller\begin{tabular}{@{}c@{}}%
-19976\\-9240\\-880
\end{tabular}\endgroup%
\kern3pt%
\begingroup \smaller\smaller\smaller\begin{tabular}{@{}c@{}}%
-9240\\-4274\\-407
\end{tabular}\endgroup%
\kern3pt%
\begingroup \smaller\smaller\smaller\begin{tabular}{@{}c@{}}%
-880\\-407\\-38
\end{tabular}\endgroup%
{$\left.\llap{\phantom{%
\begingroup \smaller\smaller\smaller\begin{tabular}{@{}c@{}}%
0\\0\\0
\end{tabular}\endgroup%
}}\!\right]$}%
\hfil\penalty500%
{$\left[\!\llap{\phantom{%
\begingroup \smaller\smaller\smaller\begin{tabular}{@{}c@{}}%
0\\0\\0
\end{tabular}\endgroup%
}}\right.$}%
\begingroup \smaller\smaller\smaller\begin{tabular}{@{}c@{}}%
4223\\-9152\\352
\end{tabular}\endgroup%
\kern3pt%
\begingroup \smaller\smaller\smaller\begin{tabular}{@{}c@{}}%
1956\\-4239\\163
\end{tabular}\endgroup%
\kern3pt%
\begingroup \smaller\smaller\smaller\begin{tabular}{@{}c@{}}%
192\\-416\\15
\end{tabular}\endgroup%
{$\left.\llap{\phantom{%
\begingroup \smaller\smaller\smaller\begin{tabular}{@{}c@{}}%
0\\0\\0
\end{tabular}\endgroup%
}}\!\right]$}%
\EasyButWeakLineBreak%
{$\left[\!\llap{\phantom{%
\begingroup \smaller\smaller\smaller\begin{tabular}{@{}c@{}}%
0\\0\\0
\end{tabular}\endgroup%
}}\right.$}%
\begingroup \smaller\smaller\smaller\begin{tabular}{@{}c@{}}%
-25\\55\\-11
\end{tabular}\endgroup%
\HardButStrongLineBreak\kern3pt%
\begingroup \smaller\smaller\smaller\begin{tabular}{@{}c@{}}%
61\\-132\\0
\end{tabular}\endgroup%
\HardButStrongLineBreak\kern3pt%
\begingroup \smaller\smaller\smaller\begin{tabular}{@{}c@{}}%
5\\-11\\2
\end{tabular}\endgroup%
{$\left.\llap{\phantom{%
\begingroup \smaller\smaller\smaller\begin{tabular}{@{}c@{}}%
0\\0\\0
\end{tabular}\endgroup%
}}\!\right]$}%

\medskip%
%
\leavevmode\llap{}%
$W_{27\phantom{0}}$%
\qquad\llap{6} lattices, $\chi=22$%
\hfill%
$423423\rtimes C_{2}$%
\nopagebreak\smallskip\hrule\nopagebreak\medskip%
%
%
\leavevmode%
${L_{27.1}}$%
{} : {$1\above{1pt}{1pt}{-2}{{\rm II}}4\above{1pt}{1pt}{-}{3}{\cdot}1\above{1pt}{1pt}{2}{}23\above{1pt}{1pt}{-}{}$}\spacer%
\instructions{2\rightarrow N_{27}}%
\EasyButWeakLineBreak%
{${2}\above{1pt}{1pt}{*}{4}{4}\above{1pt}{1pt}{b}{2}{2}\above{1pt}{1pt}{-}{3}$}\relax$\,(\times2)$%
\nopagebreak\par%
\nopagebreak\par\leavevmode%
{$\left[\!\llap{\phantom{%
\begingroup \smaller\smaller\smaller\begin{tabular}{@{}c@{}}%
0\\0\\0
\end{tabular}\endgroup%
}}\right.$}%
\begingroup \smaller\smaller\smaller\begin{tabular}{@{}c@{}}%
-182068\\2024\\4048
\end{tabular}\endgroup%
\kern3pt%
\begingroup \smaller\smaller\smaller\begin{tabular}{@{}c@{}}%
2024\\-22\\-45
\end{tabular}\endgroup%
\kern3pt%
\begingroup \smaller\smaller\smaller\begin{tabular}{@{}c@{}}%
4048\\-45\\-90
\end{tabular}\endgroup%
{$\left.\llap{\phantom{%
\begingroup \smaller\smaller\smaller\begin{tabular}{@{}c@{}}%
0\\0\\0
\end{tabular}\endgroup%
}}\!\right]$}%
\hfil\penalty500%
{$\left[\!\llap{\phantom{%
\begingroup \smaller\smaller\smaller\begin{tabular}{@{}c@{}}%
0\\0\\0
\end{tabular}\endgroup%
}}\right.$}%
\begingroup \smaller\smaller\smaller\begin{tabular}{@{}c@{}}%
17939\\-21528\\814476
\end{tabular}\endgroup%
\kern3pt%
\begingroup \smaller\smaller\smaller\begin{tabular}{@{}c@{}}%
-185\\221\\-8399
\end{tabular}\endgroup%
\kern3pt%
\begingroup \smaller\smaller\smaller\begin{tabular}{@{}c@{}}%
-400\\480\\-18161
\end{tabular}\endgroup%
{$\left.\llap{\phantom{%
\begingroup \smaller\smaller\smaller\begin{tabular}{@{}c@{}}%
0\\0\\0
\end{tabular}\endgroup%
}}\!\right]$}%
\EasyButWeakLineBreak%
{$\left[\!\llap{\phantom{%
\begingroup \smaller\smaller\smaller\begin{tabular}{@{}c@{}}%
0\\0\\0
\end{tabular}\endgroup%
}}\right.$}%
\begingroup \smaller\smaller\smaller\begin{tabular}{@{}c@{}}%
0\\2\\-1
\end{tabular}\endgroup%
\HardButStrongLineBreak\kern3pt%
\begingroup \smaller\smaller\smaller\begin{tabular}{@{}c@{}}%
-1\\-2\\-44
\end{tabular}\endgroup%
\HardButStrongLineBreak\kern3pt%
\begingroup \smaller\smaller\smaller\begin{tabular}{@{}c@{}}%
2\\-5\\92
\end{tabular}\endgroup%
{$\left.\llap{\phantom{%
\begingroup \smaller\smaller\smaller\begin{tabular}{@{}c@{}}%
0\\0\\0
\end{tabular}\endgroup%
}}\!\right]$}%

\medskip%
%
\leavevmode\llap{}%
$W_{28\phantom{0}}$%
\qquad\llap{8} lattices, $\chi=14$%
\hfill%
$\infty2322$%
\nopagebreak\smallskip\hrule\nopagebreak\medskip%
%
%
\leavevmode%
${L_{28.1}}$%
{} : {$1\above{1pt}{1pt}{-2}{{\rm II}}8\above{1pt}{1pt}{1}{1}{\cdot}1\above{1pt}{1pt}{2}{}13\above{1pt}{1pt}{-}{}$}\spacer%
\instructions{2\rightarrow N_{28}}%
\EasyButWeakLineBreak%
{${26}\above{1pt}{1pt}{4,1}{\infty a}{104}\above{1pt}{1pt}{b}{2}{2}\above{1pt}{1pt}{+}{3}{2}\above{1pt}{1pt}{l}{2}{8}\above{1pt}{1pt}{r}{2}$}%
\nopagebreak\par%
\nopagebreak\par\leavevmode%
{$\left[\!\llap{\phantom{%
\begingroup \smaller\smaller\smaller\begin{tabular}{@{}c@{}}%
0\\0\\0
\end{tabular}\endgroup%
}}\right.$}%
\begingroup \smaller\smaller\smaller\begin{tabular}{@{}c@{}}%
-255855288\\3601312\\-238784
\end{tabular}\endgroup%
\kern3pt%
\begingroup \smaller\smaller\smaller\begin{tabular}{@{}c@{}}%
3601312\\-50690\\3357
\end{tabular}\endgroup%
\kern3pt%
\begingroup \smaller\smaller\smaller\begin{tabular}{@{}c@{}}%
-238784\\3357\\-194
\end{tabular}\endgroup%
{$\left.\llap{\phantom{%
\begingroup \smaller\smaller\smaller\begin{tabular}{@{}c@{}}%
0\\0\\0
\end{tabular}\endgroup%
}}\!\right]$}%
\EasyButWeakLineBreak%
{$\left[\!\llap{\phantom{%
\begingroup \smaller\smaller\smaller\begin{tabular}{@{}c@{}}%
0\\0\\0
\end{tabular}\endgroup%
}}\right.$}%
\begingroup \smaller\smaller\smaller\begin{tabular}{@{}c@{}}%
-1062\\-76154\\-10621
\end{tabular}\endgroup%
\HardButStrongLineBreak\kern3pt%
\begingroup \smaller\smaller\smaller\begin{tabular}{@{}c@{}}%
343\\24596\\3432
\end{tabular}\endgroup%
\HardButStrongLineBreak\kern3pt%
\begingroup \smaller\smaller\smaller\begin{tabular}{@{}c@{}}%
161\\11545\\1610
\end{tabular}\endgroup%
\HardButStrongLineBreak\kern3pt%
\begingroup \smaller\smaller\smaller\begin{tabular}{@{}c@{}}%
-168\\-12047\\-1681
\end{tabular}\endgroup%
\HardButStrongLineBreak\kern3pt%
\begingroup \smaller\smaller\smaller\begin{tabular}{@{}c@{}}%
-795\\-57008\\-7952
\end{tabular}\endgroup%
{$\left.\llap{\phantom{%
\begingroup \smaller\smaller\smaller\begin{tabular}{@{}c@{}}%
0\\0\\0
\end{tabular}\endgroup%
}}\!\right]$}%

\medskip%
%
\leavevmode\llap{}%
$W_{29\phantom{0}}$%
\qquad\llap{6} lattices, $\chi=28$%
\hfill%
$22232223\rtimes C_{2}$%
\nopagebreak\smallskip\hrule\nopagebreak\medskip%
%
%
\leavevmode%
${L_{29.1}}$%
{} : {$1\above{1pt}{1pt}{-2}{{\rm II}}4\above{1pt}{1pt}{1}{1}{\cdot}1\above{1pt}{1pt}{-2}{}29\above{1pt}{1pt}{-}{}$}\spacer%
\instructions{2\rightarrow N_{29}}%
\EasyButWeakLineBreak%
{${2}\above{1pt}{1pt}{b}{2}{58}\above{1pt}{1pt}{l}{2}{4}\above{1pt}{1pt}{r}{2}{2}\above{1pt}{1pt}{-}{3}$}\relax$\,(\times2)$%
\nopagebreak\par%
\nopagebreak\par\leavevmode%
{$\left[\!\llap{\phantom{%
\begingroup \smaller\smaller\smaller\begin{tabular}{@{}c@{}}%
0\\0\\0
\end{tabular}\endgroup%
}}\right.$}%
\begingroup \smaller\smaller\smaller\begin{tabular}{@{}c@{}}%
-42108\\1508\\812
\end{tabular}\endgroup%
\kern3pt%
\begingroup \smaller\smaller\smaller\begin{tabular}{@{}c@{}}%
1508\\-54\\-29
\end{tabular}\endgroup%
\kern3pt%
\begingroup \smaller\smaller\smaller\begin{tabular}{@{}c@{}}%
812\\-29\\-14
\end{tabular}\endgroup%
{$\left.\llap{\phantom{%
\begingroup \smaller\smaller\smaller\begin{tabular}{@{}c@{}}%
0\\0\\0
\end{tabular}\endgroup%
}}\!\right]$}%
\hfil\penalty500%
{$\left[\!\llap{\phantom{%
\begingroup \smaller\smaller\smaller\begin{tabular}{@{}c@{}}%
0\\0\\0
\end{tabular}\endgroup%
}}\right.$}%
\begingroup \smaller\smaller\smaller\begin{tabular}{@{}c@{}}%
3189\\90596\\-3828
\end{tabular}\endgroup%
\kern3pt%
\begingroup \smaller\smaller\smaller\begin{tabular}{@{}c@{}}%
-115\\-3267\\138
\end{tabular}\endgroup%
\kern3pt%
\begingroup \smaller\smaller\smaller\begin{tabular}{@{}c@{}}%
-65\\-1846\\77
\end{tabular}\endgroup%
{$\left.\llap{\phantom{%
\begingroup \smaller\smaller\smaller\begin{tabular}{@{}c@{}}%
0\\0\\0
\end{tabular}\endgroup%
}}\!\right]$}%
\EasyButWeakLineBreak%
{$\left[\!\llap{\phantom{%
\begingroup \smaller\smaller\smaller\begin{tabular}{@{}c@{}}%
0\\0\\0
\end{tabular}\endgroup%
}}\right.$}%
\begingroup \smaller\smaller\smaller\begin{tabular}{@{}c@{}}%
2\\56\\-1
\end{tabular}\endgroup%
\HardButStrongLineBreak\kern3pt%
\begingroup \smaller\smaller\smaller\begin{tabular}{@{}c@{}}%
133\\3770\\-145
\end{tabular}\endgroup%
\HardButStrongLineBreak\kern3pt%
\begingroup \smaller\smaller\smaller\begin{tabular}{@{}c@{}}%
31\\880\\-36
\end{tabular}\endgroup%
\HardButStrongLineBreak\kern3pt%
\begingroup \smaller\smaller\smaller\begin{tabular}{@{}c@{}}%
16\\455\\-20
\end{tabular}\endgroup%
{$\left.\llap{\phantom{%
\begingroup \smaller\smaller\smaller\begin{tabular}{@{}c@{}}%
0\\0\\0
\end{tabular}\endgroup%
}}\!\right]$}%

\medskip%
%
\leavevmode\llap{}%
$W_{30\phantom{0}}$%
\qquad\llap{48} lattices, $\chi=24$%
\hfill%
$\infty22\infty22\rtimes C_{2}$%
\nopagebreak\smallskip\hrule\nopagebreak\medskip%
%
%
\leavevmode%
${L_{30.1}}$%
{} : {$1\above{1pt}{1pt}{-2}{{\rm II}}8\above{1pt}{1pt}{-}{3}{\cdot}1\above{1pt}{1pt}{1}{}3\above{1pt}{1pt}{1}{}9\above{1pt}{1pt}{-}{}{\cdot}1\above{1pt}{1pt}{-}{}5\above{1pt}{1pt}{1}{}25\above{1pt}{1pt}{-}{}$}\spacer%
\instructions{235\rightarrow N_{30},35,25,23,5,3,2}%
\EasyButWeakLineBreak%
{${30}\above{1pt}{1pt}{60,31}{\infty b}{120}\above{1pt}{1pt}{b}{2}{18}\above{1pt}{1pt}{b}{2}{30}\above{1pt}{1pt}{60,19}{\infty a}{120}\above{1pt}{1pt}{b}{2}{450}\above{1pt}{1pt}{b}{2}$}%
\nopagebreak\par%
\nopagebreak\par\leavevmode%
{$\left[\!\llap{\phantom{%
\begingroup \smaller\smaller\smaller\begin{tabular}{@{}c@{}}%
0\\0\\0
\end{tabular}\endgroup%
}}\right.$}%
\begingroup \smaller\smaller\smaller\begin{tabular}{@{}c@{}}%
-1348200\\21600\\-55800
\end{tabular}\endgroup%
\kern3pt%
\begingroup \smaller\smaller\smaller\begin{tabular}{@{}c@{}}%
21600\\-330\\915
\end{tabular}\endgroup%
\kern3pt%
\begingroup \smaller\smaller\smaller\begin{tabular}{@{}c@{}}%
-55800\\915\\-2282
\end{tabular}\endgroup%
{$\left.\llap{\phantom{%
\begingroup \smaller\smaller\smaller\begin{tabular}{@{}c@{}}%
0\\0\\0
\end{tabular}\endgroup%
}}\!\right]$}%
\EasyButWeakLineBreak%
{$\left[\!\llap{\phantom{%
\begingroup \smaller\smaller\smaller\begin{tabular}{@{}c@{}}%
0\\0\\0
\end{tabular}\endgroup%
}}\right.$}%
\begingroup \smaller\smaller\smaller\begin{tabular}{@{}c@{}}%
-15\\-317\\240
\end{tabular}\endgroup%
\HardButStrongLineBreak\kern3pt%
\begingroup \smaller\smaller\smaller\begin{tabular}{@{}c@{}}%
-131\\-2756\\2100
\end{tabular}\endgroup%
\HardButStrongLineBreak\kern3pt%
\begingroup \smaller\smaller\smaller\begin{tabular}{@{}c@{}}%
-23\\-483\\369
\end{tabular}\endgroup%
\HardButStrongLineBreak\kern3pt%
\begingroup \smaller\smaller\smaller\begin{tabular}{@{}c@{}}%
-14\\-293\\225
\end{tabular}\endgroup%
\HardButStrongLineBreak\kern3pt%
\begingroup \smaller\smaller\smaller\begin{tabular}{@{}c@{}}%
15\\316\\-240
\end{tabular}\endgroup%
\HardButStrongLineBreak\kern3pt%
\begingroup \smaller\smaller\smaller\begin{tabular}{@{}c@{}}%
28\\585\\-450
\end{tabular}\endgroup%
{$\left.\llap{\phantom{%
\begingroup \smaller\smaller\smaller\begin{tabular}{@{}c@{}}%
0\\0\\0
\end{tabular}\endgroup%
}}\!\right]$}%

\medskip%
%
\leavevmode\llap{}%
$W_{31\phantom{0}}$%
\qquad\llap{16} lattices, $\chi=4$%
\hfill%
$6222$%
\nopagebreak\smallskip\hrule\nopagebreak\medskip%
%
%
\leavevmode%
${L_{31.1}}$%
{} : {$1\above{1pt}{1pt}{-2}{{\rm II}}8\above{1pt}{1pt}{-}{3}{\cdot}1\above{1pt}{1pt}{2}{}3\above{1pt}{1pt}{-}{}{\cdot}1\above{1pt}{1pt}{-2}{}5\above{1pt}{1pt}{-}{}$}\spacer%
\instructions{2\rightarrow N_{31}}%
\EasyButWeakLineBreak%
{${2}\above{1pt}{1pt}{}{6}{6}\above{1pt}{1pt}{b}{2}{10}\above{1pt}{1pt}{l}{2}{24}\above{1pt}{1pt}{r}{2}$}%
\nopagebreak\par%
\nopagebreak\par\leavevmode%
{$\left[\!\llap{\phantom{%
\begingroup \smaller\smaller\smaller\begin{tabular}{@{}c@{}}%
0\\0\\0
\end{tabular}\endgroup%
}}\right.$}%
\begingroup \smaller\smaller\smaller\begin{tabular}{@{}c@{}}%
-221160\\960\\1800
\end{tabular}\endgroup%
\kern3pt%
\begingroup \smaller\smaller\smaller\begin{tabular}{@{}c@{}}%
960\\-2\\-9
\end{tabular}\endgroup%
\kern3pt%
\begingroup \smaller\smaller\smaller\begin{tabular}{@{}c@{}}%
1800\\-9\\-14
\end{tabular}\endgroup%
{$\left.\llap{\phantom{%
\begingroup \smaller\smaller\smaller\begin{tabular}{@{}c@{}}%
0\\0\\0
\end{tabular}\endgroup%
}}\!\right]$}%
\EasyButWeakLineBreak%
{$\left[\!\llap{\phantom{%
\begingroup \smaller\smaller\smaller\begin{tabular}{@{}c@{}}%
0\\0\\0
\end{tabular}\endgroup%
}}\right.$}%
\begingroup \smaller\smaller\smaller\begin{tabular}{@{}c@{}}%
1\\52\\95
\end{tabular}\endgroup%
\HardButStrongLineBreak\kern3pt%
\begingroup \smaller\smaller\smaller\begin{tabular}{@{}c@{}}%
-1\\-51\\-96
\end{tabular}\endgroup%
\HardButStrongLineBreak\kern3pt%
\begingroup \smaller\smaller\smaller\begin{tabular}{@{}c@{}}%
-2\\-105\\-190
\end{tabular}\endgroup%
\HardButStrongLineBreak\kern3pt%
\begingroup \smaller\smaller\smaller\begin{tabular}{@{}c@{}}%
1\\48\\96
\end{tabular}\endgroup%
{$\left.\llap{\phantom{%
\begingroup \smaller\smaller\smaller\begin{tabular}{@{}c@{}}%
0\\0\\0
\end{tabular}\endgroup%
}}\!\right]$}%

\medskip%
%
\leavevmode\llap{}%
$W_{32\phantom{0}}$%
\qquad\llap{24} lattices, $\chi=20$%
\hfill%
$226226\rtimes C_{2}$%
\nopagebreak\smallskip\hrule\nopagebreak\medskip%
%
%
\leavevmode%
${L_{32.1}}$%
{} : {$1\above{1pt}{1pt}{-2}{{\rm II}}4\above{1pt}{1pt}{-}{5}{\cdot}1\above{1pt}{1pt}{1}{}3\above{1pt}{1pt}{-}{}9\above{1pt}{1pt}{-}{}{\cdot}1\above{1pt}{1pt}{2}{}11\above{1pt}{1pt}{-}{}$}\spacer%
\instructions{23\rightarrow N_{32},3,2}%
\EasyButWeakLineBreak%
{${18}\above{1pt}{1pt}{s}{2}{22}\above{1pt}{1pt}{b}{2}{6}\above{1pt}{1pt}{}{6}$}\relax$\,(\times2)$%
\nopagebreak\par%
\nopagebreak\par\leavevmode%
{$\left[\!\llap{\phantom{%
\begingroup \smaller\smaller\smaller\begin{tabular}{@{}c@{}}%
0\\0\\0
\end{tabular}\endgroup%
}}\right.$}%
\begingroup \smaller\smaller\smaller\begin{tabular}{@{}c@{}}%
-1983564\\-559152\\9108
\end{tabular}\endgroup%
\kern3pt%
\begingroup \smaller\smaller\smaller\begin{tabular}{@{}c@{}}%
-559152\\-157602\\2559
\end{tabular}\endgroup%
\kern3pt%
\begingroup \smaller\smaller\smaller\begin{tabular}{@{}c@{}}%
9108\\2559\\-38
\end{tabular}\endgroup%
{$\left.\llap{\phantom{%
\begingroup \smaller\smaller\smaller\begin{tabular}{@{}c@{}}%
0\\0\\0
\end{tabular}\endgroup%
}}\!\right]$}%
\hfil\penalty500%
{$\left[\!\llap{\phantom{%
\begingroup \smaller\smaller\smaller\begin{tabular}{@{}c@{}}%
0\\0\\0
\end{tabular}\endgroup%
}}\right.$}%
\begingroup \smaller\smaller\smaller\begin{tabular}{@{}c@{}}%
-99001\\364320\\803880
\end{tabular}\endgroup%
\kern3pt%
\begingroup \smaller\smaller\smaller\begin{tabular}{@{}c@{}}%
-27675\\101843\\224721
\end{tabular}\endgroup%
\kern3pt%
\begingroup \smaller\smaller\smaller\begin{tabular}{@{}c@{}}%
350\\-1288\\-2843
\end{tabular}\endgroup%
{$\left.\llap{\phantom{%
\begingroup \smaller\smaller\smaller\begin{tabular}{@{}c@{}}%
0\\0\\0
\end{tabular}\endgroup%
}}\!\right]$}%
\EasyButWeakLineBreak%
{$\left[\!\llap{\phantom{%
\begingroup \smaller\smaller\smaller\begin{tabular}{@{}c@{}}%
0\\0\\0
\end{tabular}\endgroup%
}}\right.$}%
\begingroup \smaller\smaller\smaller\begin{tabular}{@{}c@{}}%
-322\\1185\\2619
\end{tabular}\endgroup%
\HardButStrongLineBreak\kern3pt%
\begingroup \smaller\smaller\smaller\begin{tabular}{@{}c@{}}%
-278\\1023\\2255
\end{tabular}\endgroup%
\HardButStrongLineBreak\kern3pt%
\begingroup \smaller\smaller\smaller\begin{tabular}{@{}c@{}}%
-3\\11\\21
\end{tabular}\endgroup%
{$\left.\llap{\phantom{%
\begingroup \smaller\smaller\smaller\begin{tabular}{@{}c@{}}%
0\\0\\0
\end{tabular}\endgroup%
}}\!\right]$}%

\medskip%
%
\leavevmode\llap{}%
$W_{33\phantom{0}}$%
\qquad\llap{44} lattices, $\chi=15$%
\hfill%
$222224$%
\nopagebreak\smallskip\hrule\nopagebreak\medskip%
%
%
\leavevmode%
${L_{33.1}}$%
{} : {$1\above{1pt}{1pt}{2}{{\rm II}}4\above{1pt}{1pt}{1}{1}{\cdot}1\above{1pt}{1pt}{2}{}3\above{1pt}{1pt}{1}{}{\cdot}1\above{1pt}{1pt}{2}{}11\above{1pt}{1pt}{-}{}$}\spacer%
\instructions{2\rightarrow N_{33}}%
\EasyButWeakLineBreak%
{${4}\above{1pt}{1pt}{*}{2}{12}\above{1pt}{1pt}{b}{2}{22}\above{1pt}{1pt}{l}{2}{4}\above{1pt}{1pt}{r}{2}{66}\above{1pt}{1pt}{b}{2}{2}\above{1pt}{1pt}{*}{4}$}%
\nopagebreak\par%
\nopagebreak\par\leavevmode%
{$\left[\!\llap{\phantom{%
\begingroup \smaller\smaller\smaller\begin{tabular}{@{}c@{}}%
0\\0\\0
\end{tabular}\endgroup%
}}\right.$}%
\begingroup \smaller\smaller\smaller\begin{tabular}{@{}c@{}}%
-292380\\1584\\1848
\end{tabular}\endgroup%
\kern3pt%
\begingroup \smaller\smaller\smaller\begin{tabular}{@{}c@{}}%
1584\\-8\\-11
\end{tabular}\endgroup%
\kern3pt%
\begingroup \smaller\smaller\smaller\begin{tabular}{@{}c@{}}%
1848\\-11\\-10
\end{tabular}\endgroup%
{$\left.\llap{\phantom{%
\begingroup \smaller\smaller\smaller\begin{tabular}{@{}c@{}}%
0\\0\\0
\end{tabular}\endgroup%
}}\!\right]$}%
\EasyButWeakLineBreak%
{$\left[\!\llap{\phantom{%
\begingroup \smaller\smaller\smaller\begin{tabular}{@{}c@{}}%
0\\0\\0
\end{tabular}\endgroup%
}}\right.$}%
\begingroup \smaller\smaller\smaller\begin{tabular}{@{}c@{}}%
-1\\-110\\-64
\end{tabular}\endgroup%
\HardButStrongLineBreak\kern3pt%
\begingroup \smaller\smaller\smaller\begin{tabular}{@{}c@{}}%
-1\\-108\\-66
\end{tabular}\endgroup%
\HardButStrongLineBreak\kern3pt%
\begingroup \smaller\smaller\smaller\begin{tabular}{@{}c@{}}%
4\\440\\253
\end{tabular}\endgroup%
\HardButStrongLineBreak\kern3pt%
\begingroup \smaller\smaller\smaller\begin{tabular}{@{}c@{}}%
3\\328\\192
\end{tabular}\endgroup%
\HardButStrongLineBreak\kern3pt%
\begingroup \smaller\smaller\smaller\begin{tabular}{@{}c@{}}%
20\\2178\\1287
\end{tabular}\endgroup%
\HardButStrongLineBreak\kern3pt%
\begingroup \smaller\smaller\smaller\begin{tabular}{@{}c@{}}%
1\\108\\65
\end{tabular}\endgroup%
{$\left.\llap{\phantom{%
\begingroup \smaller\smaller\smaller\begin{tabular}{@{}c@{}}%
0\\0\\0
\end{tabular}\endgroup%
}}\!\right]$}%
%
%
\hbox{}\par\smallskip%
%
%
\leavevmode%
${L_{33.2}}$%
{} : {$1\above{1pt}{1pt}{2}{2}8\above{1pt}{1pt}{1}{7}{\cdot}1\above{1pt}{1pt}{2}{}3\above{1pt}{1pt}{-}{}{\cdot}1\above{1pt}{1pt}{2}{}11\above{1pt}{1pt}{1}{}$}\spacer%
\instructions{2\rightarrow N'_{22}}%
\EasyButWeakLineBreak%
{${2}\above{1pt}{1pt}{b}{2}{24}\above{1pt}{1pt}{*}{2}{44}\above{1pt}{1pt}{s}{2}{8}\above{1pt}{1pt}{l}{2}{33}\above{1pt}{1pt}{}{2}{1}\above{1pt}{1pt}{}{4}$}%
\nopagebreak\par%
\nopagebreak\par\leavevmode%
{$\left[\!\llap{\phantom{%
\begingroup \smaller\smaller\smaller\begin{tabular}{@{}c@{}}%
0\\0\\0
\end{tabular}\endgroup%
}}\right.$}%
\begingroup \smaller\smaller\smaller\begin{tabular}{@{}c@{}}%
-1159752\\-577368\\7128
\end{tabular}\endgroup%
\kern3pt%
\begingroup \smaller\smaller\smaller\begin{tabular}{@{}c@{}}%
-577368\\-287435\\3548
\end{tabular}\endgroup%
\kern3pt%
\begingroup \smaller\smaller\smaller\begin{tabular}{@{}c@{}}%
7128\\3548\\-43
\end{tabular}\endgroup%
{$\left.\llap{\phantom{%
\begingroup \smaller\smaller\smaller\begin{tabular}{@{}c@{}}%
0\\0\\0
\end{tabular}\endgroup%
}}\!\right]$}%
\EasyButWeakLineBreak%
{$\left[\!\llap{\phantom{%
\begingroup \smaller\smaller\smaller\begin{tabular}{@{}c@{}}%
0\\0\\0
\end{tabular}\endgroup%
}}\right.$}%
\begingroup \smaller\smaller\smaller\begin{tabular}{@{}c@{}}%
37\\-75\\-55
\end{tabular}\endgroup%
\HardButStrongLineBreak\kern3pt%
\begingroup \smaller\smaller\smaller\begin{tabular}{@{}c@{}}%
77\\-156\\-108
\end{tabular}\endgroup%
\HardButStrongLineBreak\kern3pt%
\begingroup \smaller\smaller\smaller\begin{tabular}{@{}c@{}}%
-293\\594\\440
\end{tabular}\endgroup%
\HardButStrongLineBreak\kern3pt%
\begingroup \smaller\smaller\smaller\begin{tabular}{@{}c@{}}%
-223\\452\\328
\end{tabular}\endgroup%
\HardButStrongLineBreak\kern3pt%
\begingroup \smaller\smaller\smaller\begin{tabular}{@{}c@{}}%
-749\\1518\\1089
\end{tabular}\endgroup%
\HardButStrongLineBreak\kern3pt%
\begingroup \smaller\smaller\smaller\begin{tabular}{@{}c@{}}%
-38\\77\\54
\end{tabular}\endgroup%
{$\left.\llap{\phantom{%
\begingroup \smaller\smaller\smaller\begin{tabular}{@{}c@{}}%
0\\0\\0
\end{tabular}\endgroup%
}}\!\right]$}%
%
%
\hbox{}\par\smallskip%
%
%
\leavevmode%
${L_{33.3}}$%
{} : {$1\above{1pt}{1pt}{-2}{2}8\above{1pt}{1pt}{-}{3}{\cdot}1\above{1pt}{1pt}{2}{}3\above{1pt}{1pt}{-}{}{\cdot}1\above{1pt}{1pt}{2}{}11\above{1pt}{1pt}{1}{}$}\spacer%
\instructions{m}%
\EasyButWeakLineBreak%
{${2}\above{1pt}{1pt}{l}{2}{24}\above{1pt}{1pt}{}{2}{11}\above{1pt}{1pt}{r}{2}{8}\above{1pt}{1pt}{s}{2}{132}\above{1pt}{1pt}{*}{2}{4}\above{1pt}{1pt}{*}{4}$}%
\nopagebreak\par%
\nopagebreak\par\leavevmode%
{$\left[\!\llap{\phantom{%
\begingroup \smaller\smaller\smaller\begin{tabular}{@{}c@{}}%
0\\0\\0
\end{tabular}\endgroup%
}}\right.$}%
\begingroup \smaller\smaller\smaller\begin{tabular}{@{}c@{}}%
-24359016\\35904\\78936
\end{tabular}\endgroup%
\kern3pt%
\begingroup \smaller\smaller\smaller\begin{tabular}{@{}c@{}}%
35904\\-49\\-119
\end{tabular}\endgroup%
\kern3pt%
\begingroup \smaller\smaller\smaller\begin{tabular}{@{}c@{}}%
78936\\-119\\-254
\end{tabular}\endgroup%
{$\left.\llap{\phantom{%
\begingroup \smaller\smaller\smaller\begin{tabular}{@{}c@{}}%
0\\0\\0
\end{tabular}\endgroup%
}}\!\right]$}%
\EasyButWeakLineBreak%
{$\left[\!\llap{\phantom{%
\begingroup \smaller\smaller\smaller\begin{tabular}{@{}c@{}}%
0\\0\\0
\end{tabular}\endgroup%
}}\right.$}%
\begingroup \smaller\smaller\smaller\begin{tabular}{@{}c@{}}%
4\\640\\943
\end{tabular}\endgroup%
\HardButStrongLineBreak\kern3pt%
\begingroup \smaller\smaller\smaller\begin{tabular}{@{}c@{}}%
35\\5592\\8256
\end{tabular}\endgroup%
\HardButStrongLineBreak\kern3pt%
\begingroup \smaller\smaller\smaller\begin{tabular}{@{}c@{}}%
29\\4631\\6842
\end{tabular}\endgroup%
\HardButStrongLineBreak\kern3pt%
\begingroup \smaller\smaller\smaller\begin{tabular}{@{}c@{}}%
9\\1436\\2124
\end{tabular}\endgroup%
\HardButStrongLineBreak\kern3pt%
\begingroup \smaller\smaller\smaller\begin{tabular}{@{}c@{}}%
-7\\-1122\\-1650
\end{tabular}\endgroup%
\HardButStrongLineBreak\kern3pt%
\begingroup \smaller\smaller\smaller\begin{tabular}{@{}c@{}}%
-5\\-798\\-1180
\end{tabular}\endgroup%
{$\left.\llap{\phantom{%
\begingroup \smaller\smaller\smaller\begin{tabular}{@{}c@{}}%
0\\0\\0
\end{tabular}\endgroup%
}}\!\right]$}%

\medskip%
%
\leavevmode\llap{}%
$W_{34\phantom{0}}$%
\qquad\llap{88} lattices, $\chi=72$%
\hfill%
$\infty222222\infty222222\rtimes C_{2}$%
\nopagebreak\smallskip\hrule\nopagebreak\medskip%
%
%
\leavevmode%
${L_{34.1}}$%
{} : {$1\above{1pt}{1pt}{2}{{\rm II}}4\above{1pt}{1pt}{1}{1}{\cdot}1\above{1pt}{1pt}{1}{}3\above{1pt}{1pt}{-}{}9\above{1pt}{1pt}{-}{}{\cdot}1\above{1pt}{1pt}{-2}{}11\above{1pt}{1pt}{1}{}$}\spacer%
\instructions{23\rightarrow N_{34},3,2}%
\EasyButWeakLineBreak%
{${132}\above{1pt}{1pt}{3,2}{\infty b}{132}\above{1pt}{1pt}{r}{2}{18}\above{1pt}{1pt}{l}{2}{4}\above{1pt}{1pt}{r}{2}{6}\above{1pt}{1pt}{b}{2}{396}\above{1pt}{1pt}{*}{2}{4}\above{1pt}{1pt}{*}{2}$}\relax$\,(\times2)$%
\nopagebreak\par%
\nopagebreak\par\leavevmode%
{$\left[\!\llap{\phantom{%
\begingroup \smaller\smaller\smaller
\begin{tabular}{@{}c@{}}%
0\\0\\0
\end{tabular}\endgroup%
}}\right.$}%
\begingroup \smaller\smaller\smaller\begin{tabular}{@{}c@{}}%
-6574916700\\52718688\\1769328
\end{tabular}\endgroup%
\kern3pt%
\begingroup \smaller\smaller\smaller\begin{tabular}{@{}c@{}}%
52718688\\-422706\\-14187
\end{tabular}\endgroup%
\kern3pt%
\begingroup \smaller\smaller\smaller\begin{tabular}{@{}c@{}}%
1769328\\-14187\\-476
\end{tabular}\endgroup%
{$\left.\llap{\phantom{%
\begingroup \smaller\smaller\smaller\begin{tabular}{@{}c@{}}%
0\\0\\0
\end{tabular}\endgroup%
}}\!\right]$}%
\hfil\penalty500%
{$\left[\!\llap{\phantom{%
\begingroup \smaller\smaller\smaller\begin{tabular}{@{}c@{}}%
0\\0\\0
\end{tabular}\endgroup%
}}\right.$}%
\begingroup \smaller\smaller\smaller\begin{tabular}{@{}c@{}}%
4245449\\496887468\\971019324
\end{tabular}\endgroup%
\kern3pt%
\begingroup \smaller\smaller\smaller\begin{tabular}{@{}c@{}}%
-34075\\-3988139\\-7793634
\end{tabular}\endgroup%
\kern3pt%
\begingroup \smaller\smaller\smaller\begin{tabular}{@{}c@{}}%
-1125\\-131670\\-257311
\end{tabular}\endgroup%
{$\left.\llap{\phantom{%
\begingroup \smaller\smaller\smaller\begin{tabular}{@{}c@{}}%
0\\0\\0
\end{tabular}\endgroup%
}}\!\right]$}%
\EasyButWeakLineBreak%
{$\left[\!\llap{\phantom{%
\begingroup \smaller\smaller\smaller\begin{tabular}{@{}c@{}}%
0\\0\\0
\end{tabular}\endgroup%
}}\right.$}%
\begingroup \smaller\smaller\smaller\begin{tabular}{@{}c@{}}%
-311\\-36388\\-71478
\end{tabular}\endgroup%
\HardButStrongLineBreak\kern3pt%
\begingroup \smaller\smaller\smaller\begin{tabular}{@{}c@{}}%
3\\352\\660
\end{tabular}\endgroup%
\HardButStrongLineBreak\kern3pt%
\begingroup \smaller\smaller\smaller\begin{tabular}{@{}c@{}}%
31\\3627\\7128
\end{tabular}\endgroup%
\HardButStrongLineBreak\kern3pt%
\begingroup \smaller\smaller\smaller\begin{tabular}{@{}c@{}}%
5\\584\\1180
\end{tabular}\endgroup%
\HardButStrongLineBreak\kern3pt%
\begingroup \smaller\smaller\smaller\begin{tabular}{@{}c@{}}%
-41\\-4801\\-9306
\end{tabular}\endgroup%
\HardButStrongLineBreak\kern3pt%
\begingroup \smaller\smaller\smaller\begin{tabular}{@{}c@{}}%
-3505\\-410322\\-798732
\end{tabular}\endgroup%
\HardButStrongLineBreak\kern3pt%
\begingroup \smaller\smaller\smaller\begin{tabular}{@{}c@{}}%
-209\\-24466\\-47662
\end{tabular}\endgroup%
{$\left.\llap{\phantom{%
\begingroup \smaller\smaller\smaller\begin{tabular}{@{}c@{}}%
0\\0\\0
\end{tabular}\endgroup%
}}\!\right]$}%
%
%
\hbox{}\par\smallskip%
%
%
\leavevmode%
${L_{34.2}}$%
{} : {$1\above{1pt}{1pt}{-2}{2}8\above{1pt}{1pt}{-}{3}{\cdot}1\above{1pt}{1pt}{-}{}3\above{1pt}{1pt}{1}{}9\above{1pt}{1pt}{1}{}{\cdot}1\above{1pt}{1pt}{-2}{}11\above{1pt}{1pt}{-}{}$}\spacer%
\instructions{3m,3,2}%
\EasyButWeakLineBreak%
{${66}\above{1pt}{1pt}{12,11}{\infty a}{264}\above{1pt}{1pt}{s}{2}{36}\above{1pt}{1pt}{s}{2}{8}\above{1pt}{1pt}{l}{2}{3}\above{1pt}{1pt}{}{2}{792}\above{1pt}{1pt}{r}{2}{2}\above{1pt}{1pt}{b}{2}$}\relax$\,(\times2)$%
\nopagebreak\par%
\nopagebreak\par\leavevmode%
{$\left[\!\llap{\phantom{%
\begingroup \smaller\smaller\smaller
\begin{tabular}{@{}c@{}}%
0\\0\\0
\end{tabular}\endgroup%
}}\right.$}%
\begingroup \smaller\smaller\smaller\begin{tabular}{@{}c@{}}%
-2375208\\-49104\\-50688
\end{tabular}\endgroup%
\kern3pt%
\begingroup \smaller\smaller\smaller\begin{tabular}{@{}c@{}}%
-49104\\-1014\\-1047
\end{tabular}\endgroup%
\kern3pt%
\begingroup \smaller\smaller\smaller\begin{tabular}{@{}c@{}}%
-50688\\-1047\\-1081
\end{tabular}\endgroup%
{$\left.\llap{\phantom{%
\begingroup \smaller\smaller\smaller\begin{tabular}{@{}c@{}}%
0\\0\\0
\end{tabular}\endgroup%
}}\!\right]$}%
\hfil\penalty500%
{$\left[\!\llap{\phantom{%
\begingroup \smaller\smaller\smaller\begin{tabular}{@{}c@{}}%
0\\0\\0
\end{tabular}\endgroup%
}}\right.$}%
\begingroup \smaller\smaller\smaller\begin{tabular}{@{}c@{}}%
-1\\92400\\-88704
\end{tabular}\endgroup%
\kern3pt%
\begingroup \smaller\smaller\smaller\begin{tabular}{@{}c@{}}%
0\\1849\\-1776
\end{tabular}\endgroup%
\kern3pt%
\begingroup \smaller\smaller\smaller\begin{tabular}{@{}c@{}}%
0\\1925\\-1849
\end{tabular}\endgroup%
{$\left.\llap{\phantom{%
\begingroup \smaller\smaller\smaller\begin{tabular}{@{}c@{}}%
0\\0\\0
\end{tabular}\endgroup%
}}\!\right]$}%
\EasyButWeakLineBreak%
{$\left[\!\llap{\phantom{%
\begingroup \smaller\smaller\smaller\begin{tabular}{@{}c@{}}%
0\\0\\0
\end{tabular}\endgroup%
}}\right.$}%
\begingroup \smaller\smaller\smaller\begin{tabular}{@{}c@{}}%
-5\\-781\\990
\end{tabular}\endgroup%
\HardButStrongLineBreak\kern3pt%
\begingroup \smaller\smaller\smaller\begin{tabular}{@{}c@{}}%
-1\\-88\\132
\end{tabular}\endgroup%
\HardButStrongLineBreak\kern3pt%
\begingroup \smaller\smaller\smaller\begin{tabular}{@{}c@{}}%
1\\156\\-198
\end{tabular}\endgroup%
\HardButStrongLineBreak\kern3pt%
\begingroup \smaller\smaller\smaller\begin{tabular}{@{}c@{}}%
1\\104\\-148
\end{tabular}\endgroup%
\HardButStrongLineBreak\kern3pt%
\begingroup \smaller\smaller\smaller\begin{tabular}{@{}c@{}}%
1\\53\\-99
\end{tabular}\endgroup%
\HardButStrongLineBreak\kern3pt%
\begingroup \smaller\smaller\smaller\begin{tabular}{@{}c@{}}%
83\\792\\-4752
\end{tabular}\endgroup%
\HardButStrongLineBreak\kern3pt%
\begingroup \smaller\smaller\smaller\begin{tabular}{@{}c@{}}%
2\\-21\\-76
\end{tabular}\endgroup%
{$\left.\llap{\phantom{%
\begingroup \smaller\smaller\smaller\begin{tabular}{@{}c@{}}%
0\\0\\0
\end{tabular}\endgroup%
}}\!\right]$}%
%
%
\hbox{}\par\smallskip%
%
%
\leavevmode%
${L_{34.3}}$%
{} : {$1\above{1pt}{1pt}{2}{2}8\above{1pt}{1pt}{1}{7}{\cdot}1\above{1pt}{1pt}{-}{}3\above{1pt}{1pt}{1}{}9\above{1pt}{1pt}{1}{}{\cdot}1\above{1pt}{1pt}{-2}{}11\above{1pt}{1pt}{-}{}$}\spacer%
\instructions{32\rightarrow N'_{21},3,m}%
\EasyButWeakLineBreak%
{${66}\above{1pt}{1pt}{12,11}{\infty b}{264}\above{1pt}{1pt}{l}{2}{9}\above{1pt}{1pt}{r}{2}{8}\above{1pt}{1pt}{s}{2}{12}\above{1pt}{1pt}{*}{2}{792}\above{1pt}{1pt}{b}{2}{2}\above{1pt}{1pt}{s}{2}$}\relax$\,(\times2)$%
\nopagebreak\par%
\nopagebreak\par\leavevmode%
{$\left[\!\llap{\phantom{%
\begingroup \smaller\smaller\smaller
\begin{tabular}{@{}c@{}}%
0\\0\\0
\end{tabular}\endgroup%
}}\right.$}%
\begingroup \smaller\smaller\smaller\begin{tabular}{@{}c@{}}%
-35211528\\6384312\\-1161864
\end{tabular}\endgroup%
\kern3pt%
\begingroup \smaller\smaller\smaller\begin{tabular}{@{}c@{}}%
6384312\\-1157559\\210660
\end{tabular}\endgroup%
\kern3pt%
\begingroup \smaller\smaller\smaller\begin{tabular}{@{}c@{}}%
-1161864\\210660\\-38335
\end{tabular}\endgroup%
{$\left.\llap{\phantom{%
\begingroup \smaller\smaller\smaller\begin{tabular}{@{}c@{}}%
0\\0\\0
\end{tabular}\endgroup%
}}\!\right]$}%
\hfil\penalty500%
{$\left[\!\llap{\phantom{%
\begingroup \smaller\smaller\smaller\begin{tabular}{@{}c@{}}%
0\\0\\0
\end{tabular}\endgroup%
}}\right.$}%
\begingroup \smaller\smaller\smaller\begin{tabular}{@{}c@{}}%
-2299969\\-13817232\\-6220368
\end{tabular}\endgroup%
\kern3pt%
\begingroup \smaller\smaller\smaller\begin{tabular}{@{}c@{}}%
417252\\2506672\\1128477
\end{tabular}\endgroup%
\kern3pt%
\begingroup \smaller\smaller\smaller\begin{tabular}{@{}c@{}}%
-76428\\-459147\\-206704
\end{tabular}\endgroup%
{$\left.\llap{\phantom{%
\begingroup \smaller\smaller\smaller\begin{tabular}{@{}c@{}}%
0\\0\\0
\end{tabular}\endgroup%
}}\!\right]$}%
\EasyButWeakLineBreak%
{$\left[\!\llap{\phantom{%
\begingroup \smaller\smaller\smaller\begin{tabular}{@{}c@{}}%
0\\0\\0
\end{tabular}\endgroup%
}}\right.$}%
\begingroup \smaller\smaller\smaller\begin{tabular}{@{}c@{}}%
1076\\6457\\2871
\end{tabular}\endgroup%
\HardButStrongLineBreak\kern3pt%
\begingroup \smaller\smaller\smaller\begin{tabular}{@{}c@{}}%
-95\\-572\\-264
\end{tabular}\endgroup%
\HardButStrongLineBreak\kern3pt%
\begingroup \smaller\smaller\smaller\begin{tabular}{@{}c@{}}%
-98\\-588\\-261
\end{tabular}\endgroup%
\HardButStrongLineBreak\kern3pt%
\begingroup \smaller\smaller\smaller\begin{tabular}{@{}c@{}}%
71\\428\\200
\end{tabular}\endgroup%
\HardButStrongLineBreak\kern3pt%
\begingroup \smaller\smaller\smaller\begin{tabular}{@{}c@{}}%
673\\4046\\1836
\end{tabular}\endgroup%
\HardButStrongLineBreak\kern3pt%
\begingroup \smaller\smaller\smaller\begin{tabular}{@{}c@{}}%
46759\\281028\\127116
\end{tabular}\endgroup%
\HardButStrongLineBreak\kern3pt%
\begingroup \smaller\smaller\smaller\begin{tabular}{@{}c@{}}%
1336\\8029\\3629
\end{tabular}\endgroup%
{$\left.\llap{\phantom{%
\begingroup \smaller\smaller\smaller\begin{tabular}{@{}c@{}}%
0\\0\\0
\end{tabular}\endgroup%
}}\!\right]$}%

\medskip%
%
\leavevmode\llap{}%
$W_{35\phantom{0}}$%
\qquad\llap{18} lattices, $\chi=12$%
\hfill%
$222222\rtimes C_{2}$%
\nopagebreak\smallskip\hrule\nopagebreak\medskip%
%
%
\leavevmode%
${L_{35.1}}$%
{} : {$1\above{1pt}{1pt}{-2}{{\rm II}}4\above{1pt}{1pt}{-}{5}{\cdot}1\above{1pt}{1pt}{-}{}3\above{1pt}{1pt}{1}{}9\above{1pt}{1pt}{-}{}{\cdot}1\above{1pt}{1pt}{-2}{}11\above{1pt}{1pt}{1}{}$}\spacer%
\instructions{23\rightarrow N_{35},3,2}%
\EasyButWeakLineBreak%
{${12}\above{1pt}{1pt}{*}{2}{396}\above{1pt}{1pt}{b}{2}{2}\above{1pt}{1pt}{b}{2}{12}\above{1pt}{1pt}{*}{2}{44}\above{1pt}{1pt}{b}{2}{18}\above{1pt}{1pt}{b}{2}$}%
\nopagebreak\par%
\nopagebreak\par\leavevmode%
{$\left[\!\llap{\phantom{%
\begingroup \smaller\smaller\smaller\begin{tabular}{@{}c@{}}%
0\\0\\0
\end{tabular}\endgroup%
}}\right.$}%
\begingroup \smaller\smaller\smaller\begin{tabular}{@{}c@{}}%
-16236\\1584\\0
\end{tabular}\endgroup%
\kern3pt%
\begingroup \smaller\smaller\smaller\begin{tabular}{@{}c@{}}%
1584\\-150\\-3
\end{tabular}\endgroup%
\kern3pt%
\begingroup \smaller\smaller\smaller\begin{tabular}{@{}c@{}}%
0\\-3\\2
\end{tabular}\endgroup%
{$\left.\llap{\phantom{%
\begingroup \smaller\smaller\smaller\begin{tabular}{@{}c@{}}%
0\\0\\0
\end{tabular}\endgroup%
}}\!\right]$}%
\EasyButWeakLineBreak%
{$\left[\!\llap{\phantom{%
\begingroup \smaller\smaller\smaller\begin{tabular}{@{}c@{}}%
0\\0\\0
\end{tabular}\endgroup%
}}\right.$}%
\begingroup \smaller\smaller\smaller\begin{tabular}{@{}c@{}}%
1\\10\\18
\end{tabular}\endgroup%
\HardButStrongLineBreak\kern3pt%
\begingroup \smaller\smaller\smaller\begin{tabular}{@{}c@{}}%
13\\132\\198
\end{tabular}\endgroup%
\HardButStrongLineBreak\kern3pt%
\begingroup \smaller\smaller\smaller\begin{tabular}{@{}c@{}}%
0\\0\\-1
\end{tabular}\endgroup%
\HardButStrongLineBreak\kern3pt%
\begingroup \smaller\smaller\smaller\begin{tabular}{@{}c@{}}%
-5\\-52\\-78
\end{tabular}\endgroup%
\HardButStrongLineBreak\kern3pt%
\begingroup \smaller\smaller\smaller\begin{tabular}{@{}c@{}}%
-19\\-198\\-286
\end{tabular}\endgroup%
\HardButStrongLineBreak\kern3pt%
\begingroup \smaller\smaller\smaller\begin{tabular}{@{}c@{}}%
-2\\-21\\-27
\end{tabular}\endgroup%
{$\left.\llap{\phantom{%
\begingroup \smaller\smaller\smaller\begin{tabular}{@{}c@{}}%
0\\0\\0
\end{tabular}\endgroup%
}}\!\right]$}%

\medskip%
%
\leavevmode\llap{}%
$W_{36\phantom{0}}$%
\qquad\llap{8} lattices, $\chi=36$%
\hfill%
$\infty222\infty222\rtimes C_{2}$%
\nopagebreak\smallskip\hrule\nopagebreak\medskip%
%
%
\leavevmode%
${L_{36.1}}$%
{} : {$1\above{1pt}{1pt}{-2}{{\rm II}}8\above{1pt}{1pt}{-}{5}{\cdot}1\above{1pt}{1pt}{2}{}17\above{1pt}{1pt}{1}{}$}\spacer%
\instructions{2\rightarrow N_{36}}%
\EasyButWeakLineBreak%
{${34}\above{1pt}{1pt}{4,1}{\infty b}{136}\above{1pt}{1pt}{b}{2}{2}\above{1pt}{1pt}{b}{2}{8}\above{1pt}{1pt}{b}{2}$}\relax$\,(\times2)$%
\nopagebreak\par%
\nopagebreak\par\leavevmode%
{$\left[\!\llap{\phantom{%
\begingroup \smaller\smaller\smaller\begin{tabular}{@{}c@{}}%
0\\0\\0
\end{tabular}\endgroup%
}}\right.$}%
\begingroup \smaller\smaller\smaller\begin{tabular}{@{}c@{}}%
-171224\\-40800\\-4352
\end{tabular}\endgroup%
\kern3pt%
\begingroup \smaller\smaller\smaller\begin{tabular}{@{}c@{}}%
-40800\\-9722\\-1037
\end{tabular}\endgroup%
\kern3pt%
\begingroup \smaller\smaller\smaller\begin{tabular}{@{}c@{}}%
-4352\\-1037\\-110
\end{tabular}\endgroup%
{$\left.\llap{\phantom{%
\begingroup \smaller\smaller\smaller\begin{tabular}{@{}c@{}}%
0\\0\\0
\end{tabular}\endgroup%
}}\!\right]$}%
\hfil\penalty500%
{$\left[\!\llap{\phantom{%
\begingroup \smaller\smaller\smaller\begin{tabular}{@{}c@{}}%
0\\0\\0
\end{tabular}\endgroup%
}}\right.$}%
\begingroup \smaller\smaller\smaller\begin{tabular}{@{}c@{}}%
39167\\-164832\\4896
\end{tabular}\endgroup%
\kern3pt%
\begingroup \smaller\smaller\smaller\begin{tabular}{@{}c@{}}%
9336\\-39290\\1167
\end{tabular}\endgroup%
\kern3pt%
\begingroup \smaller\smaller\smaller\begin{tabular}{@{}c@{}}%
984\\-4141\\122
\end{tabular}\endgroup%
{$\left.\llap{\phantom{%
\begingroup \smaller\smaller\smaller\begin{tabular}{@{}c@{}}%
0\\0\\0
\end{tabular}\endgroup%
}}\!\right]$}%
\EasyButWeakLineBreak%
{$\left[\!\llap{\phantom{%
\begingroup \smaller\smaller\smaller\begin{tabular}{@{}c@{}}%
0\\0\\0
\end{tabular}\endgroup%
}}\right.$}%
\begingroup \smaller\smaller\smaller\begin{tabular}{@{}c@{}}%
-32\\136\\-17
\end{tabular}\endgroup%
\HardButStrongLineBreak\kern3pt%
\begingroup \smaller\smaller\smaller\begin{tabular}{@{}c@{}}%
81\\-340\\0
\end{tabular}\endgroup%
\HardButStrongLineBreak\kern3pt%
\begingroup \smaller\smaller\smaller\begin{tabular}{@{}c@{}}%
4\\-17\\2
\end{tabular}\endgroup%
\HardButStrongLineBreak\kern3pt%
\begingroup \smaller\smaller\smaller\begin{tabular}{@{}c@{}}%
-41\\172\\0
\end{tabular}\endgroup%
{$\left.\llap{\phantom{%
\begingroup \smaller\smaller\smaller\begin{tabular}{@{}c@{}}%
0\\0\\0
\end{tabular}\endgroup%
}}\!\right]$}%

\medskip%
%
\leavevmode\llap{}%
$W_{37\phantom{0}}$%
\qquad\llap{12} lattices, $\chi=8$%
\hfill%
$22322$%
\nopagebreak\smallskip\hrule\nopagebreak\medskip%
%
%
\leavevmode%
${L_{37.1}}$%
{} : {$1\above{1pt}{1pt}{-2}{{\rm II}}4\above{1pt}{1pt}{1}{7}{\cdot}1\above{1pt}{1pt}{-2}{}5\above{1pt}{1pt}{1}{}{\cdot}1\above{1pt}{1pt}{-2}{}7\above{1pt}{1pt}{-}{}$}\spacer%
\instructions{2\rightarrow N_{37}}%
\EasyButWeakLineBreak%
{${4}\above{1pt}{1pt}{*}{2}{20}\above{1pt}{1pt}{b}{2}{2}\above{1pt}{1pt}{+}{3}{2}\above{1pt}{1pt}{s}{2}{70}\above{1pt}{1pt}{b}{2}$}%
\nopagebreak\par%
\nopagebreak\par\leavevmode%
{$\left[\!\llap{\phantom{%
\begingroup \smaller\smaller\smaller\begin{tabular}{@{}c@{}}%
0\\0\\0
\end{tabular}\endgroup%
}}\right.$}%
\begingroup \smaller\smaller\smaller\begin{tabular}{@{}c@{}}%
-246820\\3080\\-3360
\end{tabular}\endgroup%
\kern3pt%
\begingroup \smaller\smaller\smaller\begin{tabular}{@{}c@{}}%
3080\\-38\\39
\end{tabular}\endgroup%
\kern3pt%
\begingroup \smaller\smaller\smaller\begin{tabular}{@{}c@{}}%
-3360\\39\\-26
\end{tabular}\endgroup%
{$\left.\llap{\phantom{%
\begingroup \smaller\smaller\smaller\begin{tabular}{@{}c@{}}%
0\\0\\0
\end{tabular}\endgroup%
}}\!\right]$}%
\EasyButWeakLineBreak%
{$\left[\!\llap{\phantom{%
\begingroup \smaller\smaller\smaller\begin{tabular}{@{}c@{}}%
0\\0\\0
\end{tabular}\endgroup%
}}\right.$}%
\begingroup \smaller\smaller\smaller\begin{tabular}{@{}c@{}}%
-3\\-286\\-42
\end{tabular}\endgroup%
\HardButStrongLineBreak\kern3pt%
\begingroup \smaller\smaller\smaller\begin{tabular}{@{}c@{}}%
7\\670\\100
\end{tabular}\endgroup%
\HardButStrongLineBreak\kern3pt%
\begingroup \smaller\smaller\smaller\begin{tabular}{@{}c@{}}%
2\\191\\28
\end{tabular}\endgroup%
\HardButStrongLineBreak\kern3pt%
\begingroup \smaller\smaller\smaller\begin{tabular}{@{}c@{}}%
-3\\-287\\-43
\end{tabular}\endgroup%
\HardButStrongLineBreak\kern3pt%
\begingroup \smaller\smaller\smaller\begin{tabular}{@{}c@{}}%
-37\\-3535\\-525
\end{tabular}\endgroup%
{$\left.\llap{\phantom{%
\begingroup \smaller\smaller\smaller\begin{tabular}{@{}c@{}}%
0\\0\\0
\end{tabular}\endgroup%
}}\!\right]$}%

\medskip%
%
\leavevmode\llap{}%
$W_{38\phantom{0}}$%
\qquad\llap{44} lattices, $\chi=36$%
\hfill%
$2222222222\rtimes C_{2}$%
\nopagebreak\smallskip\hrule\nopagebreak\medskip%
%
%
\leavevmode%
${L_{38.1}}$%
{} : {$1\above{1pt}{1pt}{2}{{\rm II}}4\above{1pt}{1pt}{-}{3}{\cdot}1\above{1pt}{1pt}{-2}{}5\above{1pt}{1pt}{1}{}{\cdot}1\above{1pt}{1pt}{2}{}7\above{1pt}{1pt}{1}{}$}\spacer%
\instructions{2\rightarrow N_{38}}%
\EasyButWeakLineBreak%
{${20}\above{1pt}{1pt}{*}{2}{28}\above{1pt}{1pt}{*}{2}{4}\above{1pt}{1pt}{b}{2}{2}\above{1pt}{1pt}{s}{2}{14}\above{1pt}{1pt}{b}{2}$}\relax$\,(\times2)$%
\nopagebreak\par%
\nopagebreak\par\leavevmode%
{$\left[\!\llap{\phantom{%
\begingroup \smaller\smaller\smaller
\begin{tabular}{@{}c@{}}%
0\\0\\0
\end{tabular}\endgroup%
}}\right.$}%
\begingroup \smaller\smaller\smaller\begin{tabular}{@{}c@{}}%
2380\\140\\-140
\end{tabular}\endgroup%
\kern3pt%
\begingroup \smaller\smaller\smaller\begin{tabular}{@{}c@{}}%
140\\-4\\-3
\end{tabular}\endgroup%
\kern3pt%
\begingroup \smaller\smaller\smaller\begin{tabular}{@{}c@{}}%
-140\\-3\\6
\end{tabular}\endgroup%
{$\left.\llap{\phantom{%
\begingroup \smaller\smaller\smaller\begin{tabular}{@{}c@{}}%
0\\0\\0
\end{tabular}\endgroup%
}}\!\right]$}%
\hfil\penalty500%
{$\left[\!\llap{\phantom{%
\begingroup \smaller\smaller\smaller\begin{tabular}{@{}c@{}}%
0\\0\\0
\end{tabular}\endgroup%
}}\right.$}%
\begingroup \smaller\smaller\smaller\begin{tabular}{@{}c@{}}%
-71\\-700\\-1820
\end{tabular}\endgroup%
\kern3pt%
\begingroup \smaller\smaller\smaller\begin{tabular}{@{}c@{}}%
-11\\-111\\-286
\end{tabular}\endgroup%
\kern3pt%
\begingroup \smaller\smaller\smaller\begin{tabular}{@{}c@{}}%
7\\70\\181
\end{tabular}\endgroup%
{$\left.\llap{\phantom{%
\begingroup \smaller\smaller\smaller\begin{tabular}{@{}c@{}}%
0\\0\\0
\end{tabular}\endgroup%
}}\!\right]$}%
\EasyButWeakLineBreak%
{$\left[\!\llap{\phantom{%
\begingroup \smaller\smaller\smaller\begin{tabular}{@{}c@{}}%
0\\0\\0
\end{tabular}\endgroup%
}}\right.$}%
\begingroup \smaller\smaller\smaller\begin{tabular}{@{}c@{}}%
23\\240\\610
\end{tabular}\endgroup%
\HardButStrongLineBreak\kern3pt%
\begingroup \smaller\smaller\smaller\begin{tabular}{@{}c@{}}%
29\\294\\756
\end{tabular}\endgroup%
\HardButStrongLineBreak\kern3pt%
\begingroup \smaller\smaller\smaller\begin{tabular}{@{}c@{}}%
7\\68\\178
\end{tabular}\endgroup%
\HardButStrongLineBreak\kern3pt%
\begingroup \smaller\smaller\smaller\begin{tabular}{@{}c@{}}%
0\\-2\\-3
\end{tabular}\endgroup%
\HardButStrongLineBreak\kern3pt%
\begingroup \smaller\smaller\smaller\begin{tabular}{@{}c@{}}%
-2\\-28\\-63
\end{tabular}\endgroup%
{$\left.\llap{\phantom{%
\begingroup \smaller\smaller\smaller\begin{tabular}{@{}c@{}}%
0\\0\\0
\end{tabular}\endgroup%
}}\!\right]$}%
%
%
\hbox{}\par\smallskip%
%
%
\leavevmode%
${L_{38.2}}$%
{} : {$1\above{1pt}{1pt}{-2}{6}8\above{1pt}{1pt}{1}{1}{\cdot}1\above{1pt}{1pt}{-2}{}5\above{1pt}{1pt}{-}{}{\cdot}1\above{1pt}{1pt}{2}{}7\above{1pt}{1pt}{1}{}$}\spacer%
\instructions{2\rightarrow N'_{25}}%
\EasyButWeakLineBreak%
{${40}\above{1pt}{1pt}{b}{2}{14}\above{1pt}{1pt}{l}{2}{8}\above{1pt}{1pt}{}{2}{1}\above{1pt}{1pt}{r}{2}{28}\above{1pt}{1pt}{*}{2}$}\relax$\,(\times2)$%
\nopagebreak\par%
\nopagebreak\par\leavevmode%
{$\left[\!\llap{\phantom{%
\begingroup \smaller\smaller\smaller
\begin{tabular}{@{}c@{}}%
0\\0\\0
\end{tabular}\endgroup%
}}\right.$}%
\begingroup \smaller\smaller\smaller\begin{tabular}{@{}c@{}}%
-133560\\560\\840
\end{tabular}\endgroup%
\kern3pt%
\begingroup \smaller\smaller\smaller\begin{tabular}{@{}c@{}}%
560\\-2\\-5
\end{tabular}\endgroup%
\kern3pt%
\begingroup \smaller\smaller\smaller\begin{tabular}{@{}c@{}}%
840\\-5\\1
\end{tabular}\endgroup%
{$\left.\llap{\phantom{%
\begingroup \smaller\smaller\smaller\begin{tabular}{@{}c@{}}%
0\\0\\0
\end{tabular}\endgroup%
}}\!\right]$}%
\hfil\penalty500%
{$\left[\!\llap{\phantom{%
\begingroup \smaller\smaller\smaller\begin{tabular}{@{}c@{}}%
0\\0\\0
\end{tabular}\endgroup%
}}\right.$}%
\begingroup \smaller\smaller\smaller\begin{tabular}{@{}c@{}}%
-2521\\-446040\\-110880
\end{tabular}\endgroup%
\kern3pt%
\begingroup \smaller\smaller\smaller\begin{tabular}{@{}c@{}}%
14\\2477\\616
\end{tabular}\endgroup%
\kern3pt%
\begingroup \smaller\smaller\smaller\begin{tabular}{@{}c@{}}%
1\\177\\43
\end{tabular}\endgroup%
{$\left.\llap{\phantom{%
\begingroup \smaller\smaller\smaller\begin{tabular}{@{}c@{}}%
0\\0\\0
\end{tabular}\endgroup%
}}\!\right]$}%
\EasyButWeakLineBreak%
{$\left[\!\llap{\phantom{%
\begingroup \smaller\smaller\smaller\begin{tabular}{@{}c@{}}%
0\\0\\0
\end{tabular}\endgroup%
}}\right.$}%
\begingroup \smaller\smaller\smaller\begin{tabular}{@{}c@{}}%
-39\\-6900\\-1720
\end{tabular}\endgroup%
\HardButStrongLineBreak\kern3pt%
\begingroup \smaller\smaller\smaller\begin{tabular}{@{}c@{}}%
-29\\-5131\\-1274
\end{tabular}\endgroup%
\HardButStrongLineBreak\kern3pt%
\begingroup \smaller\smaller\smaller\begin{tabular}{@{}c@{}}%
-17\\-3008\\-744
\end{tabular}\endgroup%
\HardButStrongLineBreak\kern3pt%
\begingroup \smaller\smaller\smaller\begin{tabular}{@{}c@{}}%
-1\\-177\\-43
\end{tabular}\endgroup%
\HardButStrongLineBreak\kern3pt%
\begingroup \smaller\smaller\smaller\begin{tabular}{@{}c@{}}%
-3\\-532\\-126
\end{tabular}\endgroup%
{$\left.\llap{\phantom{%
\begingroup \smaller\smaller\smaller\begin{tabular}{@{}c@{}}%
0\\0\\0
\end{tabular}\endgroup%
}}\!\right]$}%
%
%
\hbox{}\par\smallskip%
%
%
\leavevmode%
${L_{38.3}}$%
{} : {$1\above{1pt}{1pt}{2}{6}8\above{1pt}{1pt}{-}{5}{\cdot}1\above{1pt}{1pt}{-2}{}5\above{1pt}{1pt}{-}{}{\cdot}1\above{1pt}{1pt}{2}{}7\above{1pt}{1pt}{1}{}$}\spacer%
\instructions{m}%
\EasyButWeakLineBreak%
{${40}\above{1pt}{1pt}{r}{2}{14}\above{1pt}{1pt}{b}{2}{8}\above{1pt}{1pt}{*}{2}{4}\above{1pt}{1pt}{l}{2}{7}\above{1pt}{1pt}{}{2}$}\relax$\,(\times2)$%
\nopagebreak\par%
\nopagebreak\par\leavevmode%
{$\left[\!\llap{\phantom{%
\begingroup \smaller\smaller\smaller
\begin{tabular}{@{}c@{}}%
0\\0\\0
\end{tabular}\endgroup%
}}\right.$}%
\begingroup \smaller\smaller\smaller\begin{tabular}{@{}c@{}}%
264040\\-2240\\-840
\end{tabular}\endgroup%
\kern3pt%
\begingroup \smaller\smaller\smaller\begin{tabular}{@{}c@{}}%
-2240\\19\\7
\end{tabular}\endgroup%
\kern3pt%
\begingroup \smaller\smaller\smaller\begin{tabular}{@{}c@{}}%
-840\\7\\-2
\end{tabular}\endgroup%
{$\left.\llap{\phantom{%
\begingroup \smaller\smaller\smaller\begin{tabular}{@{}c@{}}%
0\\0\\0
\end{tabular}\endgroup%
}}\!\right]$}%
\hfil\penalty500%
{$\left[\!\llap{\phantom{%
\begingroup \smaller\smaller\smaller\begin{tabular}{@{}c@{}}%
0\\0\\0
\end{tabular}\endgroup%
}}\right.$}%
\begingroup \smaller\smaller\smaller\begin{tabular}{@{}c@{}}%
-36961\\-4385920\\258720
\end{tabular}\endgroup%
\kern3pt%
\begingroup \smaller\smaller\smaller\begin{tabular}{@{}c@{}}%
309\\36667\\-2163
\end{tabular}\endgroup%
\kern3pt%
\begingroup \smaller\smaller\smaller\begin{tabular}{@{}c@{}}%
-42\\-4984\\293
\end{tabular}\endgroup%
{$\left.\llap{\phantom{%
\begingroup \smaller\smaller\smaller\begin{tabular}{@{}c@{}}%
0\\0\\0
\end{tabular}\endgroup%
}}\!\right]$}%
\EasyButWeakLineBreak%
{$\left[\!\llap{\phantom{%
\begingroup \smaller\smaller\smaller\begin{tabular}{@{}c@{}}%
0\\0\\0
\end{tabular}\endgroup%
}}\right.$}%
\begingroup \smaller\smaller\smaller\begin{tabular}{@{}c@{}}%
119\\14120\\-840
\end{tabular}\endgroup%
\HardButStrongLineBreak\kern3pt%
\begingroup \smaller\smaller\smaller\begin{tabular}{@{}c@{}}%
86\\10206\\-595
\end{tabular}\endgroup%
\HardButStrongLineBreak\kern3pt%
\begingroup \smaller\smaller\smaller\begin{tabular}{@{}c@{}}%
49\\5816\\-332
\end{tabular}\endgroup%
\HardButStrongLineBreak\kern3pt%
\begingroup \smaller\smaller\smaller\begin{tabular}{@{}c@{}}%
5\\594\\-30
\end{tabular}\endgroup%
\HardButStrongLineBreak\kern3pt%
\begingroup \smaller\smaller\smaller\begin{tabular}{@{}c@{}}%
3\\357\\-14
\end{tabular}\endgroup%
{$\left.\llap{\phantom{%
\begingroup \smaller\smaller\smaller\begin{tabular}{@{}c@{}}%
0\\0\\0
\end{tabular}\endgroup%
}}\!\right]$}%

\medskip%
%
\leavevmode\llap{}%
$W_{39\phantom{0}}$%
\qquad\llap{12} lattices, $\chi=9$%
\hfill%
$42222$%
\nopagebreak\smallskip\hrule\nopagebreak\medskip%
%
%
\leavevmode%
${L_{39.1}}$%
{} : {$1\above{1pt}{1pt}{-2}{{\rm II}}4\above{1pt}{1pt}{1}{7}{\cdot}1\above{1pt}{1pt}{2}{}5\above{1pt}{1pt}{-}{}{\cdot}1\above{1pt}{1pt}{2}{}7\above{1pt}{1pt}{1}{}$}\spacer%
\instructions{2\rightarrow N_{39}}%
\EasyButWeakLineBreak%
{${2}\above{1pt}{1pt}{*}{4}{4}\above{1pt}{1pt}{b}{2}{14}\above{1pt}{1pt}{s}{2}{10}\above{1pt}{1pt}{l}{2}{28}\above{1pt}{1pt}{r}{2}$}%
\nopagebreak\par%
\nopagebreak\par\leavevmode%
{$\left[\!\llap{\phantom{%
\begingroup \smaller\smaller\smaller\begin{tabular}{@{}c@{}}%
0\\0\\0
\end{tabular}\endgroup%
}}\right.$}%
\begingroup \smaller\smaller\smaller\begin{tabular}{@{}c@{}}%
-691460\\2240\\3640
\end{tabular}\endgroup%
\kern3pt%
\begingroup \smaller\smaller\smaller\begin{tabular}{@{}c@{}}%
2240\\-6\\-13
\end{tabular}\endgroup%
\kern3pt%
\begingroup \smaller\smaller\smaller\begin{tabular}{@{}c@{}}%
3640\\-13\\-18
\end{tabular}\endgroup%
{$\left.\llap{\phantom{%
\begingroup \smaller\smaller\smaller\begin{tabular}{@{}c@{}}%
0\\0\\0
\end{tabular}\endgroup%
}}\!\right]$}%
\EasyButWeakLineBreak%
{$\left[\!\llap{\phantom{%
\begingroup \smaller\smaller\smaller\begin{tabular}{@{}c@{}}%
0\\0\\0
\end{tabular}\endgroup%
}}\right.$}%
\begingroup \smaller\smaller\smaller\begin{tabular}{@{}c@{}}%
1\\115\\119
\end{tabular}\endgroup%
\HardButStrongLineBreak\kern3pt%
\begingroup \smaller\smaller\smaller\begin{tabular}{@{}c@{}}%
-1\\-114\\-120
\end{tabular}\endgroup%
\HardButStrongLineBreak\kern3pt%
\begingroup \smaller\smaller\smaller\begin{tabular}{@{}c@{}}%
-2\\-231\\-238
\end{tabular}\endgroup%
\HardButStrongLineBreak\kern3pt%
\begingroup \smaller\smaller\smaller\begin{tabular}{@{}c@{}}%
2\\225\\240
\end{tabular}\endgroup%
\HardButStrongLineBreak\kern3pt%
\begingroup \smaller\smaller\smaller\begin{tabular}{@{}c@{}}%
15\\1708\\1792
\end{tabular}\endgroup%
{$\left.\llap{\phantom{%
\begingroup \smaller\smaller\smaller\begin{tabular}{@{}c@{}}%
0\\0\\0
\end{tabular}\endgroup%
}}\!\right]$}%

\medskip%
%
\leavevmode\llap{}%
$W_{40\phantom{0}}$%
\qquad\llap{8} lattices, $\chi=40$%
\hfill%
$\infty232\infty232\rtimes C_{2}$%
\nopagebreak\smallskip\hrule\nopagebreak\medskip%
%
%
\leavevmode%
${L_{40.1}}$%
{} : {$1\above{1pt}{1pt}{-2}{{\rm II}}8\above{1pt}{1pt}{1}{7}{\cdot}1\above{1pt}{1pt}{-2}{}19\above{1pt}{1pt}{-}{}$}\spacer%
\instructions{2\rightarrow N_{40}}%
\EasyButWeakLineBreak%
{${38}\above{1pt}{1pt}{4,3}{\infty b}{152}\above{1pt}{1pt}{b}{2}{2}\above{1pt}{1pt}{+}{3}{2}\above{1pt}{1pt}{s}{2}$}\relax$\,(\times2)$%
\nopagebreak\par%
\nopagebreak\par\leavevmode%
{$\left[\!\llap{\phantom{%
\begingroup \smaller\smaller\smaller\begin{tabular}{@{}c@{}}%
0\\0\\0
\end{tabular}\endgroup%
}}\right.$}%
\begingroup \smaller\smaller\smaller\begin{tabular}{@{}c@{}}%
-1205512\\6080\\8968
\end{tabular}\endgroup%
\kern3pt%
\begingroup \smaller\smaller\smaller\begin{tabular}{@{}c@{}}%
6080\\-30\\-47
\end{tabular}\endgroup%
\kern3pt%
\begingroup \smaller\smaller\smaller\begin{tabular}{@{}c@{}}%
8968\\-47\\-62
\end{tabular}\endgroup%
{$\left.\llap{\phantom{%
\begingroup \smaller\smaller\smaller\begin{tabular}{@{}c@{}}%
0\\0\\0
\end{tabular}\endgroup%
}}\!\right]$}%
\hfil\penalty500%
{$\left[\!\llap{\phantom{%
\begingroup \smaller\smaller\smaller\begin{tabular}{@{}c@{}}%
0\\0\\0
\end{tabular}\endgroup%
}}\right.$}%
\begingroup \smaller\smaller\smaller\begin{tabular}{@{}c@{}}%
67031\\8557752\\3202640
\end{tabular}\endgroup%
\kern3pt%
\begingroup \smaller\smaller\smaller\begin{tabular}{@{}c@{}}%
-360\\-45961\\-17200
\end{tabular}\endgroup%
\kern3pt%
\begingroup \smaller\smaller\smaller\begin{tabular}{@{}c@{}}%
-441\\-56301\\-21071
\end{tabular}\endgroup%
{$\left.\llap{\phantom{%
\begingroup \smaller\smaller\smaller\begin{tabular}{@{}c@{}}%
0\\0\\0
\end{tabular}\endgroup%
}}\!\right]$}%
\EasyButWeakLineBreak%
{$\left[\!\llap{\phantom{%
\begingroup \smaller\smaller\smaller\begin{tabular}{@{}c@{}}%
0\\0\\0
\end{tabular}\endgroup%
}}\right.$}%
\begingroup \smaller\smaller\smaller\begin{tabular}{@{}c@{}}%
346\\44175\\16530
\end{tabular}\endgroup%
\HardButStrongLineBreak\kern3pt%
\begingroup \smaller\smaller\smaller\begin{tabular}{@{}c@{}}%
619\\79040\\29564
\end{tabular}\endgroup%
\HardButStrongLineBreak\kern3pt%
\begingroup \smaller\smaller\smaller\begin{tabular}{@{}c@{}}%
2\\256\\95
\end{tabular}\endgroup%
\HardButStrongLineBreak\kern3pt%
\begingroup \smaller\smaller\smaller\begin{tabular}{@{}c@{}}%
-2\\-255\\-96
\end{tabular}\endgroup%
{$\left.\llap{\phantom{%
\begingroup \smaller\smaller\smaller\begin{tabular}{@{}c@{}}%
0\\0\\0
\end{tabular}\endgroup%
}}\!\right]$}%

\medskip%
%
\leavevmode\llap{}%
$W_{41\phantom{0}}$%
\qquad\llap{12} lattices, $\chi=7$%
\hfill%
$2264$%
\nopagebreak\smallskip\hrule\nopagebreak\medskip%
%
%
\leavevmode%
${L_{41.1}}$%
{} : {$1\above{1pt}{1pt}{-2}{{\rm II}}4\above{1pt}{1pt}{-}{3}{\cdot}1\above{1pt}{1pt}{2}{}3\above{1pt}{1pt}{-}{}{\cdot}1\above{1pt}{1pt}{2}{}13\above{1pt}{1pt}{1}{}$}\spacer%
\instructions{2\rightarrow N_{41}}%
\EasyButWeakLineBreak%
{${4}\above{1pt}{1pt}{*}{2}{52}\above{1pt}{1pt}{b}{2}{6}\above{1pt}{1pt}{}{6}{2}\above{1pt}{1pt}{*}{4}$}%
\nopagebreak\par%
\nopagebreak\par\leavevmode%
{$\left[\!\llap{\phantom{%
\begingroup \smaller\smaller\smaller\begin{tabular}{@{}c@{}}%
0\\0\\0
\end{tabular}\endgroup%
}}\right.$}%
\begingroup \smaller\smaller\smaller\begin{tabular}{@{}c@{}}%
-235092\\2496\\4836
\end{tabular}\endgroup%
\kern3pt%
\begingroup \smaller\smaller\smaller\begin{tabular}{@{}c@{}}%
2496\\-26\\-53
\end{tabular}\endgroup%
\kern3pt%
\begingroup \smaller\smaller\smaller\begin{tabular}{@{}c@{}}%
4836\\-53\\-94
\end{tabular}\endgroup%
{$\left.\llap{\phantom{%
\begingroup \smaller\smaller\smaller\begin{tabular}{@{}c@{}}%
0\\0\\0
\end{tabular}\endgroup%
}}\!\right]$}%
\EasyButWeakLineBreak%
{$\left[\!\llap{\phantom{%
\begingroup \smaller\smaller\smaller\begin{tabular}{@{}c@{}}%
0\\0\\0
\end{tabular}\endgroup%
}}\right.$}%
\begingroup \smaller\smaller\smaller\begin{tabular}{@{}c@{}}%
-3\\-178\\-54
\end{tabular}\endgroup%
\HardButStrongLineBreak\kern3pt%
\begingroup \smaller\smaller\smaller\begin{tabular}{@{}c@{}}%
-3\\-182\\-52
\end{tabular}\endgroup%
\HardButStrongLineBreak\kern3pt%
\begingroup \smaller\smaller\smaller\begin{tabular}{@{}c@{}}%
4\\237\\72
\end{tabular}\endgroup%
\HardButStrongLineBreak\kern3pt%
\begingroup \smaller\smaller\smaller\begin{tabular}{@{}c@{}}%
2\\120\\35
\end{tabular}\endgroup%
{$\left.\llap{\phantom{%
\begingroup \smaller\smaller\smaller\begin{tabular}{@{}c@{}}%
0\\0\\0
\end{tabular}\endgroup%
}}\!\right]$}%

\medskip%
%
\leavevmode\llap{}%
$W_{42\phantom{0}}$%
\qquad\llap{104} lattices, $\chi=36$%
\hfill%
$2222222222\rtimes C_{2}$%
\nopagebreak\smallskip\hrule\nopagebreak\medskip%
%
%
\leavevmode%
${L_{42.1}}$%
{} : {$1\above{1pt}{1pt}{-2}{6}8\above{1pt}{1pt}{-}{5}{\cdot}1\above{1pt}{1pt}{2}{}3\above{1pt}{1pt}{1}{}{\cdot}1\above{1pt}{1pt}{-2}{}13\above{1pt}{1pt}{1}{}$}\spacer%
\instructions{2m\rightarrow N_{42},2\rightarrow N'_{26}}%
\EasyButWeakLineBreak%
{${1}\above{1pt}{1pt}{}{2}{13}\above{1pt}{1pt}{r}{2}{12}\above{1pt}{1pt}{*}{2}{8}\above{1pt}{1pt}{*}{2}{156}\above{1pt}{1pt}{l}{2}$}\relax$\,(\times2)$%
\nopagebreak\par%
\nopagebreak\par\leavevmode%
{$\left[\!\llap{\phantom{%
\begingroup \smaller\smaller\smaller
\begin{tabular}{@{}c@{}}%
0\\0\\0
\end{tabular}\endgroup%
}}\right.$}%
\begingroup \smaller\smaller\smaller\begin{tabular}{@{}c@{}}%
-2834520\\5616\\5928
\end{tabular}\endgroup%
\kern3pt%
\begingroup \smaller\smaller\smaller\begin{tabular}{@{}c@{}}%
5616\\-10\\-13
\end{tabular}\endgroup%
\kern3pt%
\begingroup \smaller\smaller\smaller\begin{tabular}{@{}c@{}}%
5928\\-13\\-11
\end{tabular}\endgroup%
{$\left.\llap{\phantom{%
\begingroup \smaller\smaller\smaller\begin{tabular}{@{}c@{}}%
0\\0\\0
\end{tabular}\endgroup%
}}\!\right]$}%
\hfil\penalty500%
{$\left[\!\llap{\phantom{%
\begingroup \smaller\smaller\smaller\begin{tabular}{@{}c@{}}%
0\\0\\0
\end{tabular}\endgroup%
}}\right.$}%
\begingroup \smaller\smaller\smaller\begin{tabular}{@{}c@{}}%
4679\\1212120\\1085760
\end{tabular}\endgroup%
\kern3pt%
\begingroup \smaller\smaller\smaller\begin{tabular}{@{}c@{}}%
-10\\-2591\\-2320
\end{tabular}\endgroup%
\kern3pt%
\begingroup \smaller\smaller\smaller\begin{tabular}{@{}c@{}}%
-9\\-2331\\-2089
\end{tabular}\endgroup%
{$\left.\llap{\phantom{%
\begingroup \smaller\smaller\smaller\begin{tabular}{@{}c@{}}%
0\\0\\0
\end{tabular}\endgroup%
}}\!\right]$}%
\EasyButWeakLineBreak%
{$\left[\!\llap{\phantom{%
\begingroup \smaller\smaller\smaller\begin{tabular}{@{}c@{}}%
0\\0\\0
\end{tabular}\endgroup%
}}\right.$}%
\begingroup \smaller\smaller\smaller\begin{tabular}{@{}c@{}}%
1\\258\\233
\end{tabular}\endgroup%
\HardButStrongLineBreak\kern3pt%
\begingroup \smaller\smaller\smaller\begin{tabular}{@{}c@{}}%
9\\2327\\2093
\end{tabular}\endgroup%
\HardButStrongLineBreak\kern3pt%
\begingroup \smaller\smaller\smaller\begin{tabular}{@{}c@{}}%
7\\1812\\1626
\end{tabular}\endgroup%
\HardButStrongLineBreak\kern3pt%
\begingroup \smaller\smaller\smaller\begin{tabular}{@{}c@{}}%
9\\2332\\2088
\end{tabular}\endgroup%
\HardButStrongLineBreak\kern3pt%
\begingroup \smaller\smaller\smaller\begin{tabular}{@{}c@{}}%
71\\18408\\16458
\end{tabular}\endgroup%
{$\left.\llap{\phantom{%
\begingroup \smaller\smaller\smaller\begin{tabular}{@{}c@{}}%
0\\0\\0
\end{tabular}\endgroup%
}}\!\right]$}%
%
%
\hbox{}\par\smallskip%
%
%
\leavevmode%
${L_{42.2}}$%
{} : {$1\above{1pt}{1pt}{2}{0}8\above{1pt}{1pt}{1}{7}{\cdot}1\above{1pt}{1pt}{2}{}3\above{1pt}{1pt}{1}{}{\cdot}1\above{1pt}{1pt}{-2}{}13\above{1pt}{1pt}{1}{}$}\spacer%
\instructions{m}%
\EasyButWeakLineBreak%
{${1}\above{1pt}{1pt}{r}{2}{52}\above{1pt}{1pt}{*}{2}{12}\above{1pt}{1pt}{s}{2}{8}\above{1pt}{1pt}{s}{2}{156}\above{1pt}{1pt}{*}{2}{4}\above{1pt}{1pt}{l}{2}{13}\above{1pt}{1pt}{}{2}{3}\above{1pt}{1pt}{r}{2}{8}\above{1pt}{1pt}{l}{2}{39}\above{1pt}{1pt}{}{2}$}%
\nopagebreak\par%
\nopagebreak\par\leavevmode%
{$\left[\!\llap{\phantom{%
\begingroup \smaller\smaller\smaller
\begin{tabular}{@{}c@{}}%
0\\0\\0
\end{tabular}\endgroup%
}}\right.$}%
\begingroup \smaller\smaller\smaller\begin{tabular}{@{}c@{}}%
-3004872\\12168\\6240
\end{tabular}\endgroup%
\kern3pt%
\begingroup \smaller\smaller\smaller\begin{tabular}{@{}c@{}}%
12168\\-49\\-26
\end{tabular}\endgroup%
\kern3pt%
\begingroup \smaller\smaller\smaller\begin{tabular}{@{}c@{}}%
6240\\-26\\-11
\end{tabular}\endgroup%
{$\left.\llap{\phantom{%
\begingroup \smaller\smaller\smaller\begin{tabular}{@{}c@{}}%
0\\0\\0
\end{tabular}\endgroup%
}}\!\right]$}%
\EasyButWeakLineBreak%
{$\left[\!\llap{\phantom{%
\begingroup \smaller\smaller\smaller\begin{tabular}{@{}c@{}}%
0\\0\\0
\end{tabular}\endgroup%
}}\right.$}%
\begingroup \smaller\smaller\smaller\begin{tabular}{@{}c@{}}%
2\\414\\155
\end{tabular}\endgroup%
\HardButStrongLineBreak\kern3pt%
\begingroup \smaller\smaller\smaller\begin{tabular}{@{}c@{}}%
31\\6422\\2392
\end{tabular}\endgroup%
\HardButStrongLineBreak\kern3pt%
\begingroup \smaller\smaller\smaller\begin{tabular}{@{}c@{}}%
11\\2280\\846
\end{tabular}\endgroup%
\HardButStrongLineBreak\kern3pt%
\begingroup \smaller\smaller\smaller\begin{tabular}{@{}c@{}}%
13\\2696\\996
\end{tabular}\endgroup%
\HardButStrongLineBreak\kern3pt%
\begingroup \smaller\smaller\smaller\begin{tabular}{@{}c@{}}%
97\\20124\\7410
\end{tabular}\endgroup%
\HardButStrongLineBreak\kern3pt%
\begingroup \smaller\smaller\smaller\begin{tabular}{@{}c@{}}%
5\\1038\\380
\end{tabular}\endgroup%
\HardButStrongLineBreak\kern3pt%
\begingroup \smaller\smaller\smaller\begin{tabular}{@{}c@{}}%
4\\832\\299
\end{tabular}\endgroup%
\HardButStrongLineBreak\kern3pt%
\begingroup \smaller\smaller\smaller\begin{tabular}{@{}c@{}}%
-1\\-207\\-78
\end{tabular}\endgroup%
\HardButStrongLineBreak\kern3pt%
\begingroup \smaller\smaller\smaller\begin{tabular}{@{}c@{}}%
-1\\-208\\-76
\end{tabular}\endgroup%
\HardButStrongLineBreak\kern3pt%
\begingroup \smaller\smaller\smaller\begin{tabular}{@{}c@{}}%
10\\2067\\780
\end{tabular}\endgroup%
{$\left.\llap{\phantom{%
\begingroup \smaller\smaller\smaller\begin{tabular}{@{}c@{}}%
0\\0\\0
\end{tabular}\endgroup%
}}\!\right]$}%
%
%
\hbox{}\par\smallskip%
%
%
\leavevmode%
${L_{42.3}}$%
{} : {$1\above{1pt}{1pt}{2}{6}8\above{1pt}{1pt}{1}{1}{\cdot}1\above{1pt}{1pt}{2}{}3\above{1pt}{1pt}{1}{}{\cdot}1\above{1pt}{1pt}{-2}{}13\above{1pt}{1pt}{1}{}$}\spacer%
\instructions{m}%
\EasyButWeakLineBreak%
{${4}\above{1pt}{1pt}{*}{2}{52}\above{1pt}{1pt}{l}{2}{3}\above{1pt}{1pt}{}{2}{8}\above{1pt}{1pt}{}{2}{39}\above{1pt}{1pt}{r}{2}$}\relax$\,(\times2)$%
\nopagebreak\par%
\nopagebreak\par\leavevmode%
{$\left[\!\llap{\phantom{%
\begingroup \smaller\smaller\smaller
\begin{tabular}{@{}c@{}}%
0\\0\\0
\end{tabular}\endgroup%
}}\right.$}%
\begingroup \smaller\smaller\smaller\begin{tabular}{@{}c@{}}%
-50843832\\-1271088\\49920
\end{tabular}\endgroup%
\kern3pt%
\begingroup \smaller\smaller\smaller\begin{tabular}{@{}c@{}}%
-1271088\\-31777\\1248
\end{tabular}\endgroup%
\kern3pt%
\begingroup \smaller\smaller\smaller\begin{tabular}{@{}c@{}}%
49920\\1248\\-49
\end{tabular}\endgroup%
{$\left.\llap{\phantom{%
\begingroup \smaller\smaller\smaller\begin{tabular}{@{}c@{}}%
0\\0\\0
\end{tabular}\endgroup%
}}\!\right]$}%
\hfil\penalty500%
{$\left[\!\llap{\phantom{%
\begingroup \smaller\smaller\smaller\begin{tabular}{@{}c@{}}%
0\\0\\0
\end{tabular}\endgroup%
}}\right.$}%
\begingroup \smaller\smaller\smaller\begin{tabular}{@{}c@{}}%
47423\\-1873248\\592800
\end{tabular}\endgroup%
\kern3pt%
\begingroup \smaller\smaller\smaller\begin{tabular}{@{}c@{}}%
1186\\-46848\\14825
\end{tabular}\endgroup%
\kern3pt%
\begingroup \smaller\smaller\smaller\begin{tabular}{@{}c@{}}%
-46\\1817\\-576
\end{tabular}\endgroup%
{$\left.\llap{\phantom{%
\begingroup \smaller\smaller\smaller\begin{tabular}{@{}c@{}}%
0\\0\\0
\end{tabular}\endgroup%
}}\!\right]$}%
\EasyButWeakLineBreak%
{$\left[\!\llap{\phantom{%
\begingroup \smaller\smaller\smaller\begin{tabular}{@{}c@{}}%
0\\0\\0
\end{tabular}\endgroup%
}}\right.$}%
\begingroup \smaller\smaller\smaller\begin{tabular}{@{}c@{}}%
7\\-274\\152
\end{tabular}\endgroup%
\HardButStrongLineBreak\kern3pt%
\begingroup \smaller\smaller\smaller\begin{tabular}{@{}c@{}}%
49\\-1924\\910
\end{tabular}\endgroup%
\HardButStrongLineBreak\kern3pt%
\begingroup \smaller\smaller\smaller\begin{tabular}{@{}c@{}}%
8\\-315\\126
\end{tabular}\endgroup%
\HardButStrongLineBreak\kern3pt%
\begingroup \smaller\smaller\smaller\begin{tabular}{@{}c@{}}%
17\\-672\\200
\end{tabular}\endgroup%
\HardButStrongLineBreak\kern3pt%
\begingroup \smaller\smaller\smaller\begin{tabular}{@{}c@{}}%
58\\-2301\\468
\end{tabular}\endgroup%
{$\left.\llap{\phantom{%
\begingroup \smaller\smaller\smaller\begin{tabular}{@{}c@{}}%
0\\0\\0
\end{tabular}\endgroup%
}}\!\right]$}%
%
%
\hbox{}\par\smallskip%
%
%
\leavevmode%
${L_{42.4}}$%
{} : {$[1\above{1pt}{1pt}{-}{}2\above{1pt}{1pt}{1}{}]\above{1pt}{1pt}{}{2}16\above{1pt}{1pt}{-}{5}{\cdot}1\above{1pt}{1pt}{2}{}3\above{1pt}{1pt}{1}{}{\cdot}1\above{1pt}{1pt}{-2}{}13\above{1pt}{1pt}{1}{}$}\spacer%
\instructions{2}%
\EasyButWeakLineBreak%
{${16}\above{1pt}{1pt}{s}{2}{52}\above{1pt}{1pt}{*}{2}{48}\above{1pt}{1pt}{s}{2}{8}\above{1pt}{1pt}{*}{2}{624}\above{1pt}{1pt}{l}{2}{1}\above{1pt}{1pt}{}{2}{208}\above{1pt}{1pt}{}{2}{3}\above{1pt}{1pt}{r}{2}{8}\above{1pt}{1pt}{*}{2}{156}\above{1pt}{1pt}{s}{2}$}%
\nopagebreak\par%
\nopagebreak\par\leavevmode%
{$\left[\!\llap{\phantom{%
\begingroup \smaller\smaller\smaller
\begin{tabular}{@{}c@{}}%
0\\0\\0
\end{tabular}\endgroup%
}}\right.$}%
\begingroup \smaller\smaller\smaller\begin{tabular}{@{}c@{}}%
341328\\169728\\-624
\end{tabular}\endgroup%
\kern3pt%
\begingroup \smaller\smaller\smaller\begin{tabular}{@{}c@{}}%
169728\\84398\\-310
\end{tabular}\endgroup%
\kern3pt%
\begingroup \smaller\smaller\smaller\begin{tabular}{@{}c@{}}%
-624\\-310\\1
\end{tabular}\endgroup%
{$\left.\llap{\phantom{%
\begingroup \smaller\smaller\smaller\begin{tabular}{@{}c@{}}%
0\\0\\0
\end{tabular}\endgroup%
}}\!\right]$}%
\EasyButWeakLineBreak%
{$\left[\!\llap{\phantom{%
\begingroup \smaller\smaller\smaller\begin{tabular}{@{}c@{}}%
0\\0\\0
\end{tabular}\endgroup%
}}\right.$}%
\begingroup \smaller\smaller\smaller\begin{tabular}{@{}c@{}}%
201\\-408\\-1016
\end{tabular}\endgroup%
\HardButStrongLineBreak\kern3pt%
\begingroup \smaller\smaller\smaller\begin{tabular}{@{}c@{}}%
333\\-676\\-1690
\end{tabular}\endgroup%
\HardButStrongLineBreak\kern3pt%
\begingroup \smaller\smaller\smaller\begin{tabular}{@{}c@{}}%
65\\-132\\-336
\end{tabular}\endgroup%
\HardButStrongLineBreak\kern3pt%
\begingroup \smaller\smaller\smaller\begin{tabular}{@{}c@{}}%
-1\\2\\0
\end{tabular}\endgroup%
\HardButStrongLineBreak\kern3pt%
\begingroup \smaller\smaller\smaller\begin{tabular}{@{}c@{}}%
-77\\156\\312
\end{tabular}\endgroup%
\HardButStrongLineBreak\kern3pt%
\begingroup \smaller\smaller\smaller\begin{tabular}{@{}c@{}}%
0\\0\\-1
\end{tabular}\endgroup%
\HardButStrongLineBreak\kern3pt%
\begingroup \smaller\smaller\smaller\begin{tabular}{@{}c@{}}%
205\\-416\\-1040
\end{tabular}\endgroup%
\HardButStrongLineBreak\kern3pt%
\begingroup \smaller\smaller\smaller\begin{tabular}{@{}c@{}}%
34\\-69\\-171
\end{tabular}\endgroup%
\HardButStrongLineBreak\kern3pt%
\begingroup \smaller\smaller\smaller\begin{tabular}{@{}c@{}}%
135\\-274\\-680
\end{tabular}\endgroup%
\HardButStrongLineBreak\kern3pt%
\begingroup \smaller\smaller\smaller\begin{tabular}{@{}c@{}}%
1345\\-2730\\-6786
\end{tabular}\endgroup%
{$\left.\llap{\phantom{%
\begingroup \smaller\smaller\smaller\begin{tabular}{@{}c@{}}%
0\\0\\0
\end{tabular}\endgroup%
}}\!\right]$}%
%
%
\hbox{}\par\smallskip%
%
%
\leavevmode%
${L_{42.5}}$%
{} : {$[1\above{1pt}{1pt}{1}{}2\above{1pt}{1pt}{1}{}]\above{1pt}{1pt}{}{0}16\above{1pt}{1pt}{1}{7}{\cdot}1\above{1pt}{1pt}{2}{}3\above{1pt}{1pt}{1}{}{\cdot}1\above{1pt}{1pt}{-2}{}13\above{1pt}{1pt}{1}{}$}\spacer%
\instructions{m}%
\EasyButWeakLineBreak%
{${16}\above{1pt}{1pt}{*}{2}{52}\above{1pt}{1pt}{s}{2}{48}\above{1pt}{1pt}{l}{2}{2}\above{1pt}{1pt}{}{2}{624}\above{1pt}{1pt}{}{2}{1}\above{1pt}{1pt}{r}{2}{208}\above{1pt}{1pt}{*}{2}{12}\above{1pt}{1pt}{l}{2}{2}\above{1pt}{1pt}{}{2}{39}\above{1pt}{1pt}{r}{2}$}%
\nopagebreak\par%
\nopagebreak\par\leavevmode%
{$\left[\!\llap{\phantom{%
\begingroup \smaller\smaller\smaller
\begin{tabular}{@{}c@{}}%
0\\0\\0
\end{tabular}\endgroup%
}}\right.$}%
\begingroup \smaller\smaller\smaller\begin{tabular}{@{}c@{}}%
-7951632\\7488\\18096
\end{tabular}\endgroup%
\kern3pt%
\begingroup \smaller\smaller\smaller\begin{tabular}{@{}c@{}}%
7488\\-2\\-18
\end{tabular}\endgroup%
\kern3pt%
\begingroup \smaller\smaller\smaller\begin{tabular}{@{}c@{}}%
18096\\-18\\-41
\end{tabular}\endgroup%
{$\left.\llap{\phantom{%
\begingroup \smaller\smaller\smaller\begin{tabular}{@{}c@{}}%
0\\0\\0
\end{tabular}\endgroup%
}}\!\right]$}%
\EasyButWeakLineBreak%
{$\left[\!\llap{\phantom{%
\begingroup \smaller\smaller\smaller\begin{tabular}{@{}c@{}}%
0\\0\\0
\end{tabular}\endgroup%
}}\right.$}%
\begingroup \smaller\smaller\smaller\begin{tabular}{@{}c@{}}%
-1\\-76\\-408
\end{tabular}\endgroup%
\HardButStrongLineBreak\kern3pt%
\begingroup \smaller\smaller\smaller\begin{tabular}{@{}c@{}}%
-3\\-234\\-1222
\end{tabular}\endgroup%
\HardButStrongLineBreak\kern3pt%
\begingroup \smaller\smaller\smaller\begin{tabular}{@{}c@{}}%
1\\72\\408
\end{tabular}\endgroup%
\HardButStrongLineBreak\kern3pt%
\begingroup \smaller\smaller\smaller\begin{tabular}{@{}c@{}}%
3\\229\\1222
\end{tabular}\endgroup%
\HardButStrongLineBreak\kern3pt%
\begingroup \smaller\smaller\smaller\begin{tabular}{@{}c@{}}%
167\\12792\\68016
\end{tabular}\endgroup%
\HardButStrongLineBreak\kern3pt%
\begingroup \smaller\smaller\smaller\begin{tabular}{@{}c@{}}%
4\\307\\1629
\end{tabular}\endgroup%
\HardButStrongLineBreak\kern3pt%
\begingroup \smaller\smaller\smaller\begin{tabular}{@{}c@{}}%
71\\5460\\28912
\end{tabular}\endgroup%
\HardButStrongLineBreak\kern3pt%
\begingroup \smaller\smaller\smaller\begin{tabular}{@{}c@{}}%
7\\540\\2850
\end{tabular}\endgroup%
\HardButStrongLineBreak\kern3pt%
\begingroup \smaller\smaller\smaller\begin{tabular}{@{}c@{}}%
2\\155\\814
\end{tabular}\endgroup%
\HardButStrongLineBreak\kern3pt%
\begingroup \smaller\smaller\smaller\begin{tabular}{@{}c@{}}%
7\\546\\2847
\end{tabular}\endgroup%
{$\left.\llap{\phantom{%
\begingroup \smaller\smaller\smaller\begin{tabular}{@{}c@{}}%
0\\0\\0
\end{tabular}\endgroup%
}}\!\right]$}%
%
%
\hbox{}\par\smallskip%
%
%
\leavevmode%
${L_{42.6}}$%
{} : {$[1\above{1pt}{1pt}{-}{}2\above{1pt}{1pt}{1}{}]\above{1pt}{1pt}{}{4}16\above{1pt}{1pt}{-}{3}{\cdot}1\above{1pt}{1pt}{2}{}3\above{1pt}{1pt}{1}{}{\cdot}1\above{1pt}{1pt}{-2}{}13\above{1pt}{1pt}{1}{}$}\EasyButWeakLineBreak%
{${16}\above{1pt}{1pt}{l}{2}{13}\above{1pt}{1pt}{}{2}{48}\above{1pt}{1pt}{}{2}{2}\above{1pt}{1pt}{r}{2}{624}\above{1pt}{1pt}{s}{2}{4}\above{1pt}{1pt}{*}{2}{208}\above{1pt}{1pt}{l}{2}{3}\above{1pt}{1pt}{}{2}{2}\above{1pt}{1pt}{r}{2}{156}\above{1pt}{1pt}{*}{2}$}%
\nopagebreak\par%
\nopagebreak\par\leavevmode%
{$\left[\!\llap{\phantom{%
\begingroup \smaller\smaller\smaller
\begin{tabular}{@{}c@{}}%
0\\0\\0
\end{tabular}\endgroup%
}}\right.$}%
\begingroup \smaller\smaller\smaller\begin{tabular}{@{}c@{}}%
-501072\\1872\\2496
\end{tabular}\endgroup%
\kern3pt%
\begingroup \smaller\smaller\smaller\begin{tabular}{@{}c@{}}%
1872\\-2\\-12
\end{tabular}\endgroup%
\kern3pt%
\begingroup \smaller\smaller\smaller\begin{tabular}{@{}c@{}}%
2496\\-12\\-11
\end{tabular}\endgroup%
{$\left.\llap{\phantom{%
\begingroup \smaller\smaller\smaller\begin{tabular}{@{}c@{}}%
0\\0\\0
\end{tabular}\endgroup%
}}\!\right]$}%
\EasyButWeakLineBreak%
{$\left[\!\llap{\phantom{%
\begingroup \smaller\smaller\smaller\begin{tabular}{@{}c@{}}%
0\\0\\0
\end{tabular}\endgroup%
}}\right.$}%
\begingroup \smaller\smaller\smaller\begin{tabular}{@{}c@{}}%
-1\\-76\\-144
\end{tabular}\endgroup%
\HardButStrongLineBreak\kern3pt%
\begingroup \smaller\smaller\smaller\begin{tabular}{@{}c@{}}%
-1\\-78\\-143
\end{tabular}\endgroup%
\HardButStrongLineBreak\kern3pt%
\begingroup \smaller\smaller\smaller\begin{tabular}{@{}c@{}}%
1\\72\\144
\end{tabular}\endgroup%
\HardButStrongLineBreak\kern3pt%
\begingroup \smaller\smaller\smaller\begin{tabular}{@{}c@{}}%
2\\151\\286
\end{tabular}\endgroup%
\HardButStrongLineBreak\kern3pt%
\begingroup \smaller\smaller\smaller\begin{tabular}{@{}c@{}}%
107\\8112\\15288
\end{tabular}\endgroup%
\HardButStrongLineBreak\kern3pt%
\begingroup \smaller\smaller\smaller\begin{tabular}{@{}c@{}}%
5\\380\\714
\end{tabular}\endgroup%
\HardButStrongLineBreak\kern3pt%
\begingroup \smaller\smaller\smaller\begin{tabular}{@{}c@{}}%
43\\3276\\6136
\end{tabular}\endgroup%
\HardButStrongLineBreak\kern3pt%
\begingroup \smaller\smaller\smaller\begin{tabular}{@{}c@{}}%
2\\153\\285
\end{tabular}\endgroup%
\HardButStrongLineBreak\kern3pt%
\begingroup \smaller\smaller\smaller\begin{tabular}{@{}c@{}}%
1\\77\\142
\end{tabular}\endgroup%
\HardButStrongLineBreak\kern3pt%
\begingroup \smaller\smaller\smaller\begin{tabular}{@{}c@{}}%
5\\390\\702
\end{tabular}\endgroup%
{$\left.\llap{\phantom{%
\begingroup \smaller\smaller\smaller\begin{tabular}{@{}c@{}}%
0\\0\\0
\end{tabular}\endgroup%
}}\!\right]$}%
%
%
\hbox{}\par\smallskip%
%
%
\leavevmode%
${L_{42.7}}$%
{} : {$[1\above{1pt}{1pt}{1}{}2\above{1pt}{1pt}{1}{}]\above{1pt}{1pt}{}{6}16\above{1pt}{1pt}{1}{1}{\cdot}1\above{1pt}{1pt}{2}{}3\above{1pt}{1pt}{1}{}{\cdot}1\above{1pt}{1pt}{-2}{}13\above{1pt}{1pt}{1}{}$}\spacer%
\instructions{m}%
\EasyButWeakLineBreak%
{${16}\above{1pt}{1pt}{}{2}{13}\above{1pt}{1pt}{r}{2}{48}\above{1pt}{1pt}{*}{2}{8}\above{1pt}{1pt}{s}{2}{624}\above{1pt}{1pt}{*}{2}{4}\above{1pt}{1pt}{s}{2}{208}\above{1pt}{1pt}{s}{2}{12}\above{1pt}{1pt}{*}{2}{8}\above{1pt}{1pt}{l}{2}{39}\above{1pt}{1pt}{}{2}$}%
\nopagebreak\par%
\nopagebreak\par\leavevmode%
{$\left[\!\llap{\phantom{%
\begingroup \smaller\smaller\smaller
\begin{tabular}{@{}c@{}}%
0\\0\\0
\end{tabular}\endgroup%
}}\right.$}%
\begingroup \smaller\smaller\smaller\begin{tabular}{@{}c@{}}%
243984\\-1872\\1872
\end{tabular}\endgroup%
\kern3pt%
\begingroup \smaller\smaller\smaller\begin{tabular}{@{}c@{}}%
-1872\\14\\-12
\end{tabular}\endgroup%
\kern3pt%
\begingroup \smaller\smaller\smaller\begin{tabular}{@{}c@{}}%
1872\\-12\\-1
\end{tabular}\endgroup%
{$\left.\llap{\phantom{%
\begingroup \smaller\smaller\smaller\begin{tabular}{@{}c@{}}%
0\\0\\0
\end{tabular}\endgroup%
}}\!\right]$}%
\EasyButWeakLineBreak%
{$\left[\!\llap{\phantom{%
\begingroup \smaller\smaller\smaller\begin{tabular}{@{}c@{}}%
0\\0\\0
\end{tabular}\endgroup%
}}\right.$}%
\begingroup \smaller\smaller\smaller\begin{tabular}{@{}c@{}}%
7\\1072\\160
\end{tabular}\endgroup%
\HardButStrongLineBreak\kern3pt%
\begingroup \smaller\smaller\smaller\begin{tabular}{@{}c@{}}%
4\\611\\91
\end{tabular}\endgroup%
\HardButStrongLineBreak\kern3pt%
\begingroup \smaller\smaller\smaller\begin{tabular}{@{}c@{}}%
-1\\-156\\-24
\end{tabular}\endgroup%
\HardButStrongLineBreak\kern3pt%
\begingroup \smaller\smaller\smaller\begin{tabular}{@{}c@{}}%
-1\\-154\\-24
\end{tabular}\endgroup%
\HardButStrongLineBreak\kern3pt%
\begingroup \smaller\smaller\smaller\begin{tabular}{@{}c@{}}%
1\\156\\0
\end{tabular}\endgroup%
\HardButStrongLineBreak\kern3pt%
\begingroup \smaller\smaller\smaller\begin{tabular}{@{}c@{}}%
1\\154\\22
\end{tabular}\endgroup%
\HardButStrongLineBreak\kern3pt%
\begingroup \smaller\smaller\smaller\begin{tabular}{@{}c@{}}%
23\\3536\\520
\end{tabular}\endgroup%
\HardButStrongLineBreak\kern3pt%
\begingroup \smaller\smaller\smaller\begin{tabular}{@{}c@{}}%
5\\768\\114
\end{tabular}\endgroup%
\HardButStrongLineBreak\kern3pt%
\begingroup \smaller\smaller\smaller\begin{tabular}{@{}c@{}}%
7\\1074\\160
\end{tabular}\endgroup%
\HardButStrongLineBreak\kern3pt%
\begingroup \smaller\smaller\smaller\begin{tabular}{@{}c@{}}%
29\\4446\\663
\end{tabular}\endgroup%
{$\left.\llap{\phantom{%
\begingroup \smaller\smaller\smaller\begin{tabular}{@{}c@{}}%
0\\0\\0
\end{tabular}\endgroup%
}}\!\right]$}%

\medskip%
%
\leavevmode\llap{}%
$W_{43\phantom{0}}$%
\qquad\llap{24} lattices, $\chi=24$%
\hfill%
$22222222\rtimes C_{2}$%
\nopagebreak\smallskip\hrule\nopagebreak\medskip%
%
%
\leavevmode%
${L_{43.1}}$%
{} : {$1\above{1pt}{1pt}{-2}{{\rm II}}4\above{1pt}{1pt}{-}{3}{\cdot}1\above{1pt}{1pt}{-}{}3\above{1pt}{1pt}{1}{}9\above{1pt}{1pt}{1}{}{\cdot}1\above{1pt}{1pt}{-2}{}13\above{1pt}{1pt}{-}{}$}\spacer%
\instructions{23\rightarrow N_{43},3,2}%
\EasyButWeakLineBreak%
{${12}\above{1pt}{1pt}{r}{2}{26}\above{1pt}{1pt}{b}{2}{36}\above{1pt}{1pt}{b}{2}{2}\above{1pt}{1pt}{l}{2}$}\relax$\,(\times2)$%
\nopagebreak\par%
\nopagebreak\par\leavevmode%
{$\left[\!\llap{\phantom{%
\begingroup \smaller\smaller\smaller\begin{tabular}{@{}c@{}}%
0\\0\\0
\end{tabular}\endgroup%
}}\right.$}%
\begingroup \smaller\smaller\smaller\begin{tabular}{@{}c@{}}%
-805428\\2808\\4212
\end{tabular}\endgroup%
\kern3pt%
\begingroup \smaller\smaller\smaller\begin{tabular}{@{}c@{}}%
2808\\-6\\-15
\end{tabular}\endgroup%
\kern3pt%
\begingroup \smaller\smaller\smaller\begin{tabular}{@{}c@{}}%
4212\\-15\\-22
\end{tabular}\endgroup%
{$\left.\llap{\phantom{%
\begingroup \smaller\smaller\smaller\begin{tabular}{@{}c@{}}%
0\\0\\0
\end{tabular}\endgroup%
}}\!\right]$}%
\hfil\penalty500%
{$\left[\!\llap{\phantom{%
\begingroup \smaller\smaller\smaller\begin{tabular}{@{}c@{}}%
0\\0\\0
\end{tabular}\endgroup%
}}\right.$}%
\begingroup \smaller\smaller\smaller\begin{tabular}{@{}c@{}}%
-1561\\-21840\\-285480
\end{tabular}\endgroup%
\kern3pt%
\begingroup \smaller\smaller\smaller\begin{tabular}{@{}c@{}}%
7\\97\\1281
\end{tabular}\endgroup%
\kern3pt%
\begingroup \smaller\smaller\smaller\begin{tabular}{@{}c@{}}%
8\\112\\1463
\end{tabular}\endgroup%
{$\left.\llap{\phantom{%
\begingroup \smaller\smaller\smaller\begin{tabular}{@{}c@{}}%
0\\0\\0
\end{tabular}\endgroup%
}}\!\right]$}%
\EasyButWeakLineBreak%
{$\left[\!\llap{\phantom{%
\begingroup \smaller\smaller\smaller\begin{tabular}{@{}c@{}}%
0\\0\\0
\end{tabular}\endgroup%
}}\right.$}%
\begingroup \smaller\smaller\smaller\begin{tabular}{@{}c@{}}%
-9\\-128\\-1644
\end{tabular}\endgroup%
\HardButStrongLineBreak\kern3pt%
\begingroup \smaller\smaller\smaller\begin{tabular}{@{}c@{}}%
-14\\-195\\-2561
\end{tabular}\endgroup%
\HardButStrongLineBreak\kern3pt%
\begingroup \smaller\smaller\smaller\begin{tabular}{@{}c@{}}%
-5\\-66\\-918
\end{tabular}\endgroup%
\HardButStrongLineBreak\kern3pt%
\begingroup \smaller\smaller\smaller\begin{tabular}{@{}c@{}}%
0\\1\\-1
\end{tabular}\endgroup%
{$\left.\llap{\phantom{%
\begingroup \smaller\smaller\smaller\begin{tabular}{@{}c@{}}%
0\\0\\0
\end{tabular}\endgroup%
}}\!\right]$}%

\medskip%
%
\leavevmode\llap{}%
$W_{44\phantom{0}}$%
\qquad\llap{32} lattices, $\chi=16$%
\hfill%
$\infty2622$%
\nopagebreak\smallskip\hrule\nopagebreak\medskip%
%
%
\leavevmode%
${L_{44.1}}$%
{} : {$1\above{1pt}{1pt}{-2}{{\rm II}}8\above{1pt}{1pt}{1}{1}{\cdot}1\above{1pt}{1pt}{1}{}3\above{1pt}{1pt}{-}{}9\above{1pt}{1pt}{-}{}{\cdot}1\above{1pt}{1pt}{-2}{}7\above{1pt}{1pt}{-}{}$}\spacer%
\instructions{23\rightarrow N_{44},3,2}%
\EasyButWeakLineBreak%
{${42}\above{1pt}{1pt}{12,1}{\infty a}{168}\above{1pt}{1pt}{b}{2}{18}\above{1pt}{1pt}{}{6}{6}\above{1pt}{1pt}{l}{2}{72}\above{1pt}{1pt}{r}{2}$}%
\nopagebreak\par%
\nopagebreak\par\leavevmode%
{$\left[\!\llap{\phantom{%
\begingroup \smaller\smaller\smaller\begin{tabular}{@{}c@{}}%
0\\0\\0
\end{tabular}\endgroup%
}}\right.$}%
\begingroup \smaller\smaller\smaller\begin{tabular}{@{}c@{}}%
-1218168\\-337176\\-519624
\end{tabular}\endgroup%
\kern3pt%
\begingroup \smaller\smaller\smaller\begin{tabular}{@{}c@{}}%
-337176\\-91686\\-159513
\end{tabular}\endgroup%
\kern3pt%
\begingroup \smaller\smaller\smaller\begin{tabular}{@{}c@{}}%
-519624\\-159513\\-71678
\end{tabular}\endgroup%
{$\left.\llap{\phantom{%
\begingroup \smaller\smaller\smaller\begin{tabular}{@{}c@{}}%
0\\0\\0
\end{tabular}\endgroup%
}}\!\right]$}%
\EasyButWeakLineBreak%
{$\left[\!\llap{\phantom{%
\begingroup \smaller\smaller\smaller\begin{tabular}{@{}c@{}}%
0\\0\\0
\end{tabular}\endgroup%
}}\right.$}%
\begingroup \smaller\smaller\smaller\begin{tabular}{@{}c@{}}%
42396\\-131908\\-13797
\end{tabular}\endgroup%
\HardButStrongLineBreak\kern3pt%
\begingroup \smaller\smaller\smaller\begin{tabular}{@{}c@{}}%
-28393\\88340\\9240
\end{tabular}\endgroup%
\HardButStrongLineBreak\kern3pt%
\begingroup \smaller\smaller\smaller\begin{tabular}{@{}c@{}}%
-12611\\39237\\4104
\end{tabular}\endgroup%
\HardButStrongLineBreak\kern3pt%
\begingroup \smaller\smaller\smaller\begin{tabular}{@{}c@{}}%
13671\\-42535\\-4449
\end{tabular}\endgroup%
\HardButStrongLineBreak\kern3pt%
\begingroup \smaller\smaller\smaller\begin{tabular}{@{}c@{}}%
131419\\-408888\\-42768
\end{tabular}\endgroup%
{$\left.\llap{\phantom{%
\begingroup \smaller\smaller\smaller\begin{tabular}{@{}c@{}}%
0\\0\\0
\end{tabular}\endgroup%
}}\!\right]$}%

\medskip%
%
\leavevmode\llap{}%
$W_{45\phantom{0}}$%
\qquad\llap{24} lattices, $\chi=12$%
\hfill%
$222222\rtimes C_{2}$%
\nopagebreak\smallskip\hrule\nopagebreak\medskip%
%
%
\leavevmode%
${L_{45.1}}$%
{} : {$1\above{1pt}{1pt}{-2}{{\rm II}}8\above{1pt}{1pt}{1}{1}{\cdot}1\above{1pt}{1pt}{-}{}3\above{1pt}{1pt}{1}{}9\above{1pt}{1pt}{-}{}{\cdot}1\above{1pt}{1pt}{2}{}7\above{1pt}{1pt}{1}{}$}\spacer%
\instructions{23\rightarrow N_{45},3,2}%
\EasyButWeakLineBreak%
{${72}\above{1pt}{1pt}{r}{2}{14}\above{1pt}{1pt}{b}{2}{18}\above{1pt}{1pt}{l}{2}{8}\above{1pt}{1pt}{r}{2}{126}\above{1pt}{1pt}{b}{2}{2}\above{1pt}{1pt}{l}{2}$}%
\nopagebreak\par%
\nopagebreak\par\leavevmode%
{$\left[\!\llap{\phantom{%
\begingroup \smaller\smaller\smaller\begin{tabular}{@{}c@{}}%
0\\0\\0
\end{tabular}\endgroup%
}}\right.$}%
\begingroup \smaller\smaller\smaller\begin{tabular}{@{}c@{}}%
-20664\\3024\\0
\end{tabular}\endgroup%
\kern3pt%
\begingroup \smaller\smaller\smaller\begin{tabular}{@{}c@{}}%
3024\\-438\\-3
\end{tabular}\endgroup%
\kern3pt%
\begingroup \smaller\smaller\smaller\begin{tabular}{@{}c@{}}%
0\\-3\\2
\end{tabular}\endgroup%
{$\left.\llap{\phantom{%
\begingroup \smaller\smaller\smaller\begin{tabular}{@{}c@{}}%
0\\0\\0
\end{tabular}\endgroup%
}}\!\right]$}%
\EasyButWeakLineBreak%
{$\left[\!\llap{\phantom{%
\begingroup \smaller\smaller\smaller\begin{tabular}{@{}c@{}}%
0\\0\\0
\end{tabular}\endgroup%
}}\right.$}%
\begingroup \smaller\smaller\smaller\begin{tabular}{@{}c@{}}%
7\\48\\72
\end{tabular}\endgroup%
\HardButStrongLineBreak\kern3pt%
\begingroup \smaller\smaller\smaller\begin{tabular}{@{}c@{}}%
1\\7\\14
\end{tabular}\endgroup%
\HardButStrongLineBreak\kern3pt%
\begingroup \smaller\smaller\smaller\begin{tabular}{@{}c@{}}%
-4\\-27\\-36
\end{tabular}\endgroup%
\HardButStrongLineBreak\kern3pt%
\begingroup \smaller\smaller\smaller\begin{tabular}{@{}c@{}}%
-13\\-88\\-128
\end{tabular}\endgroup%
\HardButStrongLineBreak\kern3pt%
\begingroup \smaller\smaller\smaller\begin{tabular}{@{}c@{}}%
-31\\-210\\-315
\end{tabular}\endgroup%
\HardButStrongLineBreak\kern3pt%
\begingroup \smaller\smaller\smaller\begin{tabular}{@{}c@{}}%
0\\0\\-1
\end{tabular}\endgroup%
{$\left.\llap{\phantom{%
\begingroup \smaller\smaller\smaller\begin{tabular}{@{}c@{}}%
0\\0\\0
\end{tabular}\endgroup%
}}\!\right]$}%

\medskip%
%
\leavevmode\llap{}%
$W_{46\phantom{0}}$%
\qquad\llap{24} lattices, $\chi=16$%
\hfill%
$222262$%
\nopagebreak\smallskip\hrule\nopagebreak\medskip%
%
%
\leavevmode%
${L_{46.1}}$%
{} : {$1\above{1pt}{1pt}{-2}{{\rm II}}4\above{1pt}{1pt}{1}{7}{\cdot}1\above{1pt}{1pt}{-}{}3\above{1pt}{1pt}{-}{}9\above{1pt}{1pt}{1}{}{\cdot}1\above{1pt}{1pt}{-2}{}17\above{1pt}{1pt}{1}{}$}\spacer%
\instructions{23\rightarrow N_{46},3,2}%
\EasyButWeakLineBreak%
{${36}\above{1pt}{1pt}{*}{2}{68}\above{1pt}{1pt}{b}{2}{6}\above{1pt}{1pt}{s}{2}{306}\above{1pt}{1pt}{b}{2}{2}\above{1pt}{1pt}{}{6}{6}\above{1pt}{1pt}{b}{2}$}%
\nopagebreak\par%
\nopagebreak\par\leavevmode%
{$\left[\!\llap{\phantom{%
\begingroup \smaller\smaller\smaller\begin{tabular}{@{}c@{}}%
0\\0\\0
\end{tabular}\endgroup%
}}\right.$}%
\begingroup \smaller\smaller\smaller\begin{tabular}{@{}c@{}}%
-85866660\\-89964\\332928
\end{tabular}\endgroup%
\kern3pt%
\begingroup \smaller\smaller\smaller\begin{tabular}{@{}c@{}}%
-89964\\-66\\333
\end{tabular}\endgroup%
\kern3pt%
\begingroup \smaller\smaller\smaller\begin{tabular}{@{}c@{}}%
332928\\333\\-1282
\end{tabular}\endgroup%
{$\left.\llap{\phantom{%
\begingroup \smaller\smaller\smaller\begin{tabular}{@{}c@{}}%
0\\0\\0
\end{tabular}\endgroup%
}}\!\right]$}%
\EasyButWeakLineBreak%
{$\left[\!\llap{\phantom{%
\begingroup \smaller\smaller\smaller\begin{tabular}{@{}c@{}}%
0\\0\\0
\end{tabular}\endgroup%
}}\right.$}%
\begingroup \smaller\smaller\smaller\begin{tabular}{@{}c@{}}%
17\\2892\\5166
\end{tabular}\endgroup%
\HardButStrongLineBreak\kern3pt%
\begingroup \smaller\smaller\smaller\begin{tabular}{@{}c@{}}%
-63\\-10710\\-19142
\end{tabular}\endgroup%
\HardButStrongLineBreak\kern3pt%
\begingroup \smaller\smaller\smaller\begin{tabular}{@{}c@{}}%
-21\\-3571\\-6381
\end{tabular}\endgroup%
\HardButStrongLineBreak\kern3pt%
\begingroup \smaller\smaller\smaller\begin{tabular}{@{}c@{}}%
-215\\-36567\\-65331
\end{tabular}\endgroup%
\HardButStrongLineBreak\kern3pt%
\begingroup \smaller\smaller\smaller\begin{tabular}{@{}c@{}}%
-8\\-1361\\-2431
\end{tabular}\endgroup%
\HardButStrongLineBreak\kern3pt%
\begingroup \smaller\smaller\smaller\begin{tabular}{@{}c@{}}%
7\\1190\\2127
\end{tabular}\endgroup%
{$\left.\llap{\phantom{%
\begingroup \smaller\smaller\smaller\begin{tabular}{@{}c@{}}%
0\\0\\0
\end{tabular}\endgroup%
}}\!\right]$}%

\medskip%
%
\leavevmode\llap{}%
$W_{47\phantom{0}}$%
\qquad\llap{44} lattices, $\chi=48$%
\hfill%
$222222222222\rtimes C_{2}$%
\nopagebreak\smallskip\hrule\nopagebreak\medskip%
%
%
\leavevmode%
${L_{47.1}}$%
{} : {$1\above{1pt}{1pt}{2}{{\rm II}}4\above{1pt}{1pt}{-}{3}{\cdot}1\above{1pt}{1pt}{2}{}3\above{1pt}{1pt}{1}{}{\cdot}1\above{1pt}{1pt}{-2}{}17\above{1pt}{1pt}{1}{}$}\spacer%
\instructions{2\rightarrow N_{47}}%
\EasyButWeakLineBreak%
{${4}\above{1pt}{1pt}{*}{2}{12}\above{1pt}{1pt}{*}{2}{68}\above{1pt}{1pt}{b}{2}{2}\above{1pt}{1pt}{l}{2}{12}\above{1pt}{1pt}{r}{2}{34}\above{1pt}{1pt}{b}{2}$}\relax$\,(\times2)$%
\nopagebreak\par%
\nopagebreak\par\leavevmode%
{$\left[\!\llap{\phantom{%
\begingroup \smaller\smaller\smaller
\begin{tabular}{@{}c@{}}%
0\\0\\0
\end{tabular}\endgroup%
}}\right.$}%
\begingroup \smaller\smaller\smaller\begin{tabular}{@{}c@{}}%
-1142196\\4488\\4080
\end{tabular}\endgroup%
\kern3pt%
\begingroup \smaller\smaller\smaller\begin{tabular}{@{}c@{}}%
4488\\-16\\-17
\end{tabular}\endgroup%
\kern3pt%
\begingroup \smaller\smaller\smaller\begin{tabular}{@{}c@{}}%
4080\\-17\\-14
\end{tabular}\endgroup%
{$\left.\llap{\phantom{%
\begingroup \smaller\smaller\smaller\begin{tabular}{@{}c@{}}%
0\\0\\0
\end{tabular}\endgroup%
}}\!\right]$}%
\hfil\penalty500%
{$\left[\!\llap{\phantom{%
\begingroup \smaller\smaller\smaller\begin{tabular}{@{}c@{}}%
0\\0\\0
\end{tabular}\endgroup%
}}\right.$}%
\begingroup \smaller\smaller\smaller\begin{tabular}{@{}c@{}}%
10097\\1016532\\1703196
\end{tabular}\endgroup%
\kern3pt%
\begingroup \smaller\smaller\smaller\begin{tabular}{@{}c@{}}%
-45\\-4531\\-7590
\end{tabular}\endgroup%
\kern3pt%
\begingroup \smaller\smaller\smaller\begin{tabular}{@{}c@{}}%
-33\\-3322\\-5567
\end{tabular}\endgroup%
{$\left.\llap{\phantom{%
\begingroup \smaller\smaller\smaller\begin{tabular}{@{}c@{}}%
0\\0\\0
\end{tabular}\endgroup%
}}\!\right]$}%
\EasyButWeakLineBreak%
{$\left[\!\llap{\phantom{%
\begingroup \smaller\smaller\smaller\begin{tabular}{@{}c@{}}%
0\\0\\0
\end{tabular}\endgroup%
}}\right.$}%
\begingroup \smaller\smaller\smaller\begin{tabular}{@{}c@{}}%
13\\1308\\2194
\end{tabular}\endgroup%
\HardButStrongLineBreak\kern3pt%
\begingroup \smaller\smaller\smaller\begin{tabular}{@{}c@{}}%
37\\3726\\6240
\end{tabular}\endgroup%
\HardButStrongLineBreak\kern3pt%
\begingroup \smaller\smaller\smaller\begin{tabular}{@{}c@{}}%
191\\19244\\32198
\end{tabular}\endgroup%
\HardButStrongLineBreak\kern3pt%
\begingroup \smaller\smaller\smaller\begin{tabular}{@{}c@{}}%
7\\706\\1179
\end{tabular}\endgroup%
\HardButStrongLineBreak\kern3pt%
\begingroup \smaller\smaller\smaller\begin{tabular}{@{}c@{}}%
7\\708\\1176
\end{tabular}\endgroup%
\HardButStrongLineBreak\kern3pt%
\begingroup \smaller\smaller\smaller\begin{tabular}{@{}c@{}}%
5\\510\\833
\end{tabular}\endgroup%
{$\left.\llap{\phantom{%
\begingroup \smaller\smaller\smaller\begin{tabular}{@{}c@{}}%
0\\0\\0
\end{tabular}\endgroup%
}}\!\right]$}%
%
%
\hbox{}\par\smallskip%
%
%
\leavevmode%
${L_{47.2}}$%
{} : {$1\above{1pt}{1pt}{-2}{6}8\above{1pt}{1pt}{1}{1}{\cdot}1\above{1pt}{1pt}{2}{}3\above{1pt}{1pt}{-}{}{\cdot}1\above{1pt}{1pt}{-2}{}17\above{1pt}{1pt}{1}{}$}\spacer%
\instructions{2\rightarrow N'_{28}}%
\EasyButWeakLineBreak%
{${8}\above{1pt}{1pt}{r}{2}{6}\above{1pt}{1pt}{l}{2}{136}\above{1pt}{1pt}{}{2}{1}\above{1pt}{1pt}{r}{2}{24}\above{1pt}{1pt}{l}{2}{17}\above{1pt}{1pt}{}{2}$}\relax$\,(\times2)$%
\nopagebreak\par%
\nopagebreak\par\leavevmode%
{$\left[\!\llap{\phantom{%
\begingroup \smaller\smaller\smaller
\begin{tabular}{@{}c@{}}%
0\\0\\0
\end{tabular}\endgroup%
}}\right.$}%
\begingroup \smaller\smaller\smaller\begin{tabular}{@{}c@{}}%
-2521848\\8568\\4488
\end{tabular}\endgroup%
\kern3pt%
\begingroup \smaller\smaller\smaller\begin{tabular}{@{}c@{}}%
8568\\-26\\-17
\end{tabular}\endgroup%
\kern3pt%
\begingroup \smaller\smaller\smaller\begin{tabular}{@{}c@{}}%
4488\\-17\\-7
\end{tabular}\endgroup%
{$\left.\llap{\phantom{%
\begingroup \smaller\smaller\smaller\begin{tabular}{@{}c@{}}%
0\\0\\0
\end{tabular}\endgroup%
}}\!\right]$}%
\hfil\penalty500%
{$\left[\!\llap{\phantom{%
\begingroup \smaller\smaller\smaller\begin{tabular}{@{}c@{}}%
0\\0\\0
\end{tabular}\endgroup%
}}\right.$}%
\begingroup \smaller\smaller\smaller\begin{tabular}{@{}c@{}}%
13871\\2115480\\3745440
\end{tabular}\endgroup%
\kern3pt%
\begingroup \smaller\smaller\smaller\begin{tabular}{@{}c@{}}%
-52\\-7931\\-14040
\end{tabular}\endgroup%
\kern3pt%
\begingroup \smaller\smaller\smaller\begin{tabular}{@{}c@{}}%
-22\\-3355\\-5941
\end{tabular}\endgroup%
{$\left.\llap{\phantom{%
\begingroup \smaller\smaller\smaller\begin{tabular}{@{}c@{}}%
0\\0\\0
\end{tabular}\endgroup%
}}\!\right]$}%
\EasyButWeakLineBreak%
{$\left[\!\llap{\phantom{%
\begingroup \smaller\smaller\smaller\begin{tabular}{@{}c@{}}%
0\\0\\0
\end{tabular}\endgroup%
}}\right.$}%
\begingroup \smaller\smaller\smaller\begin{tabular}{@{}c@{}}%
17\\2592\\4592
\end{tabular}\endgroup%
\HardButStrongLineBreak\kern3pt%
\begingroup \smaller\smaller\smaller\begin{tabular}{@{}c@{}}%
25\\3813\\6750
\end{tabular}\endgroup%
\HardButStrongLineBreak\kern3pt%
\begingroup \smaller\smaller\smaller\begin{tabular}{@{}c@{}}%
263\\40120\\70992
\end{tabular}\endgroup%
\HardButStrongLineBreak\kern3pt%
\begingroup \smaller\smaller\smaller\begin{tabular}{@{}c@{}}%
5\\763\\1349
\end{tabular}\endgroup%
\HardButStrongLineBreak\kern3pt%
\begingroup \smaller\smaller\smaller\begin{tabular}{@{}c@{}}%
11\\1680\\2964
\end{tabular}\endgroup%
\HardButStrongLineBreak\kern3pt%
\begingroup \smaller\smaller\smaller\begin{tabular}{@{}c@{}}%
5\\765\\1343
\end{tabular}\endgroup%
{$\left.\llap{\phantom{%
\begingroup \smaller\smaller\smaller\begin{tabular}{@{}c@{}}%
0\\0\\0
\end{tabular}\endgroup%
}}\!\right]$}%
%
%
\hbox{}\par\smallskip%
%
%
\leavevmode%
${L_{47.3}}$%
{} : {$1\above{1pt}{1pt}{2}{6}8\above{1pt}{1pt}{-}{5}{\cdot}1\above{1pt}{1pt}{2}{}3\above{1pt}{1pt}{-}{}{\cdot}1\above{1pt}{1pt}{-2}{}17\above{1pt}{1pt}{1}{}$}\spacer%
\instructions{m}%
\EasyButWeakLineBreak%
{${8}\above{1pt}{1pt}{b}{2}{6}\above{1pt}{1pt}{b}{2}{136}\above{1pt}{1pt}{*}{2}{4}\above{1pt}{1pt}{s}{2}{24}\above{1pt}{1pt}{s}{2}{68}\above{1pt}{1pt}{*}{2}$}\relax$\,(\times2)$%
\nopagebreak\par%
\nopagebreak\par\leavevmode%
{$\left[\!\llap{\phantom{%
\begingroup \smaller\smaller\smaller
\begin{tabular}{@{}c@{}}%
0\\0\\0
\end{tabular}\endgroup%
}}\right.$}%
\begingroup \smaller\smaller\smaller\begin{tabular}{@{}c@{}}%
-6404376\\3014712\\21624
\end{tabular}\endgroup%
\kern3pt%
\begingroup \smaller\smaller\smaller\begin{tabular}{@{}c@{}}%
3014712\\-1419106\\-10179
\end{tabular}\endgroup%
\kern3pt%
\begingroup \smaller\smaller\smaller\begin{tabular}{@{}c@{}}%
21624\\-10179\\-73
\end{tabular}\endgroup%
{$\left.\llap{\phantom{%
\begingroup \smaller\smaller\smaller\begin{tabular}{@{}c@{}}%
0\\0\\0
\end{tabular}\endgroup%
}}\!\right]$}%
\hfil\penalty500%
{$\left[\!\llap{\phantom{%
\begingroup \smaller\smaller\smaller\begin{tabular}{@{}c@{}}%
0\\0\\0
\end{tabular}\endgroup%
}}\right.$}%
\begingroup \smaller\smaller\smaller\begin{tabular}{@{}c@{}}%
281519\\599760\\-244800
\end{tabular}\endgroup%
\kern3pt%
\begingroup \smaller\smaller\smaller\begin{tabular}{@{}c@{}}%
-132526\\-282339\\115240
\end{tabular}\endgroup%
\kern3pt%
\begingroup \smaller\smaller\smaller\begin{tabular}{@{}c@{}}%
-943\\-2009\\819
\end{tabular}\endgroup%
{$\left.\llap{\phantom{%
\begingroup \smaller\smaller\smaller\begin{tabular}{@{}c@{}}%
0\\0\\0
\end{tabular}\endgroup%
}}\!\right]$}%
\EasyButWeakLineBreak%
{$\left[\!\llap{\phantom{%
\begingroup \smaller\smaller\smaller\begin{tabular}{@{}c@{}}%
0\\0\\0
\end{tabular}\endgroup%
}}\right.$}%
\begingroup \smaller\smaller\smaller\begin{tabular}{@{}c@{}}%
201\\428\\-144
\end{tabular}\endgroup%
\HardButStrongLineBreak\kern3pt%
\begingroup \smaller\smaller\smaller\begin{tabular}{@{}c@{}}%
283\\603\\-258
\end{tabular}\endgroup%
\HardButStrongLineBreak\kern3pt%
\begingroup \smaller\smaller\smaller\begin{tabular}{@{}c@{}}%
2903\\6188\\-2992
\end{tabular}\endgroup%
\HardButStrongLineBreak\kern3pt%
\begingroup \smaller\smaller\smaller\begin{tabular}{@{}c@{}}%
105\\224\\-134
\end{tabular}\endgroup%
\HardButStrongLineBreak\kern3pt%
\begingroup \smaller\smaller\smaller\begin{tabular}{@{}c@{}}%
101\\216\\-204
\end{tabular}\endgroup%
\HardButStrongLineBreak\kern3pt%
\begingroup \smaller\smaller\smaller\begin{tabular}{@{}c@{}}%
63\\136\\-306
\end{tabular}\endgroup%
{$\left.\llap{\phantom{%
\begingroup \smaller\smaller\smaller\begin{tabular}{@{}c@{}}%
0\\0\\0
\end{tabular}\endgroup%
}}\!\right]$}%

\medskip%
%
\leavevmode\llap{}%
$W_{48\phantom{0}}$%
\qquad\llap{12} lattices, $\chi=18$%
\hfill%
$422422\rtimes C_{2}$%
\nopagebreak\smallskip\hrule\nopagebreak\medskip%
%
%
\leavevmode%
${L_{48.1}}$%
{} : {$1\above{1pt}{1pt}{-2}{{\rm II}}4\above{1pt}{1pt}{1}{7}{\cdot}1\above{1pt}{1pt}{2}{}3\above{1pt}{1pt}{1}{}{\cdot}1\above{1pt}{1pt}{2}{}17\above{1pt}{1pt}{-}{}$}\spacer%
\instructions{2\rightarrow N_{48}}%
\EasyButWeakLineBreak%
{${2}\above{1pt}{1pt}{*}{4}{4}\above{1pt}{1pt}{b}{2}{102}\above{1pt}{1pt}{s}{2}$}\relax$\,(\times2)$%
\nopagebreak\par%
\nopagebreak\par\leavevmode%
{$\left[\!\llap{\phantom{%
\begingroup \smaller\smaller\smaller\begin{tabular}{@{}c@{}}%
0\\0\\0
\end{tabular}\endgroup%
}}\right.$}%
\begingroup \smaller\smaller\smaller\begin{tabular}{@{}c@{}}%
-13414020\\36108\\52020
\end{tabular}\endgroup%
\kern3pt%
\begingroup \smaller\smaller\smaller\begin{tabular}{@{}c@{}}%
36108\\-94\\-145
\end{tabular}\endgroup%
\kern3pt%
\begingroup \smaller\smaller\smaller\begin{tabular}{@{}c@{}}%
52020\\-145\\-194
\end{tabular}\endgroup%
{$\left.\llap{\phantom{%
\begingroup \smaller\smaller\smaller\begin{tabular}{@{}c@{}}%
0\\0\\0
\end{tabular}\endgroup%
}}\!\right]$}%
\hfil\penalty500%
{$\left[\!\llap{\phantom{%
\begingroup \smaller\smaller\smaller\begin{tabular}{@{}c@{}}%
0\\0\\0
\end{tabular}\endgroup%
}}\right.$}%
\begingroup \smaller\smaller\smaller\begin{tabular}{@{}c@{}}%
83095\\16031136\\10297512
\end{tabular}\endgroup%
\kern3pt%
\begingroup \smaller\smaller\smaller\begin{tabular}{@{}c@{}}%
-247\\-47653\\-30609
\end{tabular}\endgroup%
\kern3pt%
\begingroup \smaller\smaller\smaller\begin{tabular}{@{}c@{}}%
-286\\-55176\\-35443
\end{tabular}\endgroup%
{$\left.\llap{\phantom{%
\begingroup \smaller\smaller\smaller\begin{tabular}{@{}c@{}}%
0\\0\\0
\end{tabular}\endgroup%
}}\!\right]$}%
\EasyButWeakLineBreak%
{$\left[\!\llap{\phantom{%
\begingroup \smaller\smaller\smaller\begin{tabular}{@{}c@{}}%
0\\0\\0
\end{tabular}\endgroup%
}}\right.$}%
\begingroup \smaller\smaller\smaller\begin{tabular}{@{}c@{}}%
6\\1157\\744
\end{tabular}\endgroup%
\HardButStrongLineBreak\kern3pt%
\begingroup \smaller\smaller\smaller\begin{tabular}{@{}c@{}}%
33\\6366\\4090
\end{tabular}\endgroup%
\HardButStrongLineBreak\kern3pt%
\begingroup \smaller\smaller\smaller\begin{tabular}{@{}c@{}}%
235\\45339\\29121
\end{tabular}\endgroup%
{$\left.\llap{\phantom{%
\begingroup \smaller\smaller\smaller\begin{tabular}{@{}c@{}}%
0\\0\\0
\end{tabular}\endgroup%
}}\!\right]$}%

\medskip%
%
\leavevmode\llap{}%
$W_{49\phantom{0}}$%
\qquad\llap{12} lattices, $\chi=30$%
\hfill%
$22242224\rtimes C_{2}$%
\nopagebreak\smallskip\hrule\nopagebreak\medskip%
%
%
\leavevmode%
${L_{49.1}}$%
{} : {$1\above{1pt}{1pt}{-2}{{\rm II}}4\above{1pt}{1pt}{-}{3}{\cdot}1\above{1pt}{1pt}{2}{}5\above{1pt}{1pt}{1}{}{\cdot}1\above{1pt}{1pt}{2}{}11\above{1pt}{1pt}{-}{}$}\spacer%
\instructions{2\rightarrow N_{49}}%
\EasyButWeakLineBreak%
{${4}\above{1pt}{1pt}{*}{2}{20}\above{1pt}{1pt}{b}{2}{22}\above{1pt}{1pt}{s}{2}{2}\above{1pt}{1pt}{*}{4}$}\relax$\,(\times2)$%
\nopagebreak\par%
\nopagebreak\par\leavevmode%
{$\left[\!\llap{\phantom{%
\begingroup \smaller\smaller\smaller\begin{tabular}{@{}c@{}}%
0\\0\\0
\end{tabular}\endgroup%
}}\right.$}%
\begingroup \smaller\smaller\smaller\begin{tabular}{@{}c@{}}%
-1472020\\-81400\\9240
\end{tabular}\endgroup%
\kern3pt%
\begingroup \smaller\smaller\smaller\begin{tabular}{@{}c@{}}%
-81400\\-4494\\511
\end{tabular}\endgroup%
\kern3pt%
\begingroup \smaller\smaller\smaller\begin{tabular}{@{}c@{}}%
9240\\511\\-58
\end{tabular}\endgroup%
{$\left.\llap{\phantom{%
\begingroup \smaller\smaller\smaller\begin{tabular}{@{}c@{}}%
0\\0\\0
\end{tabular}\endgroup%
}}\!\right]$}%
\hfil\penalty500%
{$\left[\!\llap{\phantom{%
\begingroup \smaller\smaller\smaller\begin{tabular}{@{}c@{}}%
0\\0\\0
\end{tabular}\endgroup%
}}\right.$}%
\begingroup \smaller\smaller\smaller\begin{tabular}{@{}c@{}}%
23099\\-23100\\3474240
\end{tabular}\endgroup%
\kern3pt%
\begingroup \smaller\smaller\smaller\begin{tabular}{@{}c@{}}%
1290\\-1291\\194016
\end{tabular}\endgroup%
\kern3pt%
\begingroup \smaller\smaller\smaller\begin{tabular}{@{}c@{}}%
-145\\145\\-21809
\end{tabular}\endgroup%
{$\left.\llap{\phantom{%
\begingroup \smaller\smaller\smaller\begin{tabular}{@{}c@{}}%
0\\0\\0
\end{tabular}\endgroup%
}}\!\right]$}%
\EasyButWeakLineBreak%
{$\left[\!\llap{\phantom{%
\begingroup \smaller\smaller\smaller\begin{tabular}{@{}c@{}}%
0\\0\\0
\end{tabular}\endgroup%
}}\right.$}%
\begingroup \smaller\smaller\smaller\begin{tabular}{@{}c@{}}%
3\\-2\\460
\end{tabular}\endgroup%
\HardButStrongLineBreak\kern3pt%
\begingroup \smaller\smaller\smaller\begin{tabular}{@{}c@{}}%
33\\-30\\4990
\end{tabular}\endgroup%
\HardButStrongLineBreak\kern3pt%
\begingroup \smaller\smaller\smaller\begin{tabular}{@{}c@{}}%
57\\-55\\8591
\end{tabular}\endgroup%
\HardButStrongLineBreak\kern3pt%
\begingroup \smaller\smaller\smaller\begin{tabular}{@{}c@{}}%
27\\-27\\4061
\end{tabular}\endgroup%
{$\left.\llap{\phantom{%
\begingroup \smaller\smaller\smaller\begin{tabular}{@{}c@{}}%
0\\0\\0
\end{tabular}\endgroup%
}}\!\right]$}%

\medskip%
%
\leavevmode\llap{}%
$W_{50\phantom{0}}$%
\qquad\llap{12} lattices, $\chi=24$%
\hfill%
$22222222\rtimes C_{2}$%
\nopagebreak\smallskip\hrule\nopagebreak\medskip%
%
%
\leavevmode%
${L_{50.1}}$%
{} : {$1\above{1pt}{1pt}{-2}{{\rm II}}4\above{1pt}{1pt}{-}{3}{\cdot}1\above{1pt}{1pt}{-2}{}5\above{1pt}{1pt}{-}{}{\cdot}1\above{1pt}{1pt}{-2}{}11\above{1pt}{1pt}{1}{}$}\spacer%
\instructions{2\rightarrow N_{50}}%
\EasyButWeakLineBreak%
{${2}\above{1pt}{1pt}{b}{2}{4}\above{1pt}{1pt}{b}{2}{10}\above{1pt}{1pt}{l}{2}{44}\above{1pt}{1pt}{r}{2}$}\relax$\,(\times2)$%
\nopagebreak\par%
\nopagebreak\par\leavevmode%
{$\left[\!\llap{\phantom{%
\begingroup \smaller\smaller\smaller\begin{tabular}{@{}c@{}}%
0\\0\\0
\end{tabular}\endgroup%
}}\right.$}%
\begingroup \smaller\smaller\smaller\begin{tabular}{@{}c@{}}%
-604340\\110220\\5280
\end{tabular}\endgroup%
\kern3pt%
\begingroup \smaller\smaller\smaller\begin{tabular}{@{}c@{}}%
110220\\-20102\\-963
\end{tabular}\endgroup%
\kern3pt%
\begingroup \smaller\smaller\smaller\begin{tabular}{@{}c@{}}%
5280\\-963\\-46
\end{tabular}\endgroup%
{$\left.\llap{\phantom{%
\begingroup \smaller\smaller\smaller\begin{tabular}{@{}c@{}}%
0\\0\\0
\end{tabular}\endgroup%
}}\!\right]$}%
\hfil\penalty500%
{$\left[\!\llap{\phantom{%
\begingroup \smaller\smaller\smaller\begin{tabular}{@{}c@{}}%
0\\0\\0
\end{tabular}\endgroup%
}}\right.$}%
\begingroup \smaller\smaller\smaller\begin{tabular}{@{}c@{}}%
3079\\16940\\-1540
\end{tabular}\endgroup%
\kern3pt%
\begingroup \smaller\smaller\smaller\begin{tabular}{@{}c@{}}%
-562\\-3092\\281
\end{tabular}\endgroup%
\kern3pt%
\begingroup \smaller\smaller\smaller\begin{tabular}{@{}c@{}}%
-26\\-143\\12
\end{tabular}\endgroup%
{$\left.\llap{\phantom{%
\begingroup \smaller\smaller\smaller\begin{tabular}{@{}c@{}}%
0\\0\\0
\end{tabular}\endgroup%
}}\!\right]$}%
\EasyButWeakLineBreak%
{$\left[\!\llap{\phantom{%
\begingroup \smaller\smaller\smaller\begin{tabular}{@{}c@{}}%
0\\0\\0
\end{tabular}\endgroup%
}}\right.$}%
\begingroup \smaller\smaller\smaller\begin{tabular}{@{}c@{}}%
21\\114\\23
\end{tabular}\endgroup%
\HardButStrongLineBreak\kern3pt%
\begingroup \smaller\smaller\smaller\begin{tabular}{@{}c@{}}%
11\\60\\6
\end{tabular}\endgroup%
\HardButStrongLineBreak\kern3pt%
\begingroup \smaller\smaller\smaller\begin{tabular}{@{}c@{}}%
9\\50\\-15
\end{tabular}\endgroup%
\HardButStrongLineBreak\kern3pt%
\begingroup \smaller\smaller\smaller\begin{tabular}{@{}c@{}}%
-27\\-132\\-352
\end{tabular}\endgroup%
{$\left.\llap{\phantom{%
\begingroup \smaller\smaller\smaller\begin{tabular}{@{}c@{}}%
0\\0\\0
\end{tabular}\endgroup%
}}\!\right]$}%

\medskip%
%
\leavevmode\llap{}%
$W_{51\phantom{0}}$%
\qquad\llap{18} lattices, $\chi=20$%
\hfill%
$226226\rtimes C_{2}$%
\nopagebreak\smallskip\hrule\nopagebreak\medskip%
%
%
\leavevmode%
${L_{51.1}}$%
{} : {$1\above{1pt}{1pt}{-2}{{\rm II}}4\above{1pt}{1pt}{-}{5}{\cdot}1\above{1pt}{1pt}{-}{}3\above{1pt}{1pt}{-}{}9\above{1pt}{1pt}{-}{}{\cdot}1\above{1pt}{1pt}{-2}{}19\above{1pt}{1pt}{-}{}$}\spacer%
\instructions{23\rightarrow N_{51},3,2}%
\EasyButWeakLineBreak%
{${6}\above{1pt}{1pt}{b}{2}{38}\above{1pt}{1pt}{s}{2}{18}\above{1pt}{1pt}{}{6}{6}\above{1pt}{1pt}{b}{2}{342}\above{1pt}{1pt}{s}{2}{2}\above{1pt}{1pt}{}{6}$}%
\nopagebreak\par%
\nopagebreak\par\leavevmode%
{$\left[\!\llap{\phantom{%
\begingroup \smaller\smaller\smaller\begin{tabular}{@{}c@{}}%
0\\0\\0
\end{tabular}\endgroup%
}}\right.$}%
\begingroup \smaller\smaller\smaller\begin{tabular}{@{}c@{}}%
-9620460\\18468\\25992
\end{tabular}\endgroup%
\kern3pt%
\begingroup \smaller\smaller\smaller\begin{tabular}{@{}c@{}}%
18468\\-30\\-51
\end{tabular}\endgroup%
\kern3pt%
\begingroup \smaller\smaller\smaller\begin{tabular}{@{}c@{}}%
25992\\-51\\-70
\end{tabular}\endgroup%
{$\left.\llap{\phantom{%
\begingroup \smaller\smaller\smaller\begin{tabular}{@{}c@{}}%
0\\0\\0
\end{tabular}\endgroup%
}}\!\right]$}%
\EasyButWeakLineBreak%
{$\left[\!\llap{\phantom{%
\begingroup \smaller\smaller\smaller\begin{tabular}{@{}c@{}}%
0\\0\\0
\end{tabular}\endgroup%
}}\right.$}%
\begingroup \smaller\smaller\smaller\begin{tabular}{@{}c@{}}%
-1\\-65\\-324
\end{tabular}\endgroup%
\HardButStrongLineBreak\kern3pt%
\begingroup \smaller\smaller\smaller\begin{tabular}{@{}c@{}}%
-2\\-133\\-646
\end{tabular}\endgroup%
\HardButStrongLineBreak\kern3pt%
\begingroup \smaller\smaller\smaller\begin{tabular}{@{}c@{}}%
2\\129\\648
\end{tabular}\endgroup%
\HardButStrongLineBreak\kern3pt%
\begingroup \smaller\smaller\smaller\begin{tabular}{@{}c@{}}%
16\\1046\\5175
\end{tabular}\endgroup%
\HardButStrongLineBreak\kern3pt%
\begingroup \smaller\smaller\smaller\begin{tabular}{@{}c@{}}%
101\\6612\\32661
\end{tabular}\endgroup%
\HardButStrongLineBreak\kern3pt%
\begingroup \smaller\smaller\smaller\begin{tabular}{@{}c@{}}%
1\\66\\323
\end{tabular}\endgroup%
{$\left.\llap{\phantom{%
\begingroup \smaller\smaller\smaller\begin{tabular}{@{}c@{}}%
0\\0\\0
\end{tabular}\endgroup%
}}\!\right]$}%

\medskip%
%
\leavevmode\llap{}%
$W_{52\phantom{0}}$%
\qquad\llap{44} lattices, $\chi=54$%
\hfill%
$222422222422\rtimes C_{2}$%
\nopagebreak\smallskip\hrule\nopagebreak\medskip%
%
%
\leavevmode%
${L_{52.1}}$%
{} : {$1\above{1pt}{1pt}{2}{{\rm II}}4\above{1pt}{1pt}{1}{1}{\cdot}1\above{1pt}{1pt}{2}{}3\above{1pt}{1pt}{-}{}{\cdot}1\above{1pt}{1pt}{2}{}19\above{1pt}{1pt}{1}{}$}\spacer%
\instructions{2\rightarrow N_{52}}%
\EasyButWeakLineBreak%
{${6}\above{1pt}{1pt}{l}{2}{4}\above{1pt}{1pt}{r}{2}{114}\above{1pt}{1pt}{b}{2}{2}\above{1pt}{1pt}{*}{4}{4}\above{1pt}{1pt}{*}{2}{76}\above{1pt}{1pt}{b}{2}$}\relax$\,(\times2)$%
\nopagebreak\par%
\nopagebreak\par\leavevmode%
{$\left[\!\llap{\phantom{%
\begingroup \smaller\smaller\smaller
\begin{tabular}{@{}c@{}}%
0\\0\\0
\end{tabular}\endgroup%
}}\right.$}%
\begingroup \smaller\smaller\smaller\begin{tabular}{@{}c@{}}%
-833340\\-411084\\2736
\end{tabular}\endgroup%
\kern3pt%
\begingroup \smaller\smaller\smaller\begin{tabular}{@{}c@{}}%
-411084\\-202786\\1349
\end{tabular}\endgroup%
\kern3pt%
\begingroup \smaller\smaller\smaller\begin{tabular}{@{}c@{}}%
2736\\1349\\-8
\end{tabular}\endgroup%
{$\left.\llap{\phantom{%
\begingroup \smaller\smaller\smaller\begin{tabular}{@{}c@{}}%
0\\0\\0
\end{tabular}\endgroup%
}}\!\right]$}%
\hfil\penalty500%
{$\left[\!\llap{\phantom{%
\begingroup \smaller\smaller\smaller\begin{tabular}{@{}c@{}}%
0\\0\\0
\end{tabular}\endgroup%
}}\right.$}%
\begingroup \smaller\smaller\smaller\begin{tabular}{@{}c@{}}%
-684191\\1393080\\904020
\end{tabular}\endgroup%
\kern3pt%
\begingroup \smaller\smaller\smaller\begin{tabular}{@{}c@{}}%
-337109\\686387\\445422
\end{tabular}\endgroup%
\kern3pt%
\begingroup \smaller\smaller\smaller\begin{tabular}{@{}c@{}}%
1662\\-3384\\-2197
\end{tabular}\endgroup%
{$\left.\llap{\phantom{%
\begingroup \smaller\smaller\smaller\begin{tabular}{@{}c@{}}%
0\\0\\0
\end{tabular}\endgroup%
}}\!\right]$}%
\EasyButWeakLineBreak%
{$\left[\!\llap{\phantom{%
\begingroup \smaller\smaller\smaller\begin{tabular}{@{}c@{}}%
0\\0\\0
\end{tabular}\endgroup%
}}\right.$}%
\begingroup \smaller\smaller\smaller\begin{tabular}{@{}c@{}}%
-1967\\4005\\2598
\end{tabular}\endgroup%
\HardButStrongLineBreak\kern3pt%
\begingroup \smaller\smaller\smaller\begin{tabular}{@{}c@{}}%
-609\\1240\\808
\end{tabular}\endgroup%
\HardButStrongLineBreak\kern3pt%
\begingroup \smaller\smaller\smaller\begin{tabular}{@{}c@{}}%
-4619\\9405\\6156
\end{tabular}\endgroup%
\HardButStrongLineBreak\kern3pt%
\begingroup \smaller\smaller\smaller\begin{tabular}{@{}c@{}}%
-138\\281\\186
\end{tabular}\endgroup%
\HardButStrongLineBreak\kern3pt%
\begingroup \smaller\smaller\smaller\begin{tabular}{@{}c@{}}%
-165\\336\\226
\end{tabular}\endgroup%
\HardButStrongLineBreak\kern3pt%
\begingroup \smaller\smaller\smaller\begin{tabular}{@{}c@{}}%
-877\\1786\\1216
\end{tabular}\endgroup%
{$\left.\llap{\phantom{%
\begingroup \smaller\smaller\smaller\begin{tabular}{@{}c@{}}%
0\\0\\0
\end{tabular}\endgroup%
}}\!\right]$}%
%
%
\hbox{}\par\smallskip%
%
%
\leavevmode%
${L_{52.2}}$%
{} : {$1\above{1pt}{1pt}{2}{2}8\above{1pt}{1pt}{1}{7}{\cdot}1\above{1pt}{1pt}{2}{}3\above{1pt}{1pt}{1}{}{\cdot}1\above{1pt}{1pt}{2}{}19\above{1pt}{1pt}{-}{}$}\spacer%
\instructions{2\rightarrow N'_{29}}%
\EasyButWeakLineBreak%
{${12}\above{1pt}{1pt}{s}{2}{8}\above{1pt}{1pt}{l}{2}{57}\above{1pt}{1pt}{}{2}{1}\above{1pt}{1pt}{}{4}{2}\above{1pt}{1pt}{b}{2}{152}\above{1pt}{1pt}{*}{2}$}\relax$\,(\times2)$%
\nopagebreak\par%
\nopagebreak\par\leavevmode%
{$\left[\!\llap{\phantom{%
\begingroup \smaller\smaller\smaller
\begin{tabular}{@{}c@{}}%
0\\0\\0
\end{tabular}\endgroup%
}}\right.$}%
\begingroup \smaller\smaller\smaller\begin{tabular}{@{}c@{}}%
-1477896\\-246696\\4560
\end{tabular}\endgroup%
\kern3pt%
\begingroup \smaller\smaller\smaller\begin{tabular}{@{}c@{}}%
-246696\\-41179\\761
\end{tabular}\endgroup%
\kern3pt%
\begingroup \smaller\smaller\smaller\begin{tabular}{@{}c@{}}%
4560\\761\\-14
\end{tabular}\endgroup%
{$\left.\llap{\phantom{%
\begingroup \smaller\smaller\smaller\begin{tabular}{@{}c@{}}%
0\\0\\0
\end{tabular}\endgroup%
}}\!\right]$}%
\hfil\penalty500%
{$\left[\!\llap{\phantom{%
\begingroup \smaller\smaller\smaller\begin{tabular}{@{}c@{}}%
0\\0\\0
\end{tabular}\endgroup%
}}\right.$}%
\begingroup \smaller\smaller\smaller\begin{tabular}{@{}c@{}}%
6079\\-36480\\-9120
\end{tabular}\endgroup%
\kern3pt%
\begingroup \smaller\smaller\smaller\begin{tabular}{@{}c@{}}%
1018\\-6109\\-1527
\end{tabular}\endgroup%
\kern3pt%
\begingroup \smaller\smaller\smaller\begin{tabular}{@{}c@{}}%
-20\\120\\29
\end{tabular}\endgroup%
{$\left.\llap{\phantom{%
\begingroup \smaller\smaller\smaller\begin{tabular}{@{}c@{}}%
0\\0\\0
\end{tabular}\endgroup%
}}\!\right]$}%
\EasyButWeakLineBreak%
{$\left[\!\llap{\phantom{%
\begingroup \smaller\smaller\smaller\begin{tabular}{@{}c@{}}%
0\\0\\0
\end{tabular}\endgroup%
}}\right.$}%
\begingroup \smaller\smaller\smaller\begin{tabular}{@{}c@{}}%
29\\-174\\-42
\end{tabular}\endgroup%
\HardButStrongLineBreak\kern3pt%
\begingroup \smaller\smaller\smaller\begin{tabular}{@{}c@{}}%
1\\-4\\100
\end{tabular}\endgroup%
\HardButStrongLineBreak\kern3pt%
\begingroup \smaller\smaller\smaller\begin{tabular}{@{}c@{}}%
-26\\171\\798
\end{tabular}\endgroup%
\HardButStrongLineBreak\kern3pt%
\begingroup \smaller\smaller\smaller\begin{tabular}{@{}c@{}}%
-3\\19\\55
\end{tabular}\endgroup%
\HardButStrongLineBreak\kern3pt%
\begingroup \smaller\smaller\smaller\begin{tabular}{@{}c@{}}%
-7\\44\\111
\end{tabular}\endgroup%
\HardButStrongLineBreak\kern3pt%
\begingroup \smaller\smaller\smaller\begin{tabular}{@{}c@{}}%
-97\\608\\1444
\end{tabular}\endgroup%
{$\left.\llap{\phantom{%
\begingroup \smaller\smaller\smaller\begin{tabular}{@{}c@{}}%
0\\0\\0
\end{tabular}\endgroup%
}}\!\right]$}%
%
%
\hbox{}\par\smallskip%
%
%
\leavevmode%
${L_{52.3}}$%
{} : {$1\above{1pt}{1pt}{-2}{2}8\above{1pt}{1pt}{-}{3}{\cdot}1\above{1pt}{1pt}{2}{}3\above{1pt}{1pt}{1}{}{\cdot}1\above{1pt}{1pt}{2}{}19\above{1pt}{1pt}{-}{}$}\spacer%
\instructions{m}%
\EasyButWeakLineBreak%
{${3}\above{1pt}{1pt}{r}{2}{8}\above{1pt}{1pt}{s}{2}{228}\above{1pt}{1pt}{*}{2}{4}\above{1pt}{1pt}{*}{4}{2}\above{1pt}{1pt}{l}{2}{152}\above{1pt}{1pt}{}{2}$}\relax$\,(\times2)$%
\nopagebreak\par%
\nopagebreak\par\leavevmode%
{$\left[\!\llap{\phantom{%
\begingroup \smaller\smaller\smaller
\begin{tabular}{@{}c@{}}%
0\\0\\0
\end{tabular}\endgroup%
}}\right.$}%
\begingroup \smaller\smaller\smaller\begin{tabular}{@{}c@{}}%
-508748712\\-3325152\\93024
\end{tabular}\endgroup%
\kern3pt%
\begingroup \smaller\smaller\smaller\begin{tabular}{@{}c@{}}%
-3325152\\-21733\\608
\end{tabular}\endgroup%
\kern3pt%
\begingroup \smaller\smaller\smaller\begin{tabular}{@{}c@{}}%
93024\\608\\-17
\end{tabular}\endgroup%
{$\left.\llap{\phantom{%
\begingroup \smaller\smaller\smaller\begin{tabular}{@{}c@{}}%
0\\0\\0
\end{tabular}\endgroup%
}}\!\right]$}%
\hfil\penalty500%
{$\left[\!\llap{\phantom{%
\begingroup \smaller\smaller\smaller\begin{tabular}{@{}c@{}}%
0\\0\\0
\end{tabular}\endgroup%
}}\right.$}%
\begingroup \smaller\smaller\smaller\begin{tabular}{@{}c@{}}%
658615\\-100768248\\-282264
\end{tabular}\endgroup%
\kern3pt%
\begingroup \smaller\smaller\smaller\begin{tabular}{@{}c@{}}%
4305\\-658666\\-1845
\end{tabular}\endgroup%
\kern3pt%
\begingroup \smaller\smaller\smaller\begin{tabular}{@{}c@{}}%
-119\\18207\\50
\end{tabular}\endgroup%
{$\left.\llap{\phantom{%
\begingroup \smaller\smaller\smaller\begin{tabular}{@{}c@{}}%
0\\0\\0
\end{tabular}\endgroup%
}}\!\right]$}%
\EasyButWeakLineBreak%
{$\left[\!\llap{\phantom{%
\begingroup \smaller\smaller\smaller\begin{tabular}{@{}c@{}}%
0\\0\\0
\end{tabular}\endgroup%
}}\right.$}%
\begingroup \smaller\smaller\smaller\begin{tabular}{@{}c@{}}%
50\\-7650\\-21
\end{tabular}\endgroup%
\HardButStrongLineBreak\kern3pt%
\begingroup \smaller\smaller\smaller\begin{tabular}{@{}c@{}}%
29\\-4436\\24
\end{tabular}\endgroup%
\HardButStrongLineBreak\kern3pt%
\begingroup \smaller\smaller\smaller\begin{tabular}{@{}c@{}}%
205\\-31350\\456
\end{tabular}\endgroup%
\HardButStrongLineBreak\kern3pt%
\begingroup \smaller\smaller\smaller\begin{tabular}{@{}c@{}}%
5\\-764\\34
\end{tabular}\endgroup%
\HardButStrongLineBreak\kern3pt%
\begingroup \smaller\smaller\smaller\begin{tabular}{@{}c@{}}%
2\\-305\\35
\end{tabular}\endgroup%
\HardButStrongLineBreak\kern3pt%
\begingroup \smaller\smaller\smaller\begin{tabular}{@{}c@{}}%
13\\-1976\\456
\end{tabular}\endgroup%
{$\left.\llap{\phantom{%
\begingroup \smaller\smaller\smaller\begin{tabular}{@{}c@{}}%
0\\0\\0
\end{tabular}\endgroup%
}}\!\right]$}%

\medskip%
%
\leavevmode\llap{}%
$W_{53\phantom{0}}$%
\qquad\llap{12} lattices, $\chi=28$%
\hfill%
$22232223\rtimes C_{2}$%
\nopagebreak\smallskip\hrule\nopagebreak\medskip%
%
%
\leavevmode%
${L_{53.1}}$%
{} : {$1\above{1pt}{1pt}{-2}{{\rm II}}4\above{1pt}{1pt}{-}{5}{\cdot}1\above{1pt}{1pt}{-2}{}5\above{1pt}{1pt}{1}{}{\cdot}1\above{1pt}{1pt}{2}{}13\above{1pt}{1pt}{-}{}$}\spacer%
\instructions{2\rightarrow N_{53}}%
\EasyButWeakLineBreak%
{${2}\above{1pt}{1pt}{b}{2}{26}\above{1pt}{1pt}{l}{2}{20}\above{1pt}{1pt}{r}{2}{2}\above{1pt}{1pt}{-}{3}$}\relax$\,(\times2)$%
\nopagebreak\par%
\nopagebreak\par\leavevmode%
{$\left[\!\llap{\phantom{%
\begingroup \smaller\smaller\smaller\begin{tabular}{@{}c@{}}%
0\\0\\0
\end{tabular}\endgroup%
}}\right.$}%
\begingroup \smaller\smaller\smaller\begin{tabular}{@{}c@{}}%
-527020\\3120\\4680
\end{tabular}\endgroup%
\kern3pt%
\begingroup \smaller\smaller\smaller\begin{tabular}{@{}c@{}}%
3120\\-18\\-29
\end{tabular}\endgroup%
\kern3pt%
\begingroup \smaller\smaller\smaller\begin{tabular}{@{}c@{}}%
4680\\-29\\-38
\end{tabular}\endgroup%
{$\left.\llap{\phantom{%
\begingroup \smaller\smaller\smaller\begin{tabular}{@{}c@{}}%
0\\0\\0
\end{tabular}\endgroup%
}}\!\right]$}%
\hfil\penalty500%
{$\left[\!\llap{\phantom{%
\begingroup \smaller\smaller\smaller\begin{tabular}{@{}c@{}}%
0\\0\\0
\end{tabular}\endgroup%
}}\right.$}%
\begingroup \smaller\smaller\smaller\begin{tabular}{@{}c@{}}%
12349\\1348620\\489060
\end{tabular}\endgroup%
\kern3pt%
\begingroup \smaller\smaller\smaller\begin{tabular}{@{}c@{}}%
-75\\-8191\\-2970
\end{tabular}\endgroup%
\kern3pt%
\begingroup \smaller\smaller\smaller\begin{tabular}{@{}c@{}}%
-105\\-11466\\-4159
\end{tabular}\endgroup%
{$\left.\llap{\phantom{%
\begingroup \smaller\smaller\smaller\begin{tabular}{@{}c@{}}%
0\\0\\0
\end{tabular}\endgroup%
}}\!\right]$}%
\EasyButWeakLineBreak%
{$\left[\!\llap{\phantom{%
\begingroup \smaller\smaller\smaller\begin{tabular}{@{}c@{}}%
0\\0\\0
\end{tabular}\endgroup%
}}\right.$}%
\begingroup \smaller\smaller\smaller\begin{tabular}{@{}c@{}}%
-1\\-109\\-40
\end{tabular}\endgroup%
\HardButStrongLineBreak\kern3pt%
\begingroup \smaller\smaller\smaller\begin{tabular}{@{}c@{}}%
-2\\-221\\-78
\end{tabular}\endgroup%
\HardButStrongLineBreak\kern3pt%
\begingroup \smaller\smaller\smaller\begin{tabular}{@{}c@{}}%
7\\760\\280
\end{tabular}\endgroup%
\HardButStrongLineBreak\kern3pt%
\begingroup \smaller\smaller\smaller\begin{tabular}{@{}c@{}}%
4\\436\\159
\end{tabular}\endgroup%
{$\left.\llap{\phantom{%
\begingroup \smaller\smaller\smaller\begin{tabular}{@{}c@{}}%
0\\0\\0
\end{tabular}\endgroup%
}}\!\right]$}%

\medskip%
%
\leavevmode\llap{}%
$W_{54\phantom{0}}$%
\qquad\llap{32} lattices, $\chi=48$%
\hfill%
$\infty2222\infty2222\rtimes C_{2}$%
\nopagebreak\smallskip\hrule\nopagebreak\medskip%
%
%
\leavevmode%
${L_{54.1}}$%
{} : {$1\above{1pt}{1pt}{-2}{{\rm II}}8\above{1pt}{1pt}{-}{5}{\cdot}1\above{1pt}{1pt}{1}{}3\above{1pt}{1pt}{1}{}9\above{1pt}{1pt}{-}{}{\cdot}1\above{1pt}{1pt}{-2}{}11\above{1pt}{1pt}{-}{}$}\spacer%
\instructions{23\rightarrow N_{54},3,2}%
\EasyButWeakLineBreak%
{${66}\above{1pt}{1pt}{12,1}{\infty b}{264}\above{1pt}{1pt}{b}{2}{18}\above{1pt}{1pt}{s}{2}{22}\above{1pt}{1pt}{b}{2}{72}\above{1pt}{1pt}{b}{2}$}\relax$\,(\times2)$%
\nopagebreak\par%
\nopagebreak\par\leavevmode%
{$\left[\!\llap{\phantom{%
\begingroup \smaller\smaller\smaller
\begin{tabular}{@{}c@{}}%
0\\0\\0
\end{tabular}\endgroup%
}}\right.$}%
\begingroup \smaller\smaller\smaller\begin{tabular}{@{}c@{}}%
-114840\\24552\\9504
\end{tabular}\endgroup%
\kern3pt%
\begingroup \smaller\smaller\smaller\begin{tabular}{@{}c@{}}%
24552\\-5190\\-2181
\end{tabular}\endgroup%
\kern3pt%
\begingroup \smaller\smaller\smaller\begin{tabular}{@{}c@{}}%
9504\\-2181\\-410
\end{tabular}\endgroup%
{$\left.\llap{\phantom{%
\begingroup \smaller\smaller\smaller\begin{tabular}{@{}c@{}}%
0\\0\\0
\end{tabular}\endgroup%
}}\!\right]$}%
\hfil\penalty500%
{$\left[\!\llap{\phantom{%
\begingroup \smaller\smaller\smaller\begin{tabular}{@{}c@{}}%
0\\0\\0
\end{tabular}\endgroup%
}}\right.$}%
\begingroup \smaller\smaller\smaller\begin{tabular}{@{}c@{}}%
375935\\1524864\\601920
\end{tabular}\endgroup%
\kern3pt%
\begingroup \smaller\smaller\smaller\begin{tabular}{@{}c@{}}%
-88466\\-358835\\-141645
\end{tabular}\endgroup%
\kern3pt%
\begingroup \smaller\smaller\smaller\begin{tabular}{@{}c@{}}%
-10680\\-43320\\-17101
\end{tabular}\endgroup%
{$\left.\llap{\phantom{%
\begingroup \smaller\smaller\smaller\begin{tabular}{@{}c@{}}%
0\\0\\0
\end{tabular}\endgroup%
}}\!\right]$}%
\EasyButWeakLineBreak%
{$\left[\!\llap{\phantom{%
\begingroup \smaller\smaller\smaller\begin{tabular}{@{}c@{}}%
0\\0\\0
\end{tabular}\endgroup%
}}\right.$}%
\begingroup \smaller\smaller\smaller\begin{tabular}{@{}c@{}}%
-1874\\-7601\\-3003
\end{tabular}\endgroup%
\HardButStrongLineBreak\kern3pt%
\begingroup \smaller\smaller\smaller\begin{tabular}{@{}c@{}}%
-8407\\-34100\\-13464
\end{tabular}\endgroup%
\HardButStrongLineBreak\kern3pt%
\begingroup \smaller\smaller\smaller\begin{tabular}{@{}c@{}}%
-1619\\-6567\\-2592
\end{tabular}\endgroup%
\HardButStrongLineBreak\kern3pt%
\begingroup \smaller\smaller\smaller\begin{tabular}{@{}c@{}}%
-1581\\-6413\\-2530
\end{tabular}\endgroup%
\HardButStrongLineBreak\kern3pt%
\begingroup \smaller\smaller\smaller\begin{tabular}{@{}c@{}}%
-991\\-4020\\-1584
\end{tabular}\endgroup%
{$\left.\llap{\phantom{%
\begingroup \smaller\smaller\smaller\begin{tabular}{@{}c@{}}%
0\\0\\0
\end{tabular}\endgroup%
}}\!\right]$}%

\medskip%
%
\leavevmode\llap{}%
$W_{55\phantom{0}}$%
\qquad\llap{24} lattices, $\chi=20$%
\hfill%
$622622\rtimes C_{2}$%
\nopagebreak\smallskip\hrule\nopagebreak\medskip%
%
%
\leavevmode%
${L_{55.1}}$%
{} : {$1\above{1pt}{1pt}{-2}{{\rm II}}8\above{1pt}{1pt}{-}{5}{\cdot}1\above{1pt}{1pt}{-}{}3\above{1pt}{1pt}{-}{}9\above{1pt}{1pt}{-}{}{\cdot}1\above{1pt}{1pt}{2}{}11\above{1pt}{1pt}{1}{}$}\spacer%
\instructions{23\rightarrow N_{55},3,2}%
\EasyButWeakLineBreak%
{${18}\above{1pt}{1pt}{}{6}{6}\above{1pt}{1pt}{b}{2}{72}\above{1pt}{1pt}{b}{2}{2}\above{1pt}{1pt}{}{6}{6}\above{1pt}{1pt}{b}{2}{8}\above{1pt}{1pt}{b}{2}$}%
\nopagebreak\par%
\nopagebreak\par\leavevmode%
{$\left[\!\llap{\phantom{%
\begingroup \smaller\smaller\smaller\begin{tabular}{@{}c@{}}%
0\\0\\0
\end{tabular}\endgroup%
}}\right.$}%
\begingroup \smaller\smaller\smaller\begin{tabular}{@{}c@{}}%
-16860888\\101376\\-15048
\end{tabular}\endgroup%
\kern3pt%
\begingroup \smaller\smaller\smaller\begin{tabular}{@{}c@{}}%
101376\\-606\\87
\end{tabular}\endgroup%
\kern3pt%
\begingroup \smaller\smaller\smaller\begin{tabular}{@{}c@{}}%
-15048\\87\\-10
\end{tabular}\endgroup%
{$\left.\llap{\phantom{%
\begingroup \smaller\smaller\smaller\begin{tabular}{@{}c@{}}%
0\\0\\0
\end{tabular}\endgroup%
}}\!\right]$}%
\EasyButWeakLineBreak%
{$\left[\!\llap{\phantom{%
\begingroup \smaller\smaller\smaller\begin{tabular}{@{}c@{}}%
0\\0\\0
\end{tabular}\endgroup%
}}\right.$}%
\begingroup \smaller\smaller\smaller\begin{tabular}{@{}c@{}}%
4\\783\\792
\end{tabular}\endgroup%
\HardButStrongLineBreak\kern3pt%
\begingroup \smaller\smaller\smaller\begin{tabular}{@{}c@{}}%
-1\\-196\\-201
\end{tabular}\endgroup%
\HardButStrongLineBreak\kern3pt%
\begingroup \smaller\smaller\smaller\begin{tabular}{@{}c@{}}%
-19\\-3720\\-3780
\end{tabular}\endgroup%
\HardButStrongLineBreak\kern3pt%
\begingroup \smaller\smaller\smaller\begin{tabular}{@{}c@{}}%
-8\\-1566\\-1589
\end{tabular}\endgroup%
\HardButStrongLineBreak\kern3pt%
\begingroup \smaller\smaller\smaller\begin{tabular}{@{}c@{}}%
-3\\-587\\-594
\end{tabular}\endgroup%
\HardButStrongLineBreak\kern3pt%
\begingroup \smaller\smaller\smaller\begin{tabular}{@{}c@{}}%
1\\196\\200
\end{tabular}\endgroup%
{$\left.\llap{\phantom{%
\begingroup \smaller\smaller\smaller\begin{tabular}{@{}c@{}}%
0\\0\\0
\end{tabular}\endgroup%
}}\!\right]$}%

\medskip%
%
\leavevmode\llap{}%
$W_{56\phantom{0}}$%
\qquad\llap{12} lattices, $\chi=12$%
\hfill%
$222222$%
\nopagebreak\smallskip\hrule\nopagebreak\medskip%
%
%
\leavevmode%
${L_{56.1}}$%
{} : {$1\above{1pt}{1pt}{-2}{{\rm II}}4\above{1pt}{1pt}{1}{1}{\cdot}1\above{1pt}{1pt}{2}{}3\above{1pt}{1pt}{1}{}{\cdot}1\above{1pt}{1pt}{-2}{}23\above{1pt}{1pt}{1}{}$}\spacer%
\instructions{2\rightarrow N_{56}}%
\EasyButWeakLineBreak%
{${12}\above{1pt}{1pt}{*}{2}{92}\above{1pt}{1pt}{b}{2}{2}\above{1pt}{1pt}{b}{2}{138}\above{1pt}{1pt}{l}{2}{4}\above{1pt}{1pt}{r}{2}{46}\above{1pt}{1pt}{b}{2}$}%
\nopagebreak\par%
\nopagebreak\par\leavevmode%
{$\left[\!\llap{\phantom{%
\begingroup \smaller\smaller\smaller\begin{tabular}{@{}c@{}}%
0\\0\\0
\end{tabular}\endgroup%
}}\right.$}%
\begingroup \smaller\smaller\smaller\begin{tabular}{@{}c@{}}%
-5244\\276\\0
\end{tabular}\endgroup%
\kern3pt%
\begingroup \smaller\smaller\smaller\begin{tabular}{@{}c@{}}%
276\\-2\\-5
\end{tabular}\endgroup%
\kern3pt%
\begingroup \smaller\smaller\smaller\begin{tabular}{@{}c@{}}%
0\\-5\\2
\end{tabular}\endgroup%
{$\left.\llap{\phantom{%
\begingroup \smaller\smaller\smaller\begin{tabular}{@{}c@{}}%
0\\0\\0
\end{tabular}\endgroup%
}}\!\right]$}%
\EasyButWeakLineBreak%
{$\left[\!\llap{\phantom{%
\begingroup \smaller\smaller\smaller\begin{tabular}{@{}c@{}}%
0\\0\\0
\end{tabular}\endgroup%
}}\right.$}%
\begingroup \smaller\smaller\smaller\begin{tabular}{@{}c@{}}%
1\\18\\48
\end{tabular}\endgroup%
\HardButStrongLineBreak\kern3pt%
\begingroup \smaller\smaller\smaller\begin{tabular}{@{}c@{}}%
5\\92\\230
\end{tabular}\endgroup%
\HardButStrongLineBreak\kern3pt%
\begingroup \smaller\smaller\smaller\begin{tabular}{@{}c@{}}%
0\\0\\-1
\end{tabular}\endgroup%
\HardButStrongLineBreak\kern3pt%
\begingroup \smaller\smaller\smaller\begin{tabular}{@{}c@{}}%
-7\\-138\\-345
\end{tabular}\endgroup%
\HardButStrongLineBreak\kern3pt%
\begingroup \smaller\smaller\smaller\begin{tabular}{@{}c@{}}%
-1\\-20\\-48
\end{tabular}\endgroup%
\HardButStrongLineBreak\kern3pt%
\begingroup \smaller\smaller\smaller\begin{tabular}{@{}c@{}}%
-1\\-23\\-46
\end{tabular}\endgroup%
{$\left.\llap{\phantom{%
\begingroup \smaller\smaller\smaller\begin{tabular}{@{}c@{}}%
0\\0\\0
\end{tabular}\endgroup%
}}\!\right]$}%

\medskip%
%
\leavevmode\llap{}%
$W_{57\phantom{0}}$%
\qquad\llap{16} lattices, $\chi=24$%
\hfill%
$22222222\rtimes C_{2}$%
\nopagebreak\smallskip\hrule\nopagebreak\medskip%
%
%
\leavevmode%
${L_{57.1}}$%
{} : {$1\above{1pt}{1pt}{-2}{{\rm II}}8\above{1pt}{1pt}{1}{7}{\cdot}1\above{1pt}{1pt}{-2}{}5\above{1pt}{1pt}{-}{}{\cdot}1\above{1pt}{1pt}{2}{}7\above{1pt}{1pt}{1}{}$}\spacer%
\instructions{2\rightarrow N_{57}}%
\EasyButWeakLineBreak%
{${10}\above{1pt}{1pt}{s}{2}{14}\above{1pt}{1pt}{b}{2}{2}\above{1pt}{1pt}{l}{2}{56}\above{1pt}{1pt}{r}{2}$}\relax$\,(\times2)$%
\nopagebreak\par%
\nopagebreak\par\leavevmode%
{$\left[\!\llap{\phantom{%
\begingroup \smaller\smaller\smaller\begin{tabular}{@{}c@{}}%
0\\0\\0
\end{tabular}\endgroup%
}}\right.$}%
\begingroup \smaller\smaller\smaller\begin{tabular}{@{}c@{}}%
3640\\1400\\0
\end{tabular}\endgroup%
\kern3pt%
\begingroup \smaller\smaller\smaller\begin{tabular}{@{}c@{}}%
1400\\538\\1
\end{tabular}\endgroup%
\kern3pt%
\begingroup \smaller\smaller\smaller\begin{tabular}{@{}c@{}}%
0\\1\\-2
\end{tabular}\endgroup%
{$\left.\llap{\phantom{%
\begingroup \smaller\smaller\smaller\begin{tabular}{@{}c@{}}%
0\\0\\0
\end{tabular}\endgroup%
}}\!\right]$}%
\hfil\penalty500%
{$\left[\!\llap{\phantom{%
\begingroup \smaller\smaller\smaller\begin{tabular}{@{}c@{}}%
0\\0\\0
\end{tabular}\endgroup%
}}\right.$}%
\begingroup \smaller\smaller\smaller\begin{tabular}{@{}c@{}}%
-1\\0\\0
\end{tabular}\endgroup%
\kern3pt%
\begingroup \smaller\smaller\smaller\begin{tabular}{@{}c@{}}%
0\\-1\\-1
\end{tabular}\endgroup%
\kern3pt%
\begingroup \smaller\smaller\smaller\begin{tabular}{@{}c@{}}%
0\\0\\1
\end{tabular}\endgroup%
{$\left.\llap{\phantom{%
\begingroup \smaller\smaller\smaller\begin{tabular}{@{}c@{}}%
0\\0\\0
\end{tabular}\endgroup%
}}\!\right]$}%
\EasyButWeakLineBreak%
{$\left[\!\llap{\phantom{%
\begingroup \smaller\smaller\smaller\begin{tabular}{@{}c@{}}%
0\\0\\0
\end{tabular}\endgroup%
}}\right.$}%
\begingroup \smaller\smaller\smaller\begin{tabular}{@{}c@{}}%
-2\\5\\0
\end{tabular}\endgroup%
\HardButStrongLineBreak\kern3pt%
\begingroup \smaller\smaller\smaller\begin{tabular}{@{}c@{}}%
8\\-21\\-14
\end{tabular}\endgroup%
\HardButStrongLineBreak\kern3pt%
\begingroup \smaller\smaller\smaller\begin{tabular}{@{}c@{}}%
5\\-13\\-8
\end{tabular}\endgroup%
\HardButStrongLineBreak\kern3pt%
\begingroup \smaller\smaller\smaller\begin{tabular}{@{}c@{}}%
65\\-168\\-112
\end{tabular}\endgroup%
{$\left.\llap{\phantom{%
\begingroup \smaller\smaller\smaller\begin{tabular}{@{}c@{}}%
0\\0\\0
\end{tabular}\endgroup%
}}\!\right]$}%

\medskip%
%
\leavevmode\llap{}%
$W_{58\phantom{0}}$%
\qquad\llap{12} lattices, $\chi=20$%
\hfill%
$2222232$%
\nopagebreak\smallskip\hrule\nopagebreak\medskip%
%
%
\leavevmode%
${L_{58.1}}$%
{} : {$1\above{1pt}{1pt}{-2}{{\rm II}}4\above{1pt}{1pt}{1}{1}{\cdot}1\above{1pt}{1pt}{-2}{}7\above{1pt}{1pt}{1}{}{\cdot}1\above{1pt}{1pt}{2}{}11\above{1pt}{1pt}{1}{}$}\spacer%
\instructions{2\rightarrow N_{58}}%
\EasyButWeakLineBreak%
{${28}\above{1pt}{1pt}{*}{2}{44}\above{1pt}{1pt}{b}{2}{14}\above{1pt}{1pt}{l}{2}{4}\above{1pt}{1pt}{r}{2}{154}\above{1pt}{1pt}{b}{2}{2}\above{1pt}{1pt}{+}{3}{2}\above{1pt}{1pt}{b}{2}$}%
\nopagebreak\par%
\nopagebreak\par\leavevmode%
{$\left[\!\llap{\phantom{%
\begingroup \smaller\smaller\smaller\begin{tabular}{@{}c@{}}%
0\\0\\0
\end{tabular}\endgroup%
}}\right.$}%
\begingroup \smaller\smaller\smaller\begin{tabular}{@{}c@{}}%
-19740028\\77308\\-33264
\end{tabular}\endgroup%
\kern3pt%
\begingroup \smaller\smaller\smaller\begin{tabular}{@{}c@{}}%
77308\\-302\\127
\end{tabular}\endgroup%
\kern3pt%
\begingroup \smaller\smaller\smaller\begin{tabular}{@{}c@{}}%
-33264\\127\\-42
\end{tabular}\endgroup%
{$\left.\llap{\phantom{%
\begingroup \smaller\smaller\smaller\begin{tabular}{@{}c@{}}%
0\\0\\0
\end{tabular}\endgroup%
}}\!\right]$}%
\EasyButWeakLineBreak%
{$\left[\!\llap{\phantom{%
\begingroup \smaller\smaller\smaller\begin{tabular}{@{}c@{}}%
0\\0\\0
\end{tabular}\endgroup%
}}\right.$}%
\begingroup \smaller\smaller\smaller\begin{tabular}{@{}c@{}}%
11\\3122\\728
\end{tabular}\endgroup%
\HardButStrongLineBreak\kern3pt%
\begingroup \smaller\smaller\smaller\begin{tabular}{@{}c@{}}%
-29\\-8228\\-1914
\end{tabular}\endgroup%
\HardButStrongLineBreak\kern3pt%
\begingroup \smaller\smaller\smaller\begin{tabular}{@{}c@{}}%
-41\\-11634\\-2709
\end{tabular}\endgroup%
\HardButStrongLineBreak\kern3pt%
\begingroup \smaller\smaller\smaller\begin{tabular}{@{}c@{}}%
-137\\-38876\\-9056
\end{tabular}\endgroup%
\HardButStrongLineBreak\kern3pt%
\begingroup \smaller\smaller\smaller\begin{tabular}{@{}c@{}}%
-912\\-258797\\-60291
\end{tabular}\endgroup%
\HardButStrongLineBreak\kern3pt%
\begingroup \smaller\smaller\smaller\begin{tabular}{@{}c@{}}%
-5\\-1419\\-331
\end{tabular}\endgroup%
\HardButStrongLineBreak\kern3pt%
\begingroup \smaller\smaller\smaller\begin{tabular}{@{}c@{}}%
4\\1135\\264
\end{tabular}\endgroup%
{$\left.\llap{\phantom{%
\begingroup \smaller\smaller\smaller\begin{tabular}{@{}c@{}}%
0\\0\\0
\end{tabular}\endgroup%
}}\!\right]$}%

\medskip%
%
\leavevmode\llap{}%
$W_{59\phantom{0}}$%
\qquad\llap{32} lattices, $\chi=28$%
\hfill%
$\infty222262$%
\nopagebreak\smallskip\hrule\nopagebreak\medskip%
%
%
\leavevmode%
${L_{59.1}}$%
{} : {$1\above{1pt}{1pt}{-2}{{\rm II}}8\above{1pt}{1pt}{-}{3}{\cdot}1\above{1pt}{1pt}{1}{}3\above{1pt}{1pt}{-}{}9\above{1pt}{1pt}{-}{}{\cdot}1\above{1pt}{1pt}{2}{}13\above{1pt}{1pt}{-}{}$}\spacer%
\instructions{23\rightarrow N_{59},3,2}%
\EasyButWeakLineBreak%
{${78}\above{1pt}{1pt}{12,7}{\infty b}{312}\above{1pt}{1pt}{b}{2}{18}\above{1pt}{1pt}{l}{2}{24}\above{1pt}{1pt}{r}{2}{234}\above{1pt}{1pt}{b}{2}{6}\above{1pt}{1pt}{}{6}{18}\above{1pt}{1pt}{b}{2}$}%
\nopagebreak\par%
\nopagebreak\par\leavevmode%
{$\left[\!\llap{\phantom{%
\begingroup \smaller\smaller\smaller\begin{tabular}{@{}c@{}}%
0\\0\\0
\end{tabular}\endgroup%
}}\right.$}%
\begingroup \smaller\smaller\smaller\begin{tabular}{@{}c@{}}%
-15718248\\-2834208\\-2395224
\end{tabular}\endgroup%
\kern3pt%
\begingroup \smaller\smaller\smaller\begin{tabular}{@{}c@{}}%
-2834208\\-485598\\-457431
\end{tabular}\endgroup%
\kern3pt%
\begingroup \smaller\smaller\smaller\begin{tabular}{@{}c@{}}%
-2395224\\-457431\\-339362
\end{tabular}\endgroup%
{$\left.\llap{\phantom{%
\begingroup \smaller\smaller\smaller\begin{tabular}{@{}c@{}}%
0\\0\\0
\end{tabular}\endgroup%
}}\!\right]$}%
\EasyButWeakLineBreak%
{$\left[\!\llap{\phantom{%
\begingroup \smaller\smaller\smaller\begin{tabular}{@{}c@{}}%
0\\0\\0
\end{tabular}\endgroup%
}}\right.$}%
\begingroup \smaller\smaller\smaller\begin{tabular}{@{}c@{}}%
-26015\\78325\\78039
\end{tabular}\endgroup%
\HardButStrongLineBreak\kern3pt%
\begingroup \smaller\smaller\smaller\begin{tabular}{@{}c@{}}%
222369\\-669500\\-667056
\end{tabular}\endgroup%
\HardButStrongLineBreak\kern3pt%
\begingroup \smaller\smaller\smaller\begin{tabular}{@{}c@{}}%
91435\\-275289\\-274284
\end{tabular}\endgroup%
\HardButStrongLineBreak\kern3pt%
\begingroup \smaller\smaller\smaller\begin{tabular}{@{}c@{}}%
195783\\-589456\\-587304
\end{tabular}\endgroup%
\HardButStrongLineBreak\kern3pt%
\begingroup \smaller\smaller\smaller\begin{tabular}{@{}c@{}}%
443503\\-1335282\\-1330407
\end{tabular}\endgroup%
\HardButStrongLineBreak\kern3pt%
\begingroup \smaller\smaller\smaller\begin{tabular}{@{}c@{}}%
13008\\-39164\\-39021
\end{tabular}\endgroup%
\HardButStrongLineBreak\kern3pt%
\begingroup \smaller\smaller\smaller\begin{tabular}{@{}c@{}}%
-13102\\39447\\39303
\end{tabular}\endgroup%
{$\left.\llap{\phantom{%
\begingroup \smaller\smaller\smaller\begin{tabular}{@{}c@{}}%
0\\0\\0
\end{tabular}\endgroup%
}}\!\right]$}%

\medskip%
%
\leavevmode\llap{}%
$W_{60\phantom{0}}$%
\qquad\llap{12} lattices, $\chi=48$%
\hfill%
$222222222222\rtimes C_{2}$%
\nopagebreak\smallskip\hrule\nopagebreak\medskip%
%
%
\leavevmode%
${L_{60.1}}$%
{} : {$1\above{1pt}{1pt}{-2}{{\rm II}}4\above{1pt}{1pt}{1}{1}{\cdot}1\above{1pt}{1pt}{2}{}5\above{1pt}{1pt}{-}{}{\cdot}1\above{1pt}{1pt}{-2}{}17\above{1pt}{1pt}{1}{}$}\spacer%
\instructions{2\rightarrow N_{60}}%
\EasyButWeakLineBreak%
{${34}\above{1pt}{1pt}{l}{2}{4}\above{1pt}{1pt}{r}{2}{10}\above{1pt}{1pt}{b}{2}{2}\above{1pt}{1pt}{l}{2}{68}\above{1pt}{1pt}{r}{2}{10}\above{1pt}{1pt}{b}{2}$}\relax$\,(\times2)$%
\nopagebreak\par%
\nopagebreak\par\leavevmode%
{$\left[\!\llap{\phantom{%
\begingroup \smaller\smaller\smaller
\begin{tabular}{@{}c@{}}%
0\\0\\0
\end{tabular}\endgroup%
}}\right.$}%
\begingroup \smaller\smaller\smaller\begin{tabular}{@{}c@{}}%
47940\\19040\\-340
\end{tabular}\endgroup%
\kern3pt%
\begingroup \smaller\smaller\smaller\begin{tabular}{@{}c@{}}%
19040\\7562\\-135
\end{tabular}\endgroup%
\kern3pt%
\begingroup \smaller\smaller\smaller\begin{tabular}{@{}c@{}}%
-340\\-135\\2
\end{tabular}\endgroup%
{$\left.\llap{\phantom{%
\begingroup \smaller\smaller\smaller\begin{tabular}{@{}c@{}}%
0\\0\\0
\end{tabular}%
\endgroup%
}}\!\right]$}%

\medskip%
%
\leavevmode\llap{}%
$W_{61\phantom{0}}$%
\qquad\llap{12} lattices, $\chi=30$%
\hfill%
$22422242\rtimes C_{2}$%
\nopagebreak\smallskip\hrule\nopagebreak\medskip%
%
%
\leavevmode%
${L_{61.1}}$%
{} : {$1\above{1pt}{1pt}{-2}{{\rm II}}4\above{1pt}{1pt}{-}{3}{\cdot}1\above{1pt}{1pt}{2}{}3\above{1pt}{1pt}{1}{}{\cdot}1\above{1pt}{1pt}{2}{}29\above{1pt}{1pt}{-}{}$}\spacer%
\instructions{2\rightarrow N_{61}}%
\EasyButWeakLineBreak%
{${12}\above{1pt}{1pt}{r}{2}{58}\above{1pt}{1pt}{b}{2}{4}\above{1pt}{1pt}{*}{4}{2}\above{1pt}{1pt}{l}{2}$}\relax$\,(\times2)$%
\nopagebreak\par%
\nopagebreak\par\leavevmode%
{$\left[\!\llap{\phantom{%
\begingroup \smaller\smaller\smaller\begin{tabular}{@{}c@{}}%
0\\0\\0
\end{tabular}\endgroup%
}}\right.$}%
\begingroup \smaller\smaller\smaller\begin{tabular}{@{}c@{}}%
-27083796\\37932\\54288
\end{tabular}\endgroup%
\kern3pt%
\begingroup \smaller\smaller\smaller\begin{tabular}{@{}c@{}}%
37932\\-50\\-79
\end{tabular}\endgroup%
\kern3pt%
\begingroup \smaller\smaller\smaller\begin{tabular}{@{}c@{}}%
54288\\-79\\-106
\end{tabular}\endgroup%
{$\left.\llap{\phantom{%
\begingroup \smaller\smaller\smaller\begin{tabular}{@{}c@{}}%
0\\0\\0
\end{tabular}\endgroup%
}}\!\right]$}%
\hfil\penalty500%
{$\left[\!\llap{\phantom{%
\begingroup \smaller\smaller\smaller\begin{tabular}{@{}c@{}}%
0\\0\\0
\end{tabular}\endgroup%
}}\right.$}%
\begingroup \smaller\smaller\smaller\begin{tabular}{@{}c@{}}%
428735\\122111808\\128542848
\end{tabular}\endgroup%
\kern3pt%
\begingroup \smaller\smaller\smaller\begin{tabular}{@{}c@{}}%
-660\\-187981\\-197880
\end{tabular}\endgroup%
\kern3pt%
\begingroup \smaller\smaller\smaller\begin{tabular}{@{}c@{}}%
-803\\-228709\\-240755
\end{tabular}\endgroup%
{$\left.\llap{\phantom{%
\begingroup \smaller\smaller\smaller\begin{tabular}{@{}c@{}}%
0\\0\\0
\end{tabular}\endgroup%
}}\!\right]$}%
\EasyButWeakLineBreak%
{$\left[\!\llap{\phantom{%
\begingroup \smaller\smaller\smaller\begin{tabular}{@{}c@{}}%
0\\0\\0
\end{tabular}\endgroup%
}}\right.$}%
\begingroup \smaller\smaller\smaller\begin{tabular}{@{}c@{}}%
125\\35604\\37476
\end{tabular}\endgroup%
\HardButStrongLineBreak\kern3pt%
\begingroup \smaller\smaller\smaller\begin{tabular}{@{}c@{}}%
644\\183425\\193082
\end{tabular}\endgroup%
\HardButStrongLineBreak\kern3pt%
\begingroup \smaller\smaller\smaller\begin{tabular}{@{}c@{}}%
135\\38450\\40476
\end{tabular}\endgroup%
\HardButStrongLineBreak\kern3pt%
\begingroup \smaller\smaller\smaller\begin{tabular}{@{}c@{}}%
8\\2278\\2399
\end{tabular}\endgroup%
{$\left.\llap{\phantom{%
\begingroup \smaller\smaller\smaller\begin{tabular}{@{}c@{}}%
0\\0\\0
\end{tabular}\endgroup%
}}\!\right]$}%

\medskip%
%
\leavevmode\llap{}%
$W_{62\phantom{0}}$%
\qquad\llap{12} lattices, $\chi=42$%
\hfill%
$4222242222\rtimes C_{2}$%
\nopagebreak\smallskip\hrule\nopagebreak\medskip%
%
%
\leavevmode%
${L_{62.1}}$%
{} : {$1\above{1pt}{1pt}{-2}{{\rm II}}4\above{1pt}{1pt}{1}{7}{\cdot}1\above{1pt}{1pt}{2}{}7\above{1pt}{1pt}{1}{}{\cdot}1\above{1pt}{1pt}{2}{}13\above{1pt}{1pt}{-}{}$}\spacer%
\instructions{2\rightarrow N_{62}}%
\EasyButWeakLineBreak%
{${2}\above{1pt}{1pt}{*}{4}{4}\above{1pt}{1pt}{b}{2}{14}\above{1pt}{1pt}{s}{2}{26}\above{1pt}{1pt}{l}{2}{28}\above{1pt}{1pt}{r}{2}$}\relax$\,(\times2)$%
\nopagebreak\par%
\nopagebreak\par\leavevmode%
{$\left[\!\llap{\phantom{%
\begingroup \smaller\smaller\smaller
\begin{tabular}{@{}c@{}}%
0\\0\\0
\end{tabular}%
\endgroup%
}}\!\right]$}%

\medskip%
%
\leavevmode\llap{}%
$W_{63\phantom{0}}$%
\qquad\llap{12} lattices, $\chi=32$%
\hfill%
$62226222\rtimes C_{2}$%
\nopagebreak\smallskip\hrule\nopagebreak\medskip%
%
%
\leavevmode%
${L_{63.1}}$%
{} : {$1\above{1pt}{1pt}{-2}{{\rm II}}4\above{1pt}{1pt}{1}{1}{\cdot}1\above{1pt}{1pt}{2}{}3\above{1pt}{1pt}{-}{}{\cdot}1\above{1pt}{1pt}{-2}{}31\above{1pt}{1pt}{-}{}$}\spacer%
\instructions{2\rightarrow N_{63}}%
\EasyButWeakLineBreak%
{${2}\above{1pt}{1pt}{}{6}{6}\above{1pt}{1pt}{l}{2}{4}\above{1pt}{1pt}{r}{2}{186}\above{1pt}{1pt}{b}{2}$}\relax$\,(\times2)$%
\nopagebreak\par%
\nopagebreak\par\leavevmode%
{$\left[\!\llap{\phantom{%
\begingroup \smaller\smaller\smaller\begin{tabular}{@{}c@{}}%
0\\0\\0
\end{tabular}\endgroup%
}}\right.$}%
\begingroup \smaller\smaller\smaller\begin{tabular}{@{}c@{}}%
-2765820\\13020\\6324
\end{tabular}\endgroup%
\kern3pt%
\begingroup \smaller\smaller\smaller\begin{tabular}{@{}c@{}}%
13020\\-58\\-31
\end{tabular}\endgroup%
\kern3pt%
\begingroup \smaller\smaller\smaller\begin{tabular}{@{}c@{}}%
6324\\-31\\-14
\end{tabular}\endgroup%
{$\left.\llap{\phantom{%
\begingroup \smaller\smaller\smaller\begin{tabular}{@{}c@{}}%
0\\0\\0
\end{tabular}\endgroup%
}}\!\right]$}%
\hfil\penalty500%
{$\left[\!\llap{\phantom{%
\begingroup \smaller\smaller\smaller\begin{tabular}{@{}c@{}}%
0\\0\\0
\end{tabular}\endgroup%
}}\right.$}%
\begingroup \smaller\smaller\smaller\begin{tabular}{@{}c@{}}%
43151\\3980772\\10658544
\end{tabular}\endgroup%
\kern3pt%
\begingroup \smaller\smaller\smaller\begin{tabular}{@{}c@{}}%
-200\\-18451\\-49400
\end{tabular}\endgroup%
\kern3pt%
\begingroup \smaller\smaller\smaller\begin{tabular}{@{}c@{}}%
-100\\-9225\\-24701
\end{tabular}\endgroup%
{$\left.\llap{\phantom{%
\begingroup \smaller\smaller\smaller\begin{tabular}{@{}c@{}}%
0\\0\\0
\end{tabular}\endgroup%
}}\!\right]$}%
\EasyButWeakLineBreak%
{$\left[\!\llap{\phantom{%
\begingroup \smaller\smaller\smaller\begin{tabular}{@{}c@{}}%
0\\0\\0
\end{tabular}\endgroup%
}}\right.$}%
\begingroup \smaller\smaller\smaller\begin{tabular}{@{}c@{}}%
2\\184\\495
\end{tabular}\endgroup%
\HardButStrongLineBreak\kern3pt%
\begingroup \smaller\smaller\smaller\begin{tabular}{@{}c@{}}%
-1\\-93\\-246
\end{tabular}\endgroup%
\HardButStrongLineBreak\kern3pt%
\begingroup \smaller\smaller\smaller\begin{tabular}{@{}c@{}}%
-1\\-92\\-248
\end{tabular}\endgroup%
\HardButStrongLineBreak\kern3pt%
\begingroup \smaller\smaller\smaller\begin{tabular}{@{}c@{}}%
11\\1023\\2697
\end{tabular}\endgroup%
{$\left.\llap{\phantom{%
\begingroup \smaller\smaller\smaller\begin{tabular}{@{}c@{}}%
0\\0\\0
\end{tabular}\endgroup%
}}\!\right]$}%

\medskip%
%
\leavevmode\llap{}%
$W_{64\phantom{0}}$%
\qquad\llap{12} lattices, $\chi=40$%
\hfill%
$3222232222\rtimes C_{2}$%
\nopagebreak\smallskip\hrule\nopagebreak\medskip%
%
%
\leavevmode%
${L_{64.1}}$%
{} : {$1\above{1pt}{1pt}{-2}{{\rm II}}4\above{1pt}{1pt}{-}{3}{\cdot}1\above{1pt}{1pt}{-2}{}5\above{1pt}{1pt}{-}{}{\cdot}1\above{1pt}{1pt}{-2}{}19\above{1pt}{1pt}{1}{}$}\spacer%
\instructions{2\rightarrow N_{64}}%
\EasyButWeakLineBreak%
{${2}\above{1pt}{1pt}{+}{3}{2}\above{1pt}{1pt}{b}{2}{4}\above{1pt}{1pt}{b}{2}{10}\above{1pt}{1pt}{l}{2}{76}\above{1pt}{1pt}{r}{2}$}\relax$\,(\times2)$%
\nopagebreak\par%
\nopagebreak\par\leavevmode%
{$\left[\!\llap{\phantom{%
\begingroup \smaller\smaller\smaller
\endgroup%
}}\!\right]$}%

\medskip%
%
\leavevmode\llap{}%
$W_{65\phantom{0}}$%
\qquad\llap{16} lattices, $\chi=32$%
\hfill%
$22262226\rtimes C_{2}$%
\nopagebreak\smallskip\hrule\nopagebreak\medskip%
%
%
\leavevmode%
${L_{65.1}}$%
{} : {$1\above{1pt}{1pt}{-2}{{\rm II}}8\above{1pt}{1pt}{1}{7}{\cdot}1\above{1pt}{1pt}{2}{}3\above{1pt}{1pt}{-}{}{\cdot}1\above{1pt}{1pt}{-2}{}17\above{1pt}{1pt}{1}{}$}\spacer%
\instructions{2\rightarrow N_{65}}%
\EasyButWeakLineBreak%
{${6}\above{1pt}{1pt}{s}{2}{34}\above{1pt}{1pt}{b}{2}{24}\above{1pt}{1pt}{b}{2}{2}\above{1pt}{1pt}{}{6}$}\relax$\,(\times2)$%
\nopagebreak\par%
\nopagebreak\par\leavevmode%
{$\left[\!\llap{\phantom{%
\begingroup \smaller\smaller\smaller\begin{tabular}{@{}c@{}}%
0\\0\\0
\end{tabular}\endgroup%
}}\right.$}%
\begingroup \smaller\smaller\smaller\begin{tabular}{@{}c@{}}%
-15547656\\20808\\31824
\end{tabular}\endgroup%
\kern3pt%
\begingroup \smaller\smaller\smaller\begin{tabular}{@{}c@{}}%
20808\\-26\\-45
\end{tabular}\endgroup%
\kern3pt%
\begingroup \smaller\smaller\smaller\begin{tabular}{@{}c@{}}%
31824\\-45\\-62
\end{tabular}\endgroup%
{$\left.\llap{\phantom{%
\begingroup \smaller\smaller\smaller\begin{tabular}{@{}c@{}}%
0\\0\\0
\end{tabular}\endgroup%
}}\!\right]$}%
\hfil\penalty500%
{$\left[\!\llap{\phantom{%
\begingroup \smaller\smaller\smaller\begin{tabular}{@{}c@{}}%
0\\0\\0
\end{tabular}\endgroup%
}}\right.$}%
\begingroup \smaller\smaller\smaller\begin{tabular}{@{}c@{}}%
395759\\136022712\\104361912
\end{tabular}\endgroup%
\kern3pt%
\begingroup \smaller\smaller\smaller\begin{tabular}{@{}c@{}}%
-530\\-182162\\-139761
\end{tabular}\endgroup%
\kern3pt%
\begingroup \smaller\smaller\smaller\begin{tabular}{@{}c@{}}%
-810\\-278397\\-213598
\end{tabular}\endgroup%
{$\left.\llap{\phantom{%
\begingroup \smaller\smaller\smaller\begin{tabular}{@{}c@{}}%
0\\0\\0
\end{tabular}\endgroup%
}}\!\right]$}%
\EasyButWeakLineBreak%
{$\left[\!\llap{\phantom{%
\begingroup \smaller\smaller\smaller\begin{tabular}{@{}c@{}}%
0\\0\\0
\end{tabular}\endgroup%
}}\right.$}%
\begingroup \smaller\smaller\smaller\begin{tabular}{@{}c@{}}%
272\\93486\\71727
\end{tabular}\endgroup%
\HardButStrongLineBreak\kern3pt%
\begingroup \smaller\smaller\smaller\begin{tabular}{@{}c@{}}%
514\\176664\\135541
\end{tabular}\endgroup%
\HardButStrongLineBreak\kern3pt%
\begingroup \smaller\smaller\smaller\begin{tabular}{@{}c@{}}%
115\\39528\\30324
\end{tabular}\endgroup%
\HardButStrongLineBreak\kern3pt%
\begingroup \smaller\smaller\smaller\begin{tabular}{@{}c@{}}%
2\\688\\527
\end{tabular}\endgroup%
{$\left.\llap{\phantom{%
\begingroup \smaller\smaller\smaller\begin{tabular}{@{}c@{}}%
0\\0\\0
\end{tabular}\endgroup%
}}\!\right]$}%

\medskip%
%
\leavevmode\llap{}%
$W_{66\phantom{0}}$%
\qquad\llap{48} lattices, $\chi=36$%
\hfill%
$2222222222\rtimes C_{2}$%
\nopagebreak\smallskip\hrule\nopagebreak\medskip%
%
%
\leavevmode%
${L_{66.1}}$%
{} : {$1\above{1pt}{1pt}{-2}{{\rm II}}4\above{1pt}{1pt}{-}{5}{\cdot}1\above{1pt}{1pt}{1}{}3\above{1pt}{1pt}{-}{}9\above{1pt}{1pt}{-}{}{\cdot}1\above{1pt}{1pt}{2}{}5\above{1pt}{1pt}{1}{}{\cdot}1\above{1pt}{1pt}{2}{}7\above{1pt}{1pt}{-}{}$}\spacer%
\instructions{23\rightarrow N_{66},3,2}%
\EasyButWeakLineBreak%
{${180}\above{1pt}{1pt}{r}{2}{42}\above{1pt}{1pt}{b}{2}{18}\above{1pt}{1pt}{s}{2}{70}\above{1pt}{1pt}{b}{2}{6}\above{1pt}{1pt}{l}{2}$}\relax$\,(\times2)$%
\nopagebreak\par%
\nopagebreak\par\leavevmode%
{$\left[\!\llap{\phantom{%
\begingroup \smaller\smaller\smaller
\endgroup%
}}\!\right]$}%

\medskip%
%
\leavevmode\llap{}%
$W_{67\phantom{0}}$%
\qquad\llap{88} lattices, $\chi=27$%
\hfill%
$22422222$%
\nopagebreak\smallskip\hrule\nopagebreak\medskip%
%
%
\leavevmode%
${L_{67.1}}$%
{} : {$1\above{1pt}{1pt}{2}{{\rm II}}4\above{1pt}{1pt}{1}{1}{\cdot}1\above{1pt}{1pt}{2}{}3\above{1pt}{1pt}{1}{}{\cdot}1\above{1pt}{1pt}{2}{}5\above{1pt}{1pt}{1}{}{\cdot}1\above{1pt}{1pt}{2}{}7\above{1pt}{1pt}{-}{}$}\spacer%
\instructions{2\rightarrow N_{67}}%
\EasyButWeakLineBreak%
{${20}\above{1pt}{1pt}{*}{2}{84}\above{1pt}{1pt}{*}{2}{4}\above{1pt}{1pt}{*}{4}{2}\above{1pt}{1pt}{s}{2}{30}\above{1pt}{1pt}{l}{2}{4}\above{1pt}{1pt}{r}{2}{70}\above{1pt}{1pt}{b}{2}{12}\above{1pt}{1pt}{*}{2}$}%
\nopagebreak\par%
\nopagebreak\par\leavevmode%
{$\left[\!\llap{\phantom{%
\begingroup \smaller\smaller\smaller\begin{tabular}{@{}c@{}}%
0\\0\\0
\end{tabular}\endgroup%
}}\right.$}%
\begingroup \smaller\smaller\smaller\begin{tabular}{@{}c@{}}%
-5778780\\-1930740\\30240
\end{tabular}\endgroup%
\kern3pt%
\begingroup \smaller\smaller\smaller\begin{tabular}{@{}c@{}}%
-1930740\\-645076\\10103
\end{tabular}\endgroup%
\kern3pt%
\begingroup \smaller\smaller\smaller\begin{tabular}{@{}c@{}}%
30240\\10103\\-158
\end{tabular}\endgroup%
{$\left.\llap{\phantom{%
\begingroup \smaller\smaller\smaller\begin{tabular}{@{}c@{}}%
0\\0\\0
\end{tabular}\endgroup%
}}\!\right]$}%
\EasyButWeakLineBreak%
{$\left[\!\llap{\phantom{%
\begingroup \smaller\smaller\smaller\begin{tabular}{@{}c@{}}%
0\\0\\0
\end{tabular}\endgroup%
}}\right.$}%
\begingroup \smaller\smaller\smaller\begin{tabular}{@{}c@{}}%
-159\\490\\900
\end{tabular}\endgroup%
\HardButStrongLineBreak\kern3pt%
\begingroup \smaller\smaller\smaller\begin{tabular}{@{}c@{}}%
41\\-126\\-210
\end{tabular}\endgroup%
\HardButStrongLineBreak\kern3pt%
\begingroup \smaller\smaller\smaller\begin{tabular}{@{}c@{}}%
37\\-114\\-208
\end{tabular}\endgroup%
\HardButStrongLineBreak\kern3pt%
\begingroup \smaller\smaller\smaller\begin{tabular}{@{}c@{}}%
-26\\80\\139
\end{tabular}\endgroup%
\HardButStrongLineBreak\kern3pt%
\begingroup \smaller\smaller\smaller\begin{tabular}{@{}c@{}}%
-302\\930\\1665
\end{tabular}\endgroup%
\HardButStrongLineBreak\kern3pt%
\begingroup \smaller\smaller\smaller\begin{tabular}{@{}c@{}}%
-187\\576\\1040
\end{tabular}\endgroup%
\HardButStrongLineBreak\kern3pt%
\begingroup \smaller\smaller\smaller\begin{tabular}{@{}c@{}}%
-1886\\5810\\10535
\end{tabular}\endgroup%
\HardButStrongLineBreak\kern3pt%
\begingroup \smaller\smaller\smaller\begin{tabular}{@{}c@{}}%
-631\\1944\\3534
\end{tabular}\endgroup%
{$\left.\llap{\phantom{%
\begingroup \smaller\smaller\smaller\begin{tabular}{@{}c@{}}%
0\\0\\0
\end{tabular}\endgroup%
}}\!\right]$}%
%
%
\hbox{}\par\smallskip%
%
%
\leavevmode%
${L_{67.2}}$%
{} : {$1\above{1pt}{1pt}{2}{2}8\above{1pt}{1pt}{1}{7}{\cdot}1\above{1pt}{1pt}{2}{}3\above{1pt}{1pt}{-}{}{\cdot}1\above{1pt}{1pt}{2}{}5\above{1pt}{1pt}{-}{}{\cdot}1\above{1pt}{1pt}{2}{}7\above{1pt}{1pt}{-}{}$}\spacer%
\instructions{2\rightarrow N'_{33}}%
\EasyButWeakLineBreak%
{${10}\above{1pt}{1pt}{s}{2}{42}\above{1pt}{1pt}{b}{2}{2}\above{1pt}{1pt}{}{4}{1}\above{1pt}{1pt}{r}{2}{60}\above{1pt}{1pt}{s}{2}{8}\above{1pt}{1pt}{s}{2}{140}\above{1pt}{1pt}{*}{2}{24}\above{1pt}{1pt}{b}{2}$}%
\nopagebreak\par%
\nopagebreak\par\leavevmode%
{$\left[\!\llap{\phantom{%
\begingroup \smaller\smaller\smaller\begin{tabular}{@{}c@{}}%
0\\0\\0
\end{tabular}\endgroup%
}}\right.$}%
\begingroup \smaller\smaller\smaller\begin{tabular}{@{}c@{}}%
-24797640\\31920\\44520
\end{tabular}\endgroup%
\kern3pt%
\begingroup \smaller\smaller\smaller\begin{tabular}{@{}c@{}}%
31920\\-38\\-59
\end{tabular}\endgroup%
\kern3pt%
\begingroup \smaller\smaller\smaller\begin{tabular}{@{}c@{}}%
44520\\-59\\-79
\end{tabular}\endgroup%
{$\left.\llap{\phantom{%
\begingroup \smaller\smaller\smaller\begin{tabular}{@{}c@{}}%
0\\0\\0
\end{tabular}\endgroup%
}}\!\right]$}%
\EasyButWeakLineBreak%
{$\left[\!\llap{\phantom{%
\begingroup \smaller\smaller\smaller\begin{tabular}{@{}c@{}}%
0\\0\\0
\end{tabular}\endgroup%
}}\right.$}%
\begingroup \smaller\smaller\smaller\begin{tabular}{@{}c@{}}%
4\\875\\1600
\end{tabular}\endgroup%
\HardButStrongLineBreak\kern3pt%
\begingroup \smaller\smaller\smaller\begin{tabular}{@{}c@{}}%
-2\\-441\\-798
\end{tabular}\endgroup%
\HardButStrongLineBreak\kern3pt%
\begingroup \smaller\smaller\smaller\begin{tabular}{@{}c@{}}%
-1\\-219\\-400
\end{tabular}\endgroup%
\HardButStrongLineBreak\kern3pt%
\begingroup \smaller\smaller\smaller\begin{tabular}{@{}c@{}}%
1\\220\\399
\end{tabular}\endgroup%
\HardButStrongLineBreak\kern3pt%
\begingroup \smaller\smaller\smaller\begin{tabular}{@{}c@{}}%
19\\4170\\7590
\end{tabular}\endgroup%
\HardButStrongLineBreak\kern3pt%
\begingroup \smaller\smaller\smaller\begin{tabular}{@{}c@{}}%
11\\2412\\4396
\end{tabular}\endgroup%
\HardButStrongLineBreak\kern3pt%
\begingroup \smaller\smaller\smaller\begin{tabular}{@{}c@{}}%
107\\23450\\42770
\end{tabular}\endgroup%
\HardButStrongLineBreak\kern3pt%
\begingroup \smaller\smaller\smaller\begin{tabular}{@{}c@{}}%
35\\7668\\13992
\end{tabular}\endgroup%
{$\left.\llap{\phantom{%
\begingroup \smaller\smaller\smaller\begin{tabular}{@{}c@{}}%
0\\0\\0
\end{tabular}\endgroup%
}}\!\right]$}%
%
%
\hbox{}\par\smallskip%
%
%
\leavevmode%
${L_{67.3}}$%
{} : {$1\above{1pt}{1pt}{-2}{2}8\above{1pt}{1pt}{-}{3}{\cdot}1\above{1pt}{1pt}{2}{}3\above{1pt}{1pt}{-}{}{\cdot}1\above{1pt}{1pt}{2}{}5\above{1pt}{1pt}{-}{}{\cdot}1\above{1pt}{1pt}{2}{}7\above{1pt}{1pt}{-}{}$}\spacer%
\instructions{m}%
\EasyButWeakLineBreak%
{${10}\above{1pt}{1pt}{b}{2}{42}\above{1pt}{1pt}{s}{2}{2}\above{1pt}{1pt}{*}{4}{4}\above{1pt}{1pt}{l}{2}{15}\above{1pt}{1pt}{r}{2}{8}\above{1pt}{1pt}{l}{2}{35}\above{1pt}{1pt}{}{2}{24}\above{1pt}{1pt}{r}{2}$}%
\nopagebreak\par%
\nopagebreak\par\leavevmode%
{$\left[\!\llap{\phantom{%
\begingroup \smaller\smaller\smaller\begin{tabular}{@{}c@{}}%
0\\0\\0
\end{tabular}\endgroup%
}}\right.$}%
\begingroup \smaller\smaller\smaller\begin{tabular}{@{}c@{}}%
-32246760\\36960\\22680
\end{tabular}\endgroup%
\kern3pt%
\begingroup \smaller\smaller\smaller\begin{tabular}{@{}c@{}}%
36960\\-41\\-28
\end{tabular}\endgroup%
\kern3pt%
\begingroup \smaller\smaller\smaller\begin{tabular}{@{}c@{}}%
22680\\-28\\-13
\end{tabular}\endgroup%
{$\left.\llap{\phantom{%
\begingroup \smaller\smaller\smaller\begin{tabular}{@{}c@{}}%
0\\0\\0
\end{tabular}\endgroup%
}}\!\right]$}%
\EasyButWeakLineBreak%
{$\left[\!\llap{\phantom{%
\begingroup \smaller\smaller\smaller\begin{tabular}{@{}c@{}}%
0\\0\\0
\end{tabular}\endgroup%
}}\right.$}%
\begingroup \smaller\smaller\smaller\begin{tabular}{@{}c@{}}%
8\\4925\\3345
\end{tabular}\endgroup%
\HardButStrongLineBreak\kern3pt%
\begingroup \smaller\smaller\smaller\begin{tabular}{@{}c@{}}%
13\\8001\\5439
\end{tabular}\endgroup%
\HardButStrongLineBreak\kern3pt%
\begingroup \smaller\smaller\smaller\begin{tabular}{@{}c@{}}%
1\\615\\419
\end{tabular}\endgroup%
\HardButStrongLineBreak\kern3pt%
\begingroup \smaller\smaller\smaller\begin{tabular}{@{}c@{}}%
-1\\-616\\-418
\end{tabular}\endgroup%
\HardButStrongLineBreak\kern3pt%
\begingroup \smaller\smaller\smaller\begin{tabular}{@{}c@{}}%
-1\\-615\\-420
\end{tabular}\endgroup%
\HardButStrongLineBreak\kern3pt%
\begingroup \smaller\smaller\smaller\begin{tabular}{@{}c@{}}%
3\\1848\\1252
\end{tabular}\endgroup%
\HardButStrongLineBreak\kern3pt%
\begingroup \smaller\smaller\smaller\begin{tabular}{@{}c@{}}%
32\\19705\\13370
\end{tabular}\endgroup%
\HardButStrongLineBreak\kern3pt%
\begingroup \smaller\smaller\smaller\begin{tabular}{@{}c@{}}%
29\\17856\\12120
\end{tabular}\endgroup%
{$\left.\llap{\phantom{%
\begingroup \smaller\smaller\smaller\begin{tabular}{@{}c@{}}%
0\\0\\0
\end{tabular}\endgroup%
}}\!\right]$}%

\medskip%
%
\leavevmode\llap{}%
$W_{68\phantom{0}}$%
\qquad\llap{176} lattices, $\chi=72$%
\hfill%
$2222222222222222\rtimes C_{2}$%
\nopagebreak\smallskip\hrule\nopagebreak\medskip%
%
%
\leavevmode%
${L_{68.1}}$%
{} : {$1\above{1pt}{1pt}{2}{{\rm II}}4\above{1pt}{1pt}{1}{1}{\cdot}1\above{1pt}{1pt}{-}{}3\above{1pt}{1pt}{-}{}9\above{1pt}{1pt}{1}{}{\cdot}1\above{1pt}{1pt}{-2}{}5\above{1pt}{1pt}{-}{}{\cdot}1\above{1pt}{1pt}{2}{}7\above{1pt}{1pt}{-}{}$}\spacer%
\instructions{23\rightarrow N_{68},3,2}%
\EasyButWeakLineBreak%
{${36}\above{1pt}{1pt}{*}{2}{60}\above{1pt}{1pt}{b}{2}{2}\above{1pt}{1pt}{l}{2}{36}\above{1pt}{1pt}{r}{2}{42}\above{1pt}{1pt}{b}{2}{90}\above{1pt}{1pt}{s}{2}{6}\above{1pt}{1pt}{b}{2}{140}\above{1pt}{1pt}{*}{2}$}\relax$\,(\times2)$%
\nopagebreak\par%
\nopagebreak\par\leavevmode%
{$\left[\!\llap{\phantom{%
\begingroup \smaller\smaller\smaller
\endgroup%
}}\!\right]$}%
%
%
\hbox{}\par\smallskip%
%
%
\leavevmode%
${L_{68.2}}$%
{} : {$1\above{1pt}{1pt}{-2}{2}8\above{1pt}{1pt}{-}{3}{\cdot}1\above{1pt}{1pt}{1}{}3\above{1pt}{1pt}{1}{}9\above{1pt}{1pt}{-}{}{\cdot}1\above{1pt}{1pt}{-2}{}5\above{1pt}{1pt}{1}{}{\cdot}1\above{1pt}{1pt}{2}{}7\above{1pt}{1pt}{-}{}$}\spacer%
\instructions{3m,3,2}%
\EasyButWeakLineBreak%
{${18}\above{1pt}{1pt}{b}{2}{120}\above{1pt}{1pt}{*}{2}{4}\above{1pt}{1pt}{s}{2}{72}\above{1pt}{1pt}{s}{2}{84}\above{1pt}{1pt}{*}{2}{180}\above{1pt}{1pt}{l}{2}{3}\above{1pt}{1pt}{}{2}{280}\above{1pt}{1pt}{r}{2}$}\relax$\,(\times2)$%
\nopagebreak\par%
\nopagebreak\par\leavevmode%
{$\left[\!\llap{\phantom{%
\begingroup \smaller\smaller\smaller
\endgroup%
}}\!\right]$}%
%
%
\hbox{}\par\smallskip%
%
%
\leavevmode%
${L_{68.3}}$%
{} : {$1\above{1pt}{1pt}{2}{2}8\above{1pt}{1pt}{1}{7}{\cdot}1\above{1pt}{1pt}{1}{}3\above{1pt}{1pt}{1}{}9\above{1pt}{1pt}{-}{}{\cdot}1\above{1pt}{1pt}{-2}{}5\above{1pt}{1pt}{1}{}{\cdot}1\above{1pt}{1pt}{2}{}7\above{1pt}{1pt}{-}{}$}\spacer%
\instructions{32\rightarrow N'_{32},3,m}%
\EasyButWeakLineBreak%
{${18}\above{1pt}{1pt}{l}{2}{120}\above{1pt}{1pt}{}{2}{1}\above{1pt}{1pt}{r}{2}{72}\above{1pt}{1pt}{l}{2}{21}\above{1pt}{1pt}{}{2}{45}\above{1pt}{1pt}{r}{2}{12}\above{1pt}{1pt}{*}{2}{280}\above{1pt}{1pt}{b}{2}$}\relax$\,(\times2)$%
\nopagebreak\par%
\nopagebreak\par\leavevmode%
{$\left[\!\llap{\phantom{%
\begingroup \smaller\smaller\smaller
\endgroup%
}}\!\right]$}%

\medskip%
%
\leavevmode\llap{}%
$W_{69\phantom{0}}$%
\qquad\llap{24} lattices, $\chi=6$%
\hfill%
$22222$%
\nopagebreak\smallskip\hrule\nopagebreak\medskip%
%
%
\leavevmode%
${L_{69.1}}$%
{} : {$1\above{1pt}{1pt}{-2}{{\rm II}}4\above{1pt}{1pt}{-}{5}{\cdot}1\above{1pt}{1pt}{2}{}3\above{1pt}{1pt}{1}{}{\cdot}1\above{1pt}{1pt}{-2}{}5\above{1pt}{1pt}{-}{}{\cdot}1\above{1pt}{1pt}{2}{}7\above{1pt}{1pt}{-}{}$}\spacer%
\instructions{2\rightarrow N_{69}}%
\EasyButWeakLineBreak%
{${12}\above{1pt}{1pt}{b}{2}{10}\above{1pt}{1pt}{l}{2}{84}\above{1pt}{1pt}{r}{2}{2}\above{1pt}{1pt}{b}{2}{140}\above{1pt}{1pt}{*}{2}$}%
\nopagebreak\par%
\nopagebreak\par\leavevmode%
{$\left[\!\llap{\phantom{%
\begingroup \smaller\smaller\smaller\begin{tabular}{@{}c@{}}%
0\\0\\0
\end{tabular}\endgroup%
}}\right.$}%
\begingroup \smaller\smaller\smaller\begin{tabular}{@{}c@{}}%
-3898860\\1115100\\10080
\end{tabular}\endgroup%
\kern3pt%
\begingroup \smaller\smaller\smaller\begin{tabular}{@{}c@{}}%
1115100\\-318926\\-2883
\end{tabular}\endgroup%
\kern3pt%
\begingroup \smaller\smaller\smaller\begin{tabular}{@{}c@{}}%
10080\\-2883\\-26
\end{tabular}\endgroup%
{$\left.\llap{\phantom{%
\begingroup \smaller\smaller\smaller\begin{tabular}{@{}c@{}}%
0\\0\\0
\end{tabular}\endgroup%
}}\!\right]$}%
\EasyButWeakLineBreak%
{$\left[\!\llap{\phantom{%
\begingroup \smaller\smaller\smaller\begin{tabular}{@{}c@{}}%
0\\0\\0
\end{tabular}\endgroup%
}}\right.$}%
\begingroup \smaller\smaller\smaller\begin{tabular}{@{}c@{}}%
-19\\-66\\-48
\end{tabular}\endgroup%
\HardButStrongLineBreak\kern3pt%
\begingroup \smaller\smaller\smaller\begin{tabular}{@{}c@{}}%
13\\45\\50
\end{tabular}\endgroup%
\HardButStrongLineBreak\kern3pt%
\begingroup \smaller\smaller\smaller\begin{tabular}{@{}c@{}}%
121\\420\\336
\end{tabular}\endgroup%
\HardButStrongLineBreak\kern3pt%
\begingroup \smaller\smaller\smaller\begin{tabular}{@{}c@{}}%
2\\7\\-1
\end{tabular}\endgroup%
\HardButStrongLineBreak\kern3pt%
\begingroup \smaller\smaller\smaller\begin{tabular}{@{}c@{}}%
-101\\-350\\-350
\end{tabular}\endgroup%
{$\left.\llap{\phantom{%
\begingroup \smaller\smaller\smaller\begin{tabular}{@{}c@{}}%
0\\0\\0
\end{tabular}\endgroup%
}}\!\right]$}%

\medskip%
%
\leavevmode\llap{}%
$W_{70\phantom{0}}$%
\qquad\llap{24} lattices, $\chi=12$%
\hfill%
$222222$%
\nopagebreak\smallskip\hrule\nopagebreak\medskip%
%
%
\leavevmode%
${L_{70.1}}$%
{} : {$1\above{1pt}{1pt}{-2}{{\rm II}}4\above{1pt}{1pt}{-}{5}{\cdot}1\above{1pt}{1pt}{2}{}3\above{1pt}{1pt}{1}{}{\cdot}1\above{1pt}{1pt}{2}{}5\above{1pt}{1pt}{1}{}{\cdot}1\above{1pt}{1pt}{-2}{}7\above{1pt}{1pt}{1}{}$}\spacer%
\instructions{2\rightarrow N_{70}}%
\EasyButWeakLineBreak%
{${12}\above{1pt}{1pt}{*}{2}{28}\above{1pt}{1pt}{b}{2}{30}\above{1pt}{1pt}{b}{2}{14}\above{1pt}{1pt}{l}{2}{20}\above{1pt}{1pt}{r}{2}{2}\above{1pt}{1pt}{b}{2}$}%
\nopagebreak\par%
\nopagebreak\par\leavevmode%
{$\left[\!\llap{\phantom{%
\begingroup \smaller\smaller\smaller\begin{tabular}{@{}c@{}}%
0\\0\\0
\end{tabular}\endgroup%
}}\right.$}%
\begingroup \smaller\smaller\smaller\begin{tabular}{@{}c@{}}%
-713580\\-355320\\420
\end{tabular}\endgroup%
\kern3pt%
\begingroup \smaller\smaller\smaller\begin{tabular}{@{}c@{}}%
-355320\\-176926\\207
\end{tabular}\endgroup%
\kern3pt%
\begingroup \smaller\smaller\smaller\begin{tabular}{@{}c@{}}%
420\\207\\2
\end{tabular}\endgroup%
{$\left.\llap{\phantom{%
\begingroup \smaller\smaller\smaller\begin{tabular}{@{}c@{}}%
0\\0\\0
\end{tabular}\endgroup%
}}\!\right]$}%
\EasyButWeakLineBreak%
{$\left[\!\llap{\phantom{%
\begingroup \smaller\smaller\smaller\begin{tabular}{@{}c@{}}%
0\\0\\0
\end{tabular}\endgroup%
}}\right.$}%
\begingroup \smaller\smaller\smaller\begin{tabular}{@{}c@{}}%
191\\-384\\-366
\end{tabular}\endgroup%
\HardButStrongLineBreak\kern3pt%
\begingroup \smaller\smaller\smaller\begin{tabular}{@{}c@{}}%
383\\-770\\-728
\end{tabular}\endgroup%
\HardButStrongLineBreak\kern3pt%
\begingroup \smaller\smaller\smaller\begin{tabular}{@{}c@{}}%
97\\-195\\-180
\end{tabular}\endgroup%
\HardButStrongLineBreak\kern3pt%
\begingroup \smaller\smaller\smaller\begin{tabular}{@{}c@{}}%
-94\\189\\182
\end{tabular}\endgroup%
\HardButStrongLineBreak\kern3pt%
\begingroup \smaller\smaller\smaller\begin{tabular}{@{}c@{}}%
-189\\380\\360
\end{tabular}\endgroup%
\HardButStrongLineBreak\kern3pt%
\begingroup \smaller\smaller\smaller\begin{tabular}{@{}c@{}}%
0\\0\\-1
\end{tabular}\endgroup%
{$\left.\llap{\phantom{%
\begingroup \smaller\smaller\smaller\begin{tabular}{@{}c@{}}%
0\\0\\0
\end{tabular}\endgroup%
}}\!\right]$}%

\medskip%
%
\leavevmode\llap{}%
$W_{71\phantom{0}}$%
\qquad\llap{48} lattices, $\chi=16$%
\hfill%
$222262$%
\nopagebreak\smallskip\hrule\nopagebreak\medskip%
%
%
\leavevmode%
${L_{71.1}}$%
{} : {$1\above{1pt}{1pt}{-2}{{\rm II}}4\above{1pt}{1pt}{-}{5}{\cdot}1\above{1pt}{1pt}{1}{}3\above{1pt}{1pt}{-}{}9\above{1pt}{1pt}{-}{}{\cdot}1\above{1pt}{1pt}{-2}{}5\above{1pt}{1pt}{-}{}{\cdot}1\above{1pt}{1pt}{-2}{}7\above{1pt}{1pt}{1}{}$}\spacer%
\instructions{23\rightarrow N_{71},3,2}%
\EasyButWeakLineBreak%
{${28}\above{1pt}{1pt}{*}{2}{60}\above{1pt}{1pt}{b}{2}{126}\above{1pt}{1pt}{s}{2}{10}\above{1pt}{1pt}{b}{2}{18}\above{1pt}{1pt}{}{6}{6}\above{1pt}{1pt}{b}{2}$}%
\nopagebreak\par%
\nopagebreak\par\leavevmode%
{$\left[\!\llap{\phantom{%
\begingroup \smaller\smaller\smaller\begin{tabular}{@{}c@{}}%
0\\0\\0
\end{tabular}\endgroup%
}}\right.$}%
\begingroup \smaller\smaller\smaller\begin{tabular}{@{}c@{}}%
-7158060\\-2882880\\10080
\end{tabular}\endgroup%
\kern3pt%
\begingroup \smaller\smaller\smaller\begin{tabular}{@{}c@{}}%
-2882880\\-1161066\\4059
\end{tabular}\endgroup%
\kern3pt%
\begingroup \smaller\smaller\smaller\begin{tabular}{@{}c@{}}%
10080\\4059\\-14
\end{tabular}\endgroup%
{$\left.\llap{\phantom{%
\begingroup \smaller\smaller\smaller\begin{tabular}{@{}c@{}}%
0\\0\\0
\end{tabular}\endgroup%
}}\!\right]$}%
\EasyButWeakLineBreak%
{$\left[\!\llap{\phantom{%
\begingroup \smaller\smaller\smaller\begin{tabular}{@{}c@{}}%
0\\0\\0
\end{tabular}\endgroup%
}}\right.$}%
\begingroup \smaller\smaller\smaller\begin{tabular}{@{}c@{}}%
-117\\294\\994
\end{tabular}\endgroup%
\HardButStrongLineBreak\kern3pt%
\begingroup \smaller\smaller\smaller\begin{tabular}{@{}c@{}}%
-187\\470\\1620
\end{tabular}\endgroup%
\HardButStrongLineBreak\kern3pt%
\begingroup \smaller\smaller\smaller\begin{tabular}{@{}c@{}}%
-142\\357\\1260
\end{tabular}\endgroup%
\HardButStrongLineBreak\kern3pt%
\begingroup \smaller\smaller\smaller\begin{tabular}{@{}c@{}}%
2\\-5\\-10
\end{tabular}\endgroup%
\HardButStrongLineBreak\kern3pt%
\begingroup \smaller\smaller\smaller\begin{tabular}{@{}c@{}}%
37\\-93\\-324
\end{tabular}\endgroup%
\HardButStrongLineBreak\kern3pt%
\begingroup \smaller\smaller\smaller\begin{tabular}{@{}c@{}}%
-2\\5\\9
\end{tabular}\endgroup%
{$\left.\llap{\phantom{%
\begingroup \smaller\smaller\smaller\begin{tabular}{@{}c@{}}%
0\\0\\0
\end{tabular}\endgroup%
}}\!\right]$}%

\medskip%
%
\leavevmode\llap{}%
$W_{72\phantom{0}}$%
\qquad\llap{88} lattices, $\chi=24$%
\hfill%
$22222222$%
\nopagebreak\smallskip\hrule\nopagebreak\medskip%
%
%
\leavevmode%
${L_{72.1}}$%
{} : {$1\above{1pt}{1pt}{2}{{\rm II}}4\above{1pt}{1pt}{1}{1}{\cdot}1\above{1pt}{1pt}{2}{}3\above{1pt}{1pt}{1}{}{\cdot}1\above{1pt}{1pt}{-2}{}5\above{1pt}{1pt}{-}{}{\cdot}1\above{1pt}{1pt}{-2}{}7\above{1pt}{1pt}{1}{}$}\spacer%
\instructions{2\rightarrow N_{72}}%
\EasyButWeakLineBreak%
{${12}\above{1pt}{1pt}{b}{2}{14}\above{1pt}{1pt}{s}{2}{2}\above{1pt}{1pt}{b}{2}{210}\above{1pt}{1pt}{l}{2}{4}\above{1pt}{1pt}{r}{2}{10}\above{1pt}{1pt}{b}{2}{28}\above{1pt}{1pt}{*}{2}{4}\above{1pt}{1pt}{*}{2}$}%
\nopagebreak\par%
\nopagebreak\par\leavevmode%
{$\left[\!\llap{\phantom{%
\begingroup \smaller\smaller\smaller\begin{tabular}{@{}c@{}}%
0\\0\\0
\end{tabular}\endgroup%
}}\right.$}%
\begingroup \smaller\smaller\smaller\begin{tabular}{@{}c@{}}%
-133966140\\1288140\\-32760
\end{tabular}\endgroup%
\kern3pt%
\begingroup \smaller\smaller\smaller\begin{tabular}{@{}c@{}}%
1288140\\-12386\\315
\end{tabular}\endgroup%
\kern3pt%
\begingroup \smaller\smaller\smaller\begin{tabular}{@{}c@{}}%
-32760\\315\\-8
\end{tabular}\endgroup%
{$\left.\llap{\phantom{%
\begingroup \smaller\smaller\smaller\begin{tabular}{@{}c@{}}%
0\\0\\0
\end{tabular}\endgroup%
}}\!\right]$}%
\EasyButWeakLineBreak%
{$\left[\!\llap{\phantom{%
\begingroup \smaller\smaller\smaller\begin{tabular}{@{}c@{}}%
0\\0\\0
\end{tabular}\endgroup%
}}\right.$}%
\begingroup \smaller\smaller\smaller\begin{tabular}{@{}c@{}}%
11\\1146\\72
\end{tabular}\endgroup%
\HardButStrongLineBreak\kern3pt%
\begingroup \smaller\smaller\smaller\begin{tabular}{@{}c@{}}%
10\\1043\\112
\end{tabular}\endgroup%
\HardButStrongLineBreak\kern3pt%
\begingroup \smaller\smaller\smaller\begin{tabular}{@{}c@{}}%
2\\209\\38
\end{tabular}\endgroup%
\HardButStrongLineBreak\kern3pt%
\begingroup \smaller\smaller\smaller\begin{tabular}{@{}c@{}}%
11\\1155\\420
\end{tabular}\endgroup%
\HardButStrongLineBreak\kern3pt%
\begingroup \smaller\smaller\smaller\begin{tabular}{@{}c@{}}%
-1\\-104\\0
\end{tabular}\endgroup%
\HardButStrongLineBreak\kern3pt%
\begingroup \smaller\smaller\smaller\begin{tabular}{@{}c@{}}%
-1\\-105\\-40
\end{tabular}\endgroup%
\HardButStrongLineBreak\kern3pt%
\begingroup \smaller\smaller\smaller\begin{tabular}{@{}c@{}}%
5\\518\\-84
\end{tabular}\endgroup%
\HardButStrongLineBreak\kern3pt%
\begingroup \smaller\smaller\smaller\begin{tabular}{@{}c@{}}%
3\\312\\-2
\end{tabular}\endgroup%
{$\left.\llap{\phantom{%
\begingroup \smaller\smaller\smaller\begin{tabular}{@{}c@{}}%
0\\0\\0
\end{tabular}\endgroup%
}}\!\right]$}%
%
%
\hbox{}\par\smallskip%
%
%
\leavevmode%
${L_{72.2}}$%
{} : {$1\above{1pt}{1pt}{2}{2}8\above{1pt}{1pt}{1}{7}{\cdot}1\above{1pt}{1pt}{2}{}3\above{1pt}{1pt}{-}{}{\cdot}1\above{1pt}{1pt}{-2}{}5\above{1pt}{1pt}{1}{}{\cdot}1\above{1pt}{1pt}{-2}{}7\above{1pt}{1pt}{1}{}$}\spacer%
\instructions{2\rightarrow N'_{34}}%
\EasyButWeakLineBreak%
{${24}\above{1pt}{1pt}{*}{2}{28}\above{1pt}{1pt}{l}{2}{1}\above{1pt}{1pt}{}{2}{105}\above{1pt}{1pt}{r}{2}{8}\above{1pt}{1pt}{l}{2}{5}\above{1pt}{1pt}{}{2}{56}\above{1pt}{1pt}{r}{2}{2}\above{1pt}{1pt}{b}{2}$}%
\nopagebreak\par%
\nopagebreak\par\leavevmode%
{$\left[\!\llap{\phantom{%
\begingroup \smaller\smaller\smaller\begin{tabular}{@{}c@{}}%
0\\0\\0
\end{tabular}\endgroup%
}}\right.$}%
\begingroup \smaller\smaller\smaller\begin{tabular}{@{}c@{}}%
-34131720\\52920\\26040
\end{tabular}\endgroup%
\kern3pt%
\begingroup \smaller\smaller\smaller\begin{tabular}{@{}c@{}}%
52920\\-79\\-42
\end{tabular}\endgroup%
\kern3pt%
\begingroup \smaller\smaller\smaller\begin{tabular}{@{}c@{}}%
26040\\-42\\-19
\end{tabular}\endgroup%
{$\left.\llap{\phantom{%
\begingroup \smaller\smaller\smaller\begin{tabular}{@{}c@{}}%
0\\0\\0
\end{tabular}\endgroup%
}}\!\right]$}%
\EasyButWeakLineBreak%
{$\left[\!\llap{\phantom{%
\begingroup \smaller\smaller\smaller\begin{tabular}{@{}c@{}}%
0\\0\\0
\end{tabular}\endgroup%
}}\right.$}%
\begingroup \smaller\smaller\smaller\begin{tabular}{@{}c@{}}%
17\\5700\\10692
\end{tabular}\endgroup%
\HardButStrongLineBreak\kern3pt%
\begingroup \smaller\smaller\smaller\begin{tabular}{@{}c@{}}%
17\\5698\\10696
\end{tabular}\endgroup%
\HardButStrongLineBreak\kern3pt%
\begingroup \smaller\smaller\smaller\begin{tabular}{@{}c@{}}%
2\\670\\1259
\end{tabular}\endgroup%
\HardButStrongLineBreak\kern3pt%
\begingroup \smaller\smaller\smaller\begin{tabular}{@{}c@{}}%
16\\5355\\10080
\end{tabular}\endgroup%
\HardButStrongLineBreak\kern3pt%
\begingroup \smaller\smaller\smaller\begin{tabular}{@{}c@{}}%
-1\\-336\\-628
\end{tabular}\endgroup%
\HardButStrongLineBreak\kern3pt%
\begingroup \smaller\smaller\smaller\begin{tabular}{@{}c@{}}%
-1\\-335\\-630
\end{tabular}\endgroup%
\HardButStrongLineBreak\kern3pt%
\begingroup \smaller\smaller\smaller\begin{tabular}{@{}c@{}}%
5\\1680\\3136
\end{tabular}\endgroup%
\HardButStrongLineBreak\kern3pt%
\begingroup \smaller\smaller\smaller\begin{tabular}{@{}c@{}}%
2\\671\\1257
\end{tabular}\endgroup%
{$\left.\llap{\phantom{%
\begingroup \smaller\smaller\smaller\begin{tabular}{@{}c@{}}%
0\\0\\0
\end{tabular}\endgroup%
}}\!\right]$}%
%
%
\hbox{}\par\smallskip%
%
%
\leavevmode%
${L_{72.3}}$%
{} : {$1\above{1pt}{1pt}{-2}{2}8\above{1pt}{1pt}{-}{3}{\cdot}1\above{1pt}{1pt}{2}{}3\above{1pt}{1pt}{-}{}{\cdot}1\above{1pt}{1pt}{-2}{}5\above{1pt}{1pt}{1}{}{\cdot}1\above{1pt}{1pt}{-2}{}7\above{1pt}{1pt}{1}{}$}\spacer%
\instructions{m}%
\EasyButWeakLineBreak%
{${24}\above{1pt}{1pt}{}{2}{7}\above{1pt}{1pt}{r}{2}{4}\above{1pt}{1pt}{*}{2}{420}\above{1pt}{1pt}{s}{2}{8}\above{1pt}{1pt}{s}{2}{20}\above{1pt}{1pt}{*}{2}{56}\above{1pt}{1pt}{b}{2}{2}\above{1pt}{1pt}{l}{2}$}%
\nopagebreak\par%
\nopagebreak\par\leavevmode%
{$\left[\!\llap{\phantom{%
\begingroup \smaller\smaller\smaller\begin{tabular}{@{}c@{}}%
0\\0\\0
\end{tabular}\endgroup%
}}\right.$}%
\begingroup \smaller\smaller\smaller\begin{tabular}{@{}c@{}}%
-434280\\-43680\\2520
\end{tabular}\endgroup%
\kern3pt%
\begingroup \smaller\smaller\smaller\begin{tabular}{@{}c@{}}%
-43680\\-4393\\253
\end{tabular}\endgroup%
\kern3pt%
\begingroup \smaller\smaller\smaller\begin{tabular}{@{}c@{}}%
2520\\253\\-14
\end{tabular}\endgroup%
{$\left.\llap{\phantom{%
\begingroup \smaller\smaller\smaller\begin{tabular}{@{}c@{}}%
0\\0\\0
\end{tabular}\endgroup%
}}\!\right]$}%
\EasyButWeakLineBreak%
{$\left[\!\llap{\phantom{%
\begingroup \smaller\smaller\smaller\begin{tabular}{@{}c@{}}%
0\\0\\0
\end{tabular}\endgroup%
}}\right.$}%
\begingroup \smaller\smaller\smaller\begin{tabular}{@{}c@{}}%
-23\\240\\192
\end{tabular}\endgroup%
\HardButStrongLineBreak\kern3pt%
\begingroup \smaller\smaller\smaller\begin{tabular}{@{}c@{}}%
-6\\63\\56
\end{tabular}\endgroup%
\HardButStrongLineBreak\kern3pt%
\begingroup \smaller\smaller\smaller\begin{tabular}{@{}c@{}}%
1\\-10\\-2
\end{tabular}\endgroup%
\HardButStrongLineBreak\kern3pt%
\begingroup \smaller\smaller\smaller\begin{tabular}{@{}c@{}}%
61\\-630\\-420
\end{tabular}\endgroup%
\HardButStrongLineBreak\kern3pt%
\begingroup \smaller\smaller\smaller\begin{tabular}{@{}c@{}}%
5\\-52\\-40
\end{tabular}\endgroup%
\HardButStrongLineBreak\kern3pt%
\begingroup \smaller\smaller\smaller\begin{tabular}{@{}c@{}}%
-1\\10\\0
\end{tabular}\endgroup%
\HardButStrongLineBreak\kern3pt%
\begingroup \smaller\smaller\smaller\begin{tabular}{@{}c@{}}%
-27\\280\\196
\end{tabular}\endgroup%
\HardButStrongLineBreak\kern3pt%
\begingroup \smaller\smaller\smaller\begin{tabular}{@{}c@{}}%
-5\\52\\39
\end{tabular}\endgroup%
{$\left.\llap{\phantom{%
\begingroup \smaller\smaller\smaller\begin{tabular}{@{}c@{}}%
0\\0\\0
\end{tabular}\endgroup%
}}\!\right]$}%

\medskip%
%
\leavevmode\llap{}%
$W_{73\phantom{0}}$%
\qquad\llap{32} lattices, $\chi=72$%
\hfill%
$\infty222222\infty222222\rtimes C_{2}$%
\nopagebreak\smallskip\hrule\nopagebreak\medskip%
%
%
\leavevmode%
${L_{73.1}}$%
{} : {$1\above{1pt}{1pt}{-2}{{\rm II}}8\above{1pt}{1pt}{-}{3}{\cdot}1\above{1pt}{1pt}{-}{}5\above{1pt}{1pt}{-}{}25\above{1pt}{1pt}{-}{}{\cdot}1\above{1pt}{1pt}{-2}{}11\above{1pt}{1pt}{-}{}$}\spacer%
\instructions{25\rightarrow N_{73},5,2*}%
\EasyButWeakLineBreak%
{${110}\above{1pt}{1pt}{20,19}{\infty a}{440}\above{1pt}{1pt}{b}{2}{50}\above{1pt}{1pt}{l}{2}{88}\above{1pt}{1pt}{r}{2}{10}\above{1pt}{1pt}{b}{2}{550}\above{1pt}{1pt}{s}{2}{2}\above{1pt}{1pt}{b}{2}$}\relax$\,(\times2)$%
\nopagebreak\par%
shares genus with 5-dual\nopagebreak\par%
\nopagebreak\par\leavevmode%
{$\left[\!\llap{\phantom{%
\begingroup \smaller\smaller\smaller
\endgroup%
}}\!\right]$}%

\medskip%
%
\leavevmode\llap{}%
$W_{74\phantom{0}}$%
\qquad\llap{12} lattices, $\chi=38$%
\hfill%
$22642264\rtimes C_{2}$%
\nopagebreak\smallskip\hrule\nopagebreak\medskip%
%
%
\leavevmode%
${L_{74.1}}$%
{} : {$1\above{1pt}{1pt}{-2}{{\rm II}}4\above{1pt}{1pt}{-}{3}{\cdot}1\above{1pt}{1pt}{2}{}3\above{1pt}{1pt}{-}{}{\cdot}1\above{1pt}{1pt}{2}{}37\above{1pt}{1pt}{1}{}$}\spacer%
\instructions{2\rightarrow N_{74}}%
\EasyButWeakLineBreak%
{${4}\above{1pt}{1pt}{*}{2}{148}\above{1pt}{1pt}{b}{2}{6}\above{1pt}{1pt}{}{6}{2}\above{1pt}{1pt}{*}{4}$}\relax$\,(\times2)$%
\nopagebreak\par%
\nopagebreak\par\leavevmode%
{$\left[\!\llap{\phantom{%
\begingroup \smaller\smaller\smaller\begin{tabular}{@{}c@{}}%
0\\0\\0
\end{tabular}\endgroup%
}}\right.$}%
\begingroup \smaller\smaller\smaller\begin{tabular}{@{}c@{}}%
-3066708\\21312\\20868
\end{tabular}\endgroup%
\kern3pt%
\begingroup \smaller\smaller\smaller\begin{tabular}{@{}c@{}}%
21312\\-146\\-145
\end{tabular}\endgroup%
\kern3pt%
\begingroup \smaller\smaller\smaller\begin{tabular}{@{}c@{}}%
20868\\-145\\-142
\end{tabular}\endgroup%
{$\left.\llap{\phantom{%
\begingroup \smaller\smaller\smaller\begin{tabular}{@{}c@{}}%
0\\0\\0
\end{tabular}\endgroup%
}}\!\right]$}%
\hfil\penalty500%
{$\left[\!\llap{\phantom{%
\begingroup \smaller\smaller\smaller\begin{tabular}{@{}c@{}}%
0\\0\\0
\end{tabular}\endgroup%
}}\right.$}%
\begingroup \smaller\smaller\smaller\begin{tabular}{@{}c@{}}%
124319\\-124320\\18381600
\end{tabular}\endgroup%
\kern3pt%
\begingroup \smaller\smaller\smaller\begin{tabular}{@{}c@{}}%
-917\\916\\-135585
\end{tabular}\endgroup%
\kern3pt%
\begingroup \smaller\smaller\smaller\begin{tabular}{@{}c@{}}%
-847\\847\\-125236
\end{tabular}\endgroup%
{$\left.\llap{\phantom{%
\begingroup \smaller\smaller\smaller\begin{tabular}{@{}c@{}}%
0\\0\\0
\end{tabular}\endgroup%
}}\!\right]$}%
\EasyButWeakLineBreak%
{$\left[\!\llap{\phantom{%
\begingroup \smaller\smaller\smaller\begin{tabular}{@{}c@{}}%
0\\0\\0
\end{tabular}\endgroup%
}}\right.$}%
\begingroup \smaller\smaller\smaller\begin{tabular}{@{}c@{}}%
37\\-40\\5474
\end{tabular}\endgroup%
\HardButStrongLineBreak\kern3pt%
\begingroup \smaller\smaller\smaller\begin{tabular}{@{}c@{}}%
327\\-370\\48396
\end{tabular}\endgroup%
\HardButStrongLineBreak\kern3pt%
\begingroup \smaller\smaller\smaller\begin{tabular}{@{}c@{}}%
1\\-3\\150
\end{tabular}\endgroup%
\HardButStrongLineBreak\kern3pt%
\begingroup \smaller\smaller\smaller\begin{tabular}{@{}c@{}}%
-1\\2\\-149
\end{tabular}\endgroup%
{$\left.\llap{\phantom{%
\begingroup \smaller\smaller\smaller\begin{tabular}{@{}c@{}}%
0\\0\\0
\end{tabular}\endgroup%
}}\!\right]$}%

\medskip%
%
\leavevmode\llap{}%
$W_{75\phantom{0}}$%
\qquad\llap{16} lattices, $\chi=48$%
\hfill%
$222222222222\rtimes C_{2}$%
\nopagebreak\smallskip\hrule\nopagebreak\medskip%
%
%
\leavevmode%
${L_{75.1}}$%
{} : {$1\above{1pt}{1pt}{-2}{{\rm II}}8\above{1pt}{1pt}{-}{5}{\cdot}1\above{1pt}{1pt}{-2}{}5\above{1pt}{1pt}{-}{}{\cdot}1\above{1pt}{1pt}{-2}{}13\above{1pt}{1pt}{-}{}$}\spacer%
\instructions{2\rightarrow N_{75}}%
\EasyButWeakLineBreak%
{${8}\above{1pt}{1pt}{b}{2}{26}\above{1pt}{1pt}{l}{2}{40}\above{1pt}{1pt}{r}{2}{2}\above{1pt}{1pt}{l}{2}{104}\above{1pt}{1pt}{r}{2}{10}\above{1pt}{1pt}{b}{2}$}\relax$\,(\times2)$%
\nopagebreak\par%
\nopagebreak\par\leavevmode%
{$\left[\!\llap{\phantom{%
\begingroup \smaller\smaller\smaller
\endgroup%
}}\!\right]$}%

\medskip%
%
\leavevmode\llap{}%
$W_{76\phantom{0}}$%
\qquad\llap{12} lattices, $\chi=48$%
\hfill%
$222222222222\rtimes C_{2}$%
\nopagebreak\smallskip\hrule\nopagebreak\medskip%
%
%
\leavevmode%
${L_{76.1}}$%
{} : {$1\above{1pt}{1pt}{-2}{{\rm II}}4\above{1pt}{1pt}{1}{1}{\cdot}1\above{1pt}{1pt}{2}{}3\above{1pt}{1pt}{1}{}{\cdot}1\above{1pt}{1pt}{-2}{}47\above{1pt}{1pt}{1}{}$}\spacer%
\instructions{2\rightarrow N_{76}}%
\EasyButWeakLineBreak%
{${12}\above{1pt}{1pt}{*}{2}{188}\above{1pt}{1pt}{b}{2}{2}\above{1pt}{1pt}{b}{2}{282}\above{1pt}{1pt}{l}{2}{4}\above{1pt}{1pt}{r}{2}{94}\above{1pt}{1pt}{b}{2}$}\relax$\,(\times2)$%
\nopagebreak\par%
\nopagebreak\par\leavevmode%
{$\left[\!\llap{\phantom{%
\begingroup \smaller\smaller\smaller
\endgroup%
}}\!\right]$}%

\medskip%
%
\leavevmode\llap{}%
$W_{77\phantom{0}}$%
\qquad\llap{12} lattices, $\chi=64$%
\hfill%
$22222232222223\rtimes C_{2}$%
\nopagebreak\smallskip\hrule\nopagebreak\medskip%
%
%
\leavevmode%
${L_{77.1}}$%
{} : {$1\above{1pt}{1pt}{-2}{{\rm II}}4\above{1pt}{1pt}{1}{7}{\cdot}1\above{1pt}{1pt}{-2}{}5\above{1pt}{1pt}{-}{}{\cdot}1\above{1pt}{1pt}{-2}{}31\above{1pt}{1pt}{1}{}$}\spacer%
\instructions{2\rightarrow N_{77}}%
\EasyButWeakLineBreak%
{${2}\above{1pt}{1pt}{s}{2}{310}\above{1pt}{1pt}{b}{2}{4}\above{1pt}{1pt}{b}{2}{62}\above{1pt}{1pt}{s}{2}{10}\above{1pt}{1pt}{l}{2}{124}\above{1pt}{1pt}{r}{2}{2}\above{1pt}{1pt}{+}{3}$}\relax$\,(\times2)$%
\nopagebreak\par%
\nopagebreak\par\leavevmode%
{$\left[\!\llap{\phantom{%
\begingroup \smaller\smaller\smaller
\endgroup%
}}\!\right]$}%

\medskip%
%
\leavevmode\llap{}%
$W_{78\phantom{0}}$%
\qquad\llap{24} lattices, $\chi=10$%
\hfill%
$26222$%
\nopagebreak\smallskip\hrule\nopagebreak\medskip%
%
%
\leavevmode%
${L_{78.1}}$%
{} : {$1\above{1pt}{1pt}{-2}{{\rm II}}4\above{1pt}{1pt}{1}{1}{\cdot}1\above{1pt}{1pt}{2}{}3\above{1pt}{1pt}{-}{}{\cdot}1\above{1pt}{1pt}{-2}{}5\above{1pt}{1pt}{1}{}{\cdot}1\above{1pt}{1pt}{2}{}11\above{1pt}{1pt}{-}{}$}\spacer%
\instructions{2\rightarrow N_{78}}%
\EasyButWeakLineBreak%
{${22}\above{1pt}{1pt}{b}{2}{6}\above{1pt}{1pt}{}{6}{2}\above{1pt}{1pt}{b}{2}{330}\above{1pt}{1pt}{l}{2}{4}\above{1pt}{1pt}{r}{2}$}%
\nopagebreak\par%
\nopagebreak\par\leavevmode%
{$\left[\!\llap{\phantom{%
\begingroup \smaller\smaller\smaller\begin{tabular}{@{}c@{}}%
0\\0\\0
\end{tabular}\endgroup%
}}\right.$}%
\begingroup \smaller\smaller\smaller\begin{tabular}{@{}c@{}}%
-34021020\\23100\\36960
\end{tabular}\endgroup%
\kern3pt%
\begingroup \smaller\smaller\smaller\begin{tabular}{@{}c@{}}%
23100\\-14\\-27
\end{tabular}\endgroup%
\kern3pt%
\begingroup \smaller\smaller\smaller\begin{tabular}{@{}c@{}}%
36960\\-27\\-38
\end{tabular}\endgroup%
{$\left.\llap{\phantom{%
\begingroup \smaller\smaller\smaller\begin{tabular}{@{}c@{}}%
0\\0\\0
\end{tabular}\endgroup%
}}\!\right]$}%
\EasyButWeakLineBreak%
{$\left[\!\llap{\phantom{%
\begingroup \smaller\smaller\smaller\begin{tabular}{@{}c@{}}%
0\\0\\0
\end{tabular}\endgroup%
}}\right.$}%
\begingroup \smaller\smaller\smaller\begin{tabular}{@{}c@{}}%
-2\\-1221\\-1078
\end{tabular}\endgroup%
\HardButStrongLineBreak\kern3pt%
\begingroup \smaller\smaller\smaller\begin{tabular}{@{}c@{}}%
-1\\-609\\-540
\end{tabular}\endgroup%
\HardButStrongLineBreak\kern3pt%
\begingroup \smaller\smaller\smaller\begin{tabular}{@{}c@{}}%
1\\610\\539
\end{tabular}\endgroup%
\HardButStrongLineBreak\kern3pt%
\begingroup \smaller\smaller\smaller\begin{tabular}{@{}c@{}}%
26\\15840\\14025
\end{tabular}\endgroup%
\HardButStrongLineBreak\kern3pt%
\begingroup \smaller\smaller\smaller\begin{tabular}{@{}c@{}}%
1\\608\\540
\end{tabular}\endgroup%
{$\left.\llap{\phantom{%
\begingroup \smaller\smaller\smaller\begin{tabular}{@{}c@{}}%
0\\0\\0
\end{tabular}\endgroup%
}}\!\right]$}%

\medskip%
%
\leavevmode\llap{}%
$W_{79\phantom{0}}$%
\qquad\llap{48} lattices, $\chi=60$%
\hfill%
$22222222222222\rtimes C_{2}$%
\nopagebreak\smallskip\hrule\nopagebreak\medskip%
%
%
\leavevmode%
${L_{79.1}}$%
{} : {$1\above{1pt}{1pt}{-2}{{\rm II}}4\above{1pt}{1pt}{1}{1}{\cdot}1\above{1pt}{1pt}{1}{}3\above{1pt}{1pt}{1}{}9\above{1pt}{1pt}{-}{}{\cdot}1\above{1pt}{1pt}{2}{}5\above{1pt}{1pt}{-}{}{\cdot}1\above{1pt}{1pt}{2}{}11\above{1pt}{1pt}{-}{}$}\spacer%
\instructions{23\rightarrow N_{79},3,2}%
\EasyButWeakLineBreak%
{${18}\above{1pt}{1pt}{s}{2}{22}\above{1pt}{1pt}{b}{2}{12}\above{1pt}{1pt}{b}{2}{990}\above{1pt}{1pt}{l}{2}{4}\above{1pt}{1pt}{r}{2}{66}\above{1pt}{1pt}{b}{2}{10}\above{1pt}{1pt}{b}{2}$}\relax$\,(\times2)$%
\nopagebreak\par%
\nopagebreak\par\leavevmode%
{$\left[\!\llap{\phantom{%
\begingroup \smaller\smaller\smaller
\endgroup%
}}\!\right]$}%

\medskip%
%
\leavevmode\llap{}%
$W_{80\phantom{0}}$%
\qquad\llap{88} lattices, $\chi=72$%
\hfill%
$2222222222222222\rtimes C_{2}$%
\nopagebreak\smallskip\hrule\nopagebreak\medskip%
%
%
\leavevmode%
${L_{80.1}}$%
{} : {$1\above{1pt}{1pt}{2}{{\rm II}}4\above{1pt}{1pt}{-}{5}{\cdot}1\above{1pt}{1pt}{2}{}3\above{1pt}{1pt}{-}{}{\cdot}1\above{1pt}{1pt}{-2}{}5\above{1pt}{1pt}{1}{}{\cdot}1\above{1pt}{1pt}{-2}{}11\above{1pt}{1pt}{1}{}$}\spacer%
\instructions{2\rightarrow N_{80}}%
\EasyButWeakLineBreak%
{${44}\above{1pt}{1pt}{*}{2}{20}\above{1pt}{1pt}{*}{2}{132}\above{1pt}{1pt}{*}{2}{4}\above{1pt}{1pt}{*}{2}{220}\above{1pt}{1pt}{b}{2}{6}\above{1pt}{1pt}{l}{2}{20}\above{1pt}{1pt}{r}{2}{2}\above{1pt}{1pt}{b}{2}$}\relax$\,(\times2)$%
\nopagebreak\par%
\nopagebreak\par\leavevmode%
{$\left[\!\llap{\phantom{%
\begingroup \smaller\smaller\smaller
\endgroup%
}}\!\right]$}%
%
%
\hbox{}\par\smallskip%
%
%
\leavevmode%
${L_{80.2}}$%
{} : {$1\above{1pt}{1pt}{2}{2}8\above{1pt}{1pt}{-}{3}{\cdot}1\above{1pt}{1pt}{2}{}3\above{1pt}{1pt}{1}{}{\cdot}1\above{1pt}{1pt}{-2}{}5\above{1pt}{1pt}{-}{}{\cdot}1\above{1pt}{1pt}{-2}{}11\above{1pt}{1pt}{-}{}$}\spacer%
\instructions{2\rightarrow N'_{36}}%
\EasyButWeakLineBreak%
{${88}\above{1pt}{1pt}{r}{2}{10}\above{1pt}{1pt}{b}{2}{66}\above{1pt}{1pt}{s}{2}{2}\above{1pt}{1pt}{b}{2}{440}\above{1pt}{1pt}{*}{2}{12}\above{1pt}{1pt}{s}{2}{40}\above{1pt}{1pt}{l}{2}{1}\above{1pt}{1pt}{}{2}$}\relax$\,(\times2)$%
\nopagebreak\par%
\nopagebreak\par\leavevmode%
{$\left[\!\llap{\phantom{%
\begingroup \smaller\smaller\smaller
\endgroup%
}}\!\right]$}%
%
%
\hbox{}\par\smallskip%
%
%
\leavevmode%
${L_{80.3}}$%
{} : {$1\above{1pt}{1pt}{-2}{2}8\above{1pt}{1pt}{1}{7}{\cdot}1\above{1pt}{1pt}{2}{}3\above{1pt}{1pt}{1}{}{\cdot}1\above{1pt}{1pt}{-2}{}5\above{1pt}{1pt}{-}{}{\cdot}1\above{1pt}{1pt}{-2}{}11\above{1pt}{1pt}{-}{}$}\spacer%
\instructions{m}%
\EasyButWeakLineBreak%
{${88}\above{1pt}{1pt}{b}{2}{10}\above{1pt}{1pt}{s}{2}{66}\above{1pt}{1pt}{b}{2}{2}\above{1pt}{1pt}{l}{2}{440}\above{1pt}{1pt}{}{2}{3}\above{1pt}{1pt}{r}{2}{40}\above{1pt}{1pt}{s}{2}{4}\above{1pt}{1pt}{*}{2}$}\relax$\,(\times2)$%
\nopagebreak\par%
\nopagebreak\par\leavevmode%
{$\left[\!\llap{\phantom{%
\begingroup \smaller\smaller\smaller
\endgroup%
}}\!\right]$}%

\medskip%
%
\leavevmode\llap{}%
$W_{81\phantom{0}}$%
\qquad\llap{48} lattices, $\chi=48$%
\hfill%
$222222222222\rtimes C_{2}$%
\nopagebreak\smallskip\hrule\nopagebreak\medskip%
%
%
\leavevmode%
${L_{81.1}}$%
{} : {$1\above{1pt}{1pt}{-2}{{\rm II}}4\above{1pt}{1pt}{1}{1}{\cdot}1\above{1pt}{1pt}{-}{}3\above{1pt}{1pt}{1}{}9\above{1pt}{1pt}{1}{}{\cdot}1\above{1pt}{1pt}{-2}{}5\above{1pt}{1pt}{1}{}{\cdot}1\above{1pt}{1pt}{-2}{}11\above{1pt}{1pt}{1}{}$}\spacer%
\instructions{23\rightarrow N_{81},3,2}%
\EasyButWeakLineBreak%
{${12}\above{1pt}{1pt}{*}{2}{1980}\above{1pt}{1pt}{b}{2}{2}\above{1pt}{1pt}{l}{2}{36}\above{1pt}{1pt}{r}{2}{30}\above{1pt}{1pt}{b}{2}{44}\above{1pt}{1pt}{*}{2}$}\relax$\,(\times2)$%
\nopagebreak\par%
\nopagebreak\par\leavevmode%
{$\left[\!\llap{\phantom{%
\begingroup \smaller\smaller\smaller
\endgroup%
}}\!\right]$}%

\medskip%
%
\leavevmode\llap{}%
$W_{82\phantom{0}}$%
\qquad\llap{24} lattices, $\chi=18$%
\hfill%
$2222222$%
\nopagebreak\smallskip\hrule\nopagebreak\medskip%
%
%
\leavevmode%
${L_{82.1}}$%
{} : {$1\above{1pt}{1pt}{-2}{{\rm II}}4\above{1pt}{1pt}{1}{1}{\cdot}1\above{1pt}{1pt}{2}{}3\above{1pt}{1pt}{-}{}{\cdot}1\above{1pt}{1pt}{2}{}5\above{1pt}{1pt}{-}{}{\cdot}1\above{1pt}{1pt}{-2}{}11\above{1pt}{1pt}{1}{}$}\spacer%
\instructions{2\rightarrow N_{82}}%
\EasyButWeakLineBreak%
{${44}\above{1pt}{1pt}{*}{2}{60}\above{1pt}{1pt}{b}{2}{2}\above{1pt}{1pt}{l}{2}{132}\above{1pt}{1pt}{r}{2}{10}\above{1pt}{1pt}{l}{2}{4}\above{1pt}{1pt}{r}{2}{6}\above{1pt}{1pt}{b}{2}$}%
\nopagebreak\par%
\nopagebreak\par\leavevmode%
{$\left[\!\llap{\phantom{%
\begingroup \smaller\smaller\smaller\begin{tabular}{@{}c@{}}%
0\\0\\0
\end{tabular}\endgroup%
}}\right.$}%
\begingroup \smaller\smaller\smaller\begin{tabular}{@{}c@{}}%
-13682460\\-6811200\\14520
\end{tabular}\endgroup%
\kern3pt%
\begingroup \smaller\smaller\smaller\begin{tabular}{@{}c@{}}%
-6811200\\-3390650\\7227
\end{tabular}\endgroup%
\kern3pt%
\begingroup \smaller\smaller\smaller\begin{tabular}{@{}c@{}}%
14520\\7227\\-14
\end{tabular}\endgroup%
{$\left.\llap{\phantom{%
\begingroup \smaller\smaller\smaller\begin{tabular}{@{}c@{}}%
0\\0\\0
\end{tabular}\endgroup%
}}\!\right]$}%
\EasyButWeakLineBreak%
{$\left[\!\llap{\phantom{%
\begingroup \smaller\smaller\smaller\begin{tabular}{@{}c@{}}%
0\\0\\0
\end{tabular}\endgroup%
}}\right.$}%
\begingroup \smaller\smaller\smaller\begin{tabular}{@{}c@{}}%
-2285\\4598\\3674
\end{tabular}\endgroup%
\HardButStrongLineBreak\kern3pt%
\begingroup \smaller\smaller\smaller\begin{tabular}{@{}c@{}}%
-2773\\5580\\4470
\end{tabular}\endgroup%
\HardButStrongLineBreak\kern3pt%
\begingroup \smaller\smaller\smaller\begin{tabular}{@{}c@{}}%
-326\\656\\527
\end{tabular}\endgroup%
\HardButStrongLineBreak\kern3pt%
\begingroup \smaller\smaller\smaller\begin{tabular}{@{}c@{}}%
-4723\\9504\\7656
\end{tabular}\endgroup%
\HardButStrongLineBreak\kern3pt%
\begingroup \smaller\smaller\smaller\begin{tabular}{@{}c@{}}%
82\\-165\\-130
\end{tabular}\endgroup%
\HardButStrongLineBreak\kern3pt%
\begingroup \smaller\smaller\smaller\begin{tabular}{@{}c@{}}%
163\\-328\\-264
\end{tabular}\endgroup%
\HardButStrongLineBreak\kern3pt%
\begingroup \smaller\smaller\smaller\begin{tabular}{@{}c@{}}%
-82\\165\\129
\end{tabular}\endgroup%
{$\left.\llap{\phantom{%
\begingroup \smaller\smaller\smaller\begin{tabular}{@{}c@{}}%
0\\0\\0
\end{tabular}\endgroup%
}}\!\right]$}%

\medskip%
%
\leavevmode\llap{}%
$W_{83\phantom{0}}$%
\qquad\llap{16} lattices, $\chi=64$%
\hfill%
$23222222322222\rtimes C_{2}$%
\nopagebreak\smallskip\hrule\nopagebreak\medskip%
%
%
\leavevmode%
${L_{83.1}}$%
{} : {$1\above{1pt}{1pt}{-2}{{\rm II}}8\above{1pt}{1pt}{1}{1}{\cdot}1\above{1pt}{1pt}{-2}{}5\above{1pt}{1pt}{-}{}{\cdot}1\above{1pt}{1pt}{-2}{}17\above{1pt}{1pt}{1}{}$}\spacer%
\instructions{2\rightarrow N_{83}}%
\EasyButWeakLineBreak%
{${40}\above{1pt}{1pt}{b}{2}{2}\above{1pt}{1pt}{-}{3}{2}\above{1pt}{1pt}{l}{2}{136}\above{1pt}{1pt}{r}{2}{10}\above{1pt}{1pt}{l}{2}{8}\above{1pt}{1pt}{r}{2}{34}\above{1pt}{1pt}{b}{2}$}\relax$\,(\times2)$%
\nopagebreak\par%
\nopagebreak\par\leavevmode%
{$\left[\!\llap{\phantom{%
\begingroup \smaller\smaller\smaller
\endgroup%
}}\!\right]$}%

\medskip%
%
\leavevmode\llap{}%
$W_{84\phantom{0}}$%
\qquad\llap{24} lattices, $\chi=12$%
\hfill%
$222222$%
\nopagebreak\smallskip\hrule\nopagebreak\medskip%
%
%
\leavevmode%
${L_{84.1}}$%
{} : {$1\above{1pt}{1pt}{-2}{{\rm II}}4\above{1pt}{1pt}{1}{7}{\cdot}1\above{1pt}{1pt}{2}{}3\above{1pt}{1pt}{1}{}{\cdot}1\above{1pt}{1pt}{-2}{}5\above{1pt}{1pt}{-}{}{\cdot}1\above{1pt}{1pt}{-2}{}13\above{1pt}{1pt}{1}{}$}\spacer%
\instructions{2\rightarrow N_{84}}%
\EasyButWeakLineBreak%
{${4}\above{1pt}{1pt}{*}{2}{52}\above{1pt}{1pt}{b}{2}{10}\above{1pt}{1pt}{l}{2}{156}\above{1pt}{1pt}{r}{2}{2}\above{1pt}{1pt}{s}{2}{390}\above{1pt}{1pt}{b}{2}$}%
\nopagebreak\par%
\nopagebreak\par\leavevmode%
{$\left[\!\llap{\phantom{%
\begingroup \smaller\smaller\smaller\begin{tabular}{@{}c@{}}%
0\\0\\0
\end{tabular}\endgroup%
}}\right.$}%
\begingroup \smaller\smaller\smaller\begin{tabular}{@{}c@{}}%
-11003460\\-3664440\\15600
\end{tabular}\endgroup%
\kern3pt%
\begingroup \smaller\smaller\smaller\begin{tabular}{@{}c@{}}%
-3664440\\-1220354\\5195
\end{tabular}\endgroup%
\kern3pt%
\begingroup \smaller\smaller\smaller\begin{tabular}{@{}c@{}}%
15600\\5195\\-22
\end{tabular}\endgroup%
{$\left.\llap{\phantom{%
\begingroup \smaller\smaller\smaller\begin{tabular}{@{}c@{}}%
0\\0\\0
\end{tabular}\endgroup%
}}\!\right]$}%
\EasyButWeakLineBreak%
{$\left[\!\llap{\phantom{%
\begingroup \smaller\smaller\smaller\begin{tabular}{@{}c@{}}%
0\\0\\0
\end{tabular}\endgroup%
}}\right.$}%
\begingroup \smaller\smaller\smaller\begin{tabular}{@{}c@{}}%
39\\-118\\-210
\end{tabular}\endgroup%
\HardButStrongLineBreak\kern3pt%
\begingroup \smaller\smaller\smaller\begin{tabular}{@{}c@{}}%
43\\-130\\-208
\end{tabular}\endgroup%
\HardButStrongLineBreak\kern3pt%
\begingroup \smaller\smaller\smaller\begin{tabular}{@{}c@{}}%
-38\\115\\210
\end{tabular}\endgroup%
\HardButStrongLineBreak\kern3pt%
\begingroup \smaller\smaller\smaller\begin{tabular}{@{}c@{}}%
-361\\1092\\1872
\end{tabular}\endgroup%
\HardButStrongLineBreak\kern3pt%
\begingroup \smaller\smaller\smaller\begin{tabular}{@{}c@{}}%
-1\\3\\-1
\end{tabular}\endgroup%
\HardButStrongLineBreak\kern3pt%
\begingroup \smaller\smaller\smaller\begin{tabular}{@{}c@{}}%
451\\-1365\\-2535
\end{tabular}\endgroup%
{$\left.\llap{\phantom{%
\begingroup \smaller\smaller\smaller\begin{tabular}{@{}c@{}}%
0\\0\\0
\end{tabular}\endgroup%
}}\!\right]$}%

\medskip%
%
\leavevmode\llap{}%
$W_{85\phantom{0}}$%
\qquad\llap{24} lattices, $\chi=42$%
\hfill%
$2224222242\rtimes C_{2}$%
\nopagebreak\smallskip\hrule\nopagebreak\medskip%
%
%
\leavevmode%
${L_{85.1}}$%
{} : {$1\above{1pt}{1pt}{-2}{{\rm II}}4\above{1pt}{1pt}{1}{7}{\cdot}1\above{1pt}{1pt}{2}{}3\above{1pt}{1pt}{1}{}{\cdot}1\above{1pt}{1pt}{2}{}5\above{1pt}{1pt}{1}{}{\cdot}1\above{1pt}{1pt}{2}{}13\above{1pt}{1pt}{-}{}$}\spacer%
\instructions{2\rightarrow N_{85}}%
\EasyButWeakLineBreak%
{${20}\above{1pt}{1pt}{b}{2}{26}\above{1pt}{1pt}{s}{2}{30}\above{1pt}{1pt}{b}{2}{4}\above{1pt}{1pt}{*}{4}{2}\above{1pt}{1pt}{b}{2}$}\relax$\,(\times2)$%
\nopagebreak\par%
\nopagebreak\par\leavevmode%
{$\left[\!\llap{\phantom{%
\begingroup \smaller\smaller\smaller
\endgroup%
}}\!\right]$}%

\medskip%
%
\leavevmode\llap{}%
$W_{86\phantom{0}}$%
\qquad\llap{48} lattices, $\chi=28$%
\hfill%
$22622222$%
\nopagebreak\smallskip\hrule\nopagebreak\medskip%
%
%
\leavevmode%
${L_{86.1}}$%
{} : {$1\above{1pt}{1pt}{-2}{{\rm II}}4\above{1pt}{1pt}{1}{7}{\cdot}1\above{1pt}{1pt}{1}{}3\above{1pt}{1pt}{-}{}9\above{1pt}{1pt}{-}{}{\cdot}1\above{1pt}{1pt}{-2}{}5\above{1pt}{1pt}{-}{}{\cdot}1\above{1pt}{1pt}{2}{}13\above{1pt}{1pt}{-}{}$}\spacer%
\instructions{23\rightarrow N_{86},3,2}%
\EasyButWeakLineBreak%
{${4}\above{1pt}{1pt}{*}{2}{2340}\above{1pt}{1pt}{b}{2}{6}\above{1pt}{1pt}{}{6}{18}\above{1pt}{1pt}{l}{2}{60}\above{1pt}{1pt}{r}{2}{234}\above{1pt}{1pt}{b}{2}{10}\above{1pt}{1pt}{s}{2}{78}\above{1pt}{1pt}{b}{2}$}%
\nopagebreak\par%
\nopagebreak\par\leavevmode%
{$\left[\!\llap{\phantom{%
\begingroup \smaller\smaller\smaller\begin{tabular}{@{}c@{}}%
0\\0\\0
\end{tabular}\endgroup%
}}\right.$}%
\begingroup \smaller\smaller\smaller\begin{tabular}{@{}c@{}}%
-143865540\\159120\\10754640
\end{tabular}\endgroup%
\kern3pt%
\begingroup \smaller\smaller\smaller\begin{tabular}{@{}c@{}}%
159120\\-174\\-12363
\end{tabular}\endgroup%
\kern3pt%
\begingroup \smaller\smaller\smaller\begin{tabular}{@{}c@{}}%
10754640\\-12363\\-693998
\end{tabular}\endgroup%
{$\left.\llap{\phantom{%
\begingroup \smaller\smaller\smaller\begin{tabular}{@{}c@{}}%
0\\0\\0
\end{tabular}\endgroup%
}}\!\right]$}%
\EasyButWeakLineBreak%
{$\left[\!\llap{\phantom{%
\begingroup \smaller\smaller\smaller\begin{tabular}{@{}c@{}}%
0\\0\\0
\end{tabular}\endgroup%
}}\right.$}%
\begingroup \smaller\smaller\smaller\begin{tabular}{@{}c@{}}%
3027\\2125402\\9046
\end{tabular}\endgroup%
\HardButStrongLineBreak\kern3pt%
\begingroup \smaller\smaller\smaller\begin{tabular}{@{}c@{}}%
112363\\78895440\\335790
\end{tabular}\endgroup%
\HardButStrongLineBreak\kern3pt%
\begingroup \smaller\smaller\smaller\begin{tabular}{@{}c@{}}%
260\\182558\\777
\end{tabular}\endgroup%
\HardButStrongLineBreak\kern3pt%
\begingroup \smaller\smaller\smaller\begin{tabular}{@{}c@{}}%
-259\\-181857\\-774
\end{tabular}\endgroup%
\HardButStrongLineBreak\kern3pt%
\begingroup \smaller\smaller\smaller\begin{tabular}{@{}c@{}}%
-261\\-183260\\-780
\end{tabular}\endgroup%
\HardButStrongLineBreak\kern3pt%
\begingroup \smaller\smaller\smaller\begin{tabular}{@{}c@{}}%
1292\\907179\\3861
\end{tabular}\endgroup%
\HardButStrongLineBreak\kern3pt%
\begingroup \smaller\smaller\smaller\begin{tabular}{@{}c@{}}%
691\\485185\\2065
\end{tabular}\endgroup%
\HardButStrongLineBreak\kern3pt%
\begingroup \smaller\smaller\smaller\begin{tabular}{@{}c@{}}%
6225\\4370873\\18603
\end{tabular}\endgroup%
{$\left.\llap{\phantom{%
\begingroup \smaller\smaller\smaller\begin{tabular}{@{}c@{}}%
0\\0\\0
\end{tabular}\endgroup%
}}\!\right]$}%

\medskip%
%
\leavevmode\llap{}%
$W_{87\phantom{0}}$%
\qquad\llap{12} lattices, $\chi=84$%
\hfill%
$222222222222222222\rtimes C_{2}$%
\nopagebreak\smallskip\hrule\nopagebreak\medskip%
%
%
\leavevmode%
${L_{87.1}}$%
{} : {$1\above{1pt}{1pt}{-2}{{\rm II}}4\above{1pt}{1pt}{1}{1}{\cdot}1\above{1pt}{1pt}{-2}{}5\above{1pt}{1pt}{-}{}{\cdot}1\above{1pt}{1pt}{2}{}41\above{1pt}{1pt}{1}{}$}\spacer%
\instructions{2\rightarrow N_{87}}%
\EasyButWeakLineBreak%
{${4}\above{1pt}{1pt}{r}{2}{10}\above{1pt}{1pt}{b}{2}{82}\above{1pt}{1pt}{l}{2}{4}\above{1pt}{1pt}{r}{2}{410}\above{1pt}{1pt}{b}{2}{2}\above{1pt}{1pt}{b}{2}{10}\above{1pt}{1pt}{l}{2}{164}\above{1pt}{1pt}{r}{2}{2}\above{1pt}{1pt}{l}{2}$}\relax$\,(\times2)$%
\nopagebreak\par%
\nopagebreak\par\leavevmode%
{$\left[\!\llap{\phantom{%
\begingroup \smaller\smaller\smaller
\endgroup%
}}\!\right]$}%

\medskip%
%
\leavevmode\llap{}%
$W_{88\phantom{0}}$%
\qquad\llap{32} lattices, $\chi=18$%
\hfill%
$2222222$%
\nopagebreak\smallskip\hrule\nopagebreak\medskip%
%
%
\leavevmode%
${L_{88.1}}$%
{} : {$1\above{1pt}{1pt}{-2}{{\rm II}}8\above{1pt}{1pt}{-}{5}{\cdot}1\above{1pt}{1pt}{2}{}3\above{1pt}{1pt}{-}{}{\cdot}1\above{1pt}{1pt}{2}{}5\above{1pt}{1pt}{-}{}{\cdot}1\above{1pt}{1pt}{2}{}7\above{1pt}{1pt}{-}{}$}\spacer%
\instructions{2\rightarrow N_{88}}%
\EasyButWeakLineBreak%
{${168}\above{1pt}{1pt}{r}{2}{10}\above{1pt}{1pt}{b}{2}{6}\above{1pt}{1pt}{l}{2}{40}\above{1pt}{1pt}{r}{2}{42}\above{1pt}{1pt}{b}{2}{8}\above{1pt}{1pt}{b}{2}{2}\above{1pt}{1pt}{l}{2}$}%
\nopagebreak\par%
\nopagebreak\par\leavevmode%
{$\left[\!\llap{\phantom{%
\begingroup \smaller\smaller\smaller\begin{tabular}{@{}c@{}}%
0\\0\\0
\end{tabular}\endgroup%
}}\right.$}%
\begingroup \smaller\smaller\smaller\begin{tabular}{@{}c@{}}%
-620760\\-248640\\2520
\end{tabular}\endgroup%
\kern3pt%
\begingroup \smaller\smaller\smaller\begin{tabular}{@{}c@{}}%
-248640\\-99590\\1009
\end{tabular}\endgroup%
\kern3pt%
\begingroup \smaller\smaller\smaller\begin{tabular}{@{}c@{}}%
2520\\1009\\-10
\end{tabular}\endgroup%
{$\left.\llap{\phantom{%
\begingroup \smaller\smaller\smaller\begin{tabular}{@{}c@{}}%
0\\0\\0
\end{tabular}\endgroup%
}}\!\right]$}%
\EasyButWeakLineBreak%
{$\left[\!\llap{\phantom{%
\begingroup \smaller\smaller\smaller\begin{tabular}{@{}c@{}}%
0\\0\\0
\end{tabular}\endgroup%
}}\right.$}%
\begingroup \smaller\smaller\smaller\begin{tabular}{@{}c@{}}%
-265\\672\\1008
\end{tabular}\endgroup%
\HardButStrongLineBreak\kern3pt%
\begingroup \smaller\smaller\smaller\begin{tabular}{@{}c@{}}%
-2\\5\\0
\end{tabular}\endgroup%
\HardButStrongLineBreak\kern3pt%
\begingroup \smaller\smaller\smaller\begin{tabular}{@{}c@{}}%
13\\-33\\-54
\end{tabular}\endgroup%
\HardButStrongLineBreak\kern3pt%
\begingroup \smaller\smaller\smaller\begin{tabular}{@{}c@{}}%
79\\-200\\-280
\end{tabular}\endgroup%
\HardButStrongLineBreak\kern3pt%
\begingroup \smaller\smaller\smaller\begin{tabular}{@{}c@{}}%
25\\-63\\-63
\end{tabular}\endgroup%
\HardButStrongLineBreak\kern3pt%
\begingroup \smaller\smaller\smaller\begin{tabular}{@{}c@{}}%
-11\\28\\52
\end{tabular}\endgroup%
\HardButStrongLineBreak\kern3pt%
\begingroup \smaller\smaller\smaller\begin{tabular}{@{}c@{}}%
-13\\33\\53
\end{tabular}\endgroup%
{$\left.\llap{\phantom{%
\begingroup \smaller\smaller\smaller\begin{tabular}{@{}c@{}}%
0\\0\\0
\end{tabular}\endgroup%
}}\!\right]$}%

\medskip%
%
\leavevmode\llap{}%
$W_{89\phantom{0}}$%
\qquad\llap{48} lattices, $\chi=32$%
\hfill%
$22262226\rtimes C_{2}$%
\nopagebreak\smallskip\hrule\nopagebreak\medskip%
%
%
\leavevmode%
${L_{89.1}}$%
{} : {$1\above{1pt}{1pt}{-2}{{\rm II}}8\above{1pt}{1pt}{-}{5}{\cdot}1\above{1pt}{1pt}{-}{}3\above{1pt}{1pt}{-}{}9\above{1pt}{1pt}{-}{}{\cdot}1\above{1pt}{1pt}{-2}{}5\above{1pt}{1pt}{1}{}{\cdot}1\above{1pt}{1pt}{-2}{}7\above{1pt}{1pt}{1}{}$}\spacer%
\instructions{23\rightarrow N_{89},3,2}%
\EasyButWeakLineBreak%
{${18}\above{1pt}{1pt}{b}{2}{14}\above{1pt}{1pt}{b}{2}{72}\above{1pt}{1pt}{b}{2}{6}\above{1pt}{1pt}{}{6}{2}\above{1pt}{1pt}{b}{2}{126}\above{1pt}{1pt}{b}{2}{8}\above{1pt}{1pt}{b}{2}{6}\above{1pt}{1pt}{}{6}$}%
\nopagebreak\par%
\nopagebreak\par\leavevmode%
{$\left[\!\llap{\phantom{%
\begingroup \smaller\smaller\smaller\begin{tabular}{@{}c@{}}%
0\\0\\0
\end{tabular}\endgroup%
}}\right.$}%
\begingroup \smaller\smaller\smaller\begin{tabular}{@{}c@{}}%
-22538902680\\-1942920\\19106640
\end{tabular}\endgroup%
\kern3pt%
\begingroup \smaller\smaller\smaller\begin{tabular}{@{}c@{}}%
-1942920\\-138\\1629
\end{tabular}\endgroup%
\kern3pt%
\begingroup \smaller\smaller\smaller\begin{tabular}{@{}c@{}}%
19106640\\1629\\-16186
\end{tabular}\endgroup%
{$\left.\llap{\phantom{%
\begingroup \smaller\smaller\smaller\begin{tabular}{@{}c@{}}%
0\\0\\0
\end{tabular}\endgroup%
}}\!\right]$}%
\EasyButWeakLineBreak%
{$\left[\!\llap{\phantom{%
\begingroup \smaller\smaller\smaller\begin{tabular}{@{}c@{}}%
0\\0\\0
\end{tabular}\endgroup%
}}\right.$}%
\begingroup \smaller\smaller\smaller\begin{tabular}{@{}c@{}}%
-397\\-305694\\-499401
\end{tabular}\endgroup%
\HardButStrongLineBreak\kern3pt%
\begingroup \smaller\smaller\smaller\begin{tabular}{@{}c@{}}%
-78\\-60060\\-98119
\end{tabular}\endgroup%
\HardButStrongLineBreak\kern3pt%
\begingroup \smaller\smaller\smaller\begin{tabular}{@{}c@{}}%
35\\26952\\44028
\end{tabular}\endgroup%
\HardButStrongLineBreak\kern3pt%
\begingroup \smaller\smaller\smaller\begin{tabular}{@{}c@{}}%
16\\12320\\20127
\end{tabular}\endgroup%
\HardButStrongLineBreak\kern3pt%
\begingroup \smaller\smaller\smaller\begin{tabular}{@{}c@{}}%
-17\\-13091\\-21385
\end{tabular}\endgroup%
\HardButStrongLineBreak\kern3pt%
\begingroup \smaller\smaller\smaller\begin{tabular}{@{}c@{}}%
-214\\-164787\\-269199
\end{tabular}\endgroup%
\HardButStrongLineBreak\kern3pt%
\begingroup \smaller\smaller\smaller\begin{tabular}{@{}c@{}}%
-97\\-74692\\-122020
\end{tabular}\endgroup%
\HardButStrongLineBreak\kern3pt%
\begingroup \smaller\smaller\smaller\begin{tabular}{@{}c@{}}%
-208\\-160163\\-261651
\end{tabular}\endgroup%
{$\left.\llap{\phantom{%
\begingroup \smaller\smaller\smaller\begin{tabular}{@{}c@{}}%
0\\0\\0
\end{tabular}\endgroup%
}}\!\right]$}%

\medskip%
%
\leavevmode\llap{}%
$W_{90\phantom{0}}$%
\qquad\llap{12} lattices, $\chi=37$%
\hfill%
$222462222$%
\nopagebreak\smallskip\hrule\nopagebreak\medskip%
%
%
\leavevmode%
${L_{90.1}}$%
{} : {$1\above{1pt}{1pt}{-2}{{\rm II}}4\above{1pt}{1pt}{1}{7}{\cdot}1\above{1pt}{1pt}{2}{}3\above{1pt}{1pt}{-}{}{\cdot}1\above{1pt}{1pt}{2}{}73\above{1pt}{1pt}{1}{}$}\spacer%
\instructions{2\rightarrow N_{90}}%
\EasyButWeakLineBreak%
{${4}\above{1pt}{1pt}{*}{2}{292}\above{1pt}{1pt}{b}{2}{6}\above{1pt}{1pt}{b}{2}{4}\above{1pt}{1pt}{*}{4}{2}\above{1pt}{1pt}{}{6}{6}\above{1pt}{1pt}{s}{2}{146}\above{1pt}{1pt}{b}{2}{2}\above{1pt}{1pt}{s}{2}{438}\above{1pt}{1pt}{b}{2}$}%
\nopagebreak\par%
\nopagebreak\par\leavevmode%
{$\left[\!\llap{\phantom{%
\begingroup \smaller\smaller\smaller
\endgroup%
}}\!\right]$}%

\medskip%
%
\leavevmode\llap{}%
$W_{91\phantom{0}}$%
\qquad\llap{24} lattices, $\chi=48$%
\hfill%
$222222222222\rtimes C_{2}$%
\nopagebreak\smallskip\hrule\nopagebreak\medskip%
%
%
\leavevmode%
${L_{91.1}}$%
{} : {$1\above{1pt}{1pt}{-2}{{\rm II}}4\above{1pt}{1pt}{-}{3}{\cdot}1\above{1pt}{1pt}{2}{}3\above{1pt}{1pt}{1}{}{\cdot}1\above{1pt}{1pt}{-2}{}7\above{1pt}{1pt}{-}{}{\cdot}1\above{1pt}{1pt}{-2}{}11\above{1pt}{1pt}{-}{}$}\spacer%
\instructions{2\rightarrow N_{91}}%
\EasyButWeakLineBreak%
{${12}\above{1pt}{1pt}{r}{2}{22}\above{1pt}{1pt}{b}{2}{84}\above{1pt}{1pt}{*}{2}{4}\above{1pt}{1pt}{b}{2}{66}\above{1pt}{1pt}{b}{2}{2}\above{1pt}{1pt}{l}{2}$}\relax$\,(\times2)$%
\nopagebreak\par%
\nopagebreak\par\leavevmode%
{$\left[\!\llap{\phantom{%
\begingroup \smaller\smaller\smaller
\endgroup%
}}\!\right]$}%

\medskip%
%
\leavevmode\llap{}%
$W_{92\phantom{0}}$%
\qquad\llap{24} lattices, $\chi=15$%
\hfill%
$222242$%
\nopagebreak\smallskip\hrule\nopagebreak\medskip%
%
%
\leavevmode%
${L_{92.1}}$%
{} : {$1\above{1pt}{1pt}{-2}{{\rm II}}4\above{1pt}{1pt}{-}{3}{\cdot}1\above{1pt}{1pt}{2}{}3\above{1pt}{1pt}{1}{}{\cdot}1\above{1pt}{1pt}{2}{}7\above{1pt}{1pt}{1}{}{\cdot}1\above{1pt}{1pt}{2}{}11\above{1pt}{1pt}{1}{}$}\spacer%
\instructions{2\rightarrow N_{92}}%
\EasyButWeakLineBreak%
{${44}\above{1pt}{1pt}{r}{2}{14}\above{1pt}{1pt}{l}{2}{12}\above{1pt}{1pt}{r}{2}{154}\above{1pt}{1pt}{b}{2}{4}\above{1pt}{1pt}{*}{4}{2}\above{1pt}{1pt}{l}{2}$}%
\nopagebreak\par%
\nopagebreak\par\leavevmode%
{$\left[\!\llap{\phantom{%
\begingroup \smaller\smaller\smaller\begin{tabular}{@{}c@{}}%
0\\0\\0
\end{tabular}\endgroup%
}}\right.$}%
\begingroup \smaller\smaller\smaller\begin{tabular}{@{}c@{}}%
-3694273044\\3152688\\6305376
\end{tabular}\endgroup%
\kern3pt%
\begingroup \smaller\smaller\smaller\begin{tabular}{@{}c@{}}%
3152688\\-2690\\-5381
\end{tabular}\endgroup%
\kern3pt%
\begingroup \smaller\smaller\smaller\begin{tabular}{@{}c@{}}%
6305376\\-5381\\-10762
\end{tabular}\endgroup%
{$\left.\llap{\phantom{%
\begingroup \smaller\smaller\smaller\begin{tabular}{@{}c@{}}%
0\\0\\0
\end{tabular}\endgroup%
}}\!\right]$}%
\EasyButWeakLineBreak%
{$\left[\!\llap{\phantom{%
\begingroup \smaller\smaller\smaller\begin{tabular}{@{}c@{}}%
0\\0\\0
\end{tabular}\endgroup%
}}\right.$}%
\begingroup \smaller\smaller\smaller\begin{tabular}{@{}c@{}}%
19\\0\\11132
\end{tabular}\endgroup%
\HardButStrongLineBreak\kern3pt%
\begingroup \smaller\smaller\smaller\begin{tabular}{@{}c@{}}%
5\\-7\\2933
\end{tabular}\endgroup%
\HardButStrongLineBreak\kern3pt%
\begingroup \smaller\smaller\smaller\begin{tabular}{@{}c@{}}%
-17\\-24\\-9948
\end{tabular}\endgroup%
\HardButStrongLineBreak\kern3pt%
\begingroup \smaller\smaller\smaller\begin{tabular}{@{}c@{}}%
-101\\-77\\-59136
\end{tabular}\endgroup%
\HardButStrongLineBreak\kern3pt%
\begingroup \smaller\smaller\smaller\begin{tabular}{@{}c@{}}%
-9\\-2\\-5272
\end{tabular}\endgroup%
\HardButStrongLineBreak\kern3pt%
\begingroup \smaller\smaller\smaller\begin{tabular}{@{}c@{}}%
0\\2\\-1
\end{tabular}\endgroup%
{$\left.\llap{\phantom{%
\begingroup \smaller\smaller\smaller\begin{tabular}{@{}c@{}}%
0\\0\\0
\end{tabular}\endgroup%
}}\!\right]$}%

\medskip%
%
\leavevmode\llap{}%
$W_{93\phantom{0}}$%
\qquad\llap{24} lattices, $\chi=54$%
\hfill%
$224222224222\rtimes C_{2}$%
\nopagebreak\smallskip\hrule\nopagebreak\medskip%
%
%
\leavevmode%
${L_{93.1}}$%
{} : {$1\above{1pt}{1pt}{-2}{{\rm II}}4\above{1pt}{1pt}{-}{3}{\cdot}1\above{1pt}{1pt}{2}{}3\above{1pt}{1pt}{-}{}{\cdot}1\above{1pt}{1pt}{2}{}5\above{1pt}{1pt}{1}{}{\cdot}1\above{1pt}{1pt}{2}{}17\above{1pt}{1pt}{1}{}$}\spacer%
\instructions{2\rightarrow N_{93}}%
\EasyButWeakLineBreak%
{${20}\above{1pt}{1pt}{*}{2}{68}\above{1pt}{1pt}{b}{2}{2}\above{1pt}{1pt}{*}{4}{4}\above{1pt}{1pt}{b}{2}{34}\above{1pt}{1pt}{s}{2}{6}\above{1pt}{1pt}{b}{2}$}\relax$\,(\times2)$%
\nopagebreak\par%
\nopagebreak\par\leavevmode%
{$\left[\!\llap{\phantom{%
\begingroup \smaller\smaller\smaller
\endgroup%
}}\!\right]$}%

\medskip%
%
\leavevmode\llap{}%
$W_{94\phantom{0}}$%
\qquad\llap{24} lattices, $\chi=16$%
\hfill%
$226222$%
\nopagebreak\smallskip\hrule\nopagebreak\medskip%
%
%
\leavevmode%
${L_{94.1}}$%
{} : {$1\above{1pt}{1pt}{-2}{{\rm II}}4\above{1pt}{1pt}{-}{3}{\cdot}1\above{1pt}{1pt}{2}{}3\above{1pt}{1pt}{-}{}{\cdot}1\above{1pt}{1pt}{-2}{}5\above{1pt}{1pt}{-}{}{\cdot}1\above{1pt}{1pt}{-2}{}17\above{1pt}{1pt}{-}{}$}\spacer%
\instructions{2\rightarrow N_{94}}%
\EasyButWeakLineBreak%
{${4}\above{1pt}{1pt}{*}{2}{340}\above{1pt}{1pt}{b}{2}{6}\above{1pt}{1pt}{}{6}{2}\above{1pt}{1pt}{l}{2}{204}\above{1pt}{1pt}{r}{2}{10}\above{1pt}{1pt}{b}{2}$}%
\nopagebreak\par%
\nopagebreak\par\leavevmode%
{$\left[\!\llap{\phantom{%
\begingroup \smaller\smaller\smaller\begin{tabular}{@{}c@{}}%
0\\0\\0
\end{tabular}\endgroup%
}}\right.$}%
\begingroup \smaller\smaller\smaller\begin{tabular}{@{}c@{}}%
-16365070740\\-63936660\\1970640
\end{tabular}\endgroup%
\kern3pt%
\begingroup \smaller\smaller\smaller\begin{tabular}{@{}c@{}}%
-63936660\\-249794\\7699
\end{tabular}\endgroup%
\kern3pt%
\begingroup \smaller\smaller\smaller\begin{tabular}{@{}c@{}}%
1970640\\7699\\-226
\end{tabular}\endgroup%
{$\left.\llap{\phantom{%
\begingroup \smaller\smaller\smaller\begin{tabular}{@{}c@{}}%
0\\0\\0
\end{tabular}\endgroup%
}}\!\right]$}%
\EasyButWeakLineBreak%
{$\left[\!\llap{\phantom{%
\begingroup \smaller\smaller\smaller\begin{tabular}{@{}c@{}}%
0\\0\\0
\end{tabular}\endgroup%
}}\right.$}%
\begingroup \smaller\smaller\smaller\begin{tabular}{@{}c@{}}%
101\\-25858\\-204
\end{tabular}\endgroup%
\HardButStrongLineBreak\kern3pt%
\begingroup \smaller\smaller\smaller\begin{tabular}{@{}c@{}}%
-583\\149260\\1190
\end{tabular}\endgroup%
\HardButStrongLineBreak\kern3pt%
\begingroup \smaller\smaller\smaller\begin{tabular}{@{}c@{}}%
-200\\51204\\405
\end{tabular}\endgroup%
\HardButStrongLineBreak\kern3pt%
\begingroup \smaller\smaller\smaller\begin{tabular}{@{}c@{}}%
-252\\64517\\509
\end{tabular}\endgroup%
\HardButStrongLineBreak\kern3pt%
\begingroup \smaller\smaller\smaller\begin{tabular}{@{}c@{}}%
-4651\\1190748\\9384
\end{tabular}\endgroup%
\HardButStrongLineBreak\kern3pt%
\begingroup \smaller\smaller\smaller\begin{tabular}{@{}c@{}}%
49\\-12545\\-100
\end{tabular}\endgroup%
{$\left.\llap{\phantom{%
\begingroup \smaller\smaller\smaller\begin{tabular}{@{}c@{}}%
0\\0\\0
\end{tabular}\endgroup%
}}\!\right]$}%

\medskip%
%
\leavevmode\llap{}%
$W_{95\phantom{0}}$%
\qquad\llap{24} lattices, $\chi=36$%
\hfill%
$2222222222\rtimes C_{2}$%
\nopagebreak\smallskip\hrule\nopagebreak\medskip%
%
%
\leavevmode%
${L_{95.1}}$%
{} : {$1\above{1pt}{1pt}{-2}{{\rm II}}4\above{1pt}{1pt}{-}{5}{\cdot}1\above{1pt}{1pt}{2}{}3\above{1pt}{1pt}{-}{}{\cdot}1\above{1pt}{1pt}{2}{}7\above{1pt}{1pt}{-}{}{\cdot}1\above{1pt}{1pt}{-2}{}13\above{1pt}{1pt}{1}{}$}\spacer%
\instructions{2\rightarrow N_{95}}%
\EasyButWeakLineBreak%
{${42}\above{1pt}{1pt}{b}{2}{2}\above{1pt}{1pt}{s}{2}{182}\above{1pt}{1pt}{b}{2}{6}\above{1pt}{1pt}{l}{2}{52}\above{1pt}{1pt}{r}{2}$}\relax$\,(\times2)$%
\nopagebreak\par%
\nopagebreak\par\leavevmode%
{$\left[\!\llap{\phantom{%
\begingroup \smaller\smaller\smaller
\endgroup%
}}\!\right]$}%

\medskip%
%
\leavevmode\llap{}%
$W_{96\phantom{0}}$%
\qquad\llap{24} lattices, $\chi=56$%
\hfill%
$222262222262\rtimes C_{2}$%
\nopagebreak\smallskip\hrule\nopagebreak\medskip%
%
%
\leavevmode%
${L_{96.1}}$%
{} : {$1\above{1pt}{1pt}{-2}{{\rm II}}4\above{1pt}{1pt}{-}{5}{\cdot}1\above{1pt}{1pt}{2}{}3\above{1pt}{1pt}{-}{}{\cdot}1\above{1pt}{1pt}{-2}{}7\above{1pt}{1pt}{1}{}{\cdot}1\above{1pt}{1pt}{2}{}13\above{1pt}{1pt}{-}{}$}\spacer%
\instructions{2\rightarrow N_{96}}%
\EasyButWeakLineBreak%
{${28}\above{1pt}{1pt}{b}{2}{78}\above{1pt}{1pt}{b}{2}{14}\above{1pt}{1pt}{s}{2}{26}\above{1pt}{1pt}{b}{2}{2}\above{1pt}{1pt}{}{6}{6}\above{1pt}{1pt}{b}{2}$}\relax$\,(\times2)$%
\nopagebreak\par%
\nopagebreak\par\leavevmode%
{$\left[\!\llap{\phantom{%
\begingroup \smaller\smaller\smaller
\endgroup%
}}\!\right]$}%

\medskip%
%
\leavevmode\llap{}%
$W_{97\phantom{0}}$%
\qquad\llap{24} lattices, $\chi=60$%
\hfill%
$22222222222222\rtimes C_{2}$%
\nopagebreak\smallskip\hrule\nopagebreak\medskip%
%
%
\leavevmode%
${L_{97.1}}$%
{} : {$1\above{1pt}{1pt}{-2}{{\rm II}}4\above{1pt}{1pt}{1}{1}{\cdot}1\above{1pt}{1pt}{2}{}3\above{1pt}{1pt}{1}{}{\cdot}1\above{1pt}{1pt}{2}{}5\above{1pt}{1pt}{-}{}{\cdot}1\above{1pt}{1pt}{-2}{}19\above{1pt}{1pt}{-}{}$}\spacer%
\instructions{2\rightarrow N_{97}}%
\EasyButWeakLineBreak%
{${12}\above{1pt}{1pt}{b}{2}{38}\above{1pt}{1pt}{s}{2}{2}\above{1pt}{1pt}{l}{2}{228}\above{1pt}{1pt}{r}{2}{10}\above{1pt}{1pt}{l}{2}{4}\above{1pt}{1pt}{r}{2}{190}\above{1pt}{1pt}{b}{2}$}\relax$\,(\times2)$%
\nopagebreak\par%
\nopagebreak\par\leavevmode%
{$\left[\!\llap{\phantom{%
\begingroup \smaller\smaller\smaller
\endgroup%
}}\!\right]$}%

\medskip%
%
\leavevmode\llap{}%
$W_{98\phantom{0}}$%
\qquad\llap{24} lattices, $\chi=18$%
\hfill%
$2222222$%
\nopagebreak\smallskip\hrule\nopagebreak\medskip%
%
%
\leavevmode%
${L_{98.1}}$%
{} : {$1\above{1pt}{1pt}{-2}{{\rm II}}4\above{1pt}{1pt}{1}{1}{\cdot}1\above{1pt}{1pt}{2}{}3\above{1pt}{1pt}{1}{}{\cdot}1\above{1pt}{1pt}{-2}{}5\above{1pt}{1pt}{1}{}{\cdot}1\above{1pt}{1pt}{2}{}19\above{1pt}{1pt}{1}{}$}\spacer%
\instructions{2\rightarrow N_{98}}%
\EasyButWeakLineBreak%
{${12}\above{1pt}{1pt}{*}{2}{380}\above{1pt}{1pt}{b}{2}{2}\above{1pt}{1pt}{b}{2}{570}\above{1pt}{1pt}{l}{2}{4}\above{1pt}{1pt}{r}{2}{30}\above{1pt}{1pt}{b}{2}{76}\above{1pt}{1pt}{*}{2}$}%
\nopagebreak\par%
\nopagebreak\par\leavevmode%
{$\left[\!\llap{\phantom{%
\begingroup \smaller\smaller\smaller\begin{tabular}{@{}c@{}}%
0\\0\\0
\end{tabular}\endgroup%
}}\right.$}%
\begingroup \smaller\smaller\smaller\begin{tabular}{@{}c@{}}%
-15644220\\9120\\5700
\end{tabular}\endgroup%
\kern3pt%
\begingroup \smaller\smaller\smaller\begin{tabular}{@{}c@{}}%
9120\\-2\\-7
\end{tabular}\endgroup%
\kern3pt%
\begingroup \smaller\smaller\smaller\begin{tabular}{@{}c@{}}%
5700\\-7\\2
\end{tabular}\endgroup%
{$\left.\llap{\phantom{%
\begingroup \smaller\smaller\smaller\begin{tabular}{@{}c@{}}%
0\\0\\0
\end{tabular}\endgroup%
}}\!\right]$}%
\EasyButWeakLineBreak%
{$\left[\!\llap{\phantom{%
\begingroup \smaller\smaller\smaller\begin{tabular}{@{}c@{}}%
0\\0\\0
\end{tabular}\endgroup%
}}\right.$}%
\begingroup \smaller\smaller\smaller\begin{tabular}{@{}c@{}}%
-1\\-1098\\-990
\end{tabular}\endgroup%
\HardButStrongLineBreak\kern3pt%
\begingroup \smaller\smaller\smaller\begin{tabular}{@{}c@{}}%
-9\\-9880\\-8930
\end{tabular}\endgroup%
\HardButStrongLineBreak\kern3pt%
\begingroup \smaller\smaller\smaller\begin{tabular}{@{}c@{}}%
0\\0\\-1
\end{tabular}\endgroup%
\HardButStrongLineBreak\kern3pt%
\begingroup \smaller\smaller\smaller\begin{tabular}{@{}c@{}}%
13\\14250\\12825
\end{tabular}\endgroup%
\HardButStrongLineBreak\kern3pt%
\begingroup \smaller\smaller\smaller\begin{tabular}{@{}c@{}}%
1\\1096\\988
\end{tabular}\endgroup%
\HardButStrongLineBreak\kern3pt%
\begingroup \smaller\smaller\smaller\begin{tabular}{@{}c@{}}%
1\\1095\\990
\end{tabular}\endgroup%
\HardButStrongLineBreak\kern3pt%
\begingroup \smaller\smaller\smaller\begin{tabular}{@{}c@{}}%
-1\\-1102\\-988
\end{tabular}\endgroup%
{$\left.\llap{\phantom{%
\begingroup \smaller\smaller\smaller\begin{tabular}{@{}c@{}}%
0\\0\\0
\end{tabular}\endgroup%
}}\!\right]$}%

\medskip%
%
\leavevmode\llap{}%
$W_{99\phantom{0}}$%
\qquad\llap{12} lattices, $\chi=98$%
\hfill%
$222242262222242262\rtimes C_{2}$%
\nopagebreak\smallskip\hrule\nopagebreak\medskip%
%
%
\leavevmode%
${L_{99.1}}$%
{} : {$1\above{1pt}{1pt}{-2}{{\rm II}}4\above{1pt}{1pt}{1}{7}{\cdot}1\above{1pt}{1pt}{2}{}3\above{1pt}{1pt}{-}{}{\cdot}1\above{1pt}{1pt}{2}{}97\above{1pt}{1pt}{1}{}$}\spacer%
\instructions{2\rightarrow N_{99}}%
\EasyButWeakLineBreak%
{${4}\above{1pt}{1pt}{*}{2}{388}\above{1pt}{1pt}{b}{2}{6}\above{1pt}{1pt}{s}{2}{194}\above{1pt}{1pt}{b}{2}{2}\above{1pt}{1pt}{*}{4}{4}\above{1pt}{1pt}{b}{2}{582}\above{1pt}{1pt}{s}{2}{2}\above{1pt}{1pt}{}{6}{6}\above{1pt}{1pt}{b}{2}$}\relax$\,(\times2)$%
\nopagebreak\par%
\nopagebreak\par\leavevmode%
{$\left[\!\llap{\phantom{%
\begingroup \smaller\smaller\smaller
\endgroup%
}}\!\right]$}%

\medskip%
%
\leavevmode\llap{}%
$W_{100}$%
\qquad\llap{32} lattices, $\chi=24$%
\hfill%
$22222222$%
\nopagebreak\smallskip\hrule\nopagebreak\medskip%
%
%
\leavevmode%
${L_{100.1}}$%
{} : {$1\above{1pt}{1pt}{-2}{{\rm II}}8\above{1pt}{1pt}{1}{1}{\cdot}1\above{1pt}{1pt}{2}{}3\above{1pt}{1pt}{1}{}{\cdot}1\above{1pt}{1pt}{-2}{}5\above{1pt}{1pt}{-}{}{\cdot}1\above{1pt}{1pt}{-2}{}11\above{1pt}{1pt}{-}{}$}\spacer%
\instructions{2\rightarrow N_{100}}%
\EasyButWeakLineBreak%
{${8}\above{1pt}{1pt}{r}{2}{110}\above{1pt}{1pt}{b}{2}{2}\above{1pt}{1pt}{l}{2}{264}\above{1pt}{1pt}{r}{2}{10}\above{1pt}{1pt}{b}{2}{22}\above{1pt}{1pt}{b}{2}{40}\above{1pt}{1pt}{b}{2}{66}\above{1pt}{1pt}{l}{2}$}%
\nopagebreak\par%
\nopagebreak\par\leavevmode%
{$\left[\!\llap{\phantom{%
\begingroup \smaller\smaller\smaller\begin{tabular}{@{}c@{}}%
0\\0\\0
\end{tabular}\endgroup%
}}\right.$}%
\begingroup \smaller\smaller\smaller\begin{tabular}{@{}c@{}}%
703560\\-141240\\-1320
\end{tabular}\endgroup%
\kern3pt%
\begingroup \smaller\smaller\smaller\begin{tabular}{@{}c@{}}%
-141240\\28354\\265
\end{tabular}\endgroup%
\kern3pt%
\begingroup \smaller\smaller\smaller\begin{tabular}{@{}c@{}}%
-1320\\265\\2
\end{tabular}\endgroup%
{$\left.\llap{\phantom{%
\begingroup \smaller\smaller\smaller\begin{tabular}{@{}c@{}}%
0\\0\\0
\end{tabular}\endgroup%
}}\!\right]$}%
\EasyButWeakLineBreak%
{$\left[\!\llap{\phantom{%
\begingroup \smaller\smaller\smaller\begin{tabular}{@{}c@{}}%
0\\0\\0
\end{tabular}\endgroup%
}}\right.$}%
\begingroup \smaller\smaller\smaller\begin{tabular}{@{}c@{}}%
-299\\-1488\\-176
\end{tabular}\endgroup%
\HardButStrongLineBreak\kern3pt%
\begingroup \smaller\smaller\smaller\begin{tabular}{@{}c@{}}%
-641\\-3190\\-385
\end{tabular}\endgroup%
\HardButStrongLineBreak\kern3pt%
\begingroup \smaller\smaller\smaller\begin{tabular}{@{}c@{}}%
0\\0\\-1
\end{tabular}\endgroup%
\HardButStrongLineBreak\kern3pt%
\begingroup \smaller\smaller\smaller\begin{tabular}{@{}c@{}}%
53\\264\\0
\end{tabular}\endgroup%
\HardButStrongLineBreak\kern3pt%
\begingroup \smaller\smaller\smaller\begin{tabular}{@{}c@{}}%
1\\5\\0
\end{tabular}\endgroup%
\HardButStrongLineBreak\kern3pt%
\begingroup \smaller\smaller\smaller\begin{tabular}{@{}c@{}}%
-42\\-209\\-22
\end{tabular}\endgroup%
\HardButStrongLineBreak\kern3pt%
\begingroup \smaller\smaller\smaller\begin{tabular}{@{}c@{}}%
-213\\-1060\\-120
\end{tabular}\endgroup%
\HardButStrongLineBreak\kern3pt%
\begingroup \smaller\smaller\smaller\begin{tabular}{@{}c@{}}%
-683\\-3399\\-396
\end{tabular}\endgroup%
{$\left.\llap{\phantom{%
\begingroup \smaller\smaller\smaller\begin{tabular}{@{}c@{}}%
0\\0\\0
\end{tabular}\endgroup%
}}\!\right]$}%

\medskip%
%
\leavevmode\llap{}%
$W_{101}$%
\qquad\llap{24} lattices, $\chi=22$%
\hfill%
$2262222$%
\nopagebreak\smallskip\hrule\nopagebreak\medskip%
%
%
\leavevmode%
${L_{101.1}}$%
{} : {$1\above{1pt}{1pt}{-2}{{\rm II}}4\above{1pt}{1pt}{-}{5}{\cdot}1\above{1pt}{1pt}{2}{}3\above{1pt}{1pt}{-}{}{\cdot}1\above{1pt}{1pt}{-2}{}5\above{1pt}{1pt}{-}{}{\cdot}1\above{1pt}{1pt}{2}{}23\above{1pt}{1pt}{1}{}$}\spacer%
\instructions{2\rightarrow N_{101}}%
\EasyButWeakLineBreak%
{${60}\above{1pt}{1pt}{*}{2}{92}\above{1pt}{1pt}{b}{2}{6}\above{1pt}{1pt}{}{6}{2}\above{1pt}{1pt}{l}{2}{276}\above{1pt}{1pt}{r}{2}{10}\above{1pt}{1pt}{s}{2}{46}\above{1pt}{1pt}{b}{2}$}%
\nopagebreak\par%
\nopagebreak\par\leavevmode%
{$\left[\!\llap{\phantom{%
\begingroup \smaller\smaller\smaller\begin{tabular}{@{}c@{}}%
0\\0\\0
\end{tabular}\endgroup%
}}\right.$}%
\begingroup \smaller\smaller\smaller\begin{tabular}{@{}c@{}}%
-579195660\\194669700\\-401580
\end{tabular}\endgroup%
\kern3pt%
\begingroup \smaller\smaller\smaller\begin{tabular}{@{}c@{}}%
194669700\\-65428814\\134783
\end{tabular}\endgroup%
\kern3pt%
\begingroup \smaller\smaller\smaller\begin{tabular}{@{}c@{}}%
-401580\\134783\\-178
\end{tabular}\endgroup%
{$\left.\llap{\phantom{%
\begingroup \smaller\smaller\smaller\begin{tabular}{@{}c@{}}%
0\\0\\0
\end{tabular}\endgroup%
}}\!\right]$}%
\EasyButWeakLineBreak%
{$\left[\!\llap{\phantom{%
\begingroup \smaller\smaller\smaller\begin{tabular}{@{}c@{}}%
0\\0\\0
\end{tabular}\endgroup%
}}\right.$}%
\begingroup \smaller\smaller\smaller\begin{tabular}{@{}c@{}}%
-19003\\-56760\\-107070
\end{tabular}\endgroup%
\HardButStrongLineBreak\kern3pt%
\begingroup \smaller\smaller\smaller\begin{tabular}{@{}c@{}}%
10303\\30774\\58052
\end{tabular}\endgroup%
\HardButStrongLineBreak\kern3pt%
\begingroup \smaller\smaller\smaller\begin{tabular}{@{}c@{}}%
4808\\14361\\27090
\end{tabular}\endgroup%
\HardButStrongLineBreak\kern3pt%
\begingroup \smaller\smaller\smaller\begin{tabular}{@{}c@{}}%
-5037\\-15045\\-28381
\end{tabular}\endgroup%
\HardButStrongLineBreak\kern3pt%
\begingroup \smaller\smaller\smaller\begin{tabular}{@{}c@{}}%
-612451\\-1829328\\-3450828
\end{tabular}\endgroup%
\HardButStrongLineBreak\kern3pt%
\begingroup \smaller\smaller\smaller\begin{tabular}{@{}c@{}}%
-45867\\-137000\\-258435
\end{tabular}\endgroup%
\HardButStrongLineBreak\kern3pt%
\begingroup \smaller\smaller\smaller\begin{tabular}{@{}c@{}}%
-61741\\-184414\\-347875
\end{tabular}\endgroup%
{$\left.\llap{\phantom{%
\begingroup \smaller\smaller\smaller\begin{tabular}{@{}c@{}}%
0\\0\\0
\end{tabular}\endgroup%
}}\!\right]$}%

\medskip%
%
\leavevmode\llap{}%
$W_{102}$%
\qquad\llap{24} lattices, $\chi=72$%
\hfill%
$2222222222222222\rtimes C_{2}$%
\nopagebreak\smallskip\hrule\nopagebreak\medskip%
%
%
\leavevmode%
${L_{102.1}}$%
{} : {$1\above{1pt}{1pt}{-2}{{\rm II}}4\above{1pt}{1pt}{1}{1}{\cdot}1\above{1pt}{1pt}{2}{}3\above{1pt}{1pt}{1}{}{\cdot}1\above{1pt}{1pt}{-2}{}7\above{1pt}{1pt}{1}{}{\cdot}1\above{1pt}{1pt}{2}{}17\above{1pt}{1pt}{1}{}$}\spacer%
\instructions{2\rightarrow N_{102}}%
\EasyButWeakLineBreak%
{${2}\above{1pt}{1pt}{l}{2}{68}\above{1pt}{1pt}{r}{2}{14}\above{1pt}{1pt}{b}{2}{12}\above{1pt}{1pt}{*}{2}{28}\above{1pt}{1pt}{b}{2}{34}\above{1pt}{1pt}{l}{2}{4}\above{1pt}{1pt}{r}{2}{714}\above{1pt}{1pt}{b}{2}$}\relax$\,(\times2)$%
\nopagebreak\par%
\nopagebreak\par\leavevmode%
{$\left[\!\llap{\phantom{%
\begingroup \smaller\smaller\smaller
\endgroup%
}}\!\right]$}%

\medskip%
%
\leavevmode\llap{}%
$W_{103}$%
\qquad\llap{24} lattices, $\chi=48$%
\hfill%
$222222222222\rtimes C_{2}$%
\nopagebreak\smallskip\hrule\nopagebreak\medskip%
%
%
\leavevmode%
${L_{103.1}}$%
{} : {$1\above{1pt}{1pt}{-2}{{\rm II}}4\above{1pt}{1pt}{1}{1}{\cdot}1\above{1pt}{1pt}{2}{}3\above{1pt}{1pt}{1}{}{\cdot}1\above{1pt}{1pt}{2}{}7\above{1pt}{1pt}{-}{}{\cdot}1\above{1pt}{1pt}{-2}{}17\above{1pt}{1pt}{-}{}$}\spacer%
\instructions{2\rightarrow N_{103}}%
\EasyButWeakLineBreak%
{${12}\above{1pt}{1pt}{*}{2}{476}\above{1pt}{1pt}{b}{2}{2}\above{1pt}{1pt}{s}{2}{102}\above{1pt}{1pt}{l}{2}{4}\above{1pt}{1pt}{r}{2}{238}\above{1pt}{1pt}{b}{2}$}\relax$\,(\times2)$%
\nopagebreak\par%
\nopagebreak\par\leavevmode%
{$\left[\!\llap{\phantom{%
\begingroup \smaller\smaller\smaller
\endgroup%
}}\!\right]$}%

\medskip%
%
\leavevmode\llap{}%
$W_{104}$%
\qquad\llap{24} lattices, $\chi=60$%
\hfill%
$22222222222222\rtimes C_{2}$%
\nopagebreak\smallskip\hrule\nopagebreak\medskip%
%
%
\leavevmode%
${L_{104.1}}$%
{} : {$1\above{1pt}{1pt}{-2}{{\rm II}}4\above{1pt}{1pt}{-}{5}{\cdot}1\above{1pt}{1pt}{-2}{}5\above{1pt}{1pt}{1}{}{\cdot}1\above{1pt}{1pt}{2}{}7\above{1pt}{1pt}{1}{}{\cdot}1\above{1pt}{1pt}{2}{}11\above{1pt}{1pt}{1}{}$}\spacer%
\instructions{2\rightarrow N_{104}}%
\EasyButWeakLineBreak%
{${28}\above{1pt}{1pt}{*}{2}{44}\above{1pt}{1pt}{b}{2}{2}\above{1pt}{1pt}{b}{2}{154}\above{1pt}{1pt}{l}{2}{20}\above{1pt}{1pt}{r}{2}{14}\above{1pt}{1pt}{b}{2}{220}\above{1pt}{1pt}{*}{2}$}\relax$\,(\times2)$%
\nopagebreak\par%
\nopagebreak\par\leavevmode%
{$\left[\!\llap{\phantom{%
\begingroup \smaller\smaller\smaller
\endgroup%
}}\!\right]$}%

\medskip%
%
\leavevmode\llap{}%
$W_{105}$%
\qquad\llap{32} lattices, $\chi=72$%
\hfill%
$2222222222222222\rtimes C_{2}$%
\nopagebreak\smallskip\hrule\nopagebreak\medskip%
%
%
\leavevmode%
${L_{105.1}}$%
{} : {$1\above{1pt}{1pt}{-2}{{\rm II}}8\above{1pt}{1pt}{1}{7}{\cdot}1\above{1pt}{1pt}{2}{}3\above{1pt}{1pt}{-}{}{\cdot}1\above{1pt}{1pt}{2}{}5\above{1pt}{1pt}{-}{}{\cdot}1\above{1pt}{1pt}{-2}{}13\above{1pt}{1pt}{-}{}$}\spacer%
\instructions{2\rightarrow N_{105}}%
\EasyButWeakLineBreak%
{${24}\above{1pt}{1pt}{b}{2}{26}\above{1pt}{1pt}{b}{2}{6}\above{1pt}{1pt}{b}{2}{10}\above{1pt}{1pt}{s}{2}{78}\above{1pt}{1pt}{b}{2}{2}\above{1pt}{1pt}{l}{2}{312}\above{1pt}{1pt}{r}{2}{10}\above{1pt}{1pt}{b}{2}$}\relax$\,(\times2)$%
\nopagebreak\par%
\nopagebreak\par\leavevmode%
{$\left[\!\llap{\phantom{%
\begingroup \smaller\smaller\smaller
\endgroup%
}}\!\right]$}%

\medskip%
%
\leavevmode\llap{}%
$W_{106}$%
\qquad\llap{24} lattices, $\chi=27$%
\hfill%
$22242222$%
\nopagebreak\smallskip\hrule\nopagebreak\medskip%
%
%
\leavevmode%
${L_{106.1}}$%
{} : {$1\above{1pt}{1pt}{-2}{{\rm II}}4\above{1pt}{1pt}{-}{3}{\cdot}1\above{1pt}{1pt}{2}{}3\above{1pt}{1pt}{-}{}{\cdot}1\above{1pt}{1pt}{2}{}7\above{1pt}{1pt}{-}{}{\cdot}1\above{1pt}{1pt}{2}{}19\above{1pt}{1pt}{1}{}$}\spacer%
\instructions{2\rightarrow N_{106}}%
\EasyButWeakLineBreak%
{${6}\above{1pt}{1pt}{l}{2}{76}\above{1pt}{1pt}{r}{2}{42}\above{1pt}{1pt}{b}{2}{4}\above{1pt}{1pt}{*}{4}{2}\above{1pt}{1pt}{b}{2}{114}\above{1pt}{1pt}{b}{2}{4}\above{1pt}{1pt}{*}{2}{532}\above{1pt}{1pt}{b}{2}$}%
\nopagebreak\par%
\nopagebreak\par\leavevmode%
{$\left[\!\llap{\phantom{%
\begingroup \smaller\smaller\smaller\begin{tabular}{@{}c@{}}%
0\\0\\0
\end{tabular}\endgroup%
}}\right.$}%
\begingroup \smaller\smaller\smaller\begin{tabular}{@{}c@{}}%
-165273780\\78204\\111720
\end{tabular}\endgroup%
\kern3pt%
\begingroup \smaller\smaller\smaller\begin{tabular}{@{}c@{}}%
78204\\-34\\-55
\end{tabular}\endgroup%
\kern3pt%
\begingroup \smaller\smaller\smaller\begin{tabular}{@{}c@{}}%
111720\\-55\\-74
\end{tabular}\endgroup%
{$\left.\llap{\phantom{%
\begingroup \smaller\smaller\smaller\begin{tabular}{@{}c@{}}%
0\\0\\0
\end{tabular}\endgroup%
}}\!\right]$}%
\EasyButWeakLineBreak%
{$\left[\!\llap{\phantom{%
\begingroup \smaller\smaller\smaller\begin{tabular}{@{}c@{}}%
0\\0\\0
\end{tabular}\endgroup%
}}\right.$}%
\begingroup \smaller\smaller\smaller\begin{tabular}{@{}c@{}}%
4\\2808\\3951
\end{tabular}\endgroup%
\HardButStrongLineBreak\kern3pt%
\begingroup \smaller\smaller\smaller\begin{tabular}{@{}c@{}}%
13\\9120\\12844
\end{tabular}\endgroup%
\HardButStrongLineBreak\kern3pt%
\begingroup \smaller\smaller\smaller\begin{tabular}{@{}c@{}}%
-2\\-1407\\-1974
\end{tabular}\endgroup%
\HardButStrongLineBreak\kern3pt%
\begingroup \smaller\smaller\smaller\begin{tabular}{@{}c@{}}%
-1\\-702\\-988
\end{tabular}\endgroup%
\HardButStrongLineBreak\kern3pt%
\begingroup \smaller\smaller\smaller\begin{tabular}{@{}c@{}}%
1\\703\\987
\end{tabular}\endgroup%
\HardButStrongLineBreak\kern3pt%
\begingroup \smaller\smaller\smaller\begin{tabular}{@{}c@{}}%
34\\23883\\33573
\end{tabular}\endgroup%
\HardButStrongLineBreak\kern3pt%
\begingroup \smaller\smaller\smaller\begin{tabular}{@{}c@{}}%
13\\9130\\12838
\end{tabular}\endgroup%
\HardButStrongLineBreak\kern3pt%
\begingroup \smaller\smaller\smaller\begin{tabular}{@{}c@{}}%
275\\193116\\271586
\end{tabular}\endgroup%
{$\left.\llap{\phantom{%
\begingroup \smaller\smaller\smaller\begin{tabular}{@{}c@{}}%
0\\0\\0
\end{tabular}\endgroup%
}}\!\right]$}%

\medskip%
%
\leavevmode\llap{}%
$W_{107}$%
\qquad\llap{24} lattices, $\chi=60$%
\hfill%
$22222222222222\rtimes C_{2}$%
\nopagebreak\smallskip\hrule\nopagebreak\medskip%
%
%
\leavevmode%
${L_{107.1}}$%
{} : {$1\above{1pt}{1pt}{-2}{{\rm II}}4\above{1pt}{1pt}{1}{1}{\cdot}1\above{1pt}{1pt}{2}{}3\above{1pt}{1pt}{1}{}{\cdot}1\above{1pt}{1pt}{2}{}11\above{1pt}{1pt}{1}{}{\cdot}1\above{1pt}{1pt}{-2}{}13\above{1pt}{1pt}{1}{}$}\spacer%
\instructions{2\rightarrow N_{107}}%
\EasyButWeakLineBreak%
{${44}\above{1pt}{1pt}{*}{2}{156}\above{1pt}{1pt}{b}{2}{2}\above{1pt}{1pt}{b}{2}{858}\above{1pt}{1pt}{l}{2}{4}\above{1pt}{1pt}{r}{2}{286}\above{1pt}{1pt}{b}{2}{12}\above{1pt}{1pt}{*}{2}$}\relax$\,(\times2)$%
\nopagebreak\par%
\nopagebreak\par\leavevmode%
{$\left[\!\llap{\phantom{%
\begingroup \smaller\smaller\smaller
\endgroup%
}}\!\right]$}%

\medskip%
%
\leavevmode\llap{}%
$W_{108}$%
\qquad\llap{24} lattices, $\chi=28$%
\hfill%
$22222262$%
\nopagebreak\smallskip\hrule\nopagebreak\medskip%
%
%
\leavevmode%
${L_{108.1}}$%
{} : {$1\above{1pt}{1pt}{-2}{{\rm II}}4\above{1pt}{1pt}{1}{7}{\cdot}1\above{1pt}{1pt}{2}{}3\above{1pt}{1pt}{-}{}{\cdot}1\above{1pt}{1pt}{-2}{}5\above{1pt}{1pt}{1}{}{\cdot}1\above{1pt}{1pt}{-2}{}29\above{1pt}{1pt}{1}{}$}\spacer%
\instructions{2\rightarrow N_{108}}%
\EasyButWeakLineBreak%
{${20}\above{1pt}{1pt}{*}{2}{116}\above{1pt}{1pt}{b}{2}{2}\above{1pt}{1pt}{s}{2}{870}\above{1pt}{1pt}{b}{2}{4}\above{1pt}{1pt}{*}{2}{580}\above{1pt}{1pt}{b}{2}{6}\above{1pt}{1pt}{}{6}{2}\above{1pt}{1pt}{b}{2}$}%
\nopagebreak\par%
\nopagebreak\par\leavevmode%
{$\left[\!\llap{\phantom{%
\begingroup \smaller\smaller\smaller\begin{tabular}{@{}c@{}}%
0\\0\\0
\end{tabular}\endgroup%
}}\right.$}%
\begingroup \smaller\smaller\smaller\begin{tabular}{@{}c@{}}%
-964563780\\5310480\\-210540
\end{tabular}\endgroup%
\kern3pt%
\begingroup \smaller\smaller\smaller\begin{tabular}{@{}c@{}}%
5310480\\-29014\\1117
\end{tabular}\endgroup%
\kern3pt%
\begingroup \smaller\smaller\smaller\begin{tabular}{@{}c@{}}%
-210540\\1117\\-38
\end{tabular}\endgroup%
{$\left.\llap{\phantom{%
\begingroup \smaller\smaller\smaller\begin{tabular}{@{}c@{}}%
0\\0\\0
\end{tabular}\endgroup%
}}\!\right]$}%
\EasyButWeakLineBreak%
{$\left[\!\llap{\phantom{%
\begingroup \smaller\smaller\smaller\begin{tabular}{@{}c@{}}%
0\\0\\0
\end{tabular}\endgroup%
}}\right.$}%
\begingroup \smaller\smaller\smaller\begin{tabular}{@{}c@{}}%
1\\230\\1220
\end{tabular}\endgroup%
\HardButStrongLineBreak\kern3pt%
\begingroup \smaller\smaller\smaller\begin{tabular}{@{}c@{}}%
-335\\-77024\\-408030
\end{tabular}\endgroup%
\HardButStrongLineBreak\kern3pt%
\begingroup \smaller\smaller\smaller\begin{tabular}{@{}c@{}}%
-52\\-11956\\-63337
\end{tabular}\endgroup%
\HardButStrongLineBreak\kern3pt%
\begingroup \smaller\smaller\smaller\begin{tabular}{@{}c@{}}%
-2486\\-571590\\-3028035
\end{tabular}\endgroup%
\HardButStrongLineBreak\kern3pt%
\begingroup \smaller\smaller\smaller\begin{tabular}{@{}c@{}}%
-79\\-18164\\-96226
\end{tabular}\endgroup%
\HardButStrongLineBreak\kern3pt%
\begingroup \smaller\smaller\smaller\begin{tabular}{@{}c@{}}%
-1593\\-366270\\-1940390
\end{tabular}\endgroup%
\HardButStrongLineBreak\kern3pt%
\begingroup \smaller\smaller\smaller\begin{tabular}{@{}c@{}}%
-14\\-3219\\-17055
\end{tabular}\endgroup%
\HardButStrongLineBreak\kern3pt%
\begingroup \smaller\smaller\smaller\begin{tabular}{@{}c@{}}%
13\\2989\\15834
\end{tabular}\endgroup%
{$\left.\llap{\phantom{%
\begingroup \smaller\smaller\smaller\begin{tabular}{@{}c@{}}%
0\\0\\0
\end{tabular}\endgroup%
}}\!\right]$}%

\medskip%
%
\leavevmode\llap{}%
$W_{109}$%
\qquad\llap{24} lattices, $\chi=90$%
\hfill%
$222422222222422222\rtimes C_{2}$%
\nopagebreak\smallskip\hrule\nopagebreak\medskip%
%
%
\leavevmode%
${L_{109.1}}$%
{} : {$1\above{1pt}{1pt}{-2}{{\rm II}}4\above{1pt}{1pt}{1}{7}{\cdot}1\above{1pt}{1pt}{2}{}3\above{1pt}{1pt}{-}{}{\cdot}1\above{1pt}{1pt}{2}{}5\above{1pt}{1pt}{-}{}{\cdot}1\above{1pt}{1pt}{2}{}29\above{1pt}{1pt}{-}{}$}\spacer%
\instructions{2\rightarrow N_{109}}%
\EasyButWeakLineBreak%
{${10}\above{1pt}{1pt}{b}{2}{58}\above{1pt}{1pt}{l}{2}{60}\above{1pt}{1pt}{r}{2}{2}\above{1pt}{1pt}{*}{4}{4}\above{1pt}{1pt}{b}{2}{6}\above{1pt}{1pt}{s}{2}{290}\above{1pt}{1pt}{b}{2}{2}\above{1pt}{1pt}{l}{2}{348}\above{1pt}{1pt}{r}{2}$}\relax$\,(\times2)$%
\nopagebreak\par%
\nopagebreak\par\leavevmode%
{$\left[\!\llap{\phantom{%
\begingroup \smaller\smaller\smaller
\endgroup%
}}\!\right]$}%

\medskip%
%
\leavevmode\llap{}%
$W_{110}$%
\qquad\llap{24} lattices, $\chi=72$%
\hfill%
$2222222222222222\rtimes C_{2}$%
\nopagebreak\smallskip\hrule\nopagebreak\medskip%
%
%
\leavevmode%
${L_{110.1}}$%
{} : {$1\above{1pt}{1pt}{-2}{{\rm II}}4\above{1pt}{1pt}{-}{3}{\cdot}1\above{1pt}{1pt}{-2}{}5\above{1pt}{1pt}{-}{}{\cdot}1\above{1pt}{1pt}{2}{}7\above{1pt}{1pt}{-}{}{\cdot}1\above{1pt}{1pt}{-2}{}13\above{1pt}{1pt}{-}{}$}\spacer%
\instructions{2\rightarrow N_{110}}%
\EasyButWeakLineBreak%
{${26}\above{1pt}{1pt}{l}{2}{140}\above{1pt}{1pt}{r}{2}{2}\above{1pt}{1pt}{b}{2}{260}\above{1pt}{1pt}{*}{2}{4}\above{1pt}{1pt}{b}{2}{2}\above{1pt}{1pt}{l}{2}{364}\above{1pt}{1pt}{r}{2}{10}\above{1pt}{1pt}{b}{2}$}\relax$\,(\times2)$%
\nopagebreak\par%
\nopagebreak\par\leavevmode%
{$\left[\!\llap{\phantom{%
\begingroup \smaller\smaller\smaller
\endgroup%
}}\!\right]$}%

\medskip%
%
\leavevmode\llap{}%
$W_{111}$%
\qquad\llap{24} lattices, $\chi=60$%
\hfill%
$22222222222222\rtimes C_{2}$%
\nopagebreak\smallskip\hrule\nopagebreak\medskip%
%
%
\leavevmode%
${L_{111.1}}$%
{} : {$1\above{1pt}{1pt}{-2}{{\rm II}}4\above{1pt}{1pt}{-}{5}{\cdot}1\above{1pt}{1pt}{2}{}3\above{1pt}{1pt}{1}{}{\cdot}1\above{1pt}{1pt}{-2}{}5\above{1pt}{1pt}{1}{}{\cdot}1\above{1pt}{1pt}{2}{}31\above{1pt}{1pt}{1}{}$}\spacer%
\instructions{2\rightarrow N_{111}}%
\EasyButWeakLineBreak%
{${12}\above{1pt}{1pt}{*}{2}{124}\above{1pt}{1pt}{b}{2}{30}\above{1pt}{1pt}{b}{2}{62}\above{1pt}{1pt}{l}{2}{20}\above{1pt}{1pt}{r}{2}{2}\above{1pt}{1pt}{b}{2}{620}\above{1pt}{1pt}{*}{2}$}\relax$\,(\times2)$%
\nopagebreak\par%
\nopagebreak\par\leavevmode%
{$\left[\!\llap{\phantom{%
\begingroup \smaller\smaller\smaller
\endgroup%
}}\!\right]$}%

\medskip%
%
\leavevmode\llap{}%
$W_{112}$%
\qquad\llap{24} lattices, $\chi=66$%
\hfill%
$22222242222224\rtimes C_{2}$%
\nopagebreak\smallskip\hrule\nopagebreak\medskip%
%
%
\leavevmode%
${L_{112.1}}$%
{} : {$1\above{1pt}{1pt}{-2}{{\rm II}}4\above{1pt}{1pt}{1}{7}{\cdot}1\above{1pt}{1pt}{2}{}3\above{1pt}{1pt}{1}{}{\cdot}1\above{1pt}{1pt}{2}{}7\above{1pt}{1pt}{1}{}{\cdot}1\above{1pt}{1pt}{2}{}23\above{1pt}{1pt}{1}{}$}\spacer%
\instructions{2\rightarrow N_{112}}%
\EasyButWeakLineBreak%
{${2}\above{1pt}{1pt}{l}{2}{92}\above{1pt}{1pt}{r}{2}{14}\above{1pt}{1pt}{s}{2}{138}\above{1pt}{1pt}{l}{2}{28}\above{1pt}{1pt}{r}{2}{46}\above{1pt}{1pt}{b}{2}{4}\above{1pt}{1pt}{*}{4}$}\relax$\,(\times2)$%
\nopagebreak\par%
\nopagebreak\par\leavevmode%
{$\left[\!\llap{\phantom{%
\begingroup \smaller\smaller\smaller
\endgroup%
}}\!\right]$}%

\medskip%
%
\leavevmode\llap{}%
$W_{113}$%
\qquad\llap{32} lattices, $\chi=80$%
\hfill%
$2226222222262222\rtimes C_{2}$%
\nopagebreak\smallskip\hrule\nopagebreak\medskip%
%
%
\leavevmode%
${L_{113.1}}$%
{} : {$1\above{1pt}{1pt}{-2}{{\rm II}}8\above{1pt}{1pt}{1}{1}{\cdot}1\above{1pt}{1pt}{2}{}3\above{1pt}{1pt}{-}{}{\cdot}1\above{1pt}{1pt}{-2}{}5\above{1pt}{1pt}{-}{}{\cdot}1\above{1pt}{1pt}{-2}{}19\above{1pt}{1pt}{1}{}$}\spacer%
\instructions{2\rightarrow N_{113}}%
\EasyButWeakLineBreak%
{${6}\above{1pt}{1pt}{b}{2}{10}\above{1pt}{1pt}{l}{2}{456}\above{1pt}{1pt}{r}{2}{2}\above{1pt}{1pt}{}{6}{6}\above{1pt}{1pt}{b}{2}{40}\above{1pt}{1pt}{b}{2}{114}\above{1pt}{1pt}{l}{2}{8}\above{1pt}{1pt}{r}{2}$}\relax$\,(\times2)$%
\nopagebreak\par%
\nopagebreak\par\leavevmode%
{$\left[\!\llap{\phantom{%
\begingroup \smaller\smaller\smaller
\endgroup%
}}\!\right]$}%

\medskip%
%
\leavevmode\llap{}%
$W_{114}$%
\qquad\llap{24} lattices, $\chi=80$%
\hfill%
$2262222222622222\rtimes C_{2}$%
\nopagebreak\smallskip\hrule\nopagebreak\medskip%
%
%
\leavevmode%
${L_{114.1}}$%
{} : {$1\above{1pt}{1pt}{-2}{{\rm II}}4\above{1pt}{1pt}{-}{3}{\cdot}1\above{1pt}{1pt}{2}{}3\above{1pt}{1pt}{-}{}{\cdot}1\above{1pt}{1pt}{-2}{}5\above{1pt}{1pt}{1}{}{\cdot}1\above{1pt}{1pt}{-2}{}41\above{1pt}{1pt}{1}{}$}\spacer%
\instructions{2\rightarrow N_{114}}%
\EasyButWeakLineBreak%
{${20}\above{1pt}{1pt}{*}{2}{164}\above{1pt}{1pt}{b}{2}{2}\above{1pt}{1pt}{}{6}{6}\above{1pt}{1pt}{s}{2}{82}\above{1pt}{1pt}{b}{2}{4}\above{1pt}{1pt}{*}{2}{820}\above{1pt}{1pt}{b}{2}{6}\above{1pt}{1pt}{b}{2}$}\relax$\,(\times2)$%
\nopagebreak\par%
\nopagebreak\par\leavevmode%
{$\left[\!\llap{\phantom{%
\begingroup \smaller\smaller\smaller
\endgroup%
}}\!\right]$}%

\medskip%
%
\leavevmode\llap{}%
$W_{115}$%
\qquad\llap{24} lattices, $\chi=84$%
\hfill%
$222222222222222222\rtimes C_{2}$%
\nopagebreak\smallskip\hrule\nopagebreak\medskip%
%
%
\leavevmode%
${L_{115.1}}$%
{} : {$1\above{1pt}{1pt}{-2}{{\rm II}}4\above{1pt}{1pt}{1}{1}{\cdot}1\above{1pt}{1pt}{2}{}3\above{1pt}{1pt}{1}{}{\cdot}1\above{1pt}{1pt}{-2}{}5\above{1pt}{1pt}{-}{}{\cdot}1\above{1pt}{1pt}{2}{}43\above{1pt}{1pt}{-}{}$}\spacer%
\instructions{2\rightarrow N_{115}}%
\EasyButWeakLineBreak%
{${12}\above{1pt}{1pt}{b}{2}{86}\above{1pt}{1pt}{s}{2}{2}\above{1pt}{1pt}{l}{2}{516}\above{1pt}{1pt}{r}{2}{10}\above{1pt}{1pt}{l}{2}{4}\above{1pt}{1pt}{r}{2}{1290}\above{1pt}{1pt}{b}{2}{2}\above{1pt}{1pt}{b}{2}{860}\above{1pt}{1pt}{*}{2}$}\relax$\,(\times2)$%
\nopagebreak\par%
\nopagebreak\par\leavevmode%
{$\left[\!\llap{\phantom{%
\begingroup \smaller\smaller\smaller
\endgroup%
}}\!\right]$}%

\medskip%
%
\leavevmode\llap{}%
$W_{116}$%
\qquad\llap{24} lattices, $\chi=90$%
\hfill%
$224222222224222222\rtimes C_{2}$%
\nopagebreak\smallskip\hrule\nopagebreak\medskip%
%
%
\leavevmode%
${L_{116.1}}$%
{} : {$1\above{1pt}{1pt}{-2}{{\rm II}}4\above{1pt}{1pt}{1}{7}{\cdot}1\above{1pt}{1pt}{2}{}3\above{1pt}{1pt}{-}{}{\cdot}1\above{1pt}{1pt}{2}{}7\above{1pt}{1pt}{-}{}{\cdot}1\above{1pt}{1pt}{2}{}31\above{1pt}{1pt}{1}{}$}\spacer%
\instructions{2\rightarrow N_{116}}%
\EasyButWeakLineBreak%
{${6}\above{1pt}{1pt}{s}{2}{434}\above{1pt}{1pt}{b}{2}{2}\above{1pt}{1pt}{*}{4}{4}\above{1pt}{1pt}{*}{2}{868}\above{1pt}{1pt}{b}{2}{6}\above{1pt}{1pt}{l}{2}{124}\above{1pt}{1pt}{r}{2}{42}\above{1pt}{1pt}{s}{2}{62}\above{1pt}{1pt}{b}{2}$}\relax$\,(\times2)$%
\nopagebreak\par%
\nopagebreak\par\leavevmode%
{$\left[\!\llap{\phantom{%
\begingroup \smaller\smaller\smaller
\endgroup%
}}\!\right]$}%

\medskip%
%
\leavevmode\llap{}%
$W_{117}$%
\qquad\llap{24} lattices, $\chi=104$%
\hfill%
$62222222226222222222\rtimes C_{2}$%
\nopagebreak\smallskip\hrule\nopagebreak\medskip%
%
%
\leavevmode%
${L_{117.1}}$%
{} : {$1\above{1pt}{1pt}{-2}{{\rm II}}4\above{1pt}{1pt}{1}{7}{\cdot}1\above{1pt}{1pt}{2}{}3\above{1pt}{1pt}{-}{}{\cdot}1\above{1pt}{1pt}{-2}{}5\above{1pt}{1pt}{-}{}{\cdot}1\above{1pt}{1pt}{-2}{}53\above{1pt}{1pt}{-}{}$}\spacer%
\instructions{2\rightarrow N_{117}}%
\EasyButWeakLineBreak%
{${2}\above{1pt}{1pt}{}{6}{6}\above{1pt}{1pt}{b}{2}{1060}\above{1pt}{1pt}{*}{2}{4}\above{1pt}{1pt}{b}{2}{1590}\above{1pt}{1pt}{s}{2}{2}\above{1pt}{1pt}{l}{2}{636}\above{1pt}{1pt}{r}{2}{10}\above{1pt}{1pt}{b}{2}{106}\above{1pt}{1pt}{l}{2}{60}\above{1pt}{1pt}{r}{2}$}\relax$\,(\times2)$%
\nopagebreak\par%
\nopagebreak\par\leavevmode%
{$\left[\!\llap{\phantom{%
\begingroup \smaller\smaller\smaller
\endgroup%
}}\!\right]$}%

\medskip%
%
\leavevmode\llap{}%
$W_{118}$%
\qquad\llap{48} lattices, $\chi=36$%
\hfill%
$2222222222$%
\nopagebreak\smallskip\hrule\nopagebreak\medskip%
%
%
\leavevmode%
${L_{118.1}}$%
{} : {$1\above{1pt}{1pt}{-2}{{\rm II}}4\above{1pt}{1pt}{1}{7}{\cdot}1\above{1pt}{1pt}{2}{}3\above{1pt}{1pt}{-}{}{\cdot}1\above{1pt}{1pt}{-2}{}5\above{1pt}{1pt}{-}{}{\cdot}1\above{1pt}{1pt}{2}{}7\above{1pt}{1pt}{-}{}{\cdot}1\above{1pt}{1pt}{-2}{}11\above{1pt}{1pt}{-}{}$}\spacer%
\instructions{2\rightarrow N_{118}}%
\EasyButWeakLineBreak%
{${60}\above{1pt}{1pt}{r}{2}{22}\above{1pt}{1pt}{b}{2}{4}\above{1pt}{1pt}{*}{2}{1540}\above{1pt}{1pt}{b}{2}{6}\above{1pt}{1pt}{b}{2}{110}\above{1pt}{1pt}{s}{2}{42}\above{1pt}{1pt}{b}{2}{10}\above{1pt}{1pt}{l}{2}{924}\above{1pt}{1pt}{r}{2}{2}\above{1pt}{1pt}{l}{2}$}%
\nopagebreak\par%
\nopagebreak\par\leavevmode%
{$\left[\!\llap{\phantom{%
\begingroup \smaller\smaller\smaller
\endgroup%
}}\!\right]$}%

\medskip%
%
\leavevmode\llap{}%
$W_{119}$%
\qquad\llap{48} lattices, $\chi=90$%
\hfill%
$422222222422222222\rtimes C_{2}$%
\nopagebreak\smallskip\hrule\nopagebreak\medskip%
%
%
\leavevmode%
${L_{119.1}}$%
{} : {$1\above{1pt}{1pt}{-2}{{\rm II}}4\above{1pt}{1pt}{1}{7}{\cdot}1\above{1pt}{1pt}{2}{}3\above{1pt}{1pt}{-}{}{\cdot}1\above{1pt}{1pt}{2}{}5\above{1pt}{1pt}{1}{}{\cdot}1\above{1pt}{1pt}{2}{}7\above{1pt}{1pt}{-}{}{\cdot}1\above{1pt}{1pt}{2}{}11\above{1pt}{1pt}{1}{}$}\spacer%
\instructions{2\rightarrow N_{119}}%
\EasyButWeakLineBreak%
{${2}\above{1pt}{1pt}{*}{4}{4}\above{1pt}{1pt}{b}{2}{70}\above{1pt}{1pt}{b}{2}{132}\above{1pt}{1pt}{*}{2}{20}\above{1pt}{1pt}{b}{2}{42}\above{1pt}{1pt}{l}{2}{220}\above{1pt}{1pt}{r}{2}{6}\above{1pt}{1pt}{s}{2}{770}\above{1pt}{1pt}{b}{2}$}\relax$\,(\times2)$%
\nopagebreak\par%
\nopagebreak\par\leavevmode%
{$\left[\!\llap{\phantom{%
\begingroup \smaller\smaller\smaller
\endgroup%
}}\!\right]$}%

\medskip%
%
\leavevmode\llap{}%
$W_{120}$%
\qquad\llap{48} lattices, $\chi=40$%
\hfill%
$2222262222$%
\nopagebreak\smallskip\hrule\nopagebreak\medskip%
%
%
\leavevmode%
${L_{120.1}}$%
{} : {$1\above{1pt}{1pt}{-2}{{\rm II}}4\above{1pt}{1pt}{1}{7}{\cdot}1\above{1pt}{1pt}{2}{}3\above{1pt}{1pt}{-}{}{\cdot}1\above{1pt}{1pt}{-2}{}5\above{1pt}{1pt}{-}{}{\cdot}1\above{1pt}{1pt}{-2}{}7\above{1pt}{1pt}{1}{}{\cdot}1\above{1pt}{1pt}{2}{}11\above{1pt}{1pt}{1}{}$}\spacer%
\instructions{2\rightarrow N_{120}}%
\EasyButWeakLineBreak%
{${14}\above{1pt}{1pt}{l}{2}{60}\above{1pt}{1pt}{r}{2}{154}\above{1pt}{1pt}{b}{2}{10}\above{1pt}{1pt}{l}{2}{28}\above{1pt}{1pt}{r}{2}{6}\above{1pt}{1pt}{}{6}{2}\above{1pt}{1pt}{s}{2}{2310}\above{1pt}{1pt}{b}{2}{4}\above{1pt}{1pt}{*}{2}{132}\above{1pt}{1pt}{b}{2}$}%
\nopagebreak\par%
\nopagebreak\par\leavevmode%
{$\left[\!\llap{\phantom{%
\begingroup \smaller\smaller\smaller
\endgroup%
}}\!\right]$}%

\medskip%
%
\leavevmode\llap{}%
$W_{121}$%
\qquad\llap{48} lattices, $\chi=84$%
\hfill%
$222222222222222222\rtimes C_{2}$%
\nopagebreak\smallskip\hrule\nopagebreak\medskip%
%
%
\leavevmode%
${L_{121.1}}$%
{} : {$1\above{1pt}{1pt}{-2}{{\rm II}}4\above{1pt}{1pt}{1}{1}{\cdot}1\above{1pt}{1pt}{2}{}3\above{1pt}{1pt}{1}{}{\cdot}1\above{1pt}{1pt}{-2}{}5\above{1pt}{1pt}{1}{}{\cdot}1\above{1pt}{1pt}{2}{}7\above{1pt}{1pt}{1}{}{\cdot}1\above{1pt}{1pt}{2}{}13\above{1pt}{1pt}{1}{}$}\spacer%
\instructions{2\rightarrow N_{121}}%
\EasyButWeakLineBreak%
{${12}\above{1pt}{1pt}{b}{2}{14}\above{1pt}{1pt}{b}{2}{30}\above{1pt}{1pt}{l}{2}{4}\above{1pt}{1pt}{r}{2}{130}\above{1pt}{1pt}{b}{2}{28}\above{1pt}{1pt}{*}{2}{156}\above{1pt}{1pt}{b}{2}{2}\above{1pt}{1pt}{b}{2}{1820}\above{1pt}{1pt}{*}{2}$}\relax$\,(\times2)$%
\nopagebreak\par%
\nopagebreak\par\leavevmode%
{$\left[\!\llap{\phantom{%
\begingroup \smaller\smaller\smaller
\endgroup%
}}\!\right]$}%

\medskip%
%
\leavevmode\llap{}%
$W_{122}$%
\qquad\llap{48} lattices, $\chi=108$%
\hfill%
$2222222222222222222222\rtimes C_{2}$%
\nopagebreak\smallskip\hrule\nopagebreak\medskip%
%
%
\leavevmode%
${L_{122.1}}$%
{} : {$1\above{1pt}{1pt}{-2}{{\rm II}}4\above{1pt}{1pt}{1}{1}{\cdot}1\above{1pt}{1pt}{2}{}3\above{1pt}{1pt}{1}{}{\cdot}1\above{1pt}{1pt}{2}{}5\above{1pt}{1pt}{-}{}{\cdot}1\above{1pt}{1pt}{2}{}7\above{1pt}{1pt}{1}{}{\cdot}1\above{1pt}{1pt}{-2}{}13\above{1pt}{1pt}{-}{}$}\spacer%
\instructions{2\rightarrow N_{122}}%
\EasyButWeakLineBreak%
{${12}\above{1pt}{1pt}{b}{2}{910}\above{1pt}{1pt}{l}{2}{4}\above{1pt}{1pt}{r}{2}{210}\above{1pt}{1pt}{b}{2}{26}\above{1pt}{1pt}{s}{2}{14}\above{1pt}{1pt}{l}{2}{260}\above{1pt}{1pt}{r}{2}{2}\above{1pt}{1pt}{l}{2}{1092}\above{1pt}{1pt}{r}{2}{10}\above{1pt}{1pt}{b}{2}{28}\above{1pt}{1pt}{*}{2}$}\relax$\,(\times2)$%
\nopagebreak\par%
\nopagebreak\par\leavevmode%
{$\left[\!\llap{\phantom{%
\begingroup \smaller\smaller\smaller
\endgroup%
}}\!\right]$}%

\medskip%
%
\leavevmode\llap{}%
$W_{123}$%
\qquad\llap{16} lattices, $\chi=3$%
\hfill%
$4222$%
\nopagebreak\smallskip\hrule\nopagebreak\medskip%
%
%
\leavevmode%
${L_{123.1}}$%
{} : {$1\above{1pt}{1pt}{2}{2}8\above{1pt}{1pt}{-}{5}{\cdot}1\above{1pt}{1pt}{2}{}3\above{1pt}{1pt}{1}{}$}\spacer%
\instructions{2\rightarrow N'_{3}}%
\EasyButWeakLineBreak%
{${2}\above{1pt}{1pt}{}{4}{1}\above{1pt}{1pt}{r}{2}{12}\above{1pt}{1pt}{*}{2}{8}\above{1pt}{1pt}{b}{2}$}%
\nopagebreak\par%
\nopagebreak\par\leavevmode%
{$\left[\!\llap{\phantom{%
\begingroup \smaller\smaller\smaller\begin{tabular}{@{}c@{}}%
0\\0\\0
\end{tabular}\endgroup%
}}\right.$}%
\begingroup \smaller\smaller\smaller\begin{tabular}{@{}c@{}}%
-1176\\-552\\120
\end{tabular}\endgroup%
\kern3pt%
\begingroup \smaller\smaller\smaller\begin{tabular}{@{}c@{}}%
-552\\-259\\56
\end{tabular}\endgroup%
\kern3pt%
\begingroup \smaller\smaller\smaller\begin{tabular}{@{}c@{}}%
120\\56\\-11
\end{tabular}\endgroup%
{$\left.\llap{\phantom{%
\begingroup \smaller\smaller\smaller\begin{tabular}{@{}c@{}}%
0\\0\\0
\end{tabular}\endgroup%
}}\!\right]$}%
\EasyButWeakLineBreak%
{$\left[\!\llap{\phantom{%
\begingroup \smaller\smaller\smaller\begin{tabular}{@{}c@{}}%
0\\0\\0
\end{tabular}\endgroup%
}}\right.$}%
\begingroup \smaller\smaller\smaller\begin{tabular}{@{}c@{}}%
3\\-7\\-3
\end{tabular}\endgroup%
\HardButStrongLineBreak\kern3pt%
\begingroup \smaller\smaller\smaller\begin{tabular}{@{}c@{}}%
-4\\9\\2
\end{tabular}\endgroup%
\HardButStrongLineBreak\kern3pt%
\begingroup \smaller\smaller\smaller\begin{tabular}{@{}c@{}}%
-5\\12\\6
\end{tabular}\endgroup%
\HardButStrongLineBreak\kern3pt%
\begingroup \smaller\smaller\smaller\begin{tabular}{@{}c@{}}%
9\\-20\\-4
\end{tabular}\endgroup%
{$\left.\llap{\phantom{%
\begingroup \smaller\smaller\smaller\begin{tabular}{@{}c@{}}%
0\\0\\0
\end{tabular}\endgroup%
}}\!\right]$}%
%
%
\hbox{}\par\smallskip%
%
%
\leavevmode%
${L_{123.2}}$%
{} : {$1\above{1pt}{1pt}{-2}{2}8\above{1pt}{1pt}{1}{1}{\cdot}1\above{1pt}{1pt}{2}{}3\above{1pt}{1pt}{1}{}$}\spacer%
\instructions{m}%
\EasyButWeakLineBreak%
{${2}\above{1pt}{1pt}{*}{4}{4}\above{1pt}{1pt}{l}{2}{3}\above{1pt}{1pt}{}{2}{8}\above{1pt}{1pt}{r}{2}$}%
\nopagebreak\par%
\nopagebreak\par\leavevmode%
{$\left[\!\llap{\phantom{%
\begingroup \smaller\smaller\smaller\begin{tabular}{@{}c@{}}%
0\\0\\0
\end{tabular}\endgroup%
}}\right.$}%
\begingroup \smaller\smaller\smaller\begin{tabular}{@{}c@{}}%
-2616\\120\\72
\end{tabular}\endgroup%
\kern3pt%
\begingroup \smaller\smaller\smaller\begin{tabular}{@{}c@{}}%
120\\-5\\-4
\end{tabular}\endgroup%
\kern3pt%
\begingroup \smaller\smaller\smaller\begin{tabular}{@{}c@{}}%
72\\-4\\-1
\end{tabular}\endgroup%
{$\left.\llap{\phantom{%
\begingroup \smaller\smaller\smaller\begin{tabular}{@{}c@{}}%
0\\0\\0
\end{tabular}\endgroup%
}}\!\right]$}%
\EasyButWeakLineBreak%
{$\left[\!\llap{\phantom{%
\begingroup \smaller\smaller\smaller\begin{tabular}{@{}c@{}}%
0\\0\\0
\end{tabular}\endgroup%
}}\right.$}%
\begingroup \smaller\smaller\smaller\begin{tabular}{@{}c@{}}%
1\\15\\11
\end{tabular}\endgroup%
\HardButStrongLineBreak\kern3pt%
\begingroup \smaller\smaller\smaller\begin{tabular}{@{}c@{}}%
-1\\-16\\-10
\end{tabular}\endgroup%
\HardButStrongLineBreak\kern3pt%
\begingroup \smaller\smaller\smaller\begin{tabular}{@{}c@{}}%
-1\\-15\\-12
\end{tabular}\endgroup%
\HardButStrongLineBreak\kern3pt%
\begingroup \smaller\smaller\smaller\begin{tabular}{@{}c@{}}%
1\\16\\8
\end{tabular}\endgroup%
{$\left.\llap{\phantom{%
\begingroup \smaller\smaller\smaller\begin{tabular}{@{}c@{}}%
0\\0\\0
\end{tabular}\endgroup%
}}\!\right]$}%

\medskip%
%
\leavevmode\llap{}%
$W_{124}$%
\qquad\llap{46} lattices, $\chi=12$%
\hfill%
$222222\rtimes C_{2}$%
\nopagebreak\smallskip\hrule\nopagebreak\medskip%
%
%
\leavevmode%
${L_{124.1}}$%
{} : {$1\above{1pt}{1pt}{-2}{4}8\above{1pt}{1pt}{1}{1}{\cdot}1\above{1pt}{1pt}{-2}{}5\above{1pt}{1pt}{1}{}$}\spacer%
\instructions{2\rightarrow N'_{4}}%
\EasyButWeakLineBreak%
{${4}\above{1pt}{1pt}{l}{2}{5}\above{1pt}{1pt}{}{2}{8}\above{1pt}{1pt}{}{2}{1}\above{1pt}{1pt}{r}{2}{20}\above{1pt}{1pt}{*}{2}{8}\above{1pt}{1pt}{*}{2}$}%
\nopagebreak\par%
\nopagebreak\par\leavevmode%
{$\left[\!\llap{\phantom{%
\begingroup \smaller\smaller\smaller\begin{tabular}{@{}c@{}}%
0\\0\\0
\end{tabular}\endgroup%
}}\right.$}%
\begingroup \smaller\smaller\smaller\begin{tabular}{@{}c@{}}%
-3320\\120\\120
\end{tabular}\endgroup%
\kern3pt%
\begingroup \smaller\smaller\smaller\begin{tabular}{@{}c@{}}%
120\\-4\\-5
\end{tabular}\endgroup%
\kern3pt%
\begingroup \smaller\smaller\smaller\begin{tabular}{@{}c@{}}%
120\\-5\\-3
\end{tabular}\endgroup%
{$\left.\llap{\phantom{%
\begingroup \smaller\smaller\smaller\begin{tabular}{@{}c@{}}%
0\\0\\0
\end{tabular}\endgroup%
}}\!\right]$}%
\EasyButWeakLineBreak%
{$\left[\!\llap{\phantom{%
\begingroup \smaller\smaller\smaller\begin{tabular}{@{}c@{}}%
0\\0\\0
\end{tabular}\endgroup%
}}\right.$}%
\begingroup \smaller\smaller\smaller\begin{tabular}{@{}c@{}}%
-1\\-18\\-10
\end{tabular}\endgroup%
\HardButStrongLineBreak\kern3pt%
\begingroup \smaller\smaller\smaller\begin{tabular}{@{}c@{}}%
1\\20\\5
\end{tabular}\endgroup%
\HardButStrongLineBreak\kern3pt%
\begingroup \smaller\smaller\smaller\begin{tabular}{@{}c@{}}%
3\\56\\24
\end{tabular}\endgroup%
\HardButStrongLineBreak\kern3pt%
\begingroup \smaller\smaller\smaller\begin{tabular}{@{}c@{}}%
1\\18\\9
\end{tabular}\endgroup%
\HardButStrongLineBreak\kern3pt%
\begingroup \smaller\smaller\smaller\begin{tabular}{@{}c@{}}%
3\\50\\30
\end{tabular}\endgroup%
\HardButStrongLineBreak\kern3pt%
\begingroup \smaller\smaller\smaller\begin{tabular}{@{}c@{}}%
-1\\-20\\-8
\end{tabular}\endgroup%
{$\left.\llap{\phantom{%
\begingroup \smaller\smaller\smaller\begin{tabular}{@{}c@{}}%
0\\0\\0
\end{tabular}\endgroup%
}}\!\right]$}%
%
%
\hbox{}\par\smallskip%
%
%
\leavevmode%
${L_{124.2}}$%
{} : {$1\above{1pt}{1pt}{-2}{6}8\above{1pt}{1pt}{1}{7}{\cdot}1\above{1pt}{1pt}{-2}{}5\above{1pt}{1pt}{1}{}$}\spacer%
\instructions{m}%
\EasyButWeakLineBreak%
{${1}\above{1pt}{1pt}{}{2}{5}\above{1pt}{1pt}{r}{2}{8}\above{1pt}{1pt}{l}{2}$}\relax$\,(\times2)$%
\nopagebreak\par%
\nopagebreak\par\leavevmode%
{$\left[\!\llap{\phantom{%
\begingroup \smaller\smaller\smaller\begin{tabular}{@{}c@{}}%
0\\0\\0
\end{tabular}\endgroup%
}}\right.$}%
\begingroup \smaller\smaller\smaller\begin{tabular}{@{}c@{}}%
-2120\\280\\80
\end{tabular}\endgroup%
\kern3pt%
\begingroup \smaller\smaller\smaller\begin{tabular}{@{}c@{}}%
280\\-19\\-10
\end{tabular}\endgroup%
\kern3pt%
\begingroup \smaller\smaller\smaller\begin{tabular}{@{}c@{}}%
80\\-10\\-3
\end{tabular}\endgroup%
{$\left.\llap{\phantom{%
\begingroup \smaller\smaller\smaller\begin{tabular}{@{}c@{}}%
0\\0\\0
\end{tabular}\endgroup%
}}\!\right]$}%
\hfil\penalty500%
{$\left[\!\llap{\phantom{%
\begingroup \smaller\smaller\smaller\begin{tabular}{@{}c@{}}%
0\\0\\0
\end{tabular}\endgroup%
}}\right.$}%
\begingroup \smaller\smaller\smaller\begin{tabular}{@{}c@{}}%
79\\-80\\2320
\end{tabular}\endgroup%
\kern3pt%
\begingroup \smaller\smaller\smaller\begin{tabular}{@{}c@{}}%
-9\\8\\-261
\end{tabular}\endgroup%
\kern3pt%
\begingroup \smaller\smaller\smaller\begin{tabular}{@{}c@{}}%
-3\\3\\-88
\end{tabular}\endgroup%
{$\left.\llap{\phantom{%
\begingroup \smaller\smaller\smaller\begin{tabular}{@{}c@{}}%
0\\0\\0
\end{tabular}\endgroup%
}}\!\right]$}%
\EasyButWeakLineBreak%
{$\left[\!\llap{\phantom{%
\begingroup \smaller\smaller\smaller\begin{tabular}{@{}c@{}}%
0\\0\\0
\end{tabular}\endgroup%
}}\right.$}%
\begingroup \smaller\smaller\smaller\begin{tabular}{@{}c@{}}%
-1\\1\\-30
\end{tabular}\endgroup%
\HardButStrongLineBreak\kern3pt%
\begingroup \smaller\smaller\smaller\begin{tabular}{@{}c@{}}%
1\\0\\25
\end{tabular}\endgroup%
\HardButStrongLineBreak\kern3pt%
\begingroup \smaller\smaller\smaller\begin{tabular}{@{}c@{}}%
5\\-4\\144
\end{tabular}\endgroup%
{$\left.\llap{\phantom{%
\begingroup \smaller\smaller\smaller\begin{tabular}{@{}c@{}}%
0\\0\\0
\end{tabular}\endgroup%
}}\!\right]$}%
%
%
\hbox{}\par\smallskip%
%
%
\leavevmode%
${L_{124.3}}$%
{} : {$1\above{1pt}{1pt}{2}{6}8\above{1pt}{1pt}{-}{3}{\cdot}1\above{1pt}{1pt}{-2}{}5\above{1pt}{1pt}{1}{}$}\EasyButWeakLineBreak%
{${4}\above{1pt}{1pt}{*}{2}{20}\above{1pt}{1pt}{s}{2}{8}\above{1pt}{1pt}{s}{2}$}\relax$\,(\times2)$%
\nopagebreak\par%
\nopagebreak\par\leavevmode%
{$\left[\!\llap{\phantom{%
\begingroup \smaller\smaller\smaller\begin{tabular}{@{}c@{}}%
0\\0\\0
\end{tabular}\endgroup%
}}\right.$}%
\begingroup \smaller\smaller\smaller\begin{tabular}{@{}c@{}}%
920\\200\\0
\end{tabular}\endgroup%
\kern3pt%
\begingroup \smaller\smaller\smaller\begin{tabular}{@{}c@{}}%
200\\43\\1
\end{tabular}\endgroup%
\kern3pt%
\begingroup \smaller\smaller\smaller\begin{tabular}{@{}c@{}}%
0\\1\\-2
\end{tabular}\endgroup%
{$\left.\llap{\phantom{%
\begingroup \smaller\smaller\smaller\begin{tabular}{@{}c@{}}%
0\\0\\0
\end{tabular}\endgroup%
}}\!\right]$}%
\hfil\penalty500%
{$\left[\!\llap{\phantom{%
\begingroup \smaller\smaller\smaller\begin{tabular}{@{}c@{}}%
0\\0\\0
\end{tabular}\endgroup%
}}\right.$}%
\begingroup \smaller\smaller\smaller\begin{tabular}{@{}c@{}}%
-1\\0\\0
\end{tabular}\endgroup%
\kern3pt%
\begingroup \smaller\smaller\smaller\begin{tabular}{@{}c@{}}%
0\\-1\\-1
\end{tabular}\endgroup%
\kern3pt%
\begingroup \smaller\smaller\smaller\begin{tabular}{@{}c@{}}%
0\\0\\1
\end{tabular}\endgroup%
{$\left.\llap{\phantom{%
\begingroup \smaller\smaller\smaller\begin{tabular}{@{}c@{}}%
0\\0\\0
\end{tabular}\endgroup%
}}\!\right]$}%
\EasyButWeakLineBreak%
{$\left[\!\llap{\phantom{%
\begingroup \smaller\smaller\smaller\begin{tabular}{@{}c@{}}%
0\\0\\0
\end{tabular}\endgroup%
}}\right.$}%
\begingroup \smaller\smaller\smaller\begin{tabular}{@{}c@{}}%
3\\-14\\-8
\end{tabular}\endgroup%
\HardButStrongLineBreak\kern3pt%
\begingroup \smaller\smaller\smaller\begin{tabular}{@{}c@{}}%
11\\-50\\-30
\end{tabular}\endgroup%
\HardButStrongLineBreak\kern3pt%
\begingroup \smaller\smaller\smaller\begin{tabular}{@{}c@{}}%
1\\-4\\-4
\end{tabular}\endgroup%
{$\left.\llap{\phantom{%
\begingroup \smaller\smaller\smaller\begin{tabular}{@{}c@{}}%
0\\0\\0
\end{tabular}\endgroup%
}}\!\right]$}%
%
%
\hbox{}\par\smallskip%
%
%
\leavevmode%
${L_{124.4}}$%
{} : {$[1\above{1pt}{1pt}{1}{}2\above{1pt}{1pt}{1}{}]\above{1pt}{1pt}{}{2}16\above{1pt}{1pt}{-}{3}{\cdot}1\above{1pt}{1pt}{-2}{}5\above{1pt}{1pt}{1}{}$}\spacer%
\instructions{2}%
\EasyButWeakLineBreak%
{${1}\above{1pt}{1pt}{r}{2}{80}\above{1pt}{1pt}{*}{2}{8}\above{1pt}{1pt}{s}{2}{16}\above{1pt}{1pt}{*}{2}{20}\above{1pt}{1pt}{l}{2}{2}\above{1pt}{1pt}{}{2}$}%
\nopagebreak\par%
\nopagebreak\par\leavevmode%
{$\left[\!\llap{\phantom{%
\begingroup \smaller\smaller\smaller\begin{tabular}{@{}c@{}}%
0\\0\\0
\end{tabular}\endgroup%
}}\right.$}%
\begingroup \smaller\smaller\smaller\begin{tabular}{@{}c@{}}%
5680\\2720\\-80
\end{tabular}\endgroup%
\kern3pt%
\begingroup \smaller\smaller\smaller\begin{tabular}{@{}c@{}}%
2720\\1302\\-38
\end{tabular}\endgroup%
\kern3pt%
\begingroup \smaller\smaller\smaller\begin{tabular}{@{}c@{}}%
-80\\-38\\1
\end{tabular}\endgroup%
{$\left.\llap{\phantom{%
\begingroup \smaller\smaller\smaller\begin{tabular}{@{}c@{}}%
0\\0\\0
\end{tabular}\endgroup%
}}\!\right]$}%
\EasyButWeakLineBreak%
{$\left[\!\llap{\phantom{%
\begingroup \smaller\smaller\smaller\begin{tabular}{@{}c@{}}%
0\\0\\0
\end{tabular}\endgroup%
}}\right.$}%
\begingroup \smaller\smaller\smaller\begin{tabular}{@{}c@{}}%
0\\0\\-1
\end{tabular}\endgroup%
\HardButStrongLineBreak\kern3pt%
\begingroup \smaller\smaller\smaller\begin{tabular}{@{}c@{}}%
-9\\20\\40
\end{tabular}\endgroup%
\HardButStrongLineBreak\kern3pt%
\begingroup \smaller\smaller\smaller\begin{tabular}{@{}c@{}}%
-1\\2\\0
\end{tabular}\endgroup%
\HardButStrongLineBreak\kern3pt%
\begingroup \smaller\smaller\smaller\begin{tabular}{@{}c@{}}%
5\\-12\\-48
\end{tabular}\endgroup%
\HardButStrongLineBreak\kern3pt%
\begingroup \smaller\smaller\smaller\begin{tabular}{@{}c@{}}%
17\\-40\\-150
\end{tabular}\endgroup%
\HardButStrongLineBreak\kern3pt%
\begingroup \smaller\smaller\smaller\begin{tabular}{@{}c@{}}%
3\\-7\\-26
\end{tabular}\endgroup%
{$\left.\llap{\phantom{%
\begingroup \smaller\smaller\smaller\begin{tabular}{@{}c@{}}%
0\\0\\0
\end{tabular}\endgroup%
}}\!\right]$}%
%
%
\hbox{}\par\smallskip%
%
%
\leavevmode%
${L_{124.5}}$%
{} : {$[1\above{1pt}{1pt}{-}{}2\above{1pt}{1pt}{1}{}]\above{1pt}{1pt}{}{6}16\above{1pt}{1pt}{1}{7}{\cdot}1\above{1pt}{1pt}{-2}{}5\above{1pt}{1pt}{1}{}$}\spacer%
\instructions{m}%
\EasyButWeakLineBreak%
{${4}\above{1pt}{1pt}{*}{2}{80}\above{1pt}{1pt}{s}{2}{8}\above{1pt}{1pt}{*}{2}{16}\above{1pt}{1pt}{l}{2}{5}\above{1pt}{1pt}{}{2}{2}\above{1pt}{1pt}{r}{2}$}%
\nopagebreak\par%
\nopagebreak\par\leavevmode%
{$\left[\!\llap{\phantom{%
\begingroup \smaller\smaller\smaller\begin{tabular}{@{}c@{}}%
0\\0\\0
\end{tabular}\endgroup%
}}\right.$}%
\begingroup \smaller\smaller\smaller\begin{tabular}{@{}c@{}}%
4080\\-160\\80
\end{tabular}\endgroup%
\kern3pt%
\begingroup \smaller\smaller\smaller\begin{tabular}{@{}c@{}}%
-160\\6\\-2
\end{tabular}\endgroup%
\kern3pt%
\begingroup \smaller\smaller\smaller\begin{tabular}{@{}c@{}}%
80\\-2\\-3
\end{tabular}\endgroup%
{$\left.\llap{\phantom{%
\begingroup \smaller\smaller\smaller\begin{tabular}{@{}c@{}}%
0\\0\\0
\end{tabular}\endgroup%
}}\!\right]$}%
\EasyButWeakLineBreak%
{$\left[\!\llap{\phantom{%
\begingroup \smaller\smaller\smaller\begin{tabular}{@{}c@{}}%
0\\0\\0
\end{tabular}\endgroup%
}}\right.$}%
\begingroup \smaller\smaller\smaller\begin{tabular}{@{}c@{}}%
1\\28\\6
\end{tabular}\endgroup%
\HardButStrongLineBreak\kern3pt%
\begingroup \smaller\smaller\smaller\begin{tabular}{@{}c@{}}%
1\\20\\0
\end{tabular}\endgroup%
\HardButStrongLineBreak\kern3pt%
\begingroup \smaller\smaller\smaller\begin{tabular}{@{}c@{}}%
-1\\-30\\-8
\end{tabular}\endgroup%
\HardButStrongLineBreak\kern3pt%
\begingroup \smaller\smaller\smaller\begin{tabular}{@{}c@{}}%
-1\\-28\\-8
\end{tabular}\endgroup%
\HardButStrongLineBreak\kern3pt%
\begingroup \smaller\smaller\smaller\begin{tabular}{@{}c@{}}%
1\\30\\5
\end{tabular}\endgroup%
\HardButStrongLineBreak\kern3pt%
\begingroup \smaller\smaller\smaller\begin{tabular}{@{}c@{}}%
1\\29\\6
\end{tabular}\endgroup%
{$\left.\llap{\phantom{%
\begingroup \smaller\smaller\smaller\begin{tabular}{@{}c@{}}%
0\\0\\0
\end{tabular}\endgroup%
}}\!\right]$}%
%
%
\hbox{}\par\smallskip%
%
%
\leavevmode%
${L_{124.6}}$%
{} : {$[1\above{1pt}{1pt}{1}{}2\above{1pt}{1pt}{1}{}]\above{1pt}{1pt}{}{0}16\above{1pt}{1pt}{-}{5}{\cdot}1\above{1pt}{1pt}{-2}{}5\above{1pt}{1pt}{1}{}$}\spacer%
\instructions{m}%
\EasyButWeakLineBreak%
{${1}\above{1pt}{1pt}{}{2}{80}\above{1pt}{1pt}{}{2}{2}\above{1pt}{1pt}{r}{2}{16}\above{1pt}{1pt}{s}{2}{20}\above{1pt}{1pt}{*}{2}{8}\above{1pt}{1pt}{l}{2}$}%
\nopagebreak\par%
\nopagebreak\par\leavevmode%
{$\left[\!\llap{\phantom{%
\begingroup \smaller\smaller\smaller\begin{tabular}{@{}c@{}}%
0\\0\\0
\end{tabular}\endgroup%
}}\right.$}%
\begingroup \smaller\smaller\smaller\begin{tabular}{@{}c@{}}%
-2480\\80\\80
\end{tabular}\endgroup%
\kern3pt%
\begingroup \smaller\smaller\smaller\begin{tabular}{@{}c@{}}%
80\\-2\\-4
\end{tabular}\endgroup%
\kern3pt%
\begingroup \smaller\smaller\smaller\begin{tabular}{@{}c@{}}%
80\\-4\\1
\end{tabular}\endgroup%
{$\left.\llap{\phantom{%
\begingroup \smaller\smaller\smaller\begin{tabular}{@{}c@{}}%
0\\0\\0
\end{tabular}\endgroup%
}}\!\right]$}%
\EasyButWeakLineBreak%
{$\left[\!\llap{\phantom{%
\begingroup \smaller\smaller\smaller\begin{tabular}{@{}c@{}}%
0\\0\\0
\end{tabular}\endgroup%
}}\right.$}%
\begingroup \smaller\smaller\smaller\begin{tabular}{@{}c@{}}%
0\\0\\-1
\end{tabular}\endgroup%
\HardButStrongLineBreak\kern3pt%
\begingroup \smaller\smaller\smaller\begin{tabular}{@{}c@{}}%
-7\\-160\\-80
\end{tabular}\endgroup%
\HardButStrongLineBreak\kern3pt%
\begingroup \smaller\smaller\smaller\begin{tabular}{@{}c@{}}%
-1\\-23\\-10
\end{tabular}\endgroup%
\HardButStrongLineBreak\kern3pt%
\begingroup \smaller\smaller\smaller\begin{tabular}{@{}c@{}}%
-1\\-24\\-8
\end{tabular}\endgroup%
\HardButStrongLineBreak\kern3pt%
\begingroup \smaller\smaller\smaller\begin{tabular}{@{}c@{}}%
1\\20\\10
\end{tabular}\endgroup%
\HardButStrongLineBreak\kern3pt%
\begingroup \smaller\smaller\smaller\begin{tabular}{@{}c@{}}%
1\\22\\8
\end{tabular}\endgroup%
{$\left.\llap{\phantom{%
\begingroup \smaller\smaller\smaller\begin{tabular}{@{}c@{}}%
0\\0\\0
\end{tabular}\endgroup%
}}\!\right]$}%
%
%
\hbox{}\par\smallskip%
%
%
\leavevmode%
${L_{124.7}}$%
{} : {$[1\above{1pt}{1pt}{-}{}2\above{1pt}{1pt}{1}{}]\above{1pt}{1pt}{}{4}16\above{1pt}{1pt}{1}{1}{\cdot}1\above{1pt}{1pt}{-2}{}5\above{1pt}{1pt}{1}{}$}\EasyButWeakLineBreak%
{${4}\above{1pt}{1pt}{s}{2}{80}\above{1pt}{1pt}{l}{2}{2}\above{1pt}{1pt}{}{2}{16}\above{1pt}{1pt}{}{2}{5}\above{1pt}{1pt}{r}{2}{8}\above{1pt}{1pt}{*}{2}$}%
\nopagebreak\par%
\nopagebreak\par\leavevmode%
{$\left[\!\llap{\phantom{%
\begingroup \smaller\smaller\smaller\begin{tabular}{@{}c@{}}%
0\\0\\0
\end{tabular}\endgroup%
}}\right.$}%
\begingroup \smaller\smaller\smaller\begin{tabular}{@{}c@{}}%
-64880\\1120\\1120
\end{tabular}\endgroup%
\kern3pt%
\begingroup \smaller\smaller\smaller\begin{tabular}{@{}c@{}}%
1120\\-18\\-20
\end{tabular}\endgroup%
\kern3pt%
\begingroup \smaller\smaller\smaller\begin{tabular}{@{}c@{}}%
1120\\-20\\-19
\end{tabular}\endgroup%
{$\left.\llap{\phantom{%
\begingroup \smaller\smaller\smaller\begin{tabular}{@{}c@{}}%
0\\0\\0
\end{tabular}\endgroup%
}}\!\right]$}%
\EasyButWeakLineBreak%
{$\left[\!\llap{\phantom{%
\begingroup \smaller\smaller\smaller\begin{tabular}{@{}c@{}}%
0\\0\\0
\end{tabular}\endgroup%
}}\right.$}%
\begingroup \smaller\smaller\smaller\begin{tabular}{@{}c@{}}%
-1\\-20\\-38
\end{tabular}\endgroup%
\HardButStrongLineBreak\kern3pt%
\begingroup \smaller\smaller\smaller\begin{tabular}{@{}c@{}}%
7\\120\\280
\end{tabular}\endgroup%
\HardButStrongLineBreak\kern3pt%
\begingroup \smaller\smaller\smaller\begin{tabular}{@{}c@{}}%
2\\37\\78
\end{tabular}\endgroup%
\HardButStrongLineBreak\kern3pt%
\begingroup \smaller\smaller\smaller\begin{tabular}{@{}c@{}}%
5\\96\\192
\end{tabular}\endgroup%
\HardButStrongLineBreak\kern3pt%
\begingroup \smaller\smaller\smaller\begin{tabular}{@{}c@{}}%
2\\40\\75
\end{tabular}\endgroup%
\HardButStrongLineBreak\kern3pt%
\begingroup \smaller\smaller\smaller\begin{tabular}{@{}c@{}}%
-1\\-18\\-40
\end{tabular}\endgroup%
{$\left.\llap{\phantom{%
\begingroup \smaller\smaller\smaller\begin{tabular}{@{}c@{}}%
0\\0\\0
\end{tabular}\endgroup%
}}\!\right]$}%

\medskip%
%
\leavevmode\llap{}%
$W_{125}$%
\qquad\llap{20} lattices, $\chi=18$%
\hfill%
$422422\rtimes C_{2}$%
\nopagebreak\smallskip\hrule\nopagebreak\medskip%
%
%
\leavevmode%
${L_{125.1}}$%
{} : {$1\above{1pt}{1pt}{-2}{2}8\above{1pt}{1pt}{-}{5}{\cdot}1\above{1pt}{1pt}{2}{}7\above{1pt}{1pt}{-}{}$}\spacer%
\instructions{2\rightarrow N'_{6}}%
\EasyButWeakLineBreak%
{${2}\above{1pt}{1pt}{*}{4}{4}\above{1pt}{1pt}{*}{2}{8}\above{1pt}{1pt}{b}{2}$}\relax$\,(\times2)$%
\nopagebreak\par%
\nopagebreak\par\leavevmode%
{$\left[\!\llap{\phantom{%
\begingroup \smaller\smaller\smaller\begin{tabular}{@{}c@{}}%
0\\0\\0
\end{tabular}\endgroup%
}}\right.$}%
\begingroup \smaller\smaller\smaller\begin{tabular}{@{}c@{}}%
-28056\\392\\448
\end{tabular}\endgroup%
\kern3pt%
\begingroup \smaller\smaller\smaller\begin{tabular}{@{}c@{}}%
392\\-5\\-7
\end{tabular}\endgroup%
\kern3pt%
\begingroup \smaller\smaller\smaller\begin{tabular}{@{}c@{}}%
448\\-7\\-6
\end{tabular}\endgroup%
{$\left.\llap{\phantom{%
\begingroup \smaller\smaller\smaller\begin{tabular}{@{}c@{}}%
0\\0\\0
\end{tabular}\endgroup%
}}\!\right]$}%
\hfil\penalty500%
{$\left[\!\llap{\phantom{%
\begingroup \smaller\smaller\smaller\begin{tabular}{@{}c@{}}%
0\\0\\0
\end{tabular}\endgroup%
}}\right.$}%
\begingroup \smaller\smaller\smaller\begin{tabular}{@{}c@{}}%
671\\27552\\17472
\end{tabular}\endgroup%
\kern3pt%
\begingroup \smaller\smaller\smaller\begin{tabular}{@{}c@{}}%
-10\\-411\\-260
\end{tabular}\endgroup%
\kern3pt%
\begingroup \smaller\smaller\smaller\begin{tabular}{@{}c@{}}%
-10\\-410\\-261
\end{tabular}\endgroup%
{$\left.\llap{\phantom{%
\begingroup \smaller\smaller\smaller\begin{tabular}{@{}c@{}}%
0\\0\\0
\end{tabular}\endgroup%
}}\!\right]$}%
\EasyButWeakLineBreak%
{$\left[\!\llap{\phantom{%
\begingroup \smaller\smaller\smaller\begin{tabular}{@{}c@{}}%
0\\0\\0
\end{tabular}\endgroup%
}}\right.$}%
\begingroup \smaller\smaller\smaller\begin{tabular}{@{}c@{}}%
1\\40\\27
\end{tabular}\endgroup%
\HardButStrongLineBreak\kern3pt%
\begingroup \smaller\smaller\smaller\begin{tabular}{@{}c@{}}%
-1\\-42\\-26
\end{tabular}\endgroup%
\HardButStrongLineBreak\kern3pt%
\begingroup \smaller\smaller\smaller\begin{tabular}{@{}c@{}}%
-1\\-40\\-28
\end{tabular}\endgroup%
{$\left.\llap{\phantom{%
\begingroup \smaller\smaller\smaller\begin{tabular}{@{}c@{}}%
0\\0\\0
\end{tabular}\endgroup%
}}\!\right]$}%
%
%
\hbox{}\par\smallskip%
%
%
\leavevmode%
${L_{125.2}}$%
{} : {$1\above{1pt}{1pt}{2}{2}16\above{1pt}{1pt}{1}{1}{\cdot}1\above{1pt}{1pt}{2}{}7\above{1pt}{1pt}{-}{}$}\spacer%
\instructions{2,m}%
\EasyButWeakLineBreak%
{${1}\above{1pt}{1pt}{}{4}{2}\above{1pt}{1pt}{b}{2}{16}\above{1pt}{1pt}{*}{2}{4}\above{1pt}{1pt}{*}{4}{2}\above{1pt}{1pt}{l}{2}{16}\above{1pt}{1pt}{}{2}$}%
\nopagebreak\par%
\nopagebreak\par\leavevmode%
{$\left[\!\llap{\phantom{%
\begingroup \smaller\smaller\smaller\begin{tabular}{@{}c@{}}%
0\\0\\0
\end{tabular}\endgroup%
}}\right.$}%
\begingroup \smaller\smaller\smaller\begin{tabular}{@{}c@{}}%
-115696\\-57232\\1792
\end{tabular}\endgroup%
\kern3pt%
\begingroup \smaller\smaller\smaller\begin{tabular}{@{}c@{}}%
-57232\\-28311\\886
\end{tabular}\endgroup%
\kern3pt%
\begingroup \smaller\smaller\smaller\begin{tabular}{@{}c@{}}%
1792\\886\\-27
\end{tabular}\endgroup%
{$\left.\llap{\phantom{%
\begingroup \smaller\smaller\smaller\begin{tabular}{@{}c@{}}%
0\\0\\0
\end{tabular}\endgroup%
}}\!\right]$}%
\EasyButWeakLineBreak%
{$\left[\!\llap{\phantom{%
\begingroup \smaller\smaller\smaller\begin{tabular}{@{}c@{}}%
0\\0\\0
\end{tabular}\endgroup%
}}\right.$}%
\begingroup \smaller\smaller\smaller\begin{tabular}{@{}c@{}}%
-17\\35\\20
\end{tabular}\endgroup%
\HardButStrongLineBreak\kern3pt%
\begingroup \smaller\smaller\smaller\begin{tabular}{@{}c@{}}%
16\\-33\\-21
\end{tabular}\endgroup%
\HardButStrongLineBreak\kern3pt%
\begingroup \smaller\smaller\smaller\begin{tabular}{@{}c@{}}%
35\\-72\\-40
\end{tabular}\endgroup%
\HardButStrongLineBreak\kern3pt%
\begingroup \smaller\smaller\smaller\begin{tabular}{@{}c@{}}%
-31\\64\\42
\end{tabular}\endgroup%
\HardButStrongLineBreak\kern3pt%
\begingroup \smaller\smaller\smaller\begin{tabular}{@{}c@{}}%
-146\\301\\185
\end{tabular}\endgroup%
\HardButStrongLineBreak\kern3pt%
\begingroup \smaller\smaller\smaller\begin{tabular}{@{}c@{}}%
-295\\608\\368
\end{tabular}\endgroup%
{$\left.\llap{\phantom{%
\begingroup \smaller\smaller\smaller\begin{tabular}{@{}c@{}}%
0\\0\\0
\end{tabular}\endgroup%
}}\!\right]$}%

\medskip%
%
\leavevmode\llap{}%
$W_{126}$%
\qquad\llap{16} lattices, $\chi=30$%
\hfill%
$42224222\rtimes C_{2}$%
\nopagebreak\smallskip\hrule\nopagebreak\medskip%
%
%
\leavevmode%
${L_{126.1}}$%
{} : {$1\above{1pt}{1pt}{2}{2}8\above{1pt}{1pt}{-}{5}{\cdot}1\above{1pt}{1pt}{2}{}11\above{1pt}{1pt}{1}{}$}\spacer%
\instructions{2\rightarrow N'_{9}}%
\EasyButWeakLineBreak%
{${2}\above{1pt}{1pt}{}{4}{1}\above{1pt}{1pt}{r}{2}{44}\above{1pt}{1pt}{*}{2}{8}\above{1pt}{1pt}{b}{2}$}\relax$\,(\times2)$%
\nopagebreak\par%
\nopagebreak\par\leavevmode%
{$\left[\!\llap{\phantom{%
\begingroup \smaller\smaller\smaller\begin{tabular}{@{}c@{}}%
0\\0\\0
\end{tabular}\endgroup%
}}\right.$}%
\begingroup \smaller\smaller\smaller\begin{tabular}{@{}c@{}}%
616\\-176\\-88
\end{tabular}\endgroup%
\kern3pt%
\begingroup \smaller\smaller\smaller\begin{tabular}{@{}c@{}}%
-176\\50\\23
\end{tabular}\endgroup%
\kern3pt%
\begingroup \smaller\smaller\smaller\begin{tabular}{@{}c@{}}%
-88\\23\\-3
\end{tabular}\endgroup%
{$\left.\llap{\phantom{%
\begingroup \smaller\smaller\smaller\begin{tabular}{@{}c@{}}%
0\\0\\0
\end{tabular}\endgroup%
}}\!\right]$}%
\hfil\penalty500%
{$\left[\!\llap{\phantom{%
\begingroup \smaller\smaller\smaller\begin{tabular}{@{}c@{}}%
0\\0\\0
\end{tabular}\endgroup%
}}\right.$}%
\begingroup \smaller\smaller\smaller\begin{tabular}{@{}c@{}}%
109\\440\\440
\end{tabular}\endgroup%
\kern3pt%
\begingroup \smaller\smaller\smaller\begin{tabular}{@{}c@{}}%
-29\\-117\\-116
\end{tabular}\endgroup%
\kern3pt%
\begingroup \smaller\smaller\smaller\begin{tabular}{@{}c@{}}%
2\\8\\7
\end{tabular}\endgroup%
{$\left.\llap{\phantom{%
\begingroup \smaller\smaller\smaller\begin{tabular}{@{}c@{}}%
0\\0\\0
\end{tabular}\endgroup%
}}\!\right]$}%
\EasyButWeakLineBreak%
{$\left[\!\llap{\phantom{%
\begingroup \smaller\smaller\smaller\begin{tabular}{@{}c@{}}%
0\\0\\0
\end{tabular}\endgroup%
}}\right.$}%
\begingroup \smaller\smaller\smaller\begin{tabular}{@{}c@{}}%
3\\11\\-6
\end{tabular}\endgroup%
\HardButStrongLineBreak\kern3pt%
\begingroup \smaller\smaller\smaller\begin{tabular}{@{}c@{}}%
4\\15\\-3
\end{tabular}\endgroup%
\HardButStrongLineBreak\kern3pt%
\begingroup \smaller\smaller\smaller\begin{tabular}{@{}c@{}}%
35\\132\\-22
\end{tabular}\endgroup%
\HardButStrongLineBreak\kern3pt%
\begingroup \smaller\smaller\smaller\begin{tabular}{@{}c@{}}%
1\\4\\0
\end{tabular}\endgroup%
{$\left.\llap{\phantom{%
\begingroup \smaller\smaller\smaller\begin{tabular}{@{}c@{}}%
0\\0\\0
\end{tabular}\endgroup%
}}\!\right]$}%
%
%
\hbox{}\par\smallskip%
%
%
\leavevmode%
${L_{126.2}}$%
{} : {$1\above{1pt}{1pt}{-2}{2}8\above{1pt}{1pt}{1}{1}{\cdot}1\above{1pt}{1pt}{2}{}11\above{1pt}{1pt}{1}{}$}\spacer%
\instructions{m}%
\EasyButWeakLineBreak%
{${2}\above{1pt}{1pt}{*}{4}{4}\above{1pt}{1pt}{l}{2}{11}\above{1pt}{1pt}{}{2}{8}\above{1pt}{1pt}{r}{2}$}\relax$\,(\times2)$%
\nopagebreak\par%
\nopagebreak\par\leavevmode%
{$\left[\!\llap{\phantom{%
\begingroup \smaller\smaller\smaller\begin{tabular}{@{}c@{}}%
0\\0\\0
\end{tabular}\endgroup%
}}\right.$}%
\begingroup \smaller\smaller\smaller\begin{tabular}{@{}c@{}}%
-13816\\440\\352
\end{tabular}\endgroup%
\kern3pt%
\begingroup \smaller\smaller\smaller\begin{tabular}{@{}c@{}}%
440\\-14\\-11
\end{tabular}\endgroup%
\kern3pt%
\begingroup \smaller\smaller\smaller\begin{tabular}{@{}c@{}}%
352\\-11\\-5
\end{tabular}\endgroup%
{$\left.\llap{\phantom{%
\begingroup \smaller\smaller\smaller\begin{tabular}{@{}c@{}}%
0\\0\\0
\end{tabular}\endgroup%
}}\!\right]$}%
\hfil\penalty500%
{$\left[\!\llap{\phantom{%
\begingroup \smaller\smaller\smaller\begin{tabular}{@{}c@{}}%
0\\0\\0
\end{tabular}\endgroup%
}}\right.$}%
\begingroup \smaller\smaller\smaller\begin{tabular}{@{}c@{}}%
461\\14784\\-616
\end{tabular}\endgroup%
\kern3pt%
\begingroup \smaller\smaller\smaller\begin{tabular}{@{}c@{}}%
-15\\-481\\20
\end{tabular}\endgroup%
\kern3pt%
\begingroup \smaller\smaller\smaller\begin{tabular}{@{}c@{}}%
-15\\-480\\19
\end{tabular}\endgroup%
{$\left.\llap{\phantom{%
\begingroup \smaller\smaller\smaller\begin{tabular}{@{}c@{}}%
0\\0\\0
\end{tabular}\endgroup%
}}\!\right]$}%
\EasyButWeakLineBreak%
{$\left[\!\llap{\phantom{%
\begingroup \smaller\smaller\smaller\begin{tabular}{@{}c@{}}%
0\\0\\0
\end{tabular}\endgroup%
}}\right.$}%
\begingroup \smaller\smaller\smaller\begin{tabular}{@{}c@{}}%
4\\129\\-6
\end{tabular}\endgroup%
\HardButStrongLineBreak\kern3pt%
\begingroup \smaller\smaller\smaller\begin{tabular}{@{}c@{}}%
3\\98\\-6
\end{tabular}\endgroup%
\HardButStrongLineBreak\kern3pt%
\begingroup \smaller\smaller\smaller\begin{tabular}{@{}c@{}}%
4\\132\\-11
\end{tabular}\endgroup%
\HardButStrongLineBreak\kern3pt%
\begingroup \smaller\smaller\smaller\begin{tabular}{@{}c@{}}%
-1\\-32\\0
\end{tabular}\endgroup%
{$\left.\llap{\phantom{%
\begingroup \smaller\smaller\smaller\begin{tabular}{@{}c@{}}%
0\\0\\0
\end{tabular}\endgroup%
}}\!\right]$}%

\medskip%
%
\leavevmode\llap{}%
$W_{127}$%
\qquad\llap{32} lattices, $\chi=9$%
\hfill%
$42222$%
\nopagebreak\smallskip\hrule\nopagebreak\medskip%
%
%
\leavevmode%
${L_{127.1}}$%
{} : {$1\above{1pt}{1pt}{-2}{2}8\above{1pt}{1pt}{-}{5}{\cdot}1\above{1pt}{1pt}{2}{}3\above{1pt}{1pt}{-}{}{\cdot}1\above{1pt}{1pt}{2}{}5\above{1pt}{1pt}{1}{}$}\spacer%
\instructions{2\rightarrow N'_{11}}%
\EasyButWeakLineBreak%
{${2}\above{1pt}{1pt}{*}{4}{4}\above{1pt}{1pt}{s}{2}{24}\above{1pt}{1pt}{s}{2}{20}\above{1pt}{1pt}{*}{2}{8}\above{1pt}{1pt}{b}{2}$}%
\nopagebreak\par%
\nopagebreak\par\leavevmode%
{$\left[\!\llap{\phantom{%
\begingroup \smaller\smaller\smaller\begin{tabular}{@{}c@{}}%
0\\0\\0
\end{tabular}\endgroup%
}}\right.$}%
\begingroup \smaller\smaller\smaller\begin{tabular}{@{}c@{}}%
-27480\\600\\1200
\end{tabular}\endgroup%
\kern3pt%
\begingroup \smaller\smaller\smaller\begin{tabular}{@{}c@{}}%
600\\-13\\-25
\end{tabular}\endgroup%
\kern3pt%
\begingroup \smaller\smaller\smaller\begin{tabular}{@{}c@{}}%
1200\\-25\\-38
\end{tabular}\endgroup%
{$\left.\llap{\phantom{%
\begingroup \smaller\smaller\smaller\begin{tabular}{@{}c@{}}%
0\\0\\0
\end{tabular}\endgroup%
}}\!\right]$}%
\EasyButWeakLineBreak%
{$\left[\!\llap{\phantom{%
\begingroup \smaller\smaller\smaller\begin{tabular}{@{}c@{}}%
0\\0\\0
\end{tabular}\endgroup%
}}\right.$}%
\begingroup \smaller\smaller\smaller\begin{tabular}{@{}c@{}}%
-2\\-110\\9
\end{tabular}\endgroup%
\HardButStrongLineBreak\kern3pt%
\begingroup \smaller\smaller\smaller\begin{tabular}{@{}c@{}}%
1\\54\\-4
\end{tabular}\endgroup%
\HardButStrongLineBreak\kern3pt%
\begingroup \smaller\smaller\smaller\begin{tabular}{@{}c@{}}%
5\\276\\-24
\end{tabular}\endgroup%
\HardButStrongLineBreak\kern3pt%
\begingroup \smaller\smaller\smaller\begin{tabular}{@{}c@{}}%
-1\\-50\\0
\end{tabular}\endgroup%
\HardButStrongLineBreak\kern3pt%
\begingroup \smaller\smaller\smaller\begin{tabular}{@{}c@{}}%
-5\\-272\\20
\end{tabular}\endgroup%
{$\left.\llap{\phantom{%
\begingroup \smaller\smaller\smaller\begin{tabular}{@{}c@{}}%
0\\0\\0
\end{tabular}\endgroup%
}}\!\right]$}%
%
%
\hbox{}\par\smallskip%
%
%
\leavevmode%
${L_{127.2}}$%
{} : {$1\above{1pt}{1pt}{2}{2}8\above{1pt}{1pt}{1}{1}{\cdot}1\above{1pt}{1pt}{2}{}3\above{1pt}{1pt}{-}{}{\cdot}1\above{1pt}{1pt}{2}{}5\above{1pt}{1pt}{1}{}$}\spacer%
\instructions{m}%
\EasyButWeakLineBreak%
{${2}\above{1pt}{1pt}{}{4}{1}\above{1pt}{1pt}{r}{2}{24}\above{1pt}{1pt}{l}{2}{5}\above{1pt}{1pt}{}{2}{8}\above{1pt}{1pt}{r}{2}$}%
\nopagebreak\par%
\nopagebreak\par\leavevmode%
{$\left[\!\llap{\phantom{%
\begingroup \smaller\smaller\smaller\begin{tabular}{@{}c@{}}%
0\\0\\0
\end{tabular}\endgroup%
}}\right.$}%
\begingroup \smaller\smaller\smaller\begin{tabular}{@{}c@{}}%
-236280\\2280\\1560
\end{tabular}\endgroup%
\kern3pt%
\begingroup \smaller\smaller\smaller\begin{tabular}{@{}c@{}}%
2280\\-22\\-15
\end{tabular}\endgroup%
\kern3pt%
\begingroup \smaller\smaller\smaller\begin{tabular}{@{}c@{}}%
1560\\-15\\-7
\end{tabular}\endgroup%
{$\left.\llap{\phantom{%
\begingroup \smaller\smaller\smaller\begin{tabular}{@{}c@{}}%
0\\0\\0
\end{tabular}\endgroup%
}}\!\right]$}%
\EasyButWeakLineBreak%
{$\left[\!\llap{\phantom{%
\begingroup \smaller\smaller\smaller\begin{tabular}{@{}c@{}}%
0\\0\\0
\end{tabular}\endgroup%
}}\right.$}%
\begingroup \smaller\smaller\smaller\begin{tabular}{@{}c@{}}%
-1\\-105\\2
\end{tabular}\endgroup%
\HardButStrongLineBreak\kern3pt%
\begingroup \smaller\smaller\smaller\begin{tabular}{@{}c@{}}%
1\\104\\-1
\end{tabular}\endgroup%
\HardButStrongLineBreak\kern3pt%
\begingroup \smaller\smaller\smaller\begin{tabular}{@{}c@{}}%
7\\732\\-12
\end{tabular}\endgroup%
\HardButStrongLineBreak\kern3pt%
\begingroup \smaller\smaller\smaller\begin{tabular}{@{}c@{}}%
2\\210\\-5
\end{tabular}\endgroup%
\HardButStrongLineBreak\kern3pt%
\begingroup \smaller\smaller\smaller\begin{tabular}{@{}c@{}}%
-1\\-104\\0
\end{tabular}\endgroup%
{$\left.\llap{\phantom{%
\begingroup \smaller\smaller\smaller\begin{tabular}{@{}c@{}}%
0\\0\\0
\end{tabular}\endgroup%
}}\!\right]$}%

\medskip%
%
\leavevmode\llap{}%
$W_{128}$%
\qquad\llap{64} lattices, $\chi=12$%
\hfill%
$222222$%
\nopagebreak\smallskip\hrule\nopagebreak\medskip%
%
%
\leavevmode%
${L_{128.1}}$%
{} : {$1\above{1pt}{1pt}{2}{2}8\above{1pt}{1pt}{1}{1}{\cdot}1\above{1pt}{1pt}{-}{}3\above{1pt}{1pt}{1}{}9\above{1pt}{1pt}{1}{}{\cdot}1\above{1pt}{1pt}{-2}{}5\above{1pt}{1pt}{-}{}$}\spacer%
\instructions{3m,3,2}%
\EasyButWeakLineBreak%
{${8}\above{1pt}{1pt}{r}{2}{90}\above{1pt}{1pt}{b}{2}{2}\above{1pt}{1pt}{b}{2}{360}\above{1pt}{1pt}{*}{2}{12}\above{1pt}{1pt}{l}{2}{9}\above{1pt}{1pt}{}{2}$}%
\nopagebreak\par%
\nopagebreak\par\leavevmode%
{$\left[\!\llap{\phantom{%
\begingroup \smaller\smaller\smaller\begin{tabular}{@{}c@{}}%
0\\0\\0
\end{tabular}\endgroup%
}}\right.$}%
\begingroup \smaller\smaller\smaller\begin{tabular}{@{}c@{}}%
-1187640\\-8640\\6120
\end{tabular}\endgroup%
\kern3pt%
\begingroup \smaller\smaller\smaller\begin{tabular}{@{}c@{}}%
-8640\\-51\\42
\end{tabular}\endgroup%
\kern3pt%
\begingroup \smaller\smaller\smaller\begin{tabular}{@{}c@{}}%
6120\\42\\-31
\end{tabular}\endgroup%
{$\left.\llap{\phantom{%
\begingroup \smaller\smaller\smaller\begin{tabular}{@{}c@{}}%
0\\0\\0
\end{tabular}\endgroup%
}}\!\right]$}%
\EasyButWeakLineBreak%
{$\left[\!\llap{\phantom{%
\begingroup \smaller\smaller\smaller\begin{tabular}{@{}c@{}}%
0\\0\\0
\end{tabular}\endgroup%
}}\right.$}%
\begingroup \smaller\smaller\smaller\begin{tabular}{@{}c@{}}%
1\\56\\272
\end{tabular}\endgroup%
\HardButStrongLineBreak\kern3pt%
\begingroup \smaller\smaller\smaller\begin{tabular}{@{}c@{}}%
8\\465\\2205
\end{tabular}\endgroup%
\HardButStrongLineBreak\kern3pt%
\begingroup \smaller\smaller\smaller\begin{tabular}{@{}c@{}}%
1\\59\\277
\end{tabular}\endgroup%
\HardButStrongLineBreak\kern3pt%
\begingroup \smaller\smaller\smaller\begin{tabular}{@{}c@{}}%
11\\660\\3060
\end{tabular}\endgroup%
\HardButStrongLineBreak\kern3pt%
\begingroup \smaller\smaller\smaller\begin{tabular}{@{}c@{}}%
-1\\-58\\-276
\end{tabular}\endgroup%
\HardButStrongLineBreak\kern3pt%
\begingroup \smaller\smaller\smaller\begin{tabular}{@{}c@{}}%
-1\\-60\\-279
\end{tabular}\endgroup%
{$\left.\llap{\phantom{%
\begingroup \smaller\smaller\smaller\begin{tabular}{@{}c@{}}%
0\\0\\0
\end{tabular}\endgroup%
}}\!\right]$}%
%
%
\hbox{}\par\smallskip%
%
%
\leavevmode%
${L_{128.2}}$%
{} : {$1\above{1pt}{1pt}{-2}{2}8\above{1pt}{1pt}{-}{5}{\cdot}1\above{1pt}{1pt}{1}{}3\above{1pt}{1pt}{1}{}9\above{1pt}{1pt}{-}{}{\cdot}1\above{1pt}{1pt}{-2}{}5\above{1pt}{1pt}{-}{}$}\spacer%
\instructions{32\rightarrow N'_{14},3,m}%
\EasyButWeakLineBreak%
{${72}\above{1pt}{1pt}{b}{2}{10}\above{1pt}{1pt}{s}{2}{18}\above{1pt}{1pt}{l}{2}{40}\above{1pt}{1pt}{}{2}{3}\above{1pt}{1pt}{r}{2}{4}\above{1pt}{1pt}{*}{2}$}%
\nopagebreak\par%
\nopagebreak\par\leavevmode%
{$\left[\!\llap{\phantom{%
\begingroup \smaller\smaller\smaller\begin{tabular}{@{}c@{}}%
0\\0\\0
\end{tabular}\endgroup%
}}\right.$}%
\begingroup \smaller\smaller\smaller\begin{tabular}{@{}c@{}}%
-1281240\\-428040\\4680
\end{tabular}\endgroup%
\kern3pt%
\begingroup \smaller\smaller\smaller\begin{tabular}{@{}c@{}}%
-428040\\-142998\\1563
\end{tabular}\endgroup%
\kern3pt%
\begingroup \smaller\smaller\smaller\begin{tabular}{@{}c@{}}%
4680\\1563\\-17
\end{tabular}\endgroup%
{$\left.\llap{\phantom{%
\begingroup \smaller\smaller\smaller\begin{tabular}{@{}c@{}}%
0\\0\\0
\end{tabular}\endgroup%
}}\!\right]$}%
\EasyButWeakLineBreak%
{$\left[\!\llap{\phantom{%
\begingroup \smaller\smaller\smaller\begin{tabular}{@{}c@{}}%
0\\0\\0
\end{tabular}\endgroup%
}}\right.$}%
\begingroup \smaller\smaller\smaller\begin{tabular}{@{}c@{}}%
19\\-60\\-288
\end{tabular}\endgroup%
\HardButStrongLineBreak\kern3pt%
\begingroup \smaller\smaller\smaller\begin{tabular}{@{}c@{}}%
11\\-35\\-190
\end{tabular}\endgroup%
\HardButStrongLineBreak\kern3pt%
\begingroup \smaller\smaller\smaller\begin{tabular}{@{}c@{}}%
-1\\3\\0
\end{tabular}\endgroup%
\HardButStrongLineBreak\kern3pt%
\begingroup \smaller\smaller\smaller\begin{tabular}{@{}c@{}}%
-63\\200\\1040
\end{tabular}\endgroup%
\HardButStrongLineBreak\kern3pt%
\begingroup \smaller\smaller\smaller\begin{tabular}{@{}c@{}}%
-11\\35\\189
\end{tabular}\endgroup%
\HardButStrongLineBreak\kern3pt%
\begingroup \smaller\smaller\smaller\begin{tabular}{@{}c@{}}%
-5\\16\\94
\end{tabular}\endgroup%
{$\left.\llap{\phantom{%
\begingroup \smaller\smaller\smaller\begin{tabular}{@{}c@{}}%
0\\0\\0
\end{tabular}\endgroup%
}}\!\right]$}%

\medskip%
%
\leavevmode\llap{}%
$W_{129}$%
\qquad\llap{92} lattices, $\chi=24$%
\hfill%
$22222222\rtimes C_{2}$%
\nopagebreak\smallskip\hrule\nopagebreak\medskip%
%
%
\leavevmode%
${L_{129.1}}$%
{} : {$1\above{1pt}{1pt}{-2}{4}8\above{1pt}{1pt}{1}{1}{\cdot}1\above{1pt}{1pt}{2}{}3\above{1pt}{1pt}{1}{}{\cdot}1\above{1pt}{1pt}{-2}{}7\above{1pt}{1pt}{-}{}$}\spacer%
\instructions{2\rightarrow N'_{17}}%
\EasyButWeakLineBreak%
{${8}\above{1pt}{1pt}{*}{2}{12}\above{1pt}{1pt}{*}{2}{4}\above{1pt}{1pt}{l}{2}{21}\above{1pt}{1pt}{}{2}{8}\above{1pt}{1pt}{}{2}{3}\above{1pt}{1pt}{}{2}{1}\above{1pt}{1pt}{r}{2}{84}\above{1pt}{1pt}{*}{2}$}%
\nopagebreak\par%
\nopagebreak\par\leavevmode%
{$\left[\!\llap{\phantom{%
\begingroup \smaller\smaller\smaller\begin{tabular}{@{}c@{}}%
0\\0\\0
\end{tabular}\endgroup%
}}\right.$}%
\begingroup \smaller\smaller\smaller\begin{tabular}{@{}c@{}}%
-700728\\2016\\4032
\end{tabular}\endgroup%
\kern3pt%
\begingroup \smaller\smaller\smaller\begin{tabular}{@{}c@{}}%
2016\\-5\\-12
\end{tabular}\endgroup%
\kern3pt%
\begingroup \smaller\smaller\smaller\begin{tabular}{@{}c@{}}%
4032\\-12\\-23
\end{tabular}\endgroup%
{$\left.\llap{\phantom{%
\begingroup \smaller\smaller\smaller\begin{tabular}{@{}c@{}}%
0\\0\\0
\end{tabular}\endgroup%
}}\!\right]$}%
\EasyButWeakLineBreak%
{$\left[\!\llap{\phantom{%
\begingroup \smaller\smaller\smaller\begin{tabular}{@{}c@{}}%
0\\0\\0
\end{tabular}\endgroup%
}}\right.$}%
\begingroup \smaller\smaller\smaller\begin{tabular}{@{}c@{}}%
-1\\-68\\-140
\end{tabular}\endgroup%
\HardButStrongLineBreak\kern3pt%
\begingroup \smaller\smaller\smaller\begin{tabular}{@{}c@{}}%
-1\\-72\\-138
\end{tabular}\endgroup%
\HardButStrongLineBreak\kern3pt%
\begingroup \smaller\smaller\smaller\begin{tabular}{@{}c@{}}%
1\\66\\140
\end{tabular}\endgroup%
\HardButStrongLineBreak\kern3pt%
\begingroup \smaller\smaller\smaller\begin{tabular}{@{}c@{}}%
8\\546\\1113
\end{tabular}\endgroup%
\HardButStrongLineBreak\kern3pt%
\begingroup \smaller\smaller\smaller\begin{tabular}{@{}c@{}}%
3\\208\\416
\end{tabular}\endgroup%
\HardButStrongLineBreak\kern3pt%
\begingroup \smaller\smaller\smaller\begin{tabular}{@{}c@{}}%
2\\141\\276
\end{tabular}\endgroup%
\HardButStrongLineBreak\kern3pt%
\begingroup \smaller\smaller\smaller\begin{tabular}{@{}c@{}}%
1\\72\\137
\end{tabular}\endgroup%
\HardButStrongLineBreak\kern3pt%
\begingroup \smaller\smaller\smaller\begin{tabular}{@{}c@{}}%
5\\378\\672
\end{tabular}\endgroup%
{$\left.\llap{\phantom{%
\begingroup \smaller\smaller\smaller\begin{tabular}{@{}c@{}}%
0\\0\\0
\end{tabular}\endgroup%
}}\!\right]$}%
%
%
\hbox{}\par\smallskip%
%
%
\leavevmode%
${L_{129.2}}$%
{} : {$1\above{1pt}{1pt}{-2}{6}8\above{1pt}{1pt}{1}{7}{\cdot}1\above{1pt}{1pt}{2}{}3\above{1pt}{1pt}{1}{}{\cdot}1\above{1pt}{1pt}{-2}{}7\above{1pt}{1pt}{-}{}$}\spacer%
\instructions{m}%
\EasyButWeakLineBreak%
{${8}\above{1pt}{1pt}{s}{2}{12}\above{1pt}{1pt}{l}{2}{1}\above{1pt}{1pt}{}{2}{21}\above{1pt}{1pt}{r}{2}$}\relax$\,(\times2)$%
\nopagebreak\par%
\nopagebreak\par\leavevmode%
{$\left[\!\llap{\phantom{%
\begingroup \smaller\smaller\smaller\begin{tabular}{@{}c@{}}%
0\\0\\0
\end{tabular}\endgroup%
}}\right.$}%
\begingroup \smaller\smaller\smaller\begin{tabular}{@{}c@{}}%
1848\\-504\\0
\end{tabular}\endgroup%
\kern3pt%
\begingroup \smaller\smaller\smaller\begin{tabular}{@{}c@{}}%
-504\\137\\1
\end{tabular}\endgroup%
\kern3pt%
\begingroup \smaller\smaller\smaller\begin{tabular}{@{}c@{}}%
0\\1\\-2
\end{tabular}\endgroup%
{$\left.\llap{\phantom{%
\begingroup \smaller\smaller\smaller\begin{tabular}{@{}c@{}}%
0\\0\\0
\end{tabular}\endgroup%
}}\!\right]$}%
\hfil\penalty500%
{$\left[\!\llap{\phantom{%
\begingroup \smaller\smaller\smaller\begin{tabular}{@{}c@{}}%
0\\0\\0
\end{tabular}\endgroup%
}}\right.$}%
\begingroup \smaller\smaller\smaller\begin{tabular}{@{}c@{}}%
-1\\0\\0
\end{tabular}\endgroup%
\kern3pt%
\begingroup \smaller\smaller\smaller\begin{tabular}{@{}c@{}}%
0\\-1\\-1
\end{tabular}\endgroup%
\kern3pt%
\begingroup \smaller\smaller\smaller\begin{tabular}{@{}c@{}}%
0\\0\\1
\end{tabular}\endgroup%
{$\left.\llap{\phantom{%
\begingroup \smaller\smaller\smaller\begin{tabular}{@{}c@{}}%
0\\0\\0
\end{tabular}\endgroup%
}}\!\right]$}%
\EasyButWeakLineBreak%
{$\left[\!\llap{\phantom{%
\begingroup \smaller\smaller\smaller\begin{tabular}{@{}c@{}}%
0\\0\\0
\end{tabular}\endgroup%
}}\right.$}%
\begingroup \smaller\smaller\smaller\begin{tabular}{@{}c@{}}%
1\\4\\0
\end{tabular}\endgroup%
\HardButStrongLineBreak\kern3pt%
\begingroup \smaller\smaller\smaller\begin{tabular}{@{}c@{}}%
-5\\-18\\-12
\end{tabular}\endgroup%
\HardButStrongLineBreak\kern3pt%
\begingroup \smaller\smaller\smaller\begin{tabular}{@{}c@{}}%
-3\\-11\\-7
\end{tabular}\endgroup%
\HardButStrongLineBreak\kern3pt%
\begingroup \smaller\smaller\smaller\begin{tabular}{@{}c@{}}%
-17\\-63\\-42
\end{tabular}\endgroup%
{$\left.\llap{\phantom{%
\begingroup \smaller\smaller\smaller\begin{tabular}{@{}c@{}}%
0\\0\\0
\end{tabular}\endgroup%
}}\!\right]$}%
%
%
\hbox{}\par\smallskip%
%
%
\leavevmode%
${L_{129.3}}$%
{} : {$1\above{1pt}{1pt}{2}{6}8\above{1pt}{1pt}{-}{3}{\cdot}1\above{1pt}{1pt}{2}{}3\above{1pt}{1pt}{1}{}{\cdot}1\above{1pt}{1pt}{-2}{}7\above{1pt}{1pt}{-}{}$}\EasyButWeakLineBreak%
{${8}\above{1pt}{1pt}{l}{2}{3}\above{1pt}{1pt}{r}{2}{4}\above{1pt}{1pt}{*}{2}{84}\above{1pt}{1pt}{s}{2}$}\relax$\,(\times2)$%
\nopagebreak\par%
\nopagebreak\par\leavevmode%
{$\left[\!\llap{\phantom{%
\begingroup \smaller\smaller\smaller\begin{tabular}{@{}c@{}}%
0\\0\\0
\end{tabular}\endgroup%
}}\right.$}%
\begingroup \smaller\smaller\smaller\begin{tabular}{@{}c@{}}%
241752\\336\\-2184
\end{tabular}\endgroup%
\kern3pt%
\begingroup \smaller\smaller\smaller\begin{tabular}{@{}c@{}}%
336\\-1\\-2
\end{tabular}\endgroup%
\kern3pt%
\begingroup \smaller\smaller\smaller\begin{tabular}{@{}c@{}}%
-2184\\-2\\19
\end{tabular}\endgroup%
{$\left.\llap{\phantom{%
\begingroup \smaller\smaller\smaller\begin{tabular}{@{}c@{}}%
0\\0\\0
\end{tabular}\endgroup%
}}\!\right]$}%
\hfil\penalty500%
{$\left[\!\llap{\phantom{%
\begingroup \smaller\smaller\smaller\begin{tabular}{@{}c@{}}%
0\\0\\0
\end{tabular}\endgroup%
}}\right.$}%
\begingroup \smaller\smaller\smaller\begin{tabular}{@{}c@{}}%
-1\\3528\\1176
\end{tabular}\endgroup%
\kern3pt%
\begingroup \smaller\smaller\smaller\begin{tabular}{@{}c@{}}%
0\\14\\5
\end{tabular}\endgroup%
\kern3pt%
\begingroup \smaller\smaller\smaller\begin{tabular}{@{}c@{}}%
0\\-39\\-14
\end{tabular}\endgroup%
{$\left.\llap{\phantom{%
\begingroup \smaller\smaller\smaller\begin{tabular}{@{}c@{}}%
0\\0\\0
\end{tabular}\endgroup%
}}\!\right]$}%
\EasyButWeakLineBreak%
{$\left[\!\llap{\phantom{%
\begingroup \smaller\smaller\smaller\begin{tabular}{@{}c@{}}%
0\\0\\0
\end{tabular}\endgroup%
}}\right.$}%
\begingroup \smaller\smaller\smaller\begin{tabular}{@{}c@{}}%
1\\76\\120
\end{tabular}\endgroup%
\HardButStrongLineBreak\kern3pt%
\begingroup \smaller\smaller\smaller\begin{tabular}{@{}c@{}}%
1\\84\\123
\end{tabular}\endgroup%
\HardButStrongLineBreak\kern3pt%
\begingroup \smaller\smaller\smaller\begin{tabular}{@{}c@{}}%
1\\86\\124
\end{tabular}\endgroup%
\HardButStrongLineBreak\kern3pt%
\begingroup \smaller\smaller\smaller\begin{tabular}{@{}c@{}}%
1\\84\\126
\end{tabular}\endgroup%
{$\left.\llap{\phantom{%
\begingroup \smaller\smaller\smaller\begin{tabular}{@{}c@{}}%
0\\0\\0
\end{tabular}\endgroup%
}}\!\right]$}%
%
%
\hbox{}\par\smallskip%
%
%
\leavevmode%
${L_{129.4}}$%
{} : {$[1\above{1pt}{1pt}{-}{}2\above{1pt}{1pt}{1}{}]\above{1pt}{1pt}{}{6}16\above{1pt}{1pt}{1}{7}{\cdot}1\above{1pt}{1pt}{2}{}3\above{1pt}{1pt}{1}{}{\cdot}1\above{1pt}{1pt}{-2}{}7\above{1pt}{1pt}{-}{}$}\spacer%
\instructions{2}%
\EasyButWeakLineBreak%
{${2}\above{1pt}{1pt}{r}{2}{48}\above{1pt}{1pt}{s}{2}{4}\above{1pt}{1pt}{*}{2}{336}\above{1pt}{1pt}{s}{2}{8}\above{1pt}{1pt}{*}{2}{12}\above{1pt}{1pt}{*}{2}{16}\above{1pt}{1pt}{l}{2}{21}\above{1pt}{1pt}{}{2}$}%
\nopagebreak\par%
\nopagebreak\par\leavevmode%
{$\left[\!\llap{\phantom{%
\begingroup \smaller\smaller\smaller\begin{tabular}{@{}c@{}}%
0\\0\\0
\end{tabular}\endgroup%
}}\right.$}%
\begingroup \smaller\smaller\smaller\begin{tabular}{@{}c@{}}%
-380688\\2016\\2352
\end{tabular}\endgroup%
\kern3pt%
\begingroup \smaller\smaller\smaller\begin{tabular}{@{}c@{}}%
2016\\-10\\-14
\end{tabular}\endgroup%
\kern3pt%
\begingroup \smaller\smaller\smaller\begin{tabular}{@{}c@{}}%
2352\\-14\\-11
\end{tabular}\endgroup%
{$\left.\llap{\phantom{%
\begingroup \smaller\smaller\smaller\begin{tabular}{@{}c@{}}%
0\\0\\0
\end{tabular}\endgroup%
}}\!\right]$}%
\EasyButWeakLineBreak%
{$\left[\!\llap{\phantom{%
\begingroup \smaller\smaller\smaller\begin{tabular}{@{}c@{}}%
0\\0\\0
\end{tabular}\endgroup%
}}\right.$}%
\begingroup \smaller\smaller\smaller\begin{tabular}{@{}c@{}}%
1\\125\\54
\end{tabular}\endgroup%
\HardButStrongLineBreak\kern3pt%
\begingroup \smaller\smaller\smaller\begin{tabular}{@{}c@{}}%
11\\1368\\600
\end{tabular}\endgroup%
\HardButStrongLineBreak\kern3pt%
\begingroup \smaller\smaller\smaller\begin{tabular}{@{}c@{}}%
5\\620\\274
\end{tabular}\endgroup%
\HardButStrongLineBreak\kern3pt%
\begingroup \smaller\smaller\smaller\begin{tabular}{@{}c@{}}%
55\\6804\\3024
\end{tabular}\endgroup%
\HardButStrongLineBreak\kern3pt%
\begingroup \smaller\smaller\smaller\begin{tabular}{@{}c@{}}%
1\\122\\56
\end{tabular}\endgroup%
\HardButStrongLineBreak\kern3pt%
\begingroup \smaller\smaller\smaller\begin{tabular}{@{}c@{}}%
-1\\-126\\-54
\end{tabular}\endgroup%
\HardButStrongLineBreak\kern3pt%
\begingroup \smaller\smaller\smaller\begin{tabular}{@{}c@{}}%
-1\\-124\\-56
\end{tabular}\endgroup%
\HardButStrongLineBreak\kern3pt%
\begingroup \smaller\smaller\smaller\begin{tabular}{@{}c@{}}%
2\\252\\105
\end{tabular}\endgroup%
{$\left.\llap{\phantom{%
\begingroup \smaller\smaller\smaller\begin{tabular}{@{}c@{}}%
0\\0\\0
\end{tabular}\endgroup%
}}\!\right]$}%
%
%
\hbox{}\par\smallskip%
%
%
\leavevmode%
${L_{129.5}}$%
{} : {$[1\above{1pt}{1pt}{1}{}2\above{1pt}{1pt}{1}{}]\above{1pt}{1pt}{}{2}16\above{1pt}{1pt}{-}{3}{\cdot}1\above{1pt}{1pt}{2}{}3\above{1pt}{1pt}{1}{}{\cdot}1\above{1pt}{1pt}{-2}{}7\above{1pt}{1pt}{-}{}$}\spacer%
\instructions{m}%
\EasyButWeakLineBreak%
{${2}\above{1pt}{1pt}{}{2}{48}\above{1pt}{1pt}{}{2}{1}\above{1pt}{1pt}{r}{2}{336}\above{1pt}{1pt}{*}{2}{8}\above{1pt}{1pt}{l}{2}{3}\above{1pt}{1pt}{r}{2}{16}\above{1pt}{1pt}{*}{2}{84}\above{1pt}{1pt}{l}{2}$}%
\nopagebreak\par%
\nopagebreak\par\leavevmode%
{$\left[\!\llap{\phantom{%
\begingroup \smaller\smaller\smaller\begin{tabular}{@{}c@{}}%
0\\0\\0
\end{tabular}\endgroup%
}}\right.$}%
\begingroup \smaller\smaller\smaller\begin{tabular}{@{}c@{}}%
-1594320\\4368\\9408
\end{tabular}\endgroup%
\kern3pt%
\begingroup \smaller\smaller\smaller\begin{tabular}{@{}c@{}}%
4368\\-10\\-28
\end{tabular}\endgroup%
\kern3pt%
\begingroup \smaller\smaller\smaller\begin{tabular}{@{}c@{}}%
9408\\-28\\-53
\end{tabular}\endgroup%
{$\left.\llap{\phantom{%
\begingroup \smaller\smaller\smaller\begin{tabular}{@{}c@{}}%
0\\0\\0
\end{tabular}\endgroup%
}}\!\right]$}%
\EasyButWeakLineBreak%
{$\left[\!\llap{\phantom{%
\begingroup \smaller\smaller\smaller\begin{tabular}{@{}c@{}}%
0\\0\\0
\end{tabular}\endgroup%
}}\right.$}%
\begingroup \smaller\smaller\smaller\begin{tabular}{@{}c@{}}%
2\\251\\222
\end{tabular}\endgroup%
\HardButStrongLineBreak\kern3pt%
\begingroup \smaller\smaller\smaller\begin{tabular}{@{}c@{}}%
19\\2376\\2112
\end{tabular}\endgroup%
\HardButStrongLineBreak\kern3pt%
\begingroup \smaller\smaller\smaller\begin{tabular}{@{}c@{}}%
4\\499\\445
\end{tabular}\endgroup%
\HardButStrongLineBreak\kern3pt%
\begingroup \smaller\smaller\smaller\begin{tabular}{@{}c@{}}%
83\\10332\\9240
\end{tabular}\endgroup%
\HardButStrongLineBreak\kern3pt%
\begingroup \smaller\smaller\smaller\begin{tabular}{@{}c@{}}%
1\\122\\112
\end{tabular}\endgroup%
\HardButStrongLineBreak\kern3pt%
\begingroup \smaller\smaller\smaller\begin{tabular}{@{}c@{}}%
-1\\-126\\-111
\end{tabular}\endgroup%
\HardButStrongLineBreak\kern3pt%
\begingroup \smaller\smaller\smaller\begin{tabular}{@{}c@{}}%
-1\\-124\\-112
\end{tabular}\endgroup%
\HardButStrongLineBreak\kern3pt%
\begingroup \smaller\smaller\smaller\begin{tabular}{@{}c@{}}%
11\\1386\\1218
\end{tabular}\endgroup%
{$\left.\llap{\phantom{%
\begingroup \smaller\smaller\smaller\begin{tabular}{@{}c@{}}%
0\\0\\0
\end{tabular}\endgroup%
}}\!\right]$}%
%
%
\hbox{}\par\smallskip%
%
%
\leavevmode%
${L_{129.6}}$%
{} : {$[1\above{1pt}{1pt}{-}{}2\above{1pt}{1pt}{1}{}]\above{1pt}{1pt}{}{4}16\above{1pt}{1pt}{1}{1}{\cdot}1\above{1pt}{1pt}{2}{}3\above{1pt}{1pt}{1}{}{\cdot}1\above{1pt}{1pt}{-2}{}7\above{1pt}{1pt}{-}{}$}\spacer%
\instructions{m}%
\EasyButWeakLineBreak%
{${8}\above{1pt}{1pt}{*}{2}{48}\above{1pt}{1pt}{*}{2}{4}\above{1pt}{1pt}{s}{2}{336}\above{1pt}{1pt}{l}{2}{2}\above{1pt}{1pt}{}{2}{3}\above{1pt}{1pt}{}{2}{16}\above{1pt}{1pt}{}{2}{21}\above{1pt}{1pt}{r}{2}$}%
\nopagebreak\par%
\nopagebreak\par\leavevmode%
{$\left[\!\llap{\phantom{%
\begingroup \smaller\smaller\smaller\begin{tabular}{@{}c@{}}%
0\\0\\0
\end{tabular}\endgroup%
}}\right.$}%
\begingroup \smaller\smaller\smaller\begin{tabular}{@{}c@{}}%
-423024\\2352\\2352
\end{tabular}\endgroup%
\kern3pt%
\begingroup \smaller\smaller\smaller\begin{tabular}{@{}c@{}}%
2352\\-2\\-14
\end{tabular}\endgroup%
\kern3pt%
\begingroup \smaller\smaller\smaller\begin{tabular}{@{}c@{}}%
2352\\-14\\-13
\end{tabular}\endgroup%
{$\left.\llap{\phantom{%
\begingroup \smaller\smaller\smaller\begin{tabular}{@{}c@{}}%
0\\0\\0
\end{tabular}\endgroup%
}}\!\right]$}%
\EasyButWeakLineBreak%
{$\left[\!\llap{\phantom{%
\begingroup \smaller\smaller\smaller\begin{tabular}{@{}c@{}}%
0\\0\\0
\end{tabular}\endgroup%
}}\right.$}%
\begingroup \smaller\smaller\smaller\begin{tabular}{@{}c@{}}%
3\\42\\496
\end{tabular}\endgroup%
\HardButStrongLineBreak\kern3pt%
\begingroup \smaller\smaller\smaller\begin{tabular}{@{}c@{}}%
-1\\-12\\-168
\end{tabular}\endgroup%
\HardButStrongLineBreak\kern3pt%
\begingroup \smaller\smaller\smaller\begin{tabular}{@{}c@{}}%
-1\\-14\\-166
\end{tabular}\endgroup%
\HardButStrongLineBreak\kern3pt%
\begingroup \smaller\smaller\smaller\begin{tabular}{@{}c@{}}%
1\\0\\168
\end{tabular}\endgroup%
\HardButStrongLineBreak\kern3pt%
\begingroup \smaller\smaller\smaller\begin{tabular}{@{}c@{}}%
1\\13\\166
\end{tabular}\endgroup%
\HardButStrongLineBreak\kern3pt%
\begingroup \smaller\smaller\smaller\begin{tabular}{@{}c@{}}%
4\\54\\663
\end{tabular}\endgroup%
\HardButStrongLineBreak\kern3pt%
\begingroup \smaller\smaller\smaller\begin{tabular}{@{}c@{}}%
17\\232\\2816
\end{tabular}\endgroup%
\HardButStrongLineBreak\kern3pt%
\begingroup \smaller\smaller\smaller\begin{tabular}{@{}c@{}}%
26\\357\\4305
\end{tabular}\endgroup%
{$\left.\llap{\phantom{%
\begingroup \smaller\smaller\smaller\begin{tabular}{@{}c@{}}%
0\\0\\0
\end{tabular}\endgroup%
}}\!\right]$}%
%
%
\hbox{}\par\smallskip%
%
%
\leavevmode%
${L_{129.7}}$%
{} : {$[1\above{1pt}{1pt}{1}{}2\above{1pt}{1pt}{1}{}]\above{1pt}{1pt}{}{0}16\above{1pt}{1pt}{-}{5}{\cdot}1\above{1pt}{1pt}{2}{}3\above{1pt}{1pt}{1}{}{\cdot}1\above{1pt}{1pt}{-2}{}7\above{1pt}{1pt}{-}{}$}\EasyButWeakLineBreak%
{${8}\above{1pt}{1pt}{s}{2}{48}\above{1pt}{1pt}{l}{2}{1}\above{1pt}{1pt}{}{2}{336}\above{1pt}{1pt}{}{2}{2}\above{1pt}{1pt}{r}{2}{12}\above{1pt}{1pt}{s}{2}{16}\above{1pt}{1pt}{s}{2}{84}\above{1pt}{1pt}{*}{2}$}%
\nopagebreak\par%
\nopagebreak\par\leavevmode%
{$\left[\!\llap{\phantom{%
\begingroup \smaller\smaller\smaller\begin{tabular}{@{}c@{}}%
0\\0\\0
\end{tabular}\endgroup%
}}\right.$}%
\begingroup \smaller\smaller\smaller\begin{tabular}{@{}c@{}}%
-72240\\0\\1344
\end{tabular}\endgroup%
\kern3pt%
\begingroup \smaller\smaller\smaller\begin{tabular}{@{}c@{}}%
0\\2\\0
\end{tabular}\endgroup%
\kern3pt%
\begingroup \smaller\smaller\smaller\begin{tabular}{@{}c@{}}%
1344\\0\\-25
\end{tabular}\endgroup%
{$\left.\llap{\phantom{%
\begingroup \smaller\smaller\smaller\begin{tabular}{@{}c@{}}%
0\\0\\0
\end{tabular}\endgroup%
}}\!\right]$}%
\EasyButWeakLineBreak%
{$\left[\!\llap{\phantom{%
\begingroup \smaller\smaller\smaller\begin{tabular}{@{}c@{}}%
0\\0\\0
\end{tabular}\endgroup%
}}\right.$}%
\begingroup \smaller\smaller\smaller\begin{tabular}{@{}c@{}}%
1\\-6\\52
\end{tabular}\endgroup%
\HardButStrongLineBreak\kern3pt%
\begingroup \smaller\smaller\smaller\begin{tabular}{@{}c@{}}%
5\\-12\\264
\end{tabular}\endgroup%
\HardButStrongLineBreak\kern3pt%
\begingroup \smaller\smaller\smaller\begin{tabular}{@{}c@{}}%
1\\-1\\53
\end{tabular}\endgroup%
\HardButStrongLineBreak\kern3pt%
\begingroup \smaller\smaller\smaller\begin{tabular}{@{}c@{}}%
19\\0\\1008
\end{tabular}\endgroup%
\HardButStrongLineBreak\kern3pt%
\begingroup \smaller\smaller\smaller\begin{tabular}{@{}c@{}}%
0\\1\\0
\end{tabular}\endgroup%
\HardButStrongLineBreak\kern3pt%
\begingroup \smaller\smaller\smaller\begin{tabular}{@{}c@{}}%
-1\\0\\-54
\end{tabular}\endgroup%
\HardButStrongLineBreak\kern3pt%
\begingroup \smaller\smaller\smaller\begin{tabular}{@{}c@{}}%
-1\\-8\\-56
\end{tabular}\endgroup%
\HardButStrongLineBreak\kern3pt%
\begingroup \smaller\smaller\smaller\begin{tabular}{@{}c@{}}%
1\\-42\\42
\end{tabular}\endgroup%
{$\left.\llap{\phantom{%
\begingroup \smaller\smaller\smaller\begin{tabular}{@{}c@{}}%
0\\0\\0
\end{tabular}\endgroup%
}}\!\right]$}%

\medskip%
%
\leavevmode\llap{}%
$W_{130}$%
\qquad\llap{64} lattices, $\chi=18$%
\hfill%
$2222222$%
\nopagebreak\smallskip\hrule\nopagebreak\medskip%
%
%
\leavevmode%
${L_{130.1}}$%
{} : {$1\above{1pt}{1pt}{2}{6}8\above{1pt}{1pt}{-}{3}{\cdot}1\above{1pt}{1pt}{-}{}3\above{1pt}{1pt}{-}{}9\above{1pt}{1pt}{1}{}{\cdot}1\above{1pt}{1pt}{2}{}7\above{1pt}{1pt}{1}{}$}\spacer%
\instructions{3,2}%
\EasyButWeakLineBreak%
{${6}\above{1pt}{1pt}{b}{2}{14}\above{1pt}{1pt}{l}{2}{24}\above{1pt}{1pt}{}{2}{63}\above{1pt}{1pt}{r}{2}{8}\above{1pt}{1pt}{s}{2}{36}\above{1pt}{1pt}{*}{2}{56}\above{1pt}{1pt}{b}{2}$}%
\nopagebreak\par%
\nopagebreak\par\leavevmode%
{$\left[\!\llap{\phantom{%
\begingroup \smaller\smaller\smaller\begin{tabular}{@{}c@{}}%
0\\0\\0
\end{tabular}\endgroup%
}}\right.$}%
\begingroup \smaller\smaller\smaller\begin{tabular}{@{}c@{}}%
-94248\\-16128\\1008
\end{tabular}\endgroup%
\kern3pt%
\begingroup \smaller\smaller\smaller\begin{tabular}{@{}c@{}}%
-16128\\-2757\\171
\end{tabular}\endgroup%
\kern3pt%
\begingroup \smaller\smaller\smaller\begin{tabular}{@{}c@{}}%
1008\\171\\-10
\end{tabular}\endgroup%
{$\left.\llap{\phantom{%
\begingroup \smaller\smaller\smaller\begin{tabular}{@{}c@{}}%
0\\0\\0
\end{tabular}\endgroup%
}}\!\right]$}%
\EasyButWeakLineBreak%
{$\left[\!\llap{\phantom{%
\begingroup \smaller\smaller\smaller\begin{tabular}{@{}c@{}}%
0\\0\\0
\end{tabular}\endgroup%
}}\right.$}%
\begingroup \smaller\smaller\smaller\begin{tabular}{@{}c@{}}%
-3\\20\\39
\end{tabular}\endgroup%
\HardButStrongLineBreak\kern3pt%
\begingroup \smaller\smaller\smaller\begin{tabular}{@{}c@{}}%
-2\\14\\35
\end{tabular}\endgroup%
\HardButStrongLineBreak\kern3pt%
\begingroup \smaller\smaller\smaller\begin{tabular}{@{}c@{}}%
5\\-32\\-48
\end{tabular}\endgroup%
\HardButStrongLineBreak\kern3pt%
\begingroup \smaller\smaller\smaller\begin{tabular}{@{}c@{}}%
16\\-105\\-189
\end{tabular}\endgroup%
\HardButStrongLineBreak\kern3pt%
\begingroup \smaller\smaller\smaller\begin{tabular}{@{}c@{}}%
3\\-20\\-40
\end{tabular}\endgroup%
\HardButStrongLineBreak\kern3pt%
\begingroup \smaller\smaller\smaller\begin{tabular}{@{}c@{}}%
-1\\6\\0
\end{tabular}\endgroup%
\HardButStrongLineBreak\kern3pt%
\begingroup \smaller\smaller\smaller\begin{tabular}{@{}c@{}}%
-17\\112\\196
\end{tabular}\endgroup%
{$\left.\llap{\phantom{%
\begingroup \smaller\smaller\smaller\begin{tabular}{@{}c@{}}%
0\\0\\0
\end{tabular}\endgroup%
}}\!\right]$}%
%
%
\hbox{}\par\smallskip%
%
%
\leavevmode%
${L_{130.2}}$%
{} : {$1\above{1pt}{1pt}{-2}{6}8\above{1pt}{1pt}{1}{7}{\cdot}1\above{1pt}{1pt}{-}{}3\above{1pt}{1pt}{-}{}9\above{1pt}{1pt}{1}{}{\cdot}1\above{1pt}{1pt}{2}{}7\above{1pt}{1pt}{1}{}$}\spacer%
\instructions{32\rightarrow N'_{18},3m,3,m}%
\EasyButWeakLineBreak%
{${6}\above{1pt}{1pt}{s}{2}{14}\above{1pt}{1pt}{b}{2}{24}\above{1pt}{1pt}{*}{2}{252}\above{1pt}{1pt}{s}{2}{8}\above{1pt}{1pt}{l}{2}{9}\above{1pt}{1pt}{}{2}{56}\above{1pt}{1pt}{r}{2}$}%
\nopagebreak\par%
\nopagebreak\par\leavevmode%
{$\left[\!\llap{\phantom{%
\begingroup \smaller\smaller\smaller\begin{tabular}{@{}c@{}}%
0\\0\\0
\end{tabular}\endgroup%
}}\right.$}%
\begingroup \smaller\smaller\smaller\begin{tabular}{@{}c@{}}%
915768\\504\\-7560
\end{tabular}\endgroup%
\kern3pt%
\begingroup \smaller\smaller\smaller\begin{tabular}{@{}c@{}}%
504\\-3\\-3
\end{tabular}\endgroup%
\kern3pt%
\begingroup \smaller\smaller\smaller\begin{tabular}{@{}c@{}}%
-7560\\-3\\62
\end{tabular}\endgroup%
{$\left.\llap{\phantom{%
\begingroup \smaller\smaller\smaller\begin{tabular}{@{}c@{}}%
0\\0\\0
\end{tabular}\endgroup%
}}\!\right]$}%
\EasyButWeakLineBreak%
{$\left[\!\llap{\phantom{%
\begingroup \smaller\smaller\smaller\begin{tabular}{@{}c@{}}%
0\\0\\0
\end{tabular}\endgroup%
}}\right.$}%
\begingroup \smaller\smaller\smaller\begin{tabular}{@{}c@{}}%
1\\40\\123
\end{tabular}\endgroup%
\HardButStrongLineBreak\kern3pt%
\begingroup \smaller\smaller\smaller\begin{tabular}{@{}c@{}}%
3\\126\\371
\end{tabular}\endgroup%
\HardButStrongLineBreak\kern3pt%
\begingroup \smaller\smaller\smaller\begin{tabular}{@{}c@{}}%
3\\128\\372
\end{tabular}\endgroup%
\HardButStrongLineBreak\kern3pt%
\begingroup \smaller\smaller\smaller\begin{tabular}{@{}c@{}}%
1\\42\\126
\end{tabular}\endgroup%
\HardButStrongLineBreak\kern3pt%
\begingroup \smaller\smaller\smaller\begin{tabular}{@{}c@{}}%
-1\\-44\\-124
\end{tabular}\endgroup%
\HardButStrongLineBreak\kern3pt%
\begingroup \smaller\smaller\smaller\begin{tabular}{@{}c@{}}%
-1\\-51\\-126
\end{tabular}\endgroup%
\HardButStrongLineBreak\kern3pt%
\begingroup \smaller\smaller\smaller\begin{tabular}{@{}c@{}}%
1\\0\\112
\end{tabular}\endgroup%
{$\left.\llap{\phantom{%
\begingroup \smaller\smaller\smaller\begin{tabular}{@{}c@{}}%
0\\0\\0
\end{tabular}\endgroup%
}}\!\right]$}%

\medskip%
%
\leavevmode\llap{}%
$W_{131}$%
\qquad\llap{64} lattices, $\chi=36$%
\hfill%
$2222222222\rtimes C_{2}$%
\nopagebreak\smallskip\hrule\nopagebreak\medskip%
%
%
\leavevmode%
${L_{131.1}}$%
{} : {$1\above{1pt}{1pt}{2}{6}8\above{1pt}{1pt}{1}{7}{\cdot}1\above{1pt}{1pt}{1}{}3\above{1pt}{1pt}{-}{}9\above{1pt}{1pt}{1}{}{\cdot}1\above{1pt}{1pt}{-2}{}11\above{1pt}{1pt}{-}{}$}\spacer%
\instructions{23\rightarrow N'_{20},3m,3,2,m}%
\EasyButWeakLineBreak%
{${24}\above{1pt}{1pt}{b}{2}{198}\above{1pt}{1pt}{s}{2}{6}\above{1pt}{1pt}{b}{2}{792}\above{1pt}{1pt}{*}{2}{4}\above{1pt}{1pt}{*}{2}{24}\above{1pt}{1pt}{b}{2}{22}\above{1pt}{1pt}{s}{2}{6}\above{1pt}{1pt}{b}{2}{88}\above{1pt}{1pt}{*}{2}{36}\above{1pt}{1pt}{*}{2}$}%
\nopagebreak\par%
\nopagebreak\par\leavevmode%
{$\left[\!\llap{\phantom{%
\begingroup \smaller\smaller\smaller
\endgroup%
}}\!\right]$}%
%
%
\hbox{}\par\smallskip%
%
%
\leavevmode%
${L_{131.2}}$%
{} : {$1\above{1pt}{1pt}{-2}{2}16\above{1pt}{1pt}{-}{3}{\cdot}1\above{1pt}{1pt}{-}{}3\above{1pt}{1pt}{1}{}9\above{1pt}{1pt}{-}{}{\cdot}1\above{1pt}{1pt}{-2}{}11\above{1pt}{1pt}{1}{}$}\spacer%
\instructions{3m,3,m}%
\EasyButWeakLineBreak%
{${48}\above{1pt}{1pt}{}{2}{99}\above{1pt}{1pt}{r}{2}{12}\above{1pt}{1pt}{*}{2}{1584}\above{1pt}{1pt}{b}{2}{2}\above{1pt}{1pt}{b}{2}{48}\above{1pt}{1pt}{*}{2}{44}\above{1pt}{1pt}{l}{2}{3}\above{1pt}{1pt}{}{2}{176}\above{1pt}{1pt}{r}{2}{18}\above{1pt}{1pt}{l}{2}$}%
\nopagebreak\par%
shares genus with 3-dual\nopagebreak\par%
\nopagebreak\par\leavevmode%
{$\left[\!\llap{\phantom{%
\begingroup \smaller\smaller\smaller
\endgroup%
}}\!\right]$}%

\medskip%
%
\leavevmode\llap{}%
$W_{132}$%
\qquad\llap{32} lattices, $\chi=48$%
\hfill%
$222222222222\rtimes C_{2}$%
\nopagebreak\smallskip\hrule\nopagebreak\medskip%
%
%
\leavevmode%
${L_{132.1}}$%
{} : {$1\above{1pt}{1pt}{2}{2}8\above{1pt}{1pt}{-}{5}{\cdot}1\above{1pt}{1pt}{-2}{}5\above{1pt}{1pt}{-}{}{\cdot}1\above{1pt}{1pt}{-2}{}7\above{1pt}{1pt}{-}{}$}\spacer%
\instructions{2\rightarrow N'_{23}}%
\EasyButWeakLineBreak%
{${10}\above{1pt}{1pt}{b}{2}{2}\above{1pt}{1pt}{l}{2}{40}\above{1pt}{1pt}{}{2}{1}\above{1pt}{1pt}{r}{2}{140}\above{1pt}{1pt}{*}{2}{8}\above{1pt}{1pt}{b}{2}$}\relax$\,(\times2)$%
\nopagebreak\par%
\nopagebreak\par\leavevmode%
{$\left[\!\llap{\phantom{%
\begingroup \smaller\smaller\smaller
\endgroup%
}}\!\right]$}%
%
%
\hbox{}\par\smallskip%
%
%
\leavevmode%
${L_{132.2}}$%
{} : {$1\above{1pt}{1pt}{-2}{2}8\above{1pt}{1pt}{1}{1}{\cdot}1\above{1pt}{1pt}{-2}{}5\above{1pt}{1pt}{-}{}{\cdot}1\above{1pt}{1pt}{-2}{}7\above{1pt}{1pt}{-}{}$}\spacer%
\instructions{m}%
\EasyButWeakLineBreak%
{${10}\above{1pt}{1pt}{s}{2}{2}\above{1pt}{1pt}{b}{2}{40}\above{1pt}{1pt}{*}{2}{4}\above{1pt}{1pt}{l}{2}{35}\above{1pt}{1pt}{}{2}{8}\above{1pt}{1pt}{r}{2}$}\relax$\,(\times2)$%
\nopagebreak\par%
\nopagebreak\par\leavevmode%
{$\left[\!\llap{\phantom{%
\begingroup \smaller\smaller\smaller
\endgroup%
}}\!\right]$}%

\medskip%
%
\leavevmode\llap{}%
$W_{133}$%
\qquad\llap{32} lattices, $\chi=54$%
\hfill%
$422222422222\rtimes C_{2}$%
\nopagebreak\smallskip\hrule\nopagebreak\medskip%
%
%
\leavevmode%
${L_{133.1}}$%
{} : {$1\above{1pt}{1pt}{2}{2}8\above{1pt}{1pt}{-}{5}{\cdot}1\above{1pt}{1pt}{2}{}5\above{1pt}{1pt}{1}{}{\cdot}1\above{1pt}{1pt}{2}{}7\above{1pt}{1pt}{1}{}$}\spacer%
\instructions{2\rightarrow N'_{24}}%
\EasyButWeakLineBreak%
{${2}\above{1pt}{1pt}{}{4}{1}\above{1pt}{1pt}{r}{2}{56}\above{1pt}{1pt}{l}{2}{5}\above{1pt}{1pt}{r}{2}{28}\above{1pt}{1pt}{*}{2}{8}\above{1pt}{1pt}{b}{2}$}\relax$\,(\times2)$%
\nopagebreak\par%
\nopagebreak\par\leavevmode%
{$\left[\!\llap{\phantom{%
\begingroup \smaller\smaller\smaller
\endgroup%
}}\!\right]$}%
%
%
\hbox{}\par\smallskip%
%
%
\leavevmode%
${L_{133.2}}$%
{} : {$1\above{1pt}{1pt}{-2}{2}8\above{1pt}{1pt}{1}{1}{\cdot}1\above{1pt}{1pt}{2}{}5\above{1pt}{1pt}{1}{}{\cdot}1\above{1pt}{1pt}{2}{}7\above{1pt}{1pt}{1}{}$}\spacer%
\instructions{m}%
\EasyButWeakLineBreak%
{${2}\above{1pt}{1pt}{*}{4}{4}\above{1pt}{1pt}{s}{2}{56}\above{1pt}{1pt}{s}{2}{20}\above{1pt}{1pt}{l}{2}{7}\above{1pt}{1pt}{}{2}{8}\above{1pt}{1pt}{r}{2}$}\relax$\,(\times2)$%
\nopagebreak\par%
\nopagebreak\par\leavevmode%
{$\left[\!\llap{\phantom{%
\begingroup \smaller\smaller\smaller
\endgroup%
}}\!\right]$}%

\medskip%
%
\leavevmode\llap{}%
$W_{134}$%
\qquad\llap{32} lattices, $\chi=21$%
\hfill%
$2222224$%
\nopagebreak\smallskip\hrule\nopagebreak\medskip%
%
%
\leavevmode%
${L_{134.1}}$%
{} : {$1\above{1pt}{1pt}{2}{2}8\above{1pt}{1pt}{1}{1}{\cdot}1\above{1pt}{1pt}{2}{}3\above{1pt}{1pt}{1}{}{\cdot}1\above{1pt}{1pt}{2}{}13\above{1pt}{1pt}{-}{}$}\spacer%
\instructions{2\rightarrow N'_{27}}%
\EasyButWeakLineBreak%
{${2}\above{1pt}{1pt}{b}{2}{26}\above{1pt}{1pt}{l}{2}{8}\above{1pt}{1pt}{r}{2}{2}\above{1pt}{1pt}{b}{2}{104}\above{1pt}{1pt}{*}{2}{12}\above{1pt}{1pt}{l}{2}{1}\above{1pt}{1pt}{}{4}$}%
\nopagebreak\par%
\nopagebreak\par\leavevmode%
{$\left[\!\llap{\phantom{%
\begingroup \smaller\smaller\smaller\begin{tabular}{@{}c@{}}%
0\\0\\0
\end{tabular}\endgroup%
}}\right.$}%
\begingroup \smaller\smaller\smaller\begin{tabular}{@{}c@{}}%
-1962168\\6864\\9672
\end{tabular}\endgroup%
\kern3pt%
\begingroup \smaller\smaller\smaller\begin{tabular}{@{}c@{}}%
6864\\-22\\-35
\end{tabular}\endgroup%
\kern3pt%
\begingroup \smaller\smaller\smaller\begin{tabular}{@{}c@{}}%
9672\\-35\\-47
\end{tabular}\endgroup%
{$\left.\llap{\phantom{%
\begingroup \smaller\smaller\smaller\begin{tabular}{@{}c@{}}%
0\\0\\0
\end{tabular}\endgroup%
}}\!\right]$}%
\EasyButWeakLineBreak%
{$\left[\!\llap{\phantom{%
\begingroup \smaller\smaller\smaller\begin{tabular}{@{}c@{}}%
0\\0\\0
\end{tabular}\endgroup%
}}\right.$}%
\begingroup \smaller\smaller\smaller\begin{tabular}{@{}c@{}}%
-1\\-83\\-144
\end{tabular}\endgroup%
\HardButStrongLineBreak\kern3pt%
\begingroup \smaller\smaller\smaller\begin{tabular}{@{}c@{}}%
-2\\-169\\-286
\end{tabular}\endgroup%
\HardButStrongLineBreak\kern3pt%
\begingroup \smaller\smaller\smaller\begin{tabular}{@{}c@{}}%
3\\248\\432
\end{tabular}\endgroup%
\HardButStrongLineBreak\kern3pt%
\begingroup \smaller\smaller\smaller\begin{tabular}{@{}c@{}}%
6\\499\\862
\end{tabular}\endgroup%
\HardButStrongLineBreak\kern3pt%
\begingroup \smaller\smaller\smaller\begin{tabular}{@{}c@{}}%
113\\9412\\16224
\end{tabular}\endgroup%
\HardButStrongLineBreak\kern3pt%
\begingroup \smaller\smaller\smaller\begin{tabular}{@{}c@{}}%
11\\918\\1578
\end{tabular}\endgroup%
\HardButStrongLineBreak\kern3pt%
\begingroup \smaller\smaller\smaller\begin{tabular}{@{}c@{}}%
1\\84\\143
\end{tabular}\endgroup%
{$\left.\llap{\phantom{%
\begingroup \smaller\smaller\smaller\begin{tabular}{@{}c@{}}%
0\\0\\0
\end{tabular}\endgroup%
}}\!\right]$}%
%
%
\hbox{}\par\smallskip%
%
%
\leavevmode%
${L_{134.2}}$%
{} : {$1\above{1pt}{1pt}{-2}{2}8\above{1pt}{1pt}{-}{5}{\cdot}1\above{1pt}{1pt}{2}{}3\above{1pt}{1pt}{1}{}{\cdot}1\above{1pt}{1pt}{2}{}13\above{1pt}{1pt}{-}{}$}\spacer%
\instructions{m}%
\EasyButWeakLineBreak%
{${2}\above{1pt}{1pt}{s}{2}{26}\above{1pt}{1pt}{b}{2}{8}\above{1pt}{1pt}{b}{2}{2}\above{1pt}{1pt}{l}{2}{104}\above{1pt}{1pt}{}{2}{3}\above{1pt}{1pt}{r}{2}{4}\above{1pt}{1pt}{*}{4}$}%
\nopagebreak\par%
\nopagebreak\par\leavevmode%
{$\left[\!\llap{\phantom{%
\begingroup \smaller\smaller\smaller\begin{tabular}{@{}c@{}}%
0\\0\\0
\end{tabular}\endgroup%
}}\right.$}%
\begingroup \smaller\smaller\smaller\begin{tabular}{@{}c@{}}%
-8293272\\29016\\14664
\end{tabular}\endgroup%
\kern3pt%
\begingroup \smaller\smaller\smaller\begin{tabular}{@{}c@{}}%
29016\\-101\\-52
\end{tabular}\endgroup%
\kern3pt%
\begingroup \smaller\smaller\smaller\begin{tabular}{@{}c@{}}%
14664\\-52\\-25
\end{tabular}\endgroup%
{$\left.\llap{\phantom{%
\begingroup \smaller\smaller\smaller\begin{tabular}{@{}c@{}}%
0\\0\\0
\end{tabular}\endgroup%
}}\!\right]$}%
\EasyButWeakLineBreak%
{$\left[\!\llap{\phantom{%
\begingroup \smaller\smaller\smaller\begin{tabular}{@{}c@{}}%
0\\0\\0
\end{tabular}\endgroup%
}}\right.$}%
\begingroup \smaller\smaller\smaller\begin{tabular}{@{}c@{}}%
3\\621\\467
\end{tabular}\endgroup%
\HardButStrongLineBreak\kern3pt%
\begingroup \smaller\smaller\smaller\begin{tabular}{@{}c@{}}%
17\\3523\\2639
\end{tabular}\endgroup%
\HardButStrongLineBreak\kern3pt%
\begingroup \smaller\smaller\smaller\begin{tabular}{@{}c@{}}%
7\\1452\\1084
\end{tabular}\endgroup%
\HardButStrongLineBreak\kern3pt%
\begingroup \smaller\smaller\smaller\begin{tabular}{@{}c@{}}%
3\\623\\463
\end{tabular}\endgroup%
\HardButStrongLineBreak\kern3pt%
\begingroup \smaller\smaller\smaller\begin{tabular}{@{}c@{}}%
21\\4368\\3224
\end{tabular}\endgroup%
\HardButStrongLineBreak\kern3pt%
\begingroup \smaller\smaller\smaller\begin{tabular}{@{}c@{}}%
-1\\-207\\-156
\end{tabular}\endgroup%
\HardButStrongLineBreak\kern3pt%
\begingroup \smaller\smaller\smaller\begin{tabular}{@{}c@{}}%
-1\\-208\\-154
\end{tabular}\endgroup%
{$\left.\llap{\phantom{%
\begingroup \smaller\smaller\smaller\begin{tabular}{@{}c@{}}%
0\\0\\0
\end{tabular}\endgroup%
}}\!\right]$}%

\medskip%
%
\leavevmode\llap{}%
$W_{135}$%
\qquad\llap{32} lattices, $\chi=72$%
\hfill%
$2222222222222222\rtimes C_{2}$%
\nopagebreak\smallskip\hrule\nopagebreak\medskip%
%
%
\leavevmode%
${L_{135.1}}$%
{} : {$1\above{1pt}{1pt}{-2}{6}8\above{1pt}{1pt}{1}{7}{\cdot}1\above{1pt}{1pt}{2}{}3\above{1pt}{1pt}{-}{}{\cdot}1\above{1pt}{1pt}{-2}{}23\above{1pt}{1pt}{1}{}$}\spacer%
\instructions{2\rightarrow N'_{30}}%
\EasyButWeakLineBreak%
{${6}\above{1pt}{1pt}{s}{2}{46}\above{1pt}{1pt}{b}{2}{24}\above{1pt}{1pt}{*}{2}{92}\above{1pt}{1pt}{s}{2}{8}\above{1pt}{1pt}{l}{2}{69}\above{1pt}{1pt}{}{2}{1}\above{1pt}{1pt}{}{2}{184}\above{1pt}{1pt}{r}{2}$}\relax$\,(\times2)$%
\nopagebreak\par%
\nopagebreak\par\leavevmode%
{$\left[\!\llap{\phantom{%
\begingroup \smaller\smaller\smaller
\endgroup%
}}\!\right]$}%
%
%
\hbox{}\par\smallskip%
%
%
\leavevmode%
${L_{135.2}}$%
{} : {$1\above{1pt}{1pt}{2}{6}8\above{1pt}{1pt}{-}{3}{\cdot}1\above{1pt}{1pt}{2}{}3\above{1pt}{1pt}{-}{}{\cdot}1\above{1pt}{1pt}{-2}{}23\above{1pt}{1pt}{1}{}$}\spacer%
\instructions{m}%
\EasyButWeakLineBreak%
{${6}\above{1pt}{1pt}{b}{2}{46}\above{1pt}{1pt}{l}{2}{24}\above{1pt}{1pt}{}{2}{23}\above{1pt}{1pt}{r}{2}{8}\above{1pt}{1pt}{s}{2}{276}\above{1pt}{1pt}{*}{2}{4}\above{1pt}{1pt}{*}{2}{184}\above{1pt}{1pt}{b}{2}$}\relax$\,(\times2)$%
\nopagebreak\par%
\nopagebreak\par\leavevmode%
{$\left[\!\llap{\phantom{%
\begingroup \smaller\smaller\smaller
\endgroup%
}}\!\right]$}%

\medskip%
%
\leavevmode\llap{}%
$W_{136}$%
\qquad\llap{64} lattices, $\chi=18$%
\hfill%
$2222222$%
\nopagebreak\smallskip\hrule\nopagebreak\medskip%
%
%
\leavevmode%
${L_{136.1}}$%
{} : {$1\above{1pt}{1pt}{-2}{6}8\above{1pt}{1pt}{-}{3}{\cdot}1\above{1pt}{1pt}{2}{}3\above{1pt}{1pt}{-}{}{\cdot}1\above{1pt}{1pt}{-2}{}5\above{1pt}{1pt}{1}{}{\cdot}1\above{1pt}{1pt}{2}{}7\above{1pt}{1pt}{-}{}$}\spacer%
\instructions{2\rightarrow N'_{31}}%
\EasyButWeakLineBreak%
{${24}\above{1pt}{1pt}{}{2}{5}\above{1pt}{1pt}{r}{2}{168}\above{1pt}{1pt}{l}{2}{1}\above{1pt}{1pt}{}{2}{280}\above{1pt}{1pt}{r}{2}{6}\above{1pt}{1pt}{b}{2}{70}\above{1pt}{1pt}{l}{2}$}%
\nopagebreak\par%
\nopagebreak\par\leavevmode%
{$\left[\!\llap{\phantom{%
\begingroup \smaller\smaller\smaller\begin{tabular}{@{}c@{}}%
0\\0\\0
\end{tabular}\endgroup%
}}\right.$}%
\begingroup \smaller\smaller\smaller\begin{tabular}{@{}c@{}}%
472920\\157080\\-840
\end{tabular}\endgroup%
\kern3pt%
\begingroup \smaller\smaller\smaller\begin{tabular}{@{}c@{}}%
157080\\52174\\-279
\end{tabular}\endgroup%
\kern3pt%
\begingroup \smaller\smaller\smaller\begin{tabular}{@{}c@{}}%
-840\\-279\\1
\end{tabular}\endgroup%
{$\left.\llap{\phantom{%
\begingroup \smaller\smaller\smaller\begin{tabular}{@{}c@{}}%
0\\0\\0
\end{tabular}\endgroup%
}}\!\right]$}%
\EasyButWeakLineBreak%
{$\left[\!\llap{\phantom{%
\begingroup \smaller\smaller\smaller\begin{tabular}{@{}c@{}}%
0\\0\\0
\end{tabular}\endgroup%
}}\right.$}%
\begingroup \smaller\smaller\smaller\begin{tabular}{@{}c@{}}%
247\\-744\\-72
\end{tabular}\endgroup%
\HardButStrongLineBreak\kern3pt%
\begingroup \smaller\smaller\smaller\begin{tabular}{@{}c@{}}%
83\\-250\\-25
\end{tabular}\endgroup%
\HardButStrongLineBreak\kern3pt%
\begingroup \smaller\smaller\smaller\begin{tabular}{@{}c@{}}%
251\\-756\\-84
\end{tabular}\endgroup%
\HardButStrongLineBreak\kern3pt%
\begingroup \smaller\smaller\smaller\begin{tabular}{@{}c@{}}%
0\\0\\-1
\end{tabular}\endgroup%
\HardButStrongLineBreak\kern3pt%
\begingroup \smaller\smaller\smaller\begin{tabular}{@{}c@{}}%
-93\\280\\0
\end{tabular}\endgroup%
\HardButStrongLineBreak\kern3pt%
\begingroup \smaller\smaller\smaller\begin{tabular}{@{}c@{}}%
-1\\3\\0
\end{tabular}\endgroup%
\HardButStrongLineBreak\kern3pt%
\begingroup \smaller\smaller\smaller\begin{tabular}{@{}c@{}}%
244\\-735\\-70
\end{tabular}\endgroup%
{$\left.\llap{\phantom{%
\begingroup \smaller\smaller\smaller\begin{tabular}{@{}c@{}}%
0\\0\\0
\end{tabular}\endgroup%
}}\!\right]$}%
%
%
\hbox{}\par\smallskip%
%
%
\leavevmode%
${L_{136.2}}$%
{} : {$1\above{1pt}{1pt}{2}{6}8\above{1pt}{1pt}{1}{7}{\cdot}1\above{1pt}{1pt}{2}{}3\above{1pt}{1pt}{-}{}{\cdot}1\above{1pt}{1pt}{-2}{}5\above{1pt}{1pt}{1}{}{\cdot}1\above{1pt}{1pt}{2}{}7\above{1pt}{1pt}{-}{}$}\spacer%
\instructions{m}%
\EasyButWeakLineBreak%
{${24}\above{1pt}{1pt}{*}{2}{20}\above{1pt}{1pt}{s}{2}{168}\above{1pt}{1pt}{s}{2}{4}\above{1pt}{1pt}{*}{2}{280}\above{1pt}{1pt}{b}{2}{6}\above{1pt}{1pt}{s}{2}{70}\above{1pt}{1pt}{b}{2}$}%
\nopagebreak\par%
\nopagebreak\par\leavevmode%
{$\left[\!\llap{\phantom{%
\begingroup \smaller\smaller\smaller\begin{tabular}{@{}c@{}}%
0\\0\\0
\end{tabular}\endgroup%
}}\right.$}%
\begingroup \smaller\smaller\smaller\begin{tabular}{@{}c@{}}%
-867720\\288960\\3360
\end{tabular}\endgroup%
\kern3pt%
\begingroup \smaller\smaller\smaller\begin{tabular}{@{}c@{}}%
288960\\-96226\\-1119
\end{tabular}\endgroup%
\kern3pt%
\begingroup \smaller\smaller\smaller\begin{tabular}{@{}c@{}}%
3360\\-1119\\-13
\end{tabular}\endgroup%
{$\left.\llap{\phantom{%
\begingroup \smaller\smaller\smaller\begin{tabular}{@{}c@{}}%
0\\0\\0
\end{tabular}\endgroup%
}}\!\right]$}%
\EasyButWeakLineBreak%
{$\left[\!\llap{\phantom{%
\begingroup \smaller\smaller\smaller\begin{tabular}{@{}c@{}}%
0\\0\\0
\end{tabular}\endgroup%
}}\right.$}%
\begingroup \smaller\smaller\smaller\begin{tabular}{@{}c@{}}%
-13\\-36\\-264
\end{tabular}\endgroup%
\HardButStrongLineBreak\kern3pt%
\begingroup \smaller\smaller\smaller\begin{tabular}{@{}c@{}}%
-7\\-20\\-90
\end{tabular}\endgroup%
\HardButStrongLineBreak\kern3pt%
\begingroup \smaller\smaller\smaller\begin{tabular}{@{}c@{}}%
1\\0\\252
\end{tabular}\endgroup%
\HardButStrongLineBreak\kern3pt%
\begingroup \smaller\smaller\smaller\begin{tabular}{@{}c@{}}%
3\\8\\86
\end{tabular}\endgroup%
\HardButStrongLineBreak\kern3pt%
\begingroup \smaller\smaller\smaller\begin{tabular}{@{}c@{}}%
51\\140\\1120
\end{tabular}\endgroup%
\HardButStrongLineBreak\kern3pt%
\begingroup \smaller\smaller\smaller\begin{tabular}{@{}c@{}}%
1\\3\\0
\end{tabular}\endgroup%
\HardButStrongLineBreak\kern3pt%
\begingroup \smaller\smaller\smaller\begin{tabular}{@{}c@{}}%
-13\\-35\\-350
\end{tabular}\endgroup%
{$\left.\llap{\phantom{%
\begingroup \smaller\smaller\smaller\begin{tabular}{@{}c@{}}%
0\\0\\0
\end{tabular}\endgroup%
}}\!\right]$}%

\medskip%
%
\leavevmode\llap{}%
$W_{137}$%
\qquad\llap{64} lattices, $\chi=72$%
\hfill%
$2222222222222222\rtimes C_{2}$%
\nopagebreak\smallskip\hrule\nopagebreak\medskip%
%
%
\leavevmode%
${L_{137.1}}$%
{} : {$1\above{1pt}{1pt}{-2}{6}8\above{1pt}{1pt}{-}{3}{\cdot}1\above{1pt}{1pt}{2}{}3\above{1pt}{1pt}{-}{}{\cdot}1\above{1pt}{1pt}{2}{}5\above{1pt}{1pt}{-}{}{\cdot}1\above{1pt}{1pt}{-2}{}7\above{1pt}{1pt}{1}{}$}\spacer%
\instructions{2\rightarrow N'_{35}}%
\EasyButWeakLineBreak%
{${6}\above{1pt}{1pt}{s}{2}{14}\above{1pt}{1pt}{l}{2}{24}\above{1pt}{1pt}{}{2}{1}\above{1pt}{1pt}{r}{2}{40}\above{1pt}{1pt}{s}{2}{28}\above{1pt}{1pt}{*}{2}{60}\above{1pt}{1pt}{*}{2}{56}\above{1pt}{1pt}{b}{2}$}\relax$\,(\times2)$%
\nopagebreak\par%
\nopagebreak\par\leavevmode%
{$\left[\!\llap{\phantom{%
\begingroup \smaller\smaller\smaller
\endgroup%
}}\!\right]$}%
%
%
\hbox{}\par\smallskip%
%
%
\leavevmode%
${L_{137.2}}$%
{} : {$1\above{1pt}{1pt}{2}{6}8\above{1pt}{1pt}{1}{7}{\cdot}1\above{1pt}{1pt}{2}{}3\above{1pt}{1pt}{-}{}{\cdot}1\above{1pt}{1pt}{2}{}5\above{1pt}{1pt}{-}{}{\cdot}1\above{1pt}{1pt}{-2}{}7\above{1pt}{1pt}{1}{}$}\spacer%
\instructions{m}%
\EasyButWeakLineBreak%
{${6}\above{1pt}{1pt}{b}{2}{14}\above{1pt}{1pt}{b}{2}{24}\above{1pt}{1pt}{*}{2}{4}\above{1pt}{1pt}{s}{2}{40}\above{1pt}{1pt}{l}{2}{7}\above{1pt}{1pt}{}{2}{15}\above{1pt}{1pt}{}{2}{56}\above{1pt}{1pt}{r}{2}$}\relax$\,(\times2)$%
\nopagebreak\par%
\nopagebreak\par\leavevmode%
{$\left[\!\llap{\phantom{%
\begingroup \smaller\smaller\smaller
\endgroup%
}}\!\right]$}%

\medskip%
%
\leavevmode\llap{}%
$W_{138}$%
\qquad\llap{64} lattices, $\chi=72$%
\hfill%
$2222222222222222\rtimes C_{2}$%
\nopagebreak\smallskip\hrule\nopagebreak\medskip%
%
%
\leavevmode%
${L_{138.1}}$%
{} : {$1\above{1pt}{1pt}{2}{2}8\above{1pt}{1pt}{-}{5}{\cdot}1\above{1pt}{1pt}{2}{}3\above{1pt}{1pt}{-}{}{\cdot}1\above{1pt}{1pt}{-2}{}5\above{1pt}{1pt}{1}{}{\cdot}1\above{1pt}{1pt}{-2}{}13\above{1pt}{1pt}{-}{}$}\spacer%
\instructions{2\rightarrow N'_{37}}%
\EasyButWeakLineBreak%
{${2}\above{1pt}{1pt}{b}{2}{26}\above{1pt}{1pt}{b}{2}{8}\above{1pt}{1pt}{*}{2}{780}\above{1pt}{1pt}{l}{2}{1}\above{1pt}{1pt}{r}{2}{312}\above{1pt}{1pt}{l}{2}{5}\above{1pt}{1pt}{}{2}{104}\above{1pt}{1pt}{r}{2}$}\relax$\,(\times2)$%
\nopagebreak\par%
\nopagebreak\par\leavevmode%
{$\left[\!\llap{\phantom{%
\begingroup \smaller\smaller\smaller
\endgroup%
}}\!\right]$}%
%
%
\hbox{}\par\smallskip%
%
%
\leavevmode%
${L_{138.2}}$%
{} : {$1\above{1pt}{1pt}{-2}{2}8\above{1pt}{1pt}{1}{1}{\cdot}1\above{1pt}{1pt}{2}{}3\above{1pt}{1pt}{-}{}{\cdot}1\above{1pt}{1pt}{-2}{}5\above{1pt}{1pt}{1}{}{\cdot}1\above{1pt}{1pt}{-2}{}13\above{1pt}{1pt}{-}{}$}\spacer%
\instructions{m}%
\EasyButWeakLineBreak%
{${2}\above{1pt}{1pt}{s}{2}{26}\above{1pt}{1pt}{l}{2}{8}\above{1pt}{1pt}{}{2}{195}\above{1pt}{1pt}{r}{2}{4}\above{1pt}{1pt}{s}{2}{312}\above{1pt}{1pt}{s}{2}{20}\above{1pt}{1pt}{*}{2}{104}\above{1pt}{1pt}{b}{2}$}\relax$\,(\times2)$%
\nopagebreak\par%
\nopagebreak\par\leavevmode%
{$\left[\!\llap{\phantom{%
\begingroup \smaller\smaller\smaller
\endgroup%
}}\!\right]$}%

\medskip%
%
\leavevmode\llap{}%
$W_{139}$%
\qquad\llap{64} lattices, $\chi=108$%
\hfill%
$2222222222222222222222\rtimes C_{2}$%
\nopagebreak\smallskip\hrule\nopagebreak\medskip%
%
%
\leavevmode%
${L_{139.1}}$%
{} : {$1\above{1pt}{1pt}{-2}{6}8\above{1pt}{1pt}{1}{7}{\cdot}1\above{1pt}{1pt}{2}{}3\above{1pt}{1pt}{-}{}{\cdot}1\above{1pt}{1pt}{-2}{}5\above{1pt}{1pt}{-}{}{\cdot}1\above{1pt}{1pt}{2}{}19\above{1pt}{1pt}{-}{}$}\spacer%
\instructions{2\rightarrow N'_{38}}%
\EasyButWeakLineBreak%
{${1}\above{1pt}{1pt}{}{2}{285}\above{1pt}{1pt}{r}{2}{8}\above{1pt}{1pt}{s}{2}{60}\above{1pt}{1pt}{*}{2}{152}\above{1pt}{1pt}{b}{2}{6}\above{1pt}{1pt}{s}{2}{190}\above{1pt}{1pt}{b}{2}{24}\above{1pt}{1pt}{b}{2}{38}\above{1pt}{1pt}{b}{2}{6}\above{1pt}{1pt}{l}{2}{760}\above{1pt}{1pt}{}{2}$}\relax$\,(\times2)$%
\nopagebreak\par%
\nopagebreak\par\leavevmode%
{$\left[\!\llap{\phantom{%
\begingroup \smaller\smaller\smaller
\endgroup%
}}\!\right]$}%
%
%
\hbox{}\par\smallskip%
%
%
\leavevmode%
${L_{139.2}}$%
{} : {$1\above{1pt}{1pt}{2}{6}8\above{1pt}{1pt}{-}{3}{\cdot}1\above{1pt}{1pt}{2}{}3\above{1pt}{1pt}{-}{}{\cdot}1\above{1pt}{1pt}{-2}{}5\above{1pt}{1pt}{-}{}{\cdot}1\above{1pt}{1pt}{2}{}19\above{1pt}{1pt}{-}{}$}\spacer%
\instructions{m}%
\EasyButWeakLineBreak%
{${4}\above{1pt}{1pt}{*}{2}{1140}\above{1pt}{1pt}{s}{2}{8}\above{1pt}{1pt}{l}{2}{15}\above{1pt}{1pt}{}{2}{152}\above{1pt}{1pt}{r}{2}{6}\above{1pt}{1pt}{b}{2}{190}\above{1pt}{1pt}{l}{2}{24}\above{1pt}{1pt}{r}{2}{38}\above{1pt}{1pt}{s}{2}{6}\above{1pt}{1pt}{b}{2}{760}\above{1pt}{1pt}{*}{2}$}\relax$\,(\times2)$%
\nopagebreak\par%
\nopagebreak\par\leavevmode%
{$\left[\!\llap{\phantom{%
\begingroup \smaller\smaller\smaller
\endgroup%
}}\!\right]$}%

\medskip%
%
\leavevmode\llap{}%
$W_{140}$%
\qquad\llap{4} lattices, $\chi=6$%
\hfill%
$\infty\slashtwo\infty|\rtimes D_{2}$%
\nopagebreak\smallskip\hrule\nopagebreak\medskip%
%
%
\leavevmode%
${L_{140.1}}$%
{} : {$1\above{1pt}{1pt}{2}{0}4\above{1pt}{1pt}{1}{1}$}\EasyButWeakLineBreak%
{${4}\above{1pt}{1pt}{1,0}{\infty a}{4}\above{1pt}{1pt}{}{2}{1}\above{1pt}{1pt}{4,1}{\infty}$}%
\nopagebreak\par%
\nopagebreak\par\leavevmode%
{$\left[\!\llap{\phantom{%
\begingroup \smaller\smaller\smaller\begin{tabular}{@{}c@{}}%
0\\0\\0
\end{tabular}\endgroup%
}}\right.$}%
\begingroup \smaller\smaller\smaller\begin{tabular}{@{}c@{}}%
-44\\-16\\12
\end{tabular}\endgroup%
\kern3pt%
\begingroup \smaller\smaller\smaller\begin{tabular}{@{}c@{}}%
-16\\-5\\4
\end{tabular}\endgroup%
\kern3pt%
\begingroup \smaller\smaller\smaller\begin{tabular}{@{}c@{}}%
12\\4\\-3
\end{tabular}\endgroup%
{$\left.\llap{\phantom{%
\begingroup \smaller\smaller\smaller\begin{tabular}{@{}c@{}}%
0\\0\\0
\end{tabular}\endgroup%
}}\!\right]$}%
\EasyButWeakLineBreak%
{$\left[\!\llap{\phantom{%
\begingroup \smaller\smaller\smaller\begin{tabular}{@{}c@{}}%
0\\0\\0
\end{tabular}\endgroup%
}}\right.$}%
\begingroup \smaller\smaller\smaller\begin{tabular}{@{}c@{}}%
1\\2\\6
\end{tabular}\endgroup%
\HardButStrongLineBreak\kern3pt%
\begingroup \smaller\smaller\smaller\begin{tabular}{@{}c@{}}%
-1\\0\\-4
\end{tabular}\endgroup%
\HardButStrongLineBreak\kern3pt%
\begingroup \smaller\smaller\smaller\begin{tabular}{@{}c@{}}%
0\\-2\\-3
\end{tabular}\endgroup%
{$\left.\llap{\phantom{%
\begingroup \smaller\smaller\smaller\begin{tabular}{@{}c@{}}%
0\\0\\0
\end{tabular}\endgroup%
}}\!\right]$}%
%
%
\hbox{}\par\smallskip%
%
%
\leavevmode%
${L_{140.2}}$%
{} : {$1\above{1pt}{1pt}{2}{{\rm II}}8\above{1pt}{1pt}{1}{1}$}\EasyButWeakLineBreak%
{${8}\above{1pt}{1pt}{1,0}{\infty a}{8}\above{1pt}{1pt}{r}{2}{2}\above{1pt}{1pt}{4,1}{\infty a}$}%
\nopagebreak\par%
\nopagebreak\par\leavevmode%
{$\left[\!\llap{\phantom{%
\begingroup \smaller\smaller\smaller\begin{tabular}{@{}c@{}}%
0\\0\\0
\end{tabular}\endgroup%
}}\right.$}%
\begingroup \smaller\smaller\smaller\begin{tabular}{@{}c@{}}%
8\\0\\0
\end{tabular}\endgroup%
\kern3pt%
\begingroup \smaller\smaller\smaller\begin{tabular}{@{}c@{}}%
0\\0\\1
\end{tabular}\endgroup%
\kern3pt%
\begingroup \smaller\smaller\smaller\begin{tabular}{@{}c@{}}%
0\\1\\-6
\end{tabular}\endgroup%
{$\left.\llap{\phantom{%
\begingroup \smaller\smaller\smaller\begin{tabular}{@{}c@{}}%
0\\0\\0
\end{tabular}\endgroup%
}}\!\right]$}%
\EasyButWeakLineBreak%
{$\left[\!\llap{\phantom{%
\begingroup \smaller\smaller\smaller\begin{tabular}{@{}c@{}}%
0\\0\\0
\end{tabular}\endgroup%
}}\right.$}%
\begingroup \smaller\smaller\smaller\begin{tabular}{@{}c@{}}%
-1\\4\\0
\end{tabular}\endgroup%
\HardButStrongLineBreak\kern3pt%
\begingroup \smaller\smaller\smaller\begin{tabular}{@{}c@{}}%
1\\0\\0
\end{tabular}\endgroup%
\HardButStrongLineBreak\kern3pt%
\begingroup \smaller\smaller\smaller\begin{tabular}{@{}c@{}}%
0\\-4\\-1
\end{tabular}\endgroup%
{$\left.\llap{\phantom{%
\begingroup \smaller\smaller\smaller\begin{tabular}{@{}c@{}}%
0\\0\\0
\end{tabular}\endgroup%
}}\!\right]$}%
%
%
%
%
%
%
%
%
%
%

\medskip%
%
\leavevmode\llap{}%
$W_{141}$%
\qquad\llap{17} lattices, $\chi=6$%
\hfill%
$\slashtwo2\slashinfty2\rtimes D_{2}$%
\nopagebreak\smallskip\hrule\nopagebreak\medskip%
%
%
\leavevmode%
${L_{141.1}}$%
{} : {$1\above{1pt}{1pt}{2}{0}8\above{1pt}{1pt}{1}{1}$}\EasyButWeakLineBreak%
{${1}\above{1pt}{1pt}{r}{2}{4}\above{1pt}{1pt}{*}{2}{8}\above{1pt}{1pt}{1,0}{\infty b}{8}\above{1pt}{1pt}{}{2}$}%
\nopagebreak\par%
\nopagebreak\par\leavevmode%
{$\left[\!\llap{\phantom{%
\begingroup \smaller\smaller\smaller\begin{tabular}{@{}c@{}}%
0\\0\\0
\end{tabular}\endgroup%
}}\right.$}%
\begingroup \smaller\smaller\smaller\begin{tabular}{@{}c@{}}%
-696\\72\\80
\end{tabular}\endgroup%
\kern3pt%
\begingroup \smaller\smaller\smaller\begin{tabular}{@{}c@{}}%
72\\-7\\-9
\end{tabular}\endgroup%
\kern3pt%
\begingroup \smaller\smaller\smaller\begin{tabular}{@{}c@{}}%
80\\-9\\-8
\end{tabular}\endgroup%
{$\left.\llap{\phantom{%
\begingroup \smaller\smaller\smaller\begin{tabular}{@{}c@{}}%
0\\0\\0
\end{tabular}\endgroup%
}}\!\right]$}%
\EasyButWeakLineBreak%
{$\left[\!\llap{\phantom{%
\begingroup \smaller\smaller\smaller\begin{tabular}{@{}c@{}}%
0\\0\\0
\end{tabular}\endgroup%
}}\right.$}%
\begingroup \smaller\smaller\smaller\begin{tabular}{@{}c@{}}%
1\\5\\4
\end{tabular}\endgroup%
\HardButStrongLineBreak\kern3pt%
\begingroup \smaller\smaller\smaller\begin{tabular}{@{}c@{}}%
-3\\-18\\-10
\end{tabular}\endgroup%
\HardButStrongLineBreak\kern3pt%
\begingroup \smaller\smaller\smaller\begin{tabular}{@{}c@{}}%
-3\\-16\\-12
\end{tabular}\endgroup%
\HardButStrongLineBreak\kern3pt%
\begingroup \smaller\smaller\smaller\begin{tabular}{@{}c@{}}%
7\\40\\24
\end{tabular}\endgroup%
{$\left.\llap{\phantom{%
\begingroup \smaller\smaller\smaller\begin{tabular}{@{}c@{}}%
0\\0\\0
\end{tabular}\endgroup%
}}\!\right]$}%
%
%
\hbox{}\par\smallskip%
%
%
\leavevmode%
${L_{141.2}}$%
{} : {$[1\above{1pt}{1pt}{-}{}2\above{1pt}{1pt}{1}{}]\above{1pt}{1pt}{}{6}16\above{1pt}{1pt}{-}{3}$}\spacer%
\instructions{2}%
\EasyButWeakLineBreak%
{${16}\above{1pt}{1pt}{*}{2}{4}\above{1pt}{1pt}{l}{2}{2}\above{1pt}{1pt}{8,3}{\infty}{8}\above{1pt}{1pt}{s}{2}$}%
\nopagebreak\par%
\nopagebreak\par\leavevmode%
{$\left[\!\llap{\phantom{%
\begingroup \smaller\smaller\smaller\begin{tabular}{@{}c@{}}%
0\\0\\0
\end{tabular}\endgroup%
}}\right.$}%
\begingroup \smaller\smaller\smaller\begin{tabular}{@{}c@{}}%
-848\\160\\32
\end{tabular}\endgroup%
\kern3pt%
\begingroup \smaller\smaller\smaller\begin{tabular}{@{}c@{}}%
160\\-30\\-6
\end{tabular}\endgroup%
\kern3pt%
\begingroup \smaller\smaller\smaller\begin{tabular}{@{}c@{}}%
32\\-6\\-1
\end{tabular}\endgroup%
{$\left.\llap{\phantom{%
\begingroup \smaller\smaller\smaller\begin{tabular}{@{}c@{}}%
0\\0\\0
\end{tabular}\endgroup%
}}\!\right]$}%
\EasyButWeakLineBreak%
{$\left[\!\llap{\phantom{%
\begingroup \smaller\smaller\smaller\begin{tabular}{@{}c@{}}%
0\\0\\0
\end{tabular}\endgroup%
}}\right.$}%
\begingroup \smaller\smaller\smaller\begin{tabular}{@{}c@{}}%
-1\\-4\\-8
\end{tabular}\endgroup%
\HardButStrongLineBreak\kern3pt%
\begingroup \smaller\smaller\smaller\begin{tabular}{@{}c@{}}%
-1\\-6\\2
\end{tabular}\endgroup%
\HardButStrongLineBreak\kern3pt%
\begingroup \smaller\smaller\smaller\begin{tabular}{@{}c@{}}%
0\\-1\\4
\end{tabular}\endgroup%
\HardButStrongLineBreak\kern3pt%
\begingroup \smaller\smaller\smaller\begin{tabular}{@{}c@{}}%
1\\6\\-4
\end{tabular}\endgroup%
{$\left.\llap{\phantom{%
\begingroup \smaller\smaller\smaller\begin{tabular}{@{}c@{}}%
0\\0\\0
\end{tabular}\endgroup%
}}\!\right]$}%
%
%
\hbox{}\par\smallskip%
%
%
\leavevmode%
${L_{141.3}}$%
{} : {$[1\above{1pt}{1pt}{1}{}2\above{1pt}{1pt}{1}{}]\above{1pt}{1pt}{}{2}16\above{1pt}{1pt}{1}{7}$}\spacer%
\instructions{m}%
\EasyButWeakLineBreak%
{${16}\above{1pt}{1pt}{l}{2}{1}\above{1pt}{1pt}{}{2}{2}\above{1pt}{1pt}{8,7}{\infty}{8}\above{1pt}{1pt}{*}{2}$}%
\nopagebreak\par%
\nopagebreak\par\leavevmode%
{$\left[\!\llap{\phantom{%
\begingroup \smaller\smaller\smaller\begin{tabular}{@{}c@{}}%
0\\0\\0
\end{tabular}\endgroup%
}}\right.$}%
\begingroup \smaller\smaller\smaller\begin{tabular}{@{}c@{}}%
-1296\\208\\96
\end{tabular}\endgroup%
\kern3pt%
\begingroup \smaller\smaller\smaller\begin{tabular}{@{}c@{}}%
208\\-30\\-16
\end{tabular}\endgroup%
\kern3pt%
\begingroup \smaller\smaller\smaller\begin{tabular}{@{}c@{}}%
96\\-16\\-7
\end{tabular}\endgroup%
{$\left.\llap{\phantom{%
\begingroup \smaller\smaller\smaller\begin{tabular}{@{}c@{}}%
0\\0\\0
\end{tabular}\endgroup%
}}\!\right]$}%
\EasyButWeakLineBreak%
{$\left[\!\llap{\phantom{%
\begingroup \smaller\smaller\smaller\begin{tabular}{@{}c@{}}%
0\\0\\0
\end{tabular}\endgroup%
}}\right.$}%
\begingroup \smaller\smaller\smaller\begin{tabular}{@{}c@{}}%
-3\\-4\\-32
\end{tabular}\endgroup%
\HardButStrongLineBreak\kern3pt%
\begingroup \smaller\smaller\smaller\begin{tabular}{@{}c@{}}%
1\\2\\9
\end{tabular}\endgroup%
\HardButStrongLineBreak\kern3pt%
\begingroup \smaller\smaller\smaller\begin{tabular}{@{}c@{}}%
2\\3\\20
\end{tabular}\endgroup%
\HardButStrongLineBreak\kern3pt%
\begingroup \smaller\smaller\smaller\begin{tabular}{@{}c@{}}%
-3\\-6\\-28
\end{tabular}\endgroup%
{$\left.\llap{\phantom{%
\begingroup \smaller\smaller\smaller\begin{tabular}{@{}c@{}}%
0\\0\\0
\end{tabular}\endgroup%
}}\!\right]$}%
%
%
\hbox{}\par\smallskip%
%
%
\leavevmode%
${L_{141.4}}$%
{} : {$[1\above{1pt}{1pt}{-}{}2\above{1pt}{1pt}{1}{}]\above{1pt}{1pt}{}{4}16\above{1pt}{1pt}{-}{5}$}\spacer%
\instructions{m}%
\EasyButWeakLineBreak%
{${16}\above{1pt}{1pt}{s}{2}{4}\above{1pt}{1pt}{*}{2}{8}\above{1pt}{1pt}{8,3}{\infty z}{2}\above{1pt}{1pt}{r}{2}$}%
\nopagebreak\par%
\nopagebreak\par\leavevmode%
{$\left[\!\llap{\phantom{%
\begingroup \smaller\smaller\smaller\begin{tabular}{@{}c@{}}%
0\\0\\0
\end{tabular}\endgroup%
}}\right.$}%
\begingroup \smaller\smaller\smaller\begin{tabular}{@{}c@{}}%
-304\\32\\32
\end{tabular}\endgroup%
\kern3pt%
\begingroup \smaller\smaller\smaller\begin{tabular}{@{}c@{}}%
32\\-2\\-4
\end{tabular}\endgroup%
\kern3pt%
\begingroup \smaller\smaller\smaller\begin{tabular}{@{}c@{}}%
32\\-4\\-3
\end{tabular}\endgroup%
{$\left.\llap{\phantom{%
\begingroup \smaller\smaller\smaller\begin{tabular}{@{}c@{}}%
0\\0\\0
\end{tabular}\endgroup%
}}\!\right]$}%
\EasyButWeakLineBreak%
{$\left[\!\llap{\phantom{%
\begingroup \smaller\smaller\smaller\begin{tabular}{@{}c@{}}%
0\\0\\0
\end{tabular}\endgroup%
}}\right.$}%
\begingroup \smaller\smaller\smaller\begin{tabular}{@{}c@{}}%
1\\0\\8
\end{tabular}\endgroup%
\HardButStrongLineBreak\kern3pt%
\begingroup \smaller\smaller\smaller\begin{tabular}{@{}c@{}}%
-1\\-4\\-6
\end{tabular}\endgroup%
\HardButStrongLineBreak\kern3pt%
\begingroup \smaller\smaller\smaller\begin{tabular}{@{}c@{}}%
-1\\-2\\-8
\end{tabular}\endgroup%
\HardButStrongLineBreak\kern3pt%
\begingroup \smaller\smaller\smaller\begin{tabular}{@{}c@{}}%
1\\3\\6
\end{tabular}\endgroup%
{$\left.\llap{\phantom{%
\begingroup \smaller\smaller\smaller\begin{tabular}{@{}c@{}}%
0\\0\\0
\end{tabular}\endgroup%
}}\!\right]$}%
%
%
\hbox{}\par\smallskip%
%
%
\leavevmode%
${L_{141.5}}$%
{} : {$1\above{1pt}{1pt}{1}{1}8\above{1pt}{1pt}{1}{7}64\above{1pt}{1pt}{1}{1}$}\spacer%
\instructions{2}%
\EasyButWeakLineBreak%
{${64}\above{1pt}{1pt}{s}{2}{4}\above{1pt}{1pt}{*}{2}{32}\above{1pt}{1pt}{16,7}{\infty z}{8}\above{1pt}{1pt}{b}{2}$}%
\nopagebreak\par%
shares genus with {$ {L_{141.6}}$}%
; isometric to own %
2-dual\nopagebreak\par%
\nopagebreak\par\leavevmode%
{$\left[\!\llap{\phantom{%
\begingroup \smaller\smaller\smaller\begin{tabular}{@{}c@{}}%
0\\0\\0
\end{tabular}\endgroup%
}}\right.$}%
\begingroup \smaller\smaller\smaller\begin{tabular}{@{}c@{}}%
-17344\\512\\512
\end{tabular}\endgroup%
\kern3pt%
\begingroup \smaller\smaller\smaller\begin{tabular}{@{}c@{}}%
512\\-8\\-16
\end{tabular}\endgroup%
\kern3pt%
\begingroup \smaller\smaller\smaller\begin{tabular}{@{}c@{}}%
512\\-16\\-15
\end{tabular}\endgroup%
{$\left.\llap{\phantom{%
\begingroup \smaller\smaller\smaller\begin{tabular}{@{}c@{}}%
0\\0\\0
\end{tabular}\endgroup%
}}\!\right]$}%
\EasyButWeakLineBreak%
{$\left[\!\llap{\phantom{%
\begingroup \smaller\smaller\smaller\begin{tabular}{@{}c@{}}%
0\\0\\0
\end{tabular}\endgroup%
}}\right.$}%
\begingroup \smaller\smaller\smaller\begin{tabular}{@{}c@{}}%
1\\0\\32
\end{tabular}\endgroup%
\HardButStrongLineBreak\kern3pt%
\begingroup \smaller\smaller\smaller\begin{tabular}{@{}c@{}}%
-1\\-4\\-30
\end{tabular}\endgroup%
\HardButStrongLineBreak\kern3pt%
\begingroup \smaller\smaller\smaller\begin{tabular}{@{}c@{}}%
-1\\-2\\-32
\end{tabular}\endgroup%
\HardButStrongLineBreak\kern3pt%
\begingroup \smaller\smaller\smaller\begin{tabular}{@{}c@{}}%
2\\7\\60
\end{tabular}\endgroup%
{$\left.\llap{\phantom{%
\begingroup \smaller\smaller\smaller\begin{tabular}{@{}c@{}}%
0\\0\\0
\end{tabular}\endgroup%
}}\!\right]$}%
%
%
\hbox{}\par\smallskip%
%
%
\leavevmode%
${L_{141.6}}$%
{} : {$1\above{1pt}{1pt}{1}{1}8\above{1pt}{1pt}{1}{7}64\above{1pt}{1pt}{1}{1}$}\EasyButWeakLineBreak%
{${64}\above{1pt}{1pt}{}{2}{1}\above{1pt}{1pt}{r}{2}{32}\above{1pt}{1pt}{16,15}{\infty z}{8}\above{1pt}{1pt}{l}{2}$}%
\nopagebreak\par%
shares genus with {$ {L_{141.5}}$}%
; isometric to own %
2-dual\nopagebreak\par%
\nopagebreak\par\leavevmode%
{$\left[\!\llap{\phantom{%
\begingroup \smaller\smaller\smaller\begin{tabular}{@{}c@{}}%
0\\0\\0
\end{tabular}\endgroup%
}}\right.$}%
\begingroup \smaller\smaller\smaller\begin{tabular}{@{}c@{}}%
-266176\\2048\\4096
\end{tabular}\endgroup%
\kern3pt%
\begingroup \smaller\smaller\smaller\begin{tabular}{@{}c@{}}%
2048\\-8\\-32
\end{tabular}\endgroup%
\kern3pt%
\begingroup \smaller\smaller\smaller\begin{tabular}{@{}c@{}}%
4096\\-32\\-63
\end{tabular}\endgroup%
{$\left.\llap{\phantom{%
\begingroup \smaller\smaller\smaller\begin{tabular}{@{}c@{}}%
0\\0\\0
\end{tabular}\endgroup%
}}\!\right]$}%
\EasyButWeakLineBreak%
{$\left[\!\llap{\phantom{%
\begingroup \smaller\smaller\smaller\begin{tabular}{@{}c@{}}%
0\\0\\0
\end{tabular}\endgroup%
}}\right.$}%
\begingroup \smaller\smaller\smaller\begin{tabular}{@{}c@{}}%
1\\0\\64
\end{tabular}\endgroup%
\HardButStrongLineBreak\kern3pt%
\begingroup \smaller\smaller\smaller\begin{tabular}{@{}c@{}}%
-1\\-4\\-63
\end{tabular}\endgroup%
\HardButStrongLineBreak\kern3pt%
\begingroup \smaller\smaller\smaller\begin{tabular}{@{}c@{}}%
-1\\-2\\-64
\end{tabular}\endgroup%
\HardButStrongLineBreak\kern3pt%
\begingroup \smaller\smaller\smaller\begin{tabular}{@{}c@{}}%
4\\15\\252
\end{tabular}\endgroup%
{$\left.\llap{\phantom{%
\begingroup \smaller\smaller\smaller\begin{tabular}{@{}c@{}}%
0\\0\\0
\end{tabular}\endgroup%
}}\!\right]$}%

\medskip%
%
\leavevmode\llap{}%
$W_{142}$%
\qquad\llap{22} lattices, $\chi=9$%
\hfill%
$4\infty22$%
\nopagebreak\smallskip\hrule\nopagebreak\medskip%
%
%
\leavevmode%
${L_{142.1}}$%
{} : {$1\above{1pt}{1pt}{2}{{\rm II}}4\above{1pt}{1pt}{1}{1}{\cdot}1\above{1pt}{1pt}{2}{}9\above{1pt}{1pt}{-}{}$}\spacer%
\instructions{2}%
\EasyButWeakLineBreak%
{${2}\above{1pt}{1pt}{*}{4}{4}\above{1pt}{1pt}{3,2}{\infty b}{4}\above{1pt}{1pt}{r}{2}{18}\above{1pt}{1pt}{b}{2}$}%
\nopagebreak\par%
\nopagebreak\par\leavevmode%
{$\left[\!\llap{\phantom{%
\begingroup \smaller\smaller\smaller\begin{tabular}{@{}c@{}}%
0\\0\\0
\end{tabular}\endgroup%
}}\right.$}%
\begingroup \smaller\smaller\smaller\begin{tabular}{@{}c@{}}%
-1116\\288\\108
\end{tabular}\endgroup%
\kern3pt%
\begingroup \smaller\smaller\smaller\begin{tabular}{@{}c@{}}%
288\\-56\\-25
\end{tabular}\endgroup%
\kern3pt%
\begingroup \smaller\smaller\smaller\begin{tabular}{@{}c@{}}%
108\\-25\\-10
\end{tabular}\endgroup%
{$\left.\llap{\phantom{%
\begingroup \smaller\smaller\smaller\begin{tabular}{@{}c@{}}%
0\\0\\0
\end{tabular}\endgroup%
}}\!\right]$}%
\EasyButWeakLineBreak%
{$\left[\!\llap{\phantom{%
\begingroup \smaller\smaller\smaller\begin{tabular}{@{}c@{}}%
0\\0\\0
\end{tabular}\endgroup%
}}\right.$}%
\begingroup \smaller\smaller\smaller\begin{tabular}{@{}c@{}}%
-2\\6\\-37
\end{tabular}\endgroup%
\HardButStrongLineBreak\kern3pt%
\begingroup \smaller\smaller\smaller\begin{tabular}{@{}c@{}}%
-1\\2\\-16
\end{tabular}\endgroup%
\HardButStrongLineBreak\kern3pt%
\begingroup \smaller\smaller\smaller\begin{tabular}{@{}c@{}}%
3\\-8\\52
\end{tabular}\endgroup%
\HardButStrongLineBreak\kern3pt%
\begingroup \smaller\smaller\smaller\begin{tabular}{@{}c@{}}%
1\\0\\9
\end{tabular}\endgroup%
{$\left.\llap{\phantom{%
\begingroup \smaller\smaller\smaller\begin{tabular}{@{}c@{}}%
0\\0\\0
\end{tabular}\endgroup%
}}\!\right]$}%
%
%
\hbox{}\par\smallskip%
%
%
\leavevmode%
${L_{142.2}}$%
{} : {$1\above{1pt}{1pt}{-2}{2}8\above{1pt}{1pt}{-}{3}{\cdot}1\above{1pt}{1pt}{2}{}9\above{1pt}{1pt}{1}{}$}\spacer%
\instructions{2}%
\EasyButWeakLineBreak%
{${4}\above{1pt}{1pt}{*}{4}{2}\above{1pt}{1pt}{12,11}{\infty a}{8}\above{1pt}{1pt}{s}{2}{36}\above{1pt}{1pt}{*}{2}$}%
\nopagebreak\par%
\nopagebreak\par\leavevmode%
{$\left[\!\llap{\phantom{%
\begingroup \smaller\smaller\smaller\begin{tabular}{@{}c@{}}%
0\\0\\0
\end{tabular}\endgroup%
}}\right.$}%
\begingroup \smaller\smaller\smaller\begin{tabular}{@{}c@{}}%
-47592\\-16128\\1728
\end{tabular}\endgroup%
\kern3pt%
\begingroup \smaller\smaller\smaller\begin{tabular}{@{}c@{}}%
-16128\\-5465\\585
\end{tabular}\endgroup%
\kern3pt%
\begingroup \smaller\smaller\smaller\begin{tabular}{@{}c@{}}%
1728\\585\\-62
\end{tabular}\endgroup%
{$\left.\llap{\phantom{%
\begingroup \smaller\smaller\smaller\begin{tabular}{@{}c@{}}%
0\\0\\0
\end{tabular}\endgroup%
}}\!\right]$}%
\EasyButWeakLineBreak%
{$\left[\!\llap{\phantom{%
\begingroup \smaller\smaller\smaller\begin{tabular}{@{}c@{}}%
0\\0\\0
\end{tabular}\endgroup%
}}\right.$}%
\begingroup \smaller\smaller\smaller\begin{tabular}{@{}c@{}}%
13\\-42\\-34
\end{tabular}\endgroup%
\HardButStrongLineBreak\kern3pt%
\begingroup \smaller\smaller\smaller\begin{tabular}{@{}c@{}}%
-10\\32\\23
\end{tabular}\endgroup%
\HardButStrongLineBreak\kern3pt%
\begingroup \smaller\smaller\smaller\begin{tabular}{@{}c@{}}%
-21\\68\\56
\end{tabular}\endgroup%
\HardButStrongLineBreak\kern3pt%
\begingroup \smaller\smaller\smaller\begin{tabular}{@{}c@{}}%
17\\-54\\-36
\end{tabular}\endgroup%
{$\left.\llap{\phantom{%
\begingroup \smaller\smaller\smaller\begin{tabular}{@{}c@{}}%
0\\0\\0
\end{tabular}\endgroup%
}}\!\right]$}%
%
%
\hbox{}\par\smallskip%
%
%
\leavevmode%
${L_{142.3}}$%
{} : {$1\above{1pt}{1pt}{2}{2}8\above{1pt}{1pt}{1}{7}{\cdot}1\above{1pt}{1pt}{2}{}9\above{1pt}{1pt}{1}{}$}\spacer%
\instructions{m}%
\EasyButWeakLineBreak%
{${1}\above{1pt}{1pt}{}{4}{2}\above{1pt}{1pt}{12,11}{\infty b}{8}\above{1pt}{1pt}{l}{2}{9}\above{1pt}{1pt}{}{2}$}%
\nopagebreak\par%
\nopagebreak\par\leavevmode%
{$\left[\!\llap{\phantom{%
\begingroup \smaller\smaller\smaller\begin{tabular}{@{}c@{}}%
0\\0\\0
\end{tabular}\endgroup%
}}\right.$}%
\begingroup \smaller\smaller\smaller\begin{tabular}{@{}c@{}}%
-38088\\576\\648
\end{tabular}\endgroup%
\kern3pt%
\begingroup \smaller\smaller\smaller\begin{tabular}{@{}c@{}}%
576\\-7\\-10
\end{tabular}\endgroup%
\kern3pt%
\begingroup \smaller\smaller\smaller\begin{tabular}{@{}c@{}}%
648\\-10\\-11
\end{tabular}\endgroup%
{$\left.\llap{\phantom{%
\begingroup \smaller\smaller\smaller\begin{tabular}{@{}c@{}}%
0\\0\\0
\end{tabular}\endgroup%
}}\!\right]$}%
\EasyButWeakLineBreak%
{$\left[\!\llap{\phantom{%
\begingroup \smaller\smaller\smaller\begin{tabular}{@{}c@{}}%
0\\0\\0
\end{tabular}\endgroup%
}}\right.$}%
\begingroup \smaller\smaller\smaller\begin{tabular}{@{}c@{}}%
1\\6\\53
\end{tabular}\endgroup%
\HardButStrongLineBreak\kern3pt%
\begingroup \smaller\smaller\smaller\begin{tabular}{@{}c@{}}%
0\\1\\-1
\end{tabular}\endgroup%
\HardButStrongLineBreak\kern3pt%
\begingroup \smaller\smaller\smaller\begin{tabular}{@{}c@{}}%
-1\\-8\\-52
\end{tabular}\endgroup%
\HardButStrongLineBreak\kern3pt%
\begingroup \smaller\smaller\smaller\begin{tabular}{@{}c@{}}%
2\\9\\108
\end{tabular}\endgroup%
{$\left.\llap{\phantom{%
\begingroup \smaller\smaller\smaller\begin{tabular}{@{}c@{}}%
0\\0\\0
\end{tabular}\endgroup%
}}\!\right]$}%

\medskip%
%
\leavevmode\llap{}%
$W_{143}$%
\qquad\llap{22} lattices, $\chi=18$%
\hfill%
$2\infty\slashtwo\infty2|\rtimes D_{2}$%
\nopagebreak\smallskip\hrule\nopagebreak\medskip%
%
%
\leavevmode%
${L_{143.1}}$%
{} : {$1\above{1pt}{1pt}{2}{{\rm II}}4\above{1pt}{1pt}{1}{1}{\cdot}1\above{1pt}{1pt}{-2}{}9\above{1pt}{1pt}{1}{}$}\spacer%
\instructions{2}%
\EasyButWeakLineBreak%
{${2}\above{1pt}{1pt}{l}{2}{36}\above{1pt}{1pt}{1,0}{\infty}{36}\above{1pt}{1pt}{*}{2}{4}\above{1pt}{1pt}{3,1}{\infty a}{4}\above{1pt}{1pt}{r}{2}$}%
\nopagebreak\par%
\nopagebreak\par\leavevmode%
{$\left[\!\llap{\phantom{%
\begingroup \smaller\smaller\smaller\begin{tabular}{@{}c@{}}%
0\\0\\0
\end{tabular}\endgroup%
}}\right.$}%
\begingroup \smaller\smaller\smaller\begin{tabular}{@{}c@{}}%
-2556\\900\\144
\end{tabular}\endgroup%
\kern3pt%
\begingroup \smaller\smaller\smaller\begin{tabular}{@{}c@{}}%
900\\-316\\-51
\end{tabular}\endgroup%
\kern3pt%
\begingroup \smaller\smaller\smaller\begin{tabular}{@{}c@{}}%
144\\-51\\-8
\end{tabular}\endgroup%
{$\left.\llap{\phantom{%
\begingroup \smaller\smaller\smaller\begin{tabular}{@{}c@{}}%
0\\0\\0
\end{tabular}\endgroup%
}}\!\right]$}%
\EasyButWeakLineBreak%
{$\left[\!\llap{\phantom{%
\begingroup \smaller\smaller\smaller\begin{tabular}{@{}c@{}}%
0\\0\\0
\end{tabular}\endgroup%
}}\right.$}%
\begingroup \smaller\smaller\smaller\begin{tabular}{@{}c@{}}%
3\\5\\21
\end{tabular}\endgroup%
\HardButStrongLineBreak\kern3pt%
\begingroup \smaller\smaller\smaller\begin{tabular}{@{}c@{}}%
19\\36\\108
\end{tabular}\endgroup%
\HardButStrongLineBreak\kern3pt%
\begingroup \smaller\smaller\smaller\begin{tabular}{@{}c@{}}%
-1\\0\\-18
\end{tabular}\endgroup%
\HardButStrongLineBreak\kern3pt%
\begingroup \smaller\smaller\smaller\begin{tabular}{@{}c@{}}%
-3\\-6\\-16
\end{tabular}\endgroup%
\HardButStrongLineBreak\kern3pt%
\begingroup \smaller\smaller\smaller\begin{tabular}{@{}c@{}}%
1\\0\\16
\end{tabular}\endgroup%
{$\left.\llap{\phantom{%
\begingroup \smaller\smaller\smaller\begin{tabular}{@{}c@{}}%
0\\0\\0
\end{tabular}\endgroup%
}}\!\right]$}%
%
%
\hbox{}\par\smallskip%
%
%
\leavevmode%
${L_{143.2}}$%
{} : {$1\above{1pt}{1pt}{-2}{2}8\above{1pt}{1pt}{-}{3}{\cdot}1\above{1pt}{1pt}{-2}{}9\above{1pt}{1pt}{-}{}$}\spacer%
\instructions{2}%
\EasyButWeakLineBreak%
{${4}\above{1pt}{1pt}{s}{2}{72}\above{1pt}{1pt}{4,1}{\infty z}{18}\above{1pt}{1pt}{b}{2}{2}\above{1pt}{1pt}{12,7}{\infty b}{8}\above{1pt}{1pt}{s}{2}$}%
\nopagebreak\par%
\nopagebreak\par\leavevmode%
{$\left[\!\llap{\phantom{%
\begingroup \smaller\smaller\smaller\begin{tabular}{@{}c@{}}%
0\\0\\0
\end{tabular}\endgroup%
}}\right.$}%
\begingroup \smaller\smaller\smaller\begin{tabular}{@{}c@{}}%
-52200\\-25488\\576
\end{tabular}\endgroup%
\kern3pt%
\begingroup \smaller\smaller\smaller\begin{tabular}{@{}c@{}}%
-25488\\-12445\\281
\end{tabular}\endgroup%
\kern3pt%
\begingroup \smaller\smaller\smaller\begin{tabular}{@{}c@{}}%
576\\281\\-6
\end{tabular}\endgroup%
{$\left.\llap{\phantom{%
\begingroup \smaller\smaller\smaller\begin{tabular}{@{}c@{}}%
0\\0\\0
\end{tabular}\endgroup%
}}\!\right]$}%
\EasyButWeakLineBreak%
{$\left[\!\llap{\phantom{%
\begingroup \smaller\smaller\smaller\begin{tabular}{@{}c@{}}%
0\\0\\0
\end{tabular}\endgroup%
}}\right.$}%
\begingroup \smaller\smaller\smaller\begin{tabular}{@{}c@{}}%
1\\-2\\2
\end{tabular}\endgroup%
\HardButStrongLineBreak\kern3pt%
\begingroup \smaller\smaller\smaller\begin{tabular}{@{}c@{}}%
-121\\252\\180
\end{tabular}\endgroup%
\HardButStrongLineBreak\kern3pt%
\begingroup \smaller\smaller\smaller\begin{tabular}{@{}c@{}}%
-52\\108\\63
\end{tabular}\endgroup%
\HardButStrongLineBreak\kern3pt%
\begingroup \smaller\smaller\smaller\begin{tabular}{@{}c@{}}%
-1\\2\\-3
\end{tabular}\endgroup%
\HardButStrongLineBreak\kern3pt%
\begingroup \smaller\smaller\smaller\begin{tabular}{@{}c@{}}%
25\\-52\\-36
\end{tabular}\endgroup%
{$\left.\llap{\phantom{%
\begingroup \smaller\smaller\smaller\begin{tabular}{@{}c@{}}%
0\\0\\0
\end{tabular}\endgroup%
}}\!\right]$}%
%
%
\hbox{}\par\smallskip%
%
%
\leavevmode%
${L_{143.3}}$%
{} : {$1\above{1pt}{1pt}{2}{2}8\above{1pt}{1pt}{1}{7}{\cdot}1\above{1pt}{1pt}{-2}{}9\above{1pt}{1pt}{-}{}$}\spacer%
\instructions{m}%
\EasyButWeakLineBreak%
{${1}\above{1pt}{1pt}{r}{2}{72}\above{1pt}{1pt}{4,3}{\infty z}{18}\above{1pt}{1pt}{s}{2}{2}\above{1pt}{1pt}{12,7}{\infty a}{8}\above{1pt}{1pt}{l}{2}$}%
\nopagebreak\par%
\nopagebreak\par\leavevmode%
{$\left[\!\llap{\phantom{%
\begingroup \smaller\smaller\smaller\begin{tabular}{@{}c@{}}%
0\\0\\0
\end{tabular}\endgroup%
}}\right.$}%
\begingroup \smaller\smaller\smaller\begin{tabular}{@{}c@{}}%
-334728\\3456\\1656
\end{tabular}\endgroup%
\kern3pt%
\begingroup \smaller\smaller\smaller\begin{tabular}{@{}c@{}}%
3456\\-35\\-18
\end{tabular}\endgroup%
\kern3pt%
\begingroup \smaller\smaller\smaller\begin{tabular}{@{}c@{}}%
1656\\-18\\-7
\end{tabular}\endgroup%
{$\left.\llap{\phantom{%
\begingroup \smaller\smaller\smaller\begin{tabular}{@{}c@{}}%
0\\0\\0
\end{tabular}\endgroup%
}}\!\right]$}%
\EasyButWeakLineBreak%
{$\left[\!\llap{\phantom{%
\begingroup \smaller\smaller\smaller\begin{tabular}{@{}c@{}}%
0\\0\\0
\end{tabular}\endgroup%
}}\right.$}%
\begingroup \smaller\smaller\smaller\begin{tabular}{@{}c@{}}%
-1\\-71\\-54
\end{tabular}\endgroup%
\HardButStrongLineBreak\kern3pt%
\begingroup \smaller\smaller\smaller\begin{tabular}{@{}c@{}}%
5\\360\\252
\end{tabular}\endgroup%
\HardButStrongLineBreak\kern3pt%
\begingroup \smaller\smaller\smaller\begin{tabular}{@{}c@{}}%
11\\783\\585
\end{tabular}\endgroup%
\HardButStrongLineBreak\kern3pt%
\begingroup \smaller\smaller\smaller\begin{tabular}{@{}c@{}}%
3\\213\\161
\end{tabular}\endgroup%
\HardButStrongLineBreak\kern3pt%
\begingroup \smaller\smaller\smaller\begin{tabular}{@{}c@{}}%
-1\\-72\\-52
\end{tabular}\endgroup%
{$\left.\llap{\phantom{%
\begingroup \smaller\smaller\smaller\begin{tabular}{@{}c@{}}%
0\\0\\0
\end{tabular}\endgroup%
}}\!\right]$}%

\medskip%
%
\leavevmode\llap{}%
$W_{144}$%
\qquad\llap{16} lattices, $\chi=12$%
\hfill%
$2|2\infty|\infty\rtimes D_{2}$%
\nopagebreak\smallskip\hrule\nopagebreak\medskip%
%
%
\leavevmode%
${L_{144.1}}$%
{} : {$1\above{1pt}{1pt}{-2}{4}4\above{1pt}{1pt}{1}{7}{\cdot}1\above{1pt}{1pt}{-}{}3\above{1pt}{1pt}{1}{}9\above{1pt}{1pt}{1}{}$}\spacer%
\instructions{3}%
\EasyButWeakLineBreak%
{${12}\above{1pt}{1pt}{}{2}{9}\above{1pt}{1pt}{}{2}{3}\above{1pt}{1pt}{12,7}{\infty}{12}\above{1pt}{1pt}{3,1}{\infty b}$}%
\nopagebreak\par%
\nopagebreak\par\leavevmode%
{$\left[\!\llap{\phantom{%
\begingroup \smaller\smaller\smaller\begin{tabular}{@{}c@{}}%
0\\0\\0
\end{tabular}\endgroup%
}}\right.$}%
\begingroup \smaller\smaller\smaller\begin{tabular}{@{}c@{}}%
-612\\-180\\-792
\end{tabular}\endgroup%
\kern3pt%
\begingroup \smaller\smaller\smaller\begin{tabular}{@{}c@{}}%
-180\\-33\\-177
\end{tabular}\endgroup%
\kern3pt%
\begingroup \smaller\smaller\smaller\begin{tabular}{@{}c@{}}%
-792\\-177\\-868
\end{tabular}\endgroup%
{$\left.\llap{\phantom{%
\begingroup \smaller\smaller\smaller\begin{tabular}{@{}c@{}}%
0\\0\\0
\end{tabular}\endgroup%
}}\!\right]$}%
\EasyButWeakLineBreak%
{$\left[\!\llap{\phantom{%
\begingroup \smaller\smaller\smaller\begin{tabular}{@{}c@{}}%
0\\0\\0
\end{tabular}\endgroup%
}}\right.$}%
\begingroup \smaller\smaller\smaller\begin{tabular}{@{}c@{}}%
17\\100\\-36
\end{tabular}\endgroup%
\HardButStrongLineBreak\kern3pt%
\begingroup \smaller\smaller\smaller\begin{tabular}{@{}c@{}}%
-4\\-27\\9
\end{tabular}\endgroup%
\HardButStrongLineBreak\kern3pt%
\begingroup \smaller\smaller\smaller\begin{tabular}{@{}c@{}}%
-14\\-85\\30
\end{tabular}\endgroup%
\HardButStrongLineBreak\kern3pt%
\begingroup \smaller\smaller\smaller\begin{tabular}{@{}c@{}}%
-3\\-16\\6
\end{tabular}\endgroup%
{$\left.\llap{\phantom{%
\begingroup \smaller\smaller\smaller\begin{tabular}{@{}c@{}}%
0\\0\\0
\end{tabular}\endgroup%
}}\!\right]$}%
%
%
\hbox{}\par\smallskip%
%
%
\leavevmode%
${L_{144.2}}$%
{} : {$1\above{1pt}{1pt}{2}{{\rm II}}8\above{1pt}{1pt}{-}{3}{\cdot}1\above{1pt}{1pt}{1}{}3\above{1pt}{1pt}{-}{}9\above{1pt}{1pt}{-}{}$}\spacer%
\instructions{3}%
\EasyButWeakLineBreak%
{${24}\above{1pt}{1pt}{r}{2}{18}\above{1pt}{1pt}{b}{2}{6}\above{1pt}{1pt}{12,7}{\infty b}{24}\above{1pt}{1pt}{3,2}{\infty a}$}%
\nopagebreak\par%
\nopagebreak\par\leavevmode%
{$\left[\!\llap{\phantom{%
\begingroup \smaller\smaller\smaller\begin{tabular}{@{}c@{}}%
0\\0\\0
\end{tabular}\endgroup%
}}\right.$}%
\begingroup \smaller\smaller\smaller\begin{tabular}{@{}c@{}}%
-19368\\-3672\\-1944
\end{tabular}\endgroup%
\kern3pt%
\begingroup \smaller\smaller\smaller\begin{tabular}{@{}c@{}}%
-3672\\-696\\-369
\end{tabular}\endgroup%
\kern3pt%
\begingroup \smaller\smaller\smaller\begin{tabular}{@{}c@{}}%
-1944\\-369\\-194
\end{tabular}\endgroup%
{$\left.\llap{\phantom{%
\begingroup \smaller\smaller\smaller\begin{tabular}{@{}c@{}}%
0\\0\\0
\end{tabular}\endgroup%
}}\!\right]$}%
\EasyButWeakLineBreak%
{$\left[\!\llap{\phantom{%
\begingroup \smaller\smaller\smaller\begin{tabular}{@{}c@{}}%
0\\0\\0
\end{tabular}\endgroup%
}}\right.$}%
\begingroup \smaller\smaller\smaller\begin{tabular}{@{}c@{}}%
13\\-56\\-24
\end{tabular}\endgroup%
\HardButStrongLineBreak\kern3pt%
\begingroup \smaller\smaller\smaller\begin{tabular}{@{}c@{}}%
2\\-6\\-9
\end{tabular}\endgroup%
\HardButStrongLineBreak\kern3pt%
\begingroup \smaller\smaller\smaller\begin{tabular}{@{}c@{}}%
-7\\32\\9
\end{tabular}\endgroup%
\HardButStrongLineBreak\kern3pt%
\begingroup \smaller\smaller\smaller\begin{tabular}{@{}c@{}}%
-5\\20\\12
\end{tabular}\endgroup%
{$\left.\llap{\phantom{%
\begingroup \smaller\smaller\smaller\begin{tabular}{@{}c@{}}%
0\\0\\0
\end{tabular}\endgroup%
}}\!\right]$}%

\medskip%
%
\leavevmode\llap{}%
$W_{145}$%
\qquad\llap{2} lattices, $\chi=6$%
\hfill%
$\slashinfty4|4\rtimes D_{2}$%
\nopagebreak\smallskip\hrule\nopagebreak\medskip%
%
%
\leavevmode%
${L_{145.1}}$%
{} : {$1\above{1pt}{1pt}{2}{2}16\above{1pt}{1pt}{1}{7}$}\EasyButWeakLineBreak%
{${1}\above{1pt}{1pt}{8,7}{\infty}{4}\above{1pt}{1pt}{*}{4}{2}\above{1pt}{1pt}{}{4}$}%
\nopagebreak\par%
\nopagebreak\par\leavevmode%
{$\left[\!\llap{\phantom{%
\begingroup \smaller\smaller\smaller\begin{tabular}{@{}c@{}}%
0\\0\\0
\end{tabular}\endgroup%
}}\right.$}%
\begingroup \smaller\smaller\smaller\begin{tabular}{@{}c@{}}%
-656\\-48\\64
\end{tabular}\endgroup%
\kern3pt%
\begingroup \smaller\smaller\smaller\begin{tabular}{@{}c@{}}%
-48\\-3\\5
\end{tabular}\endgroup%
\kern3pt%
\begingroup \smaller\smaller\smaller\begin{tabular}{@{}c@{}}%
64\\5\\-6
\end{tabular}\endgroup%
{$\left.\llap{\phantom{%
\begingroup \smaller\smaller\smaller\begin{tabular}{@{}c@{}}%
0\\0\\0
\end{tabular}\endgroup%
}}\!\right]$}%
\EasyButWeakLineBreak%
{$\left[\!\llap{\phantom{%
\begingroup \smaller\smaller\smaller\begin{tabular}{@{}c@{}}%
0\\0\\0
\end{tabular}\endgroup%
}}\right.$}%
\begingroup \smaller\smaller\smaller\begin{tabular}{@{}c@{}}%
0\\-1\\-1
\end{tabular}\endgroup%
\HardButStrongLineBreak\kern3pt%
\begingroup \smaller\smaller\smaller\begin{tabular}{@{}c@{}}%
-1\\6\\-6
\end{tabular}\endgroup%
\HardButStrongLineBreak\kern3pt%
\begingroup \smaller\smaller\smaller\begin{tabular}{@{}c@{}}%
1\\-4\\7
\end{tabular}\endgroup%
{$\left.\llap{\phantom{%
\begingroup \smaller\smaller\smaller\begin{tabular}{@{}c@{}}%
0\\0\\0
\end{tabular}\endgroup%
}}\!\right]$}%
%
%
%
%
%
%

\medskip%
%
\leavevmode\llap{}%
$W_{146}$%
\qquad\llap{9} lattices, $\chi=12$%
\hfill%
$\slashtwo\slashinfty\slashtwo\slashinfty\rtimes D_{4}$%
\nopagebreak\smallskip\hrule\nopagebreak\medskip%
%
%
\leavevmode%
${L_{146.1}}$%
{} : {$1\above{1pt}{1pt}{-2}{4}16\above{1pt}{1pt}{-}{5}$}\EasyButWeakLineBreak%
{${16}\above{1pt}{1pt}{l}{2}{1}\above{1pt}{1pt}{8,5}{\infty}{4}\above{1pt}{1pt}{s}{2}{16}\above{1pt}{1pt}{4,1}{\infty z}$}%
\nopagebreak\par%
\nopagebreak\par\leavevmode%
{$\left[\!\llap{\phantom{%
\begingroup \smaller\smaller\smaller\begin{tabular}{@{}c@{}}%
0\\0\\0
\end{tabular}\endgroup%
}}\right.$}%
\begingroup \smaller\smaller\smaller\begin{tabular}{@{}c@{}}%
-3632\\64\\336
\end{tabular}\endgroup%
\kern3pt%
\begingroup \smaller\smaller\smaller\begin{tabular}{@{}c@{}}%
64\\-1\\-6
\end{tabular}\endgroup%
\kern3pt%
\begingroup \smaller\smaller\smaller\begin{tabular}{@{}c@{}}%
336\\-6\\-31
\end{tabular}\endgroup%
{$\left.\llap{\phantom{%
\begingroup \smaller\smaller\smaller\begin{tabular}{@{}c@{}}%
0\\0\\0
\end{tabular}\endgroup%
}}\!\right]$}%
\EasyButWeakLineBreak%
{$\left[\!\llap{\phantom{%
\begingroup \smaller\smaller\smaller\begin{tabular}{@{}c@{}}%
0\\0\\0
\end{tabular}\endgroup%
}}\right.$}%
\begingroup \smaller\smaller\smaller\begin{tabular}{@{}c@{}}%
-3\\-8\\-32
\end{tabular}\endgroup%
\HardButStrongLineBreak\kern3pt%
\begingroup \smaller\smaller\smaller\begin{tabular}{@{}c@{}}%
0\\4\\-1
\end{tabular}\endgroup%
\HardButStrongLineBreak\kern3pt%
\begingroup \smaller\smaller\smaller\begin{tabular}{@{}c@{}}%
1\\4\\10
\end{tabular}\endgroup%
\HardButStrongLineBreak\kern3pt%
\begingroup \smaller\smaller\smaller\begin{tabular}{@{}c@{}}%
-1\\-16\\-8
\end{tabular}\endgroup%
{$\left.\llap{\phantom{%
\begingroup \smaller\smaller\smaller\begin{tabular}{@{}c@{}}%
0\\0\\0
\end{tabular}\endgroup%
}}\!\right]$}%
%
%
\hbox{}\par\smallskip%
%
%
\leavevmode%
${L_{146.2}}$%
{} : {$1\above{1pt}{1pt}{2}{0}16\above{1pt}{1pt}{1}{1}$}\EasyButWeakLineBreak%
{${16}\above{1pt}{1pt}{}{2}{1}\above{1pt}{1pt}{8,1}{\infty}{4}\above{1pt}{1pt}{*}{2}{16}\above{1pt}{1pt}{1,0}{\infty b}$}%
\nopagebreak\par%
\nopagebreak\par\leavevmode%
{$\left[\!\llap{\phantom{%
\begingroup \smaller\smaller\smaller\begin{tabular}{@{}c@{}}%
0\\0\\0
\end{tabular}\endgroup%
}}\right.$}%
\begingroup \smaller\smaller\smaller\begin{tabular}{@{}c@{}}%
-1392\\112\\208
\end{tabular}\endgroup%
\kern3pt%
\begingroup \smaller\smaller\smaller\begin{tabular}{@{}c@{}}%
112\\-8\\-17
\end{tabular}\endgroup%
\kern3pt%
\begingroup \smaller\smaller\smaller\begin{tabular}{@{}c@{}}%
208\\-17\\-31
\end{tabular}\endgroup%
{$\left.\llap{\phantom{%
\begingroup \smaller\smaller\smaller\begin{tabular}{@{}c@{}}%
0\\0\\0
\end{tabular}\endgroup%
}}\!\right]$}%
\EasyButWeakLineBreak%
{$\left[\!\llap{\phantom{%
\begingroup \smaller\smaller\smaller\begin{tabular}{@{}c@{}}%
0\\0\\0
\end{tabular}\endgroup%
}}\right.$}%
\begingroup \smaller\smaller\smaller\begin{tabular}{@{}c@{}}%
11\\16\\64
\end{tabular}\endgroup%
\HardButStrongLineBreak\kern3pt%
\begingroup \smaller\smaller\smaller\begin{tabular}{@{}c@{}}%
2\\4\\11
\end{tabular}\endgroup%
\HardButStrongLineBreak\kern3pt%
\begingroup \smaller\smaller\smaller\begin{tabular}{@{}c@{}}%
-3\\-4\\-18
\end{tabular}\endgroup%
\HardButStrongLineBreak\kern3pt%
\begingroup \smaller\smaller\smaller\begin{tabular}{@{}c@{}}%
-3\\-8\\-16
\end{tabular}\endgroup%
{$\left.\llap{\phantom{%
\begingroup \smaller\smaller\smaller\begin{tabular}{@{}c@{}}%
0\\0\\0
\end{tabular}\endgroup%
}}\!\right]$}%
%
%
\hbox{}\par\smallskip%
%
%
\leavevmode%
${L_{146.3}}$%
{} : {$1\above{1pt}{1pt}{-2}{6}16\above{1pt}{1pt}{-}{3}$}\EasyButWeakLineBreak%
{${4}\above{1pt}{1pt}{l}{2}{1}\above{1pt}{1pt}{8,3}{\infty}$}\relax$\,(\times2)$%
\nopagebreak\par%
\nopagebreak\par\leavevmode%
{$\left[\!\llap{\phantom{%
\begingroup \smaller\smaller\smaller\begin{tabular}{@{}c@{}}%
0\\0\\0
\end{tabular}\endgroup%
}}\right.$}%
\begingroup \smaller\smaller\smaller\begin{tabular}{@{}c@{}}%
-80\\48\\16
\end{tabular}\endgroup%
\kern3pt%
\begingroup \smaller\smaller\smaller\begin{tabular}{@{}c@{}}%
48\\-15\\-8
\end{tabular}\endgroup%
\kern3pt%
\begingroup \smaller\smaller\smaller\begin{tabular}{@{}c@{}}%
16\\-8\\-3
\end{tabular}\endgroup%
{$\left.\llap{\phantom{%
\begingroup \smaller\smaller\smaller\begin{tabular}{@{}c@{}}%
0\\0\\0
\end{tabular}\endgroup%
}}\!\right]$}%
\hfil\penalty500%
{$\left[\!\llap{\phantom{%
\begingroup \smaller\smaller\smaller\begin{tabular}{@{}c@{}}%
0\\0\\0
\end{tabular}\endgroup%
}}\right.$}%
\begingroup \smaller\smaller\smaller\begin{tabular}{@{}c@{}}%
15\\-16\\112
\end{tabular}\endgroup%
\kern3pt%
\begingroup \smaller\smaller\smaller\begin{tabular}{@{}c@{}}%
-7\\6\\-49
\end{tabular}\endgroup%
\kern3pt%
\begingroup \smaller\smaller\smaller\begin{tabular}{@{}c@{}}%
-3\\3\\-22
\end{tabular}\endgroup%
{$\left.\llap{\phantom{%
\begingroup \smaller\smaller\smaller\begin{tabular}{@{}c@{}}%
0\\0\\0
\end{tabular}\endgroup%
}}\!\right]$}%
\EasyButWeakLineBreak%
{$\left[\!\llap{\phantom{%
\begingroup \smaller\smaller\smaller\begin{tabular}{@{}c@{}}%
0\\0\\0
\end{tabular}\endgroup%
}}\right.$}%
\begingroup \smaller\smaller\smaller\begin{tabular}{@{}c@{}}%
-1\\0\\-6
\end{tabular}\endgroup%
\HardButStrongLineBreak\kern3pt%
\begingroup \smaller\smaller\smaller\begin{tabular}{@{}c@{}}%
-1\\1\\-8
\end{tabular}\endgroup%
{$\left.\llap{\phantom{%
\begingroup \smaller\smaller\smaller\begin{tabular}{@{}c@{}}%
0\\0\\0
\end{tabular}\endgroup%
}}\!\right]$}%
%
%
\hbox{}\par\smallskip%
%
%
\leavevmode%
${L_{146.4}}$%
{} : {$1\above{1pt}{1pt}{1}{1}4\above{1pt}{1pt}{1}{7}16\above{1pt}{1pt}{1}{1}$}\EasyButWeakLineBreak%
{${16}\above{1pt}{1pt}{s}{2}{4}\above{1pt}{1pt}{4,1}{\infty z}{1}\above{1pt}{1pt}{}{2}{16}\above{1pt}{1pt}{4,3}{\infty}$}%
\nopagebreak\par%
\nopagebreak\par\leavevmode%
{$\left[\!\llap{\phantom{%
\begingroup \smaller\smaller\smaller\begin{tabular}{@{}c@{}}%
0\\0\\0
\end{tabular}\endgroup%
}}\right.$}%
\begingroup \smaller\smaller\smaller\begin{tabular}{@{}c@{}}%
-2672\\112\\288
\end{tabular}\endgroup%
\kern3pt%
\begingroup \smaller\smaller\smaller\begin{tabular}{@{}c@{}}%
112\\-4\\-12
\end{tabular}\endgroup%
\kern3pt%
\begingroup \smaller\smaller\smaller\begin{tabular}{@{}c@{}}%
288\\-12\\-31
\end{tabular}\endgroup%
{$\left.\llap{\phantom{%
\begingroup \smaller\smaller\smaller\begin{tabular}{@{}c@{}}%
0\\0\\0
\end{tabular}\endgroup%
}}\!\right]$}%
\EasyButWeakLineBreak%
{$\left[\!\llap{\phantom{%
\begingroup \smaller\smaller\smaller\begin{tabular}{@{}c@{}}%
0\\0\\0
\end{tabular}\endgroup%
}}\right.$}%
\begingroup \smaller\smaller\smaller\begin{tabular}{@{}c@{}}%
-1\\-4\\-8
\end{tabular}\endgroup%
\HardButStrongLineBreak\kern3pt%
\begingroup \smaller\smaller\smaller\begin{tabular}{@{}c@{}}%
1\\-2\\10
\end{tabular}\endgroup%
\HardButStrongLineBreak\kern3pt%
\begingroup \smaller\smaller\smaller\begin{tabular}{@{}c@{}}%
0\\2\\-1
\end{tabular}\endgroup%
\HardButStrongLineBreak\kern3pt%
\begingroup \smaller\smaller\smaller\begin{tabular}{@{}c@{}}%
-3\\8\\-32
\end{tabular}\endgroup%
{$\left.\llap{\phantom{%
\begingroup \smaller\smaller\smaller\begin{tabular}{@{}c@{}}%
0\\0\\0
\end{tabular}\endgroup%
}}\!\right]$}%
%
%
\hbox{}\par\smallskip%
%
%
\leavevmode%
${L_{146.5}}$%
{} : {$1\above{1pt}{1pt}{1}{1}4\above{1pt}{1pt}{1}{1}16\above{1pt}{1pt}{1}{7}$}\EasyButWeakLineBreak%
{${4}\above{1pt}{1pt}{r}{2}{4}\above{1pt}{1pt}{4,3}{\infty z}{1}\above{1pt}{1pt}{}{2}{4}\above{1pt}{1pt}{8,7}{\infty}$}%
\nopagebreak\par%
\nopagebreak\par\leavevmode%
{$\left[\!\llap{\phantom{%
\begingroup \smaller\smaller\smaller\begin{tabular}{@{}c@{}}%
0\\0\\0
\end{tabular}\endgroup%
}}\right.$}%
\begingroup \smaller\smaller\smaller\begin{tabular}{@{}c@{}}%
-80\\48\\16
\end{tabular}\endgroup%
\kern3pt%
\begingroup \smaller\smaller\smaller\begin{tabular}{@{}c@{}}%
48\\-12\\-8
\end{tabular}\endgroup%
\kern3pt%
\begingroup \smaller\smaller\smaller\begin{tabular}{@{}c@{}}%
16\\-8\\-3
\end{tabular}\endgroup%
{$\left.\llap{\phantom{%
\begingroup \smaller\smaller\smaller\begin{tabular}{@{}c@{}}%
0\\0\\0
\end{tabular}\endgroup%
}}\!\right]$}%
\EasyButWeakLineBreak%
{$\left[\!\llap{\phantom{%
\begingroup \smaller\smaller\smaller\begin{tabular}{@{}c@{}}%
0\\0\\0
\end{tabular}\endgroup%
}}\right.$}%
\begingroup \smaller\smaller\smaller\begin{tabular}{@{}c@{}}%
-1\\1\\-8
\end{tabular}\endgroup%
\HardButStrongLineBreak\kern3pt%
\begingroup \smaller\smaller\smaller\begin{tabular}{@{}c@{}}%
-1\\0\\-6
\end{tabular}\endgroup%
\HardButStrongLineBreak\kern3pt%
\begingroup \smaller\smaller\smaller\begin{tabular}{@{}c@{}}%
1\\-1\\7
\end{tabular}\endgroup%
\HardButStrongLineBreak\kern3pt%
\begingroup \smaller\smaller\smaller\begin{tabular}{@{}c@{}}%
2\\-1\\12
\end{tabular}\endgroup%
{$\left.\llap{\phantom{%
\begingroup \smaller\smaller\smaller\begin{tabular}{@{}c@{}}%
0\\0\\0
\end{tabular}\endgroup%
}}\!\right]$}%
%
%
%
%
%
%
%
%
%
%
%
%
%
%
%
%
%
%

\medskip%
%
\leavevmode\llap{}%
$W_{147}$%
\qquad\llap{4} lattices, $\chi=12$%
\hfill%
$2|2\infty|\infty\rtimes D_{2}$%
\nopagebreak\smallskip\hrule\nopagebreak\medskip%
%
%
\leavevmode%
${L_{147.1}}$%
{} : {$1\above{1pt}{1pt}{2}{{\rm II}}16\above{1pt}{1pt}{1}{1}$}\EasyButWeakLineBreak%
{${16}\above{1pt}{1pt}{r}{2}{2}\above{1pt}{1pt}{b}{2}{16}\above{1pt}{1pt}{2,1}{\infty z}{16}\above{1pt}{1pt}{1,0}{\infty a}$}%
\nopagebreak\par%
\nopagebreak\par\leavevmode%
{$\left[\!\llap{\phantom{%
\begingroup \smaller\smaller\smaller\begin{tabular}{@{}c@{}}%
0\\0\\0
\end{tabular}\endgroup%
}}\right.$}%
\begingroup \smaller\smaller\smaller\begin{tabular}{@{}c@{}}%
16\\0\\0
\end{tabular}\endgroup%
\kern3pt%
\begingroup \smaller\smaller\smaller\begin{tabular}{@{}c@{}}%
0\\-4\\11
\end{tabular}\endgroup%
\kern3pt%
\begingroup \smaller\smaller\smaller\begin{tabular}{@{}c@{}}%
0\\11\\-30
\end{tabular}\endgroup%
{$\left.\llap{\phantom{%
\begingroup \smaller\smaller\smaller\begin{tabular}{@{}c@{}}%
0\\0\\0
\end{tabular}\endgroup%
}}\!\right]$}%
\EasyButWeakLineBreak%
{$\left[\!\llap{\phantom{%
\begingroup \smaller\smaller\smaller\begin{tabular}{@{}c@{}}%
0\\0\\0
\end{tabular}\endgroup%
}}\right.$}%
\begingroup \smaller\smaller\smaller\begin{tabular}{@{}c@{}}%
1\\0\\0
\end{tabular}\endgroup%
\HardButStrongLineBreak\kern3pt%
\begingroup \smaller\smaller\smaller\begin{tabular}{@{}c@{}}%
0\\-8\\-3
\end{tabular}\endgroup%
\HardButStrongLineBreak\kern3pt%
\begingroup \smaller\smaller\smaller\begin{tabular}{@{}c@{}}%
-3\\-16\\-8
\end{tabular}\endgroup%
\HardButStrongLineBreak\kern3pt%
\begingroup \smaller\smaller\smaller\begin{tabular}{@{}c@{}}%
-1\\24\\8
\end{tabular}\endgroup%
{$\left.\llap{\phantom{%
\begingroup \smaller\smaller\smaller\begin{tabular}{@{}c@{}}%
0\\0\\0
\end{tabular}\endgroup%
}}\!\right]$}%
%
%
\hbox{}\par\smallskip%
%
%
\leavevmode%
${L_{147.2}}$%
{} : {$1\above{1pt}{1pt}{1}{7}8\above{1pt}{1pt}{-}{3}64\above{1pt}{1pt}{-}{3}$}\EasyButWeakLineBreak%
{${32}\above{1pt}{1pt}{*}{2}{64}\above{1pt}{1pt}{s}{2}{32}\above{1pt}{1pt}{16,1}{\infty z}{8}\above{1pt}{1pt}{8,3}{\infty b}$}%
\nopagebreak\par%
\nopagebreak\par\leavevmode%
{$\left[\!\llap{\phantom{%
\begingroup \smaller\smaller\smaller\begin{tabular}{@{}c@{}}%
0\\0\\0
\end{tabular}\endgroup%
}}\right.$}%
\begingroup \smaller\smaller\smaller\begin{tabular}{@{}c@{}}%
-3392\\384\\64
\end{tabular}\endgroup%
\kern3pt%
\begingroup \smaller\smaller\smaller\begin{tabular}{@{}c@{}}%
384\\-40\\-8
\end{tabular}\endgroup%
\kern3pt%
\begingroup \smaller\smaller\smaller\begin{tabular}{@{}c@{}}%
64\\-8\\-1
\end{tabular}\endgroup%
{$\left.\llap{\phantom{%
\begingroup \smaller\smaller\smaller\begin{tabular}{@{}c@{}}%
0\\0\\0
\end{tabular}\endgroup%
}}\!\right]$}%
\EasyButWeakLineBreak%
{$\left[\!\llap{\phantom{%
\begingroup \smaller\smaller\smaller\begin{tabular}{@{}c@{}}%
0\\0\\0
\end{tabular}\endgroup%
}}\right.$}%
\begingroup \smaller\smaller\smaller\begin{tabular}{@{}c@{}}%
-3\\-18\\-64
\end{tabular}\endgroup%
\HardButStrongLineBreak\kern3pt%
\begingroup \smaller\smaller\smaller\begin{tabular}{@{}c@{}}%
-1\\-4\\-32
\end{tabular}\endgroup%
\HardButStrongLineBreak\kern3pt%
\begingroup \smaller\smaller\smaller\begin{tabular}{@{}c@{}}%
1\\6\\16
\end{tabular}\endgroup%
\HardButStrongLineBreak\kern3pt%
\begingroup \smaller\smaller\smaller\begin{tabular}{@{}c@{}}%
0\\-1\\4
\end{tabular}\endgroup%
{$\left.\llap{\phantom{%
\begingroup \smaller\smaller\smaller\begin{tabular}{@{}c@{}}%
0\\0\\0
\end{tabular}\endgroup%
}}\!\right]$}%
%
%
%
%
%
%
%
%
%
%

\medskip%
%
\leavevmode\llap{}%
$W_{148}$%
\qquad\llap{16} lattices, $\chi=12$%
\hfill%
$2\infty2\infty\rtimes C_{2}$%
\nopagebreak\smallskip\hrule\nopagebreak\medskip%
%
%
\leavevmode%
${L_{148.1}}$%
{} : {$1\above{1pt}{1pt}{-2}{{\rm II}}8\above{1pt}{1pt}{-}{5}{\cdot}1\above{1pt}{1pt}{-2}{}9\above{1pt}{1pt}{-}{}$}\spacer%
\instructions{2}%
\EasyButWeakLineBreak%
{${18}\above{1pt}{1pt}{b}{2}{8}\above{1pt}{1pt}{6,1}{\infty z}{2}\above{1pt}{1pt}{b}{2}{72}\above{1pt}{1pt}{2,1}{\infty z}$}%
\nopagebreak\par%
\nopagebreak\par\leavevmode%
{$\left[\!\llap{\phantom{%
\begingroup \smaller\smaller\smaller\begin{tabular}{@{}c@{}}%
0\\0\\0
\end{tabular}\endgroup%
}}\right.$}%
\begingroup \smaller\smaller\smaller\begin{tabular}{@{}c@{}}%
936\\288\\0
\end{tabular}\endgroup%
\kern3pt%
\begingroup \smaller\smaller\smaller\begin{tabular}{@{}c@{}}%
288\\82\\15
\end{tabular}\endgroup%
\kern3pt%
\begingroup \smaller\smaller\smaller\begin{tabular}{@{}c@{}}%
0\\15\\-34
\end{tabular}\endgroup%
{$\left.\llap{\phantom{%
\begingroup \smaller\smaller\smaller\begin{tabular}{@{}c@{}}%
0\\0\\0
\end{tabular}\endgroup%
}}\!\right]$}%
\EasyButWeakLineBreak%
{$\left[\!\llap{\phantom{%
\begingroup \smaller\smaller\smaller\begin{tabular}{@{}c@{}}%
0\\0\\0
\end{tabular}\endgroup%
}}\right.$}%
\begingroup \smaller\smaller\smaller\begin{tabular}{@{}c@{}}%
25\\-81\\-36
\end{tabular}\endgroup%
\HardButStrongLineBreak\kern3pt%
\begingroup \smaller\smaller\smaller\begin{tabular}{@{}c@{}}%
11\\-36\\-16
\end{tabular}\endgroup%
\HardButStrongLineBreak\kern3pt%
\begingroup \smaller\smaller\smaller\begin{tabular}{@{}c@{}}%
-12\\39\\17
\end{tabular}\endgroup%
\HardButStrongLineBreak\kern3pt%
\begingroup \smaller\smaller\smaller\begin{tabular}{@{}c@{}}%
-77\\252\\108
\end{tabular}\endgroup%
{$\left.\llap{\phantom{%
\begingroup \smaller\smaller\smaller\begin{tabular}{@{}c@{}}%
0\\0\\0
\end{tabular}\endgroup%
}}\!\right]$}%
%
%
\hbox{}\par\smallskip%
%
%
\leavevmode%
${L_{148.2}}$%
{} : {$1\above{1pt}{1pt}{-2}{{\rm II}}8\above{1pt}{1pt}{-}{5}{\cdot}1\above{1pt}{1pt}{2}{}9\above{1pt}{1pt}{1}{}$}\spacer%
\instructions{2}%
\EasyButWeakLineBreak%
{${2}\above{1pt}{1pt}{b}{2}{8}\above{1pt}{1pt}{6,5}{\infty z}$}\relax$\,(\times2)$%
\nopagebreak\par%
\nopagebreak\par\leavevmode%
{$\left[\!\llap{\phantom{%
\begingroup \smaller\smaller\smaller\begin{tabular}{@{}c@{}}%
0\\0\\0
\end{tabular}\endgroup%
}}\right.$}%
\begingroup \smaller\smaller\smaller\begin{tabular}{@{}c@{}}%
-18648\\432\\216
\end{tabular}\endgroup%
\kern3pt%
\begingroup \smaller\smaller\smaller\begin{tabular}{@{}c@{}}%
432\\-10\\-5
\end{tabular}\endgroup%
\kern3pt%
\begingroup \smaller\smaller\smaller\begin{tabular}{@{}c@{}}%
216\\-5\\-2
\end{tabular}\endgroup%
{$\left.\llap{\phantom{%
\begingroup \smaller\smaller\smaller\begin{tabular}{@{}c@{}}%
0\\0\\0
\end{tabular}\endgroup%
}}\!\right]$}%
\hfil\penalty500%
{$\left[\!\llap{\phantom{%
\begingroup \smaller\smaller\smaller\begin{tabular}{@{}c@{}}%
0\\0\\0
\end{tabular}\endgroup%
}}\right.$}%
\begingroup \smaller\smaller\smaller\begin{tabular}{@{}c@{}}%
-1\\0\\-216
\end{tabular}\endgroup%
\kern3pt%
\begingroup \smaller\smaller\smaller\begin{tabular}{@{}c@{}}%
0\\-1\\5
\end{tabular}\endgroup%
\kern3pt%
\begingroup \smaller\smaller\smaller\begin{tabular}{@{}c@{}}%
0\\0\\1
\end{tabular}\endgroup%
{$\left.\llap{\phantom{%
\begingroup \smaller\smaller\smaller\begin{tabular}{@{}c@{}}%
0\\0\\0
\end{tabular}\endgroup%
}}\!\right]$}%
\EasyButWeakLineBreak%
{$\left[\!\llap{\phantom{%
\begingroup \smaller\smaller\smaller\begin{tabular}{@{}c@{}}%
0\\0\\0
\end{tabular}\endgroup%
}}\right.$}%
\begingroup \smaller\smaller\smaller\begin{tabular}{@{}c@{}}%
0\\1\\-3
\end{tabular}\endgroup%
\HardButStrongLineBreak\kern3pt%
\begingroup \smaller\smaller\smaller\begin{tabular}{@{}c@{}}%
1\\44\\-4
\end{tabular}\endgroup%
{$\left.\llap{\phantom{%
\begingroup \smaller\smaller\smaller\begin{tabular}{@{}c@{}}%
0\\0\\0
\end{tabular}\endgroup%
}}\!\right]$}%

\medskip%
%
\leavevmode\llap{}%
$W_{149}$%
\qquad\llap{27} lattices, $\chi=18$%
\hfill%
$\slashtwo2\infty|\infty2\rtimes D_{2}$%
\nopagebreak\smallskip\hrule\nopagebreak\medskip%
%
%
\leavevmode%
${L_{149.1}}$%
{} : {$1\above{1pt}{1pt}{-2}{4}4\above{1pt}{1pt}{1}{1}{\cdot}1\above{1pt}{1pt}{1}{}5\above{1pt}{1pt}{1}{}25\above{1pt}{1pt}{1}{}$}\spacer%
\instructions{5}%
\EasyButWeakLineBreak%
{${100}\above{1pt}{1pt}{}{2}{1}\above{1pt}{1pt}{}{2}{20}\above{1pt}{1pt}{10,9}{\infty}{20}\above{1pt}{1pt}{5,4}{\infty z}{5}\above{1pt}{1pt}{}{2}$}%
\nopagebreak\par%
shares genus with 5-dual\nopagebreak\par%
\nopagebreak\par\leavevmode%
{$\left[\!\llap{\phantom{%
\begingroup \smaller\smaller\smaller\begin{tabular}{@{}c@{}}%
0\\0\\0
\end{tabular}\endgroup%
}}\right.$}%
\begingroup \smaller\smaller\smaller\begin{tabular}{@{}c@{}}%
-11100\\-800\\-5200
\end{tabular}\endgroup%
\kern3pt%
\begingroup \smaller\smaller\smaller\begin{tabular}{@{}c@{}}%
-800\\-55\\-375
\end{tabular}\endgroup%
\kern3pt%
\begingroup \smaller\smaller\smaller\begin{tabular}{@{}c@{}}%
-5200\\-375\\-2436
\end{tabular}\endgroup%
{$\left.\llap{\phantom{%
\begingroup \smaller\smaller\smaller\begin{tabular}{@{}c@{}}%
0\\0\\0
\end{tabular}\endgroup%
}}\!\right]$}%
\EasyButWeakLineBreak%
{$\left[\!\llap{\phantom{%
\begingroup \smaller\smaller\smaller\begin{tabular}{@{}c@{}}%
0\\0\\0
\end{tabular}\endgroup%
}}\right.$}%
\begingroup \smaller\smaller\smaller\begin{tabular}{@{}c@{}}%
-47\\0\\100
\end{tabular}\endgroup%
\HardButStrongLineBreak\kern3pt%
\begingroup \smaller\smaller\smaller\begin{tabular}{@{}c@{}}%
1\\-1\\-2
\end{tabular}\endgroup%
\HardButStrongLineBreak\kern3pt%
\begingroup \smaller\smaller\smaller\begin{tabular}{@{}c@{}}%
19\\-4\\-40
\end{tabular}\endgroup%
\HardButStrongLineBreak\kern3pt%
\begingroup \smaller\smaller\smaller\begin{tabular}{@{}c@{}}%
-5\\4\\10
\end{tabular}\endgroup%
\HardButStrongLineBreak\kern3pt%
\begingroup \smaller\smaller\smaller\begin{tabular}{@{}c@{}}%
-19\\3\\40
\end{tabular}\endgroup%
{$\left.\llap{\phantom{%
\begingroup \smaller\smaller\smaller\begin{tabular}{@{}c@{}}%
0\\0\\0
\end{tabular}\endgroup%
}}\!\right]$}%
%
%
\hbox{}\par\smallskip%
%
%
\leavevmode%
${L_{149.2}}$%
{} : {$1\above{1pt}{1pt}{2}{{\rm II}}4\above{1pt}{1pt}{-}{5}{\cdot}1\above{1pt}{1pt}{1}{}5\above{1pt}{1pt}{1}{}25\above{1pt}{1pt}{1}{}$}\spacer%
\instructions{2}%
\EasyButWeakLineBreak%
{${100}\above{1pt}{1pt}{*}{2}{4}\above{1pt}{1pt}{*}{2}{20}\above{1pt}{1pt}{5,4}{\infty b}{20}\above{1pt}{1pt}{5,3}{\infty}{20}\above{1pt}{1pt}{*}{2}$}%
\nopagebreak\par%
\nopagebreak\par\leavevmode%
{$\left[\!\llap{\phantom{%
\begingroup \smaller\smaller\smaller\begin{tabular}{@{}c@{}}%
0\\0\\0
\end{tabular}\endgroup%
}}\right.$}%
\begingroup \smaller\smaller\smaller\begin{tabular}{@{}c@{}}%
-1900\\200\\0
\end{tabular}\endgroup%
\kern3pt%
\begingroup \smaller\smaller\smaller\begin{tabular}{@{}c@{}}%
200\\-20\\-5
\end{tabular}\endgroup%
\kern3pt%
\begingroup \smaller\smaller\smaller\begin{tabular}{@{}c@{}}%
0\\-5\\24
\end{tabular}\endgroup%
{$\left.\llap{\phantom{%
\begingroup \smaller\smaller\smaller\begin{tabular}{@{}c@{}}%
0\\0\\0
\end{tabular}\endgroup%
}}\!\right]$}%
\EasyButWeakLineBreak%
{$\left[\!\llap{\phantom{%
\begingroup \smaller\smaller\smaller\begin{tabular}{@{}c@{}}%
0\\0\\0
\end{tabular}\endgroup%
}}\right.$}%
\begingroup \smaller\smaller\smaller\begin{tabular}{@{}c@{}}%
-1\\-10\\0
\end{tabular}\endgroup%
\HardButStrongLineBreak\kern3pt%
\begingroup \smaller\smaller\smaller\begin{tabular}{@{}c@{}}%
3\\28\\6
\end{tabular}\endgroup%
\HardButStrongLineBreak\kern3pt%
\begingroup \smaller\smaller\smaller\begin{tabular}{@{}c@{}}%
11\\102\\20
\end{tabular}\endgroup%
\HardButStrongLineBreak\kern3pt%
\begingroup \smaller\smaller\smaller\begin{tabular}{@{}c@{}}%
1\\8\\0
\end{tabular}\endgroup%
\HardButStrongLineBreak\kern3pt%
\begingroup \smaller\smaller\smaller\begin{tabular}{@{}c@{}}%
-5\\-48\\-10
\end{tabular}\endgroup%
{$\left.\llap{\phantom{%
\begingroup \smaller\smaller\smaller\begin{tabular}{@{}c@{}}%
0\\0\\0
\end{tabular}\endgroup%
}}\!\right]$}%
%
%
\hbox{}\par\smallskip%
%
%
\leavevmode%
${L_{149.3}}$%
{} : {$1\above{1pt}{1pt}{-2}{2}8\above{1pt}{1pt}{1}{7}{\cdot}1\above{1pt}{1pt}{-}{}5\above{1pt}{1pt}{-}{}25\above{1pt}{1pt}{-}{}$}\spacer%
\instructions{2}%
\EasyButWeakLineBreak%
{${50}\above{1pt}{1pt}{b}{2}{2}\above{1pt}{1pt}{s}{2}{10}\above{1pt}{1pt}{20,19}{\infty b}{40}\above{1pt}{1pt}{20,13}{\infty z}{10}\above{1pt}{1pt}{s}{2}$}%
\nopagebreak\par%
\nopagebreak\par\leavevmode%
{$\left[\!\llap{\phantom{%
\begingroup \smaller\smaller\smaller\begin{tabular}{@{}c@{}}%
0\\0\\0
\end{tabular}\endgroup%
}}\right.$}%
\begingroup \smaller\smaller\smaller\begin{tabular}{@{}c@{}}%
143800\\-3600\\-1000
\end{tabular}\endgroup%
\kern3pt%
\begingroup \smaller\smaller\smaller\begin{tabular}{@{}c@{}}%
-3600\\90\\25
\end{tabular}\endgroup%
\kern3pt%
\begingroup \smaller\smaller\smaller\begin{tabular}{@{}c@{}}%
-1000\\25\\7
\end{tabular}\endgroup%
{$\left.\llap{\phantom{%
\begingroup \smaller\smaller\smaller\begin{tabular}{@{}c@{}}%
0\\0\\0
\end{tabular}\endgroup%
}}\!\right]$}%
\EasyButWeakLineBreak%
{$\left[\!\llap{\phantom{%
\begingroup \smaller\smaller\smaller\begin{tabular}{@{}c@{}}%
0\\0\\0
\end{tabular}\endgroup%
}}\right.$}%
\begingroup \smaller\smaller\smaller\begin{tabular}{@{}c@{}}%
-1\\-25\\-50
\end{tabular}\endgroup%
\HardButStrongLineBreak\kern3pt%
\begingroup \smaller\smaller\smaller\begin{tabular}{@{}c@{}}%
0\\-1\\4
\end{tabular}\endgroup%
\HardButStrongLineBreak\kern3pt%
\begingroup \smaller\smaller\smaller\begin{tabular}{@{}c@{}}%
0\\-3\\10
\end{tabular}\endgroup%
\HardButStrongLineBreak\kern3pt%
\begingroup \smaller\smaller\smaller\begin{tabular}{@{}c@{}}%
-1\\-24\\-60
\end{tabular}\endgroup%
\HardButStrongLineBreak\kern3pt%
\begingroup \smaller\smaller\smaller\begin{tabular}{@{}c@{}}%
-1\\-23\\-60
\end{tabular}\endgroup%
{$\left.\llap{\phantom{%
\begingroup \smaller\smaller\smaller\begin{tabular}{@{}c@{}}%
0\\0\\0
\end{tabular}\endgroup%
}}\!\right]$}%
%
%
\hbox{}\par\smallskip%
%
%
\leavevmode%
${L_{149.4}}$%
{} : {$1\above{1pt}{1pt}{2}{2}8\above{1pt}{1pt}{-}{3}{\cdot}1\above{1pt}{1pt}{-}{}5\above{1pt}{1pt}{-}{}25\above{1pt}{1pt}{-}{}$}\spacer%
\instructions{m}%
\EasyButWeakLineBreak%
{${50}\above{1pt}{1pt}{s}{2}{2}\above{1pt}{1pt}{b}{2}{10}\above{1pt}{1pt}{20,19}{\infty a}{40}\above{1pt}{1pt}{20,3}{\infty z}{10}\above{1pt}{1pt}{b}{2}$}%
\nopagebreak\par%
\nopagebreak\par\leavevmode%
{$\left[\!\llap{\phantom{%
\begingroup \smaller\smaller\smaller\begin{tabular}{@{}c@{}}%
0\\0\\0
\end{tabular}\endgroup%
}}\right.$}%
\begingroup \smaller\smaller\smaller\begin{tabular}{@{}c@{}}%
442200\\-6400\\-2400
\end{tabular}\endgroup%
\kern3pt%
\begingroup \smaller\smaller\smaller\begin{tabular}{@{}c@{}}%
-6400\\90\\35
\end{tabular}\endgroup%
\kern3pt%
\begingroup \smaller\smaller\smaller\begin{tabular}{@{}c@{}}%
-2400\\35\\13
\end{tabular}\endgroup%
{$\left.\llap{\phantom{%
\begingroup \smaller\smaller\smaller\begin{tabular}{@{}c@{}}%
0\\0\\0
\end{tabular}\endgroup%
}}\!\right]$}%
\EasyButWeakLineBreak%
{$\left[\!\llap{\phantom{%
\begingroup \smaller\smaller\smaller\begin{tabular}{@{}c@{}}%
0\\0\\0
\end{tabular}\endgroup%
}}\right.$}%
\begingroup \smaller\smaller\smaller\begin{tabular}{@{}c@{}}%
-2\\-25\\-300
\end{tabular}\endgroup%
\HardButStrongLineBreak\kern3pt%
\begingroup \smaller\smaller\smaller\begin{tabular}{@{}c@{}}%
0\\1\\-2
\end{tabular}\endgroup%
\HardButStrongLineBreak\kern3pt%
\begingroup \smaller\smaller\smaller\begin{tabular}{@{}c@{}}%
1\\17\\140
\end{tabular}\endgroup%
\HardButStrongLineBreak\kern3pt%
\begingroup \smaller\smaller\smaller\begin{tabular}{@{}c@{}}%
1\\16\\140
\end{tabular}\endgroup%
\HardButStrongLineBreak\kern3pt%
\begingroup \smaller\smaller\smaller\begin{tabular}{@{}c@{}}%
-1\\-13\\-150
\end{tabular}\endgroup%
{$\left.\llap{\phantom{%
\begingroup \smaller\smaller\smaller\begin{tabular}{@{}c@{}}%
0\\0\\0
\end{tabular}\endgroup%
}}\!\right]$}%
%
%
\hbox{}\par\smallskip%
%
%
\leavevmode%
${L_{149.5}}$%
{} : {$1\above{1pt}{1pt}{2}{{\rm II}}8\above{1pt}{1pt}{-}{5}{\cdot}1\above{1pt}{1pt}{-}{}5\above{1pt}{1pt}{-}{}25\above{1pt}{1pt}{-}{}$}\spacer%
\instructions{5}%
\EasyButWeakLineBreak%
{${200}\above{1pt}{1pt}{b}{2}{2}\above{1pt}{1pt}{l}{2}{40}\above{1pt}{1pt}{5,4}{\infty}{40}\above{1pt}{1pt}{5,2}{\infty z}{10}\above{1pt}{1pt}{b}{2}$}%
\nopagebreak\par%
shares genus with 5-dual\nopagebreak\par%
\nopagebreak\par\leavevmode%
{$\left[\!\llap{\phantom{%
\begingroup \smaller\smaller\smaller\begin{tabular}{@{}c@{}}%
0\\0\\0
\end{tabular}\endgroup%
}}\right.$}%
\begingroup \smaller\smaller\smaller\begin{tabular}{@{}c@{}}%
-2200\\600\\-2600
\end{tabular}\endgroup%
\kern3pt%
\begingroup \smaller\smaller\smaller\begin{tabular}{@{}c@{}}%
600\\-110\\425
\end{tabular}\endgroup%
\kern3pt%
\begingroup \smaller\smaller\smaller\begin{tabular}{@{}c@{}}%
-2600\\425\\-1568
\end{tabular}\endgroup%
{$\left.\llap{\phantom{%
\begingroup \smaller\smaller\smaller\begin{tabular}{@{}c@{}}%
0\\0\\0
\end{tabular}\endgroup%
}}\!\right]$}%
\EasyButWeakLineBreak%
{$\left[\!\llap{\phantom{%
\begingroup \smaller\smaller\smaller\begin{tabular}{@{}c@{}}%
0\\0\\0
\end{tabular}\endgroup%
}}\right.$}%
\begingroup \smaller\smaller\smaller\begin{tabular}{@{}c@{}}%
53\\1060\\200
\end{tabular}\endgroup%
\HardButStrongLineBreak\kern3pt%
\begingroup \smaller\smaller\smaller\begin{tabular}{@{}c@{}}%
-1\\-21\\-4
\end{tabular}\endgroup%
\HardButStrongLineBreak\kern3pt%
\begingroup \smaller\smaller\smaller\begin{tabular}{@{}c@{}}%
-21\\-424\\-80
\end{tabular}\endgroup%
\HardButStrongLineBreak\kern3pt%
\begingroup \smaller\smaller\smaller\begin{tabular}{@{}c@{}}%
5\\104\\20
\end{tabular}\endgroup%
\HardButStrongLineBreak\kern3pt%
\begingroup \smaller\smaller\smaller\begin{tabular}{@{}c@{}}%
21\\423\\80
\end{tabular}\endgroup%
{$\left.\llap{\phantom{%
\begingroup \smaller\smaller\smaller\begin{tabular}{@{}c@{}}%
0\\0\\0
\end{tabular}\endgroup%
}}\!\right]$}%

\medskip%
%
\leavevmode\llap{}%
$W_{150}$%
\qquad\llap{60} lattices, $\chi=12$%
\hfill%
$22|22\slashinfty\rtimes D_{2}$%
\nopagebreak\smallskip\hrule\nopagebreak\medskip%
%
%
\leavevmode%
${L_{150.1}}$%
{} : {$1\above{1pt}{1pt}{2}{0}8\above{1pt}{1pt}{-}{3}{\cdot}1\above{1pt}{1pt}{-}{}3\above{1pt}{1pt}{-}{}9\above{1pt}{1pt}{1}{}$}\spacer%
\instructions{3}%
\EasyButWeakLineBreak%
{${24}\above{1pt}{1pt}{}{2}{9}\above{1pt}{1pt}{r}{2}{8}\above{1pt}{1pt}{s}{2}{36}\above{1pt}{1pt}{*}{2}{24}\above{1pt}{1pt}{3,1}{\infty a}$}%
\nopagebreak\par%
\nopagebreak\par\leavevmode%
{$\left[\!\llap{\phantom{%
\begingroup \smaller\smaller\smaller\begin{tabular}{@{}c@{}}%
0\\0\\0
\end{tabular}\endgroup%
}}\right.$}%
\begingroup \smaller\smaller\smaller\begin{tabular}{@{}c@{}}%
-28584\\144\\1080
\end{tabular}\endgroup%
\kern3pt%
\begingroup \smaller\smaller\smaller\begin{tabular}{@{}c@{}}%
144\\15\\-9
\end{tabular}\endgroup%
\kern3pt%
\begingroup \smaller\smaller\smaller\begin{tabular}{@{}c@{}}%
1080\\-9\\-40
\end{tabular}\endgroup%
{$\left.\llap{\phantom{%
\begingroup \smaller\smaller\smaller\begin{tabular}{@{}c@{}}%
0\\0\\0
\end{tabular}\endgroup%
}}\!\right]$}%
\EasyButWeakLineBreak%
{$\left[\!\llap{\phantom{%
\begingroup \smaller\smaller\smaller\begin{tabular}{@{}c@{}}%
0\\0\\0
\end{tabular}\endgroup%
}}\right.$}%
\begingroup \smaller\smaller\smaller\begin{tabular}{@{}c@{}}%
15\\88\\384
\end{tabular}\endgroup%
\HardButStrongLineBreak\kern3pt%
\begingroup \smaller\smaller\smaller\begin{tabular}{@{}c@{}}%
13\\75\\333
\end{tabular}\endgroup%
\HardButStrongLineBreak\kern3pt%
\begingroup \smaller\smaller\smaller\begin{tabular}{@{}c@{}}%
5\\28\\128
\end{tabular}\endgroup%
\HardButStrongLineBreak\kern3pt%
\begingroup \smaller\smaller\smaller\begin{tabular}{@{}c@{}}%
-7\\-42\\-180
\end{tabular}\endgroup%
\HardButStrongLineBreak\kern3pt%
\begingroup \smaller\smaller\smaller\begin{tabular}{@{}c@{}}%
-7\\-40\\-180
\end{tabular}\endgroup%
{$\left.\llap{\phantom{%
\begingroup \smaller\smaller\smaller\begin{tabular}{@{}c@{}}%
0\\0\\0
\end{tabular}\endgroup%
}}\!\right]$}%
%
%
\hbox{}\par\smallskip%
%
%
\leavevmode%
${L_{150.2}}$%
{} : {$[1\above{1pt}{1pt}{1}{}2\above{1pt}{1pt}{-}{}]\above{1pt}{1pt}{}{4}16\above{1pt}{1pt}{1}{7}{\cdot}1\above{1pt}{1pt}{-}{}3\above{1pt}{1pt}{-}{}9\above{1pt}{1pt}{1}{}$}\spacer%
\instructions{3m,3,2}%
\EasyButWeakLineBreak%
{${6}\above{1pt}{1pt}{}{2}{9}\above{1pt}{1pt}{r}{2}{8}\above{1pt}{1pt}{*}{2}{144}\above{1pt}{1pt}{s}{2}{24}\above{1pt}{1pt}{24,1}{\infty z}$}%
\nopagebreak\par%
\nopagebreak\par\leavevmode%
{$\left[\!\llap{\phantom{%
\begingroup \smaller\smaller\smaller\begin{tabular}{@{}c@{}}%
0\\0\\0
\end{tabular}\endgroup%
}}\right.$}%
\begingroup \smaller\smaller\smaller\begin{tabular}{@{}c@{}}%
-75024\\4608\\288
\end{tabular}\endgroup%
\kern3pt%
\begingroup \smaller\smaller\smaller\begin{tabular}{@{}c@{}}%
4608\\-282\\-18
\end{tabular}\endgroup%
\kern3pt%
\begingroup \smaller\smaller\smaller\begin{tabular}{@{}c@{}}%
288\\-18\\-1
\end{tabular}\endgroup%
{$\left.\llap{\phantom{%
\begingroup \smaller\smaller\smaller\begin{tabular}{@{}c@{}}%
0\\0\\0
\end{tabular}\endgroup%
}}\!\right]$}%
\EasyButWeakLineBreak%
{$\left[\!\llap{\phantom{%
\begingroup \smaller\smaller\smaller\begin{tabular}{@{}c@{}}%
0\\0\\0
\end{tabular}\endgroup%
}}\right.$}%
\begingroup \smaller\smaller\smaller\begin{tabular}{@{}c@{}}%
0\\-1\\12
\end{tabular}\endgroup%
\HardButStrongLineBreak\kern3pt%
\begingroup \smaller\smaller\smaller\begin{tabular}{@{}c@{}}%
-1\\-15\\-27
\end{tabular}\endgroup%
\HardButStrongLineBreak\kern3pt%
\begingroup \smaller\smaller\smaller\begin{tabular}{@{}c@{}}%
-1\\-14\\-40
\end{tabular}\endgroup%
\HardButStrongLineBreak\kern3pt%
\begingroup \smaller\smaller\smaller\begin{tabular}{@{}c@{}}%
-1\\-12\\-72
\end{tabular}\endgroup%
\HardButStrongLineBreak\kern3pt%
\begingroup \smaller\smaller\smaller\begin{tabular}{@{}c@{}}%
1\\14\\36
\end{tabular}\endgroup%
{$\left.\llap{\phantom{%
\begingroup \smaller\smaller\smaller\begin{tabular}{@{}c@{}}%
0\\0\\0
\end{tabular}\endgroup%
}}\!\right]$}%
%
%
\hbox{}\par\smallskip%
%
%
\leavevmode%
${L_{150.3}}$%
{} : {$[1\above{1pt}{1pt}{-}{}2\above{1pt}{1pt}{-}{}]\above{1pt}{1pt}{}{0}16\above{1pt}{1pt}{-}{3}{\cdot}1\above{1pt}{1pt}{-}{}3\above{1pt}{1pt}{-}{}9\above{1pt}{1pt}{1}{}$}\spacer%
\instructions{32,3,m}%
\EasyButWeakLineBreak%
{${6}\above{1pt}{1pt}{r}{2}{36}\above{1pt}{1pt}{*}{2}{8}\above{1pt}{1pt}{s}{2}{144}\above{1pt}{1pt}{*}{2}{24}\above{1pt}{1pt}{24,13}{\infty z}$}%
\nopagebreak\par%
\nopagebreak\par\leavevmode%
{$\left[\!\llap{\phantom{%
\begingroup \smaller\smaller\smaller\begin{tabular}{@{}c@{}}%
0\\0\\0
\end{tabular}\endgroup%
}}\right.$}%
\begingroup \smaller\smaller\smaller\begin{tabular}{@{}c@{}}%
-118224\\5904\\2304
\end{tabular}\endgroup%
\kern3pt%
\begingroup \smaller\smaller\smaller\begin{tabular}{@{}c@{}}%
5904\\-282\\-120
\end{tabular}\endgroup%
\kern3pt%
\begingroup \smaller\smaller\smaller\begin{tabular}{@{}c@{}}%
2304\\-120\\-43
\end{tabular}\endgroup%
{$\left.\llap{\phantom{%
\begingroup \smaller\smaller\smaller\begin{tabular}{@{}c@{}}%
0\\0\\0
\end{tabular}\endgroup%
}}\!\right]$}%
\EasyButWeakLineBreak%
{$\left[\!\llap{\phantom{%
\begingroup \smaller\smaller\smaller\begin{tabular}{@{}c@{}}%
0\\0\\0
\end{tabular}\endgroup%
}}\right.$}%
\begingroup \smaller\smaller\smaller\begin{tabular}{@{}c@{}}%
6\\59\\156
\end{tabular}\endgroup%
\HardButStrongLineBreak\kern3pt%
\begingroup \smaller\smaller\smaller\begin{tabular}{@{}c@{}}%
23\\228\\594
\end{tabular}\endgroup%
\HardButStrongLineBreak\kern3pt%
\begingroup \smaller\smaller\smaller\begin{tabular}{@{}c@{}}%
5\\50\\128
\end{tabular}\endgroup%
\HardButStrongLineBreak\kern3pt%
\begingroup \smaller\smaller\smaller\begin{tabular}{@{}c@{}}%
-11\\-108\\-288
\end{tabular}\endgroup%
\HardButStrongLineBreak\kern3pt%
\begingroup \smaller\smaller\smaller\begin{tabular}{@{}c@{}}%
-7\\-70\\-180
\end{tabular}\endgroup%
{$\left.\llap{\phantom{%
\begingroup \smaller\smaller\smaller\begin{tabular}{@{}c@{}}%
0\\0\\0
\end{tabular}\endgroup%
}}\!\right]$}%
%
%
\hbox{}\par\smallskip%
%
%
\leavevmode%
${L_{150.4}}$%
{} : {$[1\above{1pt}{1pt}{-}{}2\above{1pt}{1pt}{1}{}]\above{1pt}{1pt}{}{6}16\above{1pt}{1pt}{1}{1}{\cdot}1\above{1pt}{1pt}{-}{}3\above{1pt}{1pt}{-}{}9\above{1pt}{1pt}{1}{}$}\spacer%
\instructions{3m,3,m}%
\EasyButWeakLineBreak%
{${24}\above{1pt}{1pt}{*}{2}{36}\above{1pt}{1pt}{l}{2}{2}\above{1pt}{1pt}{}{2}{144}\above{1pt}{1pt}{}{2}{6}\above{1pt}{1pt}{24,1}{\infty}$}%
\nopagebreak\par%
\nopagebreak\par\leavevmode%
{$\left[\!\llap{\phantom{%
\begingroup \smaller\smaller\smaller\begin{tabular}{@{}c@{}}%
0\\0\\0
\end{tabular}\endgroup%
}}\right.$}%
\begingroup \smaller\smaller\smaller\begin{tabular}{@{}c@{}}%
-10224\\576\\1296
\end{tabular}\endgroup%
\kern3pt%
\begingroup \smaller\smaller\smaller\begin{tabular}{@{}c@{}}%
576\\-30\\-84
\end{tabular}\endgroup%
\kern3pt%
\begingroup \smaller\smaller\smaller\begin{tabular}{@{}c@{}}%
1296\\-84\\-115
\end{tabular}\endgroup%
{$\left.\llap{\phantom{%
\begingroup \smaller\smaller\smaller\begin{tabular}{@{}c@{}}%
0\\0\\0
\end{tabular}\endgroup%
}}\!\right]$}%
\EasyButWeakLineBreak%
{$\left[\!\llap{\phantom{%
\begingroup \smaller\smaller\smaller\begin{tabular}{@{}c@{}}%
0\\0\\0
\end{tabular}\endgroup%
}}\right.$}%
\begingroup \smaller\smaller\smaller\begin{tabular}{@{}c@{}}%
-9\\-106\\-24
\end{tabular}\endgroup%
\HardButStrongLineBreak\kern3pt%
\begingroup \smaller\smaller\smaller\begin{tabular}{@{}c@{}}%
-7\\-84\\-18
\end{tabular}\endgroup%
\HardButStrongLineBreak\kern3pt%
\begingroup \smaller\smaller\smaller\begin{tabular}{@{}c@{}}%
3\\35\\8
\end{tabular}\endgroup%
\HardButStrongLineBreak\kern3pt%
\begingroup \smaller\smaller\smaller\begin{tabular}{@{}c@{}}%
55\\648\\144
\end{tabular}\endgroup%
\HardButStrongLineBreak\kern3pt%
\begingroup \smaller\smaller\smaller\begin{tabular}{@{}c@{}}%
7\\83\\18
\end{tabular}\endgroup%
{$\left.\llap{\phantom{%
\begingroup \smaller\smaller\smaller\begin{tabular}{@{}c@{}}%
0\\0\\0
\end{tabular}\endgroup%
}}\!\right]$}%
%
%
\hbox{}\par\smallskip%
%
%
\leavevmode%
${L_{150.5}}$%
{} : {$[1\above{1pt}{1pt}{1}{}2\above{1pt}{1pt}{1}{}]\above{1pt}{1pt}{}{2}16\above{1pt}{1pt}{-}{5}{\cdot}1\above{1pt}{1pt}{-}{}3\above{1pt}{1pt}{-}{}9\above{1pt}{1pt}{1}{}$}\spacer%
\instructions{3}%
\EasyButWeakLineBreak%
{${24}\above{1pt}{1pt}{l}{2}{9}\above{1pt}{1pt}{}{2}{2}\above{1pt}{1pt}{r}{2}{144}\above{1pt}{1pt}{l}{2}{6}\above{1pt}{1pt}{24,13}{\infty}$}%
\nopagebreak\par%
\nopagebreak\par\leavevmode%
{$\left[\!\llap{\phantom{%
\begingroup \smaller\smaller\smaller\begin{tabular}{@{}c@{}}%
0\\0\\0
\end{tabular}\endgroup%
}}\right.$}%
\begingroup \smaller\smaller\smaller\begin{tabular}{@{}c@{}}%
-60336\\1440\\-13392
\end{tabular}\endgroup%
\kern3pt%
\begingroup \smaller\smaller\smaller\begin{tabular}{@{}c@{}}%
1440\\-30\\336
\end{tabular}\endgroup%
\kern3pt%
\begingroup \smaller\smaller\smaller\begin{tabular}{@{}c@{}}%
-13392\\336\\-2911
\end{tabular}\endgroup%
{$\left.\llap{\phantom{%
\begingroup \smaller\smaller\smaller\begin{tabular}{@{}c@{}}%
0\\0\\0
\end{tabular}\endgroup%
}}\!\right]$}%
\EasyButWeakLineBreak%
{$\left[\!\llap{\phantom{%
\begingroup \smaller\smaller\smaller\begin{tabular}{@{}c@{}}%
0\\0\\0
\end{tabular}\endgroup%
}}\right.$}%
\begingroup \smaller\smaller\smaller\begin{tabular}{@{}c@{}}%
15\\182\\-48
\end{tabular}\endgroup%
\HardButStrongLineBreak\kern3pt%
\begingroup \smaller\smaller\smaller\begin{tabular}{@{}c@{}}%
14\\168\\-45
\end{tabular}\endgroup%
\HardButStrongLineBreak\kern3pt%
\begingroup \smaller\smaller\smaller\begin{tabular}{@{}c@{}}%
-5\\-61\\16
\end{tabular}\endgroup%
\HardButStrongLineBreak\kern3pt%
\begingroup \smaller\smaller\smaller\begin{tabular}{@{}c@{}}%
-157\\-1896\\504
\end{tabular}\endgroup%
\HardButStrongLineBreak\kern3pt%
\begingroup \smaller\smaller\smaller\begin{tabular}{@{}c@{}}%
-28\\-337\\90
\end{tabular}\endgroup%
{$\left.\llap{\phantom{%
\begingroup \smaller\smaller\smaller\begin{tabular}{@{}c@{}}%
0\\0\\0
\end{tabular}\endgroup%
}}\!\right]$}%

\medskip%
%
\leavevmode\llap{}%
$W_{151}$%
\qquad\llap{30} lattices, $\chi=6$%
\hfill%
$\slashtwo22|22\rtimes D_{2}$%
\nopagebreak\smallskip\hrule\nopagebreak\medskip%
%
%
\leavevmode%
${L_{151.1}}$%
{} : {$1\above{1pt}{1pt}{-2}{4}8\above{1pt}{1pt}{1}{7}{\cdot}1\above{1pt}{1pt}{2}{}3\above{1pt}{1pt}{1}{}$}\EasyButWeakLineBreak%
{${1}\above{1pt}{1pt}{r}{2}{4}\above{1pt}{1pt}{*}{2}{12}\above{1pt}{1pt}{s}{2}{8}\above{1pt}{1pt}{l}{2}{3}\above{1pt}{1pt}{}{2}$}%
\nopagebreak\par%
\nopagebreak\par\leavevmode%
{$\left[\!\llap{\phantom{%
\begingroup \smaller\smaller\smaller\begin{tabular}{@{}c@{}}%
0\\0\\0
\end{tabular}\endgroup%
}}\right.$}%
\begingroup \smaller\smaller\smaller\begin{tabular}{@{}c@{}}%
-456\\48\\0
\end{tabular}\endgroup%
\kern3pt%
\begingroup \smaller\smaller\smaller\begin{tabular}{@{}c@{}}%
48\\-1\\-2
\end{tabular}\endgroup%
\kern3pt%
\begingroup \smaller\smaller\smaller\begin{tabular}{@{}c@{}}%
0\\-2\\1
\end{tabular}\endgroup%
{$\left.\llap{\phantom{%
\begingroup \smaller\smaller\smaller\begin{tabular}{@{}c@{}}%
0\\0\\0
\end{tabular}\endgroup%
}}\!\right]$}%
\EasyButWeakLineBreak%
{$\left[\!\llap{\phantom{%
\begingroup \smaller\smaller\smaller\begin{tabular}{@{}c@{}}%
0\\0\\0
\end{tabular}\endgroup%
}}\right.$}%
\begingroup \smaller\smaller\smaller\begin{tabular}{@{}c@{}}%
0\\0\\-1
\end{tabular}\endgroup%
\HardButStrongLineBreak\kern3pt%
\begingroup \smaller\smaller\smaller\begin{tabular}{@{}c@{}}%
-1\\-10\\-20
\end{tabular}\endgroup%
\HardButStrongLineBreak\kern3pt%
\begingroup \smaller\smaller\smaller\begin{tabular}{@{}c@{}}%
-1\\-12\\-18
\end{tabular}\endgroup%
\HardButStrongLineBreak\kern3pt%
\begingroup \smaller\smaller\smaller\begin{tabular}{@{}c@{}}%
1\\8\\20
\end{tabular}\endgroup%
\HardButStrongLineBreak\kern3pt%
\begingroup \smaller\smaller\smaller\begin{tabular}{@{}c@{}}%
1\\9\\18
\end{tabular}\endgroup%
{$\left.\llap{\phantom{%
\begingroup \smaller\smaller\smaller\begin{tabular}{@{}c@{}}%
0\\0\\0
\end{tabular}\endgroup%
}}\!\right]$}%
%
%
\hbox{}\par\smallskip%
%
%
\leavevmode%
${L_{151.2}}$%
{} : {$[1\above{1pt}{1pt}{-}{}2\above{1pt}{1pt}{1}{}]\above{1pt}{1pt}{}{4}16\above{1pt}{1pt}{1}{7}{\cdot}1\above{1pt}{1pt}{2}{}3\above{1pt}{1pt}{1}{}$}\spacer%
\instructions{2}%
\EasyButWeakLineBreak%
{${16}\above{1pt}{1pt}{*}{2}{4}\above{1pt}{1pt}{s}{2}{48}\above{1pt}{1pt}{l}{2}{2}\above{1pt}{1pt}{}{2}{3}\above{1pt}{1pt}{r}{2}$}%
\nopagebreak\par%
\nopagebreak\par\leavevmode%
{$\left[\!\llap{\phantom{%
\begingroup \smaller\smaller\smaller\begin{tabular}{@{}c@{}}%
0\\0\\0
\end{tabular}\endgroup%
}}\right.$}%
\begingroup \smaller\smaller\smaller\begin{tabular}{@{}c@{}}%
-3984\\144\\144
\end{tabular}\endgroup%
\kern3pt%
\begingroup \smaller\smaller\smaller\begin{tabular}{@{}c@{}}%
144\\-2\\-6
\end{tabular}\endgroup%
\kern3pt%
\begingroup \smaller\smaller\smaller\begin{tabular}{@{}c@{}}%
144\\-6\\-5
\end{tabular}\endgroup%
{$\left.\llap{\phantom{%
\begingroup \smaller\smaller\smaller\begin{tabular}{@{}c@{}}%
0\\0\\0
\end{tabular}\endgroup%
}}\!\right]$}%
\EasyButWeakLineBreak%
{$\left[\!\llap{\phantom{%
\begingroup \smaller\smaller\smaller\begin{tabular}{@{}c@{}}%
0\\0\\0
\end{tabular}\endgroup%
}}\right.$}%
\begingroup \smaller\smaller\smaller\begin{tabular}{@{}c@{}}%
-1\\-4\\-24
\end{tabular}\endgroup%
\HardButStrongLineBreak\kern3pt%
\begingroup \smaller\smaller\smaller\begin{tabular}{@{}c@{}}%
-1\\-6\\-22
\end{tabular}\endgroup%
\HardButStrongLineBreak\kern3pt%
\begingroup \smaller\smaller\smaller\begin{tabular}{@{}c@{}}%
1\\0\\24
\end{tabular}\endgroup%
\HardButStrongLineBreak\kern3pt%
\begingroup \smaller\smaller\smaller\begin{tabular}{@{}c@{}}%
1\\5\\22
\end{tabular}\endgroup%
\HardButStrongLineBreak\kern3pt%
\begingroup \smaller\smaller\smaller\begin{tabular}{@{}c@{}}%
1\\6\\21
\end{tabular}\endgroup%
{$\left.\llap{\phantom{%
\begingroup \smaller\smaller\smaller\begin{tabular}{@{}c@{}}%
0\\0\\0
\end{tabular}\endgroup%
}}\!\right]$}%
%
%
\hbox{}\par\smallskip%
%
%
\leavevmode%
${L_{151.3}}$%
{} : {$[1\above{1pt}{1pt}{-}{}2\above{1pt}{1pt}{1}{}]\above{1pt}{1pt}{}{2}16\above{1pt}{1pt}{1}{1}{\cdot}1\above{1pt}{1pt}{2}{}3\above{1pt}{1pt}{1}{}$}\spacer%
\instructions{m}%
\EasyButWeakLineBreak%
{${16}\above{1pt}{1pt}{}{2}{1}\above{1pt}{1pt}{r}{2}{48}\above{1pt}{1pt}{*}{2}{8}\above{1pt}{1pt}{l}{2}{3}\above{1pt}{1pt}{}{2}$}%
\nopagebreak\par%
\nopagebreak\par\leavevmode%
{$\left[\!\llap{\phantom{%
\begingroup \smaller\smaller\smaller\begin{tabular}{@{}c@{}}%
0\\0\\0
\end{tabular}\endgroup%
}}\right.$}%
\begingroup \smaller\smaller\smaller\begin{tabular}{@{}c@{}}%
-6000\\144\\288
\end{tabular}\endgroup%
\kern3pt%
\begingroup \smaller\smaller\smaller\begin{tabular}{@{}c@{}}%
144\\-2\\-8
\end{tabular}\endgroup%
\kern3pt%
\begingroup \smaller\smaller\smaller\begin{tabular}{@{}c@{}}%
288\\-8\\-13
\end{tabular}\endgroup%
{$\left.\llap{\phantom{%
\begingroup \smaller\smaller\smaller\begin{tabular}{@{}c@{}}%
0\\0\\0
\end{tabular}\endgroup%
}}\!\right]$}%
\EasyButWeakLineBreak%
{$\left[\!\llap{\phantom{%
\begingroup \smaller\smaller\smaller\begin{tabular}{@{}c@{}}%
0\\0\\0
\end{tabular}\endgroup%
}}\right.$}%
\begingroup \smaller\smaller\smaller\begin{tabular}{@{}c@{}}%
1\\8\\16
\end{tabular}\endgroup%
\HardButStrongLineBreak\kern3pt%
\begingroup \smaller\smaller\smaller\begin{tabular}{@{}c@{}}%
1\\11\\15
\end{tabular}\endgroup%
\HardButStrongLineBreak\kern3pt%
\begingroup \smaller\smaller\smaller\begin{tabular}{@{}c@{}}%
5\\60\\72
\end{tabular}\endgroup%
\HardButStrongLineBreak\kern3pt%
\begingroup \smaller\smaller\smaller\begin{tabular}{@{}c@{}}%
-1\\-10\\-16
\end{tabular}\endgroup%
\HardButStrongLineBreak\kern3pt%
\begingroup \smaller\smaller\smaller\begin{tabular}{@{}c@{}}%
-1\\-12\\-15
\end{tabular}\endgroup%
{$\left.\llap{\phantom{%
\begingroup \smaller\smaller\smaller\begin{tabular}{@{}c@{}}%
0\\0\\0
\end{tabular}\endgroup%
}}\!\right]$}%
%
%
\hbox{}\par\smallskip%
%
%
\leavevmode%
${L_{151.4}}$%
{} : {$[1\above{1pt}{1pt}{1}{}2\above{1pt}{1pt}{1}{}]\above{1pt}{1pt}{}{0}16\above{1pt}{1pt}{-}{3}{\cdot}1\above{1pt}{1pt}{2}{}3\above{1pt}{1pt}{1}{}$}\spacer%
\instructions{m}%
\EasyButWeakLineBreak%
{${16}\above{1pt}{1pt}{l}{2}{1}\above{1pt}{1pt}{}{2}{48}\above{1pt}{1pt}{}{2}{2}\above{1pt}{1pt}{r}{2}{12}\above{1pt}{1pt}{*}{2}$}%
\nopagebreak\par%
\nopagebreak\par\leavevmode%
{$\left[\!\llap{\phantom{%
\begingroup \smaller\smaller\smaller\begin{tabular}{@{}c@{}}%
0\\0\\0
\end{tabular}\endgroup%
}}\right.$}%
\begingroup \smaller\smaller\smaller\begin{tabular}{@{}c@{}}%
48\\0\\48
\end{tabular}\endgroup%
\kern3pt%
\begingroup \smaller\smaller\smaller\begin{tabular}{@{}c@{}}%
0\\-2\\-12
\end{tabular}\endgroup%
\kern3pt%
\begingroup \smaller\smaller\smaller\begin{tabular}{@{}c@{}}%
48\\-12\\-23
\end{tabular}\endgroup%
{$\left.\llap{\phantom{%
\begingroup \smaller\smaller\smaller\begin{tabular}{@{}c@{}}%
0\\0\\0
\end{tabular}\endgroup%
}}\!\right]$}%
\EasyButWeakLineBreak%
{$\left[\!\llap{\phantom{%
\begingroup \smaller\smaller\smaller\begin{tabular}{@{}c@{}}%
0\\0\\0
\end{tabular}\endgroup%
}}\right.$}%
\begingroup \smaller\smaller\smaller\begin{tabular}{@{}c@{}}%
-1\\-4\\0
\end{tabular}\endgroup%
\HardButStrongLineBreak\kern3pt%
\begingroup \smaller\smaller\smaller\begin{tabular}{@{}c@{}}%
-1\\-6\\1
\end{tabular}\endgroup%
\HardButStrongLineBreak\kern3pt%
\begingroup \smaller\smaller\smaller\begin{tabular}{@{}c@{}}%
1\\0\\0
\end{tabular}\endgroup%
\HardButStrongLineBreak\kern3pt%
\begingroup \smaller\smaller\smaller\begin{tabular}{@{}c@{}}%
2\\11\\-2
\end{tabular}\endgroup%
\HardButStrongLineBreak\kern3pt%
\begingroup \smaller\smaller\smaller\begin{tabular}{@{}c@{}}%
5\\30\\-6
\end{tabular}\endgroup%
{$\left.\llap{\phantom{%
\begingroup \smaller\smaller\smaller\begin{tabular}{@{}c@{}}%
0\\0\\0
\end{tabular}\endgroup%
}}\!\right]$}%
%
%
\hbox{}\par\smallskip%
%
%
\leavevmode%
${L_{151.5}}$%
{} : {$[1\above{1pt}{1pt}{1}{}2\above{1pt}{1pt}{1}{}]\above{1pt}{1pt}{}{6}16\above{1pt}{1pt}{-}{5}{\cdot}1\above{1pt}{1pt}{2}{}3\above{1pt}{1pt}{1}{}$}\EasyButWeakLineBreak%
{${16}\above{1pt}{1pt}{s}{2}{4}\above{1pt}{1pt}{*}{2}{48}\above{1pt}{1pt}{s}{2}{8}\above{1pt}{1pt}{*}{2}{12}\above{1pt}{1pt}{s}{2}$}%
\nopagebreak\par%
\nopagebreak\par\leavevmode%
{$\left[\!\llap{\phantom{%
\begingroup \smaller\smaller\smaller\begin{tabular}{@{}c@{}}%
0\\0\\0
\end{tabular}\endgroup%
}}\right.$}%
\begingroup \smaller\smaller\smaller\begin{tabular}{@{}c@{}}%
-816\\48\\48
\end{tabular}\endgroup%
\kern3pt%
\begingroup \smaller\smaller\smaller\begin{tabular}{@{}c@{}}%
48\\-2\\-4
\end{tabular}\endgroup%
\kern3pt%
\begingroup \smaller\smaller\smaller\begin{tabular}{@{}c@{}}%
48\\-4\\-1
\end{tabular}\endgroup%
{$\left.\llap{\phantom{%
\begingroup \smaller\smaller\smaller\begin{tabular}{@{}c@{}}%
0\\0\\0
\end{tabular}\endgroup%
}}\!\right]$}%
\EasyButWeakLineBreak%
{$\left[\!\llap{\phantom{%
\begingroup \smaller\smaller\smaller\begin{tabular}{@{}c@{}}%
0\\0\\0
\end{tabular}\endgroup%
}}\right.$}%
\begingroup \smaller\smaller\smaller\begin{tabular}{@{}c@{}}%
1\\8\\8
\end{tabular}\endgroup%
\HardButStrongLineBreak\kern3pt%
\begingroup \smaller\smaller\smaller\begin{tabular}{@{}c@{}}%
1\\10\\6
\end{tabular}\endgroup%
\HardButStrongLineBreak\kern3pt%
\begingroup \smaller\smaller\smaller\begin{tabular}{@{}c@{}}%
1\\12\\0
\end{tabular}\endgroup%
\HardButStrongLineBreak\kern3pt%
\begingroup \smaller\smaller\smaller\begin{tabular}{@{}c@{}}%
-1\\-10\\-8
\end{tabular}\endgroup%
\HardButStrongLineBreak\kern3pt%
\begingroup \smaller\smaller\smaller\begin{tabular}{@{}c@{}}%
-1\\-12\\-6
\end{tabular}\endgroup%
{$\left.\llap{\phantom{%
\begingroup \smaller\smaller\smaller\begin{tabular}{@{}c@{}}%
0\\0\\0
\end{tabular}\endgroup%
}}\!\right]$}%

\medskip%
%
\leavevmode\llap{}%
$W_{152}$%
\qquad\llap{22} lattices, $\chi=90$%
\hfill%
$4\infty\infty2\infty24\infty\infty2\infty2\rtimes C_{2}$%
\nopagebreak\smallskip\hrule\nopagebreak\medskip%
%
%
\leavevmode%
${L_{152.1}}$%
{} : {$1\above{1pt}{1pt}{2}{{\rm II}}4\above{1pt}{1pt}{1}{1}{\cdot}1\above{1pt}{1pt}{2}{}25\above{1pt}{1pt}{1}{}$}\spacer%
\instructions{2}%
\EasyButWeakLineBreak%
{${2}\above{1pt}{1pt}{*}{4}{4}\above{1pt}{1pt}{5,4}{\infty b}{4}\above{1pt}{1pt}{5,3}{\infty}{4}\above{1pt}{1pt}{*}{2}{100}\above{1pt}{1pt}{1,0}{\infty b}{100}\above{1pt}{1pt}{r}{2}$}\relax$\,(\times2)$%
\nopagebreak\par%
\nopagebreak\par\leavevmode%
{$\left[\!\llap{\phantom{%
\begingroup \smaller\smaller\smaller
\endgroup%
}}\!\right]$}%
%
%
\hbox{}\par\smallskip%
%
%
\leavevmode%
${L_{152.2}}$%
{} : {$1\above{1pt}{1pt}{2}{2}8\above{1pt}{1pt}{1}{7}{\cdot}1\above{1pt}{1pt}{2}{}25\above{1pt}{1pt}{-}{}$}\spacer%
\instructions{2}%
\EasyButWeakLineBreak%
{${1}\above{1pt}{1pt}{}{4}{2}\above{1pt}{1pt}{20,19}{\infty b}{8}\above{1pt}{1pt}{20,3}{\infty z}{2}\above{1pt}{1pt}{s}{2}{50}\above{1pt}{1pt}{4,3}{\infty b}{200}\above{1pt}{1pt}{l}{2}$}\relax$\,(\times2)$%
\nopagebreak\par%
\nopagebreak\par\leavevmode%
{$\left[\!\llap{\phantom{%
\begingroup \smaller\smaller\smaller
\endgroup%
}}\!\right]$}%
%
%
\hbox{}\par\smallskip%
%
%
\leavevmode%
${L_{152.3}}$%
{} : {$1\above{1pt}{1pt}{-2}{2}8\above{1pt}{1pt}{-}{3}{\cdot}1\above{1pt}{1pt}{2}{}25\above{1pt}{1pt}{-}{}$}\spacer%
\instructions{m}%
\EasyButWeakLineBreak%
{${4}\above{1pt}{1pt}{*}{4}{2}\above{1pt}{1pt}{20,19}{\infty a}{8}\above{1pt}{1pt}{20,13}{\infty z}{2}\above{1pt}{1pt}{b}{2}{50}\above{1pt}{1pt}{4,3}{\infty a}{200}\above{1pt}{1pt}{s}{2}$}\relax$\,(\times2)$%
\nopagebreak\par%
\nopagebreak\par\leavevmode%
{$\left[\!\llap{\phantom{%
\begingroup \smaller\smaller\smaller
\endgroup%
}}\!\right]$}%

\medskip%
%
\leavevmode\llap{}%
$W_{153}$%
\qquad\llap{22} lattices, $\chi=60$%
\hfill%
$\infty222\infty\infty222\infty\rtimes C_{2}$%
\nopagebreak\smallskip\hrule\nopagebreak\medskip%
%
%
\leavevmode%
${L_{153.1}}$%
{} : {$1\above{1pt}{1pt}{2}{{\rm II}}4\above{1pt}{1pt}{1}{1}{\cdot}1\above{1pt}{1pt}{-2}{}25\above{1pt}{1pt}{-}{}$}\spacer%
\instructions{2}%
\EasyButWeakLineBreak%
{${4}\above{1pt}{1pt}{5,2}{\infty b}{4}\above{1pt}{1pt}{r}{2}{50}\above{1pt}{1pt}{b}{2}{2}\above{1pt}{1pt}{l}{2}{4}\above{1pt}{1pt}{5,4}{\infty}$}\relax$\,(\times2)$%
\nopagebreak\par%
\nopagebreak\par\leavevmode%
{$\left[\!\llap{\phantom{%
\begingroup \smaller\smaller\smaller
\endgroup%
}}\!\right]$}%
%
%
\hbox{}\par\smallskip%
%
%
\leavevmode%
${L_{153.2}}$%
{} : {$1\above{1pt}{1pt}{-2}{2}8\above{1pt}{1pt}{-}{3}{\cdot}1\above{1pt}{1pt}{-2}{}25\above{1pt}{1pt}{1}{}$}\spacer%
\instructions{2}%
\EasyButWeakLineBreak%
{${2}\above{1pt}{1pt}{20,7}{\infty b}{8}\above{1pt}{1pt}{s}{2}{100}\above{1pt}{1pt}{*}{2}{4}\above{1pt}{1pt}{s}{2}{8}\above{1pt}{1pt}{20,9}{\infty z}$}\relax$\,(\times2)$%
\nopagebreak\par%
\nopagebreak\par\leavevmode%
{$\left[\!\llap{\phantom{%
\begingroup \smaller\smaller\smaller
\endgroup%
}}\!\right]$}%
%
%
\hbox{}\par\smallskip%
%
%
\leavevmode%
${L_{153.3}}$%
{} : {$1\above{1pt}{1pt}{2}{2}8\above{1pt}{1pt}{1}{7}{\cdot}1\above{1pt}{1pt}{-2}{}25\above{1pt}{1pt}{1}{}$}\spacer%
\instructions{m}%
\EasyButWeakLineBreak%
{${2}\above{1pt}{1pt}{20,7}{\infty a}{8}\above{1pt}{1pt}{l}{2}{25}\above{1pt}{1pt}{}{2}{1}\above{1pt}{1pt}{r}{2}{8}\above{1pt}{1pt}{20,19}{\infty z}$}\relax$\,(\times2)$%
\nopagebreak\par%
\nopagebreak\par\leavevmode%
{$\left[\!\llap{\phantom{%
\begingroup \smaller\smaller\smaller
\endgroup%
}}\!\right]$}%

\medskip%
%
\leavevmode\llap{}%
$W_{154}$%
\qquad\llap{3} lattices, $\chi=2$%
\hfill%
$\slashthree2\slashtwo2\rtimes D_{2}$%
\nopagebreak\smallskip\hrule\nopagebreak\medskip%
%
%
\leavevmode%
${L_{154.1}}$%
{} : {$1\above{1pt}{1pt}{-2}{{\rm II}}4\above{1pt}{1pt}{1}{7}{\cdot}1\above{1pt}{1pt}{1}{}3\above{1pt}{1pt}{-}{}9\above{1pt}{1pt}{1}{}$}\spacer%
\instructions{2}%
\EasyButWeakLineBreak%
{${6}\above{1pt}{1pt}{+}{3}{6}\above{1pt}{1pt}{b}{2}{4}\above{1pt}{1pt}{*}{2}{36}\above{1pt}{1pt}{b}{2}$}%
\nopagebreak\par%
\nopagebreak\par\leavevmode%
{$\left[\!\llap{\phantom{%
\begingroup \smaller\smaller\smaller\begin{tabular}{@{}c@{}}%
0\\0\\0
\end{tabular}\endgroup%
}}\right.$}%
\begingroup \smaller\smaller\smaller\begin{tabular}{@{}c@{}}%
252\\72\\-36
\end{tabular}\endgroup%
\kern3pt%
\begingroup \smaller\smaller\smaller\begin{tabular}{@{}c@{}}%
72\\6\\-27
\end{tabular}\endgroup%
\kern3pt%
\begingroup \smaller\smaller\smaller\begin{tabular}{@{}c@{}}%
-36\\-27\\-14
\end{tabular}\endgroup%
{$\left.\llap{\phantom{%
\begingroup \smaller\smaller\smaller\begin{tabular}{@{}c@{}}%
0\\0\\0
\end{tabular}\endgroup%
}}\!\right]$}%
\EasyButWeakLineBreak%
{$\left[\!\llap{\phantom{%
\begingroup \smaller\smaller\smaller\begin{tabular}{@{}c@{}}%
0\\0\\0
\end{tabular}\endgroup%
}}\right.$}%
\begingroup \smaller\smaller\smaller\begin{tabular}{@{}c@{}}%
3\\-7\\6
\end{tabular}\endgroup%
\HardButStrongLineBreak\kern3pt%
\begingroup \smaller\smaller\smaller\begin{tabular}{@{}c@{}}%
-7\\17\\-15
\end{tabular}\endgroup%
\HardButStrongLineBreak\kern3pt%
\begingroup \smaller\smaller\smaller\begin{tabular}{@{}c@{}}%
-1\\2\\-2
\end{tabular}\endgroup%
\HardButStrongLineBreak\kern3pt%
\begingroup \smaller\smaller\smaller\begin{tabular}{@{}c@{}}%
17\\-42\\36
\end{tabular}\endgroup%
{$\left.\llap{\phantom{%
\begingroup \smaller\smaller\smaller\begin{tabular}{@{}c@{}}%
0\\0\\0
\end{tabular}\endgroup%
}}\!\right]$}%
%
%
%
%
%
%
%
%
%
%

\medskip%
%
\leavevmode\llap{}%
$W_{155}$%
\qquad\llap{3} lattices, $\chi=2$%
\hfill%
$6\slashtwo6|\rtimes D_{2}$%
\nopagebreak\smallskip\hrule\nopagebreak\medskip%
%
%
\leavevmode%
${L_{155.1}}$%
{} : {$1\above{1pt}{1pt}{-2}{{\rm II}}4\above{1pt}{1pt}{1}{7}{\cdot}1\above{1pt}{1pt}{-}{}3\above{1pt}{1pt}{-}{}9\above{1pt}{1pt}{-}{}$}\spacer%
\instructions{2}%
\EasyButWeakLineBreak%
{${6}\above{1pt}{1pt}{}{6}{18}\above{1pt}{1pt}{b}{2}{2}\above{1pt}{1pt}{}{6}$}%
\nopagebreak\par%
\nopagebreak\par\leavevmode%
{$\left[\!\llap{\phantom{%
\begingroup \smaller\smaller\smaller\begin{tabular}{@{}c@{}}%
0\\0\\0
\end{tabular}\endgroup%
}}\right.$}%
\begingroup \smaller\smaller\smaller\begin{tabular}{@{}c@{}}%
-36\\36\\36
\end{tabular}\endgroup%
\kern3pt%
\begingroup \smaller\smaller\smaller\begin{tabular}{@{}c@{}}%
36\\-30\\-27
\end{tabular}\endgroup%
\kern3pt%
\begingroup \smaller\smaller\smaller\begin{tabular}{@{}c@{}}%
36\\-27\\-22
\end{tabular}\endgroup%
{$\left.\llap{\phantom{%
\begingroup \smaller\smaller\smaller\begin{tabular}{@{}c@{}}%
0\\0\\0
\end{tabular}\endgroup%
}}\!\right]$}%
\EasyButWeakLineBreak%
{$\left[\!\llap{\phantom{%
\begingroup \smaller\smaller\smaller\begin{tabular}{@{}c@{}}%
0\\0\\0
\end{tabular}\endgroup%
}}\right.$}%
\begingroup \smaller\smaller\smaller\begin{tabular}{@{}c@{}}%
2\\5\\-3
\end{tabular}\endgroup%
\HardButStrongLineBreak\kern3pt%
\begingroup \smaller\smaller\smaller\begin{tabular}{@{}c@{}}%
-2\\-3\\0
\end{tabular}\endgroup%
\HardButStrongLineBreak\kern3pt%
\begingroup \smaller\smaller\smaller\begin{tabular}{@{}c@{}}%
-1\\-3\\2
\end{tabular}\endgroup%
{$\left.\llap{\phantom{%
\begingroup \smaller\smaller\smaller\begin{tabular}{@{}c@{}}%
0\\0\\0
\end{tabular}\endgroup%
}}\!\right]$}%
%
%
%
%
%
%
%
%
%
%

\medskip%
%
\leavevmode\llap{}%
$W_{156}$%
\qquad\llap{6} lattices, $\chi=18$%
\hfill%
$242242\rtimes C_{2}$%
\nopagebreak\smallskip\hrule\nopagebreak\medskip%
%
%
\leavevmode%
${L_{156.1}}$%
{} : {$1\above{1pt}{1pt}{-2}{{\rm II}}4\above{1pt}{1pt}{1}{7}{\cdot}1\above{1pt}{1pt}{2}{}27\above{1pt}{1pt}{-}{}$}\spacer%
\instructions{2}%
\EasyButWeakLineBreak%
{${54}\above{1pt}{1pt}{b}{2}{4}\above{1pt}{1pt}{*}{4}{2}\above{1pt}{1pt}{s}{2}$}\relax$\,(\times2)$%
\nopagebreak\par%
\nopagebreak\par\leavevmode%
{$\left[\!\llap{\phantom{%
\begingroup \smaller\smaller\smaller\begin{tabular}{@{}c@{}}%
0\\0\\0
\end{tabular}\endgroup%
}}\right.$}%
\begingroup \smaller\smaller\smaller\begin{tabular}{@{}c@{}}%
-1352484\\8316\\11556
\end{tabular}\endgroup%
\kern3pt%
\begingroup \smaller\smaller\smaller\begin{tabular}{@{}c@{}}%
8316\\-46\\-73
\end{tabular}\endgroup%
\kern3pt%
\begingroup \smaller\smaller\smaller\begin{tabular}{@{}c@{}}%
11556\\-73\\-98
\end{tabular}\endgroup%
{$\left.\llap{\phantom{%
\begingroup \smaller\smaller\smaller\begin{tabular}{@{}c@{}}%
0\\0\\0
\end{tabular}\endgroup%
}}\!\right]$}%
\hfil\penalty500%
{$\left[\!\llap{\phantom{%
\begingroup \smaller\smaller\smaller\begin{tabular}{@{}c@{}}%
0\\0\\0
\end{tabular}\endgroup%
}}\right.$}%
\begingroup \smaller\smaller\smaller\begin{tabular}{@{}c@{}}%
38879\\1356480\\3572640
\end{tabular}\endgroup%
\kern3pt%
\begingroup \smaller\smaller\smaller\begin{tabular}{@{}c@{}}%
-261\\-9107\\-23983
\end{tabular}\endgroup%
\kern3pt%
\begingroup \smaller\smaller\smaller\begin{tabular}{@{}c@{}}%
-324\\-11304\\-29773
\end{tabular}\endgroup%
{$\left.\llap{\phantom{%
\begingroup \smaller\smaller\smaller\begin{tabular}{@{}c@{}}%
0\\0\\0
\end{tabular}\endgroup%
}}\!\right]$}%
\EasyButWeakLineBreak%
{$\left[\!\llap{\phantom{%
\begingroup \smaller\smaller\smaller\begin{tabular}{@{}c@{}}%
0\\0\\0
\end{tabular}\endgroup%
}}\right.$}%
\begingroup \smaller\smaller\smaller\begin{tabular}{@{}c@{}}%
-10\\-351\\-918
\end{tabular}\endgroup%
\HardButStrongLineBreak\kern3pt%
\begingroup \smaller\smaller\smaller\begin{tabular}{@{}c@{}}%
-5\\-174\\-460
\end{tabular}\endgroup%
\HardButStrongLineBreak\kern3pt%
\begingroup \smaller\smaller\smaller\begin{tabular}{@{}c@{}}%
5\\175\\459
\end{tabular}\endgroup%
{$\left.\llap{\phantom{%
\begingroup \smaller\smaller\smaller\begin{tabular}{@{}c@{}}%
0\\0\\0
\end{tabular}\endgroup%
}}\!\right]$}%

\medskip%
%
\leavevmode\llap{}%
$W_{157}$%
\qquad\llap{8} lattices, $\chi=48$%
\hfill%
$\infty2|2\infty|\infty2|2\infty|\rtimes D_{4}$%
\nopagebreak\smallskip\hrule\nopagebreak\medskip%
%
%
\leavevmode%
${L_{157.1}}$%
{} : {$1\above{1pt}{1pt}{2}{0}4\above{1pt}{1pt}{1}{7}{\cdot}1\above{1pt}{1pt}{-2}{}7\above{1pt}{1pt}{1}{}$}\EasyButWeakLineBreak%
{${28}\above{1pt}{1pt}{1,0}{\infty z}{7}\above{1pt}{1pt}{}{2}{1}\above{1pt}{1pt}{}{2}{28}\above{1pt}{1pt}{2,1}{\infty}$}\relax$\,(\times2)$%
\nopagebreak\par%
\nopagebreak\par\leavevmode%
{$\left[\!\llap{\phantom{%
\begingroup \smaller\smaller\smaller\begin{tabular}{@{}c@{}}%
0\\0\\0
\end{tabular}\endgroup%
}}\right.$}%
\begingroup \smaller\smaller\smaller\begin{tabular}{@{}c@{}}%
-868\\-56\\-728
\end{tabular}\endgroup%
\kern3pt%
\begingroup \smaller\smaller\smaller\begin{tabular}{@{}c@{}}%
-56\\0\\-25
\end{tabular}\endgroup%
\kern3pt%
\begingroup \smaller\smaller\smaller\begin{tabular}{@{}c@{}}%
-728\\-25\\-477
\end{tabular}\endgroup%
{$\left.\llap{\phantom{%
\begingroup \smaller\smaller\smaller\begin{tabular}{@{}c@{}}%
0\\0\\0
\end{tabular}\endgroup%
}}\!\right]$}%
\hfil\penalty500%
{$\left[\!\llap{\phantom{%
\begingroup \smaller\smaller\smaller\begin{tabular}{@{}c@{}}%
0\\0\\0
\end{tabular}\endgroup%
}}\right.$}%
\begingroup \smaller\smaller\smaller\begin{tabular}{@{}c@{}}%
377\\5124\\-840
\end{tabular}\endgroup%
\kern3pt%
\begingroup \smaller\smaller\smaller\begin{tabular}{@{}c@{}}%
18\\243\\-40
\end{tabular}\endgroup%
\kern3pt%
\begingroup \smaller\smaller\smaller\begin{tabular}{@{}c@{}}%
279\\3782\\-621
\end{tabular}\endgroup%
{$\left.\llap{\phantom{%
\begingroup \smaller\smaller\smaller\begin{tabular}{@{}c@{}}%
0\\0\\0
\end{tabular}\endgroup%
}}\!\right]$}%
\EasyButWeakLineBreak%
{$\left[\!\llap{\phantom{%
\begingroup \smaller\smaller\smaller\begin{tabular}{@{}c@{}}%
0\\0\\0
\end{tabular}\endgroup%
}}\right.$}%
\begingroup \smaller\smaller\smaller\begin{tabular}{@{}c@{}}%
-25\\-350\\56
\end{tabular}\endgroup%
\HardButStrongLineBreak\kern3pt%
\begingroup \smaller\smaller\smaller\begin{tabular}{@{}c@{}}%
-47\\-644\\105
\end{tabular}\endgroup%
\HardButStrongLineBreak\kern3pt%
\begingroup \smaller\smaller\smaller\begin{tabular}{@{}c@{}}%
-13\\-177\\29
\end{tabular}\endgroup%
\HardButStrongLineBreak\kern3pt%
\begingroup \smaller\smaller\smaller\begin{tabular}{@{}c@{}}%
-151\\-2044\\336
\end{tabular}\endgroup%
{$\left.\llap{\phantom{%
\begingroup \smaller\smaller\smaller\begin{tabular}{@{}c@{}}%
0\\0\\0
\end{tabular}\endgroup%
}}\!\right]$}%
%
%
\hbox{}\par\smallskip%
%
%
\leavevmode%
${L_{157.2}}$%
{} : {$1\above{1pt}{1pt}{2}{{\rm II}}8\above{1pt}{1pt}{1}{7}{\cdot}1\above{1pt}{1pt}{-2}{}7\above{1pt}{1pt}{1}{}$}\EasyButWeakLineBreak%
{${56}\above{1pt}{1pt}{1,0}{\infty a}{56}\above{1pt}{1pt}{r}{2}{2}\above{1pt}{1pt}{s}{2}{14}\above{1pt}{1pt}{4,3}{\infty b}$}\relax$\,(\times2)$%
\nopagebreak\par%
\nopagebreak\par\leavevmode%
{$\left[\!\llap{\phantom{%
\begingroup \smaller\smaller\smaller\begin{tabular}{@{}c@{}}%
0\\0\\0
\end{tabular}\endgroup%
}}\right.$}%
\begingroup \smaller\smaller\smaller\begin{tabular}{@{}c@{}}%
-1736\\-672\\-56
\end{tabular}\endgroup%
\kern3pt%
\begingroup \smaller\smaller\smaller\begin{tabular}{@{}c@{}}%
-672\\-254\\-25
\end{tabular}\endgroup%
\kern3pt%
\begingroup \smaller\smaller\smaller\begin{tabular}{@{}c@{}}%
-56\\-25\\0
\end{tabular}\endgroup%
{$\left.\llap{\phantom{%
\begingroup \smaller\smaller\smaller\begin{tabular}{@{}c@{}}%
0\\0\\0
\end{tabular}\endgroup%
}}\!\right]$}%
\hfil\penalty500%
{$\left[\!\llap{\phantom{%
\begingroup \smaller\smaller\smaller\begin{tabular}{@{}c@{}}%
0\\0\\0
\end{tabular}\endgroup%
}}\right.$}%
\begingroup \smaller\smaller\smaller\begin{tabular}{@{}c@{}}%
377\\-840\\-1512
\end{tabular}\endgroup%
\kern3pt%
\begingroup \smaller\smaller\smaller\begin{tabular}{@{}c@{}}%
153\\-341\\-612
\end{tabular}\endgroup%
\kern3pt%
\begingroup \smaller\smaller\smaller\begin{tabular}{@{}c@{}}%
9\\-20\\-37
\end{tabular}\endgroup%
{$\left.\llap{\phantom{%
\begingroup \smaller\smaller\smaller\begin{tabular}{@{}c@{}}%
0\\0\\0
\end{tabular}\endgroup%
}}\!\right]$}%
\EasyButWeakLineBreak%
{$\left[\!\llap{\phantom{%
\begingroup \smaller\smaller\smaller\begin{tabular}{@{}c@{}}%
0\\0\\0
\end{tabular}\endgroup%
}}\right.$}%
\begingroup \smaller\smaller\smaller\begin{tabular}{@{}c@{}}%
-101\\224\\420
\end{tabular}\endgroup%
\HardButStrongLineBreak\kern3pt%
\begingroup \smaller\smaller\smaller\begin{tabular}{@{}c@{}}%
-151\\336\\616
\end{tabular}\endgroup%
\HardButStrongLineBreak\kern3pt%
\begingroup \smaller\smaller\smaller\begin{tabular}{@{}c@{}}%
-13\\29\\52
\end{tabular}\endgroup%
\HardButStrongLineBreak\kern3pt%
\begingroup \smaller\smaller\smaller\begin{tabular}{@{}c@{}}%
-47\\105\\182
\end{tabular}\endgroup%
{$\left.\llap{\phantom{%
\begingroup \smaller\smaller\smaller\begin{tabular}{@{}c@{}}%
0\\0\\0
\end{tabular}\endgroup%
}}\!\right]$}%

\medskip%
%
\leavevmode\llap{}%
$W_{158}$%
\qquad\llap{2} lattices, $\chi=16$%
\hfill%
$3\infty3\infty\rtimes C_{2}$%
\nopagebreak\smallskip\hrule\nopagebreak\medskip%
%
%
\leavevmode%
${L_{158.1}}$%
{} : {$1\above{1pt}{1pt}{-2}{{\rm II}}32\above{1pt}{1pt}{-}{5}$}\EasyButWeakLineBreak%
{${2}\above{1pt}{1pt}{+}{3}{2}\above{1pt}{1pt}{8,5}{\infty a}$}\relax$\,(\times2)$%
\nopagebreak\par%
\nopagebreak\par\leavevmode%
{$\left[\!\llap{\phantom{%
\begingroup \smaller\smaller\smaller\begin{tabular}{@{}c@{}}%
0\\0\\0
\end{tabular}\endgroup%
}}\right.$}%
\begingroup \smaller\smaller\smaller\begin{tabular}{@{}c@{}}%
160\\32\\128
\end{tabular}\endgroup%
\kern3pt%
\begingroup \smaller\smaller\smaller\begin{tabular}{@{}c@{}}%
32\\6\\25
\end{tabular}\endgroup%
\kern3pt%
\begingroup \smaller\smaller\smaller\begin{tabular}{@{}c@{}}%
128\\25\\102
\end{tabular}\endgroup%
{$\left.\llap{\phantom{%
\begingroup \smaller\smaller\smaller\begin{tabular}{@{}c@{}}%
0\\0\\0
\end{tabular}\endgroup%
}}\!\right]$}%
\hfil\penalty500%
{$\left[\!\llap{\phantom{%
\begingroup \smaller\smaller\smaller\begin{tabular}{@{}c@{}}%
0\\0\\0
\end{tabular}\endgroup%
}}\right.$}%
\begingroup \smaller\smaller\smaller\begin{tabular}{@{}c@{}}%
-1\\0\\0
\end{tabular}\endgroup%
\kern3pt%
\begingroup \smaller\smaller\smaller\begin{tabular}{@{}c@{}}%
-1\\0\\1
\end{tabular}\endgroup%
\kern3pt%
\begingroup \smaller\smaller\smaller\begin{tabular}{@{}c@{}}%
-1\\1\\0
\end{tabular}\endgroup%
{$\left.\llap{\phantom{%
\begingroup \smaller\smaller\smaller\begin{tabular}{@{}c@{}}%
0\\0\\0
\end{tabular}\endgroup%
}}\!\right]$}%
\EasyButWeakLineBreak%
{$\left[\!\llap{\phantom{%
\begingroup \smaller\smaller\smaller\begin{tabular}{@{}c@{}}%
0\\0\\0
\end{tabular}\endgroup%
}}\right.$}%
\begingroup \smaller\smaller\smaller\begin{tabular}{@{}c@{}}%
1\\-2\\-1
\end{tabular}\endgroup%
\HardButStrongLineBreak\kern3pt%
\begingroup \smaller\smaller\smaller\begin{tabular}{@{}c@{}}%
-1\\-3\\2
\end{tabular}\endgroup%
{$\left.\llap{\phantom{%
\begingroup \smaller\smaller\smaller\begin{tabular}{@{}c@{}}%
0\\0\\0
\end{tabular}\endgroup%
}}\!\right]$}%
%
%
%
%
%
%

\medskip%
%
\leavevmode\llap{}%
$W_{159}$%
\qquad\llap{2} lattices, $\chi=12$%
\hfill%
$\slashtwo4\slashinfty4\rtimes D_{2}$%
\nopagebreak\smallskip\hrule\nopagebreak\medskip%
%
%
\leavevmode%
${L_{159.1}}$%
{} : {$1\above{1pt}{1pt}{2}{2}32\above{1pt}{1pt}{1}{7}$}\EasyButWeakLineBreak%
{${1}\above{1pt}{1pt}{r}{2}{4}\above{1pt}{1pt}{*}{4}{2}\above{1pt}{1pt}{8,7}{\infty a}{2}\above{1pt}{1pt}{}{4}$}%
\nopagebreak\par%
\nopagebreak\par\leavevmode%
{$\left[\!\llap{\phantom{%
\begingroup \smaller\smaller\smaller\begin{tabular}{@{}c@{}}%
0\\0\\0
\end{tabular}\endgroup%
}}\right.$}%
\begingroup \smaller\smaller\smaller\begin{tabular}{@{}c@{}}%
-47392\\1216\\608
\end{tabular}\endgroup%
\kern3pt%
\begingroup \smaller\smaller\smaller\begin{tabular}{@{}c@{}}%
1216\\-31\\-16
\end{tabular}\endgroup%
\kern3pt%
\begingroup \smaller\smaller\smaller\begin{tabular}{@{}c@{}}%
608\\-16\\-7
\end{tabular}\endgroup%
{$\left.\llap{\phantom{%
\begingroup \smaller\smaller\smaller\begin{tabular}{@{}c@{}}%
0\\0\\0
\end{tabular}\endgroup%
}}\!\right]$}%
\EasyButWeakLineBreak%
{$\left[\!\llap{\phantom{%
\begingroup \smaller\smaller\smaller\begin{tabular}{@{}c@{}}%
0\\0\\0
\end{tabular}\endgroup%
}}\right.$}%
\begingroup \smaller\smaller\smaller\begin{tabular}{@{}c@{}}%
-1\\-31\\-16
\end{tabular}\endgroup%
\HardButStrongLineBreak\kern3pt%
\begingroup \smaller\smaller\smaller\begin{tabular}{@{}c@{}}%
-1\\-32\\-14
\end{tabular}\endgroup%
\HardButStrongLineBreak\kern3pt%
\begingroup \smaller\smaller\smaller\begin{tabular}{@{}c@{}}%
3\\93\\47
\end{tabular}\endgroup%
\HardButStrongLineBreak\kern3pt%
\begingroup \smaller\smaller\smaller\begin{tabular}{@{}c@{}}%
2\\63\\29
\end{tabular}\endgroup%
{$\left.\llap{\phantom{%
\begingroup \smaller\smaller\smaller\begin{tabular}{@{}c@{}}%
0\\0\\0
\end{tabular}\endgroup%
}}\!\right]$}%
%
%
%
%
%
%

\medskip%
%
\leavevmode\llap{}%
$W_{160}$%
\qquad\llap{10} lattices, $\chi=24$%
\hfill%
$2\slashinfty2\infty|\infty\rtimes D_{2}$%
\nopagebreak\smallskip\hrule\nopagebreak\medskip%
%
%
\leavevmode%
${L_{160.1}}$%
{} : {$1\above{1pt}{1pt}{2}{{\rm II}}32\above{1pt}{1pt}{1}{1}$}\EasyButWeakLineBreak%
{${32}\above{1pt}{1pt}{b}{2}{2}\above{1pt}{1pt}{8,1}{\infty b}{2}\above{1pt}{1pt}{l}{2}{32}\above{1pt}{1pt}{1,0}{\infty}{32}\above{1pt}{1pt}{1,0}{\infty z}$}%
\nopagebreak\par%
\nopagebreak\par\leavevmode%
{$\left[\!\llap{\phantom{%
\begingroup \smaller\smaller\smaller\begin{tabular}{@{}c@{}}%
0\\0\\0
\end{tabular}\endgroup%
}}\right.$}%
\begingroup \smaller\smaller\smaller\begin{tabular}{@{}c@{}}%
-26848\\992\\928
\end{tabular}\endgroup%
\kern3pt%
\begingroup \smaller\smaller\smaller\begin{tabular}{@{}c@{}}%
992\\-30\\-35
\end{tabular}\endgroup%
\kern3pt%
\begingroup \smaller\smaller\smaller\begin{tabular}{@{}c@{}}%
928\\-35\\-32
\end{tabular}\endgroup%
{$\left.\llap{\phantom{%
\begingroup \smaller\smaller\smaller\begin{tabular}{@{}c@{}}%
0\\0\\0
\end{tabular}\endgroup%
}}\!\right]$}%
\EasyButWeakLineBreak%
{$\left[\!\llap{\phantom{%
\begingroup \smaller\smaller\smaller\begin{tabular}{@{}c@{}}%
0\\0\\0
\end{tabular}\endgroup%
}}\right.$}%
\begingroup \smaller\smaller\smaller\begin{tabular}{@{}c@{}}%
-5\\-16\\-128
\end{tabular}\endgroup%
\HardButStrongLineBreak\kern3pt%
\begingroup \smaller\smaller\smaller\begin{tabular}{@{}c@{}}%
-4\\-11\\-104
\end{tabular}\endgroup%
\HardButStrongLineBreak\kern3pt%
\begingroup \smaller\smaller\smaller\begin{tabular}{@{}c@{}}%
5\\15\\128
\end{tabular}\endgroup%
\HardButStrongLineBreak\kern3pt%
\begingroup \smaller\smaller\smaller\begin{tabular}{@{}c@{}}%
67\\192\\1728
\end{tabular}\endgroup%
\HardButStrongLineBreak\kern3pt%
\begingroup \smaller\smaller\smaller\begin{tabular}{@{}c@{}}%
29\\80\\752
\end{tabular}\endgroup%
{$\left.\llap{\phantom{%
\begingroup \smaller\smaller\smaller\begin{tabular}{@{}c@{}}%
0\\0\\0
\end{tabular}\endgroup%
}}\!\right]$}%
%
%
\hbox{}\par\smallskip%
%
%
\leavevmode%
${L_{160.2}}$%
{} : {$1\above{1pt}{1pt}{-}{3}4\above{1pt}{1pt}{1}{1}64\above{1pt}{1pt}{-}{5}$}\EasyButWeakLineBreak%
{${16}\above{1pt}{1pt}{s}{2}{64}\above{1pt}{1pt}{8,1}{\infty z}{64}\above{1pt}{1pt}{*}{2}{16}\above{1pt}{1pt}{16,9}{\infty z}{4}\above{1pt}{1pt}{16,5}{\infty}$}%
\nopagebreak\par%
\nopagebreak\par\leavevmode%
{$\left[\!\llap{\phantom{%
\begingroup \smaller\smaller\smaller\begin{tabular}{@{}c@{}}%
0\\0\\0
\end{tabular}\endgroup%
}}\right.$}%
\begingroup \smaller\smaller\smaller\begin{tabular}{@{}c@{}}%
-14528\\-1344\\-1472
\end{tabular}\endgroup%
\kern3pt%
\begingroup \smaller\smaller\smaller\begin{tabular}{@{}c@{}}%
-1344\\-124\\-136
\end{tabular}\endgroup%
\kern3pt%
\begingroup \smaller\smaller\smaller\begin{tabular}{@{}c@{}}%
-1472\\-136\\-149
\end{tabular}\endgroup%
{$\left.\llap{\phantom{%
\begingroup \smaller\smaller\smaller\begin{tabular}{@{}c@{}}%
0\\0\\0
\end{tabular}\endgroup%
}}\!\right]$}%
\EasyButWeakLineBreak%
{$\left[\!\llap{\phantom{%
\begingroup \smaller\smaller\smaller\begin{tabular}{@{}c@{}}%
0\\0\\0
\end{tabular}\endgroup%
}}\right.$}%
\begingroup \smaller\smaller\smaller\begin{tabular}{@{}c@{}}%
-1\\2\\8
\end{tabular}\endgroup%
\HardButStrongLineBreak\kern3pt%
\begingroup \smaller\smaller\smaller\begin{tabular}{@{}c@{}}%
1\\24\\-32
\end{tabular}\endgroup%
\HardButStrongLineBreak\kern3pt%
\begingroup \smaller\smaller\smaller\begin{tabular}{@{}c@{}}%
7\\-8\\-64
\end{tabular}\endgroup%
\HardButStrongLineBreak\kern3pt%
\begingroup \smaller\smaller\smaller\begin{tabular}{@{}c@{}}%
5\\-30\\-24
\end{tabular}\endgroup%
\HardButStrongLineBreak\kern3pt%
\begingroup \smaller\smaller\smaller\begin{tabular}{@{}c@{}}%
0\\-9\\8
\end{tabular}\endgroup%
{$\left.\llap{\phantom{%
\begingroup \smaller\smaller\smaller\begin{tabular}{@{}c@{}}%
0\\0\\0
\end{tabular}\endgroup%
}}\!\right]$}%
%
%
\hbox{}\par\smallskip%
%
%
\leavevmode%
${L_{160.3}}$%
{} : {$1\above{1pt}{1pt}{-}{5}4\above{1pt}{1pt}{1}{1}64\above{1pt}{1pt}{-}{3}$}\EasyButWeakLineBreak%
{${16}\above{1pt}{1pt}{*}{2}{4}\above{1pt}{1pt}{8,3}{\infty z}{1}\above{1pt}{1pt}{r}{2}{16}\above{1pt}{1pt}{16,15}{\infty z}{4}\above{1pt}{1pt}{16,11}{\infty}$}%
\nopagebreak\par%
\nopagebreak\par\leavevmode%
{$\left[\!\llap{\phantom{%
\begingroup \smaller\smaller\smaller\begin{tabular}{@{}c@{}}%
0\\0\\0
\end{tabular}\endgroup%
}}\right.$}%
\begingroup \smaller\smaller\smaller\begin{tabular}{@{}c@{}}%
-21824\\576\\576
\end{tabular}\endgroup%
\kern3pt%
\begingroup \smaller\smaller\smaller\begin{tabular}{@{}c@{}}%
576\\-12\\-16
\end{tabular}\endgroup%
\kern3pt%
\begingroup \smaller\smaller\smaller\begin{tabular}{@{}c@{}}%
576\\-16\\-15
\end{tabular}\endgroup%
{$\left.\llap{\phantom{%
\begingroup \smaller\smaller\smaller\begin{tabular}{@{}c@{}}%
0\\0\\0
\end{tabular}\endgroup%
}}\!\right]$}%
\EasyButWeakLineBreak%
{$\left[\!\llap{\phantom{%
\begingroup \smaller\smaller\smaller\begin{tabular}{@{}c@{}}%
0\\0\\0
\end{tabular}\endgroup%
}}\right.$}%
\begingroup \smaller\smaller\smaller\begin{tabular}{@{}c@{}}%
-1\\-6\\-32
\end{tabular}\endgroup%
\HardButStrongLineBreak\kern3pt%
\begingroup \smaller\smaller\smaller\begin{tabular}{@{}c@{}}%
-1\\-8\\-30
\end{tabular}\endgroup%
\HardButStrongLineBreak\kern3pt%
\begingroup \smaller\smaller\smaller\begin{tabular}{@{}c@{}}%
1\\6\\31
\end{tabular}\endgroup%
\HardButStrongLineBreak\kern3pt%
\begingroup \smaller\smaller\smaller\begin{tabular}{@{}c@{}}%
11\\74\\336
\end{tabular}\endgroup%
\HardButStrongLineBreak\kern3pt%
\begingroup \smaller\smaller\smaller\begin{tabular}{@{}c@{}}%
2\\15\\60
\end{tabular}\endgroup%
{$\left.\llap{\phantom{%
\begingroup \smaller\smaller\smaller\begin{tabular}{@{}c@{}}%
0\\0\\0
\end{tabular}\endgroup%
}}\!\right]$}%
%
%
\hbox{}\par\smallskip%
%
%
\leavevmode%
${L_{160.4}}$%
{} : {$1\above{1pt}{1pt}{1}{1}4\above{1pt}{1pt}{1}{1}64\above{1pt}{1pt}{1}{7}$}\EasyButWeakLineBreak%
{${4}\above{1pt}{1pt}{r}{2}{4}\above{1pt}{1pt}{8,7}{\infty z}{1}\above{1pt}{1pt}{}{2}{4}\above{1pt}{1pt}{16,15}{\infty}{16}\above{1pt}{1pt}{16,3}{\infty z}$}%
\nopagebreak\par%
\nopagebreak\par\leavevmode%
{$\left[\!\llap{\phantom{%
\begingroup \smaller\smaller\smaller\begin{tabular}{@{}c@{}}%
0\\0\\0
\end{tabular}\endgroup%
}}\right.$}%
\begingroup \smaller\smaller\smaller\begin{tabular}{@{}c@{}}%
-13376\\960\\448
\end{tabular}\endgroup%
\kern3pt%
\begingroup \smaller\smaller\smaller\begin{tabular}{@{}c@{}}%
960\\-60\\-32
\end{tabular}\endgroup%
\kern3pt%
\begingroup \smaller\smaller\smaller\begin{tabular}{@{}c@{}}%
448\\-32\\-15
\end{tabular}\endgroup%
{$\left.\llap{\phantom{%
\begingroup \smaller\smaller\smaller\begin{tabular}{@{}c@{}}%
0\\0\\0
\end{tabular}\endgroup%
}}\!\right]$}%
\EasyButWeakLineBreak%
{$\left[\!\llap{\phantom{%
\begingroup \smaller\smaller\smaller\begin{tabular}{@{}c@{}}%
0\\0\\0
\end{tabular}\endgroup%
}}\right.$}%
\begingroup \smaller\smaller\smaller\begin{tabular}{@{}c@{}}%
-1\\1\\-32
\end{tabular}\endgroup%
\HardButStrongLineBreak\kern3pt%
\begingroup \smaller\smaller\smaller\begin{tabular}{@{}c@{}}%
-1\\0\\-30
\end{tabular}\endgroup%
\HardButStrongLineBreak\kern3pt%
\begingroup \smaller\smaller\smaller\begin{tabular}{@{}c@{}}%
2\\-2\\63
\end{tabular}\endgroup%
\HardButStrongLineBreak\kern3pt%
\begingroup \smaller\smaller\smaller\begin{tabular}{@{}c@{}}%
9\\-7\\280
\end{tabular}\endgroup%
\HardButStrongLineBreak\kern3pt%
\begingroup \smaller\smaller\smaller\begin{tabular}{@{}c@{}}%
5\\-2\\152
\end{tabular}\endgroup%
{$\left.\llap{\phantom{%
\begingroup \smaller\smaller\smaller\begin{tabular}{@{}c@{}}%
0\\0\\0
\end{tabular}\endgroup%
}}\!\right]$}%
%
%
\hbox{}\par\smallskip%
%
%
\leavevmode%
${L_{160.5}}$%
{} : {$1\above{1pt}{1pt}{1}{7}4\above{1pt}{1pt}{1}{1}64\above{1pt}{1pt}{1}{1}$}\EasyButWeakLineBreak%
{${4}\above{1pt}{1pt}{r}{2}{64}\above{1pt}{1pt}{2,1}{\infty a}{64}\above{1pt}{1pt}{}{2}{4}\above{1pt}{1pt}{16,1}{\infty}{16}\above{1pt}{1pt}{16,5}{\infty z}$}%
\nopagebreak\par%
\nopagebreak\par\leavevmode%
{$\left[\!\llap{\phantom{%
\begingroup \smaller\smaller\smaller\begin{tabular}{@{}c@{}}%
0\\0\\0
\end{tabular}\endgroup%
}}\right.$}%
\begingroup \smaller\smaller\smaller\begin{tabular}{@{}c@{}}%
64\\0\\0
\end{tabular}\endgroup%
\kern3pt%
\begingroup \smaller\smaller\smaller\begin{tabular}{@{}c@{}}%
0\\-1020\\-32
\end{tabular}\endgroup%
\kern3pt%
\begingroup \smaller\smaller\smaller\begin{tabular}{@{}c@{}}%
0\\-32\\-1
\end{tabular}\endgroup%
{$\left.\llap{\phantom{%
\begingroup \smaller\smaller\smaller\begin{tabular}{@{}c@{}}%
0\\0\\0
\end{tabular}\endgroup%
}}\!\right]$}%
\EasyButWeakLineBreak%
{$\left[\!\llap{\phantom{%
\begingroup \smaller\smaller\smaller\begin{tabular}{@{}c@{}}%
0\\0\\0
\end{tabular}\endgroup%
}}\right.$}%
\begingroup \smaller\smaller\smaller\begin{tabular}{@{}c@{}}%
-1\\-7\\208
\end{tabular}\endgroup%
\HardButStrongLineBreak\kern3pt%
\begingroup \smaller\smaller\smaller\begin{tabular}{@{}c@{}}%
-1\\-16\\480
\end{tabular}\endgroup%
\HardButStrongLineBreak\kern3pt%
\begingroup \smaller\smaller\smaller\begin{tabular}{@{}c@{}}%
1\\0\\0
\end{tabular}\endgroup%
\HardButStrongLineBreak\kern3pt%
\begingroup \smaller\smaller\smaller\begin{tabular}{@{}c@{}}%
0\\1\\-32
\end{tabular}\endgroup%
\HardButStrongLineBreak\kern3pt%
\begingroup \smaller\smaller\smaller\begin{tabular}{@{}c@{}}%
-1\\-2\\56
\end{tabular}\endgroup%
{$\left.\llap{\phantom{%
\begingroup \smaller\smaller\smaller\begin{tabular}{@{}c@{}}%
0\\0\\0
\end{tabular}\endgroup%
}}\!\right]$}%

\medskip%
%
\leavevmode\llap{}%
$W_{161}$%
\qquad\llap{12} lattices, $\chi=24$%
\hfill%
$\slashinfty2|2\slashinfty2|2\rtimes D_{4}$%
\nopagebreak\smallskip\hrule\nopagebreak\medskip%
%
%
\leavevmode%
${L_{161.1}}$%
{} : {$1\above{1pt}{1pt}{2}{0}32\above{1pt}{1pt}{1}{1}$}\EasyButWeakLineBreak%
{${32}\above{1pt}{1pt}{1,0}{\infty b}{32}\above{1pt}{1pt}{}{2}{1}\above{1pt}{1pt}{r}{2}{32}\above{1pt}{1pt}{4,3}{\infty z}{32}\above{1pt}{1pt}{s}{2}{4}\above{1pt}{1pt}{*}{2}$}%
\nopagebreak\par%
\nopagebreak\par\leavevmode%
{$\left[\!\llap{\phantom{%
\begingroup \smaller\smaller\smaller\begin{tabular}{@{}c@{}}%
0\\0\\0
\end{tabular}\endgroup%
}}\right.$}%
\begingroup \smaller\smaller\smaller\begin{tabular}{@{}c@{}}%
-63712\\256\\3968
\end{tabular}\endgroup%
\kern3pt%
\begingroup \smaller\smaller\smaller\begin{tabular}{@{}c@{}}%
256\\-1\\-16
\end{tabular}\endgroup%
\kern3pt%
\begingroup \smaller\smaller\smaller\begin{tabular}{@{}c@{}}%
3968\\-16\\-247
\end{tabular}\endgroup%
{$\left.\llap{\phantom{%
\begingroup \smaller\smaller\smaller\begin{tabular}{@{}c@{}}%
0\\0\\0
\end{tabular}\endgroup%
}}\!\right]$}%
\EasyButWeakLineBreak%
{$\left[\!\llap{\phantom{%
\begingroup \smaller\smaller\smaller\begin{tabular}{@{}c@{}}%
0\\0\\0
\end{tabular}\endgroup%
}}\right.$}%
\begingroup \smaller\smaller\smaller\begin{tabular}{@{}c@{}}%
-3\\-16\\-48
\end{tabular}\endgroup%
\HardButStrongLineBreak\kern3pt%
\begingroup \smaller\smaller\smaller\begin{tabular}{@{}c@{}}%
3\\224\\32
\end{tabular}\endgroup%
\HardButStrongLineBreak\kern3pt%
\begingroup \smaller\smaller\smaller\begin{tabular}{@{}c@{}}%
1\\44\\13
\end{tabular}\endgroup%
\HardButStrongLineBreak\kern3pt%
\begingroup \smaller\smaller\smaller\begin{tabular}{@{}c@{}}%
7\\240\\96
\end{tabular}\endgroup%
\HardButStrongLineBreak\kern3pt%
\begingroup \smaller\smaller\smaller\begin{tabular}{@{}c@{}}%
1\\0\\16
\end{tabular}\endgroup%
\HardButStrongLineBreak\kern3pt%
\begingroup \smaller\smaller\smaller\begin{tabular}{@{}c@{}}%
-1\\-32\\-14
\end{tabular}\endgroup%
{$\left.\llap{\phantom{%
\begingroup \smaller\smaller\smaller\begin{tabular}{@{}c@{}}%
0\\0\\0
\end{tabular}\endgroup%
}}\!\right]$}%
%
%
\hbox{}\par\smallskip%
%
%
\leavevmode%
${L_{161.2}}$%
{} : {$[1\above{1pt}{1pt}{-}{}2\above{1pt}{1pt}{1}{}]\above{1pt}{1pt}{}{2}32\above{1pt}{1pt}{-}{3}$}\EasyButWeakLineBreak%
{${4}\above{1pt}{1pt}{8,3}{\infty z}{1}\above{1pt}{1pt}{r}{2}{8}\above{1pt}{1pt}{*}{2}$}\relax$\,(\times2)$%
\nopagebreak\par%
\nopagebreak\par\leavevmode%
{$\left[\!\llap{\phantom{%
\begingroup \smaller\smaller\smaller\begin{tabular}{@{}c@{}}%
0\\0\\0
\end{tabular}\endgroup%
}}\right.$}%
\begingroup \smaller\smaller\smaller\begin{tabular}{@{}c@{}}%
-3488\\160\\160
\end{tabular}\endgroup%
\kern3pt%
\begingroup \smaller\smaller\smaller\begin{tabular}{@{}c@{}}%
160\\-6\\-8
\end{tabular}\endgroup%
\kern3pt%
\begingroup \smaller\smaller\smaller\begin{tabular}{@{}c@{}}%
160\\-8\\-7
\end{tabular}\endgroup%
{$\left.\llap{\phantom{%
\begingroup \smaller\smaller\smaller\begin{tabular}{@{}c@{}}%
0\\0\\0
\end{tabular}\endgroup%
}}\!\right]$}%
\hfil\penalty500%
{$\left[\!\llap{\phantom{%
\begingroup \smaller\smaller\smaller\begin{tabular}{@{}c@{}}%
0\\0\\0
\end{tabular}\endgroup%
}}\right.$}%
\begingroup \smaller\smaller\smaller\begin{tabular}{@{}c@{}}%
127\\896\\1792
\end{tabular}\endgroup%
\kern3pt%
\begingroup \smaller\smaller\smaller\begin{tabular}{@{}c@{}}%
-6\\-43\\-84
\end{tabular}\endgroup%
\kern3pt%
\begingroup \smaller\smaller\smaller\begin{tabular}{@{}c@{}}%
-6\\-42\\-85
\end{tabular}\endgroup%
{$\left.\llap{\phantom{%
\begingroup \smaller\smaller\smaller\begin{tabular}{@{}c@{}}%
0\\0\\0
\end{tabular}\endgroup%
}}\!\right]$}%
\EasyButWeakLineBreak%
{$\left[\!\llap{\phantom{%
\begingroup \smaller\smaller\smaller\begin{tabular}{@{}c@{}}%
0\\0\\0
\end{tabular}\endgroup%
}}\right.$}%
\begingroup \smaller\smaller\smaller\begin{tabular}{@{}c@{}}%
5\\36\\70
\end{tabular}\endgroup%
\HardButStrongLineBreak\kern3pt%
\begingroup \smaller\smaller\smaller\begin{tabular}{@{}c@{}}%
1\\8\\13
\end{tabular}\endgroup%
\HardButStrongLineBreak\kern3pt%
\begingroup \smaller\smaller\smaller\begin{tabular}{@{}c@{}}%
-1\\-6\\-16
\end{tabular}\endgroup%
{$\left.\llap{\phantom{%
\begingroup \smaller\smaller\smaller\begin{tabular}{@{}c@{}}%
0\\0\\0
\end{tabular}\endgroup%
}}\!\right]$}%
%
%
\hbox{}\par\smallskip%
%
%
\leavevmode%
${L_{161.3}}$%
{} : {$[1\above{1pt}{1pt}{1}{}2\above{1pt}{1pt}{1}{}]\above{1pt}{1pt}{}{2}64\above{1pt}{1pt}{1}{7}$}\spacer%
\instructions{m}%
\EasyButWeakLineBreak%
{${2}\above{1pt}{1pt}{16,7}{\infty}{8}\above{1pt}{1pt}{l}{2}{1}\above{1pt}{1pt}{}{2}{2}\above{1pt}{1pt}{16,15}{\infty}{8}\above{1pt}{1pt}{*}{2}{4}\above{1pt}{1pt}{l}{2}$}%
\nopagebreak\par%
\nopagebreak\par\leavevmode%
{$\left[\!\llap{\phantom{%
\begingroup \smaller\smaller\smaller\begin{tabular}{@{}c@{}}%
0\\0\\0
\end{tabular}\endgroup%
}}\right.$}%
\begingroup \smaller\smaller\smaller\begin{tabular}{@{}c@{}}%
-32320\\704\\704
\end{tabular}\endgroup%
\kern3pt%
\begingroup \smaller\smaller\smaller\begin{tabular}{@{}c@{}}%
704\\-14\\-16
\end{tabular}\endgroup%
\kern3pt%
\begingroup \smaller\smaller\smaller\begin{tabular}{@{}c@{}}%
704\\-16\\-15
\end{tabular}\endgroup%
{$\left.\llap{\phantom{%
\begingroup \smaller\smaller\smaller\begin{tabular}{@{}c@{}}%
0\\0\\0
\end{tabular}\endgroup%
}}\!\right]$}%
\EasyButWeakLineBreak%
{$\left[\!\llap{\phantom{%
\begingroup \smaller\smaller\smaller\begin{tabular}{@{}c@{}}%
0\\0\\0
\end{tabular}\endgroup%
}}\right.$}%
\begingroup \smaller\smaller\smaller\begin{tabular}{@{}c@{}}%
1\\13\\32
\end{tabular}\endgroup%
\HardButStrongLineBreak\kern3pt%
\begingroup \smaller\smaller\smaller\begin{tabular}{@{}c@{}}%
7\\102\\216
\end{tabular}\endgroup%
\HardButStrongLineBreak\kern3pt%
\begingroup \smaller\smaller\smaller\begin{tabular}{@{}c@{}}%
2\\30\\61
\end{tabular}\endgroup%
\HardButStrongLineBreak\kern3pt%
\begingroup \smaller\smaller\smaller\begin{tabular}{@{}c@{}}%
2\\31\\60
\end{tabular}\endgroup%
\HardButStrongLineBreak\kern3pt%
\begingroup \smaller\smaller\smaller\begin{tabular}{@{}c@{}}%
-1\\-14\\-32
\end{tabular}\endgroup%
\HardButStrongLineBreak\kern3pt%
\begingroup \smaller\smaller\smaller\begin{tabular}{@{}c@{}}%
-1\\-16\\-30
\end{tabular}\endgroup%
{$\left.\llap{\phantom{%
\begingroup \smaller\smaller\smaller\begin{tabular}{@{}c@{}}%
0\\0\\0
\end{tabular}\endgroup%
}}\!\right]$}%
%
%
\hbox{}\par\smallskip%
%
%
\leavevmode%
${L_{161.4}}$%
{} : {$1\above{1pt}{1pt}{-}{3}4\above{1pt}{1pt}{1}{7}32\above{1pt}{1pt}{-}{3}$}\EasyButWeakLineBreak%
{${32}\above{1pt}{1pt}{8,5}{\infty z}{32}\above{1pt}{1pt}{s}{2}{16}\above{1pt}{1pt}{*}{2}$}\relax$\,(\times2)$%
\nopagebreak\par%
\nopagebreak\par\leavevmode%
{$\left[\!\llap{\phantom{%
\begingroup \smaller\smaller\smaller\begin{tabular}{@{}c@{}}%
0\\0\\0
\end{tabular}\endgroup%
}}\right.$}%
\begingroup \smaller\smaller\smaller\begin{tabular}{@{}c@{}}%
-4000\\992\\64
\end{tabular}\endgroup%
\kern3pt%
\begingroup \smaller\smaller\smaller\begin{tabular}{@{}c@{}}%
992\\-244\\-16
\end{tabular}\endgroup%
\kern3pt%
\begingroup \smaller\smaller\smaller\begin{tabular}{@{}c@{}}%
64\\-16\\-1
\end{tabular}\endgroup%
{$\left.\llap{\phantom{%
\begingroup \smaller\smaller\smaller\begin{tabular}{@{}c@{}}%
0\\0\\0
\end{tabular}\endgroup%
}}\!\right]$}%
\hfil\penalty500%
{$\left[\!\llap{\phantom{%
\begingroup \smaller\smaller\smaller\begin{tabular}{@{}c@{}}%
0\\0\\0
\end{tabular}\endgroup%
}}\right.$}%
\begingroup \smaller\smaller\smaller\begin{tabular}{@{}c@{}}%
-1\\-112\\1344
\end{tabular}\endgroup%
\kern3pt%
\begingroup \smaller\smaller\smaller\begin{tabular}{@{}c@{}}%
0\\25\\-312
\end{tabular}\endgroup%
\kern3pt%
\begingroup \smaller\smaller\smaller\begin{tabular}{@{}c@{}}%
0\\2\\-25
\end{tabular}\endgroup%
{$\left.\llap{\phantom{%
\begingroup \smaller\smaller\smaller\begin{tabular}{@{}c@{}}%
0\\0\\0
\end{tabular}\endgroup%
}}\!\right]$}%
\EasyButWeakLineBreak%
{$\left[\!\llap{\phantom{%
\begingroup \smaller\smaller\smaller\begin{tabular}{@{}c@{}}%
0\\0\\0
\end{tabular}\endgroup%
}}\right.$}%
\begingroup \smaller\smaller\smaller\begin{tabular}{@{}c@{}}%
-3\\-12\\-16
\end{tabular}\endgroup%
\HardButStrongLineBreak\kern3pt%
\begingroup \smaller\smaller\smaller\begin{tabular}{@{}c@{}}%
-1\\-12\\96
\end{tabular}\endgroup%
\HardButStrongLineBreak\kern3pt%
\begingroup \smaller\smaller\smaller\begin{tabular}{@{}c@{}}%
1\\-2\\80
\end{tabular}\endgroup%
{$\left.\llap{\phantom{%
\begingroup \smaller\smaller\smaller\begin{tabular}{@{}c@{}}%
0\\0\\0
\end{tabular}\endgroup%
}}\!\right]$}%
%
%
\hbox{}\par\smallskip%
%
%
\leavevmode%
${L_{161.5}}$%
{} : {$1\above{1pt}{1pt}{1}{7}4\above{1pt}{1pt}{1}{1}32\above{1pt}{1pt}{1}{1}$}\EasyButWeakLineBreak%
{${32}\above{1pt}{1pt}{2,1}{\infty a}{32}\above{1pt}{1pt}{}{2}{4}\above{1pt}{1pt}{r}{2}$}\relax$\,(\times2)$%
\nopagebreak\par%
\nopagebreak\par\leavevmode%
{$\left[\!\llap{\phantom{%
\begingroup \smaller\smaller\smaller\begin{tabular}{@{}c@{}}%
0\\0\\0
\end{tabular}\endgroup%
}}\right.$}%
\begingroup \smaller\smaller\smaller\begin{tabular}{@{}c@{}}%
32\\0\\0
\end{tabular}\endgroup%
\kern3pt%
\begingroup \smaller\smaller\smaller\begin{tabular}{@{}c@{}}%
0\\-252\\-16
\end{tabular}\endgroup%
\kern3pt%
\begingroup \smaller\smaller\smaller\begin{tabular}{@{}c@{}}%
0\\-16\\-1
\end{tabular}\endgroup%
{$\left.\llap{\phantom{%
\begingroup \smaller\smaller\smaller\begin{tabular}{@{}c@{}}%
0\\0\\0
\end{tabular}\endgroup%
}}\!\right]$}%
\hfil\penalty500%
{$\left[\!\llap{\phantom{%
\begingroup \smaller\smaller\smaller\begin{tabular}{@{}c@{}}%
0\\0\\0
\end{tabular}\endgroup%
}}\right.$}%
\begingroup \smaller\smaller\smaller\begin{tabular}{@{}c@{}}%
-5\\-8\\96
\end{tabular}\endgroup%
\kern3pt%
\begingroup \smaller\smaller\smaller\begin{tabular}{@{}c@{}}%
15\\29\\-360
\end{tabular}\endgroup%
\kern3pt%
\begingroup \smaller\smaller\smaller\begin{tabular}{@{}c@{}}%
1\\2\\-25
\end{tabular}\endgroup%
{$\left.\llap{\phantom{%
\begingroup \smaller\smaller\smaller\begin{tabular}{@{}c@{}}%
0\\0\\0
\end{tabular}\endgroup%
}}\!\right]$}%
\EasyButWeakLineBreak%
{$\left[\!\llap{\phantom{%
\begingroup \smaller\smaller\smaller\begin{tabular}{@{}c@{}}%
0\\0\\0
\end{tabular}\endgroup%
}}\right.$}%
\begingroup \smaller\smaller\smaller\begin{tabular}{@{}c@{}}%
-3\\0\\-16
\end{tabular}\endgroup%
\HardButStrongLineBreak\kern3pt%
\begingroup \smaller\smaller\smaller\begin{tabular}{@{}c@{}}%
-5\\-8\\96
\end{tabular}\endgroup%
\HardButStrongLineBreak\kern3pt%
\begingroup \smaller\smaller\smaller\begin{tabular}{@{}c@{}}%
-1\\-3\\40
\end{tabular}\endgroup%
{$\left.\llap{\phantom{%
\begingroup \smaller\smaller\smaller\begin{tabular}{@{}c@{}}%
0\\0\\0
\end{tabular}\endgroup%
}}\!\right]$}%

\medskip%
%
\leavevmode\llap{}%
$W_{162}$%
\qquad\llap{16} lattices, $\chi=36$%
\hfill%
$\slashtwo\infty|\infty\slashtwo\infty|\infty\rtimes D_{4}$%
\nopagebreak\smallskip\hrule\nopagebreak\medskip%
%
%
\leavevmode%
${L_{162.1}}$%
{} : {$1\above{1pt}{1pt}{2}{0}4\above{1pt}{1pt}{1}{1}{\cdot}1\above{1pt}{1pt}{-2}{}9\above{1pt}{1pt}{1}{}$}\EasyButWeakLineBreak%
{${9}\above{1pt}{1pt}{}{2}{4}\above{1pt}{1pt}{6,1}{\infty}{4}\above{1pt}{1pt}{3,1}{\infty z}{1}\above{1pt}{1pt}{}{2}{36}\above{1pt}{1pt}{2,1}{\infty}{36}\above{1pt}{1pt}{1,0}{\infty z}$}%
\nopagebreak\par%
\nopagebreak\par\leavevmode%
{$\left[\!\llap{\phantom{%
\begingroup \smaller\smaller\smaller\begin{tabular}{@{}c@{}}%
0\\0\\0
\end{tabular}\endgroup%
}}\right.$}%
\begingroup \smaller\smaller\smaller\begin{tabular}{@{}c@{}}%
-298476\\-100440\\-33696
\end{tabular}\endgroup%
\kern3pt%
\begingroup \smaller\smaller\smaller\begin{tabular}{@{}c@{}}%
-100440\\-33799\\-11339
\end{tabular}\endgroup%
\kern3pt%
\begingroup \smaller\smaller\smaller\begin{tabular}{@{}c@{}}%
-33696\\-11339\\-3804
\end{tabular}\endgroup%
{$\left.\llap{\phantom{%
\begingroup \smaller\smaller\smaller\begin{tabular}{@{}c@{}}%
0\\0\\0
\end{tabular}\endgroup%
}}\!\right]$}%
\EasyButWeakLineBreak%
{$\left[\!\llap{\phantom{%
\begingroup \smaller\smaller\smaller\begin{tabular}{@{}c@{}}%
0\\0\\0
\end{tabular}\endgroup%
}}\right.$}%
\begingroup \smaller\smaller\smaller\begin{tabular}{@{}c@{}}%
-1\\9\\-18
\end{tabular}\endgroup%
\HardButStrongLineBreak\kern3pt%
\begingroup \smaller\smaller\smaller\begin{tabular}{@{}c@{}}%
-13\\44\\-16
\end{tabular}\endgroup%
\HardButStrongLineBreak\kern3pt%
\begingroup \smaller\smaller\smaller\begin{tabular}{@{}c@{}}%
15\\-56\\34
\end{tabular}\endgroup%
\HardButStrongLineBreak\kern3pt%
\begingroup \smaller\smaller\smaller\begin{tabular}{@{}c@{}}%
55\\-191\\82
\end{tabular}\endgroup%
\HardButStrongLineBreak\kern3pt%
\begingroup \smaller\smaller\smaller\begin{tabular}{@{}c@{}}%
367\\-1260\\504
\end{tabular}\endgroup%
\HardButStrongLineBreak\kern3pt%
\begingroup \smaller\smaller\smaller\begin{tabular}{@{}c@{}}%
119\\-396\\126
\end{tabular}\endgroup%
{$\left.\llap{\phantom{%
\begingroup \smaller\smaller\smaller\begin{tabular}{@{}c@{}}%
0\\0\\0
\end{tabular}\endgroup%
}}\!\right]$}%
%
%
\hbox{}\par\smallskip%
%
%
\leavevmode%
${L_{162.2}}$%
{} : {$1\above{1pt}{1pt}{2}{0}4\above{1pt}{1pt}{1}{1}{\cdot}1\above{1pt}{1pt}{2}{}9\above{1pt}{1pt}{-}{}$}\EasyButWeakLineBreak%
{${4}\above{1pt}{1pt}{}{2}{1}\above{1pt}{1pt}{12,5}{\infty}{4}\above{1pt}{1pt}{3,2}{\infty a}$}\relax$\,(\times2)$%
\nopagebreak\par%
\nopagebreak\par\leavevmode%
{$\left[\!\llap{\phantom{%
\begingroup \smaller\smaller\smaller\begin{tabular}{@{}c@{}}%
0\\0\\0
\end{tabular}\endgroup%
}}\right.$}%
\begingroup \smaller\smaller\smaller\begin{tabular}{@{}c@{}}%
-348012\\-147888\\-3492
\end{tabular}\endgroup%
\kern3pt%
\begingroup \smaller\smaller\smaller\begin{tabular}{@{}c@{}}%
-147888\\-62845\\-1484
\end{tabular}\endgroup%
\kern3pt%
\begingroup \smaller\smaller\smaller\begin{tabular}{@{}c@{}}%
-3492\\-1484\\-35
\end{tabular}\endgroup%
{$\left.\llap{\phantom{%
\begingroup \smaller\smaller\smaller\begin{tabular}{@{}c@{}}%
0\\0\\0
\end{tabular}\endgroup%
}}\!\right]$}%
\hfil\penalty500%
{$\left[\!\llap{\phantom{%
\begingroup \smaller\smaller\smaller\begin{tabular}{@{}c@{}}%
0\\0\\0
\end{tabular}\endgroup%
}}\right.$}%
\begingroup \smaller\smaller\smaller\begin{tabular}{@{}c@{}}%
61559\\-139320\\-233280
\end{tabular}\endgroup%
\kern3pt%
\begingroup \smaller\smaller\smaller\begin{tabular}{@{}c@{}}%
26182\\-59255\\-99216
\end{tabular}\endgroup%
\kern3pt%
\begingroup \smaller\smaller\smaller\begin{tabular}{@{}c@{}}%
608\\-1376\\-2305
\end{tabular}\endgroup%
{$\left.\llap{\phantom{%
\begingroup \smaller\smaller\smaller\begin{tabular}{@{}c@{}}%
0\\0\\0
\end{tabular}\endgroup%
}}\!\right]$}%
\EasyButWeakLineBreak%
{$\left[\!\llap{\phantom{%
\begingroup \smaller\smaller\smaller\begin{tabular}{@{}c@{}}%
0\\0\\0
\end{tabular}\endgroup%
}}\right.$}%
\begingroup \smaller\smaller\smaller\begin{tabular}{@{}c@{}}%
-159\\360\\596
\end{tabular}\endgroup%
\HardButStrongLineBreak\kern3pt%
\begingroup \smaller\smaller\smaller\begin{tabular}{@{}c@{}}%
-84\\190\\323
\end{tabular}\endgroup%
\HardButStrongLineBreak\kern3pt%
\begingroup \smaller\smaller\smaller\begin{tabular}{@{}c@{}}%
-47\\106\\194
\end{tabular}\endgroup%
{$\left.\llap{\phantom{%
\begingroup \smaller\smaller\smaller\begin{tabular}{@{}c@{}}%
0\\0\\0
\end{tabular}\endgroup%
}}\!\right]$}%
%
%
\hbox{}\par\smallskip%
%
%
\leavevmode%
${L_{162.3}}$%
{} : {$1\above{1pt}{1pt}{2}{{\rm II}}8\above{1pt}{1pt}{1}{1}{\cdot}1\above{1pt}{1pt}{-2}{}9\above{1pt}{1pt}{-}{}$}\EasyButWeakLineBreak%
{${8}\above{1pt}{1pt}{r}{2}{18}\above{1pt}{1pt}{4,1}{\infty a}{72}\above{1pt}{1pt}{1,0}{\infty a}{72}\above{1pt}{1pt}{r}{2}{2}\above{1pt}{1pt}{12,1}{\infty a}{8}\above{1pt}{1pt}{3,2}{\infty a}$}%
\nopagebreak\par%
\nopagebreak\par\leavevmode%
{$\left[\!\llap{\phantom{%
\begingroup \smaller\smaller\smaller\begin{tabular}{@{}c@{}}%
0\\0\\0
\end{tabular}\endgroup%
}}\right.$}%
\begingroup \smaller\smaller\smaller\begin{tabular}{@{}c@{}}%
-741240\\333072\\-99000
\end{tabular}\endgroup%
\kern3pt%
\begingroup \smaller\smaller\smaller\begin{tabular}{@{}c@{}}%
333072\\-149664\\44485
\end{tabular}\endgroup%
\kern3pt%
\begingroup \smaller\smaller\smaller\begin{tabular}{@{}c@{}}%
-99000\\44485\\-13222
\end{tabular}\endgroup%
{$\left.\llap{\phantom{%
\begingroup \smaller\smaller\smaller\begin{tabular}{@{}c@{}}%
0\\0\\0
\end{tabular}\endgroup%
}}\!\right]$}%
\EasyButWeakLineBreak%
{$\left[\!\llap{\phantom{%
\begingroup \smaller\smaller\smaller\begin{tabular}{@{}c@{}}%
0\\0\\0
\end{tabular}\endgroup%
}}\right.$}%
\begingroup \smaller\smaller\smaller\begin{tabular}{@{}c@{}}%
-41\\-96\\-16
\end{tabular}\endgroup%
\HardButStrongLineBreak\kern3pt%
\begingroup \smaller\smaller\smaller\begin{tabular}{@{}c@{}}%
-121\\-288\\-63
\end{tabular}\endgroup%
\HardButStrongLineBreak\kern3pt%
\begingroup \smaller\smaller\smaller\begin{tabular}{@{}c@{}}%
233\\540\\72
\end{tabular}\endgroup%
\HardButStrongLineBreak\kern3pt%
\begingroup \smaller\smaller\smaller\begin{tabular}{@{}c@{}}%
1531\\3600\\648
\end{tabular}\endgroup%
\HardButStrongLineBreak\kern3pt%
\begingroup \smaller\smaller\smaller\begin{tabular}{@{}c@{}}%
275\\648\\121
\end{tabular}\endgroup%
\HardButStrongLineBreak\kern3pt%
\begingroup \smaller\smaller\smaller\begin{tabular}{@{}c@{}}%
157\\372\\76
\end{tabular}\endgroup%
{$\left.\llap{\phantom{%
\begingroup \smaller\smaller\smaller\begin{tabular}{@{}c@{}}%
0\\0\\0
\end{tabular}\endgroup%
}}\!\right]$}%
%
%
\hbox{}\par\smallskip%
%
%
\leavevmode%
${L_{162.4}}$%
{} : {$1\above{1pt}{1pt}{2}{{\rm II}}8\above{1pt}{1pt}{1}{1}{\cdot}1\above{1pt}{1pt}{2}{}9\above{1pt}{1pt}{1}{}$}\EasyButWeakLineBreak%
{${2}\above{1pt}{1pt}{l}{2}{8}\above{1pt}{1pt}{3,2}{\infty}{8}\above{1pt}{1pt}{3,1}{\infty z}$}\relax$\,(\times2)$%
\nopagebreak\par%
\nopagebreak\par\leavevmode%
{$\left[\!\llap{\phantom{%
\begingroup \smaller\smaller\smaller\begin{tabular}{@{}c@{}}%
0\\0\\0
\end{tabular}\endgroup%
}}\right.$}%
\begingroup \smaller\smaller\smaller\begin{tabular}{@{}c@{}}%
-1289592\\-282096\\-9504
\end{tabular}\endgroup%
\kern3pt%
\begingroup \smaller\smaller\smaller\begin{tabular}{@{}c@{}}%
-282096\\-61708\\-2079
\end{tabular}\endgroup%
\kern3pt%
\begingroup \smaller\smaller\smaller\begin{tabular}{@{}c@{}}%
-9504\\-2079\\-70
\end{tabular}\endgroup%
{$\left.\llap{\phantom{%
\begingroup \smaller\smaller\smaller\begin{tabular}{@{}c@{}}%
0\\0\\0
\end{tabular}\endgroup%
}}\!\right]$}%
\hfil\penalty500%
{$\left[\!\llap{\phantom{%
\begingroup \smaller\smaller\smaller\begin{tabular}{@{}c@{}}%
0\\0\\0
\end{tabular}\endgroup%
}}\right.$}%
\begingroup \smaller\smaller\smaller\begin{tabular}{@{}c@{}}%
83447\\-377712\\-109800
\end{tabular}\endgroup%
\kern3pt%
\begingroup \smaller\smaller\smaller\begin{tabular}{@{}c@{}}%
18259\\-82647\\-24025
\end{tabular}\endgroup%
\kern3pt%
\begingroup \smaller\smaller\smaller\begin{tabular}{@{}c@{}}%
608\\-2752\\-801
\end{tabular}\endgroup%
{$\left.\llap{\phantom{%
\begingroup \smaller\smaller\smaller\begin{tabular}{@{}c@{}}%
0\\0\\0
\end{tabular}\endgroup%
}}\!\right]$}%
\EasyButWeakLineBreak%
{$\left[\!\llap{\phantom{%
\begingroup \smaller\smaller\smaller\begin{tabular}{@{}c@{}}%
0\\0\\0
\end{tabular}\endgroup%
}}\right.$}%
\begingroup \smaller\smaller\smaller\begin{tabular}{@{}c@{}}%
8\\-36\\-17
\end{tabular}\endgroup%
\HardButStrongLineBreak\kern3pt%
\begingroup \smaller\smaller\smaller\begin{tabular}{@{}c@{}}%
7\\-32\\0
\end{tabular}\endgroup%
\HardButStrongLineBreak\kern3pt%
\begingroup \smaller\smaller\smaller\begin{tabular}{@{}c@{}}%
-47\\212\\84
\end{tabular}\endgroup%
{$\left.\llap{\phantom{%
\begingroup \smaller\smaller\smaller\begin{tabular}{@{}c@{}}%
0\\0\\0
\end{tabular}\endgroup%
}}\!\right]$}%

\medskip%
%
\leavevmode\llap{}%
$W_{163}$%
\qquad\llap{49} lattices, $\chi=36$%
\hfill%
$\slashtwo2\slashinfty2\slashtwo2\slashinfty2\rtimes D_{4}$%
\nopagebreak\smallskip\hrule\nopagebreak\medskip%
%
%
\leavevmode%
${L_{163.1}}$%
{} : {$1\above{1pt}{1pt}{2}{0}8\above{1pt}{1pt}{-}{5}{\cdot}1\above{1pt}{1pt}{1}{}5\above{1pt}{1pt}{-}{}25\above{1pt}{1pt}{1}{}$}\spacer%
\instructions{5}%
\EasyButWeakLineBreak%
{${25}\above{1pt}{1pt}{r}{2}{4}\above{1pt}{1pt}{*}{2}{40}\above{1pt}{1pt}{5,4}{\infty b}{40}\above{1pt}{1pt}{}{2}{1}\above{1pt}{1pt}{r}{2}{100}\above{1pt}{1pt}{*}{2}{40}\above{1pt}{1pt}{5,1}{\infty a}{40}\above{1pt}{1pt}{}{2}$}%
\nopagebreak\par%
\nopagebreak\par\leavevmode%
{$\left[\!\llap{\phantom{%
\begingroup \smaller\smaller\smaller\begin{tabular}{@{}c@{}}%
0\\0\\0
\end{tabular}\endgroup%
}}\right.$}%
\begingroup \smaller\smaller\smaller\begin{tabular}{@{}c@{}}%
-733400\\1000\\8600
\end{tabular}\endgroup%
\kern3pt%
\begingroup \smaller\smaller\smaller\begin{tabular}{@{}c@{}}%
1000\\35\\-25
\end{tabular}\endgroup%
\kern3pt%
\begingroup \smaller\smaller\smaller\begin{tabular}{@{}c@{}}%
8600\\-25\\-96
\end{tabular}\endgroup%
{$\left.\llap{\phantom{%
\begingroup \smaller\smaller\smaller\begin{tabular}{@{}c@{}}%
0\\0\\0
\end{tabular}\endgroup%
}}\!\right]$}%
\EasyButWeakLineBreak%
{$\left[\!\llap{\phantom{%
\begingroup \smaller\smaller\smaller\begin{tabular}{@{}c@{}}%
0\\0\\0
\end{tabular}\endgroup%
}}\right.$}%
\begingroup \smaller\smaller\smaller\begin{tabular}{@{}c@{}}%
37\\1105\\3025
\end{tabular}\endgroup%
\HardButStrongLineBreak\kern3pt%
\begingroup \smaller\smaller\smaller\begin{tabular}{@{}c@{}}%
25\\746\\2044
\end{tabular}\endgroup%
\HardButStrongLineBreak\kern3pt%
\begingroup \smaller\smaller\smaller\begin{tabular}{@{}c@{}}%
125\\3728\\10220
\end{tabular}\endgroup%
\HardButStrongLineBreak\kern3pt%
\begingroup \smaller\smaller\smaller\begin{tabular}{@{}c@{}}%
91\\2712\\7440
\end{tabular}\endgroup%
\HardButStrongLineBreak\kern3pt%
\begingroup \smaller\smaller\smaller\begin{tabular}{@{}c@{}}%
4\\119\\327
\end{tabular}\endgroup%
\HardButStrongLineBreak\kern3pt%
\begingroup \smaller\smaller\smaller\begin{tabular}{@{}c@{}}%
-11\\-330\\-900
\end{tabular}\endgroup%
\HardButStrongLineBreak\kern3pt%
\begingroup \smaller\smaller\smaller\begin{tabular}{@{}c@{}}%
-11\\-328\\-900
\end{tabular}\endgroup%
\HardButStrongLineBreak\kern3pt%
\begingroup \smaller\smaller\smaller\begin{tabular}{@{}c@{}}%
23\\688\\1880
\end{tabular}\endgroup%
{$\left.\llap{\phantom{%
\begingroup \smaller\smaller\smaller\begin{tabular}{@{}c@{}}%
0\\0\\0
\end{tabular}\endgroup%
}}\!\right]$}%
%
%
\hbox{}\par\smallskip%
%
%
\leavevmode%
${L_{163.2}}$%
{} : {$[1\above{1pt}{1pt}{1}{}2\above{1pt}{1pt}{1}{}]\above{1pt}{1pt}{}{6}16\above{1pt}{1pt}{-}{3}{\cdot}1\above{1pt}{1pt}{1}{}5\above{1pt}{1pt}{-}{}25\above{1pt}{1pt}{1}{}$}\spacer%
\instructions{5m,5,2}%
\EasyButWeakLineBreak%
{${400}\above{1pt}{1pt}{*}{2}{4}\above{1pt}{1pt}{l}{2}{10}\above{1pt}{1pt}{40,19}{\infty}{40}\above{1pt}{1pt}{s}{2}{16}\above{1pt}{1pt}{*}{2}{100}\above{1pt}{1pt}{l}{2}{10}\above{1pt}{1pt}{40,11}{\infty}{40}\above{1pt}{1pt}{s}{2}$}%
\nopagebreak\par%
\nopagebreak\par\leavevmode%
{$\left[\!\llap{\phantom{%
\begingroup \smaller\smaller\smaller\begin{tabular}{@{}c@{}}%
0\\0\\0
\end{tabular}\endgroup%
}}\right.$}%
\begingroup \smaller\smaller\smaller\begin{tabular}{@{}c@{}}%
-587600\\21600\\800
\end{tabular}\endgroup%
\kern3pt%
\begingroup \smaller\smaller\smaller\begin{tabular}{@{}c@{}}%
21600\\-790\\-30
\end{tabular}\endgroup%
\kern3pt%
\begingroup \smaller\smaller\smaller\begin{tabular}{@{}c@{}}%
800\\-30\\-1
\end{tabular}\endgroup%
{$\left.\llap{\phantom{%
\begingroup \smaller\smaller\smaller\begin{tabular}{@{}c@{}}%
0\\0\\0
\end{tabular}\endgroup%
}}\!\right]$}%
\EasyButWeakLineBreak%
{$\left[\!\llap{\phantom{%
\begingroup \smaller\smaller\smaller\begin{tabular}{@{}c@{}}%
0\\0\\0
\end{tabular}\endgroup%
}}\right.$}%
\begingroup \smaller\smaller\smaller\begin{tabular}{@{}c@{}}%
-1\\-20\\-200
\end{tabular}\endgroup%
\HardButStrongLineBreak\kern3pt%
\begingroup \smaller\smaller\smaller\begin{tabular}{@{}c@{}}%
-1\\-22\\-146
\end{tabular}\endgroup%
\HardButStrongLineBreak\kern3pt%
\begingroup \smaller\smaller\smaller\begin{tabular}{@{}c@{}}%
-4\\-89\\-560
\end{tabular}\endgroup%
\HardButStrongLineBreak\kern3pt%
\begingroup \smaller\smaller\smaller\begin{tabular}{@{}c@{}}%
-9\\-202\\-1220
\end{tabular}\endgroup%
\HardButStrongLineBreak\kern3pt%
\begingroup \smaller\smaller\smaller\begin{tabular}{@{}c@{}}%
-3\\-68\\-392
\end{tabular}\endgroup%
\HardButStrongLineBreak\kern3pt%
\begingroup \smaller\smaller\smaller\begin{tabular}{@{}c@{}}%
-3\\-70\\-350
\end{tabular}\endgroup%
\HardButStrongLineBreak\kern3pt%
\begingroup \smaller\smaller\smaller\begin{tabular}{@{}c@{}}%
0\\-1\\20
\end{tabular}\endgroup%
\HardButStrongLineBreak\kern3pt%
\begingroup \smaller\smaller\smaller\begin{tabular}{@{}c@{}}%
1\\22\\140
\end{tabular}\endgroup%
{$\left.\llap{\phantom{%
\begingroup \smaller\smaller\smaller\begin{tabular}{@{}c@{}}%
0\\0\\0
\end{tabular}\endgroup%
}}\!\right]$}%
%
%
\hbox{}\par\smallskip%
%
%
\leavevmode%
${L_{163.3}}$%
{} : {$[1\above{1pt}{1pt}{-}{}2\above{1pt}{1pt}{1}{}]\above{1pt}{1pt}{}{2}16\above{1pt}{1pt}{1}{7}{\cdot}1\above{1pt}{1pt}{1}{}5\above{1pt}{1pt}{-}{}25\above{1pt}{1pt}{1}{}$}\spacer%
\instructions{52,5,m}%
\EasyButWeakLineBreak%
{${400}\above{1pt}{1pt}{l}{2}{1}\above{1pt}{1pt}{}{2}{10}\above{1pt}{1pt}{40,39}{\infty}{40}\above{1pt}{1pt}{*}{2}{16}\above{1pt}{1pt}{l}{2}{25}\above{1pt}{1pt}{}{2}{10}\above{1pt}{1pt}{40,31}{\infty}{40}\above{1pt}{1pt}{*}{2}$}%
\nopagebreak\par%
\nopagebreak\par\leavevmode%
{$\left[\!\llap{\phantom{%
\begingroup \smaller\smaller\smaller\begin{tabular}{@{}c@{}}%
0\\0\\0
\end{tabular}\endgroup%
}}\right.$}%
\begingroup \smaller\smaller\smaller\begin{tabular}{@{}c@{}}%
-931600\\27600\\10400
\end{tabular}\endgroup%
\kern3pt%
\begingroup \smaller\smaller\smaller\begin{tabular}{@{}c@{}}%
27600\\-790\\-320
\end{tabular}\endgroup%
\kern3pt%
\begingroup \smaller\smaller\smaller\begin{tabular}{@{}c@{}}%
10400\\-320\\-111
\end{tabular}\endgroup%
{$\left.\llap{\phantom{%
\begingroup \smaller\smaller\smaller\begin{tabular}{@{}c@{}}%
0\\0\\0
\end{tabular}\endgroup%
}}\!\right]$}%
\EasyButWeakLineBreak%
{$\left[\!\llap{\phantom{%
\begingroup \smaller\smaller\smaller\begin{tabular}{@{}c@{}}%
0\\0\\0
\end{tabular}\endgroup%
}}\right.$}%
\begingroup \smaller\smaller\smaller\begin{tabular}{@{}c@{}}%
-19\\-340\\-800
\end{tabular}\endgroup%
\HardButStrongLineBreak\kern3pt%
\begingroup \smaller\smaller\smaller\begin{tabular}{@{}c@{}}%
4\\72\\167
\end{tabular}\endgroup%
\HardButStrongLineBreak\kern3pt%
\begingroup \smaller\smaller\smaller\begin{tabular}{@{}c@{}}%
44\\791\\1840
\end{tabular}\endgroup%
\HardButStrongLineBreak\kern3pt%
\begingroup \smaller\smaller\smaller\begin{tabular}{@{}c@{}}%
119\\2138\\4980
\end{tabular}\endgroup%
\HardButStrongLineBreak\kern3pt%
\begingroup \smaller\smaller\smaller\begin{tabular}{@{}c@{}}%
47\\844\\1968
\end{tabular}\endgroup%
\HardButStrongLineBreak\kern3pt%
\begingroup \smaller\smaller\smaller\begin{tabular}{@{}c@{}}%
34\\610\\1425
\end{tabular}\endgroup%
\HardButStrongLineBreak\kern3pt%
\begingroup \smaller\smaller\smaller\begin{tabular}{@{}c@{}}%
10\\179\\420
\end{tabular}\endgroup%
\HardButStrongLineBreak\kern3pt%
\begingroup \smaller\smaller\smaller\begin{tabular}{@{}c@{}}%
-11\\-198\\-460
\end{tabular}\endgroup%
{$\left.\llap{\phantom{%
\begingroup \smaller\smaller\smaller\begin{tabular}{@{}c@{}}%
0\\0\\0
\end{tabular}\endgroup%
}}\!\right]$}%
%
%
\hbox{}\par\smallskip%
%
%
\leavevmode%
${L_{163.4}}$%
{} : {$[1\above{1pt}{1pt}{-}{}2\above{1pt}{1pt}{-}{}]\above{1pt}{1pt}{}{0}16\above{1pt}{1pt}{-}{5}{\cdot}1\above{1pt}{1pt}{1}{}5\above{1pt}{1pt}{-}{}25\above{1pt}{1pt}{1}{}$}\spacer%
\instructions{5m,5,m}%
\EasyButWeakLineBreak%
{${400}\above{1pt}{1pt}{s}{2}{4}\above{1pt}{1pt}{*}{2}{40}\above{1pt}{1pt}{40,19}{\infty z}{10}\above{1pt}{1pt}{r}{2}{16}\above{1pt}{1pt}{s}{2}{100}\above{1pt}{1pt}{*}{2}{40}\above{1pt}{1pt}{40,11}{\infty z}{10}\above{1pt}{1pt}{r}{2}$}%
\nopagebreak\par%
\nopagebreak\par\leavevmode%
{$\left[\!\llap{\phantom{%
\begingroup \smaller\smaller\smaller\begin{tabular}{@{}c@{}}%
0\\0\\0
\end{tabular}\endgroup%
}}\right.$}%
\begingroup \smaller\smaller\smaller\begin{tabular}{@{}c@{}}%
-676400\\8800\\201600
\end{tabular}\endgroup%
\kern3pt%
\begingroup \smaller\smaller\smaller\begin{tabular}{@{}c@{}}%
8800\\-90\\-2820
\end{tabular}\endgroup%
\kern3pt%
\begingroup \smaller\smaller\smaller\begin{tabular}{@{}c@{}}%
201600\\-2820\\-58499
\end{tabular}\endgroup%
{$\left.\llap{\phantom{%
\begingroup \smaller\smaller\smaller\begin{tabular}{@{}c@{}}%
0\\0\\0
\end{tabular}\endgroup%
}}\!\right]$}%
\EasyButWeakLineBreak%
{$\left[\!\llap{\phantom{%
\begingroup \smaller\smaller\smaller\begin{tabular}{@{}c@{}}%
0\\0\\0
\end{tabular}\endgroup%
}}\right.$}%
\begingroup \smaller\smaller\smaller\begin{tabular}{@{}c@{}}%
1853\\37040\\4600
\end{tabular}\endgroup%
\HardButStrongLineBreak\kern3pt%
\begingroup \smaller\smaller\smaller\begin{tabular}{@{}c@{}}%
315\\6296\\782
\end{tabular}\endgroup%
\HardButStrongLineBreak\kern3pt%
\begingroup \smaller\smaller\smaller\begin{tabular}{@{}c@{}}%
1579\\31558\\3920
\end{tabular}\endgroup%
\HardButStrongLineBreak\kern3pt%
\begingroup \smaller\smaller\smaller\begin{tabular}{@{}c@{}}%
576\\11511\\1430
\end{tabular}\endgroup%
\HardButStrongLineBreak\kern3pt%
\begingroup \smaller\smaller\smaller\begin{tabular}{@{}c@{}}%
203\\4056\\504
\end{tabular}\endgroup%
\HardButStrongLineBreak\kern3pt%
\begingroup \smaller\smaller\smaller\begin{tabular}{@{}c@{}}%
-141\\-2820\\-350
\end{tabular}\endgroup%
\HardButStrongLineBreak\kern3pt%
\begingroup \smaller\smaller\smaller\begin{tabular}{@{}c@{}}%
-145\\-2898\\-360
\end{tabular}\endgroup%
\HardButStrongLineBreak\kern3pt%
\begingroup \smaller\smaller\smaller\begin{tabular}{@{}c@{}}%
141\\2819\\350
\end{tabular}\endgroup%
{$\left.\llap{\phantom{%
\begingroup \smaller\smaller\smaller\begin{tabular}{@{}c@{}}%
0\\0\\0
\end{tabular}\endgroup%
}}\!\right]$}%
%
%
\hbox{}\par\smallskip%
%
%
\leavevmode%
${L_{163.5}}$%
{} : {$1\above{1pt}{1pt}{1}{1}8\above{1pt}{1pt}{-}{3}64\above{1pt}{1pt}{1}{1}{\cdot}1\above{1pt}{1pt}{1}{}5\above{1pt}{1pt}{-}{}25\above{1pt}{1pt}{1}{}$}\spacer%
\instructions{25,5,2}%
\EasyButWeakLineBreak%
{${1600}\above{1pt}{1pt}{}{2}{1}\above{1pt}{1pt}{r}{2}{160}\above{1pt}{1pt}{80,79}{\infty z}{40}\above{1pt}{1pt}{b}{2}{64}\above{1pt}{1pt}{s}{2}{100}\above{1pt}{1pt}{*}{2}{160}\above{1pt}{1pt}{80,71}{\infty z}{40}\above{1pt}{1pt}{l}{2}$}%
\nopagebreak\par%
shares genus with 2-dual${}\iso{}$5-dual; isometric to own %
2.5-dual\nopagebreak\par%
\nopagebreak\par\leavevmode%
{$\left[\!\llap{\phantom{%
\begingroup \smaller\smaller\smaller\begin{tabular}{@{}c@{}}%
0\\0\\0
\end{tabular}\endgroup%
}}\right.$}%
\begingroup \smaller\smaller\smaller\begin{tabular}{@{}c@{}}%
-958400\\19200\\233600
\end{tabular}\endgroup%
\kern3pt%
\begingroup \smaller\smaller\smaller\begin{tabular}{@{}c@{}}%
19200\\-360\\-4880
\end{tabular}\endgroup%
\kern3pt%
\begingroup \smaller\smaller\smaller\begin{tabular}{@{}c@{}}%
233600\\-4880\\-55311
\end{tabular}\endgroup%
{$\left.\llap{\phantom{%
\begingroup \smaller\smaller\smaller\begin{tabular}{@{}c@{}}%
0\\0\\0
\end{tabular}\endgroup%
}}\!\right]$}%
\EasyButWeakLineBreak%
{$\left[\!\llap{\phantom{%
\begingroup \smaller\smaller\smaller\begin{tabular}{@{}c@{}}%
0\\0\\0
\end{tabular}\endgroup%
}}\right.$}%
\begingroup \smaller\smaller\smaller\begin{tabular}{@{}c@{}}%
1301\\26000\\3200
\end{tabular}\endgroup%
\HardButStrongLineBreak\kern3pt%
\begingroup \smaller\smaller\smaller\begin{tabular}{@{}c@{}}%
50\\999\\123
\end{tabular}\endgroup%
\HardButStrongLineBreak\kern3pt%
\begingroup \smaller\smaller\smaller\begin{tabular}{@{}c@{}}%
943\\18838\\2320
\end{tabular}\endgroup%
\HardButStrongLineBreak\kern3pt%
\begingroup \smaller\smaller\smaller\begin{tabular}{@{}c@{}}%
317\\6331\\780
\end{tabular}\endgroup%
\HardButStrongLineBreak\kern3pt%
\begingroup \smaller\smaller\smaller\begin{tabular}{@{}c@{}}%
91\\1816\\224
\end{tabular}\endgroup%
\HardButStrongLineBreak\kern3pt%
\begingroup \smaller\smaller\smaller\begin{tabular}{@{}c@{}}%
-61\\-1220\\-150
\end{tabular}\endgroup%
\HardButStrongLineBreak\kern3pt%
\begingroup \smaller\smaller\smaller\begin{tabular}{@{}c@{}}%
-65\\-1298\\-160
\end{tabular}\endgroup%
\HardButStrongLineBreak\kern3pt%
\begingroup \smaller\smaller\smaller\begin{tabular}{@{}c@{}}%
122\\2439\\300
\end{tabular}\endgroup%
{$\left.\llap{\phantom{%
\begingroup \smaller\smaller\smaller\begin{tabular}{@{}c@{}}%
0\\0\\0
\end{tabular}\endgroup%
}}\!\right]$}%

\medskip%
%
\leavevmode\llap{}%
$W_{164}$%
\qquad\llap{44} lattices, $\chi=54$%
\hfill%
$42\infty2242\infty22\rtimes C_{2}$%
\nopagebreak\smallskip\hrule\nopagebreak\medskip%
%
%
\leavevmode%
${L_{164.1}}$%
{} : {$1\above{1pt}{1pt}{2}{{\rm II}}4\above{1pt}{1pt}{-}{5}{\cdot}1\above{1pt}{1pt}{2}{}9\above{1pt}{1pt}{1}{}{\cdot}1\above{1pt}{1pt}{2}{}5\above{1pt}{1pt}{1}{}$}\spacer%
\instructions{2}%
\EasyButWeakLineBreak%
{${4}\above{1pt}{1pt}{*}{4}{2}\above{1pt}{1pt}{l}{2}{20}\above{1pt}{1pt}{3,2}{\infty}{20}\above{1pt}{1pt}{*}{2}{36}\above{1pt}{1pt}{*}{2}$}\relax$\,(\times2)$%
\nopagebreak\par%
\nopagebreak\par\leavevmode%
{$\left[\!\llap{\phantom{%
\begingroup \smaller\smaller\smaller
\endgroup%
}}\!\right]$}%
%
%
\hbox{}\par\smallskip%
%
%
\leavevmode%
${L_{164.2}}$%
{} : {$1\above{1pt}{1pt}{-2}{2}8\above{1pt}{1pt}{1}{7}{\cdot}1\above{1pt}{1pt}{2}{}9\above{1pt}{1pt}{-}{}{\cdot}1\above{1pt}{1pt}{2}{}5\above{1pt}{1pt}{-}{}$}\spacer%
\instructions{2}%
\EasyButWeakLineBreak%
{${2}\above{1pt}{1pt}{*}{4}{4}\above{1pt}{1pt}{s}{2}{40}\above{1pt}{1pt}{12,5}{\infty z}{10}\above{1pt}{1pt}{s}{2}{18}\above{1pt}{1pt}{b}{2}$}\relax$\,(\times2)$%
\nopagebreak\par%
\nopagebreak\par\leavevmode%
{$\left[\!\llap{\phantom{%
\begingroup \smaller\smaller\smaller
\endgroup%
}}\!\right]$}%
%
%
\hbox{}\par\smallskip%
%
%
\leavevmode%
${L_{164.3}}$%
{} : {$1\above{1pt}{1pt}{2}{2}8\above{1pt}{1pt}{-}{3}{\cdot}1\above{1pt}{1pt}{2}{}9\above{1pt}{1pt}{-}{}{\cdot}1\above{1pt}{1pt}{2}{}5\above{1pt}{1pt}{-}{}$}\spacer%
\instructions{m}%
\EasyButWeakLineBreak%
{${2}\above{1pt}{1pt}{}{4}{1}\above{1pt}{1pt}{r}{2}{40}\above{1pt}{1pt}{12,11}{\infty z}{10}\above{1pt}{1pt}{b}{2}{18}\above{1pt}{1pt}{s}{2}$}\relax$\,(\times2)$%
\nopagebreak\par%
\nopagebreak\par\leavevmode%
{$\left[\!\llap{\phantom{%
\begingroup \smaller\smaller\smaller
\endgroup%
}}\!\right]$}%

\medskip%
%
\leavevmode\llap{}%
$W_{165}$%
\qquad\llap{12} lattices, $\chi=6$%
\hfill%
$22222$%
\nopagebreak\smallskip\hrule\nopagebreak\medskip%
%
%
\leavevmode%
${L_{165.1}}$%
{} : {$1\above{1pt}{1pt}{-2}{{\rm II}}4\above{1pt}{1pt}{1}{1}{\cdot}1\above{1pt}{1pt}{2}{}9\above{1pt}{1pt}{1}{}{\cdot}1\above{1pt}{1pt}{-2}{}5\above{1pt}{1pt}{-}{}$}\spacer%
\instructions{2}%
\EasyButWeakLineBreak%
{${90}\above{1pt}{1pt}{l}{2}{4}\above{1pt}{1pt}{r}{2}{10}\above{1pt}{1pt}{l}{2}{36}\above{1pt}{1pt}{r}{2}{2}\above{1pt}{1pt}{b}{2}$}%
\nopagebreak\par%
\nopagebreak\par\leavevmode%
{$\left[\!\llap{\phantom{%
\begingroup \smaller\smaller\smaller\begin{tabular}{@{}c@{}}%
0\\0\\0
\end{tabular}\endgroup%
}}\right.$}%
\begingroup \smaller\smaller\smaller\begin{tabular}{@{}c@{}}%
-400860\\-197640\\1080
\end{tabular}\endgroup%
\kern3pt%
\begingroup \smaller\smaller\smaller\begin{tabular}{@{}c@{}}%
-197640\\-97442\\531
\end{tabular}\endgroup%
\kern3pt%
\begingroup \smaller\smaller\smaller\begin{tabular}{@{}c@{}}%
1080\\531\\-2
\end{tabular}\endgroup%
{$\left.\llap{\phantom{%
\begingroup \smaller\smaller\smaller\begin{tabular}{@{}c@{}}%
0\\0\\0
\end{tabular}\endgroup%
}}\!\right]$}%
\EasyButWeakLineBreak%
{$\left[\!\llap{\phantom{%
\begingroup \smaller\smaller\smaller\begin{tabular}{@{}c@{}}%
0\\0\\0
\end{tabular}\endgroup%
}}\right.$}%
\begingroup \smaller\smaller\smaller\begin{tabular}{@{}c@{}}%
-44\\90\\135
\end{tabular}\endgroup%
\HardButStrongLineBreak\kern3pt%
\begingroup \smaller\smaller\smaller\begin{tabular}{@{}c@{}}%
43\\-88\\-144
\end{tabular}\endgroup%
\HardButStrongLineBreak\kern3pt%
\begingroup \smaller\smaller\smaller\begin{tabular}{@{}c@{}}%
22\\-45\\-70
\end{tabular}\endgroup%
\HardButStrongLineBreak\kern3pt%
\begingroup \smaller\smaller\smaller\begin{tabular}{@{}c@{}}%
-299\\612\\1008
\end{tabular}\endgroup%
\HardButStrongLineBreak\kern3pt%
\begingroup \smaller\smaller\smaller\begin{tabular}{@{}c@{}}%
-43\\88\\143
\end{tabular}\endgroup%
{$\left.\llap{\phantom{%
\begingroup \smaller\smaller\smaller\begin{tabular}{@{}c@{}}%
0\\0\\0
\end{tabular}\endgroup%
}}\!\right]$}%

\medskip%
%
\leavevmode\llap{}%
$W_{166}$%
\qquad\llap{12} lattices, $\chi=24$%
\hfill%
$22|22|22|22|\rtimes D_{4}$%
\nopagebreak\smallskip\hrule\nopagebreak\medskip%
%
%
\leavevmode%
${L_{166.1}}$%
{} : {$1\above{1pt}{1pt}{-2}{{\rm II}}4\above{1pt}{1pt}{1}{1}{\cdot}1\above{1pt}{1pt}{-2}{}9\above{1pt}{1pt}{-}{}{\cdot}1\above{1pt}{1pt}{-2}{}5\above{1pt}{1pt}{-}{}$}\spacer%
\instructions{2}%
\EasyButWeakLineBreak%
{${4}\above{1pt}{1pt}{r}{2}{18}\above{1pt}{1pt}{b}{2}{10}\above{1pt}{1pt}{b}{2}{2}\above{1pt}{1pt}{l}{2}$}\relax$\,(\times2)$%
\nopagebreak\par%
\nopagebreak\par\leavevmode%
{$\left[\!\llap{\phantom{%
\begingroup \smaller\smaller\smaller\begin{tabular}{@{}c@{}}%
0\\0\\0
\end{tabular}\endgroup%
}}\right.$}%
\begingroup \smaller\smaller\smaller\begin{tabular}{@{}c@{}}%
2340\\-1080\\0
\end{tabular}\endgroup%
\kern3pt%
\begingroup \smaller\smaller\smaller\begin{tabular}{@{}c@{}}%
-1080\\498\\1
\end{tabular}\endgroup%
\kern3pt%
\begingroup \smaller\smaller\smaller\begin{tabular}{@{}c@{}}%
0\\1\\-2
\end{tabular}\endgroup%
{$\left.\llap{\phantom{%
\begingroup \smaller\smaller\smaller\begin{tabular}{@{}c@{}}%
0\\0\\0
\end{tabular}\endgroup%
}}\!\right]$}%
\hfil\penalty500%
{$\left[\!\llap{\phantom{%
\begingroup \smaller\smaller\smaller\begin{tabular}{@{}c@{}}%
0\\0\\0
\end{tabular}\endgroup%
}}\right.$}%
\begingroup \smaller\smaller\smaller\begin{tabular}{@{}c@{}}%
-1\\0\\0
\end{tabular}\endgroup%
\kern3pt%
\begingroup \smaller\smaller\smaller\begin{tabular}{@{}c@{}}%
0\\-1\\-1
\end{tabular}\endgroup%
\kern3pt%
\begingroup \smaller\smaller\smaller\begin{tabular}{@{}c@{}}%
0\\0\\1
\end{tabular}\endgroup%
{$\left.\llap{\phantom{%
\begingroup \smaller\smaller\smaller\begin{tabular}{@{}c@{}}%
0\\0\\0
\end{tabular}\endgroup%
}}\!\right]$}%
\EasyButWeakLineBreak%
{$\left[\!\llap{\phantom{%
\begingroup \smaller\smaller\smaller\begin{tabular}{@{}c@{}}%
0\\0\\0
\end{tabular}\endgroup%
}}\right.$}%
\begingroup \smaller\smaller\smaller\begin{tabular}{@{}c@{}}%
11\\24\\8
\end{tabular}\endgroup%
\HardButStrongLineBreak\kern3pt%
\begingroup \smaller\smaller\smaller\begin{tabular}{@{}c@{}}%
4\\9\\0
\end{tabular}\endgroup%
\HardButStrongLineBreak\kern3pt%
\begingroup \smaller\smaller\smaller\begin{tabular}{@{}c@{}}%
-7\\-15\\-10
\end{tabular}\endgroup%
\HardButStrongLineBreak\kern3pt%
\begingroup \smaller\smaller\smaller\begin{tabular}{@{}c@{}}%
-6\\-13\\-8
\end{tabular}\endgroup%
{$\left.\llap{\phantom{%
\begingroup \smaller\smaller\smaller\begin{tabular}{@{}c@{}}%
0\\0\\0
\end{tabular}\endgroup%
}}\!\right]$}%

\medskip%
%
\leavevmode\llap{}%
$W_{167}$%
\qquad\llap{4} lattices, $\chi=6$%
\hfill%
$4|42|2\rtimes D_{2}$%
\nopagebreak\smallskip\hrule\nopagebreak\medskip%
%
%
\leavevmode%
${L_{167.1}}$%
{} : {$1\above{1pt}{1pt}{2}{2}16\above{1pt}{1pt}{-}{5}{\cdot}1\above{1pt}{1pt}{2}{}3\above{1pt}{1pt}{-}{}$}\EasyButWeakLineBreak%
{${1}\above{1pt}{1pt}{}{4}{2}\above{1pt}{1pt}{*}{4}{4}\above{1pt}{1pt}{s}{2}{16}\above{1pt}{1pt}{l}{2}$}%
\nopagebreak\par%
\nopagebreak\par\leavevmode%
{$\left[\!\llap{\phantom{%
\begingroup \smaller\smaller\smaller\begin{tabular}{@{}c@{}}%
0\\0\\0
\end{tabular}\endgroup%
}}\right.$}%
\begingroup \smaller\smaller\smaller\begin{tabular}{@{}c@{}}%
-2352\\-672\\240
\end{tabular}\endgroup%
\kern3pt%
\begingroup \smaller\smaller\smaller\begin{tabular}{@{}c@{}}%
-672\\-191\\67
\end{tabular}\endgroup%
\kern3pt%
\begingroup \smaller\smaller\smaller\begin{tabular}{@{}c@{}}%
240\\67\\-22
\end{tabular}\endgroup%
{$\left.\llap{\phantom{%
\begingroup \smaller\smaller\smaller\begin{tabular}{@{}c@{}}%
0\\0\\0
\end{tabular}\endgroup%
}}\!\right]$}%
\EasyButWeakLineBreak%
{$\left[\!\llap{\phantom{%
\begingroup \smaller\smaller\smaller\begin{tabular}{@{}c@{}}%
0\\0\\0
\end{tabular}\endgroup%
}}\right.$}%
\begingroup \smaller\smaller\smaller\begin{tabular}{@{}c@{}}%
-1\\5\\4
\end{tabular}\endgroup%
\HardButStrongLineBreak\kern3pt%
\begingroup \smaller\smaller\smaller\begin{tabular}{@{}c@{}}%
-4\\18\\11
\end{tabular}\endgroup%
\HardButStrongLineBreak\kern3pt%
\begingroup \smaller\smaller\smaller\begin{tabular}{@{}c@{}}%
3\\-14\\-10
\end{tabular}\endgroup%
\HardButStrongLineBreak\kern3pt%
\begingroup \smaller\smaller\smaller\begin{tabular}{@{}c@{}}%
9\\-40\\-24
\end{tabular}\endgroup%
{$\left.\llap{\phantom{%
\begingroup \smaller\smaller\smaller\begin{tabular}{@{}c@{}}%
0\\0\\0
\end{tabular}\endgroup%
}}\!\right]$}%
%
%
%
%
%
%
%
%
%
%
%
%
%
%

\medskip%
%
\leavevmode\llap{}%
$W_{168}$%
\qquad\llap{6} lattices, $\chi=8$%
\hfill%
$6262\rtimes C_{2}$%
\nopagebreak\smallskip\hrule\nopagebreak\medskip%
%
%
\leavevmode%
${L_{168.1}}$%
{} : {$1\above{1pt}{1pt}{-2}{{\rm II}}16\above{1pt}{1pt}{1}{7}{\cdot}1\above{1pt}{1pt}{-}{}3\above{1pt}{1pt}{-}{}9\above{1pt}{1pt}{-}{}$}\spacer%
\instructions{3}%
\EasyButWeakLineBreak%
{${6}\above{1pt}{1pt}{}{6}{18}\above{1pt}{1pt}{b}{2}{6}\above{1pt}{1pt}{}{6}{2}\above{1pt}{1pt}{b}{2}$}%
\nopagebreak\par%
\nopagebreak\par\leavevmode%
{$\left[\!\llap{\phantom{%
\begingroup \smaller\smaller\smaller\begin{tabular}{@{}c@{}}%
0\\0\\0
\end{tabular}\endgroup%
}}\right.$}%
\begingroup \smaller\smaller\smaller\begin{tabular}{@{}c@{}}%
-48528\\1296\\720
\end{tabular}\endgroup%
\kern3pt%
\begingroup \smaller\smaller\smaller\begin{tabular}{@{}c@{}}%
1296\\-30\\-21
\end{tabular}\endgroup%
\kern3pt%
\begingroup \smaller\smaller\smaller\begin{tabular}{@{}c@{}}%
720\\-21\\-10
\end{tabular}\endgroup%
{$\left.\llap{\phantom{%
\begingroup \smaller\smaller\smaller\begin{tabular}{@{}c@{}}%
0\\0\\0
\end{tabular}\endgroup%
}}\!\right]$}%
\EasyButWeakLineBreak%
{$\left[\!\llap{\phantom{%
\begingroup \smaller\smaller\smaller\begin{tabular}{@{}c@{}}%
0\\0\\0
\end{tabular}\endgroup%
}}\right.$}%
\begingroup \smaller\smaller\smaller\begin{tabular}{@{}c@{}}%
2\\31\\78
\end{tabular}\endgroup%
\HardButStrongLineBreak\kern3pt%
\begingroup \smaller\smaller\smaller\begin{tabular}{@{}c@{}}%
-2\\-30\\-81
\end{tabular}\endgroup%
\HardButStrongLineBreak\kern3pt%
\begingroup \smaller\smaller\smaller\begin{tabular}{@{}c@{}}%
-1\\-16\\-39
\end{tabular}\endgroup%
\HardButStrongLineBreak\kern3pt%
\begingroup \smaller\smaller\smaller\begin{tabular}{@{}c@{}}%
1\\15\\40
\end{tabular}\endgroup%
{$\left.\llap{\phantom{%
\begingroup \smaller\smaller\smaller\begin{tabular}{@{}c@{}}%
0\\0\\0
\end{tabular}\endgroup%
}}\!\right]$}%

\medskip%
%
\leavevmode\llap{}%
$W_{169}$%
\qquad\llap{36} lattices, $\chi=24$%
\hfill%
$2|2\slashinfty2|2\slashinfty\rtimes D_{4}$%
\nopagebreak\smallskip\hrule\nopagebreak\medskip%
%
%
\leavevmode%
${L_{169.1}}$%
{} : {$1\above{1pt}{1pt}{2}{6}16\above{1pt}{1pt}{-}{5}{\cdot}1\above{1pt}{1pt}{-}{}3\above{1pt}{1pt}{1}{}9\above{1pt}{1pt}{1}{}$}\spacer%
\instructions{3}%
\EasyButWeakLineBreak%
{${3}\above{1pt}{1pt}{r}{2}{144}\above{1pt}{1pt}{s}{2}{12}\above{1pt}{1pt}{24,1}{\infty z}$}\relax$\,(\times2)$%
\nopagebreak\par%
\nopagebreak\par\leavevmode%
{$\left[\!\llap{\phantom{%
\begingroup \smaller\smaller\smaller\begin{tabular}{@{}c@{}}%
0\\0\\0
\end{tabular}\endgroup%
}}\right.$}%
\begingroup \smaller\smaller\smaller\begin{tabular}{@{}c@{}}%
-305712\\-9360\\-9936
\end{tabular}\endgroup%
\kern3pt%
\begingroup \smaller\smaller\smaller\begin{tabular}{@{}c@{}}%
-9360\\-285\\-303
\end{tabular}\endgroup%
\kern3pt%
\begingroup \smaller\smaller\smaller\begin{tabular}{@{}c@{}}%
-9936\\-303\\-322
\end{tabular}\endgroup%
{$\left.\llap{\phantom{%
\begingroup \smaller\smaller\smaller\begin{tabular}{@{}c@{}}%
0\\0\\0
\end{tabular}\endgroup%
}}\!\right]$}%
\hfil\penalty500%
{$\left[\!\llap{\phantom{%
\begingroup \smaller\smaller\smaller\begin{tabular}{@{}c@{}}%
0\\0\\0
\end{tabular}\endgroup%
}}\right.$}%
\begingroup \smaller\smaller\smaller\begin{tabular}{@{}c@{}}%
-449\\-32256\\44352
\end{tabular}\endgroup%
\kern3pt%
\begingroup \smaller\smaller\smaller\begin{tabular}{@{}c@{}}%
-13\\-937\\1287
\end{tabular}\endgroup%
\kern3pt%
\begingroup \smaller\smaller\smaller\begin{tabular}{@{}c@{}}%
-14\\-1008\\1385
\end{tabular}\endgroup%
{$\left.\llap{\phantom{%
\begingroup \smaller\smaller\smaller\begin{tabular}{@{}c@{}}%
0\\0\\0
\end{tabular}\endgroup%
}}\!\right]$}%
\EasyButWeakLineBreak%
{$\left[\!\llap{\phantom{%
\begingroup \smaller\smaller\smaller\begin{tabular}{@{}c@{}}%
0\\0\\0
\end{tabular}\endgroup%
}}\right.$}%
\begingroup \smaller\smaller\smaller\begin{tabular}{@{}c@{}}%
1\\85\\-111
\end{tabular}\endgroup%
\HardButStrongLineBreak\kern3pt%
\begingroup \smaller\smaller\smaller\begin{tabular}{@{}c@{}}%
1\\120\\-144
\end{tabular}\endgroup%
\HardButStrongLineBreak\kern3pt%
\begingroup \smaller\smaller\smaller\begin{tabular}{@{}c@{}}%
-1\\-82\\108
\end{tabular}\endgroup%
{$\left.\llap{\phantom{%
\begingroup \smaller\smaller\smaller\begin{tabular}{@{}c@{}}%
0\\0\\0
\end{tabular}\endgroup%
}}\!\right]$}%
%
%
\hbox{}\par\smallskip%
%
%
\leavevmode%
${L_{169.2}}$%
{} : {$1\above{1pt}{1pt}{-2}{4}16\above{1pt}{1pt}{1}{7}{\cdot}1\above{1pt}{1pt}{-}{}3\above{1pt}{1pt}{1}{}9\above{1pt}{1pt}{1}{}$}\spacer%
\instructions{3}%
\EasyButWeakLineBreak%
{${12}\above{1pt}{1pt}{*}{2}{36}\above{1pt}{1pt}{s}{2}{48}\above{1pt}{1pt}{12,1}{\infty z}{48}\above{1pt}{1pt}{l}{2}{9}\above{1pt}{1pt}{}{2}{3}\above{1pt}{1pt}{24,7}{\infty}$}%
\nopagebreak\par%
\nopagebreak\par\leavevmode%
{$\left[\!\llap{\phantom{%
\begingroup \smaller\smaller\smaller\begin{tabular}{@{}c@{}}%
0\\0\\0
\end{tabular}\endgroup%
}}\right.$}%
\begingroup \smaller\smaller\smaller\begin{tabular}{@{}c@{}}%
-30096\\4032\\-6336
\end{tabular}\endgroup%
\kern3pt%
\begingroup \smaller\smaller\smaller\begin{tabular}{@{}c@{}}%
4032\\-492\\891
\end{tabular}\endgroup%
\kern3pt%
\begingroup \smaller\smaller\smaller\begin{tabular}{@{}c@{}}%
-6336\\891\\-1297
\end{tabular}\endgroup%
{$\left.\llap{\phantom{%
\begingroup \smaller\smaller\smaller\begin{tabular}{@{}c@{}}%
0\\0\\0
\end{tabular}\endgroup%
}}\!\right]$}%
\EasyButWeakLineBreak%
{$\left[\!\llap{\phantom{%
\begingroup \smaller\smaller\smaller\begin{tabular}{@{}c@{}}%
0\\0\\0
\end{tabular}\endgroup%
}}\right.$}%
\begingroup \smaller\smaller\smaller\begin{tabular}{@{}c@{}}%
57\\152\\-174
\end{tabular}\endgroup%
\HardButStrongLineBreak\kern3pt%
\begingroup \smaller\smaller\smaller\begin{tabular}{@{}c@{}}%
65\\174\\-198
\end{tabular}\endgroup%
\HardButStrongLineBreak\kern3pt%
\begingroup \smaller\smaller\smaller\begin{tabular}{@{}c@{}}%
-165\\-440\\504
\end{tabular}\endgroup%
\HardButStrongLineBreak\kern3pt%
\begingroup \smaller\smaller\smaller\begin{tabular}{@{}c@{}}%
-527\\-1408\\1608
\end{tabular}\endgroup%
\HardButStrongLineBreak\kern3pt%
\begingroup \smaller\smaller\smaller\begin{tabular}{@{}c@{}}%
-239\\-639\\729
\end{tabular}\endgroup%
\HardButStrongLineBreak\kern3pt%
\begingroup \smaller\smaller\smaller\begin{tabular}{@{}c@{}}%
-62\\-166\\189
\end{tabular}\endgroup%
{$\left.\llap{\phantom{%
\begingroup \smaller\smaller\smaller\begin{tabular}{@{}c@{}}%
0\\0\\0
\end{tabular}\endgroup%
}}\!\right]$}%
%
%
\hbox{}\par\smallskip%
%
%
\leavevmode%
${L_{169.3}}$%
{} : {$1\above{1pt}{1pt}{2}{0}16\above{1pt}{1pt}{-}{3}{\cdot}1\above{1pt}{1pt}{-}{}3\above{1pt}{1pt}{1}{}9\above{1pt}{1pt}{1}{}$}\spacer%
\instructions{3}%
\EasyButWeakLineBreak%
{${12}\above{1pt}{1pt}{l}{2}{9}\above{1pt}{1pt}{}{2}{48}\above{1pt}{1pt}{6,1}{\infty}{48}\above{1pt}{1pt}{*}{2}{36}\above{1pt}{1pt}{l}{2}{3}\above{1pt}{1pt}{24,19}{\infty}$}%
\nopagebreak\par%
\nopagebreak\par\leavevmode%
{$\left[\!\llap{\phantom{%
\begingroup \smaller\smaller\smaller\begin{tabular}{@{}c@{}}%
0\\0\\0
\end{tabular}\endgroup%
}}\right.$}%
\begingroup \smaller\smaller\smaller\begin{tabular}{@{}c@{}}%
-100944\\-20592\\-11808
\end{tabular}\endgroup%
\kern3pt%
\begingroup \smaller\smaller\smaller\begin{tabular}{@{}c@{}}%
-20592\\-4128\\-2487
\end{tabular}\endgroup%
\kern3pt%
\begingroup \smaller\smaller\smaller\begin{tabular}{@{}c@{}}%
-11808\\-2487\\-1297
\end{tabular}\endgroup%
{$\left.\llap{\phantom{%
\begingroup \smaller\smaller\smaller\begin{tabular}{@{}c@{}}%
0\\0\\0
\end{tabular}\endgroup%
}}\!\right]$}%
\EasyButWeakLineBreak%
{$\left[\!\llap{\phantom{%
\begingroup \smaller\smaller\smaller\begin{tabular}{@{}c@{}}%
0\\0\\0
\end{tabular}\endgroup%
}}\right.$}%
\begingroup \smaller\smaller\smaller\begin{tabular}{@{}c@{}}%
-95\\304\\282
\end{tabular}\endgroup%
\HardButStrongLineBreak\kern3pt%
\begingroup \smaller\smaller\smaller\begin{tabular}{@{}c@{}}%
-106\\339\\315
\end{tabular}\endgroup%
\HardButStrongLineBreak\kern3pt%
\begingroup \smaller\smaller\smaller\begin{tabular}{@{}c@{}}%
275\\-880\\-816
\end{tabular}\endgroup%
\HardButStrongLineBreak\kern3pt%
\begingroup \smaller\smaller\smaller\begin{tabular}{@{}c@{}}%
1293\\-4136\\-3840
\end{tabular}\endgroup%
\HardButStrongLineBreak\kern3pt%
\begingroup \smaller\smaller\smaller\begin{tabular}{@{}c@{}}%
1315\\-4206\\-3906
\end{tabular}\endgroup%
\HardButStrongLineBreak\kern3pt%
\begingroup \smaller\smaller\smaller\begin{tabular}{@{}c@{}}%
207\\-662\\-615
\end{tabular}\endgroup%
{$\left.\llap{\phantom{%
\begingroup \smaller\smaller\smaller\begin{tabular}{@{}c@{}}%
0\\0\\0
\end{tabular}\endgroup%
}}\!\right]$}%
%
%
\hbox{}\par\smallskip%
%
%
\leavevmode%
${L_{169.4}}$%
{} : {$1\above{1pt}{1pt}{-}{3}4\above{1pt}{1pt}{1}{7}16\above{1pt}{1pt}{1}{1}{\cdot}1\above{1pt}{1pt}{-}{}3\above{1pt}{1pt}{1}{}9\above{1pt}{1pt}{1}{}$}\spacer%
\instructions{3}%
\EasyButWeakLineBreak%
{${12}\above{1pt}{1pt}{s}{2}{144}\above{1pt}{1pt}{l}{2}{12}\above{1pt}{1pt}{24,13}{\infty}{12}\above{1pt}{1pt}{}{2}{144}\above{1pt}{1pt}{}{2}{3}\above{1pt}{1pt}{12,1}{\infty}$}%
\nopagebreak\par%
\nopagebreak\par\leavevmode%
{$\left[\!\llap{\phantom{%
\begingroup \smaller\smaller\smaller\begin{tabular}{@{}c@{}}%
0\\0\\0
\end{tabular}\endgroup%
}}\right.$}%
\begingroup \smaller\smaller\smaller\begin{tabular}{@{}c@{}}%
-227952\\-35280\\9072
\end{tabular}\endgroup%
\kern3pt%
\begingroup \smaller\smaller\smaller\begin{tabular}{@{}c@{}}%
-35280\\-5460\\1404
\end{tabular}\endgroup%
\kern3pt%
\begingroup \smaller\smaller\smaller\begin{tabular}{@{}c@{}}%
9072\\1404\\-361
\end{tabular}\endgroup%
{$\left.\llap{\phantom{%
\begingroup \smaller\smaller\smaller\begin{tabular}{@{}c@{}}%
0\\0\\0
\end{tabular}\endgroup%
}}\!\right]$}%
\EasyButWeakLineBreak%
{$\left[\!\llap{\phantom{%
\begingroup \smaller\smaller\smaller\begin{tabular}{@{}c@{}}%
0\\0\\0
\end{tabular}\endgroup%
}}\right.$}%
\begingroup \smaller\smaller\smaller\begin{tabular}{@{}c@{}}%
1\\-8\\-6
\end{tabular}\endgroup%
\HardButStrongLineBreak\kern3pt%
\begingroup \smaller\smaller\smaller\begin{tabular}{@{}c@{}}%
-1\\-12\\-72
\end{tabular}\endgroup%
\HardButStrongLineBreak\kern3pt%
\begingroup \smaller\smaller\smaller\begin{tabular}{@{}c@{}}%
-2\\13\\0
\end{tabular}\endgroup%
\HardButStrongLineBreak\kern3pt%
\begingroup \smaller\smaller\smaller\begin{tabular}{@{}c@{}}%
-3\\35\\60
\end{tabular}\endgroup%
\HardButStrongLineBreak\kern3pt%
\begingroup \smaller\smaller\smaller\begin{tabular}{@{}c@{}}%
-7\\120\\288
\end{tabular}\endgroup%
\HardButStrongLineBreak\kern3pt%
\begingroup \smaller\smaller\smaller\begin{tabular}{@{}c@{}}%
0\\7\\27
\end{tabular}\endgroup%
{$\left.\llap{\phantom{%
\begingroup \smaller\smaller\smaller\begin{tabular}{@{}c@{}}%
0\\0\\0
\end{tabular}\endgroup%
}}\!\right]$}%
%
%
\hbox{}\par\smallskip%
%
%
\leavevmode%
${L_{169.5}}$%
{} : {$1\above{1pt}{1pt}{-}{3}4\above{1pt}{1pt}{1}{1}16\above{1pt}{1pt}{1}{7}{\cdot}1\above{1pt}{1pt}{-}{}3\above{1pt}{1pt}{1}{}9\above{1pt}{1pt}{1}{}$}\spacer%
\instructions{3}%
\EasyButWeakLineBreak%
{${12}\above{1pt}{1pt}{l}{2}{36}\above{1pt}{1pt}{}{2}{48}\above{1pt}{1pt}{12,1}{\infty}{48}\above{1pt}{1pt}{l}{2}{36}\above{1pt}{1pt}{}{2}{3}\above{1pt}{1pt}{12,7}{\infty}$}%
\nopagebreak\par%
\nopagebreak\par\leavevmode%
{$\left[\!\llap{\phantom{%
\begingroup \smaller\smaller\smaller\begin{tabular}{@{}c@{}}%
0\\0\\0
\end{tabular}\endgroup%
}}\right.$}%
\begingroup \smaller\smaller\smaller\begin{tabular}{@{}c@{}}%
-30096\\15984\\-6336
\end{tabular}\endgroup%
\kern3pt%
\begingroup \smaller\smaller\smaller\begin{tabular}{@{}c@{}}%
15984\\-7980\\3228
\end{tabular}\endgroup%
\kern3pt%
\begingroup \smaller\smaller\smaller\begin{tabular}{@{}c@{}}%
-6336\\3228\\-1297
\end{tabular}\endgroup%
{$\left.\llap{\phantom{%
\begingroup \smaller\smaller\smaller\begin{tabular}{@{}c@{}}%
0\\0\\0
\end{tabular}\endgroup%
}}\!\right]$}%
\EasyButWeakLineBreak%
{$\left[\!\llap{\phantom{%
\begingroup \smaller\smaller\smaller\begin{tabular}{@{}c@{}}%
0\\0\\0
\end{tabular}\endgroup%
}}\right.$}%
\begingroup \smaller\smaller\smaller\begin{tabular}{@{}c@{}}%
-19\\76\\282
\end{tabular}\endgroup%
\HardButStrongLineBreak\kern3pt%
\begingroup \smaller\smaller\smaller\begin{tabular}{@{}c@{}}%
-22\\87\\324
\end{tabular}\endgroup%
\HardButStrongLineBreak\kern3pt%
\begingroup \smaller\smaller\smaller\begin{tabular}{@{}c@{}}%
55\\-220\\-816
\end{tabular}\endgroup%
\HardButStrongLineBreak\kern3pt%
\begingroup \smaller\smaller\smaller\begin{tabular}{@{}c@{}}%
177\\-704\\-2616
\end{tabular}\endgroup%
\HardButStrongLineBreak\kern3pt%
\begingroup \smaller\smaller\smaller\begin{tabular}{@{}c@{}}%
161\\-639\\-2376
\end{tabular}\endgroup%
\HardButStrongLineBreak\kern3pt%
\begingroup \smaller\smaller\smaller\begin{tabular}{@{}c@{}}%
21\\-83\\-309
\end{tabular}\endgroup%
{$\left.\llap{\phantom{%
\begingroup \smaller\smaller\smaller\begin{tabular}{@{}c@{}}%
0\\0\\0
\end{tabular}\endgroup%
}}\!\right]$}%

\medskip%
%
\leavevmode\llap{}%
$W_{170}$%
\qquad\llap{8} lattices, $\chi=12$%
\hfill%
$22|22\slashinfty\rtimes D_{2}$%
\nopagebreak\smallskip\hrule\nopagebreak\medskip%
%
%
\leavevmode%
${L_{170.1}}$%
{} : {$1\above{1pt}{1pt}{-2}{2}16\above{1pt}{1pt}{1}{1}{\cdot}1\above{1pt}{1pt}{-}{}3\above{1pt}{1pt}{1}{}9\above{1pt}{1pt}{1}{}$}\spacer%
\instructions{3}%
\EasyButWeakLineBreak%
{${12}\above{1pt}{1pt}{*}{2}{144}\above{1pt}{1pt}{b}{2}{2}\above{1pt}{1pt}{l}{2}{144}\above{1pt}{1pt}{}{2}{3}\above{1pt}{1pt}{24,1}{\infty}$}%
\nopagebreak\par%
\nopagebreak\par\leavevmode%
{$\left[\!\llap{\phantom{%
\begingroup \smaller\smaller\smaller\begin{tabular}{@{}c@{}}%
0\\0\\0
\end{tabular}\endgroup%
}}\right.$}%
\begingroup \smaller\smaller\smaller\begin{tabular}{@{}c@{}}%
-1105776\\18288\\7200
\end{tabular}\endgroup%
\kern3pt%
\begingroup \smaller\smaller\smaller\begin{tabular}{@{}c@{}}%
18288\\-285\\-123
\end{tabular}\endgroup%
\kern3pt%
\begingroup \smaller\smaller\smaller\begin{tabular}{@{}c@{}}%
7200\\-123\\-46
\end{tabular}\endgroup%
{$\left.\llap{\phantom{%
\begingroup \smaller\smaller\smaller\begin{tabular}{@{}c@{}}%
0\\0\\0
\end{tabular}\endgroup%
}}\!\right]$}%
\EasyButWeakLineBreak%
{$\left[\!\llap{\phantom{%
\begingroup \smaller\smaller\smaller\begin{tabular}{@{}c@{}}%
0\\0\\0
\end{tabular}\endgroup%
}}\right.$}%
\begingroup \smaller\smaller\smaller\begin{tabular}{@{}c@{}}%
-7\\-154\\-684
\end{tabular}\endgroup%
\HardButStrongLineBreak\kern3pt%
\begingroup \smaller\smaller\smaller\begin{tabular}{@{}c@{}}%
-11\\-240\\-1080
\end{tabular}\endgroup%
\HardButStrongLineBreak\kern3pt%
\begingroup \smaller\smaller\smaller\begin{tabular}{@{}c@{}}%
3\\66\\293
\end{tabular}\endgroup%
\HardButStrongLineBreak\kern3pt%
\begingroup \smaller\smaller\smaller\begin{tabular}{@{}c@{}}%
103\\2256\\10080
\end{tabular}\endgroup%
\HardButStrongLineBreak\kern3pt%
\begingroup \smaller\smaller\smaller\begin{tabular}{@{}c@{}}%
6\\131\\588
\end{tabular}\endgroup%
{$\left.\llap{\phantom{%
\begingroup \smaller\smaller\smaller\begin{tabular}{@{}c@{}}%
0\\0\\0
\end{tabular}\endgroup%
}}\!\right]$}%

\medskip%
%
\leavevmode\llap{}%
$W_{171}$%
\qquad\llap{8} lattices, $\chi=48$%
\hfill%
$\infty2|2\infty|\infty2|2\infty|\rtimes D_{4}$%
\nopagebreak\smallskip\hrule\nopagebreak\medskip%
%
%
\leavevmode%
${L_{171.1}}$%
{} : {$1\above{1pt}{1pt}{2}{{\rm II}}16\above{1pt}{1pt}{-}{3}{\cdot}1\above{1pt}{1pt}{1}{}3\above{1pt}{1pt}{1}{}9\above{1pt}{1pt}{-}{}$}\spacer%
\instructions{3}%
\EasyButWeakLineBreak%
{${48}\above{1pt}{1pt}{3,2}{\infty a}{48}\above{1pt}{1pt}{r}{2}{18}\above{1pt}{1pt}{b}{2}{48}\above{1pt}{1pt}{6,1}{\infty z}$}\relax$\,(\times2)$%
\nopagebreak\par%
\nopagebreak\par\leavevmode%
{$\left[\!\llap{\phantom{%
\begingroup \smaller\smaller\smaller\begin{tabular}{@{}c@{}}%
0\\0\\0
\end{tabular}\endgroup%
}}\right.$}%
\begingroup \smaller\smaller\smaller\begin{tabular}{@{}c@{}}%
-8784\\1152\\2304
\end{tabular}\endgroup%
\kern3pt%
\begingroup \smaller\smaller\smaller\begin{tabular}{@{}c@{}}%
1152\\-132\\-327
\end{tabular}\endgroup%
\kern3pt%
\begingroup \smaller\smaller\smaller\begin{tabular}{@{}c@{}}%
2304\\-327\\-572
\end{tabular}\endgroup%
{$\left.\llap{\phantom{%
\begingroup \smaller\smaller\smaller\begin{tabular}{@{}c@{}}%
0\\0\\0
\end{tabular}\endgroup%
}}\!\right]$}%
\hfil\penalty500%
{$\left[\!\llap{\phantom{%
\begingroup \smaller\smaller\smaller\begin{tabular}{@{}c@{}}%
0\\0\\0
\end{tabular}\endgroup%
}}\right.$}%
\begingroup \smaller\smaller\smaller\begin{tabular}{@{}c@{}}%
10943\\33120\\25056
\end{tabular}\endgroup%
\kern3pt%
\begingroup \smaller\smaller\smaller\begin{tabular}{@{}c@{}}%
-1862\\-5636\\-4263
\end{tabular}\endgroup%
\kern3pt%
\begingroup \smaller\smaller\smaller\begin{tabular}{@{}c@{}}%
-2318\\-7015\\-5308
\end{tabular}\endgroup%
{$\left.\llap{\phantom{%
\begingroup \smaller\smaller\smaller\begin{tabular}{@{}c@{}}%
0\\0\\0
\end{tabular}\endgroup%
}}\!\right]$}%
\EasyButWeakLineBreak%
{$\left[\!\llap{\phantom{%
\begingroup \smaller\smaller\smaller\begin{tabular}{@{}c@{}}%
0\\0\\0
\end{tabular}\endgroup%
}}\right.$}%
\begingroup \smaller\smaller\smaller\begin{tabular}{@{}c@{}}%
179\\544\\408
\end{tabular}\endgroup%
\HardButStrongLineBreak\kern3pt%
\begingroup \smaller\smaller\smaller\begin{tabular}{@{}c@{}}%
629\\1904\\1440
\end{tabular}\endgroup%
\HardButStrongLineBreak\kern3pt%
\begingroup \smaller\smaller\smaller\begin{tabular}{@{}c@{}}%
263\\795\\603
\end{tabular}\endgroup%
\HardButStrongLineBreak\kern3pt%
\begingroup \smaller\smaller\smaller\begin{tabular}{@{}c@{}}%
575\\1736\\1320
\end{tabular}\endgroup%
{$\left.\llap{\phantom{%
\begingroup \smaller\smaller\smaller\begin{tabular}{@{}c@{}}%
0\\0\\0
\end{tabular}\endgroup%
}}\!\right]$}%

\medskip%
%
\leavevmode\llap{}%
$W_{172}$%
\qquad\llap{16} lattices, $\chi=6$%
\hfill%
$2|22\slashtwo2\rtimes D_{2}$%
\nopagebreak\smallskip\hrule\nopagebreak\medskip%
%
%
\leavevmode%
${L_{172.1}}$%
{} : {$1\above{1pt}{1pt}{2}{2}8\above{1pt}{1pt}{-}{5}{\cdot}1\above{1pt}{1pt}{-}{}3\above{1pt}{1pt}{1}{}9\above{1pt}{1pt}{-}{}$}\spacer%
\instructions{2,m}%
\EasyButWeakLineBreak%
{${8}\above{1pt}{1pt}{*}{2}{12}\above{1pt}{1pt}{*}{2}{72}\above{1pt}{1pt}{b}{2}{2}\above{1pt}{1pt}{s}{2}{18}\above{1pt}{1pt}{b}{2}$}%
\nopagebreak\par%
\nopagebreak\par\leavevmode%
{$\left[\!\llap{\phantom{%
\begingroup \smaller\smaller\smaller\begin{tabular}{@{}c@{}}%
0\\0\\0
\end{tabular}\endgroup%
}}\right.$}%
\begingroup \smaller\smaller\smaller\begin{tabular}{@{}c@{}}%
-5976\\216\\216
\end{tabular}\endgroup%
\kern3pt%
\begingroup \smaller\smaller\smaller\begin{tabular}{@{}c@{}}%
216\\-6\\-9
\end{tabular}\endgroup%
\kern3pt%
\begingroup \smaller\smaller\smaller\begin{tabular}{@{}c@{}}%
216\\-9\\-7
\end{tabular}\endgroup%
{$\left.\llap{\phantom{%
\begingroup \smaller\smaller\smaller\begin{tabular}{@{}c@{}}%
0\\0\\0
\end{tabular}\endgroup%
}}\!\right]$}%
\EasyButWeakLineBreak%
{$\left[\!\llap{\phantom{%
\begingroup \smaller\smaller\smaller\begin{tabular}{@{}c@{}}%
0\\0\\0
\end{tabular}\endgroup%
}}\right.$}%
\begingroup \smaller\smaller\smaller\begin{tabular}{@{}c@{}}%
-1\\-12\\-16
\end{tabular}\endgroup%
\HardButStrongLineBreak\kern3pt%
\begingroup \smaller\smaller\smaller\begin{tabular}{@{}c@{}}%
-1\\-10\\-18
\end{tabular}\endgroup%
\HardButStrongLineBreak\kern3pt%
\begingroup \smaller\smaller\smaller\begin{tabular}{@{}c@{}}%
5\\60\\72
\end{tabular}\endgroup%
\HardButStrongLineBreak\kern3pt%
\begingroup \smaller\smaller\smaller\begin{tabular}{@{}c@{}}%
1\\11\\16
\end{tabular}\endgroup%
\HardButStrongLineBreak\kern3pt%
\begingroup \smaller\smaller\smaller\begin{tabular}{@{}c@{}}%
1\\9\\18
\end{tabular}\endgroup%
{$\left.\llap{\phantom{%
\begingroup \smaller\smaller\smaller\begin{tabular}{@{}c@{}}%
0\\0\\0
\end{tabular}\endgroup%
}}\!\right]$}%
%
%
\hbox{}\par\smallskip%
%
%
\leavevmode%
${L_{172.2}}$%
{} : {$1\above{1pt}{1pt}{-2}{6}16\above{1pt}{1pt}{1}{1}{\cdot}1\above{1pt}{1pt}{1}{}3\above{1pt}{1pt}{-}{}9\above{1pt}{1pt}{1}{}$}\spacer%
\instructions{3,m}%
\EasyButWeakLineBreak%
{${16}\above{1pt}{1pt}{b}{2}{6}\above{1pt}{1pt}{l}{2}{144}\above{1pt}{1pt}{}{2}{1}\above{1pt}{1pt}{r}{2}{36}\above{1pt}{1pt}{*}{2}$}%
\nopagebreak\par%
shares genus with 3-dual\nopagebreak\par%
\nopagebreak\par\leavevmode%
{$\left[\!\llap{\phantom{%
\begingroup \smaller\smaller\smaller\begin{tabular}{@{}c@{}}%
0\\0\\0
\end{tabular}\endgroup%
}}\right.$}%
\begingroup \smaller\smaller\smaller\begin{tabular}{@{}c@{}}%
13968\\-288\\-144
\end{tabular}\endgroup%
\kern3pt%
\begingroup \smaller\smaller\smaller\begin{tabular}{@{}c@{}}%
-288\\6\\3
\end{tabular}\endgroup%
\kern3pt%
\begingroup \smaller\smaller\smaller\begin{tabular}{@{}c@{}}%
-144\\3\\1
\end{tabular}\endgroup%
{$\left.\llap{\phantom{%
\begingroup \smaller\smaller\smaller\begin{tabular}{@{}c@{}}%
0\\0\\0
\end{tabular}\endgroup%
}}\!\right]$}%
\EasyButWeakLineBreak%
{$\left[\!\llap{\phantom{%
\begingroup \smaller\smaller\smaller\begin{tabular}{@{}c@{}}%
0\\0\\0
\end{tabular}\endgroup%
}}\right.$}%
\begingroup \smaller\smaller\smaller\begin{tabular}{@{}c@{}}%
-1\\-40\\-16
\end{tabular}\endgroup%
\HardButStrongLineBreak\kern3pt%
\begingroup \smaller\smaller\smaller\begin{tabular}{@{}c@{}}%
0\\1\\0
\end{tabular}\endgroup%
\HardButStrongLineBreak\kern3pt%
\begingroup \smaller\smaller\smaller\begin{tabular}{@{}c@{}}%
1\\48\\0
\end{tabular}\endgroup%
\HardButStrongLineBreak\kern3pt%
\begingroup \smaller\smaller\smaller\begin{tabular}{@{}c@{}}%
0\\0\\-1
\end{tabular}\endgroup%
\HardButStrongLineBreak\kern3pt%
\begingroup \smaller\smaller\smaller\begin{tabular}{@{}c@{}}%
-1\\-42\\-18
\end{tabular}\endgroup%
{$\left.\llap{\phantom{%
\begingroup \smaller\smaller\smaller\begin{tabular}{@{}c@{}}%
0\\0\\0
\end{tabular}\endgroup%
}}\!\right]$}%

\medskip%
%
\leavevmode\llap{}%
$W_{173}$%
\qquad\llap{8} lattices, $\chi=40$%
\hfill%
$3\infty\infty3\infty\infty\rtimes C_{2}$%
\nopagebreak\smallskip\hrule\nopagebreak\medskip%
%
%
\leavevmode%
${L_{173.1}}$%
{} : {$1\above{1pt}{1pt}{-2}{{\rm II}}8\above{1pt}{1pt}{-}{5}{\cdot}1\above{1pt}{1pt}{-2}{}25\above{1pt}{1pt}{1}{}$}\spacer%
\instructions{2}%
\EasyButWeakLineBreak%
{${2}\above{1pt}{1pt}{-}{3}{2}\above{1pt}{1pt}{20,17}{\infty b}{8}\above{1pt}{1pt}{10,7}{\infty z}$}\relax$\,(\times2)$%
\nopagebreak\par%
\nopagebreak\par\leavevmode%
{$\left[\!\llap{\phantom{%
\begingroup \smaller\smaller\smaller\begin{tabular}{@{}c@{}}%
0\\0\\0
\end{tabular}\endgroup%
}}\right.$}%
\begingroup \smaller\smaller\smaller\begin{tabular}{@{}c@{}}%
-45400\\600\\800
\end{tabular}\endgroup%
\kern3pt%
\begingroup \smaller\smaller\smaller\begin{tabular}{@{}c@{}}%
600\\-6\\-11
\end{tabular}\endgroup%
\kern3pt%
\begingroup \smaller\smaller\smaller\begin{tabular}{@{}c@{}}%
800\\-11\\-14
\end{tabular}\endgroup%
{$\left.\llap{\phantom{%
\begingroup \smaller\smaller\smaller\begin{tabular}{@{}c@{}}%
0\\0\\0
\end{tabular}\endgroup%
}}\!\right]$}%
\hfil\penalty500%
{$\left[\!\llap{\phantom{%
\begingroup \smaller\smaller\smaller\begin{tabular}{@{}c@{}}%
0\\0\\0
\end{tabular}\endgroup%
}}\right.$}%
\begingroup \smaller\smaller\smaller\begin{tabular}{@{}c@{}}%
799\\7200\\39200
\end{tabular}\endgroup%
\kern3pt%
\begingroup \smaller\smaller\smaller\begin{tabular}{@{}c@{}}%
-7\\-64\\-343
\end{tabular}\endgroup%
\kern3pt%
\begingroup \smaller\smaller\smaller\begin{tabular}{@{}c@{}}%
-15\\-135\\-736
\end{tabular}\endgroup%
{$\left.\llap{\phantom{%
\begingroup \smaller\smaller\smaller\begin{tabular}{@{}c@{}}%
0\\0\\0
\end{tabular}\endgroup%
}}\!\right]$}%
\EasyButWeakLineBreak%
{$\left[\!\llap{\phantom{%
\begingroup \smaller\smaller\smaller\begin{tabular}{@{}c@{}}%
0\\0\\0
\end{tabular}\endgroup%
}}\right.$}%
\begingroup \smaller\smaller\smaller\begin{tabular}{@{}c@{}}%
2\\16\\99
\end{tabular}\endgroup%
\HardButStrongLineBreak\kern3pt%
\begingroup \smaller\smaller\smaller\begin{tabular}{@{}c@{}}%
8\\71\\393
\end{tabular}\endgroup%
\HardButStrongLineBreak\kern3pt%
\begingroup \smaller\smaller\smaller\begin{tabular}{@{}c@{}}%
5\\48\\244
\end{tabular}\endgroup%
{$\left.\llap{\phantom{%
\begingroup \smaller\smaller\smaller\begin{tabular}{@{}c@{}}%
0\\0\\0
\end{tabular}\endgroup%
}}\!\right]$}%

\medskip%
%
\leavevmode\llap{}%
$W_{174}$%
\qquad\llap{8} lattices, $\chi=30$%
\hfill%
$2\infty2\infty2\infty\rtimes C_{3}$%
\nopagebreak\smallskip\hrule\nopagebreak\medskip%
%
%
\leavevmode%
${L_{174.1}}$%
{} : {$1\above{1pt}{1pt}{-2}{{\rm II}}8\above{1pt}{1pt}{-}{5}{\cdot}1\above{1pt}{1pt}{2}{}25\above{1pt}{1pt}{-}{}$}\spacer%
\instructions{2}%
\EasyButWeakLineBreak%
{${50}\above{1pt}{1pt}{b}{2}{8}\above{1pt}{1pt}{10,1}{\infty z}{2}\above{1pt}{1pt}{b}{2}{8}\above{1pt}{1pt}{10,9}{\infty z}{2}\above{1pt}{1pt}{b}{2}{200}\above{1pt}{1pt}{2,1}{\infty z}$}%
\nopagebreak\par%
\nopagebreak\par\leavevmode%
{$\left[\!\llap{\phantom{%
\begingroup \smaller\smaller\smaller\begin{tabular}{@{}c@{}}%
0\\0\\0
\end{tabular}\endgroup%
}}\right.$}%
\begingroup \smaller\smaller\smaller\begin{tabular}{@{}c@{}}%
12200\\-5000\\200
\end{tabular}\endgroup%
\kern3pt%
\begingroup \smaller\smaller\smaller\begin{tabular}{@{}c@{}}%
-5000\\2042\\-55
\end{tabular}\endgroup%
\kern3pt%
\begingroup \smaller\smaller\smaller\begin{tabular}{@{}c@{}}%
200\\-55\\-98
\end{tabular}\endgroup%
{$\left.\llap{\phantom{%
\begingroup \smaller\smaller\smaller\begin{tabular}{@{}c@{}}%
0\\0\\0
\end{tabular}\endgroup%
}}\!\right]$}%
\EasyButWeakLineBreak%
{$\left[\!\llap{\phantom{%
\begingroup \smaller\smaller\smaller\begin{tabular}{@{}c@{}}%
0\\0\\0
\end{tabular}\endgroup%
}}\right.$}%
\begingroup \smaller\smaller\smaller\begin{tabular}{@{}c@{}}%
-152\\-375\\-100
\end{tabular}\endgroup%
\HardButStrongLineBreak\kern3pt%
\begingroup \smaller\smaller\smaller\begin{tabular}{@{}c@{}}%
-73\\-180\\-48
\end{tabular}\endgroup%
\HardButStrongLineBreak\kern3pt%
\begingroup \smaller\smaller\smaller\begin{tabular}{@{}c@{}}%
75\\185\\49
\end{tabular}\endgroup%
\HardButStrongLineBreak\kern3pt%
\begingroup \smaller\smaller\smaller\begin{tabular}{@{}c@{}}%
459\\1132\\300
\end{tabular}\endgroup%
\HardButStrongLineBreak\kern3pt%
\begingroup \smaller\smaller\smaller\begin{tabular}{@{}c@{}}%
176\\434\\115
\end{tabular}\endgroup%
\HardButStrongLineBreak\kern3pt%
\begingroup \smaller\smaller\smaller\begin{tabular}{@{}c@{}}%
1379\\3400\\900
\end{tabular}\endgroup%
{$\left.\llap{\phantom{%
\begingroup \smaller\smaller\smaller\begin{tabular}{@{}c@{}}%
0\\0\\0
\end{tabular}\endgroup%
}}\!\right]$}%

\medskip%
%
\leavevmode\llap{}%
$W_{175}$%
\qquad\llap{27} lattices, $\chi=12$%
\hfill%
$\slashtwo2|2\slashtwo2|2\rtimes D_{4}$%
\nopagebreak\smallskip\hrule\nopagebreak\medskip%
%
%
\leavevmode%
${L_{175.1}}$%
{} : {$1\above{1pt}{1pt}{-2}{4}8\above{1pt}{1pt}{1}{7}{\cdot}1\above{1pt}{1pt}{1}{}3\above{1pt}{1pt}{1}{}9\above{1pt}{1pt}{1}{}$}\spacer%
\instructions{2}%
\EasyButWeakLineBreak%
{${36}\above{1pt}{1pt}{l}{2}{1}\above{1pt}{1pt}{}{2}{3}\above{1pt}{1pt}{}{2}{9}\above{1pt}{1pt}{r}{2}{4}\above{1pt}{1pt}{*}{2}{12}\above{1pt}{1pt}{*}{2}$}%
\nopagebreak\par%
\nopagebreak\par\leavevmode%
{$\left[\!\llap{\phantom{%
\begingroup \smaller\smaller\smaller\begin{tabular}{@{}c@{}}%
0\\0\\0
\end{tabular}\endgroup%
}}\right.$}%
\begingroup \smaller\smaller\smaller\begin{tabular}{@{}c@{}}%
-1800\\288\\-72
\end{tabular}\endgroup%
\kern3pt%
\begingroup \smaller\smaller\smaller\begin{tabular}{@{}c@{}}%
288\\-15\\-6
\end{tabular}\endgroup%
\kern3pt%
\begingroup \smaller\smaller\smaller\begin{tabular}{@{}c@{}}%
-72\\-6\\7
\end{tabular}\endgroup%
{$\left.\llap{\phantom{%
\begingroup \smaller\smaller\smaller\begin{tabular}{@{}c@{}}%
0\\0\\0
\end{tabular}\endgroup%
}}\!\right]$}%
\EasyButWeakLineBreak%
{$\left[\!\llap{\phantom{%
\begingroup \smaller\smaller\smaller\begin{tabular}{@{}c@{}}%
0\\0\\0
\end{tabular}\endgroup%
}}\right.$}%
\begingroup \smaller\smaller\smaller\begin{tabular}{@{}c@{}}%
-1\\-12\\-18
\end{tabular}\endgroup%
\HardButStrongLineBreak\kern3pt%
\begingroup \smaller\smaller\smaller\begin{tabular}{@{}c@{}}%
1\\11\\20
\end{tabular}\endgroup%
\HardButStrongLineBreak\kern3pt%
\begingroup \smaller\smaller\smaller\begin{tabular}{@{}c@{}}%
2\\22\\39
\end{tabular}\endgroup%
\HardButStrongLineBreak\kern3pt%
\begingroup \smaller\smaller\smaller\begin{tabular}{@{}c@{}}%
2\\21\\36
\end{tabular}\endgroup%
\HardButStrongLineBreak\kern3pt%
\begingroup \smaller\smaller\smaller\begin{tabular}{@{}c@{}}%
-1\\-12\\-22
\end{tabular}\endgroup%
\HardButStrongLineBreak\kern3pt%
\begingroup \smaller\smaller\smaller\begin{tabular}{@{}c@{}}%
-3\\-34\\-60
\end{tabular}\endgroup%
{$\left.\llap{\phantom{%
\begingroup \smaller\smaller\smaller\begin{tabular}{@{}c@{}}%
0\\0\\0
\end{tabular}\endgroup%
}}\!\right]$}%
%
%
\hbox{}\par\smallskip%
%
%
\leavevmode%
${L_{175.2}}$%
{} : {$1\above{1pt}{1pt}{-2}{2}8\above{1pt}{1pt}{1}{1}{\cdot}1\above{1pt}{1pt}{1}{}3\above{1pt}{1pt}{1}{}9\above{1pt}{1pt}{1}{}$}\spacer%
\instructions{m}%
\EasyButWeakLineBreak%
{${36}\above{1pt}{1pt}{*}{2}{4}\above{1pt}{1pt}{l}{2}{3}\above{1pt}{1pt}{r}{2}$}\relax$\,(\times2)$%
\nopagebreak\par%
\nopagebreak\par\leavevmode%
{$\left[\!\llap{\phantom{%
\begingroup \smaller\smaller\smaller\begin{tabular}{@{}c@{}}%
0\\0\\0
\end{tabular}\endgroup%
}}\right.$}%
\begingroup \smaller\smaller\smaller\begin{tabular}{@{}c@{}}%
8136\\-576\\288
\end{tabular}\endgroup%
\kern3pt%
\begingroup \smaller\smaller\smaller\begin{tabular}{@{}c@{}}%
-576\\39\\-18
\end{tabular}\endgroup%
\kern3pt%
\begingroup \smaller\smaller\smaller\begin{tabular}{@{}c@{}}%
288\\-18\\7
\end{tabular}\endgroup%
{$\left.\llap{\phantom{%
\begingroup \smaller\smaller\smaller\begin{tabular}{@{}c@{}}%
0\\0\\0
\end{tabular}\endgroup%
}}\!\right]$}%
\hfil\penalty500%
{$\left[\!\llap{\phantom{%
\begingroup \smaller\smaller\smaller\begin{tabular}{@{}c@{}}%
0\\0\\0
\end{tabular}\endgroup%
}}\right.$}%
\begingroup \smaller\smaller\smaller\begin{tabular}{@{}c@{}}%
-1\\-96\\-288
\end{tabular}\endgroup%
\kern3pt%
\begingroup \smaller\smaller\smaller\begin{tabular}{@{}c@{}}%
0\\4\\15
\end{tabular}\endgroup%
\kern3pt%
\begingroup \smaller\smaller\smaller\begin{tabular}{@{}c@{}}%
0\\-1\\-4
\end{tabular}\endgroup%
{$\left.\llap{\phantom{%
\begingroup \smaller\smaller\smaller\begin{tabular}{@{}c@{}}%
0\\0\\0
\end{tabular}\endgroup%
}}\!\right]$}%
\EasyButWeakLineBreak%
{$\left[\!\llap{\phantom{%
\begingroup \smaller\smaller\smaller\begin{tabular}{@{}c@{}}%
0\\0\\0
\end{tabular}\endgroup%
}}\right.$}%
\begingroup \smaller\smaller\smaller\begin{tabular}{@{}c@{}}%
-1\\-24\\-18
\end{tabular}\endgroup%
\HardButStrongLineBreak\kern3pt%
\begingroup \smaller\smaller\smaller\begin{tabular}{@{}c@{}}%
1\\22\\16
\end{tabular}\endgroup%
\HardButStrongLineBreak\kern3pt%
\begingroup \smaller\smaller\smaller\begin{tabular}{@{}c@{}}%
1\\22\\15
\end{tabular}\endgroup%
{$\left.\llap{\phantom{%
\begingroup \smaller\smaller\smaller\begin{tabular}{@{}c@{}}%
0\\0\\0
\end{tabular}\endgroup%
}}\!\right]$}%
%
%
\hbox{}\par\smallskip%
%
%
\leavevmode%
${L_{175.3}}$%
{} : {$1\above{1pt}{1pt}{2}{2}8\above{1pt}{1pt}{-}{5}{\cdot}1\above{1pt}{1pt}{1}{}3\above{1pt}{1pt}{1}{}9\above{1pt}{1pt}{1}{}$}\EasyButWeakLineBreak%
{${9}\above{1pt}{1pt}{}{2}{1}\above{1pt}{1pt}{r}{2}{12}\above{1pt}{1pt}{l}{2}$}\relax$\,(\times2)$%
\nopagebreak\par%
\nopagebreak\par\leavevmode%
{$\left[\!\llap{\phantom{%
\begingroup \smaller\smaller\smaller\begin{tabular}{@{}c@{}}%
0\\0\\0
\end{tabular}\endgroup%
}}\right.$}%
\begingroup \smaller\smaller\smaller\begin{tabular}{@{}c@{}}%
-1368\\144\\360
\end{tabular}\endgroup%
\kern3pt%
\begingroup \smaller\smaller\smaller\begin{tabular}{@{}c@{}}%
144\\-15\\-36
\end{tabular}\endgroup%
\kern3pt%
\begingroup \smaller\smaller\smaller\begin{tabular}{@{}c@{}}%
360\\-36\\-71
\end{tabular}\endgroup%
{$\left.\llap{\phantom{%
\begingroup \smaller\smaller\smaller\begin{tabular}{@{}c@{}}%
0\\0\\0
\end{tabular}\endgroup%
}}\!\right]$}%
\hfil\penalty500%
{$\left[\!\llap{\phantom{%
\begingroup \smaller\smaller\smaller\begin{tabular}{@{}c@{}}%
0\\0\\0
\end{tabular}\endgroup%
}}\right.$}%
\begingroup \smaller\smaller\smaller\begin{tabular}{@{}c@{}}%
143\\1680\\-144
\end{tabular}\endgroup%
\kern3pt%
\begingroup \smaller\smaller\smaller\begin{tabular}{@{}c@{}}%
-15\\-176\\15
\end{tabular}\endgroup%
\kern3pt%
\begingroup \smaller\smaller\smaller\begin{tabular}{@{}c@{}}%
-33\\-385\\32
\end{tabular}\endgroup%
{$\left.\llap{\phantom{%
\begingroup \smaller\smaller\smaller\begin{tabular}{@{}c@{}}%
0\\0\\0
\end{tabular}\endgroup%
}}\!\right]$}%
\EasyButWeakLineBreak%
{$\left[\!\llap{\phantom{%
\begingroup \smaller\smaller\smaller\begin{tabular}{@{}c@{}}%
0\\0\\0
\end{tabular}\endgroup%
}}\right.$}%
\begingroup \smaller\smaller\smaller\begin{tabular}{@{}c@{}}%
1\\9\\0
\end{tabular}\endgroup%
\HardButStrongLineBreak\kern3pt%
\begingroup \smaller\smaller\smaller\begin{tabular}{@{}c@{}}%
-1\\-12\\1
\end{tabular}\endgroup%
\HardButStrongLineBreak\kern3pt%
\begingroup \smaller\smaller\smaller\begin{tabular}{@{}c@{}}%
-1\\-10\\0
\end{tabular}\endgroup%
{$\left.\llap{\phantom{%
\begingroup \smaller\smaller\smaller\begin{tabular}{@{}c@{}}%
0\\0\\0
\end{tabular}\endgroup%
}}\!\right]$}%
%
%
\hbox{}\par\smallskip%
%
%
\leavevmode%
${L_{175.4}}$%
{} : {$[1\above{1pt}{1pt}{-}{}2\above{1pt}{1pt}{1}{}]\above{1pt}{1pt}{}{4}16\above{1pt}{1pt}{1}{7}{\cdot}1\above{1pt}{1pt}{1}{}3\above{1pt}{1pt}{1}{}9\above{1pt}{1pt}{1}{}$}\spacer%
\instructions{2}%
\EasyButWeakLineBreak%
{${4}\above{1pt}{1pt}{*}{2}{144}\above{1pt}{1pt}{l}{2}{3}\above{1pt}{1pt}{r}{2}{16}\above{1pt}{1pt}{*}{2}{36}\above{1pt}{1pt}{s}{2}{48}\above{1pt}{1pt}{s}{2}$}%
\nopagebreak\par%
\nopagebreak\par\leavevmode%
{$\left[\!\llap{\phantom{%
\begingroup \smaller\smaller\smaller\begin{tabular}{@{}c@{}}%
0\\0\\0
\end{tabular}\endgroup%
}}\right.$}%
\begingroup \smaller\smaller\smaller\begin{tabular}{@{}c@{}}%
-5904\\864\\288
\end{tabular}\endgroup%
\kern3pt%
\begingroup \smaller\smaller\smaller\begin{tabular}{@{}c@{}}%
864\\-78\\-30
\end{tabular}\endgroup%
\kern3pt%
\begingroup \smaller\smaller\smaller\begin{tabular}{@{}c@{}}%
288\\-30\\-11
\end{tabular}\endgroup%
{$\left.\llap{\phantom{%
\begingroup \smaller\smaller\smaller\begin{tabular}{@{}c@{}}%
0\\0\\0
\end{tabular}\endgroup%
}}\!\right]$}%
\EasyButWeakLineBreak%
{$\left[\!\llap{\phantom{%
\begingroup \smaller\smaller\smaller\begin{tabular}{@{}c@{}}%
0\\0\\0
\end{tabular}\endgroup%
}}\right.$}%
\begingroup \smaller\smaller\smaller\begin{tabular}{@{}c@{}}%
-1\\20\\-82
\end{tabular}\endgroup%
\HardButStrongLineBreak\kern3pt%
\begingroup \smaller\smaller\smaller\begin{tabular}{@{}c@{}}%
-5\\108\\-432
\end{tabular}\endgroup%
\HardButStrongLineBreak\kern3pt%
\begingroup \smaller\smaller\smaller\begin{tabular}{@{}c@{}}%
0\\1\\-3
\end{tabular}\endgroup%
\HardButStrongLineBreak\kern3pt%
\begingroup \smaller\smaller\smaller\begin{tabular}{@{}c@{}}%
1\\-20\\80
\end{tabular}\endgroup%
\HardButStrongLineBreak\kern3pt%
\begingroup \smaller\smaller\smaller\begin{tabular}{@{}c@{}}%
1\\-24\\90
\end{tabular}\endgroup%
\HardButStrongLineBreak\kern3pt%
\begingroup \smaller\smaller\smaller\begin{tabular}{@{}c@{}}%
-1\\16\\-72
\end{tabular}\endgroup%
{$\left.\llap{\phantom{%
\begingroup \smaller\smaller\smaller\begin{tabular}{@{}c@{}}%
0\\0\\0
\end{tabular}\endgroup%
}}\!\right]$}%
%
%
\hbox{}\par\smallskip%
%
%
\leavevmode%
${L_{175.5}}$%
{} : {$[1\above{1pt}{1pt}{1}{}2\above{1pt}{1pt}{1}{}]\above{1pt}{1pt}{}{0}16\above{1pt}{1pt}{-}{3}{\cdot}1\above{1pt}{1pt}{1}{}3\above{1pt}{1pt}{1}{}9\above{1pt}{1pt}{1}{}$}\spacer%
\instructions{m}%
\EasyButWeakLineBreak%
{${1}\above{1pt}{1pt}{r}{2}{144}\above{1pt}{1pt}{*}{2}{12}\above{1pt}{1pt}{*}{2}{16}\above{1pt}{1pt}{l}{2}{9}\above{1pt}{1pt}{}{2}{48}\above{1pt}{1pt}{}{2}$}%
\nopagebreak\par%
\nopagebreak\par\leavevmode%
{$\left[\!\llap{\phantom{%
\begingroup \smaller\smaller\smaller\begin{tabular}{@{}c@{}}%
0\\0\\0
\end{tabular}\endgroup%
}}\right.$}%
\begingroup \smaller\smaller\smaller\begin{tabular}{@{}c@{}}%
-90576\\-6048\\2448
\end{tabular}\endgroup%
\kern3pt%
\begingroup \smaller\smaller\smaller\begin{tabular}{@{}c@{}}%
-6048\\-402\\162
\end{tabular}\endgroup%
\kern3pt%
\begingroup \smaller\smaller\smaller\begin{tabular}{@{}c@{}}%
2448\\162\\-65
\end{tabular}\endgroup%
{$\left.\llap{\phantom{%
\begingroup \smaller\smaller\smaller\begin{tabular}{@{}c@{}}%
0\\0\\0
\end{tabular}\endgroup%
}}\!\right]$}%
\EasyButWeakLineBreak%
{$\left[\!\llap{\phantom{%
\begingroup \smaller\smaller\smaller\begin{tabular}{@{}c@{}}%
0\\0\\0
\end{tabular}\endgroup%
}}\right.$}%
\begingroup \smaller\smaller\smaller\begin{tabular}{@{}c@{}}%
1\\-29\\-35
\end{tabular}\endgroup%
\HardButStrongLineBreak\kern3pt%
\begingroup \smaller\smaller\smaller\begin{tabular}{@{}c@{}}%
5\\-132\\-144
\end{tabular}\endgroup%
\HardButStrongLineBreak\kern3pt%
\begingroup \smaller\smaller\smaller\begin{tabular}{@{}c@{}}%
-1\\32\\42
\end{tabular}\endgroup%
\HardButStrongLineBreak\kern3pt%
\begingroup \smaller\smaller\smaller\begin{tabular}{@{}c@{}}%
-1\\28\\32
\end{tabular}\endgroup%
\HardButStrongLineBreak\kern3pt%
\begingroup \smaller\smaller\smaller\begin{tabular}{@{}c@{}}%
1\\-33\\-45
\end{tabular}\endgroup%
\HardButStrongLineBreak\kern3pt%
\begingroup \smaller\smaller\smaller\begin{tabular}{@{}c@{}}%
5\\-152\\-192
\end{tabular}\endgroup%
{$\left.\llap{\phantom{%
\begingroup \smaller\smaller\smaller\begin{tabular}{@{}c@{}}%
0\\0\\0
\end{tabular}\endgroup%
}}\!\right]$}%
%
%
\hbox{}\par\smallskip%
%
%
\leavevmode%
${L_{175.6}}$%
{} : {$[1\above{1pt}{1pt}{1}{}2\above{1pt}{1pt}{1}{}]\above{1pt}{1pt}{}{6}16\above{1pt}{1pt}{-}{5}{\cdot}1\above{1pt}{1pt}{1}{}3\above{1pt}{1pt}{1}{}9\above{1pt}{1pt}{1}{}$}\spacer%
\instructions{m}%
\EasyButWeakLineBreak%
{${4}\above{1pt}{1pt}{s}{2}{144}\above{1pt}{1pt}{s}{2}{12}\above{1pt}{1pt}{s}{2}{16}\above{1pt}{1pt}{s}{2}{36}\above{1pt}{1pt}{*}{2}{48}\above{1pt}{1pt}{*}{2}$}%
\nopagebreak\par%
\nopagebreak\par\leavevmode%
{$\left[\!\llap{\phantom{%
\begingroup \smaller\smaller\smaller\begin{tabular}{@{}c@{}}%
0\\0\\0
\end{tabular}\endgroup%
}}\right.$}%
\begingroup \smaller\smaller\smaller\begin{tabular}{@{}c@{}}%
8784\\432\\-288
\end{tabular}\endgroup%
\kern3pt%
\begingroup \smaller\smaller\smaller\begin{tabular}{@{}c@{}}%
432\\-6\\-6
\end{tabular}\endgroup%
\kern3pt%
\begingroup \smaller\smaller\smaller\begin{tabular}{@{}c@{}}%
-288\\-6\\7
\end{tabular}\endgroup%
{$\left.\llap{\phantom{%
\begingroup \smaller\smaller\smaller\begin{tabular}{@{}c@{}}%
0\\0\\0
\end{tabular}\endgroup%
}}\!\right]$}%
\EasyButWeakLineBreak%
{$\left[\!\llap{\phantom{%
\begingroup \smaller\smaller\smaller\begin{tabular}{@{}c@{}}%
0\\0\\0
\end{tabular}\endgroup%
}}\right.$}%
\begingroup \smaller\smaller\smaller\begin{tabular}{@{}c@{}}%
-1\\-18\\-58
\end{tabular}\endgroup%
\HardButStrongLineBreak\kern3pt%
\begingroup \smaller\smaller\smaller\begin{tabular}{@{}c@{}}%
-1\\-24\\-72
\end{tabular}\endgroup%
\HardButStrongLineBreak\kern3pt%
\begingroup \smaller\smaller\smaller\begin{tabular}{@{}c@{}}%
1\\16\\54
\end{tabular}\endgroup%
\HardButStrongLineBreak\kern3pt%
\begingroup \smaller\smaller\smaller\begin{tabular}{@{}c@{}}%
1\\16\\56
\end{tabular}\endgroup%
\HardButStrongLineBreak\kern3pt%
\begingroup \smaller\smaller\smaller\begin{tabular}{@{}c@{}}%
-1\\-18\\-54
\end{tabular}\endgroup%
\HardButStrongLineBreak\kern3pt%
\begingroup \smaller\smaller\smaller\begin{tabular}{@{}c@{}}%
-3\\-52\\-168
\end{tabular}\endgroup%
{$\left.\llap{\phantom{%
\begingroup \smaller\smaller\smaller\begin{tabular}{@{}c@{}}%
0\\0\\0
\end{tabular}\endgroup%
}}\!\right]$}%
%
%
\hbox{}\par\smallskip%
%
%
\leavevmode%
${L_{175.7}}$%
{} : {$1\above{1pt}{1pt}{1}{1}8\above{1pt}{1pt}{-}{5}64\above{1pt}{1pt}{1}{1}{\cdot}1\above{1pt}{1pt}{1}{}3\above{1pt}{1pt}{1}{}9\above{1pt}{1pt}{1}{}$}\spacer%
\instructions{3,2}%
\EasyButWeakLineBreak%
{${64}\above{1pt}{1pt}{s}{2}{36}\above{1pt}{1pt}{*}{2}{192}\above{1pt}{1pt}{l}{2}{1}\above{1pt}{1pt}{}{2}{576}\above{1pt}{1pt}{r}{2}{12}\above{1pt}{1pt}{b}{2}$}%
\nopagebreak\par%
shares genus with 2-dual${}\iso{}$3-dual; isometric to own %
2.3-dual\nopagebreak\par%
\nopagebreak\par\leavevmode%
{$\left[\!\llap{\phantom{%
\begingroup \smaller\smaller\smaller\begin{tabular}{@{}c@{}}%
0\\0\\0
\end{tabular}\endgroup%
}}\right.$}%
\begingroup \smaller\smaller\smaller\begin{tabular}{@{}c@{}}%
318528\\576\\-576
\end{tabular}\endgroup%
\kern3pt%
\begingroup \smaller\smaller\smaller\begin{tabular}{@{}c@{}}%
576\\-24\\0
\end{tabular}\endgroup%
\kern3pt%
\begingroup \smaller\smaller\smaller\begin{tabular}{@{}c@{}}%
-576\\0\\1
\end{tabular}\endgroup%
{$\left.\llap{\phantom{%
\begingroup \smaller\smaller\smaller\begin{tabular}{@{}c@{}}%
0\\0\\0
\end{tabular}\endgroup%
}}\!\right]$}%
\EasyButWeakLineBreak%
{$\left[\!\llap{\phantom{%
\begingroup \smaller\smaller\smaller\begin{tabular}{@{}c@{}}%
0\\0\\0
\end{tabular}\endgroup%
}}\right.$}%
\begingroup \smaller\smaller\smaller\begin{tabular}{@{}c@{}}%
-1\\-16\\-544
\end{tabular}\endgroup%
\HardButStrongLineBreak\kern3pt%
\begingroup \smaller\smaller\smaller\begin{tabular}{@{}c@{}}%
-1\\-18\\-558
\end{tabular}\endgroup%
\HardButStrongLineBreak\kern3pt%
\begingroup \smaller\smaller\smaller\begin{tabular}{@{}c@{}}%
-1\\-20\\-576
\end{tabular}\endgroup%
\HardButStrongLineBreak\kern3pt%
\begingroup \smaller\smaller\smaller\begin{tabular}{@{}c@{}}%
0\\0\\-1
\end{tabular}\endgroup%
\HardButStrongLineBreak\kern3pt%
\begingroup \smaller\smaller\smaller\begin{tabular}{@{}c@{}}%
1\\24\\576
\end{tabular}\endgroup%
\HardButStrongLineBreak\kern3pt%
\begingroup \smaller\smaller\smaller\begin{tabular}{@{}c@{}}%
0\\1\\6
\end{tabular}\endgroup%
{$\left.\llap{\phantom{%
\begingroup \smaller\smaller\smaller\begin{tabular}{@{}c@{}}%
0\\0\\0
\end{tabular}\endgroup%
}}\!\right]$}%

\medskip%
%
\leavevmode\llap{}%
$W_{176}$%
\qquad\llap{34} lattices, $\chi=24$%
\hfill%
$222|222\slashinfty\rtimes D_{2}$%
\nopagebreak\smallskip\hrule\nopagebreak\medskip%
%
%
\leavevmode%
${L_{176.1}}$%
{} : {$1\above{1pt}{1pt}{2}{0}8\above{1pt}{1pt}{1}{7}{\cdot}1\above{1pt}{1pt}{-2}{}7\above{1pt}{1pt}{1}{}$}\EasyButWeakLineBreak%
{${56}\above{1pt}{1pt}{}{2}{1}\above{1pt}{1pt}{}{2}{7}\above{1pt}{1pt}{r}{2}{8}\above{1pt}{1pt}{s}{2}{28}\above{1pt}{1pt}{*}{2}{4}\above{1pt}{1pt}{*}{2}{56}\above{1pt}{1pt}{1,0}{\infty b}$}%
\nopagebreak\par%
\nopagebreak\par\leavevmode%
{$\left[\!\llap{\phantom{%
\begingroup \smaller\smaller\smaller\begin{tabular}{@{}c@{}}%
0\\0\\0
\end{tabular}\endgroup%
}}\right.$}%
\begingroup \smaller\smaller\smaller\begin{tabular}{@{}c@{}}%
-740040\\12600\\12992
\end{tabular}\endgroup%
\kern3pt%
\begingroup \smaller\smaller\smaller\begin{tabular}{@{}c@{}}%
12600\\-211\\-225
\end{tabular}\endgroup%
\kern3pt%
\begingroup \smaller\smaller\smaller\begin{tabular}{@{}c@{}}%
12992\\-225\\-224
\end{tabular}\endgroup%
{$\left.\llap{\phantom{%
\begingroup \smaller\smaller\smaller\begin{tabular}{@{}c@{}}%
0\\0\\0
\end{tabular}\endgroup%
}}\!\right]$}%
\EasyButWeakLineBreak%
{$\left[\!\llap{\phantom{%
\begingroup \smaller\smaller\smaller\begin{tabular}{@{}c@{}}%
0\\0\\0
\end{tabular}\endgroup%
}}\right.$}%
\begingroup \smaller\smaller\smaller\begin{tabular}{@{}c@{}}%
211\\6328\\5880
\end{tabular}\endgroup%
\HardButStrongLineBreak\kern3pt%
\begingroup \smaller\smaller\smaller\begin{tabular}{@{}c@{}}%
49\\1469\\1366
\end{tabular}\endgroup%
\HardButStrongLineBreak\kern3pt%
\begingroup \smaller\smaller\smaller\begin{tabular}{@{}c@{}}%
181\\5425\\5047
\end{tabular}\endgroup%
\HardButStrongLineBreak\kern3pt%
\begingroup \smaller\smaller\smaller\begin{tabular}{@{}c@{}}%
41\\1228\\1144
\end{tabular}\endgroup%
\HardButStrongLineBreak\kern3pt%
\begingroup \smaller\smaller\smaller\begin{tabular}{@{}c@{}}%
23\\686\\644
\end{tabular}\endgroup%
\HardButStrongLineBreak\kern3pt%
\begingroup \smaller\smaller\smaller\begin{tabular}{@{}c@{}}%
-15\\-450\\-418
\end{tabular}\endgroup%
\HardButStrongLineBreak\kern3pt%
\begingroup \smaller\smaller\smaller\begin{tabular}{@{}c@{}}%
-15\\-448\\-420
\end{tabular}\endgroup%
{$\left.\llap{\phantom{%
\begingroup \smaller\smaller\smaller\begin{tabular}{@{}c@{}}%
0\\0\\0
\end{tabular}\endgroup%
}}\!\right]$}%
%
%
\hbox{}\par\smallskip%
%
%
\leavevmode%
${L_{176.2}}$%
{} : {$[1\above{1pt}{1pt}{-}{}2\above{1pt}{1pt}{1}{}]\above{1pt}{1pt}{}{4}16\above{1pt}{1pt}{-}{3}{\cdot}1\above{1pt}{1pt}{-2}{}7\above{1pt}{1pt}{1}{}$}\spacer%
\instructions{2,m}%
\EasyButWeakLineBreak%
{${14}\above{1pt}{1pt}{r}{2}{4}\above{1pt}{1pt}{s}{2}{112}\above{1pt}{1pt}{l}{2}{2}\above{1pt}{1pt}{r}{2}{28}\above{1pt}{1pt}{*}{2}{16}\above{1pt}{1pt}{*}{2}{56}\above{1pt}{1pt}{8,5}{\infty z}$}%
\nopagebreak\par%
\nopagebreak\par\leavevmode%
{$\left[\!\llap{\phantom{%
\begingroup \smaller\smaller\smaller\begin{tabular}{@{}c@{}}%
0\\0\\0
\end{tabular}\endgroup%
}}\right.$}%
\begingroup \smaller\smaller\smaller\begin{tabular}{@{}c@{}}%
-345296\\7840\\8624
\end{tabular}\endgroup%
\kern3pt%
\begingroup \smaller\smaller\smaller\begin{tabular}{@{}c@{}}%
7840\\-178\\-196
\end{tabular}\endgroup%
\kern3pt%
\begingroup \smaller\smaller\smaller\begin{tabular}{@{}c@{}}%
8624\\-196\\-211
\end{tabular}\endgroup%
{$\left.\llap{\phantom{%
\begingroup \smaller\smaller\smaller\begin{tabular}{@{}c@{}}%
0\\0\\0
\end{tabular}\endgroup%
}}\!\right]$}%
\EasyButWeakLineBreak%
{$\left[\!\llap{\phantom{%
\begingroup \smaller\smaller\smaller\begin{tabular}{@{}c@{}}%
0\\0\\0
\end{tabular}\endgroup%
}}\right.$}%
\begingroup \smaller\smaller\smaller\begin{tabular}{@{}c@{}}%
-3\\-133\\0
\end{tabular}\endgroup%
\HardButStrongLineBreak\kern3pt%
\begingroup \smaller\smaller\smaller\begin{tabular}{@{}c@{}}%
5\\204\\14
\end{tabular}\endgroup%
\HardButStrongLineBreak\kern3pt%
\begingroup \smaller\smaller\smaller\begin{tabular}{@{}c@{}}%
73\\3024\\168
\end{tabular}\endgroup%
\HardButStrongLineBreak\kern3pt%
\begingroup \smaller\smaller\smaller\begin{tabular}{@{}c@{}}%
5\\209\\10
\end{tabular}\endgroup%
\HardButStrongLineBreak\kern3pt%
\begingroup \smaller\smaller\smaller\begin{tabular}{@{}c@{}}%
23\\966\\42
\end{tabular}\endgroup%
\HardButStrongLineBreak\kern3pt%
\begingroup \smaller\smaller\smaller\begin{tabular}{@{}c@{}}%
1\\44\\0
\end{tabular}\endgroup%
\HardButStrongLineBreak\kern3pt%
\begingroup \smaller\smaller\smaller\begin{tabular}{@{}c@{}}%
-15\\-630\\-28
\end{tabular}\endgroup%
{$\left.\llap{\phantom{%
\begingroup \smaller\smaller\smaller\begin{tabular}{@{}c@{}}%
0\\0\\0
\end{tabular}\endgroup%
}}\!\right]$}%
%
%
\hbox{}\par\smallskip%
%
%
\leavevmode%
${L_{176.3}}$%
{} : {$[1\above{1pt}{1pt}{1}{}2\above{1pt}{1pt}{1}{}]\above{1pt}{1pt}{}{6}16\above{1pt}{1pt}{1}{1}{\cdot}1\above{1pt}{1pt}{-2}{}7\above{1pt}{1pt}{1}{}$}\spacer%
\instructions{m}%
\EasyButWeakLineBreak%
{${56}\above{1pt}{1pt}{*}{2}{4}\above{1pt}{1pt}{*}{2}{112}\above{1pt}{1pt}{s}{2}{8}\above{1pt}{1pt}{l}{2}{7}\above{1pt}{1pt}{}{2}{16}\above{1pt}{1pt}{}{2}{14}\above{1pt}{1pt}{8,1}{\infty}$}%
\nopagebreak\par%
\nopagebreak\par\leavevmode%
{$\left[\!\llap{\phantom{%
\begingroup \smaller\smaller\smaller\begin{tabular}{@{}c@{}}%
0\\0\\0
\end{tabular}\endgroup%
}}\right.$}%
\begingroup \smaller\smaller\smaller\begin{tabular}{@{}c@{}}%
-655984\\-185808\\11312
\end{tabular}\endgroup%
\kern3pt%
\begingroup \smaller\smaller\smaller\begin{tabular}{@{}c@{}}%
-185808\\-52630\\3204
\end{tabular}\endgroup%
\kern3pt%
\begingroup \smaller\smaller\smaller\begin{tabular}{@{}c@{}}%
11312\\3204\\-195
\end{tabular}\endgroup%
{$\left.\llap{\phantom{%
\begingroup \smaller\smaller\smaller\begin{tabular}{@{}c@{}}%
0\\0\\0
\end{tabular}\endgroup%
}}\!\right]$}%
\EasyButWeakLineBreak%
{$\left[\!\llap{\phantom{%
\begingroup \smaller\smaller\smaller\begin{tabular}{@{}c@{}}%
0\\0\\0
\end{tabular}\endgroup%
}}\right.$}%
\begingroup \smaller\smaller\smaller\begin{tabular}{@{}c@{}}%
3\\-14\\-56
\end{tabular}\endgroup%
\HardButStrongLineBreak\kern3pt%
\begingroup \smaller\smaller\smaller\begin{tabular}{@{}c@{}}%
7\\-28\\-54
\end{tabular}\endgroup%
\HardButStrongLineBreak\kern3pt%
\begingroup \smaller\smaller\smaller\begin{tabular}{@{}c@{}}%
-5\\28\\168
\end{tabular}\endgroup%
\HardButStrongLineBreak\kern3pt%
\begingroup \smaller\smaller\smaller\begin{tabular}{@{}c@{}}%
-15\\62\\148
\end{tabular}\endgroup%
\HardButStrongLineBreak\kern3pt%
\begingroup \smaller\smaller\smaller\begin{tabular}{@{}c@{}}%
-74\\301\\651
\end{tabular}\endgroup%
\HardButStrongLineBreak\kern3pt%
\begingroup \smaller\smaller\smaller\begin{tabular}{@{}c@{}}%
-83\\336\\704
\end{tabular}\endgroup%
\HardButStrongLineBreak\kern3pt%
\begingroup \smaller\smaller\smaller\begin{tabular}{@{}c@{}}%
-47\\189\\378
\end{tabular}\endgroup%
{$\left.\llap{\phantom{%
\begingroup \smaller\smaller\smaller\begin{tabular}{@{}c@{}}%
0\\0\\0
\end{tabular}\endgroup%
}}\!\right]$}%
%
%
\hbox{}\par\smallskip%
%
%
\leavevmode%
${L_{176.4}}$%
{} : {$[1\above{1pt}{1pt}{-}{}2\above{1pt}{1pt}{1}{}]\above{1pt}{1pt}{}{2}16\above{1pt}{1pt}{-}{5}{\cdot}1\above{1pt}{1pt}{-2}{}7\above{1pt}{1pt}{1}{}$}\EasyButWeakLineBreak%
{${56}\above{1pt}{1pt}{l}{2}{1}\above{1pt}{1pt}{r}{2}{112}\above{1pt}{1pt}{*}{2}{8}\above{1pt}{1pt}{*}{2}{28}\above{1pt}{1pt}{s}{2}{16}\above{1pt}{1pt}{l}{2}{14}\above{1pt}{1pt}{8,5}{\infty}$}%
\nopagebreak\par%
\nopagebreak\par\leavevmode%
{$\left[\!\llap{\phantom{%
\begingroup \smaller\smaller\smaller\begin{tabular}{@{}c@{}}%
0\\0\\0
\end{tabular}\endgroup%
}}\right.$}%
\begingroup \smaller\smaller\smaller\begin{tabular}{@{}c@{}}%
-560\\-112\\784
\end{tabular}\endgroup%
\kern3pt%
\begingroup \smaller\smaller\smaller\begin{tabular}{@{}c@{}}%
-112\\-22\\168
\end{tabular}\endgroup%
\kern3pt%
\begingroup \smaller\smaller\smaller\begin{tabular}{@{}c@{}}%
784\\168\\-783
\end{tabular}\endgroup%
{$\left.\llap{\phantom{%
\begingroup \smaller\smaller\smaller\begin{tabular}{@{}c@{}}%
0\\0\\0
\end{tabular}\endgroup%
}}\!\right]$}%
\EasyButWeakLineBreak%
{$\left[\!\llap{\phantom{%
\begingroup \smaller\smaller\smaller\begin{tabular}{@{}c@{}}%
0\\0\\0
\end{tabular}\endgroup%
}}\right.$}%
\begingroup \smaller\smaller\smaller\begin{tabular}{@{}c@{}}%
3\\-14\\0
\end{tabular}\endgroup%
\HardButStrongLineBreak\kern3pt%
\begingroup \smaller\smaller\smaller\begin{tabular}{@{}c@{}}%
7\\-28\\1
\end{tabular}\endgroup%
\HardButStrongLineBreak\kern3pt%
\begingroup \smaller\smaller\smaller\begin{tabular}{@{}c@{}}%
-5\\28\\0
\end{tabular}\endgroup%
\HardButStrongLineBreak\kern3pt%
\begingroup \smaller\smaller\smaller\begin{tabular}{@{}c@{}}%
-29\\118\\-4
\end{tabular}\endgroup%
\HardButStrongLineBreak\kern3pt%
\begingroup \smaller\smaller\smaller\begin{tabular}{@{}c@{}}%
-295\\1190\\-42
\end{tabular}\endgroup%
\HardButStrongLineBreak\kern3pt%
\begingroup \smaller\smaller\smaller\begin{tabular}{@{}c@{}}%
-167\\672\\-24
\end{tabular}\endgroup%
\HardButStrongLineBreak\kern3pt%
\begingroup \smaller\smaller\smaller\begin{tabular}{@{}c@{}}%
-96\\385\\-14
\end{tabular}\endgroup%
{$\left.\llap{\phantom{%
\begingroup \smaller\smaller\smaller\begin{tabular}{@{}c@{}}%
0\\0\\0
\end{tabular}\endgroup%
}}\!\right]$}%
%
%
\hbox{}\par\smallskip%
%
%
\leavevmode%
${L_{176.5}}$%
{} : {$1\above{1pt}{1pt}{1}{7}8\above{1pt}{1pt}{1}{1}64\above{1pt}{1pt}{1}{7}{\cdot}1\above{1pt}{1pt}{-2}{}7\above{1pt}{1pt}{1}{}$}\spacer%
\instructions{2}%
\EasyButWeakLineBreak%
{${56}\above{1pt}{1pt}{s}{2}{4}\above{1pt}{1pt}{b}{2}{448}\above{1pt}{1pt}{l}{2}{8}\above{1pt}{1pt}{r}{2}{28}\above{1pt}{1pt}{*}{2}{64}\above{1pt}{1pt}{s}{2}{224}\above{1pt}{1pt}{16,1}{\infty z}$}%
\nopagebreak\par%
shares genus with {$ {L_{176.6}}$}%
; isometric to own %
2-dual\nopagebreak\par%
\nopagebreak\par\leavevmode%
{$\left[\!\llap{\phantom{%
\begingroup \smaller\smaller\smaller\begin{tabular}{@{}c@{}}%
0\\0\\0
\end{tabular}\endgroup%
}}\right.$}%
\begingroup \smaller\smaller\smaller\begin{tabular}{@{}c@{}}%
-189504\\-9856\\448
\end{tabular}\endgroup%
\kern3pt%
\begingroup \smaller\smaller\smaller\begin{tabular}{@{}c@{}}%
-9856\\-504\\24
\end{tabular}\endgroup%
\kern3pt%
\begingroup \smaller\smaller\smaller\begin{tabular}{@{}c@{}}%
448\\24\\-1
\end{tabular}\endgroup%
{$\left.\llap{\phantom{%
\begingroup \smaller\smaller\smaller\begin{tabular}{@{}c@{}}%
0\\0\\0
\end{tabular}\endgroup%
}}\!\right]$}%
\EasyButWeakLineBreak%
{$\left[\!\llap{\phantom{%
\begingroup \smaller\smaller\smaller\begin{tabular}{@{}c@{}}%
0\\0\\0
\end{tabular}\endgroup%
}}\right.$}%
\begingroup \smaller\smaller\smaller\begin{tabular}{@{}c@{}}%
3\\-35\\476
\end{tabular}\endgroup%
\HardButStrongLineBreak\kern3pt%
\begingroup \smaller\smaller\smaller\begin{tabular}{@{}c@{}}%
1\\-11\\170
\end{tabular}\endgroup%
\HardButStrongLineBreak\kern3pt%
\begingroup \smaller\smaller\smaller\begin{tabular}{@{}c@{}}%
11\\-112\\2016
\end{tabular}\endgroup%
\HardButStrongLineBreak\kern3pt%
\begingroup \smaller\smaller\smaller\begin{tabular}{@{}c@{}}%
0\\1\\16
\end{tabular}\endgroup%
\HardButStrongLineBreak\kern3pt%
\begingroup \smaller\smaller\smaller\begin{tabular}{@{}c@{}}%
-1\\14\\-126
\end{tabular}\endgroup%
\HardButStrongLineBreak\kern3pt%
\begingroup \smaller\smaller\smaller\begin{tabular}{@{}c@{}}%
-1\\12\\-160
\end{tabular}\endgroup%
\HardButStrongLineBreak\kern3pt%
\begingroup \smaller\smaller\smaller\begin{tabular}{@{}c@{}}%
1\\-14\\112
\end{tabular}\endgroup%
{$\left.\llap{\phantom{%
\begingroup \smaller\smaller\smaller\begin{tabular}{@{}c@{}}%
0\\0\\0
\end{tabular}\endgroup%
}}\!\right]$}%
%
%
\hbox{}\par\smallskip%
%
%
\leavevmode%
${L_{176.6}}$%
{} : {$1\above{1pt}{1pt}{1}{7}8\above{1pt}{1pt}{1}{1}64\above{1pt}{1pt}{1}{7}{\cdot}1\above{1pt}{1pt}{-2}{}7\above{1pt}{1pt}{1}{}$}\EasyButWeakLineBreak%
{${56}\above{1pt}{1pt}{b}{2}{4}\above{1pt}{1pt}{l}{2}{448}\above{1pt}{1pt}{}{2}{8}\above{1pt}{1pt}{}{2}{7}\above{1pt}{1pt}{r}{2}{64}\above{1pt}{1pt}{*}{2}{224}\above{1pt}{1pt}{16,9}{\infty z}$}%
\nopagebreak\par%
shares genus with {$ {L_{176.5}}$}%
; isometric to own %
2-dual\nopagebreak\par%
\nopagebreak\par\leavevmode%
{$\left[\!\llap{\phantom{%
\begingroup \smaller\smaller\smaller\begin{tabular}{@{}c@{}}%
0\\0\\0
\end{tabular}\endgroup%
}}\right.$}%
\begingroup \smaller\smaller\smaller\begin{tabular}{@{}c@{}}%
-3777088\\165312\\-54656
\end{tabular}\endgroup%
\kern3pt%
\begingroup \smaller\smaller\smaller\begin{tabular}{@{}c@{}}%
165312\\-7224\\2384
\end{tabular}\endgroup%
\kern3pt%
\begingroup \smaller\smaller\smaller\begin{tabular}{@{}c@{}}%
-54656\\2384\\-785
\end{tabular}\endgroup%
{$\left.\llap{\phantom{%
\begingroup \smaller\smaller\smaller\begin{tabular}{@{}c@{}}%
0\\0\\0
\end{tabular}\endgroup%
}}\!\right]$}%
\EasyButWeakLineBreak%
{$\left[\!\llap{\phantom{%
\begingroup \smaller\smaller\smaller\begin{tabular}{@{}c@{}}%
0\\0\\0
\end{tabular}\endgroup%
}}\right.$}%
\begingroup \smaller\smaller\smaller\begin{tabular}{@{}c@{}}%
-3\\-133\\-196
\end{tabular}\endgroup%
\HardButStrongLineBreak\kern3pt%
\begingroup \smaller\smaller\smaller\begin{tabular}{@{}c@{}}%
6\\249\\338
\end{tabular}\endgroup%
\HardButStrongLineBreak\kern3pt%
\begingroup \smaller\smaller\smaller\begin{tabular}{@{}c@{}}%
157\\6552\\8960
\end{tabular}\endgroup%
\HardButStrongLineBreak\kern3pt%
\begingroup \smaller\smaller\smaller\begin{tabular}{@{}c@{}}%
10\\419\\576
\end{tabular}\endgroup%
\HardButStrongLineBreak\kern3pt%
\begingroup \smaller\smaller\smaller\begin{tabular}{@{}c@{}}%
11\\462\\637
\end{tabular}\endgroup%
\HardButStrongLineBreak\kern3pt%
\begingroup \smaller\smaller\smaller\begin{tabular}{@{}c@{}}%
1\\44\\64
\end{tabular}\endgroup%
\HardButStrongLineBreak\kern3pt%
\begingroup \smaller\smaller\smaller\begin{tabular}{@{}c@{}}%
-29\\-1218\\-1680
\end{tabular}\endgroup%
{$\left.\llap{\phantom{%
\begingroup \smaller\smaller\smaller\begin{tabular}{@{}c@{}}%
0\\0\\0
\end{tabular}\endgroup%
}}\!\right]$}%

\medskip%
%
\leavevmode\llap{}%
$W_{177}$%
\qquad\llap{44} lattices, $\chi=24$%
\hfill%
$2|22|22|22|2\rtimes D_{4}$%
\nopagebreak\smallskip\hrule\nopagebreak\medskip%
%
%
\leavevmode%
${L_{177.1}}$%
{} : {$1\above{1pt}{1pt}{-2}{4}16\above{1pt}{1pt}{-}{3}{\cdot}1\above{1pt}{1pt}{2}{}3\above{1pt}{1pt}{1}{}{\cdot}1\above{1pt}{1pt}{-2}{}5\above{1pt}{1pt}{1}{}$}\spacer%
\instructions{2}%
\EasyButWeakLineBreak%
{${48}\above{1pt}{1pt}{*}{2}{4}\above{1pt}{1pt}{*}{2}{12}\above{1pt}{1pt}{l}{2}{5}\above{1pt}{1pt}{}{2}{48}\above{1pt}{1pt}{}{2}{1}\above{1pt}{1pt}{}{2}{3}\above{1pt}{1pt}{r}{2}{20}\above{1pt}{1pt}{*}{2}$}%
\nopagebreak\par%
\nopagebreak\par\leavevmode%
{$\left[\!\llap{\phantom{%
\begingroup \smaller\smaller\smaller\begin{tabular}{@{}c@{}}%
0\\0\\0
\end{tabular}\endgroup%
}}\right.$}%
\begingroup \smaller\smaller\smaller\begin{tabular}{@{}c@{}}%
-96720\\720\\240
\end{tabular}\endgroup%
\kern3pt%
\begingroup \smaller\smaller\smaller\begin{tabular}{@{}c@{}}%
720\\-4\\-5
\end{tabular}\endgroup%
\kern3pt%
\begingroup \smaller\smaller\smaller\begin{tabular}{@{}c@{}}%
240\\-5\\7
\end{tabular}\endgroup%
{$\left.\llap{\phantom{%
\begingroup \smaller\smaller\smaller\begin{tabular}{@{}c@{}}%
0\\0\\0
\end{tabular}\endgroup%
}}\!\right]$}%
\EasyButWeakLineBreak%
{$\left[\!\llap{\phantom{%
\begingroup \smaller\smaller\smaller\begin{tabular}{@{}c@{}}%
0\\0\\0
\end{tabular}\endgroup%
}}\right.$}%
\begingroup \smaller\smaller\smaller\begin{tabular}{@{}c@{}}%
-1\\-120\\-48
\end{tabular}\endgroup%
\HardButStrongLineBreak\kern3pt%
\begingroup \smaller\smaller\smaller\begin{tabular}{@{}c@{}}%
-1\\-118\\-50
\end{tabular}\endgroup%
\HardButStrongLineBreak\kern3pt%
\begingroup \smaller\smaller\smaller\begin{tabular}{@{}c@{}}%
-1\\-120\\-54
\end{tabular}\endgroup%
\HardButStrongLineBreak\kern3pt%
\begingroup \smaller\smaller\smaller\begin{tabular}{@{}c@{}}%
1\\115\\45
\end{tabular}\endgroup%
\HardButStrongLineBreak\kern3pt%
\begingroup \smaller\smaller\smaller\begin{tabular}{@{}c@{}}%
7\\816\\336
\end{tabular}\endgroup%
\HardButStrongLineBreak\kern3pt%
\begingroup \smaller\smaller\smaller\begin{tabular}{@{}c@{}}%
1\\117\\49
\end{tabular}\endgroup%
\HardButStrongLineBreak\kern3pt%
\begingroup \smaller\smaller\smaller\begin{tabular}{@{}c@{}}%
2\\234\\99
\end{tabular}\endgroup%
\HardButStrongLineBreak\kern3pt%
\begingroup \smaller\smaller\smaller\begin{tabular}{@{}c@{}}%
3\\350\\150
\end{tabular}\endgroup%
{$\left.\llap{\phantom{%
\begingroup \smaller\smaller\smaller\begin{tabular}{@{}c@{}}%
0\\0\\0
\end{tabular}\endgroup%
}}\!\right]$}%
%
%
\hbox{}\par\smallskip%
%
%
\leavevmode%
${L_{177.2}}$%
{} : {$1\above{1pt}{1pt}{2}{0}16\above{1pt}{1pt}{1}{7}{\cdot}1\above{1pt}{1pt}{2}{}3\above{1pt}{1pt}{1}{}{\cdot}1\above{1pt}{1pt}{-2}{}5\above{1pt}{1pt}{1}{}$}\spacer%
\instructions{m}%
\EasyButWeakLineBreak%
{${48}\above{1pt}{1pt}{l}{2}{1}\above{1pt}{1pt}{r}{2}{12}\above{1pt}{1pt}{*}{2}{20}\above{1pt}{1pt}{s}{2}{48}\above{1pt}{1pt}{s}{2}{4}\above{1pt}{1pt}{l}{2}{3}\above{1pt}{1pt}{}{2}{5}\above{1pt}{1pt}{r}{2}$}%
\nopagebreak\par%
\nopagebreak\par\leavevmode%
{$\left[\!\llap{\phantom{%
\begingroup \smaller\smaller\smaller\begin{tabular}{@{}c@{}}%
0\\0\\0
\end{tabular}\endgroup%
}}\right.$}%
\begingroup \smaller\smaller\smaller\begin{tabular}{@{}c@{}}%
205680\\1680\\-1200
\end{tabular}\endgroup%
\kern3pt%
\begingroup \smaller\smaller\smaller\begin{tabular}{@{}c@{}}%
1680\\-19\\-10
\end{tabular}\endgroup%
\kern3pt%
\begingroup \smaller\smaller\smaller\begin{tabular}{@{}c@{}}%
-1200\\-10\\7
\end{tabular}\endgroup%
{$\left.\llap{\phantom{%
\begingroup \smaller\smaller\smaller\begin{tabular}{@{}c@{}}%
0\\0\\0
\end{tabular}\endgroup%
}}\!\right]$}%
\EasyButWeakLineBreak%
{$\left[\!\llap{\phantom{%
\begingroup \smaller\smaller\smaller\begin{tabular}{@{}c@{}}%
0\\0\\0
\end{tabular}\endgroup%
}}\right.$}%
\begingroup \smaller\smaller\smaller\begin{tabular}{@{}c@{}}%
-1\\0\\-168
\end{tabular}\endgroup%
\HardButStrongLineBreak\kern3pt%
\begingroup \smaller\smaller\smaller\begin{tabular}{@{}c@{}}%
-1\\1\\-170
\end{tabular}\endgroup%
\HardButStrongLineBreak\kern3pt%
\begingroup \smaller\smaller\smaller\begin{tabular}{@{}c@{}}%
-1\\0\\-174
\end{tabular}\endgroup%
\HardButStrongLineBreak\kern3pt%
\begingroup \smaller\smaller\smaller\begin{tabular}{@{}c@{}}%
7\\-10\\1180
\end{tabular}\endgroup%
\HardButStrongLineBreak\kern3pt%
\begingroup \smaller\smaller\smaller\begin{tabular}{@{}c@{}}%
19\\-24\\3216
\end{tabular}\endgroup%
\HardButStrongLineBreak\kern3pt%
\begingroup \smaller\smaller\smaller\begin{tabular}{@{}c@{}}%
5\\-6\\848
\end{tabular}\endgroup%
\HardButStrongLineBreak\kern3pt%
\begingroup \smaller\smaller\smaller\begin{tabular}{@{}c@{}}%
5\\-6\\849
\end{tabular}\endgroup%
\HardButStrongLineBreak\kern3pt%
\begingroup \smaller\smaller\smaller\begin{tabular}{@{}c@{}}%
4\\-5\\680
\end{tabular}\endgroup%
{$\left.\llap{\phantom{%
\begingroup \smaller\smaller\smaller\begin{tabular}{@{}c@{}}%
0\\0\\0
\end{tabular}\endgroup%
}}\!\right]$}%
%
%
\hbox{}\par\smallskip%
%
%
\leavevmode%
${L_{177.3}}$%
{} : {$1\above{1pt}{1pt}{1}{7}4\above{1pt}{1pt}{1}{1}16\above{1pt}{1pt}{1}{7}{\cdot}1\above{1pt}{1pt}{2}{}3\above{1pt}{1pt}{1}{}{\cdot}1\above{1pt}{1pt}{-2}{}5\above{1pt}{1pt}{1}{}$}\EasyButWeakLineBreak%
{${48}\above{1pt}{1pt}{l}{2}{4}\above{1pt}{1pt}{r}{2}{12}\above{1pt}{1pt}{l}{2}{20}\above{1pt}{1pt}{}{2}{48}\above{1pt}{1pt}{}{2}{4}\above{1pt}{1pt}{}{2}{3}\above{1pt}{1pt}{}{2}{20}\above{1pt}{1pt}{r}{2}$}%
\nopagebreak\par%
\nopagebreak\par\leavevmode%
{$\left[\!\llap{\phantom{%
\begingroup \smaller\smaller\smaller\begin{tabular}{@{}c@{}}%
0\\0\\0
\end{tabular}\endgroup%
}}\right.$}%
\begingroup \smaller\smaller\smaller\begin{tabular}{@{}c@{}}%
205680\\4080\\-1200
\end{tabular}\endgroup%
\kern3pt%
\begingroup \smaller\smaller\smaller\begin{tabular}{@{}c@{}}%
4080\\52\\-24
\end{tabular}\endgroup%
\kern3pt%
\begingroup \smaller\smaller\smaller\begin{tabular}{@{}c@{}}%
-1200\\-24\\7
\end{tabular}\endgroup%
{$\left.\llap{\phantom{%
\begingroup \smaller\smaller\smaller\begin{tabular}{@{}c@{}}%
0\\0\\0
\end{tabular}\endgroup%
}}\!\right]$}%
\EasyButWeakLineBreak%
{$\left[\!\llap{\phantom{%
\begingroup \smaller\smaller\smaller\begin{tabular}{@{}c@{}}%
0\\0\\0
\end{tabular}\endgroup%
}}\right.$}%
\begingroup \smaller\smaller\smaller\begin{tabular}{@{}c@{}}%
-1\\0\\-168
\end{tabular}\endgroup%
\HardButStrongLineBreak\kern3pt%
\begingroup \smaller\smaller\smaller\begin{tabular}{@{}c@{}}%
-1\\1\\-168
\end{tabular}\endgroup%
\HardButStrongLineBreak\kern3pt%
\begingroup \smaller\smaller\smaller\begin{tabular}{@{}c@{}}%
-1\\0\\-174
\end{tabular}\endgroup%
\HardButStrongLineBreak\kern3pt%
\begingroup \smaller\smaller\smaller\begin{tabular}{@{}c@{}}%
2\\-5\\320
\end{tabular}\endgroup%
\HardButStrongLineBreak\kern3pt%
\begingroup \smaller\smaller\smaller\begin{tabular}{@{}c@{}}%
7\\-12\\1152
\end{tabular}\endgroup%
\HardButStrongLineBreak\kern3pt%
\begingroup \smaller\smaller\smaller\begin{tabular}{@{}c@{}}%
2\\-3\\332
\end{tabular}\endgroup%
\HardButStrongLineBreak\kern3pt%
\begingroup \smaller\smaller\smaller\begin{tabular}{@{}c@{}}%
2\\-3\\333
\end{tabular}\endgroup%
\HardButStrongLineBreak\kern3pt%
\begingroup \smaller\smaller\smaller\begin{tabular}{@{}c@{}}%
3\\-5\\500
\end{tabular}\endgroup%
{$\left.\llap{\phantom{%
\begingroup \smaller\smaller\smaller\begin{tabular}{@{}c@{}}%
0\\0\\0
\end{tabular}\endgroup%
}}\!\right]$}%
%
%
\hbox{}\par\smallskip%
%
%
\leavevmode%
${L_{177.4}}$%
{} : {$1\above{1pt}{1pt}{1}{1}4\above{1pt}{1pt}{1}{7}16\above{1pt}{1pt}{1}{7}{\cdot}1\above{1pt}{1pt}{2}{}3\above{1pt}{1pt}{1}{}{\cdot}1\above{1pt}{1pt}{-2}{}5\above{1pt}{1pt}{1}{}$}\EasyButWeakLineBreak%
{${12}\above{1pt}{1pt}{}{2}{1}\above{1pt}{1pt}{}{2}{48}\above{1pt}{1pt}{}{2}{5}\above{1pt}{1pt}{}{2}{12}\above{1pt}{1pt}{r}{2}{4}\above{1pt}{1pt}{s}{2}{48}\above{1pt}{1pt}{s}{2}{20}\above{1pt}{1pt}{l}{2}$}%
\nopagebreak\par%
\nopagebreak\par\leavevmode%
{$\left[\!\llap{\phantom{%
\begingroup \smaller\smaller\smaller\begin{tabular}{@{}c@{}}%
0\\0\\0
\end{tabular}\endgroup%
}}\right.$}%
\begingroup \smaller\smaller\smaller\begin{tabular}{@{}c@{}}%
-30480\\1440\\0
\end{tabular}\endgroup%
\kern3pt%
\begingroup \smaller\smaller\smaller\begin{tabular}{@{}c@{}}%
1440\\-4\\-8
\end{tabular}\endgroup%
\kern3pt%
\begingroup \smaller\smaller\smaller\begin{tabular}{@{}c@{}}%
0\\-8\\1
\end{tabular}\endgroup%
{$\left.\llap{\phantom{%
\begingroup \smaller\smaller\smaller\begin{tabular}{@{}c@{}}%
0\\0\\0
\end{tabular}\endgroup%
}}\!\right]$}%
\EasyButWeakLineBreak%
{$\left[\!\llap{\phantom{%
\begingroup \smaller\smaller\smaller\begin{tabular}{@{}c@{}}%
0\\0\\0
\end{tabular}\endgroup%
}}\right.$}%
\begingroup \smaller\smaller\smaller\begin{tabular}{@{}c@{}}%
1\\21\\168
\end{tabular}\endgroup%
\HardButStrongLineBreak\kern3pt%
\begingroup \smaller\smaller\smaller\begin{tabular}{@{}c@{}}%
0\\0\\-1
\end{tabular}\endgroup%
\HardButStrongLineBreak\kern3pt%
\begingroup \smaller\smaller\smaller\begin{tabular}{@{}c@{}}%
-5\\-108\\-864
\end{tabular}\endgroup%
\HardButStrongLineBreak\kern3pt%
\begingroup \smaller\smaller\smaller\begin{tabular}{@{}c@{}}%
-3\\-65\\-515
\end{tabular}\endgroup%
\HardButStrongLineBreak\kern3pt%
\begingroup \smaller\smaller\smaller\begin{tabular}{@{}c@{}}%
-4\\-87\\-684
\end{tabular}\endgroup%
\HardButStrongLineBreak\kern3pt%
\begingroup \smaller\smaller\smaller\begin{tabular}{@{}c@{}}%
-1\\-22\\-170
\end{tabular}\endgroup%
\HardButStrongLineBreak\kern3pt%
\begingroup \smaller\smaller\smaller\begin{tabular}{@{}c@{}}%
-1\\-24\\-168
\end{tabular}\endgroup%
\HardButStrongLineBreak\kern3pt%
\begingroup \smaller\smaller\smaller\begin{tabular}{@{}c@{}}%
1\\20\\170
\end{tabular}\endgroup%
{$\left.\llap{\phantom{%
\begingroup \smaller\smaller\smaller\begin{tabular}{@{}c@{}}%
0\\0\\0
\end{tabular}\endgroup%
}}\!\right]$}%

\medskip%
%
\leavevmode\llap{}%
$W_{178}$%
\qquad\llap{12} lattices, $\chi=18$%
\hfill%
$422422\rtimes C_{2}$%
\nopagebreak\smallskip\hrule\nopagebreak\medskip%
%
%
\leavevmode%
${L_{178.1}}$%
{} : {$1\above{1pt}{1pt}{-2}{{\rm II}}4\above{1pt}{1pt}{-}{3}{\cdot}1\above{1pt}{1pt}{2}{}9\above{1pt}{1pt}{-}{}{\cdot}1\above{1pt}{1pt}{2}{}7\above{1pt}{1pt}{-}{}$}\spacer%
\instructions{2}%
\EasyButWeakLineBreak%
{${2}\above{1pt}{1pt}{*}{4}{4}\above{1pt}{1pt}{b}{2}{18}\above{1pt}{1pt}{b}{2}$}\relax$\,(\times2)$%
\nopagebreak\par%
\nopagebreak\par\leavevmode%
{$\left[\!\llap{\phantom{%
\begingroup \smaller\smaller\smaller\begin{tabular}{@{}c@{}}%
0\\0\\0
\end{tabular}\endgroup%
}}\right.$}%
\begingroup \smaller\smaller\smaller\begin{tabular}{@{}c@{}}%
-506772\\2268\\4032
\end{tabular}\endgroup%
\kern3pt%
\begingroup \smaller\smaller\smaller\begin{tabular}{@{}c@{}}%
2268\\-10\\-19
\end{tabular}\endgroup%
\kern3pt%
\begingroup \smaller\smaller\smaller\begin{tabular}{@{}c@{}}%
4032\\-19\\-26
\end{tabular}\endgroup%
{$\left.\llap{\phantom{%
\begingroup \smaller\smaller\smaller\begin{tabular}{@{}c@{}}%
0\\0\\0
\end{tabular}\endgroup%
}}\!\right]$}%
\hfil\penalty500%
{$\left[\!\llap{\phantom{%
\begingroup \smaller\smaller\smaller\begin{tabular}{@{}c@{}}%
0\\0\\0
\end{tabular}\endgroup%
}}\right.$}%
\begingroup \smaller\smaller\smaller\begin{tabular}{@{}c@{}}%
5039\\876960\\138600
\end{tabular}\endgroup%
\kern3pt%
\begingroup \smaller\smaller\smaller\begin{tabular}{@{}c@{}}%
-22\\-3829\\-605
\end{tabular}\endgroup%
\kern3pt%
\begingroup \smaller\smaller\smaller\begin{tabular}{@{}c@{}}%
-44\\-7656\\-1211
\end{tabular}\endgroup%
{$\left.\llap{\phantom{%
\begingroup \smaller\smaller\smaller\begin{tabular}{@{}c@{}}%
0\\0\\0
\end{tabular}\endgroup%
}}\!\right]$}%
\EasyButWeakLineBreak%
{$\left[\!\llap{\phantom{%
\begingroup \smaller\smaller\smaller\begin{tabular}{@{}c@{}}%
0\\0\\0
\end{tabular}\endgroup%
}}\right.$}%
\begingroup \smaller\smaller\smaller\begin{tabular}{@{}c@{}}%
1\\173\\28
\end{tabular}\endgroup%
\HardButStrongLineBreak\kern3pt%
\begingroup \smaller\smaller\smaller\begin{tabular}{@{}c@{}}%
21\\3654\\578
\end{tabular}\endgroup%
\HardButStrongLineBreak\kern3pt%
\begingroup \smaller\smaller\smaller\begin{tabular}{@{}c@{}}%
20\\3483\\549
\end{tabular}\endgroup%
{$\left.\llap{\phantom{%
\begingroup \smaller\smaller\smaller\begin{tabular}{@{}c@{}}%
0\\0\\0
\end{tabular}\endgroup%
}}\!\right]$}%

\medskip%
%
\leavevmode\llap{}%
$W_{179}$%
\qquad\llap{44} lattices, $\chi=72$%
\hfill%
$222\infty222222\infty222\rtimes C_{2}$%
\nopagebreak\smallskip\hrule\nopagebreak\medskip%
%
%
\leavevmode%
${L_{179.1}}$%
{} : {$1\above{1pt}{1pt}{2}{{\rm II}}4\above{1pt}{1pt}{1}{7}{\cdot}1\above{1pt}{1pt}{2}{}9\above{1pt}{1pt}{-}{}{\cdot}1\above{1pt}{1pt}{-2}{}7\above{1pt}{1pt}{1}{}$}\spacer%
\instructions{2}%
\EasyButWeakLineBreak%
{${2}\above{1pt}{1pt}{s}{2}{126}\above{1pt}{1pt}{b}{2}{4}\above{1pt}{1pt}{*}{2}{28}\above{1pt}{1pt}{3,2}{\infty b}{28}\above{1pt}{1pt}{r}{2}{18}\above{1pt}{1pt}{s}{2}{14}\above{1pt}{1pt}{s}{2}$}\relax$\,(\times2)$%
\nopagebreak\par%
\nopagebreak\par\leavevmode%
{$\left[\!\llap{\phantom{%
\begingroup \smaller\smaller\smaller
\endgroup%
}}\!\right]$}%
%
%
\hbox{}\par\smallskip%
%
%
\leavevmode%
${L_{179.2}}$%
{} : {$1\above{1pt}{1pt}{-2}{6}8\above{1pt}{1pt}{-}{5}{\cdot}1\above{1pt}{1pt}{2}{}9\above{1pt}{1pt}{1}{}{\cdot}1\above{1pt}{1pt}{-2}{}7\above{1pt}{1pt}{1}{}$}\spacer%
\instructions{2}%
\EasyButWeakLineBreak%
{${1}\above{1pt}{1pt}{r}{2}{252}\above{1pt}{1pt}{*}{2}{8}\above{1pt}{1pt}{b}{2}{14}\above{1pt}{1pt}{12,5}{\infty a}{56}\above{1pt}{1pt}{l}{2}{9}\above{1pt}{1pt}{r}{2}{28}\above{1pt}{1pt}{l}{2}$}\relax$\,(\times2)$%
\nopagebreak\par%
\nopagebreak\par\leavevmode%
{$\left[\!\llap{\phantom{%
\begingroup \smaller\smaller\smaller
\endgroup%
}}\!\right]$}%
%
%
\hbox{}\par\smallskip%
%
%
\leavevmode%
${L_{179.3}}$%
{} : {$1\above{1pt}{1pt}{2}{6}8\above{1pt}{1pt}{1}{1}{\cdot}1\above{1pt}{1pt}{2}{}9\above{1pt}{1pt}{1}{}{\cdot}1\above{1pt}{1pt}{-2}{}7\above{1pt}{1pt}{1}{}$}\spacer%
\instructions{m}%
\EasyButWeakLineBreak%
{${4}\above{1pt}{1pt}{l}{2}{63}\above{1pt}{1pt}{}{2}{8}\above{1pt}{1pt}{r}{2}{14}\above{1pt}{1pt}{12,5}{\infty b}{56}\above{1pt}{1pt}{s}{2}{36}\above{1pt}{1pt}{l}{2}{7}\above{1pt}{1pt}{r}{2}$}\relax$\,(\times2)$%
\nopagebreak\par%
\nopagebreak\par\leavevmode%
{$\left[\!\llap{\phantom{%
\begingroup \smaller\smaller\smaller
\endgroup%
}}\!\right]$}%

\medskip%
%
\leavevmode\llap{}%
$W_{180}$%
\qquad\llap{2} lattices, $\chi=24$%
\hfill%
$\slashtwo\infty4|4\infty\rtimes D_{2}$%
\nopagebreak\smallskip\hrule\nopagebreak\medskip%
%
%
\leavevmode%
${L_{180.1}}$%
{} : {$1\above{1pt}{1pt}{2}{2}64\above{1pt}{1pt}{1}{7}$}\EasyButWeakLineBreak%
{${1}\above{1pt}{1pt}{r}{2}{4}\above{1pt}{1pt}{16,15}{\infty z}{1}\above{1pt}{1pt}{}{4}{2}\above{1pt}{1pt}{*}{4}{4}\above{1pt}{1pt}{16,7}{\infty z}$}%
\nopagebreak\par%
\nopagebreak\par\leavevmode%
{$\left[\!\llap{\phantom{%
\begingroup \smaller\smaller\smaller\begin{tabular}{@{}c@{}}%
0\\0\\0
\end{tabular}\endgroup%
}}\right.$}%
\begingroup \smaller\smaller\smaller\begin{tabular}{@{}c@{}}%
-16960\\-896\\320
\end{tabular}\endgroup%
\kern3pt%
\begingroup \smaller\smaller\smaller\begin{tabular}{@{}c@{}}%
-896\\-47\\17
\end{tabular}\endgroup%
\kern3pt%
\begingroup \smaller\smaller\smaller\begin{tabular}{@{}c@{}}%
320\\17\\-6
\end{tabular}\endgroup%
{$\left.\llap{\phantom{%
\begingroup \smaller\smaller\smaller\begin{tabular}{@{}c@{}}%
0\\0\\0
\end{tabular}\endgroup%
}}\!\right]$}%
\EasyButWeakLineBreak%
{$\left[\!\llap{\phantom{%
\begingroup \smaller\smaller\smaller\begin{tabular}{@{}c@{}}%
0\\0\\0
\end{tabular}\endgroup%
}}\right.$}%
\begingroup \smaller\smaller\smaller\begin{tabular}{@{}c@{}}%
0\\-1\\-3
\end{tabular}\endgroup%
\HardButStrongLineBreak\kern3pt%
\begingroup \smaller\smaller\smaller\begin{tabular}{@{}c@{}}%
1\\-10\\24
\end{tabular}\endgroup%
\HardButStrongLineBreak\kern3pt%
\begingroup \smaller\smaller\smaller\begin{tabular}{@{}c@{}}%
1\\-7\\32
\end{tabular}\endgroup%
\HardButStrongLineBreak\kern3pt%
\begingroup \smaller\smaller\smaller\begin{tabular}{@{}c@{}}%
0\\2\\5
\end{tabular}\endgroup%
\HardButStrongLineBreak\kern3pt%
\begingroup \smaller\smaller\smaller\begin{tabular}{@{}c@{}}%
-1\\10\\-26
\end{tabular}\endgroup%
{$\left.\llap{\phantom{%
\begingroup \smaller\smaller\smaller\begin{tabular}{@{}c@{}}%
0\\0\\0
\end{tabular}\endgroup%
}}\!\right]$}%
%
%
%
%
%
%

\medskip%
%
\leavevmode\llap{}%
$W_{181}$%
\qquad\llap{14} lattices, $\chi=12$%
\hfill%
$2|22\slashinfty2\rtimes D_{2}$%
\nopagebreak\smallskip\hrule\nopagebreak\medskip%
%
%
\leavevmode%
${L_{181.1}}$%
{} : {$[1\above{1pt}{1pt}{1}{}2\above{1pt}{1pt}{-}{}]\above{1pt}{1pt}{}{4}32\above{1pt}{1pt}{-}{5}$}\EasyButWeakLineBreak%
{${32}\above{1pt}{1pt}{s}{2}{8}\above{1pt}{1pt}{*}{2}{32}\above{1pt}{1pt}{l}{2}{1}\above{1pt}{1pt}{8,5}{\infty}{4}\above{1pt}{1pt}{*}{2}$}%
\nopagebreak\par%
\nopagebreak\par\leavevmode%
{$\left[\!\llap{\phantom{%
\begingroup \smaller\smaller\smaller\begin{tabular}{@{}c@{}}%
0\\0\\0
\end{tabular}\endgroup%
}}\right.$}%
\begingroup \smaller\smaller\smaller\begin{tabular}{@{}c@{}}%
-13664\\-416\\928
\end{tabular}\endgroup%
\kern3pt%
\begingroup \smaller\smaller\smaller\begin{tabular}{@{}c@{}}%
-416\\-10\\28
\end{tabular}\endgroup%
\kern3pt%
\begingroup \smaller\smaller\smaller\begin{tabular}{@{}c@{}}%
928\\28\\-63
\end{tabular}\endgroup%
{$\left.\llap{\phantom{%
\begingroup \smaller\smaller\smaller\begin{tabular}{@{}c@{}}%
0\\0\\0
\end{tabular}\endgroup%
}}\!\right]$}%
\EasyButWeakLineBreak%
{$\left[\!\llap{\phantom{%
\begingroup \smaller\smaller\smaller\begin{tabular}{@{}c@{}}%
0\\0\\0
\end{tabular}\endgroup%
}}\right.$}%
\begingroup \smaller\smaller\smaller\begin{tabular}{@{}c@{}}%
-3\\-8\\-48
\end{tabular}\endgroup%
\HardButStrongLineBreak\kern3pt%
\begingroup \smaller\smaller\smaller\begin{tabular}{@{}c@{}}%
5\\6\\76
\end{tabular}\endgroup%
\HardButStrongLineBreak\kern3pt%
\begingroup \smaller\smaller\smaller\begin{tabular}{@{}c@{}}%
25\\40\\384
\end{tabular}\endgroup%
\HardButStrongLineBreak\kern3pt%
\begingroup \smaller\smaller\smaller\begin{tabular}{@{}c@{}}%
2\\4\\31
\end{tabular}\endgroup%
\HardButStrongLineBreak\kern3pt%
\begingroup \smaller\smaller\smaller\begin{tabular}{@{}c@{}}%
-3\\-4\\-46
\end{tabular}\endgroup%
{$\left.\llap{\phantom{%
\begingroup \smaller\smaller\smaller\begin{tabular}{@{}c@{}}%
0\\0\\0
\end{tabular}\endgroup%
}}\!\right]$}%
%
%
\hbox{}\par\smallskip%
%
%
\leavevmode%
${L_{181.2}}$%
{} : {$[1\above{1pt}{1pt}{1}{}2\above{1pt}{1pt}{1}{}]\above{1pt}{1pt}{}{0}64\above{1pt}{1pt}{1}{1}$}\spacer%
\instructions{m}%
\EasyButWeakLineBreak%
{${64}\above{1pt}{1pt}{*}{2}{4}\above{1pt}{1pt}{s}{2}{64}\above{1pt}{1pt}{l}{2}{2}\above{1pt}{1pt}{16,9}{\infty}{8}\above{1pt}{1pt}{s}{2}$}%
\nopagebreak\par%
shares genus with {$ {L_{181.3}}$}%
\nopagebreak\par%
\nopagebreak\par\leavevmode%
{$\left[\!\llap{\phantom{%
\begingroup \smaller\smaller\smaller\begin{tabular}{@{}c@{}}%
0\\0\\0
\end{tabular}\endgroup%
}}\right.$}%
\begingroup \smaller\smaller\smaller\begin{tabular}{@{}c@{}}%
-59840\\-2752\\-3008
\end{tabular}\endgroup%
\kern3pt%
\begingroup \smaller\smaller\smaller\begin{tabular}{@{}c@{}}%
-2752\\-126\\-138
\end{tabular}\endgroup%
\kern3pt%
\begingroup \smaller\smaller\smaller\begin{tabular}{@{}c@{}}%
-3008\\-138\\-151
\end{tabular}\endgroup%
{$\left.\llap{\phantom{%
\begingroup \smaller\smaller\smaller\begin{tabular}{@{}c@{}}%
0\\0\\0
\end{tabular}\endgroup%
}}\!\right]$}%
\EasyButWeakLineBreak%
{$\left[\!\llap{\phantom{%
\begingroup \smaller\smaller\smaller\begin{tabular}{@{}c@{}}%
0\\0\\0
\end{tabular}\endgroup%
}}\right.$}%
\begingroup \smaller\smaller\smaller\begin{tabular}{@{}c@{}}%
1\\48\\-64
\end{tabular}\endgroup%
\HardButStrongLineBreak\kern3pt%
\begingroup \smaller\smaller\smaller\begin{tabular}{@{}c@{}}%
1\\24\\-42
\end{tabular}\endgroup%
\HardButStrongLineBreak\kern3pt%
\begingroup \smaller\smaller\smaller\begin{tabular}{@{}c@{}}%
5\\64\\-160
\end{tabular}\endgroup%
\HardButStrongLineBreak\kern3pt%
\begingroup \smaller\smaller\smaller\begin{tabular}{@{}c@{}}%
0\\-9\\8
\end{tabular}\endgroup%
\HardButStrongLineBreak\kern3pt%
\begingroup \smaller\smaller\smaller\begin{tabular}{@{}c@{}}%
-1\\-22\\40
\end{tabular}\endgroup%
{$\left.\llap{\phantom{%
\begingroup \smaller\smaller\smaller\begin{tabular}{@{}c@{}}%
0\\0\\0
\end{tabular}\endgroup%
}}\!\right]$}%
%
%
\hbox{}\par\smallskip%
%
%
\leavevmode%
${L_{181.3}}$%
{} : {$[1\above{1pt}{1pt}{1}{}2\above{1pt}{1pt}{1}{}]\above{1pt}{1pt}{}{0}64\above{1pt}{1pt}{1}{1}$}\EasyButWeakLineBreak%
{${64}\above{1pt}{1pt}{l}{2}{1}\above{1pt}{1pt}{}{2}{64}\above{1pt}{1pt}{}{2}{2}\above{1pt}{1pt}{16,1}{\infty}{8}\above{1pt}{1pt}{*}{2}$}%
\nopagebreak\par%
shares genus with {$ {L_{181.2}}$}%
\nopagebreak\par%
\nopagebreak\par\leavevmode%
{$\left[\!\llap{\phantom{%
\begingroup \smaller\smaller\smaller\begin{tabular}{@{}c@{}}%
0\\0\\0
\end{tabular}\endgroup%
}}\right.$}%
\begingroup \smaller\smaller\smaller\begin{tabular}{@{}c@{}}%
-342464\\6656\\2880
\end{tabular}\endgroup%
\kern3pt%
\begingroup \smaller\smaller\smaller\begin{tabular}{@{}c@{}}%
6656\\-126\\-58
\end{tabular}\endgroup%
\kern3pt%
\begingroup \smaller\smaller\smaller\begin{tabular}{@{}c@{}}%
2880\\-58\\-23
\end{tabular}\endgroup%
{$\left.\llap{\phantom{%
\begingroup \smaller\smaller\smaller\begin{tabular}{@{}c@{}}%
0\\0\\0
\end{tabular}\endgroup%
}}\!\right]$}%
\EasyButWeakLineBreak%
{$\left[\!\llap{\phantom{%
\begingroup \smaller\smaller\smaller\begin{tabular}{@{}c@{}}%
0\\0\\0
\end{tabular}\endgroup%
}}\right.$}%
\begingroup \smaller\smaller\smaller\begin{tabular}{@{}c@{}}%
-7\\-208\\-352
\end{tabular}\endgroup%
\HardButStrongLineBreak\kern3pt%
\begingroup \smaller\smaller\smaller\begin{tabular}{@{}c@{}}%
2\\60\\99
\end{tabular}\endgroup%
\HardButStrongLineBreak\kern3pt%
\begingroup \smaller\smaller\smaller\begin{tabular}{@{}c@{}}%
45\\1344\\2240
\end{tabular}\endgroup%
\HardButStrongLineBreak\kern3pt%
\begingroup \smaller\smaller\smaller\begin{tabular}{@{}c@{}}%
4\\119\\200
\end{tabular}\endgroup%
\HardButStrongLineBreak\kern3pt%
\begingroup \smaller\smaller\smaller\begin{tabular}{@{}c@{}}%
-5\\-150\\-248
\end{tabular}\endgroup%
{$\left.\llap{\phantom{%
\begingroup \smaller\smaller\smaller\begin{tabular}{@{}c@{}}%
0\\0\\0
\end{tabular}\endgroup%
}}\!\right]$}%
%
%
\hbox{}\par\smallskip%
%
%
\leavevmode%
${L_{181.4}}$%
{} : {$1\above{1pt}{1pt}{1}{1}4\above{1pt}{1pt}{1}{7}32\above{1pt}{1pt}{1}{1}$}\EasyButWeakLineBreak%
{${4}\above{1pt}{1pt}{*}{2}{16}\above{1pt}{1pt}{l}{2}{1}\above{1pt}{1pt}{}{2}{32}\above{1pt}{1pt}{4,3}{\infty}{32}\above{1pt}{1pt}{s}{2}$}%
\nopagebreak\par%
\nopagebreak\par\leavevmode%
{$\left[\!\llap{\phantom{%
\begingroup \smaller\smaller\smaller\begin{tabular}{@{}c@{}}%
0\\0\\0
\end{tabular}\endgroup%
}}\right.$}%
\begingroup \smaller\smaller\smaller\begin{tabular}{@{}c@{}}%
-55776\\480\\3712
\end{tabular}\endgroup%
\kern3pt%
\begingroup \smaller\smaller\smaller\begin{tabular}{@{}c@{}}%
480\\-4\\-32
\end{tabular}\endgroup%
\kern3pt%
\begingroup \smaller\smaller\smaller\begin{tabular}{@{}c@{}}%
3712\\-32\\-247
\end{tabular}\endgroup%
{$\left.\llap{\phantom{%
\begingroup \smaller\smaller\smaller\begin{tabular}{@{}c@{}}%
0\\0\\0
\end{tabular}\endgroup%
}}\!\right]$}%
\EasyButWeakLineBreak%
{$\left[\!\llap{\phantom{%
\begingroup \smaller\smaller\smaller\begin{tabular}{@{}c@{}}%
0\\0\\0
\end{tabular}\endgroup%
}}\right.$}%
\begingroup \smaller\smaller\smaller\begin{tabular}{@{}c@{}}%
-1\\-8\\-14
\end{tabular}\endgroup%
\HardButStrongLineBreak\kern3pt%
\begingroup \smaller\smaller\smaller\begin{tabular}{@{}c@{}}%
-1\\6\\-16
\end{tabular}\endgroup%
\HardButStrongLineBreak\kern3pt%
\begingroup \smaller\smaller\smaller\begin{tabular}{@{}c@{}}%
1\\14\\13
\end{tabular}\endgroup%
\HardButStrongLineBreak\kern3pt%
\begingroup \smaller\smaller\smaller\begin{tabular}{@{}c@{}}%
7\\64\\96
\end{tabular}\endgroup%
\HardButStrongLineBreak\kern3pt%
\begingroup \smaller\smaller\smaller\begin{tabular}{@{}c@{}}%
1\\-8\\16
\end{tabular}\endgroup%
{$\left.\llap{\phantom{%
\begingroup \smaller\smaller\smaller\begin{tabular}{@{}c@{}}%
0\\0\\0
\end{tabular}\endgroup%
}}\!\right]$}%
%
%
\hbox{}\par\smallskip%
%
%
\leavevmode%
${L_{181.5}}$%
{} : {$1\above{1pt}{1pt}{1}{1}4\above{1pt}{1pt}{1}{1}32\above{1pt}{1pt}{1}{7}$}\EasyButWeakLineBreak%
{${4}\above{1pt}{1pt}{l}{2}{4}\above{1pt}{1pt}{}{2}{1}\above{1pt}{1pt}{r}{2}{32}\above{1pt}{1pt}{8,7}{\infty z}{32}\above{1pt}{1pt}{*}{2}$}%
\nopagebreak\par%
\nopagebreak\par\leavevmode%
{$\left[\!\llap{\phantom{%
\begingroup \smaller\smaller\smaller\begin{tabular}{@{}c@{}}%
0\\0\\0
\end{tabular}\endgroup%
}}\right.$}%
\begingroup \smaller\smaller\smaller\begin{tabular}{@{}c@{}}%
-46112\\2112\\2144
\end{tabular}\endgroup%
\kern3pt%
\begingroup \smaller\smaller\smaller\begin{tabular}{@{}c@{}}%
2112\\-92\\-100
\end{tabular}\endgroup%
\kern3pt%
\begingroup \smaller\smaller\smaller\begin{tabular}{@{}c@{}}%
2144\\-100\\-99
\end{tabular}\endgroup%
{$\left.\llap{\phantom{%
\begingroup \smaller\smaller\smaller\begin{tabular}{@{}c@{}}%
0\\0\\0
\end{tabular}\endgroup%
}}\!\right]$}%
\EasyButWeakLineBreak%
{$\left[\!\llap{\phantom{%
\begingroup \smaller\smaller\smaller\begin{tabular}{@{}c@{}}%
0\\0\\0
\end{tabular}\endgroup%
}}\right.$}%
\begingroup \smaller\smaller\smaller\begin{tabular}{@{}c@{}}%
-5\\-30\\-78
\end{tabular}\endgroup%
\HardButStrongLineBreak\kern3pt%
\begingroup \smaller\smaller\smaller\begin{tabular}{@{}c@{}}%
2\\11\\32
\end{tabular}\endgroup%
\HardButStrongLineBreak\kern3pt%
\begingroup \smaller\smaller\smaller\begin{tabular}{@{}c@{}}%
9\\53\\141
\end{tabular}\endgroup%
\HardButStrongLineBreak\kern3pt%
\begingroup \smaller\smaller\smaller\begin{tabular}{@{}c@{}}%
41\\244\\640
\end{tabular}\endgroup%
\HardButStrongLineBreak\kern3pt%
\begingroup \smaller\smaller\smaller\begin{tabular}{@{}c@{}}%
-5\\-28\\-80
\end{tabular}\endgroup%
{$\left.\llap{\phantom{%
\begingroup \smaller\smaller\smaller\begin{tabular}{@{}c@{}}%
0\\0\\0
\end{tabular}\endgroup%
}}\!\right]$}%
%
%
\hbox{}\par\smallskip%
%
%
\leavevmode%
${L_{181.6}}$%
{} : {$1\above{1pt}{1pt}{1}{7}8\above{1pt}{1pt}{1}{1}64\above{1pt}{1pt}{1}{1}$}\EasyButWeakLineBreak%
{${64}\above{1pt}{1pt}{r}{2}{4}\above{1pt}{1pt}{b}{2}{64}\above{1pt}{1pt}{l}{2}{8}\above{1pt}{1pt}{16,9}{\infty}{8}\above{1pt}{1pt}{}{2}$}%
\nopagebreak\par%
\nopagebreak\par\leavevmode%
{$\left[\!\llap{\phantom{%
\begingroup \smaller\smaller\smaller\begin{tabular}{@{}c@{}}%
0\\0\\0
\end{tabular}\endgroup%
}}\right.$}%
\begingroup \smaller\smaller\smaller\begin{tabular}{@{}c@{}}%
64\\0\\0
\end{tabular}\endgroup%
\kern3pt%
\begingroup \smaller\smaller\smaller\begin{tabular}{@{}c@{}}%
0\\-56\\-8
\end{tabular}\endgroup%
\kern3pt%
\begingroup \smaller\smaller\smaller\begin{tabular}{@{}c@{}}%
0\\-8\\-1
\end{tabular}\endgroup%
{$\left.\llap{\phantom{%
\begingroup \smaller\smaller\smaller\begin{tabular}{@{}c@{}}%
0\\0\\0
\end{tabular}\endgroup%
}}\!\right]$}%
\EasyButWeakLineBreak%
{$\left[\!\llap{\phantom{%
\begingroup \smaller\smaller\smaller\begin{tabular}{@{}c@{}}%
0\\0\\0
\end{tabular}\endgroup%
}}\right.$}%
\begingroup \smaller\smaller\smaller\begin{tabular}{@{}c@{}}%
-1\\0\\0
\end{tabular}\endgroup%
\HardButStrongLineBreak\kern3pt%
\begingroup \smaller\smaller\smaller\begin{tabular}{@{}c@{}}%
0\\-1\\6
\end{tabular}\endgroup%
\HardButStrongLineBreak\kern3pt%
\begingroup \smaller\smaller\smaller\begin{tabular}{@{}c@{}}%
3\\-8\\32
\end{tabular}\endgroup%
\HardButStrongLineBreak\kern3pt%
\begingroup \smaller\smaller\smaller\begin{tabular}{@{}c@{}}%
1\\-1\\0
\end{tabular}\endgroup%
\HardButStrongLineBreak\kern3pt%
\begingroup \smaller\smaller\smaller\begin{tabular}{@{}c@{}}%
0\\1\\-8
\end{tabular}\endgroup%
{$\left.\llap{\phantom{%
\begingroup \smaller\smaller\smaller\begin{tabular}{@{}c@{}}%
0\\0\\0
\end{tabular}\endgroup%
}}\!\right]$}%

\medskip%
%
\leavevmode\llap{}%
$W_{182}$%
\qquad\llap{6} lattices, $\chi=48$%
\hfill%
$\slashinfty\slashtwo\slashinfty\slashtwo\slashinfty\slashtwo\slashinfty\slashtwo\rtimes D_{8}$%
\nopagebreak\smallskip\hrule\nopagebreak\medskip%
%
%
\leavevmode%
${L_{182.1}}$%
{} : {$1\above{1pt}{1pt}{2}{0}64\above{1pt}{1pt}{1}{1}$}\EasyButWeakLineBreak%
{${64}\above{1pt}{1pt}{4,3}{\infty z}{64}\above{1pt}{1pt}{s}{2}{4}\above{1pt}{1pt}{16,1}{\infty z}{1}\above{1pt}{1pt}{}{2}{64}\above{1pt}{1pt}{2,1}{\infty}{64}\above{1pt}{1pt}{*}{2}{4}\above{1pt}{1pt}{16,9}{\infty z}{1}\above{1pt}{1pt}{r}{2}$}%
\nopagebreak\par%
\nopagebreak\par\leavevmode%
{$\left[\!\llap{\phantom{%
\begingroup \smaller\smaller\smaller\begin{tabular}{@{}c@{}}%
0\\0\\0
\end{tabular}\endgroup%
}}\right.$}%
\begingroup \smaller\smaller\smaller\begin{tabular}{@{}c@{}}%
-181184\\448\\4800
\end{tabular}\endgroup%
\kern3pt%
\begingroup \smaller\smaller\smaller\begin{tabular}{@{}c@{}}%
448\\-1\\-12
\end{tabular}\endgroup%
\kern3pt%
\begingroup \smaller\smaller\smaller\begin{tabular}{@{}c@{}}%
4800\\-12\\-127
\end{tabular}\endgroup%
{$\left.\llap{\phantom{%
\begingroup \smaller\smaller\smaller\begin{tabular}{@{}c@{}}%
0\\0\\0
\end{tabular}\endgroup%
}}\!\right]$}%
\EasyButWeakLineBreak%
{$\left[\!\llap{\phantom{%
\begingroup \smaller\smaller\smaller\begin{tabular}{@{}c@{}}%
0\\0\\0
\end{tabular}\endgroup%
}}\right.$}%
\begingroup \smaller\smaller\smaller\begin{tabular}{@{}c@{}}%
-15\\-608\\-512
\end{tabular}\endgroup%
\HardButStrongLineBreak\kern3pt%
\begingroup \smaller\smaller\smaller\begin{tabular}{@{}c@{}}%
-1\\-64\\-32
\end{tabular}\endgroup%
\HardButStrongLineBreak\kern3pt%
\begingroup \smaller\smaller\smaller\begin{tabular}{@{}c@{}}%
1\\40\\34
\end{tabular}\endgroup%
\HardButStrongLineBreak\kern3pt%
\begingroup \smaller\smaller\smaller\begin{tabular}{@{}c@{}}%
0\\8\\-1
\end{tabular}\endgroup%
\HardButStrongLineBreak\kern3pt%
\begingroup \smaller\smaller\smaller\begin{tabular}{@{}c@{}}%
-9\\-256\\-320
\end{tabular}\endgroup%
\HardButStrongLineBreak\kern3pt%
\begingroup \smaller\smaller\smaller\begin{tabular}{@{}c@{}}%
-23\\-800\\-800
\end{tabular}\endgroup%
\HardButStrongLineBreak\kern3pt%
\begingroup \smaller\smaller\smaller\begin{tabular}{@{}c@{}}%
-7\\-256\\-242
\end{tabular}\endgroup%
\HardButStrongLineBreak\kern3pt%
\begingroup \smaller\smaller\smaller\begin{tabular}{@{}c@{}}%
-3\\-116\\-103
\end{tabular}\endgroup%
{$\left.\llap{\phantom{%
\begingroup \smaller\smaller\smaller\begin{tabular}{@{}c@{}}%
0\\0\\0
\end{tabular}\endgroup%
}}\!\right]$}%
%
%
\hbox{}\par\smallskip%
%
%
\leavevmode%
${L_{182.2}}$%
{} : {$1\above{1pt}{1pt}{1}{1}4\above{1pt}{1pt}{1}{7}64\above{1pt}{1pt}{1}{1}$}\EasyButWeakLineBreak%
{${64}\above{1pt}{1pt}{2,1}{\infty b}{64}\above{1pt}{1pt}{}{2}{1}\above{1pt}{1pt}{8,1}{\infty}{4}\above{1pt}{1pt}{s}{2}$}\relax$\,(\times2)$%
\nopagebreak\par%
\nopagebreak\par\leavevmode%
{$\left[\!\llap{\phantom{%
\begingroup \smaller\smaller\smaller\begin{tabular}{@{}c@{}}%
0\\0\\0
\end{tabular}\endgroup%
}}\right.$}%
\begingroup \smaller\smaller\smaller\begin{tabular}{@{}c@{}}%
-153536\\832\\4416
\end{tabular}\endgroup%
\kern3pt%
\begingroup \smaller\smaller\smaller\begin{tabular}{@{}c@{}}%
832\\-4\\-24
\end{tabular}\endgroup%
\kern3pt%
\begingroup \smaller\smaller\smaller\begin{tabular}{@{}c@{}}%
4416\\-24\\-127
\end{tabular}\endgroup%
{$\left.\llap{\phantom{%
\begingroup \smaller\smaller\smaller\begin{tabular}{@{}c@{}}%
0\\0\\0
\end{tabular}\endgroup%
}}\!\right]$}%
\hfil\penalty500%
{$\left[\!\llap{\phantom{%
\begingroup \smaller\smaller\smaller\begin{tabular}{@{}c@{}}%
0\\0\\0
\end{tabular}\endgroup%
}}\right.$}%
\begingroup \smaller\smaller\smaller\begin{tabular}{@{}c@{}}%
-1369\\-2736\\-47424
\end{tabular}\endgroup%
\kern3pt%
\begingroup \smaller\smaller\smaller\begin{tabular}{@{}c@{}}%
9\\17\\312
\end{tabular}\endgroup%
\kern3pt%
\begingroup \smaller\smaller\smaller\begin{tabular}{@{}c@{}}%
39\\78\\1351
\end{tabular}\endgroup%
{$\left.\llap{\phantom{%
\begingroup \smaller\smaller\smaller\begin{tabular}{@{}c@{}}%
0\\0\\0
\end{tabular}\endgroup%
}}\!\right]$}%
\EasyButWeakLineBreak%
{$\left[\!\llap{\phantom{%
\begingroup \smaller\smaller\smaller\begin{tabular}{@{}c@{}}%
0\\0\\0
\end{tabular}\endgroup%
}}\right.$}%
\begingroup \smaller\smaller\smaller\begin{tabular}{@{}c@{}}%
-23\\-32\\-800
\end{tabular}\endgroup%
\HardButStrongLineBreak\kern3pt%
\begingroup \smaller\smaller\smaller\begin{tabular}{@{}c@{}}%
-9\\16\\-320
\end{tabular}\endgroup%
\HardButStrongLineBreak\kern3pt%
\begingroup \smaller\smaller\smaller\begin{tabular}{@{}c@{}}%
0\\4\\-1
\end{tabular}\endgroup%
\HardButStrongLineBreak\kern3pt%
\begingroup \smaller\smaller\smaller\begin{tabular}{@{}c@{}}%
1\\4\\34
\end{tabular}\endgroup%
{$\left.\llap{\phantom{%
\begingroup \smaller\smaller\smaller\begin{tabular}{@{}c@{}}%
0\\0\\0
\end{tabular}\endgroup%
}}\!\right]$}%
%
%
\hbox{}\par\smallskip%
%
%
\leavevmode%
${L_{182.3}}$%
{} : {$1\above{1pt}{1pt}{-}{5}4\above{1pt}{1pt}{1}{7}64\above{1pt}{1pt}{-}{5}$}\EasyButWeakLineBreak%
{${64}\above{1pt}{1pt}{8,3}{\infty z}{64}\above{1pt}{1pt}{l}{2}{1}\above{1pt}{1pt}{8,5}{\infty}{4}\above{1pt}{1pt}{*}{2}$}\relax$\,(\times2)$%
\nopagebreak\par%
\nopagebreak\par\leavevmode%
{$\left[\!\llap{\phantom{%
\begingroup \smaller\smaller\smaller\begin{tabular}{@{}c@{}}%
0\\0\\0
\end{tabular}\endgroup%
}}\right.$}%
\begingroup \smaller\smaller\smaller\begin{tabular}{@{}c@{}}%
-120000\\-1728\\3904
\end{tabular}\endgroup%
\kern3pt%
\begingroup \smaller\smaller\smaller\begin{tabular}{@{}c@{}}%
-1728\\-20\\56
\end{tabular}\endgroup%
\kern3pt%
\begingroup \smaller\smaller\smaller\begin{tabular}{@{}c@{}}%
3904\\56\\-127
\end{tabular}\endgroup%
{$\left.\llap{\phantom{%
\begingroup \smaller\smaller\smaller\begin{tabular}{@{}c@{}}%
0\\0\\0
\end{tabular}\endgroup%
}}\!\right]$}%
\hfil\penalty500%
{$\left[\!\llap{\phantom{%
\begingroup \smaller\smaller\smaller\begin{tabular}{@{}c@{}}%
0\\0\\0
\end{tabular}\endgroup%
}}\right.$}%
\begingroup \smaller\smaller\smaller\begin{tabular}{@{}c@{}}%
4799\\7200\\150400
\end{tabular}\endgroup%
\kern3pt%
\begingroup \smaller\smaller\smaller\begin{tabular}{@{}c@{}}%
60\\89\\1880
\end{tabular}\endgroup%
\kern3pt%
\begingroup \smaller\smaller\smaller\begin{tabular}{@{}c@{}}%
-156\\-234\\-4889
\end{tabular}\endgroup%
{$\left.\llap{\phantom{%
\begingroup \smaller\smaller\smaller\begin{tabular}{@{}c@{}}%
0\\0\\0
\end{tabular}\endgroup%
}}\!\right]$}%
\EasyButWeakLineBreak%
{$\left[\!\llap{\phantom{%
\begingroup \smaller\smaller\smaller\begin{tabular}{@{}c@{}}%
0\\0\\0
\end{tabular}\endgroup%
}}\right.$}%
\begingroup \smaller\smaller\smaller\begin{tabular}{@{}c@{}}%
99\\152\\3104
\end{tabular}\endgroup%
\HardButStrongLineBreak\kern3pt%
\begingroup \smaller\smaller\smaller\begin{tabular}{@{}c@{}}%
53\\88\\1664
\end{tabular}\endgroup%
\HardButStrongLineBreak\kern3pt%
\begingroup \smaller\smaller\smaller\begin{tabular}{@{}c@{}}%
2\\4\\63
\end{tabular}\endgroup%
\HardButStrongLineBreak\kern3pt%
\begingroup \smaller\smaller\smaller\begin{tabular}{@{}c@{}}%
-3\\-4\\-94
\end{tabular}\endgroup%
{$\left.\llap{\phantom{%
\begingroup \smaller\smaller\smaller\begin{tabular}{@{}c@{}}%
0\\0\\0
\end{tabular}\endgroup%
}}\!\right]$}%
%
%
%
%
%
%
%
%
%
%
%
%
%
%

\medskip%
%
\leavevmode\llap{}%
$W_{183}$%
\qquad\llap{34} lattices, $\chi=18$%
\hfill%
$\slashinfty22\slashtwo22\rtimes D_{2}$%
\nopagebreak\smallskip\hrule\nopagebreak\medskip%
%
%
\leavevmode%
${L_{183.1}}$%
{} : {$1\above{1pt}{1pt}{2}{0}8\above{1pt}{1pt}{1}{1}{\cdot}1\above{1pt}{1pt}{2}{}9\above{1pt}{1pt}{1}{}$}\EasyButWeakLineBreak%
{${8}\above{1pt}{1pt}{6,1}{\infty}{8}\above{1pt}{1pt}{*}{2}{36}\above{1pt}{1pt}{l}{2}{1}\above{1pt}{1pt}{r}{2}{4}\above{1pt}{1pt}{l}{2}{9}\above{1pt}{1pt}{}{2}$}%
\nopagebreak\par%
\nopagebreak\par\leavevmode%
{$\left[\!\llap{\phantom{%
\begingroup \smaller\smaller\smaller\begin{tabular}{@{}c@{}}%
0\\0\\0
\end{tabular}\endgroup%
}}\right.$}%
\begingroup \smaller\smaller\smaller\begin{tabular}{@{}c@{}}%
-444600\\4680\\4392
\end{tabular}\endgroup%
\kern3pt%
\begingroup \smaller\smaller\smaller\begin{tabular}{@{}c@{}}%
4680\\-47\\-49
\end{tabular}\endgroup%
\kern3pt%
\begingroup \smaller\smaller\smaller\begin{tabular}{@{}c@{}}%
4392\\-49\\-40
\end{tabular}\endgroup%
{$\left.\llap{\phantom{%
\begingroup \smaller\smaller\smaller\begin{tabular}{@{}c@{}}%
0\\0\\0
\end{tabular}\endgroup%
}}\!\right]$}%
\EasyButWeakLineBreak%
{$\left[\!\llap{\phantom{%
\begingroup \smaller\smaller\smaller\begin{tabular}{@{}c@{}}%
0\\0\\0
\end{tabular}\endgroup%
}}\right.$}%
\begingroup \smaller\smaller\smaller\begin{tabular}{@{}c@{}}%
15\\808\\656
\end{tabular}\endgroup%
\HardButStrongLineBreak\kern3pt%
\begingroup \smaller\smaller\smaller\begin{tabular}{@{}c@{}}%
-7\\-376\\-308
\end{tabular}\endgroup%
\HardButStrongLineBreak\kern3pt%
\begingroup \smaller\smaller\smaller\begin{tabular}{@{}c@{}}%
-7\\-378\\-306
\end{tabular}\endgroup%
\HardButStrongLineBreak\kern3pt%
\begingroup \smaller\smaller\smaller\begin{tabular}{@{}c@{}}%
3\\161\\132
\end{tabular}\endgroup%
\HardButStrongLineBreak\kern3pt%
\begingroup \smaller\smaller\smaller\begin{tabular}{@{}c@{}}%
17\\914\\746
\end{tabular}\endgroup%
\HardButStrongLineBreak\kern3pt%
\begingroup \smaller\smaller\smaller\begin{tabular}{@{}c@{}}%
46\\2475\\2016
\end{tabular}\endgroup%
{$\left.\llap{\phantom{%
\begingroup \smaller\smaller\smaller\begin{tabular}{@{}c@{}}%
0\\0\\0
\end{tabular}\endgroup%
}}\!\right]$}%
%
%
\hbox{}\par\smallskip%
%
%
\leavevmode%
${L_{183.2}}$%
{} : {$[1\above{1pt}{1pt}{-}{}2\above{1pt}{1pt}{1}{}]\above{1pt}{1pt}{}{6}16\above{1pt}{1pt}{-}{3}{\cdot}1\above{1pt}{1pt}{2}{}9\above{1pt}{1pt}{1}{}$}\spacer%
\instructions{2}%
\EasyButWeakLineBreak%
{${8}\above{1pt}{1pt}{24,1}{\infty z}{2}\above{1pt}{1pt}{r}{2}{36}\above{1pt}{1pt}{*}{2}{16}\above{1pt}{1pt}{*}{2}{4}\above{1pt}{1pt}{*}{2}{144}\above{1pt}{1pt}{s}{2}$}%
\nopagebreak\par%
\nopagebreak\par\leavevmode%
{$\left[\!\llap{\phantom{%
\begingroup \smaller\smaller\smaller\begin{tabular}{@{}c@{}}%
0\\0\\0
\end{tabular}\endgroup%
}}\right.$}%
\begingroup \smaller\smaller\smaller\begin{tabular}{@{}c@{}}%
-463824\\-11520\\-12240
\end{tabular}\endgroup%
\kern3pt%
\begingroup \smaller\smaller\smaller\begin{tabular}{@{}c@{}}%
-11520\\-286\\-304
\end{tabular}\endgroup%
\kern3pt%
\begingroup \smaller\smaller\smaller\begin{tabular}{@{}c@{}}%
-12240\\-304\\-323
\end{tabular}\endgroup%
{$\left.\llap{\phantom{%
\begingroup \smaller\smaller\smaller\begin{tabular}{@{}c@{}}%
0\\0\\0
\end{tabular}\endgroup%
}}\!\right]$}%
\EasyButWeakLineBreak%
{$\left[\!\llap{\phantom{%
\begingroup \smaller\smaller\smaller\begin{tabular}{@{}c@{}}%
0\\0\\0
\end{tabular}\endgroup%
}}\right.$}%
\begingroup \smaller\smaller\smaller\begin{tabular}{@{}c@{}}%
-1\\2\\36
\end{tabular}\endgroup%
\HardButStrongLineBreak\kern3pt%
\begingroup \smaller\smaller\smaller\begin{tabular}{@{}c@{}}%
0\\-13\\12
\end{tabular}\endgroup%
\HardButStrongLineBreak\kern3pt%
\begingroup \smaller\smaller\smaller\begin{tabular}{@{}c@{}}%
5\\-108\\-90
\end{tabular}\endgroup%
\HardButStrongLineBreak\kern3pt%
\begingroup \smaller\smaller\smaller\begin{tabular}{@{}c@{}}%
3\\-28\\-88
\end{tabular}\endgroup%
\HardButStrongLineBreak\kern3pt%
\begingroup \smaller\smaller\smaller\begin{tabular}{@{}c@{}}%
1\\0\\-38
\end{tabular}\endgroup%
\HardButStrongLineBreak\kern3pt%
\begingroup \smaller\smaller\smaller\begin{tabular}{@{}c@{}}%
1\\36\\-72
\end{tabular}\endgroup%
{$\left.\llap{\phantom{%
\begingroup \smaller\smaller\smaller\begin{tabular}{@{}c@{}}%
0\\0\\0
\end{tabular}\endgroup%
}}\!\right]$}%
%
%
\hbox{}\par\smallskip%
%
%
\leavevmode%
${L_{183.3}}$%
{} : {$[1\above{1pt}{1pt}{1}{}2\above{1pt}{1pt}{1}{}]\above{1pt}{1pt}{}{2}16\above{1pt}{1pt}{1}{7}{\cdot}1\above{1pt}{1pt}{2}{}9\above{1pt}{1pt}{1}{}$}\spacer%
\instructions{m}%
\EasyButWeakLineBreak%
{${8}\above{1pt}{1pt}{24,13}{\infty z}{2}\above{1pt}{1pt}{}{2}{9}\above{1pt}{1pt}{r}{2}{16}\above{1pt}{1pt}{l}{2}{1}\above{1pt}{1pt}{r}{2}{144}\above{1pt}{1pt}{*}{2}$}%
\nopagebreak\par%
\nopagebreak\par\leavevmode%
{$\left[\!\llap{\phantom{%
\begingroup \smaller\smaller\smaller\begin{tabular}{@{}c@{}}%
0\\0\\0
\end{tabular}\endgroup%
}}\right.$}%
\begingroup \smaller\smaller\smaller\begin{tabular}{@{}c@{}}%
-123408\\6048\\2448
\end{tabular}\endgroup%
\kern3pt%
\begingroup \smaller\smaller\smaller\begin{tabular}{@{}c@{}}%
6048\\-286\\-124
\end{tabular}\endgroup%
\kern3pt%
\begingroup \smaller\smaller\smaller\begin{tabular}{@{}c@{}}%
2448\\-124\\-47
\end{tabular}\endgroup%
{$\left.\llap{\phantom{%
\begingroup \smaller\smaller\smaller\begin{tabular}{@{}c@{}}%
0\\0\\0
\end{tabular}\endgroup%
}}\!\right]$}%
\EasyButWeakLineBreak%
{$\left[\!\llap{\phantom{%
\begingroup \smaller\smaller\smaller\begin{tabular}{@{}c@{}}%
0\\0\\0
\end{tabular}\endgroup%
}}\right.$}%
\begingroup \smaller\smaller\smaller\begin{tabular}{@{}c@{}}%
-7\\-70\\-180
\end{tabular}\endgroup%
\HardButStrongLineBreak\kern3pt%
\begingroup \smaller\smaller\smaller\begin{tabular}{@{}c@{}}%
6\\59\\156
\end{tabular}\endgroup%
\HardButStrongLineBreak\kern3pt%
\begingroup \smaller\smaller\smaller\begin{tabular}{@{}c@{}}%
40\\396\\1035
\end{tabular}\endgroup%
\HardButStrongLineBreak\kern3pt%
\begingroup \smaller\smaller\smaller\begin{tabular}{@{}c@{}}%
31\\308\\800
\end{tabular}\endgroup%
\HardButStrongLineBreak\kern3pt%
\begingroup \smaller\smaller\smaller\begin{tabular}{@{}c@{}}%
3\\30\\77
\end{tabular}\endgroup%
\HardButStrongLineBreak\kern3pt%
\begingroup \smaller\smaller\smaller\begin{tabular}{@{}c@{}}%
-11\\-108\\-288
\end{tabular}\endgroup%
{$\left.\llap{\phantom{%
\begingroup \smaller\smaller\smaller\begin{tabular}{@{}c@{}}%
0\\0\\0
\end{tabular}\endgroup%
}}\!\right]$}%
%
%
\hbox{}\par\smallskip%
%
%
\leavevmode%
${L_{183.4}}$%
{} : {$[1\above{1pt}{1pt}{-}{}2\above{1pt}{1pt}{1}{}]\above{1pt}{1pt}{}{4}16\above{1pt}{1pt}{-}{5}{\cdot}1\above{1pt}{1pt}{2}{}9\above{1pt}{1pt}{1}{}$}\spacer%
\instructions{m}%
\EasyButWeakLineBreak%
{${2}\above{1pt}{1pt}{24,13}{\infty}{8}\above{1pt}{1pt}{*}{2}{36}\above{1pt}{1pt}{s}{2}{16}\above{1pt}{1pt}{s}{2}{4}\above{1pt}{1pt}{s}{2}{144}\above{1pt}{1pt}{l}{2}$}%
\nopagebreak\par%
\nopagebreak\par\leavevmode%
{$\left[\!\llap{\phantom{%
\begingroup \smaller\smaller\smaller\begin{tabular}{@{}c@{}}%
0\\0\\0
\end{tabular}\endgroup%
}}\right.$}%
\begingroup \smaller\smaller\smaller\begin{tabular}{@{}c@{}}%
-60336\\17712\\2880
\end{tabular}\endgroup%
\kern3pt%
\begingroup \smaller\smaller\smaller\begin{tabular}{@{}c@{}}%
17712\\-5198\\-848
\end{tabular}\endgroup%
\kern3pt%
\begingroup \smaller\smaller\smaller\begin{tabular}{@{}c@{}}%
2880\\-848\\-133
\end{tabular}\endgroup%
{$\left.\llap{\phantom{%
\begingroup \smaller\smaller\smaller\begin{tabular}{@{}c@{}}%
0\\0\\0
\end{tabular}\endgroup%
}}\!\right]$}%
\EasyButWeakLineBreak%
{$\left[\!\llap{\phantom{%
\begingroup \smaller\smaller\smaller\begin{tabular}{@{}c@{}}%
0\\0\\0
\end{tabular}\endgroup%
}}\right.$}%
\begingroup \smaller\smaller\smaller\begin{tabular}{@{}c@{}}%
-28\\-87\\-52
\end{tabular}\endgroup%
\HardButStrongLineBreak\kern3pt%
\begingroup \smaller\smaller\smaller\begin{tabular}{@{}c@{}}%
25\\78\\44
\end{tabular}\endgroup%
\HardButStrongLineBreak\kern3pt%
\begingroup \smaller\smaller\smaller\begin{tabular}{@{}c@{}}%
29\\90\\54
\end{tabular}\endgroup%
\HardButStrongLineBreak\kern3pt%
\begingroup \smaller\smaller\smaller\begin{tabular}{@{}c@{}}%
-41\\-128\\-72
\end{tabular}\endgroup%
\HardButStrongLineBreak\kern3pt%
\begingroup \smaller\smaller\smaller\begin{tabular}{@{}c@{}}%
-61\\-190\\-110
\end{tabular}\endgroup%
\HardButStrongLineBreak\kern3pt%
\begingroup \smaller\smaller\smaller\begin{tabular}{@{}c@{}}%
-671\\-2088\\-1224
\end{tabular}\endgroup%
{$\left.\llap{\phantom{%
\begingroup \smaller\smaller\smaller\begin{tabular}{@{}c@{}}%
0\\0\\0
\end{tabular}\endgroup%
}}\!\right]$}%
%
%
\hbox{}\par\smallskip%
%
%
\leavevmode%
${L_{183.5}}$%
{} : {$1\above{1pt}{1pt}{1}{1}8\above{1pt}{1pt}{1}{7}64\above{1pt}{1pt}{1}{1}{\cdot}1\above{1pt}{1pt}{2}{}9\above{1pt}{1pt}{1}{}$}\spacer%
\instructions{2}%
\EasyButWeakLineBreak%
{${8}\above{1pt}{1pt}{24,1}{\infty b}{32}\above{1pt}{1pt}{*}{2}{36}\above{1pt}{1pt}{s}{2}{64}\above{1pt}{1pt}{s}{2}{4}\above{1pt}{1pt}{s}{2}{576}\above{1pt}{1pt}{b}{2}$}%
\nopagebreak\par%
shares genus with {$ {L_{183.6}}$}%
; isometric to own %
2-dual\nopagebreak\par%
\nopagebreak\par\leavevmode%
{$\left[\!\llap{\phantom{%
\begingroup \smaller\smaller\smaller\begin{tabular}{@{}c@{}}%
0\\0\\0
\end{tabular}\endgroup%
}}\right.$}%
\begingroup \smaller\smaller\smaller\begin{tabular}{@{}c@{}}%
-607680\\9216\\95040
\end{tabular}\endgroup%
\kern3pt%
\begingroup \smaller\smaller\smaller\begin{tabular}{@{}c@{}}%
9216\\-136\\-1488
\end{tabular}\endgroup%
\kern3pt%
\begingroup \smaller\smaller\smaller\begin{tabular}{@{}c@{}}%
95040\\-1488\\-14287
\end{tabular}\endgroup%
{$\left.\llap{\phantom{%
\begingroup \smaller\smaller\smaller\begin{tabular}{@{}c@{}}%
0\\0\\0
\end{tabular}\endgroup%
}}\!\right]$}%
\EasyButWeakLineBreak%
{$\left[\!\llap{\phantom{%
\begingroup \smaller\smaller\smaller\begin{tabular}{@{}c@{}}%
0\\0\\0
\end{tabular}\endgroup%
}}\right.$}%
\begingroup \smaller\smaller\smaller\begin{tabular}{@{}c@{}}%
62\\2231\\180
\end{tabular}\endgroup%
\HardButStrongLineBreak\kern3pt%
\begingroup \smaller\smaller\smaller\begin{tabular}{@{}c@{}}%
-33\\-1186\\-96
\end{tabular}\endgroup%
\HardButStrongLineBreak\kern3pt%
\begingroup \smaller\smaller\smaller\begin{tabular}{@{}c@{}}%
-31\\-1116\\-90
\end{tabular}\endgroup%
\HardButStrongLineBreak\kern3pt%
\begingroup \smaller\smaller\smaller\begin{tabular}{@{}c@{}}%
55\\1976\\160
\end{tabular}\endgroup%
\HardButStrongLineBreak\kern3pt%
\begingroup \smaller\smaller\smaller\begin{tabular}{@{}c@{}}%
53\\1906\\154
\end{tabular}\endgroup%
\HardButStrongLineBreak\kern3pt%
\begingroup \smaller\smaller\smaller\begin{tabular}{@{}c@{}}%
1289\\46368\\3744
\end{tabular}\endgroup%
{$\left.\llap{\phantom{%
\begingroup \smaller\smaller\smaller\begin{tabular}{@{}c@{}}%
0\\0\\0
\end{tabular}\endgroup%
}}\!\right]$}%
%
%
\hbox{}\par\smallskip%
%
%
\leavevmode%
${L_{183.6}}$%
{} : {$1\above{1pt}{1pt}{1}{1}8\above{1pt}{1pt}{1}{7}64\above{1pt}{1pt}{1}{1}{\cdot}1\above{1pt}{1pt}{2}{}9\above{1pt}{1pt}{1}{}$}\EasyButWeakLineBreak%
{${8}\above{1pt}{1pt}{24,1}{\infty a}{32}\above{1pt}{1pt}{l}{2}{9}\above{1pt}{1pt}{}{2}{64}\above{1pt}{1pt}{}{2}{1}\above{1pt}{1pt}{}{2}{576}\above{1pt}{1pt}{r}{2}$}%
\nopagebreak\par%
shares genus with {$ {L_{183.5}}$}%
; isometric to own %
2-dual\nopagebreak\par%
\nopagebreak\par\leavevmode%
{$\left[\!\llap{\phantom{%
\begingroup \smaller\smaller\smaller\begin{tabular}{@{}c@{}}%
0\\0\\0
\end{tabular}\endgroup%
}}\right.$}%
\begingroup \smaller\smaller\smaller\begin{tabular}{@{}c@{}}%
-42328512\\78336\\-13117824
\end{tabular}\endgroup%
\kern3pt%
\begingroup \smaller\smaller\smaller\begin{tabular}{@{}c@{}}%
78336\\-136\\24384
\end{tabular}\endgroup%
\kern3pt%
\begingroup \smaller\smaller\smaller\begin{tabular}{@{}c@{}}%
-13117824\\24384\\-4063999
\end{tabular}\endgroup%
{$\left.\llap{\phantom{%
\begingroup \smaller\smaller\smaller\begin{tabular}{@{}c@{}}%
0\\0\\0
\end{tabular}\endgroup%
}}\!\right]$}%
\EasyButWeakLineBreak%
{$\left[\!\llap{\phantom{%
\begingroup \smaller\smaller\smaller\begin{tabular}{@{}c@{}}%
0\\0\\0
\end{tabular}\endgroup%
}}\right.$}%
\begingroup \smaller\smaller\smaller\begin{tabular}{@{}c@{}}%
-1016\\-36577\\3060
\end{tabular}\endgroup%
\HardButStrongLineBreak\kern3pt%
\begingroup \smaller\smaller\smaller\begin{tabular}{@{}c@{}}%
255\\9182\\-768
\end{tabular}\endgroup%
\HardButStrongLineBreak\kern3pt%
\begingroup \smaller\smaller\smaller\begin{tabular}{@{}c@{}}%
254\\9144\\-765
\end{tabular}\endgroup%
\HardButStrongLineBreak\kern3pt%
\begingroup \smaller\smaller\smaller\begin{tabular}{@{}c@{}}%
-425\\-15304\\1280
\end{tabular}\endgroup%
\HardButStrongLineBreak\kern3pt%
\begingroup \smaller\smaller\smaller\begin{tabular}{@{}c@{}}%
-339\\-12205\\1021
\end{tabular}\endgroup%
\HardButStrongLineBreak\kern3pt%
\begingroup \smaller\smaller\smaller\begin{tabular}{@{}c@{}}%
-18551\\-667872\\55872
\end{tabular}\endgroup%
{$\left.\llap{\phantom{%
\begingroup \smaller\smaller\smaller\begin{tabular}{@{}c@{}}%
0\\0\\0
\end{tabular}\endgroup%
}}\!\right]$}%

\medskip%
%
\leavevmode\llap{}%
$W_{184}$%
\qquad\llap{34} lattices, $\chi=72$%
\hfill%
$\slashinfty2|2\slashinfty2|2\slashinfty2|2\slashinfty2|2\rtimes D_{8}$%
\nopagebreak\smallskip\hrule\nopagebreak\medskip%
%
%
\leavevmode%
${L_{184.1}}$%
{} : {$1\above{1pt}{1pt}{2}{0}8\above{1pt}{1pt}{1}{1}{\cdot}1\above{1pt}{1pt}{-2}{}9\above{1pt}{1pt}{-}{}$}\EasyButWeakLineBreak%
{${72}\above{1pt}{1pt}{2,1}{\infty}{72}\above{1pt}{1pt}{*}{2}{4}\above{1pt}{1pt}{*}{2}{8}\above{1pt}{1pt}{3,2}{\infty b}{8}\above{1pt}{1pt}{}{2}{1}\above{1pt}{1pt}{}{2}$}\relax$\,(\times2)$%
\nopagebreak\par%
\nopagebreak\par\leavevmode%
{$\left[\!\llap{\phantom{%
\begingroup \smaller\smaller\smaller
\endgroup%
}}\!\right]$}%
%
%
\hbox{}\par\smallskip%
%
%
\leavevmode%
${L_{184.2}}$%
{} : {$[1\above{1pt}{1pt}{-}{}2\above{1pt}{1pt}{1}{}]\above{1pt}{1pt}{}{6}16\above{1pt}{1pt}{-}{3}{\cdot}1\above{1pt}{1pt}{-2}{}9\above{1pt}{1pt}{-}{}$}\spacer%
\instructions{2}%
\EasyButWeakLineBreak%
{${72}\above{1pt}{1pt}{8,1}{\infty z}{18}\above{1pt}{1pt}{r}{2}{4}\above{1pt}{1pt}{l}{2}{2}\above{1pt}{1pt}{24,11}{\infty}{8}\above{1pt}{1pt}{s}{2}{16}\above{1pt}{1pt}{s}{2}$}\relax$\,(\times2)$%
\nopagebreak\par%
\nopagebreak\par\leavevmode%
{$\left[\!\llap{\phantom{%
\begingroup \smaller\smaller\smaller
\endgroup%
}}\!\right]$}%
%
%
\hbox{}\par\smallskip%
%
%
\leavevmode%
${L_{184.3}}$%
{} : {$[1\above{1pt}{1pt}{1}{}2\above{1pt}{1pt}{1}{}]\above{1pt}{1pt}{}{2}16\above{1pt}{1pt}{1}{7}{\cdot}1\above{1pt}{1pt}{-2}{}9\above{1pt}{1pt}{-}{}$}\spacer%
\instructions{m}%
\EasyButWeakLineBreak%
{${72}\above{1pt}{1pt}{8,5}{\infty z}{18}\above{1pt}{1pt}{}{2}{1}\above{1pt}{1pt}{}{2}{2}\above{1pt}{1pt}{24,23}{\infty}{8}\above{1pt}{1pt}{*}{2}{16}\above{1pt}{1pt}{*}{2}$}\relax$\,(\times2)$%
\nopagebreak\par%
\nopagebreak\par\leavevmode%
{$\left[\!\llap{\phantom{%
\begingroup \smaller\smaller\smaller
\endgroup%
}}\!\right]$}%
%
%
\hbox{}\par\smallskip%
%
%
\leavevmode%
${L_{184.4}}$%
{} : {$[1\above{1pt}{1pt}{-}{}2\above{1pt}{1pt}{1}{}]\above{1pt}{1pt}{}{4}16\above{1pt}{1pt}{-}{5}{\cdot}1\above{1pt}{1pt}{-2}{}9\above{1pt}{1pt}{-}{}$}\spacer%
\instructions{m}%
\EasyButWeakLineBreak%
{${18}\above{1pt}{1pt}{8,5}{\infty}{72}\above{1pt}{1pt}{*}{2}{4}\above{1pt}{1pt}{*}{2}{8}\above{1pt}{1pt}{24,11}{\infty z}{2}\above{1pt}{1pt}{r}{2}{16}\above{1pt}{1pt}{l}{2}$}\relax$\,(\times2)$%
\nopagebreak\par%
\nopagebreak\par\leavevmode%
{$\left[\!\llap{\phantom{%
\begingroup \smaller\smaller\smaller
\endgroup%
}}\!\right]$}%
%
%
\hbox{}\par\smallskip%
%
%
\leavevmode%
${L_{184.5}}$%
{} : {$1\above{1pt}{1pt}{1}{1}8\above{1pt}{1pt}{1}{7}64\above{1pt}{1pt}{1}{1}{\cdot}1\above{1pt}{1pt}{-2}{}9\above{1pt}{1pt}{-}{}$}\spacer%
\instructions{2}%
\EasyButWeakLineBreak%
{${72}\above{1pt}{1pt}{8,1}{\infty b}{288}\above{1pt}{1pt}{*}{2}{4}\above{1pt}{1pt}{*}{2}{32}\above{1pt}{1pt}{48,23}{\infty z}{8}\above{1pt}{1pt}{b}{2}{64}\above{1pt}{1pt}{b}{2}$}\relax$\,(\times2)$%
\nopagebreak\par%
shares genus with {$ {L_{184.6}}$}%
; isometric to own %
2-dual\nopagebreak\par%
\nopagebreak\par\leavevmode%
{$\left[\!\llap{\phantom{%
\begingroup \smaller\smaller\smaller
\endgroup%
}}\!\right]$}%
%
%
\hbox{}\par\smallskip%
%
%
\leavevmode%
${L_{184.6}}$%
{} : {$1\above{1pt}{1pt}{1}{1}8\above{1pt}{1pt}{1}{7}64\above{1pt}{1pt}{1}{1}{\cdot}1\above{1pt}{1pt}{-2}{}9\above{1pt}{1pt}{-}{}$}\EasyButWeakLineBreak%
{${72}\above{1pt}{1pt}{8,1}{\infty a}{288}\above{1pt}{1pt}{l}{2}{1}\above{1pt}{1pt}{r}{2}{32}\above{1pt}{1pt}{48,47}{\infty z}{8}\above{1pt}{1pt}{l}{2}{64}\above{1pt}{1pt}{r}{2}$}\relax$\,(\times2)$%
\nopagebreak\par%
shares genus with {$ {L_{184.5}}$}%
; isometric to own %
2-dual\nopagebreak\par%
\nopagebreak\par\leavevmode%
{$\left[\!\llap{\phantom{%
\begingroup \smaller\smaller\smaller
\endgroup%
}}\!\right]$}%

\medskip%
%
\leavevmode\llap{}%
$W_{185}$%
\qquad\llap{12} lattices, $\chi=30$%
\hfill%
$22242224\rtimes C_{2}$%
\nopagebreak\smallskip\hrule\nopagebreak\medskip%
%
%
\leavevmode%
${L_{185.1}}$%
{} : {$1\above{1pt}{1pt}{-2}{{\rm II}}4\above{1pt}{1pt}{1}{7}{\cdot}1\above{1pt}{1pt}{2}{}3\above{1pt}{1pt}{-}{}{\cdot}1\above{1pt}{1pt}{2}{}25\above{1pt}{1pt}{-}{}$}\spacer%
\instructions{2}%
\EasyButWeakLineBreak%
{${2}\above{1pt}{1pt}{b}{2}{50}\above{1pt}{1pt}{s}{2}{6}\above{1pt}{1pt}{b}{2}{4}\above{1pt}{1pt}{*}{4}$}\relax$\,(\times2)$%
\nopagebreak\par%
\nopagebreak\par\leavevmode%
{$\left[\!\llap{\phantom{%
\begingroup \smaller\smaller\smaller\begin{tabular}{@{}c@{}}%
0\\0\\0
\end{tabular}\endgroup%
}}\right.$}%
\begingroup \smaller\smaller\smaller\begin{tabular}{@{}c@{}}%
-221700\\73800\\1800
\end{tabular}\endgroup%
\kern3pt%
\begingroup \smaller\smaller\smaller\begin{tabular}{@{}c@{}}%
73800\\-24566\\-601
\end{tabular}\endgroup%
\kern3pt%
\begingroup \smaller\smaller\smaller\begin{tabular}{@{}c@{}}%
1800\\-601\\-10
\end{tabular}\endgroup%
{$\left.\llap{\phantom{%
\begingroup \smaller\smaller\smaller\begin{tabular}{@{}c@{}}%
0\\0\\0
\end{tabular}\endgroup%
}}\!\right]$}%
\hfil\penalty500%
{$\left[\!\llap{\phantom{%
\begingroup \smaller\smaller\smaller\begin{tabular}{@{}c@{}}%
0\\0\\0
\end{tabular}\endgroup%
}}\right.$}%
\begingroup \smaller\smaller\smaller\begin{tabular}{@{}c@{}}%
61499\\183000\\70500
\end{tabular}\endgroup%
\kern3pt%
\begingroup \smaller\smaller\smaller\begin{tabular}{@{}c@{}}%
-20541\\-61123\\-23547
\end{tabular}\endgroup%
\kern3pt%
\begingroup \smaller\smaller\smaller\begin{tabular}{@{}c@{}}%
-328\\-976\\-377
\end{tabular}\endgroup%
{$\left.\llap{\phantom{%
\begingroup \smaller\smaller\smaller\begin{tabular}{@{}c@{}}%
0\\0\\0
\end{tabular}\endgroup%
}}\!\right]$}%
\EasyButWeakLineBreak%
{$\left[\!\llap{\phantom{%
\begingroup \smaller\smaller\smaller\begin{tabular}{@{}c@{}}%
0\\0\\0
\end{tabular}\endgroup%
}}\right.$}%
\begingroup \smaller\smaller\smaller\begin{tabular}{@{}c@{}}%
1\\3\\-1
\end{tabular}\endgroup%
\HardButStrongLineBreak\kern3pt%
\begingroup \smaller\smaller\smaller\begin{tabular}{@{}c@{}}%
294\\875\\325
\end{tabular}\endgroup%
\HardButStrongLineBreak\kern3pt%
\begingroup \smaller\smaller\smaller\begin{tabular}{@{}c@{}}%
124\\369\\141
\end{tabular}\endgroup%
\HardButStrongLineBreak\kern3pt%
\begingroup \smaller\smaller\smaller\begin{tabular}{@{}c@{}}%
287\\854\\330
\end{tabular}\endgroup%
{$\left.\llap{\phantom{%
\begingroup \smaller\smaller\smaller\begin{tabular}{@{}c@{}}%
0\\0\\0
\end{tabular}\endgroup%
}}\!\right]$}%

\medskip%
%
\leavevmode\llap{}%
$W_{186}$%
\qquad\llap{12} lattices, $\chi=10$%
\hfill%
$22226$%
\nopagebreak\smallskip\hrule\nopagebreak\medskip%
%
%
\leavevmode%
${L_{186.1}}$%
{} : {$1\above{1pt}{1pt}{-2}{{\rm II}}4\above{1pt}{1pt}{1}{7}{\cdot}1\above{1pt}{1pt}{2}{}3\above{1pt}{1pt}{-}{}{\cdot}1\above{1pt}{1pt}{-2}{}25\above{1pt}{1pt}{1}{}$}\spacer%
\instructions{2}%
\EasyButWeakLineBreak%
{${6}\above{1pt}{1pt}{b}{2}{100}\above{1pt}{1pt}{*}{2}{4}\above{1pt}{1pt}{b}{2}{150}\above{1pt}{1pt}{s}{2}{2}\above{1pt}{1pt}{}{6}$}%
\nopagebreak\par%
\nopagebreak\par\leavevmode%
{$\left[\!\llap{\phantom{%
\begingroup \smaller\smaller\smaller\begin{tabular}{@{}c@{}}%
0\\0\\0
\end{tabular}\endgroup%
}}\right.$}%
\begingroup \smaller\smaller\smaller\begin{tabular}{@{}c@{}}%
-2716400100\\11466900\\-790200
\end{tabular}\endgroup%
\kern3pt%
\begingroup \smaller\smaller\smaller\begin{tabular}{@{}c@{}}%
11466900\\-48386\\3309
\end{tabular}\endgroup%
\kern3pt%
\begingroup \smaller\smaller\smaller\begin{tabular}{@{}c@{}}%
-790200\\3309\\-194
\end{tabular}\endgroup%
{$\left.\llap{\phantom{%
\begingroup \smaller\smaller\smaller\begin{tabular}{@{}c@{}}%
0\\0\\0
\end{tabular}\endgroup%
}}\!\right]$}%
\EasyButWeakLineBreak%
{$\left[\!\llap{\phantom{%
\begingroup \smaller\smaller\smaller\begin{tabular}{@{}c@{}}%
0\\0\\0
\end{tabular}\endgroup%
}}\right.$}%
\begingroup \smaller\smaller\smaller\begin{tabular}{@{}c@{}}%
161\\40203\\29946
\end{tabular}\endgroup%
\HardButStrongLineBreak\kern3pt%
\begingroup \smaller\smaller\smaller\begin{tabular}{@{}c@{}}%
343\\85650\\63800
\end{tabular}\endgroup%
\HardButStrongLineBreak\kern3pt%
\begingroup \smaller\smaller\smaller\begin{tabular}{@{}c@{}}%
-137\\-34210\\-25482
\end{tabular}\endgroup%
\HardButStrongLineBreak\kern3pt%
\begingroup \smaller\smaller\smaller\begin{tabular}{@{}c@{}}%
-2656\\-663225\\-494025
\end{tabular}\endgroup%
\HardButStrongLineBreak\kern3pt%
\begingroup \smaller\smaller\smaller\begin{tabular}{@{}c@{}}%
-168\\-41951\\-31249
\end{tabular}\endgroup%
{$\left.\llap{\phantom{%
\begingroup \smaller\smaller\smaller\begin{tabular}{@{}c@{}}%
0\\0\\0
\end{tabular}\endgroup%
}}\!\right]$}%

\medskip%
%
\leavevmode\llap{}%
$W_{187}$%
\qquad\llap{4} lattices, $\chi=8$%
\hfill%
$22232$%
\nopagebreak\smallskip\hrule\nopagebreak\medskip%
%
%
\leavevmode%
${L_{187.1}}$%
{} : {$1\above{1pt}{1pt}{-2}{{\rm II}}16\above{1pt}{1pt}{1}{1}{\cdot}1\above{1pt}{1pt}{-2}{}5\above{1pt}{1pt}{-}{}$}\EasyButWeakLineBreak%
{${16}\above{1pt}{1pt}{r}{2}{10}\above{1pt}{1pt}{b}{2}{16}\above{1pt}{1pt}{b}{2}{2}\above{1pt}{1pt}{-}{3}{2}\above{1pt}{1pt}{l}{2}$}%
\nopagebreak\par%
\nopagebreak\par\leavevmode%
{$\left[\!\llap{\phantom{%
\begingroup \smaller\smaller\smaller\begin{tabular}{@{}c@{}}%
0\\0\\0
\end{tabular}\endgroup%
}}\right.$}%
\begingroup \smaller\smaller\smaller\begin{tabular}{@{}c@{}}%
-57840\\-880\\1040
\end{tabular}\endgroup%
\kern3pt%
\begingroup \smaller\smaller\smaller\begin{tabular}{@{}c@{}}%
-880\\-2\\13
\end{tabular}\endgroup%
\kern3pt%
\begingroup \smaller\smaller\smaller\begin{tabular}{@{}c@{}}%
1040\\13\\-18
\end{tabular}\endgroup%
{$\left.\llap{\phantom{%
\begingroup \smaller\smaller\smaller\begin{tabular}{@{}c@{}}%
0\\0\\0
\end{tabular}\endgroup%
}}\!\right]$}%
\EasyButWeakLineBreak%
{$\left[\!\llap{\phantom{%
\begingroup \smaller\smaller\smaller\begin{tabular}{@{}c@{}}%
0\\0\\0
\end{tabular}\endgroup%
}}\right.$}%
\begingroup \smaller\smaller\smaller\begin{tabular}{@{}c@{}}%
-15\\-256\\-1056
\end{tabular}\endgroup%
\HardButStrongLineBreak\kern3pt%
\begingroup \smaller\smaller\smaller\begin{tabular}{@{}c@{}}%
-1\\-15\\-70
\end{tabular}\endgroup%
\HardButStrongLineBreak\kern3pt%
\begingroup \smaller\smaller\smaller\begin{tabular}{@{}c@{}}%
5\\88\\352
\end{tabular}\endgroup%
\HardButStrongLineBreak\kern3pt%
\begingroup \smaller\smaller\smaller\begin{tabular}{@{}c@{}}%
1\\17\\70
\end{tabular}\endgroup%
\HardButStrongLineBreak\kern3pt%
\begingroup \smaller\smaller\smaller\begin{tabular}{@{}c@{}}%
-2\\-35\\-141
\end{tabular}\endgroup%
{$\left.\llap{\phantom{%
\begingroup \smaller\smaller\smaller\begin{tabular}{@{}c@{}}%
0\\0\\0
\end{tabular}\endgroup%
}}\!\right]$}%
%
%
%
%
%
%
%
%
%
%
%
%
%
%

\medskip%
%
\leavevmode\llap{}%
$W_{188}$%
\qquad\llap{4} lattices, $\chi=18$%
\hfill%
$42\slashinfty24|\rtimes D_{2}$%
\nopagebreak\smallskip\hrule\nopagebreak\medskip%
%
%
\leavevmode%
${L_{188.1}}$%
{} : {$1\above{1pt}{1pt}{2}{2}16\above{1pt}{1pt}{-}{3}{\cdot}1\above{1pt}{1pt}{2}{}5\above{1pt}{1pt}{1}{}$}\EasyButWeakLineBreak%
{${2}\above{1pt}{1pt}{*}{4}{4}\above{1pt}{1pt}{l}{2}{5}\above{1pt}{1pt}{8,3}{\infty}{20}\above{1pt}{1pt}{l}{2}{1}\above{1pt}{1pt}{}{4}$}%
\nopagebreak\par%
\nopagebreak\par\leavevmode%
{$\left[\!\llap{\phantom{%
\begingroup \smaller\smaller\smaller\begin{tabular}{@{}c@{}}%
0\\0\\0
\end{tabular}\endgroup%
}}\right.$}%
\begingroup \smaller\smaller\smaller\begin{tabular}{@{}c@{}}%
-16720\\480\\240
\end{tabular}\endgroup%
\kern3pt%
\begingroup \smaller\smaller\smaller\begin{tabular}{@{}c@{}}%
480\\-11\\-8
\end{tabular}\endgroup%
\kern3pt%
\begingroup \smaller\smaller\smaller\begin{tabular}{@{}c@{}}%
240\\-8\\-3
\end{tabular}\endgroup%
{$\left.\llap{\phantom{%
\begingroup \smaller\smaller\smaller\begin{tabular}{@{}c@{}}%
0\\0\\0
\end{tabular}\endgroup%
}}\!\right]$}%
\EasyButWeakLineBreak%
{$\left[\!\llap{\phantom{%
\begingroup \smaller\smaller\smaller\begin{tabular}{@{}c@{}}%
0\\0\\0
\end{tabular}\endgroup%
}}\right.$}%
\begingroup \smaller\smaller\smaller\begin{tabular}{@{}c@{}}%
1\\15\\39
\end{tabular}\endgroup%
\HardButStrongLineBreak\kern3pt%
\begingroup \smaller\smaller\smaller\begin{tabular}{@{}c@{}}%
-1\\-16\\-38
\end{tabular}\endgroup%
\HardButStrongLineBreak\kern3pt%
\begingroup \smaller\smaller\smaller\begin{tabular}{@{}c@{}}%
-1\\-15\\-40
\end{tabular}\endgroup%
\HardButStrongLineBreak\kern3pt%
\begingroup \smaller\smaller\smaller\begin{tabular}{@{}c@{}}%
7\\110\\260
\end{tabular}\endgroup%
\HardButStrongLineBreak\kern3pt%
\begingroup \smaller\smaller\smaller\begin{tabular}{@{}c@{}}%
4\\62\\151
\end{tabular}\endgroup%
{$\left.\llap{\phantom{%
\begingroup \smaller\smaller\smaller\begin{tabular}{@{}c@{}}%
0\\0\\0
\end{tabular}\endgroup%
}}\!\right]$}%
%
%
%
%
%
%
%
%
%
%
%
%
%
%

\medskip%
%
\leavevmode\llap{}%
$W_{189}$%
\qquad\llap{36} lattices, $\chi=36$%
\hfill%
$2\slashtwo2\slashinfty2\slashtwo2\slashinfty\rtimes D_{4}$%
\nopagebreak\smallskip\hrule\nopagebreak\medskip%
%
%
\leavevmode%
${L_{189.1}}$%
{} : {$1\above{1pt}{1pt}{2}{0}16\above{1pt}{1pt}{-}{5}{\cdot}1\above{1pt}{1pt}{1}{}5\above{1pt}{1pt}{1}{}25\above{1pt}{1pt}{1}{}$}\spacer%
\instructions{5}%
\EasyButWeakLineBreak%
{${20}\above{1pt}{1pt}{s}{2}{400}\above{1pt}{1pt}{s}{2}{4}\above{1pt}{1pt}{*}{2}{80}\above{1pt}{1pt}{5,4}{\infty b}{80}\above{1pt}{1pt}{}{2}{1}\above{1pt}{1pt}{r}{2}{400}\above{1pt}{1pt}{l}{2}{5}\above{1pt}{1pt}{40,21}{\infty}$}%
\nopagebreak\par%
shares genus with 5-dual\nopagebreak\par%
\nopagebreak\par\leavevmode%
{$\left[\!\llap{\phantom{%
\begingroup \smaller\smaller\smaller\begin{tabular}{@{}c@{}}%
0\\0\\0
\end{tabular}\endgroup%
}}\right.$}%
\begingroup \smaller\smaller\smaller\begin{tabular}{@{}c@{}}%
-2375600\\-43600\\-45200
\end{tabular}\endgroup%
\kern3pt%
\begingroup \smaller\smaller\smaller\begin{tabular}{@{}c@{}}%
-43600\\-795\\-825
\end{tabular}\endgroup%
\kern3pt%
\begingroup \smaller\smaller\smaller\begin{tabular}{@{}c@{}}%
-45200\\-825\\-856
\end{tabular}\endgroup%
{$\left.\llap{\phantom{%
\begingroup \smaller\smaller\smaller\begin{tabular}{@{}c@{}}%
0\\0\\0
\end{tabular}\endgroup%
}}\!\right]$}%
\EasyButWeakLineBreak%
{$\left[\!\llap{\phantom{%
\begingroup \smaller\smaller\smaller\begin{tabular}{@{}c@{}}%
0\\0\\0
\end{tabular}\endgroup%
}}\right.$}%
\begingroup \smaller\smaller\smaller\begin{tabular}{@{}c@{}}%
-1\\-298\\340
\end{tabular}\endgroup%
\HardButStrongLineBreak\kern3pt%
\begingroup \smaller\smaller\smaller\begin{tabular}{@{}c@{}}%
1\\360\\-400
\end{tabular}\endgroup%
\HardButStrongLineBreak\kern3pt%
\begingroup \smaller\smaller\smaller\begin{tabular}{@{}c@{}}%
1\\302\\-344
\end{tabular}\endgroup%
\HardButStrongLineBreak\kern3pt%
\begingroup \smaller\smaller\smaller\begin{tabular}{@{}c@{}}%
7\\2064\\-2360
\end{tabular}\endgroup%
\HardButStrongLineBreak\kern3pt%
\begingroup \smaller\smaller\smaller\begin{tabular}{@{}c@{}}%
9\\2576\\-2960
\end{tabular}\endgroup%
\HardButStrongLineBreak\kern3pt%
\begingroup \smaller\smaller\smaller\begin{tabular}{@{}c@{}}%
1\\279\\-322
\end{tabular}\endgroup%
\HardButStrongLineBreak\kern3pt%
\begingroup \smaller\smaller\smaller\begin{tabular}{@{}c@{}}%
11\\2920\\-3400
\end{tabular}\endgroup%
\HardButStrongLineBreak\kern3pt%
\begingroup \smaller\smaller\smaller\begin{tabular}{@{}c@{}}%
0\\-21\\20
\end{tabular}\endgroup%
{$\left.\llap{\phantom{%
\begingroup \smaller\smaller\smaller\begin{tabular}{@{}c@{}}%
0\\0\\0
\end{tabular}\endgroup%
}}\!\right]$}%
%
%
\hbox{}\par\smallskip%
%
%
\leavevmode%
${L_{189.2}}$%
{} : {$1\above{1pt}{1pt}{-2}{4}16\above{1pt}{1pt}{1}{1}{\cdot}1\above{1pt}{1pt}{1}{}5\above{1pt}{1pt}{1}{}25\above{1pt}{1pt}{1}{}$}\spacer%
\instructions{5}%
\EasyButWeakLineBreak%
{${20}\above{1pt}{1pt}{*}{2}{400}\above{1pt}{1pt}{*}{2}{4}\above{1pt}{1pt}{s}{2}{80}\above{1pt}{1pt}{20,9}{\infty z}{80}\above{1pt}{1pt}{l}{2}{1}\above{1pt}{1pt}{}{2}{400}\above{1pt}{1pt}{}{2}{5}\above{1pt}{1pt}{40,1}{\infty}$}%
\nopagebreak\par%
shares genus with 5-dual\nopagebreak\par%
\nopagebreak\par\leavevmode%
{$\left[\!\llap{\phantom{%
\begingroup \smaller\smaller\smaller\begin{tabular}{@{}c@{}}%
0\\0\\0
\end{tabular}\endgroup%
}}\right.$}%
\begingroup \smaller\smaller\smaller\begin{tabular}{@{}c@{}}%
-13903600\\108800\\40400
\end{tabular}\endgroup%
\kern3pt%
\begingroup \smaller\smaller\smaller\begin{tabular}{@{}c@{}}%
108800\\-795\\-325
\end{tabular}\endgroup%
\kern3pt%
\begingroup \smaller\smaller\smaller\begin{tabular}{@{}c@{}}%
40400\\-325\\-116
\end{tabular}\endgroup%
{$\left.\llap{\phantom{%
\begingroup \smaller\smaller\smaller\begin{tabular}{@{}c@{}}%
0\\0\\0
\end{tabular}\endgroup%
}}\!\right]$}%
\EasyButWeakLineBreak%
{$\left[\!\llap{\phantom{%
\begingroup \smaller\smaller\smaller\begin{tabular}{@{}c@{}}%
0\\0\\0
\end{tabular}\endgroup%
}}\right.$}%
\begingroup \smaller\smaller\smaller\begin{tabular}{@{}c@{}}%
-11\\-418\\-2660
\end{tabular}\endgroup%
\HardButStrongLineBreak\kern3pt%
\begingroup \smaller\smaller\smaller\begin{tabular}{@{}c@{}}%
-19\\-720\\-4600
\end{tabular}\endgroup%
\HardButStrongLineBreak\kern3pt%
\begingroup \smaller\smaller\smaller\begin{tabular}{@{}c@{}}%
9\\342\\2176
\end{tabular}\endgroup%
\HardButStrongLineBreak\kern3pt%
\begingroup \smaller\smaller\smaller\begin{tabular}{@{}c@{}}%
87\\3304\\21040
\end{tabular}\endgroup%
\HardButStrongLineBreak\kern3pt%
\begingroup \smaller\smaller\smaller\begin{tabular}{@{}c@{}}%
149\\5656\\36040
\end{tabular}\endgroup%
\HardButStrongLineBreak\kern3pt%
\begingroup \smaller\smaller\smaller\begin{tabular}{@{}c@{}}%
20\\759\\4838
\end{tabular}\endgroup%
\HardButStrongLineBreak\kern3pt%
\begingroup \smaller\smaller\smaller\begin{tabular}{@{}c@{}}%
291\\11040\\70400
\end{tabular}\endgroup%
\HardButStrongLineBreak\kern3pt%
\begingroup \smaller\smaller\smaller\begin{tabular}{@{}c@{}}%
10\\379\\2420
\end{tabular}\endgroup%
{$\left.\llap{\phantom{%
\begingroup \smaller\smaller\smaller\begin{tabular}{@{}c@{}}%
0\\0\\0
\end{tabular}\endgroup%
}}\!\right]$}%
%
%
\hbox{}\par\smallskip%
%
%
\leavevmode%
${L_{189.3}}$%
{} : {$1\above{1pt}{1pt}{-2}{6}16\above{1pt}{1pt}{1}{7}{\cdot}1\above{1pt}{1pt}{1}{}5\above{1pt}{1pt}{1}{}25\above{1pt}{1pt}{1}{}$}\spacer%
\instructions{5}%
\EasyButWeakLineBreak%
{${20}\above{1pt}{1pt}{*}{2}{100}\above{1pt}{1pt}{l}{2}{1}\above{1pt}{1pt}{}{2}{5}\above{1pt}{1pt}{40,39}{\infty}{20}\above{1pt}{1pt}{*}{2}{4}\above{1pt}{1pt}{l}{2}{25}\above{1pt}{1pt}{}{2}{5}\above{1pt}{1pt}{40,31}{\infty}$}%
\nopagebreak\par%
\nopagebreak\par\leavevmode%
{$\left[\!\llap{\phantom{%
\begingroup \smaller\smaller\smaller\begin{tabular}{@{}c@{}}%
0\\0\\0
\end{tabular}\endgroup%
}}\right.$}%
\begingroup \smaller\smaller\smaller\begin{tabular}{@{}c@{}}%
-1072400\\288400\\-71200
\end{tabular}\endgroup%
\kern3pt%
\begingroup \smaller\smaller\smaller\begin{tabular}{@{}c@{}}%
288400\\-75055\\20390
\end{tabular}\endgroup%
\kern3pt%
\begingroup \smaller\smaller\smaller\begin{tabular}{@{}c@{}}%
-71200\\20390\\-4111
\end{tabular}\endgroup%
{$\left.\llap{\phantom{%
\begingroup \smaller\smaller\smaller\begin{tabular}{@{}c@{}}%
0\\0\\0
\end{tabular}\endgroup%
}}\!\right]$}%
\EasyButWeakLineBreak%
{$\left[\!\llap{\phantom{%
\begingroup \smaller\smaller\smaller\begin{tabular}{@{}c@{}}%
0\\0\\0
\end{tabular}\endgroup%
}}\right.$}%
\begingroup \smaller\smaller\smaller\begin{tabular}{@{}c@{}}%
957\\2376\\-4790
\end{tabular}\endgroup%
\HardButStrongLineBreak\kern3pt%
\begingroup \smaller\smaller\smaller\begin{tabular}{@{}c@{}}%
1019\\2530\\-5100
\end{tabular}\endgroup%
\HardButStrongLineBreak\kern3pt%
\begingroup \smaller\smaller\smaller\begin{tabular}{@{}c@{}}%
-377\\-936\\1887
\end{tabular}\endgroup%
\HardButStrongLineBreak\kern3pt%
\begingroup \smaller\smaller\smaller\begin{tabular}{@{}c@{}}%
-1916\\-4757\\9590
\end{tabular}\endgroup%
\HardButStrongLineBreak\kern3pt%
\begingroup \smaller\smaller\smaller\begin{tabular}{@{}c@{}}%
-6769\\-16806\\33880
\end{tabular}\endgroup%
\HardButStrongLineBreak\kern3pt%
\begingroup \smaller\smaller\smaller\begin{tabular}{@{}c@{}}%
-3691\\-9164\\18474
\end{tabular}\endgroup%
\HardButStrongLineBreak\kern3pt%
\begingroup \smaller\smaller\smaller\begin{tabular}{@{}c@{}}%
-6833\\-16965\\34200
\end{tabular}\endgroup%
\HardButStrongLineBreak\kern3pt%
\begingroup \smaller\smaller\smaller\begin{tabular}{@{}c@{}}%
-990\\-2458\\4955
\end{tabular}\endgroup%
{$\left.\llap{\phantom{%
\begingroup \smaller\smaller\smaller\begin{tabular}{@{}c@{}}%
0\\0\\0
\end{tabular}\endgroup%
}}\!\right]$}%
%
%
\hbox{}\par\smallskip%
%
%
\leavevmode%
${L_{189.4}}$%
{} : {$1\above{1pt}{1pt}{2}{2}16\above{1pt}{1pt}{-}{3}{\cdot}1\above{1pt}{1pt}{1}{}5\above{1pt}{1pt}{1}{}25\above{1pt}{1pt}{1}{}$}\EasyButWeakLineBreak%
{${20}\above{1pt}{1pt}{l}{2}{25}\above{1pt}{1pt}{}{2}{1}\above{1pt}{1pt}{r}{2}{20}\above{1pt}{1pt}{40,39}{\infty z}{5}\above{1pt}{1pt}{r}{2}{4}\above{1pt}{1pt}{*}{2}{100}\above{1pt}{1pt}{l}{2}{5}\above{1pt}{1pt}{40,11}{\infty}$}%
\nopagebreak\par%
\nopagebreak\par\leavevmode%
{$\left[\!\llap{\phantom{%
\begingroup \smaller\smaller\smaller\begin{tabular}{@{}c@{}}%
0\\0\\0
\end{tabular}\endgroup%
}}\right.$}%
\begingroup \smaller\smaller\smaller\begin{tabular}{@{}c@{}}%
-475600\\-47600\\-44800
\end{tabular}\endgroup%
\kern3pt%
\begingroup \smaller\smaller\smaller\begin{tabular}{@{}c@{}}%
-47600\\-4670\\-4585
\end{tabular}\endgroup%
\kern3pt%
\begingroup \smaller\smaller\smaller\begin{tabular}{@{}c@{}}%
-44800\\-4585\\-4111
\end{tabular}\endgroup%
{$\left.\llap{\phantom{%
\begingroup \smaller\smaller\smaller\begin{tabular}{@{}c@{}}%
0\\0\\0
\end{tabular}\endgroup%
}}\!\right]$}%
\EasyButWeakLineBreak%
{$\left[\!\llap{\phantom{%
\begingroup \smaller\smaller\smaller\begin{tabular}{@{}c@{}}%
0\\0\\0
\end{tabular}\endgroup%
}}\right.$}%
\begingroup \smaller\smaller\smaller\begin{tabular}{@{}c@{}}%
-99\\528\\490
\end{tabular}\endgroup%
\HardButStrongLineBreak\kern3pt%
\begingroup \smaller\smaller\smaller\begin{tabular}{@{}c@{}}%
-106\\565\\525
\end{tabular}\endgroup%
\HardButStrongLineBreak\kern3pt%
\begingroup \smaller\smaller\smaller\begin{tabular}{@{}c@{}}%
39\\-208\\-193
\end{tabular}\endgroup%
\HardButStrongLineBreak\kern3pt%
\begingroup \smaller\smaller\smaller\begin{tabular}{@{}c@{}}%
503\\-2682\\-2490
\end{tabular}\endgroup%
\HardButStrongLineBreak\kern3pt%
\begingroup \smaller\smaller\smaller\begin{tabular}{@{}c@{}}%
510\\-2719\\-2525
\end{tabular}\endgroup%
\HardButStrongLineBreak\kern3pt%
\begingroup \smaller\smaller\smaller\begin{tabular}{@{}c@{}}%
595\\-3172\\-2946
\end{tabular}\endgroup%
\HardButStrongLineBreak\kern3pt%
\begingroup \smaller\smaller\smaller\begin{tabular}{@{}c@{}}%
2373\\-12650\\-11750
\end{tabular}\endgroup%
\HardButStrongLineBreak\kern3pt%
\begingroup \smaller\smaller\smaller\begin{tabular}{@{}c@{}}%
209\\-1114\\-1035
\end{tabular}\endgroup%
{$\left.\llap{\phantom{%
\begingroup \smaller\smaller\smaller\begin{tabular}{@{}c@{}}%
0\\0\\0
\end{tabular}\endgroup%
}}\!\right]$}%
%
%
\hbox{}\par\smallskip%
%
%
\leavevmode%
${L_{189.5}}$%
{} : {$1\above{1pt}{1pt}{-}{5}4\above{1pt}{1pt}{1}{7}16\above{1pt}{1pt}{1}{1}{\cdot}1\above{1pt}{1pt}{1}{}5\above{1pt}{1pt}{1}{}25\above{1pt}{1pt}{1}{}$}\spacer%
\instructions{5}%
\EasyButWeakLineBreak%
{${20}\above{1pt}{1pt}{s}{2}{400}\above{1pt}{1pt}{s}{2}{4}\above{1pt}{1pt}{s}{2}{80}\above{1pt}{1pt}{10,9}{\infty b}{80}\above{1pt}{1pt}{}{2}{1}\above{1pt}{1pt}{}{2}{400}\above{1pt}{1pt}{}{2}{5}\above{1pt}{1pt}{20,1}{\infty}$}%
\nopagebreak\par%
shares genus with 2-dual${}\iso{}$5-dual; isometric to own %
2.5-dual\nopagebreak\par%
\nopagebreak\par\leavevmode%
{$\left[\!\llap{\phantom{%
\begingroup \smaller\smaller\smaller\begin{tabular}{@{}c@{}}%
0\\0\\0
\end{tabular}\endgroup%
}}\right.$}%
\begingroup \smaller\smaller\smaller\begin{tabular}{@{}c@{}}%
-1775600\\164400\\40400
\end{tabular}\endgroup%
\kern3pt%
\begingroup \smaller\smaller\smaller\begin{tabular}{@{}c@{}}%
164400\\-15220\\-3740
\end{tabular}\endgroup%
\kern3pt%
\begingroup \smaller\smaller\smaller\begin{tabular}{@{}c@{}}%
40400\\-3740\\-919
\end{tabular}\endgroup%
{$\left.\llap{\phantom{%
\begingroup \smaller\smaller\smaller\begin{tabular}{@{}c@{}}%
0\\0\\0
\end{tabular}\endgroup%
}}\!\right]$}%
\EasyButWeakLineBreak%
{$\left[\!\llap{\phantom{%
\begingroup \smaller\smaller\smaller\begin{tabular}{@{}c@{}}%
0\\0\\0
\end{tabular}\endgroup%
}}\right.$}%
\begingroup \smaller\smaller\smaller\begin{tabular}{@{}c@{}}%
1\\28\\-70
\end{tabular}\endgroup%
\HardButStrongLineBreak\kern3pt%
\begingroup \smaller\smaller\smaller\begin{tabular}{@{}c@{}}%
-1\\-60\\200
\end{tabular}\endgroup%
\HardButStrongLineBreak\kern3pt%
\begingroup \smaller\smaller\smaller\begin{tabular}{@{}c@{}}%
-1\\-30\\78
\end{tabular}\endgroup%
\HardButStrongLineBreak\kern3pt%
\begingroup \smaller\smaller\smaller\begin{tabular}{@{}c@{}}%
-7\\-184\\440
\end{tabular}\endgroup%
\HardButStrongLineBreak\kern3pt%
\begingroup \smaller\smaller\smaller\begin{tabular}{@{}c@{}}%
-9\\-196\\400
\end{tabular}\endgroup%
\HardButStrongLineBreak\kern3pt%
\begingroup \smaller\smaller\smaller\begin{tabular}{@{}c@{}}%
-1\\-18\\29
\end{tabular}\endgroup%
\HardButStrongLineBreak\kern3pt%
\begingroup \smaller\smaller\smaller\begin{tabular}{@{}c@{}}%
-11\\-120\\0
\end{tabular}\endgroup%
\HardButStrongLineBreak\kern3pt%
\begingroup \smaller\smaller\smaller\begin{tabular}{@{}c@{}}%
0\\11\\-45
\end{tabular}\endgroup%
{$\left.\llap{\phantom{%
\begingroup \smaller\smaller\smaller\begin{tabular}{@{}c@{}}%
0\\0\\0
\end{tabular}\endgroup%
}}\!\right]$}%
%
%
\hbox{}\par\smallskip%
%
%
\leavevmode%
${L_{189.6}}$%
{} : {$1\above{1pt}{1pt}{-}{5}4\above{1pt}{1pt}{1}{1}16\above{1pt}{1pt}{1}{7}{\cdot}1\above{1pt}{1pt}{1}{}5\above{1pt}{1pt}{1}{}25\above{1pt}{1pt}{1}{}$}\spacer%
\instructions{5}%
\EasyButWeakLineBreak%
{${20}\above{1pt}{1pt}{l}{2}{100}\above{1pt}{1pt}{}{2}{1}\above{1pt}{1pt}{}{2}{20}\above{1pt}{1pt}{40,39}{\infty}{20}\above{1pt}{1pt}{r}{2}{4}\above{1pt}{1pt}{l}{2}{100}\above{1pt}{1pt}{}{2}{5}\above{1pt}{1pt}{20,11}{\infty}$}%
\nopagebreak\par%
shares genus with 5-dual\nopagebreak\par%
\nopagebreak\par\leavevmode%
{$\left[\!\llap{\phantom{%
\begingroup \smaller\smaller\smaller\begin{tabular}{@{}c@{}}%
0\\0\\0
\end{tabular}\endgroup%
}}\right.$}%
\begingroup \smaller\smaller\smaller\begin{tabular}{@{}c@{}}%
-456400\\-268000\\-44800
\end{tabular}\endgroup%
\kern3pt%
\begingroup \smaller\smaller\smaller\begin{tabular}{@{}c@{}}%
-268000\\-152620\\-25140
\end{tabular}\endgroup%
\kern3pt%
\begingroup \smaller\smaller\smaller\begin{tabular}{@{}c@{}}%
-44800\\-25140\\-4111
\end{tabular}\endgroup%
{$\left.\llap{\phantom{%
\begingroup \smaller\smaller\smaller\begin{tabular}{@{}c@{}}%
0\\0\\0
\end{tabular}\endgroup%
}}\!\right]$}%
\EasyButWeakLineBreak%
{$\left[\!\llap{\phantom{%
\begingroup \smaller\smaller\smaller\begin{tabular}{@{}c@{}}%
0\\0\\0
\end{tabular}\endgroup%
}}\right.$}%
\begingroup \smaller\smaller\smaller\begin{tabular}{@{}c@{}}%
-99\\528\\-2150
\end{tabular}\endgroup%
\HardButStrongLineBreak\kern3pt%
\begingroup \smaller\smaller\smaller\begin{tabular}{@{}c@{}}%
-106\\565\\-2300
\end{tabular}\endgroup%
\HardButStrongLineBreak\kern3pt%
\begingroup \smaller\smaller\smaller\begin{tabular}{@{}c@{}}%
39\\-208\\847
\end{tabular}\endgroup%
\HardButStrongLineBreak\kern3pt%
\begingroup \smaller\smaller\smaller\begin{tabular}{@{}c@{}}%
397\\-2117\\8620
\end{tabular}\endgroup%
\HardButStrongLineBreak\kern3pt%
\begingroup \smaller\smaller\smaller\begin{tabular}{@{}c@{}}%
702\\-3743\\15240
\end{tabular}\endgroup%
\HardButStrongLineBreak\kern3pt%
\begingroup \smaller\smaller\smaller\begin{tabular}{@{}c@{}}%
383\\-2042\\8314
\end{tabular}\endgroup%
\HardButStrongLineBreak\kern3pt%
\begingroup \smaller\smaller\smaller\begin{tabular}{@{}c@{}}%
1419\\-7565\\30800
\end{tabular}\endgroup%
\HardButStrongLineBreak\kern3pt%
\begingroup \smaller\smaller\smaller\begin{tabular}{@{}c@{}}%
103\\-549\\2235
\end{tabular}\endgroup%
{$\left.\llap{\phantom{%
\begingroup \smaller\smaller\smaller\begin{tabular}{@{}c@{}}%
0\\0\\0
\end{tabular}\endgroup%
}}\!\right]$}%

\medskip%
%
\leavevmode\llap{}%
$W_{190}$%
\qquad\llap{6} lattices, $\chi=6$%
\hfill%
$6223$%
\nopagebreak\smallskip\hrule\nopagebreak\medskip%
%
%
\leavevmode%
${L_{190.1}}$%
{} : {$1\above{1pt}{1pt}{-2}{{\rm II}}4\above{1pt}{1pt}{-}{5}{\cdot}1\above{1pt}{1pt}{-}{}3\above{1pt}{1pt}{-}{}27\above{1pt}{1pt}{-}{}$}\spacer%
\instructions{2}%
\EasyButWeakLineBreak%
{${6}\above{1pt}{1pt}{}{6}{2}\above{1pt}{1pt}{s}{2}{54}\above{1pt}{1pt}{b}{2}{6}\above{1pt}{1pt}{-}{3}$}%
\nopagebreak\par%
\nopagebreak\par\leavevmode%
{$\left[\!\llap{\phantom{%
\begingroup \smaller\smaller\smaller\begin{tabular}{@{}c@{}}%
0\\0\\0
\end{tabular}\endgroup%
}}\right.$}%
\begingroup \smaller\smaller\smaller\begin{tabular}{@{}c@{}}%
-124524\\2052\\-29592
\end{tabular}\endgroup%
\kern3pt%
\begingroup \smaller\smaller\smaller\begin{tabular}{@{}c@{}}%
2052\\-30\\507
\end{tabular}\endgroup%
\kern3pt%
\begingroup \smaller\smaller\smaller\begin{tabular}{@{}c@{}}%
-29592\\507\\-6934
\end{tabular}\endgroup%
{$\left.\llap{\phantom{%
\begingroup \smaller\smaller\smaller\begin{tabular}{@{}c@{}}%
0\\0\\0
\end{tabular}\endgroup%
}}\!\right]$}%
\EasyButWeakLineBreak%
{$\left[\!\llap{\phantom{%
\begingroup \smaller\smaller\smaller\begin{tabular}{@{}c@{}}%
0\\0\\0
\end{tabular}\endgroup%
}}\right.$}%
\begingroup \smaller\smaller\smaller\begin{tabular}{@{}c@{}}%
-26\\-410\\81
\end{tabular}\endgroup%
\HardButStrongLineBreak\kern3pt%
\begingroup \smaller\smaller\smaller\begin{tabular}{@{}c@{}}%
-18\\-285\\56
\end{tabular}\endgroup%
\HardButStrongLineBreak\kern3pt%
\begingroup \smaller\smaller\smaller\begin{tabular}{@{}c@{}}%
52\\819\\-162
\end{tabular}\endgroup%
\HardButStrongLineBreak\kern3pt%
\begingroup \smaller\smaller\smaller\begin{tabular}{@{}c@{}}%
27\\427\\-84
\end{tabular}\endgroup%
{$\left.\llap{\phantom{%
\begingroup \smaller\smaller\smaller\begin{tabular}{@{}c@{}}%
0\\0\\0
\end{tabular}\endgroup%
}}\!\right]$}%

\medskip%
%
\leavevmode\llap{}%
$W_{191}$%
\qquad\llap{22} lattices, $\chi=36$%
\hfill%
$2\infty222\infty22\rtimes C_{2}$%
\nopagebreak\smallskip\hrule\nopagebreak\medskip%
%
%
\leavevmode%
${L_{191.1}}$%
{} : {$1\above{1pt}{1pt}{2}{{\rm II}}4\above{1pt}{1pt}{1}{1}{\cdot}1\above{1pt}{1pt}{1}{}3\above{1pt}{1pt}{1}{}27\above{1pt}{1pt}{-}{}$}\spacer%
\instructions{2}%
\EasyButWeakLineBreak%
{${12}\above{1pt}{1pt}{*}{2}{4}\above{1pt}{1pt}{3,2}{\infty b}{4}\above{1pt}{1pt}{r}{2}{54}\above{1pt}{1pt}{b}{2}$}\relax$\,(\times2)$%
\nopagebreak\par%
\nopagebreak\par\leavevmode%
{$\left[\!\llap{\phantom{%
\begingroup \smaller\smaller\smaller\begin{tabular}{@{}c@{}}%
0\\0\\0
\end{tabular}\endgroup%
}}\right.$}%
\begingroup \smaller\smaller\smaller\begin{tabular}{@{}c@{}}%
9739332\\-38988\\-110160
\end{tabular}\endgroup%
\kern3pt%
\begingroup \smaller\smaller\smaller\begin{tabular}{@{}c@{}}%
-38988\\156\\441
\end{tabular}\endgroup%
\kern3pt%
\begingroup \smaller\smaller\smaller\begin{tabular}{@{}c@{}}%
-110160\\441\\1246
\end{tabular}\endgroup%
{$\left.\llap{\phantom{%
\begingroup \smaller\smaller\smaller\begin{tabular}{@{}c@{}}%
0\\0\\0
\end{tabular}\endgroup%
}}\!\right]$}%
\hfil\penalty500%
{$\left[\!\llap{\phantom{%
\begingroup \smaller\smaller\smaller\begin{tabular}{@{}c@{}}%
0\\0\\0
\end{tabular}\endgroup%
}}\right.$}%
\begingroup \smaller\smaller\smaller\begin{tabular}{@{}c@{}}%
-9217\\-230400\\-732672
\end{tabular}\endgroup%
\kern3pt%
\begingroup \smaller\smaller\smaller\begin{tabular}{@{}c@{}}%
38\\949\\3021
\end{tabular}\endgroup%
\kern3pt%
\begingroup \smaller\smaller\smaller\begin{tabular}{@{}c@{}}%
104\\2600\\8267
\end{tabular}\endgroup%
{$\left.\llap{\phantom{%
\begingroup \smaller\smaller\smaller\begin{tabular}{@{}c@{}}%
0\\0\\0
\end{tabular}\endgroup%
}}\!\right]$}%
\EasyButWeakLineBreak%
{$\left[\!\llap{\phantom{%
\begingroup \smaller\smaller\smaller\begin{tabular}{@{}c@{}}%
0\\0\\0
\end{tabular}\endgroup%
}}\right.$}%
\begingroup \smaller\smaller\smaller\begin{tabular}{@{}c@{}}%
-1\\4\\-90
\end{tabular}\endgroup%
\HardButStrongLineBreak\kern3pt%
\begingroup \smaller\smaller\smaller\begin{tabular}{@{}c@{}}%
-1\\-18\\-82
\end{tabular}\endgroup%
\HardButStrongLineBreak\kern3pt%
\begingroup \smaller\smaller\smaller\begin{tabular}{@{}c@{}}%
3\\60\\244
\end{tabular}\endgroup%
\HardButStrongLineBreak\kern3pt%
\begingroup \smaller\smaller\smaller\begin{tabular}{@{}c@{}}%
25\\522\\2025
\end{tabular}\endgroup%
{$\left.\llap{\phantom{%
\begingroup \smaller\smaller\smaller\begin{tabular}{@{}c@{}}%
0\\0\\0
\end{tabular}\endgroup%
}}\!\right]$}%
%
%
\hbox{}\par\smallskip%
%
%
\leavevmode%
${L_{191.2}}$%
{} : {$1\above{1pt}{1pt}{2}{2}8\above{1pt}{1pt}{1}{7}{\cdot}1\above{1pt}{1pt}{-}{}3\above{1pt}{1pt}{-}{}27\above{1pt}{1pt}{1}{}$}\spacer%
\instructions{2}%
\EasyButWeakLineBreak%
{${24}\above{1pt}{1pt}{b}{2}{2}\above{1pt}{1pt}{12,5}{\infty a}{8}\above{1pt}{1pt}{s}{2}{108}\above{1pt}{1pt}{*}{2}$}\relax$\,(\times2)$%
\nopagebreak\par%
\nopagebreak\par\leavevmode%
{$\left[\!\llap{\phantom{%
\begingroup \smaller\smaller\smaller\begin{tabular}{@{}c@{}}%
0\\0\\0
\end{tabular}\endgroup%
}}\right.$}%
\begingroup \smaller\smaller\smaller\begin{tabular}{@{}c@{}}%
-4970376\\-835272\\32616
\end{tabular}\endgroup%
\kern3pt%
\begingroup \smaller\smaller\smaller\begin{tabular}{@{}c@{}}%
-835272\\-140367\\5481
\end{tabular}\endgroup%
\kern3pt%
\begingroup \smaller\smaller\smaller\begin{tabular}{@{}c@{}}%
32616\\5481\\-214
\end{tabular}\endgroup%
{$\left.\llap{\phantom{%
\begingroup \smaller\smaller\smaller\begin{tabular}{@{}c@{}}%
0\\0\\0
\end{tabular}\endgroup%
}}\!\right]$}%
\hfil\penalty500%
{$\left[\!\llap{\phantom{%
\begingroup \smaller\smaller\smaller\begin{tabular}{@{}c@{}}%
0\\0\\0
\end{tabular}\endgroup%
}}\right.$}%
\begingroup \smaller\smaller\smaller\begin{tabular}{@{}c@{}}%
-76321\\541872\\2243808
\end{tabular}\endgroup%
\kern3pt%
\begingroup \smaller\smaller\smaller\begin{tabular}{@{}c@{}}%
-12820\\91021\\376908
\end{tabular}\endgroup%
\kern3pt%
\begingroup \smaller\smaller\smaller\begin{tabular}{@{}c@{}}%
500\\-3550\\-14701
\end{tabular}\endgroup%
{$\left.\llap{\phantom{%
\begingroup \smaller\smaller\smaller\begin{tabular}{@{}c@{}}%
0\\0\\0
\end{tabular}\endgroup%
}}\!\right]$}%
\EasyButWeakLineBreak%
{$\left[\!\llap{\phantom{%
\begingroup \smaller\smaller\smaller\begin{tabular}{@{}c@{}}%
0\\0\\0
\end{tabular}\endgroup%
}}\right.$}%
\begingroup \smaller\smaller\smaller\begin{tabular}{@{}c@{}}%
-11\\80\\372
\end{tabular}\endgroup%
\HardButStrongLineBreak\kern3pt%
\begingroup \smaller\smaller\smaller\begin{tabular}{@{}c@{}}%
-14\\100\\427
\end{tabular}\endgroup%
\HardButStrongLineBreak\kern3pt%
\begingroup \smaller\smaller\smaller\begin{tabular}{@{}c@{}}%
-85\\604\\2512
\end{tabular}\endgroup%
\HardButStrongLineBreak\kern3pt%
\begingroup \smaller\smaller\smaller\begin{tabular}{@{}c@{}}%
-373\\2646\\10908
\end{tabular}\endgroup%
{$\left.\llap{\phantom{%
\begingroup \smaller\smaller\smaller\begin{tabular}{@{}c@{}}%
0\\0\\0
\end{tabular}\endgroup%
}}\!\right]$}%
%
%
\hbox{}\par\smallskip%
%
%
\leavevmode%
${L_{191.3}}$%
{} : {$1\above{1pt}{1pt}{-2}{2}8\above{1pt}{1pt}{-}{3}{\cdot}1\above{1pt}{1pt}{-}{}3\above{1pt}{1pt}{-}{}27\above{1pt}{1pt}{1}{}$}\spacer%
\instructions{m}%
\EasyButWeakLineBreak%
{${24}\above{1pt}{1pt}{r}{2}{2}\above{1pt}{1pt}{12,5}{\infty b}{8}\above{1pt}{1pt}{l}{2}{27}\above{1pt}{1pt}{}{2}$}\relax$\,(\times2)$%
\nopagebreak\par%
\nopagebreak\par\leavevmode%
{$\left[\!\llap{\phantom{%
\begingroup \smaller\smaller\smaller\begin{tabular}{@{}c@{}}%
0\\0\\0
\end{tabular}\endgroup%
}}\right.$}%
\begingroup \smaller\smaller\smaller\begin{tabular}{@{}c@{}}%
-758376\\-250128\\4752
\end{tabular}\endgroup%
\kern3pt%
\begingroup \smaller\smaller\smaller\begin{tabular}{@{}c@{}}%
-250128\\-82497\\1566
\end{tabular}\endgroup%
\kern3pt%
\begingroup \smaller\smaller\smaller\begin{tabular}{@{}c@{}}%
4752\\1566\\-25
\end{tabular}\endgroup%
{$\left.\llap{\phantom{%
\begingroup \smaller\smaller\smaller\begin{tabular}{@{}c@{}}%
0\\0\\0
\end{tabular}\endgroup%
}}\!\right]$}%
\hfil\penalty500%
{$\left[\!\llap{\phantom{%
\begingroup \smaller\smaller\smaller\begin{tabular}{@{}c@{}}%
0\\0\\0
\end{tabular}\endgroup%
}}\right.$}%
\begingroup \smaller\smaller\smaller\begin{tabular}{@{}c@{}}%
-338005\\1030248\\283392
\end{tabular}\endgroup%
\kern3pt%
\begingroup \smaller\smaller\smaller\begin{tabular}{@{}c@{}}%
-111523\\339925\\93504
\end{tabular}\endgroup%
\kern3pt%
\begingroup \smaller\smaller\smaller\begin{tabular}{@{}c@{}}%
2290\\-6980\\-1921
\end{tabular}\endgroup%
{$\left.\llap{\phantom{%
\begingroup \smaller\smaller\smaller\begin{tabular}{@{}c@{}}%
0\\0\\0
\end{tabular}\endgroup%
}}\!\right]$}%
\EasyButWeakLineBreak%
{$\left[\!\llap{\phantom{%
\begingroup \smaller\smaller\smaller\begin{tabular}{@{}c@{}}%
0\\0\\0
\end{tabular}\endgroup%
}}\right.$}%
\begingroup \smaller\smaller\smaller\begin{tabular}{@{}c@{}}%
-147\\448\\120
\end{tabular}\endgroup%
\HardButStrongLineBreak\kern3pt%
\begingroup \smaller\smaller\smaller\begin{tabular}{@{}c@{}}%
-167\\509\\139
\end{tabular}\endgroup%
\HardButStrongLineBreak\kern3pt%
\begingroup \smaller\smaller\smaller\begin{tabular}{@{}c@{}}%
-979\\2984\\820
\end{tabular}\endgroup%
\HardButStrongLineBreak\kern3pt%
\begingroup \smaller\smaller\smaller\begin{tabular}{@{}c@{}}%
-2123\\6471\\1782
\end{tabular}\endgroup%
{$\left.\llap{\phantom{%
\begingroup \smaller\smaller\smaller\begin{tabular}{@{}c@{}}%
0\\0\\0
\end{tabular}\endgroup%
}}\!\right]$}%

\medskip%
%
\leavevmode\llap{}%
$W_{192}$%
\qquad\llap{22} lattices, $\chi=36$%
\hfill%
$2\infty222\infty22\rtimes C_{2}$%
\nopagebreak\smallskip\hrule\nopagebreak\medskip%
%
%
\leavevmode%
${L_{192.1}}$%
{} : {$1\above{1pt}{1pt}{2}{{\rm II}}4\above{1pt}{1pt}{1}{1}{\cdot}1\above{1pt}{1pt}{1}{}3\above{1pt}{1pt}{-}{}27\above{1pt}{1pt}{1}{}$}\spacer%
\instructions{2}%
\EasyButWeakLineBreak%
{${108}\above{1pt}{1pt}{*}{2}{4}\above{1pt}{1pt}{3,1}{\infty a}{4}\above{1pt}{1pt}{r}{2}{6}\above{1pt}{1pt}{b}{2}$}\relax$\,(\times2)$%
\nopagebreak\par%
\nopagebreak\par\leavevmode%
{$\left[\!\llap{\phantom{%
\begingroup \smaller\smaller\smaller\begin{tabular}{@{}c@{}}%
0\\0\\0
\end{tabular}\endgroup%
}}\right.$}%
\begingroup \smaller\smaller\smaller\begin{tabular}{@{}c@{}}%
-1255068\\-800172\\84348
\end{tabular}\endgroup%
\kern3pt%
\begingroup \smaller\smaller\smaller\begin{tabular}{@{}c@{}}%
-800172\\-510132\\53769
\end{tabular}\endgroup%
\kern3pt%
\begingroup \smaller\smaller\smaller\begin{tabular}{@{}c@{}}%
84348\\53769\\-5666
\end{tabular}\endgroup%
{$\left.\llap{\phantom{%
\begingroup \smaller\smaller\smaller\begin{tabular}{@{}c@{}}%
0\\0\\0
\end{tabular}\endgroup%
}}\!\right]$}%
\hfil\penalty500%
{$\left[\!\llap{\phantom{%
\begingroup \smaller\smaller\smaller\begin{tabular}{@{}c@{}}%
0\\0\\0
\end{tabular}\endgroup%
}}\right.$}%
\begingroup \smaller\smaller\smaller\begin{tabular}{@{}c@{}}%
-31051\\67500\\178200
\end{tabular}\endgroup%
\kern3pt%
\begingroup \smaller\smaller\smaller\begin{tabular}{@{}c@{}}%
-19688\\42799\\112992
\end{tabular}\endgroup%
\kern3pt%
\begingroup \smaller\smaller\smaller\begin{tabular}{@{}c@{}}%
2047\\-4450\\-11749
\end{tabular}\endgroup%
{$\left.\llap{\phantom{%
\begingroup \smaller\smaller\smaller\begin{tabular}{@{}c@{}}%
0\\0\\0
\end{tabular}\endgroup%
}}\!\right]$}%
\EasyButWeakLineBreak%
{$\left[\!\llap{\phantom{%
\begingroup \smaller\smaller\smaller\begin{tabular}{@{}c@{}}%
0\\0\\0
\end{tabular}\endgroup%
}}\right.$}%
\begingroup \smaller\smaller\smaller\begin{tabular}{@{}c@{}}%
2819\\-6192\\-16794
\end{tabular}\endgroup%
\HardButStrongLineBreak\kern3pt%
\begingroup \smaller\smaller\smaller\begin{tabular}{@{}c@{}}%
243\\-534\\-1450
\end{tabular}\endgroup%
\HardButStrongLineBreak\kern3pt%
\begingroup \smaller\smaller\smaller\begin{tabular}{@{}c@{}}%
-11\\24\\64
\end{tabular}\endgroup%
\HardButStrongLineBreak\kern3pt%
\begingroup \smaller\smaller\smaller\begin{tabular}{@{}c@{}}%
-121\\266\\723
\end{tabular}\endgroup%
{$\left.\llap{\phantom{%
\begingroup \smaller\smaller\smaller\begin{tabular}{@{}c@{}}%
0\\0\\0
\end{tabular}\endgroup%
}}\!\right]$}%
%
%
\hbox{}\par\smallskip%
%
%
\leavevmode%
${L_{192.2}}$%
{} : {$1\above{1pt}{1pt}{2}{2}8\above{1pt}{1pt}{1}{7}{\cdot}1\above{1pt}{1pt}{-}{}3\above{1pt}{1pt}{1}{}27\above{1pt}{1pt}{-}{}$}\spacer%
\instructions{2}%
\EasyButWeakLineBreak%
{${216}\above{1pt}{1pt}{b}{2}{2}\above{1pt}{1pt}{12,1}{\infty b}{8}\above{1pt}{1pt}{s}{2}{12}\above{1pt}{1pt}{*}{2}$}\relax$\,(\times2)$%
\nopagebreak\par%
\nopagebreak\par\leavevmode%
{$\left[\!\llap{\phantom{%
\begingroup \smaller\smaller\smaller\begin{tabular}{@{}c@{}}%
0\\0\\0
\end{tabular}\endgroup%
}}\right.$}%
\begingroup \smaller\smaller\smaller\begin{tabular}{@{}c@{}}%
-603720\\-301536\\3672
\end{tabular}\endgroup%
\kern3pt%
\begingroup \smaller\smaller\smaller\begin{tabular}{@{}c@{}}%
-301536\\-150603\\1833
\end{tabular}\endgroup%
\kern3pt%
\begingroup \smaller\smaller\smaller\begin{tabular}{@{}c@{}}%
3672\\1833\\-22
\end{tabular}\endgroup%
{$\left.\llap{\phantom{%
\begingroup \smaller\smaller\smaller\begin{tabular}{@{}c@{}}%
0\\0\\0
\end{tabular}\endgroup%
}}\!\right]$}%
\hfil\penalty500%
{$\left[\!\llap{\phantom{%
\begingroup \smaller\smaller\smaller\begin{tabular}{@{}c@{}}%
0\\0\\0
\end{tabular}\endgroup%
}}\right.$}%
\begingroup \smaller\smaller\smaller\begin{tabular}{@{}c@{}}%
-25201\\52416\\160272
\end{tabular}\endgroup%
\kern3pt%
\begingroup \smaller\smaller\smaller\begin{tabular}{@{}c@{}}%
-12575\\26155\\79977
\end{tabular}\endgroup%
\kern3pt%
\begingroup \smaller\smaller\smaller\begin{tabular}{@{}c@{}}%
150\\-312\\-955
\end{tabular}\endgroup%
{$\left.\llap{\phantom{%
\begingroup \smaller\smaller\smaller\begin{tabular}{@{}c@{}}%
0\\0\\0
\end{tabular}\endgroup%
}}\!\right]$}%
\EasyButWeakLineBreak%
{$\left[\!\llap{\phantom{%
\begingroup \smaller\smaller\smaller\begin{tabular}{@{}c@{}}%
0\\0\\0
\end{tabular}\endgroup%
}}\right.$}%
\begingroup \smaller\smaller\smaller\begin{tabular}{@{}c@{}}%
-1453\\3024\\9396
\end{tabular}\endgroup%
\HardButStrongLineBreak\kern3pt%
\begingroup \smaller\smaller\smaller\begin{tabular}{@{}c@{}}%
-99\\206\\637
\end{tabular}\endgroup%
\HardButStrongLineBreak\kern3pt%
\begingroup \smaller\smaller\smaller\begin{tabular}{@{}c@{}}%
-125\\260\\796
\end{tabular}\endgroup%
\HardButStrongLineBreak\kern3pt%
\begingroup \smaller\smaller\smaller\begin{tabular}{@{}c@{}}%
-51\\106\\318
\end{tabular}\endgroup%
{$\left.\llap{\phantom{%
\begingroup \smaller\smaller\smaller\begin{tabular}{@{}c@{}}%
0\\0\\0
\end{tabular}\endgroup%
}}\!\right]$}%
%
%
\hbox{}\par\smallskip%
%
%
\leavevmode%
${L_{192.3}}$%
{} : {$1\above{1pt}{1pt}{-2}{2}8\above{1pt}{1pt}{-}{3}{\cdot}1\above{1pt}{1pt}{-}{}3\above{1pt}{1pt}{1}{}27\above{1pt}{1pt}{-}{}$}\spacer%
\instructions{m}%
\EasyButWeakLineBreak%
{${216}\above{1pt}{1pt}{r}{2}{2}\above{1pt}{1pt}{12,1}{\infty a}{8}\above{1pt}{1pt}{l}{2}{3}\above{1pt}{1pt}{}{2}$}\relax$\,(\times2)$%
\nopagebreak\par%
\nopagebreak\par\leavevmode%
{$\left[\!\llap{\phantom{%
\begingroup \smaller\smaller\smaller\begin{tabular}{@{}c@{}}%
0\\0\\0
\end{tabular}\endgroup%
}}\right.$}%
\begingroup \smaller\smaller\smaller\begin{tabular}{@{}c@{}}%
-1015848\\10368\\5184
\end{tabular}\endgroup%
\kern3pt%
\begingroup \smaller\smaller\smaller\begin{tabular}{@{}c@{}}%
10368\\-105\\-54
\end{tabular}\endgroup%
\kern3pt%
\begingroup \smaller\smaller\smaller\begin{tabular}{@{}c@{}}%
5184\\-54\\-25
\end{tabular}\endgroup%
{$\left.\llap{\phantom{%
\begingroup \smaller\smaller\smaller\begin{tabular}{@{}c@{}}%
0\\0\\0
\end{tabular}\endgroup%
}}\!\right]$}%
\hfil\penalty500%
{$\left[\!\llap{\phantom{%
\begingroup \smaller\smaller\smaller\begin{tabular}{@{}c@{}}%
0\\0\\0
\end{tabular}\endgroup%
}}\right.$}%
\begingroup \smaller\smaller\smaller\begin{tabular}{@{}c@{}}%
3167\\226512\\166320
\end{tabular}\endgroup%
\kern3pt%
\begingroup \smaller\smaller\smaller\begin{tabular}{@{}c@{}}%
-34\\-2432\\-1785
\end{tabular}\endgroup%
\kern3pt%
\begingroup \smaller\smaller\smaller\begin{tabular}{@{}c@{}}%
-14\\-1001\\-736
\end{tabular}\endgroup%
{$\left.\llap{\phantom{%
\begingroup \smaller\smaller\smaller\begin{tabular}{@{}c@{}}%
0\\0\\0
\end{tabular}\endgroup%
}}\!\right]$}%
\EasyButWeakLineBreak%
{$\left[\!\llap{\phantom{%
\begingroup \smaller\smaller\smaller\begin{tabular}{@{}c@{}}%
0\\0\\0
\end{tabular}\endgroup%
}}\right.$}%
\begingroup \smaller\smaller\smaller\begin{tabular}{@{}c@{}}%
59\\4248\\3024
\end{tabular}\endgroup%
\HardButStrongLineBreak\kern3pt%
\begingroup \smaller\smaller\smaller\begin{tabular}{@{}c@{}}%
5\\359\\259
\end{tabular}\endgroup%
\HardButStrongLineBreak\kern3pt%
\begingroup \smaller\smaller\smaller\begin{tabular}{@{}c@{}}%
9\\644\\472
\end{tabular}\endgroup%
\HardButStrongLineBreak\kern3pt%
\begingroup \smaller\smaller\smaller\begin{tabular}{@{}c@{}}%
3\\214\\159
\end{tabular}\endgroup%
{$\left.\llap{\phantom{%
\begingroup \smaller\smaller\smaller\begin{tabular}{@{}c@{}}%
0\\0\\0
\end{tabular}\endgroup%
}}\!\right]$}%

\medskip%
%
\leavevmode\llap{}%
$W_{193}$%
\qquad\llap{6} lattices, $\chi=12$%
\hfill%
$222222\rtimes C_{2}$%
\nopagebreak\smallskip\hrule\nopagebreak\medskip%
%
%
\leavevmode%
${L_{193.1}}$%
{} : {$1\above{1pt}{1pt}{-2}{{\rm II}}4\above{1pt}{1pt}{-}{5}{\cdot}1\above{1pt}{1pt}{-}{}3\above{1pt}{1pt}{1}{}27\above{1pt}{1pt}{1}{}$}\spacer%
\instructions{2}%
\EasyButWeakLineBreak%
{${12}\above{1pt}{1pt}{*}{2}{108}\above{1pt}{1pt}{b}{2}{2}\above{1pt}{1pt}{b}{2}$}\relax$\,(\times2)$%
\nopagebreak\par%
\nopagebreak\par\leavevmode%
{$\left[\!\llap{\phantom{%
\begingroup \smaller\smaller\smaller\begin{tabular}{@{}c@{}}%
0\\0\\0
\end{tabular}\endgroup%
}}\right.$}%
\begingroup \smaller\smaller\smaller\begin{tabular}{@{}c@{}}%
756\\324\\-108
\end{tabular}\endgroup%
\kern3pt%
\begingroup \smaller\smaller\smaller\begin{tabular}{@{}c@{}}%
324\\138\\-51
\end{tabular}\endgroup%
\kern3pt%
\begingroup \smaller\smaller\smaller\begin{tabular}{@{}c@{}}%
-108\\-51\\-10
\end{tabular}\endgroup%
{$\left.\llap{\phantom{%
\begingroup \smaller\smaller\smaller\begin{tabular}{@{}c@{}}%
0\\0\\0
\end{tabular}\endgroup%
}}\!\right]$}%
\hfil\penalty500%
{$\left[\!\llap{\phantom{%
\begingroup \smaller\smaller\smaller\begin{tabular}{@{}c@{}}%
0\\0\\0
\end{tabular}\endgroup%
}}\right.$}%
\begingroup \smaller\smaller\smaller\begin{tabular}{@{}c@{}}%
179\\-396\\108
\end{tabular}\endgroup%
\kern3pt%
\begingroup \smaller\smaller\smaller\begin{tabular}{@{}c@{}}%
85\\-188\\51
\end{tabular}\endgroup%
\kern3pt%
\begingroup \smaller\smaller\smaller\begin{tabular}{@{}c@{}}%
15\\-33\\8
\end{tabular}\endgroup%
{$\left.\llap{\phantom{%
\begingroup \smaller\smaller\smaller\begin{tabular}{@{}c@{}}%
0\\0\\0
\end{tabular}\endgroup%
}}\!\right]$}%
\EasyButWeakLineBreak%
{$\left[\!\llap{\phantom{%
\begingroup \smaller\smaller\smaller\begin{tabular}{@{}c@{}}%
0\\0\\0
\end{tabular}\endgroup%
}}\right.$}%
\begingroup \smaller\smaller\smaller\begin{tabular}{@{}c@{}}%
-1\\2\\0
\end{tabular}\endgroup%
\HardButStrongLineBreak\kern3pt%
\begingroup \smaller\smaller\smaller\begin{tabular}{@{}c@{}}%
-115\\252\\-54
\end{tabular}\endgroup%
\HardButStrongLineBreak\kern3pt%
\begingroup \smaller\smaller\smaller\begin{tabular}{@{}c@{}}%
-10\\22\\-5
\end{tabular}\endgroup%
{$\left.\llap{\phantom{%
\begingroup \smaller\smaller\smaller\begin{tabular}{@{}c@{}}%
0\\0\\0
\end{tabular}\endgroup%
}}\!\right]$}%

\medskip%
%
\leavevmode\llap{}%
$W_{194}$%
\qquad\llap{32} lattices, $\chi=18$%
\hfill%
$2\infty2222$%
\nopagebreak\smallskip\hrule\nopagebreak\medskip%
%
%
\leavevmode%
${L_{194.1}}$%
{} : {$1\above{1pt}{1pt}{-2}{{\rm II}}8\above{1pt}{1pt}{1}{1}{\cdot}1\above{1pt}{1pt}{2}{}9\above{1pt}{1pt}{-}{}{\cdot}1\above{1pt}{1pt}{-}{}5\above{1pt}{1pt}{-}{}25\above{1pt}{1pt}{-}{}$}\spacer%
\instructions{25,5,2*}%
\EasyButWeakLineBreak%
{${72}\above{1pt}{1pt}{r}{2}{10}\above{1pt}{1pt}{60,49}{\infty a}{40}\above{1pt}{1pt}{b}{2}{450}\above{1pt}{1pt}{l}{2}{8}\above{1pt}{1pt}{r}{2}{50}\above{1pt}{1pt}{l}{2}$}%
\nopagebreak\par%
shares genus with 5-dual\nopagebreak\par%
\nopagebreak\par\leavevmode%
{$\left[\!\llap{\phantom{%
\begingroup \smaller\smaller\smaller\begin{tabular}{@{}c@{}}%
0\\0\\0
\end{tabular}\endgroup%
}}\right.$}%
\begingroup \smaller\smaller\smaller\begin{tabular}{@{}c@{}}%
-8854200\\1479600\\145800
\end{tabular}\endgroup%
\kern3pt%
\begingroup \smaller\smaller\smaller\begin{tabular}{@{}c@{}}%
1479600\\-239810\\-20665
\end{tabular}\endgroup%
\kern3pt%
\begingroup \smaller\smaller\smaller\begin{tabular}{@{}c@{}}%
145800\\-20665\\-562
\end{tabular}\endgroup%
{$\left.\llap{\phantom{%
\begingroup \smaller\smaller\smaller\begin{tabular}{@{}c@{}}%
0\\0\\0
\end{tabular}\endgroup%
}}\!\right]$}%
\EasyButWeakLineBreak%
{$\left[\!\llap{\phantom{%
\begingroup \smaller\smaller\smaller\begin{tabular}{@{}c@{}}%
0\\0\\0
\end{tabular}\endgroup%
}}\right.$}%
\begingroup \smaller\smaller\smaller\begin{tabular}{@{}c@{}}%
-12367\\-92304\\185688
\end{tabular}\endgroup%
\HardButStrongLineBreak\kern3pt%
\begingroup \smaller\smaller\smaller\begin{tabular}{@{}c@{}}%
-1020\\-7613\\15315
\end{tabular}\endgroup%
\HardButStrongLineBreak\kern3pt%
\begingroup \smaller\smaller\smaller\begin{tabular}{@{}c@{}}%
1007\\7516\\-15120
\end{tabular}\endgroup%
\HardButStrongLineBreak\kern3pt%
\begingroup \smaller\smaller\smaller\begin{tabular}{@{}c@{}}%
2068\\15435\\-31050
\end{tabular}\endgroup%
\HardButStrongLineBreak\kern3pt%
\begingroup \smaller\smaller\smaller\begin{tabular}{@{}c@{}}%
-731\\-5456\\10976
\end{tabular}\endgroup%
\HardButStrongLineBreak\kern3pt%
\begingroup \smaller\smaller\smaller\begin{tabular}{@{}c@{}}%
-2689\\-20070\\40375
\end{tabular}\endgroup%
{$\left.\llap{\phantom{%
\begingroup \smaller\smaller\smaller\begin{tabular}{@{}c@{}}%
0\\0\\0
\end{tabular}\endgroup%
}}\!\right]$}%

\medskip%
%
\leavevmode\llap{}%
$W_{195}$%
\qquad\llap{32} lattices, $\chi=72$%
\hfill%
$2\infty2|2\infty2|2\infty2|2\infty2|\rtimes D_{4}$%
\nopagebreak\smallskip\hrule\nopagebreak\medskip%
%
%
\leavevmode%
${L_{195.1}}$%
{} : {$1\above{1pt}{1pt}{-2}{{\rm II}}8\above{1pt}{1pt}{1}{1}{\cdot}1\above{1pt}{1pt}{-2}{}9\above{1pt}{1pt}{1}{}{\cdot}1\above{1pt}{1pt}{-}{}5\above{1pt}{1pt}{-}{}25\above{1pt}{1pt}{-}{}$}\spacer%
\instructions{25,5,2*}%
\EasyButWeakLineBreak%
{${8}\above{1pt}{1pt}{r}{2}{10}\above{1pt}{1pt}{60,29}{\infty b}{40}\above{1pt}{1pt}{b}{2}{50}\above{1pt}{1pt}{b}{2}{360}\above{1pt}{1pt}{10,1}{\infty z}{90}\above{1pt}{1pt}{l}{2}$}\relax$\,(\times2)$%
\nopagebreak\par%
shares genus with 5-dual\nopagebreak\par%
\nopagebreak\par\leavevmode%
{$\left[\!\llap{\phantom{%
\begingroup \smaller\smaller\smaller
\endgroup%
}}\!\right]$}%

\medskip%
%
\leavevmode\llap{}%
$W_{196}$%
\qquad\llap{92} lattices, $\chi=36$%
\hfill%
$2222222222\rtimes C_{2}$%
\nopagebreak\smallskip\hrule\nopagebreak\medskip%
%
%
\leavevmode%
${L_{196.1}}$%
{} : {$1\above{1pt}{1pt}{-2}{4}8\above{1pt}{1pt}{1}{1}{\cdot}1\above{1pt}{1pt}{2}{}9\above{1pt}{1pt}{-}{}{\cdot}1\above{1pt}{1pt}{-2}{}5\above{1pt}{1pt}{1}{}$}\spacer%
\instructions{2}%
\EasyButWeakLineBreak%
{${45}\above{1pt}{1pt}{r}{2}{4}\above{1pt}{1pt}{*}{2}{72}\above{1pt}{1pt}{*}{2}{20}\above{1pt}{1pt}{*}{2}{8}\above{1pt}{1pt}{*}{2}{180}\above{1pt}{1pt}{l}{2}{1}\above{1pt}{1pt}{}{2}{72}\above{1pt}{1pt}{}{2}{5}\above{1pt}{1pt}{}{2}{8}\above{1pt}{1pt}{}{2}$}%
\nopagebreak\par%
\nopagebreak\par\leavevmode%
{$\left[\!\llap{\phantom{%
\begingroup \smaller\smaller\smaller
\endgroup%
}}\!\right]$}%
%
%
\hbox{}\par\smallskip%
%
%
\leavevmode%
${L_{196.2}}$%
{} : {$1\above{1pt}{1pt}{-2}{6}8\above{1pt}{1pt}{1}{7}{\cdot}1\above{1pt}{1pt}{2}{}9\above{1pt}{1pt}{-}{}{\cdot}1\above{1pt}{1pt}{-2}{}5\above{1pt}{1pt}{1}{}$}\spacer%
\instructions{m}%
\EasyButWeakLineBreak%
{${45}\above{1pt}{1pt}{}{2}{1}\above{1pt}{1pt}{r}{2}{72}\above{1pt}{1pt}{l}{2}{5}\above{1pt}{1pt}{r}{2}{8}\above{1pt}{1pt}{l}{2}$}\relax$\,(\times2)$%
\nopagebreak\par%
\nopagebreak\par\leavevmode%
{$\left[\!\llap{\phantom{%
\begingroup \smaller\smaller\smaller
\endgroup%
}}\!\right]$}%
%
%
\hbox{}\par\smallskip%
%
%
\leavevmode%
${L_{196.3}}$%
{} : {$1\above{1pt}{1pt}{2}{6}8\above{1pt}{1pt}{-}{3}{\cdot}1\above{1pt}{1pt}{2}{}9\above{1pt}{1pt}{-}{}{\cdot}1\above{1pt}{1pt}{-2}{}5\above{1pt}{1pt}{1}{}$}\EasyButWeakLineBreak%
{${180}\above{1pt}{1pt}{*}{2}{4}\above{1pt}{1pt}{s}{2}{72}\above{1pt}{1pt}{s}{2}{20}\above{1pt}{1pt}{s}{2}{8}\above{1pt}{1pt}{s}{2}$}\relax$\,(\times2)$%
\nopagebreak\par%
\nopagebreak\par\leavevmode%
{$\left[\!\llap{\phantom{%
\begingroup \smaller\smaller\smaller
\endgroup%
}}\!\right]$}%
%
%
\hbox{}\par\smallskip%
%
%
\leavevmode%
${L_{196.4}}$%
{} : {$[1\above{1pt}{1pt}{1}{}2\above{1pt}{1pt}{1}{}]\above{1pt}{1pt}{}{2}16\above{1pt}{1pt}{-}{3}{\cdot}1\above{1pt}{1pt}{2}{}9\above{1pt}{1pt}{-}{}{\cdot}1\above{1pt}{1pt}{-2}{}5\above{1pt}{1pt}{1}{}$}\spacer%
\instructions{2}%
\EasyButWeakLineBreak%
{${180}\above{1pt}{1pt}{*}{2}{16}\above{1pt}{1pt}{s}{2}{72}\above{1pt}{1pt}{*}{2}{80}\above{1pt}{1pt}{*}{2}{8}\above{1pt}{1pt}{*}{2}{720}\above{1pt}{1pt}{l}{2}{1}\above{1pt}{1pt}{}{2}{18}\above{1pt}{1pt}{r}{2}{20}\above{1pt}{1pt}{l}{2}{2}\above{1pt}{1pt}{r}{2}$}%
\nopagebreak\par%
\nopagebreak\par\leavevmode%
{$\left[\!\llap{\phantom{%
\begingroup \smaller\smaller\smaller
\endgroup%
}}\!\right]$}%
%
%
\hbox{}\par\smallskip%
%
%
\leavevmode%
${L_{196.5}}$%
{} : {$[1\above{1pt}{1pt}{-}{}2\above{1pt}{1pt}{1}{}]\above{1pt}{1pt}{}{6}16\above{1pt}{1pt}{1}{7}{\cdot}1\above{1pt}{1pt}{2}{}9\above{1pt}{1pt}{-}{}{\cdot}1\above{1pt}{1pt}{-2}{}5\above{1pt}{1pt}{1}{}$}\spacer%
\instructions{m}%
\EasyButWeakLineBreak%
{${45}\above{1pt}{1pt}{r}{2}{16}\above{1pt}{1pt}{*}{2}{72}\above{1pt}{1pt}{s}{2}{80}\above{1pt}{1pt}{s}{2}{8}\above{1pt}{1pt}{s}{2}{720}\above{1pt}{1pt}{*}{2}{4}\above{1pt}{1pt}{l}{2}{18}\above{1pt}{1pt}{}{2}{5}\above{1pt}{1pt}{}{2}{2}\above{1pt}{1pt}{}{2}$}%
\nopagebreak\par%
\nopagebreak\par\leavevmode%
{$\left[\!\llap{\phantom{%
\begingroup \smaller\smaller\smaller
\endgroup%
}}\!\right]$}%
%
%
\hbox{}\par\smallskip%
%
%
\leavevmode%
${L_{196.6}}$%
{} : {$[1\above{1pt}{1pt}{1}{}2\above{1pt}{1pt}{1}{}]\above{1pt}{1pt}{}{0}16\above{1pt}{1pt}{-}{5}{\cdot}1\above{1pt}{1pt}{2}{}9\above{1pt}{1pt}{-}{}{\cdot}1\above{1pt}{1pt}{-2}{}5\above{1pt}{1pt}{1}{}$}\spacer%
\instructions{m}%
\EasyButWeakLineBreak%
{${180}\above{1pt}{1pt}{s}{2}{16}\above{1pt}{1pt}{l}{2}{18}\above{1pt}{1pt}{}{2}{80}\above{1pt}{1pt}{}{2}{2}\above{1pt}{1pt}{}{2}{720}\above{1pt}{1pt}{}{2}{1}\above{1pt}{1pt}{r}{2}{72}\above{1pt}{1pt}{*}{2}{20}\above{1pt}{1pt}{*}{2}{8}\above{1pt}{1pt}{*}{2}$}%
\nopagebreak\par%
\nopagebreak\par\leavevmode%
{$\left[\!\llap{\phantom{%
\begingroup \smaller\smaller\smaller
\endgroup%
}}\!\right]$}%
%
%
\hbox{}\par\smallskip%
%
%
\leavevmode%
${L_{196.7}}$%
{} : {$[1\above{1pt}{1pt}{-}{}2\above{1pt}{1pt}{1}{}]\above{1pt}{1pt}{}{4}16\above{1pt}{1pt}{1}{1}{\cdot}1\above{1pt}{1pt}{2}{}9\above{1pt}{1pt}{-}{}{\cdot}1\above{1pt}{1pt}{-2}{}5\above{1pt}{1pt}{1}{}$}\EasyButWeakLineBreak%
{${45}\above{1pt}{1pt}{}{2}{16}\above{1pt}{1pt}{}{2}{18}\above{1pt}{1pt}{r}{2}{80}\above{1pt}{1pt}{l}{2}{2}\above{1pt}{1pt}{r}{2}{720}\above{1pt}{1pt}{s}{2}{4}\above{1pt}{1pt}{*}{2}{72}\above{1pt}{1pt}{l}{2}{5}\above{1pt}{1pt}{r}{2}{8}\above{1pt}{1pt}{l}{2}$}%
\nopagebreak\par%
\nopagebreak\par\leavevmode%
{$\left[\!\llap{\phantom{%
\begingroup \smaller\smaller\smaller
\endgroup%
}}\!\right]$}%

\medskip%
%
\leavevmode\llap{}%
$W_{197}$%
\qquad\llap{4} lattices, $\chi=24$%
\hfill%
$4|4\slashtwo4|4\slashtwo\rtimes D_{4}$%
\nopagebreak\smallskip\hrule\nopagebreak\medskip%
%
%
\leavevmode%
${L_{197.1}}$%
{} : {$1\above{1pt}{1pt}{2}{2}32\above{1pt}{1pt}{-}{5}{\cdot}1\above{1pt}{1pt}{2}{}3\above{1pt}{1pt}{1}{}$}\EasyButWeakLineBreak%
{${1}\above{1pt}{1pt}{}{4}{2}\above{1pt}{1pt}{*}{4}{4}\above{1pt}{1pt}{l}{2}$}\relax$\,(\times2)$%
\nopagebreak\par%
\nopagebreak\par\leavevmode%
{$\left[\!\llap{\phantom{%
\begingroup \smaller\smaller\smaller\begin{tabular}{@{}c@{}}%
0\\0\\0
\end{tabular}\endgroup%
}}\right.$}%
\begingroup \smaller\smaller\smaller\begin{tabular}{@{}c@{}}%
-1347936\\11328\\5664
\end{tabular}\endgroup%
\kern3pt%
\begingroup \smaller\smaller\smaller\begin{tabular}{@{}c@{}}%
11328\\-95\\-48
\end{tabular}\endgroup%
\kern3pt%
\begingroup \smaller\smaller\smaller\begin{tabular}{@{}c@{}}%
5664\\-48\\-23
\end{tabular}\endgroup%
{$\left.\llap{\phantom{%
\begingroup \smaller\smaller\smaller\begin{tabular}{@{}c@{}}%
0\\0\\0
\end{tabular}\endgroup%
}}\!\right]$}%
\hfil\penalty500%
{$\left[\!\llap{\phantom{%
\begingroup \smaller\smaller\smaller\begin{tabular}{@{}c@{}}%
0\\0\\0
\end{tabular}\endgroup%
}}\right.$}%
\begingroup \smaller\smaller\smaller\begin{tabular}{@{}c@{}}%
7039\\670560\\332640
\end{tabular}\endgroup%
\kern3pt%
\begingroup \smaller\smaller\smaller\begin{tabular}{@{}c@{}}%
-60\\-5716\\-2835
\end{tabular}\endgroup%
\kern3pt%
\begingroup \smaller\smaller\smaller\begin{tabular}{@{}c@{}}%
-28\\-2667\\-1324
\end{tabular}\endgroup%
{$\left.\llap{\phantom{%
\begingroup \smaller\smaller\smaller\begin{tabular}{@{}c@{}}%
0\\0\\0
\end{tabular}\endgroup%
}}\!\right]$}%
\EasyButWeakLineBreak%
{$\left[\!\llap{\phantom{%
\begingroup \smaller\smaller\smaller\begin{tabular}{@{}c@{}}%
0\\0\\0
\end{tabular}\endgroup%
}}\right.$}%
\begingroup \smaller\smaller\smaller\begin{tabular}{@{}c@{}}%
-1\\-95\\-48
\end{tabular}\endgroup%
\HardButStrongLineBreak\kern3pt%
\begingroup \smaller\smaller\smaller\begin{tabular}{@{}c@{}}%
3\\287\\139
\end{tabular}\endgroup%
\HardButStrongLineBreak\kern3pt%
\begingroup \smaller\smaller\smaller\begin{tabular}{@{}c@{}}%
9\\858\\424
\end{tabular}\endgroup%
{$\left.\llap{\phantom{%
\begingroup \smaller\smaller\smaller\begin{tabular}{@{}c@{}}%
0\\0\\0
\end{tabular}\endgroup%
}}\!\right]$}%
%
%
%
%
%
%
%
%
%
%
%
%
%
%

\medskip%
%
\leavevmode\llap{}%
$W_{198}$%
\qquad\llap{8} lattices, $\chi=32$%
\hfill%
$62\infty62\infty\rtimes C_{2}$%
\nopagebreak\smallskip\hrule\nopagebreak\medskip%
%
%
\leavevmode%
${L_{198.1}}$%
{} : {$1\above{1pt}{1pt}{-2}{{\rm II}}32\above{1pt}{1pt}{1}{7}{\cdot}1\above{1pt}{1pt}{1}{}3\above{1pt}{1pt}{-}{}9\above{1pt}{1pt}{-}{}$}\spacer%
\instructions{3}%
\EasyButWeakLineBreak%
{${6}\above{1pt}{1pt}{}{6}{18}\above{1pt}{1pt}{b}{2}{6}\above{1pt}{1pt}{24,7}{\infty a}$}\relax$\,(\times2)$%
\nopagebreak\par%
\nopagebreak\par\leavevmode%
{$\left[\!\llap{\phantom{%
\begingroup \smaller\smaller\smaller\begin{tabular}{@{}c@{}}%
0\\0\\0
\end{tabular}\endgroup%
}}\right.$}%
\begingroup \smaller\smaller\smaller\begin{tabular}{@{}c@{}}%
-3322656\\-643968\\-453600
\end{tabular}\endgroup%
\kern3pt%
\begingroup \smaller\smaller\smaller\begin{tabular}{@{}c@{}}%
-643968\\-124806\\-87915
\end{tabular}\endgroup%
\kern3pt%
\begingroup \smaller\smaller\smaller\begin{tabular}{@{}c@{}}%
-453600\\-87915\\-61922
\end{tabular}\endgroup%
{$\left.\llap{\phantom{%
\begingroup \smaller\smaller\smaller\begin{tabular}{@{}c@{}}%
0\\0\\0
\end{tabular}\endgroup%
}}\!\right]$}%
\hfil\penalty500%
{$\left[\!\llap{\phantom{%
\begingroup \smaller\smaller\smaller\begin{tabular}{@{}c@{}}%
0\\0\\0
\end{tabular}\endgroup%
}}\right.$}%
\begingroup \smaller\smaller\smaller\begin{tabular}{@{}c@{}}%
-250561\\760320\\756000
\end{tabular}\endgroup%
\kern3pt%
\begingroup \smaller\smaller\smaller\begin{tabular}{@{}c@{}}%
-48546\\147311\\146475
\end{tabular}\endgroup%
\kern3pt%
\begingroup \smaller\smaller\smaller\begin{tabular}{@{}c@{}}%
-34220\\103840\\103249
\end{tabular}\endgroup%
{$\left.\llap{\phantom{%
\begingroup \smaller\smaller\smaller\begin{tabular}{@{}c@{}}%
0\\0\\0
\end{tabular}\endgroup%
}}\!\right]$}%
\EasyButWeakLineBreak%
{$\left[\!\llap{\phantom{%
\begingroup \smaller\smaller\smaller\begin{tabular}{@{}c@{}}%
0\\0\\0
\end{tabular}\endgroup%
}}\right.$}%
\begingroup \smaller\smaller\smaller\begin{tabular}{@{}c@{}}%
-150\\457\\450
\end{tabular}\endgroup%
\HardButStrongLineBreak\kern3pt%
\begingroup \smaller\smaller\smaller\begin{tabular}{@{}c@{}}%
110\\-333\\-333
\end{tabular}\endgroup%
\HardButStrongLineBreak\kern3pt%
\begingroup \smaller\smaller\smaller\begin{tabular}{@{}c@{}}%
43\\-131\\-129
\end{tabular}\endgroup%
{$\left.\llap{\phantom{%
\begingroup \smaller\smaller\smaller\begin{tabular}{@{}c@{}}%
0\\0\\0
\end{tabular}\endgroup%
}}\!\right]$}%

\medskip%
%
\leavevmode\llap{}%
$W_{199}$%
\qquad\llap{8} lattices, $\chi=12$%
\hfill%
$222|222|\rtimes D_{2}$%
\nopagebreak\smallskip\hrule\nopagebreak\medskip%
%
%
\leavevmode%
${L_{199.1}}$%
{} : {$1\above{1pt}{1pt}{-2}{2}32\above{1pt}{1pt}{1}{1}{\cdot}1\above{1pt}{1pt}{-}{}3\above{1pt}{1pt}{1}{}9\above{1pt}{1pt}{-}{}$}\spacer%
\instructions{3}%
\EasyButWeakLineBreak%
{${32}\above{1pt}{1pt}{l}{2}{3}\above{1pt}{1pt}{}{2}{288}\above{1pt}{1pt}{r}{2}{2}\above{1pt}{1pt}{b}{2}{288}\above{1pt}{1pt}{*}{2}{12}\above{1pt}{1pt}{s}{2}$}%
\nopagebreak\par%
shares genus with 3-dual\nopagebreak\par%
\nopagebreak\par\leavevmode%
{$\left[\!\llap{\phantom{%
\begingroup \smaller\smaller\smaller\begin{tabular}{@{}c@{}}%
0\\0\\0
\end{tabular}\endgroup%
}}\right.$}%
\begingroup \smaller\smaller\smaller\begin{tabular}{@{}c@{}}%
152352\\75456\\-576
\end{tabular}\endgroup%
\kern3pt%
\begingroup \smaller\smaller\smaller\begin{tabular}{@{}c@{}}%
75456\\37371\\-285
\end{tabular}\endgroup%
\kern3pt%
\begingroup \smaller\smaller\smaller\begin{tabular}{@{}c@{}}%
-576\\-285\\2
\end{tabular}\endgroup%
{$\left.\llap{\phantom{%
\begingroup \smaller\smaller\smaller\begin{tabular}{@{}c@{}}%
0\\0\\0
\end{tabular}\endgroup%
}}\!\right]$}%
\EasyButWeakLineBreak%
{$\left[\!\llap{\phantom{%
\begingroup \smaller\smaller\smaller\begin{tabular}{@{}c@{}}%
0\\0\\0
\end{tabular}\endgroup%
}}\right.$}%
\begingroup \smaller\smaller\smaller\begin{tabular}{@{}c@{}}%
39\\-80\\-160
\end{tabular}\endgroup%
\HardButStrongLineBreak\kern3pt%
\begingroup \smaller\smaller\smaller\begin{tabular}{@{}c@{}}%
20\\-41\\-81
\end{tabular}\endgroup%
\HardButStrongLineBreak\kern3pt%
\begingroup \smaller\smaller\smaller\begin{tabular}{@{}c@{}}%
281\\-576\\-1152
\end{tabular}\endgroup%
\HardButStrongLineBreak\kern3pt%
\begingroup \smaller\smaller\smaller\begin{tabular}{@{}c@{}}%
0\\0\\-1
\end{tabular}\endgroup%
\HardButStrongLineBreak\kern3pt%
\begingroup \smaller\smaller\smaller\begin{tabular}{@{}c@{}}%
-47\\96\\144
\end{tabular}\endgroup%
\HardButStrongLineBreak\kern3pt%
\begingroup \smaller\smaller\smaller\begin{tabular}{@{}c@{}}%
-1\\2\\0
\end{tabular}\endgroup%
{$\left.\llap{\phantom{%
\begingroup \smaller\smaller\smaller\begin{tabular}{@{}c@{}}%
0\\0\\0
\end{tabular}\endgroup%
}}\!\right]$}%

\medskip%
%
\leavevmode\llap{}%
$W_{200}$%
\qquad\llap{8} lattices, $\chi=24$%
\hfill%
$2\slashinfty222|22\rtimes D_{2}$%
\nopagebreak\smallskip\hrule\nopagebreak\medskip%
%
%
\leavevmode%
${L_{200.1}}$%
{} : {$1\above{1pt}{1pt}{-2}{6}32\above{1pt}{1pt}{1}{1}{\cdot}1\above{1pt}{1pt}{-}{}3\above{1pt}{1pt}{-}{}9\above{1pt}{1pt}{1}{}$}\spacer%
\instructions{3}%
\EasyButWeakLineBreak%
{${32}\above{1pt}{1pt}{b}{2}{6}\above{1pt}{1pt}{24,17}{\infty b}{6}\above{1pt}{1pt}{l}{2}{32}\above{1pt}{1pt}{}{2}{9}\above{1pt}{1pt}{r}{2}{32}\above{1pt}{1pt}{s}{2}{36}\above{1pt}{1pt}{*}{2}$}%
\nopagebreak\par%
\nopagebreak\par\leavevmode%
{$\left[\!\llap{\phantom{%
\begingroup \smaller\smaller\smaller\begin{tabular}{@{}c@{}}%
0\\0\\0
\end{tabular}\endgroup%
}}\right.$}%
\begingroup \smaller\smaller\smaller\begin{tabular}{@{}c@{}}%
-202464\\-34272\\1440
\end{tabular}\endgroup%
\kern3pt%
\begingroup \smaller\smaller\smaller\begin{tabular}{@{}c@{}}%
-34272\\-5799\\243
\end{tabular}\endgroup%
\kern3pt%
\begingroup \smaller\smaller\smaller\begin{tabular}{@{}c@{}}%
1440\\243\\-10
\end{tabular}\endgroup%
{$\left.\llap{\phantom{%
\begingroup \smaller\smaller\smaller\begin{tabular}{@{}c@{}}%
0\\0\\0
\end{tabular}\endgroup%
}}\!\right]$}%
\EasyButWeakLineBreak%
{$\left[\!\llap{\phantom{%
\begingroup \smaller\smaller\smaller\begin{tabular}{@{}c@{}}%
0\\0\\0
\end{tabular}\endgroup%
}}\right.$}%
\begingroup \smaller\smaller\smaller\begin{tabular}{@{}c@{}}%
-19\\128\\368
\end{tabular}\endgroup%
\HardButStrongLineBreak\kern3pt%
\begingroup \smaller\smaller\smaller\begin{tabular}{@{}c@{}}%
-5\\34\\105
\end{tabular}\endgroup%
\HardButStrongLineBreak\kern3pt%
\begingroup \smaller\smaller\smaller\begin{tabular}{@{}c@{}}%
-2\\14\\51
\end{tabular}\endgroup%
\HardButStrongLineBreak\kern3pt%
\begingroup \smaller\smaller\smaller\begin{tabular}{@{}c@{}}%
5\\-32\\-64
\end{tabular}\endgroup%
\HardButStrongLineBreak\kern3pt%
\begingroup \smaller\smaller\smaller\begin{tabular}{@{}c@{}}%
4\\-27\\-81
\end{tabular}\endgroup%
\HardButStrongLineBreak\kern3pt%
\begingroup \smaller\smaller\smaller\begin{tabular}{@{}c@{}}%
7\\-48\\-160
\end{tabular}\endgroup%
\HardButStrongLineBreak\kern3pt%
\begingroup \smaller\smaller\smaller\begin{tabular}{@{}c@{}}%
-1\\6\\0
\end{tabular}\endgroup%
{$\left.\llap{\phantom{%
\begingroup \smaller\smaller\smaller\begin{tabular}{@{}c@{}}%
0\\0\\0
\end{tabular}\endgroup%
}}\!\right]$}%

\medskip%
%
\leavevmode\llap{}%
$W_{201}$%
\qquad\llap{8} lattices, $\chi=112$%
\hfill%
$3\infty2\infty2\infty\infty3\infty2\infty2\infty\infty\rtimes C_{2}$%
\nopagebreak\smallskip\hrule\nopagebreak\medskip%
%
%
\leavevmode%
${L_{201.1}}$%
{} : {$1\above{1pt}{1pt}{-2}{{\rm II}}8\above{1pt}{1pt}{-}{5}{\cdot}1\above{1pt}{1pt}{-2}{}49\above{1pt}{1pt}{1}{}$}\spacer%
\instructions{2}%
\EasyButWeakLineBreak%
{${2}\above{1pt}{1pt}{+}{3}{2}\above{1pt}{1pt}{28,25}{\infty b}{8}\above{1pt}{1pt}{b}{2}{98}\above{1pt}{1pt}{4,1}{\infty b}{392}\above{1pt}{1pt}{b}{2}{2}\above{1pt}{1pt}{28,1}{\infty b}{8}\above{1pt}{1pt}{14,11}{\infty z}$}\relax$\,(\times2)$%
\nopagebreak\par%
\nopagebreak\par\leavevmode%
{$\left[\!\llap{\phantom{%
\begingroup \smaller\smaller\smaller
\endgroup%
}}\!\right]$}%

\medskip%
%
\leavevmode\llap{}%
$W_{202}$%
\qquad\llap{12} lattices, $\chi=15$%
\hfill%
$222224$%
\nopagebreak\smallskip\hrule\nopagebreak\medskip%
%
%
\leavevmode%
${L_{202.1}}$%
{} : {$1\above{1pt}{1pt}{-2}{{\rm II}}4\above{1pt}{1pt}{1}{7}{\cdot}1\above{1pt}{1pt}{2}{}9\above{1pt}{1pt}{1}{}{\cdot}1\above{1pt}{1pt}{2}{}11\above{1pt}{1pt}{-}{}$}\spacer%
\instructions{2}%
\EasyButWeakLineBreak%
{${2}\above{1pt}{1pt}{s}{2}{198}\above{1pt}{1pt}{b}{2}{4}\above{1pt}{1pt}{*}{2}{36}\above{1pt}{1pt}{b}{2}{22}\above{1pt}{1pt}{b}{2}{4}\above{1pt}{1pt}{*}{4}$}%
\nopagebreak\par%
\nopagebreak\par\leavevmode%
{$\left[\!\llap{\phantom{%
\begingroup \smaller\smaller\smaller\begin{tabular}{@{}c@{}}%
0\\0\\0
\end{tabular}\endgroup%
}}\right.$}%
\begingroup \smaller\smaller\smaller\begin{tabular}{@{}c@{}}%
-900900\\3960\\5544
\end{tabular}\endgroup%
\kern3pt%
\begingroup \smaller\smaller\smaller\begin{tabular}{@{}c@{}}%
3960\\-14\\-25
\end{tabular}\endgroup%
\kern3pt%
\begingroup \smaller\smaller\smaller\begin{tabular}{@{}c@{}}%
5544\\-25\\-34
\end{tabular}\endgroup%
{$\left.\llap{\phantom{%
\begingroup \smaller\smaller\smaller\begin{tabular}{@{}c@{}}%
0\\0\\0
\end{tabular}\endgroup%
}}\!\right]$}%
\EasyButWeakLineBreak%
{$\left[\!\llap{\phantom{%
\begingroup \smaller\smaller\smaller\begin{tabular}{@{}c@{}}%
0\\0\\0
\end{tabular}\endgroup%
}}\right.$}%
\begingroup \smaller\smaller\smaller\begin{tabular}{@{}c@{}}%
1\\27\\143
\end{tabular}\endgroup%
\HardButStrongLineBreak\kern3pt%
\begingroup \smaller\smaller\smaller\begin{tabular}{@{}c@{}}%
49\\1287\\7029
\end{tabular}\endgroup%
\HardButStrongLineBreak\kern3pt%
\begingroup \smaller\smaller\smaller\begin{tabular}{@{}c@{}}%
5\\130\\718
\end{tabular}\endgroup%
\HardButStrongLineBreak\kern3pt%
\begingroup \smaller\smaller\smaller\begin{tabular}{@{}c@{}}%
5\\126\\720
\end{tabular}\endgroup%
\HardButStrongLineBreak\kern3pt%
\begingroup \smaller\smaller\smaller\begin{tabular}{@{}c@{}}%
-2\\-55\\-286
\end{tabular}\endgroup%
\HardButStrongLineBreak\kern3pt%
\begingroup \smaller\smaller\smaller\begin{tabular}{@{}c@{}}%
-1\\-26\\-144
\end{tabular}\endgroup%
{$\left.\llap{\phantom{%
\begingroup \smaller\smaller\smaller\begin{tabular}{@{}c@{}}%
0\\0\\0
\end{tabular}\endgroup%
}}\!\right]$}%

\medskip%
%
\leavevmode\llap{}%
$W_{203}$%
\qquad\llap{4} lattices, $\chi=48$%
\hfill%
$\slashinfty22|22\slashinfty22|22\rtimes D_{4}$%
\nopagebreak\smallskip\hrule\nopagebreak\medskip%
%
%
\leavevmode%
${L_{203.1}}$%
{} : {$1\above{1pt}{1pt}{2}{6}16\above{1pt}{1pt}{1}{1}{\cdot}1\above{1pt}{1pt}{-2}{}7\above{1pt}{1pt}{1}{}$}\EasyButWeakLineBreak%
{${28}\above{1pt}{1pt}{8,5}{\infty z}{7}\above{1pt}{1pt}{}{2}{16}\above{1pt}{1pt}{r}{2}{14}\above{1pt}{1pt}{b}{2}{16}\above{1pt}{1pt}{*}{2}$}\relax$\,(\times2)$%
\nopagebreak\par%
\nopagebreak\par\leavevmode%
{$\left[\!\llap{\phantom{%
\begingroup \smaller\smaller\smaller
\endgroup%
}}\!\right]$}%
%
%
%
%
%
%
%
%
%
%
%
%
%
%

\medskip%
%
\leavevmode\llap{}%
$W_{204}$%
\qquad\llap{12} lattices, $\chi=36$%
\hfill%
$2222222222\rtimes C_{2}$%
\nopagebreak\smallskip\hrule\nopagebreak\medskip%
%
%
\leavevmode%
${L_{204.1}}$%
{} : {$1\above{1pt}{1pt}{-2}{{\rm II}}4\above{1pt}{1pt}{1}{1}{\cdot}1\above{1pt}{1pt}{2}{}9\above{1pt}{1pt}{-}{}{\cdot}1\above{1pt}{1pt}{-2}{}13\above{1pt}{1pt}{-}{}$}\spacer%
\instructions{2}%
\EasyButWeakLineBreak%
{${234}\above{1pt}{1pt}{l}{2}{4}\above{1pt}{1pt}{r}{2}{18}\above{1pt}{1pt}{b}{2}{26}\above{1pt}{1pt}{b}{2}{2}\above{1pt}{1pt}{b}{2}$}\relax$\,(\times2)$%
\nopagebreak\par%
\nopagebreak\par\leavevmode%
{$\left[\!\llap{\phantom{%
\begingroup \smaller\smaller\smaller
\endgroup%
}}\!\right]$}%

\medskip%
%
\leavevmode\llap{}%
$W_{205}$%
\qquad\llap{60} lattices, $\chi=12$%
\hfill%
$22|222|2\rtimes D_{2}$%
\nopagebreak\smallskip\hrule\nopagebreak\medskip%
%
%
\leavevmode%
${L_{205.1}}$%
{} : {$1\above{1pt}{1pt}{2}{0}8\above{1pt}{1pt}{1}{7}{\cdot}1\above{1pt}{1pt}{2}{}3\above{1pt}{1pt}{-}{}{\cdot}1\above{1pt}{1pt}{-2}{}5\above{1pt}{1pt}{-}{}$}\EasyButWeakLineBreak%
{${60}\above{1pt}{1pt}{*}{2}{4}\above{1pt}{1pt}{s}{2}{40}\above{1pt}{1pt}{l}{2}{1}\above{1pt}{1pt}{}{2}{15}\above{1pt}{1pt}{r}{2}{8}\above{1pt}{1pt}{s}{2}$}%
\nopagebreak\par%
\nopagebreak\par\leavevmode%
{$\left[\!\llap{\phantom{%
\begingroup \smaller\smaller\smaller\begin{tabular}{@{}c@{}}%
0\\0\\0
\end{tabular}\endgroup%
}}\right.$}%
\begingroup \smaller\smaller\smaller\begin{tabular}{@{}c@{}}%
347640\\6240\\840
\end{tabular}\endgroup%
\kern3pt%
\begingroup \smaller\smaller\smaller\begin{tabular}{@{}c@{}}%
6240\\112\\15
\end{tabular}\endgroup%
\kern3pt%
\begingroup \smaller\smaller\smaller\begin{tabular}{@{}c@{}}%
840\\15\\1
\end{tabular}\endgroup%
{$\left.\llap{\phantom{%
\begingroup \smaller\smaller\smaller\begin{tabular}{@{}c@{}}%
0\\0\\0
\end{tabular}\endgroup%
}}\!\right]$}%
\EasyButWeakLineBreak%
{$\left[\!\llap{\phantom{%
\begingroup \smaller\smaller\smaller\begin{tabular}{@{}c@{}}%
0\\0\\0
\end{tabular}\endgroup%
}}\right.$}%
\begingroup \smaller\smaller\smaller\begin{tabular}{@{}c@{}}%
-17\\960\\-90
\end{tabular}\endgroup%
\HardButStrongLineBreak\kern3pt%
\begingroup \smaller\smaller\smaller\begin{tabular}{@{}c@{}}%
-5\\282\\-26
\end{tabular}\endgroup%
\HardButStrongLineBreak\kern3pt%
\begingroup \smaller\smaller\smaller\begin{tabular}{@{}c@{}}%
-11\\620\\-60
\end{tabular}\endgroup%
\HardButStrongLineBreak\kern3pt%
\begingroup \smaller\smaller\smaller\begin{tabular}{@{}c@{}}%
0\\0\\-1
\end{tabular}\endgroup%
\HardButStrongLineBreak\kern3pt%
\begingroup \smaller\smaller\smaller\begin{tabular}{@{}c@{}}%
4\\-225\\15
\end{tabular}\endgroup%
\HardButStrongLineBreak\kern3pt%
\begingroup \smaller\smaller\smaller\begin{tabular}{@{}c@{}}%
1\\-56\\4
\end{tabular}\endgroup%
{$\left.\llap{\phantom{%
\begingroup \smaller\smaller\smaller\begin{tabular}{@{}c@{}}%
0\\0\\0
\end{tabular}\endgroup%
}}\!\right]$}%
%
%
\hbox{}\par\smallskip%
%
%
\leavevmode%
${L_{205.2}}$%
{} : {$[1\above{1pt}{1pt}{1}{}2\above{1pt}{1pt}{1}{}]\above{1pt}{1pt}{}{0}16\above{1pt}{1pt}{1}{7}{\cdot}1\above{1pt}{1pt}{2}{}3\above{1pt}{1pt}{-}{}{\cdot}1\above{1pt}{1pt}{-2}{}5\above{1pt}{1pt}{-}{}$}\spacer%
\instructions{2}%
\EasyButWeakLineBreak%
{${15}\above{1pt}{1pt}{r}{2}{16}\above{1pt}{1pt}{*}{2}{40}\above{1pt}{1pt}{l}{2}{1}\above{1pt}{1pt}{}{2}{240}\above{1pt}{1pt}{}{2}{2}\above{1pt}{1pt}{}{2}$}%
\nopagebreak\par%
\nopagebreak\par\leavevmode%
{$\left[\!\llap{\phantom{%
\begingroup \smaller\smaller\smaller\begin{tabular}{@{}c@{}}%
0\\0\\0
\end{tabular}\endgroup%
}}\right.$}%
\begingroup \smaller\smaller\smaller\begin{tabular}{@{}c@{}}%
-161040\\-960\\1680
\end{tabular}\endgroup%
\kern3pt%
\begingroup \smaller\smaller\smaller\begin{tabular}{@{}c@{}}%
-960\\2\\8
\end{tabular}\endgroup%
\kern3pt%
\begingroup \smaller\smaller\smaller\begin{tabular}{@{}c@{}}%
1680\\8\\-17
\end{tabular}\endgroup%
{$\left.\llap{\phantom{%
\begingroup \smaller\smaller\smaller\begin{tabular}{@{}c@{}}%
0\\0\\0
\end{tabular}\endgroup%
}}\!\right]$}%
\EasyButWeakLineBreak%
{$\left[\!\llap{\phantom{%
\begingroup \smaller\smaller\smaller\begin{tabular}{@{}c@{}}%
0\\0\\0
\end{tabular}\endgroup%
}}\right.$}%
\begingroup \smaller\smaller\smaller\begin{tabular}{@{}c@{}}%
2\\60\\225
\end{tabular}\endgroup%
\HardButStrongLineBreak\kern3pt%
\begingroup \smaller\smaller\smaller\begin{tabular}{@{}c@{}}%
-1\\-28\\-112
\end{tabular}\endgroup%
\HardButStrongLineBreak\kern3pt%
\begingroup \smaller\smaller\smaller\begin{tabular}{@{}c@{}}%
-3\\-90\\-340
\end{tabular}\endgroup%
\HardButStrongLineBreak\kern3pt%
\begingroup \smaller\smaller\smaller\begin{tabular}{@{}c@{}}%
0\\-1\\-1
\end{tabular}\endgroup%
\HardButStrongLineBreak\kern3pt%
\begingroup \smaller\smaller\smaller\begin{tabular}{@{}c@{}}%
13\\360\\1440
\end{tabular}\endgroup%
\HardButStrongLineBreak\kern3pt%
\begingroup \smaller\smaller\smaller\begin{tabular}{@{}c@{}}%
1\\29\\112
\end{tabular}\endgroup%
{$\left.\llap{\phantom{%
\begingroup \smaller\smaller\smaller\begin{tabular}{@{}c@{}}%
0\\0\\0
\end{tabular}\endgroup%
}}\!\right]$}%
%
%
\hbox{}\par\smallskip%
%
%
\leavevmode%
${L_{205.3}}$%
{} : {$[1\above{1pt}{1pt}{-}{}2\above{1pt}{1pt}{1}{}]\above{1pt}{1pt}{}{4}16\above{1pt}{1pt}{-}{3}{\cdot}1\above{1pt}{1pt}{2}{}3\above{1pt}{1pt}{-}{}{\cdot}1\above{1pt}{1pt}{-2}{}5\above{1pt}{1pt}{-}{}$}\spacer%
\instructions{m}%
\EasyButWeakLineBreak%
{${60}\above{1pt}{1pt}{*}{2}{16}\above{1pt}{1pt}{s}{2}{40}\above{1pt}{1pt}{*}{2}{4}\above{1pt}{1pt}{s}{2}{240}\above{1pt}{1pt}{l}{2}{2}\above{1pt}{1pt}{r}{2}$}%
\nopagebreak\par%
\nopagebreak\par\leavevmode%
{$\left[\!\llap{\phantom{%
\begingroup \smaller\smaller\smaller\begin{tabular}{@{}c@{}}%
0\\0\\0
\end{tabular}\endgroup%
}}\right.$}%
\begingroup \smaller\smaller\smaller\begin{tabular}{@{}c@{}}%
26160\\-6720\\2880
\end{tabular}\endgroup%
\kern3pt%
\begingroup \smaller\smaller\smaller\begin{tabular}{@{}c@{}}%
-6720\\1726\\-740
\end{tabular}\endgroup%
\kern3pt%
\begingroup \smaller\smaller\smaller\begin{tabular}{@{}c@{}}%
2880\\-740\\317
\end{tabular}\endgroup%
{$\left.\llap{\phantom{%
\begingroup \smaller\smaller\smaller\begin{tabular}{@{}c@{}}%
0\\0\\0
\end{tabular}\endgroup%
}}\!\right]$}%
\EasyButWeakLineBreak%
{$\left[\!\llap{\phantom{%
\begingroup \smaller\smaller\smaller\begin{tabular}{@{}c@{}}%
0\\0\\0
\end{tabular}\endgroup%
}}\right.$}%
\begingroup \smaller\smaller\smaller\begin{tabular}{@{}c@{}}%
-11\\-30\\30
\end{tabular}\endgroup%
\HardButStrongLineBreak\kern3pt%
\begingroup \smaller\smaller\smaller\begin{tabular}{@{}c@{}}%
1\\4\\0
\end{tabular}\endgroup%
\HardButStrongLineBreak\kern3pt%
\begingroup \smaller\smaller\smaller\begin{tabular}{@{}c@{}}%
37\\110\\-80
\end{tabular}\endgroup%
\HardButStrongLineBreak\kern3pt%
\begingroup \smaller\smaller\smaller\begin{tabular}{@{}c@{}}%
19\\56\\-42
\end{tabular}\endgroup%
\HardButStrongLineBreak\kern3pt%
\begingroup \smaller\smaller\smaller\begin{tabular}{@{}c@{}}%
163\\480\\-360
\end{tabular}\endgroup%
\HardButStrongLineBreak\kern3pt%
\begingroup \smaller\smaller\smaller\begin{tabular}{@{}c@{}}%
1\\3\\-2
\end{tabular}\endgroup%
{$\left.\llap{\phantom{%
\begingroup \smaller\smaller\smaller\begin{tabular}{@{}c@{}}%
0\\0\\0
\end{tabular}\endgroup%
}}\!\right]$}%
%
%
\hbox{}\par\smallskip%
%
%
\leavevmode%
${L_{205.4}}$%
{} : {$[1\above{1pt}{1pt}{-}{}2\above{1pt}{1pt}{1}{}]\above{1pt}{1pt}{}{2}16\above{1pt}{1pt}{-}{5}{\cdot}1\above{1pt}{1pt}{2}{}3\above{1pt}{1pt}{-}{}{\cdot}1\above{1pt}{1pt}{-2}{}5\above{1pt}{1pt}{-}{}$}\spacer%
\instructions{m}%
\EasyButWeakLineBreak%
{${60}\above{1pt}{1pt}{s}{2}{16}\above{1pt}{1pt}{l}{2}{10}\above{1pt}{1pt}{}{2}{1}\above{1pt}{1pt}{r}{2}{240}\above{1pt}{1pt}{*}{2}{8}\above{1pt}{1pt}{*}{2}$}%
\nopagebreak\par%
\nopagebreak\par\leavevmode%
{$\left[\!\llap{\phantom{%
\begingroup \smaller\smaller\smaller\begin{tabular}{@{}c@{}}%
0\\0\\0
\end{tabular}\endgroup%
}}\right.$}%
\begingroup \smaller\smaller\smaller\begin{tabular}{@{}c@{}}%
-24240\\960\\0
\end{tabular}\endgroup%
\kern3pt%
\begingroup \smaller\smaller\smaller\begin{tabular}{@{}c@{}}%
960\\-2\\-6
\end{tabular}\endgroup%
\kern3pt%
\begingroup \smaller\smaller\smaller\begin{tabular}{@{}c@{}}%
0\\-6\\1
\end{tabular}\endgroup%
{$\left.\llap{\phantom{%
\begingroup \smaller\smaller\smaller\begin{tabular}{@{}c@{}}%
0\\0\\0
\end{tabular}\endgroup%
}}\!\right]$}%
\EasyButWeakLineBreak%
{$\left[\!\llap{\phantom{%
\begingroup \smaller\smaller\smaller\begin{tabular}{@{}c@{}}%
0\\0\\0
\end{tabular}\endgroup%
}}\right.$}%
\begingroup \smaller\smaller\smaller\begin{tabular}{@{}c@{}}%
-1\\-30\\-150
\end{tabular}\endgroup%
\HardButStrongLineBreak\kern3pt%
\begingroup \smaller\smaller\smaller\begin{tabular}{@{}c@{}}%
1\\24\\152
\end{tabular}\endgroup%
\HardButStrongLineBreak\kern3pt%
\begingroup \smaller\smaller\smaller\begin{tabular}{@{}c@{}}%
1\\25\\150
\end{tabular}\endgroup%
\HardButStrongLineBreak\kern3pt%
\begingroup \smaller\smaller\smaller\begin{tabular}{@{}c@{}}%
0\\0\\-1
\end{tabular}\endgroup%
\HardButStrongLineBreak\kern3pt%
\begingroup \smaller\smaller\smaller\begin{tabular}{@{}c@{}}%
-7\\-180\\-1080
\end{tabular}\endgroup%
\HardButStrongLineBreak\kern3pt%
\begingroup \smaller\smaller\smaller\begin{tabular}{@{}c@{}}%
-1\\-26\\-152
\end{tabular}\endgroup%
{$\left.\llap{\phantom{%
\begingroup \smaller\smaller\smaller\begin{tabular}{@{}c@{}}%
0\\0\\0
\end{tabular}\endgroup%
}}\!\right]$}%
%
%
\hbox{}\par\smallskip%
%
%
\leavevmode%
${L_{205.5}}$%
{} : {$[1\above{1pt}{1pt}{1}{}2\above{1pt}{1pt}{1}{}]\above{1pt}{1pt}{}{6}16\above{1pt}{1pt}{1}{1}{\cdot}1\above{1pt}{1pt}{2}{}3\above{1pt}{1pt}{-}{}{\cdot}1\above{1pt}{1pt}{-2}{}5\above{1pt}{1pt}{-}{}$}\EasyButWeakLineBreak%
{${15}\above{1pt}{1pt}{}{2}{16}\above{1pt}{1pt}{}{2}{10}\above{1pt}{1pt}{r}{2}{4}\above{1pt}{1pt}{*}{2}{240}\above{1pt}{1pt}{s}{2}{8}\above{1pt}{1pt}{l}{2}$}%
\nopagebreak\par%
\nopagebreak\par\leavevmode%
{$\left[\!\llap{\phantom{%
\begingroup \smaller\smaller\smaller\begin{tabular}{@{}c@{}}%
0\\0\\0
\end{tabular}\endgroup%
}}\right.$}%
\begingroup \smaller\smaller\smaller\begin{tabular}{@{}c@{}}%
151440\\3120\\2880
\end{tabular}\endgroup%
\kern3pt%
\begingroup \smaller\smaller\smaller\begin{tabular}{@{}c@{}}%
3120\\58\\56
\end{tabular}\endgroup%
\kern3pt%
\begingroup \smaller\smaller\smaller\begin{tabular}{@{}c@{}}%
2880\\56\\53
\end{tabular}\endgroup%
{$\left.\llap{\phantom{%
\begingroup \smaller\smaller\smaller\begin{tabular}{@{}c@{}}%
0\\0\\0
\end{tabular}\endgroup%
}}\!\right]$}%
\EasyButWeakLineBreak%
{$\left[\!\llap{\phantom{%
\begingroup \smaller\smaller\smaller\begin{tabular}{@{}c@{}}%
0\\0\\0
\end{tabular}\endgroup%
}}\right.$}%
\begingroup \smaller\smaller\smaller\begin{tabular}{@{}c@{}}%
1\\75\\-135
\end{tabular}\endgroup%
\HardButStrongLineBreak\kern3pt%
\begingroup \smaller\smaller\smaller\begin{tabular}{@{}c@{}}%
-1\\-56\\112
\end{tabular}\endgroup%
\HardButStrongLineBreak\kern3pt%
\begingroup \smaller\smaller\smaller\begin{tabular}{@{}c@{}}%
-2\\-125\\240
\end{tabular}\endgroup%
\HardButStrongLineBreak\kern3pt%
\begingroup \smaller\smaller\smaller\begin{tabular}{@{}c@{}}%
-1\\-64\\122
\end{tabular}\endgroup%
\HardButStrongLineBreak\kern3pt%
\begingroup \smaller\smaller\smaller\begin{tabular}{@{}c@{}}%
-1\\-60\\120
\end{tabular}\endgroup%
\HardButStrongLineBreak\kern3pt%
\begingroup \smaller\smaller\smaller\begin{tabular}{@{}c@{}}%
1\\66\\-124
\end{tabular}\endgroup%
{$\left.\llap{\phantom{%
\begingroup \smaller\smaller\smaller\begin{tabular}{@{}c@{}}%
0\\0\\0
\end{tabular}\endgroup%
}}\!\right]$}%

\medskip%
%
\leavevmode\llap{}%
$W_{206}$%
\qquad\llap{120} lattices, $\chi=36$%
\hfill%
$2\slashinfty2222|222\rtimes D_{2}$%
\nopagebreak\smallskip\hrule\nopagebreak\medskip%
%
%
\leavevmode%
${L_{206.1}}$%
{} : {$1\above{1pt}{1pt}{2}{0}8\above{1pt}{1pt}{1}{7}{\cdot}1\above{1pt}{1pt}{-}{}3\above{1pt}{1pt}{1}{}9\above{1pt}{1pt}{1}{}{\cdot}1\above{1pt}{1pt}{2}{}5\above{1pt}{1pt}{1}{}$}\spacer%
\instructions{3}%
\EasyButWeakLineBreak%
{${9}\above{1pt}{1pt}{}{2}{120}\above{1pt}{1pt}{6,1}{\infty}{120}\above{1pt}{1pt}{*}{2}{36}\above{1pt}{1pt}{l}{2}{5}\above{1pt}{1pt}{}{2}{3}\above{1pt}{1pt}{r}{2}{8}\above{1pt}{1pt}{s}{2}{12}\above{1pt}{1pt}{*}{2}{20}\above{1pt}{1pt}{l}{2}$}%
\nopagebreak\par%
\nopagebreak\par\leavevmode%
{$\left[\!\llap{\phantom{%
\begingroup \smaller\smaller\smaller\begin{tabular}{@{}c@{}}%
0\\0\\0
\end{tabular}\endgroup%
}}\right.$}%
\begingroup \smaller\smaller\smaller\begin{tabular}{@{}c@{}}%
-330120\\-1800\\2520
\end{tabular}\endgroup%
\kern3pt%
\begingroup \smaller\smaller\smaller\begin{tabular}{@{}c@{}}%
-1800\\3\\12
\end{tabular}\endgroup%
\kern3pt%
\begingroup \smaller\smaller\smaller\begin{tabular}{@{}c@{}}%
2520\\12\\-19
\end{tabular}\endgroup%
{$\left.\llap{\phantom{%
\begingroup \smaller\smaller\smaller\begin{tabular}{@{}c@{}}%
0\\0\\0
\end{tabular}\endgroup%
}}\!\right]$}%
\EasyButWeakLineBreak%
{$\left[\!\llap{\phantom{%
\begingroup \smaller\smaller\smaller\begin{tabular}{@{}c@{}}%
0\\0\\0
\end{tabular}\endgroup%
}}\right.$}%
\begingroup \smaller\smaller\smaller\begin{tabular}{@{}c@{}}%
-5\\-114\\-747
\end{tabular}\endgroup%
\HardButStrongLineBreak\kern3pt%
\begingroup \smaller\smaller\smaller\begin{tabular}{@{}c@{}}%
-13\\-280\\-1920
\end{tabular}\endgroup%
\HardButStrongLineBreak\kern3pt%
\begingroup \smaller\smaller\smaller\begin{tabular}{@{}c@{}}%
-7\\-140\\-1020
\end{tabular}\endgroup%
\HardButStrongLineBreak\kern3pt%
\begingroup \smaller\smaller\smaller\begin{tabular}{@{}c@{}}%
-1\\-18\\-144
\end{tabular}\endgroup%
\HardButStrongLineBreak\kern3pt%
\begingroup \smaller\smaller\smaller\begin{tabular}{@{}c@{}}%
1\\20\\145
\end{tabular}\endgroup%
\HardButStrongLineBreak\kern3pt%
\begingroup \smaller\smaller\smaller\begin{tabular}{@{}c@{}}%
1\\19\\144
\end{tabular}\endgroup%
\HardButStrongLineBreak\kern3pt%
\begingroup \smaller\smaller\smaller\begin{tabular}{@{}c@{}}%
1\\16\\140
\end{tabular}\endgroup%
\HardButStrongLineBreak\kern3pt%
\begingroup \smaller\smaller\smaller\begin{tabular}{@{}c@{}}%
-1\\-32\\-162
\end{tabular}\endgroup%
\HardButStrongLineBreak\kern3pt%
\begingroup \smaller\smaller\smaller\begin{tabular}{@{}c@{}}%
-7\\-170\\-1060
\end{tabular}\endgroup%
{$\left.\llap{\phantom{%
\begingroup \smaller\smaller\smaller\begin{tabular}{@{}c@{}}%
0\\0\\0
\end{tabular}\endgroup%
}}\!\right]$}%
%
%
\hbox{}\par\smallskip%
%
%
\leavevmode%
${L_{206.2}}$%
{} : {$[1\above{1pt}{1pt}{-}{}2\above{1pt}{1pt}{1}{}]\above{1pt}{1pt}{}{4}16\above{1pt}{1pt}{-}{3}{\cdot}1\above{1pt}{1pt}{-}{}3\above{1pt}{1pt}{1}{}9\above{1pt}{1pt}{1}{}{\cdot}1\above{1pt}{1pt}{2}{}5\above{1pt}{1pt}{1}{}$}\spacer%
\instructions{3,2}%
\EasyButWeakLineBreak%
{${36}\above{1pt}{1pt}{l}{2}{30}\above{1pt}{1pt}{24,19}{\infty}{120}\above{1pt}{1pt}{*}{2}{144}\above{1pt}{1pt}{l}{2}{5}\above{1pt}{1pt}{}{2}{48}\above{1pt}{1pt}{}{2}{2}\above{1pt}{1pt}{}{2}{3}\above{1pt}{1pt}{r}{2}{80}\above{1pt}{1pt}{*}{2}$}%
\nopagebreak\par%
\nopagebreak\par\leavevmode%
{$\left[\!\llap{\phantom{%
\begingroup \smaller\smaller\smaller\begin{tabular}{@{}c@{}}%
0\\0\\0
\end{tabular}\endgroup%
}}\right.$}%
\begingroup \smaller\smaller\smaller\begin{tabular}{@{}c@{}}%
223920\\14400\\-1440
\end{tabular}\endgroup%
\kern3pt%
\begingroup \smaller\smaller\smaller\begin{tabular}{@{}c@{}}%
14400\\174\\-36
\end{tabular}\endgroup%
\kern3pt%
\begingroup \smaller\smaller\smaller\begin{tabular}{@{}c@{}}%
-1440\\-36\\5
\end{tabular}\endgroup%
{$\left.\llap{\phantom{%
\begingroup \smaller\smaller\smaller\begin{tabular}{@{}c@{}}%
0\\0\\0
\end{tabular}\endgroup%
}}\!\right]$}%
\EasyButWeakLineBreak%
{$\left[\!\llap{\phantom{%
\begingroup \smaller\smaller\smaller\begin{tabular}{@{}c@{}}%
0\\0\\0
\end{tabular}\endgroup%
}}\right.$}%
\begingroup \smaller\smaller\smaller\begin{tabular}{@{}c@{}}%
-1\\-48\\-630
\end{tabular}\endgroup%
\HardButStrongLineBreak\kern3pt%
\begingroup \smaller\smaller\smaller\begin{tabular}{@{}c@{}}%
-2\\-95\\-1260
\end{tabular}\endgroup%
\HardButStrongLineBreak\kern3pt%
\begingroup \smaller\smaller\smaller\begin{tabular}{@{}c@{}}%
-1\\-50\\-660
\end{tabular}\endgroup%
\HardButStrongLineBreak\kern3pt%
\begingroup \smaller\smaller\smaller\begin{tabular}{@{}c@{}}%
7\\324\\4320
\end{tabular}\endgroup%
\HardButStrongLineBreak\kern3pt%
\begingroup \smaller\smaller\smaller\begin{tabular}{@{}c@{}}%
3\\140\\1865
\end{tabular}\endgroup%
\HardButStrongLineBreak\kern3pt%
\begingroup \smaller\smaller\smaller\begin{tabular}{@{}c@{}}%
7\\328\\4368
\end{tabular}\endgroup%
\HardButStrongLineBreak\kern3pt%
\begingroup \smaller\smaller\smaller\begin{tabular}{@{}c@{}}%
1\\47\\626
\end{tabular}\endgroup%
\HardButStrongLineBreak\kern3pt%
\begingroup \smaller\smaller\smaller\begin{tabular}{@{}c@{}}%
1\\47\\627
\end{tabular}\endgroup%
\HardButStrongLineBreak\kern3pt%
\begingroup \smaller\smaller\smaller\begin{tabular}{@{}c@{}}%
3\\140\\1880
\end{tabular}\endgroup%
{$\left.\llap{\phantom{%
\begingroup \smaller\smaller\smaller\begin{tabular}{@{}c@{}}%
0\\0\\0
\end{tabular}\endgroup%
}}\!\right]$}%
%
%
\hbox{}\par\smallskip%
%
%
\leavevmode%
${L_{206.3}}$%
{} : {$[1\above{1pt}{1pt}{1}{}2\above{1pt}{1pt}{1}{}]\above{1pt}{1pt}{}{0}16\above{1pt}{1pt}{1}{7}{\cdot}1\above{1pt}{1pt}{1}{}3\above{1pt}{1pt}{1}{}9\above{1pt}{1pt}{-}{}{\cdot}1\above{1pt}{1pt}{2}{}5\above{1pt}{1pt}{1}{}$}\spacer%
\instructions{32,3m,3,m}%
\EasyButWeakLineBreak%
{${16}\above{1pt}{1pt}{s}{2}{120}\above{1pt}{1pt}{24,17}{\infty z}{30}\above{1pt}{1pt}{}{2}{1}\above{1pt}{1pt}{r}{2}{720}\above{1pt}{1pt}{*}{2}{12}\above{1pt}{1pt}{l}{2}{18}\above{1pt}{1pt}{r}{2}{48}\above{1pt}{1pt}{s}{2}{180}\above{1pt}{1pt}{*}{2}$}%
\nopagebreak\par%
\nopagebreak\par\leavevmode%
{$\left[\!\llap{\phantom{%
\begingroup \smaller\smaller\smaller\begin{tabular}{@{}c@{}}%
0\\0\\0
\end{tabular}\endgroup%
}}\right.$}%
\begingroup \smaller\smaller\smaller\begin{tabular}{@{}c@{}}%
-712080\\0\\3600
\end{tabular}\endgroup%
\kern3pt%
\begingroup \smaller\smaller\smaller\begin{tabular}{@{}c@{}}%
0\\30\\-6
\end{tabular}\endgroup%
\kern3pt%
\begingroup \smaller\smaller\smaller\begin{tabular}{@{}c@{}}%
3600\\-6\\-17
\end{tabular}\endgroup%
{$\left.\llap{\phantom{%
\begingroup \smaller\smaller\smaller\begin{tabular}{@{}c@{}}%
0\\0\\0
\end{tabular}\endgroup%
}}\!\right]$}%
\EasyButWeakLineBreak%
{$\left[\!\llap{\phantom{%
\begingroup \smaller\smaller\smaller\begin{tabular}{@{}c@{}}%
0\\0\\0
\end{tabular}\endgroup%
}}\right.$}%
\begingroup \smaller\smaller\smaller\begin{tabular}{@{}c@{}}%
11\\436\\2168
\end{tabular}\endgroup%
\HardButStrongLineBreak\kern3pt%
\begingroup \smaller\smaller\smaller\begin{tabular}{@{}c@{}}%
21\\830\\4140
\end{tabular}\endgroup%
\HardButStrongLineBreak\kern3pt%
\begingroup \smaller\smaller\smaller\begin{tabular}{@{}c@{}}%
7\\275\\1380
\end{tabular}\endgroup%
\HardButStrongLineBreak\kern3pt%
\begingroup \smaller\smaller\smaller\begin{tabular}{@{}c@{}}%
1\\39\\197
\end{tabular}\endgroup%
\HardButStrongLineBreak\kern3pt%
\begingroup \smaller\smaller\smaller\begin{tabular}{@{}c@{}}%
11\\420\\2160
\end{tabular}\endgroup%
\HardButStrongLineBreak\kern3pt%
\begingroup \smaller\smaller\smaller\begin{tabular}{@{}c@{}}%
-1\\-40\\-198
\end{tabular}\endgroup%
\HardButStrongLineBreak\kern3pt%
\begingroup \smaller\smaller\smaller\begin{tabular}{@{}c@{}}%
-1\\-39\\-198
\end{tabular}\endgroup%
\HardButStrongLineBreak\kern3pt%
\begingroup \smaller\smaller\smaller\begin{tabular}{@{}c@{}}%
5\\200\\984
\end{tabular}\endgroup%
\HardButStrongLineBreak\kern3pt%
\begingroup \smaller\smaller\smaller\begin{tabular}{@{}c@{}}%
37\\1470\\7290
\end{tabular}\endgroup%
{$\left.\llap{\phantom{%
\begingroup \smaller\smaller\smaller\begin{tabular}{@{}c@{}}%
0\\0\\0
\end{tabular}\endgroup%
}}\!\right]$}%
%
%
\hbox{}\par\smallskip%
%
%
\leavevmode%
${L_{206.4}}$%
{} : {$[1\above{1pt}{1pt}{-}{}2\above{1pt}{1pt}{1}{}]\above{1pt}{1pt}{}{2}16\above{1pt}{1pt}{-}{5}{\cdot}1\above{1pt}{1pt}{-}{}3\above{1pt}{1pt}{1}{}9\above{1pt}{1pt}{1}{}{\cdot}1\above{1pt}{1pt}{2}{}5\above{1pt}{1pt}{1}{}$}\spacer%
\instructions{3m,3,m}%
\EasyButWeakLineBreak%
{${9}\above{1pt}{1pt}{r}{2}{120}\above{1pt}{1pt}{24,7}{\infty z}{30}\above{1pt}{1pt}{r}{2}{144}\above{1pt}{1pt}{s}{2}{20}\above{1pt}{1pt}{*}{2}{48}\above{1pt}{1pt}{s}{2}{8}\above{1pt}{1pt}{l}{2}{3}\above{1pt}{1pt}{}{2}{80}\above{1pt}{1pt}{}{2}$}%
\nopagebreak\par%
\nopagebreak\par\leavevmode%
{$\left[\!\llap{\phantom{%
\begingroup \smaller\smaller\smaller
\endgroup%
}}\!\right]$}%
%
%
\hbox{}\par\smallskip%
%
%
\leavevmode%
${L_{206.5}}$%
{} : {$[1\above{1pt}{1pt}{1}{}2\above{1pt}{1pt}{1}{}]\above{1pt}{1pt}{}{6}16\above{1pt}{1pt}{1}{1}{\cdot}1\above{1pt}{1pt}{-}{}3\above{1pt}{1pt}{1}{}9\above{1pt}{1pt}{1}{}{\cdot}1\above{1pt}{1pt}{2}{}5\above{1pt}{1pt}{1}{}$}\spacer%
\instructions{3}%
\EasyButWeakLineBreak%
{${36}\above{1pt}{1pt}{*}{2}{120}\above{1pt}{1pt}{24,19}{\infty z}{30}\above{1pt}{1pt}{}{2}{144}\above{1pt}{1pt}{}{2}{5}\above{1pt}{1pt}{r}{2}{48}\above{1pt}{1pt}{*}{2}{8}\above{1pt}{1pt}{*}{2}{12}\above{1pt}{1pt}{s}{2}{80}\above{1pt}{1pt}{s}{2}$}%
\nopagebreak\par%
\nopagebreak\par\leavevmode%
{$\left[\!\llap{\phantom{%
\begingroup \smaller\smaller\smaller\begin{tabular}{@{}c@{}}%
0\\0\\0
\end{tabular}\endgroup%
}}\right.$}%
\begingroup \smaller\smaller\smaller\begin{tabular}{@{}c@{}}%
-238320\\2880\\0
\end{tabular}\endgroup%
\kern3pt%
\begingroup \smaller\smaller\smaller\begin{tabular}{@{}c@{}}%
2880\\-6\\-12
\end{tabular}\endgroup%
\kern3pt%
\begingroup \smaller\smaller\smaller\begin{tabular}{@{}c@{}}%
0\\-12\\5
\end{tabular}\endgroup%
{$\left.\llap{\phantom{%
\begingroup \smaller\smaller\smaller\begin{tabular}{@{}c@{}}%
0\\0\\0
\end{tabular}\endgroup%
}}\!\right]$}%
\EasyButWeakLineBreak%
{$\left[\!\llap{\phantom{%
\begingroup \smaller\smaller\smaller\begin{tabular}{@{}c@{}}%
0\\0\\0
\end{tabular}\endgroup%
}}\right.$}%
\begingroup \smaller\smaller\smaller\begin{tabular}{@{}c@{}}%
-1\\-84\\-198
\end{tabular}\endgroup%
\HardButStrongLineBreak\kern3pt%
\begingroup \smaller\smaller\smaller\begin{tabular}{@{}c@{}}%
-3\\-250\\-600
\end{tabular}\endgroup%
\HardButStrongLineBreak\kern3pt%
\begingroup \smaller\smaller\smaller\begin{tabular}{@{}c@{}}%
-1\\-85\\-210
\end{tabular}\endgroup%
\HardButStrongLineBreak\kern3pt%
\begingroup \smaller\smaller\smaller\begin{tabular}{@{}c@{}}%
1\\72\\144
\end{tabular}\endgroup%
\HardButStrongLineBreak\kern3pt%
\begingroup \smaller\smaller\smaller\begin{tabular}{@{}c@{}}%
1\\80\\185
\end{tabular}\endgroup%
\HardButStrongLineBreak\kern3pt%
\begingroup \smaller\smaller\smaller\begin{tabular}{@{}c@{}}%
3\\244\\576
\end{tabular}\endgroup%
\HardButStrongLineBreak\kern3pt%
\begingroup \smaller\smaller\smaller\begin{tabular}{@{}c@{}}%
1\\82\\196
\end{tabular}\endgroup%
\HardButStrongLineBreak\kern3pt%
\begingroup \smaller\smaller\smaller\begin{tabular}{@{}c@{}}%
1\\82\\198
\end{tabular}\endgroup%
\HardButStrongLineBreak\kern3pt%
\begingroup \smaller\smaller\smaller\begin{tabular}{@{}c@{}}%
1\\80\\200
\end{tabular}\endgroup%
{$\left.\llap{\phantom{%
\begingroup \smaller\smaller\smaller\begin{tabular}{@{}c@{}}%
0\\0\\0
\end{tabular}\endgroup%
}}\!\right]$}%

\medskip%
%
\leavevmode\llap{}%
$W_{207}$%
\qquad\llap{60} lattices, $\chi=18$%
\hfill%
$22|222\slashtwo2\rtimes D_{2}$%
\nopagebreak\smallskip\hrule\nopagebreak\medskip%
%
%
\leavevmode%
${L_{207.1}}$%
{} : {$1\above{1pt}{1pt}{-2}{4}8\above{1pt}{1pt}{-}{3}{\cdot}1\above{1pt}{1pt}{2}{}3\above{1pt}{1pt}{-}{}{\cdot}1\above{1pt}{1pt}{2}{}5\above{1pt}{1pt}{1}{}$}\EasyButWeakLineBreak%
{${24}\above{1pt}{1pt}{}{2}{5}\above{1pt}{1pt}{r}{2}{8}\above{1pt}{1pt}{s}{2}{20}\above{1pt}{1pt}{*}{2}{24}\above{1pt}{1pt}{*}{2}{4}\above{1pt}{1pt}{l}{2}{1}\above{1pt}{1pt}{}{2}$}%
\nopagebreak\par%
\nopagebreak\par\leavevmode%
{$\left[\!\llap{\phantom{%
\begingroup \smaller\smaller\smaller\begin{tabular}{@{}c@{}}%
0\\0\\0
\end{tabular}\endgroup%
}}\right.$}%
\begingroup \smaller\smaller\smaller\begin{tabular}{@{}c@{}}%
-39720\\480\\0
\end{tabular}\endgroup%
\kern3pt%
\begingroup \smaller\smaller\smaller\begin{tabular}{@{}c@{}}%
480\\-4\\-3
\end{tabular}\endgroup%
\kern3pt%
\begingroup \smaller\smaller\smaller\begin{tabular}{@{}c@{}}%
0\\-3\\5
\end{tabular}\endgroup%
{$\left.\llap{\phantom{%
\begingroup \smaller\smaller\smaller\begin{tabular}{@{}c@{}}%
0\\0\\0
\end{tabular}\endgroup%
}}\!\right]$}%
\EasyButWeakLineBreak%
{$\left[\!\llap{\phantom{%
\begingroup \smaller\smaller\smaller\begin{tabular}{@{}c@{}}%
0\\0\\0
\end{tabular}\endgroup%
}}\right.$}%
\begingroup \smaller\smaller\smaller\begin{tabular}{@{}c@{}}%
5\\408\\240
\end{tabular}\endgroup%
\HardButStrongLineBreak\kern3pt%
\begingroup \smaller\smaller\smaller\begin{tabular}{@{}c@{}}%
1\\80\\45
\end{tabular}\endgroup%
\HardButStrongLineBreak\kern3pt%
\begingroup \smaller\smaller\smaller\begin{tabular}{@{}c@{}}%
-1\\-84\\-52
\end{tabular}\endgroup%
\HardButStrongLineBreak\kern3pt%
\begingroup \smaller\smaller\smaller\begin{tabular}{@{}c@{}}%
-3\\-250\\-150
\end{tabular}\endgroup%
\HardButStrongLineBreak\kern3pt%
\begingroup \smaller\smaller\smaller\begin{tabular}{@{}c@{}}%
-1\\-84\\-48
\end{tabular}\endgroup%
\HardButStrongLineBreak\kern3pt%
\begingroup \smaller\smaller\smaller\begin{tabular}{@{}c@{}}%
1\\82\\50
\end{tabular}\endgroup%
\HardButStrongLineBreak\kern3pt%
\begingroup \smaller\smaller\smaller\begin{tabular}{@{}c@{}}%
1\\82\\49
\end{tabular}\endgroup%
{$\left.\llap{\phantom{%
\begingroup \smaller\smaller\smaller\begin{tabular}{@{}c@{}}%
0\\0\\0
\end{tabular}\endgroup%
}}\!\right]$}%
%
%
\hbox{}\par\smallskip%
%
%
\leavevmode%
${L_{207.2}}$%
{} : {$[1\above{1pt}{1pt}{-}{}2\above{1pt}{1pt}{1}{}]\above{1pt}{1pt}{}{6}16\above{1pt}{1pt}{-}{5}{\cdot}1\above{1pt}{1pt}{2}{}3\above{1pt}{1pt}{-}{}{\cdot}1\above{1pt}{1pt}{2}{}5\above{1pt}{1pt}{1}{}$}\spacer%
\instructions{2}%
\EasyButWeakLineBreak%
{${24}\above{1pt}{1pt}{l}{2}{5}\above{1pt}{1pt}{}{2}{2}\above{1pt}{1pt}{}{2}{80}\above{1pt}{1pt}{}{2}{6}\above{1pt}{1pt}{r}{2}{16}\above{1pt}{1pt}{s}{2}{4}\above{1pt}{1pt}{*}{2}$}%
\nopagebreak\par%
\nopagebreak\par\leavevmode%
{$\left[\!\llap{\phantom{%
\begingroup \smaller\smaller\smaller\begin{tabular}{@{}c@{}}%
0\\0\\0
\end{tabular}\endgroup%
}}\right.$}%
\begingroup \smaller\smaller\smaller\begin{tabular}{@{}c@{}}%
-364080\\2640\\2640
\end{tabular}\endgroup%
\kern3pt%
\begingroup \smaller\smaller\smaller\begin{tabular}{@{}c@{}}%
2640\\-14\\-20
\end{tabular}\endgroup%
\kern3pt%
\begingroup \smaller\smaller\smaller\begin{tabular}{@{}c@{}}%
2640\\-20\\-19
\end{tabular}\endgroup%
{$\left.\llap{\phantom{%
\begingroup \smaller\smaller\smaller\begin{tabular}{@{}c@{}}%
0\\0\\0
\end{tabular}\endgroup%
}}\!\right]$}%
\EasyButWeakLineBreak%
{$\left[\!\llap{\phantom{%
\begingroup \smaller\smaller\smaller\begin{tabular}{@{}c@{}}%
0\\0\\0
\end{tabular}\endgroup%
}}\right.$}%
\begingroup \smaller\smaller\smaller\begin{tabular}{@{}c@{}}%
-1\\-18\\-120
\end{tabular}\endgroup%
\HardButStrongLineBreak\kern3pt%
\begingroup \smaller\smaller\smaller\begin{tabular}{@{}c@{}}%
2\\40\\235
\end{tabular}\endgroup%
\HardButStrongLineBreak\kern3pt%
\begingroup \smaller\smaller\smaller\begin{tabular}{@{}c@{}}%
2\\39\\236
\end{tabular}\endgroup%
\HardButStrongLineBreak\kern3pt%
\begingroup \smaller\smaller\smaller\begin{tabular}{@{}c@{}}%
23\\440\\2720
\end{tabular}\endgroup%
\HardButStrongLineBreak\kern3pt%
\begingroup \smaller\smaller\smaller\begin{tabular}{@{}c@{}}%
4\\75\\474
\end{tabular}\endgroup%
\HardButStrongLineBreak\kern3pt%
\begingroup \smaller\smaller\smaller\begin{tabular}{@{}c@{}}%
1\\16\\120
\end{tabular}\endgroup%
\HardButStrongLineBreak\kern3pt%
\begingroup \smaller\smaller\smaller\begin{tabular}{@{}c@{}}%
-1\\-20\\-118
\end{tabular}\endgroup%
{$\left.\llap{\phantom{%
\begingroup \smaller\smaller\smaller\begin{tabular}{@{}c@{}}%
0\\0\\0
\end{tabular}\endgroup%
}}\!\right]$}%
%
%
\hbox{}\par\smallskip%
%
%
\leavevmode%
${L_{207.3}}$%
{} : {$[1\above{1pt}{1pt}{1}{}2\above{1pt}{1pt}{1}{}]\above{1pt}{1pt}{}{2}16\above{1pt}{1pt}{1}{1}{\cdot}1\above{1pt}{1pt}{2}{}3\above{1pt}{1pt}{-}{}{\cdot}1\above{1pt}{1pt}{2}{}5\above{1pt}{1pt}{1}{}$}\spacer%
\instructions{m}%
\EasyButWeakLineBreak%
{${24}\above{1pt}{1pt}{*}{2}{20}\above{1pt}{1pt}{l}{2}{2}\above{1pt}{1pt}{r}{2}{80}\above{1pt}{1pt}{l}{2}{6}\above{1pt}{1pt}{}{2}{16}\above{1pt}{1pt}{}{2}{1}\above{1pt}{1pt}{r}{2}$}%
\nopagebreak\par%
\nopagebreak\par\leavevmode%
{$\left[\!\llap{\phantom{%
\begingroup \smaller\smaller\smaller\begin{tabular}{@{}c@{}}%
0\\0\\0
\end{tabular}\endgroup%
}}\right.$}%
\begingroup \smaller\smaller\smaller\begin{tabular}{@{}c@{}}%
-4080\\240\\720
\end{tabular}\endgroup%
\kern3pt%
\begingroup \smaller\smaller\smaller\begin{tabular}{@{}c@{}}%
240\\-14\\-40
\end{tabular}\endgroup%
\kern3pt%
\begingroup \smaller\smaller\smaller\begin{tabular}{@{}c@{}}%
720\\-40\\-79
\end{tabular}\endgroup%
{$\left.\llap{\phantom{%
\begingroup \smaller\smaller\smaller\begin{tabular}{@{}c@{}}%
0\\0\\0
\end{tabular}\endgroup%
}}\!\right]$}%
\EasyButWeakLineBreak%
{$\left[\!\llap{\phantom{%
\begingroup \smaller\smaller\smaller\begin{tabular}{@{}c@{}}%
0\\0\\0
\end{tabular}\endgroup%
}}\right.$}%
\begingroup \smaller\smaller\smaller\begin{tabular}{@{}c@{}}%
-1\\-18\\0
\end{tabular}\endgroup%
\HardButStrongLineBreak\kern3pt%
\begingroup \smaller\smaller\smaller\begin{tabular}{@{}c@{}}%
9\\180\\-10
\end{tabular}\endgroup%
\HardButStrongLineBreak\kern3pt%
\begingroup \smaller\smaller\smaller\begin{tabular}{@{}c@{}}%
4\\79\\-4
\end{tabular}\endgroup%
\HardButStrongLineBreak\kern3pt%
\begingroup \smaller\smaller\smaller\begin{tabular}{@{}c@{}}%
43\\840\\-40
\end{tabular}\endgroup%
\HardButStrongLineBreak\kern3pt%
\begingroup \smaller\smaller\smaller\begin{tabular}{@{}c@{}}%
7\\135\\-6
\end{tabular}\endgroup%
\HardButStrongLineBreak\kern3pt%
\begingroup \smaller\smaller\smaller\begin{tabular}{@{}c@{}}%
1\\16\\0
\end{tabular}\endgroup%
\HardButStrongLineBreak\kern3pt%
\begingroup \smaller\smaller\smaller\begin{tabular}{@{}c@{}}%
-1\\-20\\1
\end{tabular}\endgroup%
{$\left.\llap{\phantom{%
\begingroup \smaller\smaller\smaller\begin{tabular}{@{}c@{}}%
0\\0\\0
\end{tabular}\endgroup%
}}\!\right]$}%
%
%
\hbox{}\par\smallskip%
%
%
\leavevmode%
${L_{207.4}}$%
{} : {$[1\above{1pt}{1pt}{1}{}2\above{1pt}{1pt}{-}{}]\above{1pt}{1pt}{}{4}16\above{1pt}{1pt}{-}{3}{\cdot}1\above{1pt}{1pt}{2}{}3\above{1pt}{1pt}{-}{}{\cdot}1\above{1pt}{1pt}{2}{}5\above{1pt}{1pt}{1}{}$}\spacer%
\instructions{m}%
\EasyButWeakLineBreak%
{${6}\above{1pt}{1pt}{r}{2}{20}\above{1pt}{1pt}{*}{2}{8}\above{1pt}{1pt}{*}{2}{80}\above{1pt}{1pt}{s}{2}{24}\above{1pt}{1pt}{*}{2}{16}\above{1pt}{1pt}{l}{2}{1}\above{1pt}{1pt}{}{2}$}%
\nopagebreak\par%
\nopagebreak\par\leavevmode%
{$\left[\!\llap{\phantom{%
\begingroup \smaller\smaller\smaller\begin{tabular}{@{}c@{}}%
0\\0\\0
\end{tabular}\endgroup%
}}\right.$}%
\begingroup \smaller\smaller\smaller\begin{tabular}{@{}c@{}}%
-169680\\1440\\1200
\end{tabular}\endgroup%
\kern3pt%
\begingroup \smaller\smaller\smaller\begin{tabular}{@{}c@{}}%
1440\\-10\\-12
\end{tabular}\endgroup%
\kern3pt%
\begingroup \smaller\smaller\smaller\begin{tabular}{@{}c@{}}%
1200\\-12\\-7
\end{tabular}\endgroup%
{$\left.\llap{\phantom{%
\begingroup \smaller\smaller\smaller\begin{tabular}{@{}c@{}}%
0\\0\\0
\end{tabular}\endgroup%
}}\!\right]$}%
\EasyButWeakLineBreak%
{$\left[\!\llap{\phantom{%
\begingroup \smaller\smaller\smaller\begin{tabular}{@{}c@{}}%
0\\0\\0
\end{tabular}\endgroup%
}}\right.$}%
\begingroup \smaller\smaller\smaller\begin{tabular}{@{}c@{}}%
1\\57\\72
\end{tabular}\endgroup%
\HardButStrongLineBreak\kern3pt%
\begingroup \smaller\smaller\smaller\begin{tabular}{@{}c@{}}%
-1\\-60\\-70
\end{tabular}\endgroup%
\HardButStrongLineBreak\kern3pt%
\begingroup \smaller\smaller\smaller\begin{tabular}{@{}c@{}}%
-1\\-58\\-72
\end{tabular}\endgroup%
\HardButStrongLineBreak\kern3pt%
\begingroup \smaller\smaller\smaller\begin{tabular}{@{}c@{}}%
3\\180\\200
\end{tabular}\endgroup%
\HardButStrongLineBreak\kern3pt%
\begingroup \smaller\smaller\smaller\begin{tabular}{@{}c@{}}%
5\\294\\348
\end{tabular}\endgroup%
\HardButStrongLineBreak\kern3pt%
\begingroup \smaller\smaller\smaller\begin{tabular}{@{}c@{}}%
5\\292\\352
\end{tabular}\endgroup%
\HardButStrongLineBreak\kern3pt%
\begingroup \smaller\smaller\smaller\begin{tabular}{@{}c@{}}%
1\\58\\71
\end{tabular}\endgroup%
{$\left.\llap{\phantom{%
\begingroup \smaller\smaller\smaller\begin{tabular}{@{}c@{}}%
0\\0\\0
\end{tabular}\endgroup%
}}\!\right]$}%
%
%
\hbox{}\par\smallskip%
%
%
\leavevmode%
${L_{207.5}}$%
{} : {$[1\above{1pt}{1pt}{-}{}2\above{1pt}{1pt}{-}{}]\above{1pt}{1pt}{}{0}16\above{1pt}{1pt}{1}{7}{\cdot}1\above{1pt}{1pt}{2}{}3\above{1pt}{1pt}{-}{}{\cdot}1\above{1pt}{1pt}{2}{}5\above{1pt}{1pt}{1}{}$}\EasyButWeakLineBreak%
{${6}\above{1pt}{1pt}{}{2}{5}\above{1pt}{1pt}{r}{2}{8}\above{1pt}{1pt}{s}{2}{80}\above{1pt}{1pt}{*}{2}{24}\above{1pt}{1pt}{s}{2}{16}\above{1pt}{1pt}{*}{2}{4}\above{1pt}{1pt}{l}{2}$}%
\nopagebreak\par%
\nopagebreak\par\leavevmode%
{$\left[\!\llap{\phantom{%
\begingroup \smaller\smaller\smaller\begin{tabular}{@{}c@{}}%
0\\0\\0
\end{tabular}\endgroup%
}}\right.$}%
\begingroup \smaller\smaller\smaller\begin{tabular}{@{}c@{}}%
478320\\11280\\1200
\end{tabular}\endgroup%
\kern3pt%
\begingroup \smaller\smaller\smaller\begin{tabular}{@{}c@{}}%
11280\\266\\28
\end{tabular}\endgroup%
\kern3pt%
\begingroup \smaller\smaller\smaller\begin{tabular}{@{}c@{}}%
1200\\28\\-5
\end{tabular}\endgroup%
{$\left.\llap{\phantom{%
\begingroup \smaller\smaller\smaller\begin{tabular}{@{}c@{}}%
0\\0\\0
\end{tabular}\endgroup%
}}\!\right]$}%
\EasyButWeakLineBreak%
{$\left[\!\llap{\phantom{%
\begingroup \smaller\smaller\smaller\begin{tabular}{@{}c@{}}%
0\\0\\0
\end{tabular}\endgroup%
}}\right.$}%
\begingroup \smaller\smaller\smaller\begin{tabular}{@{}c@{}}%
-5\\213\\-6
\end{tabular}\endgroup%
\HardButStrongLineBreak\kern3pt%
\begingroup \smaller\smaller\smaller\begin{tabular}{@{}c@{}}%
2\\-85\\5
\end{tabular}\endgroup%
\HardButStrongLineBreak\kern3pt%
\begingroup \smaller\smaller\smaller\begin{tabular}{@{}c@{}}%
7\\-298\\12
\end{tabular}\endgroup%
\HardButStrongLineBreak\kern3pt%
\begingroup \smaller\smaller\smaller\begin{tabular}{@{}c@{}}%
23\\-980\\40
\end{tabular}\endgroup%
\HardButStrongLineBreak\kern3pt%
\begingroup \smaller\smaller\smaller\begin{tabular}{@{}c@{}}%
-1\\42\\0
\end{tabular}\endgroup%
\HardButStrongLineBreak\kern3pt%
\begingroup \smaller\smaller\smaller\begin{tabular}{@{}c@{}}%
-11\\468\\-16
\end{tabular}\endgroup%
\HardButStrongLineBreak\kern3pt%
\begingroup \smaller\smaller\smaller\begin{tabular}{@{}c@{}}%
-7\\298\\-10
\end{tabular}\endgroup%
{$\left.\llap{\phantom{%
\begingroup \smaller\smaller\smaller\begin{tabular}{@{}c@{}}%
0\\0\\0
\end{tabular}\endgroup%
}}\!\right]$}%

\medskip%
%
\leavevmode\llap{}%
$W_{208}$%
\qquad\llap{120} lattices, $\chi=48$%
\hfill%
$222|222|222|222|\rtimes D_{4}$%
\nopagebreak\smallskip\hrule\nopagebreak\medskip%
%
%
\leavevmode%
${L_{208.1}}$%
{} : {$1\above{1pt}{1pt}{-2}{4}8\above{1pt}{1pt}{-}{3}{\cdot}1\above{1pt}{1pt}{-}{}3\above{1pt}{1pt}{1}{}9\above{1pt}{1pt}{1}{}{\cdot}1\above{1pt}{1pt}{-2}{}5\above{1pt}{1pt}{-}{}$}\spacer%
\instructions{3}%
\EasyButWeakLineBreak%
{${8}\above{1pt}{1pt}{s}{2}{36}\above{1pt}{1pt}{*}{2}{12}\above{1pt}{1pt}{s}{2}{360}\above{1pt}{1pt}{l}{2}{3}\above{1pt}{1pt}{}{2}{9}\above{1pt}{1pt}{r}{2}$}\relax$\,(\times2)$%
\nopagebreak\par%
\nopagebreak\par\leavevmode%
{$\left[\!\llap{\phantom{%
\begingroup \smaller\smaller\smaller
\endgroup%
}}\!\right]$}%
%
%
\hbox{}\par\smallskip%
%
%
\leavevmode%
${L_{208.2}}$%
{} : {$[1\above{1pt}{1pt}{-}{}2\above{1pt}{1pt}{-}{}]\above{1pt}{1pt}{}{0}16\above{1pt}{1pt}{1}{7}{\cdot}1\above{1pt}{1pt}{-}{}3\above{1pt}{1pt}{1}{}9\above{1pt}{1pt}{1}{}{\cdot}1\above{1pt}{1pt}{-2}{}5\above{1pt}{1pt}{-}{}$}\spacer%
\instructions{3,2}%
\EasyButWeakLineBreak%
{${8}\above{1pt}{1pt}{*}{2}{144}\above{1pt}{1pt}{l}{2}{3}\above{1pt}{1pt}{}{2}{90}\above{1pt}{1pt}{r}{2}{48}\above{1pt}{1pt}{s}{2}{36}\above{1pt}{1pt}{*}{2}$}\relax$\,(\times2)$%
\nopagebreak\par%
\nopagebreak\par\leavevmode%
{$\left[\!\llap{\phantom{%
\begingroup \smaller\smaller\smaller
\endgroup%
}}\!\right]$}%
%
%
\hbox{}\par\smallskip%
%
%
\leavevmode%
${L_{208.3}}$%
{} : {$[1\above{1pt}{1pt}{1}{}2\above{1pt}{1pt}{1}{}]\above{1pt}{1pt}{}{2}16\above{1pt}{1pt}{1}{1}{\cdot}1\above{1pt}{1pt}{-}{}3\above{1pt}{1pt}{1}{}9\above{1pt}{1pt}{1}{}{\cdot}1\above{1pt}{1pt}{-2}{}5\above{1pt}{1pt}{-}{}$}\spacer%
\instructions{32,3,m}%
\EasyButWeakLineBreak%
{${2}\above{1pt}{1pt}{}{2}{144}\above{1pt}{1pt}{}{2}{3}\above{1pt}{1pt}{r}{2}{360}\above{1pt}{1pt}{*}{2}{48}\above{1pt}{1pt}{l}{2}{9}\above{1pt}{1pt}{}{2}$}\relax$\,(\times2)$%
\nopagebreak\par%
\nopagebreak\par\leavevmode%
{$\left[\!\llap{\phantom{%
\begingroup \smaller\smaller\smaller
\endgroup%
}}\!\right]$}%
%
%
\hbox{}\par\smallskip%
%
%
\leavevmode%
${L_{208.4}}$%
{} : {$[1\above{1pt}{1pt}{-}{}2\above{1pt}{1pt}{1}{}]\above{1pt}{1pt}{}{6}16\above{1pt}{1pt}{-}{5}{\cdot}1\above{1pt}{1pt}{-}{}3\above{1pt}{1pt}{1}{}9\above{1pt}{1pt}{1}{}{\cdot}1\above{1pt}{1pt}{-2}{}5\above{1pt}{1pt}{-}{}$}\spacer%
\instructions{3m,3}%
\EasyButWeakLineBreak%
{${2}\above{1pt}{1pt}{r}{2}{144}\above{1pt}{1pt}{s}{2}{12}\above{1pt}{1pt}{*}{2}{360}\above{1pt}{1pt}{s}{2}{48}\above{1pt}{1pt}{*}{2}{36}\above{1pt}{1pt}{l}{2}$}\relax$\,(\times2)$%
\nopagebreak\par%
\nopagebreak\par\leavevmode%
{$\left[\!\llap{\phantom{%
\begingroup \smaller\smaller\smaller
\endgroup%
}}\!\right]$}%
%
%
\hbox{}\par\smallskip%
%
%
\leavevmode%
${L_{208.5}}$%
{} : {$[1\above{1pt}{1pt}{1}{}2\above{1pt}{1pt}{-}{}]\above{1pt}{1pt}{}{4}16\above{1pt}{1pt}{-}{3}{\cdot}1\above{1pt}{1pt}{-}{}3\above{1pt}{1pt}{1}{}9\above{1pt}{1pt}{1}{}{\cdot}1\above{1pt}{1pt}{-2}{}5\above{1pt}{1pt}{-}{}$}\spacer%
\instructions{3m,3,m}%
\EasyButWeakLineBreak%
{${8}\above{1pt}{1pt}{s}{2}{144}\above{1pt}{1pt}{*}{2}{12}\above{1pt}{1pt}{l}{2}{90}\above{1pt}{1pt}{}{2}{48}\above{1pt}{1pt}{}{2}{9}\above{1pt}{1pt}{r}{2}$}\relax$\,(\times2)$%
\nopagebreak\par%
\nopagebreak\par\leavevmode%
{$\left[\!\llap{\phantom{%
\begingroup \smaller\smaller\smaller
\endgroup%
}}\!\right]$}%

\medskip%
%
\leavevmode\llap{}%
$W_{209}$%
\qquad\llap{11} lattices, $\chi=18$%
\hfill%
$\slashtwo2\infty|\infty2\rtimes D_{2}$%
\nopagebreak\smallskip\hrule\nopagebreak\medskip%
%
%
\leavevmode%
${L_{209.1}}$%
{} : {$1\above{1pt}{1pt}{2}{{\rm II}}4\above{1pt}{1pt}{-}{5}{\cdot}1\above{1pt}{1pt}{-}{}5\above{1pt}{1pt}{1}{}25\above{1pt}{1pt}{-}{}$}\spacer%
\instructions{2}%
\EasyButWeakLineBreak%
{${50}\above{1pt}{1pt}{b}{2}{2}\above{1pt}{1pt}{l}{2}{20}\above{1pt}{1pt}{5,4}{\infty}{20}\above{1pt}{1pt}{5,2}{\infty b}{20}\above{1pt}{1pt}{r}{2}$}%
\nopagebreak\par%
\nopagebreak\par\leavevmode%
{$\left[\!\llap{\phantom{%
\begingroup \smaller\smaller\smaller\begin{tabular}{@{}c@{}}%
0\\0\\0
\end{tabular}\endgroup%
}}\right.$}%
\begingroup \smaller\smaller\smaller\begin{tabular}{@{}c@{}}%
9300\\4800\\-1300
\end{tabular}\endgroup%
\kern3pt%
\begingroup \smaller\smaller\smaller\begin{tabular}{@{}c@{}}%
4800\\2480\\-675
\end{tabular}\endgroup%
\kern3pt%
\begingroup \smaller\smaller\smaller\begin{tabular}{@{}c@{}}%
-1300\\-675\\188
\end{tabular}\endgroup%
{$\left.\llap{\phantom{%
\begingroup \smaller\smaller\smaller\begin{tabular}{@{}c@{}}%
0\\0\\0
\end{tabular}\endgroup%
}}\!\right]$}%
\EasyButWeakLineBreak%
{$\left[\!\llap{\phantom{%
\begingroup \smaller\smaller\smaller\begin{tabular}{@{}c@{}}%
0\\0\\0
\end{tabular}\endgroup%
}}\right.$}%
\begingroup \smaller\smaller\smaller\begin{tabular}{@{}c@{}}%
-54\\125\\75
\end{tabular}\endgroup%
\HardButStrongLineBreak\kern3pt%
\begingroup \smaller\smaller\smaller\begin{tabular}{@{}c@{}}%
-9\\21\\13
\end{tabular}\endgroup%
\HardButStrongLineBreak\kern3pt%
\begingroup \smaller\smaller\smaller\begin{tabular}{@{}c@{}}%
-39\\92\\60
\end{tabular}\endgroup%
\HardButStrongLineBreak\kern3pt%
\begingroup \smaller\smaller\smaller\begin{tabular}{@{}c@{}}%
1\\-2\\0
\end{tabular}\endgroup%
\HardButStrongLineBreak\kern3pt%
\begingroup \smaller\smaller\smaller\begin{tabular}{@{}c@{}}%
-57\\132\\80
\end{tabular}\endgroup%
{$\left.\llap{\phantom{%
\begingroup \smaller\smaller\smaller\begin{tabular}{@{}c@{}}%
0\\0\\0
\end{tabular}\endgroup%
}}\!\right]$}%
%
%
\hbox{}\par\smallskip%
%
%
\leavevmode%
${L_{209.2}}$%
{} : {$1\above{1pt}{1pt}{-2}{2}8\above{1pt}{1pt}{1}{7}{\cdot}1\above{1pt}{1pt}{1}{}5\above{1pt}{1pt}{-}{}25\above{1pt}{1pt}{1}{}$}\spacer%
\instructions{2}%
\EasyButWeakLineBreak%
{${100}\above{1pt}{1pt}{*}{2}{4}\above{1pt}{1pt}{s}{2}{40}\above{1pt}{1pt}{20,9}{\infty z}{10}\above{1pt}{1pt}{20,7}{\infty a}{40}\above{1pt}{1pt}{s}{2}$}%
\nopagebreak\par%
\nopagebreak\par\leavevmode%
{$\left[\!\llap{\phantom{%
\begingroup \smaller\smaller\smaller\begin{tabular}{@{}c@{}}%
0\\0\\0
\end{tabular}\endgroup%
}}\right.$}%
\begingroup \smaller\smaller\smaller\begin{tabular}{@{}c@{}}%
-302600\\-148800\\-52000
\end{tabular}\endgroup%
\kern3pt%
\begingroup \smaller\smaller\smaller\begin{tabular}{@{}c@{}}%
-148800\\-73165\\-25555
\end{tabular}\endgroup%
\kern3pt%
\begingroup \smaller\smaller\smaller\begin{tabular}{@{}c@{}}%
-52000\\-25555\\-8894
\end{tabular}\endgroup%
{$\left.\llap{\phantom{%
\begingroup \smaller\smaller\smaller\begin{tabular}{@{}c@{}}%
0\\0\\0
\end{tabular}\endgroup%
}}\!\right]$}%
\EasyButWeakLineBreak%
{$\left[\!\llap{\phantom{%
\begingroup \smaller\smaller\smaller\begin{tabular}{@{}c@{}}%
0\\0\\0
\end{tabular}\endgroup%
}}\right.$}%
\begingroup \smaller\smaller\smaller\begin{tabular}{@{}c@{}}%
-407\\950\\-350
\end{tabular}\endgroup%
\HardButStrongLineBreak\kern3pt%
\begingroup \smaller\smaller\smaller\begin{tabular}{@{}c@{}}%
85\\-198\\72
\end{tabular}\endgroup%
\HardButStrongLineBreak\kern3pt%
\begingroup \smaller\smaller\smaller\begin{tabular}{@{}c@{}}%
1453\\-3388\\1240
\end{tabular}\endgroup%
\HardButStrongLineBreak\kern3pt%
\begingroup \smaller\smaller\smaller\begin{tabular}{@{}c@{}}%
204\\-476\\175
\end{tabular}\endgroup%
\HardButStrongLineBreak\kern3pt%
\begingroup \smaller\smaller\smaller\begin{tabular}{@{}c@{}}%
-211\\492\\-180
\end{tabular}\endgroup%
{$\left.\llap{\phantom{%
\begingroup \smaller\smaller\smaller\begin{tabular}{@{}c@{}}%
0\\0\\0
\end{tabular}\endgroup%
}}\!\right]$}%
%
%
\hbox{}\par\smallskip%
%
%
\leavevmode%
${L_{209.3}}$%
{} : {$1\above{1pt}{1pt}{2}{2}8\above{1pt}{1pt}{-}{3}{\cdot}1\above{1pt}{1pt}{1}{}5\above{1pt}{1pt}{-}{}25\above{1pt}{1pt}{1}{}$}\spacer%
\instructions{m}%
\EasyButWeakLineBreak%
{${25}\above{1pt}{1pt}{}{2}{1}\above{1pt}{1pt}{r}{2}{40}\above{1pt}{1pt}{20,19}{\infty z}{10}\above{1pt}{1pt}{20,7}{\infty b}{40}\above{1pt}{1pt}{l}{2}$}%
\nopagebreak\par%
\nopagebreak\par\leavevmode%
{$\left[\!\llap{\phantom{%
\begingroup \smaller\smaller\smaller\begin{tabular}{@{}c@{}}%
0\\0\\0
\end{tabular}\endgroup%
}}\right.$}%
\begingroup \smaller\smaller\smaller\begin{tabular}{@{}c@{}}%
-74600\\-3800\\-10600
\end{tabular}\endgroup%
\kern3pt%
\begingroup \smaller\smaller\smaller\begin{tabular}{@{}c@{}}%
-3800\\-190\\-545
\end{tabular}\endgroup%
\kern3pt%
\begingroup \smaller\smaller\smaller\begin{tabular}{@{}c@{}}%
-10600\\-545\\-1499
\end{tabular}\endgroup%
{$\left.\llap{\phantom{%
\begingroup \smaller\smaller\smaller\begin{tabular}{@{}c@{}}%
0\\0\\0
\end{tabular}\endgroup%
}}\!\right]$}%
\EasyButWeakLineBreak%
{$\left[\!\llap{\phantom{%
\begingroup \smaller\smaller\smaller\begin{tabular}{@{}c@{}}%
0\\0\\0
\end{tabular}\endgroup%
}}\right.$}%
\begingroup \smaller\smaller\smaller\begin{tabular}{@{}c@{}}%
-16\\105\\75
\end{tabular}\endgroup%
\HardButStrongLineBreak\kern3pt%
\begingroup \smaller\smaller\smaller\begin{tabular}{@{}c@{}}%
2\\-14\\-9
\end{tabular}\endgroup%
\HardButStrongLineBreak\kern3pt%
\begingroup \smaller\smaller\smaller\begin{tabular}{@{}c@{}}%
91\\-612\\-420
\end{tabular}\endgroup%
\HardButStrongLineBreak\kern3pt%
\begingroup \smaller\smaller\smaller\begin{tabular}{@{}c@{}}%
15\\-99\\-70
\end{tabular}\endgroup%
\HardButStrongLineBreak\kern3pt%
\begingroup \smaller\smaller\smaller\begin{tabular}{@{}c@{}}%
-13\\88\\60
\end{tabular}\endgroup%
{$\left.\llap{\phantom{%
\begingroup \smaller\smaller\smaller\begin{tabular}{@{}c@{}}%
0\\0\\0
\end{tabular}\endgroup%
}}\!\right]$}%

\medskip%
%
\leavevmode\llap{}%
$W_{210}$%
\qquad\llap{6} lattices, $\chi=12$%
\hfill%
$222222\rtimes C_{2}$%
\nopagebreak\smallskip\hrule\nopagebreak\medskip%
%
%
\leavevmode%
${L_{210.1}}$%
{} : {$1\above{1pt}{1pt}{-2}{{\rm II}}4\above{1pt}{1pt}{1}{1}{\cdot}1\above{1pt}{1pt}{1}{}5\above{1pt}{1pt}{-}{}25\above{1pt}{1pt}{-}{}$}\spacer%
\instructions{2}%
\EasyButWeakLineBreak%
{${4}\above{1pt}{1pt}{r}{2}{50}\above{1pt}{1pt}{b}{2}{10}\above{1pt}{1pt}{l}{2}$}\relax$\,(\times2)$%
\nopagebreak\par%
\nopagebreak\par\leavevmode%
{$\left[\!\llap{\phantom{%
\begingroup \smaller\smaller\smaller\begin{tabular}{@{}c@{}}%
0\\0\\0
\end{tabular}\endgroup%
}}\right.$}%
\begingroup \smaller\smaller\smaller\begin{tabular}{@{}c@{}}%
5700\\2400\\200
\end{tabular}\endgroup%
\kern3pt%
\begingroup \smaller\smaller\smaller\begin{tabular}{@{}c@{}}%
2400\\1010\\85
\end{tabular}\endgroup%
\kern3pt%
\begingroup \smaller\smaller\smaller\begin{tabular}{@{}c@{}}%
200\\85\\6
\end{tabular}\endgroup%
{$\left.\llap{\phantom{%
\begingroup \smaller\smaller\smaller\begin{tabular}{@{}c@{}}%
0\\0\\0
\end{tabular}\endgroup%
}}\!\right]$}%
\hfil\penalty500%
{$\left[\!\llap{\phantom{%
\begingroup \smaller\smaller\smaller\begin{tabular}{@{}c@{}}%
0\\0\\0
\end{tabular}\endgroup%
}}\right.$}%
\begingroup \smaller\smaller\smaller\begin{tabular}{@{}c@{}}%
19\\-40\\-100
\end{tabular}\endgroup%
\kern3pt%
\begingroup \smaller\smaller\smaller\begin{tabular}{@{}c@{}}%
9\\-19\\-45
\end{tabular}\endgroup%
\kern3pt%
\begingroup \smaller\smaller\smaller\begin{tabular}{@{}c@{}}%
0\\0\\-1
\end{tabular}\endgroup%
{$\left.\llap{\phantom{%
\begingroup \smaller\smaller\smaller\begin{tabular}{@{}c@{}}%
0\\0\\0
\end{tabular}\endgroup%
}}\!\right]$}%
\EasyButWeakLineBreak%
{$\left[\!\llap{\phantom{%
\begingroup \smaller\smaller\smaller\begin{tabular}{@{}c@{}}%
0\\0\\0
\end{tabular}\endgroup%
}}\right.$}%
\begingroup \smaller\smaller\smaller\begin{tabular}{@{}c@{}}%
11\\-24\\-28
\end{tabular}\endgroup%
\HardButStrongLineBreak\kern3pt%
\begingroup \smaller\smaller\smaller\begin{tabular}{@{}c@{}}%
7\\-15\\-25
\end{tabular}\endgroup%
\HardButStrongLineBreak\kern3pt%
\begingroup \smaller\smaller\smaller\begin{tabular}{@{}c@{}}%
-4\\9\\5
\end{tabular}\endgroup%
{$\left.\llap{\phantom{%
\begingroup \smaller\smaller\smaller\begin{tabular}{@{}c@{}}%
0\\0\\0
\end{tabular}\endgroup%
}}\!\right]$}%

\medskip%
%
\leavevmode\llap{}%
$W_{211}$%
\qquad\llap{22} lattices, $\chi=12$%
\hfill%
$222|222|\rtimes D_{2}$%
\nopagebreak\smallskip\hrule\nopagebreak\medskip%
%
%
\leavevmode%
${L_{211.1}}$%
{} : {$1\above{1pt}{1pt}{2}{{\rm II}}4\above{1pt}{1pt}{1}{7}{\cdot}1\above{1pt}{1pt}{-}{}3\above{1pt}{1pt}{1}{}9\above{1pt}{1pt}{-}{}{\cdot}1\above{1pt}{1pt}{-2}{}5\above{1pt}{1pt}{1}{}$}\spacer%
\instructions{2}%
\EasyButWeakLineBreak%
{${12}\above{1pt}{1pt}{*}{2}{20}\above{1pt}{1pt}{b}{2}{18}\above{1pt}{1pt}{s}{2}{30}\above{1pt}{1pt}{s}{2}{2}\above{1pt}{1pt}{b}{2}{180}\above{1pt}{1pt}{*}{2}$}%
\nopagebreak\par%
\nopagebreak\par\leavevmode%
{$\left[\!\llap{\phantom{%
\begingroup \smaller\smaller\smaller\begin{tabular}{@{}c@{}}%
0\\0\\0
\end{tabular}\endgroup%
}}\right.$}%
\begingroup \smaller\smaller\smaller\begin{tabular}{@{}c@{}}%
3420\\-1260\\0
\end{tabular}\endgroup%
\kern3pt%
\begingroup \smaller\smaller\smaller\begin{tabular}{@{}c@{}}%
-1260\\462\\3
\end{tabular}\endgroup%
\kern3pt%
\begingroup \smaller\smaller\smaller\begin{tabular}{@{}c@{}}%
0\\3\\-4
\end{tabular}\endgroup%
{$\left.\llap{\phantom{%
\begingroup \smaller\smaller\smaller\begin{tabular}{@{}c@{}}%
0\\0\\0
\end{tabular}\endgroup%
}}\!\right]$}%
\EasyButWeakLineBreak%
{$\left[\!\llap{\phantom{%
\begingroup \smaller\smaller\smaller\begin{tabular}{@{}c@{}}%
0\\0\\0
\end{tabular}\endgroup%
}}\right.$}%
\begingroup \smaller\smaller\smaller\begin{tabular}{@{}c@{}}%
3\\8\\6
\end{tabular}\endgroup%
\HardButStrongLineBreak\kern3pt%
\begingroup \smaller\smaller\smaller\begin{tabular}{@{}c@{}}%
11\\30\\20
\end{tabular}\endgroup%
\HardButStrongLineBreak\kern3pt%
\begingroup \smaller\smaller\smaller\begin{tabular}{@{}c@{}}%
1\\3\\0
\end{tabular}\endgroup%
\HardButStrongLineBreak\kern3pt%
\begingroup \smaller\smaller\smaller\begin{tabular}{@{}c@{}}%
-13\\-35\\-30
\end{tabular}\endgroup%
\HardButStrongLineBreak\kern3pt%
\begingroup \smaller\smaller\smaller\begin{tabular}{@{}c@{}}%
-7\\-19\\-16
\end{tabular}\endgroup%
\HardButStrongLineBreak\kern3pt%
\begingroup \smaller\smaller\smaller\begin{tabular}{@{}c@{}}%
-77\\-210\\-180
\end{tabular}\endgroup%
{$\left.\llap{\phantom{%
\begingroup \smaller\smaller\smaller\begin{tabular}{@{}c@{}}%
0\\0\\0
\end{tabular}\endgroup%
}}\!\right]$}%
%
%
\hbox{}\par\smallskip%
%
%
\leavevmode%
${L_{211.2}}$%
{} : {$1\above{1pt}{1pt}{2}{6}8\above{1pt}{1pt}{1}{1}{\cdot}1\above{1pt}{1pt}{1}{}3\above{1pt}{1pt}{-}{}9\above{1pt}{1pt}{1}{}{\cdot}1\above{1pt}{1pt}{-2}{}5\above{1pt}{1pt}{-}{}$}\spacer%
\instructions{2}%
\EasyButWeakLineBreak%
{${6}\above{1pt}{1pt}{b}{2}{40}\above{1pt}{1pt}{*}{2}{36}\above{1pt}{1pt}{l}{2}{15}\above{1pt}{1pt}{r}{2}{4}\above{1pt}{1pt}{*}{2}{360}\above{1pt}{1pt}{b}{2}$}%
\nopagebreak\par%
\nopagebreak\par\leavevmode%
{$\left[\!\llap{\phantom{%
\begingroup \smaller\smaller\smaller\begin{tabular}{@{}c@{}}%
0\\0\\0
\end{tabular}\endgroup%
}}\right.$}%
\begingroup \smaller\smaller\smaller\begin{tabular}{@{}c@{}}%
-1490040\\-2880\\5040
\end{tabular}\endgroup%
\kern3pt%
\begingroup \smaller\smaller\smaller\begin{tabular}{@{}c@{}}%
-2880\\6\\9
\end{tabular}\endgroup%
\kern3pt%
\begingroup \smaller\smaller\smaller\begin{tabular}{@{}c@{}}%
5040\\9\\-17
\end{tabular}\endgroup%
{$\left.\llap{\phantom{%
\begingroup \smaller\smaller\smaller\begin{tabular}{@{}c@{}}%
0\\0\\0
\end{tabular}\endgroup%
}}\!\right]$}%
\EasyButWeakLineBreak%
{$\left[\!\llap{\phantom{%
\begingroup \smaller\smaller\smaller\begin{tabular}{@{}c@{}}%
0\\0\\0
\end{tabular}\endgroup%
}}\right.$}%
\begingroup \smaller\smaller\smaller\begin{tabular}{@{}c@{}}%
1\\19\\306
\end{tabular}\endgroup%
\HardButStrongLineBreak\kern3pt%
\begingroup \smaller\smaller\smaller\begin{tabular}{@{}c@{}}%
3\\60\\920
\end{tabular}\endgroup%
\HardButStrongLineBreak\kern3pt%
\begingroup \smaller\smaller\smaller\begin{tabular}{@{}c@{}}%
-1\\-18\\-306
\end{tabular}\endgroup%
\HardButStrongLineBreak\kern3pt%
\begingroup \smaller\smaller\smaller\begin{tabular}{@{}c@{}}%
-2\\-40\\-615
\end{tabular}\endgroup%
\HardButStrongLineBreak\kern3pt%
\begingroup \smaller\smaller\smaller\begin{tabular}{@{}c@{}}%
-1\\-22\\-310
\end{tabular}\endgroup%
\HardButStrongLineBreak\kern3pt%
\begingroup \smaller\smaller\smaller\begin{tabular}{@{}c@{}}%
-1\\-60\\-360
\end{tabular}\endgroup%
{$\left.\llap{\phantom{%
\begingroup \smaller\smaller\smaller\begin{tabular}{@{}c@{}}%
0\\0\\0
\end{tabular}\endgroup%
}}\!\right]$}%
%
%
\hbox{}\par\smallskip%
%
%
\leavevmode%
${L_{211.3}}$%
{} : {$1\above{1pt}{1pt}{-2}{6}8\above{1pt}{1pt}{-}{5}{\cdot}1\above{1pt}{1pt}{1}{}3\above{1pt}{1pt}{-}{}9\above{1pt}{1pt}{1}{}{\cdot}1\above{1pt}{1pt}{-2}{}5\above{1pt}{1pt}{-}{}$}\spacer%
\instructions{m}%
\EasyButWeakLineBreak%
{${6}\above{1pt}{1pt}{l}{2}{40}\above{1pt}{1pt}{}{2}{9}\above{1pt}{1pt}{r}{2}{60}\above{1pt}{1pt}{l}{2}{1}\above{1pt}{1pt}{}{2}{360}\above{1pt}{1pt}{r}{2}$}%
\nopagebreak\par%
\nopagebreak\par\leavevmode%
{$\left[\!\llap{\phantom{%
\begingroup \smaller\smaller\smaller\begin{tabular}{@{}c@{}}%
0\\0\\0
\end{tabular}\endgroup%
}}\right.$}%
\begingroup \smaller\smaller\smaller\begin{tabular}{@{}c@{}}%
85251240\\-2185920\\56160
\end{tabular}\endgroup%
\kern3pt%
\begingroup \smaller\smaller\smaller\begin{tabular}{@{}c@{}}%
-2185920\\56049\\-1440
\end{tabular}\endgroup%
\kern3pt%
\begingroup \smaller\smaller\smaller\begin{tabular}{@{}c@{}}%
56160\\-1440\\37
\end{tabular}\endgroup%
{$\left.\llap{\phantom{%
\begingroup \smaller\smaller\smaller\begin{tabular}{@{}c@{}}%
0\\0\\0
\end{tabular}\endgroup%
}}\!\right]$}%
\EasyButWeakLineBreak%
{$\left[\!\llap{\phantom{%
\begingroup \smaller\smaller\smaller\begin{tabular}{@{}c@{}}%
0\\0\\0
\end{tabular}\endgroup%
}}\right.$}%
\begingroup \smaller\smaller\smaller\begin{tabular}{@{}c@{}}%
2\\79\\39
\end{tabular}\endgroup%
\HardButStrongLineBreak\kern3pt%
\begingroup \smaller\smaller\smaller\begin{tabular}{@{}c@{}}%
1\\40\\40
\end{tabular}\endgroup%
\HardButStrongLineBreak\kern3pt%
\begingroup \smaller\smaller\smaller\begin{tabular}{@{}c@{}}%
-1\\-39\\0
\end{tabular}\endgroup%
\HardButStrongLineBreak\kern3pt%
\begingroup \smaller\smaller\smaller\begin{tabular}{@{}c@{}}%
7\\280\\270
\end{tabular}\endgroup%
\HardButStrongLineBreak\kern3pt%
\begingroup \smaller\smaller\smaller\begin{tabular}{@{}c@{}}%
4\\159\\116
\end{tabular}\endgroup%
\HardButStrongLineBreak\kern3pt%
\begingroup \smaller\smaller\smaller\begin{tabular}{@{}c@{}}%
133\\5280\\3600
\end{tabular}\endgroup%
{$\left.\llap{\phantom{%
\begingroup \smaller\smaller\smaller\begin{tabular}{@{}c@{}}%
0\\0\\0
\end{tabular}\endgroup%
}}\!\right]$}%

\medskip%
%
\leavevmode\llap{}%
$W_{212}$%
\qquad\llap{9} lattices, $\chi=12$%
\hfill%
$2\slashtwo2|2\slashtwo2|\rtimes D_{4}$%
\nopagebreak\smallskip\hrule\nopagebreak\medskip%
%
%
\leavevmode%
${L_{212.1}}$%
{} : {$1\above{1pt}{1pt}{-2}{{\rm II}}4\above{1pt}{1pt}{-}{3}{\cdot}1\above{1pt}{1pt}{-}{}3\above{1pt}{1pt}{1}{}9\above{1pt}{1pt}{-}{}{\cdot}1\above{1pt}{1pt}{-}{}5\above{1pt}{1pt}{-}{}25\above{1pt}{1pt}{-}{}$}\spacer%
\instructions{25,5,2}%
\EasyButWeakLineBreak%
{${12}\above{1pt}{1pt}{r}{2}{50}\above{1pt}{1pt}{b}{2}{18}\above{1pt}{1pt}{l}{2}{300}\above{1pt}{1pt}{r}{2}{2}\above{1pt}{1pt}{b}{2}{450}\above{1pt}{1pt}{l}{2}$}%
\nopagebreak\par%
\nopagebreak\par\leavevmode%
{$\left[\!\llap{\phantom{%
\begingroup \smaller\smaller\smaller\begin{tabular}{@{}c@{}}%
0\\0\\0
\end{tabular}\endgroup%
}}\right.$}%
\begingroup \smaller\smaller\smaller\begin{tabular}{@{}c@{}}%
398700\\-126000\\-900
\end{tabular}\endgroup%
\kern3pt%
\begingroup \smaller\smaller\smaller\begin{tabular}{@{}c@{}}%
-126000\\39810\\285
\end{tabular}\endgroup%
\kern3pt%
\begingroup \smaller\smaller\smaller\begin{tabular}{@{}c@{}}%
-900\\285\\2
\end{tabular}\endgroup%
{$\left.\llap{\phantom{%
\begingroup \smaller\smaller\smaller\begin{tabular}{@{}c@{}}%
0\\0\\0
\end{tabular}\endgroup%
}}\!\right]$}%
\EasyButWeakLineBreak%
{$\left[\!\llap{\phantom{%
\begingroup \smaller\smaller\smaller\begin{tabular}{@{}c@{}}%
0\\0\\0
\end{tabular}\endgroup%
}}\right.$}%
\begingroup \smaller\smaller\smaller\begin{tabular}{@{}c@{}}%
-3\\-8\\-204
\end{tabular}\endgroup%
\HardButStrongLineBreak\kern3pt%
\begingroup \smaller\smaller\smaller\begin{tabular}{@{}c@{}}%
-2\\-5\\-175
\end{tabular}\endgroup%
\HardButStrongLineBreak\kern3pt%
\begingroup \smaller\smaller\smaller\begin{tabular}{@{}c@{}}%
1\\3\\27
\end{tabular}\endgroup%
\HardButStrongLineBreak\kern3pt%
\begingroup \smaller\smaller\smaller\begin{tabular}{@{}c@{}}%
7\\20\\300
\end{tabular}\endgroup%
\HardButStrongLineBreak\kern3pt%
\begingroup \smaller\smaller\smaller\begin{tabular}{@{}c@{}}%
0\\0\\-1
\end{tabular}\endgroup%
\HardButStrongLineBreak\kern3pt%
\begingroup \smaller\smaller\smaller\begin{tabular}{@{}c@{}}%
-11\\-30\\-675
\end{tabular}\endgroup%
{$\left.\llap{\phantom{%
\begingroup \smaller\smaller\smaller\begin{tabular}{@{}c@{}}%
0\\0\\0
\end{tabular}\endgroup%
}}\!\right]$}%

\medskip%
%
\leavevmode\llap{}%
$W_{213}$%
\qquad\llap{6} lattices, $\chi=6$%
\hfill%
$2|22\slashtwo2\rtimes D_{2}$%
\nopagebreak\smallskip\hrule\nopagebreak\medskip%
%
%
\leavevmode%
${L_{213.1}}$%
{} : {$1\above{1pt}{1pt}{-2}{{\rm II}}4\above{1pt}{1pt}{-}{3}{\cdot}1\above{1pt}{1pt}{1}{}3\above{1pt}{1pt}{1}{}9\above{1pt}{1pt}{1}{}{\cdot}1\above{1pt}{1pt}{2}{}5\above{1pt}{1pt}{-}{}$}\spacer%
\instructions{2}%
\EasyButWeakLineBreak%
{${90}\above{1pt}{1pt}{l}{2}{12}\above{1pt}{1pt}{r}{2}{10}\above{1pt}{1pt}{b}{2}{36}\above{1pt}{1pt}{*}{2}{4}\above{1pt}{1pt}{b}{2}$}%
\nopagebreak\par%
\nopagebreak\par\leavevmode%
{$\left[\!\llap{\phantom{%
\begingroup \smaller\smaller\smaller\begin{tabular}{@{}c@{}}%
0\\0\\0
\end{tabular}\endgroup%
}}\right.$}%
\begingroup \smaller\smaller\smaller\begin{tabular}{@{}c@{}}%
47340\\900\\-720
\end{tabular}\endgroup%
\kern3pt%
\begingroup \smaller\smaller\smaller\begin{tabular}{@{}c@{}}%
900\\-6\\-9
\end{tabular}\endgroup%
\kern3pt%
\begingroup \smaller\smaller\smaller\begin{tabular}{@{}c@{}}%
-720\\-9\\10
\end{tabular}\endgroup%
{$\left.\llap{\phantom{%
\begingroup \smaller\smaller\smaller\begin{tabular}{@{}c@{}}%
0\\0\\0
\end{tabular}\endgroup%
}}\!\right]$}%
\EasyButWeakLineBreak%
{$\left[\!\llap{\phantom{%
\begingroup \smaller\smaller\smaller\begin{tabular}{@{}c@{}}%
0\\0\\0
\end{tabular}\endgroup%
}}\right.$}%
\begingroup \smaller\smaller\smaller\begin{tabular}{@{}c@{}}%
1\\15\\90
\end{tabular}\endgroup%
\HardButStrongLineBreak\kern3pt%
\begingroup \smaller\smaller\smaller\begin{tabular}{@{}c@{}}%
3\\52\\264
\end{tabular}\endgroup%
\HardButStrongLineBreak\kern3pt%
\begingroup \smaller\smaller\smaller\begin{tabular}{@{}c@{}}%
2\\35\\175
\end{tabular}\endgroup%
\HardButStrongLineBreak\kern3pt%
\begingroup \smaller\smaller\smaller\begin{tabular}{@{}c@{}}%
-1\\-18\\-90
\end{tabular}\endgroup%
\HardButStrongLineBreak\kern3pt%
\begingroup \smaller\smaller\smaller\begin{tabular}{@{}c@{}}%
-1\\-18\\-88
\end{tabular}\endgroup%
{$\left.\llap{\phantom{%
\begingroup \smaller\smaller\smaller\begin{tabular}{@{}c@{}}%
0\\0\\0
\end{tabular}\endgroup%
}}\!\right]$}%

\medskip%
%
\leavevmode\llap{}%
$W_{214}$%
\qquad\llap{26} lattices, $\chi=24$%
\hfill%
$22|22|22|22|\rtimes D_{4}$%
\nopagebreak\smallskip\hrule\nopagebreak\medskip%
%
%
\leavevmode%
${L_{214.1}}$%
{} : {$1\above{1pt}{1pt}{2}{{\rm II}}4\above{1pt}{1pt}{1}{7}{\cdot}1\above{1pt}{1pt}{1}{}3\above{1pt}{1pt}{1}{}9\above{1pt}{1pt}{1}{}{\cdot}1\above{1pt}{1pt}{-2}{}5\above{1pt}{1pt}{1}{}$}\spacer%
\instructions{2}%
\EasyButWeakLineBreak%
{${12}\above{1pt}{1pt}{*}{2}{4}\above{1pt}{1pt}{b}{2}{30}\above{1pt}{1pt}{b}{2}{36}\above{1pt}{1pt}{*}{2}$}\relax$\,(\times2)$%
\nopagebreak\par%
\nopagebreak\par\leavevmode%
{$\left[\!\llap{\phantom{%
\begingroup \smaller\smaller\smaller\begin{tabular}{@{}c@{}}%
0\\0\\0
\end{tabular}\endgroup%
}}\right.$}%
\begingroup \smaller\smaller\smaller\begin{tabular}{@{}c@{}}%
23580\\-540\\180
\end{tabular}\endgroup%
\kern3pt%
\begingroup \smaller\smaller\smaller\begin{tabular}{@{}c@{}}%
-540\\12\\-3
\end{tabular}\endgroup%
\kern3pt%
\begingroup \smaller\smaller\smaller\begin{tabular}{@{}c@{}}%
180\\-3\\-2
\end{tabular}\endgroup%
{$\left.\llap{\phantom{%
\begingroup \smaller\smaller\smaller\begin{tabular}{@{}c@{}}%
0\\0\\0
\end{tabular}\endgroup%
}}\!\right]$}%
\hfil\penalty500%
{$\left[\!\llap{\phantom{%
\begingroup \smaller\smaller\smaller\begin{tabular}{@{}c@{}}%
0\\0\\0
\end{tabular}\endgroup%
}}\right.$}%
\begingroup \smaller\smaller\smaller\begin{tabular}{@{}c@{}}%
-1\\360\\540
\end{tabular}\endgroup%
\kern3pt%
\begingroup \smaller\smaller\smaller\begin{tabular}{@{}c@{}}%
0\\-11\\-15
\end{tabular}\endgroup%
\kern3pt%
\begingroup \smaller\smaller\smaller\begin{tabular}{@{}c@{}}%
0\\8\\11
\end{tabular}\endgroup%
{$\left.\llap{\phantom{%
\begingroup \smaller\smaller\smaller\begin{tabular}{@{}c@{}}%
0\\0\\0
\end{tabular}\endgroup%
}}\!\right]$}%
\EasyButWeakLineBreak%
{$\left[\!\llap{\phantom{%
\begingroup \smaller\smaller\smaller\begin{tabular}{@{}c@{}}%
0\\0\\0
\end{tabular}\endgroup%
}}\right.$}%
\begingroup \smaller\smaller\smaller\begin{tabular}{@{}c@{}}%
-1\\-50\\-18
\end{tabular}\endgroup%
\HardButStrongLineBreak\kern3pt%
\begingroup \smaller\smaller\smaller\begin{tabular}{@{}c@{}}%
-1\\-56\\-26
\end{tabular}\endgroup%
\HardButStrongLineBreak\kern3pt%
\begingroup \smaller\smaller\smaller\begin{tabular}{@{}c@{}}%
-1\\-70\\-45
\end{tabular}\endgroup%
\HardButStrongLineBreak\kern3pt%
\begingroup \smaller\smaller\smaller\begin{tabular}{@{}c@{}}%
1\\24\\-18
\end{tabular}\endgroup%
{$\left.\llap{\phantom{%
\begingroup \smaller\smaller\smaller\begin{tabular}{@{}c@{}}%
0\\0\\0
\end{tabular}\endgroup%
}}\!\right]$}%
%
%
\hbox{}\par\smallskip%
%
%
\leavevmode%
${L_{214.2}}$%
{} : {$1\above{1pt}{1pt}{-2}{6}8\above{1pt}{1pt}{-}{5}{\cdot}1\above{1pt}{1pt}{-}{}3\above{1pt}{1pt}{-}{}9\above{1pt}{1pt}{-}{}{\cdot}1\above{1pt}{1pt}{-2}{}5\above{1pt}{1pt}{-}{}$}\spacer%
\instructions{2}%
\EasyButWeakLineBreak%
{${6}\above{1pt}{1pt}{b}{2}{8}\above{1pt}{1pt}{*}{2}{60}\above{1pt}{1pt}{*}{2}{72}\above{1pt}{1pt}{b}{2}$}\relax$\,(\times2)$%
\nopagebreak\par%
\nopagebreak\par\leavevmode%
{$\left[\!\llap{\phantom{%
\begingroup \smaller\smaller\smaller\begin{tabular}{@{}c@{}}%
0\\0\\0
\end{tabular}\endgroup%
}}\right.$}%
\begingroup \smaller\smaller\smaller\begin{tabular}{@{}c@{}}%
230760\\4320\\-1800
\end{tabular}\endgroup%
\kern3pt%
\begingroup \smaller\smaller\smaller\begin{tabular}{@{}c@{}}%
4320\\69\\-33
\end{tabular}\endgroup%
\kern3pt%
\begingroup \smaller\smaller\smaller\begin{tabular}{@{}c@{}}%
-1800\\-33\\14
\end{tabular}\endgroup%
{$\left.\llap{\phantom{%
\begingroup \smaller\smaller\smaller\begin{tabular}{@{}c@{}}%
0\\0\\0
\end{tabular}\endgroup%
}}\!\right]$}%
\hfil\penalty500%
{$\left[\!\llap{\phantom{%
\begingroup \smaller\smaller\smaller\begin{tabular}{@{}c@{}}%
0\\0\\0
\end{tabular}\endgroup%
}}\right.$}%
\begingroup \smaller\smaller\smaller\begin{tabular}{@{}c@{}}%
-1\\720\\2160
\end{tabular}\endgroup%
\kern3pt%
\begingroup \smaller\smaller\smaller\begin{tabular}{@{}c@{}}%
0\\19\\60
\end{tabular}\endgroup%
\kern3pt%
\begingroup \smaller\smaller\smaller\begin{tabular}{@{}c@{}}%
0\\-6\\-19
\end{tabular}\endgroup%
{$\left.\llap{\phantom{%
\begingroup \smaller\smaller\smaller\begin{tabular}{@{}c@{}}%
0\\0\\0
\end{tabular}\endgroup%
}}\!\right]$}%
\EasyButWeakLineBreak%
{$\left[\!\llap{\phantom{%
\begingroup \smaller\smaller\smaller\begin{tabular}{@{}c@{}}%
0\\0\\0
\end{tabular}\endgroup%
}}\right.$}%
\begingroup \smaller\smaller\smaller\begin{tabular}{@{}c@{}}%
-1\\-10\\-153
\end{tabular}\endgroup%
\HardButStrongLineBreak\kern3pt%
\begingroup \smaller\smaller\smaller\begin{tabular}{@{}c@{}}%
-1\\-16\\-172
\end{tabular}\endgroup%
\HardButStrongLineBreak\kern3pt%
\begingroup \smaller\smaller\smaller\begin{tabular}{@{}c@{}}%
1\\-10\\90
\end{tabular}\endgroup%
\HardButStrongLineBreak\kern3pt%
\begingroup \smaller\smaller\smaller\begin{tabular}{@{}c@{}}%
5\\24\\684
\end{tabular}\endgroup%
{$\left.\llap{\phantom{%
\begingroup \smaller\smaller\smaller\begin{tabular}{@{}c@{}}%
0\\0\\0
\end{tabular}\endgroup%
}}\!\right]$}%
%
%
\hbox{}\par\smallskip%
%
%
\leavevmode%
${L_{214.3}}$%
{} : {$1\above{1pt}{1pt}{2}{6}16\above{1pt}{1pt}{1}{1}{\cdot}1\above{1pt}{1pt}{1}{}3\above{1pt}{1pt}{1}{}9\above{1pt}{1pt}{1}{}{\cdot}1\above{1pt}{1pt}{-2}{}5\above{1pt}{1pt}{1}{}$}\spacer%
\instructions{2,m}%
\EasyButWeakLineBreak%
{${3}\above{1pt}{1pt}{}{2}{16}\above{1pt}{1pt}{r}{2}{30}\above{1pt}{1pt}{b}{2}{144}\above{1pt}{1pt}{*}{2}{12}\above{1pt}{1pt}{*}{2}{16}\above{1pt}{1pt}{b}{2}{30}\above{1pt}{1pt}{l}{2}{144}\above{1pt}{1pt}{}{2}$}%
\nopagebreak\par%
\nopagebreak\par\leavevmode%
{$\left[\!\llap{\phantom{%
\begingroup \smaller\smaller\smaller\begin{tabular}{@{}c@{}}%
0\\0\\0
\end{tabular}\endgroup%
}}\right.$}%
\begingroup \smaller\smaller\smaller\begin{tabular}{@{}c@{}}%
-929520\\-4320\\-5760
\end{tabular}\endgroup%
\kern3pt%
\begingroup \smaller\smaller\smaller\begin{tabular}{@{}c@{}}%
-4320\\3\\-6
\end{tabular}\endgroup%
\kern3pt%
\begingroup \smaller\smaller\smaller\begin{tabular}{@{}c@{}}%
-5760\\-6\\-17
\end{tabular}\endgroup%
{$\left.\llap{\phantom{%
\begingroup \smaller\smaller\smaller\begin{tabular}{@{}c@{}}%
0\\0\\0
\end{tabular}\endgroup%
}}\!\right]$}%
\EasyButWeakLineBreak%
{$\left[\!\llap{\phantom{%
\begingroup \smaller\smaller\smaller\begin{tabular}{@{}c@{}}%
0\\0\\0
\end{tabular}\endgroup%
}}\right.$}%
\begingroup \smaller\smaller\smaller\begin{tabular}{@{}c@{}}%
0\\-1\\0
\end{tabular}\endgroup%
\HardButStrongLineBreak\kern3pt%
\begingroup \smaller\smaller\smaller\begin{tabular}{@{}c@{}}%
-1\\-448\\496
\end{tabular}\endgroup%
\HardButStrongLineBreak\kern3pt%
\begingroup \smaller\smaller\smaller\begin{tabular}{@{}c@{}}%
-1\\-445\\495
\end{tabular}\endgroup%
\HardButStrongLineBreak\kern3pt%
\begingroup \smaller\smaller\smaller\begin{tabular}{@{}c@{}}%
1\\456\\-504
\end{tabular}\endgroup%
\HardButStrongLineBreak\kern3pt%
\begingroup \smaller\smaller\smaller\begin{tabular}{@{}c@{}}%
1\\448\\-498
\end{tabular}\endgroup%
\HardButStrongLineBreak\kern3pt%
\begingroup \smaller\smaller\smaller\begin{tabular}{@{}c@{}}%
5\\2232\\-2488
\end{tabular}\endgroup%
\HardButStrongLineBreak\kern3pt%
\begingroup \smaller\smaller\smaller\begin{tabular}{@{}c@{}}%
6\\2675\\-2985
\end{tabular}\endgroup%
\HardButStrongLineBreak\kern3pt%
\begingroup \smaller\smaller\smaller\begin{tabular}{@{}c@{}}%
11\\4896\\-5472
\end{tabular}\endgroup%
{$\left.\llap{\phantom{%
\begingroup \smaller\smaller\smaller\begin{tabular}{@{}c@{}}%
0\\0\\0
\end{tabular}\endgroup%
}}\!\right]$}%

\medskip%
%
\leavevmode\llap{}%
$W_{215}$%
\qquad\llap{34} lattices, $\chi=108$%
\hfill%
$22\slashtwo2222\slashinfty2222\slashtwo2222\slashinfty22\rtimes D_{4}$%
\nopagebreak\smallskip\hrule\nopagebreak\medskip%
%
%
\leavevmode%
${L_{215.1}}$%
{} : {$1\above{1pt}{1pt}{2}{0}8\above{1pt}{1pt}{1}{1}{\cdot}1\above{1pt}{1pt}{2}{}17\above{1pt}{1pt}{1}{}$}\EasyButWeakLineBreak%
{${68}\above{1pt}{1pt}{*}{2}{8}\above{1pt}{1pt}{*}{2}{4}\above{1pt}{1pt}{l}{2}{1}\above{1pt}{1pt}{}{2}{8}\above{1pt}{1pt}{}{2}{17}\above{1pt}{1pt}{r}{2}{4}\above{1pt}{1pt}{*}{2}{136}\above{1pt}{1pt}{1,0}{\infty b}{136}\above{1pt}{1pt}{}{2}{1}\above{1pt}{1pt}{r}{2}$}\relax$\,(\times2)$%
\nopagebreak\par%
\nopagebreak\par\leavevmode%
{$\left[\!\llap{\phantom{%
\begingroup \smaller\smaller\smaller
\endgroup%
}}\!\right]$}%
%
%
\hbox{}\par\smallskip%
%
%
\leavevmode%
${L_{215.2}}$%
{} : {$[1\above{1pt}{1pt}{1}{}2\above{1pt}{1pt}{1}{}]\above{1pt}{1pt}{}{2}16\above{1pt}{1pt}{1}{7}{\cdot}1\above{1pt}{1pt}{2}{}17\above{1pt}{1pt}{1}{}$}\spacer%
\instructions{2}%
\EasyButWeakLineBreak%
{${17}\above{1pt}{1pt}{}{2}{2}\above{1pt}{1pt}{}{2}{1}\above{1pt}{1pt}{r}{2}{16}\above{1pt}{1pt}{*}{2}{8}\above{1pt}{1pt}{*}{2}{272}\above{1pt}{1pt}{l}{2}{1}\above{1pt}{1pt}{}{2}{34}\above{1pt}{1pt}{8,7}{\infty}{136}\above{1pt}{1pt}{*}{2}{16}\above{1pt}{1pt}{l}{2}$}\relax$\,(\times2)$%
\nopagebreak\par%
\nopagebreak\par\leavevmode%
{$\left[\!\llap{\phantom{%
\begingroup \smaller\smaller\smaller
\endgroup%
}}\!\right]$}%
%
%
\hbox{}\par\smallskip%
%
%
\leavevmode%
${L_{215.3}}$%
{} : {$[1\above{1pt}{1pt}{-}{}2\above{1pt}{1pt}{1}{}]\above{1pt}{1pt}{}{6}16\above{1pt}{1pt}{-}{3}{\cdot}1\above{1pt}{1pt}{2}{}17\above{1pt}{1pt}{1}{}$}\spacer%
\instructions{m}%
\EasyButWeakLineBreak%
{${68}\above{1pt}{1pt}{l}{2}{2}\above{1pt}{1pt}{r}{2}{4}\above{1pt}{1pt}{*}{2}{16}\above{1pt}{1pt}{s}{2}{8}\above{1pt}{1pt}{s}{2}{272}\above{1pt}{1pt}{*}{2}{4}\above{1pt}{1pt}{l}{2}{34}\above{1pt}{1pt}{8,3}{\infty}{136}\above{1pt}{1pt}{s}{2}{16}\above{1pt}{1pt}{*}{2}$}\relax$\,(\times2)$%
\nopagebreak\par%
\nopagebreak\par\leavevmode%
{$\left[\!\llap{\phantom{%
\begingroup \smaller\smaller\smaller
\endgroup%
}}\!\right]$}%
%
%
\hbox{}\par\smallskip%
%
%
\leavevmode%
${L_{215.4}}$%
{} : {$[1\above{1pt}{1pt}{-}{}2\above{1pt}{1pt}{1}{}]\above{1pt}{1pt}{}{4}16\above{1pt}{1pt}{-}{5}{\cdot}1\above{1pt}{1pt}{2}{}17\above{1pt}{1pt}{1}{}$}\spacer%
\instructions{m}%
\EasyButWeakLineBreak%
{${68}\above{1pt}{1pt}{*}{2}{8}\above{1pt}{1pt}{*}{2}{4}\above{1pt}{1pt}{s}{2}{16}\above{1pt}{1pt}{l}{2}{2}\above{1pt}{1pt}{r}{2}{272}\above{1pt}{1pt}{s}{2}{4}\above{1pt}{1pt}{*}{2}{136}\above{1pt}{1pt}{8,3}{\infty z}{34}\above{1pt}{1pt}{r}{2}{16}\above{1pt}{1pt}{s}{2}$}\relax$\,(\times2)$%
\nopagebreak\par%
\nopagebreak\par\leavevmode%
{$\left[\!\llap{\phantom{%
\begingroup \smaller\smaller\smaller
\endgroup%
}}\!\right]$}%
%
%
\hbox{}\par\smallskip%
%
%
\leavevmode%
${L_{215.5}}$%
{} : {$1\above{1pt}{1pt}{1}{1}8\above{1pt}{1pt}{1}{7}64\above{1pt}{1pt}{1}{1}{\cdot}1\above{1pt}{1pt}{2}{}17\above{1pt}{1pt}{1}{}$}\spacer%
\instructions{2}%
\EasyButWeakLineBreak%
{${17}\above{1pt}{1pt}{r}{2}{32}\above{1pt}{1pt}{*}{2}{4}\above{1pt}{1pt}{s}{2}{64}\above{1pt}{1pt}{b}{2}{8}\above{1pt}{1pt}{l}{2}{1088}\above{1pt}{1pt}{}{2}{1}\above{1pt}{1pt}{r}{2}{544}\above{1pt}{1pt}{16,15}{\infty z}{136}\above{1pt}{1pt}{l}{2}{64}\above{1pt}{1pt}{}{2}$}\relax$\,(\times2)$%
\nopagebreak\par%
shares genus with {$ {L_{215.6}}$}%
; isometric to own %
2-dual\nopagebreak\par%
\nopagebreak\par\leavevmode%
{$\left[\!\llap{\phantom{%
\begingroup \smaller\smaller\smaller
\endgroup%
}}\!\right]$}%
%
%
\hbox{}\par\smallskip%
%
%
\leavevmode%
${L_{215.6}}$%
{} : {$1\above{1pt}{1pt}{1}{1}8\above{1pt}{1pt}{1}{7}64\above{1pt}{1pt}{1}{1}{\cdot}1\above{1pt}{1pt}{2}{}17\above{1pt}{1pt}{1}{}$}\EasyButWeakLineBreak%
{${68}\above{1pt}{1pt}{*}{2}{32}\above{1pt}{1pt}{l}{2}{1}\above{1pt}{1pt}{}{2}{64}\above{1pt}{1pt}{r}{2}{8}\above{1pt}{1pt}{b}{2}{1088}\above{1pt}{1pt}{s}{2}{4}\above{1pt}{1pt}{*}{2}{544}\above{1pt}{1pt}{16,7}{\infty z}{136}\above{1pt}{1pt}{b}{2}{64}\above{1pt}{1pt}{s}{2}$}\relax$\,(\times2)$%
\nopagebreak\par%
shares genus with {$ {L_{215.5}}$}%
; isometric to own %
2-dual\nopagebreak\par%
\nopagebreak\par\leavevmode%
{$\left[\!\llap{\phantom{%
\begingroup \smaller\smaller\smaller
\endgroup%
}}\!\right]$}%

\medskip%
%
\leavevmode\llap{}%
$W_{216}$%
\qquad\llap{4} lattices, $\chi=18$%
\hfill%
$24\slashinfty42|\rtimes D_{2}$%
\nopagebreak\smallskip\hrule\nopagebreak\medskip%
%
%
\leavevmode%
${L_{216.1}}$%
{} : {$1\above{1pt}{1pt}{2}{2}16\above{1pt}{1pt}{1}{7}{\cdot}1\above{1pt}{1pt}{2}{}9\above{1pt}{1pt}{-}{}$}\EasyButWeakLineBreak%
{${18}\above{1pt}{1pt}{b}{2}{2}\above{1pt}{1pt}{}{4}{1}\above{1pt}{1pt}{24,23}{\infty}{4}\above{1pt}{1pt}{*}{4}{2}\above{1pt}{1pt}{s}{2}$}%
\nopagebreak\par%
\nopagebreak\par\leavevmode%
{$\left[\!\llap{\phantom{%
\begingroup \smaller\smaller\smaller\begin{tabular}{@{}c@{}}%
0\\0\\0
\end{tabular}\endgroup%
}}\right.$}%
\begingroup \smaller\smaller\smaller\begin{tabular}{@{}c@{}}%
-262800\\2016\\2880
\end{tabular}\endgroup%
\kern3pt%
\begingroup \smaller\smaller\smaller\begin{tabular}{@{}c@{}}%
2016\\-14\\-23
\end{tabular}\endgroup%
\kern3pt%
\begingroup \smaller\smaller\smaller\begin{tabular}{@{}c@{}}%
2880\\-23\\-31
\end{tabular}\endgroup%
{$\left.\llap{\phantom{%
\begingroup \smaller\smaller\smaller\begin{tabular}{@{}c@{}}%
0\\0\\0
\end{tabular}\endgroup%
}}\!\right]$}%
\EasyButWeakLineBreak%
{$\left[\!\llap{\phantom{%
\begingroup \smaller\smaller\smaller\begin{tabular}{@{}c@{}}%
0\\0\\0
\end{tabular}\endgroup%
}}\right.$}%
\begingroup \smaller\smaller\smaller\begin{tabular}{@{}c@{}}%
-2\\-81\\-126
\end{tabular}\endgroup%
\HardButStrongLineBreak\kern3pt%
\begingroup \smaller\smaller\smaller\begin{tabular}{@{}c@{}}%
-1\\-39\\-64
\end{tabular}\endgroup%
\HardButStrongLineBreak\kern3pt%
\begingroup \smaller\smaller\smaller\begin{tabular}{@{}c@{}}%
1\\40\\63
\end{tabular}\endgroup%
\HardButStrongLineBreak\kern3pt%
\begingroup \smaller\smaller\smaller\begin{tabular}{@{}c@{}}%
5\\196\\318
\end{tabular}\endgroup%
\HardButStrongLineBreak\kern3pt%
\begingroup \smaller\smaller\smaller\begin{tabular}{@{}c@{}}%
2\\77\\128
\end{tabular}\endgroup%
{$\left.\llap{\phantom{%
\begingroup \smaller\smaller\smaller\begin{tabular}{@{}c@{}}%
0\\0\\0
\end{tabular}\endgroup%
}}\!\right]$}%
%
%
%
%
%
%
%
%
%
%
%
%
%
%

\medskip%
%
\leavevmode\llap{}%
$W_{217}$%
\qquad\llap{12} lattices, $\chi=42$%
\hfill%
$2224222242\rtimes C_{2}$%
\nopagebreak\smallskip\hrule\nopagebreak\medskip%
%
%
\leavevmode%
${L_{217.1}}$%
{} : {$1\above{1pt}{1pt}{-2}{{\rm II}}4\above{1pt}{1pt}{1}{7}{\cdot}1\above{1pt}{1pt}{2}{}3\above{1pt}{1pt}{-}{}{\cdot}1\above{1pt}{1pt}{2}{}49\above{1pt}{1pt}{1}{}$}\spacer%
\instructions{2}%
\EasyButWeakLineBreak%
{${196}\above{1pt}{1pt}{b}{2}{6}\above{1pt}{1pt}{s}{2}{98}\above{1pt}{1pt}{b}{2}{2}\above{1pt}{1pt}{*}{4}{4}\above{1pt}{1pt}{*}{2}$}\relax$\,(\times2)$%
\nopagebreak\par%
\nopagebreak\par\leavevmode%
{$\left[\!\llap{\phantom{%
\begingroup \smaller\smaller\smaller
\endgroup%
}}\!\right]$}%

\medskip%
%
\leavevmode\llap{}%
$W_{218}$%
\qquad\llap{32} lattices, $\chi=80$%
\hfill%
$62\infty22\infty62\infty22\infty\rtimes C_{2}$%
\nopagebreak\smallskip\hrule\nopagebreak\medskip%
%
%
\leavevmode%
${L_{218.1}}$%
{} : {$1\above{1pt}{1pt}{-2}{{\rm II}}8\above{1pt}{1pt}{1}{7}{\cdot}1\above{1pt}{1pt}{-}{}3\above{1pt}{1pt}{-}{}9\above{1pt}{1pt}{1}{}{\cdot}1\above{1pt}{1pt}{-2}{}25\above{1pt}{1pt}{-}{}$}\spacer%
\instructions{23,3,2}%
\EasyButWeakLineBreak%
{${6}\above{1pt}{1pt}{}{6}{2}\above{1pt}{1pt}{b}{2}{24}\above{1pt}{1pt}{30,29}{\infty z}{6}\above{1pt}{1pt}{s}{2}{50}\above{1pt}{1pt}{b}{2}{24}\above{1pt}{1pt}{30,11}{\infty z}$}\relax$\,(\times2)$%
\nopagebreak\par%
\nopagebreak\par\leavevmode%
{$\left[\!\llap{\phantom{%
\begingroup \smaller\smaller\smaller
\endgroup%
}}\!\right]$}%

\medskip%
%
\leavevmode\llap{}%
$W_{219}$%
\qquad\llap{32} lattices, $\chi=60$%
\hfill%
$22222222222222\rtimes C_{2}$%
\nopagebreak\smallskip\hrule\nopagebreak\medskip%
%
%
\leavevmode%
${L_{219.1}}$%
{} : {$1\above{1pt}{1pt}{2}{2}8\above{1pt}{1pt}{-}{5}{\cdot}1\above{1pt}{1pt}{2}{}3\above{1pt}{1pt}{1}{}{\cdot}1\above{1pt}{1pt}{-2}{}25\above{1pt}{1pt}{-}{}$}\spacer%
\instructions{2}%
\EasyButWeakLineBreak%
{${2}\above{1pt}{1pt}{b}{2}{200}\above{1pt}{1pt}{*}{2}{12}\above{1pt}{1pt}{l}{2}{1}\above{1pt}{1pt}{r}{2}{300}\above{1pt}{1pt}{*}{2}{8}\above{1pt}{1pt}{b}{2}{50}\above{1pt}{1pt}{s}{2}$}\relax$\,(\times2)$%
\nopagebreak\par%
\nopagebreak\par\leavevmode%
{$\left[\!\llap{\phantom{%
\begingroup \smaller\smaller\smaller
\endgroup%
}}\!\right]$}%
%
%
\hbox{}\par\smallskip%
%
%
\leavevmode%
${L_{219.2}}$%
{} : {$1\above{1pt}{1pt}{-2}{2}8\above{1pt}{1pt}{1}{1}{\cdot}1\above{1pt}{1pt}{2}{}3\above{1pt}{1pt}{1}{}{\cdot}1\above{1pt}{1pt}{-2}{}25\above{1pt}{1pt}{-}{}$}\spacer%
\instructions{m}%
\EasyButWeakLineBreak%
{${2}\above{1pt}{1pt}{l}{2}{200}\above{1pt}{1pt}{}{2}{3}\above{1pt}{1pt}{r}{2}{4}\above{1pt}{1pt}{l}{2}{75}\above{1pt}{1pt}{}{2}{8}\above{1pt}{1pt}{r}{2}{50}\above{1pt}{1pt}{b}{2}$}\relax$\,(\times2)$%
\nopagebreak\par%
\nopagebreak\par\leavevmode%
{$\left[\!\llap{\phantom{%
\begingroup \smaller\smaller\smaller
\endgroup%
}}\!\right]$}%

\medskip%
%
\leavevmode\llap{}%
$W_{220}$%
\qquad\llap{8} lattices, $\chi=48$%
\hfill%
$2\infty2222\infty222\rtimes C_{2}$%
\nopagebreak\smallskip\hrule\nopagebreak\medskip%
%
%
\leavevmode%
${L_{220.1}}$%
{} : {$1\above{1pt}{1pt}{-2}{{\rm II}}32\above{1pt}{1pt}{1}{1}{\cdot}1\above{1pt}{1pt}{-}{}5\above{1pt}{1pt}{-}{}25\above{1pt}{1pt}{-}{}$}\spacer%
\instructions{5}%
\EasyButWeakLineBreak%
{${32}\above{1pt}{1pt}{r}{2}{10}\above{1pt}{1pt}{40,9}{\infty b}{10}\above{1pt}{1pt}{b}{2}{32}\above{1pt}{1pt}{b}{2}{50}\above{1pt}{1pt}{l}{2}$}\relax$\,(\times2)$%
\nopagebreak\par%
shares genus with 5-dual\nopagebreak\par%
\nopagebreak\par\leavevmode%
{$\left[\!\llap{\phantom{%
\begingroup \smaller\smaller\smaller
\endgroup%
}}\!\right]$}%

\medskip%
%
\leavevmode\llap{}%
$W_{221}$%
\qquad\llap{8} lattices, $\chi=12$%
\hfill%
$\infty632$%
\nopagebreak\smallskip\hrule\nopagebreak\medskip%
%
%
\leavevmode%
${L_{221.1}}$%
{} : {$1\above{1pt}{1pt}{-2}{{\rm II}}8\above{1pt}{1pt}{-}{5}{\cdot}1\above{1pt}{1pt}{-}{}3\above{1pt}{1pt}{-}{}27\above{1pt}{1pt}{1}{}$}\spacer%
\instructions{2}%
\EasyButWeakLineBreak%
{${8}\above{1pt}{1pt}{6,5}{\infty z}{2}\above{1pt}{1pt}{}{6}{6}\above{1pt}{1pt}{-}{3}{6}\above{1pt}{1pt}{b}{2}$}%
\nopagebreak\par%
\nopagebreak\par\leavevmode%
{$\left[\!\llap{\phantom{%
\begingroup \smaller\smaller\smaller\begin{tabular}{@{}c@{}}%
0\\0\\0
\end{tabular}\endgroup%
}}\right.$}%
\begingroup \smaller\smaller\smaller\begin{tabular}{@{}c@{}}%
-684504\\16416\\-3024
\end{tabular}\endgroup%
\kern3pt%
\begingroup \smaller\smaller\smaller\begin{tabular}{@{}c@{}}%
16416\\-390\\69
\end{tabular}\endgroup%
\kern3pt%
\begingroup \smaller\smaller\smaller\begin{tabular}{@{}c@{}}%
-3024\\69\\-10
\end{tabular}\endgroup%
{$\left.\llap{\phantom{%
\begingroup \smaller\smaller\smaller\begin{tabular}{@{}c@{}}%
0\\0\\0
\end{tabular}\endgroup%
}}\!\right]$}%
\EasyButWeakLineBreak%
{$\left[\!\llap{\phantom{%
\begingroup \smaller\smaller\smaller\begin{tabular}{@{}c@{}}%
0\\0\\0
\end{tabular}\endgroup%
}}\right.$}%
\begingroup \smaller\smaller\smaller\begin{tabular}{@{}c@{}}%
1\\52\\56
\end{tabular}\endgroup%
\HardButStrongLineBreak\kern3pt%
\begingroup \smaller\smaller\smaller\begin{tabular}{@{}c@{}}%
-3\\-155\\-163
\end{tabular}\endgroup%
\HardButStrongLineBreak\kern3pt%
\begingroup \smaller\smaller\smaller\begin{tabular}{@{}c@{}}%
-1\\-52\\-57
\end{tabular}\endgroup%
\HardButStrongLineBreak\kern3pt%
\begingroup \smaller\smaller\smaller\begin{tabular}{@{}c@{}}%
3\\155\\162
\end{tabular}\endgroup%
{$\left.\llap{\phantom{%
\begingroup \smaller\smaller\smaller\begin{tabular}{@{}c@{}}%
0\\0\\0
\end{tabular}\endgroup%
}}\!\right]$}%

\medskip%
%
\leavevmode\llap{}%
$W_{222}$%
\qquad\llap{8} lattices, $\chi=24$%
\hfill%
$2\infty22\infty2\rtimes C_{2}$%
\nopagebreak\smallskip\hrule\nopagebreak\medskip%
%
%
\leavevmode%
${L_{222.1}}$%
{} : {$1\above{1pt}{1pt}{-2}{{\rm II}}8\above{1pt}{1pt}{-}{5}{\cdot}1\above{1pt}{1pt}{-}{}3\above{1pt}{1pt}{1}{}27\above{1pt}{1pt}{-}{}$}\spacer%
\instructions{2}%
\EasyButWeakLineBreak%
{${54}\above{1pt}{1pt}{b}{2}{8}\above{1pt}{1pt}{6,1}{\infty z}{2}\above{1pt}{1pt}{s}{2}$}\relax$\,(\times2)$%
\nopagebreak\par%
\nopagebreak\par\leavevmode%
{$\left[\!\llap{\phantom{%
\begingroup \smaller\smaller\smaller\begin{tabular}{@{}c@{}}%
0\\0\\0
\end{tabular}\endgroup%
}}\right.$}%
\begingroup \smaller\smaller\smaller\begin{tabular}{@{}c@{}}%
-1367064\\-460944\\12096
\end{tabular}\endgroup%
\kern3pt%
\begingroup \smaller\smaller\smaller\begin{tabular}{@{}c@{}}%
-460944\\-155418\\4077
\end{tabular}\endgroup%
\kern3pt%
\begingroup \smaller\smaller\smaller\begin{tabular}{@{}c@{}}%
12096\\4077\\-106
\end{tabular}\endgroup%
{$\left.\llap{\phantom{%
\begingroup \smaller\smaller\smaller\begin{tabular}{@{}c@{}}%
0\\0\\0
\end{tabular}\endgroup%
}}\!\right]$}%
\hfil\penalty500%
{$\left[\!\llap{\phantom{%
\begingroup \smaller\smaller\smaller\begin{tabular}{@{}c@{}}%
0\\0\\0
\end{tabular}\endgroup%
}}\right.$}%
\begingroup \smaller\smaller\smaller\begin{tabular}{@{}c@{}}%
-106201\\327600\\480600
\end{tabular}\endgroup%
\kern3pt%
\begingroup \smaller\smaller\smaller\begin{tabular}{@{}c@{}}%
-35813\\110473\\162069
\end{tabular}\endgroup%
\kern3pt%
\begingroup \smaller\smaller\smaller\begin{tabular}{@{}c@{}}%
944\\-2912\\-4273
\end{tabular}\endgroup%
{$\left.\llap{\phantom{%
\begingroup \smaller\smaller\smaller\begin{tabular}{@{}c@{}}%
0\\0\\0
\end{tabular}\endgroup%
}}\!\right]$}%
\EasyButWeakLineBreak%
{$\left[\!\llap{\phantom{%
\begingroup \smaller\smaller\smaller\begin{tabular}{@{}c@{}}%
0\\0\\0
\end{tabular}\endgroup%
}}\right.$}%
\begingroup \smaller\smaller\smaller\begin{tabular}{@{}c@{}}%
73\\-225\\-324
\end{tabular}\endgroup%
\HardButStrongLineBreak\kern3pt%
\begingroup \smaller\smaller\smaller\begin{tabular}{@{}c@{}}%
35\\-108\\-160
\end{tabular}\endgroup%
\HardButStrongLineBreak\kern3pt%
\begingroup \smaller\smaller\smaller\begin{tabular}{@{}c@{}}%
-36\\111\\161
\end{tabular}\endgroup%
{$\left.\llap{\phantom{%
\begingroup \smaller\smaller\smaller\begin{tabular}{@{}c@{}}%
0\\0\\0
\end{tabular}\endgroup%
}}\!\right]$}%

\medskip%
%
\leavevmode\llap{}%
$W_{223}$%
\qquad\llap{24} lattices, $\chi=36$%
\hfill%
$2222222222\rtimes C_{2}$%
\nopagebreak\smallskip\hrule\nopagebreak\medskip%
%
%
\leavevmode%
${L_{223.1}}$%
{} : {$1\above{1pt}{1pt}{2}{6}8\above{1pt}{1pt}{1}{7}{\cdot}1\above{1pt}{1pt}{1}{}3\above{1pt}{1pt}{-}{}27\above{1pt}{1pt}{-}{}$}\spacer%
\instructions{2}%
\EasyButWeakLineBreak%
{${216}\above{1pt}{1pt}{*}{2}{4}\above{1pt}{1pt}{*}{2}{24}\above{1pt}{1pt}{b}{2}{54}\above{1pt}{1pt}{s}{2}{6}\above{1pt}{1pt}{b}{2}$}\relax$\,(\times2)$%
\nopagebreak\par%
\nopagebreak\par\leavevmode%
{$\left[\!\llap{\phantom{%
\begingroup \smaller\smaller\smaller
\endgroup%
}}\!\right]$}%
%
%
\hbox{}\par\smallskip%
%
%
\leavevmode%
${L_{223.2}}$%
{} : {$1\above{1pt}{1pt}{-2}{2}16\above{1pt}{1pt}{-}{3}{\cdot}1\above{1pt}{1pt}{-}{}3\above{1pt}{1pt}{1}{}27\above{1pt}{1pt}{1}{}$}\spacer%
\instructions{2}%
\EasyButWeakLineBreak%
{${432}\above{1pt}{1pt}{r}{2}{2}\above{1pt}{1pt}{b}{2}{48}\above{1pt}{1pt}{*}{2}{108}\above{1pt}{1pt}{l}{2}{3}\above{1pt}{1pt}{}{2}$}\relax$\,(\times2)$%
\nopagebreak\par%
shares genus with {$ {L_{223.3}}$}%
\nopagebreak\par%
\nopagebreak\par\leavevmode%
{$\left[\!\llap{\phantom{%
\begingroup \smaller\smaller\smaller
\endgroup%
}}\!\right]$}%
%
%
\hbox{}\par\smallskip%
%
%
\leavevmode%
${L_{223.3}}$%
{} : {$1\above{1pt}{1pt}{-2}{2}16\above{1pt}{1pt}{-}{3}{\cdot}1\above{1pt}{1pt}{-}{}3\above{1pt}{1pt}{1}{}27\above{1pt}{1pt}{1}{}$}\spacer%
\instructions{m}%
\EasyButWeakLineBreak%
{${432}\above{1pt}{1pt}{b}{2}{2}\above{1pt}{1pt}{l}{2}{48}\above{1pt}{1pt}{}{2}{27}\above{1pt}{1pt}{r}{2}{12}\above{1pt}{1pt}{*}{2}$}\relax$\,(\times2)$%
\nopagebreak\par%
shares genus with {$ {L_{223.2}}$}%
\nopagebreak\par%
\nopagebreak\par\leavevmode%
{$\left[\!\llap{\phantom{%
\begingroup \smaller\smaller\smaller
\endgroup%
}}\!\right]$}%

\medskip%
%
\leavevmode\llap{}%
$W_{224}$%
\qquad\llap{60} lattices, $\chi=18$%
\hfill%
$\slashtwo222|222\rtimes D_{2}$%
\nopagebreak\smallskip\hrule\nopagebreak\medskip%
%
%
\leavevmode%
${L_{224.1}}$%
{} : {$1\above{1pt}{1pt}{2}{0}8\above{1pt}{1pt}{-}{5}{\cdot}1\above{1pt}{1pt}{2}{}3\above{1pt}{1pt}{1}{}{\cdot}1\above{1pt}{1pt}{2}{}7\above{1pt}{1pt}{1}{}$}\EasyButWeakLineBreak%
{${1}\above{1pt}{1pt}{r}{2}{4}\above{1pt}{1pt}{*}{2}{28}\above{1pt}{1pt}{l}{2}{3}\above{1pt}{1pt}{r}{2}{56}\above{1pt}{1pt}{s}{2}{12}\above{1pt}{1pt}{l}{2}{7}\above{1pt}{1pt}{}{2}$}%
\nopagebreak\par%
\nopagebreak\par\leavevmode%
{$\left[\!\llap{\phantom{%
\begingroup \smaller\smaller\smaller\begin{tabular}{@{}c@{}}%
0\\0\\0
\end{tabular}\endgroup%
}}\right.$}%
\begingroup \smaller\smaller\smaller\begin{tabular}{@{}c@{}}%
-13272\\168\\168
\end{tabular}\endgroup%
\kern3pt%
\begingroup \smaller\smaller\smaller\begin{tabular}{@{}c@{}}%
168\\-1\\-4
\end{tabular}\endgroup%
\kern3pt%
\begingroup \smaller\smaller\smaller\begin{tabular}{@{}c@{}}%
168\\-4\\1
\end{tabular}\endgroup%
{$\left.\llap{\phantom{%
\begingroup \smaller\smaller\smaller\begin{tabular}{@{}c@{}}%
0\\0\\0
\end{tabular}\endgroup%
}}\!\right]$}%
\EasyButWeakLineBreak%
{$\left[\!\llap{\phantom{%
\begingroup \smaller\smaller\smaller\begin{tabular}{@{}c@{}}%
0\\0\\0
\end{tabular}\endgroup%
}}\right.$}%
\begingroup \smaller\smaller\smaller\begin{tabular}{@{}c@{}}%
0\\0\\-1
\end{tabular}\endgroup%
\HardButStrongLineBreak\kern3pt%
\begingroup \smaller\smaller\smaller\begin{tabular}{@{}c@{}}%
-1\\-50\\-32
\end{tabular}\endgroup%
\HardButStrongLineBreak\kern3pt%
\begingroup \smaller\smaller\smaller\begin{tabular}{@{}c@{}}%
-5\\-252\\-154
\end{tabular}\endgroup%
\HardButStrongLineBreak\kern3pt%
\begingroup \smaller\smaller\smaller\begin{tabular}{@{}c@{}}%
-1\\-51\\-30
\end{tabular}\endgroup%
\HardButStrongLineBreak\kern3pt%
\begingroup \smaller\smaller\smaller\begin{tabular}{@{}c@{}}%
-1\\-56\\-28
\end{tabular}\endgroup%
\HardButStrongLineBreak\kern3pt%
\begingroup \smaller\smaller\smaller\begin{tabular}{@{}c@{}}%
1\\48\\30
\end{tabular}\endgroup%
\HardButStrongLineBreak\kern3pt%
\begingroup \smaller\smaller\smaller\begin{tabular}{@{}c@{}}%
1\\49\\28
\end{tabular}\endgroup%
{$\left.\llap{\phantom{%
\begingroup \smaller\smaller\smaller\begin{tabular}{@{}c@{}}%
0\\0\\0
\end{tabular}\endgroup%
}}\!\right]$}%
%
%
\hbox{}\par\smallskip%
%
%
\leavevmode%
${L_{224.2}}$%
{} : {$[1\above{1pt}{1pt}{1}{}2\above{1pt}{1pt}{-}{}]\above{1pt}{1pt}{}{4}16\above{1pt}{1pt}{1}{1}{\cdot}1\above{1pt}{1pt}{2}{}3\above{1pt}{1pt}{1}{}{\cdot}1\above{1pt}{1pt}{2}{}7\above{1pt}{1pt}{1}{}$}\spacer%
\instructions{2}%
\EasyButWeakLineBreak%
{${16}\above{1pt}{1pt}{}{2}{1}\above{1pt}{1pt}{r}{2}{112}\above{1pt}{1pt}{*}{2}{12}\above{1pt}{1pt}{*}{2}{56}\above{1pt}{1pt}{s}{2}{48}\above{1pt}{1pt}{l}{2}{7}\above{1pt}{1pt}{}{2}$}%
\nopagebreak\par%
\nopagebreak\par\leavevmode%
{$\left[\!\llap{\phantom{%
\begingroup \smaller\smaller\smaller\begin{tabular}{@{}c@{}}%
0\\0\\0
\end{tabular}\endgroup%
}}\right.$}%
\begingroup \smaller\smaller\smaller\begin{tabular}{@{}c@{}}%
4368\\-336\\336
\end{tabular}\endgroup%
\kern3pt%
\begingroup \smaller\smaller\smaller\begin{tabular}{@{}c@{}}%
-336\\10\\-6
\end{tabular}\endgroup%
\kern3pt%
\begingroup \smaller\smaller\smaller\begin{tabular}{@{}c@{}}%
336\\-6\\1
\end{tabular}\endgroup%
{$\left.\llap{\phantom{%
\begingroup \smaller\smaller\smaller\begin{tabular}{@{}c@{}}%
0\\0\\0
\end{tabular}\endgroup%
}}\!\right]$}%
\EasyButWeakLineBreak%
{$\left[\!\llap{\phantom{%
\begingroup \smaller\smaller\smaller\begin{tabular}{@{}c@{}}%
0\\0\\0
\end{tabular}\endgroup%
}}\right.$}%
\begingroup \smaller\smaller\smaller\begin{tabular}{@{}c@{}}%
1\\64\\48
\end{tabular}\endgroup%
\HardButStrongLineBreak\kern3pt%
\begingroup \smaller\smaller\smaller\begin{tabular}{@{}c@{}}%
0\\0\\-1
\end{tabular}\endgroup%
\HardButStrongLineBreak\kern3pt%
\begingroup \smaller\smaller\smaller\begin{tabular}{@{}c@{}}%
-3\\-196\\-168
\end{tabular}\endgroup%
\HardButStrongLineBreak\kern3pt%
\begingroup \smaller\smaller\smaller\begin{tabular}{@{}c@{}}%
-1\\-66\\-54
\end{tabular}\endgroup%
\HardButStrongLineBreak\kern3pt%
\begingroup \smaller\smaller\smaller\begin{tabular}{@{}c@{}}%
-1\\-70\\-56
\end{tabular}\endgroup%
\HardButStrongLineBreak\kern3pt%
\begingroup \smaller\smaller\smaller\begin{tabular}{@{}c@{}}%
1\\60\\48
\end{tabular}\endgroup%
\HardButStrongLineBreak\kern3pt%
\begingroup \smaller\smaller\smaller\begin{tabular}{@{}c@{}}%
1\\63\\49
\end{tabular}\endgroup%
{$\left.\llap{\phantom{%
\begingroup \smaller\smaller\smaller\begin{tabular}{@{}c@{}}%
0\\0\\0
\end{tabular}\endgroup%
}}\!\right]$}%
%
%
\hbox{}\par\smallskip%
%
%
\leavevmode%
${L_{224.3}}$%
{} : {$[1\above{1pt}{1pt}{-}{}2\above{1pt}{1pt}{-}{}]\above{1pt}{1pt}{}{0}16\above{1pt}{1pt}{-}{5}{\cdot}1\above{1pt}{1pt}{2}{}3\above{1pt}{1pt}{1}{}{\cdot}1\above{1pt}{1pt}{2}{}7\above{1pt}{1pt}{1}{}$}\spacer%
\instructions{m}%
\EasyButWeakLineBreak%
{${16}\above{1pt}{1pt}{s}{2}{4}\above{1pt}{1pt}{*}{2}{112}\above{1pt}{1pt}{l}{2}{3}\above{1pt}{1pt}{r}{2}{56}\above{1pt}{1pt}{*}{2}{48}\above{1pt}{1pt}{*}{2}{28}\above{1pt}{1pt}{s}{2}$}%
\nopagebreak\par%
\nopagebreak\par\leavevmode%
{$\left[\!\llap{\phantom{%
\begingroup \smaller\smaller\smaller\begin{tabular}{@{}c@{}}%
0\\0\\0
\end{tabular}\endgroup%
}}\right.$}%
\begingroup \smaller\smaller\smaller\begin{tabular}{@{}c@{}}%
306768\\-1680\\-1680
\end{tabular}\endgroup%
\kern3pt%
\begingroup \smaller\smaller\smaller\begin{tabular}{@{}c@{}}%
-1680\\10\\8
\end{tabular}\endgroup%
\kern3pt%
\begingroup \smaller\smaller\smaller\begin{tabular}{@{}c@{}}%
-1680\\8\\11
\end{tabular}\endgroup%
{$\left.\llap{\phantom{%
\begingroup \smaller\smaller\smaller\begin{tabular}{@{}c@{}}%
0\\0\\0
\end{tabular}\endgroup%
}}\!\right]$}%
\EasyButWeakLineBreak%
{$\left[\!\llap{\phantom{%
\begingroup \smaller\smaller\smaller\begin{tabular}{@{}c@{}}%
0\\0\\0
\end{tabular}\endgroup%
}}\right.$}%
\begingroup \smaller\smaller\smaller\begin{tabular}{@{}c@{}}%
-1\\-112\\-72
\end{tabular}\endgroup%
\HardButStrongLineBreak\kern3pt%
\begingroup \smaller\smaller\smaller\begin{tabular}{@{}c@{}}%
-1\\-110\\-74
\end{tabular}\endgroup%
\HardButStrongLineBreak\kern3pt%
\begingroup \smaller\smaller\smaller\begin{tabular}{@{}c@{}}%
-9\\-980\\-672
\end{tabular}\endgroup%
\HardButStrongLineBreak\kern3pt%
\begingroup \smaller\smaller\smaller\begin{tabular}{@{}c@{}}%
-1\\-108\\-75
\end{tabular}\endgroup%
\HardButStrongLineBreak\kern3pt%
\begingroup \smaller\smaller\smaller\begin{tabular}{@{}c@{}}%
-3\\-322\\-224
\end{tabular}\endgroup%
\HardButStrongLineBreak\kern3pt%
\begingroup \smaller\smaller\smaller\begin{tabular}{@{}c@{}}%
-1\\-108\\-72
\end{tabular}\endgroup%
\HardButStrongLineBreak\kern3pt%
\begingroup \smaller\smaller\smaller\begin{tabular}{@{}c@{}}%
-1\\-112\\-70
\end{tabular}\endgroup%
{$\left.\llap{\phantom{%
\begingroup \smaller\smaller\smaller\begin{tabular}{@{}c@{}}%
0\\0\\0
\end{tabular}\endgroup%
}}\!\right]$}%
%
%
\hbox{}\par\smallskip%
%
%
\leavevmode%
${L_{224.4}}$%
{} : {$[1\above{1pt}{1pt}{-}{}2\above{1pt}{1pt}{1}{}]\above{1pt}{1pt}{}{2}16\above{1pt}{1pt}{1}{7}{\cdot}1\above{1pt}{1pt}{2}{}3\above{1pt}{1pt}{1}{}{\cdot}1\above{1pt}{1pt}{2}{}7\above{1pt}{1pt}{1}{}$}\spacer%
\instructions{m}%
\EasyButWeakLineBreak%
{${16}\above{1pt}{1pt}{l}{2}{1}\above{1pt}{1pt}{}{2}{112}\above{1pt}{1pt}{}{2}{3}\above{1pt}{1pt}{}{2}{14}\above{1pt}{1pt}{r}{2}{48}\above{1pt}{1pt}{s}{2}{28}\above{1pt}{1pt}{*}{2}$}%
\nopagebreak\par%
\nopagebreak\par\leavevmode%
{$\left[\!\llap{\phantom{%
\begingroup \smaller\smaller\smaller\begin{tabular}{@{}c@{}}%
0\\0\\0
\end{tabular}\endgroup%
}}\right.$}%
\begingroup \smaller\smaller\smaller\begin{tabular}{@{}c@{}}%
38640\\336\\-672
\end{tabular}\endgroup%
\kern3pt%
\begingroup \smaller\smaller\smaller\begin{tabular}{@{}c@{}}%
336\\-2\\-4
\end{tabular}\endgroup%
\kern3pt%
\begingroup \smaller\smaller\smaller\begin{tabular}{@{}c@{}}%
-672\\-4\\11
\end{tabular}\endgroup%
{$\left.\llap{\phantom{%
\begingroup \smaller\smaller\smaller\begin{tabular}{@{}c@{}}%
0\\0\\0
\end{tabular}\endgroup%
}}\!\right]$}%
\EasyButWeakLineBreak%
{$\left[\!\llap{\phantom{%
\begingroup \smaller\smaller\smaller\begin{tabular}{@{}c@{}}%
0\\0\\0
\end{tabular}\endgroup%
}}\right.$}%
\begingroup \smaller\smaller\smaller\begin{tabular}{@{}c@{}}%
-1\\-28\\-72
\end{tabular}\endgroup%
\HardButStrongLineBreak\kern3pt%
\begingroup \smaller\smaller\smaller\begin{tabular}{@{}c@{}}%
0\\-1\\-1
\end{tabular}\endgroup%
\HardButStrongLineBreak\kern3pt%
\begingroup \smaller\smaller\smaller\begin{tabular}{@{}c@{}}%
5\\112\\336
\end{tabular}\endgroup%
\HardButStrongLineBreak\kern3pt%
\begingroup \smaller\smaller\smaller\begin{tabular}{@{}c@{}}%
1\\24\\69
\end{tabular}\endgroup%
\HardButStrongLineBreak\kern3pt%
\begingroup \smaller\smaller\smaller\begin{tabular}{@{}c@{}}%
2\\49\\140
\end{tabular}\endgroup%
\HardButStrongLineBreak\kern3pt%
\begingroup \smaller\smaller\smaller\begin{tabular}{@{}c@{}}%
1\\24\\72
\end{tabular}\endgroup%
\HardButStrongLineBreak\kern3pt%
\begingroup \smaller\smaller\smaller\begin{tabular}{@{}c@{}}%
-1\\-28\\-70
\end{tabular}\endgroup%
{$\left.\llap{\phantom{%
\begingroup \smaller\smaller\smaller\begin{tabular}{@{}c@{}}%
0\\0\\0
\end{tabular}\endgroup%
}}\!\right]$}%
%
%
\hbox{}\par\smallskip%
%
%
\leavevmode%
${L_{224.5}}$%
{} : {$[1\above{1pt}{1pt}{1}{}2\above{1pt}{1pt}{1}{}]\above{1pt}{1pt}{}{6}16\above{1pt}{1pt}{-}{3}{\cdot}1\above{1pt}{1pt}{2}{}3\above{1pt}{1pt}{1}{}{\cdot}1\above{1pt}{1pt}{2}{}7\above{1pt}{1pt}{1}{}$}\EasyButWeakLineBreak%
{${16}\above{1pt}{1pt}{*}{2}{4}\above{1pt}{1pt}{s}{2}{112}\above{1pt}{1pt}{s}{2}{12}\above{1pt}{1pt}{l}{2}{14}\above{1pt}{1pt}{}{2}{48}\above{1pt}{1pt}{}{2}{7}\above{1pt}{1pt}{r}{2}$}%
\nopagebreak\par%
\nopagebreak\par\leavevmode%
{$\left[\!\llap{\phantom{%
\begingroup \smaller\smaller\smaller\begin{tabular}{@{}c@{}}%
0\\0\\0
\end{tabular}\endgroup%
}}\right.$}%
\begingroup \smaller\smaller\smaller\begin{tabular}{@{}c@{}}%
-134736\\672\\1008
\end{tabular}\endgroup%
\kern3pt%
\begingroup \smaller\smaller\smaller\begin{tabular}{@{}c@{}}%
672\\-2\\-8
\end{tabular}\endgroup%
\kern3pt%
\begingroup \smaller\smaller\smaller\begin{tabular}{@{}c@{}}%
1008\\-8\\-1
\end{tabular}\endgroup%
{$\left.\llap{\phantom{%
\begingroup \smaller\smaller\smaller\begin{tabular}{@{}c@{}}%
0\\0\\0
\end{tabular}\endgroup%
}}\!\right]$}%
\EasyButWeakLineBreak%
{$\left[\!\llap{\phantom{%
\begingroup \smaller\smaller\smaller\begin{tabular}{@{}c@{}}%
0\\0\\0
\end{tabular}\endgroup%
}}\right.$}%
\begingroup \smaller\smaller\smaller\begin{tabular}{@{}c@{}}%
3\\356\\160
\end{tabular}\endgroup%
\HardButStrongLineBreak\kern3pt%
\begingroup \smaller\smaller\smaller\begin{tabular}{@{}c@{}}%
1\\118\\54
\end{tabular}\endgroup%
\HardButStrongLineBreak\kern3pt%
\begingroup \smaller\smaller\smaller\begin{tabular}{@{}c@{}}%
1\\112\\56
\end{tabular}\endgroup%
\HardButStrongLineBreak\kern3pt%
\begingroup \smaller\smaller\smaller\begin{tabular}{@{}c@{}}%
-1\\-120\\-54
\end{tabular}\endgroup%
\HardButStrongLineBreak\kern3pt%
\begingroup \smaller\smaller\smaller\begin{tabular}{@{}c@{}}%
-1\\-119\\-56
\end{tabular}\endgroup%
\HardButStrongLineBreak\kern3pt%
\begingroup \smaller\smaller\smaller\begin{tabular}{@{}c@{}}%
1\\120\\48
\end{tabular}\endgroup%
\HardButStrongLineBreak\kern3pt%
\begingroup \smaller\smaller\smaller\begin{tabular}{@{}c@{}}%
2\\238\\105
\end{tabular}\endgroup%
{$\left.\llap{\phantom{%
\begingroup \smaller\smaller\smaller\begin{tabular}{@{}c@{}}%
0\\0\\0
\end{tabular}\endgroup%
}}\!\right]$}%

\medskip%
%
\leavevmode\llap{}%
$W_{225}$%
\qquad\llap{120} lattices, $\chi=72$%
\hfill%
$2222|2222|2222|2222|\rtimes D_{4}$%
\nopagebreak\smallskip\hrule\nopagebreak\medskip%
%
%
\leavevmode%
${L_{225.1}}$%
{} : {$1\above{1pt}{1pt}{-2}{4}8\above{1pt}{1pt}{1}{1}{\cdot}1\above{1pt}{1pt}{1}{}3\above{1pt}{1pt}{-}{}9\above{1pt}{1pt}{-}{}{\cdot}1\above{1pt}{1pt}{2}{}7\above{1pt}{1pt}{1}{}$}\spacer%
\instructions{3}%
\EasyButWeakLineBreak%
{${24}\above{1pt}{1pt}{l}{2}{7}\above{1pt}{1pt}{}{2}{72}\above{1pt}{1pt}{}{2}{1}\above{1pt}{1pt}{r}{2}{504}\above{1pt}{1pt}{s}{2}{4}\above{1pt}{1pt}{*}{2}{72}\above{1pt}{1pt}{*}{2}{28}\above{1pt}{1pt}{s}{2}$}\relax$\,(\times2)$%
\nopagebreak\par%
\nopagebreak\par\leavevmode%
{$\left[\!\llap{\phantom{%
\begingroup \smaller\smaller\smaller
\endgroup%
}}\!\right]$}%
%
%
\hbox{}\par\smallskip%
%
%
\leavevmode%
${L_{225.2}}$%
{} : {$[1\above{1pt}{1pt}{1}{}2\above{1pt}{1pt}{1}{}]\above{1pt}{1pt}{}{2}16\above{1pt}{1pt}{-}{3}{\cdot}1\above{1pt}{1pt}{1}{}3\above{1pt}{1pt}{-}{}9\above{1pt}{1pt}{-}{}{\cdot}1\above{1pt}{1pt}{2}{}7\above{1pt}{1pt}{1}{}$}\spacer%
\instructions{3,2}%
\EasyButWeakLineBreak%
{${6}\above{1pt}{1pt}{r}{2}{112}\above{1pt}{1pt}{l}{2}{18}\above{1pt}{1pt}{}{2}{1}\above{1pt}{1pt}{r}{2}{504}\above{1pt}{1pt}{*}{2}{16}\above{1pt}{1pt}{s}{2}{72}\above{1pt}{1pt}{*}{2}{28}\above{1pt}{1pt}{l}{2}$}\relax$\,(\times2)$%
\nopagebreak\par%
\nopagebreak\par\leavevmode%
{$\left[\!\llap{\phantom{%
\begingroup \smaller\smaller\smaller
\endgroup%
}}\!\right]$}%
%
%
\hbox{}\par\smallskip%
%
%
\leavevmode%
${L_{225.3}}$%
{} : {$[1\above{1pt}{1pt}{1}{}2\above{1pt}{1pt}{1}{}]\above{1pt}{1pt}{}{0}16\above{1pt}{1pt}{-}{5}{\cdot}1\above{1pt}{1pt}{1}{}3\above{1pt}{1pt}{-}{}9\above{1pt}{1pt}{-}{}{\cdot}1\above{1pt}{1pt}{2}{}7\above{1pt}{1pt}{1}{}$}\spacer%
\instructions{32,3,m}%
\EasyButWeakLineBreak%
{${24}\above{1pt}{1pt}{s}{2}{112}\above{1pt}{1pt}{*}{2}{72}\above{1pt}{1pt}{l}{2}{1}\above{1pt}{1pt}{}{2}{126}\above{1pt}{1pt}{r}{2}{16}\above{1pt}{1pt}{l}{2}{18}\above{1pt}{1pt}{}{2}{7}\above{1pt}{1pt}{r}{2}$}\relax$\,(\times2)$%
\nopagebreak\par%
\nopagebreak\par\leavevmode%
{$\left[\!\llap{\phantom{%
\begingroup \smaller\smaller\smaller
\endgroup%
}}\!\right]$}%
%
%
\hbox{}\par\smallskip%
%
%
\leavevmode%
${L_{225.4}}$%
{} : {$[1\above{1pt}{1pt}{-}{}2\above{1pt}{1pt}{1}{}]\above{1pt}{1pt}{}{6}16\above{1pt}{1pt}{1}{7}{\cdot}1\above{1pt}{1pt}{-}{}3\above{1pt}{1pt}{-}{}9\above{1pt}{1pt}{1}{}{\cdot}1\above{1pt}{1pt}{2}{}7\above{1pt}{1pt}{1}{}$}\spacer%
\instructions{3m,3,m}%
\EasyButWeakLineBreak%
{${6}\above{1pt}{1pt}{}{2}{63}\above{1pt}{1pt}{r}{2}{8}\above{1pt}{1pt}{*}{2}{144}\above{1pt}{1pt}{s}{2}{56}\above{1pt}{1pt}{*}{2}{36}\above{1pt}{1pt}{l}{2}{2}\above{1pt}{1pt}{}{2}{1008}\above{1pt}{1pt}{}{2}$}\relax$\,(\times2)$%
\nopagebreak\par%
\nopagebreak\par\leavevmode%
{$\left[\!\llap{\phantom{%
\begingroup \smaller\smaller\smaller
\endgroup%
}}\!\right]$}%
%
%
\hbox{}\par\smallskip%
%
%
\leavevmode%
${L_{225.5}}$%
{} : {$[1\above{1pt}{1pt}{-}{}2\above{1pt}{1pt}{1}{}]\above{1pt}{1pt}{}{4}16\above{1pt}{1pt}{1}{1}{\cdot}1\above{1pt}{1pt}{1}{}3\above{1pt}{1pt}{-}{}9\above{1pt}{1pt}{-}{}{\cdot}1\above{1pt}{1pt}{2}{}7\above{1pt}{1pt}{1}{}$}\spacer%
\instructions{3m,3}%
\EasyButWeakLineBreak%
{${24}\above{1pt}{1pt}{*}{2}{112}\above{1pt}{1pt}{s}{2}{72}\above{1pt}{1pt}{*}{2}{4}\above{1pt}{1pt}{l}{2}{126}\above{1pt}{1pt}{}{2}{16}\above{1pt}{1pt}{}{2}{18}\above{1pt}{1pt}{r}{2}{28}\above{1pt}{1pt}{*}{2}$}\relax$\,(\times2)$%
\nopagebreak\par%
\nopagebreak\par\leavevmode%
{$\left[\!\llap{\phantom{%
\begingroup \smaller\smaller\smaller
\endgroup%
}}\!\right]$}%

\medskip%
%
\leavevmode\llap{}%
$W_{226}$%
\qquad\llap{12} lattices, $\chi=54$%
\hfill%
$422222422222\rtimes C_{2}$%
\nopagebreak\smallskip\hrule\nopagebreak\medskip%
%
%
\leavevmode%
${L_{226.1}}$%
{} : {$1\above{1pt}{1pt}{-2}{{\rm II}}4\above{1pt}{1pt}{1}{7}{\cdot}1\above{1pt}{1pt}{2}{}9\above{1pt}{1pt}{-}{}{\cdot}1\above{1pt}{1pt}{2}{}19\above{1pt}{1pt}{-}{}$}\spacer%
\instructions{2}%
\EasyButWeakLineBreak%
{${4}\above{1pt}{1pt}{*}{4}{2}\above{1pt}{1pt}{s}{2}{38}\above{1pt}{1pt}{s}{2}{18}\above{1pt}{1pt}{b}{2}{2}\above{1pt}{1pt}{s}{2}{342}\above{1pt}{1pt}{b}{2}$}\relax$\,(\times2)$%
\nopagebreak\par%
\nopagebreak\par\leavevmode%
{$\left[\!\llap{\phantom{%
\begingroup \smaller\smaller\smaller
\endgroup%
}}\!\right]$}%

\medskip%
%
\leavevmode\llap{}%
$W_{227}$%
\qquad\llap{6} lattices, $\chi=8$%
\hfill%
$\slashthree22|22\rtimes D_{2}$%
\nopagebreak\smallskip\hrule\nopagebreak\medskip%
%
%
\leavevmode%
${L_{227.1}}$%
{} : {$1\above{1pt}{1pt}{-2}{{\rm II}}4\above{1pt}{1pt}{1}{1}{\cdot}1\above{1pt}{1pt}{1}{}3\above{1pt}{1pt}{-}{}9\above{1pt}{1pt}{1}{}{\cdot}1\above{1pt}{1pt}{-2}{}7\above{1pt}{1pt}{-}{}$}\spacer%
\instructions{2}%
\EasyButWeakLineBreak%
{${6}\above{1pt}{1pt}{+}{3}{6}\above{1pt}{1pt}{l}{2}{36}\above{1pt}{1pt}{r}{2}{42}\above{1pt}{1pt}{l}{2}{4}\above{1pt}{1pt}{r}{2}$}%
\nopagebreak\par%
\nopagebreak\par\leavevmode%
{$\left[\!\llap{\phantom{%
\begingroup \smaller\smaller\smaller\begin{tabular}{@{}c@{}}%
0\\0\\0
\end{tabular}\endgroup%
}}\right.$}%
\begingroup \smaller\smaller\smaller\begin{tabular}{@{}c@{}}%
-203868\\-41580\\-24948
\end{tabular}\endgroup%
\kern3pt%
\begingroup \smaller\smaller\smaller\begin{tabular}{@{}c@{}}%
-41580\\-8454\\-5115
\end{tabular}\endgroup%
\kern3pt%
\begingroup \smaller\smaller\smaller\begin{tabular}{@{}c@{}}%
-24948\\-5115\\-3026
\end{tabular}\endgroup%
{$\left.\llap{\phantom{%
\begingroup \smaller\smaller\smaller\begin{tabular}{@{}c@{}}%
0\\0\\0
\end{tabular}\endgroup%
}}\!\right]$}%
\EasyButWeakLineBreak%
{$\left[\!\llap{\phantom{%
\begingroup \smaller\smaller\smaller\begin{tabular}{@{}c@{}}%
0\\0\\0
\end{tabular}\endgroup%
}}\right.$}%
\begingroup \smaller\smaller\smaller\begin{tabular}{@{}c@{}}%
-56\\172\\171
\end{tabular}\endgroup%
\HardButStrongLineBreak\kern3pt%
\begingroup \smaller\smaller\smaller\begin{tabular}{@{}c@{}}%
66\\-203\\-201
\end{tabular}\endgroup%
\HardButStrongLineBreak\kern3pt%
\begingroup \smaller\smaller\smaller\begin{tabular}{@{}c@{}}%
47\\-144\\-144
\end{tabular}\endgroup%
\HardButStrongLineBreak\kern3pt%
\begingroup \smaller\smaller\smaller\begin{tabular}{@{}c@{}}%
-207\\637\\630
\end{tabular}\endgroup%
\HardButStrongLineBreak\kern3pt%
\begingroup \smaller\smaller\smaller\begin{tabular}{@{}c@{}}%
-147\\452\\448
\end{tabular}\endgroup%
{$\left.\llap{\phantom{%
\begingroup \smaller\smaller\smaller\begin{tabular}{@{}c@{}}%
0\\0\\0
\end{tabular}\endgroup%
}}\!\right]$}%

\medskip%
%
\leavevmode\llap{}%
$W_{228}$%
\qquad\llap{6} lattices, $\chi=8$%
\hfill%
$6|62|2\rtimes D_{2}$%
\nopagebreak\smallskip\hrule\nopagebreak\medskip%
%
%
\leavevmode%
${L_{228.1}}$%
{} : {$1\above{1pt}{1pt}{-2}{{\rm II}}4\above{1pt}{1pt}{1}{1}{\cdot}1\above{1pt}{1pt}{-}{}3\above{1pt}{1pt}{-}{}9\above{1pt}{1pt}{-}{}{\cdot}1\above{1pt}{1pt}{-2}{}7\above{1pt}{1pt}{-}{}$}\spacer%
\instructions{2}%
\EasyButWeakLineBreak%
{${2}\above{1pt}{1pt}{}{6}{6}\above{1pt}{1pt}{}{6}{18}\above{1pt}{1pt}{b}{2}{42}\above{1pt}{1pt}{b}{2}$}%
\nopagebreak\par%
\nopagebreak\par\leavevmode%
{$\left[\!\llap{\phantom{%
\begingroup \smaller\smaller\smaller\begin{tabular}{@{}c@{}}%
0\\0\\0
\end{tabular}\endgroup%
}}\right.$}%
\begingroup \smaller\smaller\smaller\begin{tabular}{@{}c@{}}%
-115164\\-4536\\-33264
\end{tabular}\endgroup%
\kern3pt%
\begingroup \smaller\smaller\smaller\begin{tabular}{@{}c@{}}%
-4536\\-174\\-1317
\end{tabular}\endgroup%
\kern3pt%
\begingroup \smaller\smaller\smaller\begin{tabular}{@{}c@{}}%
-33264\\-1317\\-9598
\end{tabular}\endgroup%
{$\left.\llap{\phantom{%
\begingroup \smaller\smaller\smaller\begin{tabular}{@{}c@{}}%
0\\0\\0
\end{tabular}\endgroup%
}}\!\right]$}%
\EasyButWeakLineBreak%
{$\left[\!\llap{\phantom{%
\begingroup \smaller\smaller\smaller\begin{tabular}{@{}c@{}}%
0\\0\\0
\end{tabular}\endgroup%
}}\right.$}%
\begingroup \smaller\smaller\smaller\begin{tabular}{@{}c@{}}%
10\\-41\\-29
\end{tabular}\endgroup%
\HardButStrongLineBreak\kern3pt%
\begingroup \smaller\smaller\smaller\begin{tabular}{@{}c@{}}%
26\\-110\\-75
\end{tabular}\endgroup%
\HardButStrongLineBreak\kern3pt%
\begingroup \smaller\smaller\smaller\begin{tabular}{@{}c@{}}%
-28\\117\\81
\end{tabular}\endgroup%
\HardButStrongLineBreak\kern3pt%
\begingroup \smaller\smaller\smaller\begin{tabular}{@{}c@{}}%
-51\\217\\147
\end{tabular}\endgroup%
{$\left.\llap{\phantom{%
\begingroup \smaller\smaller\smaller\begin{tabular}{@{}c@{}}%
0\\0\\0
\end{tabular}\endgroup%
}}\!\right]$}%

\medskip%
%
\leavevmode\llap{}%
$W_{229}$%
\qquad\llap{26} lattices, $\chi=36$%
\hfill%
$\slashtwo22|22\slashtwo22|22\rtimes D_{4}$%
\nopagebreak\smallskip\hrule\nopagebreak\medskip%
%
%
\leavevmode%
${L_{229.1}}$%
{} : {$1\above{1pt}{1pt}{2}{{\rm II}}4\above{1pt}{1pt}{-}{5}{\cdot}1\above{1pt}{1pt}{1}{}3\above{1pt}{1pt}{-}{}9\above{1pt}{1pt}{1}{}{\cdot}1\above{1pt}{1pt}{2}{}7\above{1pt}{1pt}{1}{}$}\spacer%
\instructions{2}%
\EasyButWeakLineBreak%
{${36}\above{1pt}{1pt}{*}{2}{4}\above{1pt}{1pt}{*}{2}{252}\above{1pt}{1pt}{b}{2}{6}\above{1pt}{1pt}{b}{2}{28}\above{1pt}{1pt}{*}{2}$}\relax$\,(\times2)$%
\nopagebreak\par%
\nopagebreak\par\leavevmode%
{$\left[\!\llap{\phantom{%
\begingroup \smaller\smaller\smaller
\endgroup%
}}\!\right]$}%
%
%
\hbox{}\par\smallskip%
%
%
\leavevmode%
${L_{229.2}}$%
{} : {$1\above{1pt}{1pt}{2}{2}8\above{1pt}{1pt}{-}{3}{\cdot}1\above{1pt}{1pt}{-}{}3\above{1pt}{1pt}{1}{}9\above{1pt}{1pt}{-}{}{\cdot}1\above{1pt}{1pt}{2}{}7\above{1pt}{1pt}{1}{}$}\spacer%
\instructions{2}%
\EasyButWeakLineBreak%
{${2}\above{1pt}{1pt}{s}{2}{18}\above{1pt}{1pt}{b}{2}{56}\above{1pt}{1pt}{*}{2}{12}\above{1pt}{1pt}{*}{2}{504}\above{1pt}{1pt}{b}{2}$}\relax$\,(\times2)$%
\nopagebreak\par%
\nopagebreak\par\leavevmode%
{$\left[\!\llap{\phantom{%
\begingroup \smaller\smaller\smaller
\endgroup%
}}\!\right]$}%
%
%
\hbox{}\par\smallskip%
%
%
\leavevmode%
${L_{229.3}}$%
{} : {$1\above{1pt}{1pt}{-2}{6}16\above{1pt}{1pt}{1}{7}{\cdot}1\above{1pt}{1pt}{1}{}3\above{1pt}{1pt}{-}{}9\above{1pt}{1pt}{1}{}{\cdot}1\above{1pt}{1pt}{2}{}7\above{1pt}{1pt}{1}{}$}\spacer%
\instructions{2,m}%
\EasyButWeakLineBreak%
{${36}\above{1pt}{1pt}{l}{2}{1}\above{1pt}{1pt}{}{2}{1008}\above{1pt}{1pt}{r}{2}{6}\above{1pt}{1pt}{l}{2}{112}\above{1pt}{1pt}{}{2}{9}\above{1pt}{1pt}{r}{2}{4}\above{1pt}{1pt}{*}{2}{1008}\above{1pt}{1pt}{b}{2}{6}\above{1pt}{1pt}{b}{2}{112}\above{1pt}{1pt}{*}{2}$}%
\nopagebreak\par%
\nopagebreak\par\leavevmode%
{$\left[\!\llap{\phantom{%
\begingroup \smaller\smaller\smaller
\endgroup%
}}\!\right]$}%

\medskip%
%
\leavevmode\llap{}%
$W_{230}$%
\qquad\llap{52} lattices, $\chi=36$%
\hfill%
$\slashtwo22|22\slashtwo22|22\rtimes D_{4}$%
\nopagebreak\smallskip\hrule\nopagebreak\medskip%
%
%
\leavevmode%
${L_{230.1}}$%
{} : {$1\above{1pt}{1pt}{-2}{2}8\above{1pt}{1pt}{1}{7}{\cdot}1\above{1pt}{1pt}{1}{}3\above{1pt}{1pt}{1}{}9\above{1pt}{1pt}{1}{}{\cdot}1\above{1pt}{1pt}{2}{}7\above{1pt}{1pt}{1}{}$}\spacer%
\instructions{2m,2}%
\EasyButWeakLineBreak%
{${4}\above{1pt}{1pt}{*}{2}{36}\above{1pt}{1pt}{l}{2}{7}\above{1pt}{1pt}{}{2}{3}\above{1pt}{1pt}{}{2}{63}\above{1pt}{1pt}{r}{2}$}\relax$\,(\times2)$%
\nopagebreak\par%
\nopagebreak\par\leavevmode%
{$\left[\!\llap{\phantom{%
\begingroup \smaller\smaller\smaller
\endgroup%
}}\!\right]$}%
%
%
\hbox{}\par\smallskip%
%
%
\leavevmode%
${L_{230.2}}$%
{} : {$1\above{1pt}{1pt}{2}{0}8\above{1pt}{1pt}{-}{5}{\cdot}1\above{1pt}{1pt}{1}{}3\above{1pt}{1pt}{1}{}9\above{1pt}{1pt}{1}{}{\cdot}1\above{1pt}{1pt}{2}{}7\above{1pt}{1pt}{1}{}$}\spacer%
\instructions{m}%
\EasyButWeakLineBreak%
{${4}\above{1pt}{1pt}{l}{2}{9}\above{1pt}{1pt}{}{2}{7}\above{1pt}{1pt}{r}{2}{12}\above{1pt}{1pt}{l}{2}{63}\above{1pt}{1pt}{}{2}{1}\above{1pt}{1pt}{r}{2}{36}\above{1pt}{1pt}{*}{2}{28}\above{1pt}{1pt}{l}{2}{3}\above{1pt}{1pt}{r}{2}{252}\above{1pt}{1pt}{*}{2}$}%
\nopagebreak\par%
\nopagebreak\par\leavevmode%
{$\left[\!\llap{\phantom{%
\begingroup \smaller\smaller\smaller
\endgroup%
}}\!\right]$}%
%
%
\hbox{}\par\smallskip%
%
%
\leavevmode%
${L_{230.3}}$%
{} : {$1\above{1pt}{1pt}{2}{2}8\above{1pt}{1pt}{-}{3}{\cdot}1\above{1pt}{1pt}{1}{}3\above{1pt}{1pt}{1}{}9\above{1pt}{1pt}{1}{}{\cdot}1\above{1pt}{1pt}{2}{}7\above{1pt}{1pt}{1}{}$}\spacer%
\instructions{m}%
\EasyButWeakLineBreak%
{${1}\above{1pt}{1pt}{}{2}{9}\above{1pt}{1pt}{r}{2}{28}\above{1pt}{1pt}{*}{2}{12}\above{1pt}{1pt}{*}{2}{252}\above{1pt}{1pt}{l}{2}$}\relax$\,(\times2)$%
\nopagebreak\par%
\nopagebreak\par\leavevmode%
{$\left[\!\llap{\phantom{%
\begingroup \smaller\smaller\smaller
\endgroup%
}}\!\right]$}%
%
%
\hbox{}\par\smallskip%
%
%
\leavevmode%
${L_{230.4}}$%
{} : {$[1\above{1pt}{1pt}{-}{}2\above{1pt}{1pt}{-}{}]\above{1pt}{1pt}{}{0}16\above{1pt}{1pt}{-}{5}{\cdot}1\above{1pt}{1pt}{1}{}3\above{1pt}{1pt}{1}{}9\above{1pt}{1pt}{1}{}{\cdot}1\above{1pt}{1pt}{2}{}7\above{1pt}{1pt}{1}{}$}\spacer%
\instructions{2}%
\EasyButWeakLineBreak%
{${16}\above{1pt}{1pt}{s}{2}{36}\above{1pt}{1pt}{*}{2}{112}\above{1pt}{1pt}{l}{2}{3}\above{1pt}{1pt}{r}{2}{1008}\above{1pt}{1pt}{*}{2}{4}\above{1pt}{1pt}{s}{2}{144}\above{1pt}{1pt}{s}{2}{28}\above{1pt}{1pt}{*}{2}{48}\above{1pt}{1pt}{*}{2}{252}\above{1pt}{1pt}{s}{2}$}%
\nopagebreak\par%
\nopagebreak\par\leavevmode%
{$\left[\!\llap{\phantom{%
\begingroup \smaller\smaller\smaller
\endgroup%
}}\!\right]$}%
%
%
\hbox{}\par\smallskip%
%
%
\leavevmode%
${L_{230.5}}$%
{} : {$[1\above{1pt}{1pt}{-}{}2\above{1pt}{1pt}{1}{}]\above{1pt}{1pt}{}{2}16\above{1pt}{1pt}{1}{7}{\cdot}1\above{1pt}{1pt}{1}{}3\above{1pt}{1pt}{1}{}9\above{1pt}{1pt}{1}{}{\cdot}1\above{1pt}{1pt}{2}{}7\above{1pt}{1pt}{1}{}$}\spacer%
\instructions{m}%
\EasyButWeakLineBreak%
{${16}\above{1pt}{1pt}{l}{2}{9}\above{1pt}{1pt}{}{2}{112}\above{1pt}{1pt}{}{2}{3}\above{1pt}{1pt}{}{2}{1008}\above{1pt}{1pt}{}{2}{1}\above{1pt}{1pt}{r}{2}{144}\above{1pt}{1pt}{*}{2}{28}\above{1pt}{1pt}{s}{2}{48}\above{1pt}{1pt}{s}{2}{252}\above{1pt}{1pt}{*}{2}$}%
\nopagebreak\par%
\nopagebreak\par\leavevmode%
{$\left[\!\llap{\phantom{%
\begingroup \smaller\smaller\smaller
\endgroup%
}}\!\right]$}%
%
%
\hbox{}\par\smallskip%
%
%
\leavevmode%
${L_{230.6}}$%
{} : {$[1\above{1pt}{1pt}{1}{}2\above{1pt}{1pt}{1}{}]\above{1pt}{1pt}{}{6}16\above{1pt}{1pt}{-}{3}{\cdot}1\above{1pt}{1pt}{1}{}3\above{1pt}{1pt}{1}{}9\above{1pt}{1pt}{1}{}{\cdot}1\above{1pt}{1pt}{2}{}7\above{1pt}{1pt}{1}{}$}\EasyButWeakLineBreak%
{${16}\above{1pt}{1pt}{*}{2}{36}\above{1pt}{1pt}{s}{2}{112}\above{1pt}{1pt}{s}{2}{12}\above{1pt}{1pt}{s}{2}{1008}\above{1pt}{1pt}{s}{2}{4}\above{1pt}{1pt}{*}{2}{144}\above{1pt}{1pt}{l}{2}{7}\above{1pt}{1pt}{}{2}{48}\above{1pt}{1pt}{}{2}{63}\above{1pt}{1pt}{r}{2}$}%
\nopagebreak\par%
\nopagebreak\par\leavevmode%
{$\left[\!\llap{\phantom{%
\begingroup \smaller\smaller\smaller
\endgroup%
}}\!\right]$}%
%
%
\hbox{}\par\smallskip%
%
%
\leavevmode%
${L_{230.7}}$%
{} : {$[1\above{1pt}{1pt}{1}{}2\above{1pt}{1pt}{-}{}]\above{1pt}{1pt}{}{4}16\above{1pt}{1pt}{1}{1}{\cdot}1\above{1pt}{1pt}{1}{}3\above{1pt}{1pt}{1}{}9\above{1pt}{1pt}{1}{}{\cdot}1\above{1pt}{1pt}{2}{}7\above{1pt}{1pt}{1}{}$}\spacer%
\instructions{m}%
\EasyButWeakLineBreak%
{${16}\above{1pt}{1pt}{}{2}{9}\above{1pt}{1pt}{r}{2}{112}\above{1pt}{1pt}{*}{2}{12}\above{1pt}{1pt}{*}{2}{1008}\above{1pt}{1pt}{l}{2}{1}\above{1pt}{1pt}{}{2}{144}\above{1pt}{1pt}{}{2}{7}\above{1pt}{1pt}{r}{2}{48}\above{1pt}{1pt}{l}{2}{63}\above{1pt}{1pt}{}{2}$}%
\nopagebreak\par%
\nopagebreak\par\leavevmode%
{$\left[\!\llap{\phantom{%
\begingroup \smaller\smaller\smaller
\endgroup%
}}\!\right]$}%

\medskip%
%
\leavevmode\llap{}%
$W_{231}$%
\qquad\llap{8} lattices, $\chi=96$%
\hfill%
$2\infty2|2\infty22|22\infty2|2\infty22|2\rtimes D_{4}$%
\nopagebreak\smallskip\hrule\nopagebreak\medskip%
%
%
\leavevmode%
${L_{231.1}}$%
{} : {$1\above{1pt}{1pt}{-2}{2}64\above{1pt}{1pt}{1}{1}{\cdot}1\above{1pt}{1pt}{-}{}3\above{1pt}{1pt}{1}{}9\above{1pt}{1pt}{1}{}$}\spacer%
\instructions{3}%
\EasyButWeakLineBreak%
{${576}\above{1pt}{1pt}{*}{2}{12}\above{1pt}{1pt}{48,25}{\infty z}{3}\above{1pt}{1pt}{r}{2}{576}\above{1pt}{1pt}{s}{2}{12}\above{1pt}{1pt}{48,1}{\infty z}{3}\above{1pt}{1pt}{}{2}{576}\above{1pt}{1pt}{r}{2}{2}\above{1pt}{1pt}{b}{2}$}\relax$\,(\times2)$%
\nopagebreak\par%
\nopagebreak\par\leavevmode%
{$\left[\!\llap{\phantom{%
\begingroup \smaller\smaller\smaller
\endgroup%
}}\!\right]$}%

\medskip%
%
\leavevmode\llap{}%
$W_{232}$%
\qquad\llap{16} lattices, $\chi=12$%
\hfill%
$2|222|22\rtimes D_{2}$%
\nopagebreak\smallskip\hrule\nopagebreak\medskip%
%
%
\leavevmode%
${L_{232.1}}$%
{} : {$[1\above{1pt}{1pt}{1}{}2\above{1pt}{1pt}{-}{}]\above{1pt}{1pt}{}{4}32\above{1pt}{1pt}{1}{7}{\cdot}1\above{1pt}{1pt}{2}{}3\above{1pt}{1pt}{-}{}$}\EasyButWeakLineBreak%
{${96}\above{1pt}{1pt}{s}{2}{8}\above{1pt}{1pt}{*}{2}{96}\above{1pt}{1pt}{l}{2}{1}\above{1pt}{1pt}{}{2}{6}\above{1pt}{1pt}{r}{2}{4}\above{1pt}{1pt}{*}{2}$}%
\nopagebreak\par%
\nopagebreak\par\leavevmode%
{$\left[\!\llap{\phantom{%
\begingroup \smaller\smaller\smaller\begin{tabular}{@{}c@{}}%
0\\0\\0
\end{tabular}\endgroup%
}}\right.$}%
\begingroup \smaller\smaller\smaller\begin{tabular}{@{}c@{}}%
8160\\3936\\-96
\end{tabular}\endgroup%
\kern3pt%
\begingroup \smaller\smaller\smaller\begin{tabular}{@{}c@{}}%
3936\\1898\\-46
\end{tabular}\endgroup%
\kern3pt%
\begingroup \smaller\smaller\smaller\begin{tabular}{@{}c@{}}%
-96\\-46\\1
\end{tabular}\endgroup%
{$\left.\llap{\phantom{%
\begingroup \smaller\smaller\smaller\begin{tabular}{@{}c@{}}%
0\\0\\0
\end{tabular}\endgroup%
}}\!\right]$}%
\EasyButWeakLineBreak%
{$\left[\!\llap{\phantom{%
\begingroup \smaller\smaller\smaller\begin{tabular}{@{}c@{}}%
0\\0\\0
\end{tabular}\endgroup%
}}\right.$}%
\begingroup \smaller\smaller\smaller\begin{tabular}{@{}c@{}}%
31\\-72\\-288
\end{tabular}\endgroup%
\HardButStrongLineBreak\kern3pt%
\begingroup \smaller\smaller\smaller\begin{tabular}{@{}c@{}}%
-1\\2\\0
\end{tabular}\endgroup%
\HardButStrongLineBreak\kern3pt%
\begingroup \smaller\smaller\smaller\begin{tabular}{@{}c@{}}%
-11\\24\\48
\end{tabular}\endgroup%
\HardButStrongLineBreak\kern3pt%
\begingroup \smaller\smaller\smaller\begin{tabular}{@{}c@{}}%
0\\0\\-1
\end{tabular}\endgroup%
\HardButStrongLineBreak\kern3pt%
\begingroup \smaller\smaller\smaller\begin{tabular}{@{}c@{}}%
4\\-9\\-30
\end{tabular}\endgroup%
\HardButStrongLineBreak\kern3pt%
\begingroup \smaller\smaller\smaller\begin{tabular}{@{}c@{}}%
7\\-16\\-58
\end{tabular}\endgroup%
{$\left.\llap{\phantom{%
\begingroup \smaller\smaller\smaller\begin{tabular}{@{}c@{}}%
0\\0\\0
\end{tabular}\endgroup%
}}\!\right]$}%
%
%
\hbox{}\par\smallskip%
%
%
\leavevmode%
${L_{232.2}}$%
{} : {$[1\above{1pt}{1pt}{1}{}2\above{1pt}{1pt}{1}{}]\above{1pt}{1pt}{}{0}32\above{1pt}{1pt}{-}{3}{\cdot}1\above{1pt}{1pt}{2}{}3\above{1pt}{1pt}{-}{}$}\EasyButWeakLineBreak%
{${96}\above{1pt}{1pt}{}{2}{2}\above{1pt}{1pt}{r}{2}{96}\above{1pt}{1pt}{s}{2}{4}\above{1pt}{1pt}{*}{2}{24}\above{1pt}{1pt}{l}{2}{1}\above{1pt}{1pt}{}{2}$}%
\nopagebreak\par%
\nopagebreak\par\leavevmode%
{$\left[\!\llap{\phantom{%
\begingroup \smaller\smaller\smaller\begin{tabular}{@{}c@{}}%
0\\0\\0
\end{tabular}\endgroup%
}}\right.$}%
\begingroup \smaller\smaller\smaller\begin{tabular}{@{}c@{}}%
-20640\\384\\384
\end{tabular}\endgroup%
\kern3pt%
\begingroup \smaller\smaller\smaller\begin{tabular}{@{}c@{}}%
384\\-2\\-8
\end{tabular}\endgroup%
\kern3pt%
\begingroup \smaller\smaller\smaller\begin{tabular}{@{}c@{}}%
384\\-8\\-7
\end{tabular}\endgroup%
{$\left.\llap{\phantom{%
\begingroup \smaller\smaller\smaller\begin{tabular}{@{}c@{}}%
0\\0\\0
\end{tabular}\endgroup%
}}\!\right]$}%
\EasyButWeakLineBreak%
{$\left[\!\llap{\phantom{%
\begingroup \smaller\smaller\smaller\begin{tabular}{@{}c@{}}%
0\\0\\0
\end{tabular}\endgroup%
}}\right.$}%
\begingroup \smaller\smaller\smaller\begin{tabular}{@{}c@{}}%
19\\144\\864
\end{tabular}\endgroup%
\HardButStrongLineBreak\kern3pt%
\begingroup \smaller\smaller\smaller\begin{tabular}{@{}c@{}}%
1\\7\\46
\end{tabular}\endgroup%
\HardButStrongLineBreak\kern3pt%
\begingroup \smaller\smaller\smaller\begin{tabular}{@{}c@{}}%
1\\0\\48
\end{tabular}\endgroup%
\HardButStrongLineBreak\kern3pt%
\begingroup \smaller\smaller\smaller\begin{tabular}{@{}c@{}}%
-1\\-8\\-46
\end{tabular}\endgroup%
\HardButStrongLineBreak\kern3pt%
\begingroup \smaller\smaller\smaller\begin{tabular}{@{}c@{}}%
-1\\-6\\-48
\end{tabular}\endgroup%
\HardButStrongLineBreak\kern3pt%
\begingroup \smaller\smaller\smaller\begin{tabular}{@{}c@{}}%
1\\8\\45
\end{tabular}\endgroup%
{$\left.\llap{\phantom{%
\begingroup \smaller\smaller\smaller\begin{tabular}{@{}c@{}}%
0\\0\\0
\end{tabular}\endgroup%
}}\!\right]$}%
%
%
\hbox{}\par\smallskip%
%
%
\leavevmode%
${L_{232.3}}$%
{} : {$1\above{1pt}{1pt}{-}{3}4\above{1pt}{1pt}{1}{1}32\above{1pt}{1pt}{1}{7}{\cdot}1\above{1pt}{1pt}{2}{}3\above{1pt}{1pt}{1}{}$}\EasyButWeakLineBreak%
{${12}\above{1pt}{1pt}{l}{2}{4}\above{1pt}{1pt}{}{2}{3}\above{1pt}{1pt}{r}{2}{32}\above{1pt}{1pt}{*}{2}{48}\above{1pt}{1pt}{s}{2}{32}\above{1pt}{1pt}{*}{2}$}%
\nopagebreak\par%
\nopagebreak\par\leavevmode%
{$\left[\!\llap{\phantom{%
\begingroup \smaller\smaller\smaller\begin{tabular}{@{}c@{}}%
0\\0\\0
\end{tabular}\endgroup%
}}\right.$}%
\begingroup \smaller\smaller\smaller\begin{tabular}{@{}c@{}}%
-11040\\-288\\384
\end{tabular}\endgroup%
\kern3pt%
\begingroup \smaller\smaller\smaller\begin{tabular}{@{}c@{}}%
-288\\4\\8
\end{tabular}\endgroup%
\kern3pt%
\begingroup \smaller\smaller\smaller\begin{tabular}{@{}c@{}}%
384\\8\\-13
\end{tabular}\endgroup%
{$\left.\llap{\phantom{%
\begingroup \smaller\smaller\smaller\begin{tabular}{@{}c@{}}%
0\\0\\0
\end{tabular}\endgroup%
}}\!\right]$}%
\EasyButWeakLineBreak%
{$\left[\!\llap{\phantom{%
\begingroup \smaller\smaller\smaller\begin{tabular}{@{}c@{}}%
0\\0\\0
\end{tabular}\endgroup%
}}\right.$}%
\begingroup \smaller\smaller\smaller\begin{tabular}{@{}c@{}}%
-1\\-12\\-42
\end{tabular}\endgroup%
\HardButStrongLineBreak\kern3pt%
\begingroup \smaller\smaller\smaller\begin{tabular}{@{}c@{}}%
1\\5\\32
\end{tabular}\endgroup%
\HardButStrongLineBreak\kern3pt%
\begingroup \smaller\smaller\smaller\begin{tabular}{@{}c@{}}%
1\\6\\33
\end{tabular}\endgroup%
\HardButStrongLineBreak\kern3pt%
\begingroup \smaller\smaller\smaller\begin{tabular}{@{}c@{}}%
-1\\-4\\-32
\end{tabular}\endgroup%
\HardButStrongLineBreak\kern3pt%
\begingroup \smaller\smaller\smaller\begin{tabular}{@{}c@{}}%
-5\\-30\\-168
\end{tabular}\endgroup%
\HardButStrongLineBreak\kern3pt%
\begingroup \smaller\smaller\smaller\begin{tabular}{@{}c@{}}%
-5\\-36\\-176
\end{tabular}\endgroup%
{$\left.\llap{\phantom{%
\begingroup \smaller\smaller\smaller\begin{tabular}{@{}c@{}}%
0\\0\\0
\end{tabular}\endgroup%
}}\!\right]$}%
%
%
\hbox{}\par\smallskip%
%
%
\leavevmode%
${L_{232.4}}$%
{} : {$1\above{1pt}{1pt}{-}{3}4\above{1pt}{1pt}{1}{7}32\above{1pt}{1pt}{1}{1}{\cdot}1\above{1pt}{1pt}{2}{}3\above{1pt}{1pt}{1}{}$}\EasyButWeakLineBreak%
{${12}\above{1pt}{1pt}{*}{2}{16}\above{1pt}{1pt}{l}{2}{3}\above{1pt}{1pt}{}{2}{32}\above{1pt}{1pt}{}{2}{12}\above{1pt}{1pt}{r}{2}{32}\above{1pt}{1pt}{s}{2}$}%
\nopagebreak\par%
\nopagebreak\par\leavevmode%
{$\left[\!\llap{\phantom{%
\begingroup \smaller\smaller\smaller\begin{tabular}{@{}c@{}}%
0\\0\\0
\end{tabular}\endgroup%
}}\right.$}%
\begingroup \smaller\smaller\smaller\begin{tabular}{@{}c@{}}%
3360\\384\\0
\end{tabular}\endgroup%
\kern3pt%
\begingroup \smaller\smaller\smaller\begin{tabular}{@{}c@{}}%
384\\44\\0
\end{tabular}\endgroup%
\kern3pt%
\begingroup \smaller\smaller\smaller\begin{tabular}{@{}c@{}}%
0\\0\\-1
\end{tabular}\endgroup%
{$\left.\llap{\phantom{%
\begingroup \smaller\smaller\smaller\begin{tabular}{@{}c@{}}%
0\\0\\0
\end{tabular}\endgroup%
}}\!\right]$}%
\EasyButWeakLineBreak%
{$\left[\!\llap{\phantom{%
\begingroup \smaller\smaller\smaller\begin{tabular}{@{}c@{}}%
0\\0\\0
\end{tabular}\endgroup%
}}\right.$}%
\begingroup \smaller\smaller\smaller\begin{tabular}{@{}c@{}}%
-1\\6\\-18
\end{tabular}\endgroup%
\HardButStrongLineBreak\kern3pt%
\begingroup \smaller\smaller\smaller\begin{tabular}{@{}c@{}}%
1\\-10\\-8
\end{tabular}\endgroup%
\HardButStrongLineBreak\kern3pt%
\begingroup \smaller\smaller\smaller\begin{tabular}{@{}c@{}}%
1\\-9\\-3
\end{tabular}\endgroup%
\HardButStrongLineBreak\kern3pt%
\begingroup \smaller\smaller\smaller\begin{tabular}{@{}c@{}}%
1\\-8\\0
\end{tabular}\endgroup%
\HardButStrongLineBreak\kern3pt%
\begingroup \smaller\smaller\smaller\begin{tabular}{@{}c@{}}%
-1\\9\\0
\end{tabular}\endgroup%
\HardButStrongLineBreak\kern3pt%
\begingroup \smaller\smaller\smaller\begin{tabular}{@{}c@{}}%
-3\\24\\-16
\end{tabular}\endgroup%
{$\left.\llap{\phantom{%
\begingroup \smaller\smaller\smaller\begin{tabular}{@{}c@{}}%
0\\0\\0
\end{tabular}\endgroup%
}}\!\right]$}%

\medskip%
%
\leavevmode\llap{}%
$W_{233}$%
\qquad\llap{32} lattices, $\chi=24$%
\hfill%
$2\slashinfty222|22\rtimes D_{2}$%
\nopagebreak\smallskip\hrule\nopagebreak\medskip%
%
%
\leavevmode%
${L_{233.1}}$%
{} : {$[1\above{1pt}{1pt}{1}{}2\above{1pt}{1pt}{1}{}]\above{1pt}{1pt}{}{2}32\above{1pt}{1pt}{-}{5}{\cdot}1\above{1pt}{1pt}{-}{}3\above{1pt}{1pt}{1}{}9\above{1pt}{1pt}{1}{}$}\spacer%
\instructions{3}%
\EasyButWeakLineBreak%
{${32}\above{1pt}{1pt}{l}{2}{3}\above{1pt}{1pt}{24,5}{\infty}{12}\above{1pt}{1pt}{*}{2}{32}\above{1pt}{1pt}{l}{2}{9}\above{1pt}{1pt}{}{2}{2}\above{1pt}{1pt}{r}{2}{36}\above{1pt}{1pt}{*}{2}$}%
\nopagebreak\par%
\nopagebreak\par\leavevmode%
{$\left[\!\llap{\phantom{%
\begingroup \smaller\smaller\smaller\begin{tabular}{@{}c@{}}%
0\\0\\0
\end{tabular}\endgroup%
}}\right.$}%
\begingroup \smaller\smaller\smaller\begin{tabular}{@{}c@{}}%
-3115872\\-67968\\25344
\end{tabular}\endgroup%
\kern3pt%
\begingroup \smaller\smaller\smaller\begin{tabular}{@{}c@{}}%
-67968\\-1482\\552
\end{tabular}\endgroup%
\kern3pt%
\begingroup \smaller\smaller\smaller\begin{tabular}{@{}c@{}}%
25344\\552\\-205
\end{tabular}\endgroup%
{$\left.\llap{\phantom{%
\begingroup \smaller\smaller\smaller\begin{tabular}{@{}c@{}}%
0\\0\\0
\end{tabular}\endgroup%
}}\!\right]$}%
\EasyButWeakLineBreak%
{$\left[\!\llap{\phantom{%
\begingroup \smaller\smaller\smaller\begin{tabular}{@{}c@{}}%
0\\0\\0
\end{tabular}\endgroup%
}}\right.$}%
\begingroup \smaller\smaller\smaller\begin{tabular}{@{}c@{}}%
7\\-440\\-320
\end{tabular}\endgroup%
\HardButStrongLineBreak\kern3pt%
\begingroup \smaller\smaller\smaller\begin{tabular}{@{}c@{}}%
3\\-190\\-141
\end{tabular}\endgroup%
\HardButStrongLineBreak\kern3pt%
\begingroup \smaller\smaller\smaller\begin{tabular}{@{}c@{}}%
7\\-448\\-342
\end{tabular}\endgroup%
\HardButStrongLineBreak\kern3pt%
\begingroup \smaller\smaller\smaller\begin{tabular}{@{}c@{}}%
11\\-712\\-560
\end{tabular}\endgroup%
\HardButStrongLineBreak\kern3pt%
\begingroup \smaller\smaller\smaller\begin{tabular}{@{}c@{}}%
1\\-69\\-63
\end{tabular}\endgroup%
\HardButStrongLineBreak\kern3pt%
\begingroup \smaller\smaller\smaller\begin{tabular}{@{}c@{}}%
-1\\63\\46
\end{tabular}\endgroup%
\HardButStrongLineBreak\kern3pt%
\begingroup \smaller\smaller\smaller\begin{tabular}{@{}c@{}}%
-1\\66\\54
\end{tabular}\endgroup%
{$\left.\llap{\phantom{%
\begingroup \smaller\smaller\smaller\begin{tabular}{@{}c@{}}%
0\\0\\0
\end{tabular}\endgroup%
}}\!\right]$}%
%
%
\hbox{}\par\smallskip%
%
%
\leavevmode%
${L_{233.2}}$%
{} : {$[1\above{1pt}{1pt}{-}{}2\above{1pt}{1pt}{1}{}]\above{1pt}{1pt}{}{2}32\above{1pt}{1pt}{1}{1}{\cdot}1\above{1pt}{1pt}{-}{}3\above{1pt}{1pt}{1}{}9\above{1pt}{1pt}{1}{}$}\spacer%
\instructions{3}%
\EasyButWeakLineBreak%
{${32}\above{1pt}{1pt}{}{2}{3}\above{1pt}{1pt}{24,17}{\infty}{12}\above{1pt}{1pt}{s}{2}{32}\above{1pt}{1pt}{s}{2}{36}\above{1pt}{1pt}{*}{2}{8}\above{1pt}{1pt}{l}{2}{9}\above{1pt}{1pt}{}{2}$}%
\nopagebreak\par%
\nopagebreak\par\leavevmode%
{$\left[\!\llap{\phantom{%
\begingroup \smaller\smaller\smaller\begin{tabular}{@{}c@{}}%
0\\0\\0
\end{tabular}\endgroup%
}}\right.$}%
\begingroup \smaller\smaller\smaller\begin{tabular}{@{}c@{}}%
-4736736\\310176\\-9216
\end{tabular}\endgroup%
\kern3pt%
\begingroup \smaller\smaller\smaller\begin{tabular}{@{}c@{}}%
310176\\-20310\\606
\end{tabular}\endgroup%
\kern3pt%
\begingroup \smaller\smaller\smaller\begin{tabular}{@{}c@{}}%
-9216\\606\\-13
\end{tabular}\endgroup%
{$\left.\llap{\phantom{%
\begingroup \smaller\smaller\smaller\begin{tabular}{@{}c@{}}%
0\\0\\0
\end{tabular}\endgroup%
}}\!\right]$}%
\EasyButWeakLineBreak%
{$\left[\!\llap{\phantom{%
\begingroup \smaller\smaller\smaller\begin{tabular}{@{}c@{}}%
0\\0\\0
\end{tabular}\endgroup%
}}\right.$}%
\begingroup \smaller\smaller\smaller\begin{tabular}{@{}c@{}}%
67\\1008\\-512
\end{tabular}\endgroup%
\HardButStrongLineBreak\kern3pt%
\begingroup \smaller\smaller\smaller\begin{tabular}{@{}c@{}}%
23\\346\\-177
\end{tabular}\endgroup%
\HardButStrongLineBreak\kern3pt%
\begingroup \smaller\smaller\smaller\begin{tabular}{@{}c@{}}%
25\\376\\-198
\end{tabular}\endgroup%
\HardButStrongLineBreak\kern3pt%
\begingroup \smaller\smaller\smaller\begin{tabular}{@{}c@{}}%
-17\\-256\\112
\end{tabular}\endgroup%
\HardButStrongLineBreak\kern3pt%
\begingroup \smaller\smaller\smaller\begin{tabular}{@{}c@{}}%
-65\\-978\\486
\end{tabular}\endgroup%
\HardButStrongLineBreak\kern3pt%
\begingroup \smaller\smaller\smaller\begin{tabular}{@{}c@{}}%
-23\\-346\\176
\end{tabular}\endgroup%
\HardButStrongLineBreak\kern3pt%
\begingroup \smaller\smaller\smaller\begin{tabular}{@{}c@{}}%
-1\\-15\\9
\end{tabular}\endgroup%
{$\left.\llap{\phantom{%
\begingroup \smaller\smaller\smaller\begin{tabular}{@{}c@{}}%
0\\0\\0
\end{tabular}\endgroup%
}}\!\right]$}%
%
%
\hbox{}\par\smallskip%
%
%
\leavevmode%
${L_{233.3}}$%
{} : {$1\above{1pt}{1pt}{1}{1}4\above{1pt}{1pt}{1}{7}32\above{1pt}{1pt}{-}{3}{\cdot}1\above{1pt}{1pt}{1}{}3\above{1pt}{1pt}{-}{}9\above{1pt}{1pt}{-}{}$}\spacer%
\instructions{3}%
\EasyButWeakLineBreak%
{${1}\above{1pt}{1pt}{}{2}{96}\above{1pt}{1pt}{12,11}{\infty}{96}\above{1pt}{1pt}{s}{2}{4}\above{1pt}{1pt}{*}{2}{288}\above{1pt}{1pt}{s}{2}{16}\above{1pt}{1pt}{*}{2}{288}\above{1pt}{1pt}{l}{2}$}%
\nopagebreak\par%
\nopagebreak\par\leavevmode%
{$\left[\!\llap{\phantom{%
\begingroup \smaller\smaller\smaller\begin{tabular}{@{}c@{}}%
0\\0\\0
\end{tabular}\endgroup%
}}\right.$}%
\begingroup \smaller\smaller\smaller\begin{tabular}{@{}c@{}}%
-277920\\-7776\\2592
\end{tabular}\endgroup%
\kern3pt%
\begingroup \smaller\smaller\smaller\begin{tabular}{@{}c@{}}%
-7776\\-84\\60
\end{tabular}\endgroup%
\kern3pt%
\begingroup \smaller\smaller\smaller\begin{tabular}{@{}c@{}}%
2592\\60\\-23
\end{tabular}\endgroup%
{$\left.\llap{\phantom{%
\begingroup \smaller\smaller\smaller\begin{tabular}{@{}c@{}}%
0\\0\\0
\end{tabular}\endgroup%
}}\!\right]$}%
\EasyButWeakLineBreak%
{$\left[\!\llap{\phantom{%
\begingroup \smaller\smaller\smaller\begin{tabular}{@{}c@{}}%
0\\0\\0
\end{tabular}\endgroup%
}}\right.$}%
\begingroup \smaller\smaller\smaller\begin{tabular}{@{}c@{}}%
1\\14\\149
\end{tabular}\endgroup%
\HardButStrongLineBreak\kern3pt%
\begingroup \smaller\smaller\smaller\begin{tabular}{@{}c@{}}%
11\\152\\1632
\end{tabular}\endgroup%
\HardButStrongLineBreak\kern3pt%
\begingroup \smaller\smaller\smaller\begin{tabular}{@{}c@{}}%
5\\64\\720
\end{tabular}\endgroup%
\HardButStrongLineBreak\kern3pt%
\begingroup \smaller\smaller\smaller\begin{tabular}{@{}c@{}}%
-1\\-16\\-158
\end{tabular}\endgroup%
\HardButStrongLineBreak\kern3pt%
\begingroup \smaller\smaller\smaller\begin{tabular}{@{}c@{}}%
-19\\-276\\-2880
\end{tabular}\endgroup%
\HardButStrongLineBreak\kern3pt%
\begingroup \smaller\smaller\smaller\begin{tabular}{@{}c@{}}%
-3\\-42\\-448
\end{tabular}\endgroup%
\HardButStrongLineBreak\kern3pt%
\begingroup \smaller\smaller\smaller\begin{tabular}{@{}c@{}}%
-1\\-12\\-144
\end{tabular}\endgroup%
{$\left.\llap{\phantom{%
\begingroup \smaller\smaller\smaller\begin{tabular}{@{}c@{}}%
0\\0\\0
\end{tabular}\endgroup%
}}\!\right]$}%
%
%
\hbox{}\par\smallskip%
%
%
\leavevmode%
${L_{233.4}}$%
{} : {$1\above{1pt}{1pt}{-}{5}4\above{1pt}{1pt}{1}{1}32\above{1pt}{1pt}{1}{1}{\cdot}1\above{1pt}{1pt}{1}{}3\above{1pt}{1pt}{-}{}9\above{1pt}{1pt}{-}{}$}\spacer%
\instructions{3}%
\EasyButWeakLineBreak%
{${1}\above{1pt}{1pt}{r}{2}{96}\above{1pt}{1pt}{24,23}{\infty z}{96}\above{1pt}{1pt}{*}{2}{4}\above{1pt}{1pt}{s}{2}{288}\above{1pt}{1pt}{l}{2}{4}\above{1pt}{1pt}{}{2}{288}\above{1pt}{1pt}{}{2}$}%
\nopagebreak\par%
\nopagebreak\par\leavevmode%
{$\left[\!\llap{\phantom{%
\begingroup \smaller\smaller\smaller\begin{tabular}{@{}c@{}}%
0\\0\\0
\end{tabular}\endgroup%
}}\right.$}%
\begingroup \smaller\smaller\smaller\begin{tabular}{@{}c@{}}%
288\\0\\0
\end{tabular}\endgroup%
\kern3pt%
\begingroup \smaller\smaller\smaller\begin{tabular}{@{}c@{}}%
0\\-156\\-60
\end{tabular}\endgroup%
\kern3pt%
\begingroup \smaller\smaller\smaller\begin{tabular}{@{}c@{}}%
0\\-60\\-23
\end{tabular}\endgroup%
{$\left.\llap{\phantom{%
\begingroup \smaller\smaller\smaller\begin{tabular}{@{}c@{}}%
0\\0\\0
\end{tabular}\endgroup%
}}\!\right]$}%
\EasyButWeakLineBreak%
{$\left[\!\llap{\phantom{%
\begingroup \smaller\smaller\smaller\begin{tabular}{@{}c@{}}%
0\\0\\0
\end{tabular}\endgroup%
}}\right.$}%
\begingroup \smaller\smaller\smaller\begin{tabular}{@{}c@{}}%
0\\-2\\5
\end{tabular}\endgroup%
\HardButStrongLineBreak\kern3pt%
\begingroup \smaller\smaller\smaller\begin{tabular}{@{}c@{}}%
-1\\-20\\48
\end{tabular}\endgroup%
\HardButStrongLineBreak\kern3pt%
\begingroup \smaller\smaller\smaller\begin{tabular}{@{}c@{}}%
-3\\-4\\0
\end{tabular}\endgroup%
\HardButStrongLineBreak\kern3pt%
\begingroup \smaller\smaller\smaller\begin{tabular}{@{}c@{}}%
-1\\4\\-14
\end{tabular}\endgroup%
\HardButStrongLineBreak\kern3pt%
\begingroup \smaller\smaller\smaller\begin{tabular}{@{}c@{}}%
-5\\48\\-144
\end{tabular}\endgroup%
\HardButStrongLineBreak\kern3pt%
\begingroup \smaller\smaller\smaller\begin{tabular}{@{}c@{}}%
0\\3\\-8
\end{tabular}\endgroup%
\HardButStrongLineBreak\kern3pt%
\begingroup \smaller\smaller\smaller\begin{tabular}{@{}c@{}}%
1\\0\\0
\end{tabular}\endgroup%
{$\left.\llap{\phantom{%
\begingroup \smaller\smaller\smaller\begin{tabular}{@{}c@{}}%
0\\0\\0
\end{tabular}\endgroup%
}}\!\right]$}%

\medskip%
%
\leavevmode\llap{}%
$W_{234}$%
\qquad\llap{6} lattices, $\chi=12$%
\hfill%
$\slashtwo|\slashtwo|\slashtwo|\slashtwo|\slashtwo|\slashtwo|\rtimes D_{12}$%
\nopagebreak\smallskip\hrule\nopagebreak\medskip%
%
%
\leavevmode%
${L_{234.1}}$%
{} : {$1\above{1pt}{1pt}{2}{2}16\above{1pt}{1pt}{-}{5}{\cdot}1\above{1pt}{1pt}{1}{}3\above{1pt}{1pt}{-}{}9\above{1pt}{1pt}{1}{}$}\EasyButWeakLineBreak%
{${36}\above{1pt}{1pt}{s}{2}{16}\above{1pt}{1pt}{l}{2}{9}\above{1pt}{1pt}{}{2}{1}\above{1pt}{1pt}{r}{2}{144}\above{1pt}{1pt}{s}{2}{4}\above{1pt}{1pt}{*}{2}$}%
\nopagebreak\par%
\nopagebreak\par\leavevmode%
{$\left[\!\llap{\phantom{%
\begingroup \smaller\smaller\smaller\begin{tabular}{@{}c@{}}%
0\\0\\0
\end{tabular}\endgroup%
}}\right.$}%
\begingroup \smaller\smaller\smaller\begin{tabular}{@{}c@{}}%
-22320\\432\\144
\end{tabular}\endgroup%
\kern3pt%
\begingroup \smaller\smaller\smaller\begin{tabular}{@{}c@{}}%
432\\-3\\-6
\end{tabular}\endgroup%
\kern3pt%
\begingroup \smaller\smaller\smaller\begin{tabular}{@{}c@{}}%
144\\-6\\1
\end{tabular}\endgroup%
{$\left.\llap{\phantom{%
\begingroup \smaller\smaller\smaller\begin{tabular}{@{}c@{}}%
0\\0\\0
\end{tabular}\endgroup%
}}\!\right]$}%
\EasyButWeakLineBreak%
{$\left[\!\llap{\phantom{%
\begingroup \smaller\smaller\smaller\begin{tabular}{@{}c@{}}%
0\\0\\0
\end{tabular}\endgroup%
}}\right.$}%
\begingroup \smaller\smaller\smaller\begin{tabular}{@{}c@{}}%
-1\\-36\\-54
\end{tabular}\endgroup%
\HardButStrongLineBreak\kern3pt%
\begingroup \smaller\smaller\smaller\begin{tabular}{@{}c@{}}%
1\\32\\56
\end{tabular}\endgroup%
\HardButStrongLineBreak\kern3pt%
\begingroup \smaller\smaller\smaller\begin{tabular}{@{}c@{}}%
1\\33\\54
\end{tabular}\endgroup%
\HardButStrongLineBreak\kern3pt%
\begingroup \smaller\smaller\smaller\begin{tabular}{@{}c@{}}%
0\\0\\-1
\end{tabular}\endgroup%
\HardButStrongLineBreak\kern3pt%
\begingroup \smaller\smaller\smaller\begin{tabular}{@{}c@{}}%
-5\\-168\\-288
\end{tabular}\endgroup%
\HardButStrongLineBreak\kern3pt%
\begingroup \smaller\smaller\smaller\begin{tabular}{@{}c@{}}%
-1\\-34\\-56
\end{tabular}\endgroup%
{$\left.\llap{\phantom{%
\begingroup \smaller\smaller\smaller\begin{tabular}{@{}c@{}}%
0\\0\\0
\end{tabular}\endgroup%
}}\!\right]$}%
%
%
\hbox{}\par\smallskip%
%
%
\leavevmode%
${L_{234.2}}$%
{} : {$1\above{1pt}{1pt}{-}{5}4\above{1pt}{1pt}{1}{1}16\above{1pt}{1pt}{1}{1}{\cdot}1\above{1pt}{1pt}{1}{}3\above{1pt}{1pt}{-}{}9\above{1pt}{1pt}{1}{}$}\spacer%
\instructions{3}%
\EasyButWeakLineBreak%
{${36}\above{1pt}{1pt}{r}{2}{16}\above{1pt}{1pt}{s}{2}{36}\above{1pt}{1pt}{l}{2}{4}\above{1pt}{1pt}{}{2}{144}\above{1pt}{1pt}{}{2}{1}\above{1pt}{1pt}{}{2}$}%
\nopagebreak\par%
shares genus with 2-dual${}\iso{}$3-dual; isometric to own %
2.3-dual\nopagebreak\par%
\nopagebreak\par\leavevmode%
{$\left[\!\llap{\phantom{%
\begingroup \smaller\smaller\smaller\begin{tabular}{@{}c@{}}%
0\\0\\0
\end{tabular}\endgroup%
}}\right.$}%
\begingroup \smaller\smaller\smaller\begin{tabular}{@{}c@{}}%
-3312\\720\\0
\end{tabular}\endgroup%
\kern3pt%
\begingroup \smaller\smaller\smaller\begin{tabular}{@{}c@{}}%
720\\-12\\-12
\end{tabular}\endgroup%
\kern3pt%
\begingroup \smaller\smaller\smaller\begin{tabular}{@{}c@{}}%
0\\-12\\1
\end{tabular}\endgroup%
{$\left.\llap{\phantom{%
\begingroup \smaller\smaller\smaller\begin{tabular}{@{}c@{}}%
0\\0\\0
\end{tabular}\endgroup%
}}\!\right]$}%
\EasyButWeakLineBreak%
{$\left[\!\llap{\phantom{%
\begingroup \smaller\smaller\smaller\begin{tabular}{@{}c@{}}%
0\\0\\0
\end{tabular}\endgroup%
}}\right.$}%
\begingroup \smaller\smaller\smaller\begin{tabular}{@{}c@{}}%
2\\9\\108
\end{tabular}\endgroup%
\HardButStrongLineBreak\kern3pt%
\begingroup \smaller\smaller\smaller\begin{tabular}{@{}c@{}}%
1\\4\\56
\end{tabular}\endgroup%
\HardButStrongLineBreak\kern3pt%
\begingroup \smaller\smaller\smaller\begin{tabular}{@{}c@{}}%
-1\\-6\\-54
\end{tabular}\endgroup%
\HardButStrongLineBreak\kern3pt%
\begingroup \smaller\smaller\smaller\begin{tabular}{@{}c@{}}%
-1\\-5\\-56
\end{tabular}\endgroup%
\HardButStrongLineBreak\kern3pt%
\begingroup \smaller\smaller\smaller\begin{tabular}{@{}c@{}}%
-5\\-24\\-288
\end{tabular}\endgroup%
\HardButStrongLineBreak\kern3pt%
\begingroup \smaller\smaller\smaller\begin{tabular}{@{}c@{}}%
0\\0\\-1
\end{tabular}\endgroup%
{$\left.\llap{\phantom{%
\begingroup \smaller\smaller\smaller\begin{tabular}{@{}c@{}}%
0\\0\\0
\end{tabular}\endgroup%
}}\!\right]$}%
%
%
%
%
%
%
%
%
%
%
%
%
%
%
%
%
%
%

\medskip%
%
\leavevmode\llap{}%
$W_{235}$%
\qquad\llap{34} lattices, $\chi=24$%
\hfill%
$2|22|22|22|2\rtimes D_{4}$%
\nopagebreak\smallskip\hrule\nopagebreak\medskip%
%
%
\leavevmode%
${L_{235.1}}$%
{} : {$[1\above{1pt}{1pt}{1}{}2\above{1pt}{1pt}{1}{}]\above{1pt}{1pt}{}{6}32\above{1pt}{1pt}{-}{5}{\cdot}1\above{1pt}{1pt}{-}{}3\above{1pt}{1pt}{-}{}9\above{1pt}{1pt}{-}{}$}\spacer%
\instructions{3}%
\EasyButWeakLineBreak%
{${288}\above{1pt}{1pt}{*}{2}{8}\above{1pt}{1pt}{s}{2}{288}\above{1pt}{1pt}{*}{2}{24}\above{1pt}{1pt}{s}{2}{32}\above{1pt}{1pt}{*}{2}{72}\above{1pt}{1pt}{s}{2}{32}\above{1pt}{1pt}{*}{2}{24}\above{1pt}{1pt}{s}{2}$}%
\nopagebreak\par%
\nopagebreak\par\leavevmode%
{$\left[\!\llap{\phantom{%
\begingroup \smaller\smaller\smaller\begin{tabular}{@{}c@{}}%
0\\0\\0
\end{tabular}\endgroup%
}}\right.$}%
\begingroup \smaller\smaller\smaller\begin{tabular}{@{}c@{}}%
-3168\\-1728\\-864
\end{tabular}\endgroup%
\kern3pt%
\begingroup \smaller\smaller\smaller\begin{tabular}{@{}c@{}}%
-1728\\-930\\-456
\end{tabular}\endgroup%
\kern3pt%
\begingroup \smaller\smaller\smaller\begin{tabular}{@{}c@{}}%
-864\\-456\\-217
\end{tabular}\endgroup%
{$\left.\llap{\phantom{%
\begingroup \smaller\smaller\smaller\begin{tabular}{@{}c@{}}%
0\\0\\0
\end{tabular}\endgroup%
}}\!\right]$}%
\EasyButWeakLineBreak%
{$\left[\!\llap{\phantom{%
\begingroup \smaller\smaller\smaller\begin{tabular}{@{}c@{}}%
0\\0\\0
\end{tabular}\endgroup%
}}\right.$}%
\begingroup \smaller\smaller\smaller\begin{tabular}{@{}c@{}}%
119\\-360\\288
\end{tabular}\endgroup%
\HardButStrongLineBreak\kern3pt%
\begingroup \smaller\smaller\smaller\begin{tabular}{@{}c@{}}%
11\\-34\\28
\end{tabular}\endgroup%
\HardButStrongLineBreak\kern3pt%
\begingroup \smaller\smaller\smaller\begin{tabular}{@{}c@{}}%
53\\-168\\144
\end{tabular}\endgroup%
\HardButStrongLineBreak\kern3pt%
\begingroup \smaller\smaller\smaller\begin{tabular}{@{}c@{}}%
-1\\2\\0
\end{tabular}\endgroup%
\HardButStrongLineBreak\kern3pt%
\begingroup \smaller\smaller\smaller\begin{tabular}{@{}c@{}}%
-13\\40\\-32
\end{tabular}\endgroup%
\HardButStrongLineBreak\kern3pt%
\begingroup \smaller\smaller\smaller\begin{tabular}{@{}c@{}}%
-13\\42\\-36
\end{tabular}\endgroup%
\HardButStrongLineBreak\kern3pt%
\begingroup \smaller\smaller\smaller\begin{tabular}{@{}c@{}}%
9\\-24\\16
\end{tabular}\endgroup%
\HardButStrongLineBreak\kern3pt%
\begingroup \smaller\smaller\smaller\begin{tabular}{@{}c@{}}%
21\\-62\\48
\end{tabular}\endgroup%
{$\left.\llap{\phantom{%
\begingroup \smaller\smaller\smaller\begin{tabular}{@{}c@{}}%
0\\0\\0
\end{tabular}\endgroup%
}}\!\right]$}%
%
%
\hbox{}\par\smallskip%
%
%
\leavevmode%
${L_{235.2}}$%
{} : {$[1\above{1pt}{1pt}{-}{}2\above{1pt}{1pt}{1}{}]\above{1pt}{1pt}{}{2}64\above{1pt}{1pt}{1}{1}{\cdot}1\above{1pt}{1pt}{1}{}3\above{1pt}{1pt}{1}{}9\above{1pt}{1pt}{1}{}$}\spacer%
\instructions{3m,3,m}%
\EasyButWeakLineBreak%
{${576}\above{1pt}{1pt}{}{2}{1}\above{1pt}{1pt}{r}{2}{576}\above{1pt}{1pt}{*}{2}{12}\above{1pt}{1pt}{s}{2}{64}\above{1pt}{1pt}{s}{2}{36}\above{1pt}{1pt}{*}{2}{64}\above{1pt}{1pt}{l}{2}{3}\above{1pt}{1pt}{}{2}$}%
\nopagebreak\par%
shares genus with 3-dual\nopagebreak\par%
\nopagebreak\par\leavevmode%
{$\left[\!\llap{\phantom{%
\begingroup \smaller\smaller\smaller\begin{tabular}{@{}c@{}}%
0\\0\\0
\end{tabular}\endgroup%
}}\right.$}%
\begingroup \smaller\smaller\smaller\begin{tabular}{@{}c@{}}%
318528\\159552\\-576
\end{tabular}\endgroup%
\kern3pt%
\begingroup \smaller\smaller\smaller\begin{tabular}{@{}c@{}}%
159552\\79914\\-288
\end{tabular}\endgroup%
\kern3pt%
\begingroup \smaller\smaller\smaller\begin{tabular}{@{}c@{}}%
-576\\-288\\1
\end{tabular}\endgroup%
{$\left.\llap{\phantom{%
\begingroup \smaller\smaller\smaller\begin{tabular}{@{}c@{}}%
0\\0\\0
\end{tabular}\endgroup%
}}\!\right]$}%
\EasyButWeakLineBreak%
{$\left[\!\llap{\phantom{%
\begingroup \smaller\smaller\smaller\begin{tabular}{@{}c@{}}%
0\\0\\0
\end{tabular}\endgroup%
}}\right.$}%
\begingroup \smaller\smaller\smaller\begin{tabular}{@{}c@{}}%
91\\-192\\-2880
\end{tabular}\endgroup%
\HardButStrongLineBreak\kern3pt%
\begingroup \smaller\smaller\smaller\begin{tabular}{@{}c@{}}%
0\\0\\-1
\end{tabular}\endgroup%
\HardButStrongLineBreak\kern3pt%
\begingroup \smaller\smaller\smaller\begin{tabular}{@{}c@{}}%
-23\\48\\576
\end{tabular}\endgroup%
\HardButStrongLineBreak\kern3pt%
\begingroup \smaller\smaller\smaller\begin{tabular}{@{}c@{}}%
-1\\2\\6
\end{tabular}\endgroup%
\HardButStrongLineBreak\kern3pt%
\begingroup \smaller\smaller\smaller\begin{tabular}{@{}c@{}}%
15\\-32\\-544
\end{tabular}\endgroup%
\HardButStrongLineBreak\kern3pt%
\begingroup \smaller\smaller\smaller\begin{tabular}{@{}c@{}}%
17\\-36\\-558
\end{tabular}\endgroup%
\HardButStrongLineBreak\kern3pt%
\begingroup \smaller\smaller\smaller\begin{tabular}{@{}c@{}}%
53\\-112\\-1696
\end{tabular}\endgroup%
\HardButStrongLineBreak\kern3pt%
\begingroup \smaller\smaller\smaller\begin{tabular}{@{}c@{}}%
9\\-19\\-285
\end{tabular}\endgroup%
{$\left.\llap{\phantom{%
\begingroup \smaller\smaller\smaller\begin{tabular}{@{}c@{}}%
0\\0\\0
\end{tabular}\endgroup%
}}\!\right]$}%
%
%
\hbox{}\par\smallskip%
%
%
\leavevmode%
${L_{235.3}}$%
{} : {$1\above{1pt}{1pt}{1}{1}4\above{1pt}{1pt}{1}{1}32\above{1pt}{1pt}{-}{5}{\cdot}1\above{1pt}{1pt}{1}{}3\above{1pt}{1pt}{1}{}9\above{1pt}{1pt}{1}{}$}\spacer%
\instructions{3}%
\EasyButWeakLineBreak%
{${36}\above{1pt}{1pt}{l}{2}{4}\above{1pt}{1pt}{}{2}{9}\above{1pt}{1pt}{r}{2}{48}\above{1pt}{1pt}{*}{2}{4}\above{1pt}{1pt}{l}{2}{36}\above{1pt}{1pt}{}{2}{1}\above{1pt}{1pt}{r}{2}{48}\above{1pt}{1pt}{*}{2}$}%
\nopagebreak\par%
\nopagebreak\par\leavevmode%
{$\left[\!\llap{\phantom{%
\begingroup \smaller\smaller\smaller\begin{tabular}{@{}c@{}}%
0\\0\\0
\end{tabular}\endgroup%
}}\right.$}%
\begingroup \smaller\smaller\smaller\begin{tabular}{@{}c@{}}%
-1685088\\10944\\10944
\end{tabular}\endgroup%
\kern3pt%
\begingroup \smaller\smaller\smaller\begin{tabular}{@{}c@{}}%
10944\\-60\\-72
\end{tabular}\endgroup%
\kern3pt%
\begingroup \smaller\smaller\smaller\begin{tabular}{@{}c@{}}%
10944\\-72\\-71
\end{tabular}\endgroup%
{$\left.\llap{\phantom{%
\begingroup \smaller\smaller\smaller\begin{tabular}{@{}c@{}}%
0\\0\\0
\end{tabular}\endgroup%
}}\!\right]$}%
\EasyButWeakLineBreak%
{$\left[\!\llap{\phantom{%
\begingroup \smaller\smaller\smaller\begin{tabular}{@{}c@{}}%
0\\0\\0
\end{tabular}\endgroup%
}}\right.$}%
\begingroup \smaller\smaller\smaller\begin{tabular}{@{}c@{}}%
29\\336\\4122
\end{tabular}\endgroup%
\HardButStrongLineBreak\kern3pt%
\begingroup \smaller\smaller\smaller\begin{tabular}{@{}c@{}}%
4\\47\\568
\end{tabular}\endgroup%
\HardButStrongLineBreak\kern3pt%
\begingroup \smaller\smaller\smaller\begin{tabular}{@{}c@{}}%
4\\48\\567
\end{tabular}\endgroup%
\HardButStrongLineBreak\kern3pt%
\begingroup \smaller\smaller\smaller\begin{tabular}{@{}c@{}}%
-1\\-10\\-144
\end{tabular}\endgroup%
\HardButStrongLineBreak\kern3pt%
\begingroup \smaller\smaller\smaller\begin{tabular}{@{}c@{}}%
-1\\-12\\-142
\end{tabular}\endgroup%
\HardButStrongLineBreak\kern3pt%
\begingroup \smaller\smaller\smaller\begin{tabular}{@{}c@{}}%
1\\9\\144
\end{tabular}\endgroup%
\HardButStrongLineBreak\kern3pt%
\begingroup \smaller\smaller\smaller\begin{tabular}{@{}c@{}}%
3\\34\\427
\end{tabular}\endgroup%
\HardButStrongLineBreak\kern3pt%
\begingroup \smaller\smaller\smaller\begin{tabular}{@{}c@{}}%
27\\310\\3840
\end{tabular}\endgroup%
{$\left.\llap{\phantom{%
\begingroup \smaller\smaller\smaller\begin{tabular}{@{}c@{}}%
0\\0\\0
\end{tabular}\endgroup%
}}\!\right]$}%
%
%
\hbox{}\par\smallskip%
%
%
\leavevmode%
${L_{235.4}}$%
{} : {$1\above{1pt}{1pt}{-}{5}4\above{1pt}{1pt}{1}{7}32\above{1pt}{1pt}{1}{7}{\cdot}1\above{1pt}{1pt}{1}{}3\above{1pt}{1pt}{1}{}9\above{1pt}{1pt}{1}{}$}\spacer%
\instructions{3}%
\EasyButWeakLineBreak%
{${36}\above{1pt}{1pt}{*}{2}{16}\above{1pt}{1pt}{l}{2}{9}\above{1pt}{1pt}{}{2}{12}\above{1pt}{1pt}{r}{2}{4}\above{1pt}{1pt}{*}{2}{144}\above{1pt}{1pt}{l}{2}{1}\above{1pt}{1pt}{}{2}{12}\above{1pt}{1pt}{r}{2}$}%
\nopagebreak\par%
\nopagebreak\par\leavevmode%
{$\left[\!\llap{\phantom{%
\begingroup \smaller\smaller\smaller\begin{tabular}{@{}c@{}}%
0\\0\\0
\end{tabular}\endgroup%
}}\right.$}%
\begingroup \smaller\smaller\smaller\begin{tabular}{@{}c@{}}%
-2571552\\27072\\13536
\end{tabular}\endgroup%
\kern3pt%
\begingroup \smaller\smaller\smaller\begin{tabular}{@{}c@{}}%
27072\\-276\\-144
\end{tabular}\endgroup%
\kern3pt%
\begingroup \smaller\smaller\smaller\begin{tabular}{@{}c@{}}%
13536\\-144\\-71
\end{tabular}\endgroup%
{$\left.\llap{\phantom{%
\begingroup \smaller\smaller\smaller\begin{tabular}{@{}c@{}}%
0\\0\\0
\end{tabular}\endgroup%
}}\!\right]$}%
\EasyButWeakLineBreak%
{$\left[\!\llap{\phantom{%
\begingroup \smaller\smaller\smaller\begin{tabular}{@{}c@{}}%
0\\0\\0
\end{tabular}\endgroup%
}}\right.$}%
\begingroup \smaller\smaller\smaller\begin{tabular}{@{}c@{}}%
35\\828\\4986
\end{tabular}\endgroup%
\HardButStrongLineBreak\kern3pt%
\begingroup \smaller\smaller\smaller\begin{tabular}{@{}c@{}}%
9\\214\\1280
\end{tabular}\endgroup%
\HardButStrongLineBreak\kern3pt%
\begingroup \smaller\smaller\smaller\begin{tabular}{@{}c@{}}%
4\\96\\567
\end{tabular}\endgroup%
\HardButStrongLineBreak\kern3pt%
\begingroup \smaller\smaller\smaller\begin{tabular}{@{}c@{}}%
-1\\-23\\-144
\end{tabular}\endgroup%
\HardButStrongLineBreak\kern3pt%
\begingroup \smaller\smaller\smaller\begin{tabular}{@{}c@{}}%
-1\\-24\\-142
\end{tabular}\endgroup%
\HardButStrongLineBreak\kern3pt%
\begingroup \smaller\smaller\smaller\begin{tabular}{@{}c@{}}%
5\\114\\720
\end{tabular}\endgroup%
\HardButStrongLineBreak\kern3pt%
\begingroup \smaller\smaller\smaller\begin{tabular}{@{}c@{}}%
4\\94\\571
\end{tabular}\endgroup%
\HardButStrongLineBreak\kern3pt%
\begingroup \smaller\smaller\smaller\begin{tabular}{@{}c@{}}%
17\\401\\2424
\end{tabular}\endgroup%
{$\left.\llap{\phantom{%
\begingroup \smaller\smaller\smaller\begin{tabular}{@{}c@{}}%
0\\0\\0
\end{tabular}\endgroup%
}}\!\right]$}%
%
%
\hbox{}\par\smallskip%
%
%
\leavevmode%
${L_{235.5}}$%
{} : {$1\above{1pt}{1pt}{-}{3}8\above{1pt}{1pt}{1}{7}64\above{1pt}{1pt}{1}{1}{\cdot}1\above{1pt}{1pt}{1}{}3\above{1pt}{1pt}{1}{}9\above{1pt}{1pt}{1}{}$}\EasyButWeakLineBreak%
{${576}\above{1pt}{1pt}{r}{2}{4}\above{1pt}{1pt}{b}{2}{576}\above{1pt}{1pt}{s}{2}{12}\above{1pt}{1pt}{s}{2}{64}\above{1pt}{1pt}{b}{2}{36}\above{1pt}{1pt}{l}{2}{64}\above{1pt}{1pt}{}{2}{3}\above{1pt}{1pt}{}{2}$}%
\nopagebreak\par%
\nopagebreak\par\leavevmode%
{$\left[\!\llap{\phantom{%
\begingroup \smaller\smaller\smaller\begin{tabular}{@{}c@{}}%
0\\0\\0
\end{tabular}\endgroup%
}}\right.$}%
\begingroup \smaller\smaller\smaller\begin{tabular}{@{}c@{}}%
-400320\\-241920\\4608
\end{tabular}\endgroup%
\kern3pt%
\begingroup \smaller\smaller\smaller\begin{tabular}{@{}c@{}}%
-241920\\-146184\\2784
\end{tabular}\endgroup%
\kern3pt%
\begingroup \smaller\smaller\smaller\begin{tabular}{@{}c@{}}%
4608\\2784\\-53
\end{tabular}\endgroup%
{$\left.\llap{\phantom{%
\begingroup \smaller\smaller\smaller\begin{tabular}{@{}c@{}}%
0\\0\\0
\end{tabular}\endgroup%
}}\!\right]$}%
\EasyButWeakLineBreak%
{$\left[\!\llap{\phantom{%
\begingroup \smaller\smaller\smaller\begin{tabular}{@{}c@{}}%
0\\0\\0
\end{tabular}\endgroup%
}}\right.$}%
\begingroup \smaller\smaller\smaller\begin{tabular}{@{}c@{}}%
-1\\24\\1152
\end{tabular}\endgroup%
\HardButStrongLineBreak\kern3pt%
\begingroup \smaller\smaller\smaller\begin{tabular}{@{}c@{}}%
-1\\3\\70
\end{tabular}\endgroup%
\HardButStrongLineBreak\kern3pt%
\begingroup \smaller\smaller\smaller\begin{tabular}{@{}c@{}}%
-19\\48\\864
\end{tabular}\endgroup%
\HardButStrongLineBreak\kern3pt%
\begingroup \smaller\smaller\smaller\begin{tabular}{@{}c@{}}%
-1\\2\\18
\end{tabular}\endgroup%
\HardButStrongLineBreak\kern3pt%
\begingroup \smaller\smaller\smaller\begin{tabular}{@{}c@{}}%
3\\-8\\-160
\end{tabular}\endgroup%
\HardButStrongLineBreak\kern3pt%
\begingroup \smaller\smaller\smaller\begin{tabular}{@{}c@{}}%
4\\-9\\-126
\end{tabular}\endgroup%
\HardButStrongLineBreak\kern3pt%
\begingroup \smaller\smaller\smaller\begin{tabular}{@{}c@{}}%
9\\-16\\-64
\end{tabular}\endgroup%
\HardButStrongLineBreak\kern3pt%
\begingroup \smaller\smaller\smaller\begin{tabular}{@{}c@{}}%
1\\-1\\33
\end{tabular}\endgroup%
{$\left.\llap{\phantom{%
\begingroup \smaller\smaller\smaller\begin{tabular}{@{}c@{}}%
0\\0\\0
\end{tabular}\endgroup%
}}\!\right]$}%

\medskip%
%
\leavevmode\llap{}%
$W_{236}$%
\qquad\llap{8} lattices, $\chi=48$%
\hfill%
$2|222|222|222|22\rtimes D_{4}$%
\nopagebreak\smallskip\hrule\nopagebreak\medskip%
%
%
\leavevmode%
${L_{236.1}}$%
{} : {$1\above{1pt}{1pt}{-2}{6}64\above{1pt}{1pt}{1}{1}{\cdot}1\above{1pt}{1pt}{1}{}3\above{1pt}{1pt}{-}{}9\above{1pt}{1pt}{1}{}$}\spacer%
\instructions{3}%
\EasyButWeakLineBreak%
{${64}\above{1pt}{1pt}{b}{2}{6}\above{1pt}{1pt}{l}{2}{64}\above{1pt}{1pt}{}{2}{9}\above{1pt}{1pt}{r}{2}{64}\above{1pt}{1pt}{s}{2}{36}\above{1pt}{1pt}{*}{2}$}\relax$\,(\times2)$%
\nopagebreak\par%
shares genus with 3-dual\nopagebreak\par%
\nopagebreak\par\leavevmode%
{$\left[\!\llap{\phantom{%
\begingroup \smaller\smaller\smaller
\endgroup%
}}\!\right]$}%

\medskip%
%
\leavevmode\llap{}%
$W_{237}$%
\qquad\llap{32} lattices, $\chi=90$%
\hfill%
$242222222242222222\rtimes C_{2}$%
\nopagebreak\smallskip\hrule\nopagebreak\medskip%
%
%
\leavevmode%
${L_{237.1}}$%
{} : {$1\above{1pt}{1pt}{2}{2}8\above{1pt}{1pt}{-}{5}{\cdot}1\above{1pt}{1pt}{2}{}9\above{1pt}{1pt}{-}{}{\cdot}1\above{1pt}{1pt}{2}{}11\above{1pt}{1pt}{1}{}$}\spacer%
\instructions{2}%
\EasyButWeakLineBreak%
{${396}\above{1pt}{1pt}{l}{2}{1}\above{1pt}{1pt}{}{4}{2}\above{1pt}{1pt}{b}{2}{8}\above{1pt}{1pt}{*}{2}{44}\above{1pt}{1pt}{*}{2}{72}\above{1pt}{1pt}{b}{2}{2}\above{1pt}{1pt}{s}{2}{18}\above{1pt}{1pt}{b}{2}{8}\above{1pt}{1pt}{*}{2}$}\relax$\,(\times2)$%
\nopagebreak\par%
\nopagebreak\par\leavevmode%
{$\left[\!\llap{\phantom{%
\begingroup \smaller\smaller\smaller
\endgroup%
}}\!\right]$}%
%
%
\hbox{}\par\smallskip%
%
%
\leavevmode%
${L_{237.2}}$%
{} : {$1\above{1pt}{1pt}{-2}{2}8\above{1pt}{1pt}{1}{1}{\cdot}1\above{1pt}{1pt}{2}{}9\above{1pt}{1pt}{-}{}{\cdot}1\above{1pt}{1pt}{2}{}11\above{1pt}{1pt}{1}{}$}\spacer%
\instructions{m}%
\EasyButWeakLineBreak%
{${99}\above{1pt}{1pt}{r}{2}{4}\above{1pt}{1pt}{*}{4}{2}\above{1pt}{1pt}{l}{2}{8}\above{1pt}{1pt}{}{2}{11}\above{1pt}{1pt}{}{2}{72}\above{1pt}{1pt}{r}{2}{2}\above{1pt}{1pt}{b}{2}{18}\above{1pt}{1pt}{l}{2}{8}\above{1pt}{1pt}{}{2}$}\relax$\,(\times2)$%
\nopagebreak\par%
\nopagebreak\par\leavevmode%
{$\left[\!\llap{\phantom{%
\begingroup \smaller\smaller\smaller
\endgroup%
}}\!\right]$}%

\medskip%
%
\leavevmode\llap{}%
$W_{238}$%
\qquad\llap{23} lattices, $\chi=24$%
\hfill%
$2|2|2|2|2|2|2|2|\rtimes D_{8}$%
\nopagebreak\smallskip\hrule\nopagebreak\medskip%
%
%
\leavevmode%
${L_{238.1}}$%
{} : {$[1\above{1pt}{1pt}{-}{}2\above{1pt}{1pt}{1}{}]\above{1pt}{1pt}{}{2}16\above{1pt}{1pt}{1}{1}{\cdot}1\above{1pt}{1pt}{-}{}3\above{1pt}{1pt}{1}{}9\above{1pt}{1pt}{-}{}$}\spacer%
\instructions{2}%
\EasyButWeakLineBreak%
{${3}\above{1pt}{1pt}{r}{2}{8}\above{1pt}{1pt}{*}{2}{48}\above{1pt}{1pt}{*}{2}{72}\above{1pt}{1pt}{l}{2}$}\relax$\,(\times2)$%
\nopagebreak\par%
\nopagebreak\par\leavevmode%
{$\left[\!\llap{\phantom{%
\begingroup \smaller\smaller\smaller\begin{tabular}{@{}c@{}}%
0\\0\\0
\end{tabular}\endgroup%
}}\right.$}%
\begingroup \smaller\smaller\smaller\begin{tabular}{@{}c@{}}%
30096\\3168\\-576
\end{tabular}\endgroup%
\kern3pt%
\begingroup \smaller\smaller\smaller\begin{tabular}{@{}c@{}}%
3168\\318\\-60
\end{tabular}\endgroup%
\kern3pt%
\begingroup \smaller\smaller\smaller\begin{tabular}{@{}c@{}}%
-576\\-60\\11
\end{tabular}\endgroup%
{$\left.\llap{\phantom{%
\begingroup \smaller\smaller\smaller\begin{tabular}{@{}c@{}}%
0\\0\\0
\end{tabular}\endgroup%
}}\!\right]$}%
\hfil\penalty500%
{$\left[\!\llap{\phantom{%
\begingroup \smaller\smaller\smaller\begin{tabular}{@{}c@{}}%
0\\0\\0
\end{tabular}\endgroup%
}}\right.$}%
\begingroup \smaller\smaller\smaller\begin{tabular}{@{}c@{}}%
-1\\96\\576
\end{tabular}\endgroup%
\kern3pt%
\begingroup \smaller\smaller\smaller\begin{tabular}{@{}c@{}}%
0\\13\\84
\end{tabular}\endgroup%
\kern3pt%
\begingroup \smaller\smaller\smaller\begin{tabular}{@{}c@{}}%
0\\-2\\-13
\end{tabular}\endgroup%
{$\left.\llap{\phantom{%
\begingroup \smaller\smaller\smaller\begin{tabular}{@{}c@{}}%
0\\0\\0
\end{tabular}\endgroup%
}}\!\right]$}%
\EasyButWeakLineBreak%
{$\left[\!\llap{\phantom{%
\begingroup \smaller\smaller\smaller\begin{tabular}{@{}c@{}}%
0\\0\\0
\end{tabular}\endgroup%
}}\right.$}%
\begingroup \smaller\smaller\smaller\begin{tabular}{@{}c@{}}%
-1\\-4\\-75
\end{tabular}\endgroup%
\HardButStrongLineBreak\kern3pt%
\begingroup \smaller\smaller\smaller\begin{tabular}{@{}c@{}}%
-1\\-6\\-88
\end{tabular}\endgroup%
\HardButStrongLineBreak\kern3pt%
\begingroup \smaller\smaller\smaller\begin{tabular}{@{}c@{}}%
1\\-4\\24
\end{tabular}\endgroup%
\HardButStrongLineBreak\kern3pt%
\begingroup \smaller\smaller\smaller\begin{tabular}{@{}c@{}}%
5\\6\\288
\end{tabular}\endgroup%
{$\left.\llap{\phantom{%
\begingroup \smaller\smaller\smaller\begin{tabular}{@{}c@{}}%
0\\0\\0
\end{tabular}\endgroup%
}}\!\right]$}%
%
%
\hbox{}\par\smallskip%
%
%
\leavevmode%
${L_{238.2}}$%
{} : {$[1\above{1pt}{1pt}{1}{}2\above{1pt}{1pt}{1}{}]\above{1pt}{1pt}{}{6}16\above{1pt}{1pt}{-}{5}{\cdot}1\above{1pt}{1pt}{-}{}3\above{1pt}{1pt}{1}{}9\above{1pt}{1pt}{-}{}$}\spacer%
\instructions{m}%
\EasyButWeakLineBreak%
{${12}\above{1pt}{1pt}{*}{2}{8}\above{1pt}{1pt}{s}{2}{48}\above{1pt}{1pt}{s}{2}{72}\above{1pt}{1pt}{*}{2}$}\relax$\,(\times2)$%
\nopagebreak\par%
\nopagebreak\par\leavevmode%
{$\left[\!\llap{\phantom{%
\begingroup \smaller\smaller\smaller\begin{tabular}{@{}c@{}}%
0\\0\\0
\end{tabular}\endgroup%
}}\right.$}%
\begingroup \smaller\smaller\smaller\begin{tabular}{@{}c@{}}%
100944\\24480\\-720
\end{tabular}\endgroup%
\kern3pt%
\begingroup \smaller\smaller\smaller\begin{tabular}{@{}c@{}}%
24480\\5934\\-174
\end{tabular}\endgroup%
\kern3pt%
\begingroup \smaller\smaller\smaller\begin{tabular}{@{}c@{}}%
-720\\-174\\5
\end{tabular}\endgroup%
{$\left.\llap{\phantom{%
\begingroup \smaller\smaller\smaller\begin{tabular}{@{}c@{}}%
0\\0\\0
\end{tabular}\endgroup%
}}\!\right]$}%
\hfil\penalty500%
{$\left[\!\llap{\phantom{%
\begingroup \smaller\smaller\smaller\begin{tabular}{@{}c@{}}%
0\\0\\0
\end{tabular}\endgroup%
}}\right.$}%
\begingroup \smaller\smaller\smaller\begin{tabular}{@{}c@{}}%
-49\\240\\1440
\end{tabular}\endgroup%
\kern3pt%
\begingroup \smaller\smaller\smaller\begin{tabular}{@{}c@{}}%
-10\\49\\300
\end{tabular}\endgroup%
\kern3pt%
\begingroup \smaller\smaller\smaller\begin{tabular}{@{}c@{}}%
0\\0\\-1
\end{tabular}\endgroup%
{$\left.\llap{\phantom{%
\begingroup \smaller\smaller\smaller\begin{tabular}{@{}c@{}}%
0\\0\\0
\end{tabular}\endgroup%
}}\!\right]$}%
\EasyButWeakLineBreak%
{$\left[\!\llap{\phantom{%
\begingroup \smaller\smaller\smaller\begin{tabular}{@{}c@{}}%
0\\0\\0
\end{tabular}\endgroup%
}}\right.$}%
\begingroup \smaller\smaller\smaller\begin{tabular}{@{}c@{}}%
7\\-34\\-174
\end{tabular}\endgroup%
\HardButStrongLineBreak\kern3pt%
\begingroup \smaller\smaller\smaller\begin{tabular}{@{}c@{}}%
7\\-34\\-176
\end{tabular}\endgroup%
\HardButStrongLineBreak\kern3pt%
\begingroup \smaller\smaller\smaller\begin{tabular}{@{}c@{}}%
9\\-44\\-240
\end{tabular}\endgroup%
\HardButStrongLineBreak\kern3pt%
\begingroup \smaller\smaller\smaller\begin{tabular}{@{}c@{}}%
1\\-6\\-72
\end{tabular}\endgroup%
{$\left.\llap{\phantom{%
\begingroup \smaller\smaller\smaller\begin{tabular}{@{}c@{}}%
0\\0\\0
\end{tabular}\endgroup%
}}\!\right]$}%
%
%
\hbox{}\par\smallskip%
%
%
\leavevmode%
${L_{238.3}}$%
{} : {$[1\above{1pt}{1pt}{1}{}2\above{1pt}{1pt}{-}{}]\above{1pt}{1pt}{}{4}32\above{1pt}{1pt}{1}{7}{\cdot}1\above{1pt}{1pt}{1}{}3\above{1pt}{1pt}{-}{}9\above{1pt}{1pt}{1}{}$}\spacer%
\instructions{2,m}%
\EasyButWeakLineBreak%
{${96}\above{1pt}{1pt}{*}{2}{36}\above{1pt}{1pt}{l}{2}{6}\above{1pt}{1pt}{}{2}{1}\above{1pt}{1pt}{r}{2}{96}\above{1pt}{1pt}{l}{2}{9}\above{1pt}{1pt}{}{2}{6}\above{1pt}{1pt}{r}{2}{4}\above{1pt}{1pt}{*}{2}$}%
\nopagebreak\par%
\nopagebreak\par\leavevmode%
{$\left[\!\llap{\phantom{%
\begingroup \smaller\smaller\smaller\begin{tabular}{@{}c@{}}%
0\\0\\0
\end{tabular}\endgroup%
}}\right.$}%
\begingroup \smaller\smaller\smaller\begin{tabular}{@{}c@{}}%
13536\\864\\1440
\end{tabular}\endgroup%
\kern3pt%
\begingroup \smaller\smaller\smaller\begin{tabular}{@{}c@{}}%
864\\-282\\222
\end{tabular}\endgroup%
\kern3pt%
\begingroup \smaller\smaller\smaller\begin{tabular}{@{}c@{}}%
1440\\222\\103
\end{tabular}\endgroup%
{$\left.\llap{\phantom{%
\begingroup \smaller\smaller\smaller\begin{tabular}{@{}c@{}}%
0\\0\\0
\end{tabular}\endgroup%
}}\!\right]$}%
\EasyButWeakLineBreak%
{$\left[\!\llap{\phantom{%
\begingroup \smaller\smaller\smaller\begin{tabular}{@{}c@{}}%
0\\0\\0
\end{tabular}\endgroup%
}}\right.$}%
\begingroup \smaller\smaller\smaller\begin{tabular}{@{}c@{}}%
-19\\56\\144
\end{tabular}\endgroup%
\HardButStrongLineBreak\kern3pt%
\begingroup \smaller\smaller\smaller\begin{tabular}{@{}c@{}}%
-59\\174\\450
\end{tabular}\endgroup%
\HardButStrongLineBreak\kern3pt%
\begingroup \smaller\smaller\smaller\begin{tabular}{@{}c@{}}%
-22\\65\\168
\end{tabular}\endgroup%
\HardButStrongLineBreak\kern3pt%
\begingroup \smaller\smaller\smaller\begin{tabular}{@{}c@{}}%
-3\\9\\23
\end{tabular}\endgroup%
\HardButStrongLineBreak\kern3pt%
\begingroup \smaller\smaller\smaller\begin{tabular}{@{}c@{}}%
63\\-184\\-480
\end{tabular}\endgroup%
\HardButStrongLineBreak\kern3pt%
\begingroup \smaller\smaller\smaller\begin{tabular}{@{}c@{}}%
46\\-135\\-351
\end{tabular}\endgroup%
\HardButStrongLineBreak\kern3pt%
\begingroup \smaller\smaller\smaller\begin{tabular}{@{}c@{}}%
33\\-97\\-252
\end{tabular}\endgroup%
\HardButStrongLineBreak\kern3pt%
\begingroup \smaller\smaller\smaller\begin{tabular}{@{}c@{}}%
17\\-50\\-130
\end{tabular}\endgroup%
{$\left.\llap{\phantom{%
\begingroup \smaller\smaller\smaller\begin{tabular}{@{}c@{}}%
0\\0\\0
\end{tabular}\endgroup%
}}\!\right]$}%
%
%
\hbox{}\par\smallskip%
%
%
\leavevmode%
${L_{238.4}}$%
{} : {$[1\above{1pt}{1pt}{1}{}2\above{1pt}{1pt}{1}{}]\above{1pt}{1pt}{}{0}32\above{1pt}{1pt}{-}{3}{\cdot}1\above{1pt}{1pt}{1}{}3\above{1pt}{1pt}{-}{}9\above{1pt}{1pt}{1}{}$}\spacer%
\instructions{m}%
\EasyButWeakLineBreak%
{${96}\above{1pt}{1pt}{s}{2}{36}\above{1pt}{1pt}{*}{2}{24}\above{1pt}{1pt}{l}{2}{1}\above{1pt}{1pt}{}{2}{96}\above{1pt}{1pt}{}{2}{9}\above{1pt}{1pt}{r}{2}{24}\above{1pt}{1pt}{*}{2}{4}\above{1pt}{1pt}{s}{2}$}%
\nopagebreak\par%
\nopagebreak\par\leavevmode%
{$\left[\!\llap{\phantom{%
\begingroup \smaller\smaller\smaller\begin{tabular}{@{}c@{}}%
0\\0\\0
\end{tabular}\endgroup%
}}\right.$}%
\begingroup \smaller\smaller\smaller\begin{tabular}{@{}c@{}}%
37728\\0\\-288
\end{tabular}\endgroup%
\kern3pt%
\begingroup \smaller\smaller\smaller\begin{tabular}{@{}c@{}}%
0\\-30\\-6
\end{tabular}\endgroup%
\kern3pt%
\begingroup \smaller\smaller\smaller\begin{tabular}{@{}c@{}}%
-288\\-6\\1
\end{tabular}\endgroup%
{$\left.\llap{\phantom{%
\begingroup \smaller\smaller\smaller\begin{tabular}{@{}c@{}}%
0\\0\\0
\end{tabular}\endgroup%
}}\!\right]$}%
\EasyButWeakLineBreak%
{$\left[\!\llap{\phantom{%
\begingroup \smaller\smaller\smaller\begin{tabular}{@{}c@{}}%
0\\0\\0
\end{tabular}\endgroup%
}}\right.$}%
\begingroup \smaller\smaller\smaller\begin{tabular}{@{}c@{}}%
1\\-32\\144
\end{tabular}\endgroup%
\HardButStrongLineBreak\kern3pt%
\begingroup \smaller\smaller\smaller\begin{tabular}{@{}c@{}}%
-1\\24\\-126
\end{tabular}\endgroup%
\HardButStrongLineBreak\kern3pt%
\begingroup \smaller\smaller\smaller\begin{tabular}{@{}c@{}}%
-1\\26\\-132
\end{tabular}\endgroup%
\HardButStrongLineBreak\kern3pt%
\begingroup \smaller\smaller\smaller\begin{tabular}{@{}c@{}}%
0\\0\\-1
\end{tabular}\endgroup%
\HardButStrongLineBreak\kern3pt%
\begingroup \smaller\smaller\smaller\begin{tabular}{@{}c@{}}%
3\\-80\\384
\end{tabular}\endgroup%
\HardButStrongLineBreak\kern3pt%
\begingroup \smaller\smaller\smaller\begin{tabular}{@{}c@{}}%
2\\-54\\261
\end{tabular}\endgroup%
\HardButStrongLineBreak\kern3pt%
\begingroup \smaller\smaller\smaller\begin{tabular}{@{}c@{}}%
3\\-82\\396
\end{tabular}\endgroup%
\HardButStrongLineBreak\kern3pt%
\begingroup \smaller\smaller\smaller\begin{tabular}{@{}c@{}}%
1\\-28\\134
\end{tabular}\endgroup%
{$\left.\llap{\phantom{%
\begingroup \smaller\smaller\smaller\begin{tabular}{@{}c@{}}%
0\\0\\0
\end{tabular}\endgroup%
}}\!\right]$}%
%
%
\hbox{}\par\smallskip%
%
%
\leavevmode%
${L_{238.5}}$%
{} : {$1\above{1pt}{1pt}{-}{3}4\above{1pt}{1pt}{1}{1}32\above{1pt}{1pt}{1}{7}{\cdot}1\above{1pt}{1pt}{-}{}3\above{1pt}{1pt}{1}{}9\above{1pt}{1pt}{-}{}$}\spacer%
\instructions{m}%
\EasyButWeakLineBreak%
{${48}\above{1pt}{1pt}{*}{2}{32}\above{1pt}{1pt}{l}{2}{3}\above{1pt}{1pt}{r}{2}{288}\above{1pt}{1pt}{*}{2}{48}\above{1pt}{1pt}{s}{2}{32}\above{1pt}{1pt}{*}{2}{12}\above{1pt}{1pt}{*}{2}{288}\above{1pt}{1pt}{s}{2}$}%
\nopagebreak\par%
\nopagebreak\par\leavevmode%
{$\left[\!\llap{\phantom{%
\begingroup \smaller\smaller\smaller\begin{tabular}{@{}c@{}}%
0\\0\\0
\end{tabular}\endgroup%
}}\right.$}%
\begingroup \smaller\smaller\smaller\begin{tabular}{@{}c@{}}%
-11808\\-1728\\576
\end{tabular}\endgroup%
\kern3pt%
\begingroup \smaller\smaller\smaller\begin{tabular}{@{}c@{}}%
-1728\\-60\\12
\end{tabular}\endgroup%
\kern3pt%
\begingroup \smaller\smaller\smaller\begin{tabular}{@{}c@{}}%
576\\12\\-1
\end{tabular}\endgroup%
{$\left.\llap{\phantom{%
\begingroup \smaller\smaller\smaller\begin{tabular}{@{}c@{}}%
0\\0\\0
\end{tabular}\endgroup%
}}\!\right]$}%
\EasyButWeakLineBreak%
{$\left[\!\llap{\phantom{%
\begingroup \smaller\smaller\smaller\begin{tabular}{@{}c@{}}%
0\\0\\0
\end{tabular}\endgroup%
}}\right.$}%
\begingroup \smaller\smaller\smaller\begin{tabular}{@{}c@{}}%
1\\-62\\-168
\end{tabular}\endgroup%
\HardButStrongLineBreak\kern3pt%
\begingroup \smaller\smaller\smaller\begin{tabular}{@{}c@{}}%
1\\-60\\-160
\end{tabular}\endgroup%
\HardButStrongLineBreak\kern3pt%
\begingroup \smaller\smaller\smaller\begin{tabular}{@{}c@{}}%
0\\1\\3
\end{tabular}\endgroup%
\HardButStrongLineBreak\kern3pt%
\begingroup \smaller\smaller\smaller\begin{tabular}{@{}c@{}}%
-5\\324\\864
\end{tabular}\endgroup%
\HardButStrongLineBreak\kern3pt%
\begingroup \smaller\smaller\smaller\begin{tabular}{@{}c@{}}%
-3\\190\\504
\end{tabular}\endgroup%
\HardButStrongLineBreak\kern3pt%
\begingroup \smaller\smaller\smaller\begin{tabular}{@{}c@{}}%
-3\\188\\496
\end{tabular}\endgroup%
\HardButStrongLineBreak\kern3pt%
\begingroup \smaller\smaller\smaller\begin{tabular}{@{}c@{}}%
-1\\62\\162
\end{tabular}\endgroup%
\HardButStrongLineBreak\kern3pt%
\begingroup \smaller\smaller\smaller\begin{tabular}{@{}c@{}}%
-1\\60\\144
\end{tabular}\endgroup%
{$\left.\llap{\phantom{%
\begingroup \smaller\smaller\smaller\begin{tabular}{@{}c@{}}%
0\\0\\0
\end{tabular}\endgroup%
}}\!\right]$}%
%
%
\hbox{}\par\smallskip%
%
%
\leavevmode%
${L_{238.6}}$%
{} : {$1\above{1pt}{1pt}{-}{3}4\above{1pt}{1pt}{1}{7}32\above{1pt}{1pt}{1}{1}{\cdot}1\above{1pt}{1pt}{-}{}3\above{1pt}{1pt}{1}{}9\above{1pt}{1pt}{-}{}$}\EasyButWeakLineBreak%
{${12}\above{1pt}{1pt}{}{2}{32}\above{1pt}{1pt}{}{2}{3}\above{1pt}{1pt}{}{2}{288}\above{1pt}{1pt}{}{2}{12}\above{1pt}{1pt}{r}{2}{32}\above{1pt}{1pt}{s}{2}{12}\above{1pt}{1pt}{s}{2}{288}\above{1pt}{1pt}{l}{2}$}%
\nopagebreak\par%
\nopagebreak\par\leavevmode%
{$\left[\!\llap{\phantom{%
\begingroup \smaller\smaller\smaller\begin{tabular}{@{}c@{}}%
0\\0\\0
\end{tabular}\endgroup%
}}\right.$}%
\begingroup \smaller\smaller\smaller\begin{tabular}{@{}c@{}}%
90144\\0\\-1440
\end{tabular}\endgroup%
\kern3pt%
\begingroup \smaller\smaller\smaller\begin{tabular}{@{}c@{}}%
0\\12\\0
\end{tabular}\endgroup%
\kern3pt%
\begingroup \smaller\smaller\smaller\begin{tabular}{@{}c@{}}%
-1440\\0\\23
\end{tabular}\endgroup%
{$\left.\llap{\phantom{%
\begingroup \smaller\smaller\smaller\begin{tabular}{@{}c@{}}%
0\\0\\0
\end{tabular}\endgroup%
}}\!\right]$}%
\EasyButWeakLineBreak%
{$\left[\!\llap{\phantom{%
\begingroup \smaller\smaller\smaller\begin{tabular}{@{}c@{}}%
0\\0\\0
\end{tabular}\endgroup%
}}\right.$}%
\begingroup \smaller\smaller\smaller\begin{tabular}{@{}c@{}}%
0\\1\\0
\end{tabular}\endgroup%
\HardButStrongLineBreak\kern3pt%
\begingroup \smaller\smaller\smaller\begin{tabular}{@{}c@{}}%
1\\0\\64
\end{tabular}\endgroup%
\HardButStrongLineBreak\kern3pt%
\begingroup \smaller\smaller\smaller\begin{tabular}{@{}c@{}}%
1\\-1\\63
\end{tabular}\endgroup%
\HardButStrongLineBreak\kern3pt%
\begingroup \smaller\smaller\smaller\begin{tabular}{@{}c@{}}%
23\\-24\\1440
\end{tabular}\endgroup%
\HardButStrongLineBreak\kern3pt%
\begingroup \smaller\smaller\smaller\begin{tabular}{@{}c@{}}%
5\\-5\\312
\end{tabular}\endgroup%
\HardButStrongLineBreak\kern3pt%
\begingroup \smaller\smaller\smaller\begin{tabular}{@{}c@{}}%
9\\-8\\560
\end{tabular}\endgroup%
\HardButStrongLineBreak\kern3pt%
\begingroup \smaller\smaller\smaller\begin{tabular}{@{}c@{}}%
3\\-2\\186
\end{tabular}\endgroup%
\HardButStrongLineBreak\kern3pt%
\begingroup \smaller\smaller\smaller\begin{tabular}{@{}c@{}}%
7\\0\\432
\end{tabular}\endgroup%
{$\left.\llap{\phantom{%
\begingroup \smaller\smaller\smaller\begin{tabular}{@{}c@{}}%
0\\0\\0
\end{tabular}\endgroup%
}}\!\right]$}%

\medskip%
%
\leavevmode\llap{}%
$W_{239}$%
\qquad\llap{4} lattices, $\chi=36$%
\hfill%
$22\slashtwo2242|24\rtimes D_{2}$%
\nopagebreak\smallskip\hrule\nopagebreak\medskip%
%
%
\leavevmode%
${L_{239.1}}$%
{} : {$1\above{1pt}{1pt}{2}{2}32\above{1pt}{1pt}{1}{1}{\cdot}1\above{1pt}{1pt}{2}{}7\above{1pt}{1pt}{-}{}$}\EasyButWeakLineBreak%
{${2}\above{1pt}{1pt}{b}{2}{32}\above{1pt}{1pt}{*}{2}{4}\above{1pt}{1pt}{l}{2}{1}\above{1pt}{1pt}{}{2}{32}\above{1pt}{1pt}{r}{2}{2}\above{1pt}{1pt}{}{4}{1}\above{1pt}{1pt}{r}{2}{32}\above{1pt}{1pt}{s}{2}{4}\above{1pt}{1pt}{*}{4}$}%
\nopagebreak\par%
\nopagebreak\par\leavevmode%
{$\left[\!\llap{\phantom{%
\begingroup \smaller\smaller\smaller\begin{tabular}{@{}c@{}}%
0\\0\\0
\end{tabular}\endgroup%
}}\right.$}%
\begingroup \smaller\smaller\smaller\begin{tabular}{@{}c@{}}%
-231392\\56672\\3584
\end{tabular}\endgroup%
\kern3pt%
\begingroup \smaller\smaller\smaller\begin{tabular}{@{}c@{}}%
56672\\-13879\\-879
\end{tabular}\endgroup%
\kern3pt%
\begingroup \smaller\smaller\smaller\begin{tabular}{@{}c@{}}%
3584\\-879\\-54
\end{tabular}\endgroup%
{$\left.\llap{\phantom{%
\begingroup \smaller\smaller\smaller\begin{tabular}{@{}c@{}}%
0\\0\\0
\end{tabular}\endgroup%
}}\!\right]$}%
\EasyButWeakLineBreak%
{$\left[\!\llap{\phantom{%
\begingroup \smaller\smaller\smaller\begin{tabular}{@{}c@{}}%
0\\0\\0
\end{tabular}\endgroup%
}}\right.$}%
\begingroup \smaller\smaller\smaller\begin{tabular}{@{}c@{}}%
18\\70\\55
\end{tabular}\endgroup%
\HardButStrongLineBreak\kern3pt%
\begingroup \smaller\smaller\smaller\begin{tabular}{@{}c@{}}%
313\\1216\\976
\end{tabular}\endgroup%
\HardButStrongLineBreak\kern3pt%
\begingroup \smaller\smaller\smaller\begin{tabular}{@{}c@{}}%
155\\602\\486
\end{tabular}\endgroup%
\HardButStrongLineBreak\kern3pt%
\begingroup \smaller\smaller\smaller\begin{tabular}{@{}c@{}}%
94\\365\\296
\end{tabular}\endgroup%
\HardButStrongLineBreak\kern3pt%
\begingroup \smaller\smaller\smaller\begin{tabular}{@{}c@{}}%
577\\2240\\1824
\end{tabular}\endgroup%
\HardButStrongLineBreak\kern3pt%
\begingroup \smaller\smaller\smaller\begin{tabular}{@{}c@{}}%
84\\326\\267
\end{tabular}\endgroup%
\HardButStrongLineBreak\kern3pt%
\begingroup \smaller\smaller\smaller\begin{tabular}{@{}c@{}}%
41\\159\\132
\end{tabular}\endgroup%
\HardButStrongLineBreak\kern3pt%
\begingroup \smaller\smaller\smaller\begin{tabular}{@{}c@{}}%
-37\\-144\\-112
\end{tabular}\endgroup%
\HardButStrongLineBreak\kern3pt%
\begingroup \smaller\smaller\smaller\begin{tabular}{@{}c@{}}%
-17\\-66\\-54
\end{tabular}\endgroup%
{$\left.\llap{\phantom{%
\begingroup \smaller\smaller\smaller\begin{tabular}{@{}c@{}}%
0\\0\\0
\end{tabular}\endgroup%
}}\!\right]$}%
%
%
%
%
%
%
%
%
%
%
%
%
%
%

\medskip%
%
\leavevmode\llap{}%
$W_{240}$%
\qquad\llap{8} lattices, $\chi=18$%
\hfill%
$24|422|2\rtimes D_{2}$%
\nopagebreak\smallskip\hrule\nopagebreak\medskip%
%
%
\leavevmode%
${L_{240.1}}$%
{} : {$1\above{1pt}{1pt}{2}{2}16\above{1pt}{1pt}{1}{1}{\cdot}1\above{1pt}{1pt}{2}{}3\above{1pt}{1pt}{1}{}{\cdot}1\above{1pt}{1pt}{2}{}5\above{1pt}{1pt}{-}{}$}\EasyButWeakLineBreak%
{${16}\above{1pt}{1pt}{*}{2}{4}\above{1pt}{1pt}{*}{4}{2}\above{1pt}{1pt}{}{4}{1}\above{1pt}{1pt}{}{2}{16}\above{1pt}{1pt}{r}{2}{10}\above{1pt}{1pt}{b}{2}$}%
\nopagebreak\par%
\nopagebreak\par\leavevmode%
{$\left[\!\llap{\phantom{%
\begingroup \smaller\smaller\smaller\begin{tabular}{@{}c@{}}%
0\\0\\0
\end{tabular}\endgroup%
}}\right.$}%
\begingroup \smaller\smaller\smaller\begin{tabular}{@{}c@{}}%
-1422960\\4080\\4800
\end{tabular}\endgroup%
\kern3pt%
\begingroup \smaller\smaller\smaller\begin{tabular}{@{}c@{}}%
4080\\-11\\-15
\end{tabular}\endgroup%
\kern3pt%
\begingroup \smaller\smaller\smaller\begin{tabular}{@{}c@{}}%
4800\\-15\\-14
\end{tabular}\endgroup%
{$\left.\llap{\phantom{%
\begingroup \smaller\smaller\smaller\begin{tabular}{@{}c@{}}%
0\\0\\0
\end{tabular}\endgroup%
}}\!\right]$}%
\EasyButWeakLineBreak%
{$\left[\!\llap{\phantom{%
\begingroup \smaller\smaller\smaller\begin{tabular}{@{}c@{}}%
0\\0\\0
\end{tabular}\endgroup%
}}\right.$}%
\begingroup \smaller\smaller\smaller\begin{tabular}{@{}c@{}}%
-1\\-208\\-120
\end{tabular}\endgroup%
\HardButStrongLineBreak\kern3pt%
\begingroup \smaller\smaller\smaller\begin{tabular}{@{}c@{}}%
-1\\-210\\-118
\end{tabular}\endgroup%
\HardButStrongLineBreak\kern3pt%
\begingroup \smaller\smaller\smaller\begin{tabular}{@{}c@{}}%
1\\208\\119
\end{tabular}\endgroup%
\HardButStrongLineBreak\kern3pt%
\begingroup \smaller\smaller\smaller\begin{tabular}{@{}c@{}}%
3\\627\\355
\end{tabular}\endgroup%
\HardButStrongLineBreak\kern3pt%
\begingroup \smaller\smaller\smaller\begin{tabular}{@{}c@{}}%
13\\2720\\1536
\end{tabular}\endgroup%
\HardButStrongLineBreak\kern3pt%
\begingroup \smaller\smaller\smaller\begin{tabular}{@{}c@{}}%
2\\420\\235
\end{tabular}\endgroup%
{$\left.\llap{\phantom{%
\begingroup \smaller\smaller\smaller\begin{tabular}{@{}c@{}}%
0\\0\\0
\end{tabular}\endgroup%
}}\!\right]$}%

\medskip%
%
\leavevmode\llap{}%
$W_{241}$%
\qquad\llap{12} lattices, $\chi=24$%
\hfill%
$22222222\rtimes C_{2}$%
\nopagebreak\smallskip\hrule\nopagebreak\medskip%
%
%
\leavevmode%
${L_{241.1}}$%
{} : {$1\above{1pt}{1pt}{-2}{{\rm II}}16\above{1pt}{1pt}{-}{3}{\cdot}1\above{1pt}{1pt}{2}{}3\above{1pt}{1pt}{1}{}{\cdot}1\above{1pt}{1pt}{-}{}5\above{1pt}{1pt}{-}{}25\above{1pt}{1pt}{-}{}$}\spacer%
\instructions{5}%
\EasyButWeakLineBreak%
{${48}\above{1pt}{1pt}{r}{2}{10}\above{1pt}{1pt}{b}{2}{1200}\above{1pt}{1pt}{b}{2}{2}\above{1pt}{1pt}{l}{2}{1200}\above{1pt}{1pt}{r}{2}{10}\above{1pt}{1pt}{b}{2}{48}\above{1pt}{1pt}{b}{2}{50}\above{1pt}{1pt}{l}{2}$}%
\nopagebreak\par%
\nopagebreak\par\leavevmode%
{$\left[\!\llap{\phantom{%
\begingroup \smaller\smaller\smaller\begin{tabular}{@{}c@{}}%
0\\0\\0
\end{tabular}\endgroup%
}}\right.$}%
\begingroup \smaller\smaller\smaller\begin{tabular}{@{}c@{}}%
-6478800\\1299600\\25200
\end{tabular}\endgroup%
\kern3pt%
\begingroup \smaller\smaller\smaller\begin{tabular}{@{}c@{}}%
1299600\\-260690\\-5055
\end{tabular}\endgroup%
\kern3pt%
\begingroup \smaller\smaller\smaller\begin{tabular}{@{}c@{}}%
25200\\-5055\\-98
\end{tabular}\endgroup%
{$\left.\llap{\phantom{%
\begingroup \smaller\smaller\smaller\begin{tabular}{@{}c@{}}%
0\\0\\0
\end{tabular}\endgroup%
}}\!\right]$}%
\EasyButWeakLineBreak%
{$\left[\!\llap{\phantom{%
\begingroup \smaller\smaller\smaller\begin{tabular}{@{}c@{}}%
0\\0\\0
\end{tabular}\endgroup%
}}\right.$}%
\begingroup \smaller\smaller\smaller\begin{tabular}{@{}c@{}}%
-65\\-288\\-1872
\end{tabular}\endgroup%
\HardButStrongLineBreak\kern3pt%
\begingroup \smaller\smaller\smaller\begin{tabular}{@{}c@{}}%
-9\\-39\\-305
\end{tabular}\endgroup%
\HardButStrongLineBreak\kern3pt%
\begingroup \smaller\smaller\smaller\begin{tabular}{@{}c@{}}%
-31\\-120\\-1800
\end{tabular}\endgroup%
\HardButStrongLineBreak\kern3pt%
\begingroup \smaller\smaller\smaller\begin{tabular}{@{}c@{}}%
1\\5\\-1
\end{tabular}\endgroup%
\HardButStrongLineBreak\kern3pt%
\begingroup \smaller\smaller\smaller\begin{tabular}{@{}c@{}}%
101\\480\\1200
\end{tabular}\endgroup%
\HardButStrongLineBreak\kern3pt%
\begingroup \smaller\smaller\smaller\begin{tabular}{@{}c@{}}%
2\\9\\50
\end{tabular}\endgroup%
\HardButStrongLineBreak\kern3pt%
\begingroup \smaller\smaller\smaller\begin{tabular}{@{}c@{}}%
-5\\-24\\-48
\end{tabular}\endgroup%
\HardButStrongLineBreak\kern3pt%
\begingroup \smaller\smaller\smaller\begin{tabular}{@{}c@{}}%
-12\\-55\\-250
\end{tabular}\endgroup%
{$\left.\llap{\phantom{%
\begingroup \smaller\smaller\smaller\begin{tabular}{@{}c@{}}%
0\\0\\0
\end{tabular}\endgroup%
}}\!\right]$}%

\medskip%
%
\leavevmode\llap{}%
$W_{242}$%
\qquad\llap{8} lattices, $\chi=12$%
\hfill%
$222|222|\rtimes D_{2}$%
\nopagebreak\smallskip\hrule\nopagebreak\medskip%
%
%
\leavevmode%
${L_{242.1}}$%
{} : {$1\above{1pt}{1pt}{-2}{2}16\above{1pt}{1pt}{-}{5}{\cdot}1\above{1pt}{1pt}{2}{}3\above{1pt}{1pt}{1}{}{\cdot}1\above{1pt}{1pt}{-2}{}5\above{1pt}{1pt}{1}{}$}\EasyButWeakLineBreak%
{${16}\above{1pt}{1pt}{l}{2}{3}\above{1pt}{1pt}{}{2}{80}\above{1pt}{1pt}{r}{2}{2}\above{1pt}{1pt}{b}{2}{80}\above{1pt}{1pt}{*}{2}{12}\above{1pt}{1pt}{s}{2}$}%
\nopagebreak\par%
\nopagebreak\par\leavevmode%
{$\left[\!\llap{\phantom{%
\begingroup \smaller\smaller\smaller\begin{tabular}{@{}c@{}}%
0\\0\\0
\end{tabular}\endgroup%
}}\right.$}%
\begingroup \smaller\smaller\smaller\begin{tabular}{@{}c@{}}%
106320\\17520\\-480
\end{tabular}\endgroup%
\kern3pt%
\begingroup \smaller\smaller\smaller\begin{tabular}{@{}c@{}}%
17520\\2887\\-79
\end{tabular}\endgroup%
\kern3pt%
\begingroup \smaller\smaller\smaller\begin{tabular}{@{}c@{}}%
-480\\-79\\2
\end{tabular}\endgroup%
{$\left.\llap{\phantom{%
\begingroup \smaller\smaller\smaller\begin{tabular}{@{}c@{}}%
0\\0\\0
\end{tabular}\endgroup%
}}\!\right]$}%
\EasyButWeakLineBreak%
{$\left[\!\llap{\phantom{%
\begingroup \smaller\smaller\smaller\begin{tabular}{@{}c@{}}%
0\\0\\0
\end{tabular}\endgroup%
}}\right.$}%
\begingroup \smaller\smaller\smaller\begin{tabular}{@{}c@{}}%
9\\-56\\-48
\end{tabular}\endgroup%
\HardButStrongLineBreak\kern3pt%
\begingroup \smaller\smaller\smaller\begin{tabular}{@{}c@{}}%
14\\-87\\-75
\end{tabular}\endgroup%
\HardButStrongLineBreak\kern3pt%
\begingroup \smaller\smaller\smaller\begin{tabular}{@{}c@{}}%
103\\-640\\-560
\end{tabular}\endgroup%
\HardButStrongLineBreak\kern3pt%
\begingroup \smaller\smaller\smaller\begin{tabular}{@{}c@{}}%
0\\0\\-1
\end{tabular}\endgroup%
\HardButStrongLineBreak\kern3pt%
\begingroup \smaller\smaller\smaller\begin{tabular}{@{}c@{}}%
-13\\80\\40
\end{tabular}\endgroup%
\HardButStrongLineBreak\kern3pt%
\begingroup \smaller\smaller\smaller\begin{tabular}{@{}c@{}}%
-1\\6\\0
\end{tabular}\endgroup%
{$\left.\llap{\phantom{%
\begingroup \smaller\smaller\smaller\begin{tabular}{@{}c@{}}%
0\\0\\0
\end{tabular}\endgroup%
}}\!\right]$}%

\medskip%
%
\leavevmode\llap{}%
$W_{243}$%
\qquad\llap{8} lattices, $\chi=12$%
\hfill%
$22|222|2\rtimes D_{2}$%
\nopagebreak\smallskip\hrule\nopagebreak\medskip%
%
%
\leavevmode%
${L_{243.1}}$%
{} : {$1\above{1pt}{1pt}{2}{6}16\above{1pt}{1pt}{1}{1}{\cdot}1\above{1pt}{1pt}{2}{}3\above{1pt}{1pt}{1}{}{\cdot}1\above{1pt}{1pt}{-2}{}5\above{1pt}{1pt}{1}{}$}\EasyButWeakLineBreak%
{${16}\above{1pt}{1pt}{}{2}{3}\above{1pt}{1pt}{r}{2}{80}\above{1pt}{1pt}{s}{2}{12}\above{1pt}{1pt}{*}{2}{16}\above{1pt}{1pt}{b}{2}{30}\above{1pt}{1pt}{l}{2}$}%
\nopagebreak\par%
\nopagebreak\par\leavevmode%
{$\left[\!\llap{\phantom{%
\begingroup \smaller\smaller\smaller\begin{tabular}{@{}c@{}}%
0\\0\\0
\end{tabular}\endgroup%
}}\right.$}%
\begingroup \smaller\smaller\smaller\begin{tabular}{@{}c@{}}%
1595280\\-21360\\3360
\end{tabular}\endgroup%
\kern3pt%
\begingroup \smaller\smaller\smaller\begin{tabular}{@{}c@{}}%
-21360\\286\\-45
\end{tabular}\endgroup%
\kern3pt%
\begingroup \smaller\smaller\smaller\begin{tabular}{@{}c@{}}%
3360\\-45\\7
\end{tabular}\endgroup%
{$\left.\llap{\phantom{%
\begingroup \smaller\smaller\smaller\begin{tabular}{@{}c@{}}%
0\\0\\0
\end{tabular}\endgroup%
}}\!\right]$}%
\EasyButWeakLineBreak%
{$\left[\!\llap{\phantom{%
\begingroup \smaller\smaller\smaller\begin{tabular}{@{}c@{}}%
0\\0\\0
\end{tabular}\endgroup%
}}\right.$}%
\begingroup \smaller\smaller\smaller\begin{tabular}{@{}c@{}}%
-3\\-224\\0
\end{tabular}\endgroup%
\HardButStrongLineBreak\kern3pt%
\begingroup \smaller\smaller\smaller\begin{tabular}{@{}c@{}}%
-1\\-75\\-3
\end{tabular}\endgroup%
\HardButStrongLineBreak\kern3pt%
\begingroup \smaller\smaller\smaller\begin{tabular}{@{}c@{}}%
-1\\-80\\-40
\end{tabular}\endgroup%
\HardButStrongLineBreak\kern3pt%
\begingroup \smaller\smaller\smaller\begin{tabular}{@{}c@{}}%
1\\72\\-18
\end{tabular}\endgroup%
\HardButStrongLineBreak\kern3pt%
\begingroup \smaller\smaller\smaller\begin{tabular}{@{}c@{}}%
1\\72\\-16
\end{tabular}\endgroup%
\HardButStrongLineBreak\kern3pt%
\begingroup \smaller\smaller\smaller\begin{tabular}{@{}c@{}}%
-1\\-75\\0
\end{tabular}\endgroup%
{$\left.\llap{\phantom{%
\begingroup \smaller\smaller\smaller\begin{tabular}{@{}c@{}}%
0\\0\\0
\end{tabular}\endgroup%
}}\!\right]$}%

\medskip%
%
\leavevmode\llap{}%
$W_{244}$%
\qquad\llap{16} lattices, $\chi=24$%
\hfill%
$2222|2222|\rtimes D_{2}$%
\nopagebreak\smallskip\hrule\nopagebreak\medskip%
%
%
\leavevmode%
${L_{244.1}}$%
{} : {$1\above{1pt}{1pt}{2}{2}16\above{1pt}{1pt}{1}{1}{\cdot}1\above{1pt}{1pt}{-}{}3\above{1pt}{1pt}{-}{}9\above{1pt}{1pt}{1}{}{\cdot}1\above{1pt}{1pt}{-2}{}5\above{1pt}{1pt}{1}{}$}\spacer%
\instructions{3}%
\EasyButWeakLineBreak%
{${80}\above{1pt}{1pt}{l}{2}{9}\above{1pt}{1pt}{r}{2}{20}\above{1pt}{1pt}{*}{2}{144}\above{1pt}{1pt}{b}{2}{2}\above{1pt}{1pt}{l}{2}{144}\above{1pt}{1pt}{}{2}{5}\above{1pt}{1pt}{r}{2}{36}\above{1pt}{1pt}{s}{2}$}%
\nopagebreak\par%
\nopagebreak\par\leavevmode%
{$\left[\!\llap{\phantom{%
\begingroup \smaller\smaller\smaller\begin{tabular}{@{}c@{}}%
0\\0\\0
\end{tabular}\endgroup%
}}\right.$}%
\begingroup \smaller\smaller\smaller\begin{tabular}{@{}c@{}}%
-1568880\\7920\\4320
\end{tabular}\endgroup%
\kern3pt%
\begingroup \smaller\smaller\smaller\begin{tabular}{@{}c@{}}%
7920\\-39\\-24
\end{tabular}\endgroup%
\kern3pt%
\begingroup \smaller\smaller\smaller\begin{tabular}{@{}c@{}}%
4320\\-24\\-7
\end{tabular}\endgroup%
{$\left.\llap{\phantom{%
\begingroup \smaller\smaller\smaller\begin{tabular}{@{}c@{}}%
0\\0\\0
\end{tabular}\endgroup%
}}\!\right]$}%
\EasyButWeakLineBreak%
{$\left[\!\llap{\phantom{%
\begingroup \smaller\smaller\smaller\begin{tabular}{@{}c@{}}%
0\\0\\0
\end{tabular}\endgroup%
}}\right.$}%
\begingroup \smaller\smaller\smaller\begin{tabular}{@{}c@{}}%
3\\480\\200
\end{tabular}\endgroup%
\HardButStrongLineBreak\kern3pt%
\begingroup \smaller\smaller\smaller\begin{tabular}{@{}c@{}}%
-1\\-159\\-72
\end{tabular}\endgroup%
\HardButStrongLineBreak\kern3pt%
\begingroup \smaller\smaller\smaller\begin{tabular}{@{}c@{}}%
-1\\-160\\-70
\end{tabular}\endgroup%
\HardButStrongLineBreak\kern3pt%
\begingroup \smaller\smaller\smaller\begin{tabular}{@{}c@{}}%
5\\792\\360
\end{tabular}\endgroup%
\HardButStrongLineBreak\kern3pt%
\begingroup \smaller\smaller\smaller\begin{tabular}{@{}c@{}}%
1\\159\\71
\end{tabular}\endgroup%
\HardButStrongLineBreak\kern3pt%
\begingroup \smaller\smaller\smaller\begin{tabular}{@{}c@{}}%
41\\6528\\2880
\end{tabular}\endgroup%
\HardButStrongLineBreak\kern3pt%
\begingroup \smaller\smaller\smaller\begin{tabular}{@{}c@{}}%
7\\1115\\490
\end{tabular}\endgroup%
\HardButStrongLineBreak\kern3pt%
\begingroup \smaller\smaller\smaller\begin{tabular}{@{}c@{}}%
7\\1116\\486
\end{tabular}\endgroup%
{$\left.\llap{\phantom{%
\begingroup \smaller\smaller\smaller\begin{tabular}{@{}c@{}}%
0\\0\\0
\end{tabular}\endgroup%
}}\!\right]$}%

\medskip%
%
\leavevmode\llap{}%
$W_{245}$%
\qquad\llap{16} lattices, $\chi=24$%
\hfill%
$222|2222|2\rtimes D_{2}$%
\nopagebreak\smallskip\hrule\nopagebreak\medskip%
%
%
\leavevmode%
${L_{245.1}}$%
{} : {$1\above{1pt}{1pt}{-2}{6}16\above{1pt}{1pt}{-}{5}{\cdot}1\above{1pt}{1pt}{-}{}3\above{1pt}{1pt}{-}{}9\above{1pt}{1pt}{1}{}{\cdot}1\above{1pt}{1pt}{-2}{}5\above{1pt}{1pt}{1}{}$}\spacer%
\instructions{3}%
\EasyButWeakLineBreak%
{${80}\above{1pt}{1pt}{}{2}{9}\above{1pt}{1pt}{}{2}{5}\above{1pt}{1pt}{r}{2}{144}\above{1pt}{1pt}{s}{2}{20}\above{1pt}{1pt}{*}{2}{36}\above{1pt}{1pt}{*}{2}{80}\above{1pt}{1pt}{b}{2}{6}\above{1pt}{1pt}{l}{2}$}%
\nopagebreak\par%
\nopagebreak\par\leavevmode%
{$\left[\!\llap{\phantom{%
\begingroup \smaller\smaller\smaller\begin{tabular}{@{}c@{}}%
0\\0\\0
\end{tabular}\endgroup%
}}\right.$}%
\begingroup \smaller\smaller\smaller\begin{tabular}{@{}c@{}}%
-488880\\1440\\4320
\end{tabular}\endgroup%
\kern3pt%
\begingroup \smaller\smaller\smaller\begin{tabular}{@{}c@{}}%
1440\\-3\\-15
\end{tabular}\endgroup%
\kern3pt%
\begingroup \smaller\smaller\smaller\begin{tabular}{@{}c@{}}%
4320\\-15\\-34
\end{tabular}\endgroup%
{$\left.\llap{\phantom{%
\begingroup \smaller\smaller\smaller\begin{tabular}{@{}c@{}}%
0\\0\\0
\end{tabular}\endgroup%
}}\!\right]$}%
\EasyButWeakLineBreak%
{$\left[\!\llap{\phantom{%
\begingroup \smaller\smaller\smaller\begin{tabular}{@{}c@{}}%
0\\0\\0
\end{tabular}\endgroup%
}}\right.$}%
\begingroup \smaller\smaller\smaller\begin{tabular}{@{}c@{}}%
23\\2960\\1600
\end{tabular}\endgroup%
\HardButStrongLineBreak\kern3pt%
\begingroup \smaller\smaller\smaller\begin{tabular}{@{}c@{}}%
4\\513\\279
\end{tabular}\endgroup%
\HardButStrongLineBreak\kern3pt%
\begingroup \smaller\smaller\smaller\begin{tabular}{@{}c@{}}%
2\\255\\140
\end{tabular}\endgroup%
\HardButStrongLineBreak\kern3pt%
\begingroup \smaller\smaller\smaller\begin{tabular}{@{}c@{}}%
1\\120\\72
\end{tabular}\endgroup%
\HardButStrongLineBreak\kern3pt%
\begingroup \smaller\smaller\smaller\begin{tabular}{@{}c@{}}%
-1\\-130\\-70
\end{tabular}\endgroup%
\HardButStrongLineBreak\kern3pt%
\begingroup \smaller\smaller\smaller\begin{tabular}{@{}c@{}}%
-1\\-126\\-72
\end{tabular}\endgroup%
\HardButStrongLineBreak\kern3pt%
\begingroup \smaller\smaller\smaller\begin{tabular}{@{}c@{}}%
3\\400\\200
\end{tabular}\endgroup%
\HardButStrongLineBreak\kern3pt%
\begingroup \smaller\smaller\smaller\begin{tabular}{@{}c@{}}%
1\\130\\69
\end{tabular}\endgroup%
{$\left.\llap{\phantom{%
\begingroup \smaller\smaller\smaller\begin{tabular}{@{}c@{}}%
0\\0\\0
\end{tabular}\endgroup%
}}\!\right]$}%

\medskip%
%
\leavevmode\llap{}%
$W_{246}$%
\qquad\llap{6} lattices, $\chi=36$%
\hfill%
$26322632\rtimes C_{2}$%
\nopagebreak\smallskip\hrule\nopagebreak\medskip%
%
%
\leavevmode%
${L_{246.1}}$%
{} : {$1\above{1pt}{1pt}{-2}{{\rm II}}4\above{1pt}{1pt}{1}{7}{\cdot}1\above{1pt}{1pt}{-}{}3\above{1pt}{1pt}{-}{}81\above{1pt}{1pt}{-}{}$}\spacer%
\instructions{2}%
\EasyButWeakLineBreak%
{${162}\above{1pt}{1pt}{b}{2}{2}\above{1pt}{1pt}{}{6}{6}\above{1pt}{1pt}{-}{3}{6}\above{1pt}{1pt}{s}{2}$}\relax$\,(\times2)$%
\nopagebreak\par%
\nopagebreak\par\leavevmode%
{$\left[\!\llap{\phantom{%
\begingroup \smaller\smaller\smaller\begin{tabular}{@{}c@{}}%
0\\0\\0
\end{tabular}\endgroup%
}}\right.$}%
\begingroup \smaller\smaller\smaller\begin{tabular}{@{}c@{}}%
-1920996\\17172\\602316
\end{tabular}\endgroup%
\kern3pt%
\begingroup \smaller\smaller\smaller\begin{tabular}{@{}c@{}}%
17172\\-138\\-5427
\end{tabular}\endgroup%
\kern3pt%
\begingroup \smaller\smaller\smaller\begin{tabular}{@{}c@{}}%
602316\\-5427\\-188734
\end{tabular}\endgroup%
{$\left.\llap{\phantom{%
\begingroup \smaller\smaller\smaller\begin{tabular}{@{}c@{}}%
0\\0\\0
\end{tabular}\endgroup%
}}\!\right]$}%
\hfil\penalty500%
{$\left[\!\llap{\phantom{%
\begingroup \smaller\smaller\smaller\begin{tabular}{@{}c@{}}%
0\\0\\0
\end{tabular}\endgroup%
}}\right.$}%
\begingroup \smaller\smaller\smaller\begin{tabular}{@{}c@{}}%
2411801\\19702224\\7130268
\end{tabular}\endgroup%
\kern3pt%
\begingroup \smaller\smaller\smaller\begin{tabular}{@{}c@{}}%
-23467\\-191705\\-69378
\end{tabular}\endgroup%
\kern3pt%
\begingroup \smaller\smaller\smaller\begin{tabular}{@{}c@{}}%
-750944\\-6134528\\-2220097
\end{tabular}\endgroup%
{$\left.\llap{\phantom{%
\begingroup \smaller\smaller\smaller\begin{tabular}{@{}c@{}}%
0\\0\\0
\end{tabular}\endgroup%
}}\!\right]$}%
\EasyButWeakLineBreak%
{$\left[\!\llap{\phantom{%
\begingroup \smaller\smaller\smaller\begin{tabular}{@{}c@{}}%
0\\0\\0
\end{tabular}\endgroup%
}}\right.$}%
\begingroup \smaller\smaller\smaller\begin{tabular}{@{}c@{}}%
-274\\-2241\\-810
\end{tabular}\endgroup%
\HardButStrongLineBreak\kern3pt%
\begingroup \smaller\smaller\smaller\begin{tabular}{@{}c@{}}%
115\\939\\340
\end{tabular}\endgroup%
\HardButStrongLineBreak\kern3pt%
\begingroup \smaller\smaller\smaller\begin{tabular}{@{}c@{}}%
1377\\11248\\4071
\end{tabular}\endgroup%
\HardButStrongLineBreak\kern3pt%
\begingroup \smaller\smaller\smaller\begin{tabular}{@{}c@{}}%
3923\\32047\\11598
\end{tabular}\endgroup%
{$\left.\llap{\phantom{%
\begingroup \smaller\smaller\smaller\begin{tabular}{@{}c@{}}%
0\\0\\0
\end{tabular}\endgroup%
}}\!\right]$}%

\medskip%
%
\leavevmode\llap{}%
$W_{247}$%
\qquad\llap{4} lattices, $\chi=24$%
\hfill%
$\infty|\infty|\infty|\infty|\rtimes D_{4}$%
\nopagebreak\smallskip\hrule\nopagebreak\medskip%
%
%
\leavevmode%
${L_{247.1}}$%
{} : {$1\above{1pt}{1pt}{-2}{{\rm II}}8\above{1pt}{1pt}{1}{1}{\cdot}1\above{1pt}{1pt}{1}{}5\above{1pt}{1pt}{-}{}25\above{1pt}{1pt}{1}{}$}\spacer%
\instructions{2}%
\EasyButWeakLineBreak%
{${40}\above{1pt}{1pt}{10,3}{\infty z}{10}\above{1pt}{1pt}{20,13}{\infty b}$}\relax$\,(\times2)$%
\nopagebreak\par%
\nopagebreak\par\leavevmode%
{$\left[\!\llap{\phantom{%
\begingroup \smaller\smaller\smaller\begin{tabular}{@{}c@{}}%
0\\0\\0
\end{tabular}\endgroup%
}}\right.$}%
\begingroup \smaller\smaller\smaller\begin{tabular}{@{}c@{}}%
-124600\\2400\\1800
\end{tabular}\endgroup%
\kern3pt%
\begingroup \smaller\smaller\smaller\begin{tabular}{@{}c@{}}%
2400\\-10\\-35
\end{tabular}\endgroup%
\kern3pt%
\begingroup \smaller\smaller\smaller\begin{tabular}{@{}c@{}}%
1800\\-35\\-26
\end{tabular}\endgroup%
{$\left.\llap{\phantom{%
\begingroup \smaller\smaller\smaller\begin{tabular}{@{}c@{}}%
0\\0\\0
\end{tabular}\endgroup%
}}\!\right]$}%
\hfil\penalty500%
{$\left[\!\llap{\phantom{%
\begingroup \smaller\smaller\smaller\begin{tabular}{@{}c@{}}%
0\\0\\0
\end{tabular}\endgroup%
}}\right.$}%
\begingroup \smaller\smaller\smaller\begin{tabular}{@{}c@{}}%
-1\\-480\\0
\end{tabular}\endgroup%
\kern3pt%
\begingroup \smaller\smaller\smaller\begin{tabular}{@{}c@{}}%
0\\1\\0
\end{tabular}\endgroup%
\kern3pt%
\begingroup \smaller\smaller\smaller\begin{tabular}{@{}c@{}}%
0\\7\\-1
\end{tabular}\endgroup%
{$\left.\llap{\phantom{%
\begingroup \smaller\smaller\smaller\begin{tabular}{@{}c@{}}%
0\\0\\0
\end{tabular}\endgroup%
}}\!\right]$}%
\EasyButWeakLineBreak%
{$\left[\!\llap{\phantom{%
\begingroup \smaller\smaller\smaller\begin{tabular}{@{}c@{}}%
0\\0\\0
\end{tabular}\endgroup%
}}\right.$}%
\begingroup \smaller\smaller\smaller\begin{tabular}{@{}c@{}}%
7\\-4\\480
\end{tabular}\endgroup%
\HardButStrongLineBreak\kern3pt%
\begingroup \smaller\smaller\smaller\begin{tabular}{@{}c@{}}%
3\\2\\205
\end{tabular}\endgroup%
{$\left.\llap{\phantom{%
\begingroup \smaller\smaller\smaller\begin{tabular}{@{}c@{}}%
0\\0\\0
\end{tabular}\endgroup%
}}\!\right]$}%
%
%
%
%
%
%
%
%
%
%
%
%
%
%

\medskip%
%
\leavevmode\llap{}%
$W_{248}$%
\qquad\llap{10} lattices, $\chi=48$%
\hfill%
$2\slashinfty22|22\slashinfty22|2\rtimes D_{4}$%
\nopagebreak\smallskip\hrule\nopagebreak\medskip%
%
%
\leavevmode%
${L_{248.1}}$%
{} : {$[1\above{1pt}{1pt}{1}{}2\above{1pt}{1pt}{1}{}]\above{1pt}{1pt}{}{0}128\above{1pt}{1pt}{1}{1}$}\EasyButWeakLineBreak%
{${128}\above{1pt}{1pt}{l}{2}{1}\above{1pt}{1pt}{16,9}{\infty}{4}\above{1pt}{1pt}{*}{2}{128}\above{1pt}{1pt}{s}{2}{8}\above{1pt}{1pt}{*}{2}$}\relax$\,(\times2)$%
\nopagebreak\par%
shares genus with {$ \main({L_{248.2}})$}%
\nopagebreak\par%
\nopagebreak\par\leavevmode%
{$\left[\!\llap{\phantom{%
\begingroup \smaller\smaller\smaller
\endgroup%
}}\!\right]$}%
%
%
\hbox{}\par\smallskip%
%
%
\leavevmode%
${L_{248.2}}$%
{} : {$[1\above{1pt}{1pt}{1}{}2\above{1pt}{1pt}{1}{}]\above{1pt}{1pt}{}{0}256\above{1pt}{1pt}{1}{1}$}\spacer%
\instructions{m*}%
\EasyButWeakLineBreak%
{${256}\above{1pt}{1pt}{s}{2}{8}\above{1pt}{1pt}{32,1}{\infty z}{2}\above{1pt}{1pt}{}{2}{256}\above{1pt}{1pt}{}{2}{1}\above{1pt}{1pt}{r}{2}{256}\above{1pt}{1pt}{*}{2}{8}\above{1pt}{1pt}{32,17}{\infty z}{2}\above{1pt}{1pt}{r}{2}{256}\above{1pt}{1pt}{s}{2}{4}\above{1pt}{1pt}{*}{2}$}%
\nopagebreak\par%
\nopagebreak\par\leavevmode%
{$\left[\!\llap{\phantom{%
\begingroup \smaller\smaller\smaller
\endgroup%
}}\!\right]$}%
%
%
\hbox{}\par\smallskip%
%
%
\leavevmode%
${L_{248.3}}$%
{} : {$1\above{1pt}{1pt}{1}{1}8\above{1pt}{1pt}{1}{7}128\above{1pt}{1pt}{1}{1}$}\EasyButWeakLineBreak%
{${128}\above{1pt}{1pt}{}{2}{1}\above{1pt}{1pt}{8,1}{\infty}{4}\above{1pt}{1pt}{s}{2}{128}\above{1pt}{1pt}{b}{2}{8}\above{1pt}{1pt}{l}{2}$}\relax$\,(\times2)$%
\nopagebreak\par%
\nopagebreak\par\leavevmode%
{$\left[\!\llap{\phantom{%
\begingroup \smaller\smaller\smaller
\endgroup%
}}\!\right]$}%
%
%
\hbox{}\par\smallskip%
%
%
\leavevmode%
${L_{248.4}}$%
{} : {$1\above{1pt}{1pt}{1}{7}8\above{1pt}{1pt}{1}{1}128\above{1pt}{1pt}{1}{1}$}\EasyButWeakLineBreak%
{${128}\above{1pt}{1pt}{r}{2}{4}\above{1pt}{1pt}{8,1}{\infty a}{4}\above{1pt}{1pt}{b}{2}{128}\above{1pt}{1pt}{l}{2}{8}\above{1pt}{1pt}{}{2}$}\relax$\,(\times2)$%
\nopagebreak\par%
\nopagebreak\par\leavevmode%
{$\left[\!\llap{\phantom{%
\begingroup \smaller\smaller\smaller
\endgroup%
}}\!\right]$}%

\medskip%
%
\leavevmode\llap{}%
$W_{249}$%
\qquad\llap{60} lattices, $\chi=30$%
\hfill%
$222|2222\slashtwo2\rtimes D_{2}$%
\nopagebreak\smallskip\hrule\nopagebreak\medskip%
%
%
\leavevmode%
${L_{249.1}}$%
{} : {$1\above{1pt}{1pt}{2}{0}8\above{1pt}{1pt}{1}{1}{\cdot}1\above{1pt}{1pt}{2}{}3\above{1pt}{1pt}{-}{}{\cdot}1\above{1pt}{1pt}{2}{}11\above{1pt}{1pt}{1}{}$}\EasyButWeakLineBreak%
{${132}\above{1pt}{1pt}{*}{2}{8}\above{1pt}{1pt}{*}{2}{44}\above{1pt}{1pt}{s}{2}{24}\above{1pt}{1pt}{l}{2}{11}\above{1pt}{1pt}{}{2}{8}\above{1pt}{1pt}{}{2}{33}\above{1pt}{1pt}{r}{2}{4}\above{1pt}{1pt}{l}{2}{1}\above{1pt}{1pt}{r}{2}$}%
\nopagebreak\par%
\nopagebreak\par\leavevmode%
{$\left[\!\llap{\phantom{%
\begingroup \smaller\smaller\smaller\begin{tabular}{@{}c@{}}%
0\\0\\0
\end{tabular}\endgroup%
}}\right.$}%
\begingroup \smaller\smaller\smaller\begin{tabular}{@{}c@{}}%
-217272\\792\\528
\end{tabular}\endgroup%
\kern3pt%
\begingroup \smaller\smaller\smaller\begin{tabular}{@{}c@{}}%
792\\-1\\-4
\end{tabular}\endgroup%
\kern3pt%
\begingroup \smaller\smaller\smaller\begin{tabular}{@{}c@{}}%
528\\-4\\1
\end{tabular}\endgroup%
{$\left.\llap{\phantom{%
\begingroup \smaller\smaller\smaller\begin{tabular}{@{}c@{}}%
0\\0\\0
\end{tabular}\endgroup%
}}\!\right]$}%
\EasyButWeakLineBreak%
{$\left[\!\llap{\phantom{%
\begingroup \smaller\smaller\smaller\begin{tabular}{@{}c@{}}%
0\\0\\0
\end{tabular}\endgroup%
}}\right.$}%
\begingroup \smaller\smaller\smaller\begin{tabular}{@{}c@{}}%
-5\\-858\\-792
\end{tabular}\endgroup%
\HardButStrongLineBreak\kern3pt%
\begingroup \smaller\smaller\smaller\begin{tabular}{@{}c@{}}%
-1\\-172\\-156
\end{tabular}\endgroup%
\HardButStrongLineBreak\kern3pt%
\begingroup \smaller\smaller\smaller\begin{tabular}{@{}c@{}}%
-1\\-176\\-154
\end{tabular}\endgroup%
\HardButStrongLineBreak\kern3pt%
\begingroup \smaller\smaller\smaller\begin{tabular}{@{}c@{}}%
1\\168\\156
\end{tabular}\endgroup%
\HardButStrongLineBreak\kern3pt%
\begingroup \smaller\smaller\smaller\begin{tabular}{@{}c@{}}%
5\\847\\770
\end{tabular}\endgroup%
\HardButStrongLineBreak\kern3pt%
\begingroup \smaller\smaller\smaller\begin{tabular}{@{}c@{}}%
5\\848\\768
\end{tabular}\endgroup%
\HardButStrongLineBreak\kern3pt%
\begingroup \smaller\smaller\smaller\begin{tabular}{@{}c@{}}%
14\\2376\\2145
\end{tabular}\endgroup%
\HardButStrongLineBreak\kern3pt%
\begingroup \smaller\smaller\smaller\begin{tabular}{@{}c@{}}%
1\\170\\152
\end{tabular}\endgroup%
\HardButStrongLineBreak\kern3pt%
\begingroup \smaller\smaller\smaller\begin{tabular}{@{}c@{}}%
0\\0\\-1
\end{tabular}\endgroup%
{$\left.\llap{\phantom{%
\begingroup \smaller\smaller\smaller\begin{tabular}{@{}c@{}}%
0\\0\\0
\end{tabular}\endgroup%
}}\!\right]$}%
%
%
\hbox{}\par\smallskip%
%
%
\leavevmode%
${L_{249.2}}$%
{} : {$[1\above{1pt}{1pt}{1}{}2\above{1pt}{1pt}{1}{}]\above{1pt}{1pt}{}{2}16\above{1pt}{1pt}{1}{7}{\cdot}1\above{1pt}{1pt}{2}{}3\above{1pt}{1pt}{-}{}{\cdot}1\above{1pt}{1pt}{2}{}11\above{1pt}{1pt}{1}{}$}\spacer%
\instructions{2}%
\EasyButWeakLineBreak%
{${33}\above{1pt}{1pt}{}{2}{2}\above{1pt}{1pt}{r}{2}{176}\above{1pt}{1pt}{l}{2}{6}\above{1pt}{1pt}{}{2}{11}\above{1pt}{1pt}{r}{2}{8}\above{1pt}{1pt}{*}{2}{528}\above{1pt}{1pt}{l}{2}{1}\above{1pt}{1pt}{r}{2}{16}\above{1pt}{1pt}{l}{2}$}%
\nopagebreak\par%
\nopagebreak\par\leavevmode%
{$\left[\!\llap{\phantom{%
\begingroup \smaller\smaller\smaller\begin{tabular}{@{}c@{}}%
0\\0\\0
\end{tabular}\endgroup%
}}\right.$}%
\begingroup \smaller\smaller\smaller\begin{tabular}{@{}c@{}}%
-1833744\\4752\\8448
\end{tabular}\endgroup%
\kern3pt%
\begingroup \smaller\smaller\smaller\begin{tabular}{@{}c@{}}%
4752\\-10\\-24
\end{tabular}\endgroup%
\kern3pt%
\begingroup \smaller\smaller\smaller\begin{tabular}{@{}c@{}}%
8448\\-24\\-37
\end{tabular}\endgroup%
{$\left.\llap{\phantom{%
\begingroup \smaller\smaller\smaller\begin{tabular}{@{}c@{}}%
0\\0\\0
\end{tabular}\endgroup%
}}\!\right]$}%
\EasyButWeakLineBreak%
{$\left[\!\llap{\phantom{%
\begingroup \smaller\smaller\smaller\begin{tabular}{@{}c@{}}%
0\\0\\0
\end{tabular}\endgroup%
}}\right.$}%
\begingroup \smaller\smaller\smaller\begin{tabular}{@{}c@{}}%
61\\7953\\8745
\end{tabular}\endgroup%
\HardButStrongLineBreak\kern3pt%
\begingroup \smaller\smaller\smaller\begin{tabular}{@{}c@{}}%
10\\1303\\1434
\end{tabular}\endgroup%
\HardButStrongLineBreak\kern3pt%
\begingroup \smaller\smaller\smaller\begin{tabular}{@{}c@{}}%
73\\9504\\10472
\end{tabular}\endgroup%
\HardButStrongLineBreak\kern3pt%
\begingroup \smaller\smaller\smaller\begin{tabular}{@{}c@{}}%
1\\129\\144
\end{tabular}\endgroup%
\HardButStrongLineBreak\kern3pt%
\begingroup \smaller\smaller\smaller\begin{tabular}{@{}c@{}}%
-1\\-132\\-143
\end{tabular}\endgroup%
\HardButStrongLineBreak\kern3pt%
\begingroup \smaller\smaller\smaller\begin{tabular}{@{}c@{}}%
-1\\-130\\-144
\end{tabular}\endgroup%
\HardButStrongLineBreak\kern3pt%
\begingroup \smaller\smaller\smaller\begin{tabular}{@{}c@{}}%
13\\1716\\1848
\end{tabular}\endgroup%
\HardButStrongLineBreak\kern3pt%
\begingroup \smaller\smaller\smaller\begin{tabular}{@{}c@{}}%
1\\131\\143
\end{tabular}\endgroup%
\HardButStrongLineBreak\kern3pt%
\begingroup \smaller\smaller\smaller\begin{tabular}{@{}c@{}}%
11\\1436\\1576
\end{tabular}\endgroup%
{$\left.\llap{\phantom{%
\begingroup \smaller\smaller\smaller\begin{tabular}{@{}c@{}}%
0\\0\\0
\end{tabular}\endgroup%
}}\!\right]$}%
%
%
\hbox{}\par\smallskip%
%
%
\leavevmode%
${L_{249.3}}$%
{} : {$[1\above{1pt}{1pt}{1}{}2\above{1pt}{1pt}{1}{}]\above{1pt}{1pt}{}{0}16\above{1pt}{1pt}{1}{1}{\cdot}1\above{1pt}{1pt}{2}{}3\above{1pt}{1pt}{-}{}{\cdot}1\above{1pt}{1pt}{2}{}11\above{1pt}{1pt}{1}{}$}\spacer%
\instructions{m}%
\EasyButWeakLineBreak%
{${33}\above{1pt}{1pt}{r}{2}{8}\above{1pt}{1pt}{*}{2}{176}\above{1pt}{1pt}{s}{2}{24}\above{1pt}{1pt}{*}{2}{44}\above{1pt}{1pt}{l}{2}{2}\above{1pt}{1pt}{}{2}{528}\above{1pt}{1pt}{}{2}{1}\above{1pt}{1pt}{}{2}{16}\above{1pt}{1pt}{}{2}$}%
\nopagebreak\par%
\nopagebreak\par\leavevmode%
{$\left[\!\llap{\phantom{%
\begingroup \smaller\smaller\smaller\begin{tabular}{@{}c@{}}%
0\\0\\0
\end{tabular}\endgroup%
}}\right.$}%
\begingroup \smaller\smaller\smaller\begin{tabular}{@{}c@{}}%
528\\0\\0
\end{tabular}\endgroup%
\kern3pt%
\begingroup \smaller\smaller\smaller\begin{tabular}{@{}c@{}}%
0\\2\\0
\end{tabular}\endgroup%
\kern3pt%
\begingroup \smaller\smaller\smaller\begin{tabular}{@{}c@{}}%
0\\0\\-1
\end{tabular}\endgroup%
{$\left.\llap{\phantom{%
\begingroup \smaller\smaller\smaller\begin{tabular}{@{}c@{}}%
0\\0\\0
\end{tabular}\endgroup%
}}\!\right]$}%
\EasyButWeakLineBreak%
{$\left[\!\llap{\phantom{%
\begingroup \smaller\smaller\smaller\begin{tabular}{@{}c@{}}%
0\\0\\0
\end{tabular}\endgroup%
}}\right.$}%
\begingroup \smaller\smaller\smaller\begin{tabular}{@{}c@{}}%
-8\\-99\\-231
\end{tabular}\endgroup%
\HardButStrongLineBreak\kern3pt%
\begingroup \smaller\smaller\smaller\begin{tabular}{@{}c@{}}%
-3\\-34\\-84
\end{tabular}\endgroup%
\HardButStrongLineBreak\kern3pt%
\begingroup \smaller\smaller\smaller\begin{tabular}{@{}c@{}}%
-13\\-132\\-352
\end{tabular}\endgroup%
\HardButStrongLineBreak\kern3pt%
\begingroup \smaller\smaller\smaller\begin{tabular}{@{}c@{}}%
-1\\-6\\-24
\end{tabular}\endgroup%
\HardButStrongLineBreak\kern3pt%
\begingroup \smaller\smaller\smaller\begin{tabular}{@{}c@{}}%
-1\\0\\-22
\end{tabular}\endgroup%
\HardButStrongLineBreak\kern3pt%
\begingroup \smaller\smaller\smaller\begin{tabular}{@{}c@{}}%
0\\1\\0
\end{tabular}\endgroup%
\HardButStrongLineBreak\kern3pt%
\begingroup \smaller\smaller\smaller\begin{tabular}{@{}c@{}}%
1\\0\\0
\end{tabular}\endgroup%
\HardButStrongLineBreak\kern3pt%
\begingroup \smaller\smaller\smaller\begin{tabular}{@{}c@{}}%
0\\-1\\-1
\end{tabular}\endgroup%
\HardButStrongLineBreak\kern3pt%
\begingroup \smaller\smaller\smaller\begin{tabular}{@{}c@{}}%
-1\\-16\\-32
\end{tabular}\endgroup%
{$\left.\llap{\phantom{%
\begingroup \smaller\smaller\smaller\begin{tabular}{@{}c@{}}%
0\\0\\0
\end{tabular}\endgroup%
}}\!\right]$}%
%
%
\hbox{}\par\smallskip%
%
%
\leavevmode%
${L_{249.4}}$%
{} : {$[1\above{1pt}{1pt}{-}{}2\above{1pt}{1pt}{1}{}]\above{1pt}{1pt}{}{6}16\above{1pt}{1pt}{-}{3}{\cdot}1\above{1pt}{1pt}{2}{}3\above{1pt}{1pt}{-}{}{\cdot}1\above{1pt}{1pt}{2}{}11\above{1pt}{1pt}{1}{}$}\spacer%
\instructions{m}%
\EasyButWeakLineBreak%
{${132}\above{1pt}{1pt}{l}{2}{2}\above{1pt}{1pt}{}{2}{176}\above{1pt}{1pt}{}{2}{6}\above{1pt}{1pt}{r}{2}{44}\above{1pt}{1pt}{*}{2}{8}\above{1pt}{1pt}{s}{2}{528}\above{1pt}{1pt}{*}{2}{4}\above{1pt}{1pt}{*}{2}{16}\above{1pt}{1pt}{*}{2}$}%
\nopagebreak\par%
\nopagebreak\par\leavevmode%
{$\left[\!\llap{\phantom{%
\begingroup \smaller\smaller\smaller\begin{tabular}{@{}c@{}}%
0\\0\\0
\end{tabular}\endgroup%
}}\right.$}%
\begingroup \smaller\smaller\smaller\begin{tabular}{@{}c@{}}%
-243408\\1584\\1584
\end{tabular}\endgroup%
\kern3pt%
\begingroup \smaller\smaller\smaller\begin{tabular}{@{}c@{}}%
1584\\-10\\-12
\end{tabular}\endgroup%
\kern3pt%
\begingroup \smaller\smaller\smaller\begin{tabular}{@{}c@{}}%
1584\\-12\\-1
\end{tabular}\endgroup%
{$\left.\llap{\phantom{%
\begingroup \smaller\smaller\smaller\begin{tabular}{@{}c@{}}%
0\\0\\0
\end{tabular}\endgroup%
}}\!\right]$}%
\EasyButWeakLineBreak%
{$\left[\!\llap{\phantom{%
\begingroup \smaller\smaller\smaller\begin{tabular}{@{}c@{}}%
0\\0\\0
\end{tabular}\endgroup%
}}\right.$}%
\begingroup \smaller\smaller\smaller\begin{tabular}{@{}c@{}}%
83\\10758\\1914
\end{tabular}\endgroup%
\HardButStrongLineBreak\kern3pt%
\begingroup \smaller\smaller\smaller\begin{tabular}{@{}c@{}}%
7\\907\\162
\end{tabular}\endgroup%
\HardButStrongLineBreak\kern3pt%
\begingroup \smaller\smaller\smaller\begin{tabular}{@{}c@{}}%
53\\6864\\1232
\end{tabular}\endgroup%
\HardButStrongLineBreak\kern3pt%
\begingroup \smaller\smaller\smaller\begin{tabular}{@{}c@{}}%
1\\129\\24
\end{tabular}\endgroup%
\HardButStrongLineBreak\kern3pt%
\begingroup \smaller\smaller\smaller\begin{tabular}{@{}c@{}}%
-1\\-132\\-22
\end{tabular}\endgroup%
\HardButStrongLineBreak\kern3pt%
\begingroup \smaller\smaller\smaller\begin{tabular}{@{}c@{}}%
-1\\-130\\-24
\end{tabular}\endgroup%
\HardButStrongLineBreak\kern3pt%
\begingroup \smaller\smaller\smaller\begin{tabular}{@{}c@{}}%
1\\132\\0
\end{tabular}\endgroup%
\HardButStrongLineBreak\kern3pt%
\begingroup \smaller\smaller\smaller\begin{tabular}{@{}c@{}}%
1\\130\\22
\end{tabular}\endgroup%
\HardButStrongLineBreak\kern3pt%
\begingroup \smaller\smaller\smaller\begin{tabular}{@{}c@{}}%
7\\908\\160
\end{tabular}\endgroup%
{$\left.\llap{\phantom{%
\begingroup \smaller\smaller\smaller\begin{tabular}{@{}c@{}}%
0\\0\\0
\end{tabular}\endgroup%
}}\!\right]$}%
%
%
\hbox{}\par\smallskip%
%
%
\leavevmode%
${L_{249.5}}$%
{} : {$[1\above{1pt}{1pt}{-}{}2\above{1pt}{1pt}{1}{}]\above{1pt}{1pt}{}{4}16\above{1pt}{1pt}{-}{5}{\cdot}1\above{1pt}{1pt}{2}{}3\above{1pt}{1pt}{-}{}{\cdot}1\above{1pt}{1pt}{2}{}11\above{1pt}{1pt}{1}{}$}\EasyButWeakLineBreak%
{${132}\above{1pt}{1pt}{*}{2}{8}\above{1pt}{1pt}{s}{2}{176}\above{1pt}{1pt}{*}{2}{24}\above{1pt}{1pt}{l}{2}{11}\above{1pt}{1pt}{}{2}{2}\above{1pt}{1pt}{r}{2}{528}\above{1pt}{1pt}{s}{2}{4}\above{1pt}{1pt}{s}{2}{16}\above{1pt}{1pt}{s}{2}$}%
\nopagebreak\par%
\nopagebreak\par\leavevmode%
{$\left[\!\llap{\phantom{%
\begingroup \smaller\smaller\smaller\begin{tabular}{@{}c@{}}%
0\\0\\0
\end{tabular}\endgroup%
}}\right.$}%
\begingroup \smaller\smaller\smaller\begin{tabular}{@{}c@{}}%
19536\\4752\\1056
\end{tabular}\endgroup%
\kern3pt%
\begingroup \smaller\smaller\smaller\begin{tabular}{@{}c@{}}%
4752\\1154\\238
\end{tabular}\endgroup%
\kern3pt%
\begingroup \smaller\smaller\smaller\begin{tabular}{@{}c@{}}%
1056\\238\\-131
\end{tabular}\endgroup%
{$\left.\llap{\phantom{%
\begingroup \smaller\smaller\smaller\begin{tabular}{@{}c@{}}%
0\\0\\0
\end{tabular}\endgroup%
}}\!\right]$}%
\EasyButWeakLineBreak%
{$\left[\!\llap{\phantom{%
\begingroup \smaller\smaller\smaller\begin{tabular}{@{}c@{}}%
0\\0\\0
\end{tabular}\endgroup%
}}\right.$}%
\begingroup \smaller\smaller\smaller\begin{tabular}{@{}c@{}}%
157\\-660\\66
\end{tabular}\endgroup%
\HardButStrongLineBreak\kern3pt%
\begingroup \smaller\smaller\smaller\begin{tabular}{@{}c@{}}%
9\\-38\\4
\end{tabular}\endgroup%
\HardButStrongLineBreak\kern3pt%
\begingroup \smaller\smaller\smaller\begin{tabular}{@{}c@{}}%
-639\\2684\\-264
\end{tabular}\endgroup%
\HardButStrongLineBreak\kern3pt%
\begingroup \smaller\smaller\smaller\begin{tabular}{@{}c@{}}%
-347\\1458\\-144
\end{tabular}\endgroup%
\HardButStrongLineBreak\kern3pt%
\begingroup \smaller\smaller\smaller\begin{tabular}{@{}c@{}}%
-1191\\5005\\-495
\end{tabular}\endgroup%
\HardButStrongLineBreak\kern3pt%
\begingroup \smaller\smaller\smaller\begin{tabular}{@{}c@{}}%
-558\\2345\\-232
\end{tabular}\endgroup%
\HardButStrongLineBreak\kern3pt%
\begingroup \smaller\smaller\smaller\begin{tabular}{@{}c@{}}%
-12061\\50688\\-5016
\end{tabular}\endgroup%
\HardButStrongLineBreak\kern3pt%
\begingroup \smaller\smaller\smaller\begin{tabular}{@{}c@{}}%
-197\\828\\-82
\end{tabular}\endgroup%
\HardButStrongLineBreak\kern3pt%
\begingroup \smaller\smaller\smaller\begin{tabular}{@{}c@{}}%
-19\\80\\-8
\end{tabular}\endgroup%
{$\left.\llap{\phantom{%
\begingroup \smaller\smaller\smaller\begin{tabular}{@{}c@{}}%
0\\0\\0
\end{tabular}\endgroup%
}}\!\right]$}%

\medskip%
%
\leavevmode\llap{}%
$W_{250}$%
\qquad\llap{8} lattices, $\chi=8$%
\hfill%
$2\slashthree22|2\rtimes D_{2}$%
\nopagebreak\smallskip\hrule\nopagebreak\medskip%
%
%
\leavevmode%
${L_{250.1}}$%
{} : {$1\above{1pt}{1pt}{-2}{{\rm II}}8\above{1pt}{1pt}{-}{3}{\cdot}1\above{1pt}{1pt}{1}{}3\above{1pt}{1pt}{-}{}9\above{1pt}{1pt}{1}{}{\cdot}1\above{1pt}{1pt}{-2}{}5\above{1pt}{1pt}{-}{}$}\spacer%
\instructions{2}%
\EasyButWeakLineBreak%
{${90}\above{1pt}{1pt}{b}{2}{6}\above{1pt}{1pt}{+}{3}{6}\above{1pt}{1pt}{b}{2}{10}\above{1pt}{1pt}{l}{2}{24}\above{1pt}{1pt}{r}{2}$}%
\nopagebreak\par%
\nopagebreak\par\leavevmode%
{$\left[\!\llap{\phantom{%
\begingroup \smaller\smaller\smaller\begin{tabular}{@{}c@{}}%
0\\0\\0
\end{tabular}\endgroup%
}}\right.$}%
\begingroup \smaller\smaller\smaller\begin{tabular}{@{}c@{}}%
-377640\\5040\\-56160
\end{tabular}\endgroup%
\kern3pt%
\begingroup \smaller\smaller\smaller\begin{tabular}{@{}c@{}}%
5040\\-66\\765
\end{tabular}\endgroup%
\kern3pt%
\begingroup \smaller\smaller\smaller\begin{tabular}{@{}c@{}}%
-56160\\765\\-8162
\end{tabular}\endgroup%
{$\left.\llap{\phantom{%
\begingroup \smaller\smaller\smaller\begin{tabular}{@{}c@{}}%
0\\0\\0
\end{tabular}\endgroup%
}}\!\right]$}%
\EasyButWeakLineBreak%
{$\left[\!\llap{\phantom{%
\begingroup \smaller\smaller\smaller\begin{tabular}{@{}c@{}}%
0\\0\\0
\end{tabular}\endgroup%
}}\right.$}%
\begingroup \smaller\smaller\smaller\begin{tabular}{@{}c@{}}%
28\\1095\\-90
\end{tabular}\endgroup%
\HardButStrongLineBreak\kern3pt%
\begingroup \smaller\smaller\smaller\begin{tabular}{@{}c@{}}%
15\\589\\-48
\end{tabular}\endgroup%
\HardButStrongLineBreak\kern3pt%
\begingroup \smaller\smaller\smaller\begin{tabular}{@{}c@{}}%
-14\\-548\\45
\end{tabular}\endgroup%
\HardButStrongLineBreak\kern3pt%
\begingroup \smaller\smaller\smaller\begin{tabular}{@{}c@{}}%
-39\\-1530\\125
\end{tabular}\endgroup%
\HardButStrongLineBreak\kern3pt%
\begingroup \smaller\smaller\smaller\begin{tabular}{@{}c@{}}%
-45\\-1768\\144
\end{tabular}\endgroup%
{$\left.\llap{\phantom{%
\begingroup \smaller\smaller\smaller\begin{tabular}{@{}c@{}}%
0\\0\\0
\end{tabular}\endgroup%
}}\!\right]$}%

\medskip%
%
\leavevmode\llap{}%
$W_{251}$%
\qquad\llap{8} lattices, $\chi=8$%
\hfill%
$26|62|\rtimes D_{2}$%
\nopagebreak\smallskip\hrule\nopagebreak\medskip%
%
%
\leavevmode%
${L_{251.1}}$%
{} : {$1\above{1pt}{1pt}{-2}{{\rm II}}8\above{1pt}{1pt}{-}{3}{\cdot}1\above{1pt}{1pt}{-}{}3\above{1pt}{1pt}{-}{}9\above{1pt}{1pt}{-}{}{\cdot}1\above{1pt}{1pt}{-2}{}5\above{1pt}{1pt}{-}{}$}\spacer%
\instructions{2}%
\EasyButWeakLineBreak%
{${24}\above{1pt}{1pt}{r}{2}{18}\above{1pt}{1pt}{}{6}{6}\above{1pt}{1pt}{}{6}{2}\above{1pt}{1pt}{l}{2}$}%
\nopagebreak\par%
\nopagebreak\par\leavevmode%
{$\left[\!\llap{\phantom{%
\begingroup \smaller\smaller\smaller\begin{tabular}{@{}c@{}}%
0\\0\\0
\end{tabular}\endgroup%
}}\right.$}%
\begingroup \smaller\smaller\smaller\begin{tabular}{@{}c@{}}%
-1572840\\10440\\213120
\end{tabular}\endgroup%
\kern3pt%
\begingroup \smaller\smaller\smaller\begin{tabular}{@{}c@{}}%
10440\\-66\\-1515
\end{tabular}\endgroup%
\kern3pt%
\begingroup \smaller\smaller\smaller\begin{tabular}{@{}c@{}}%
213120\\-1515\\-25822
\end{tabular}\endgroup%
{$\left.\llap{\phantom{%
\begingroup \smaller\smaller\smaller\begin{tabular}{@{}c@{}}%
0\\0\\0
\end{tabular}\endgroup%
}}\!\right]$}%
\EasyButWeakLineBreak%
{$\left[\!\llap{\phantom{%
\begingroup \smaller\smaller\smaller\begin{tabular}{@{}c@{}}%
0\\0\\0
\end{tabular}\endgroup%
}}\right.$}%
\begingroup \smaller\smaller\smaller\begin{tabular}{@{}c@{}}%
-81\\-7304\\-240
\end{tabular}\endgroup%
\HardButStrongLineBreak\kern3pt%
\begingroup \smaller\smaller\smaller\begin{tabular}{@{}c@{}}%
-79\\-7125\\-234
\end{tabular}\endgroup%
\HardButStrongLineBreak\kern3pt%
\begingroup \smaller\smaller\smaller\begin{tabular}{@{}c@{}}%
80\\7214\\237
\end{tabular}\endgroup%
\HardButStrongLineBreak\kern3pt%
\begingroup \smaller\smaller\smaller\begin{tabular}{@{}c@{}}%
26\\2345\\77
\end{tabular}\endgroup%
{$\left.\llap{\phantom{%
\begingroup \smaller\smaller\smaller\begin{tabular}{@{}c@{}}%
0\\0\\0
\end{tabular}\endgroup%
}}\!\right]$}%

\medskip%
%
\leavevmode\llap{}%
$W_{252}$%
\qquad\llap{16} lattices, $\chi=18$%
\hfill%
$22\slashtwo222|2\rtimes D_{2}$%
\nopagebreak\smallskip\hrule\nopagebreak\medskip%
%
%
\leavevmode%
${L_{252.1}}$%
{} : {$1\above{1pt}{1pt}{-2}{2}8\above{1pt}{1pt}{-}{5}{\cdot}1\above{1pt}{1pt}{-}{}3\above{1pt}{1pt}{-}{}9\above{1pt}{1pt}{-}{}{\cdot}1\above{1pt}{1pt}{2}{}5\above{1pt}{1pt}{1}{}$}\spacer%
\instructions{2}%
\EasyButWeakLineBreak%
{${20}\above{1pt}{1pt}{*}{2}{72}\above{1pt}{1pt}{b}{2}{2}\above{1pt}{1pt}{b}{2}{18}\above{1pt}{1pt}{b}{2}{8}\above{1pt}{1pt}{*}{2}{180}\above{1pt}{1pt}{s}{2}{24}\above{1pt}{1pt}{s}{2}$}%
\nopagebreak\par%
\nopagebreak\par\leavevmode%
{$\left[\!\llap{\phantom{%
\begingroup \smaller\smaller\smaller\begin{tabular}{@{}c@{}}%
0\\0\\0
\end{tabular}\endgroup%
}}\right.$}%
\begingroup \smaller\smaller\smaller\begin{tabular}{@{}c@{}}%
49320\\-4320\\-360
\end{tabular}\endgroup%
\kern3pt%
\begingroup \smaller\smaller\smaller\begin{tabular}{@{}c@{}}%
-4320\\375\\33
\end{tabular}\endgroup%
\kern3pt%
\begingroup \smaller\smaller\smaller\begin{tabular}{@{}c@{}}%
-360\\33\\2
\end{tabular}\endgroup%
{$\left.\llap{\phantom{%
\begingroup \smaller\smaller\smaller\begin{tabular}{@{}c@{}}%
0\\0\\0
\end{tabular}\endgroup%
}}\!\right]$}%
\EasyButWeakLineBreak%
{$\left[\!\llap{\phantom{%
\begingroup \smaller\smaller\smaller\begin{tabular}{@{}c@{}}%
0\\0\\0
\end{tabular}\endgroup%
}}\right.$}%
\begingroup \smaller\smaller\smaller\begin{tabular}{@{}c@{}}%
1\\10\\20
\end{tabular}\endgroup%
\HardButStrongLineBreak\kern3pt%
\begingroup \smaller\smaller\smaller\begin{tabular}{@{}c@{}}%
5\\48\\108
\end{tabular}\endgroup%
\HardButStrongLineBreak\kern3pt%
\begingroup \smaller\smaller\smaller\begin{tabular}{@{}c@{}}%
0\\0\\-1
\end{tabular}\endgroup%
\HardButStrongLineBreak\kern3pt%
\begingroup \smaller\smaller\smaller\begin{tabular}{@{}c@{}}%
-7\\-66\\-171
\end{tabular}\endgroup%
\HardButStrongLineBreak\kern3pt%
\begingroup \smaller\smaller\smaller\begin{tabular}{@{}c@{}}%
-17\\-160\\-412
\end{tabular}\endgroup%
\HardButStrongLineBreak\kern3pt%
\begingroup \smaller\smaller\smaller\begin{tabular}{@{}c@{}}%
-67\\-630\\-1620
\end{tabular}\endgroup%
\HardButStrongLineBreak\kern3pt%
\begingroup \smaller\smaller\smaller\begin{tabular}{@{}c@{}}%
-3\\-28\\-72
\end{tabular}\endgroup%
{$\left.\llap{\phantom{%
\begingroup \smaller\smaller\smaller\begin{tabular}{@{}c@{}}%
0\\0\\0
\end{tabular}\endgroup%
}}\!\right]$}%
%
%
\hbox{}\par\smallskip%
%
%
\leavevmode%
${L_{252.2}}$%
{} : {$1\above{1pt}{1pt}{2}{2}8\above{1pt}{1pt}{1}{1}{\cdot}1\above{1pt}{1pt}{-}{}3\above{1pt}{1pt}{-}{}9\above{1pt}{1pt}{-}{}{\cdot}1\above{1pt}{1pt}{2}{}5\above{1pt}{1pt}{1}{}$}\spacer%
\instructions{m}%
\EasyButWeakLineBreak%
{${5}\above{1pt}{1pt}{}{2}{72}\above{1pt}{1pt}{r}{2}{2}\above{1pt}{1pt}{s}{2}{18}\above{1pt}{1pt}{l}{2}{8}\above{1pt}{1pt}{}{2}{45}\above{1pt}{1pt}{r}{2}{24}\above{1pt}{1pt}{l}{2}$}%
\nopagebreak\par%
\nopagebreak\par\leavevmode%
{$\left[\!\llap{\phantom{%
\begingroup \smaller\smaller\smaller\begin{tabular}{@{}c@{}}%
0\\0\\0
\end{tabular}\endgroup%
}}\right.$}%
\begingroup \smaller\smaller\smaller\begin{tabular}{@{}c@{}}%
105480\\4320\\360
\end{tabular}\endgroup%
\kern3pt%
\begingroup \smaller\smaller\smaller\begin{tabular}{@{}c@{}}%
4320\\177\\15
\end{tabular}\endgroup%
\kern3pt%
\begingroup \smaller\smaller\smaller\begin{tabular}{@{}c@{}}%
360\\15\\2
\end{tabular}\endgroup%
{$\left.\llap{\phantom{%
\begingroup \smaller\smaller\smaller\begin{tabular}{@{}c@{}}%
0\\0\\0
\end{tabular}\endgroup%
}}\!\right]$}%
\EasyButWeakLineBreak%
{$\left[\!\llap{\phantom{%
\begingroup \smaller\smaller\smaller\begin{tabular}{@{}c@{}}%
0\\0\\0
\end{tabular}\endgroup%
}}\right.$}%
\begingroup \smaller\smaller\smaller\begin{tabular}{@{}c@{}}%
-1\\25\\-5
\end{tabular}\endgroup%
\HardButStrongLineBreak\kern3pt%
\begingroup \smaller\smaller\smaller\begin{tabular}{@{}c@{}}%
1\\-24\\0
\end{tabular}\endgroup%
\HardButStrongLineBreak\kern3pt%
\begingroup \smaller\smaller\smaller\begin{tabular}{@{}c@{}}%
0\\0\\-1
\end{tabular}\endgroup%
\HardButStrongLineBreak\kern3pt%
\begingroup \smaller\smaller\smaller\begin{tabular}{@{}c@{}}%
-8\\198\\-45
\end{tabular}\endgroup%
\HardButStrongLineBreak\kern3pt%
\begingroup \smaller\smaller\smaller\begin{tabular}{@{}c@{}}%
-21\\520\\-112
\end{tabular}\endgroup%
\HardButStrongLineBreak\kern3pt%
\begingroup \smaller\smaller\smaller\begin{tabular}{@{}c@{}}%
-43\\1065\\-225
\end{tabular}\endgroup%
\HardButStrongLineBreak\kern3pt%
\begingroup \smaller\smaller\smaller\begin{tabular}{@{}c@{}}%
-5\\124\\-24
\end{tabular}\endgroup%
{$\left.\llap{\phantom{%
\begingroup \smaller\smaller\smaller\begin{tabular}{@{}c@{}}%
0\\0\\0
\end{tabular}\endgroup%
}}\!\right]$}%

\medskip%
%
\leavevmode\llap{}%
$W_{253}$%
\qquad\llap{4} lattices, $\chi=36$%
\hfill%
$24\slashinfty422\slashtwo2\rtimes D_{2}$%
\nopagebreak\smallskip\hrule\nopagebreak\medskip%
%
%
\leavevmode%
${L_{253.1}}$%
{} : {$1\above{1pt}{1pt}{2}{2}32\above{1pt}{1pt}{1}{7}{\cdot}1\above{1pt}{1pt}{2}{}9\above{1pt}{1pt}{1}{}$}\EasyButWeakLineBreak%
{${36}\above{1pt}{1pt}{l}{2}{1}\above{1pt}{1pt}{}{4}{2}\above{1pt}{1pt}{24,23}{\infty b}{2}\above{1pt}{1pt}{*}{4}{4}\above{1pt}{1pt}{l}{2}{9}\above{1pt}{1pt}{}{2}{1}\above{1pt}{1pt}{r}{2}{4}\above{1pt}{1pt}{*}{2}$}%
\nopagebreak\par%
\nopagebreak\par\leavevmode%
{$\left[\!\llap{\phantom{%
\begingroup \smaller\smaller\smaller\begin{tabular}{@{}c@{}}%
0\\0\\0
\end{tabular}\endgroup%
}}\right.$}%
\begingroup \smaller\smaller\smaller\begin{tabular}{@{}c@{}}%
-525600\\3744\\2016
\end{tabular}\endgroup%
\kern3pt%
\begingroup \smaller\smaller\smaller\begin{tabular}{@{}c@{}}%
3744\\-23\\-16
\end{tabular}\endgroup%
\kern3pt%
\begingroup \smaller\smaller\smaller\begin{tabular}{@{}c@{}}%
2016\\-16\\-7
\end{tabular}\endgroup%
{$\left.\llap{\phantom{%
\begingroup \smaller\smaller\smaller\begin{tabular}{@{}c@{}}%
0\\0\\0
\end{tabular}\endgroup%
}}\!\right]$}%
\EasyButWeakLineBreak%
{$\left[\!\llap{\phantom{%
\begingroup \smaller\smaller\smaller\begin{tabular}{@{}c@{}}%
0\\0\\0
\end{tabular}\endgroup%
}}\right.$}%
\begingroup \smaller\smaller\smaller\begin{tabular}{@{}c@{}}%
43\\2736\\6102
\end{tabular}\endgroup%
\HardButStrongLineBreak\kern3pt%
\begingroup \smaller\smaller\smaller\begin{tabular}{@{}c@{}}%
12\\763\\1704
\end{tabular}\endgroup%
\HardButStrongLineBreak\kern3pt%
\begingroup \smaller\smaller\smaller\begin{tabular}{@{}c@{}}%
6\\381\\853
\end{tabular}\endgroup%
\HardButStrongLineBreak\kern3pt%
\begingroup \smaller\smaller\smaller\begin{tabular}{@{}c@{}}%
1\\63\\143
\end{tabular}\endgroup%
\HardButStrongLineBreak\kern3pt%
\begingroup \smaller\smaller\smaller\begin{tabular}{@{}c@{}}%
-1\\-64\\-142
\end{tabular}\endgroup%
\HardButStrongLineBreak\kern3pt%
\begingroup \smaller\smaller\smaller\begin{tabular}{@{}c@{}}%
-1\\-63\\-144
\end{tabular}\endgroup%
\HardButStrongLineBreak\kern3pt%
\begingroup \smaller\smaller\smaller\begin{tabular}{@{}c@{}}%
1\\64\\141
\end{tabular}\endgroup%
\HardButStrongLineBreak\kern3pt%
\begingroup \smaller\smaller\smaller\begin{tabular}{@{}c@{}}%
7\\446\\992
\end{tabular}\endgroup%
{$\left.\llap{\phantom{%
\begingroup \smaller\smaller\smaller\begin{tabular}{@{}c@{}}%
0\\0\\0
\end{tabular}\endgroup%
}}\!\right]$}%
%
%
%
%
%
%
%
%
%
%
%
%
%
%

\medskip%
%
\leavevmode\llap{}%
$W_{254}$%
\qquad\llap{22} lattices, $\chi=30$%
\hfill%
$22\slashtwo2222|22\rtimes D_{2}$%
\nopagebreak\smallskip\hrule\nopagebreak\medskip%
%
%
\leavevmode%
${L_{254.1}}$%
{} : {$1\above{1pt}{1pt}{2}{{\rm II}}4\above{1pt}{1pt}{1}{1}{\cdot}1\above{1pt}{1pt}{1}{}3\above{1pt}{1pt}{1}{}9\above{1pt}{1pt}{1}{}{\cdot}1\above{1pt}{1pt}{2}{}11\above{1pt}{1pt}{-}{}$}\spacer%
\instructions{2}%
\EasyButWeakLineBreak%
{${22}\above{1pt}{1pt}{b}{2}{12}\above{1pt}{1pt}{*}{2}{4}\above{1pt}{1pt}{*}{2}{36}\above{1pt}{1pt}{*}{2}{12}\above{1pt}{1pt}{b}{2}{198}\above{1pt}{1pt}{l}{2}{4}\above{1pt}{1pt}{r}{2}{66}\above{1pt}{1pt}{l}{2}{36}\above{1pt}{1pt}{r}{2}$}%
\nopagebreak\par%
\nopagebreak\par\leavevmode%
{$\left[\!\llap{\phantom{%
\begingroup \smaller\smaller\smaller\begin{tabular}{@{}c@{}}%
0\\0\\0
\end{tabular}\endgroup%
}}\right.$}%
\begingroup \smaller\smaller\smaller\begin{tabular}{@{}c@{}}%
-489852\\1980\\1188
\end{tabular}\endgroup%
\kern3pt%
\begingroup \smaller\smaller\smaller\begin{tabular}{@{}c@{}}%
1980\\12\\-9
\end{tabular}\endgroup%
\kern3pt%
\begingroup \smaller\smaller\smaller\begin{tabular}{@{}c@{}}%
1188\\-9\\-2
\end{tabular}\endgroup%
{$\left.\llap{\phantom{%
\begingroup \smaller\smaller\smaller\begin{tabular}{@{}c@{}}%
0\\0\\0
\end{tabular}\endgroup%
}}\!\right]$}%
\EasyButWeakLineBreak%
{$\left[\!\llap{\phantom{%
\begingroup \smaller\smaller\smaller\begin{tabular}{@{}c@{}}%
0\\0\\0
\end{tabular}\endgroup%
}}\right.$}%
\begingroup \smaller\smaller\smaller\begin{tabular}{@{}c@{}}%
15\\946\\4565
\end{tabular}\endgroup%
\HardButStrongLineBreak\kern3pt%
\begingroup \smaller\smaller\smaller\begin{tabular}{@{}c@{}}%
7\\440\\2130
\end{tabular}\endgroup%
\HardButStrongLineBreak\kern3pt%
\begingroup \smaller\smaller\smaller\begin{tabular}{@{}c@{}}%
1\\62\\304
\end{tabular}\endgroup%
\HardButStrongLineBreak\kern3pt%
\begingroup \smaller\smaller\smaller\begin{tabular}{@{}c@{}}%
-1\\-66\\-306
\end{tabular}\endgroup%
\HardButStrongLineBreak\kern3pt%
\begingroup \smaller\smaller\smaller\begin{tabular}{@{}c@{}}%
-1\\-64\\-306
\end{tabular}\endgroup%
\HardButStrongLineBreak\kern3pt%
\begingroup \smaller\smaller\smaller\begin{tabular}{@{}c@{}}%
1\\66\\297
\end{tabular}\endgroup%
\HardButStrongLineBreak\kern3pt%
\begingroup \smaller\smaller\smaller\begin{tabular}{@{}c@{}}%
1\\64\\304
\end{tabular}\endgroup%
\HardButStrongLineBreak\kern3pt%
\begingroup \smaller\smaller\smaller\begin{tabular}{@{}c@{}}%
9\\572\\2739
\end{tabular}\endgroup%
\HardButStrongLineBreak\kern3pt%
\begingroup \smaller\smaller\smaller\begin{tabular}{@{}c@{}}%
11\\696\\3348
\end{tabular}\endgroup%
{$\left.\llap{\phantom{%
\begingroup \smaller\smaller\smaller\begin{tabular}{@{}c@{}}%
0\\0\\0
\end{tabular}\endgroup%
}}\!\right]$}%
%
%
\hbox{}\par\smallskip%
%
%
\leavevmode%
${L_{254.2}}$%
{} : {$1\above{1pt}{1pt}{2}{2}8\above{1pt}{1pt}{1}{7}{\cdot}1\above{1pt}{1pt}{-}{}3\above{1pt}{1pt}{-}{}9\above{1pt}{1pt}{-}{}{\cdot}1\above{1pt}{1pt}{2}{}11\above{1pt}{1pt}{1}{}$}\spacer%
\instructions{2}%
\EasyButWeakLineBreak%
{${396}\above{1pt}{1pt}{*}{2}{24}\above{1pt}{1pt}{b}{2}{18}\above{1pt}{1pt}{s}{2}{2}\above{1pt}{1pt}{b}{2}{24}\above{1pt}{1pt}{*}{2}{44}\above{1pt}{1pt}{s}{2}{72}\above{1pt}{1pt}{l}{2}{33}\above{1pt}{1pt}{r}{2}{8}\above{1pt}{1pt}{s}{2}$}%
\nopagebreak\par%
\nopagebreak\par\leavevmode%
{$\left[\!\llap{\phantom{%
\begingroup \smaller\smaller\smaller\begin{tabular}{@{}c@{}}%
0\\0\\0
\end{tabular}\endgroup%
}}\right.$}%
\begingroup \smaller\smaller\smaller\begin{tabular}{@{}c@{}}%
-490248\\1584\\1584
\end{tabular}\endgroup%
\kern3pt%
\begingroup \smaller\smaller\smaller\begin{tabular}{@{}c@{}}%
1584\\-3\\-9
\end{tabular}\endgroup%
\kern3pt%
\begingroup \smaller\smaller\smaller\begin{tabular}{@{}c@{}}%
1584\\-9\\2
\end{tabular}\endgroup%
{$\left.\llap{\phantom{%
\begingroup \smaller\smaller\smaller\begin{tabular}{@{}c@{}}%
0\\0\\0
\end{tabular}\endgroup%
}}\!\right]$}%
\EasyButWeakLineBreak%
{$\left[\!\llap{\phantom{%
\begingroup \smaller\smaller\smaller\begin{tabular}{@{}c@{}}%
0\\0\\0
\end{tabular}\endgroup%
}}\right.$}%
\begingroup \smaller\smaller\smaller\begin{tabular}{@{}c@{}}%
-41\\-8250\\-4554
\end{tabular}\endgroup%
\HardButStrongLineBreak\kern3pt%
\begingroup \smaller\smaller\smaller\begin{tabular}{@{}c@{}}%
-7\\-1408\\-780
\end{tabular}\endgroup%
\HardButStrongLineBreak\kern3pt%
\begingroup \smaller\smaller\smaller\begin{tabular}{@{}c@{}}%
-2\\-402\\-225
\end{tabular}\endgroup%
\HardButStrongLineBreak\kern3pt%
\begingroup \smaller\smaller\smaller\begin{tabular}{@{}c@{}}%
0\\0\\-1
\end{tabular}\endgroup%
\HardButStrongLineBreak\kern3pt%
\begingroup \smaller\smaller\smaller\begin{tabular}{@{}c@{}}%
1\\200\\108
\end{tabular}\endgroup%
\HardButStrongLineBreak\kern3pt%
\begingroup \smaller\smaller\smaller\begin{tabular}{@{}c@{}}%
1\\198\\110
\end{tabular}\endgroup%
\HardButStrongLineBreak\kern3pt%
\begingroup \smaller\smaller\smaller\begin{tabular}{@{}c@{}}%
-1\\-204\\-108
\end{tabular}\endgroup%
\HardButStrongLineBreak\kern3pt%
\begingroup \smaller\smaller\smaller\begin{tabular}{@{}c@{}}%
-3\\-605\\-330
\end{tabular}\endgroup%
\HardButStrongLineBreak\kern3pt%
\begingroup \smaller\smaller\smaller\begin{tabular}{@{}c@{}}%
-3\\-604\\-332
\end{tabular}\endgroup%
{$\left.\llap{\phantom{%
\begingroup \smaller\smaller\smaller\begin{tabular}{@{}c@{}}%
0\\0\\0
\end{tabular}\endgroup%
}}\!\right]$}%
%
%
\hbox{}\par\smallskip%
%
%
\leavevmode%
${L_{254.3}}$%
{} : {$1\above{1pt}{1pt}{-2}{2}8\above{1pt}{1pt}{-}{3}{\cdot}1\above{1pt}{1pt}{-}{}3\above{1pt}{1pt}{-}{}9\above{1pt}{1pt}{-}{}{\cdot}1\above{1pt}{1pt}{2}{}11\above{1pt}{1pt}{1}{}$}\spacer%
\instructions{m}%
\EasyButWeakLineBreak%
{${99}\above{1pt}{1pt}{}{2}{24}\above{1pt}{1pt}{r}{2}{18}\above{1pt}{1pt}{b}{2}{2}\above{1pt}{1pt}{l}{2}{24}\above{1pt}{1pt}{}{2}{11}\above{1pt}{1pt}{r}{2}{72}\above{1pt}{1pt}{s}{2}{132}\above{1pt}{1pt}{s}{2}{8}\above{1pt}{1pt}{l}{2}$}%
\nopagebreak\par%
\nopagebreak\par\leavevmode%
{$\left[\!\llap{\phantom{%
\begingroup \smaller\smaller\smaller\begin{tabular}{@{}c@{}}%
0\\0\\0
\end{tabular}\endgroup%
}}\right.$}%
\begingroup \smaller\smaller\smaller\begin{tabular}{@{}c@{}}%
127512\\-3960\\0
\end{tabular}\endgroup%
\kern3pt%
\begingroup \smaller\smaller\smaller\begin{tabular}{@{}c@{}}%
-3960\\123\\0
\end{tabular}\endgroup%
\kern3pt%
\begingroup \smaller\smaller\smaller\begin{tabular}{@{}c@{}}%
0\\0\\-1
\end{tabular}\endgroup%
{$\left.\llap{\phantom{%
\begingroup \smaller\smaller\smaller\begin{tabular}{@{}c@{}}%
0\\0\\0
\end{tabular}\endgroup%
}}\!\right]$}%
\EasyButWeakLineBreak%
{$\left[\!\llap{\phantom{%
\begingroup \smaller\smaller\smaller\begin{tabular}{@{}c@{}}%
0\\0\\0
\end{tabular}\endgroup%
}}\right.$}%
\begingroup \smaller\smaller\smaller\begin{tabular}{@{}c@{}}%
1\\33\\0
\end{tabular}\endgroup%
\HardButStrongLineBreak\kern3pt%
\begingroup \smaller\smaller\smaller\begin{tabular}{@{}c@{}}%
-1\\-32\\0
\end{tabular}\endgroup%
\HardButStrongLineBreak\kern3pt%
\begingroup \smaller\smaller\smaller\begin{tabular}{@{}c@{}}%
-1\\-33\\-9
\end{tabular}\endgroup%
\HardButStrongLineBreak\kern3pt%
\begingroup \smaller\smaller\smaller\begin{tabular}{@{}c@{}}%
0\\-1\\-11
\end{tabular}\endgroup%
\HardButStrongLineBreak\kern3pt%
\begingroup \smaller\smaller\smaller\begin{tabular}{@{}c@{}}%
3\\88\\-96
\end{tabular}\endgroup%
\HardButStrongLineBreak\kern3pt%
\begingroup \smaller\smaller\smaller\begin{tabular}{@{}c@{}}%
4\\121\\-88
\end{tabular}\endgroup%
\HardButStrongLineBreak\kern3pt%
\begingroup \smaller\smaller\smaller\begin{tabular}{@{}c@{}}%
7\\216\\-108
\end{tabular}\endgroup%
\HardButStrongLineBreak\kern3pt%
\begingroup \smaller\smaller\smaller\begin{tabular}{@{}c@{}}%
7\\220\\-66
\end{tabular}\endgroup%
\HardButStrongLineBreak\kern3pt%
\begingroup \smaller\smaller\smaller\begin{tabular}{@{}c@{}}%
1\\32\\-4
\end{tabular}\endgroup%
{$\left.\llap{\phantom{%
\begingroup \smaller\smaller\smaller\begin{tabular}{@{}c@{}}%
0\\0\\0
\end{tabular}\endgroup%
}}\!\right]$}%

\medskip%
%
\leavevmode\llap{}%
$W_{255}$%
\qquad\llap{24} lattices, $\chi=27$%
\hfill%
$42222222$%
\nopagebreak\smallskip\hrule\nopagebreak\medskip%
%
%
\leavevmode%
${L_{255.1}}$%
{} : {$1\above{1pt}{1pt}{-2}{{\rm II}}4\above{1pt}{1pt}{1}{7}{\cdot}1\above{1pt}{1pt}{2}{}9\above{1pt}{1pt}{1}{}{\cdot}1\above{1pt}{1pt}{2}{}5\above{1pt}{1pt}{-}{}{\cdot}1\above{1pt}{1pt}{2}{}7\above{1pt}{1pt}{1}{}$}\spacer%
\instructions{2}%
\EasyButWeakLineBreak%
{${2}\above{1pt}{1pt}{*}{4}{4}\above{1pt}{1pt}{*}{2}{36}\above{1pt}{1pt}{b}{2}{14}\above{1pt}{1pt}{s}{2}{90}\above{1pt}{1pt}{l}{2}{28}\above{1pt}{1pt}{r}{2}{10}\above{1pt}{1pt}{l}{2}{252}\above{1pt}{1pt}{r}{2}$}%
\nopagebreak\par%
\nopagebreak\par\leavevmode%
{$\left[\!\llap{\phantom{%
\begingroup \smaller\smaller\smaller\begin{tabular}{@{}c@{}}%
0\\0\\0
\end{tabular}\endgroup%
}}\right.$}%
\begingroup \smaller\smaller\smaller\begin{tabular}{@{}c@{}}%
-102537540\\52920\\63000
\end{tabular}\endgroup%
\kern3pt%
\begingroup \smaller\smaller\smaller\begin{tabular}{@{}c@{}}%
52920\\-26\\-35
\end{tabular}\endgroup%
\kern3pt%
\begingroup \smaller\smaller\smaller\begin{tabular}{@{}c@{}}%
63000\\-35\\-34
\end{tabular}\endgroup%
{$\left.\llap{\phantom{%
\begingroup \smaller\smaller\smaller\begin{tabular}{@{}c@{}}%
0\\0\\0
\end{tabular}\endgroup%
}}\!\right]$}%
\EasyButWeakLineBreak%
{$\left[\!\llap{\phantom{%
\begingroup \smaller\smaller\smaller\begin{tabular}{@{}c@{}}%
0\\0\\0
\end{tabular}\endgroup%
}}\right.$}%
\begingroup \smaller\smaller\smaller\begin{tabular}{@{}c@{}}%
1\\1188\\629
\end{tabular}\endgroup%
\HardButStrongLineBreak\kern3pt%
\begingroup \smaller\smaller\smaller\begin{tabular}{@{}c@{}}%
-1\\-1190\\-628
\end{tabular}\endgroup%
\HardButStrongLineBreak\kern3pt%
\begingroup \smaller\smaller\smaller\begin{tabular}{@{}c@{}}%
-1\\-1188\\-630
\end{tabular}\endgroup%
\HardButStrongLineBreak\kern3pt%
\begingroup \smaller\smaller\smaller\begin{tabular}{@{}c@{}}%
3\\3570\\1883
\end{tabular}\endgroup%
\HardButStrongLineBreak\kern3pt%
\begingroup \smaller\smaller\smaller\begin{tabular}{@{}c@{}}%
23\\27360\\14445
\end{tabular}\endgroup%
\HardButStrongLineBreak\kern3pt%
\begingroup \smaller\smaller\smaller\begin{tabular}{@{}c@{}}%
23\\27356\\14448
\end{tabular}\endgroup%
\HardButStrongLineBreak\kern3pt%
\begingroup \smaller\smaller\smaller\begin{tabular}{@{}c@{}}%
6\\7135\\3770
\end{tabular}\endgroup%
\HardButStrongLineBreak\kern3pt%
\begingroup \smaller\smaller\smaller\begin{tabular}{@{}c@{}}%
85\\101052\\53424
\end{tabular}\endgroup%
{$\left.\llap{\phantom{%
\begingroup \smaller\smaller\smaller\begin{tabular}{@{}c@{}}%
0\\0\\0
\end{tabular}\endgroup%
}}\!\right]$}%

\medskip%
%
\leavevmode\llap{}%
$W_{256}$%
\qquad\llap{24} lattices, $\chi=48$%
\hfill%
$222222222222\rtimes C_{2}$%
\nopagebreak\smallskip\hrule\nopagebreak\medskip%
%
%
\leavevmode%
${L_{256.1}}$%
{} : {$1\above{1pt}{1pt}{-2}{{\rm II}}4\above{1pt}{1pt}{1}{7}{\cdot}1\above{1pt}{1pt}{2}{}9\above{1pt}{1pt}{1}{}{\cdot}1\above{1pt}{1pt}{-2}{}5\above{1pt}{1pt}{1}{}{\cdot}1\above{1pt}{1pt}{-2}{}7\above{1pt}{1pt}{-}{}$}\spacer%
\instructions{2}%
\EasyButWeakLineBreak%
{${2}\above{1pt}{1pt}{s}{2}{630}\above{1pt}{1pt}{b}{2}{4}\above{1pt}{1pt}{b}{2}{70}\above{1pt}{1pt}{b}{2}{36}\above{1pt}{1pt}{*}{2}{20}\above{1pt}{1pt}{b}{2}$}\relax$\,(\times2)$%
\nopagebreak\par%
\nopagebreak\par\leavevmode%
{$\left[\!\llap{\phantom{%
\begingroup \smaller\smaller\smaller
\endgroup%
}}\!\right]$}%

\medskip%
%
\leavevmode\llap{}%
$W_{257}$%
\qquad\llap{16} lattices, $\chi=24$%
\hfill%
$222|2222|2\rtimes D_{2}$%
\nopagebreak\smallskip\hrule\nopagebreak\medskip%
%
%
\leavevmode%
${L_{257.1}}$%
{} : {$[1\above{1pt}{1pt}{1}{}2\above{1pt}{1pt}{1}{}]\above{1pt}{1pt}{}{0}32\above{1pt}{1pt}{-}{5}{\cdot}1\above{1pt}{1pt}{-2}{}5\above{1pt}{1pt}{-}{}$}\EasyButWeakLineBreak%
{${160}\above{1pt}{1pt}{}{2}{1}\above{1pt}{1pt}{r}{2}{32}\above{1pt}{1pt}{*}{2}{40}\above{1pt}{1pt}{s}{2}{32}\above{1pt}{1pt}{*}{2}{4}\above{1pt}{1pt}{s}{2}{160}\above{1pt}{1pt}{l}{2}{2}\above{1pt}{1pt}{}{2}$}%
\nopagebreak\par%
\nopagebreak\par\leavevmode%
{$\left[\!\llap{\phantom{%
\begingroup \smaller\smaller\smaller\begin{tabular}{@{}c@{}}%
0\\0\\0
\end{tabular}\endgroup%
}}\right.$}%
\begingroup \smaller\smaller\smaller\begin{tabular}{@{}c@{}}%
-70240\\800\\800
\end{tabular}\endgroup%
\kern3pt%
\begingroup \smaller\smaller\smaller\begin{tabular}{@{}c@{}}%
800\\-2\\-10
\end{tabular}\endgroup%
\kern3pt%
\begingroup \smaller\smaller\smaller\begin{tabular}{@{}c@{}}%
800\\-10\\-9
\end{tabular}\endgroup%
{$\left.\llap{\phantom{%
\begingroup \smaller\smaller\smaller\begin{tabular}{@{}c@{}}%
0\\0\\0
\end{tabular}\endgroup%
}}\!\right]$}%
\EasyButWeakLineBreak%
{$\left[\!\llap{\phantom{%
\begingroup \smaller\smaller\smaller\begin{tabular}{@{}c@{}}%
0\\0\\0
\end{tabular}\endgroup%
}}\right.$}%
\begingroup \smaller\smaller\smaller\begin{tabular}{@{}c@{}}%
101\\960\\7840
\end{tabular}\endgroup%
\HardButStrongLineBreak\kern3pt%
\begingroup \smaller\smaller\smaller\begin{tabular}{@{}c@{}}%
7\\67\\543
\end{tabular}\endgroup%
\HardButStrongLineBreak\kern3pt%
\begingroup \smaller\smaller\smaller\begin{tabular}{@{}c@{}}%
19\\184\\1472
\end{tabular}\endgroup%
\HardButStrongLineBreak\kern3pt%
\begingroup \smaller\smaller\smaller\begin{tabular}{@{}c@{}}%
7\\70\\540
\end{tabular}\endgroup%
\HardButStrongLineBreak\kern3pt%
\begingroup \smaller\smaller\smaller\begin{tabular}{@{}c@{}}%
-1\\-8\\-80
\end{tabular}\endgroup%
\HardButStrongLineBreak\kern3pt%
\begingroup \smaller\smaller\smaller\begin{tabular}{@{}c@{}}%
-1\\-10\\-78
\end{tabular}\endgroup%
\HardButStrongLineBreak\kern3pt%
\begingroup \smaller\smaller\smaller\begin{tabular}{@{}c@{}}%
1\\0\\80
\end{tabular}\endgroup%
\HardButStrongLineBreak\kern3pt%
\begingroup \smaller\smaller\smaller\begin{tabular}{@{}c@{}}%
1\\9\\78
\end{tabular}\endgroup%
{$\left.\llap{\phantom{%
\begingroup \smaller\smaller\smaller\begin{tabular}{@{}c@{}}%
0\\0\\0
\end{tabular}\endgroup%
}}\!\right]$}%
%
%
\hbox{}\par\smallskip%
%
%
\leavevmode%
${L_{257.2}}$%
{} : {$[1\above{1pt}{1pt}{1}{}2\above{1pt}{1pt}{-}{}]\above{1pt}{1pt}{}{4}32\above{1pt}{1pt}{1}{1}{\cdot}1\above{1pt}{1pt}{-2}{}5\above{1pt}{1pt}{-}{}$}\EasyButWeakLineBreak%
{${160}\above{1pt}{1pt}{*}{2}{4}\above{1pt}{1pt}{s}{2}{32}\above{1pt}{1pt}{l}{2}{10}\above{1pt}{1pt}{}{2}{32}\above{1pt}{1pt}{}{2}{1}\above{1pt}{1pt}{r}{2}{160}\above{1pt}{1pt}{*}{2}{8}\above{1pt}{1pt}{s}{2}$}%
\nopagebreak\par%
\nopagebreak\par\leavevmode%
{$\left[\!\llap{\phantom{%
\begingroup \smaller\smaller\smaller\begin{tabular}{@{}c@{}}%
0\\0\\0
\end{tabular}\endgroup%
}}\right.$}%
\begingroup \smaller\smaller\smaller\begin{tabular}{@{}c@{}}%
67360\\-1600\\-160
\end{tabular}\endgroup%
\kern3pt%
\begingroup \smaller\smaller\smaller\begin{tabular}{@{}c@{}}%
-1600\\38\\4
\end{tabular}\endgroup%
\kern3pt%
\begingroup \smaller\smaller\smaller\begin{tabular}{@{}c@{}}%
-160\\4\\-7
\end{tabular}\endgroup%
{$\left.\llap{\phantom{%
\begingroup \smaller\smaller\smaller\begin{tabular}{@{}c@{}}%
0\\0\\0
\end{tabular}\endgroup%
}}\!\right]$}%
\EasyButWeakLineBreak%
{$\left[\!\llap{\phantom{%
\begingroup \smaller\smaller\smaller\begin{tabular}{@{}c@{}}%
0\\0\\0
\end{tabular}\endgroup%
}}\right.$}%
\begingroup \smaller\smaller\smaller\begin{tabular}{@{}c@{}}%
-101\\-4280\\-240
\end{tabular}\endgroup%
\HardButStrongLineBreak\kern3pt%
\begingroup \smaller\smaller\smaller\begin{tabular}{@{}c@{}}%
-13\\-552\\-34
\end{tabular}\endgroup%
\HardButStrongLineBreak\kern3pt%
\begingroup \smaller\smaller\smaller\begin{tabular}{@{}c@{}}%
-15\\-640\\-48
\end{tabular}\endgroup%
\HardButStrongLineBreak\kern3pt%
\begingroup \smaller\smaller\smaller\begin{tabular}{@{}c@{}}%
-1\\-45\\-10
\end{tabular}\endgroup%
\HardButStrongLineBreak\kern3pt%
\begingroup \smaller\smaller\smaller\begin{tabular}{@{}c@{}}%
5\\208\\0
\end{tabular}\endgroup%
\HardButStrongLineBreak\kern3pt%
\begingroup \smaller\smaller\smaller\begin{tabular}{@{}c@{}}%
1\\42\\1
\end{tabular}\endgroup%
\HardButStrongLineBreak\kern3pt%
\begingroup \smaller\smaller\smaller\begin{tabular}{@{}c@{}}%
-1\\-40\\0
\end{tabular}\endgroup%
\HardButStrongLineBreak\kern3pt%
\begingroup \smaller\smaller\smaller\begin{tabular}{@{}c@{}}%
-3\\-126\\-4
\end{tabular}\endgroup%
{$\left.\llap{\phantom{%
\begingroup \smaller\smaller\smaller\begin{tabular}{@{}c@{}}%
0\\0\\0
\end{tabular}\endgroup%
}}\!\right]$}%
%
%
\hbox{}\par\smallskip%
%
%
\leavevmode%
${L_{257.3}}$%
{} : {$1\above{1pt}{1pt}{-}{5}4\above{1pt}{1pt}{1}{1}32\above{1pt}{1pt}{1}{7}{\cdot}1\above{1pt}{1pt}{-2}{}5\above{1pt}{1pt}{1}{}$}\EasyButWeakLineBreak%
{${20}\above{1pt}{1pt}{*}{2}{32}\above{1pt}{1pt}{l}{2}{1}\above{1pt}{1pt}{}{2}{20}\above{1pt}{1pt}{r}{2}{4}\above{1pt}{1pt}{*}{2}{32}\above{1pt}{1pt}{l}{2}{5}\above{1pt}{1pt}{}{2}{4}\above{1pt}{1pt}{r}{2}$}%
\nopagebreak\par%
\nopagebreak\par\leavevmode%
{$\left[\!\llap{\phantom{%
\begingroup \smaller\smaller\smaller\begin{tabular}{@{}c@{}}%
0\\0\\0
\end{tabular}\endgroup%
}}\right.$}%
\begingroup \smaller\smaller\smaller\begin{tabular}{@{}c@{}}%
-219680\\1920\\2080
\end{tabular}\endgroup%
\kern3pt%
\begingroup \smaller\smaller\smaller\begin{tabular}{@{}c@{}}%
1920\\-12\\-20
\end{tabular}\endgroup%
\kern3pt%
\begingroup \smaller\smaller\smaller\begin{tabular}{@{}c@{}}%
2080\\-20\\-19
\end{tabular}\endgroup%
{$\left.\llap{\phantom{%
\begingroup \smaller\smaller\smaller\begin{tabular}{@{}c@{}}%
0\\0\\0
\end{tabular}\endgroup%
}}\!\right]$}%
\EasyButWeakLineBreak%
{$\left[\!\llap{\phantom{%
\begingroup \smaller\smaller\smaller\begin{tabular}{@{}c@{}}%
0\\0\\0
\end{tabular}\endgroup%
}}\right.$}%
\begingroup \smaller\smaller\smaller\begin{tabular}{@{}c@{}}%
29\\840\\2270
\end{tabular}\endgroup%
\HardButStrongLineBreak\kern3pt%
\begingroup \smaller\smaller\smaller\begin{tabular}{@{}c@{}}%
29\\836\\2272
\end{tabular}\endgroup%
\HardButStrongLineBreak\kern3pt%
\begingroup \smaller\smaller\smaller\begin{tabular}{@{}c@{}}%
2\\57\\157
\end{tabular}\endgroup%
\HardButStrongLineBreak\kern3pt%
\begingroup \smaller\smaller\smaller\begin{tabular}{@{}c@{}}%
1\\25\\80
\end{tabular}\endgroup%
\HardButStrongLineBreak\kern3pt%
\begingroup \smaller\smaller\smaller\begin{tabular}{@{}c@{}}%
-1\\-30\\-78
\end{tabular}\endgroup%
\HardButStrongLineBreak\kern3pt%
\begingroup \smaller\smaller\smaller\begin{tabular}{@{}c@{}}%
-1\\-28\\-80
\end{tabular}\endgroup%
\HardButStrongLineBreak\kern3pt%
\begingroup \smaller\smaller\smaller\begin{tabular}{@{}c@{}}%
2\\60\\155
\end{tabular}\endgroup%
\HardButStrongLineBreak\kern3pt%
\begingroup \smaller\smaller\smaller\begin{tabular}{@{}c@{}}%
2\\59\\156
\end{tabular}\endgroup%
{$\left.\llap{\phantom{%
\begingroup \smaller\smaller\smaller\begin{tabular}{@{}c@{}}%
0\\0\\0
\end{tabular}\endgroup%
}}\!\right]$}%
%
%
\hbox{}\par\smallskip%
%
%
\leavevmode%
${L_{257.4}}$%
{} : {$1\above{1pt}{1pt}{-}{5}4\above{1pt}{1pt}{1}{7}32\above{1pt}{1pt}{1}{1}{\cdot}1\above{1pt}{1pt}{-2}{}5\above{1pt}{1pt}{1}{}$}\EasyButWeakLineBreak%
{${5}\above{1pt}{1pt}{}{2}{32}\above{1pt}{1pt}{}{2}{1}\above{1pt}{1pt}{r}{2}{80}\above{1pt}{1pt}{*}{2}{4}\above{1pt}{1pt}{s}{2}{32}\above{1pt}{1pt}{s}{2}{20}\above{1pt}{1pt}{*}{2}{16}\above{1pt}{1pt}{l}{2}$}%
\nopagebreak\par%
\nopagebreak\par\leavevmode%
{$\left[\!\llap{\phantom{%
\begingroup \smaller\smaller\smaller\begin{tabular}{@{}c@{}}%
0\\0\\0
\end{tabular}\endgroup%
}}\right.$}%
\begingroup \smaller\smaller\smaller\begin{tabular}{@{}c@{}}%
-593120\\39840\\-2560
\end{tabular}\endgroup%
\kern3pt%
\begingroup \smaller\smaller\smaller\begin{tabular}{@{}c@{}}%
39840\\-2676\\172
\end{tabular}\endgroup%
\kern3pt%
\begingroup \smaller\smaller\smaller\begin{tabular}{@{}c@{}}%
-2560\\172\\-11
\end{tabular}\endgroup%
{$\left.\llap{\phantom{%
\begingroup \smaller\smaller\smaller\begin{tabular}{@{}c@{}}%
0\\0\\0
\end{tabular}\endgroup%
}}\!\right]$}%
\EasyButWeakLineBreak%
{$\left[\!\llap{\phantom{%
\begingroup \smaller\smaller\smaller\begin{tabular}{@{}c@{}}%
0\\0\\0
\end{tabular}\endgroup%
}}\right.$}%
\begingroup \smaller\smaller\smaller\begin{tabular}{@{}c@{}}%
-9\\-130\\55
\end{tabular}\endgroup%
\HardButStrongLineBreak\kern3pt%
\begingroup \smaller\smaller\smaller\begin{tabular}{@{}c@{}}%
-17\\-248\\64
\end{tabular}\endgroup%
\HardButStrongLineBreak\kern3pt%
\begingroup \smaller\smaller\smaller\begin{tabular}{@{}c@{}}%
-1\\-15\\-3
\end{tabular}\endgroup%
\HardButStrongLineBreak\kern3pt%
\begingroup \smaller\smaller\smaller\begin{tabular}{@{}c@{}}%
1\\10\\-80
\end{tabular}\endgroup%
\HardButStrongLineBreak\kern3pt%
\begingroup \smaller\smaller\smaller\begin{tabular}{@{}c@{}}%
1\\14\\-14
\end{tabular}\endgroup%
\HardButStrongLineBreak\kern3pt%
\begingroup \smaller\smaller\smaller\begin{tabular}{@{}c@{}}%
1\\16\\16
\end{tabular}\endgroup%
\HardButStrongLineBreak\kern3pt%
\begingroup \smaller\smaller\smaller\begin{tabular}{@{}c@{}}%
-3\\-40\\70
\end{tabular}\endgroup%
\HardButStrongLineBreak\kern3pt%
\begingroup \smaller\smaller\smaller\begin{tabular}{@{}c@{}}%
-3\\-42\\40
\end{tabular}\endgroup%
{$\left.\llap{\phantom{%
\begingroup \smaller\smaller\smaller\begin{tabular}{@{}c@{}}%
0\\0\\0
\end{tabular}\endgroup%
}}\!\right]$}%

\medskip%
%
\leavevmode\llap{}%
$W_{258}$%
\qquad\llap{8} lattices, $\chi=36$%
\hfill%
$24|42|24|42|\rtimes D_{4}$%
\nopagebreak\smallskip\hrule\nopagebreak\medskip%
%
%
\leavevmode%
${L_{258.1}}$%
{} : {$1\above{1pt}{1pt}{2}{2}16\above{1pt}{1pt}{-}{3}{\cdot}1\above{1pt}{1pt}{2}{}3\above{1pt}{1pt}{-}{}{\cdot}1\above{1pt}{1pt}{2}{}7\above{1pt}{1pt}{1}{}$}\EasyButWeakLineBreak%
{${112}\above{1pt}{1pt}{l}{2}{1}\above{1pt}{1pt}{}{4}{2}\above{1pt}{1pt}{*}{4}{4}\above{1pt}{1pt}{s}{2}$}\relax$\,(\times2)$%
\nopagebreak\par%
\nopagebreak\par\leavevmode%
{$\left[\!\llap{\phantom{%
\begingroup \smaller\smaller\smaller\begin{tabular}{@{}c@{}}%
0\\0\\0
\end{tabular}\endgroup%
}}\right.$}%
\begingroup \smaller\smaller\smaller\begin{tabular}{@{}c@{}}%
-2258256\\557424\\22176
\end{tabular}\endgroup%
\kern3pt%
\begingroup \smaller\smaller\smaller\begin{tabular}{@{}c@{}}%
557424\\-137591\\-5477
\end{tabular}\endgroup%
\kern3pt%
\begingroup \smaller\smaller\smaller\begin{tabular}{@{}c@{}}%
22176\\-5477\\-214
\end{tabular}\endgroup%
{$\left.\llap{\phantom{%
\begingroup \smaller\smaller\smaller\begin{tabular}{@{}c@{}}%
0\\0\\0
\end{tabular}\endgroup%
}}\!\right]$}%
\hfil\penalty500%
{$\left[\!\llap{\phantom{%
\begingroup \smaller\smaller\smaller\begin{tabular}{@{}c@{}}%
0\\0\\0
\end{tabular}\endgroup%
}}\right.$}%
\begingroup \smaller\smaller\smaller\begin{tabular}{@{}c@{}}%
742937\\2914128\\2403744
\end{tabular}\endgroup%
\kern3pt%
\begingroup \smaller\smaller\smaller\begin{tabular}{@{}c@{}}%
-183749\\-720745\\-594512
\end{tabular}\endgroup%
\kern3pt%
\begingroup \smaller\smaller\smaller\begin{tabular}{@{}c@{}}%
-6859\\-26904\\-22193
\end{tabular}\endgroup%
{$\left.\llap{\phantom{%
\begingroup \smaller\smaller\smaller\begin{tabular}{@{}c@{}}%
0\\0\\0
\end{tabular}\endgroup%
}}\!\right]$}%
\EasyButWeakLineBreak%
{$\left[\!\llap{\phantom{%
\begingroup \smaller\smaller\smaller\begin{tabular}{@{}c@{}}%
0\\0\\0
\end{tabular}\endgroup%
}}\right.$}%
\begingroup \smaller\smaller\smaller\begin{tabular}{@{}c@{}}%
5211\\20440\\16856
\end{tabular}\endgroup%
\HardButStrongLineBreak\kern3pt%
\begingroup \smaller\smaller\smaller\begin{tabular}{@{}c@{}}%
297\\1165\\960
\end{tabular}\endgroup%
\HardButStrongLineBreak\kern3pt%
\begingroup \smaller\smaller\smaller\begin{tabular}{@{}c@{}}%
78\\306\\251
\end{tabular}\endgroup%
\HardButStrongLineBreak\kern3pt%
\begingroup \smaller\smaller\smaller\begin{tabular}{@{}c@{}}%
-77\\-302\\-250
\end{tabular}\endgroup%
{$\left.\llap{\phantom{%
\begingroup \smaller\smaller\smaller\begin{tabular}{@{}c@{}}%
0\\0\\0
\end{tabular}\endgroup%
}}\!\right]$}%

\medskip%
%
\leavevmode\llap{}%
$W_{259}$%
\qquad\llap{8} lattices, $\chi=16$%
\hfill%
$222226$%
\nopagebreak\smallskip\hrule\nopagebreak\medskip%
%
%
\leavevmode%
${L_{259.1}}$%
{} : {$1\above{1pt}{1pt}{-2}{{\rm II}}16\above{1pt}{1pt}{1}{1}{\cdot}1\above{1pt}{1pt}{2}{}3\above{1pt}{1pt}{-}{}{\cdot}1\above{1pt}{1pt}{-2}{}7\above{1pt}{1pt}{-}{}$}\EasyButWeakLineBreak%
{${6}\above{1pt}{1pt}{b}{2}{16}\above{1pt}{1pt}{b}{2}{42}\above{1pt}{1pt}{l}{2}{16}\above{1pt}{1pt}{r}{2}{6}\above{1pt}{1pt}{b}{2}{2}\above{1pt}{1pt}{}{6}$}%
\nopagebreak\par%
\nopagebreak\par\leavevmode%
{$\left[\!\llap{\phantom{%
\begingroup \smaller\smaller\smaller\begin{tabular}{@{}c@{}}%
0\\0\\0
\end{tabular}\endgroup%
}}\right.$}%
\begingroup \smaller\smaller\smaller\begin{tabular}{@{}c@{}}%
2549904\\-6384\\-3696
\end{tabular}\endgroup%
\kern3pt%
\begingroup \smaller\smaller\smaller\begin{tabular}{@{}c@{}}%
-6384\\-2\\33
\end{tabular}\endgroup%
\kern3pt%
\begingroup \smaller\smaller\smaller\begin{tabular}{@{}c@{}}%
-3696\\33\\-26
\end{tabular}\endgroup%
{$\left.\llap{\phantom{%
\begingroup \smaller\smaller\smaller\begin{tabular}{@{}c@{}}%
0\\0\\0
\end{tabular}\endgroup%
}}\!\right]$}%
\EasyButWeakLineBreak%
{$\left[\!\llap{\phantom{%
\begingroup \smaller\smaller\smaller\begin{tabular}{@{}c@{}}%
0\\0\\0
\end{tabular}\endgroup%
}}\right.$}%
\begingroup \smaller\smaller\smaller\begin{tabular}{@{}c@{}}%
2\\555\\420
\end{tabular}\endgroup%
\HardButStrongLineBreak\kern3pt%
\begingroup \smaller\smaller\smaller\begin{tabular}{@{}c@{}}%
7\\1944\\1472
\end{tabular}\endgroup%
\HardButStrongLineBreak\kern3pt%
\begingroup \smaller\smaller\smaller\begin{tabular}{@{}c@{}}%
-1\\-273\\-210
\end{tabular}\endgroup%
\HardButStrongLineBreak\kern3pt%
\begingroup \smaller\smaller\smaller\begin{tabular}{@{}c@{}}%
-21\\-5824\\-4416
\end{tabular}\endgroup%
\HardButStrongLineBreak\kern3pt%
\begingroup \smaller\smaller\smaller\begin{tabular}{@{}c@{}}%
-10\\-2775\\-2103
\end{tabular}\endgroup%
\HardButStrongLineBreak\kern3pt%
\begingroup \smaller\smaller\smaller\begin{tabular}{@{}c@{}}%
-3\\-833\\-631
\end{tabular}\endgroup%
{$\left.\llap{\phantom{%
\begingroup \smaller\smaller\smaller\begin{tabular}{@{}c@{}}%
0\\0\\0
\end{tabular}\endgroup%
}}\!\right]$}%

\medskip%
%
\leavevmode\llap{}%
$W_{260}$%
\qquad\llap{8} lattices, $\chi=18$%
\hfill%
$22|222\slashtwo2\rtimes D_{2}$%
\nopagebreak\smallskip\hrule\nopagebreak\medskip%
%
%
\leavevmode%
${L_{260.1}}$%
{} : {$1\above{1pt}{1pt}{-2}{6}16\above{1pt}{1pt}{1}{7}{\cdot}1\above{1pt}{1pt}{2}{}3\above{1pt}{1pt}{-}{}{\cdot}1\above{1pt}{1pt}{2}{}7\above{1pt}{1pt}{1}{}$}\EasyButWeakLineBreak%
{${112}\above{1pt}{1pt}{b}{2}{6}\above{1pt}{1pt}{s}{2}{14}\above{1pt}{1pt}{b}{2}{6}\above{1pt}{1pt}{l}{2}{112}\above{1pt}{1pt}{}{2}{1}\above{1pt}{1pt}{r}{2}{4}\above{1pt}{1pt}{*}{2}$}%
\nopagebreak\par%
\nopagebreak\par\leavevmode%
{$\left[\!\llap{\phantom{%
\begingroup \smaller\smaller\smaller\begin{tabular}{@{}c@{}}%
0\\0\\0
\end{tabular}\endgroup%
}}\right.$}%
\begingroup \smaller\smaller\smaller\begin{tabular}{@{}c@{}}%
76272\\25200\\-336
\end{tabular}\endgroup%
\kern3pt%
\begingroup \smaller\smaller\smaller\begin{tabular}{@{}c@{}}%
25200\\8326\\-111
\end{tabular}\endgroup%
\kern3pt%
\begingroup \smaller\smaller\smaller\begin{tabular}{@{}c@{}}%
-336\\-111\\1
\end{tabular}\endgroup%
{$\left.\llap{\phantom{%
\begingroup \smaller\smaller\smaller\begin{tabular}{@{}c@{}}%
0\\0\\0
\end{tabular}\endgroup%
}}\!\right]$}%
\EasyButWeakLineBreak%
{$\left[\!\llap{\phantom{%
\begingroup \smaller\smaller\smaller\begin{tabular}{@{}c@{}}%
0\\0\\0
\end{tabular}\endgroup%
}}\right.$}%
\begingroup \smaller\smaller\smaller\begin{tabular}{@{}c@{}}%
831\\-2520\\-448
\end{tabular}\endgroup%
\HardButStrongLineBreak\kern3pt%
\begingroup \smaller\smaller\smaller\begin{tabular}{@{}c@{}}%
92\\-279\\-48
\end{tabular}\endgroup%
\HardButStrongLineBreak\kern3pt%
\begingroup \smaller\smaller\smaller\begin{tabular}{@{}c@{}}%
30\\-91\\-14
\end{tabular}\endgroup%
\HardButStrongLineBreak\kern3pt%
\begingroup \smaller\smaller\smaller\begin{tabular}{@{}c@{}}%
-1\\3\\0
\end{tabular}\endgroup%
\HardButStrongLineBreak\kern3pt%
\begingroup \smaller\smaller\smaller\begin{tabular}{@{}c@{}}%
-37\\112\\0
\end{tabular}\endgroup%
\HardButStrongLineBreak\kern3pt%
\begingroup \smaller\smaller\smaller\begin{tabular}{@{}c@{}}%
0\\0\\-1
\end{tabular}\endgroup%
\HardButStrongLineBreak\kern3pt%
\begingroup \smaller\smaller\smaller\begin{tabular}{@{}c@{}}%
31\\-94\\-18
\end{tabular}\endgroup%
{$\left.\llap{\phantom{%
\begingroup \smaller\smaller\smaller\begin{tabular}{@{}c@{}}%
0\\0\\0
\end{tabular}\endgroup%
}}\!\right]$}%

\medskip%
%
\leavevmode\llap{}%
$W_{261}$%
\qquad\llap{16} lattices, $\chi=72$%
\hfill%
$2222|2222|2222|2222|\rtimes D_{4}$%
\nopagebreak\smallskip\hrule\nopagebreak\medskip%
%
%
\leavevmode%
${L_{261.1}}$%
{} : {$1\above{1pt}{1pt}{2}{6}16\above{1pt}{1pt}{-}{3}{\cdot}1\above{1pt}{1pt}{1}{}3\above{1pt}{1pt}{1}{}9\above{1pt}{1pt}{-}{}{\cdot}1\above{1pt}{1pt}{2}{}7\above{1pt}{1pt}{1}{}$}\spacer%
\instructions{3}%
\EasyButWeakLineBreak%
{${112}\above{1pt}{1pt}{s}{2}{12}\above{1pt}{1pt}{l}{2}{7}\above{1pt}{1pt}{}{2}{48}\above{1pt}{1pt}{r}{2}{126}\above{1pt}{1pt}{b}{2}{48}\above{1pt}{1pt}{*}{2}{28}\above{1pt}{1pt}{l}{2}{3}\above{1pt}{1pt}{r}{2}$}\relax$\,(\times2)$%
\nopagebreak\par%
\nopagebreak\par\leavevmode%
{$\left[\!\llap{\phantom{%
\begingroup \smaller\smaller\smaller
\endgroup%
}}\!\right]$}%

\medskip%
%
\leavevmode\llap{}%
$W_{262}$%
\qquad\llap{16} lattices, $\chi=72$%
\hfill%
$222|2222|2222|2222|2\rtimes D_{4}$%
\nopagebreak\smallskip\hrule\nopagebreak\medskip%
%
%
\leavevmode%
${L_{262.1}}$%
{} : {$1\above{1pt}{1pt}{-2}{2}16\above{1pt}{1pt}{1}{7}{\cdot}1\above{1pt}{1pt}{1}{}3\above{1pt}{1pt}{1}{}9\above{1pt}{1pt}{-}{}{\cdot}1\above{1pt}{1pt}{2}{}7\above{1pt}{1pt}{1}{}$}\spacer%
\instructions{3}%
\EasyButWeakLineBreak%
{${112}\above{1pt}{1pt}{}{2}{3}\above{1pt}{1pt}{}{2}{7}\above{1pt}{1pt}{r}{2}{48}\above{1pt}{1pt}{s}{2}{28}\above{1pt}{1pt}{*}{2}{12}\above{1pt}{1pt}{*}{2}{112}\above{1pt}{1pt}{b}{2}{18}\above{1pt}{1pt}{l}{2}$}\relax$\,(\times2)$%
\nopagebreak\par%
\nopagebreak\par\leavevmode%
{$\left[\!\llap{\phantom{%
\begingroup \smaller\smaller\smaller
\endgroup%
}}\!\right]$}%

\medskip%
%
\leavevmode\llap{}%
$W_{263}$%
\qquad\llap{22} lattices, $\chi=36$%
\hfill%
$2\infty\infty\infty222$%
\nopagebreak\smallskip\hrule\nopagebreak\medskip%
%
%
\leavevmode%
${L_{263.1}}$%
{} : {$1\above{1pt}{1pt}{2}{{\rm II}}4\above{1pt}{1pt}{1}{7}{\cdot}1\above{1pt}{1pt}{-}{}7\above{1pt}{1pt}{1}{}49\above{1pt}{1pt}{1}{}$}\spacer%
\instructions{2}%
\EasyButWeakLineBreak%
{${98}\above{1pt}{1pt}{l}{2}{28}\above{1pt}{1pt}{7,1}{\infty}{28}\above{1pt}{1pt}{7,4}{\infty b}{28}\above{1pt}{1pt}{7,4}{\infty}{28}\above{1pt}{1pt}{*}{2}{196}\above{1pt}{1pt}{b}{2}{14}\above{1pt}{1pt}{s}{2}$}%
\nopagebreak\par%
\nopagebreak\par\leavevmode%
{$\left[\!\llap{\phantom{%
\begingroup \smaller\smaller\smaller\begin{tabular}{@{}c@{}}%
0\\0\\0
\end{tabular}\endgroup%
}}\right.$}%
\begingroup \smaller\smaller\smaller\begin{tabular}{@{}c@{}}%
2940\\1568\\-588
\end{tabular}\endgroup%
\kern3pt%
\begingroup \smaller\smaller\smaller\begin{tabular}{@{}c@{}}%
1568\\840\\-315
\end{tabular}\endgroup%
\kern3pt%
\begingroup \smaller\smaller\smaller\begin{tabular}{@{}c@{}}%
-588\\-315\\118
\end{tabular}\endgroup%
{$\left.\llap{\phantom{%
\begingroup \smaller\smaller\smaller\begin{tabular}{@{}c@{}}%
0\\0\\0
\end{tabular}\endgroup%
}}\!\right]$}%
\EasyButWeakLineBreak%
{$\left[\!\llap{\phantom{%
\begingroup \smaller\smaller\smaller\begin{tabular}{@{}c@{}}%
0\\0\\0
\end{tabular}\endgroup%
}}\right.$}%
\begingroup \smaller\smaller\smaller\begin{tabular}{@{}c@{}}%
-5\\28\\49
\end{tabular}\endgroup%
\HardButStrongLineBreak\kern3pt%
\begingroup \smaller\smaller\smaller\begin{tabular}{@{}c@{}}%
-3\\16\\28
\end{tabular}\endgroup%
\HardButStrongLineBreak\kern3pt%
\begingroup \smaller\smaller\smaller\begin{tabular}{@{}c@{}}%
1\\-2\\0
\end{tabular}\endgroup%
\HardButStrongLineBreak\kern3pt%
\begingroup \smaller\smaller\smaller\begin{tabular}{@{}c@{}}%
3\\16\\56
\end{tabular}\endgroup%
\HardButStrongLineBreak\kern3pt%
\begingroup \smaller\smaller\smaller\begin{tabular}{@{}c@{}}%
3\\54\\154
\end{tabular}\endgroup%
\HardButStrongLineBreak\kern3pt%
\begingroup \smaller\smaller\smaller\begin{tabular}{@{}c@{}}%
1\\112\\294
\end{tabular}\endgroup%
\HardButStrongLineBreak\kern3pt%
\begingroup \smaller\smaller\smaller\begin{tabular}{@{}c@{}}%
-1\\10\\21
\end{tabular}\endgroup%
{$\left.\llap{\phantom{%
\begingroup \smaller\smaller\smaller\begin{tabular}{@{}c@{}}%
0\\0\\0
\end{tabular}\endgroup%
}}\!\right]$}%
%
%
\hbox{}\par\smallskip%
%
%
\leavevmode%
${L_{263.2}}$%
{} : {$1\above{1pt}{1pt}{2}{6}8\above{1pt}{1pt}{1}{1}{\cdot}1\above{1pt}{1pt}{-}{}7\above{1pt}{1pt}{1}{}49\above{1pt}{1pt}{1}{}$}\spacer%
\instructions{2}%
\EasyButWeakLineBreak%
{${196}\above{1pt}{1pt}{s}{2}{56}\above{1pt}{1pt}{28,1}{\infty z}{14}\above{1pt}{1pt}{28,25}{\infty a}{56}\above{1pt}{1pt}{28,25}{\infty z}{14}\above{1pt}{1pt}{l}{2}{392}\above{1pt}{1pt}{}{2}{7}\above{1pt}{1pt}{r}{2}$}%
\nopagebreak\par%
\nopagebreak\par\leavevmode%
{$\left[\!\llap{\phantom{%
\begingroup \smaller\smaller\smaller\begin{tabular}{@{}c@{}}%
0\\0\\0
\end{tabular}\endgroup%
}}\right.$}%
\begingroup \smaller\smaller\smaller\begin{tabular}{@{}c@{}}%
-526456\\-258328\\-68992
\end{tabular}\endgroup%
\kern3pt%
\begingroup \smaller\smaller\smaller\begin{tabular}{@{}c@{}}%
-258328\\-126749\\-33845
\end{tabular}\endgroup%
\kern3pt%
\begingroup \smaller\smaller\smaller\begin{tabular}{@{}c@{}}%
-68992\\-33845\\-9034
\end{tabular}\endgroup%
{$\left.\llap{\phantom{%
\begingroup \smaller\smaller\smaller\begin{tabular}{@{}c@{}}%
0\\0\\0
\end{tabular}\endgroup%
}}\!\right]$}%
\EasyButWeakLineBreak%
{$\left[\!\llap{\phantom{%
\begingroup \smaller\smaller\smaller\begin{tabular}{@{}c@{}}%
0\\0\\0
\end{tabular}\endgroup%
}}\right.$}%
\begingroup \smaller\smaller\smaller\begin{tabular}{@{}c@{}}%
-135\\406\\-490
\end{tabular}\endgroup%
\HardButStrongLineBreak\kern3pt%
\begingroup \smaller\smaller\smaller\begin{tabular}{@{}c@{}}%
-71\\212\\-252
\end{tabular}\endgroup%
\HardButStrongLineBreak\kern3pt%
\begingroup \smaller\smaller\smaller\begin{tabular}{@{}c@{}}%
68\\-204\\245
\end{tabular}\endgroup%
\HardButStrongLineBreak\kern3pt%
\begingroup \smaller\smaller\smaller\begin{tabular}{@{}c@{}}%
903\\-2700\\3220
\end{tabular}\endgroup%
\HardButStrongLineBreak\kern3pt%
\begingroup \smaller\smaller\smaller\begin{tabular}{@{}c@{}}%
907\\-2710\\3227
\end{tabular}\endgroup%
\HardButStrongLineBreak\kern3pt%
\begingroup \smaller\smaller\smaller\begin{tabular}{@{}c@{}}%
2869\\-8568\\10192
\end{tabular}\endgroup%
\HardButStrongLineBreak\kern3pt%
\begingroup \smaller\smaller\smaller\begin{tabular}{@{}c@{}}%
36\\-107\\126
\end{tabular}\endgroup%
{$\left.\llap{\phantom{%
\begingroup \smaller\smaller\smaller\begin{tabular}{@{}c@{}}%
0\\0\\0
\end{tabular}\endgroup%
}}\!\right]$}%
%
%
\hbox{}\par\smallskip%
%
%
\leavevmode%
${L_{263.3}}$%
{} : {$1\above{1pt}{1pt}{-2}{6}8\above{1pt}{1pt}{-}{5}{\cdot}1\above{1pt}{1pt}{-}{}7\above{1pt}{1pt}{1}{}49\above{1pt}{1pt}{1}{}$}\spacer%
\instructions{m}%
\EasyButWeakLineBreak%
{${49}\above{1pt}{1pt}{r}{2}{56}\above{1pt}{1pt}{28,15}{\infty z}{14}\above{1pt}{1pt}{28,25}{\infty b}{56}\above{1pt}{1pt}{28,11}{\infty z}{14}\above{1pt}{1pt}{b}{2}{392}\above{1pt}{1pt}{*}{2}{28}\above{1pt}{1pt}{l}{2}$}%
\nopagebreak\par%
\nopagebreak\par\leavevmode%
{$\left[\!\llap{\phantom{%
\begingroup \smaller\smaller\smaller\begin{tabular}{@{}c@{}}%
0\\0\\0
\end{tabular}\endgroup%
}}\right.$}%
\begingroup \smaller\smaller\smaller\begin{tabular}{@{}c@{}}%
-756952\\-80752\\-25872
\end{tabular}\endgroup%
\kern3pt%
\begingroup \smaller\smaller\smaller\begin{tabular}{@{}c@{}}%
-80752\\-8610\\-2765
\end{tabular}\endgroup%
\kern3pt%
\begingroup \smaller\smaller\smaller\begin{tabular}{@{}c@{}}%
-25872\\-2765\\-879
\end{tabular}\endgroup%
{$\left.\llap{\phantom{%
\begingroup \smaller\smaller\smaller\begin{tabular}{@{}c@{}}%
0\\0\\0
\end{tabular}\endgroup%
}}\!\right]$}%
\EasyButWeakLineBreak%
{$\left[\!\llap{\phantom{%
\begingroup \smaller\smaller\smaller\begin{tabular}{@{}c@{}}%
0\\0\\0
\end{tabular}\endgroup%
}}\right.$}%
\begingroup \smaller\smaller\smaller\begin{tabular}{@{}c@{}}%
-36\\259\\245
\end{tabular}\endgroup%
\HardButStrongLineBreak\kern3pt%
\begingroup \smaller\smaller\smaller\begin{tabular}{@{}c@{}}%
-21\\152\\140
\end{tabular}\endgroup%
\HardButStrongLineBreak\kern3pt%
\begingroup \smaller\smaller\smaller\begin{tabular}{@{}c@{}}%
29\\-209\\-196
\end{tabular}\endgroup%
\HardButStrongLineBreak\kern3pt%
\begingroup \smaller\smaller\smaller\begin{tabular}{@{}c@{}}%
281\\-2032\\-1876
\end{tabular}\endgroup%
\HardButStrongLineBreak\kern3pt%
\begingroup \smaller\smaller\smaller\begin{tabular}{@{}c@{}}%
259\\-1875\\-1722
\end{tabular}\endgroup%
\HardButStrongLineBreak\kern3pt%
\begingroup \smaller\smaller\smaller\begin{tabular}{@{}c@{}}%
769\\-5572\\-5096
\end{tabular}\endgroup%
\HardButStrongLineBreak\kern3pt%
\begingroup \smaller\smaller\smaller\begin{tabular}{@{}c@{}}%
7\\-52\\-42
\end{tabular}\endgroup%
{$\left.\llap{\phantom{%
\begingroup \smaller\smaller\smaller\begin{tabular}{@{}c@{}}%
0\\0\\0
\end{tabular}\endgroup%
}}\!\right]$}%

\medskip%
%
\leavevmode\llap{}%
$W_{264}$%
\qquad\llap{3} lattices, $\chi=24$%
\hfill%
$\slashtwo2\slashtwo2\slashtwo2\slashtwo2\rtimes D_{4}$%
\nopagebreak\smallskip\hrule\nopagebreak\medskip%
%
%
\leavevmode%
${L_{264.1}}$%
{} : {$1\above{1pt}{1pt}{-2}{{\rm II}}4\above{1pt}{1pt}{-}{3}{\cdot}1\above{1pt}{1pt}{1}{}7\above{1pt}{1pt}{-}{}49\above{1pt}{1pt}{1}{}$}\spacer%
\instructions{2}%
\EasyButWeakLineBreak%
{${196}\above{1pt}{1pt}{*}{2}{4}\above{1pt}{1pt}{b}{2}{98}\above{1pt}{1pt}{b}{2}{2}\above{1pt}{1pt}{b}{2}$}\relax$\,(\times2)$%
\nopagebreak\par%
\nopagebreak\par\leavevmode%
{$\left[\!\llap{\phantom{%
\begingroup \smaller\smaller\smaller\begin{tabular}{@{}c@{}}%
0\\0\\0
\end{tabular}\endgroup%
}}\right.$}%
\begingroup \smaller\smaller\smaller\begin{tabular}{@{}c@{}}%
-140532\\2352\\980
\end{tabular}\endgroup%
\kern3pt%
\begingroup \smaller\smaller\smaller\begin{tabular}{@{}c@{}}%
2352\\-14\\-21
\end{tabular}\endgroup%
\kern3pt%
\begingroup \smaller\smaller\smaller\begin{tabular}{@{}c@{}}%
980\\-21\\-6
\end{tabular}\endgroup%
{$\left.\llap{\phantom{%
\begingroup \smaller\smaller\smaller\begin{tabular}{@{}c@{}}%
0\\0\\0
\end{tabular}\endgroup%
}}\!\right]$}%
\hfil\penalty500%
{$\left[\!\llap{\phantom{%
\begingroup \smaller\smaller\smaller\begin{tabular}{@{}c@{}}%
0\\0\\0
\end{tabular}\endgroup%
}}\right.$}%
\begingroup \smaller\smaller\smaller\begin{tabular}{@{}c@{}}%
-1177\\-21560\\-117992
\end{tabular}\endgroup%
\kern3pt%
\begingroup \smaller\smaller\smaller\begin{tabular}{@{}c@{}}%
15\\274\\1505
\end{tabular}\endgroup%
\kern3pt%
\begingroup \smaller\smaller\smaller\begin{tabular}{@{}c@{}}%
9\\165\\902
\end{tabular}\endgroup%
{$\left.\llap{\phantom{%
\begingroup \smaller\smaller\smaller\begin{tabular}{@{}c@{}}%
0\\0\\0
\end{tabular}\endgroup%
}}\!\right]$}%
\EasyButWeakLineBreak%
{$\left[\!\llap{\phantom{%
\begingroup \smaller\smaller\smaller\begin{tabular}{@{}c@{}}%
0\\0\\0
\end{tabular}\endgroup%
}}\right.$}%
\begingroup \smaller\smaller\smaller\begin{tabular}{@{}c@{}}%
-45\\-826\\-4508
\end{tabular}\endgroup%
\HardButStrongLineBreak\kern3pt%
\begingroup \smaller\smaller\smaller\begin{tabular}{@{}c@{}}%
-7\\-128\\-702
\end{tabular}\endgroup%
\HardButStrongLineBreak\kern3pt%
\begingroup \smaller\smaller\smaller\begin{tabular}{@{}c@{}}%
-20\\-364\\-2009
\end{tabular}\endgroup%
\HardButStrongLineBreak\kern3pt%
\begingroup \smaller\smaller\smaller\begin{tabular}{@{}c@{}}%
-1\\-18\\-101
\end{tabular}\endgroup%
{$\left.\llap{\phantom{%
\begingroup \smaller\smaller\smaller\begin{tabular}{@{}c@{}}%
0\\0\\0
\end{tabular}\endgroup%
}}\!\right]$}%
%
%
%
%
%
%
%
%
%
%

\medskip%
%
\leavevmode\llap{}%
$W_{265}$%
\qquad\llap{6} lattices, $\chi=28$%
\hfill%
$\slashtwo6|6\slashtwo6|6\rtimes D_{4}$%
\nopagebreak\smallskip\hrule\nopagebreak\medskip%
%
%
\leavevmode%
${L_{265.1}}$%
{} : {$1\above{1pt}{1pt}{-2}{{\rm II}}4\above{1pt}{1pt}{-}{3}{\cdot}1\above{1pt}{1pt}{-}{}3\above{1pt}{1pt}{-}{}9\above{1pt}{1pt}{-}{}{\cdot}1\above{1pt}{1pt}{2}{}13\above{1pt}{1pt}{1}{}$}\spacer%
\instructions{2}%
\EasyButWeakLineBreak%
{${18}\above{1pt}{1pt}{b}{2}{2}\above{1pt}{1pt}{}{6}{6}\above{1pt}{1pt}{}{6}$}\relax$\,(\times2)$%
\nopagebreak\par%
\nopagebreak\par\leavevmode%
{$\left[\!\llap{\phantom{%
\begingroup \smaller\smaller\smaller\begin{tabular}{@{}c@{}}%
0\\0\\0
\end{tabular}\endgroup%
}}\right.$}%
\begingroup \smaller\smaller\smaller\begin{tabular}{@{}c@{}}%
-2291796\\17784\\8892
\end{tabular}\endgroup%
\kern3pt%
\begingroup \smaller\smaller\smaller\begin{tabular}{@{}c@{}}%
17784\\-138\\-69
\end{tabular}\endgroup%
\kern3pt%
\begingroup \smaller\smaller\smaller\begin{tabular}{@{}c@{}}%
8892\\-69\\-34
\end{tabular}\endgroup%
{$\left.\llap{\phantom{%
\begingroup \smaller\smaller\smaller\begin{tabular}{@{}c@{}}%
0\\0\\0
\end{tabular}\endgroup%
}}\!\right]$}%
\hfil\penalty500%
{$\left[\!\llap{\phantom{%
\begingroup \smaller\smaller\smaller\begin{tabular}{@{}c@{}}%
0\\0\\0
\end{tabular}\endgroup%
}}\right.$}%
\begingroup \smaller\smaller\smaller\begin{tabular}{@{}c@{}}%
2573\\344916\\-30888
\end{tabular}\endgroup%
\kern3pt%
\begingroup \smaller\smaller\smaller\begin{tabular}{@{}c@{}}%
-20\\-2681\\240
\end{tabular}\endgroup%
\kern3pt%
\begingroup \smaller\smaller\smaller\begin{tabular}{@{}c@{}}%
-9\\-1206\\107
\end{tabular}\endgroup%
{$\left.\llap{\phantom{%
\begingroup \smaller\smaller\smaller\begin{tabular}{@{}c@{}}%
0\\0\\0
\end{tabular}\endgroup%
}}\!\right]$}%
\EasyButWeakLineBreak%
{$\left[\!\llap{\phantom{%
\begingroup \smaller\smaller\smaller\begin{tabular}{@{}c@{}}%
0\\0\\0
\end{tabular}\endgroup%
}}\right.$}%
\begingroup \smaller\smaller\smaller\begin{tabular}{@{}c@{}}%
-1\\-129\\0
\end{tabular}\endgroup%
\HardButStrongLineBreak\kern3pt%
\begingroup \smaller\smaller\smaller\begin{tabular}{@{}c@{}}%
0\\-1\\2
\end{tabular}\endgroup%
\HardButStrongLineBreak\kern3pt%
\begingroup \smaller\smaller\smaller\begin{tabular}{@{}c@{}}%
1\\130\\-3
\end{tabular}\endgroup%
{$\left.\llap{\phantom{%
\begingroup \smaller\smaller\smaller\begin{tabular}{@{}c@{}}%
0\\0\\0
\end{tabular}\endgroup%
}}\!\right]$}%

\medskip%
%
\leavevmode\llap{}%
$W_{266}$%
\qquad\llap{6} lattices, $\chi=14$%
\hfill%
$\slashtwo22\slashthree22\rtimes D_{2}$%
\nopagebreak\smallskip\hrule\nopagebreak\medskip%
%
%
\leavevmode%
${L_{266.1}}$%
{} : {$1\above{1pt}{1pt}{-2}{{\rm II}}4\above{1pt}{1pt}{-}{3}{\cdot}1\above{1pt}{1pt}{1}{}3\above{1pt}{1pt}{-}{}9\above{1pt}{1pt}{1}{}{\cdot}1\above{1pt}{1pt}{2}{}13\above{1pt}{1pt}{1}{}$}\spacer%
\instructions{2}%
\EasyButWeakLineBreak%
{${36}\above{1pt}{1pt}{*}{2}{4}\above{1pt}{1pt}{*}{2}{468}\above{1pt}{1pt}{b}{2}{6}\above{1pt}{1pt}{-}{3}{6}\above{1pt}{1pt}{b}{2}{52}\above{1pt}{1pt}{*}{2}$}%
\nopagebreak\par%
\nopagebreak\par\leavevmode%
{$\left[\!\llap{\phantom{%
\begingroup \smaller\smaller\smaller\begin{tabular}{@{}c@{}}%
0\\0\\0
\end{tabular}\endgroup%
}}\right.$}%
\begingroup \smaller\smaller\smaller\begin{tabular}{@{}c@{}}%
-2535156\\9828\\9828
\end{tabular}\endgroup%
\kern3pt%
\begingroup \smaller\smaller\smaller\begin{tabular}{@{}c@{}}%
9828\\-30\\-39
\end{tabular}\endgroup%
\kern3pt%
\begingroup \smaller\smaller\smaller\begin{tabular}{@{}c@{}}%
9828\\-39\\-38
\end{tabular}\endgroup%
{$\left.\llap{\phantom{%
\begingroup \smaller\smaller\smaller\begin{tabular}{@{}c@{}}%
0\\0\\0
\end{tabular}\endgroup%
}}\!\right]$}%
\EasyButWeakLineBreak%
{$\left[\!\llap{\phantom{%
\begingroup \smaller\smaller\smaller\begin{tabular}{@{}c@{}}%
0\\0\\0
\end{tabular}\endgroup%
}}\right.$}%
\begingroup \smaller\smaller\smaller\begin{tabular}{@{}c@{}}%
-1\\-24\\-234
\end{tabular}\endgroup%
\HardButStrongLineBreak\kern3pt%
\begingroup \smaller\smaller\smaller\begin{tabular}{@{}c@{}}%
-1\\-26\\-232
\end{tabular}\endgroup%
\HardButStrongLineBreak\kern3pt%
\begingroup \smaller\smaller\smaller\begin{tabular}{@{}c@{}}%
7\\156\\1638
\end{tabular}\endgroup%
\HardButStrongLineBreak\kern3pt%
\begingroup \smaller\smaller\smaller\begin{tabular}{@{}c@{}}%
2\\50\\465
\end{tabular}\endgroup%
\HardButStrongLineBreak\kern3pt%
\begingroup \smaller\smaller\smaller\begin{tabular}{@{}c@{}}%
3\\77\\696
\end{tabular}\endgroup%
\HardButStrongLineBreak\kern3pt%
\begingroup \smaller\smaller\smaller\begin{tabular}{@{}c@{}}%
11\\286\\2548
\end{tabular}\endgroup%
{$\left.\llap{\phantom{%
\begingroup \smaller\smaller\smaller\begin{tabular}{@{}c@{}}%
0\\0\\0
\end{tabular}\endgroup%
}}\!\right]$}%

\medskip%
%
\leavevmode\llap{}%
$W_{267}$%
\qquad\llap{44} lattices, $\chi=72$%
\hfill%
$2\infty|\infty22|22\infty|\infty22|2\rtimes D_{4}$%
\nopagebreak\smallskip\hrule\nopagebreak\medskip%
%
%
\leavevmode%
${L_{267.1}}$%
{} : {$1\above{1pt}{1pt}{2}{{\rm II}}4\above{1pt}{1pt}{1}{7}{\cdot}1\above{1pt}{1pt}{1}{}3\above{1pt}{1pt}{-}{}9\above{1pt}{1pt}{-}{}{\cdot}1\above{1pt}{1pt}{-}{}5\above{1pt}{1pt}{-}{}25\above{1pt}{1pt}{-}{}$}\spacer%
\instructions{23,3,2}%
\EasyButWeakLineBreak%
{${450}\above{1pt}{1pt}{l}{2}{60}\above{1pt}{1pt}{15,1}{\infty}{60}\above{1pt}{1pt}{15,8}{\infty b}{60}\above{1pt}{1pt}{r}{2}{18}\above{1pt}{1pt}{b}{2}{10}\above{1pt}{1pt}{b}{2}$}\relax$\,(\times2)$%
\nopagebreak\par%
\nopagebreak\par\leavevmode%
{$\left[\!\llap{\phantom{%
\begingroup \smaller\smaller\smaller
\endgroup%
}}\!\right]$}%
%
%
\hbox{}\par\smallskip%
%
%
\leavevmode%
${L_{267.2}}$%
{} : {$1\above{1pt}{1pt}{2}{6}8\above{1pt}{1pt}{1}{1}{\cdot}1\above{1pt}{1pt}{-}{}3\above{1pt}{1pt}{1}{}9\above{1pt}{1pt}{1}{}{\cdot}1\above{1pt}{1pt}{1}{}5\above{1pt}{1pt}{1}{}25\above{1pt}{1pt}{1}{}$}\spacer%
\instructions{3m,3,2}%
\EasyButWeakLineBreak%
{${36}\above{1pt}{1pt}{s}{2}{120}\above{1pt}{1pt}{60,49}{\infty z}{30}\above{1pt}{1pt}{60,17}{\infty a}{120}\above{1pt}{1pt}{s}{2}{900}\above{1pt}{1pt}{*}{2}{20}\above{1pt}{1pt}{*}{2}$}\relax$\,(\times2)$%
\nopagebreak\par%
\nopagebreak\par\leavevmode%
{$\left[\!\llap{\phantom{%
\begingroup \smaller\smaller\smaller
\endgroup%
}}\!\right]$}%
%
%
\hbox{}\par\smallskip%
%
%
\leavevmode%
${L_{267.3}}$%
{} : {$1\above{1pt}{1pt}{-2}{6}8\above{1pt}{1pt}{-}{5}{\cdot}1\above{1pt}{1pt}{-}{}3\above{1pt}{1pt}{1}{}9\above{1pt}{1pt}{1}{}{\cdot}1\above{1pt}{1pt}{1}{}5\above{1pt}{1pt}{1}{}25\above{1pt}{1pt}{1}{}$}\spacer%
\instructions{32,3,m}%
\EasyButWeakLineBreak%
{${9}\above{1pt}{1pt}{r}{2}{120}\above{1pt}{1pt}{60,19}{\infty z}{30}\above{1pt}{1pt}{60,17}{\infty b}{120}\above{1pt}{1pt}{l}{2}{225}\above{1pt}{1pt}{}{2}{5}\above{1pt}{1pt}{}{2}$}\relax$\,(\times2)$%
\nopagebreak\par%
\nopagebreak\par\leavevmode%
{$\left[\!\llap{\phantom{%
\begingroup \smaller\smaller\smaller
\endgroup%
}}\!\right]$}%

\medskip%
%
\leavevmode\llap{}%
$W_{268}$%
\qquad\llap{44} lattices, $\chi=36$%
\hfill%
$\infty222|222\infty|\rtimes D_{2}$%
\nopagebreak\smallskip\hrule\nopagebreak\medskip%
%
%
\leavevmode%
${L_{268.1}}$%
{} : {$1\above{1pt}{1pt}{2}{{\rm II}}4\above{1pt}{1pt}{1}{7}{\cdot}1\above{1pt}{1pt}{-}{}3\above{1pt}{1pt}{-}{}9\above{1pt}{1pt}{1}{}{\cdot}1\above{1pt}{1pt}{1}{}5\above{1pt}{1pt}{-}{}25\above{1pt}{1pt}{1}{}$}\spacer%
\instructions{23,3,2}%
\EasyButWeakLineBreak%
{${60}\above{1pt}{1pt}{15,8}{\infty}{60}\above{1pt}{1pt}{*}{2}{900}\above{1pt}{1pt}{b}{2}{6}\above{1pt}{1pt}{s}{2}{90}\above{1pt}{1pt}{s}{2}{150}\above{1pt}{1pt}{b}{2}{36}\above{1pt}{1pt}{*}{2}{60}\above{1pt}{1pt}{15,4}{\infty b}$}%
\nopagebreak\par%
\nopagebreak\par\leavevmode%
{$\left[\!\llap{\phantom{%
\begingroup \smaller\smaller\smaller\begin{tabular}{@{}c@{}}%
0\\0\\0
\end{tabular}\endgroup%
}}\right.$}%
\begingroup \smaller\smaller\smaller\begin{tabular}{@{}c@{}}%
-188100\\6300\\0
\end{tabular}\endgroup%
\kern3pt%
\begingroup \smaller\smaller\smaller\begin{tabular}{@{}c@{}}%
6300\\-210\\-15
\end{tabular}\endgroup%
\kern3pt%
\begingroup \smaller\smaller\smaller\begin{tabular}{@{}c@{}}%
0\\-15\\224
\end{tabular}\endgroup%
{$\left.\llap{\phantom{%
\begingroup \smaller\smaller\smaller\begin{tabular}{@{}c@{}}%
0\\0\\0
\end{tabular}\endgroup%
}}\!\right]$}%
\EasyButWeakLineBreak%
{$\left[\!\llap{\phantom{%
\begingroup \smaller\smaller\smaller\begin{tabular}{@{}c@{}}%
0\\0\\0
\end{tabular}\endgroup%
}}\right.$}%
\begingroup \smaller\smaller\smaller\begin{tabular}{@{}c@{}}%
1\\28\\0
\end{tabular}\endgroup%
\HardButStrongLineBreak\kern3pt%
\begingroup \smaller\smaller\smaller\begin{tabular}{@{}c@{}}%
-15\\-448\\-30
\end{tabular}\endgroup%
\HardButStrongLineBreak\kern3pt%
\begingroup \smaller\smaller\smaller\begin{tabular}{@{}c@{}}%
-1\\-30\\0
\end{tabular}\endgroup%
\HardButStrongLineBreak\kern3pt%
\begingroup \smaller\smaller\smaller\begin{tabular}{@{}c@{}}%
6\\179\\12
\end{tabular}\endgroup%
\HardButStrongLineBreak\kern3pt%
\begingroup \smaller\smaller\smaller\begin{tabular}{@{}c@{}}%
46\\1371\\90
\end{tabular}\endgroup%
\HardButStrongLineBreak\kern3pt%
\begingroup \smaller\smaller\smaller\begin{tabular}{@{}c@{}}%
156\\4645\\300
\end{tabular}\endgroup%
\HardButStrongLineBreak\kern3pt%
\begingroup \smaller\smaller\smaller\begin{tabular}{@{}c@{}}%
151\\4494\\288
\end{tabular}\endgroup%
\HardButStrongLineBreak\kern3pt%
\begingroup \smaller\smaller\smaller\begin{tabular}{@{}c@{}}%
111\\3302\\210
\end{tabular}\endgroup%
{$\left.\llap{\phantom{%
\begingroup \smaller\smaller\smaller\begin{tabular}{@{}c@{}}%
0\\0\\0
\end{tabular}\endgroup%
}}\!\right]$}%
%
%
\hbox{}\par\smallskip%
%
%
\leavevmode%
${L_{268.2}}$%
{} : {$1\above{1pt}{1pt}{2}{6}8\above{1pt}{1pt}{1}{1}{\cdot}1\above{1pt}{1pt}{1}{}3\above{1pt}{1pt}{1}{}9\above{1pt}{1pt}{-}{}{\cdot}1\above{1pt}{1pt}{-}{}5\above{1pt}{1pt}{1}{}25\above{1pt}{1pt}{-}{}$}\spacer%
\instructions{3m,3,2}%
\EasyButWeakLineBreak%
{${120}\above{1pt}{1pt}{60,53}{\infty z}{30}\above{1pt}{1pt}{l}{2}{1800}\above{1pt}{1pt}{}{2}{3}\above{1pt}{1pt}{r}{2}{180}\above{1pt}{1pt}{l}{2}{75}\above{1pt}{1pt}{}{2}{72}\above{1pt}{1pt}{r}{2}{30}\above{1pt}{1pt}{60,49}{\infty a}$}%
\nopagebreak\par%
\nopagebreak\par\leavevmode%
{$\left[\!\llap{\phantom{%
\begingroup \smaller\smaller\smaller\begin{tabular}{@{}c@{}}%
0\\0\\0
\end{tabular}\endgroup%
}}\right.$}%
\begingroup \smaller\smaller\smaller\begin{tabular}{@{}c@{}}%
28529065800\\-37254600\\-18932400
\end{tabular}\endgroup%
\kern3pt%
\begingroup \smaller\smaller\smaller\begin{tabular}{@{}c@{}}%
-37254600\\48630\\24705
\end{tabular}\endgroup%
\kern3pt%
\begingroup \smaller\smaller\smaller\begin{tabular}{@{}c@{}}%
-18932400\\24705\\12547
\end{tabular}\endgroup%
{$\left.\llap{\phantom{%
\begingroup \smaller\smaller\smaller\begin{tabular}{@{}c@{}}%
0\\0\\0
\end{tabular}\endgroup%
}}\!\right]$}%
\EasyButWeakLineBreak%
{$\left[\!\llap{\phantom{%
\begingroup \smaller\smaller\smaller\begin{tabular}{@{}c@{}}%
0\\0\\0
\end{tabular}\endgroup%
}}\right.$}%
\begingroup \smaller\smaller\smaller\begin{tabular}{@{}c@{}}%
-29\\-47912\\50580
\end{tabular}\endgroup%
\HardButStrongLineBreak\kern3pt%
\begingroup \smaller\smaller\smaller\begin{tabular}{@{}c@{}}%
14\\23131\\-24420
\end{tabular}\endgroup%
\HardButStrongLineBreak\kern3pt%
\begingroup \smaller\smaller\smaller\begin{tabular}{@{}c@{}}%
1\\1680\\-1800
\end{tabular}\endgroup%
\HardButStrongLineBreak\kern3pt%
\begingroup \smaller\smaller\smaller\begin{tabular}{@{}c@{}}%
-7\\-11564\\12207
\end{tabular}\endgroup%
\HardButStrongLineBreak\kern3pt%
\begingroup \smaller\smaller\smaller\begin{tabular}{@{}c@{}}%
-127\\-209814\\221490
\end{tabular}\endgroup%
\HardButStrongLineBreak\kern3pt%
\begingroup \smaller\smaller\smaller\begin{tabular}{@{}c@{}}%
-248\\-409720\\432525
\end{tabular}\endgroup%
\HardButStrongLineBreak\kern3pt%
\begingroup \smaller\smaller\smaller\begin{tabular}{@{}c@{}}%
-511\\-844224\\891216
\end{tabular}\endgroup%
\HardButStrongLineBreak\kern3pt%
\begingroup \smaller\smaller\smaller\begin{tabular}{@{}c@{}}%
-199\\-328769\\347070
\end{tabular}\endgroup%
{$\left.\llap{\phantom{%
\begingroup \smaller\smaller\smaller\begin{tabular}{@{}c@{}}%
0\\0\\0
\end{tabular}\endgroup%
}}\!\right]$}%
%
%
\hbox{}\par\smallskip%
%
%
\leavevmode%
${L_{268.3}}$%
{} : {$1\above{1pt}{1pt}{-2}{6}8\above{1pt}{1pt}{-}{5}{\cdot}1\above{1pt}{1pt}{1}{}3\above{1pt}{1pt}{1}{}9\above{1pt}{1pt}{-}{}{\cdot}1\above{1pt}{1pt}{-}{}5\above{1pt}{1pt}{1}{}25\above{1pt}{1pt}{-}{}$}\spacer%
\instructions{32,3,m}%
\EasyButWeakLineBreak%
{${120}\above{1pt}{1pt}{60,23}{\infty z}{30}\above{1pt}{1pt}{b}{2}{1800}\above{1pt}{1pt}{*}{2}{12}\above{1pt}{1pt}{l}{2}{45}\above{1pt}{1pt}{r}{2}{300}\above{1pt}{1pt}{*}{2}{72}\above{1pt}{1pt}{b}{2}{30}\above{1pt}{1pt}{60,49}{\infty b}$}%
\nopagebreak\par%
\nopagebreak\par\leavevmode%
{$\left[\!\llap{\phantom{%
\begingroup \smaller\smaller\smaller\begin{tabular}{@{}c@{}}%
0\\0\\0
\end{tabular}\endgroup%
}}\right.$}%
\begingroup \smaller\smaller\smaller\begin{tabular}{@{}c@{}}%
719389800\\-5914800\\-3114000
\end{tabular}\endgroup%
\kern3pt%
\begingroup \smaller\smaller\smaller\begin{tabular}{@{}c@{}}%
-5914800\\48630\\25605
\end{tabular}\endgroup%
\kern3pt%
\begingroup \smaller\smaller\smaller\begin{tabular}{@{}c@{}}%
-3114000\\25605\\13477
\end{tabular}\endgroup%
{$\left.\llap{\phantom{%
\begingroup \smaller\smaller\smaller\begin{tabular}{@{}c@{}}%
0\\0\\0
\end{tabular}\endgroup%
}}\!\right]$}%
\EasyButWeakLineBreak%
{$\left[\!\llap{\phantom{%
\begingroup \smaller\smaller\smaller\begin{tabular}{@{}c@{}}%
0\\0\\0
\end{tabular}\endgroup%
}}\right.$}%
\begingroup \smaller\smaller\smaller\begin{tabular}{@{}c@{}}%
31\\2728\\1980
\end{tabular}\endgroup%
\HardButStrongLineBreak\kern3pt%
\begingroup \smaller\smaller\smaller\begin{tabular}{@{}c@{}}%
-16\\-1409\\-1020
\end{tabular}\endgroup%
\HardButStrongLineBreak\kern3pt%
\begingroup \smaller\smaller\smaller\begin{tabular}{@{}c@{}}%
-29\\-2580\\-1800
\end{tabular}\endgroup%
\HardButStrongLineBreak\kern3pt%
\begingroup \smaller\smaller\smaller\begin{tabular}{@{}c@{}}%
13\\1142\\834
\end{tabular}\endgroup%
\HardButStrongLineBreak\kern3pt%
\begingroup \smaller\smaller\smaller\begin{tabular}{@{}c@{}}%
64\\5628\\4095
\end{tabular}\endgroup%
\HardButStrongLineBreak\kern3pt%
\begingroup \smaller\smaller\smaller\begin{tabular}{@{}c@{}}%
509\\44770\\32550
\end{tabular}\endgroup%
\HardButStrongLineBreak\kern3pt%
\begingroup \smaller\smaller\smaller\begin{tabular}{@{}c@{}}%
527\\46356\\33696
\end{tabular}\endgroup%
\HardButStrongLineBreak\kern3pt%
\begingroup \smaller\smaller\smaller\begin{tabular}{@{}c@{}}%
206\\18121\\13170
\end{tabular}\endgroup%
{$\left.\llap{\phantom{%
\begingroup \smaller\smaller\smaller\begin{tabular}{@{}c@{}}%
0\\0\\0
\end{tabular}\endgroup%
}}\!\right]$}%

\medskip%
%
\leavevmode\llap{}%
$W_{269}$%
\qquad\llap{6} lattices, $\chi=6$%
\hfill%
$2|22\slashtwo2\rtimes D_{2}$%
\nopagebreak\smallskip\hrule\nopagebreak\medskip%
%
%
\leavevmode%
${L_{269.1}}$%
{} : {$1\above{1pt}{1pt}{-2}{{\rm II}}4\above{1pt}{1pt}{-}{3}{\cdot}1\above{1pt}{1pt}{2}{}3\above{1pt}{1pt}{1}{}{\cdot}1\above{1pt}{1pt}{-}{}5\above{1pt}{1pt}{-}{}25\above{1pt}{1pt}{-}{}$}\spacer%
\instructions{2}%
\EasyButWeakLineBreak%
{${12}\above{1pt}{1pt}{r}{2}{10}\above{1pt}{1pt}{l}{2}{300}\above{1pt}{1pt}{r}{2}{2}\above{1pt}{1pt}{b}{2}{50}\above{1pt}{1pt}{l}{2}$}%
\nopagebreak\par%
\nopagebreak\par\leavevmode%
{$\left[\!\llap{\phantom{%
\begingroup \smaller\smaller\smaller\begin{tabular}{@{}c@{}}%
0\\0\\0
\end{tabular}\endgroup%
}}\right.$}%
\begingroup \smaller\smaller\smaller\begin{tabular}{@{}c@{}}%
-107700\\-40200\\5700
\end{tabular}\endgroup%
\kern3pt%
\begingroup \smaller\smaller\smaller\begin{tabular}{@{}c@{}}%
-40200\\-14990\\2135
\end{tabular}\endgroup%
\kern3pt%
\begingroup \smaller\smaller\smaller\begin{tabular}{@{}c@{}}%
5700\\2135\\-298
\end{tabular}\endgroup%
{$\left.\llap{\phantom{%
\begingroup \smaller\smaller\smaller\begin{tabular}{@{}c@{}}%
0\\0\\0
\end{tabular}\endgroup%
}}\!\right]$}%
\EasyButWeakLineBreak%
{$\left[\!\llap{\phantom{%
\begingroup \smaller\smaller\smaller\begin{tabular}{@{}c@{}}%
0\\0\\0
\end{tabular}\endgroup%
}}\right.$}%
\begingroup \smaller\smaller\smaller\begin{tabular}{@{}c@{}}%
23\\-48\\96
\end{tabular}\endgroup%
\HardButStrongLineBreak\kern3pt%
\begingroup \smaller\smaller\smaller\begin{tabular}{@{}c@{}}%
14\\-29\\60
\end{tabular}\endgroup%
\HardButStrongLineBreak\kern3pt%
\begingroup \smaller\smaller\smaller\begin{tabular}{@{}c@{}}%
-431\\900\\-1800
\end{tabular}\endgroup%
\HardButStrongLineBreak\kern3pt%
\begingroup \smaller\smaller\smaller\begin{tabular}{@{}c@{}}%
-24\\50\\-101
\end{tabular}\endgroup%
\HardButStrongLineBreak\kern3pt%
\begingroup \smaller\smaller\smaller\begin{tabular}{@{}c@{}}%
-29\\60\\-125
\end{tabular}\endgroup%
{$\left.\llap{\phantom{%
\begingroup \smaller\smaller\smaller\begin{tabular}{@{}c@{}}%
0\\0\\0
\end{tabular}\endgroup%
}}\!\right]$}%

\medskip%
%
\leavevmode\llap{}%
$W_{270}$%
\qquad\llap{24} lattices, $\chi=24$%
\hfill%
$22222222\rtimes C_{2}$%
\nopagebreak\smallskip\hrule\nopagebreak\medskip%
%
%
\leavevmode%
${L_{270.1}}$%
{} : {$1\above{1pt}{1pt}{-2}{{\rm II}}4\above{1pt}{1pt}{-}{3}{\cdot}1\above{1pt}{1pt}{-}{}3\above{1pt}{1pt}{-}{}9\above{1pt}{1pt}{1}{}{\cdot}1\above{1pt}{1pt}{-}{}5\above{1pt}{1pt}{1}{}25\above{1pt}{1pt}{1}{}$}\spacer%
\instructions{23,3,2}%
\EasyButWeakLineBreak%
{${2}\above{1pt}{1pt}{b}{2}{900}\above{1pt}{1pt}{*}{2}{20}\above{1pt}{1pt}{b}{2}{150}\above{1pt}{1pt}{s}{2}$}\relax$\,(\times2)$%
\nopagebreak\par%
\nopagebreak\par\leavevmode%
{$\left[\!\llap{\phantom{%
\begingroup \smaller\smaller\smaller\begin{tabular}{@{}c@{}}%
0\\0\\0
\end{tabular}\endgroup%
}}\right.$}%
\begingroup \smaller\smaller\smaller\begin{tabular}{@{}c@{}}%
81900\\5400\\-900
\end{tabular}\endgroup%
\kern3pt%
\begingroup \smaller\smaller\smaller\begin{tabular}{@{}c@{}}%
5400\\330\\-45
\end{tabular}\endgroup%
\kern3pt%
\begingroup \smaller\smaller\smaller\begin{tabular}{@{}c@{}}%
-900\\-45\\2
\end{tabular}\endgroup%
{$\left.\llap{\phantom{%
\begingroup \smaller\smaller\smaller\begin{tabular}{@{}c@{}}%
0\\0\\0
\end{tabular}\endgroup%
}}\!\right]$}%
\hfil\penalty500%
{$\left[\!\llap{\phantom{%
\begingroup \smaller\smaller\smaller\begin{tabular}{@{}c@{}}%
0\\0\\0
\end{tabular}\endgroup%
}}\right.$}%
\begingroup \smaller\smaller\smaller\begin{tabular}{@{}c@{}}%
1919\\-41280\\-57600
\end{tabular}\endgroup%
\kern3pt%
\begingroup \smaller\smaller\smaller\begin{tabular}{@{}c@{}}%
92\\-1979\\-2760
\end{tabular}\endgroup%
\kern3pt%
\begingroup \smaller\smaller\smaller\begin{tabular}{@{}c@{}}%
-2\\43\\59
\end{tabular}\endgroup%
{$\left.\llap{\phantom{%
\begingroup \smaller\smaller\smaller\begin{tabular}{@{}c@{}}%
0\\0\\0
\end{tabular}\endgroup%
}}\!\right]$}%
\EasyButWeakLineBreak%
{$\left[\!\llap{\phantom{%
\begingroup \smaller\smaller\smaller\begin{tabular}{@{}c@{}}%
0\\0\\0
\end{tabular}\endgroup%
}}\right.$}%
\begingroup \smaller\smaller\smaller\begin{tabular}{@{}c@{}}%
2\\-43\\-59
\end{tabular}\endgroup%
\HardButStrongLineBreak\kern3pt%
\begingroup \smaller\smaller\smaller\begin{tabular}{@{}c@{}}%
71\\-1530\\-2250
\end{tabular}\endgroup%
\HardButStrongLineBreak\kern3pt%
\begingroup \smaller\smaller\smaller\begin{tabular}{@{}c@{}}%
5\\-108\\-170
\end{tabular}\endgroup%
\HardButStrongLineBreak\kern3pt%
\begingroup \smaller\smaller\smaller\begin{tabular}{@{}c@{}}%
6\\-130\\-225
\end{tabular}\endgroup%
{$\left.\llap{\phantom{%
\begingroup \smaller\smaller\smaller\begin{tabular}{@{}c@{}}%
0\\0\\0
\end{tabular}\endgroup%
}}\!\right]$}%

\medskip%
%
\leavevmode\llap{}%
$W_{271}$%
\qquad\llap{9} lattices, $\chi=12$%
\hfill%
$\slashtwo2|2\slashtwo2|2\rtimes D_{4}$%
\nopagebreak\smallskip\hrule\nopagebreak\medskip%
%
%
\leavevmode%
${L_{271.1}}$%
{} : {$1\above{1pt}{1pt}{-2}{{\rm II}}4\above{1pt}{1pt}{-}{3}{\cdot}1\above{1pt}{1pt}{1}{}3\above{1pt}{1pt}{1}{}9\above{1pt}{1pt}{1}{}{\cdot}1\above{1pt}{1pt}{1}{}5\above{1pt}{1pt}{-}{}25\above{1pt}{1pt}{1}{}$}\spacer%
\instructions{23,3,2}%
\EasyButWeakLineBreak%
{${4}\above{1pt}{1pt}{*}{2}{900}\above{1pt}{1pt}{b}{2}{10}\above{1pt}{1pt}{b}{2}{36}\above{1pt}{1pt}{*}{2}{100}\above{1pt}{1pt}{b}{2}{90}\above{1pt}{1pt}{b}{2}$}%
\nopagebreak\par%
\nopagebreak\par\leavevmode%
{$\left[\!\llap{\phantom{%
\begingroup \smaller\smaller\smaller\begin{tabular}{@{}c@{}}%
0\\0\\0
\end{tabular}\endgroup%
}}\right.$}%
\begingroup \smaller\smaller\smaller\begin{tabular}{@{}c@{}}%
-522900\\2249100\\-449100
\end{tabular}\endgroup%
\kern3pt%
\begingroup \smaller\smaller\smaller\begin{tabular}{@{}c@{}}%
2249100\\-9335490\\1863315
\end{tabular}\endgroup%
\kern3pt%
\begingroup \smaller\smaller\smaller\begin{tabular}{@{}c@{}}%
-449100\\1863315\\-371906
\end{tabular}\endgroup%
{$\left.\llap{\phantom{%
\begingroup \smaller\smaller\smaller\begin{tabular}{@{}c@{}}%
0\\0\\0
\end{tabular}\endgroup%
}}\!\right]$}%
\EasyButWeakLineBreak%
{$\left[\!\llap{\phantom{%
\begingroup \smaller\smaller\smaller\begin{tabular}{@{}c@{}}%
0\\0\\0
\end{tabular}\endgroup%
}}\right.$}%
\begingroup \smaller\smaller\smaller\begin{tabular}{@{}c@{}}%
253\\5060\\25046
\end{tabular}\endgroup%
\HardButStrongLineBreak\kern3pt%
\begingroup \smaller\smaller\smaller\begin{tabular}{@{}c@{}}%
3287\\65730\\325350
\end{tabular}\endgroup%
\HardButStrongLineBreak\kern3pt%
\begingroup \smaller\smaller\smaller\begin{tabular}{@{}c@{}}%
0\\-1\\-5
\end{tabular}\endgroup%
\HardButStrongLineBreak\kern3pt%
\begingroup \smaller\smaller\smaller\begin{tabular}{@{}c@{}}%
-251\\-5022\\-24858
\end{tabular}\endgroup%
\HardButStrongLineBreak\kern3pt%
\begingroup \smaller\smaller\smaller\begin{tabular}{@{}c@{}}%
-249\\-4980\\-24650
\end{tabular}\endgroup%
\HardButStrongLineBreak\kern3pt%
\begingroup \smaller\smaller\smaller\begin{tabular}{@{}c@{}}%
254\\5082\\25155
\end{tabular}\endgroup%
{$\left.\llap{\phantom{%
\begingroup \smaller\smaller\smaller\begin{tabular}{@{}c@{}}%
0\\0\\0
\end{tabular}\endgroup%
}}\!\right]$}%

\medskip%
%
\leavevmode\llap{}%
$W_{272}$%
\qquad\llap{44} lattices, $\chi=36$%
\hfill%
$2222222222\rtimes C_{2}$%
\nopagebreak\smallskip\hrule\nopagebreak\medskip%
%
%
\leavevmode%
${L_{272.1}}$%
{} : {$1\above{1pt}{1pt}{2}{{\rm II}}4\above{1pt}{1pt}{1}{7}{\cdot}1\above{1pt}{1pt}{2}{}3\above{1pt}{1pt}{1}{}{\cdot}1\above{1pt}{1pt}{1}{}5\above{1pt}{1pt}{1}{}25\above{1pt}{1pt}{-}{}$}\spacer%
\instructions{2}%
\EasyButWeakLineBreak%
{${300}\above{1pt}{1pt}{*}{2}{4}\above{1pt}{1pt}{b}{2}{30}\above{1pt}{1pt}{s}{2}{50}\above{1pt}{1pt}{b}{2}{20}\above{1pt}{1pt}{*}{2}$}\relax$\,(\times2)$%
\nopagebreak\par%
\nopagebreak\par\leavevmode%
{$\left[\!\llap{\phantom{%
\begingroup \smaller\smaller\smaller
\endgroup%
}}\!\right]$}%
%
%
\hbox{}\par\smallskip%
%
%
\leavevmode%
${L_{272.2}}$%
{} : {$1\above{1pt}{1pt}{2}{6}8\above{1pt}{1pt}{1}{1}{\cdot}1\above{1pt}{1pt}{2}{}3\above{1pt}{1pt}{-}{}{\cdot}1\above{1pt}{1pt}{-}{}5\above{1pt}{1pt}{-}{}25\above{1pt}{1pt}{1}{}$}\spacer%
\instructions{2}%
\EasyButWeakLineBreak%
{${150}\above{1pt}{1pt}{l}{2}{8}\above{1pt}{1pt}{}{2}{15}\above{1pt}{1pt}{r}{2}{100}\above{1pt}{1pt}{*}{2}{40}\above{1pt}{1pt}{b}{2}$}\relax$\,(\times2)$%
\nopagebreak\par%
\nopagebreak\par\leavevmode%
{$\left[\!\llap{\phantom{%
\begingroup \smaller\smaller\smaller
\endgroup%
}}\!\right]$}%
%
%
\hbox{}\par\smallskip%
%
%
\leavevmode%
${L_{272.3}}$%
{} : {$1\above{1pt}{1pt}{-2}{6}8\above{1pt}{1pt}{-}{5}{\cdot}1\above{1pt}{1pt}{2}{}3\above{1pt}{1pt}{-}{}{\cdot}1\above{1pt}{1pt}{-}{}5\above{1pt}{1pt}{-}{}25\above{1pt}{1pt}{1}{}$}\spacer%
\instructions{m}%
\EasyButWeakLineBreak%
{${150}\above{1pt}{1pt}{b}{2}{8}\above{1pt}{1pt}{*}{2}{60}\above{1pt}{1pt}{l}{2}{25}\above{1pt}{1pt}{}{2}{40}\above{1pt}{1pt}{r}{2}$}\relax$\,(\times2)$%
\nopagebreak\par%
\nopagebreak\par\leavevmode%
{$\left[\!\llap{\phantom{%
\begingroup \smaller\smaller\smaller
\endgroup%
}}\!\right]$}%

\medskip%
%
\leavevmode\llap{}%
$W_{273}$%
\qquad\llap{8} lattices, $\chi=24$%
\hfill%
$2|22|22|22|2\rtimes D_{4}$%
\nopagebreak\smallskip\hrule\nopagebreak\medskip%
%
%
\leavevmode%
${L_{273.1}}$%
{} : {$[1\above{1pt}{1pt}{1}{}2\above{1pt}{1pt}{1}{}]\above{1pt}{1pt}{}{0}64\above{1pt}{1pt}{-}{3}{\cdot}1\above{1pt}{1pt}{2}{}3\above{1pt}{1pt}{1}{}$}\EasyButWeakLineBreak%
{${192}\above{1pt}{1pt}{*}{2}{4}\above{1pt}{1pt}{s}{2}{192}\above{1pt}{1pt}{l}{2}{2}\above{1pt}{1pt}{}{2}{192}\above{1pt}{1pt}{}{2}{1}\above{1pt}{1pt}{r}{2}{192}\above{1pt}{1pt}{*}{2}{8}\above{1pt}{1pt}{s}{2}$}%
\nopagebreak\par%
\nopagebreak\par\leavevmode%
{$\left[\!\llap{\phantom{%
\begingroup \smaller\smaller\smaller\begin{tabular}{@{}c@{}}%
0\\0\\0
\end{tabular}\endgroup%
}}\right.$}%
\begingroup \smaller\smaller\smaller\begin{tabular}{@{}c@{}}%
-27456\\960\\192
\end{tabular}\endgroup%
\kern3pt%
\begingroup \smaller\smaller\smaller\begin{tabular}{@{}c@{}}%
960\\-2\\-10
\end{tabular}\endgroup%
\kern3pt%
\begingroup \smaller\smaller\smaller\begin{tabular}{@{}c@{}}%
192\\-10\\-1
\end{tabular}\endgroup%
{$\left.\llap{\phantom{%
\begingroup \smaller\smaller\smaller\begin{tabular}{@{}c@{}}%
0\\0\\0
\end{tabular}\endgroup%
}}\!\right]$}%
\EasyButWeakLineBreak%
{$\left[\!\llap{\phantom{%
\begingroup \smaller\smaller\smaller\begin{tabular}{@{}c@{}}%
0\\0\\0
\end{tabular}\endgroup%
}}\right.$}%
\begingroup \smaller\smaller\smaller\begin{tabular}{@{}c@{}}%
-5\\-48\\-480
\end{tabular}\endgroup%
\HardButStrongLineBreak\kern3pt%
\begingroup \smaller\smaller\smaller\begin{tabular}{@{}c@{}}%
-1\\-10\\-94
\end{tabular}\endgroup%
\HardButStrongLineBreak\kern3pt%
\begingroup \smaller\smaller\smaller\begin{tabular}{@{}c@{}}%
1\\0\\96
\end{tabular}\endgroup%
\HardButStrongLineBreak\kern3pt%
\begingroup \smaller\smaller\smaller\begin{tabular}{@{}c@{}}%
1\\9\\94
\end{tabular}\endgroup%
\HardButStrongLineBreak\kern3pt%
\begingroup \smaller\smaller\smaller\begin{tabular}{@{}c@{}}%
41\\384\\3840
\end{tabular}\endgroup%
\HardButStrongLineBreak\kern3pt%
\begingroup \smaller\smaller\smaller\begin{tabular}{@{}c@{}}%
2\\19\\187
\end{tabular}\endgroup%
\HardButStrongLineBreak\kern3pt%
\begingroup \smaller\smaller\smaller\begin{tabular}{@{}c@{}}%
35\\336\\3264
\end{tabular}\endgroup%
\HardButStrongLineBreak\kern3pt%
\begingroup \smaller\smaller\smaller\begin{tabular}{@{}c@{}}%
1\\10\\92
\end{tabular}\endgroup%
{$\left.\llap{\phantom{%
\begingroup \smaller\smaller\smaller\begin{tabular}{@{}c@{}}%
0\\0\\0
\end{tabular}\endgroup%
}}\!\right]$}%
%
%
\hbox{}\par\smallskip%
%
%
\leavevmode%
${L_{273.2}}$%
{} : {$1\above{1pt}{1pt}{-}{3}8\above{1pt}{1pt}{1}{1}64\above{1pt}{1pt}{1}{7}{\cdot}1\above{1pt}{1pt}{2}{}3\above{1pt}{1pt}{1}{}$}\EasyButWeakLineBreak%
{${12}\above{1pt}{1pt}{*}{2}{64}\above{1pt}{1pt}{l}{2}{3}\above{1pt}{1pt}{}{2}{8}\above{1pt}{1pt}{r}{2}$}\relax$\,(\times2)$%
\nopagebreak\par%
\nopagebreak\par\leavevmode%
{$\left[\!\llap{\phantom{%
\begingroup \smaller\smaller\smaller\begin{tabular}{@{}c@{}}%
0\\0\\0
\end{tabular}\endgroup%
}}\right.$}%
\begingroup \smaller\smaller\smaller\begin{tabular}{@{}c@{}}%
-6720\\192\\192
\end{tabular}\endgroup%
\kern3pt%
\begingroup \smaller\smaller\smaller\begin{tabular}{@{}c@{}}%
192\\8\\-8
\end{tabular}\endgroup%
\kern3pt%
\begingroup \smaller\smaller\smaller\begin{tabular}{@{}c@{}}%
192\\-8\\-5
\end{tabular}\endgroup%
{$\left.\llap{\phantom{%
\begingroup \smaller\smaller\smaller\begin{tabular}{@{}c@{}}%
0\\0\\0
\end{tabular}\endgroup%
}}\!\right]$}%
\hfil\penalty500%
{$\left[\!\llap{\phantom{%
\begingroup \smaller\smaller\smaller\begin{tabular}{@{}c@{}}%
0\\0\\0
\end{tabular}\endgroup%
}}\right.$}%
\begingroup \smaller\smaller\smaller\begin{tabular}{@{}c@{}}%
95\\480\\2688
\end{tabular}\endgroup%
\kern3pt%
\begingroup \smaller\smaller\smaller\begin{tabular}{@{}c@{}}%
-2\\-11\\-56
\end{tabular}\endgroup%
\kern3pt%
\begingroup \smaller\smaller\smaller\begin{tabular}{@{}c@{}}%
-3\\-15\\-85
\end{tabular}\endgroup%
{$\left.\llap{\phantom{%
\begingroup \smaller\smaller\smaller\begin{tabular}{@{}c@{}}%
0\\0\\0
\end{tabular}\endgroup%
}}\!\right]$}%
\EasyButWeakLineBreak%
{$\left[\!\llap{\phantom{%
\begingroup \smaller\smaller\smaller\begin{tabular}{@{}c@{}}%
0\\0\\0
\end{tabular}\endgroup%
}}\right.$}%
\begingroup \smaller\smaller\smaller\begin{tabular}{@{}c@{}}%
-1\\-6\\-30
\end{tabular}\endgroup%
\HardButStrongLineBreak\kern3pt%
\begingroup \smaller\smaller\smaller\begin{tabular}{@{}c@{}}%
-1\\-4\\-32
\end{tabular}\endgroup%
\HardButStrongLineBreak\kern3pt%
\begingroup \smaller\smaller\smaller\begin{tabular}{@{}c@{}}%
1\\6\\27
\end{tabular}\endgroup%
\HardButStrongLineBreak\kern3pt%
\begingroup \smaller\smaller\smaller\begin{tabular}{@{}c@{}}%
2\\11\\56
\end{tabular}\endgroup%
{$\left.\llap{\phantom{%
\begingroup \smaller\smaller\smaller\begin{tabular}{@{}c@{}}%
0\\0\\0
\end{tabular}\endgroup%
}}\!\right]$}%

\medskip%
%
\leavevmode\llap{}%
$W_{274}$%
\qquad\llap{16} lattices, $\chi=48$%
\hfill%
$2\slashinfty22|22\slashinfty22|2\rtimes D_{4}$%
\nopagebreak\smallskip\hrule\nopagebreak\medskip%
%
%
\leavevmode%
${L_{274.1}}$%
{} : {$[1\above{1pt}{1pt}{-}{}2\above{1pt}{1pt}{1}{}]\above{1pt}{1pt}{}{6}64\above{1pt}{1pt}{1}{1}{\cdot}1\above{1pt}{1pt}{-}{}3\above{1pt}{1pt}{-}{}9\above{1pt}{1pt}{1}{}$}\spacer%
\instructions{3}%
\EasyButWeakLineBreak%
{${576}\above{1pt}{1pt}{}{2}{6}\above{1pt}{1pt}{48,1}{\infty}{24}\above{1pt}{1pt}{*}{2}{576}\above{1pt}{1pt}{s}{2}{8}\above{1pt}{1pt}{*}{2}{576}\above{1pt}{1pt}{s}{2}{24}\above{1pt}{1pt}{48,1}{\infty z}{6}\above{1pt}{1pt}{r}{2}{576}\above{1pt}{1pt}{l}{2}{2}\above{1pt}{1pt}{}{2}$}%
\nopagebreak\par%
\nopagebreak\par\leavevmode%
{$\left[\!\llap{\phantom{%
\begingroup \smaller\smaller\smaller
\endgroup%
}}\!\right]$}%
%
%
\hbox{}\par\smallskip%
%
%
\leavevmode%
${L_{274.2}}$%
{} : {$1\above{1pt}{1pt}{1}{7}8\above{1pt}{1pt}{-}{3}64\above{1pt}{1pt}{1}{1}{\cdot}1\above{1pt}{1pt}{-}{}3\above{1pt}{1pt}{-}{}9\above{1pt}{1pt}{1}{}$}\spacer%
\instructions{3}%
\EasyButWeakLineBreak%
{${576}\above{1pt}{1pt}{}{2}{24}\above{1pt}{1pt}{48,1}{\infty}{24}\above{1pt}{1pt}{r}{2}{576}\above{1pt}{1pt}{b}{2}{8}\above{1pt}{1pt}{l}{2}$}\relax$\,(\times2)$%
\nopagebreak\par%
\nopagebreak\par\leavevmode%
{$\left[\!\llap{\phantom{%
\begingroup \smaller\smaller\smaller
\endgroup%
}}\!\right]$}%

\medskip%
%
\leavevmode\llap{}%
$W_{275}$%
\qquad\llap{12} lattices, $\chi=12$%
\hfill%
$62223$%
\nopagebreak\smallskip\hrule\nopagebreak\medskip%
%
%
\leavevmode%
${L_{275.1}}$%
{} : {$1\above{1pt}{1pt}{-2}{{\rm II}}4\above{1pt}{1pt}{1}{1}{\cdot}1\above{1pt}{1pt}{-}{}3\above{1pt}{1pt}{-}{}27\above{1pt}{1pt}{1}{}{\cdot}1\above{1pt}{1pt}{-2}{}5\above{1pt}{1pt}{-}{}$}\spacer%
\instructions{2}%
\EasyButWeakLineBreak%
{${6}\above{1pt}{1pt}{}{6}{2}\above{1pt}{1pt}{b}{2}{60}\above{1pt}{1pt}{*}{2}{108}\above{1pt}{1pt}{b}{2}{6}\above{1pt}{1pt}{-}{3}$}%
\nopagebreak\par%
\nopagebreak\par\leavevmode%
{$\left[\!\llap{\phantom{%
\begingroup \smaller\smaller\smaller\begin{tabular}{@{}c@{}}%
0\\0\\0
\end{tabular}\endgroup%
}}\right.$}%
\begingroup \smaller\smaller\smaller\begin{tabular}{@{}c@{}}%
-194940\\-76680\\62640
\end{tabular}\endgroup%
\kern3pt%
\begingroup \smaller\smaller\smaller\begin{tabular}{@{}c@{}}%
-76680\\-29982\\24027
\end{tabular}\endgroup%
\kern3pt%
\begingroup \smaller\smaller\smaller\begin{tabular}{@{}c@{}}%
62640\\24027\\-18046
\end{tabular}\endgroup%
{$\left.\llap{\phantom{%
\begingroup \smaller\smaller\smaller\begin{tabular}{@{}c@{}}%
0\\0\\0
\end{tabular}\endgroup%
}}\!\right]$}%
\EasyButWeakLineBreak%
{$\left[\!\llap{\phantom{%
\begingroup \smaller\smaller\smaller\begin{tabular}{@{}c@{}}%
0\\0\\0
\end{tabular}\endgroup%
}}\right.$}%
\begingroup \smaller\smaller\smaller\begin{tabular}{@{}c@{}}%
384\\-1285\\-378
\end{tabular}\endgroup%
\HardButStrongLineBreak\kern3pt%
\begingroup \smaller\smaller\smaller\begin{tabular}{@{}c@{}}%
909\\-3042\\-895
\end{tabular}\endgroup%
\HardButStrongLineBreak\kern3pt%
\begingroup \smaller\smaller\smaller\begin{tabular}{@{}c@{}}%
1249\\-4180\\-1230
\end{tabular}\endgroup%
\HardButStrongLineBreak\kern3pt%
\begingroup \smaller\smaller\smaller\begin{tabular}{@{}c@{}}%
-823\\2754\\810
\end{tabular}\endgroup%
\HardButStrongLineBreak\kern3pt%
\begingroup \smaller\smaller\smaller\begin{tabular}{@{}c@{}}%
-326\\1091\\321
\end{tabular}\endgroup%
{$\left.\llap{\phantom{%
\begingroup \smaller\smaller\smaller\begin{tabular}{@{}c@{}}%
0\\0\\0
\end{tabular}\endgroup%
}}\!\right]$}%

\medskip%
%
\leavevmode\llap{}%
$W_{276}$%
\qquad\llap{12} lattices, $\chi=36$%
\hfill%
$2222222222\rtimes C_{2}$%
\nopagebreak\smallskip\hrule\nopagebreak\medskip%
%
%
\leavevmode%
${L_{276.1}}$%
{} : {$1\above{1pt}{1pt}{-2}{{\rm II}}4\above{1pt}{1pt}{1}{1}{\cdot}1\above{1pt}{1pt}{1}{}3\above{1pt}{1pt}{1}{}27\above{1pt}{1pt}{1}{}{\cdot}1\above{1pt}{1pt}{2}{}5\above{1pt}{1pt}{1}{}$}\spacer%
\instructions{2}%
\EasyButWeakLineBreak%
{${270}\above{1pt}{1pt}{l}{2}{4}\above{1pt}{1pt}{r}{2}{30}\above{1pt}{1pt}{b}{2}{108}\above{1pt}{1pt}{*}{2}{12}\above{1pt}{1pt}{b}{2}$}\relax$\,(\times2)$%
\nopagebreak\par%
\nopagebreak\par\leavevmode%
{$\left[\!\llap{\phantom{%
\begingroup \smaller\smaller\smaller
\endgroup%
}}\!\right]$}%

\medskip%
%
\leavevmode\llap{}%
$W_{277}$%
\qquad\llap{12} lattices, $\chi=24$%
\hfill%
$22222222\rtimes C_{2}$%
\nopagebreak\smallskip\hrule\nopagebreak\medskip%
%
%
\leavevmode%
${L_{277.1}}$%
{} : {$1\above{1pt}{1pt}{-2}{{\rm II}}4\above{1pt}{1pt}{1}{1}{\cdot}1\above{1pt}{1pt}{-}{}3\above{1pt}{1pt}{1}{}27\above{1pt}{1pt}{-}{}{\cdot}1\above{1pt}{1pt}{-2}{}5\above{1pt}{1pt}{-}{}$}\spacer%
\instructions{2}%
\EasyButWeakLineBreak%
{${540}\above{1pt}{1pt}{b}{2}{2}\above{1pt}{1pt}{s}{2}{54}\above{1pt}{1pt}{b}{2}{12}\above{1pt}{1pt}{*}{2}$}\relax$\,(\times2)$%
\nopagebreak\par%
\nopagebreak\par\leavevmode%
{$\left[\!\llap{\phantom{%
\begingroup \smaller\smaller\smaller\begin{tabular}{@{}c@{}}%
0\\0\\0
\end{tabular}\endgroup%
}}\right.$}%
\begingroup \smaller\smaller\smaller\begin{tabular}{@{}c@{}}%
-307260\\7560\\-3780
\end{tabular}\endgroup%
\kern3pt%
\begingroup \smaller\smaller\smaller\begin{tabular}{@{}c@{}}%
7560\\-186\\93
\end{tabular}\endgroup%
\kern3pt%
\begingroup \smaller\smaller\smaller\begin{tabular}{@{}c@{}}%
-3780\\93\\-46
\end{tabular}\endgroup%
{$\left.\llap{\phantom{%
\begingroup \smaller\smaller\smaller\begin{tabular}{@{}c@{}}%
0\\0\\0
\end{tabular}\endgroup%
}}\!\right]$}%
\hfil\penalty500%
{$\left[\!\llap{\phantom{%
\begingroup \smaller\smaller\smaller\begin{tabular}{@{}c@{}}%
0\\0\\0
\end{tabular}\endgroup%
}}\right.$}%
\begingroup \smaller\smaller\smaller\begin{tabular}{@{}c@{}}%
-2161\\-73440\\32400
\end{tabular}\endgroup%
\kern3pt%
\begingroup \smaller\smaller\smaller\begin{tabular}{@{}c@{}}%
53\\1801\\-795
\end{tabular}\endgroup%
\kern3pt%
\begingroup \smaller\smaller\smaller\begin{tabular}{@{}c@{}}%
-24\\-816\\359
\end{tabular}\endgroup%
{$\left.\llap{\phantom{%
\begingroup \smaller\smaller\smaller\begin{tabular}{@{}c@{}}%
0\\0\\0
\end{tabular}\endgroup%
}}\!\right]$}%
\EasyButWeakLineBreak%
{$\left[\!\llap{\phantom{%
\begingroup \smaller\smaller\smaller\begin{tabular}{@{}c@{}}%
0\\0\\0
\end{tabular}\endgroup%
}}\right.$}%
\begingroup \smaller\smaller\smaller\begin{tabular}{@{}c@{}}%
11\\450\\0
\end{tabular}\endgroup%
\HardButStrongLineBreak\kern3pt%
\begingroup \smaller\smaller\smaller\begin{tabular}{@{}c@{}}%
1\\39\\-4
\end{tabular}\endgroup%
\HardButStrongLineBreak\kern3pt%
\begingroup \smaller\smaller\smaller\begin{tabular}{@{}c@{}}%
7\\261\\-54
\end{tabular}\endgroup%
\HardButStrongLineBreak\kern3pt%
\begingroup \smaller\smaller\smaller\begin{tabular}{@{}c@{}}%
3\\106\\-36
\end{tabular}\endgroup%
{$\left.\llap{\phantom{%
\begingroup \smaller\smaller\smaller\begin{tabular}{@{}c@{}}%
0\\0\\0
\end{tabular}\endgroup%
}}\!\right]$}%

\medskip%
%
\leavevmode\llap{}%
$W_{278}$%
\qquad\llap{28} lattices, $\chi=96$%
\hfill%
$22\slashinfty2222|2222\slashinfty2222|22\rtimes D_{4}$%
\nopagebreak\smallskip\hrule\nopagebreak\medskip%
%
%
\leavevmode%
${L_{278.1}}$%
{} : {$[1\above{1pt}{1pt}{1}{}2\above{1pt}{1pt}{-}{}]\above{1pt}{1pt}{}{4}32\above{1pt}{1pt}{-}{3}{\cdot}1\above{1pt}{1pt}{-2}{}7\above{1pt}{1pt}{1}{}$}\EasyButWeakLineBreak%
{${224}\above{1pt}{1pt}{s}{2}{8}\above{1pt}{1pt}{*}{2}{28}\above{1pt}{1pt}{8,3}{\infty z}{7}\above{1pt}{1pt}{r}{2}{8}\above{1pt}{1pt}{*}{2}{224}\above{1pt}{1pt}{l}{2}{1}\above{1pt}{1pt}{r}{2}{56}\above{1pt}{1pt}{*}{2}{4}\above{1pt}{1pt}{*}{2}$}\relax$\,(\times2)$%
\nopagebreak\par%
\nopagebreak\par\leavevmode%
{$\left[\!\llap{\phantom{%
\begingroup \smaller\smaller\smaller
\endgroup%
}}\!\right]$}%
%
%
\hbox{}\par\smallskip%
%
%
\leavevmode%
${L_{278.2}}$%
{} : {$[1\above{1pt}{1pt}{1}{}2\above{1pt}{1pt}{1}{}]\above{1pt}{1pt}{}{0}64\above{1pt}{1pt}{1}{7}{\cdot}1\above{1pt}{1pt}{-2}{}7\above{1pt}{1pt}{1}{}$}\spacer%
\instructions{m}%
\EasyButWeakLineBreak%
{${448}\above{1pt}{1pt}{s}{2}{4}\above{1pt}{1pt}{*}{2}{56}\above{1pt}{1pt}{16,7}{\infty z}{14}\above{1pt}{1pt}{}{2}{1}\above{1pt}{1pt}{r}{2}{448}\above{1pt}{1pt}{*}{2}{8}\above{1pt}{1pt}{*}{2}{28}\above{1pt}{1pt}{l}{2}{2}\above{1pt}{1pt}{r}{2}$}\relax$\,(\times2)$%
\nopagebreak\par%
shares genus with {$ {L_{278.3}}$}%
\nopagebreak\par%
\nopagebreak\par\leavevmode%
{$\left[\!\llap{\phantom{%
\begingroup \smaller\smaller\smaller
\endgroup%
}}\!\right]$}%
%
%
\hbox{}\par\smallskip%
%
%
\leavevmode%
${L_{278.3}}$%
{} : {$[1\above{1pt}{1pt}{1}{}2\above{1pt}{1pt}{1}{}]\above{1pt}{1pt}{}{0}64\above{1pt}{1pt}{1}{7}{\cdot}1\above{1pt}{1pt}{-2}{}7\above{1pt}{1pt}{1}{}$}\EasyButWeakLineBreak%
{${448}\above{1pt}{1pt}{}{2}{1}\above{1pt}{1pt}{r}{2}{56}\above{1pt}{1pt}{16,15}{\infty z}{14}\above{1pt}{1pt}{r}{2}{4}\above{1pt}{1pt}{*}{2}{448}\above{1pt}{1pt}{s}{2}{8}\above{1pt}{1pt}{l}{2}{7}\above{1pt}{1pt}{}{2}{2}\above{1pt}{1pt}{}{2}$}\relax$\,(\times2)$%
\nopagebreak\par%
shares genus with {$ {L_{278.2}}$}%
\nopagebreak\par%
\nopagebreak\par\leavevmode%
{$\left[\!\llap{\phantom{%
\begingroup \smaller\smaller\smaller
\endgroup%
}}\!\right]$}%
%
%
\hbox{}\par\smallskip%
%
%
\leavevmode%
${L_{278.4}}$%
{} : {$1\above{1pt}{1pt}{1}{7}4\above{1pt}{1pt}{1}{7}32\above{1pt}{1pt}{1}{1}{\cdot}1\above{1pt}{1pt}{-2}{}7\above{1pt}{1pt}{1}{}$}\EasyButWeakLineBreak%
{${28}\above{1pt}{1pt}{*}{2}{16}\above{1pt}{1pt}{s}{2}{224}\above{1pt}{1pt}{8,1}{\infty z}{224}\above{1pt}{1pt}{*}{2}{16}\above{1pt}{1pt}{l}{2}{7}\above{1pt}{1pt}{}{2}{32}\above{1pt}{1pt}{}{2}{28}\above{1pt}{1pt}{r}{2}{32}\above{1pt}{1pt}{s}{2}$}\relax$\,(\times2)$%
\nopagebreak\par%
\nopagebreak\par\leavevmode%
{$\left[\!\llap{\phantom{%
\begingroup \smaller\smaller\smaller
\endgroup%
}}\!\right]$}%
%
%
\hbox{}\par\smallskip%
%
%
\leavevmode%
${L_{278.5}}$%
{} : {$1\above{1pt}{1pt}{1}{7}4\above{1pt}{1pt}{1}{1}32\above{1pt}{1pt}{1}{7}{\cdot}1\above{1pt}{1pt}{-2}{}7\above{1pt}{1pt}{1}{}$}\EasyButWeakLineBreak%
{${28}\above{1pt}{1pt}{l}{2}{4}\above{1pt}{1pt}{}{2}{224}\above{1pt}{1pt}{4,1}{\infty}{224}\above{1pt}{1pt}{l}{2}{4}\above{1pt}{1pt}{}{2}{7}\above{1pt}{1pt}{r}{2}{32}\above{1pt}{1pt}{s}{2}{112}\above{1pt}{1pt}{*}{2}{32}\above{1pt}{1pt}{*}{2}$}\relax$\,(\times2)$%
\nopagebreak\par%
\nopagebreak\par\leavevmode%
{$\left[\!\llap{\phantom{%
\begingroup \smaller\smaller\smaller
\endgroup%
}}\!\right]$}%
%
%
\hbox{}\par\smallskip%
%
%
\leavevmode%
${L_{278.6}}$%
{} : {$1\above{1pt}{1pt}{1}{7}8\above{1pt}{1pt}{1}{7}64\above{1pt}{1pt}{1}{1}{\cdot}1\above{1pt}{1pt}{-2}{}7\above{1pt}{1pt}{1}{}$}\EasyButWeakLineBreak%
{${28}\above{1pt}{1pt}{s}{2}{64}\above{1pt}{1pt}{l}{2}{56}\above{1pt}{1pt}{16,9}{\infty}{56}\above{1pt}{1pt}{}{2}{64}\above{1pt}{1pt}{}{2}{7}\above{1pt}{1pt}{r}{2}{32}\above{1pt}{1pt}{s}{2}{448}\above{1pt}{1pt}{*}{2}{32}\above{1pt}{1pt}{*}{2}$}\relax$\,(\times2)$%
\nopagebreak\par%
\nopagebreak\par\leavevmode%
{$\left[\!\llap{\phantom{%
\begingroup \smaller\smaller\smaller
\endgroup%
}}\!\right]$}%

\medskip%
%
\leavevmode\llap{}%
$W_{279}$%
\qquad\llap{8} lattices, $\chi=36$%
\hfill%
$222|2224\slashtwo4\rtimes D_{2}$%
\nopagebreak\smallskip\hrule\nopagebreak\medskip%
%
%
\leavevmode%
${L_{279.1}}$%
{} : {$1\above{1pt}{1pt}{2}{2}32\above{1pt}{1pt}{1}{1}{\cdot}1\above{1pt}{1pt}{2}{}3\above{1pt}{1pt}{-}{}{\cdot}1\above{1pt}{1pt}{2}{}5\above{1pt}{1pt}{1}{}$}\EasyButWeakLineBreak%
{${2}\above{1pt}{1pt}{b}{2}{32}\above{1pt}{1pt}{*}{2}{20}\above{1pt}{1pt}{s}{2}{32}\above{1pt}{1pt}{l}{2}{5}\above{1pt}{1pt}{}{2}{32}\above{1pt}{1pt}{r}{2}{2}\above{1pt}{1pt}{*}{4}{4}\above{1pt}{1pt}{l}{2}{1}\above{1pt}{1pt}{}{4}$}%
\nopagebreak\par%
\nopagebreak\par\leavevmode%
{$\left[\!\llap{\phantom{%
\begingroup \smaller\smaller\smaller
\endgroup%
}}\!\right]$}%

\medskip%
%
\leavevmode\llap{}%
$W_{280}$%
\qquad\llap{8} lattices, $\chi=16$%
\hfill%
$222622$%
\nopagebreak\smallskip\hrule\nopagebreak\medskip%
%
%
\leavevmode%
${L_{280.1}}$%
{} : {$1\above{1pt}{1pt}{-2}{{\rm II}}32\above{1pt}{1pt}{-}{3}{\cdot}1\above{1pt}{1pt}{2}{}3\above{1pt}{1pt}{-}{}{\cdot}1\above{1pt}{1pt}{-2}{}5\above{1pt}{1pt}{-}{}$}\EasyButWeakLineBreak%
{${96}\above{1pt}{1pt}{r}{2}{10}\above{1pt}{1pt}{b}{2}{96}\above{1pt}{1pt}{b}{2}{2}\above{1pt}{1pt}{}{6}{6}\above{1pt}{1pt}{b}{2}{2}\above{1pt}{1pt}{l}{2}$}%
\nopagebreak\par%
\nopagebreak\par\leavevmode%
{$\left[\!\llap{\phantom{%
\begingroup \smaller\smaller\smaller\begin{tabular}{@{}c@{}}%
0\\0\\0
\end{tabular}\endgroup%
}}\right.$}%
\begingroup \smaller\smaller\smaller\begin{tabular}{@{}c@{}}%
-24512160\\61440\\30720
\end{tabular}\endgroup%
\kern3pt%
\begingroup \smaller\smaller\smaller\begin{tabular}{@{}c@{}}%
61440\\-154\\-77
\end{tabular}\endgroup%
\kern3pt%
\begingroup \smaller\smaller\smaller\begin{tabular}{@{}c@{}}%
30720\\-77\\-38
\end{tabular}\endgroup%
{$\left.\llap{\phantom{%
\begingroup \smaller\smaller\smaller\begin{tabular}{@{}c@{}}%
0\\0\\0
\end{tabular}\endgroup%
}}\!\right]$}%
\EasyButWeakLineBreak%
{$\left[\!\llap{\phantom{%
\begingroup \smaller\smaller\smaller\begin{tabular}{@{}c@{}}%
0\\0\\0
\end{tabular}\endgroup%
}}\right.$}%
\begingroup \smaller\smaller\smaller\begin{tabular}{@{}c@{}}%
13\\5184\\0
\end{tabular}\endgroup%
\HardButStrongLineBreak\kern3pt%
\begingroup \smaller\smaller\smaller\begin{tabular}{@{}c@{}}%
2\\800\\-5
\end{tabular}\endgroup%
\HardButStrongLineBreak\kern3pt%
\begingroup \smaller\smaller\smaller\begin{tabular}{@{}c@{}}%
5\\2016\\-48
\end{tabular}\endgroup%
\HardButStrongLineBreak\kern3pt%
\begingroup \smaller\smaller\smaller\begin{tabular}{@{}c@{}}%
0\\4\\-9
\end{tabular}\endgroup%
\HardButStrongLineBreak\kern3pt%
\begingroup \smaller\smaller\smaller\begin{tabular}{@{}c@{}}%
-1\\-399\\0
\end{tabular}\endgroup%
\HardButStrongLineBreak\kern3pt%
\begingroup \smaller\smaller\smaller\begin{tabular}{@{}c@{}}%
0\\-1\\2
\end{tabular}\endgroup%
{$\left.\llap{\phantom{%
\begingroup \smaller\smaller\smaller\begin{tabular}{@{}c@{}}%
0\\0\\0
\end{tabular}\endgroup%
}}\!\right]$}%

\medskip%
%
\leavevmode\llap{}%
$W_{281}$%
\qquad\llap{12} lattices, $\chi=48$%
\hfill%
$2|222|222|222|22\rtimes D_{4}$%
\nopagebreak\smallskip\hrule\nopagebreak\medskip%
%
%
\leavevmode%
${L_{281.1}}$%
{} : {$1\above{1pt}{1pt}{2}{6}32\above{1pt}{1pt}{1}{1}{\cdot}1\above{1pt}{1pt}{-}{}3\above{1pt}{1pt}{-}{}9\above{1pt}{1pt}{-}{}{\cdot}1\above{1pt}{1pt}{-2}{}5\above{1pt}{1pt}{-}{}$}\spacer%
\instructions{3}%
\EasyButWeakLineBreak%
{${288}\above{1pt}{1pt}{r}{2}{6}\above{1pt}{1pt}{b}{2}{32}\above{1pt}{1pt}{*}{2}{60}\above{1pt}{1pt}{s}{2}{288}\above{1pt}{1pt}{l}{2}{15}\above{1pt}{1pt}{}{2}{32}\above{1pt}{1pt}{r}{2}{6}\above{1pt}{1pt}{b}{2}{288}\above{1pt}{1pt}{*}{2}{60}\above{1pt}{1pt}{s}{2}{32}\above{1pt}{1pt}{l}{2}{15}\above{1pt}{1pt}{}{2}$}%
\nopagebreak\par%
\nopagebreak\par\leavevmode%
{$\left[\!\llap{\phantom{%
\begingroup \smaller\smaller\smaller
\endgroup%
}}\!\right]$}%

\medskip%
%
\leavevmode\llap{}%
$W_{282}$%
\qquad\llap{12} lattices, $\chi=48$%
\hfill%
$2|222|222|222|22\rtimes D_{4}$%
\nopagebreak\smallskip\hrule\nopagebreak\medskip%
%
%
\leavevmode%
${L_{282.1}}$%
{} : {$1\above{1pt}{1pt}{-2}{6}32\above{1pt}{1pt}{-}{5}{\cdot}1\above{1pt}{1pt}{1}{}3\above{1pt}{1pt}{-}{}9\above{1pt}{1pt}{1}{}{\cdot}1\above{1pt}{1pt}{-2}{}5\above{1pt}{1pt}{-}{}$}\spacer%
\instructions{3}%
\EasyButWeakLineBreak%
{${160}\above{1pt}{1pt}{b}{2}{6}\above{1pt}{1pt}{l}{2}{1440}\above{1pt}{1pt}{}{2}{1}\above{1pt}{1pt}{r}{2}{1440}\above{1pt}{1pt}{s}{2}{4}\above{1pt}{1pt}{*}{2}{1440}\above{1pt}{1pt}{b}{2}{6}\above{1pt}{1pt}{l}{2}{160}\above{1pt}{1pt}{}{2}{9}\above{1pt}{1pt}{r}{2}{160}\above{1pt}{1pt}{s}{2}{36}\above{1pt}{1pt}{*}{2}$}%
\nopagebreak\par%
\nopagebreak\par\leavevmode%
{$\left[\!\llap{\phantom{%
\begingroup \smaller\smaller\smaller
\endgroup%
}}\!\right]$}%

\medskip%
%
\leavevmode\llap{}%
$W_{283}$%
\qquad\llap{8} lattices, $\chi=36$%
\hfill%
$\infty\infty3\infty26$%
\nopagebreak\smallskip\hrule\nopagebreak\medskip%
%
%
\leavevmode%
${L_{283.1}}$%
{} : {$1\above{1pt}{1pt}{-2}{{\rm II}}8\above{1pt}{1pt}{1}{7}{\cdot}1\above{1pt}{1pt}{-}{}3\above{1pt}{1pt}{-}{}81\above{1pt}{1pt}{1}{}$}\spacer%
\instructions{2}%
\EasyButWeakLineBreak%
{${6}\above{1pt}{1pt}{36,35}{\infty b}{24}\above{1pt}{1pt}{18,5}{\infty z}{6}\above{1pt}{1pt}{-}{3}{6}\above{1pt}{1pt}{36,23}{\infty a}{24}\above{1pt}{1pt}{b}{2}{2}\above{1pt}{1pt}{}{6}$}%
\nopagebreak\par%
\nopagebreak\par\leavevmode%
{$\left[\!\llap{\phantom{%
\begingroup \smaller\smaller\smaller\begin{tabular}{@{}c@{}}%
0\\0\\0
\end{tabular}\endgroup%
}}\right.$}%
\begingroup \smaller\smaller\smaller\begin{tabular}{@{}c@{}}%
-669384\\6480\\-648
\end{tabular}\endgroup%
\kern3pt%
\begingroup \smaller\smaller\smaller\begin{tabular}{@{}c@{}}%
6480\\-30\\-3
\end{tabular}\endgroup%
\kern3pt%
\begingroup \smaller\smaller\smaller\begin{tabular}{@{}c@{}}%
-648\\-3\\2
\end{tabular}\endgroup%
{$\left.\llap{\phantom{%
\begingroup \smaller\smaller\smaller\begin{tabular}{@{}c@{}}%
0\\0\\0
\end{tabular}\endgroup%
}}\!\right]$}%
\EasyButWeakLineBreak%
{$\left[\!\llap{\phantom{%
\begingroup \smaller\smaller\smaller\begin{tabular}{@{}c@{}}%
0\\0\\0
\end{tabular}\endgroup%
}}\right.$}%
\begingroup \smaller\smaller\smaller\begin{tabular}{@{}c@{}}%
4\\635\\2250
\end{tabular}\endgroup%
\HardButStrongLineBreak\kern3pt%
\begingroup \smaller\smaller\smaller\begin{tabular}{@{}c@{}}%
3\\476\\1692
\end{tabular}\endgroup%
\HardButStrongLineBreak\kern3pt%
\begingroup \smaller\smaller\smaller\begin{tabular}{@{}c@{}}%
1\\158\\567
\end{tabular}\endgroup%
\HardButStrongLineBreak\kern3pt%
\begingroup \smaller\smaller\smaller\begin{tabular}{@{}c@{}}%
0\\-1\\3
\end{tabular}\endgroup%
\HardButStrongLineBreak\kern3pt%
\begingroup \smaller\smaller\smaller\begin{tabular}{@{}c@{}}%
-1\\-160\\-564
\end{tabular}\endgroup%
\HardButStrongLineBreak\kern3pt%
\begingroup \smaller\smaller\smaller\begin{tabular}{@{}c@{}}%
0\\0\\-1
\end{tabular}\endgroup%
{$\left.\llap{\phantom{%
\begingroup \smaller\smaller\smaller\begin{tabular}{@{}c@{}}%
0\\0\\0
\end{tabular}\endgroup%
}}\!\right]$}%

\medskip%
%
\leavevmode\llap{}%
$W_{284}$%
\qquad\llap{8} lattices, $\chi=72$%
\hfill%
$2\infty\infty\infty22\infty\infty\infty2\rtimes C_{2}$%
\nopagebreak\smallskip\hrule\nopagebreak\medskip%
%
%
\leavevmode%
${L_{284.1}}$%
{} : {$1\above{1pt}{1pt}{-2}{{\rm II}}8\above{1pt}{1pt}{1}{7}{\cdot}1\above{1pt}{1pt}{1}{}3\above{1pt}{1pt}{-}{}81\above{1pt}{1pt}{-}{}$}\spacer%
\instructions{2}%
\EasyButWeakLineBreak%
{${162}\above{1pt}{1pt}{s}{2}{6}\above{1pt}{1pt}{36,19}{\infty b}{24}\above{1pt}{1pt}{18,13}{\infty z}{6}\above{1pt}{1pt}{36,31}{\infty a}{24}\above{1pt}{1pt}{b}{2}$}\relax$\,(\times2)$%
\nopagebreak\par%
\nopagebreak\par\leavevmode%
{$\left[\!\llap{\phantom{%
\begingroup \smaller\smaller\smaller
\endgroup%
}}\!\right]$}%

\medskip%
%
\leavevmode\llap{}%
$W_{285}$%
\qquad\llap{4} lattices, $\chi=72$%
\hfill%
$\infty|\infty|\infty|\infty|\infty|\infty|\infty|\infty|\rtimes D_{8}$%
\nopagebreak\smallskip\hrule\nopagebreak\medskip%
%
%
\leavevmode%
${L_{285.1}}$%
{} : {$1\above{1pt}{1pt}{-2}{4}4\above{1pt}{1pt}{1}{1}{\cdot}1\above{1pt}{1pt}{-}{}5\above{1pt}{1pt}{1}{}25\above{1pt}{1pt}{-}{}$}\EasyButWeakLineBreak%
{${20}\above{1pt}{1pt}{5,2}{\infty a}{20}\above{1pt}{1pt}{10,3}{\infty}{20}\above{1pt}{1pt}{5,2}{\infty z}{5}\above{1pt}{1pt}{20,17}{\infty}$}\relax$\,(\times2)$%
\nopagebreak\par%
\nopagebreak\par\leavevmode%
{$\left[\!\llap{\phantom{%
\begingroup \smaller\smaller\smaller\begin{tabular}{@{}c@{}}%
0\\0\\0
\end{tabular}\endgroup%
}}\right.$}%
\begingroup \smaller\smaller\smaller\begin{tabular}{@{}c@{}}%
-8700\\-700\\2300
\end{tabular}\endgroup%
\kern3pt%
\begingroup \smaller\smaller\smaller\begin{tabular}{@{}c@{}}%
-700\\-55\\185
\end{tabular}\endgroup%
\kern3pt%
\begingroup \smaller\smaller\smaller\begin{tabular}{@{}c@{}}%
2300\\185\\-608
\end{tabular}\endgroup%
{$\left.\llap{\phantom{%
\begingroup \smaller\smaller\smaller\begin{tabular}{@{}c@{}}%
0\\0\\0
\end{tabular}\endgroup%
}}\!\right]$}%
\hfil\penalty500%
{$\left[\!\llap{\phantom{%
\begingroup \smaller\smaller\smaller\begin{tabular}{@{}c@{}}%
0\\0\\0
\end{tabular}\endgroup%
}}\right.$}%
\begingroup \smaller\smaller\smaller\begin{tabular}{@{}c@{}}%
-641\\160\\-2400
\end{tabular}\endgroup%
\kern3pt%
\begingroup \smaller\smaller\smaller\begin{tabular}{@{}c@{}}%
-48\\11\\-180
\end{tabular}\endgroup%
\kern3pt%
\begingroup \smaller\smaller\smaller\begin{tabular}{@{}c@{}}%
168\\-42\\629
\end{tabular}\endgroup%
{$\left.\llap{\phantom{%
\begingroup \smaller\smaller\smaller\begin{tabular}{@{}c@{}}%
0\\0\\0
\end{tabular}\endgroup%
}}\!\right]$}%
\EasyButWeakLineBreak%
{$\left[\!\llap{\phantom{%
\begingroup \smaller\smaller\smaller\begin{tabular}{@{}c@{}}%
0\\0\\0
\end{tabular}\endgroup%
}}\right.$}%
\begingroup \smaller\smaller\smaller\begin{tabular}{@{}c@{}}%
39\\4\\150
\end{tabular}\endgroup%
\HardButStrongLineBreak\kern3pt%
\begingroup \smaller\smaller\smaller\begin{tabular}{@{}c@{}}%
37\\-4\\140
\end{tabular}\endgroup%
\HardButStrongLineBreak\kern3pt%
\begingroup \smaller\smaller\smaller\begin{tabular}{@{}c@{}}%
51\\-16\\190
\end{tabular}\endgroup%
\HardButStrongLineBreak\kern3pt%
\begingroup \smaller\smaller\smaller\begin{tabular}{@{}c@{}}%
11\\-7\\40
\end{tabular}\endgroup%
{$\left.\llap{\phantom{%
\begingroup \smaller\smaller\smaller\begin{tabular}{@{}c@{}}%
0\\0\\0
\end{tabular}\endgroup%
}}\!\right]$}%
%
%
\hbox{}\par\smallskip%
%
%
\leavevmode%
${L_{285.2}}$%
{} : {$1\above{1pt}{1pt}{2}{{\rm II}}8\above{1pt}{1pt}{-}{5}{\cdot}1\above{1pt}{1pt}{1}{}5\above{1pt}{1pt}{-}{}25\above{1pt}{1pt}{1}{}$}\EasyButWeakLineBreak%
{${40}\above{1pt}{1pt}{5,1}{\infty b}{40}\above{1pt}{1pt}{5,2}{\infty}{40}\above{1pt}{1pt}{5,4}{\infty z}{10}\above{1pt}{1pt}{20,13}{\infty a}$}\relax$\,(\times2)$%
\nopagebreak\par%
\nopagebreak\par\leavevmode%
{$\left[\!\llap{\phantom{%
\begingroup \smaller\smaller\smaller\begin{tabular}{@{}c@{}}%
0\\0\\0
\end{tabular}\endgroup%
}}\right.$}%
\begingroup \smaller\smaller\smaller\begin{tabular}{@{}c@{}}%
-1821400\\-14400\\10800
\end{tabular}\endgroup%
\kern3pt%
\begingroup \smaller\smaller\smaller\begin{tabular}{@{}c@{}}%
-14400\\-110\\85
\end{tabular}\endgroup%
\kern3pt%
\begingroup \smaller\smaller\smaller\begin{tabular}{@{}c@{}}%
10800\\85\\-64
\end{tabular}\endgroup%
{$\left.\llap{\phantom{%
\begingroup \smaller\smaller\smaller\begin{tabular}{@{}c@{}}%
0\\0\\0
\end{tabular}\endgroup%
}}\!\right]$}%
\hfil\penalty500%
{$\left[\!\llap{\phantom{%
\begingroup \smaller\smaller\smaller\begin{tabular}{@{}c@{}}%
0\\0\\0
\end{tabular}\endgroup%
}}\right.$}%
\begingroup \smaller\smaller\smaller\begin{tabular}{@{}c@{}}%
-3041\\-59280\\-592800
\end{tabular}\endgroup%
\kern3pt%
\begingroup \smaller\smaller\smaller\begin{tabular}{@{}c@{}}%
-24\\-469\\-4680
\end{tabular}\endgroup%
\kern3pt%
\begingroup \smaller\smaller\smaller\begin{tabular}{@{}c@{}}%
18\\351\\3509
\end{tabular}\endgroup%
{$\left.\llap{\phantom{%
\begingroup \smaller\smaller\smaller\begin{tabular}{@{}c@{}}%
0\\0\\0
\end{tabular}\endgroup%
}}\!\right]$}%
\EasyButWeakLineBreak%
{$\left[\!\llap{\phantom{%
\begingroup \smaller\smaller\smaller\begin{tabular}{@{}c@{}}%
0\\0\\0
\end{tabular}\endgroup%
}}\right.$}%
\begingroup \smaller\smaller\smaller\begin{tabular}{@{}c@{}}%
15\\304\\2940
\end{tabular}\endgroup%
\HardButStrongLineBreak\kern3pt%
\begingroup \smaller\smaller\smaller\begin{tabular}{@{}c@{}}%
17\\336\\3320
\end{tabular}\endgroup%
\HardButStrongLineBreak\kern3pt%
\begingroup \smaller\smaller\smaller\begin{tabular}{@{}c@{}}%
27\\524\\5260
\end{tabular}\endgroup%
\HardButStrongLineBreak\kern3pt%
\begingroup \smaller\smaller\smaller\begin{tabular}{@{}c@{}}%
7\\133\\1360
\end{tabular}\endgroup%
{$\left.\llap{\phantom{%
\begingroup \smaller\smaller\smaller\begin{tabular}{@{}c@{}}%
0\\0\\0
\end{tabular}\endgroup%
}}\!\right]$}%
%
%
%
%
%
%
%
%
%
%

\medskip%
%
\leavevmode\llap{}%
$W_{286}$%
\qquad\llap{1} lattice, $\chi=24$%
\hfill%
$\slashinfty|\slashinfty|\slashinfty|\slashinfty|\rtimes D_{8}$%
\nopagebreak\smallskip\hrule\nopagebreak\medskip%
%
%
\leavevmode%
${L_{286.1}}$%
{} : {$1\above{1pt}{1pt}{-}{5}8\above{1pt}{1pt}{1}{7}64\above{1pt}{1pt}{-}{5}$}\EasyButWeakLineBreak%
{${32}\above{1pt}{1pt}{16,3}{\infty z}{8}\above{1pt}{1pt}{8,5}{\infty b}$}\relax$\,(\times2)$%
\nopagebreak\par%
\nopagebreak\par\leavevmode%
{$\left[\!\llap{\phantom{%
\begingroup \smaller\smaller\smaller\begin{tabular}{@{}c@{}}%
0\\0\\0
\end{tabular}\endgroup%
}}\right.$}%
\begingroup \smaller\smaller\smaller\begin{tabular}{@{}c@{}}%
6976\\768\\-192
\end{tabular}\endgroup%
\kern3pt%
\begingroup \smaller\smaller\smaller\begin{tabular}{@{}c@{}}%
768\\-8\\-16
\end{tabular}\endgroup%
\kern3pt%
\begingroup \smaller\smaller\smaller\begin{tabular}{@{}c@{}}%
-192\\-16\\5
\end{tabular}\endgroup%
{$\left.\llap{\phantom{%
\begingroup \smaller\smaller\smaller\begin{tabular}{@{}c@{}}%
0\\0\\0
\end{tabular}\endgroup%
}}\!\right]$}%
\hfil\penalty500%
{$\left[\!\llap{\phantom{%
\begingroup \smaller\smaller\smaller\begin{tabular}{@{}c@{}}%
0\\0\\0
\end{tabular}\endgroup%
}}\right.$}%
\begingroup \smaller\smaller\smaller\begin{tabular}{@{}c@{}}%
-17\\-48\\-768
\end{tabular}\endgroup%
\kern3pt%
\begingroup \smaller\smaller\smaller\begin{tabular}{@{}c@{}}%
6\\17\\288
\end{tabular}\endgroup%
\kern3pt%
\begingroup \smaller\smaller\smaller\begin{tabular}{@{}c@{}}%
0\\0\\-1
\end{tabular}\endgroup%
{$\left.\llap{\phantom{%
\begingroup \smaller\smaller\smaller\begin{tabular}{@{}c@{}}%
0\\0\\0
\end{tabular}\endgroup%
}}\!\right]$}%
\EasyButWeakLineBreak%
{$\left[\!\llap{\phantom{%
\begingroup \smaller\smaller\smaller\begin{tabular}{@{}c@{}}%
0\\0\\0
\end{tabular}\endgroup%
}}\right.$}%
\begingroup \smaller\smaller\smaller\begin{tabular}{@{}c@{}}%
1\\2\\48
\end{tabular}\endgroup%
\HardButStrongLineBreak\kern3pt%
\begingroup \smaller\smaller\smaller\begin{tabular}{@{}c@{}}%
-4\\-11\\-188
\end{tabular}\endgroup%
{$\left.\llap{\phantom{%
\begingroup \smaller\smaller\smaller\begin{tabular}{@{}c@{}}%
0\\0\\0
\end{tabular}\endgroup%
}}\!\right]$}%
%
%

\medskip%
%
\leavevmode\llap{}%
$W_{287}$%
\qquad\llap{8} lattices, $\chi=60$%
\hfill%
$242|242|242|242|\rtimes D_{4}$%
\nopagebreak\smallskip\hrule\nopagebreak\medskip%
%
%
\leavevmode%
${L_{287.1}}$%
{} : {$1\above{1pt}{1pt}{2}{2}16\above{1pt}{1pt}{1}{7}{\cdot}1\above{1pt}{1pt}{2}{}3\above{1pt}{1pt}{1}{}{\cdot}1\above{1pt}{1pt}{2}{}11\above{1pt}{1pt}{-}{}$}\EasyButWeakLineBreak%
{${66}\above{1pt}{1pt}{b}{2}{2}\above{1pt}{1pt}{*}{4}{4}\above{1pt}{1pt}{s}{2}{48}\above{1pt}{1pt}{l}{2}{1}\above{1pt}{1pt}{}{4}{2}\above{1pt}{1pt}{s}{2}$}\relax$\,(\times2)$%
\nopagebreak\par%
\nopagebreak\par\leavevmode%
{$\left[\!\llap{\phantom{%
\begingroup \smaller\smaller\smaller
\endgroup%
}}\!\right]$}%

\medskip%
%
\leavevmode\llap{}%
$W_{288}$%
\qquad\llap{12} lattices, $\chi=48$%
\hfill%
$3222632226\rtimes C_{2}$%
\nopagebreak\smallskip\hrule\nopagebreak\medskip%
%
%
\leavevmode%
${L_{288.1}}$%
{} : {$1\above{1pt}{1pt}{-2}{{\rm II}}4\above{1pt}{1pt}{-}{3}{\cdot}1\above{1pt}{1pt}{-}{}3\above{1pt}{1pt}{-}{}27\above{1pt}{1pt}{-}{}{\cdot}1\above{1pt}{1pt}{-2}{}7\above{1pt}{1pt}{1}{}$}\spacer%
\instructions{2}%
\EasyButWeakLineBreak%
{${6}\above{1pt}{1pt}{-}{3}{6}\above{1pt}{1pt}{b}{2}{14}\above{1pt}{1pt}{b}{2}{54}\above{1pt}{1pt}{s}{2}{2}\above{1pt}{1pt}{}{6}$}\relax$\,(\times2)$%
\nopagebreak\par%
\nopagebreak\par\leavevmode%
{$\left[\!\llap{\phantom{%
\begingroup \smaller\smaller\smaller
\endgroup%
}}\!\right]$}%

\medskip%
%
\leavevmode\llap{}%
$W_{289}$%
\qquad\llap{12} lattices, $\chi=36$%
\hfill%
$2222222222\rtimes C_{2}$%
\nopagebreak\smallskip\hrule\nopagebreak\medskip%
%
%
\leavevmode%
${L_{289.1}}$%
{} : {$1\above{1pt}{1pt}{-2}{{\rm II}}4\above{1pt}{1pt}{-}{3}{\cdot}1\above{1pt}{1pt}{1}{}3\above{1pt}{1pt}{1}{}27\above{1pt}{1pt}{-}{}{\cdot}1\above{1pt}{1pt}{2}{}7\above{1pt}{1pt}{-}{}$}\spacer%
\instructions{2}%
\EasyButWeakLineBreak%
{${378}\above{1pt}{1pt}{b}{2}{4}\above{1pt}{1pt}{*}{2}{84}\above{1pt}{1pt}{b}{2}{54}\above{1pt}{1pt}{l}{2}{12}\above{1pt}{1pt}{r}{2}$}\relax$\,(\times2)$%
\nopagebreak\par%
\nopagebreak\par\leavevmode%
{$\left[\!\llap{\phantom{%
\begingroup \smaller\smaller\smaller
\endgroup%
}}\!\right]$}%

\medskip%
%
\leavevmode\llap{}%
$W_{290}$%
\qquad\llap{12} lattices, $\chi=36$%
\hfill%
$2222222222\rtimes C_{2}$%
\nopagebreak\smallskip\hrule\nopagebreak\medskip%
%
%
\leavevmode%
${L_{290.1}}$%
{} : {$1\above{1pt}{1pt}{-2}{{\rm II}}4\above{1pt}{1pt}{-}{3}{\cdot}1\above{1pt}{1pt}{1}{}3\above{1pt}{1pt}{-}{}27\above{1pt}{1pt}{1}{}{\cdot}1\above{1pt}{1pt}{2}{}7\above{1pt}{1pt}{-}{}$}\spacer%
\instructions{2}%
\EasyButWeakLineBreak%
{${756}\above{1pt}{1pt}{*}{2}{4}\above{1pt}{1pt}{b}{2}{42}\above{1pt}{1pt}{l}{2}{108}\above{1pt}{1pt}{r}{2}{6}\above{1pt}{1pt}{b}{2}$}\relax$\,(\times2)$%
\nopagebreak\par%
\nopagebreak\par\leavevmode%
{$\left[\!\llap{\phantom{%
\begingroup \smaller\smaller\smaller
\endgroup%
}}\!\right]$}%

\medskip%
%
\leavevmode\llap{}%
$W_{291}$%
\qquad\llap{28} lattices, $\chi=36$%
\hfill%
$222|2222\slashinfty2\rtimes D_{2}$%
\nopagebreak\smallskip\hrule\nopagebreak\medskip%
%
%
\leavevmode%
${L_{291.1}}$%
{} : {$[1\above{1pt}{1pt}{1}{}2\above{1pt}{1pt}{-}{}]\above{1pt}{1pt}{}{4}32\above{1pt}{1pt}{-}{5}{\cdot}1\above{1pt}{1pt}{2}{}9\above{1pt}{1pt}{-}{}$}\EasyButWeakLineBreak%
{${288}\above{1pt}{1pt}{*}{2}{8}\above{1pt}{1pt}{s}{2}{32}\above{1pt}{1pt}{s}{2}{72}\above{1pt}{1pt}{*}{2}{32}\above{1pt}{1pt}{*}{2}{8}\above{1pt}{1pt}{s}{2}{288}\above{1pt}{1pt}{*}{2}{4}\above{1pt}{1pt}{24,13}{\infty z}{1}\above{1pt}{1pt}{r}{2}$}%
\nopagebreak\par%
\nopagebreak\par\leavevmode%
{$\left[\!\llap{\phantom{%
\begingroup \smaller\smaller\smaller\begin{tabular}{@{}c@{}}%
0\\0\\0
\end{tabular}\endgroup%
}}\right.$}%
\begingroup \smaller\smaller\smaller\begin{tabular}{@{}c@{}}%
72864\\36000\\-288
\end{tabular}\endgroup%
\kern3pt%
\begingroup \smaller\smaller\smaller\begin{tabular}{@{}c@{}}%
36000\\17786\\-142
\end{tabular}\endgroup%
\kern3pt%
\begingroup \smaller\smaller\smaller\begin{tabular}{@{}c@{}}%
-288\\-142\\1
\end{tabular}\endgroup%
{$\left.\llap{\phantom{%
\begingroup \smaller\smaller\smaller\begin{tabular}{@{}c@{}}%
0\\0\\0
\end{tabular}\endgroup%
}}\!\right]$}%
\EasyButWeakLineBreak%
{$\left[\!\llap{\phantom{%
\begingroup \smaller\smaller\smaller\begin{tabular}{@{}c@{}}%
0\\0\\0
\end{tabular}\endgroup%
}}\right.$}%
\begingroup \smaller\smaller\smaller\begin{tabular}{@{}c@{}}%
-35\\72\\144
\end{tabular}\endgroup%
\HardButStrongLineBreak\kern3pt%
\begingroup \smaller\smaller\smaller\begin{tabular}{@{}c@{}}%
-1\\2\\0
\end{tabular}\endgroup%
\HardButStrongLineBreak\kern3pt%
\begingroup \smaller\smaller\smaller\begin{tabular}{@{}c@{}}%
27\\-56\\-160
\end{tabular}\endgroup%
\HardButStrongLineBreak\kern3pt%
\begingroup \smaller\smaller\smaller\begin{tabular}{@{}c@{}}%
113\\-234\\-648
\end{tabular}\endgroup%
\HardButStrongLineBreak\kern3pt%
\begingroup \smaller\smaller\smaller\begin{tabular}{@{}c@{}}%
143\\-296\\-816
\end{tabular}\endgroup%
\HardButStrongLineBreak\kern3pt%
\begingroup \smaller\smaller\smaller\begin{tabular}{@{}c@{}}%
115\\-238\\-656
\end{tabular}\endgroup%
\HardButStrongLineBreak\kern3pt%
\begingroup \smaller\smaller\smaller\begin{tabular}{@{}c@{}}%
1009\\-2088\\-5760
\end{tabular}\endgroup%
\HardButStrongLineBreak\kern3pt%
\begingroup \smaller\smaller\smaller\begin{tabular}{@{}c@{}}%
29\\-60\\-166
\end{tabular}\endgroup%
\HardButStrongLineBreak\kern3pt%
\begingroup \smaller\smaller\smaller\begin{tabular}{@{}c@{}}%
0\\0\\-1
\end{tabular}\endgroup%
{$\left.\llap{\phantom{%
\begingroup \smaller\smaller\smaller\begin{tabular}{@{}c@{}}%
0\\0\\0
\end{tabular}\endgroup%
}}\!\right]$}%
%
%
\hbox{}\par\smallskip%
%
%
\leavevmode%
${L_{291.2}}$%
{} : {$[1\above{1pt}{1pt}{1}{}2\above{1pt}{1pt}{1}{}]\above{1pt}{1pt}{}{0}64\above{1pt}{1pt}{1}{1}{\cdot}1\above{1pt}{1pt}{2}{}9\above{1pt}{1pt}{1}{}$}\spacer%
\instructions{m}%
\EasyButWeakLineBreak%
{${576}\above{1pt}{1pt}{l}{2}{1}\above{1pt}{1pt}{}{2}{64}\above{1pt}{1pt}{}{2}{9}\above{1pt}{1pt}{r}{2}{64}\above{1pt}{1pt}{l}{2}{1}\above{1pt}{1pt}{}{2}{576}\above{1pt}{1pt}{}{2}{2}\above{1pt}{1pt}{48,1}{\infty}{8}\above{1pt}{1pt}{*}{2}$}%
\nopagebreak\par%
shares genus with {$ {L_{291.4}}$}%
\nopagebreak\par%
\nopagebreak\par\leavevmode%
{$\left[\!\llap{\phantom{%
\begingroup \smaller\smaller\smaller\begin{tabular}{@{}c@{}}%
0\\0\\0
\end{tabular}\endgroup%
}}\right.$}%
\begingroup \smaller\smaller\smaller\begin{tabular}{@{}c@{}}%
28224\\0\\-1152
\end{tabular}\endgroup%
\kern3pt%
\begingroup \smaller\smaller\smaller\begin{tabular}{@{}c@{}}%
0\\2\\-2
\end{tabular}\endgroup%
\kern3pt%
\begingroup \smaller\smaller\smaller\begin{tabular}{@{}c@{}}%
-1152\\-2\\49
\end{tabular}\endgroup%
{$\left.\llap{\phantom{%
\begingroup \smaller\smaller\smaller\begin{tabular}{@{}c@{}}%
0\\0\\0
\end{tabular}\endgroup%
}}\!\right]$}%
\EasyButWeakLineBreak%
{$\left[\!\llap{\phantom{%
\begingroup \smaller\smaller\smaller\begin{tabular}{@{}c@{}}%
0\\0\\0
\end{tabular}\endgroup%
}}\right.$}%
\begingroup \smaller\smaller\smaller\begin{tabular}{@{}c@{}}%
-59\\-1296\\-1440
\end{tabular}\endgroup%
\HardButStrongLineBreak\kern3pt%
\begingroup \smaller\smaller\smaller\begin{tabular}{@{}c@{}}%
-2\\-44\\-49
\end{tabular}\endgroup%
\HardButStrongLineBreak\kern3pt%
\begingroup \smaller\smaller\smaller\begin{tabular}{@{}c@{}}%
-13\\-288\\-320
\end{tabular}\endgroup%
\HardButStrongLineBreak\kern3pt%
\begingroup \smaller\smaller\smaller\begin{tabular}{@{}c@{}}%
-4\\-90\\-99
\end{tabular}\endgroup%
\HardButStrongLineBreak\kern3pt%
\begingroup \smaller\smaller\smaller\begin{tabular}{@{}c@{}}%
-9\\-208\\-224
\end{tabular}\endgroup%
\HardButStrongLineBreak\kern3pt%
\begingroup \smaller\smaller\smaller\begin{tabular}{@{}c@{}}%
-1\\-24\\-25
\end{tabular}\endgroup%
\HardButStrongLineBreak\kern3pt%
\begingroup \smaller\smaller\smaller\begin{tabular}{@{}c@{}}%
-23\\-576\\-576
\end{tabular}\endgroup%
\HardButStrongLineBreak\kern3pt%
\begingroup \smaller\smaller\smaller\begin{tabular}{@{}c@{}}%
0\\-1\\0
\end{tabular}\endgroup%
\HardButStrongLineBreak\kern3pt%
\begingroup \smaller\smaller\smaller\begin{tabular}{@{}c@{}}%
-1\\-22\\-24
\end{tabular}\endgroup%
{$\left.\llap{\phantom{%
\begingroup \smaller\smaller\smaller\begin{tabular}{@{}c@{}}%
0\\0\\0
\end{tabular}\endgroup%
}}\!\right]$}%
%
%
\hbox{}\par\smallskip%
%
%
\leavevmode%
${L_{291.3}}$%
{} : {$1\above{1pt}{1pt}{1}{1}4\above{1pt}{1pt}{1}{1}32\above{1pt}{1pt}{1}{7}{\cdot}1\above{1pt}{1pt}{2}{}9\above{1pt}{1pt}{1}{}$}\EasyButWeakLineBreak%
{${36}\above{1pt}{1pt}{l}{2}{4}\above{1pt}{1pt}{}{2}{1}\above{1pt}{1pt}{}{2}{36}\above{1pt}{1pt}{r}{2}{4}\above{1pt}{1pt}{l}{2}{4}\above{1pt}{1pt}{}{2}{9}\above{1pt}{1pt}{r}{2}{32}\above{1pt}{1pt}{24,7}{\infty z}{32}\above{1pt}{1pt}{*}{2}$}%
\nopagebreak\par%
\nopagebreak\par\leavevmode%
{$\left[\!\llap{\phantom{%
\begingroup \smaller\smaller\smaller\begin{tabular}{@{}c@{}}%
0\\0\\0
\end{tabular}\endgroup%
}}\right.$}%
\begingroup \smaller\smaller\smaller\begin{tabular}{@{}c@{}}%
6624\\-576\\288
\end{tabular}\endgroup%
\kern3pt%
\begingroup \smaller\smaller\smaller\begin{tabular}{@{}c@{}}%
-576\\68\\-40
\end{tabular}\endgroup%
\kern3pt%
\begingroup \smaller\smaller\smaller\begin{tabular}{@{}c@{}}%
288\\-40\\25
\end{tabular}\endgroup%
{$\left.\llap{\phantom{%
\begingroup \smaller\smaller\smaller\begin{tabular}{@{}c@{}}%
0\\0\\0
\end{tabular}\endgroup%
}}\!\right]$}%
\EasyButWeakLineBreak%
{$\left[\!\llap{\phantom{%
\begingroup \smaller\smaller\smaller\begin{tabular}{@{}c@{}}%
0\\0\\0
\end{tabular}\endgroup%
}}\right.$}%
\begingroup \smaller\smaller\smaller\begin{tabular}{@{}c@{}}%
-13\\-396\\-486
\end{tabular}\endgroup%
\HardButStrongLineBreak\kern3pt%
\begingroup \smaller\smaller\smaller\begin{tabular}{@{}c@{}}%
-3\\-91\\-112
\end{tabular}\endgroup%
\HardButStrongLineBreak\kern3pt%
\begingroup \smaller\smaller\smaller\begin{tabular}{@{}c@{}}%
-1\\-30\\-37
\end{tabular}\endgroup%
\HardButStrongLineBreak\kern3pt%
\begingroup \smaller\smaller\smaller\begin{tabular}{@{}c@{}}%
-4\\-117\\-144
\end{tabular}\endgroup%
\HardButStrongLineBreak\kern3pt%
\begingroup \smaller\smaller\smaller\begin{tabular}{@{}c@{}}%
-1\\-28\\-34
\end{tabular}\endgroup%
\HardButStrongLineBreak\kern3pt%
\begingroup \smaller\smaller\smaller\begin{tabular}{@{}c@{}}%
-1\\-27\\-32
\end{tabular}\endgroup%
\HardButStrongLineBreak\kern3pt%
\begingroup \smaller\smaller\smaller\begin{tabular}{@{}c@{}}%
-2\\-54\\-63
\end{tabular}\endgroup%
\HardButStrongLineBreak\kern3pt%
\begingroup \smaller\smaller\smaller\begin{tabular}{@{}c@{}}%
-1\\-28\\-32
\end{tabular}\endgroup%
\HardButStrongLineBreak\kern3pt%
\begingroup \smaller\smaller\smaller\begin{tabular}{@{}c@{}}%
-3\\-92\\-112
\end{tabular}\endgroup%
{$\left.\llap{\phantom{%
\begingroup \smaller\smaller\smaller\begin{tabular}{@{}c@{}}%
0\\0\\0
\end{tabular}\endgroup%
}}\!\right]$}%
%
%
\hbox{}\par\smallskip%
%
%
\leavevmode%
${L_{291.4}}$%
{} : {$[1\above{1pt}{1pt}{1}{}2\above{1pt}{1pt}{1}{}]\above{1pt}{1pt}{}{0}64\above{1pt}{1pt}{1}{1}{\cdot}1\above{1pt}{1pt}{2}{}9\above{1pt}{1pt}{1}{}$}\EasyButWeakLineBreak%
{${576}\above{1pt}{1pt}{s}{2}{4}\above{1pt}{1pt}{*}{2}{64}\above{1pt}{1pt}{*}{2}{36}\above{1pt}{1pt}{s}{2}{64}\above{1pt}{1pt}{s}{2}{4}\above{1pt}{1pt}{*}{2}{576}\above{1pt}{1pt}{s}{2}{8}\above{1pt}{1pt}{48,1}{\infty z}{2}\above{1pt}{1pt}{r}{2}$}%
\nopagebreak\par%
shares genus with {$ {L_{291.2}}$}%
\nopagebreak\par%
\nopagebreak\par\leavevmode%
{$\left[\!\llap{\phantom{%
\begingroup \smaller\smaller\smaller\begin{tabular}{@{}c@{}}%
0\\0\\0
\end{tabular}\endgroup%
}}\right.$}%
\begingroup \smaller\smaller\smaller\begin{tabular}{@{}c@{}}%
-1575360\\5184\\5184
\end{tabular}\endgroup%
\kern3pt%
\begingroup \smaller\smaller\smaller\begin{tabular}{@{}c@{}}%
5184\\-2\\-18
\end{tabular}\endgroup%
\kern3pt%
\begingroup \smaller\smaller\smaller\begin{tabular}{@{}c@{}}%
5184\\-18\\-17
\end{tabular}\endgroup%
{$\left.\llap{\phantom{%
\begingroup \smaller\smaller\smaller\begin{tabular}{@{}c@{}}%
0\\0\\0
\end{tabular}\endgroup%
}}\!\right]$}%
\EasyButWeakLineBreak%
{$\left[\!\llap{\phantom{%
\begingroup \smaller\smaller\smaller\begin{tabular}{@{}c@{}}%
0\\0\\0
\end{tabular}\endgroup%
}}\right.$}%
\begingroup \smaller\smaller\smaller\begin{tabular}{@{}c@{}}%
1\\0\\288
\end{tabular}\endgroup%
\HardButStrongLineBreak\kern3pt%
\begingroup \smaller\smaller\smaller\begin{tabular}{@{}c@{}}%
-1\\-18\\-286
\end{tabular}\endgroup%
\HardButStrongLineBreak\kern3pt%
\begingroup \smaller\smaller\smaller\begin{tabular}{@{}c@{}}%
-1\\-16\\-288
\end{tabular}\endgroup%
\HardButStrongLineBreak\kern3pt%
\begingroup \smaller\smaller\smaller\begin{tabular}{@{}c@{}}%
7\\126\\1998
\end{tabular}\endgroup%
\HardButStrongLineBreak\kern3pt%
\begingroup \smaller\smaller\smaller\begin{tabular}{@{}c@{}}%
27\\480\\7712
\end{tabular}\endgroup%
\HardButStrongLineBreak\kern3pt%
\begingroup \smaller\smaller\smaller\begin{tabular}{@{}c@{}}%
13\\230\\3714
\end{tabular}\endgroup%
\HardButStrongLineBreak\kern3pt%
\begingroup \smaller\smaller\smaller\begin{tabular}{@{}c@{}}%
253\\4464\\72288
\end{tabular}\endgroup%
\HardButStrongLineBreak\kern3pt%
\begingroup \smaller\smaller\smaller\begin{tabular}{@{}c@{}}%
9\\158\\2572
\end{tabular}\endgroup%
\HardButStrongLineBreak\kern3pt%
\begingroup \smaller\smaller\smaller\begin{tabular}{@{}c@{}}%
1\\17\\286
\end{tabular}\endgroup%
{$\left.\llap{\phantom{%
\begingroup \smaller\smaller\smaller\begin{tabular}{@{}c@{}}%
0\\0\\0
\end{tabular}\endgroup%
}}\!\right]$}%
%
%
\hbox{}\par\smallskip%
%
%
\leavevmode%
${L_{291.5}}$%
{} : {$1\above{1pt}{1pt}{1}{1}4\above{1pt}{1pt}{1}{7}32\above{1pt}{1pt}{1}{1}{\cdot}1\above{1pt}{1pt}{2}{}9\above{1pt}{1pt}{1}{}$}\EasyButWeakLineBreak%
{${36}\above{1pt}{1pt}{*}{2}{16}\above{1pt}{1pt}{l}{2}{1}\above{1pt}{1pt}{r}{2}{144}\above{1pt}{1pt}{*}{2}{4}\above{1pt}{1pt}{*}{2}{16}\above{1pt}{1pt}{l}{2}{9}\above{1pt}{1pt}{}{2}{32}\above{1pt}{1pt}{12,7}{\infty}{32}\above{1pt}{1pt}{s}{2}$}%
\nopagebreak\par%
\nopagebreak\par\leavevmode%
{$\left[\!\llap{\phantom{%
\begingroup \smaller\smaller\smaller\begin{tabular}{@{}c@{}}%
0\\0\\0
\end{tabular}\endgroup%
}}\right.$}%
\begingroup \smaller\smaller\smaller\begin{tabular}{@{}c@{}}%
1435680\\44640\\6048
\end{tabular}\endgroup%
\kern3pt%
\begingroup \smaller\smaller\smaller\begin{tabular}{@{}c@{}}%
44640\\1388\\188
\end{tabular}\endgroup%
\kern3pt%
\begingroup \smaller\smaller\smaller\begin{tabular}{@{}c@{}}%
6048\\188\\25
\end{tabular}\endgroup%
{$\left.\llap{\phantom{%
\begingroup \smaller\smaller\smaller\begin{tabular}{@{}c@{}}%
0\\0\\0
\end{tabular}\endgroup%
}}\!\right]$}%
\EasyButWeakLineBreak%
{$\left[\!\llap{\phantom{%
\begingroup \smaller\smaller\smaller\begin{tabular}{@{}c@{}}%
0\\0\\0
\end{tabular}\endgroup%
}}\right.$}%
\begingroup \smaller\smaller\smaller\begin{tabular}{@{}c@{}}%
-35\\1152\\-198
\end{tabular}\endgroup%
\HardButStrongLineBreak\kern3pt%
\begingroup \smaller\smaller\smaller\begin{tabular}{@{}c@{}}%
-15\\494\\-88
\end{tabular}\endgroup%
\HardButStrongLineBreak\kern3pt%
\begingroup \smaller\smaller\smaller\begin{tabular}{@{}c@{}}%
-2\\66\\-13
\end{tabular}\endgroup%
\HardButStrongLineBreak\kern3pt%
\begingroup \smaller\smaller\smaller\begin{tabular}{@{}c@{}}%
-7\\234\\-72
\end{tabular}\endgroup%
\HardButStrongLineBreak\kern3pt%
\begingroup \smaller\smaller\smaller\begin{tabular}{@{}c@{}}%
1\\-32\\-2
\end{tabular}\endgroup%
\HardButStrongLineBreak\kern3pt%
\begingroup \smaller\smaller\smaller\begin{tabular}{@{}c@{}}%
5\\-162\\8
\end{tabular}\endgroup%
\HardButStrongLineBreak\kern3pt%
\begingroup \smaller\smaller\smaller\begin{tabular}{@{}c@{}}%
5\\-162\\9
\end{tabular}\endgroup%
\HardButStrongLineBreak\kern3pt%
\begingroup \smaller\smaller\smaller\begin{tabular}{@{}c@{}}%
1\\-32\\0
\end{tabular}\endgroup%
\HardButStrongLineBreak\kern3pt%
\begingroup \smaller\smaller\smaller\begin{tabular}{@{}c@{}}%
-9\\296\\-48
\end{tabular}\endgroup%
{$\left.\llap{\phantom{%
\begingroup \smaller\smaller\smaller\begin{tabular}{@{}c@{}}%
0\\0\\0
\end{tabular}\endgroup%
}}\!\right]$}%
%
%
\hbox{}\par\smallskip%
%
%
\leavevmode%
${L_{291.6}}$%
{} : {$1\above{1pt}{1pt}{1}{1}8\above{1pt}{1pt}{1}{1}64\above{1pt}{1pt}{1}{7}{\cdot}1\above{1pt}{1pt}{2}{}9\above{1pt}{1pt}{1}{}$}\EasyButWeakLineBreak%
{${36}\above{1pt}{1pt}{*}{2}{64}\above{1pt}{1pt}{l}{2}{1}\above{1pt}{1pt}{r}{2}{576}\above{1pt}{1pt}{*}{2}{4}\above{1pt}{1pt}{*}{2}{64}\above{1pt}{1pt}{l}{2}{9}\above{1pt}{1pt}{}{2}{8}\above{1pt}{1pt}{48,31}{\infty}{8}\above{1pt}{1pt}{r}{2}$}%
\nopagebreak\par%
\nopagebreak\par\leavevmode%
{$\left[\!\llap{\phantom{%
\begingroup \smaller\smaller\smaller\begin{tabular}{@{}c@{}}%
0\\0\\0
\end{tabular}\endgroup%
}}\right.$}%
\begingroup \smaller\smaller\smaller\begin{tabular}{@{}c@{}}%
-8193600\\23616\\24192
\end{tabular}\endgroup%
\kern3pt%
\begingroup \smaller\smaller\smaller\begin{tabular}{@{}c@{}}%
23616\\-56\\-72
\end{tabular}\endgroup%
\kern3pt%
\begingroup \smaller\smaller\smaller\begin{tabular}{@{}c@{}}%
24192\\-72\\-71
\end{tabular}\endgroup%
{$\left.\llap{\phantom{%
\begingroup \smaller\smaller\smaller\begin{tabular}{@{}c@{}}%
0\\0\\0
\end{tabular}\endgroup%
}}\!\right]$}%
\EasyButWeakLineBreak%
{$\left[\!\llap{\phantom{%
\begingroup \smaller\smaller\smaller\begin{tabular}{@{}c@{}}%
0\\0\\0
\end{tabular}\endgroup%
}}\right.$}%
\begingroup \smaller\smaller\smaller\begin{tabular}{@{}c@{}}%
53\\2826\\15174
\end{tabular}\endgroup%
\HardButStrongLineBreak\kern3pt%
\begingroup \smaller\smaller\smaller\begin{tabular}{@{}c@{}}%
39\\2076\\11168
\end{tabular}\endgroup%
\HardButStrongLineBreak\kern3pt%
\begingroup \smaller\smaller\smaller\begin{tabular}{@{}c@{}}%
2\\106\\573
\end{tabular}\endgroup%
\HardButStrongLineBreak\kern3pt%
\begingroup \smaller\smaller\smaller\begin{tabular}{@{}c@{}}%
5\\252\\1440
\end{tabular}\endgroup%
\HardButStrongLineBreak\kern3pt%
\begingroup \smaller\smaller\smaller\begin{tabular}{@{}c@{}}%
-1\\-54\\-286
\end{tabular}\endgroup%
\HardButStrongLineBreak\kern3pt%
\begingroup \smaller\smaller\smaller\begin{tabular}{@{}c@{}}%
-1\\-52\\-288
\end{tabular}\endgroup%
\HardButStrongLineBreak\kern3pt%
\begingroup \smaller\smaller\smaller\begin{tabular}{@{}c@{}}%
4\\216\\1143
\end{tabular}\endgroup%
\HardButStrongLineBreak\kern3pt%
\begingroup \smaller\smaller\smaller\begin{tabular}{@{}c@{}}%
4\\215\\1144
\end{tabular}\endgroup%
\HardButStrongLineBreak\kern3pt%
\begingroup \smaller\smaller\smaller\begin{tabular}{@{}c@{}}%
9\\481\\2576
\end{tabular}\endgroup%
{$\left.\llap{\phantom{%
\begingroup \smaller\smaller\smaller\begin{tabular}{@{}c@{}}%
0\\0\\0
\end{tabular}\endgroup%
}}\!\right]$}%

\medskip%
%
\leavevmode\llap{}%
$W_{292}$%
\qquad\llap{6} lattices, $\chi=48$%
\hfill%
$22|222|222|222|2\rtimes D_{4}$%
\nopagebreak\smallskip\hrule\nopagebreak\medskip%
%
%
\leavevmode%
${L_{292.1}}$%
{} : {$1\above{1pt}{1pt}{-2}{{\rm II}}4\above{1pt}{1pt}{1}{1}{\cdot}1\above{1pt}{1pt}{-}{}3\above{1pt}{1pt}{1}{}9\above{1pt}{1pt}{-}{}{\cdot}1\above{1pt}{1pt}{-2}{}23\above{1pt}{1pt}{1}{}$}\spacer%
\instructions{2}%
\EasyButWeakLineBreak%
{${92}\above{1pt}{1pt}{b}{2}{18}\above{1pt}{1pt}{b}{2}{138}\above{1pt}{1pt}{b}{2}{2}\above{1pt}{1pt}{b}{2}{828}\above{1pt}{1pt}{*}{2}{12}\above{1pt}{1pt}{*}{2}$}\relax$\,(\times2)$%
\nopagebreak\par%
\nopagebreak\par\leavevmode%
{$\left[\!\llap{\phantom{%
\begingroup \smaller\smaller\smaller
\endgroup%
}}\!\right]$}%

\medskip%
%
\leavevmode\llap{}%
$W_{293}$%
\qquad\llap{6} lattices, $\chi=48$%
\hfill%
$222|222|222|222|\rtimes D_{4}$%
\nopagebreak\smallskip\hrule\nopagebreak\medskip%
%
%
\leavevmode%
${L_{293.1}}$%
{} : {$1\above{1pt}{1pt}{-2}{{\rm II}}4\above{1pt}{1pt}{1}{1}{\cdot}1\above{1pt}{1pt}{1}{}3\above{1pt}{1pt}{1}{}9\above{1pt}{1pt}{1}{}{\cdot}1\above{1pt}{1pt}{-2}{}23\above{1pt}{1pt}{1}{}$}\spacer%
\instructions{2}%
\EasyButWeakLineBreak%
{${12}\above{1pt}{1pt}{b}{2}{46}\above{1pt}{1pt}{l}{2}{36}\above{1pt}{1pt}{r}{2}{138}\above{1pt}{1pt}{l}{2}{4}\above{1pt}{1pt}{r}{2}{414}\above{1pt}{1pt}{b}{2}$}\relax$\,(\times2)$%
\nopagebreak\par%
\nopagebreak\par\leavevmode%
{$\left[\!\llap{\phantom{%
\begingroup \smaller\smaller\smaller
\endgroup%
}}\!\right]$}%

\medskip%
%
\leavevmode\llap{}%
$W_{294}$%
\qquad\llap{8} lattices, $\chi=42$%
\hfill%
$42222|22224|\rtimes D_{2}$%
\nopagebreak\smallskip\hrule\nopagebreak\medskip%
%
%
\leavevmode%
${L_{294.1}}$%
{} : {$1\above{1pt}{1pt}{2}{2}16\above{1pt}{1pt}{1}{1}{\cdot}1\above{1pt}{1pt}{2}{}3\above{1pt}{1pt}{-}{}{\cdot}1\above{1pt}{1pt}{2}{}13\above{1pt}{1pt}{1}{}$}\EasyButWeakLineBreak%
{${2}\above{1pt}{1pt}{*}{4}{4}\above{1pt}{1pt}{l}{2}{13}\above{1pt}{1pt}{}{2}{16}\above{1pt}{1pt}{}{2}{1}\above{1pt}{1pt}{r}{2}{208}\above{1pt}{1pt}{s}{2}{4}\above{1pt}{1pt}{*}{2}{16}\above{1pt}{1pt}{*}{2}{52}\above{1pt}{1pt}{l}{2}{1}\above{1pt}{1pt}{}{4}$}%
\nopagebreak\par%
\nopagebreak\par\leavevmode%
{$\left[\!\llap{\phantom{%
\begingroup \smaller\smaller\smaller
\endgroup%
}}\!\right]$}%

\medskip%
%
\leavevmode\llap{}%
$W_{295}$%
\qquad\llap{32} lattices, $\chi=42$%
\hfill%
$2222\slashtwo22222|2\rtimes D_{2}$%
\nopagebreak\smallskip\hrule\nopagebreak\medskip%
%
%
\leavevmode%
${L_{295.1}}$%
{} : {$1\above{1pt}{1pt}{2}{2}8\above{1pt}{1pt}{1}{1}{\cdot}1\above{1pt}{1pt}{-}{}3\above{1pt}{1pt}{1}{}9\above{1pt}{1pt}{-}{}{\cdot}1\above{1pt}{1pt}{2}{}13\above{1pt}{1pt}{-}{}$}\spacer%
\instructions{2,m}%
\EasyButWeakLineBreak%
{${104}\above{1pt}{1pt}{b}{2}{18}\above{1pt}{1pt}{l}{2}{8}\above{1pt}{1pt}{r}{2}{234}\above{1pt}{1pt}{b}{2}{2}\above{1pt}{1pt}{s}{2}{18}\above{1pt}{1pt}{b}{2}{26}\above{1pt}{1pt}{l}{2}{72}\above{1pt}{1pt}{r}{2}{2}\above{1pt}{1pt}{b}{2}{936}\above{1pt}{1pt}{*}{2}{12}\above{1pt}{1pt}{*}{2}$}%
\nopagebreak\par%
\nopagebreak\par\leavevmode%
{$\left[\!\llap{\phantom{%
\begingroup \smaller\smaller\smaller
\endgroup%
}}\!\right]$}%
%
%
\hbox{}\par\smallskip%
%
%
\leavevmode%
${L_{295.2}}$%
{} : {$1\above{1pt}{1pt}{-2}{6}16\above{1pt}{1pt}{-}{5}{\cdot}1\above{1pt}{1pt}{1}{}3\above{1pt}{1pt}{-}{}9\above{1pt}{1pt}{1}{}{\cdot}1\above{1pt}{1pt}{2}{}13\above{1pt}{1pt}{1}{}$}\spacer%
\instructions{3,m}%
\EasyButWeakLineBreak%
{${1872}\above{1pt}{1pt}{}{2}{1}\above{1pt}{1pt}{r}{2}{144}\above{1pt}{1pt}{l}{2}{13}\above{1pt}{1pt}{}{2}{9}\above{1pt}{1pt}{r}{2}{4}\above{1pt}{1pt}{*}{2}{468}\above{1pt}{1pt}{s}{2}{16}\above{1pt}{1pt}{s}{2}{36}\above{1pt}{1pt}{*}{2}{208}\above{1pt}{1pt}{b}{2}{6}\above{1pt}{1pt}{l}{2}$}%
\nopagebreak\par%
shares genus with 3-dual\nopagebreak\par%
\nopagebreak\par\leavevmode%
{$\left[\!\llap{\phantom{%
\begingroup \smaller\smaller\smaller
\endgroup%
}}\!\right]$}%

\medskip%
%
\leavevmode\llap{}%
$W_{296}$%
\qquad\llap{32} lattices, $\chi=72$%
\hfill%
$2222222222222222\rtimes C_{2}$%
\nopagebreak\smallskip\hrule\nopagebreak\medskip%
%
%
\leavevmode%
${L_{296.1}}$%
{} : {$1\above{1pt}{1pt}{-2}{{\rm II}}8\above{1pt}{1pt}{1}{7}{\cdot}1\above{1pt}{1pt}{2}{}9\above{1pt}{1pt}{-}{}{\cdot}1\above{1pt}{1pt}{-2}{}5\above{1pt}{1pt}{-}{}{\cdot}1\above{1pt}{1pt}{2}{}7\above{1pt}{1pt}{1}{}$}\spacer%
\instructions{2}%
\EasyButWeakLineBreak%
{${504}\above{1pt}{1pt}{r}{2}{10}\above{1pt}{1pt}{s}{2}{126}\above{1pt}{1pt}{b}{2}{2}\above{1pt}{1pt}{l}{2}{56}\above{1pt}{1pt}{r}{2}{18}\above{1pt}{1pt}{b}{2}{14}\above{1pt}{1pt}{b}{2}{2}\above{1pt}{1pt}{l}{2}$}\relax$\,(\times2)$%
\nopagebreak\par%
\nopagebreak\par\leavevmode%
{$\left[\!\llap{\phantom{%
\begingroup \smaller\smaller\smaller
\endgroup%
}}\!\right]$}%

\medskip%
%
\leavevmode\llap{}%
$W_{297}$%
\qquad\llap{8} lattices, $\chi=36$%
\hfill%
$22222|22222|\rtimes D_{2}$%
\nopagebreak\smallskip\hrule\nopagebreak\medskip%
%
%
\leavevmode%
${L_{297.1}}$%
{} : {$1\above{1pt}{1pt}{-2}{2}32\above{1pt}{1pt}{1}{7}{\cdot}1\above{1pt}{1pt}{2}{}3\above{1pt}{1pt}{1}{}{\cdot}1\above{1pt}{1pt}{2}{}7\above{1pt}{1pt}{1}{}$}\EasyButWeakLineBreak%
{${224}\above{1pt}{1pt}{s}{2}{12}\above{1pt}{1pt}{*}{2}{28}\above{1pt}{1pt}{l}{2}{3}\above{1pt}{1pt}{}{2}{224}\above{1pt}{1pt}{r}{2}{2}\above{1pt}{1pt}{b}{2}{224}\above{1pt}{1pt}{*}{2}{12}\above{1pt}{1pt}{l}{2}{7}\above{1pt}{1pt}{}{2}{3}\above{1pt}{1pt}{r}{2}$}%
\nopagebreak\par%
\nopagebreak\par\leavevmode%
{$\left[\!\llap{\phantom{%
\begingroup \smaller\smaller\smaller
\endgroup%
}}\!\right]$}%

\medskip%
%
\leavevmode\llap{}%
$W_{298}$%
\qquad\llap{6} lattices, $\chi=40$%
\hfill%
$26|62|26|62|\rtimes D_{4}$%
\nopagebreak\smallskip\hrule\nopagebreak\medskip%
%
%
\leavevmode%
${L_{298.1}}$%
{} : {$1\above{1pt}{1pt}{-2}{{\rm II}}4\above{1pt}{1pt}{1}{7}{\cdot}1\above{1pt}{1pt}{-}{}3\above{1pt}{1pt}{-}{}9\above{1pt}{1pt}{-}{}{\cdot}1\above{1pt}{1pt}{-2}{}25\above{1pt}{1pt}{1}{}$}\spacer%
\instructions{2}%
\EasyButWeakLineBreak%
{${150}\above{1pt}{1pt}{s}{2}{18}\above{1pt}{1pt}{}{6}{6}\above{1pt}{1pt}{}{6}{2}\above{1pt}{1pt}{s}{2}$}\relax$\,(\times2)$%
\nopagebreak\par%
\nopagebreak\par\leavevmode%
{$\left[\!\llap{\phantom{%
\begingroup \smaller\smaller\smaller
\endgroup%
}}\!\right]$}%

\medskip%
%
\leavevmode\llap{}%
$W_{299}$%
\qquad\llap{6} lattices, $\chi=20$%
\hfill%
$222|222\slashthree\rtimes D_{2}$%
\nopagebreak\smallskip\hrule\nopagebreak\medskip%
%
%
\leavevmode%
${L_{299.1}}$%
{} : {$1\above{1pt}{1pt}{-2}{{\rm II}}4\above{1pt}{1pt}{1}{7}{\cdot}1\above{1pt}{1pt}{1}{}3\above{1pt}{1pt}{-}{}9\above{1pt}{1pt}{1}{}{\cdot}1\above{1pt}{1pt}{-2}{}25\above{1pt}{1pt}{1}{}$}\spacer%
\instructions{2}%
\EasyButWeakLineBreak%
{${6}\above{1pt}{1pt}{b}{2}{100}\above{1pt}{1pt}{*}{2}{36}\above{1pt}{1pt}{b}{2}{150}\above{1pt}{1pt}{b}{2}{4}\above{1pt}{1pt}{*}{2}{900}\above{1pt}{1pt}{b}{2}{6}\above{1pt}{1pt}{+}{3}$}%
\nopagebreak\par%
\nopagebreak\par\leavevmode%
{$\left[\!\llap{\phantom{%
\begingroup \smaller\smaller\smaller\begin{tabular}{@{}c@{}}%
0\\0\\0
\end{tabular}\endgroup%
}}\right.$}%
\begingroup \smaller\smaller\smaller\begin{tabular}{@{}c@{}}%
-3536100\\4500\\3600
\end{tabular}\endgroup%
\kern3pt%
\begingroup \smaller\smaller\smaller\begin{tabular}{@{}c@{}}%
4500\\6\\-9
\end{tabular}\endgroup%
\kern3pt%
\begingroup \smaller\smaller\smaller\begin{tabular}{@{}c@{}}%
3600\\-9\\-2
\end{tabular}\endgroup%
{$\left.\llap{\phantom{%
\begingroup \smaller\smaller\smaller\begin{tabular}{@{}c@{}}%
0\\0\\0
\end{tabular}\endgroup%
}}\!\right]$}%
\EasyButWeakLineBreak%
{$\left[\!\llap{\phantom{%
\begingroup \smaller\smaller\smaller\begin{tabular}{@{}c@{}}%
0\\0\\0
\end{tabular}\endgroup%
}}\right.$}%
\begingroup \smaller\smaller\smaller\begin{tabular}{@{}c@{}}%
0\\-1\\0
\end{tabular}\endgroup%
\HardButStrongLineBreak\kern3pt%
\begingroup \smaller\smaller\smaller\begin{tabular}{@{}c@{}}%
3\\750\\2000
\end{tabular}\endgroup%
\HardButStrongLineBreak\kern3pt%
\begingroup \smaller\smaller\smaller\begin{tabular}{@{}c@{}}%
1\\252\\666
\end{tabular}\endgroup%
\HardButStrongLineBreak\kern3pt%
\begingroup \smaller\smaller\smaller\begin{tabular}{@{}c@{}}%
-1\\-250\\-675
\end{tabular}\endgroup%
\HardButStrongLineBreak\kern3pt%
\begingroup \smaller\smaller\smaller\begin{tabular}{@{}c@{}}%
-1\\-252\\-670
\end{tabular}\endgroup%
\HardButStrongLineBreak\kern3pt%
\begingroup \smaller\smaller\smaller\begin{tabular}{@{}c@{}}%
-41\\-10350\\-27450
\end{tabular}\endgroup%
\HardButStrongLineBreak\kern3pt%
\begingroup \smaller\smaller\smaller\begin{tabular}{@{}c@{}}%
-1\\-253\\-669
\end{tabular}\endgroup%
{$\left.\llap{\phantom{%
\begingroup \smaller\smaller\smaller\begin{tabular}{@{}c@{}}%
0\\0\\0
\end{tabular}\endgroup%
}}\!\right]$}%

\medskip%
%
\leavevmode\llap{}%
$W_{300}$%
\qquad\llap{8} lattices, $\chi=24$%
\hfill%
$\infty22\infty\infty$%
\nopagebreak\smallskip\hrule\nopagebreak\medskip%
%
%
\leavevmode%
${L_{300.1}}$%
{} : {$1\above{1pt}{1pt}{-2}{{\rm II}}8\above{1pt}{1pt}{-}{3}{\cdot}1\above{1pt}{1pt}{-}{}7\above{1pt}{1pt}{1}{}49\above{1pt}{1pt}{1}{}$}\spacer%
\instructions{2}%
\EasyButWeakLineBreak%
{${56}\above{1pt}{1pt}{14,9}{\infty z}{14}\above{1pt}{1pt}{b}{2}{98}\above{1pt}{1pt}{b}{2}{56}\above{1pt}{1pt}{14,1}{\infty z}{14}\above{1pt}{1pt}{28,23}{\infty b}$}%
\nopagebreak\par%
\nopagebreak\par\leavevmode%
{$\left[\!\llap{\phantom{%
\begingroup \smaller\smaller\smaller\begin{tabular}{@{}c@{}}%
0\\0\\0
\end{tabular}\endgroup%
}}\right.$}%
\begingroup \smaller\smaller\smaller\begin{tabular}{@{}c@{}}%
-171304\\2744\\-8624
\end{tabular}\endgroup%
\kern3pt%
\begingroup \smaller\smaller\smaller\begin{tabular}{@{}c@{}}%
2744\\-42\\147
\end{tabular}\endgroup%
\kern3pt%
\begingroup \smaller\smaller\smaller\begin{tabular}{@{}c@{}}%
-8624\\147\\-394
\end{tabular}\endgroup%
{$\left.\llap{\phantom{%
\begingroup \smaller\smaller\smaller\begin{tabular}{@{}c@{}}%
0\\0\\0
\end{tabular}\endgroup%
}}\!\right]$}%
\EasyButWeakLineBreak%
{$\left[\!\llap{\phantom{%
\begingroup \smaller\smaller\smaller\begin{tabular}{@{}c@{}}%
0\\0\\0
\end{tabular}\endgroup%
}}\right.$}%
\begingroup \smaller\smaller\smaller\begin{tabular}{@{}c@{}}%
-45\\-1672\\364
\end{tabular}\endgroup%
\HardButStrongLineBreak\kern3pt%
\begingroup \smaller\smaller\smaller\begin{tabular}{@{}c@{}}%
-7\\-263\\56
\end{tabular}\endgroup%
\HardButStrongLineBreak\kern3pt%
\begingroup \smaller\smaller\smaller\begin{tabular}{@{}c@{}}%
12\\441\\-98
\end{tabular}\endgroup%
\HardButStrongLineBreak\kern3pt%
\begingroup \smaller\smaller\smaller\begin{tabular}{@{}c@{}}%
7\\260\\-56
\end{tabular}\endgroup%
\HardButStrongLineBreak\kern3pt%
\begingroup \smaller\smaller\smaller\begin{tabular}{@{}c@{}}%
-6\\-221\\49
\end{tabular}\endgroup%
{$\left.\llap{\phantom{%
\begingroup \smaller\smaller\smaller\begin{tabular}{@{}c@{}}%
0\\0\\0
\end{tabular}\endgroup%
}}\!\right]$}%

\medskip%
%
\leavevmode\llap{}%
$W_{301}$%
\qquad\llap{8} lattices, $\chi=24$%
\hfill%
$22222222$%
\nopagebreak\smallskip\hrule\nopagebreak\medskip%
%
%
\leavevmode%
${L_{301.1}}$%
{} : {$1\above{1pt}{1pt}{-2}{{\rm II}}16\above{1pt}{1pt}{1}{1}{\cdot}1\above{1pt}{1pt}{2}{}9\above{1pt}{1pt}{1}{}{\cdot}1\above{1pt}{1pt}{-2}{}5\above{1pt}{1pt}{-}{}$}\EasyButWeakLineBreak%
{${144}\above{1pt}{1pt}{b}{2}{10}\above{1pt}{1pt}{l}{2}{16}\above{1pt}{1pt}{r}{2}{90}\above{1pt}{1pt}{b}{2}{16}\above{1pt}{1pt}{b}{2}{10}\above{1pt}{1pt}{l}{2}{144}\above{1pt}{1pt}{r}{2}{2}\above{1pt}{1pt}{b}{2}$}%
\nopagebreak\par%
\nopagebreak\par\leavevmode%
{$\left[\!\llap{\phantom{%
\begingroup \smaller\smaller\smaller\begin{tabular}{@{}c@{}}%
0\\0\\0
\end{tabular}\endgroup%
}}\right.$}%
\begingroup \smaller\smaller\smaller\begin{tabular}{@{}c@{}}%
-633878640\\-7824240\\146880
\end{tabular}\endgroup%
\kern3pt%
\begingroup \smaller\smaller\smaller\begin{tabular}{@{}c@{}}%
-7824240\\-96578\\1813
\end{tabular}\endgroup%
\kern3pt%
\begingroup \smaller\smaller\smaller\begin{tabular}{@{}c@{}}%
146880\\1813\\-34
\end{tabular}\endgroup%
{$\left.\llap{\phantom{%
\begingroup \smaller\smaller\smaller\begin{tabular}{@{}c@{}}%
0\\0\\0
\end{tabular}\endgroup%
}}\!\right]$}%
\EasyButWeakLineBreak%
{$\left[\!\llap{\phantom{%
\begingroup \smaller\smaller\smaller\begin{tabular}{@{}c@{}}%
0\\0\\0
\end{tabular}\endgroup%
}}\right.$}%
\begingroup \smaller\smaller\smaller\begin{tabular}{@{}c@{}}%
109\\-8856\\-1368
\end{tabular}\endgroup%
\HardButStrongLineBreak\kern3pt%
\begingroup \smaller\smaller\smaller\begin{tabular}{@{}c@{}}%
61\\-4955\\-705
\end{tabular}\endgroup%
\HardButStrongLineBreak\kern3pt%
\begingroup \smaller\smaller\smaller\begin{tabular}{@{}c@{}}%
131\\-10640\\-1456
\end{tabular}\endgroup%
\HardButStrongLineBreak\kern3pt%
\begingroup \smaller\smaller\smaller\begin{tabular}{@{}c@{}}%
128\\-10395\\-1350
\end{tabular}\endgroup%
\HardButStrongLineBreak\kern3pt%
\begingroup \smaller\smaller\smaller\begin{tabular}{@{}c@{}}%
7\\-568\\-48
\end{tabular}\endgroup%
\HardButStrongLineBreak\kern3pt%
\begingroup \smaller\smaller\smaller\begin{tabular}{@{}c@{}}%
-4\\325\\50
\end{tabular}\endgroup%
\HardButStrongLineBreak\kern3pt%
\begingroup \smaller\smaller\smaller\begin{tabular}{@{}c@{}}%
-55\\4464\\432
\end{tabular}\endgroup%
\HardButStrongLineBreak\kern3pt%
\begingroup \smaller\smaller\smaller\begin{tabular}{@{}c@{}}%
-1\\81\\-1
\end{tabular}\endgroup%
{$\left.\llap{\phantom{%
\begingroup \smaller\smaller\smaller\begin{tabular}{@{}c@{}}%
0\\0\\0
\end{tabular}\endgroup%
}}\!\right]$}%

\medskip%
%
\leavevmode\llap{}%
$W_{302}$%
\qquad\llap{22} lattices, $\chi=108$%
\hfill%
$\infty\infty2\infty22\infty\infty\infty2\infty22\infty\rtimes C_{2}$%
\nopagebreak\smallskip\hrule\nopagebreak\medskip%
%
%
\leavevmode%
${L_{302.1}}$%
{} : {$1\above{1pt}{1pt}{2}{{\rm II}}4\above{1pt}{1pt}{1}{1}{\cdot}1\above{1pt}{1pt}{1}{}9\above{1pt}{1pt}{1}{}81\above{1pt}{1pt}{-}{}$}\spacer%
\instructions{2}%
\EasyButWeakLineBreak%
{${36}\above{1pt}{1pt}{9,5}{\infty a}{36}\above{1pt}{1pt}{9,4}{\infty}{36}\above{1pt}{1pt}{*}{2}{4}\above{1pt}{1pt}{3,2}{\infty b}{4}\above{1pt}{1pt}{r}{2}{162}\above{1pt}{1pt}{l}{2}{36}\above{1pt}{1pt}{9,1}{\infty}$}\relax$\,(\times2)$%
\nopagebreak\par%
\nopagebreak\par\leavevmode%
{$\left[\!\llap{\phantom{%
\begingroup \smaller\smaller\smaller
\endgroup%
}}\!\right]$}%
%
%
\hbox{}\par\smallskip%
%
%
\leavevmode%
${L_{302.2}}$%
{} : {$1\above{1pt}{1pt}{-2}{2}8\above{1pt}{1pt}{-}{3}{\cdot}1\above{1pt}{1pt}{-}{}9\above{1pt}{1pt}{-}{}81\above{1pt}{1pt}{1}{}$}\spacer%
\instructions{2}%
\EasyButWeakLineBreak%
{${18}\above{1pt}{1pt}{36,23}{\infty b}{72}\above{1pt}{1pt}{36,13}{\infty z}{18}\above{1pt}{1pt}{b}{2}{2}\above{1pt}{1pt}{12,11}{\infty a}{8}\above{1pt}{1pt}{s}{2}{324}\above{1pt}{1pt}{s}{2}{72}\above{1pt}{1pt}{36,1}{\infty z}$}\relax$\,(\times2)$%
\nopagebreak\par%
\nopagebreak\par\leavevmode%
{$\left[\!\llap{\phantom{%
\begingroup \smaller\smaller\smaller
\endgroup%
}}\!\right]$}%
%
%
\hbox{}\par\smallskip%
%
%
\leavevmode%
${L_{302.3}}$%
{} : {$1\above{1pt}{1pt}{2}{2}8\above{1pt}{1pt}{1}{7}{\cdot}1\above{1pt}{1pt}{-}{}9\above{1pt}{1pt}{-}{}81\above{1pt}{1pt}{1}{}$}\spacer%
\instructions{m}%
\EasyButWeakLineBreak%
{${18}\above{1pt}{1pt}{36,23}{\infty a}{72}\above{1pt}{1pt}{36,31}{\infty z}{18}\above{1pt}{1pt}{s}{2}{2}\above{1pt}{1pt}{12,11}{\infty b}{8}\above{1pt}{1pt}{l}{2}{81}\above{1pt}{1pt}{r}{2}{72}\above{1pt}{1pt}{36,19}{\infty z}$}\relax$\,(\times2)$%
\nopagebreak\par%
\nopagebreak\par\leavevmode%
{$\left[\!\llap{\phantom{%
\begingroup \smaller\smaller\smaller
\endgroup%
}}\!\right]$}%

\medskip%
%
\leavevmode\llap{}%
$W_{303}$%
\qquad\llap{16} lattices, $\chi=24$%
\hfill%
$22222222\rtimes C_{2}$%
\nopagebreak\smallskip\hrule\nopagebreak\medskip%
%
%
\leavevmode%
${L_{303.1}}$%
{} : {$1\above{1pt}{1pt}{-2}{{\rm II}}8\above{1pt}{1pt}{-}{3}{\cdot}1\above{1pt}{1pt}{2}{}3\above{1pt}{1pt}{-}{}{\cdot}1\above{1pt}{1pt}{1}{}5\above{1pt}{1pt}{-}{}25\above{1pt}{1pt}{-}{}$}\spacer%
\instructions{2}%
\EasyButWeakLineBreak%
{${24}\above{1pt}{1pt}{r}{2}{50}\above{1pt}{1pt}{s}{2}{6}\above{1pt}{1pt}{b}{2}{10}\above{1pt}{1pt}{l}{2}$}\relax$\,(\times2)$%
\nopagebreak\par%
\nopagebreak\par\leavevmode%
{$\left[\!\llap{\phantom{%
\begingroup \smaller\smaller\smaller\begin{tabular}{@{}c@{}}%
0\\0\\0
\end{tabular}\endgroup%
}}\right.$}%
\begingroup \smaller\smaller\smaller\begin{tabular}{@{}c@{}}%
-877800\\-352200\\4800
\end{tabular}\endgroup%
\kern3pt%
\begingroup \smaller\smaller\smaller\begin{tabular}{@{}c@{}}%
-352200\\-141310\\1925
\end{tabular}\endgroup%
\kern3pt%
\begingroup \smaller\smaller\smaller\begin{tabular}{@{}c@{}}%
4800\\1925\\-26
\end{tabular}\endgroup%
{$\left.\llap{\phantom{%
\begingroup \smaller\smaller\smaller\begin{tabular}{@{}c@{}}%
0\\0\\0
\end{tabular}\endgroup%
}}\!\right]$}%
\hfil\penalty500%
{$\left[\!\llap{\phantom{%
\begingroup \smaller\smaller\smaller\begin{tabular}{@{}c@{}}%
0\\0\\0
\end{tabular}\endgroup%
}}\right.$}%
\begingroup \smaller\smaller\smaller\begin{tabular}{@{}c@{}}%
-13201\\34800\\138000
\end{tabular}\endgroup%
\kern3pt%
\begingroup \smaller\smaller\smaller\begin{tabular}{@{}c@{}}%
-5313\\14006\\55545
\end{tabular}\endgroup%
\kern3pt%
\begingroup \smaller\smaller\smaller\begin{tabular}{@{}c@{}}%
77\\-203\\-806
\end{tabular}\endgroup%
{$\left.\llap{\phantom{%
\begingroup \smaller\smaller\smaller\begin{tabular}{@{}c@{}}%
0\\0\\0
\end{tabular}\endgroup%
}}\!\right]$}%
\EasyButWeakLineBreak%
{$\left[\!\llap{\phantom{%
\begingroup \smaller\smaller\smaller\begin{tabular}{@{}c@{}}%
0\\0\\0
\end{tabular}\endgroup%
}}\right.$}%
\begingroup \smaller\smaller\smaller\begin{tabular}{@{}c@{}}%
-55\\144\\504
\end{tabular}\endgroup%
\HardButStrongLineBreak\kern3pt%
\begingroup \smaller\smaller\smaller\begin{tabular}{@{}c@{}}%
-2\\5\\0
\end{tabular}\endgroup%
\HardButStrongLineBreak\kern3pt%
\begingroup \smaller\smaller\smaller\begin{tabular}{@{}c@{}}%
8\\-21\\-78
\end{tabular}\endgroup%
\HardButStrongLineBreak\kern3pt%
\begingroup \smaller\smaller\smaller\begin{tabular}{@{}c@{}}%
5\\-13\\-40
\end{tabular}\endgroup%
{$\left.\llap{\phantom{%
\begingroup \smaller\smaller\smaller\begin{tabular}{@{}c@{}}%
0\\0\\0
\end{tabular}\endgroup%
}}\!\right]$}%

\medskip%
%
\leavevmode\llap{}%
$W_{304}$%
\qquad\llap{36} lattices, $\chi=36$%
\hfill%
$2|22\slashtwo22|22\slashtwo2\rtimes D_{4}$%
\nopagebreak\smallskip\hrule\nopagebreak\medskip%
%
%
\leavevmode%
${L_{304.1}}$%
{} : {$1\above{1pt}{1pt}{-2}{2}8\above{1pt}{1pt}{-}{5}{\cdot}1\above{1pt}{1pt}{-}{}3\above{1pt}{1pt}{-}{}9\above{1pt}{1pt}{-}{}{\cdot}1\above{1pt}{1pt}{-}{}5\above{1pt}{1pt}{1}{}25\above{1pt}{1pt}{-}{}$}\spacer%
\instructions{23,3m,3,2,m}%
\EasyButWeakLineBreak%
{${200}\above{1pt}{1pt}{*}{2}{180}\above{1pt}{1pt}{*}{2}{8}\above{1pt}{1pt}{b}{2}{450}\above{1pt}{1pt}{b}{2}{2}\above{1pt}{1pt}{b}{2}{1800}\above{1pt}{1pt}{*}{2}{20}\above{1pt}{1pt}{*}{2}{72}\above{1pt}{1pt}{b}{2}{50}\above{1pt}{1pt}{b}{2}{18}\above{1pt}{1pt}{b}{2}$}%
\nopagebreak\par%
\nopagebreak\par\leavevmode%
{$\left[\!\llap{\phantom{%
\begingroup \smaller\smaller\smaller
\endgroup%
}}\!\right]$}%
%
%
\hbox{}\par\smallskip%
%
%
\leavevmode%
${L_{304.2}}$%
{} : {$1\above{1pt}{1pt}{2}{2}16\above{1pt}{1pt}{1}{1}{\cdot}1\above{1pt}{1pt}{1}{}3\above{1pt}{1pt}{1}{}9\above{1pt}{1pt}{1}{}{\cdot}1\above{1pt}{1pt}{1}{}5\above{1pt}{1pt}{-}{}25\above{1pt}{1pt}{1}{}$}\spacer%
\instructions{5,3m,3,m}%
\EasyButWeakLineBreak%
{${400}\above{1pt}{1pt}{b}{2}{90}\above{1pt}{1pt}{l}{2}{16}\above{1pt}{1pt}{}{2}{225}\above{1pt}{1pt}{}{2}{1}\above{1pt}{1pt}{}{2}{3600}\above{1pt}{1pt}{r}{2}{10}\above{1pt}{1pt}{b}{2}{144}\above{1pt}{1pt}{*}{2}{100}\above{1pt}{1pt}{*}{2}{36}\above{1pt}{1pt}{*}{2}$}%
\nopagebreak\par%
shares genus with 3-dual${}\iso{}$5-dual; isometric to own %
3.5-dual\nopagebreak\par%
\nopagebreak\par\leavevmode%
{$\left[\!\llap{\phantom{%
\begingroup \smaller\smaller\smaller
\endgroup%
}}\!\right]$}%

\medskip%
%
\leavevmode\llap{}%
$W_{305}$%
\qquad\llap{46} lattices, $\chi=36$%
\hfill%
$2\slashtwo22|22\slashtwo22|2\rtimes D_{4}$%
\nopagebreak\smallskip\hrule\nopagebreak\medskip%
%
%
\leavevmode%
${L_{305.1}}$%
{} : {$1\above{1pt}{1pt}{-2}{4}8\above{1pt}{1pt}{-}{3}{\cdot}1\above{1pt}{1pt}{2}{}3\above{1pt}{1pt}{-}{}{\cdot}1\above{1pt}{1pt}{1}{}5\above{1pt}{1pt}{1}{}25\above{1pt}{1pt}{1}{}$}\spacer%
\instructions{2}%
\EasyButWeakLineBreak%
{${24}\above{1pt}{1pt}{}{2}{25}\above{1pt}{1pt}{r}{2}{4}\above{1pt}{1pt}{*}{2}{600}\above{1pt}{1pt}{*}{2}{20}\above{1pt}{1pt}{*}{2}{24}\above{1pt}{1pt}{*}{2}{100}\above{1pt}{1pt}{l}{2}{1}\above{1pt}{1pt}{}{2}{600}\above{1pt}{1pt}{}{2}{5}\above{1pt}{1pt}{}{2}$}%
\nopagebreak\par%
\nopagebreak\par\leavevmode%
{$\left[\!\llap{\phantom{%
\begingroup \smaller\smaller\smaller
\endgroup%
}}\!\right]$}%
%
%
\hbox{}\par\smallskip%
%
%
\leavevmode%
${L_{305.2}}$%
{} : {$1\above{1pt}{1pt}{2}{2}8\above{1pt}{1pt}{1}{1}{\cdot}1\above{1pt}{1pt}{2}{}3\above{1pt}{1pt}{-}{}{\cdot}1\above{1pt}{1pt}{1}{}5\above{1pt}{1pt}{1}{}25\above{1pt}{1pt}{1}{}$}\spacer%
\instructions{m}%
\EasyButWeakLineBreak%
{${24}\above{1pt}{1pt}{l}{2}{25}\above{1pt}{1pt}{}{2}{1}\above{1pt}{1pt}{r}{2}{600}\above{1pt}{1pt}{l}{2}{5}\above{1pt}{1pt}{r}{2}$}\relax$\,(\times2)$%
\nopagebreak\par%
\nopagebreak\par\leavevmode%
{$\left[\!\llap{\phantom{%
\begingroup \smaller\smaller\smaller
\endgroup%
}}\!\right]$}%
%
%
\hbox{}\par\smallskip%
%
%
\leavevmode%
${L_{305.3}}$%
{} : {$1\above{1pt}{1pt}{-2}{2}8\above{1pt}{1pt}{-}{5}{\cdot}1\above{1pt}{1pt}{2}{}3\above{1pt}{1pt}{-}{}{\cdot}1\above{1pt}{1pt}{1}{}5\above{1pt}{1pt}{1}{}25\above{1pt}{1pt}{1}{}$}\EasyButWeakLineBreak%
{${24}\above{1pt}{1pt}{s}{2}{100}\above{1pt}{1pt}{*}{2}{4}\above{1pt}{1pt}{s}{2}{600}\above{1pt}{1pt}{s}{2}{20}\above{1pt}{1pt}{s}{2}$}\relax$\,(\times2)$%
\nopagebreak\par%
\nopagebreak\par\leavevmode%
{$\left[\!\llap{\phantom{%
\begingroup \smaller\smaller\smaller
\endgroup%
}}\!\right]$}%
%
%
\hbox{}\par\smallskip%
%
%
\leavevmode%
${L_{305.4}}$%
{} : {$[1\above{1pt}{1pt}{-}{}2\above{1pt}{1pt}{-}{}]\above{1pt}{1pt}{}{0}16\above{1pt}{1pt}{1}{7}{\cdot}1\above{1pt}{1pt}{2}{}3\above{1pt}{1pt}{-}{}{\cdot}1\above{1pt}{1pt}{1}{}5\above{1pt}{1pt}{1}{}25\above{1pt}{1pt}{1}{}$}\spacer%
\instructions{2}%
\EasyButWeakLineBreak%
{${600}\above{1pt}{1pt}{s}{2}{16}\above{1pt}{1pt}{*}{2}{100}\above{1pt}{1pt}{l}{2}{6}\above{1pt}{1pt}{}{2}{5}\above{1pt}{1pt}{}{2}{150}\above{1pt}{1pt}{r}{2}{4}\above{1pt}{1pt}{*}{2}{400}\above{1pt}{1pt}{s}{2}{24}\above{1pt}{1pt}{*}{2}{80}\above{1pt}{1pt}{*}{2}$}%
\nopagebreak\par%
\nopagebreak\par\leavevmode%
{$\left[\!\llap{\phantom{%
\begingroup \smaller\smaller\smaller
\endgroup%
}}\!\right]$}%
%
%
\hbox{}\par\smallskip%
%
%
\leavevmode%
${L_{305.5}}$%
{} : {$[1\above{1pt}{1pt}{-}{}2\above{1pt}{1pt}{1}{}]\above{1pt}{1pt}{}{6}16\above{1pt}{1pt}{-}{5}{\cdot}1\above{1pt}{1pt}{2}{}3\above{1pt}{1pt}{-}{}{\cdot}1\above{1pt}{1pt}{1}{}5\above{1pt}{1pt}{1}{}25\above{1pt}{1pt}{1}{}$}\spacer%
\instructions{m}%
\EasyButWeakLineBreak%
{${150}\above{1pt}{1pt}{r}{2}{16}\above{1pt}{1pt}{s}{2}{100}\above{1pt}{1pt}{*}{2}{24}\above{1pt}{1pt}{l}{2}{5}\above{1pt}{1pt}{r}{2}{600}\above{1pt}{1pt}{*}{2}{4}\above{1pt}{1pt}{s}{2}{400}\above{1pt}{1pt}{l}{2}{6}\above{1pt}{1pt}{}{2}{80}\above{1pt}{1pt}{}{2}$}%
\nopagebreak\par%
\nopagebreak\par\leavevmode%
{$\left[\!\llap{\phantom{%
\begingroup \smaller\smaller\smaller
\endgroup%
}}\!\right]$}%
%
%
\hbox{}\par\smallskip%
%
%
\leavevmode%
${L_{305.6}}$%
{} : {$[1\above{1pt}{1pt}{1}{}2\above{1pt}{1pt}{-}{}]\above{1pt}{1pt}{}{4}16\above{1pt}{1pt}{-}{3}{\cdot}1\above{1pt}{1pt}{2}{}3\above{1pt}{1pt}{-}{}{\cdot}1\above{1pt}{1pt}{1}{}5\above{1pt}{1pt}{1}{}25\above{1pt}{1pt}{1}{}$}\spacer%
\instructions{m}%
\EasyButWeakLineBreak%
{${600}\above{1pt}{1pt}{*}{2}{16}\above{1pt}{1pt}{l}{2}{25}\above{1pt}{1pt}{}{2}{6}\above{1pt}{1pt}{r}{2}{20}\above{1pt}{1pt}{l}{2}{150}\above{1pt}{1pt}{}{2}{1}\above{1pt}{1pt}{r}{2}{400}\above{1pt}{1pt}{*}{2}{24}\above{1pt}{1pt}{s}{2}{80}\above{1pt}{1pt}{s}{2}$}%
\nopagebreak\par%
\nopagebreak\par\leavevmode%
{$\left[\!\llap{\phantom{%
\begingroup \smaller\smaller\smaller
\endgroup%
}}\!\right]$}%
%
%
\hbox{}\par\smallskip%
%
%
\leavevmode%
${L_{305.7}}$%
{} : {$[1\above{1pt}{1pt}{1}{}2\above{1pt}{1pt}{1}{}]\above{1pt}{1pt}{}{2}16\above{1pt}{1pt}{1}{1}{\cdot}1\above{1pt}{1pt}{2}{}3\above{1pt}{1pt}{-}{}{\cdot}1\above{1pt}{1pt}{1}{}5\above{1pt}{1pt}{1}{}25\above{1pt}{1pt}{1}{}$}\EasyButWeakLineBreak%
{${150}\above{1pt}{1pt}{}{2}{16}\above{1pt}{1pt}{}{2}{25}\above{1pt}{1pt}{r}{2}{24}\above{1pt}{1pt}{*}{2}{20}\above{1pt}{1pt}{*}{2}{600}\above{1pt}{1pt}{l}{2}{1}\above{1pt}{1pt}{}{2}{400}\above{1pt}{1pt}{}{2}{6}\above{1pt}{1pt}{r}{2}{80}\above{1pt}{1pt}{l}{2}$}%
\nopagebreak\par%
\nopagebreak\par\leavevmode%
{$\left[\!\llap{\phantom{%
\begingroup \smaller\smaller\smaller
\endgroup%
}}\!\right]$}%

\medskip%
%
\leavevmode\llap{}%
$W_{306}$%
\qquad\llap{32} lattices, $\chi=96$%
\hfill%
$2\infty|\infty22\slashinfty22\infty|\infty22\slashinfty2\rtimes D_{4}$%
\nopagebreak\smallskip\hrule\nopagebreak\medskip%
%
%
\leavevmode%
${L_{306.1}}$%
{} : {$1\above{1pt}{1pt}{1}{7}4\above{1pt}{1pt}{1}{7}64\above{1pt}{1pt}{-}{5}{\cdot}1\above{1pt}{1pt}{-}{}3\above{1pt}{1pt}{1}{}9\above{1pt}{1pt}{1}{}$}\spacer%
\instructions{3}%
\EasyButWeakLineBreak%
{${576}\above{1pt}{1pt}{*}{2}{48}\above{1pt}{1pt}{48,25}{\infty z}{12}\above{1pt}{1pt}{48,13}{\infty}{48}\above{1pt}{1pt}{s}{2}{576}\above{1pt}{1pt}{*}{2}{12}\above{1pt}{1pt}{24,13}{\infty z}{3}\above{1pt}{1pt}{r}{2}$}\relax$\,(\times2)$%
\nopagebreak\par%
\nopagebreak\par\leavevmode%
{$\left[\!\llap{\phantom{%
\begingroup \smaller\smaller\smaller
\endgroup%
}}\!\right]$}%
%
%
\hbox{}\par\smallskip%
%
%
\leavevmode%
${L_{306.2}}$%
{} : {$1\above{1pt}{1pt}{1}{1}4\above{1pt}{1pt}{1}{7}64\above{1pt}{1pt}{-}{3}{\cdot}1\above{1pt}{1pt}{-}{}3\above{1pt}{1pt}{1}{}9\above{1pt}{1pt}{1}{}$}\spacer%
\instructions{3}%
\EasyButWeakLineBreak%
{${36}\above{1pt}{1pt}{*}{2}{48}\above{1pt}{1pt}{48,7}{\infty z}{12}\above{1pt}{1pt}{48,43}{\infty}{48}\above{1pt}{1pt}{l}{2}{9}\above{1pt}{1pt}{}{2}{192}\above{1pt}{1pt}{12,7}{\infty}{192}\above{1pt}{1pt}{s}{2}$}\relax$\,(\times2)$%
\nopagebreak\par%
\nopagebreak\par\leavevmode%
{$\left[\!\llap{\phantom{%
\begingroup \smaller\smaller\smaller
\endgroup%
}}\!\right]$}%
%
%
\hbox{}\par\smallskip%
%
%
\leavevmode%
${L_{306.3}}$%
{} : {$1\above{1pt}{1pt}{-}{3}4\above{1pt}{1pt}{1}{7}64\above{1pt}{1pt}{1}{1}{\cdot}1\above{1pt}{1pt}{-}{}3\above{1pt}{1pt}{1}{}9\above{1pt}{1pt}{1}{}$}\spacer%
\instructions{3}%
\EasyButWeakLineBreak%
{${576}\above{1pt}{1pt}{l}{2}{12}\above{1pt}{1pt}{48,25}{\infty}{48}\above{1pt}{1pt}{48,37}{\infty z}{12}\above{1pt}{1pt}{}{2}{576}\above{1pt}{1pt}{}{2}{3}\above{1pt}{1pt}{24,1}{\infty}{12}\above{1pt}{1pt}{s}{2}$}\relax$\,(\times2)$%
\nopagebreak\par%
\nopagebreak\par\leavevmode%
{$\left[\!\llap{\phantom{%
\begingroup \smaller\smaller\smaller
\endgroup%
}}\!\right]$}%
%
%
\hbox{}\par\smallskip%
%
%
\leavevmode%
${L_{306.4}}$%
{} : {$1\above{1pt}{1pt}{-}{5}4\above{1pt}{1pt}{1}{7}64\above{1pt}{1pt}{1}{7}{\cdot}1\above{1pt}{1pt}{-}{}3\above{1pt}{1pt}{1}{}9\above{1pt}{1pt}{1}{}$}\spacer%
\instructions{3}%
\EasyButWeakLineBreak%
{${36}\above{1pt}{1pt}{l}{2}{12}\above{1pt}{1pt}{48,7}{\infty}{48}\above{1pt}{1pt}{48,19}{\infty z}{12}\above{1pt}{1pt}{}{2}{9}\above{1pt}{1pt}{r}{2}{192}\above{1pt}{1pt}{24,7}{\infty z}{192}\above{1pt}{1pt}{*}{2}$}\relax$\,(\times2)$%
\nopagebreak\par%
\nopagebreak\par\leavevmode%
{$\left[\!\llap{\phantom{%
\begingroup \smaller\smaller\smaller
\endgroup%
}}\!\right]$}%

\medskip%
%
\leavevmode\llap{}%
$W_{307}$%
\qquad\llap{16} lattices, $\chi=48$%
\hfill%
$\slashtwo2|2\slashtwo2|2\slashtwo2|2\slashtwo2|2\rtimes D_{8}$%
\nopagebreak\smallskip\hrule\nopagebreak\medskip%
%
%
\leavevmode%
${L_{307.1}}$%
{} : {$1\above{1pt}{1pt}{1}{1}4\above{1pt}{1pt}{1}{1}64\above{1pt}{1pt}{-}{5}{\cdot}1\above{1pt}{1pt}{1}{}3\above{1pt}{1pt}{-}{}9\above{1pt}{1pt}{1}{}$}\spacer%
\instructions{3}%
\EasyButWeakLineBreak%
{${64}\above{1pt}{1pt}{*}{2}{36}\above{1pt}{1pt}{l}{2}{4}\above{1pt}{1pt}{}{2}{9}\above{1pt}{1pt}{r}{2}{64}\above{1pt}{1pt}{*}{2}{144}\above{1pt}{1pt}{s}{2}$}\relax$\,(\times2)$%
\nopagebreak\par%
shares genus with 3-dual\nopagebreak\par%
\nopagebreak\par\leavevmode%
{$\left[\!\llap{\phantom{%
\begingroup \smaller\smaller\smaller
\endgroup%
}}\!\right]$}%
%
%
\hbox{}\par\smallskip%
%
%
\leavevmode%
${L_{307.2}}$%
{} : {$1\above{1pt}{1pt}{-}{5}4\above{1pt}{1pt}{1}{1}64\above{1pt}{1pt}{1}{1}{\cdot}1\above{1pt}{1pt}{1}{}3\above{1pt}{1pt}{-}{}9\above{1pt}{1pt}{1}{}$}\spacer%
\instructions{3}%
\EasyButWeakLineBreak%
{${64}\above{1pt}{1pt}{s}{2}{36}\above{1pt}{1pt}{*}{2}{16}\above{1pt}{1pt}{l}{2}{9}\above{1pt}{1pt}{}{2}{64}\above{1pt}{1pt}{}{2}{36}\above{1pt}{1pt}{r}{2}$}\relax$\,(\times2)$%
\nopagebreak\par%
shares genus with 3-dual\nopagebreak\par%
\nopagebreak\par\leavevmode%
{$\left[\!\llap{\phantom{%
\begingroup \smaller\smaller\smaller
\endgroup%
}}\!\right]$}%

\medskip%
%
\leavevmode\llap{}%
$W_{308}$%
\qquad\llap{16} lattices, $\chi=48$%
\hfill%
$6222362223\rtimes C_{2}$%
\nopagebreak\smallskip\hrule\nopagebreak\medskip%
%
%
\leavevmode%
${L_{308.1}}$%
{} : {$1\above{1pt}{1pt}{-2}{{\rm II}}8\above{1pt}{1pt}{1}{1}{\cdot}1\above{1pt}{1pt}{-}{}3\above{1pt}{1pt}{-}{}27\above{1pt}{1pt}{-}{}{\cdot}1\above{1pt}{1pt}{-2}{}5\above{1pt}{1pt}{1}{}$}\spacer%
\instructions{2}%
\EasyButWeakLineBreak%
{${6}\above{1pt}{1pt}{}{6}{2}\above{1pt}{1pt}{s}{2}{54}\above{1pt}{1pt}{l}{2}{8}\above{1pt}{1pt}{r}{2}{6}\above{1pt}{1pt}{-}{3}$}\relax$\,(\times2)$%
\nopagebreak\par%
\nopagebreak\par\leavevmode%
{$\left[\!\llap{\phantom{%
\begingroup \smaller\smaller\smaller
\endgroup%
}}\!\right]$}%

\medskip%
%
\leavevmode\llap{}%
$W_{309}$%
\qquad\llap{120} lattices, $\chi=48$%
\hfill%
$222|222222|222\rtimes D_{2}$%
\nopagebreak\smallskip\hrule\nopagebreak\medskip%
%
%
\leavevmode%
${L_{309.1}}$%
{} : {$1\above{1pt}{1pt}{2}{0}8\above{1pt}{1pt}{1}{1}{\cdot}1\above{1pt}{1pt}{2}{}3\above{1pt}{1pt}{-}{}{\cdot}1\above{1pt}{1pt}{-2}{}5\above{1pt}{1pt}{1}{}{\cdot}1\above{1pt}{1pt}{-2}{}7\above{1pt}{1pt}{1}{}$}\EasyButWeakLineBreak%
{${105}\above{1pt}{1pt}{}{2}{8}\above{1pt}{1pt}{}{2}{5}\above{1pt}{1pt}{r}{2}{56}\above{1pt}{1pt}{s}{2}{20}\above{1pt}{1pt}{*}{2}{8}\above{1pt}{1pt}{*}{2}{420}\above{1pt}{1pt}{l}{2}{1}\above{1pt}{1pt}{}{2}{7}\above{1pt}{1pt}{r}{2}{24}\above{1pt}{1pt}{s}{2}{28}\above{1pt}{1pt}{*}{2}{4}\above{1pt}{1pt}{l}{2}$}%
\nopagebreak\par%
\nopagebreak\par\leavevmode%
{$\left[\!\llap{\phantom{%
\begingroup \smaller\smaller\smaller
\endgroup%
}}\!\right]$}%
%
%
\hbox{}\par\smallskip%
%
%
\leavevmode%
${L_{309.2}}$%
{} : {$[1\above{1pt}{1pt}{1}{}2\above{1pt}{1pt}{1}{}]\above{1pt}{1pt}{}{2}16\above{1pt}{1pt}{1}{7}{\cdot}1\above{1pt}{1pt}{2}{}3\above{1pt}{1pt}{-}{}{\cdot}1\above{1pt}{1pt}{-2}{}5\above{1pt}{1pt}{1}{}{\cdot}1\above{1pt}{1pt}{-2}{}7\above{1pt}{1pt}{1}{}$}\spacer%
\instructions{2}%
\EasyButWeakLineBreak%
{${105}\above{1pt}{1pt}{}{2}{2}\above{1pt}{1pt}{r}{2}{20}\above{1pt}{1pt}{*}{2}{56}\above{1pt}{1pt}{*}{2}{80}\above{1pt}{1pt}{s}{2}{8}\above{1pt}{1pt}{*}{2}{1680}\above{1pt}{1pt}{l}{2}{1}\above{1pt}{1pt}{}{2}{112}\above{1pt}{1pt}{}{2}{6}\above{1pt}{1pt}{r}{2}{28}\above{1pt}{1pt}{*}{2}{16}\above{1pt}{1pt}{l}{2}$}%
\nopagebreak\par%
\nopagebreak\par\leavevmode%
{$\left[\!\llap{\phantom{%
\begingroup \smaller\smaller\smaller
\endgroup%
}}\!\right]$}%
%
%
\hbox{}\par\smallskip%
%
%
\leavevmode%
${L_{309.3}}$%
{} : {$[1\above{1pt}{1pt}{-}{}2\above{1pt}{1pt}{1}{}]\above{1pt}{1pt}{}{6}16\above{1pt}{1pt}{-}{3}{\cdot}1\above{1pt}{1pt}{2}{}3\above{1pt}{1pt}{-}{}{\cdot}1\above{1pt}{1pt}{-2}{}5\above{1pt}{1pt}{1}{}{\cdot}1\above{1pt}{1pt}{-2}{}7\above{1pt}{1pt}{1}{}$}\spacer%
\instructions{m}%
\EasyButWeakLineBreak%
{${420}\above{1pt}{1pt}{l}{2}{2}\above{1pt}{1pt}{}{2}{5}\above{1pt}{1pt}{r}{2}{56}\above{1pt}{1pt}{s}{2}{80}\above{1pt}{1pt}{*}{2}{8}\above{1pt}{1pt}{s}{2}{1680}\above{1pt}{1pt}{*}{2}{4}\above{1pt}{1pt}{s}{2}{112}\above{1pt}{1pt}{l}{2}{6}\above{1pt}{1pt}{}{2}{7}\above{1pt}{1pt}{r}{2}{16}\above{1pt}{1pt}{*}{2}$}%
\nopagebreak\par%
\nopagebreak\par\leavevmode%
{$\left[\!\llap{\phantom{%
\begingroup \smaller\smaller\smaller
\endgroup%
}}\!\right]$}%
%
%
\hbox{}\par\smallskip%
%
%
\leavevmode%
${L_{309.4}}$%
{} : {$[1\above{1pt}{1pt}{-}{}2\above{1pt}{1pt}{1}{}]\above{1pt}{1pt}{}{4}16\above{1pt}{1pt}{-}{5}{\cdot}1\above{1pt}{1pt}{2}{}3\above{1pt}{1pt}{-}{}{\cdot}1\above{1pt}{1pt}{-2}{}5\above{1pt}{1pt}{1}{}{\cdot}1\above{1pt}{1pt}{-2}{}7\above{1pt}{1pt}{1}{}$}\spacer%
\instructions{m}%
\EasyButWeakLineBreak%
{${420}\above{1pt}{1pt}{*}{2}{8}\above{1pt}{1pt}{l}{2}{5}\above{1pt}{1pt}{}{2}{14}\above{1pt}{1pt}{}{2}{80}\above{1pt}{1pt}{}{2}{2}\above{1pt}{1pt}{r}{2}{1680}\above{1pt}{1pt}{s}{2}{4}\above{1pt}{1pt}{*}{2}{112}\above{1pt}{1pt}{s}{2}{24}\above{1pt}{1pt}{*}{2}{28}\above{1pt}{1pt}{s}{2}{16}\above{1pt}{1pt}{s}{2}$}%
\nopagebreak\par%
\nopagebreak\par\leavevmode%
{$\left[\!\llap{\phantom{%
\begingroup \smaller\smaller\smaller
\endgroup%
}}\!\right]$}%
%
%
\hbox{}\par\smallskip%
%
%
\leavevmode%
${L_{309.5}}$%
{} : {$[1\above{1pt}{1pt}{1}{}2\above{1pt}{1pt}{1}{}]\above{1pt}{1pt}{}{0}16\above{1pt}{1pt}{1}{1}{\cdot}1\above{1pt}{1pt}{2}{}3\above{1pt}{1pt}{-}{}{\cdot}1\above{1pt}{1pt}{-2}{}5\above{1pt}{1pt}{1}{}{\cdot}1\above{1pt}{1pt}{-2}{}7\above{1pt}{1pt}{1}{}$}\EasyButWeakLineBreak%
{${105}\above{1pt}{1pt}{r}{2}{8}\above{1pt}{1pt}{*}{2}{20}\above{1pt}{1pt}{l}{2}{14}\above{1pt}{1pt}{r}{2}{80}\above{1pt}{1pt}{l}{2}{2}\above{1pt}{1pt}{}{2}{1680}\above{1pt}{1pt}{}{2}{1}\above{1pt}{1pt}{r}{2}{112}\above{1pt}{1pt}{*}{2}{24}\above{1pt}{1pt}{l}{2}{7}\above{1pt}{1pt}{}{2}{16}\above{1pt}{1pt}{}{2}$}%
\nopagebreak\par%
\nopagebreak\par\leavevmode%
{$\left[\!\llap{\phantom{%
\begingroup \smaller\smaller\smaller
\endgroup%
}}\!\right]$}%

\medskip%
%
\leavevmode\llap{}%
$W_{310}$%
\qquad\llap{120} lattices, $\chi=72$%
\hfill%
$2222|2222|2222|2222|\rtimes D_{4}$%
\nopagebreak\smallskip\hrule\nopagebreak\medskip%
%
%
\leavevmode%
${L_{310.1}}$%
{} : {$1\above{1pt}{1pt}{-2}{4}8\above{1pt}{1pt}{-}{5}{\cdot}1\above{1pt}{1pt}{2}{}3\above{1pt}{1pt}{-}{}{\cdot}1\above{1pt}{1pt}{-2}{}5\above{1pt}{1pt}{1}{}{\cdot}1\above{1pt}{1pt}{2}{}7\above{1pt}{1pt}{-}{}$}\EasyButWeakLineBreak%
{${280}\above{1pt}{1pt}{l}{2}{1}\above{1pt}{1pt}{}{2}{168}\above{1pt}{1pt}{}{2}{5}\above{1pt}{1pt}{r}{2}{24}\above{1pt}{1pt}{s}{2}{20}\above{1pt}{1pt}{*}{2}{168}\above{1pt}{1pt}{*}{2}{4}\above{1pt}{1pt}{s}{2}$}\relax$\,(\times2)$%
\nopagebreak\par%
\nopagebreak\par\leavevmode%
{$\left[\!\llap{\phantom{%
\begingroup \smaller\smaller\smaller
\endgroup%
}}\!\right]$}%
%
%
\hbox{}\par\smallskip%
%
%
\leavevmode%
${L_{310.2}}$%
{} : {$[1\above{1pt}{1pt}{-}{}2\above{1pt}{1pt}{1}{}]\above{1pt}{1pt}{}{2}16\above{1pt}{1pt}{-}{3}{\cdot}1\above{1pt}{1pt}{2}{}3\above{1pt}{1pt}{-}{}{\cdot}1\above{1pt}{1pt}{-2}{}5\above{1pt}{1pt}{1}{}{\cdot}1\above{1pt}{1pt}{2}{}7\above{1pt}{1pt}{-}{}$}\spacer%
\instructions{2}%
\EasyButWeakLineBreak%
{${280}\above{1pt}{1pt}{*}{2}{16}\above{1pt}{1pt}{s}{2}{168}\above{1pt}{1pt}{*}{2}{80}\above{1pt}{1pt}{s}{2}{24}\above{1pt}{1pt}{*}{2}{20}\above{1pt}{1pt}{l}{2}{42}\above{1pt}{1pt}{}{2}{1}\above{1pt}{1pt}{r}{2}$}\relax$\,(\times2)$%
\nopagebreak\par%
\nopagebreak\par\leavevmode%
{$\left[\!\llap{\phantom{%
\begingroup \smaller\smaller\smaller
\endgroup%
}}\!\right]$}%
%
%
\hbox{}\par\smallskip%
%
%
\leavevmode%
${L_{310.3}}$%
{} : {$[1\above{1pt}{1pt}{1}{}2\above{1pt}{1pt}{-}{}]\above{1pt}{1pt}{}{4}16\above{1pt}{1pt}{-}{5}{\cdot}1\above{1pt}{1pt}{2}{}3\above{1pt}{1pt}{-}{}{\cdot}1\above{1pt}{1pt}{-2}{}5\above{1pt}{1pt}{1}{}{\cdot}1\above{1pt}{1pt}{2}{}7\above{1pt}{1pt}{-}{}$}\spacer%
\instructions{m}%
\EasyButWeakLineBreak%
{${70}\above{1pt}{1pt}{r}{2}{16}\above{1pt}{1pt}{l}{2}{42}\above{1pt}{1pt}{}{2}{80}\above{1pt}{1pt}{}{2}{6}\above{1pt}{1pt}{r}{2}{20}\above{1pt}{1pt}{*}{2}{168}\above{1pt}{1pt}{l}{2}{1}\above{1pt}{1pt}{}{2}$}\relax$\,(\times2)$%
\nopagebreak\par%
\nopagebreak\par\leavevmode%
{$\left[\!\llap{\phantom{%
\begingroup \smaller\smaller\smaller
\endgroup%
}}\!\right]$}%
%
%
\hbox{}\par\smallskip%
%
%
\leavevmode%
${L_{310.4}}$%
{} : {$[1\above{1pt}{1pt}{-}{}2\above{1pt}{1pt}{-}{}]\above{1pt}{1pt}{}{0}16\above{1pt}{1pt}{1}{1}{\cdot}1\above{1pt}{1pt}{2}{}3\above{1pt}{1pt}{-}{}{\cdot}1\above{1pt}{1pt}{-2}{}5\above{1pt}{1pt}{1}{}{\cdot}1\above{1pt}{1pt}{2}{}7\above{1pt}{1pt}{-}{}$}\EasyButWeakLineBreak%
{${70}\above{1pt}{1pt}{}{2}{16}\above{1pt}{1pt}{}{2}{42}\above{1pt}{1pt}{r}{2}{80}\above{1pt}{1pt}{l}{2}{6}\above{1pt}{1pt}{}{2}{5}\above{1pt}{1pt}{r}{2}{168}\above{1pt}{1pt}{*}{2}{4}\above{1pt}{1pt}{l}{2}$}\relax$\,(\times2)$%
\nopagebreak\par%
\nopagebreak\par\leavevmode%
{$\left[\!\llap{\phantom{%
\begingroup \smaller\smaller\smaller
\endgroup%
}}\!\right]$}%
%
%
\hbox{}\par\smallskip%
%
%
\leavevmode%
${L_{310.5}}$%
{} : {$[1\above{1pt}{1pt}{1}{}2\above{1pt}{1pt}{1}{}]\above{1pt}{1pt}{}{6}16\above{1pt}{1pt}{1}{7}{\cdot}1\above{1pt}{1pt}{2}{}3\above{1pt}{1pt}{-}{}{\cdot}1\above{1pt}{1pt}{-2}{}5\above{1pt}{1pt}{1}{}{\cdot}1\above{1pt}{1pt}{2}{}7\above{1pt}{1pt}{-}{}$}\spacer%
\instructions{m}%
\EasyButWeakLineBreak%
{${280}\above{1pt}{1pt}{s}{2}{16}\above{1pt}{1pt}{*}{2}{168}\above{1pt}{1pt}{s}{2}{80}\above{1pt}{1pt}{*}{2}{24}\above{1pt}{1pt}{l}{2}{5}\above{1pt}{1pt}{}{2}{42}\above{1pt}{1pt}{r}{2}{4}\above{1pt}{1pt}{*}{2}$}\relax$\,(\times2)$%
\nopagebreak\par%
\nopagebreak\par\leavevmode%
{$\left[\!\llap{\phantom{%
\begingroup \smaller\smaller\smaller
\endgroup%
}}\!\right]$}%

\medskip%
%
\leavevmode\llap{}%
$W_{311}$%
\qquad\llap{6} lattices, $\chi=36$%
\hfill%
$\slashtwo22|22\slashtwo22|22\rtimes D_{4}$%
\nopagebreak\smallskip\hrule\nopagebreak\medskip%
%
%
\leavevmode%
${L_{311.1}}$%
{} : {$1\above{1pt}{1pt}{-2}{{\rm II}}4\above{1pt}{1pt}{1}{7}{\cdot}1\above{1pt}{1pt}{1}{}5\above{1pt}{1pt}{-}{}25\above{1pt}{1pt}{1}{}{\cdot}1\above{1pt}{1pt}{2}{}7\above{1pt}{1pt}{1}{}$}\spacer%
\instructions{2}%
\EasyButWeakLineBreak%
{${100}\above{1pt}{1pt}{*}{2}{4}\above{1pt}{1pt}{b}{2}{350}\above{1pt}{1pt}{s}{2}{10}\above{1pt}{1pt}{s}{2}{14}\above{1pt}{1pt}{b}{2}$}\relax$\,(\times2)$%
\nopagebreak\par%
\nopagebreak\par\leavevmode%
{$\left[\!\llap{\phantom{%
\begingroup \smaller\smaller\smaller
\endgroup%
}}\!\right]$}%

\medskip%
%
\leavevmode\llap{}%
$W_{312}$%
\qquad\llap{6} lattices, $\chi=36$%
\hfill%
$\slashtwo22|22\slashtwo22|22\rtimes D_{4}$%
\nopagebreak\smallskip\hrule\nopagebreak\medskip%
%
%
\leavevmode%
${L_{312.1}}$%
{} : {$1\above{1pt}{1pt}{-2}{{\rm II}}4\above{1pt}{1pt}{1}{7}{\cdot}1\above{1pt}{1pt}{-}{}5\above{1pt}{1pt}{-}{}25\above{1pt}{1pt}{-}{}{\cdot}1\above{1pt}{1pt}{2}{}7\above{1pt}{1pt}{1}{}$}\spacer%
\instructions{2}%
\EasyButWeakLineBreak%
{${50}\above{1pt}{1pt}{b}{2}{2}\above{1pt}{1pt}{l}{2}{700}\above{1pt}{1pt}{r}{2}{10}\above{1pt}{1pt}{l}{2}{28}\above{1pt}{1pt}{r}{2}$}\relax$\,(\times2)$%
\nopagebreak\par%
\nopagebreak\par\leavevmode%
{$\left[\!\llap{\phantom{%
\begingroup \smaller\smaller\smaller
\endgroup%
}}\!\right]$}%

\medskip%
%
\leavevmode\llap{}%
$W_{313}$%
\qquad\llap{12} lattices, $\chi=12$%
\hfill%
$22|222|2\rtimes D_{2}$%
\nopagebreak\smallskip\hrule\nopagebreak\medskip%
%
%
\leavevmode%
${L_{313.1}}$%
{} : {$1\above{1pt}{1pt}{-2}{{\rm II}}4\above{1pt}{1pt}{-}{5}{\cdot}1\above{1pt}{1pt}{-}{}3\above{1pt}{1pt}{1}{}9\above{1pt}{1pt}{-}{}{\cdot}1\above{1pt}{1pt}{-2}{}5\above{1pt}{1pt}{-}{}{\cdot}1\above{1pt}{1pt}{2}{}7\above{1pt}{1pt}{-}{}$}\spacer%
\instructions{2}%
\EasyButWeakLineBreak%
{${140}\above{1pt}{1pt}{b}{2}{18}\above{1pt}{1pt}{l}{2}{84}\above{1pt}{1pt}{r}{2}{2}\above{1pt}{1pt}{b}{2}{1260}\above{1pt}{1pt}{*}{2}{12}\above{1pt}{1pt}{*}{2}$}%
\nopagebreak\par%
\nopagebreak\par\leavevmode%
{$\left[\!\llap{\phantom{%
\begingroup \smaller\smaller\smaller\begin{tabular}{@{}c@{}}%
0\\0\\0
\end{tabular}\endgroup%
}}\right.$}%
\begingroup \smaller\smaller\smaller\begin{tabular}{@{}c@{}}%
-7006860\\12600\\3780
\end{tabular}\endgroup%
\kern3pt%
\begingroup \smaller\smaller\smaller\begin{tabular}{@{}c@{}}%
12600\\-6\\-15
\end{tabular}\endgroup%
\kern3pt%
\begingroup \smaller\smaller\smaller\begin{tabular}{@{}c@{}}%
3780\\-15\\2
\end{tabular}\endgroup%
{$\left.\llap{\phantom{%
\begingroup \smaller\smaller\smaller\begin{tabular}{@{}c@{}}%
0\\0\\0
\end{tabular}\endgroup%
}}\!\right]$}%
\EasyButWeakLineBreak%
{$\left[\!\llap{\phantom{%
\begingroup \smaller\smaller\smaller\begin{tabular}{@{}c@{}}%
0\\0\\0
\end{tabular}\endgroup%
}}\right.$}%
\begingroup \smaller\smaller\smaller\begin{tabular}{@{}c@{}}%
-1\\-350\\-700
\end{tabular}\endgroup%
\HardButStrongLineBreak\kern3pt%
\begingroup \smaller\smaller\smaller\begin{tabular}{@{}c@{}}%
1\\345\\702
\end{tabular}\endgroup%
\HardButStrongLineBreak\kern3pt%
\begingroup \smaller\smaller\smaller\begin{tabular}{@{}c@{}}%
3\\1036\\2100
\end{tabular}\endgroup%
\HardButStrongLineBreak\kern3pt%
\begingroup \smaller\smaller\smaller\begin{tabular}{@{}c@{}}%
0\\0\\-1
\end{tabular}\endgroup%
\HardButStrongLineBreak\kern3pt%
\begingroup \smaller\smaller\smaller\begin{tabular}{@{}c@{}}%
-17\\-5880\\-11970
\end{tabular}\endgroup%
\HardButStrongLineBreak\kern3pt%
\begingroup \smaller\smaller\smaller\begin{tabular}{@{}c@{}}%
-1\\-346\\-702
\end{tabular}\endgroup%
{$\left.\llap{\phantom{%
\begingroup \smaller\smaller\smaller\begin{tabular}{@{}c@{}}%
0\\0\\0
\end{tabular}\endgroup%
}}\!\right]$}%

\medskip%
%
\leavevmode\llap{}%
$W_{314}$%
\qquad\llap{12} lattices, $\chi=24$%
\hfill%
$22|2222|22\rtimes D_{2}$%
\nopagebreak\smallskip\hrule\nopagebreak\medskip%
%
%
\leavevmode%
${L_{314.1}}$%
{} : {$1\above{1pt}{1pt}{-2}{{\rm II}}4\above{1pt}{1pt}{-}{5}{\cdot}1\above{1pt}{1pt}{-}{}3\above{1pt}{1pt}{1}{}9\above{1pt}{1pt}{-}{}{\cdot}1\above{1pt}{1pt}{2}{}5\above{1pt}{1pt}{1}{}{\cdot}1\above{1pt}{1pt}{-2}{}7\above{1pt}{1pt}{1}{}$}\spacer%
\instructions{2}%
\EasyButWeakLineBreak%
{${180}\above{1pt}{1pt}{r}{2}{14}\above{1pt}{1pt}{b}{2}{30}\above{1pt}{1pt}{b}{2}{126}\above{1pt}{1pt}{l}{2}{20}\above{1pt}{1pt}{r}{2}{18}\above{1pt}{1pt}{b}{2}{12}\above{1pt}{1pt}{b}{2}{2}\above{1pt}{1pt}{l}{2}$}%
\nopagebreak\par%
\nopagebreak\par\leavevmode%
{$\left[\!\llap{\phantom{%
\begingroup \smaller\smaller\smaller\begin{tabular}{@{}c@{}}%
0\\0\\0
\end{tabular}\endgroup%
}}\right.$}%
\begingroup \smaller\smaller\smaller\begin{tabular}{@{}c@{}}%
764820\\220500\\-1260
\end{tabular}\endgroup%
\kern3pt%
\begingroup \smaller\smaller\smaller\begin{tabular}{@{}c@{}}%
220500\\63570\\-363
\end{tabular}\endgroup%
\kern3pt%
\begingroup \smaller\smaller\smaller\begin{tabular}{@{}c@{}}%
-1260\\-363\\2
\end{tabular}\endgroup%
{$\left.\llap{\phantom{%
\begingroup \smaller\smaller\smaller\begin{tabular}{@{}c@{}}%
0\\0\\0
\end{tabular}\endgroup%
}}\!\right]$}%
\EasyButWeakLineBreak%
{$\left[\!\llap{\phantom{%
\begingroup \smaller\smaller\smaller\begin{tabular}{@{}c@{}}%
0\\0\\0
\end{tabular}\endgroup%
}}\right.$}%
\begingroup \smaller\smaller\smaller\begin{tabular}{@{}c@{}}%
-17\\60\\180
\end{tabular}\endgroup%
\HardButStrongLineBreak\kern3pt%
\begingroup \smaller\smaller\smaller\begin{tabular}{@{}c@{}}%
-2\\7\\14
\end{tabular}\endgroup%
\HardButStrongLineBreak\kern3pt%
\begingroup \smaller\smaller\smaller\begin{tabular}{@{}c@{}}%
7\\-25\\-120
\end{tabular}\endgroup%
\HardButStrongLineBreak\kern3pt%
\begingroup \smaller\smaller\smaller\begin{tabular}{@{}c@{}}%
230\\-819\\-3654
\end{tabular}\endgroup%
\HardButStrongLineBreak\kern3pt%
\begingroup \smaller\smaller\smaller\begin{tabular}{@{}c@{}}%
191\\-680\\-3020
\end{tabular}\endgroup%
\HardButStrongLineBreak\kern3pt%
\begingroup \smaller\smaller\smaller\begin{tabular}{@{}c@{}}%
59\\-210\\-927
\end{tabular}\endgroup%
\HardButStrongLineBreak\kern3pt%
\begingroup \smaller\smaller\smaller\begin{tabular}{@{}c@{}}%
9\\-32\\-138
\end{tabular}\endgroup%
\HardButStrongLineBreak\kern3pt%
\begingroup \smaller\smaller\smaller\begin{tabular}{@{}c@{}}%
0\\0\\-1
\end{tabular}\endgroup%
{$\left.\llap{\phantom{%
\begingroup \smaller\smaller\smaller\begin{tabular}{@{}c@{}}%
0\\0\\0
\end{tabular}\endgroup%
}}\!\right]$}%

\medskip%
%
\leavevmode\llap{}%
$W_{315}$%
\qquad\llap{12} lattices, $\chi=24$%
\hfill%
$22|22|22|22|\rtimes D_{4}$%
\nopagebreak\smallskip\hrule\nopagebreak\medskip%
%
%
\leavevmode%
${L_{315.1}}$%
{} : {$1\above{1pt}{1pt}{-2}{{\rm II}}4\above{1pt}{1pt}{-}{5}{\cdot}1\above{1pt}{1pt}{1}{}3\above{1pt}{1pt}{1}{}9\above{1pt}{1pt}{1}{}{\cdot}1\above{1pt}{1pt}{-2}{}5\above{1pt}{1pt}{-}{}{\cdot}1\above{1pt}{1pt}{2}{}7\above{1pt}{1pt}{-}{}$}\spacer%
\instructions{2}%
\EasyButWeakLineBreak%
{${84}\above{1pt}{1pt}{r}{2}{10}\above{1pt}{1pt}{b}{2}{12}\above{1pt}{1pt}{b}{2}{90}\above{1pt}{1pt}{l}{2}$}\relax$\,(\times2)$%
\nopagebreak\par%
\nopagebreak\par\leavevmode%
{$\left[\!\llap{\phantom{%
\begingroup \smaller\smaller\smaller\begin{tabular}{@{}c@{}}%
0\\0\\0
\end{tabular}\endgroup%
}}\right.$}%
\begingroup \smaller\smaller\smaller\begin{tabular}{@{}c@{}}%
-817740\\311220\\-76860
\end{tabular}\endgroup%
\kern3pt%
\begingroup \smaller\smaller\smaller\begin{tabular}{@{}c@{}}%
311220\\-117618\\28407
\end{tabular}\endgroup%
\kern3pt%
\begingroup \smaller\smaller\smaller\begin{tabular}{@{}c@{}}%
-76860\\28407\\-6362
\end{tabular}\endgroup%
{$\left.\llap{\phantom{%
\begingroup \smaller\smaller\smaller\begin{tabular}{@{}c@{}}%
0\\0\\0
\end{tabular}\endgroup%
}}\!\right]$}%
\hfil\penalty500%
{$\left[\!\llap{\phantom{%
\begingroup \smaller\smaller\smaller\begin{tabular}{@{}c@{}}%
0\\0\\0
\end{tabular}\endgroup%
}}\right.$}%
\begingroup \smaller\smaller\smaller\begin{tabular}{@{}c@{}}%
12599\\43680\\42840
\end{tabular}\endgroup%
\kern3pt%
\begingroup \smaller\smaller\smaller\begin{tabular}{@{}c@{}}%
-4605\\-15965\\-15657
\end{tabular}\endgroup%
\kern3pt%
\begingroup \smaller\smaller\smaller\begin{tabular}{@{}c@{}}%
990\\3432\\3365
\end{tabular}\endgroup%
{$\left.\llap{\phantom{%
\begingroup \smaller\smaller\smaller\begin{tabular}{@{}c@{}}%
0\\0\\0
\end{tabular}\endgroup%
}}\!\right]$}%
\EasyButWeakLineBreak%
{$\left[\!\llap{\phantom{%
\begingroup \smaller\smaller\smaller\begin{tabular}{@{}c@{}}%
0\\0\\0
\end{tabular}\endgroup%
}}\right.$}%
\begingroup \smaller\smaller\smaller\begin{tabular}{@{}c@{}}%
-2819\\-9772\\-9576
\end{tabular}\endgroup%
\HardButStrongLineBreak\kern3pt%
\begingroup \smaller\smaller\smaller\begin{tabular}{@{}c@{}}%
-789\\-2735\\-2680
\end{tabular}\endgroup%
\HardButStrongLineBreak\kern3pt%
\begingroup \smaller\smaller\smaller\begin{tabular}{@{}c@{}}%
-371\\-1286\\-1260
\end{tabular}\endgroup%
\HardButStrongLineBreak\kern3pt%
\begingroup \smaller\smaller\smaller\begin{tabular}{@{}c@{}}%
662\\2295\\2250
\end{tabular}\endgroup%
{$\left.\llap{\phantom{%
\begingroup \smaller\smaller\smaller\begin{tabular}{@{}c@{}}%
0\\0\\0
\end{tabular}\endgroup%
}}\!\right]$}%

\medskip%
%
\leavevmode\llap{}%
$W_{316}$%
\qquad\llap{44} lattices, $\chi=96$%
\hfill%
$22222|22222|22222|22222|\rtimes D_{4}$%
\nopagebreak\smallskip\hrule\nopagebreak\medskip%
%
%
\leavevmode%
${L_{316.1}}$%
{} : {$1\above{1pt}{1pt}{2}{{\rm II}}4\above{1pt}{1pt}{1}{1}{\cdot}1\above{1pt}{1pt}{1}{}3\above{1pt}{1pt}{1}{}9\above{1pt}{1pt}{1}{}{\cdot}1\above{1pt}{1pt}{-2}{}5\above{1pt}{1pt}{-}{}{\cdot}1\above{1pt}{1pt}{-2}{}7\above{1pt}{1pt}{1}{}$}\spacer%
\instructions{2}%
\EasyButWeakLineBreak%
{${210}\above{1pt}{1pt}{l}{2}{4}\above{1pt}{1pt}{r}{2}{90}\above{1pt}{1pt}{b}{2}{28}\above{1pt}{1pt}{*}{2}{36}\above{1pt}{1pt}{*}{2}{12}\above{1pt}{1pt}{*}{2}{4}\above{1pt}{1pt}{*}{2}{252}\above{1pt}{1pt}{b}{2}{10}\above{1pt}{1pt}{l}{2}{36}\above{1pt}{1pt}{r}{2}$}\relax$\,(\times2)$%
\nopagebreak\par%
\nopagebreak\par\leavevmode%
{$\left[\!\llap{\phantom{%
\begingroup \smaller\smaller\smaller
\endgroup%
}}\!\right]$}%
%
%
\hbox{}\par\smallskip%
%
%
\leavevmode%
${L_{316.2}}$%
{} : {$1\above{1pt}{1pt}{-2}{2}8\above{1pt}{1pt}{-}{3}{\cdot}1\above{1pt}{1pt}{-}{}3\above{1pt}{1pt}{-}{}9\above{1pt}{1pt}{-}{}{\cdot}1\above{1pt}{1pt}{-2}{}5\above{1pt}{1pt}{1}{}{\cdot}1\above{1pt}{1pt}{-2}{}7\above{1pt}{1pt}{1}{}$}\spacer%
\instructions{2}%
\EasyButWeakLineBreak%
{${420}\above{1pt}{1pt}{s}{2}{8}\above{1pt}{1pt}{s}{2}{180}\above{1pt}{1pt}{*}{2}{56}\above{1pt}{1pt}{b}{2}{18}\above{1pt}{1pt}{l}{2}{24}\above{1pt}{1pt}{r}{2}{2}\above{1pt}{1pt}{b}{2}{504}\above{1pt}{1pt}{*}{2}{20}\above{1pt}{1pt}{s}{2}{72}\above{1pt}{1pt}{s}{2}$}\relax$\,(\times2)$%
\nopagebreak\par%
\nopagebreak\par\leavevmode%
{$\left[\!\llap{\phantom{%
\begingroup \smaller\smaller\smaller
\endgroup%
}}\!\right]$}%
%
%
\hbox{}\par\smallskip%
%
%
\leavevmode%
${L_{316.3}}$%
{} : {$1\above{1pt}{1pt}{2}{2}8\above{1pt}{1pt}{1}{7}{\cdot}1\above{1pt}{1pt}{-}{}3\above{1pt}{1pt}{-}{}9\above{1pt}{1pt}{-}{}{\cdot}1\above{1pt}{1pt}{-2}{}5\above{1pt}{1pt}{1}{}{\cdot}1\above{1pt}{1pt}{-2}{}7\above{1pt}{1pt}{1}{}$}\spacer%
\instructions{m}%
\EasyButWeakLineBreak%
{${105}\above{1pt}{1pt}{r}{2}{8}\above{1pt}{1pt}{l}{2}{45}\above{1pt}{1pt}{}{2}{56}\above{1pt}{1pt}{r}{2}{18}\above{1pt}{1pt}{b}{2}{24}\above{1pt}{1pt}{b}{2}{2}\above{1pt}{1pt}{l}{2}{504}\above{1pt}{1pt}{}{2}{5}\above{1pt}{1pt}{r}{2}{72}\above{1pt}{1pt}{l}{2}$}\relax$\,(\times2)$%
\nopagebreak\par%
\nopagebreak\par\leavevmode%
{$\left[\!\llap{\phantom{%
\begingroup \smaller\smaller\smaller
\endgroup%
}}\!\right]$}%

\medskip%
%
\leavevmode\llap{}%
$W_{317}$%
\qquad\llap{32} lattices, $\chi=36$%
\hfill%
$2222|22222|2\rtimes D_{2}$%
\nopagebreak\smallskip\hrule\nopagebreak\medskip%
%
%
\leavevmode%
${L_{317.1}}$%
{} : {$[1\above{1pt}{1pt}{1}{}2\above{1pt}{1pt}{1}{}]\above{1pt}{1pt}{}{2}32\above{1pt}{1pt}{1}{1}{\cdot}1\above{1pt}{1pt}{2}{}3\above{1pt}{1pt}{1}{}{\cdot}1\above{1pt}{1pt}{2}{}5\above{1pt}{1pt}{-}{}$}\EasyButWeakLineBreak%
{${32}\above{1pt}{1pt}{}{2}{3}\above{1pt}{1pt}{r}{2}{160}\above{1pt}{1pt}{*}{2}{4}\above{1pt}{1pt}{*}{2}{40}\above{1pt}{1pt}{l}{2}{1}\above{1pt}{1pt}{r}{2}{160}\above{1pt}{1pt}{*}{2}{12}\above{1pt}{1pt}{s}{2}{32}\above{1pt}{1pt}{l}{2}{2}\above{1pt}{1pt}{}{2}$}%
\nopagebreak\par%
\nopagebreak\par\leavevmode%
{$\left[\!\llap{\phantom{%
\begingroup \smaller\smaller\smaller
\endgroup%
}}\!\right]$}%
%
%
\hbox{}\par\smallskip%
%
%
\leavevmode%
${L_{317.2}}$%
{} : {$[1\above{1pt}{1pt}{-}{}2\above{1pt}{1pt}{1}{}]\above{1pt}{1pt}{}{2}32\above{1pt}{1pt}{-}{5}{\cdot}1\above{1pt}{1pt}{2}{}3\above{1pt}{1pt}{1}{}{\cdot}1\above{1pt}{1pt}{2}{}5\above{1pt}{1pt}{-}{}$}\EasyButWeakLineBreak%
{${32}\above{1pt}{1pt}{*}{2}{12}\above{1pt}{1pt}{s}{2}{160}\above{1pt}{1pt}{s}{2}{4}\above{1pt}{1pt}{l}{2}{10}\above{1pt}{1pt}{}{2}{1}\above{1pt}{1pt}{}{2}{160}\above{1pt}{1pt}{}{2}{3}\above{1pt}{1pt}{r}{2}{32}\above{1pt}{1pt}{*}{2}{8}\above{1pt}{1pt}{s}{2}$}%
\nopagebreak\par%
\nopagebreak\par\leavevmode%
{$\left[\!\llap{\phantom{%
\begingroup \smaller\smaller\smaller
\endgroup%
}}\!\right]$}%
%
%
\hbox{}\par\smallskip%
%
%
\leavevmode%
${L_{317.3}}$%
{} : {$1\above{1pt}{1pt}{1}{1}4\above{1pt}{1pt}{1}{1}32\above{1pt}{1pt}{1}{1}{\cdot}1\above{1pt}{1pt}{2}{}3\above{1pt}{1pt}{-}{}{\cdot}1\above{1pt}{1pt}{2}{}5\above{1pt}{1pt}{1}{}$}\EasyButWeakLineBreak%
{${1}\above{1pt}{1pt}{r}{2}{96}\above{1pt}{1pt}{*}{2}{20}\above{1pt}{1pt}{s}{2}{32}\above{1pt}{1pt}{l}{2}{20}\above{1pt}{1pt}{}{2}{32}\above{1pt}{1pt}{}{2}{5}\above{1pt}{1pt}{r}{2}{96}\above{1pt}{1pt}{*}{2}{4}\above{1pt}{1pt}{l}{2}{4}\above{1pt}{1pt}{}{2}$}%
\nopagebreak\par%
\nopagebreak\par\leavevmode%
{$\left[\!\llap{\phantom{%
\begingroup \smaller\smaller\smaller
\endgroup%
}}\!\right]$}%
%
%
\hbox{}\par\smallskip%
%
%
\leavevmode%
${L_{317.4}}$%
{} : {$1\above{1pt}{1pt}{-}{5}4\above{1pt}{1pt}{1}{7}32\above{1pt}{1pt}{-}{3}{\cdot}1\above{1pt}{1pt}{2}{}3\above{1pt}{1pt}{-}{}{\cdot}1\above{1pt}{1pt}{2}{}5\above{1pt}{1pt}{1}{}$}\EasyButWeakLineBreak%
{${4}\above{1pt}{1pt}{s}{2}{96}\above{1pt}{1pt}{s}{2}{20}\above{1pt}{1pt}{*}{2}{32}\above{1pt}{1pt}{*}{2}{80}\above{1pt}{1pt}{s}{2}{32}\above{1pt}{1pt}{l}{2}{5}\above{1pt}{1pt}{}{2}{96}\above{1pt}{1pt}{}{2}{1}\above{1pt}{1pt}{r}{2}{16}\above{1pt}{1pt}{*}{2}$}%
\nopagebreak\par%
\nopagebreak\par\leavevmode%
{$\left[\!\llap{\phantom{%
\begingroup \smaller\smaller\smaller
\endgroup%
}}\!\right]$}%

\medskip%
%
\leavevmode\llap{}%
$W_{318}$%
\qquad\llap{48} lattices, $\chi=48$%
\hfill%
$22|222|222|222|2\rtimes D_{4}$%
\nopagebreak\smallskip\hrule\nopagebreak\medskip%
%
%
\leavevmode%
${L_{318.1}}$%
{} : {$[1\above{1pt}{1pt}{-}{}2\above{1pt}{1pt}{1}{}]\above{1pt}{1pt}{}{6}32\above{1pt}{1pt}{-}{5}{\cdot}1\above{1pt}{1pt}{-}{}3\above{1pt}{1pt}{1}{}9\above{1pt}{1pt}{-}{}{\cdot}1\above{1pt}{1pt}{-2}{}5\above{1pt}{1pt}{1}{}$}\spacer%
\instructions{3}%
\EasyButWeakLineBreak%
{${32}\above{1pt}{1pt}{l}{2}{45}\above{1pt}{1pt}{}{2}{2}\above{1pt}{1pt}{r}{2}{180}\above{1pt}{1pt}{*}{2}{32}\above{1pt}{1pt}{s}{2}{120}\above{1pt}{1pt}{*}{2}{288}\above{1pt}{1pt}{l}{2}{5}\above{1pt}{1pt}{}{2}{18}\above{1pt}{1pt}{r}{2}{20}\above{1pt}{1pt}{*}{2}{288}\above{1pt}{1pt}{s}{2}{120}\above{1pt}{1pt}{*}{2}$}%
\nopagebreak\par%
\nopagebreak\par\leavevmode%
{$\left[\!\llap{\phantom{%
\begingroup \smaller\smaller\smaller
\endgroup%
}}\!\right]$}%
%
%
\hbox{}\par\smallskip%
%
%
\leavevmode%
${L_{318.2}}$%
{} : {$[1\above{1pt}{1pt}{1}{}2\above{1pt}{1pt}{1}{}]\above{1pt}{1pt}{}{6}32\above{1pt}{1pt}{1}{1}{\cdot}1\above{1pt}{1pt}{-}{}3\above{1pt}{1pt}{1}{}9\above{1pt}{1pt}{-}{}{\cdot}1\above{1pt}{1pt}{-2}{}5\above{1pt}{1pt}{1}{}$}\spacer%
\instructions{3}%
\EasyButWeakLineBreak%
{${32}\above{1pt}{1pt}{s}{2}{180}\above{1pt}{1pt}{*}{2}{8}\above{1pt}{1pt}{l}{2}{45}\above{1pt}{1pt}{}{2}{32}\above{1pt}{1pt}{}{2}{30}\above{1pt}{1pt}{r}{2}{288}\above{1pt}{1pt}{s}{2}{20}\above{1pt}{1pt}{*}{2}{72}\above{1pt}{1pt}{l}{2}{5}\above{1pt}{1pt}{}{2}{288}\above{1pt}{1pt}{}{2}{30}\above{1pt}{1pt}{r}{2}$}%
\nopagebreak\par%
\nopagebreak\par\leavevmode%
{$\left[\!\llap{\phantom{%
\begingroup \smaller\smaller\smaller
\endgroup%
}}\!\right]$}%
%
%
\hbox{}\par\smallskip%
%
%
\leavevmode%
${L_{318.3}}$%
{} : {$1\above{1pt}{1pt}{1}{1}4\above{1pt}{1pt}{1}{7}32\above{1pt}{1pt}{1}{7}{\cdot}1\above{1pt}{1pt}{1}{}3\above{1pt}{1pt}{-}{}9\above{1pt}{1pt}{1}{}{\cdot}1\above{1pt}{1pt}{-2}{}5\above{1pt}{1pt}{-}{}$}\spacer%
\instructions{3}%
\EasyButWeakLineBreak%
{${4}\above{1pt}{1pt}{*}{2}{1440}\above{1pt}{1pt}{s}{2}{16}\above{1pt}{1pt}{*}{2}{1440}\above{1pt}{1pt}{l}{2}{1}\above{1pt}{1pt}{}{2}{60}\above{1pt}{1pt}{r}{2}{36}\above{1pt}{1pt}{*}{2}{160}\above{1pt}{1pt}{s}{2}{144}\above{1pt}{1pt}{*}{2}{160}\above{1pt}{1pt}{l}{2}{9}\above{1pt}{1pt}{}{2}{60}\above{1pt}{1pt}{r}{2}$}%
\nopagebreak\par%
\nopagebreak\par\leavevmode%
{$\left[\!\llap{\phantom{%
\begingroup \smaller\smaller\smaller
\endgroup%
}}\!\right]$}%
%
%
\hbox{}\par\smallskip%
%
%
\leavevmode%
${L_{318.4}}$%
{} : {$1\above{1pt}{1pt}{-}{5}4\above{1pt}{1pt}{1}{1}32\above{1pt}{1pt}{-}{5}{\cdot}1\above{1pt}{1pt}{1}{}3\above{1pt}{1pt}{-}{}9\above{1pt}{1pt}{1}{}{\cdot}1\above{1pt}{1pt}{-2}{}5\above{1pt}{1pt}{-}{}$}\spacer%
\instructions{3}%
\EasyButWeakLineBreak%
{${4}\above{1pt}{1pt}{s}{2}{1440}\above{1pt}{1pt}{l}{2}{4}\above{1pt}{1pt}{}{2}{1440}\above{1pt}{1pt}{}{2}{1}\above{1pt}{1pt}{r}{2}{240}\above{1pt}{1pt}{*}{2}{36}\above{1pt}{1pt}{s}{2}{160}\above{1pt}{1pt}{l}{2}{36}\above{1pt}{1pt}{}{2}{160}\above{1pt}{1pt}{}{2}{9}\above{1pt}{1pt}{r}{2}{240}\above{1pt}{1pt}{*}{2}$}%
\nopagebreak\par%
\nopagebreak\par\leavevmode%
{$\left[\!\llap{\phantom{%
\begingroup \smaller\smaller\smaller
\endgroup%
}}\!\right]$}%

\medskip%
%
\leavevmode\llap{}%
$W_{319}$%
\qquad\llap{40} lattices, $\chi=48$%
\hfill%
$22|222|222|222|2\rtimes D_{4}$%
\nopagebreak\smallskip\hrule\nopagebreak\medskip%
%
%
\leavevmode%
${L_{319.1}}$%
{} : {$[1\above{1pt}{1pt}{1}{}2\above{1pt}{1pt}{-}{}]\above{1pt}{1pt}{}{4}32\above{1pt}{1pt}{-}{3}{\cdot}1\above{1pt}{1pt}{2}{}3\above{1pt}{1pt}{1}{}{\cdot}1\above{1pt}{1pt}{-2}{}5\above{1pt}{1pt}{1}{}$}\EasyButWeakLineBreak%
{${480}\above{1pt}{1pt}{l}{2}{1}\above{1pt}{1pt}{r}{2}{120}\above{1pt}{1pt}{*}{2}{4}\above{1pt}{1pt}{*}{2}{480}\above{1pt}{1pt}{s}{2}{8}\above{1pt}{1pt}{*}{2}$}\relax$\,(\times2)$%
\nopagebreak\par%
\nopagebreak\par\leavevmode%
{$\left[\!\llap{\phantom{%
\begingroup \smaller\smaller\smaller
\endgroup%
}}\!\right]$}%
%
%
\hbox{}\par\smallskip%
%
%
\leavevmode%
${L_{319.2}}$%
{} : {$1\above{1pt}{1pt}{1}{7}4\above{1pt}{1pt}{1}{1}32\above{1pt}{1pt}{1}{7}{\cdot}1\above{1pt}{1pt}{2}{}3\above{1pt}{1pt}{-}{}{\cdot}1\above{1pt}{1pt}{-2}{}5\above{1pt}{1pt}{-}{}$}\EasyButWeakLineBreak%
{${15}\above{1pt}{1pt}{r}{2}{32}\above{1pt}{1pt}{s}{2}{240}\above{1pt}{1pt}{*}{2}{32}\above{1pt}{1pt}{*}{2}{60}\above{1pt}{1pt}{l}{2}{4}\above{1pt}{1pt}{}{2}$}\relax$\,(\times2)$%
\nopagebreak\par%
\nopagebreak\par\leavevmode%
{$\left[\!\llap{\phantom{%
\begingroup \smaller\smaller\smaller
\endgroup%
}}\!\right]$}%
%
%
\hbox{}\par\smallskip%
%
%
\leavevmode%
${L_{319.3}}$%
{} : {$[1\above{1pt}{1pt}{1}{}2\above{1pt}{1pt}{1}{}]\above{1pt}{1pt}{}{0}64\above{1pt}{1pt}{1}{7}{\cdot}1\above{1pt}{1pt}{2}{}3\above{1pt}{1pt}{-}{}{\cdot}1\above{1pt}{1pt}{-2}{}5\above{1pt}{1pt}{-}{}$}\spacer%
\instructions{m}%
\EasyButWeakLineBreak%
{${960}\above{1pt}{1pt}{s}{2}{8}\above{1pt}{1pt}{*}{2}{60}\above{1pt}{1pt}{l}{2}{2}\above{1pt}{1pt}{}{2}{960}\above{1pt}{1pt}{}{2}{1}\above{1pt}{1pt}{r}{2}{960}\above{1pt}{1pt}{*}{2}{8}\above{1pt}{1pt}{l}{2}{15}\above{1pt}{1pt}{}{2}{2}\above{1pt}{1pt}{r}{2}{960}\above{1pt}{1pt}{s}{2}{4}\above{1pt}{1pt}{*}{2}$}%
\nopagebreak\par%
\nopagebreak\par\leavevmode%
{$\left[\!\llap{\phantom{%
\begingroup \smaller\smaller\smaller
\endgroup%
}}\!\right]$}%
%
%
\hbox{}\par\smallskip%
%
%
\leavevmode%
${L_{319.4}}$%
{} : {$1\above{1pt}{1pt}{1}{7}4\above{1pt}{1pt}{1}{7}32\above{1pt}{1pt}{1}{1}{\cdot}1\above{1pt}{1pt}{2}{}3\above{1pt}{1pt}{-}{}{\cdot}1\above{1pt}{1pt}{-2}{}5\above{1pt}{1pt}{-}{}$}\EasyButWeakLineBreak%
{${15}\above{1pt}{1pt}{}{2}{32}\above{1pt}{1pt}{}{2}{60}\above{1pt}{1pt}{r}{2}{32}\above{1pt}{1pt}{s}{2}{60}\above{1pt}{1pt}{*}{2}{16}\above{1pt}{1pt}{l}{2}$}\relax$\,(\times2)$%
\nopagebreak\par%
\nopagebreak\par\leavevmode%
{$\left[\!\llap{\phantom{%
\begingroup \smaller\smaller\smaller
\endgroup%
}}\!\right]$}%

\medskip%
%
\leavevmode\llap{}%
$W_{320}$%
\qquad\llap{8} lattices, $\chi=48$%
\hfill%
$2|2|2|2|2|2|2|2|2|2|2|2|\rtimes D_{12}$%
\nopagebreak\smallskip\hrule\nopagebreak\medskip%
%
%
\leavevmode%
${L_{320.1}}$%
{} : {$1\above{1pt}{1pt}{1}{1}4\above{1pt}{1pt}{1}{1}16\above{1pt}{1pt}{1}{1}{\cdot}1\above{1pt}{1pt}{1}{}3\above{1pt}{1pt}{-}{}9\above{1pt}{1pt}{-}{}{\cdot}1\above{1pt}{1pt}{-2}{}5\above{1pt}{1pt}{1}{}$}\spacer%
\instructions{3}%
\EasyButWeakLineBreak%
{${1}\above{1pt}{1pt}{}{2}{720}\above{1pt}{1pt}{}{2}{4}\above{1pt}{1pt}{}{2}{45}\above{1pt}{1pt}{}{2}{16}\above{1pt}{1pt}{}{2}{180}\above{1pt}{1pt}{r}{2}{4}\above{1pt}{1pt}{s}{2}{720}\above{1pt}{1pt}{l}{2}{4}\above{1pt}{1pt}{r}{2}{180}\above{1pt}{1pt}{s}{2}{16}\above{1pt}{1pt}{l}{2}{180}\above{1pt}{1pt}{}{2}$}%
\nopagebreak\par%
\nopagebreak\par\leavevmode%
{$\left[\!\llap{\phantom{%
\begingroup \smaller\smaller\smaller
\endgroup%
}}\!\right]$}%

\medskip%
%
\leavevmode\llap{}%
$W_{321}$%
\qquad\llap{12} lattices, $\chi=48$%
\hfill%
$22|222|222|222|2\rtimes D_{4}$%
\nopagebreak\smallskip\hrule\nopagebreak\medskip%
%
%
\leavevmode%
${L_{321.1}}$%
{} : {$1\above{1pt}{1pt}{-2}{{\rm II}}4\above{1pt}{1pt}{1}{1}{\cdot}1\above{1pt}{1pt}{-}{}3\above{1pt}{1pt}{1}{}9\above{1pt}{1pt}{1}{}{\cdot}1\above{1pt}{1pt}{1}{}7\above{1pt}{1pt}{1}{}49\above{1pt}{1pt}{1}{}$}\spacer%
\instructions{23,3,2}%
\EasyButWeakLineBreak%
{${1764}\above{1pt}{1pt}{r}{2}{2}\above{1pt}{1pt}{b}{2}{252}\above{1pt}{1pt}{b}{2}{98}\above{1pt}{1pt}{l}{2}{36}\above{1pt}{1pt}{r}{2}{14}\above{1pt}{1pt}{l}{2}$}\relax$\,(\times2)$%
\nopagebreak\par%
\nopagebreak\par\leavevmode%
{$\left[\!\llap{\phantom{%
\begingroup \smaller\smaller\smaller
\endgroup%
}}\!\right]$}%

\medskip%
%
\leavevmode\llap{}%
$W_{322}$%
\qquad\llap{12} lattices, $\chi=24$%
\hfill%
$22222222\rtimes C_{2}$%
\nopagebreak\smallskip\hrule\nopagebreak\medskip%
%
%
\leavevmode%
${L_{322.1}}$%
{} : {$1\above{1pt}{1pt}{-2}{{\rm II}}4\above{1pt}{1pt}{1}{1}{\cdot}1\above{1pt}{1pt}{2}{}3\above{1pt}{1pt}{-}{}{\cdot}1\above{1pt}{1pt}{1}{}7\above{1pt}{1pt}{-}{}49\above{1pt}{1pt}{-}{}$}\spacer%
\instructions{2}%
\EasyButWeakLineBreak%
{${2}\above{1pt}{1pt}{s}{2}{294}\above{1pt}{1pt}{l}{2}{4}\above{1pt}{1pt}{r}{2}{42}\above{1pt}{1pt}{b}{2}$}\relax$\,(\times2)$%
\nopagebreak\par%
\nopagebreak\par\leavevmode%
{$\left[\!\llap{\phantom{%
\begingroup \smaller\smaller\smaller\begin{tabular}{@{}c@{}}%
0\\0\\0
\end{tabular}\endgroup%
}}\right.$}%
\begingroup \smaller\smaller\smaller\begin{tabular}{@{}c@{}}%
105252\\-5292\\-2352
\end{tabular}\endgroup%
\kern3pt%
\begingroup \smaller\smaller\smaller\begin{tabular}{@{}c@{}}%
-5292\\266\\119
\end{tabular}\endgroup%
\kern3pt%
\begingroup \smaller\smaller\smaller\begin{tabular}{@{}c@{}}%
-2352\\119\\46
\end{tabular}\endgroup%
{$\left.\llap{\phantom{%
\begingroup \smaller\smaller\smaller\begin{tabular}{@{}c@{}}%
0\\0\\0
\end{tabular}\endgroup%
}}\!\right]$}%
\hfil\penalty500%
{$\left[\!\llap{\phantom{%
\begingroup \smaller\smaller\smaller\begin{tabular}{@{}c@{}}%
0\\0\\0
\end{tabular}\endgroup%
}}\right.$}%
\begingroup \smaller\smaller\smaller\begin{tabular}{@{}c@{}}%
503\\9324\\1764
\end{tabular}\endgroup%
\kern3pt%
\begingroup \smaller\smaller\smaller\begin{tabular}{@{}c@{}}%
-26\\-482\\-91
\end{tabular}\endgroup%
\kern3pt%
\begingroup \smaller\smaller\smaller\begin{tabular}{@{}c@{}}%
-6\\-111\\-22
\end{tabular}\endgroup%
{$\left.\llap{\phantom{%
\begingroup \smaller\smaller\smaller\begin{tabular}{@{}c@{}}%
0\\0\\0
\end{tabular}\endgroup%
}}\!\right]$}%
\EasyButWeakLineBreak%
{$\left[\!\llap{\phantom{%
\begingroup \smaller\smaller\smaller\begin{tabular}{@{}c@{}}%
0\\0\\0
\end{tabular}\endgroup%
}}\right.$}%
\begingroup \smaller\smaller\smaller\begin{tabular}{@{}c@{}}%
1\\19\\2
\end{tabular}\endgroup%
\HardButStrongLineBreak\kern3pt%
\begingroup \smaller\smaller\smaller\begin{tabular}{@{}c@{}}%
1\\21\\0
\end{tabular}\endgroup%
\HardButStrongLineBreak\kern3pt%
\begingroup \smaller\smaller\smaller\begin{tabular}{@{}c@{}}%
-3\\-56\\-8
\end{tabular}\endgroup%
\HardButStrongLineBreak\kern3pt%
\begingroup \smaller\smaller\smaller\begin{tabular}{@{}c@{}}%
-8\\-150\\-21
\end{tabular}\endgroup%
{$\left.\llap{\phantom{%
\begingroup \smaller\smaller\smaller\begin{tabular}{@{}c@{}}%
0\\0\\0
\end{tabular}\endgroup%
}}\!\right]$}%

\medskip%
%
\leavevmode\llap{}%
$W_{323}$%
\qquad\llap{30} lattices, $\chi=24$%
\hfill%
$22|22|22|22|\rtimes D_{4}$%
\nopagebreak\smallskip\hrule\nopagebreak\medskip%
%
%
\leavevmode%
${L_{323.1}}$%
{} : {$1\above{1pt}{1pt}{2}{0}8\above{1pt}{1pt}{1}{7}{\cdot}1\above{1pt}{1pt}{1}{}3\above{1pt}{1pt}{-}{}9\above{1pt}{1pt}{1}{}{\cdot}1\above{1pt}{1pt}{-2}{}5\above{1pt}{1pt}{-}{}$}\EasyButWeakLineBreak%
{${360}\above{1pt}{1pt}{s}{2}{4}\above{1pt}{1pt}{*}{2}{60}\above{1pt}{1pt}{*}{2}{36}\above{1pt}{1pt}{s}{2}{40}\above{1pt}{1pt}{l}{2}{9}\above{1pt}{1pt}{}{2}{15}\above{1pt}{1pt}{}{2}{1}\above{1pt}{1pt}{r}{2}$}%
\nopagebreak\par%
\nopagebreak\par\leavevmode%
{$\left[\!\llap{\phantom{%
\begingroup \smaller\smaller\smaller\begin{tabular}{@{}c@{}}%
0\\0\\0
\end{tabular}\endgroup%
}}\right.$}%
\begingroup \smaller\smaller\smaller\begin{tabular}{@{}c@{}}%
2711160\\720\\-15480
\end{tabular}\endgroup%
\kern3pt%
\begingroup \smaller\smaller\smaller\begin{tabular}{@{}c@{}}%
720\\-3\\-3
\end{tabular}\endgroup%
\kern3pt%
\begingroup \smaller\smaller\smaller\begin{tabular}{@{}c@{}}%
-15480\\-3\\88
\end{tabular}\endgroup%
{$\left.\llap{\phantom{%
\begingroup \smaller\smaller\smaller\begin{tabular}{@{}c@{}}%
0\\0\\0
\end{tabular}\endgroup%
}}\!\right]$}%
\EasyButWeakLineBreak%
{$\left[\!\llap{\phantom{%
\begingroup \smaller\smaller\smaller\begin{tabular}{@{}c@{}}%
0\\0\\0
\end{tabular}\endgroup%
}}\right.$}%
\begingroup \smaller\smaller\smaller\begin{tabular}{@{}c@{}}%
1\\60\\180
\end{tabular}\endgroup%
\HardButStrongLineBreak\kern3pt%
\begingroup \smaller\smaller\smaller\begin{tabular}{@{}c@{}}%
-1\\-62\\-178
\end{tabular}\endgroup%
\HardButStrongLineBreak\kern3pt%
\begingroup \smaller\smaller\smaller\begin{tabular}{@{}c@{}}%
-1\\-70\\-180
\end{tabular}\endgroup%
\HardButStrongLineBreak\kern3pt%
\begingroup \smaller\smaller\smaller\begin{tabular}{@{}c@{}}%
7\\414\\1242
\end{tabular}\endgroup%
\HardButStrongLineBreak\kern3pt%
\begingroup \smaller\smaller\smaller\begin{tabular}{@{}c@{}}%
17\\1020\\3020
\end{tabular}\endgroup%
\HardButStrongLineBreak\kern3pt%
\begingroup \smaller\smaller\smaller\begin{tabular}{@{}c@{}}%
8\\483\\1422
\end{tabular}\endgroup%
\HardButStrongLineBreak\kern3pt%
\begingroup \smaller\smaller\smaller\begin{tabular}{@{}c@{}}%
7\\425\\1245
\end{tabular}\endgroup%
\HardButStrongLineBreak\kern3pt%
\begingroup \smaller\smaller\smaller\begin{tabular}{@{}c@{}}%
1\\61\\178
\end{tabular}\endgroup%
{$\left.\llap{\phantom{%
\begingroup \smaller\smaller\smaller\begin{tabular}{@{}c@{}}%
0\\0\\0
\end{tabular}\endgroup%
}}\!\right]$}%
%
%
\hbox{}\par\smallskip%
%
%
\leavevmode%
${L_{323.2}}$%
{} : {$[1\above{1pt}{1pt}{-}{}2\above{1pt}{1pt}{1}{}]\above{1pt}{1pt}{}{2}16\above{1pt}{1pt}{-}{5}{\cdot}1\above{1pt}{1pt}{1}{}3\above{1pt}{1pt}{-}{}9\above{1pt}{1pt}{1}{}{\cdot}1\above{1pt}{1pt}{-2}{}5\above{1pt}{1pt}{-}{}$}\spacer%
\instructions{2}%
\EasyButWeakLineBreak%
{${10}\above{1pt}{1pt}{}{2}{9}\above{1pt}{1pt}{r}{2}{240}\above{1pt}{1pt}{l}{2}{1}\above{1pt}{1pt}{}{2}{90}\above{1pt}{1pt}{r}{2}{16}\above{1pt}{1pt}{s}{2}{60}\above{1pt}{1pt}{s}{2}{144}\above{1pt}{1pt}{l}{2}$}%
\nopagebreak\par%
\nopagebreak\par\leavevmode%
{$\left[\!\llap{\phantom{%
\begingroup \smaller\smaller\smaller\begin{tabular}{@{}c@{}}%
0\\0\\0
\end{tabular}\endgroup%
}}\right.$}%
\begingroup \smaller\smaller\smaller\begin{tabular}{@{}c@{}}%
23253840\\415440\\5040
\end{tabular}\endgroup%
\kern3pt%
\begingroup \smaller\smaller\smaller\begin{tabular}{@{}c@{}}%
415440\\7422\\90
\end{tabular}\endgroup%
\kern3pt%
\begingroup \smaller\smaller\smaller\begin{tabular}{@{}c@{}}%
5040\\90\\1
\end{tabular}\endgroup%
{$\left.\llap{\phantom{%
\begingroup \smaller\smaller\smaller\begin{tabular}{@{}c@{}}%
0\\0\\0
\end{tabular}\endgroup%
}}\!\right]$}%
\EasyButWeakLineBreak%
{$\left[\!\llap{\phantom{%
\begingroup \smaller\smaller\smaller\begin{tabular}{@{}c@{}}%
0\\0\\0
\end{tabular}\endgroup%
}}\right.$}%
\begingroup \smaller\smaller\smaller\begin{tabular}{@{}c@{}}%
-7\\395\\-250
\end{tabular}\endgroup%
\HardButStrongLineBreak\kern3pt%
\begingroup \smaller\smaller\smaller\begin{tabular}{@{}c@{}}%
-5\\282\\-171
\end{tabular}\endgroup%
\HardButStrongLineBreak\kern3pt%
\begingroup \smaller\smaller\smaller\begin{tabular}{@{}c@{}}%
-11\\620\\-360
\end{tabular}\endgroup%
\HardButStrongLineBreak\kern3pt%
\begingroup \smaller\smaller\smaller\begin{tabular}{@{}c@{}}%
0\\0\\-1
\end{tabular}\endgroup%
\HardButStrongLineBreak\kern3pt%
\begingroup \smaller\smaller\smaller\begin{tabular}{@{}c@{}}%
4\\-225\\90
\end{tabular}\endgroup%
\HardButStrongLineBreak\kern3pt%
\begingroup \smaller\smaller\smaller\begin{tabular}{@{}c@{}}%
1\\-56\\8
\end{tabular}\endgroup%
\HardButStrongLineBreak\kern3pt%
\begingroup \smaller\smaller\smaller\begin{tabular}{@{}c@{}}%
-3\\170\\-150
\end{tabular}\endgroup%
\HardButStrongLineBreak\kern3pt%
\begingroup \smaller\smaller\smaller\begin{tabular}{@{}c@{}}%
-17\\960\\-648
\end{tabular}\endgroup%
{$\left.\llap{\phantom{%
\begingroup \smaller\smaller\smaller\begin{tabular}{@{}c@{}}%
0\\0\\0
\end{tabular}\endgroup%
}}\!\right]$}%
%
%
\hbox{}\par\smallskip%
%
%
\leavevmode%
${L_{323.3}}$%
{} : {$[1\above{1pt}{1pt}{1}{}2\above{1pt}{1pt}{1}{}]\above{1pt}{1pt}{}{0}16\above{1pt}{1pt}{1}{7}{\cdot}1\above{1pt}{1pt}{1}{}3\above{1pt}{1pt}{-}{}9\above{1pt}{1pt}{1}{}{\cdot}1\above{1pt}{1pt}{-2}{}5\above{1pt}{1pt}{-}{}$}\spacer%
\instructions{m}%
\EasyButWeakLineBreak%
{${40}\above{1pt}{1pt}{l}{2}{9}\above{1pt}{1pt}{}{2}{240}\above{1pt}{1pt}{}{2}{1}\above{1pt}{1pt}{r}{2}{360}\above{1pt}{1pt}{*}{2}{16}\above{1pt}{1pt}{l}{2}{15}\above{1pt}{1pt}{r}{2}{144}\above{1pt}{1pt}{*}{2}$}%
\nopagebreak\par%
\nopagebreak\par\leavevmode%
{$\left[\!\llap{\phantom{%
\begingroup \smaller\smaller\smaller\begin{tabular}{@{}c@{}}%
0\\0\\0
\end{tabular}\endgroup%
}}\right.$}%
\begingroup \smaller\smaller\smaller\begin{tabular}{@{}c@{}}%
-3378960\\-51840\\22320
\end{tabular}\endgroup%
\kern3pt%
\begingroup \smaller\smaller\smaller\begin{tabular}{@{}c@{}}%
-51840\\-786\\336
\end{tabular}\endgroup%
\kern3pt%
\begingroup \smaller\smaller\smaller\begin{tabular}{@{}c@{}}%
22320\\336\\-143
\end{tabular}\endgroup%
{$\left.\llap{\phantom{%
\begingroup \smaller\smaller\smaller\begin{tabular}{@{}c@{}}%
0\\0\\0
\end{tabular}\endgroup%
}}\!\right]$}%
\EasyButWeakLineBreak%
{$\left[\!\llap{\phantom{%
\begingroup \smaller\smaller\smaller\begin{tabular}{@{}c@{}}%
0\\0\\0
\end{tabular}\endgroup%
}}\right.$}%
\begingroup \smaller\smaller\smaller\begin{tabular}{@{}c@{}}%
9\\-1550\\-2240
\end{tabular}\endgroup%
\HardButStrongLineBreak\kern3pt%
\begingroup \smaller\smaller\smaller\begin{tabular}{@{}c@{}}%
5\\-864\\-1251
\end{tabular}\endgroup%
\HardButStrongLineBreak\kern3pt%
\begingroup \smaller\smaller\smaller\begin{tabular}{@{}c@{}}%
21\\-3640\\-5280
\end{tabular}\endgroup%
\HardButStrongLineBreak\kern3pt%
\begingroup \smaller\smaller\smaller\begin{tabular}{@{}c@{}}%
1\\-174\\-253
\end{tabular}\endgroup%
\HardButStrongLineBreak\kern3pt%
\begingroup \smaller\smaller\smaller\begin{tabular}{@{}c@{}}%
7\\-1230\\-1800
\end{tabular}\endgroup%
\HardButStrongLineBreak\kern3pt%
\begingroup \smaller\smaller\smaller\begin{tabular}{@{}c@{}}%
-1\\172\\248
\end{tabular}\endgroup%
\HardButStrongLineBreak\kern3pt%
\begingroup \smaller\smaller\smaller\begin{tabular}{@{}c@{}}%
-1\\175\\255
\end{tabular}\endgroup%
\HardButStrongLineBreak\kern3pt%
\begingroup \smaller\smaller\smaller\begin{tabular}{@{}c@{}}%
5\\-852\\-1224
\end{tabular}\endgroup%
{$\left.\llap{\phantom{%
\begingroup \smaller\smaller\smaller\begin{tabular}{@{}c@{}}%
0\\0\\0
\end{tabular}\endgroup%
}}\!\right]$}%
%
%
\hbox{}\par\smallskip%
%
%
\leavevmode%
${L_{323.4}}$%
{} : {$[1\above{1pt}{1pt}{1}{}2\above{1pt}{1pt}{1}{}]\above{1pt}{1pt}{}{6}16\above{1pt}{1pt}{1}{1}{\cdot}1\above{1pt}{1pt}{1}{}3\above{1pt}{1pt}{-}{}9\above{1pt}{1pt}{1}{}{\cdot}1\above{1pt}{1pt}{-2}{}5\above{1pt}{1pt}{-}{}$}\spacer%
\instructions{m}%
\EasyButWeakLineBreak%
{${10}\above{1pt}{1pt}{r}{2}{36}\above{1pt}{1pt}{*}{2}{240}\above{1pt}{1pt}{*}{2}{4}\above{1pt}{1pt}{l}{2}{90}\above{1pt}{1pt}{}{2}{16}\above{1pt}{1pt}{}{2}{15}\above{1pt}{1pt}{}{2}{144}\above{1pt}{1pt}{}{2}$}%
\nopagebreak\par%
\nopagebreak\par\leavevmode%
{$\left[\!\llap{\phantom{%
\begingroup \smaller\smaller\smaller\begin{tabular}{@{}c@{}}%
0\\0\\0
\end{tabular}\endgroup%
}}\right.$}%
\begingroup \smaller\smaller\smaller\begin{tabular}{@{}c@{}}%
-3054960\\2160\\7920
\end{tabular}\endgroup%
\kern3pt%
\begingroup \smaller\smaller\smaller\begin{tabular}{@{}c@{}}%
2160\\42\\-18
\end{tabular}\endgroup%
\kern3pt%
\begingroup \smaller\smaller\smaller\begin{tabular}{@{}c@{}}%
7920\\-18\\-17
\end{tabular}\endgroup%
{$\left.\llap{\phantom{%
\begingroup \smaller\smaller\smaller\begin{tabular}{@{}c@{}}%
0\\0\\0
\end{tabular}\endgroup%
}}\!\right]$}%
\EasyButWeakLineBreak%
{$\left[\!\llap{\phantom{%
\begingroup \smaller\smaller\smaller\begin{tabular}{@{}c@{}}%
0\\0\\0
\end{tabular}\endgroup%
}}\right.$}%
\begingroup \smaller\smaller\smaller\begin{tabular}{@{}c@{}}%
8\\815\\2860
\end{tabular}\endgroup%
\HardButStrongLineBreak\kern3pt%
\begingroup \smaller\smaller\smaller\begin{tabular}{@{}c@{}}%
7\\714\\2502
\end{tabular}\endgroup%
\HardButStrongLineBreak\kern3pt%
\begingroup \smaller\smaller\smaller\begin{tabular}{@{}c@{}}%
-1\\-100\\-360
\end{tabular}\endgroup%
\HardButStrongLineBreak\kern3pt%
\begingroup \smaller\smaller\smaller\begin{tabular}{@{}c@{}}%
-1\\-102\\-358
\end{tabular}\endgroup%
\HardButStrongLineBreak\kern3pt%
\begingroup \smaller\smaller\smaller\begin{tabular}{@{}c@{}}%
-1\\-105\\-360
\end{tabular}\endgroup%
\HardButStrongLineBreak\kern3pt%
\begingroup \smaller\smaller\smaller\begin{tabular}{@{}c@{}}%
3\\304\\1072
\end{tabular}\endgroup%
\HardButStrongLineBreak\kern3pt%
\begingroup \smaller\smaller\smaller\begin{tabular}{@{}c@{}}%
6\\610\\2145
\end{tabular}\endgroup%
\HardButStrongLineBreak\kern3pt%
\begingroup \smaller\smaller\smaller\begin{tabular}{@{}c@{}}%
29\\2952\\10368
\end{tabular}\endgroup%
{$\left.\llap{\phantom{%
\begingroup \smaller\smaller\smaller\begin{tabular}{@{}c@{}}%
0\\0\\0
\end{tabular}\endgroup%
}}\!\right]$}%
%
%
\hbox{}\par\smallskip%
%
%
\leavevmode%
${L_{323.5}}$%
{} : {$[1\above{1pt}{1pt}{-}{}2\above{1pt}{1pt}{1}{}]\above{1pt}{1pt}{}{4}16\above{1pt}{1pt}{-}{3}{\cdot}1\above{1pt}{1pt}{1}{}3\above{1pt}{1pt}{-}{}9\above{1pt}{1pt}{1}{}{\cdot}1\above{1pt}{1pt}{-2}{}5\above{1pt}{1pt}{-}{}$}\EasyButWeakLineBreak%
{${40}\above{1pt}{1pt}{*}{2}{36}\above{1pt}{1pt}{s}{2}{240}\above{1pt}{1pt}{s}{2}{4}\above{1pt}{1pt}{*}{2}{360}\above{1pt}{1pt}{s}{2}{16}\above{1pt}{1pt}{*}{2}{60}\above{1pt}{1pt}{*}{2}{144}\above{1pt}{1pt}{s}{2}$}%
\nopagebreak\par%
\nopagebreak\par\leavevmode%
{$\left[\!\llap{\phantom{%
\begingroup \smaller\smaller\smaller\begin{tabular}{@{}c@{}}%
0\\0\\0
\end{tabular}\endgroup%
}}\right.$}%
\begingroup \smaller\smaller\smaller\begin{tabular}{@{}c@{}}%
-1711440\\-25200\\10080
\end{tabular}\endgroup%
\kern3pt%
\begingroup \smaller\smaller\smaller\begin{tabular}{@{}c@{}}%
-25200\\-354\\138
\end{tabular}\endgroup%
\kern3pt%
\begingroup \smaller\smaller\smaller\begin{tabular}{@{}c@{}}%
10080\\138\\-53
\end{tabular}\endgroup%
{$\left.\llap{\phantom{%
\begingroup \smaller\smaller\smaller\begin{tabular}{@{}c@{}}%
0\\0\\0
\end{tabular}\endgroup%
}}\!\right]$}%
\EasyButWeakLineBreak%
{$\left[\!\llap{\phantom{%
\begingroup \smaller\smaller\smaller\begin{tabular}{@{}c@{}}%
0\\0\\0
\end{tabular}\endgroup%
}}\right.$}%
\begingroup \smaller\smaller\smaller\begin{tabular}{@{}c@{}}%
7\\-1370\\-2240
\end{tabular}\endgroup%
\HardButStrongLineBreak\kern3pt%
\begingroup \smaller\smaller\smaller\begin{tabular}{@{}c@{}}%
7\\-1374\\-2250
\end{tabular}\endgroup%
\HardButStrongLineBreak\kern3pt%
\begingroup \smaller\smaller\smaller\begin{tabular}{@{}c@{}}%
13\\-2560\\-4200
\end{tabular}\endgroup%
\HardButStrongLineBreak\kern3pt%
\begingroup \smaller\smaller\smaller\begin{tabular}{@{}c@{}}%
1\\-198\\-326
\end{tabular}\endgroup%
\HardButStrongLineBreak\kern3pt%
\begingroup \smaller\smaller\smaller\begin{tabular}{@{}c@{}}%
1\\-210\\-360
\end{tabular}\endgroup%
\HardButStrongLineBreak\kern3pt%
\begingroup \smaller\smaller\smaller\begin{tabular}{@{}c@{}}%
-1\\196\\320
\end{tabular}\endgroup%
\HardButStrongLineBreak\kern3pt%
\begingroup \smaller\smaller\smaller\begin{tabular}{@{}c@{}}%
-1\\200\\330
\end{tabular}\endgroup%
\HardButStrongLineBreak\kern3pt%
\begingroup \smaller\smaller\smaller\begin{tabular}{@{}c@{}}%
5\\-972\\-1584
\end{tabular}\endgroup%
{$\left.\llap{\phantom{%
\begingroup \smaller\smaller\smaller\begin{tabular}{@{}c@{}}%
0\\0\\0
\end{tabular}\endgroup%
}}\!\right]$}%

\medskip%
%
\leavevmode\llap{}%
$W_{324}$%
\qquad\llap{38} lattices, $\chi=72$%
\hfill%
$22|22|22|22|22|22|22|22|\rtimes D_{8}$%
\nopagebreak\smallskip\hrule\nopagebreak\medskip%
%
%
\leavevmode%
${L_{324.1}}$%
{} : {$1\above{1pt}{1pt}{-2}{4}8\above{1pt}{1pt}{-}{3}{\cdot}1\above{1pt}{1pt}{-}{}3\above{1pt}{1pt}{-}{}9\above{1pt}{1pt}{-}{}{\cdot}1\above{1pt}{1pt}{2}{}5\above{1pt}{1pt}{1}{}$}\EasyButWeakLineBreak%
{${24}\above{1pt}{1pt}{}{2}{5}\above{1pt}{1pt}{r}{2}{72}\above{1pt}{1pt}{s}{2}{20}\above{1pt}{1pt}{*}{2}{24}\above{1pt}{1pt}{*}{2}{180}\above{1pt}{1pt}{s}{2}{8}\above{1pt}{1pt}{l}{2}{45}\above{1pt}{1pt}{}{2}$}\relax$\,(\times2)$%
\nopagebreak\par%
\nopagebreak\par\leavevmode%
{$\left[\!\llap{\phantom{%
\begingroup \smaller\smaller\smaller
\endgroup%
}}\!\right]$}%
%
%
\hbox{}\par\smallskip%
%
%
\leavevmode%
${L_{324.2}}$%
{} : {$[1\above{1pt}{1pt}{-}{}2\above{1pt}{1pt}{-}{}]\above{1pt}{1pt}{}{0}16\above{1pt}{1pt}{1}{7}{\cdot}1\above{1pt}{1pt}{-}{}3\above{1pt}{1pt}{-}{}9\above{1pt}{1pt}{-}{}{\cdot}1\above{1pt}{1pt}{2}{}5\above{1pt}{1pt}{1}{}$}\spacer%
\instructions{2}%
\EasyButWeakLineBreak%
{${24}\above{1pt}{1pt}{*}{2}{720}\above{1pt}{1pt}{s}{2}{8}\above{1pt}{1pt}{l}{2}{45}\above{1pt}{1pt}{}{2}{6}\above{1pt}{1pt}{}{2}{5}\above{1pt}{1pt}{r}{2}{72}\above{1pt}{1pt}{s}{2}{80}\above{1pt}{1pt}{*}{2}$}\relax$\,(\times2)$%
\nopagebreak\par%
\nopagebreak\par\leavevmode%
{$\left[\!\llap{\phantom{%
\begingroup \smaller\smaller\smaller
\endgroup%
}}\!\right]$}%
%
%
\hbox{}\par\smallskip%
%
%
\leavevmode%
${L_{324.3}}$%
{} : {$[1\above{1pt}{1pt}{1}{}2\above{1pt}{1pt}{-}{}]\above{1pt}{1pt}{}{4}16\above{1pt}{1pt}{-}{3}{\cdot}1\above{1pt}{1pt}{-}{}3\above{1pt}{1pt}{-}{}9\above{1pt}{1pt}{-}{}{\cdot}1\above{1pt}{1pt}{2}{}5\above{1pt}{1pt}{1}{}$}\spacer%
\instructions{m}%
\EasyButWeakLineBreak%
{${24}\above{1pt}{1pt}{s}{2}{720}\above{1pt}{1pt}{*}{2}{8}\above{1pt}{1pt}{*}{2}{180}\above{1pt}{1pt}{l}{2}{6}\above{1pt}{1pt}{r}{2}{20}\above{1pt}{1pt}{*}{2}{72}\above{1pt}{1pt}{*}{2}{80}\above{1pt}{1pt}{s}{2}$}\relax$\,(\times2)$%
\nopagebreak\par%
\nopagebreak\par\leavevmode%
{$\left[\!\llap{\phantom{%
\begingroup \smaller\smaller\smaller
\endgroup%
}}\!\right]$}%
%
%
\hbox{}\par\smallskip%
%
%
\leavevmode%
${L_{324.4}}$%
{} : {$[1\above{1pt}{1pt}{1}{}2\above{1pt}{1pt}{1}{}]\above{1pt}{1pt}{}{2}16\above{1pt}{1pt}{1}{1}{\cdot}1\above{1pt}{1pt}{-}{}3\above{1pt}{1pt}{-}{}9\above{1pt}{1pt}{-}{}{\cdot}1\above{1pt}{1pt}{2}{}5\above{1pt}{1pt}{1}{}$}\spacer%
\instructions{m}%
\EasyButWeakLineBreak%
{${6}\above{1pt}{1pt}{r}{2}{720}\above{1pt}{1pt}{l}{2}{2}\above{1pt}{1pt}{r}{2}{180}\above{1pt}{1pt}{*}{2}{24}\above{1pt}{1pt}{*}{2}{20}\above{1pt}{1pt}{l}{2}{18}\above{1pt}{1pt}{r}{2}{80}\above{1pt}{1pt}{l}{2}$}\relax$\,(\times2)$%
\nopagebreak\par%
\nopagebreak\par\leavevmode%
{$\left[\!\llap{\phantom{%
\begingroup \smaller\smaller\smaller
\endgroup%
}}\!\right]$}%
%
%
\hbox{}\par\smallskip%
%
%
\leavevmode%
${L_{324.5}}$%
{} : {$1\above{1pt}{1pt}{-}{5}8\above{1pt}{1pt}{1}{1}64\above{1pt}{1pt}{-}{5}{\cdot}1\above{1pt}{1pt}{-}{}3\above{1pt}{1pt}{-}{}9\above{1pt}{1pt}{-}{}{\cdot}1\above{1pt}{1pt}{2}{}5\above{1pt}{1pt}{1}{}$}\spacer%
\instructions{3,2}%
\EasyButWeakLineBreak%
{${96}\above{1pt}{1pt}{*}{2}{20}\above{1pt}{1pt}{l}{2}{72}\above{1pt}{1pt}{}{2}{320}\above{1pt}{1pt}{r}{2}{24}\above{1pt}{1pt}{b}{2}{2880}\above{1pt}{1pt}{l}{2}{8}\above{1pt}{1pt}{}{2}{45}\above{1pt}{1pt}{r}{2}$}\relax$\,(\times2)$%
\nopagebreak\par%
shares genus with 2-dual${}\iso{}$3-dual; isometric to own %
2.3-dual\nopagebreak\par%
\nopagebreak\par\leavevmode%
{$\left[\!\llap{\phantom{%
\begingroup \smaller\smaller\smaller
\endgroup%
}}\!\right]$}%

\medskip%
%
\leavevmode\llap{}%
$W_{325}$%
\qquad\llap{12} lattices, $\chi=36$%
\hfill%
$2222222222\rtimes C_{2}$%
\nopagebreak\smallskip\hrule\nopagebreak\medskip%
%
%
\leavevmode%
${L_{325.1}}$%
{} : {$1\above{1pt}{1pt}{-2}{{\rm II}}4\above{1pt}{1pt}{1}{1}{\cdot}1\above{1pt}{1pt}{2}{}9\above{1pt}{1pt}{1}{}{\cdot}1\above{1pt}{1pt}{1}{}5\above{1pt}{1pt}{-}{}25\above{1pt}{1pt}{-}{}$}\spacer%
\instructions{2}%
\EasyButWeakLineBreak%
{${50}\above{1pt}{1pt}{l}{2}{36}\above{1pt}{1pt}{r}{2}{10}\above{1pt}{1pt}{l}{2}{4}\above{1pt}{1pt}{r}{2}{90}\above{1pt}{1pt}{b}{2}$}\relax$\,(\times2)$%
\nopagebreak\par%
\nopagebreak\par\leavevmode%
{$\left[\!\llap{\phantom{%
\begingroup \smaller\smaller\smaller
\endgroup%
}}\!\right]$}%

\medskip%
%
\leavevmode\llap{}%
$W_{326}$%
\qquad\llap{4} lattices, $\chi=48$%
\hfill%
$6222362223\rtimes C_{2}$%
\nopagebreak\smallskip\hrule\nopagebreak\medskip%
%
%
\leavevmode%
${L_{326.1}}$%
{} : {$1\above{1pt}{1pt}{-2}{{\rm II}}16\above{1pt}{1pt}{-}{5}{\cdot}1\above{1pt}{1pt}{-}{}3\above{1pt}{1pt}{-}{}27\above{1pt}{1pt}{-}{}$}\EasyButWeakLineBreak%
{${6}\above{1pt}{1pt}{}{6}{2}\above{1pt}{1pt}{b}{2}{54}\above{1pt}{1pt}{s}{2}{2}\above{1pt}{1pt}{b}{2}{6}\above{1pt}{1pt}{-}{3}$}\relax$\,(\times2)$%
\nopagebreak\par%
\nopagebreak\par\leavevmode%
{$\left[\!\llap{\phantom{%
\begingroup \smaller\smaller\smaller
\endgroup%
}}\!\right]$}%

\medskip%
%
\leavevmode\llap{}%
$W_{327}$%
\qquad\llap{32} lattices, $\chi=36$%
\hfill%
$222|22222|22\rtimes D_{2}$%
\nopagebreak\smallskip\hrule\nopagebreak\medskip%
%
%
\leavevmode%
${L_{327.1}}$%
{} : {$[1\above{1pt}{1pt}{-}{}2\above{1pt}{1pt}{1}{}]\above{1pt}{1pt}{}{6}32\above{1pt}{1pt}{1}{7}{\cdot}1\above{1pt}{1pt}{2}{}3\above{1pt}{1pt}{-}{}{\cdot}1\above{1pt}{1pt}{2}{}7\above{1pt}{1pt}{1}{}$}\EasyButWeakLineBreak%
{${7}\above{1pt}{1pt}{}{2}{6}\above{1pt}{1pt}{r}{2}{224}\above{1pt}{1pt}{l}{2}{2}\above{1pt}{1pt}{}{2}{224}\above{1pt}{1pt}{}{2}{6}\above{1pt}{1pt}{r}{2}{28}\above{1pt}{1pt}{*}{2}{96}\above{1pt}{1pt}{*}{2}{56}\above{1pt}{1pt}{s}{2}{96}\above{1pt}{1pt}{l}{2}$}%
\nopagebreak\par%
\nopagebreak\par\leavevmode%
{$\left[\!\llap{\phantom{%
\begingroup \smaller\smaller\smaller
\endgroup%
}}\!\right]$}%
%
%
\hbox{}\par\smallskip%
%
%
\leavevmode%
${L_{327.2}}$%
{} : {$[1\above{1pt}{1pt}{1}{}2\above{1pt}{1pt}{1}{}]\above{1pt}{1pt}{}{6}32\above{1pt}{1pt}{-}{3}{\cdot}1\above{1pt}{1pt}{2}{}3\above{1pt}{1pt}{-}{}{\cdot}1\above{1pt}{1pt}{2}{}7\above{1pt}{1pt}{1}{}$}\EasyButWeakLineBreak%
{${7}\above{1pt}{1pt}{r}{2}{24}\above{1pt}{1pt}{*}{2}{224}\above{1pt}{1pt}{s}{2}{8}\above{1pt}{1pt}{*}{2}{224}\above{1pt}{1pt}{s}{2}{24}\above{1pt}{1pt}{*}{2}{28}\above{1pt}{1pt}{s}{2}{96}\above{1pt}{1pt}{l}{2}{14}\above{1pt}{1pt}{}{2}{96}\above{1pt}{1pt}{}{2}$}%
\nopagebreak\par%
\nopagebreak\par\leavevmode%
{$\left[\!\llap{\phantom{%
\begingroup \smaller\smaller\smaller
\endgroup%
}}\!\right]$}%
%
%
\hbox{}\par\smallskip%
%
%
\leavevmode%
${L_{327.3}}$%
{} : {$1\above{1pt}{1pt}{1}{7}4\above{1pt}{1pt}{1}{1}32\above{1pt}{1pt}{-}{5}{\cdot}1\above{1pt}{1pt}{2}{}3\above{1pt}{1pt}{1}{}{\cdot}1\above{1pt}{1pt}{2}{}7\above{1pt}{1pt}{1}{}$}\EasyButWeakLineBreak%
{${224}\above{1pt}{1pt}{s}{2}{48}\above{1pt}{1pt}{*}{2}{28}\above{1pt}{1pt}{l}{2}{4}\above{1pt}{1pt}{}{2}{7}\above{1pt}{1pt}{r}{2}{48}\above{1pt}{1pt}{*}{2}{224}\above{1pt}{1pt}{l}{2}{3}\above{1pt}{1pt}{r}{2}{112}\above{1pt}{1pt}{*}{2}{12}\above{1pt}{1pt}{*}{2}$}%
\nopagebreak\par%
\nopagebreak\par\leavevmode%
{$\left[\!\llap{\phantom{%
\begingroup \smaller\smaller\smaller
\endgroup%
}}\!\right]$}%
%
%
\hbox{}\par\smallskip%
%
%
\leavevmode%
${L_{327.4}}$%
{} : {$1\above{1pt}{1pt}{-}{3}4\above{1pt}{1pt}{1}{7}32\above{1pt}{1pt}{1}{7}{\cdot}1\above{1pt}{1pt}{2}{}3\above{1pt}{1pt}{1}{}{\cdot}1\above{1pt}{1pt}{2}{}7\above{1pt}{1pt}{1}{}$}\EasyButWeakLineBreak%
{${224}\above{1pt}{1pt}{}{2}{12}\above{1pt}{1pt}{r}{2}{28}\above{1pt}{1pt}{*}{2}{16}\above{1pt}{1pt}{l}{2}{7}\above{1pt}{1pt}{}{2}{12}\above{1pt}{1pt}{r}{2}{224}\above{1pt}{1pt}{s}{2}{12}\above{1pt}{1pt}{l}{2}{28}\above{1pt}{1pt}{}{2}{3}\above{1pt}{1pt}{}{2}$}%
\nopagebreak\par%
\nopagebreak\par\leavevmode%
{$\left[\!\llap{\phantom{%
\begingroup \smaller\smaller\smaller
\endgroup%
}}\!\right]$}%

\medskip%
%
\leavevmode\llap{}%
$W_{328}$%
\qquad\llap{32} lattices, $\chi=96$%
\hfill%
$2|22222|22222|22222|2222\rtimes D_{4}$%
\nopagebreak\smallskip\hrule\nopagebreak\medskip%
%
%
\leavevmode%
${L_{328.1}}$%
{} : {$[1\above{1pt}{1pt}{1}{}2\above{1pt}{1pt}{1}{}]\above{1pt}{1pt}{}{0}32\above{1pt}{1pt}{-}{5}{\cdot}1\above{1pt}{1pt}{2}{}3\above{1pt}{1pt}{-}{}{\cdot}1\above{1pt}{1pt}{-2}{}7\above{1pt}{1pt}{-}{}$}\EasyButWeakLineBreak%
{${672}\above{1pt}{1pt}{}{2}{2}\above{1pt}{1pt}{r}{2}{672}\above{1pt}{1pt}{s}{2}{4}\above{1pt}{1pt}{*}{2}{24}\above{1pt}{1pt}{s}{2}{32}\above{1pt}{1pt}{*}{2}{168}\above{1pt}{1pt}{s}{2}{32}\above{1pt}{1pt}{*}{2}{24}\above{1pt}{1pt}{l}{2}{1}\above{1pt}{1pt}{}{2}$}\relax$\,(\times2)$%
\nopagebreak\par%
\nopagebreak\par\leavevmode%
{$\left[\!\llap{\phantom{%
\begingroup \smaller\smaller\smaller
\endgroup%
}}\!\right]$}%
%
%
\hbox{}\par\smallskip%
%
%
\leavevmode%
${L_{328.2}}$%
{} : {$[1\above{1pt}{1pt}{1}{}2\above{1pt}{1pt}{-}{}]\above{1pt}{1pt}{}{4}32\above{1pt}{1pt}{1}{1}{\cdot}1\above{1pt}{1pt}{2}{}3\above{1pt}{1pt}{-}{}{\cdot}1\above{1pt}{1pt}{-2}{}7\above{1pt}{1pt}{-}{}$}\EasyButWeakLineBreak%
{${672}\above{1pt}{1pt}{*}{2}{8}\above{1pt}{1pt}{s}{2}{672}\above{1pt}{1pt}{*}{2}{4}\above{1pt}{1pt}{l}{2}{6}\above{1pt}{1pt}{}{2}{32}\above{1pt}{1pt}{}{2}{42}\above{1pt}{1pt}{r}{2}{32}\above{1pt}{1pt}{l}{2}{6}\above{1pt}{1pt}{}{2}{1}\above{1pt}{1pt}{r}{2}$}\relax$\,(\times2)$%
\nopagebreak\par%
\nopagebreak\par\leavevmode%
{$\left[\!\llap{\phantom{%
\begingroup \smaller\smaller\smaller
\endgroup%
}}\!\right]$}%
%
%
\hbox{}\par\smallskip%
%
%
\leavevmode%
${L_{328.3}}$%
{} : {$1\above{1pt}{1pt}{-}{5}4\above{1pt}{1pt}{1}{1}32\above{1pt}{1pt}{1}{7}{\cdot}1\above{1pt}{1pt}{2}{}3\above{1pt}{1pt}{1}{}{\cdot}1\above{1pt}{1pt}{-2}{}7\above{1pt}{1pt}{-}{}$}\EasyButWeakLineBreak%
{${84}\above{1pt}{1pt}{l}{2}{4}\above{1pt}{1pt}{}{2}{21}\above{1pt}{1pt}{r}{2}{32}\above{1pt}{1pt}{*}{2}{48}\above{1pt}{1pt}{l}{2}{1}\above{1pt}{1pt}{}{2}{84}\above{1pt}{1pt}{r}{2}{4}\above{1pt}{1pt}{*}{2}{48}\above{1pt}{1pt}{s}{2}{32}\above{1pt}{1pt}{*}{2}$}\relax$\,(\times2)$%
\nopagebreak\par%
\nopagebreak\par\leavevmode%
{$\left[\!\llap{\phantom{%
\begingroup \smaller\smaller\smaller
\endgroup%
}}\!\right]$}%
%
%
\hbox{}\par\smallskip%
%
%
\leavevmode%
${L_{328.4}}$%
{} : {$1\above{1pt}{1pt}{-}{5}4\above{1pt}{1pt}{1}{7}32\above{1pt}{1pt}{1}{1}{\cdot}1\above{1pt}{1pt}{2}{}3\above{1pt}{1pt}{1}{}{\cdot}1\above{1pt}{1pt}{-2}{}7\above{1pt}{1pt}{-}{}$}\EasyButWeakLineBreak%
{${21}\above{1pt}{1pt}{r}{2}{16}\above{1pt}{1pt}{*}{2}{84}\above{1pt}{1pt}{s}{2}{32}\above{1pt}{1pt}{l}{2}{12}\above{1pt}{1pt}{}{2}{1}\above{1pt}{1pt}{r}{2}{336}\above{1pt}{1pt}{*}{2}{4}\above{1pt}{1pt}{l}{2}{12}\above{1pt}{1pt}{}{2}{32}\above{1pt}{1pt}{}{2}$}\relax$\,(\times2)$%
\nopagebreak\par%
\nopagebreak\par\leavevmode%
{$\left[\!\llap{\phantom{%
\begingroup \smaller\smaller\smaller
\endgroup%
}}\!\right]$}%

\medskip%
%
\leavevmode\llap{}%
$W_{329}$%
\qquad\llap{8} lattices, $\chi=144$%
\hfill%
$\infty\infty\infty2|2\infty\infty\infty|\infty\infty\infty2|2\infty\infty\infty|\rtimes D_{4}$%
\nopagebreak\smallskip\hrule\nopagebreak\medskip%
%
%
\leavevmode%
${L_{329.1}}$%
{} : {$1\above{1pt}{1pt}{2}{0}4\above{1pt}{1pt}{1}{7}{\cdot}1\above{1pt}{1pt}{-}{}7\above{1pt}{1pt}{1}{}49\above{1pt}{1pt}{1}{}$}\EasyButWeakLineBreak%
{${28}\above{1pt}{1pt}{7,4}{\infty z}{7}\above{1pt}{1pt}{28,11}{\infty}{28}\above{1pt}{1pt}{7,1}{\infty b}{28}\above{1pt}{1pt}{}{2}{49}\above{1pt}{1pt}{}{2}{7}\above{1pt}{1pt}{28,15}{\infty}{28}\above{1pt}{1pt}{7,4}{\infty a}{28}\above{1pt}{1pt}{14,9}{\infty}$}\relax$\,(\times2)$%
\nopagebreak\par%
\nopagebreak\par\leavevmode%
{$\left[\!\llap{\phantom{%
\begingroup \smaller\smaller\smaller
\endgroup%
}}\!\right]$}%
%
%
\hbox{}\par\smallskip%
%
%
\leavevmode%
${L_{329.2}}$%
{} : {$1\above{1pt}{1pt}{2}{{\rm II}}8\above{1pt}{1pt}{1}{7}{\cdot}1\above{1pt}{1pt}{-}{}7\above{1pt}{1pt}{1}{}49\above{1pt}{1pt}{1}{}$}\EasyButWeakLineBreak%
{${56}\above{1pt}{1pt}{7,2}{\infty z}{14}\above{1pt}{1pt}{28,11}{\infty b}{56}\above{1pt}{1pt}{7,2}{\infty a}{56}\above{1pt}{1pt}{r}{2}{98}\above{1pt}{1pt}{s}{2}{14}\above{1pt}{1pt}{28,15}{\infty a}{56}\above{1pt}{1pt}{7,1}{\infty b}{56}\above{1pt}{1pt}{7,2}{\infty}$}\relax$\,(\times2)$%
\nopagebreak\par%
\nopagebreak\par\leavevmode%
{$\left[\!\llap{\phantom{%
\begingroup \smaller\smaller\smaller
\endgroup%
}}\!\right]$}%

\medskip%
%
\leavevmode\llap{}%
$W_{330}$%
\qquad\llap{8} lattices, $\chi=36$%
\hfill%
$2\infty2\infty\infty\infty$%
\nopagebreak\smallskip\hrule\nopagebreak\medskip%
%
%
\leavevmode%
${L_{330.1}}$%
{} : {$1\above{1pt}{1pt}{-2}{{\rm II}}8\above{1pt}{1pt}{-}{5}{\cdot}1\above{1pt}{1pt}{1}{}9\above{1pt}{1pt}{-}{}81\above{1pt}{1pt}{-}{}$}\spacer%
\instructions{2}%
\EasyButWeakLineBreak%
{${72}\above{1pt}{1pt}{b}{2}{162}\above{1pt}{1pt}{12,1}{\infty b}{648}\above{1pt}{1pt}{b}{2}{18}\above{1pt}{1pt}{36,1}{\infty b}{72}\above{1pt}{1pt}{18,13}{\infty z}{18}\above{1pt}{1pt}{36,13}{\infty a}$}%
\nopagebreak\par%
\nopagebreak\par\leavevmode%
{$\left[\!\llap{\phantom{%
\begingroup \smaller\smaller\smaller\begin{tabular}{@{}c@{}}%
0\\0\\0
\end{tabular}\endgroup%
}}\right.$}%
\begingroup \smaller\smaller\smaller\begin{tabular}{@{}c@{}}%
-1577880\\14256\\-166536
\end{tabular}\endgroup%
\kern3pt%
\begingroup \smaller\smaller\smaller\begin{tabular}{@{}c@{}}%
14256\\-90\\1521
\end{tabular}\endgroup%
\kern3pt%
\begingroup \smaller\smaller\smaller\begin{tabular}{@{}c@{}}%
-166536\\1521\\-17570
\end{tabular}\endgroup%
{$\left.\llap{\phantom{%
\begingroup \smaller\smaller\smaller\begin{tabular}{@{}c@{}}%
0\\0\\0
\end{tabular}\endgroup%
}}\!\right]$}%
\EasyButWeakLineBreak%
{$\left[\!\llap{\phantom{%
\begingroup \smaller\smaller\smaller\begin{tabular}{@{}c@{}}%
0\\0\\0
\end{tabular}\endgroup%
}}\right.$}%
\begingroup \smaller\smaller\smaller\begin{tabular}{@{}c@{}}%
63\\244\\-576
\end{tabular}\endgroup%
\HardButStrongLineBreak\kern3pt%
\begingroup \smaller\smaller\smaller\begin{tabular}{@{}c@{}}%
124\\477\\-1134
\end{tabular}\endgroup%
\HardButStrongLineBreak\kern3pt%
\begingroup \smaller\smaller\smaller\begin{tabular}{@{}c@{}}%
-815\\-3168\\7452
\end{tabular}\endgroup%
\HardButStrongLineBreak\kern3pt%
\begingroup \smaller\smaller\smaller\begin{tabular}{@{}c@{}}%
-375\\-1454\\3429
\end{tabular}\endgroup%
\HardButStrongLineBreak\kern3pt%
\begingroup \smaller\smaller\smaller\begin{tabular}{@{}c@{}}%
-811\\-3140\\7416
\end{tabular}\endgroup%
\HardButStrongLineBreak\kern3pt%
\begingroup \smaller\smaller\smaller\begin{tabular}{@{}c@{}}%
-62\\-239\\567
\end{tabular}\endgroup%
{$\left.\llap{\phantom{%
\begingroup \smaller\smaller\smaller\begin{tabular}{@{}c@{}}%
0\\0\\0
\end{tabular}\endgroup%
}}\!\right]$}%

\medskip%
%
\leavevmode\llap{}%
$W_{331}$%
\qquad\llap{12} lattices, $\chi=20$%
\hfill%
$\slashthree222|222\rtimes D_{2}$%
\nopagebreak\smallskip\hrule\nopagebreak\medskip%
%
%
\leavevmode%
${L_{331.1}}$%
{} : {$1\above{1pt}{1pt}{-2}{{\rm II}}4\above{1pt}{1pt}{1}{1}{\cdot}1\above{1pt}{1pt}{1}{}3\above{1pt}{1pt}{-}{}9\above{1pt}{1pt}{1}{}{\cdot}1\above{1pt}{1pt}{-2}{}5\above{1pt}{1pt}{1}{}{\cdot}1\above{1pt}{1pt}{2}{}11\above{1pt}{1pt}{-}{}$}\spacer%
\instructions{2}%
\EasyButWeakLineBreak%
{${6}\above{1pt}{1pt}{+}{3}{6}\above{1pt}{1pt}{b}{2}{22}\above{1pt}{1pt}{l}{2}{36}\above{1pt}{1pt}{r}{2}{330}\above{1pt}{1pt}{l}{2}{4}\above{1pt}{1pt}{r}{2}{198}\above{1pt}{1pt}{b}{2}$}%
\nopagebreak\par%
\nopagebreak\par\leavevmode%
{$\left[\!\llap{\phantom{%
\begingroup \smaller\smaller\smaller\begin{tabular}{@{}c@{}}%
0\\0\\0
\end{tabular}\endgroup%
}}\right.$}%
\begingroup \smaller\smaller\smaller\begin{tabular}{@{}c@{}}%
-100585980\\132660\\-7815060
\end{tabular}\endgroup%
\kern3pt%
\begingroup \smaller\smaller\smaller\begin{tabular}{@{}c@{}}%
132660\\-174\\10491
\end{tabular}\endgroup%
\kern3pt%
\begingroup \smaller\smaller\smaller\begin{tabular}{@{}c@{}}%
-7815060\\10491\\-572006
\end{tabular}\endgroup%
{$\left.\llap{\phantom{%
\begingroup \smaller\smaller\smaller\begin{tabular}{@{}c@{}}%
0\\0\\0
\end{tabular}\endgroup%
}}\!\right]$}%
\EasyButWeakLineBreak%
{$\left[\!\llap{\phantom{%
\begingroup \smaller\smaller\smaller\begin{tabular}{@{}c@{}}%
0\\0\\0
\end{tabular}\endgroup%
}}\right.$}%
\begingroup \smaller\smaller\smaller\begin{tabular}{@{}c@{}}%
99\\57391\\-300
\end{tabular}\endgroup%
\HardButStrongLineBreak\kern3pt%
\begingroup \smaller\smaller\smaller\begin{tabular}{@{}c@{}}%
-98\\-56810\\297
\end{tabular}\endgroup%
\HardButStrongLineBreak\kern3pt%
\begingroup \smaller\smaller\smaller\begin{tabular}{@{}c@{}}%
-657\\-380864\\1991
\end{tabular}\endgroup%
\HardButStrongLineBreak\kern3pt%
\begingroup \smaller\smaller\smaller\begin{tabular}{@{}c@{}}%
-1283\\-743760\\3888
\end{tabular}\endgroup%
\HardButStrongLineBreak\kern3pt%
\begingroup \smaller\smaller\smaller\begin{tabular}{@{}c@{}}%
-2668\\-1546655\\8085
\end{tabular}\endgroup%
\HardButStrongLineBreak\kern3pt%
\begingroup \smaller\smaller\smaller\begin{tabular}{@{}c@{}}%
-165\\-95652\\500
\end{tabular}\endgroup%
\HardButStrongLineBreak\kern3pt%
\begingroup \smaller\smaller\smaller\begin{tabular}{@{}c@{}}%
196\\113619\\-594
\end{tabular}\endgroup%
{$\left.\llap{\phantom{%
\begingroup \smaller\smaller\smaller\begin{tabular}{@{}c@{}}%
0\\0\\0
\end{tabular}\endgroup%
}}\!\right]$}%

\medskip%
%
\leavevmode\llap{}%
$W_{332}$%
\qquad\llap{12} lattices, $\chi=40$%
\hfill%
$6|62|26|62|2\rtimes D_{4}$%
\nopagebreak\smallskip\hrule\nopagebreak\medskip%
%
%
\leavevmode%
${L_{332.1}}$%
{} : {$1\above{1pt}{1pt}{-2}{{\rm II}}4\above{1pt}{1pt}{1}{1}{\cdot}1\above{1pt}{1pt}{-}{}3\above{1pt}{1pt}{-}{}9\above{1pt}{1pt}{-}{}{\cdot}1\above{1pt}{1pt}{-2}{}5\above{1pt}{1pt}{1}{}{\cdot}1\above{1pt}{1pt}{2}{}11\above{1pt}{1pt}{-}{}$}\spacer%
\instructions{2}%
\EasyButWeakLineBreak%
{${18}\above{1pt}{1pt}{}{6}{6}\above{1pt}{1pt}{}{6}{2}\above{1pt}{1pt}{b}{2}{330}\above{1pt}{1pt}{b}{2}$}\relax$\,(\times2)$%
\nopagebreak\par%
\nopagebreak\par\leavevmode%
{$\left[\!\llap{\phantom{%
\begingroup \smaller\smaller\smaller
\endgroup%
}}\!\right]$}%

\medskip%
%
\leavevmode\llap{}%
$W_{333}$%
\qquad\llap{12} lattices, $\chi=72$%
\hfill%
$2222|2222|2222|2222|\rtimes D_{4}$%
\nopagebreak\smallskip\hrule\nopagebreak\medskip%
%
%
\leavevmode%
${L_{333.1}}$%
{} : {$1\above{1pt}{1pt}{-2}{{\rm II}}4\above{1pt}{1pt}{1}{1}{\cdot}1\above{1pt}{1pt}{-}{}3\above{1pt}{1pt}{-}{}9\above{1pt}{1pt}{-}{}{\cdot}1\above{1pt}{1pt}{2}{}5\above{1pt}{1pt}{-}{}{\cdot}1\above{1pt}{1pt}{-2}{}11\above{1pt}{1pt}{1}{}$}\spacer%
\instructions{2}%
\EasyButWeakLineBreak%
{${132}\above{1pt}{1pt}{r}{2}{2}\above{1pt}{1pt}{b}{2}{60}\above{1pt}{1pt}{*}{2}{44}\above{1pt}{1pt}{b}{2}{6}\above{1pt}{1pt}{b}{2}{396}\above{1pt}{1pt}{*}{2}{60}\above{1pt}{1pt}{b}{2}{18}\above{1pt}{1pt}{l}{2}$}\relax$\,(\times2)$%
\nopagebreak\par%
\nopagebreak\par\leavevmode%
{$\left[\!\llap{\phantom{%
\begingroup \smaller\smaller\smaller
\endgroup%
}}\!\right]$}%

\medskip%
%
\leavevmode\llap{}%
$W_{334}$%
\qquad\llap{4} lattices, $\chi=24$%
\hfill%
$\slashinfty2|2\slashinfty2|2\rtimes D_{4}$%
\nopagebreak\smallskip\hrule\nopagebreak\medskip%
%
%
\leavevmode%
${L_{334.1}}$%
{} : {$1\above{1pt}{1pt}{1}{7}8\above{1pt}{1pt}{-}{5}64\above{1pt}{1pt}{1}{7}{\cdot}1\above{1pt}{1pt}{-}{}3\above{1pt}{1pt}{-}{}9\above{1pt}{1pt}{1}{}$}\spacer%
\instructions{3}%
\EasyButWeakLineBreak%
{${24}\above{1pt}{1pt}{24,7}{\infty b}{96}\above{1pt}{1pt}{*}{2}{576}\above{1pt}{1pt}{s}{2}{96}\above{1pt}{1pt}{48,1}{\infty z}{24}\above{1pt}{1pt}{b}{2}{36}\above{1pt}{1pt}{s}{2}$}%
\nopagebreak\par%
\nopagebreak\par\leavevmode%
{$\left[\!\llap{\phantom{%
\begingroup \smaller\smaller\smaller\begin{tabular}{@{}c@{}}%
0\\0\\0
\end{tabular}\endgroup%
}}\right.$}%
\begingroup \smaller\smaller\smaller\begin{tabular}{@{}c@{}}%
-300096\\11520\\576
\end{tabular}\endgroup%
\kern3pt%
\begingroup \smaller\smaller\smaller\begin{tabular}{@{}c@{}}%
11520\\-408\\-24
\end{tabular}\endgroup%
\kern3pt%
\begingroup \smaller\smaller\smaller\begin{tabular}{@{}c@{}}%
576\\-24\\-1
\end{tabular}\endgroup%
{$\left.\llap{\phantom{%
\begingroup \smaller\smaller\smaller\begin{tabular}{@{}c@{}}%
0\\0\\0
\end{tabular}\endgroup%
}}\!\right]$}%
\EasyButWeakLineBreak%
{$\left[\!\llap{\phantom{%
\begingroup \smaller\smaller\smaller\begin{tabular}{@{}c@{}}%
0\\0\\0
\end{tabular}\endgroup%
}}\right.$}%
\begingroup \smaller\smaller\smaller\begin{tabular}{@{}c@{}}%
-3\\-43\\-732
\end{tabular}\endgroup%
\HardButStrongLineBreak\kern3pt%
\begingroup \smaller\smaller\smaller\begin{tabular}{@{}c@{}}%
-5\\-70\\-1248
\end{tabular}\endgroup%
\HardButStrongLineBreak\kern3pt%
\begingroup \smaller\smaller\smaller\begin{tabular}{@{}c@{}}%
-1\\-12\\-288
\end{tabular}\endgroup%
\HardButStrongLineBreak\kern3pt%
\begingroup \smaller\smaller\smaller\begin{tabular}{@{}c@{}}%
1\\14\\240
\end{tabular}\endgroup%
\HardButStrongLineBreak\kern3pt%
\begingroup \smaller\smaller\smaller\begin{tabular}{@{}c@{}}%
0\\-1\\12
\end{tabular}\endgroup%
\HardButStrongLineBreak\kern3pt%
\begingroup \smaller\smaller\smaller\begin{tabular}{@{}c@{}}%
-1\\-15\\-234
\end{tabular}\endgroup%
{$\left.\llap{\phantom{%
\begingroup \smaller\smaller\smaller\begin{tabular}{@{}c@{}}%
0\\0\\0
\end{tabular}\endgroup%
}}\!\right]$}%
%
%
%
%
%
%
%
%
%
%
%
%
%
%

\medskip%
%
\leavevmode\llap{}%
$W_{335}$%
\qquad\llap{8} lattices, $\chi=24$%
\hfill%
$22|22\infty|\infty\rtimes D_{2}$%
\nopagebreak\smallskip\hrule\nopagebreak\medskip%
%
%
\leavevmode%
${L_{335.1}}$%
{} : {$1\above{1pt}{1pt}{1}{1}8\above{1pt}{1pt}{1}{1}64\above{1pt}{1pt}{-}{5}{\cdot}1\above{1pt}{1pt}{-}{}3\above{1pt}{1pt}{-}{}9\above{1pt}{1pt}{1}{}$}\spacer%
\instructions{3}%
\EasyButWeakLineBreak%
{${96}\above{1pt}{1pt}{l}{2}{9}\above{1pt}{1pt}{}{2}{8}\above{1pt}{1pt}{r}{2}{36}\above{1pt}{1pt}{*}{2}{96}\above{1pt}{1pt}{48,7}{\infty z}{24}\above{1pt}{1pt}{24,13}{\infty b}$}%
\nopagebreak\par%
\nopagebreak\par\leavevmode%
{$\left[\!\llap{\phantom{%
\begingroup \smaller\smaller\smaller\begin{tabular}{@{}c@{}}%
0\\0\\0
\end{tabular}\endgroup%
}}\right.$}%
\begingroup \smaller\smaller\smaller\begin{tabular}{@{}c@{}}%
-241344\\5760\\59328
\end{tabular}\endgroup%
\kern3pt%
\begingroup \smaller\smaller\smaller\begin{tabular}{@{}c@{}}%
5760\\-120\\-1488
\end{tabular}\endgroup%
\kern3pt%
\begingroup \smaller\smaller\smaller\begin{tabular}{@{}c@{}}%
59328\\-1488\\-14287
\end{tabular}\endgroup%
{$\left.\llap{\phantom{%
\begingroup \smaller\smaller\smaller\begin{tabular}{@{}c@{}}%
0\\0\\0
\end{tabular}\endgroup%
}}\!\right]$}%
\EasyButWeakLineBreak%
{$\left[\!\llap{\phantom{%
\begingroup \smaller\smaller\smaller\begin{tabular}{@{}c@{}}%
0\\0\\0
\end{tabular}\endgroup%
}}\right.$}%
\begingroup \smaller\smaller\smaller\begin{tabular}{@{}c@{}}%
661\\7910\\1920
\end{tabular}\endgroup%
\HardButStrongLineBreak\kern3pt%
\begingroup \smaller\smaller\smaller\begin{tabular}{@{}c@{}}%
158\\1890\\459
\end{tabular}\endgroup%
\HardButStrongLineBreak\kern3pt%
\begingroup \smaller\smaller\smaller\begin{tabular}{@{}c@{}}%
11\\131\\32
\end{tabular}\endgroup%
\HardButStrongLineBreak\kern3pt%
\begingroup \smaller\smaller\smaller\begin{tabular}{@{}c@{}}%
-31\\-372\\-90
\end{tabular}\endgroup%
\HardButStrongLineBreak\kern3pt%
\begingroup \smaller\smaller\smaller\begin{tabular}{@{}c@{}}%
-33\\-394\\-96
\end{tabular}\endgroup%
\HardButStrongLineBreak\kern3pt%
\begingroup \smaller\smaller\smaller\begin{tabular}{@{}c@{}}%
62\\743\\180
\end{tabular}\endgroup%
{$\left.\llap{\phantom{%
\begingroup \smaller\smaller\smaller\begin{tabular}{@{}c@{}}%
0\\0\\0
\end{tabular}\endgroup%
}}\!\right]$}%

\medskip%
%
\leavevmode\llap{}%
$W_{336}$%
\qquad\llap{4} lattices, $\chi=12$%
\hfill%
$2|222|22\rtimes D_{2}$%
\nopagebreak\smallskip\hrule\nopagebreak\medskip%
%
%
\leavevmode%
${L_{336.1}}$%
{} : {$1\above{1pt}{1pt}{-}{3}8\above{1pt}{1pt}{1}{7}64\above{1pt}{1pt}{1}{1}{\cdot}1\above{1pt}{1pt}{2}{}3\above{1pt}{1pt}{1}{}$}\EasyButWeakLineBreak%
{${64}\above{1pt}{1pt}{b}{2}{4}\above{1pt}{1pt}{l}{2}{64}\above{1pt}{1pt}{}{2}{3}\above{1pt}{1pt}{r}{2}{32}\above{1pt}{1pt}{*}{2}{12}\above{1pt}{1pt}{s}{2}$}%
\nopagebreak\par%
\nopagebreak\par\leavevmode%
{$\left[\!\llap{\phantom{%
\begingroup \smaller\smaller\smaller\begin{tabular}{@{}c@{}}%
0\\0\\0
\end{tabular}\endgroup%
}}\right.$}%
\begingroup \smaller\smaller\smaller\begin{tabular}{@{}c@{}}%
1478208\\170112\\-4032
\end{tabular}\endgroup%
\kern3pt%
\begingroup \smaller\smaller\smaller\begin{tabular}{@{}c@{}}%
170112\\19576\\-464
\end{tabular}\endgroup%
\kern3pt%
\begingroup \smaller\smaller\smaller\begin{tabular}{@{}c@{}}%
-4032\\-464\\11
\end{tabular}\endgroup%
{$\left.\llap{\phantom{%
\begingroup \smaller\smaller\smaller\begin{tabular}{@{}c@{}}%
0\\0\\0
\end{tabular}\endgroup%
}}\!\right]$}%
\EasyButWeakLineBreak%
{$\left[\!\llap{\phantom{%
\begingroup \smaller\smaller\smaller\begin{tabular}{@{}c@{}}%
0\\0\\0
\end{tabular}\endgroup%
}}\right.$}%
\begingroup \smaller\smaller\smaller\begin{tabular}{@{}c@{}}%
1\\-8\\32
\end{tabular}\endgroup%
\HardButStrongLineBreak\kern3pt%
\begingroup \smaller\smaller\smaller\begin{tabular}{@{}c@{}}%
0\\-1\\-42
\end{tabular}\endgroup%
\HardButStrongLineBreak\kern3pt%
\begingroup \smaller\smaller\smaller\begin{tabular}{@{}c@{}}%
3\\-32\\-256
\end{tabular}\endgroup%
\HardButStrongLineBreak\kern3pt%
\begingroup \smaller\smaller\smaller\begin{tabular}{@{}c@{}}%
1\\-9\\-15
\end{tabular}\endgroup%
\HardButStrongLineBreak\kern3pt%
\begingroup \smaller\smaller\smaller\begin{tabular}{@{}c@{}}%
1\\-6\\112
\end{tabular}\endgroup%
\HardButStrongLineBreak\kern3pt%
\begingroup \smaller\smaller\smaller\begin{tabular}{@{}c@{}}%
1\\-6\\114
\end{tabular}\endgroup%
{$\left.\llap{\phantom{%
\begingroup \smaller\smaller\smaller\begin{tabular}{@{}c@{}}%
0\\0\\0
\end{tabular}\endgroup%
}}\!\right]$}%
%
%
%
%
%
%
%
%
%
%
%
%
%
%

\medskip%
%
\leavevmode\llap{}%
$W_{337}$%
\qquad\llap{16} lattices, $\chi=36$%
\hfill%
$2222|22222|2\rtimes D_{2}$%
\nopagebreak\smallskip\hrule\nopagebreak\medskip%
%
%
\leavevmode%
${L_{337.1}}$%
{} : {$1\above{1pt}{1pt}{-2}{2}16\above{1pt}{1pt}{-}{3}{\cdot}1\above{1pt}{1pt}{2}{}3\above{1pt}{1pt}{1}{}{\cdot}1\above{1pt}{1pt}{-2}{}5\above{1pt}{1pt}{-}{}{\cdot}1\above{1pt}{1pt}{2}{}7\above{1pt}{1pt}{-}{}$}\EasyButWeakLineBreak%
{${560}\above{1pt}{1pt}{}{2}{3}\above{1pt}{1pt}{r}{2}{140}\above{1pt}{1pt}{*}{2}{48}\above{1pt}{1pt}{b}{2}{10}\above{1pt}{1pt}{l}{2}{48}\above{1pt}{1pt}{}{2}{35}\above{1pt}{1pt}{r}{2}{12}\above{1pt}{1pt}{*}{2}{560}\above{1pt}{1pt}{b}{2}{2}\above{1pt}{1pt}{l}{2}$}%
\nopagebreak\par%
\nopagebreak\par\leavevmode%
{$\left[\!\llap{\phantom{%
\begingroup \smaller\smaller\smaller
\endgroup%
}}\!\right]$}%

\medskip%
%
\leavevmode\llap{}%
$W_{338}$%
\qquad\llap{12} lattices, $\chi=24$%
\hfill%
$2222|2222|\rtimes D_{2}$%
\nopagebreak\smallskip\hrule\nopagebreak\medskip%
%
%
\leavevmode%
${L_{338.1}}$%
{} : {$1\above{1pt}{1pt}{-2}{{\rm II}}4\above{1pt}{1pt}{1}{7}{\cdot}1\above{1pt}{1pt}{1}{}3\above{1pt}{1pt}{1}{}9\above{1pt}{1pt}{1}{}{\cdot}1\above{1pt}{1pt}{-2}{}5\above{1pt}{1pt}{-}{}{\cdot}1\above{1pt}{1pt}{-2}{}13\above{1pt}{1pt}{1}{}$}\spacer%
\instructions{2}%
\EasyButWeakLineBreak%
{${156}\above{1pt}{1pt}{r}{2}{10}\above{1pt}{1pt}{b}{2}{468}\above{1pt}{1pt}{*}{2}{4}\above{1pt}{1pt}{b}{2}{390}\above{1pt}{1pt}{b}{2}{36}\above{1pt}{1pt}{*}{2}{52}\above{1pt}{1pt}{b}{2}{90}\above{1pt}{1pt}{l}{2}$}%
\nopagebreak\par%
\nopagebreak\par\leavevmode%
{$\left[\!\llap{\phantom{%
\begingroup \smaller\smaller\smaller\begin{tabular}{@{}c@{}}%
0\\0\\0
\end{tabular}\endgroup%
}}\right.$}%
\begingroup \smaller\smaller\smaller\begin{tabular}{@{}c@{}}%
73848060\\-24605100\\-65520
\end{tabular}\endgroup%
\kern3pt%
\begingroup \smaller\smaller\smaller\begin{tabular}{@{}c@{}}%
-24605100\\8198058\\21831
\end{tabular}\endgroup%
\kern3pt%
\begingroup \smaller\smaller\smaller\begin{tabular}{@{}c@{}}%
-65520\\21831\\58
\end{tabular}\endgroup%
{$\left.\llap{\phantom{%
\begingroup \smaller\smaller\smaller\begin{tabular}{@{}c@{}}%
0\\0\\0
\end{tabular}\endgroup%
}}\!\right]$}%
\EasyButWeakLineBreak%
{$\left[\!\llap{\phantom{%
\begingroup \smaller\smaller\smaller\begin{tabular}{@{}c@{}}%
0\\0\\0
\end{tabular}\endgroup%
}}\right.$}%
\begingroup \smaller\smaller\smaller\begin{tabular}{@{}c@{}}%
861\\2548\\13572
\end{tabular}\endgroup%
\HardButStrongLineBreak\kern3pt%
\begingroup \smaller\smaller\smaller\begin{tabular}{@{}c@{}}%
49\\145\\775
\end{tabular}\endgroup%
\HardButStrongLineBreak\kern3pt%
\begingroup \smaller\smaller\smaller\begin{tabular}{@{}c@{}}%
-79\\-234\\-1170
\end{tabular}\endgroup%
\HardButStrongLineBreak\kern3pt%
\begingroup \smaller\smaller\smaller\begin{tabular}{@{}c@{}}%
-25\\-74\\-388
\end{tabular}\endgroup%
\HardButStrongLineBreak\kern3pt%
\begingroup \smaller\smaller\smaller\begin{tabular}{@{}c@{}}%
22\\65\\390
\end{tabular}\endgroup%
\HardButStrongLineBreak\kern3pt%
\begingroup \smaller\smaller\smaller\begin{tabular}{@{}c@{}}%
221\\654\\3492
\end{tabular}\endgroup%
\HardButStrongLineBreak\kern3pt%
\begingroup \smaller\smaller\smaller\begin{tabular}{@{}c@{}}%
615\\1820\\9698
\end{tabular}\endgroup%
\HardButStrongLineBreak\kern3pt%
\begingroup \smaller\smaller\smaller\begin{tabular}{@{}c@{}}%
517\\1530\\8145
\end{tabular}\endgroup%
{$\left.\llap{\phantom{%
\begingroup \smaller\smaller\smaller\begin{tabular}{@{}c@{}}%
0\\0\\0
\end{tabular}\endgroup%
}}\!\right]$}%

\medskip%
%
\leavevmode\llap{}%
$W_{339}$%
\qquad\llap{32} lattices, $\chi=36$%
\hfill%
$2222|22222|2\rtimes D_{2}$%
\nopagebreak\smallskip\hrule\nopagebreak\medskip%
%
%
\leavevmode%
${L_{339.1}}$%
{} : {$1\above{1pt}{1pt}{2}{6}8\above{1pt}{1pt}{1}{7}{\cdot}1\above{1pt}{1pt}{1}{}3\above{1pt}{1pt}{-}{}9\above{1pt}{1pt}{1}{}{\cdot}1\above{1pt}{1pt}{-2}{}5\above{1pt}{1pt}{1}{}{\cdot}1\above{1pt}{1pt}{2}{}7\above{1pt}{1pt}{-}{}$}\spacer%
\instructions{2}%
\EasyButWeakLineBreak%
{${70}\above{1pt}{1pt}{s}{2}{6}\above{1pt}{1pt}{b}{2}{2520}\above{1pt}{1pt}{*}{2}{4}\above{1pt}{1pt}{s}{2}{168}\above{1pt}{1pt}{s}{2}{36}\above{1pt}{1pt}{*}{2}{280}\above{1pt}{1pt}{b}{2}{6}\above{1pt}{1pt}{s}{2}{630}\above{1pt}{1pt}{b}{2}{24}\above{1pt}{1pt}{b}{2}$}%
\nopagebreak\par%
\nopagebreak\par\leavevmode%
{$\left[\!\llap{\phantom{%
\begingroup \smaller\smaller\smaller
\endgroup%
}}\!\right]$}%
%
%
\hbox{}\par\smallskip%
%
%
\leavevmode%
${L_{339.2}}$%
{} : {$1\above{1pt}{1pt}{-2}{6}8\above{1pt}{1pt}{-}{3}{\cdot}1\above{1pt}{1pt}{1}{}3\above{1pt}{1pt}{-}{}9\above{1pt}{1pt}{1}{}{\cdot}1\above{1pt}{1pt}{-2}{}5\above{1pt}{1pt}{1}{}{\cdot}1\above{1pt}{1pt}{2}{}7\above{1pt}{1pt}{-}{}$}\spacer%
\instructions{m}%
\EasyButWeakLineBreak%
{${70}\above{1pt}{1pt}{b}{2}{6}\above{1pt}{1pt}{l}{2}{2520}\above{1pt}{1pt}{}{2}{1}\above{1pt}{1pt}{r}{2}{168}\above{1pt}{1pt}{l}{2}{9}\above{1pt}{1pt}{}{2}{280}\above{1pt}{1pt}{r}{2}{6}\above{1pt}{1pt}{b}{2}{630}\above{1pt}{1pt}{l}{2}{24}\above{1pt}{1pt}{r}{2}$}%
\nopagebreak\par%
\nopagebreak\par\leavevmode%
{$\left[\!\llap{\phantom{%
\begingroup \smaller\smaller\smaller
\endgroup%
}}\!\right]$}%

\medskip%
%
\leavevmode\llap{}%
$W_{340}$%
\qquad\llap{16} lattices, $\chi=72$%
\hfill%
$2|2222|2222|2222|222\rtimes D_{4}$%
\nopagebreak\smallskip\hrule\nopagebreak\medskip%
%
%
\leavevmode%
${L_{340.1}}$%
{} : {$1\above{1pt}{1pt}{-2}{{\rm II}}8\above{1pt}{1pt}{-}{5}{\cdot}1\above{1pt}{1pt}{1}{}3\above{1pt}{1pt}{-}{}9\above{1pt}{1pt}{1}{}{\cdot}1\above{1pt}{1pt}{2}{}5\above{1pt}{1pt}{-}{}{\cdot}1\above{1pt}{1pt}{2}{}7\above{1pt}{1pt}{-}{}$}\spacer%
\instructions{2}%
\EasyButWeakLineBreak%
{${360}\above{1pt}{1pt}{r}{2}{42}\above{1pt}{1pt}{l}{2}{40}\above{1pt}{1pt}{r}{2}{6}\above{1pt}{1pt}{b}{2}{90}\above{1pt}{1pt}{l}{2}{168}\above{1pt}{1pt}{r}{2}{10}\above{1pt}{1pt}{b}{2}{6}\above{1pt}{1pt}{l}{2}$}\relax$\,(\times2)$%
\nopagebreak\par%
\nopagebreak\par\leavevmode%
{$\left[\!\llap{\phantom{%
\begingroup \smaller\smaller\smaller
\endgroup%
}}\!\right]$}%

\medskip%
%
\leavevmode\llap{}%
$W_{341}$%
\qquad\llap{4} lattices, $\chi=48$%
\hfill%
$2\infty|\infty2|2\infty|\infty2|\rtimes D_{4}$%
\nopagebreak\smallskip\hrule\nopagebreak\medskip%
%
%
\leavevmode%
${L_{341.1}}$%
{} : {$1\above{1pt}{1pt}{1}{7}8\above{1pt}{1pt}{1}{1}256\above{1pt}{1pt}{1}{1}$}\EasyButWeakLineBreak%
{${256}\above{1pt}{1pt}{*}{2}{32}\above{1pt}{1pt}{32,17}{\infty z}{8}\above{1pt}{1pt}{32,9}{\infty}{32}\above{1pt}{1pt}{s}{2}$}\relax$\,(\times2)$%
\nopagebreak\par%
shares genus with {$ {L_{342.1}}$}%
\nopagebreak\par%
\nopagebreak\par\leavevmode%
{$\left[\!\llap{\phantom{%
\begingroup \smaller\smaller\smaller\begin{tabular}{@{}c@{}}%
0\\0\\0
\end{tabular}\endgroup%
}}\right.$}%
\begingroup \smaller\smaller\smaller\begin{tabular}{@{}c@{}}%
-239360\\-11008\\-11520
\end{tabular}\endgroup%
\kern3pt%
\begingroup \smaller\smaller\smaller\begin{tabular}{@{}c@{}}%
-11008\\-504\\-528
\end{tabular}\endgroup%
\kern3pt%
\begingroup \smaller\smaller\smaller\begin{tabular}{@{}c@{}}%
-11520\\-528\\-553
\end{tabular}\endgroup%
{$\left.\llap{\phantom{%
\begingroup \smaller\smaller\smaller\begin{tabular}{@{}c@{}}%
0\\0\\0
\end{tabular}\endgroup%
}}\!\right]$}%
\hfil\penalty500%
{$\left[\!\llap{\phantom{%
\begingroup \smaller\smaller\smaller\begin{tabular}{@{}c@{}}%
0\\0\\0
\end{tabular}\endgroup%
}}\right.$}%
\begingroup \smaller\smaller\smaller\begin{tabular}{@{}c@{}}%
-1153\\-56448\\78336
\end{tabular}\endgroup%
\kern3pt%
\begingroup \smaller\smaller\smaller\begin{tabular}{@{}c@{}}%
-50\\-2451\\3400
\end{tabular}\endgroup%
\kern3pt%
\begingroup \smaller\smaller\smaller\begin{tabular}{@{}c@{}}%
-53\\-2597\\3603
\end{tabular}\endgroup%
{$\left.\llap{\phantom{%
\begingroup \smaller\smaller\smaller\begin{tabular}{@{}c@{}}%
0\\0\\0
\end{tabular}\endgroup%
}}\!\right]$}%
\EasyButWeakLineBreak%
{$\left[\!\llap{\phantom{%
\begingroup \smaller\smaller\smaller\begin{tabular}{@{}c@{}}%
0\\0\\0
\end{tabular}\endgroup%
}}\right.$}%
\begingroup \smaller\smaller\smaller\begin{tabular}{@{}c@{}}%
1\\112\\-128
\end{tabular}\endgroup%
\HardButStrongLineBreak\kern3pt%
\begingroup \smaller\smaller\smaller\begin{tabular}{@{}c@{}}%
3\\202\\-256
\end{tabular}\endgroup%
\HardButStrongLineBreak\kern3pt%
\begingroup \smaller\smaller\smaller\begin{tabular}{@{}c@{}}%
2\\115\\-152
\end{tabular}\endgroup%
\HardButStrongLineBreak\kern3pt%
\begingroup \smaller\smaller\smaller\begin{tabular}{@{}c@{}}%
13\\650\\-896
\end{tabular}\endgroup%
{$\left.\llap{\phantom{%
\begingroup \smaller\smaller\smaller\begin{tabular}{@{}c@{}}%
0\\0\\0
\end{tabular}\endgroup%
}}\!\right]$}%
%
%
\hbox{}\par\smallskip%
%
%
\leavevmode%
${L_{341.2}}$%
{} : {$1\above{1pt}{1pt}{1}{1}8\above{1pt}{1pt}{1}{1}256\above{1pt}{1pt}{1}{7}$}\EasyButWeakLineBreak%
{${1}\above{1pt}{1pt}{}{2}{8}\above{1pt}{1pt}{32,31}{\infty}{32}\above{1pt}{1pt}{32,7}{\infty z}{8}\above{1pt}{1pt}{r}{2}{4}\above{1pt}{1pt}{*}{2}{32}\above{1pt}{1pt}{32,15}{\infty z}{8}\above{1pt}{1pt}{32,7}{\infty}{32}\above{1pt}{1pt}{l}{2}$}%
\nopagebreak\par%
\nopagebreak\par\leavevmode%
{$\left[\!\llap{\phantom{%
\begingroup \smaller\smaller\smaller\begin{tabular}{@{}c@{}}%
0\\0\\0
\end{tabular}\endgroup%
}}\right.$}%
\begingroup \smaller\smaller\smaller\begin{tabular}{@{}c@{}}%
-1272064\\8960\\8960
\end{tabular}\endgroup%
\kern3pt%
\begingroup \smaller\smaller\smaller\begin{tabular}{@{}c@{}}%
8960\\-56\\-64
\end{tabular}\endgroup%
\kern3pt%
\begingroup \smaller\smaller\smaller\begin{tabular}{@{}c@{}}%
8960\\-64\\-63
\end{tabular}\endgroup%
{$\left.\llap{\phantom{%
\begingroup \smaller\smaller\smaller\begin{tabular}{@{}c@{}}%
0\\0\\0
\end{tabular}\endgroup%
}}\!\right]$}%
\EasyButWeakLineBreak%
{$\left[\!\llap{\phantom{%
\begingroup \smaller\smaller\smaller\begin{tabular}{@{}c@{}}%
0\\0\\0
\end{tabular}\endgroup%
}}\right.$}%
\begingroup \smaller\smaller\smaller\begin{tabular}{@{}c@{}}%
10\\154\\1263
\end{tabular}\endgroup%
\HardButStrongLineBreak\kern3pt%
\begingroup \smaller\smaller\smaller\begin{tabular}{@{}c@{}}%
29\\445\\3664
\end{tabular}\endgroup%
\HardButStrongLineBreak\kern3pt%
\begingroup \smaller\smaller\smaller\begin{tabular}{@{}c@{}}%
11\\166\\1392
\end{tabular}\endgroup%
\HardButStrongLineBreak\kern3pt%
\begingroup \smaller\smaller\smaller\begin{tabular}{@{}c@{}}%
1\\13\\128
\end{tabular}\endgroup%
\HardButStrongLineBreak\kern3pt%
\begingroup \smaller\smaller\smaller\begin{tabular}{@{}c@{}}%
-1\\-16\\-126
\end{tabular}\endgroup%
\HardButStrongLineBreak\kern3pt%
\begingroup \smaller\smaller\smaller\begin{tabular}{@{}c@{}}%
-1\\-14\\-128
\end{tabular}\endgroup%
\HardButStrongLineBreak\kern3pt%
\begingroup \smaller\smaller\smaller\begin{tabular}{@{}c@{}}%
4\\63\\504
\end{tabular}\endgroup%
\HardButStrongLineBreak\kern3pt%
\begingroup \smaller\smaller\smaller\begin{tabular}{@{}c@{}}%
55\\850\\6944
\end{tabular}\endgroup%
{$\left.\llap{\phantom{%
\begingroup \smaller\smaller\smaller\begin{tabular}{@{}c@{}}%
0\\0\\0
\end{tabular}\endgroup%
}}\!\right]$}%
%
%
%
%
%
%
%
%
%
%

\medskip%
%
\leavevmode\llap{}%
$W_{342}$%
\qquad\llap{2} lattices, $\chi=24$%
\hfill%
$2|22\infty|\infty2\rtimes D_{2}$%
\nopagebreak\smallskip\hrule\nopagebreak\medskip%
%
%
\leavevmode%
${L_{342.1}}$%
{} : {$1\above{1pt}{1pt}{1}{7}8\above{1pt}{1pt}{1}{1}256\above{1pt}{1pt}{1}{1}$}\EasyButWeakLineBreak%
{${256}\above{1pt}{1pt}{r}{2}{4}\above{1pt}{1pt}{b}{2}{256}\above{1pt}{1pt}{l}{2}{8}\above{1pt}{1pt}{32,17}{\infty}{32}\above{1pt}{1pt}{32,25}{\infty z}{8}\above{1pt}{1pt}{}{2}$}%
\nopagebreak\par%
shares genus with {$ {L_{341.1}}$}%
\nopagebreak\par%
\nopagebreak\par\leavevmode%
{$\left[\!\llap{\phantom{%
\begingroup \smaller\smaller\smaller\begin{tabular}{@{}c@{}}%
0\\0\\0
\end{tabular}\endgroup%
}}\right.$}%
\begingroup \smaller\smaller\smaller\begin{tabular}{@{}c@{}}%
13951232\\6751744\\425216
\end{tabular}\endgroup%
\kern3pt%
\begingroup \smaller\smaller\smaller\begin{tabular}{@{}c@{}}%
6751744\\3267528\\205784
\end{tabular}\endgroup%
\kern3pt%
\begingroup \smaller\smaller\smaller\begin{tabular}{@{}c@{}}%
425216\\205784\\12959
\end{tabular}\endgroup%
{$\left.\llap{\phantom{%
\begingroup \smaller\smaller\smaller\begin{tabular}{@{}c@{}}%
0\\0\\0
\end{tabular}\endgroup%
}}\!\right]$}%
\EasyButWeakLineBreak%
{$\left[\!\llap{\phantom{%
\begingroup \smaller\smaller\smaller\begin{tabular}{@{}c@{}}%
0\\0\\0
\end{tabular}\endgroup%
}}\right.$}%
\begingroup \smaller\smaller\smaller\begin{tabular}{@{}c@{}}%
-209\\448\\-256
\end{tabular}\endgroup%
\HardButStrongLineBreak\kern3pt%
\begingroup \smaller\smaller\smaller\begin{tabular}{@{}c@{}}%
84\\-181\\118
\end{tabular}\endgroup%
\HardButStrongLineBreak\kern3pt%
\begingroup \smaller\smaller\smaller\begin{tabular}{@{}c@{}}%
8929\\-19232\\12416
\end{tabular}\endgroup%
\HardButStrongLineBreak\kern3pt%
\begingroup \smaller\smaller\smaller\begin{tabular}{@{}c@{}}%
1411\\-3039\\1960
\end{tabular}\endgroup%
\HardButStrongLineBreak\kern3pt%
\begingroup \smaller\smaller\smaller\begin{tabular}{@{}c@{}}%
209\\-450\\288
\end{tabular}\endgroup%
\HardButStrongLineBreak\kern3pt%
\begingroup \smaller\smaller\smaller\begin{tabular}{@{}c@{}}%
-112\\241\\-152
\end{tabular}\endgroup%
{$\left.\llap{\phantom{%
\begingroup \smaller\smaller\smaller\begin{tabular}{@{}c@{}}%
0\\0\\0
\end{tabular}\endgroup%
}}\!\right]$}%
%
%
%
%
%
%

\medskip%
%
\leavevmode\llap{}%
$W_{343}$%
\qquad\llap{2} lattices, $\chi=24$%
\hfill%
$222\slashinfty222|\rtimes D_{2}$%
\nopagebreak\smallskip\hrule\nopagebreak\medskip%
%
%
\leavevmode%
${L_{343.1}}$%
{} : {$1\above{1pt}{1pt}{1}{1}8\above{1pt}{1pt}{1}{7}256\above{1pt}{1pt}{1}{1}$}\EasyButWeakLineBreak%
{${256}\above{1pt}{1pt}{l}{2}{1}\above{1pt}{1pt}{}{2}{256}\above{1pt}{1pt}{r}{2}{8}\above{1pt}{1pt}{16,1}{\infty a}{8}\above{1pt}{1pt}{b}{2}{256}\above{1pt}{1pt}{s}{2}{4}\above{1pt}{1pt}{*}{2}$}%
\nopagebreak\par%
\nopagebreak\par\leavevmode%
{$\left[\!\llap{\phantom{%
\begingroup \smaller\smaller\smaller\begin{tabular}{@{}c@{}}%
0\\0\\0
\end{tabular}\endgroup%
}}\right.$}%
\begingroup \smaller\smaller\smaller\begin{tabular}{@{}c@{}}%
385280\\8448\\-2560
\end{tabular}\endgroup%
\kern3pt%
\begingroup \smaller\smaller\smaller\begin{tabular}{@{}c@{}}%
8448\\184\\-56
\end{tabular}\endgroup%
\kern3pt%
\begingroup \smaller\smaller\smaller\begin{tabular}{@{}c@{}}%
-2560\\-56\\17
\end{tabular}\endgroup%
{$\left.\llap{\phantom{%
\begingroup \smaller\smaller\smaller\begin{tabular}{@{}c@{}}%
0\\0\\0
\end{tabular}\endgroup%
}}\!\right]$}%
\EasyButWeakLineBreak%
{$\left[\!\llap{\phantom{%
\begingroup \smaller\smaller\smaller\begin{tabular}{@{}c@{}}%
0\\0\\0
\end{tabular}\endgroup%
}}\right.$}%
\begingroup \smaller\smaller\smaller\begin{tabular}{@{}c@{}}%
-3\\-16\\-512
\end{tabular}\endgroup%
\HardButStrongLineBreak\kern3pt%
\begingroup \smaller\smaller\smaller\begin{tabular}{@{}c@{}}%
0\\1\\3
\end{tabular}\endgroup%
\HardButStrongLineBreak\kern3pt%
\begingroup \smaller\smaller\smaller\begin{tabular}{@{}c@{}}%
1\\32\\256
\end{tabular}\endgroup%
\HardButStrongLineBreak\kern3pt%
\begingroup \smaller\smaller\smaller\begin{tabular}{@{}c@{}}%
0\\1\\4
\end{tabular}\endgroup%
\HardButStrongLineBreak\kern3pt%
\begingroup \smaller\smaller\smaller\begin{tabular}{@{}c@{}}%
-1\\-9\\-180
\end{tabular}\endgroup%
\HardButStrongLineBreak\kern3pt%
\begingroup \smaller\smaller\smaller\begin{tabular}{@{}c@{}}%
-15\\-128\\-2688
\end{tabular}\endgroup%
\HardButStrongLineBreak\kern3pt%
\begingroup \smaller\smaller\smaller\begin{tabular}{@{}c@{}}%
-1\\-8\\-178
\end{tabular}\endgroup%
{$\left.\llap{\phantom{%
\begingroup \smaller\smaller\smaller\begin{tabular}{@{}c@{}}%
0\\0\\0
\end{tabular}\endgroup%
}}\!\right]$}%
%
%
%
%
%
%

\medskip%
%
\leavevmode\llap{}%
$W_{344}$%
\qquad\llap{32} lattices, $\chi=96$%
\hfill%
$2222\infty\infty\infty2222\infty\infty\infty\rtimes C_{2}$%
\nopagebreak\smallskip\hrule\nopagebreak\medskip%
%
%
\leavevmode%
${L_{344.1}}$%
{} : {$1\above{1pt}{1pt}{-2}{{\rm II}}8\above{1pt}{1pt}{1}{1}{\cdot}1\above{1pt}{1pt}{1}{}3\above{1pt}{1pt}{-}{}9\above{1pt}{1pt}{-}{}{\cdot}1\above{1pt}{1pt}{-}{}7\above{1pt}{1pt}{-}{}49\above{1pt}{1pt}{1}{}$}\spacer%
\instructions{23,3,2}%
\EasyButWeakLineBreak%
{${42}\above{1pt}{1pt}{l}{2}{3528}\above{1pt}{1pt}{r}{2}{6}\above{1pt}{1pt}{s}{2}{882}\above{1pt}{1pt}{b}{2}{168}\above{1pt}{1pt}{42,1}{\infty z}{42}\above{1pt}{1pt}{84,37}{\infty b}{168}\above{1pt}{1pt}{42,37}{\infty z}$}\relax$\,(\times2)$%
\nopagebreak\par%
\nopagebreak\par\leavevmode%
{$\left[\!\llap{\phantom{%
\begingroup \smaller\smaller\smaller
\endgroup%
}}\!\right]$}%

\medskip%
%
\leavevmode\llap{}%
$W_{345}$%
\qquad\llap{16} lattices, $\chi=48$%
\hfill%
$222|222|222|222|\rtimes D_{4}$%
\nopagebreak\smallskip\hrule\nopagebreak\medskip%
%
%
\leavevmode%
${L_{345.1}}$%
{} : {$1\above{1pt}{1pt}{-2}{{\rm II}}8\above{1pt}{1pt}{1}{1}{\cdot}1\above{1pt}{1pt}{-}{}3\above{1pt}{1pt}{1}{}9\above{1pt}{1pt}{-}{}{\cdot}1\above{1pt}{1pt}{1}{}7\above{1pt}{1pt}{1}{}49\above{1pt}{1pt}{1}{}$}\spacer%
\instructions{23,3,2*}%
\EasyButWeakLineBreak%
{${14}\above{1pt}{1pt}{l}{2}{72}\above{1pt}{1pt}{r}{2}{98}\above{1pt}{1pt}{b}{2}{126}\above{1pt}{1pt}{b}{2}{2}\above{1pt}{1pt}{l}{2}{3528}\above{1pt}{1pt}{r}{2}$}\relax$\,(\times2)$%
\nopagebreak\par%
shares genus with 3-dual; isometric to own %
7-dual\nopagebreak\par%
\nopagebreak\par\leavevmode%
{$\left[\!\llap{\phantom{%
\begingroup \smaller\smaller\smaller
\endgroup%
}}\!\right]$}%

\medskip%
%
\leavevmode\llap{}%
$W_{346}$%
\qquad\llap{12} lattices, $\chi=60$%
\hfill%
$22|222\slashtwo222|222\slashtwo2\rtimes D_{4}$%
\nopagebreak\smallskip\hrule\nopagebreak\medskip%
%
%
\leavevmode%
${L_{346.1}}$%
{} : {$1\above{1pt}{1pt}{-2}{{\rm II}}4\above{1pt}{1pt}{-}{3}{\cdot}1\above{1pt}{1pt}{-}{}3\above{1pt}{1pt}{1}{}9\above{1pt}{1pt}{-}{}{\cdot}1\above{1pt}{1pt}{2}{}7\above{1pt}{1pt}{1}{}{\cdot}1\above{1pt}{1pt}{2}{}11\above{1pt}{1pt}{1}{}$}\spacer%
\instructions{2}%
\EasyButWeakLineBreak%
{${396}\above{1pt}{1pt}{r}{2}{14}\above{1pt}{1pt}{l}{2}{12}\above{1pt}{1pt}{r}{2}{126}\above{1pt}{1pt}{l}{2}{44}\above{1pt}{1pt}{r}{2}{18}\above{1pt}{1pt}{b}{2}{2}\above{1pt}{1pt}{l}{2}$}\relax$\,(\times2)$%
\nopagebreak\par%
\nopagebreak\par\leavevmode%
{$\left[\!\llap{\phantom{%
\begingroup \smaller\smaller\smaller
\endgroup%
}}\!\right]$}%

\medskip%
%
\leavevmode\llap{}%
$W_{347}$%
\qquad\llap{4} lattices, $\chi=24$%
\hfill%
$22|22|22|22|\rtimes D_{4}$%
\nopagebreak\smallskip\hrule\nopagebreak\medskip%
%
%
\leavevmode%
${L_{347.1}}$%
{} : {$1\above{1pt}{1pt}{-2}{2}16\above{1pt}{1pt}{-}{5}{\cdot}1\above{1pt}{1pt}{-}{}3\above{1pt}{1pt}{1}{}9\above{1pt}{1pt}{-}{}{\cdot}1\above{1pt}{1pt}{-2}{}5\above{1pt}{1pt}{1}{}$}\EasyButWeakLineBreak%
{${18}\above{1pt}{1pt}{l}{2}{80}\above{1pt}{1pt}{}{2}{3}\above{1pt}{1pt}{}{2}{720}\above{1pt}{1pt}{r}{2}{2}\above{1pt}{1pt}{b}{2}{720}\above{1pt}{1pt}{*}{2}{12}\above{1pt}{1pt}{*}{2}{80}\above{1pt}{1pt}{b}{2}$}%
\nopagebreak\par%
\nopagebreak\par\leavevmode%
{$\left[\!\llap{\phantom{%
\begingroup \smaller\smaller\smaller\begin{tabular}{@{}c@{}}%
0\\0\\0
\end{tabular}\endgroup%
}}\right.$}%
\begingroup \smaller\smaller\smaller\begin{tabular}{@{}c@{}}%
-137520\\2160\\-720
\end{tabular}\endgroup%
\kern3pt%
\begingroup \smaller\smaller\smaller\begin{tabular}{@{}c@{}}%
2160\\-6\\-9
\end{tabular}\endgroup%
\kern3pt%
\begingroup \smaller\smaller\smaller\begin{tabular}{@{}c@{}}%
-720\\-9\\11
\end{tabular}\endgroup%
{$\left.\llap{\phantom{%
\begingroup \smaller\smaller\smaller\begin{tabular}{@{}c@{}}%
0\\0\\0
\end{tabular}\endgroup%
}}\!\right]$}%
\EasyButWeakLineBreak%
{$\left[\!\llap{\phantom{%
\begingroup \smaller\smaller\smaller\begin{tabular}{@{}c@{}}%
0\\0\\0
\end{tabular}\endgroup%
}}\right.$}%
\begingroup \smaller\smaller\smaller\begin{tabular}{@{}c@{}}%
1\\117\\162
\end{tabular}\endgroup%
\HardButStrongLineBreak\kern3pt%
\begingroup \smaller\smaller\smaller\begin{tabular}{@{}c@{}}%
11\\1280\\1760
\end{tabular}\endgroup%
\HardButStrongLineBreak\kern3pt%
\begingroup \smaller\smaller\smaller\begin{tabular}{@{}c@{}}%
1\\116\\159
\end{tabular}\endgroup%
\HardButStrongLineBreak\kern3pt%
\begingroup \smaller\smaller\smaller\begin{tabular}{@{}c@{}}%
23\\2640\\3600
\end{tabular}\endgroup%
\HardButStrongLineBreak\kern3pt%
\begingroup \smaller\smaller\smaller\begin{tabular}{@{}c@{}}%
0\\-1\\-2
\end{tabular}\endgroup%
\HardButStrongLineBreak\kern3pt%
\begingroup \smaller\smaller\smaller\begin{tabular}{@{}c@{}}%
-13\\-1560\\-2160
\end{tabular}\endgroup%
\HardButStrongLineBreak\kern3pt%
\begingroup \smaller\smaller\smaller\begin{tabular}{@{}c@{}}%
-1\\-118\\-162
\end{tabular}\endgroup%
\HardButStrongLineBreak\kern3pt%
\begingroup \smaller\smaller\smaller\begin{tabular}{@{}c@{}}%
-1\\-120\\-160
\end{tabular}\endgroup%
{$\left.\llap{\phantom{%
\begingroup \smaller\smaller\smaller\begin{tabular}{@{}c@{}}%
0\\0\\0
\end{tabular}\endgroup%
}}\!\right]$}%
%
%
%
%
%
%
%
%
%
%
%
%
%
%

\medskip%
%
\leavevmode\llap{}%
$W_{348}$%
\qquad\llap{12} lattices, $\chi=32$%
\hfill%
$222|2222\slashthree2\rtimes D_{2}$%
\nopagebreak\smallskip\hrule\nopagebreak\medskip%
%
%
\leavevmode%
${L_{348.1}}$%
{} : {$1\above{1pt}{1pt}{-2}{{\rm II}}4\above{1pt}{1pt}{-}{3}{\cdot}1\above{1pt}{1pt}{1}{}3\above{1pt}{1pt}{-}{}9\above{1pt}{1pt}{1}{}{\cdot}1\above{1pt}{1pt}{-2}{}5\above{1pt}{1pt}{-}{}{\cdot}1\above{1pt}{1pt}{-2}{}17\above{1pt}{1pt}{-}{}$}\spacer%
\instructions{2}%
\EasyButWeakLineBreak%
{${340}\above{1pt}{1pt}{*}{2}{36}\above{1pt}{1pt}{b}{2}{10}\above{1pt}{1pt}{l}{2}{204}\above{1pt}{1pt}{r}{2}{90}\above{1pt}{1pt}{b}{2}{4}\above{1pt}{1pt}{*}{2}{3060}\above{1pt}{1pt}{b}{2}{6}\above{1pt}{1pt}{+}{3}{6}\above{1pt}{1pt}{b}{2}$}%
\nopagebreak\par%
\nopagebreak\par\leavevmode%
{$\left[\!\llap{\phantom{%
\begingroup \smaller\smaller\smaller\begin{tabular}{@{}c@{}}%
0\\0\\0
\end{tabular}\endgroup%
}}\right.$}%
\begingroup \smaller\smaller\smaller\begin{tabular}{@{}c@{}}%
-23993460\\-1407600\\18360
\end{tabular}\endgroup%
\kern3pt%
\begingroup \smaller\smaller\smaller\begin{tabular}{@{}c@{}}%
-1407600\\-82578\\1077
\end{tabular}\endgroup%
\kern3pt%
\begingroup \smaller\smaller\smaller\begin{tabular}{@{}c@{}}%
18360\\1077\\-14
\end{tabular}\endgroup%
{$\left.\llap{\phantom{%
\begingroup \smaller\smaller\smaller\begin{tabular}{@{}c@{}}%
0\\0\\0
\end{tabular}\endgroup%
}}\!\right]$}%
\EasyButWeakLineBreak%
{$\left[\!\llap{\phantom{%
\begingroup \smaller\smaller\smaller\begin{tabular}{@{}c@{}}%
0\\0\\0
\end{tabular}\endgroup%
}}\right.$}%
\begingroup \smaller\smaller\smaller\begin{tabular}{@{}c@{}}%
-29\\510\\1190
\end{tabular}\endgroup%
\HardButStrongLineBreak\kern3pt%
\begingroup \smaller\smaller\smaller\begin{tabular}{@{}c@{}}%
-1\\18\\72
\end{tabular}\endgroup%
\HardButStrongLineBreak\kern3pt%
\begingroup \smaller\smaller\smaller\begin{tabular}{@{}c@{}}%
2\\-35\\-70
\end{tabular}\endgroup%
\HardButStrongLineBreak\kern3pt%
\begingroup \smaller\smaller\smaller\begin{tabular}{@{}c@{}}%
27\\-476\\-1224
\end{tabular}\endgroup%
\HardButStrongLineBreak\kern3pt%
\begingroup \smaller\smaller\smaller\begin{tabular}{@{}c@{}}%
11\\-195\\-585
\end{tabular}\endgroup%
\HardButStrongLineBreak\kern3pt%
\begingroup \smaller\smaller\smaller\begin{tabular}{@{}c@{}}%
3\\-54\\-226
\end{tabular}\endgroup%
\HardButStrongLineBreak\kern3pt%
\begingroup \smaller\smaller\smaller\begin{tabular}{@{}c@{}}%
83\\-1530\\-9180
\end{tabular}\endgroup%
\HardButStrongLineBreak\kern3pt%
\begingroup \smaller\smaller\smaller\begin{tabular}{@{}c@{}}%
-1\\17\\-6
\end{tabular}\endgroup%
\HardButStrongLineBreak\kern3pt%
\begingroup \smaller\smaller\smaller\begin{tabular}{@{}c@{}}%
-2\\35\\69
\end{tabular}\endgroup%
{$\left.\llap{\phantom{%
\begingroup \smaller\smaller\smaller\begin{tabular}{@{}c@{}}%
0\\0\\0
\end{tabular}\endgroup%
}}\!\right]$}%

\medskip%
%
\leavevmode\llap{}%
$W_{349}$%
\qquad\llap{30} lattices, $\chi=120$%
\hfill%
$22|222|222|222|222|222|222|222|2\rtimes D_{8}$%
\nopagebreak\smallskip\hrule\nopagebreak\medskip%
%
%
\leavevmode%
${L_{349.1}}$%
{} : {$1\above{1pt}{1pt}{2}{0}8\above{1pt}{1pt}{1}{1}{\cdot}1\above{1pt}{1pt}{-}{}3\above{1pt}{1pt}{-}{}9\above{1pt}{1pt}{-}{}{\cdot}1\above{1pt}{1pt}{2}{}11\above{1pt}{1pt}{1}{}$}\EasyButWeakLineBreak%
{${99}\above{1pt}{1pt}{}{2}{8}\above{1pt}{1pt}{}{2}{33}\above{1pt}{1pt}{}{2}{72}\above{1pt}{1pt}{}{2}{11}\above{1pt}{1pt}{r}{2}{24}\above{1pt}{1pt}{s}{2}{44}\above{1pt}{1pt}{*}{2}{72}\above{1pt}{1pt}{*}{2}{132}\above{1pt}{1pt}{*}{2}{8}\above{1pt}{1pt}{*}{2}{396}\above{1pt}{1pt}{s}{2}{24}\above{1pt}{1pt}{l}{2}$}\relax$\,(\times2)$%
\nopagebreak\par%
\nopagebreak\par\leavevmode%
{$\left[\!\llap{\phantom{%
\begingroup \smaller\smaller\smaller
\endgroup%
}}\!\right]$}%
%
%
\hbox{}\par\smallskip%
%
%
\leavevmode%
${L_{349.2}}$%
{} : {$[1\above{1pt}{1pt}{1}{}2\above{1pt}{1pt}{1}{}]\above{1pt}{1pt}{}{2}16\above{1pt}{1pt}{1}{7}{\cdot}1\above{1pt}{1pt}{-}{}3\above{1pt}{1pt}{-}{}9\above{1pt}{1pt}{-}{}{\cdot}1\above{1pt}{1pt}{2}{}11\above{1pt}{1pt}{1}{}$}\spacer%
\instructions{2}%
\EasyButWeakLineBreak%
{${1584}\above{1pt}{1pt}{l}{2}{2}\above{1pt}{1pt}{}{2}{33}\above{1pt}{1pt}{}{2}{18}\above{1pt}{1pt}{r}{2}{176}\above{1pt}{1pt}{l}{2}{6}\above{1pt}{1pt}{}{2}{11}\above{1pt}{1pt}{r}{2}{72}\above{1pt}{1pt}{*}{2}{528}\above{1pt}{1pt}{*}{2}{8}\above{1pt}{1pt}{l}{2}{99}\above{1pt}{1pt}{}{2}{6}\above{1pt}{1pt}{r}{2}$}\relax$\,(\times2)$%
\nopagebreak\par%
\nopagebreak\par\leavevmode%
{$\left[\!\llap{\phantom{%
\begingroup \smaller\smaller\smaller
\endgroup%
}}\!\right]$}%
%
%
\hbox{}\par\smallskip%
%
%
\leavevmode%
${L_{349.3}}$%
{} : {$[1\above{1pt}{1pt}{-}{}2\above{1pt}{1pt}{1}{}]\above{1pt}{1pt}{}{6}16\above{1pt}{1pt}{-}{3}{\cdot}1\above{1pt}{1pt}{-}{}3\above{1pt}{1pt}{-}{}9\above{1pt}{1pt}{-}{}{\cdot}1\above{1pt}{1pt}{2}{}11\above{1pt}{1pt}{1}{}$}\spacer%
\instructions{m}%
\EasyButWeakLineBreak%
{${1584}\above{1pt}{1pt}{}{2}{2}\above{1pt}{1pt}{r}{2}{132}\above{1pt}{1pt}{l}{2}{18}\above{1pt}{1pt}{}{2}{176}\above{1pt}{1pt}{}{2}{6}\above{1pt}{1pt}{r}{2}{44}\above{1pt}{1pt}{*}{2}{72}\above{1pt}{1pt}{s}{2}{528}\above{1pt}{1pt}{s}{2}{8}\above{1pt}{1pt}{*}{2}{396}\above{1pt}{1pt}{l}{2}{6}\above{1pt}{1pt}{}{2}$}\relax$\,(\times2)$%
\nopagebreak\par%
\nopagebreak\par\leavevmode%
{$\left[\!\llap{\phantom{%
\begingroup \smaller\smaller\smaller
\endgroup%
}}\!\right]$}%
%
%
\hbox{}\par\smallskip%
%
%
\leavevmode%
${L_{349.4}}$%
{} : {$[1\above{1pt}{1pt}{-}{}2\above{1pt}{1pt}{1}{}]\above{1pt}{1pt}{}{4}16\above{1pt}{1pt}{-}{5}{\cdot}1\above{1pt}{1pt}{-}{}3\above{1pt}{1pt}{-}{}9\above{1pt}{1pt}{-}{}{\cdot}1\above{1pt}{1pt}{2}{}11\above{1pt}{1pt}{1}{}$}\spacer%
\instructions{m}%
\EasyButWeakLineBreak%
{${1584}\above{1pt}{1pt}{s}{2}{8}\above{1pt}{1pt}{*}{2}{132}\above{1pt}{1pt}{*}{2}{72}\above{1pt}{1pt}{s}{2}{176}\above{1pt}{1pt}{*}{2}{24}\above{1pt}{1pt}{l}{2}{11}\above{1pt}{1pt}{}{2}{18}\above{1pt}{1pt}{r}{2}{528}\above{1pt}{1pt}{l}{2}{2}\above{1pt}{1pt}{}{2}{99}\above{1pt}{1pt}{r}{2}{24}\above{1pt}{1pt}{*}{2}$}\relax$\,(\times2)$%
\nopagebreak\par%
\nopagebreak\par\leavevmode%
{$\left[\!\llap{\phantom{%
\begingroup \smaller\smaller\smaller
\endgroup%
}}\!\right]$}%
%
%
\hbox{}\par\smallskip%
%
%
\leavevmode%
${L_{349.5}}$%
{} : {$[1\above{1pt}{1pt}{1}{}2\above{1pt}{1pt}{1}{}]\above{1pt}{1pt}{}{0}16\above{1pt}{1pt}{1}{1}{\cdot}1\above{1pt}{1pt}{-}{}3\above{1pt}{1pt}{-}{}9\above{1pt}{1pt}{-}{}{\cdot}1\above{1pt}{1pt}{2}{}11\above{1pt}{1pt}{1}{}$}\EasyButWeakLineBreak%
{${1584}\above{1pt}{1pt}{*}{2}{8}\above{1pt}{1pt}{l}{2}{33}\above{1pt}{1pt}{r}{2}{72}\above{1pt}{1pt}{*}{2}{176}\above{1pt}{1pt}{s}{2}{24}\above{1pt}{1pt}{*}{2}{44}\above{1pt}{1pt}{l}{2}{18}\above{1pt}{1pt}{}{2}{528}\above{1pt}{1pt}{}{2}{2}\above{1pt}{1pt}{r}{2}{396}\above{1pt}{1pt}{*}{2}{24}\above{1pt}{1pt}{s}{2}$}\relax$\,(\times2)$%
\nopagebreak\par%
\nopagebreak\par\leavevmode%
{$\left[\!\llap{\phantom{%
\begingroup \smaller\smaller\smaller
\endgroup%
}}\!\right]$}%

\medskip%
%
\leavevmode\llap{}%
$W_{350}$%
\qquad\llap{12} lattices, $\chi=84$%
\hfill%
$2\slashtwo2|2\slashtwo2|2\slashtwo2|2\slashtwo2|2\slashtwo2|2\slashtwo2|\rtimes D_{12}$%
\nopagebreak\smallskip\hrule\nopagebreak\medskip%
%
%
\leavevmode%
${L_{350.1}}$%
{} : {$1\above{1pt}{1pt}{2}{2}16\above{1pt}{1pt}{1}{1}{\cdot}1\above{1pt}{1pt}{1}{}3\above{1pt}{1pt}{-}{}9\above{1pt}{1pt}{1}{}{\cdot}1\above{1pt}{1pt}{2}{}13\above{1pt}{1pt}{1}{}$}\EasyButWeakLineBreak%
{${1872}\above{1pt}{1pt}{l}{2}{1}\above{1pt}{1pt}{}{2}{144}\above{1pt}{1pt}{}{2}{13}\above{1pt}{1pt}{r}{2}{36}\above{1pt}{1pt}{*}{2}{4}\above{1pt}{1pt}{l}{2}{117}\above{1pt}{1pt}{}{2}{16}\above{1pt}{1pt}{}{2}{9}\above{1pt}{1pt}{r}{2}{208}\above{1pt}{1pt}{s}{2}{36}\above{1pt}{1pt}{*}{2}{16}\above{1pt}{1pt}{*}{2}{468}\above{1pt}{1pt}{l}{2}{1}\above{1pt}{1pt}{}{2}{9}\above{1pt}{1pt}{r}{2}{52}\above{1pt}{1pt}{*}{2}{144}\above{1pt}{1pt}{*}{2}{4}\above{1pt}{1pt}{s}{2}$}%
\nopagebreak\par%
\nopagebreak\par\leavevmode%
{$\left[\!\llap{\phantom{%
\begingroup \smaller\smaller\smaller
\endgroup%
}}\!\right]$}%
%
%
\hbox{}\par\smallskip%
%
%
\leavevmode%
${L_{350.2}}$%
{} : {$1\above{1pt}{1pt}{1}{1}4\above{1pt}{1pt}{1}{1}16\above{1pt}{1pt}{1}{1}{\cdot}1\above{1pt}{1pt}{1}{}3\above{1pt}{1pt}{-}{}9\above{1pt}{1pt}{1}{}{\cdot}1\above{1pt}{1pt}{2}{}13\above{1pt}{1pt}{1}{}$}\spacer%
\instructions{3}%
\EasyButWeakLineBreak%
{${52}\above{1pt}{1pt}{s}{2}{144}\above{1pt}{1pt}{s}{2}{4}\above{1pt}{1pt}{s}{2}{1872}\above{1pt}{1pt}{l}{2}{4}\above{1pt}{1pt}{}{2}{144}\above{1pt}{1pt}{}{2}{52}\above{1pt}{1pt}{r}{2}{36}\above{1pt}{1pt}{l}{2}{4}\above{1pt}{1pt}{}{2}{117}\above{1pt}{1pt}{}{2}{16}\above{1pt}{1pt}{}{2}{9}\above{1pt}{1pt}{}{2}{208}\above{1pt}{1pt}{}{2}{36}\above{1pt}{1pt}{r}{2}{16}\above{1pt}{1pt}{l}{2}{468}\above{1pt}{1pt}{}{2}{1}\above{1pt}{1pt}{}{2}{36}\above{1pt}{1pt}{r}{2}$}%
\nopagebreak\par%
shares genus with 2-dual${}\iso{}$3-dual; isometric to own %
2.3-dual\nopagebreak\par%
\nopagebreak\par\leavevmode%
{$\left[\!\llap{\phantom{%
\begingroup \smaller\smaller\smaller
\endgroup%
}}\!\right]$}%

\medskip%
%
\leavevmode\llap{}%
$W_{351}$%
\qquad\llap{4} lattices, $\chi=48$%
\hfill%
$222|222|222|222|\rtimes D_{4}$%
\nopagebreak\smallskip\hrule\nopagebreak\medskip%
%
%
\leavevmode%
${L_{351.1}}$%
{} : {$1\above{1pt}{1pt}{-}{5}8\above{1pt}{1pt}{1}{7}64\above{1pt}{1pt}{1}{1}{\cdot}1\above{1pt}{1pt}{-2}{}5\above{1pt}{1pt}{1}{}$}\EasyButWeakLineBreak%
{${8}\above{1pt}{1pt}{b}{2}{64}\above{1pt}{1pt}{s}{2}{20}\above{1pt}{1pt}{*}{2}{32}\above{1pt}{1pt}{l}{2}{5}\above{1pt}{1pt}{}{2}{64}\above{1pt}{1pt}{r}{2}$}\relax$\,(\times2)$%
\nopagebreak\par%
\nopagebreak\par\leavevmode%
{$\left[\!\llap{\phantom{%
\begingroup \smaller\smaller\smaller
\endgroup%
}}\!\right]$}%

\medskip%
%
\leavevmode\llap{}%
$W_{352}$%
\qquad\llap{12} lattices, $\chi=72$%
\hfill%
$2|2222|2222|2222|222\rtimes D_{4}$%
\nopagebreak\smallskip\hrule\nopagebreak\medskip%
%
%
\leavevmode%
${L_{352.1}}$%
{} : {$1\above{1pt}{1pt}{-2}{{\rm II}}4\above{1pt}{1pt}{1}{1}{\cdot}1\above{1pt}{1pt}{1}{}3\above{1pt}{1pt}{1}{}9\above{1pt}{1pt}{1}{}{\cdot}1\above{1pt}{1pt}{-2}{}5\above{1pt}{1pt}{1}{}{\cdot}1\above{1pt}{1pt}{2}{}19\above{1pt}{1pt}{1}{}$}\spacer%
\instructions{2}%
\EasyButWeakLineBreak%
{${684}\above{1pt}{1pt}{*}{2}{12}\above{1pt}{1pt}{*}{2}{76}\above{1pt}{1pt}{b}{2}{30}\above{1pt}{1pt}{l}{2}{4}\above{1pt}{1pt}{r}{2}{570}\above{1pt}{1pt}{l}{2}{36}\above{1pt}{1pt}{r}{2}{30}\above{1pt}{1pt}{b}{2}$}\relax$\,(\times2)$%
\nopagebreak\par%
\nopagebreak\par\leavevmode%
{$\left[\!\llap{\phantom{%
\begingroup \smaller\smaller\smaller
\endgroup%
}}\!\right]$}%

\medskip%
%
\leavevmode\llap{}%
$W_{353}$%
\qquad\llap{24} lattices, $\chi=36$%
\hfill%
$2222222222\rtimes C_{2}$%
\nopagebreak\smallskip\hrule\nopagebreak\medskip%
%
%
\leavevmode%
${L_{353.1}}$%
{} : {$1\above{1pt}{1pt}{-2}{{\rm II}}4\above{1pt}{1pt}{-}{5}{\cdot}1\above{1pt}{1pt}{2}{}3\above{1pt}{1pt}{1}{}{\cdot}1\above{1pt}{1pt}{1}{}5\above{1pt}{1pt}{-}{}25\above{1pt}{1pt}{-}{}{\cdot}1\above{1pt}{1pt}{2}{}7\above{1pt}{1pt}{-}{}$}\spacer%
\instructions{2}%
\EasyButWeakLineBreak%
{${84}\above{1pt}{1pt}{r}{2}{50}\above{1pt}{1pt}{b}{2}{140}\above{1pt}{1pt}{*}{2}{300}\above{1pt}{1pt}{b}{2}{10}\above{1pt}{1pt}{l}{2}$}\relax$\,(\times2)$%
\nopagebreak\par%
\nopagebreak\par\leavevmode%
{$\left[\!\llap{\phantom{%
\begingroup \smaller\smaller\smaller
\endgroup%
}}\!\right]$}%

\medskip%
%
\leavevmode\llap{}%
$W_{354}$%
\qquad\llap{12} lattices, $\chi=24$%
\hfill%
$222|2222|2\rtimes D_{2}$%
\nopagebreak\smallskip\hrule\nopagebreak\medskip%
%
%
\leavevmode%
${L_{354.1}}$%
{} : {$1\above{1pt}{1pt}{-2}{{\rm II}}4\above{1pt}{1pt}{-}{5}{\cdot}1\above{1pt}{1pt}{2}{}3\above{1pt}{1pt}{1}{}{\cdot}1\above{1pt}{1pt}{-}{}5\above{1pt}{1pt}{1}{}25\above{1pt}{1pt}{-}{}{\cdot}1\above{1pt}{1pt}{-2}{}7\above{1pt}{1pt}{1}{}$}\spacer%
\instructions{2}%
\EasyButWeakLineBreak%
{${700}\above{1pt}{1pt}{*}{2}{12}\above{1pt}{1pt}{b}{2}{50}\above{1pt}{1pt}{l}{2}{20}\above{1pt}{1pt}{r}{2}{2}\above{1pt}{1pt}{b}{2}{300}\above{1pt}{1pt}{*}{2}{28}\above{1pt}{1pt}{b}{2}{30}\above{1pt}{1pt}{b}{2}$}%
\nopagebreak\par%
\nopagebreak\par\leavevmode%
{$\left[\!\llap{\phantom{%
\begingroup \smaller\smaller\smaller\begin{tabular}{@{}c@{}}%
0\\0\\0
\end{tabular}\endgroup%
}}\right.$}%
\begingroup \smaller\smaller\smaller\begin{tabular}{@{}c@{}}%
-48300\\21000\\2100
\end{tabular}\endgroup%
\kern3pt%
\begingroup \smaller\smaller\smaller\begin{tabular}{@{}c@{}}%
21000\\-9130\\-905
\end{tabular}\endgroup%
\kern3pt%
\begingroup \smaller\smaller\smaller\begin{tabular}{@{}c@{}}%
2100\\-905\\58
\end{tabular}\endgroup%
{$\left.\llap{\phantom{%
\begingroup \smaller\smaller\smaller\begin{tabular}{@{}c@{}}%
0\\0\\0
\end{tabular}\endgroup%
}}\!\right]$}%
\EasyButWeakLineBreak%
{$\left[\!\llap{\phantom{%
\begingroup \smaller\smaller\smaller\begin{tabular}{@{}c@{}}%
0\\0\\0
\end{tabular}\endgroup%
}}\right.$}%
\begingroup \smaller\smaller\smaller\begin{tabular}{@{}c@{}}%
2241\\5180\\-350
\end{tabular}\endgroup%
\HardButStrongLineBreak\kern3pt%
\begingroup \smaller\smaller\smaller\begin{tabular}{@{}c@{}}%
257\\594\\-42
\end{tabular}\endgroup%
\HardButStrongLineBreak\kern3pt%
\begingroup \smaller\smaller\smaller\begin{tabular}{@{}c@{}}%
132\\305\\-25
\end{tabular}\endgroup%
\HardButStrongLineBreak\kern3pt%
\begingroup \smaller\smaller\smaller\begin{tabular}{@{}c@{}}%
-19\\-44\\0
\end{tabular}\endgroup%
\HardButStrongLineBreak\kern3pt%
\begingroup \smaller\smaller\smaller\begin{tabular}{@{}c@{}}%
-16\\-37\\2
\end{tabular}\endgroup%
\HardButStrongLineBreak\kern3pt%
\begingroup \smaller\smaller\smaller\begin{tabular}{@{}c@{}}%
13\\30\\0
\end{tabular}\endgroup%
\HardButStrongLineBreak\kern3pt%
\begingroup \smaller\smaller\smaller\begin{tabular}{@{}c@{}}%
109\\252\\-14
\end{tabular}\endgroup%
\HardButStrongLineBreak\kern3pt%
\begingroup \smaller\smaller\smaller\begin{tabular}{@{}c@{}}%
109\\252\\-15
\end{tabular}\endgroup%
{$\left.\llap{\phantom{%
\begingroup \smaller\smaller\smaller\begin{tabular}{@{}c@{}}%
0\\0\\0
\end{tabular}\endgroup%
}}\!\right]$}%

\medskip%
%
\leavevmode\llap{}%
$W_{355}$%
\qquad\llap{8} lattices, $\chi=120$%
\hfill%
$\infty\infty\infty\infty\infty22\infty\infty\infty\infty\infty22\rtimes C_{2}$%
\nopagebreak\smallskip\hrule\nopagebreak\medskip%
%
%
\leavevmode%
${L_{355.1}}$%
{} : {$1\above{1pt}{1pt}{-2}{{\rm II}}8\above{1pt}{1pt}{1}{7}{\cdot}1\above{1pt}{1pt}{1}{}11\above{1pt}{1pt}{-}{}121\above{1pt}{1pt}{-}{}$}\spacer%
\instructions{2}%
\EasyButWeakLineBreak%
{${88}\above{1pt}{1pt}{22,1}{\infty z}{22}\above{1pt}{1pt}{44,3}{\infty b}{88}\above{1pt}{1pt}{22,5}{\infty z}{22}\above{1pt}{1pt}{44,27}{\infty b}{88}\above{1pt}{1pt}{22,3}{\infty z}{22}\above{1pt}{1pt}{s}{2}{242}\above{1pt}{1pt}{b}{2}$}\relax$\,(\times2)$%
\nopagebreak\par%
\nopagebreak\par\leavevmode%
{$\left[\!\llap{\phantom{%
\begingroup \smaller\smaller\smaller
\endgroup%
}}\!\right]$}%

\medskip%
%
\leavevmode\llap{}%
$W_{356}$%
\qquad\llap{4} lattices, $\chi=32$%
\hfill%
$222\slashthree2222|2\rtimes D_{2}$%
\nopagebreak\smallskip\hrule\nopagebreak\medskip%
%
%
\leavevmode%
${L_{356.1}}$%
{} : {$1\above{1pt}{1pt}{-2}{{\rm II}}16\above{1pt}{1pt}{1}{1}{\cdot}1\above{1pt}{1pt}{1}{}3\above{1pt}{1pt}{-}{}9\above{1pt}{1pt}{1}{}{\cdot}1\above{1pt}{1pt}{-2}{}7\above{1pt}{1pt}{-}{}$}\EasyButWeakLineBreak%
{${144}\above{1pt}{1pt}{r}{2}{42}\above{1pt}{1pt}{b}{2}{16}\above{1pt}{1pt}{b}{2}{6}\above{1pt}{1pt}{-}{3}{6}\above{1pt}{1pt}{b}{2}{144}\above{1pt}{1pt}{b}{2}{42}\above{1pt}{1pt}{l}{2}{16}\above{1pt}{1pt}{r}{2}{6}\above{1pt}{1pt}{l}{2}$}%
\nopagebreak\par%
\nopagebreak\par\leavevmode%
{$\left[\!\llap{\phantom{%
\begingroup \smaller\smaller\smaller
\endgroup%
}}\!\right]$}%

\medskip%
%
\leavevmode\llap{}%
$W_{357}$%
\qquad\llap{4} lattices, $\chi=48$%
\hfill%
$22|222|222|222|2\rtimes D_{4}$%
\nopagebreak\smallskip\hrule\nopagebreak\medskip%
%
%
\leavevmode%
${L_{357.1}}$%
{} : {$1\above{1pt}{1pt}{1}{1}8\above{1pt}{1pt}{1}{7}128\above{1pt}{1pt}{-}{3}{\cdot}1\above{1pt}{1pt}{2}{}3\above{1pt}{1pt}{-}{}$}\EasyButWeakLineBreak%
{${384}\above{1pt}{1pt}{}{2}{1}\above{1pt}{1pt}{r}{2}{96}\above{1pt}{1pt}{*}{2}{4}\above{1pt}{1pt}{s}{2}{384}\above{1pt}{1pt}{b}{2}{8}\above{1pt}{1pt}{l}{2}$}\relax$\,(\times2)$%
\nopagebreak\par%
\nopagebreak\par\leavevmode%
{$\left[\!\llap{\phantom{%
\begingroup \smaller\smaller\smaller
\endgroup%
}}\!\right]$}%

\medskip%
%
\leavevmode\llap{}%
$W_{358}$%
\qquad\llap{12} lattices, $\chi=88$%
\hfill%
$262|2622|2262|2622|2\rtimes D_{4}$%
\nopagebreak\smallskip\hrule\nopagebreak\medskip%
%
%
\leavevmode%
${L_{358.1}}$%
{} : {$1\above{1pt}{1pt}{-2}{{\rm II}}4\above{1pt}{1pt}{-}{5}{\cdot}1\above{1pt}{1pt}{-}{}3\above{1pt}{1pt}{-}{}9\above{1pt}{1pt}{-}{}{\cdot}1\above{1pt}{1pt}{-2}{}5\above{1pt}{1pt}{-}{}{\cdot}1\above{1pt}{1pt}{2}{}23\above{1pt}{1pt}{1}{}$}\spacer%
\instructions{2}%
\EasyButWeakLineBreak%
{${828}\above{1pt}{1pt}{b}{2}{6}\above{1pt}{1pt}{}{6}{18}\above{1pt}{1pt}{l}{2}{276}\above{1pt}{1pt}{r}{2}{2}\above{1pt}{1pt}{}{6}{6}\above{1pt}{1pt}{b}{2}{92}\above{1pt}{1pt}{*}{2}{60}\above{1pt}{1pt}{*}{2}$}\relax$\,(\times2)$%
\nopagebreak\par%
\nopagebreak\par\leavevmode%
{$\left[\!\llap{\phantom{%
\begingroup \smaller\smaller\smaller
\endgroup%
}}\!\right]$}%

\medskip%
%
\leavevmode\llap{}%
$W_{359}$%
\qquad\llap{3} lattices, $\chi=24$%
\hfill%
$2|2|2|2|2|2|2|2|\rtimes D_{8}$%
\nopagebreak\smallskip\hrule\nopagebreak\medskip%
%
%
\leavevmode%
${L_{359.1}}$%
{} : {$1\above{1pt}{1pt}{-2}{{\rm II}}4\above{1pt}{1pt}{-}{3}{\cdot}1\above{1pt}{1pt}{1}{}3\above{1pt}{1pt}{1}{}9\above{1pt}{1pt}{1}{}{\cdot}1\above{1pt}{1pt}{-}{}5\above{1pt}{1pt}{-}{}25\above{1pt}{1pt}{-}{}$}\spacer%
\instructions{2}%
\EasyButWeakLineBreak%
{${300}\above{1pt}{1pt}{r}{2}{10}\above{1pt}{1pt}{l}{2}{12}\above{1pt}{1pt}{r}{2}{90}\above{1pt}{1pt}{l}{2}$}\relax$\,(\times2)$%
\nopagebreak\par%
\nopagebreak\par\leavevmode%
{$\left[\!\llap{\phantom{%
\begingroup \smaller\smaller\smaller\begin{tabular}{@{}c@{}}%
0\\0\\0
\end{tabular}\endgroup%
}}\right.$}%
\begingroup \smaller\smaller\smaller\begin{tabular}{@{}c@{}}%
146700\\-47700\\-1800
\end{tabular}\endgroup%
\kern3pt%
\begingroup \smaller\smaller\smaller\begin{tabular}{@{}c@{}}%
-47700\\15510\\585
\end{tabular}\endgroup%
\kern3pt%
\begingroup \smaller\smaller\smaller\begin{tabular}{@{}c@{}}%
-1800\\585\\22
\end{tabular}\endgroup%
{$\left.\llap{\phantom{%
\begingroup \smaller\smaller\smaller\begin{tabular}{@{}c@{}}%
0\\0\\0
\end{tabular}\endgroup%
}}\!\right]$}%
\hfil\penalty500%
{$\left[\!\llap{\phantom{%
\begingroup \smaller\smaller\smaller\begin{tabular}{@{}c@{}}%
0\\0\\0
\end{tabular}\endgroup%
}}\right.$}%
\begingroup \smaller\smaller\smaller\begin{tabular}{@{}c@{}}%
89\\180\\2700
\end{tabular}\endgroup%
\kern3pt%
\begingroup \smaller\smaller\smaller\begin{tabular}{@{}c@{}}%
-29\\-59\\-870
\end{tabular}\endgroup%
\kern3pt%
\begingroup \smaller\smaller\smaller\begin{tabular}{@{}c@{}}%
-1\\-2\\-31
\end{tabular}\endgroup%
{$\left.\llap{\phantom{%
\begingroup \smaller\smaller\smaller\begin{tabular}{@{}c@{}}%
0\\0\\0
\end{tabular}\endgroup%
}}\!\right]$}%
\EasyButWeakLineBreak%
{$\left[\!\llap{\phantom{%
\begingroup \smaller\smaller\smaller\begin{tabular}{@{}c@{}}%
0\\0\\0
\end{tabular}\endgroup%
}}\right.$}%
\begingroup \smaller\smaller\smaller\begin{tabular}{@{}c@{}}%
3\\20\\-300
\end{tabular}\endgroup%
\HardButStrongLineBreak\kern3pt%
\begingroup \smaller\smaller\smaller\begin{tabular}{@{}c@{}}%
1\\4\\-25
\end{tabular}\endgroup%
\HardButStrongLineBreak\kern3pt%
\begingroup \smaller\smaller\smaller\begin{tabular}{@{}c@{}}%
1\\4\\-24
\end{tabular}\endgroup%
\HardButStrongLineBreak\kern3pt%
\begingroup \smaller\smaller\smaller\begin{tabular}{@{}c@{}}%
-1\\-3\\0
\end{tabular}\endgroup%
{$\left.\llap{\phantom{%
\begingroup \smaller\smaller\smaller\begin{tabular}{@{}c@{}}%
0\\0\\0
\end{tabular}\endgroup%
}}\!\right]$}%
%
%
%
%
%
%
%
%
%
%

\medskip%
%
\leavevmode\llap{}%
$W_{360}$%
\qquad\llap{12} lattices, $\chi=54$%
\hfill%
$2\slashtwo222222|22222\rtimes D_{2}$%
\nopagebreak\smallskip\hrule\nopagebreak\medskip%
%
%
\leavevmode%
${L_{360.1}}$%
{} : {$1\above{1pt}{1pt}{-2}{{\rm II}}4\above{1pt}{1pt}{-}{3}{\cdot}1\above{1pt}{1pt}{1}{}3\above{1pt}{1pt}{-}{}9\above{1pt}{1pt}{1}{}{\cdot}1\above{1pt}{1pt}{2}{}7\above{1pt}{1pt}{-}{}{\cdot}1\above{1pt}{1pt}{2}{}19\above{1pt}{1pt}{1}{}$}\spacer%
\instructions{2}%
\EasyButWeakLineBreak%
{${42}\above{1pt}{1pt}{b}{2}{36}\above{1pt}{1pt}{*}{2}{4}\above{1pt}{1pt}{b}{2}{42}\above{1pt}{1pt}{l}{2}{76}\above{1pt}{1pt}{r}{2}{6}\above{1pt}{1pt}{b}{2}{532}\above{1pt}{1pt}{*}{2}{36}\above{1pt}{1pt}{b}{2}{114}\above{1pt}{1pt}{b}{2}{4}\above{1pt}{1pt}{*}{2}{4788}\above{1pt}{1pt}{b}{2}{6}\above{1pt}{1pt}{l}{2}{684}\above{1pt}{1pt}{r}{2}$}%
\nopagebreak\par%
\nopagebreak\par\leavevmode%
{$\left[\!\llap{\phantom{%
\begingroup \smaller\smaller\smaller
\endgroup%
}}\!\right]$}%

\medskip%
%
\leavevmode\llap{}%
$W_{361}$%
\qquad\llap{12} lattices, $\chi=112$%
\hfill%
$22226|62222|22226|62222|\rtimes D_{4}$%
\nopagebreak\smallskip\hrule\nopagebreak\medskip%
%
%
\leavevmode%
${L_{361.1}}$%
{} : {$1\above{1pt}{1pt}{-2}{{\rm II}}4\above{1pt}{1pt}{1}{7}{\cdot}1\above{1pt}{1pt}{-}{}3\above{1pt}{1pt}{-}{}9\above{1pt}{1pt}{-}{}{\cdot}1\above{1pt}{1pt}{-2}{}5\above{1pt}{1pt}{1}{}{\cdot}1\above{1pt}{1pt}{-2}{}29\above{1pt}{1pt}{1}{}$}\spacer%
\instructions{2}%
\EasyButWeakLineBreak%
{${870}\above{1pt}{1pt}{s}{2}{2}\above{1pt}{1pt}{b}{2}{1044}\above{1pt}{1pt}{*}{2}{20}\above{1pt}{1pt}{b}{2}{18}\above{1pt}{1pt}{}{6}{6}\above{1pt}{1pt}{}{6}{2}\above{1pt}{1pt}{b}{2}{180}\above{1pt}{1pt}{*}{2}{116}\above{1pt}{1pt}{b}{2}{18}\above{1pt}{1pt}{s}{2}$}\relax$\,(\times2)$%
\nopagebreak\par%
\nopagebreak\par\leavevmode%
{$\left[\!\llap{\phantom{%
\begingroup \smaller\smaller\smaller
\endgroup%
}}\!\right]$}%

\medskip%
%
\leavevmode\llap{}%
$W_{362}$%
\qquad\llap{4} lattices, $\chi=48$%
\hfill%
$\infty\infty|\infty\infty2\slashinfty2\rtimes D_{2}$%
\nopagebreak\smallskip\hrule\nopagebreak\medskip%
%
%
\leavevmode%
${L_{362.1}}$%
{} : {$1\above{1pt}{1pt}{1}{7}16\above{1pt}{1pt}{-}{5}256\above{1pt}{1pt}{-}{5}$}\EasyButWeakLineBreak%
{${64}\above{1pt}{1pt}{32,17}{\infty z}{16}\above{1pt}{1pt}{16,13}{\infty b}{64}\above{1pt}{1pt}{32,9}{\infty z}{16}\above{1pt}{1pt}{16,5}{\infty a}{64}\above{1pt}{1pt}{s}{2}{256}\above{1pt}{1pt}{16,1}{\infty z}{256}\above{1pt}{1pt}{*}{2}$}%
\nopagebreak\par%
\nopagebreak\par\leavevmode%
{$\left[\!\llap{\phantom{%
\begingroup \smaller\smaller\smaller\begin{tabular}{@{}c@{}}%
0\\0\\0
\end{tabular}\endgroup%
}}\right.$}%
\begingroup \smaller\smaller\smaller\begin{tabular}{@{}c@{}}%
-58112\\-3328\\256
\end{tabular}\endgroup%
\kern3pt%
\begingroup \smaller\smaller\smaller\begin{tabular}{@{}c@{}}%
-3328\\-176\\16
\end{tabular}\endgroup%
\kern3pt%
\begingroup \smaller\smaller\smaller\begin{tabular}{@{}c@{}}%
256\\16\\-1
\end{tabular}\endgroup%
{$\left.\llap{\phantom{%
\begingroup \smaller\smaller\smaller\begin{tabular}{@{}c@{}}%
0\\0\\0
\end{tabular}\endgroup%
}}\!\right]$}%
\EasyButWeakLineBreak%
{$\left[\!\llap{\phantom{%
\begingroup \smaller\smaller\smaller\begin{tabular}{@{}c@{}}%
0\\0\\0
\end{tabular}\endgroup%
}}\right.$}%
\begingroup \smaller\smaller\smaller\begin{tabular}{@{}c@{}}%
-27\\278\\-2720
\end{tabular}\endgroup%
\HardButStrongLineBreak\kern3pt%
\begingroup \smaller\smaller\smaller\begin{tabular}{@{}c@{}}%
-7\\73\\-696
\end{tabular}\endgroup%
\HardButStrongLineBreak\kern3pt%
\begingroup \smaller\smaller\smaller\begin{tabular}{@{}c@{}}%
-5\\54\\-480
\end{tabular}\endgroup%
\HardButStrongLineBreak\kern3pt%
\begingroup \smaller\smaller\smaller\begin{tabular}{@{}c@{}}%
0\\1\\8
\end{tabular}\endgroup%
\HardButStrongLineBreak\kern3pt%
\begingroup \smaller\smaller\smaller\begin{tabular}{@{}c@{}}%
1\\-10\\96
\end{tabular}\endgroup%
\HardButStrongLineBreak\kern3pt%
\begingroup \smaller\smaller\smaller\begin{tabular}{@{}c@{}}%
-1\\8\\-128
\end{tabular}\endgroup%
\HardButStrongLineBreak\kern3pt%
\begingroup \smaller\smaller\smaller\begin{tabular}{@{}c@{}}%
-15\\152\\-1536
\end{tabular}\endgroup%
{$\left.\llap{\phantom{%
\begingroup \smaller\smaller\smaller\begin{tabular}{@{}c@{}}%
0\\0\\0
\end{tabular}\endgroup%
}}\!\right]$}%
%
%
\hbox{}\par\smallskip%
%
%
\leavevmode%
${L_{362.2}}$%
{} : {$1\above{1pt}{1pt}{1}{1}16\above{1pt}{1pt}{-}{5}256\above{1pt}{1pt}{-}{3}$}\EasyButWeakLineBreak%
{${64}\above{1pt}{1pt}{32,31}{\infty z}{16}\above{1pt}{1pt}{16,11}{\infty a}{64}\above{1pt}{1pt}{32,23}{\infty z}{16}\above{1pt}{1pt}{16,3}{\infty b}{64}\above{1pt}{1pt}{*}{2}{4}\above{1pt}{1pt}{8,3}{\infty z}{1}\above{1pt}{1pt}{r}{2}$}%
\nopagebreak\par%
\nopagebreak\par\leavevmode%
{$\left[\!\llap{\phantom{%
\begingroup \smaller\smaller\smaller\begin{tabular}{@{}c@{}}%
0\\0\\0
\end{tabular}\endgroup%
}}\right.$}%
\begingroup \smaller\smaller\smaller\begin{tabular}{@{}c@{}}%
-1131776\\8448\\8448
\end{tabular}\endgroup%
\kern3pt%
\begingroup \smaller\smaller\smaller\begin{tabular}{@{}c@{}}%
8448\\-48\\-64
\end{tabular}\endgroup%
\kern3pt%
\begingroup \smaller\smaller\smaller\begin{tabular}{@{}c@{}}%
8448\\-64\\-63
\end{tabular}\endgroup%
{$\left.\llap{\phantom{%
\begingroup \smaller\smaller\smaller\begin{tabular}{@{}c@{}}%
0\\0\\0
\end{tabular}\endgroup%
}}\!\right]$}%
\EasyButWeakLineBreak%
{$\left[\!\llap{\phantom{%
\begingroup \smaller\smaller\smaller\begin{tabular}{@{}c@{}}%
0\\0\\0
\end{tabular}\endgroup%
}}\right.$}%
\begingroup \smaller\smaller\smaller\begin{tabular}{@{}c@{}}%
47\\314\\5952
\end{tabular}\endgroup%
\HardButStrongLineBreak\kern3pt%
\begingroup \smaller\smaller\smaller\begin{tabular}{@{}c@{}}%
16\\111\\2024
\end{tabular}\endgroup%
\HardButStrongLineBreak\kern3pt%
\begingroup \smaller\smaller\smaller\begin{tabular}{@{}c@{}}%
19\\138\\2400
\end{tabular}\endgroup%
\HardButStrongLineBreak\kern3pt%
\begingroup \smaller\smaller\smaller\begin{tabular}{@{}c@{}}%
4\\31\\504
\end{tabular}\endgroup%
\HardButStrongLineBreak\kern3pt%
\begingroup \smaller\smaller\smaller\begin{tabular}{@{}c@{}}%
-1\\-6\\-128
\end{tabular}\endgroup%
\HardButStrongLineBreak\kern3pt%
\begingroup \smaller\smaller\smaller\begin{tabular}{@{}c@{}}%
-1\\-8\\-126
\end{tabular}\endgroup%
\HardButStrongLineBreak\kern3pt%
\begingroup \smaller\smaller\smaller\begin{tabular}{@{}c@{}}%
1\\6\\127
\end{tabular}\endgroup%
{$\left.\llap{\phantom{%
\begingroup \smaller\smaller\smaller\begin{tabular}{@{}c@{}}%
0\\0\\0
\end{tabular}\endgroup%
}}\!\right]$}%
%
%
%
%
%
%
%
%
%
%

\medskip%
%
\leavevmode\llap{}%
$W_{363}$%
\qquad\llap{12} lattices, $\chi=72$%
\hfill%
$2|2222|2222|2222|222\rtimes D_{4}$%
\nopagebreak\smallskip\hrule\nopagebreak\medskip%
%
%
\leavevmode%
${L_{363.1}}$%
{} : {$1\above{1pt}{1pt}{-2}{{\rm II}}4\above{1pt}{1pt}{1}{1}{\cdot}1\above{1pt}{1pt}{2}{}3\above{1pt}{1pt}{-}{}{\cdot}1\above{1pt}{1pt}{1}{}5\above{1pt}{1pt}{-}{}25\above{1pt}{1pt}{1}{}{\cdot}1\above{1pt}{1pt}{-2}{}11\above{1pt}{1pt}{1}{}$}\spacer%
\instructions{2}%
\EasyButWeakLineBreak%
{${4}\above{1pt}{1pt}{r}{2}{10}\above{1pt}{1pt}{l}{2}{100}\above{1pt}{1pt}{r}{2}{6}\above{1pt}{1pt}{b}{2}{1100}\above{1pt}{1pt}{*}{2}{60}\above{1pt}{1pt}{*}{2}{44}\above{1pt}{1pt}{b}{2}{150}\above{1pt}{1pt}{l}{2}$}\relax$\,(\times2)$%
\nopagebreak\par%
\nopagebreak\par\leavevmode%
{$\left[\!\llap{\phantom{%
\begingroup \smaller\smaller\smaller
\endgroup%
}}\!\right]$}%

\medskip%
%
\leavevmode\llap{}%
$W_{364}$%
\qquad\llap{4} lattices, $\chi=32$%
\hfill%
$2622|2262|\rtimes D_{2}$%
\nopagebreak\smallskip\hrule\nopagebreak\medskip%
%
%
\leavevmode%
${L_{364.1}}$%
{} : {$1\above{1pt}{1pt}{-2}{{\rm II}}32\above{1pt}{1pt}{-}{3}{\cdot}1\above{1pt}{1pt}{-}{}3\above{1pt}{1pt}{-}{}9\above{1pt}{1pt}{-}{}{\cdot}1\above{1pt}{1pt}{-2}{}5\above{1pt}{1pt}{-}{}$}\EasyButWeakLineBreak%
{${96}\above{1pt}{1pt}{b}{2}{18}\above{1pt}{1pt}{}{6}{6}\above{1pt}{1pt}{b}{2}{2}\above{1pt}{1pt}{l}{2}{96}\above{1pt}{1pt}{r}{2}{18}\above{1pt}{1pt}{b}{2}{6}\above{1pt}{1pt}{}{6}{2}\above{1pt}{1pt}{b}{2}$}%
\nopagebreak\par%
\nopagebreak\par\leavevmode%
{$\left[\!\llap{\phantom{%
\begingroup \smaller\smaller\smaller\begin{tabular}{@{}c@{}}%
0\\0\\0
\end{tabular}\endgroup%
}}\right.$}%
\begingroup \smaller\smaller\smaller\begin{tabular}{@{}c@{}}%
3627360\\-1248480\\-642240
\end{tabular}\endgroup%
\kern3pt%
\begingroup \smaller\smaller\smaller\begin{tabular}{@{}c@{}}%
-1248480\\429018\\202803
\end{tabular}\endgroup%
\kern3pt%
\begingroup \smaller\smaller\smaller\begin{tabular}{@{}c@{}}%
-642240\\202803\\-369406
\end{tabular}\endgroup%
{$\left.\llap{\phantom{%
\begingroup \smaller\smaller\smaller\begin{tabular}{@{}c@{}}%
0\\0\\0
\end{tabular}\endgroup%
}}\!\right]$}%
\EasyButWeakLineBreak%
{$\left[\!\llap{\phantom{%
\begingroup \smaller\smaller\smaller\begin{tabular}{@{}c@{}}%
0\\0\\0
\end{tabular}\endgroup%
}}\right.$}%
\begingroup \smaller\smaller\smaller\begin{tabular}{@{}c@{}}%
-52331\\-155056\\5856
\end{tabular}\endgroup%
\HardButStrongLineBreak\kern3pt%
\begingroup \smaller\smaller\smaller\begin{tabular}{@{}c@{}}%
-20911\\-61959\\2340
\end{tabular}\endgroup%
\HardButStrongLineBreak\kern3pt%
\begingroup \smaller\smaller\smaller\begin{tabular}{@{}c@{}}%
24423\\72365\\-2733
\end{tabular}\endgroup%
\HardButStrongLineBreak\kern3pt%
\begingroup \smaller\smaller\smaller\begin{tabular}{@{}c@{}}%
37202\\110229\\-4163
\end{tabular}\endgroup%
\HardButStrongLineBreak\kern3pt%
\begingroup \smaller\smaller\smaller\begin{tabular}{@{}c@{}}%
645133\\1911520\\-72192
\end{tabular}\endgroup%
\HardButStrongLineBreak\kern3pt%
\begingroup \smaller\smaller\smaller\begin{tabular}{@{}c@{}}%
167369\\495912\\-18729
\end{tabular}\endgroup%
\HardButStrongLineBreak\kern3pt%
\begingroup \smaller\smaller\smaller\begin{tabular}{@{}c@{}}%
80186\\237590\\-8973
\end{tabular}\endgroup%
\HardButStrongLineBreak\kern3pt%
\begingroup \smaller\smaller\smaller\begin{tabular}{@{}c@{}}%
30205\\89497\\-3380
\end{tabular}\endgroup%
{$\left.\llap{\phantom{%
\begingroup \smaller\smaller\smaller\begin{tabular}{@{}c@{}}%
0\\0\\0
\end{tabular}\endgroup%
}}\!\right]$}%

\medskip%
%
\leavevmode\llap{}%
$W_{365}$%
\qquad\llap{12} lattices, $\chi=144$%
\hfill%
$2\slashtwo222|222\slashtwo222|222\slashtwo222|222\slashtwo222|22\rtimes D_{8}$%
\nopagebreak\smallskip\hrule\nopagebreak\medskip%
%
%
\leavevmode%
${L_{365.1}}$%
{} : {$1\above{1pt}{1pt}{2}{2}32\above{1pt}{1pt}{1}{1}{\cdot}1\above{1pt}{1pt}{-}{}3\above{1pt}{1pt}{-}{}9\above{1pt}{1pt}{-}{}{\cdot}1\above{1pt}{1pt}{-}{}5\above{1pt}{1pt}{1}{}25\above{1pt}{1pt}{-}{}$}\spacer%
\instructions{5,3}%
\EasyButWeakLineBreak%
{${800}\above{1pt}{1pt}{r}{2}{18}\above{1pt}{1pt}{b}{2}{50}\above{1pt}{1pt}{l}{2}{288}\above{1pt}{1pt}{}{2}{5}\above{1pt}{1pt}{r}{2}{7200}\above{1pt}{1pt}{s}{2}{20}\above{1pt}{1pt}{*}{2}{288}\above{1pt}{1pt}{b}{2}{50}\above{1pt}{1pt}{s}{2}{18}\above{1pt}{1pt}{b}{2}{800}\above{1pt}{1pt}{*}{2}{180}\above{1pt}{1pt}{s}{2}{32}\above{1pt}{1pt}{l}{2}{45}\above{1pt}{1pt}{}{2}$}\relax$\,(\times2)$%
\nopagebreak\par%
shares genus with 3-dual${}\iso{}$5-dual; isometric to own %
3.5-dual\nopagebreak\par%
\nopagebreak\par\leavevmode%
{$\left[\!\llap{\phantom{%
\begingroup \smaller\smaller\smaller
\endgroup%
}}\!\right]$}%

\medskip%
%
\leavevmode\llap{}%
$W_{366}$%
\qquad\llap{8} lattices, $\chi=84$%
\hfill%
$\infty\infty\infty222\infty\infty\infty222\rtimes C_{2}$%
\nopagebreak\smallskip\hrule\nopagebreak\medskip%
%
%
\leavevmode%
${L_{366.1}}$%
{} : {$1\above{1pt}{1pt}{-2}{{\rm II}}8\above{1pt}{1pt}{1}{1}{\cdot}1\above{1pt}{1pt}{-}{}13\above{1pt}{1pt}{-}{}169\above{1pt}{1pt}{-}{}$}\spacer%
\instructions{2*}%
\EasyButWeakLineBreak%
{${26}\above{1pt}{1pt}{52,25}{\infty a}{104}\above{1pt}{1pt}{26,9}{\infty z}{26}\above{1pt}{1pt}{52,17}{\infty a}{104}\above{1pt}{1pt}{b}{2}{338}\above{1pt}{1pt}{l}{2}{8}\above{1pt}{1pt}{r}{2}$}\relax$\,(\times2)$%
\nopagebreak\par%
shares genus with 13-dual\nopagebreak\par%
\nopagebreak\par\leavevmode%
{$\left[\!\llap{\phantom{%
\begingroup \smaller\smaller\smaller
\endgroup%
}}\!\right]$}%

\medskip%
%
\leavevmode\llap{}%
$W_{367}$%
\qquad\llap{12} lattices, $\chi=48$%
\hfill%
$222|222|222|222|\rtimes D_{4}$%
\nopagebreak\smallskip\hrule\nopagebreak\medskip%
%
%
\leavevmode%
${L_{367.1}}$%
{} : {$1\above{1pt}{1pt}{-2}{{\rm II}}4\above{1pt}{1pt}{-}{5}{\cdot}1\above{1pt}{1pt}{2}{}3\above{1pt}{1pt}{1}{}{\cdot}1\above{1pt}{1pt}{-2}{}5\above{1pt}{1pt}{-}{}{\cdot}1\above{1pt}{1pt}{-}{}7\above{1pt}{1pt}{-}{}49\above{1pt}{1pt}{-}{}$}\spacer%
\instructions{2}%
\EasyButWeakLineBreak%
{${84}\above{1pt}{1pt}{r}{2}{10}\above{1pt}{1pt}{b}{2}{588}\above{1pt}{1pt}{*}{2}{140}\above{1pt}{1pt}{*}{2}{12}\above{1pt}{1pt}{b}{2}{490}\above{1pt}{1pt}{l}{2}$}\relax$\,(\times2)$%
\nopagebreak\par%
\nopagebreak\par\leavevmode%
{$\left[\!\llap{\phantom{%
\begingroup \smaller\smaller\smaller
\endgroup%
}}\!\right]$}%

\medskip%
%
\leavevmode\llap{}%
$W_{368}$%
\qquad\llap{16} lattices, $\chi=36$%
\hfill%
$2222|22222|2\rtimes D_{2}$%
\nopagebreak\smallskip\hrule\nopagebreak\medskip%
%
%
\leavevmode%
${L_{368.1}}$%
{} : {$1\above{1pt}{1pt}{-2}{{\rm II}}8\above{1pt}{1pt}{-}{5}{\cdot}1\above{1pt}{1pt}{2}{}3\above{1pt}{1pt}{-}{}{\cdot}1\above{1pt}{1pt}{-}{}5\above{1pt}{1pt}{-}{}25\above{1pt}{1pt}{-}{}{\cdot}1\above{1pt}{1pt}{2}{}7\above{1pt}{1pt}{-}{}$}\spacer%
\instructions{2}%
\EasyButWeakLineBreak%
{${168}\above{1pt}{1pt}{r}{2}{50}\above{1pt}{1pt}{b}{2}{8}\above{1pt}{1pt}{b}{2}{1050}\above{1pt}{1pt}{l}{2}{40}\above{1pt}{1pt}{r}{2}{42}\above{1pt}{1pt}{b}{2}{200}\above{1pt}{1pt}{b}{2}{2}\above{1pt}{1pt}{l}{2}{4200}\above{1pt}{1pt}{r}{2}{10}\above{1pt}{1pt}{l}{2}$}%
\nopagebreak\par%
\nopagebreak\par\leavevmode%
{$\left[\!\llap{\phantom{%
\begingroup \smaller\smaller\smaller
\endgroup%
}}\!\right]$}%

\medskip%
%
\leavevmode\llap{}%
$W_{369}$%
\qquad\llap{8} lattices, $\chi=24$%
\hfill%
$22|2222|22\rtimes D_{2}$%
\nopagebreak\smallskip\hrule\nopagebreak\medskip%
%
%
\leavevmode%
${L_{369.1}}$%
{} : {$1\above{1pt}{1pt}{1}{1}8\above{1pt}{1pt}{1}{7}64\above{1pt}{1pt}{1}{7}{\cdot}1\above{1pt}{1pt}{2}{}3\above{1pt}{1pt}{-}{}{\cdot}1\above{1pt}{1pt}{-2}{}5\above{1pt}{1pt}{-}{}$}\EasyButWeakLineBreak%
{${960}\above{1pt}{1pt}{b}{2}{8}\above{1pt}{1pt}{s}{2}{60}\above{1pt}{1pt}{b}{2}{8}\above{1pt}{1pt}{l}{2}{960}\above{1pt}{1pt}{}{2}{1}\above{1pt}{1pt}{r}{2}{160}\above{1pt}{1pt}{*}{2}{4}\above{1pt}{1pt}{s}{2}$}%
\nopagebreak\par%
\nopagebreak\par\leavevmode%
{$\left[\!\llap{\phantom{%
\begingroup \smaller\smaller\smaller\begin{tabular}{@{}c@{}}%
0\\0\\0
\end{tabular}\endgroup%
}}\right.$}%
\begingroup \smaller\smaller\smaller\begin{tabular}{@{}c@{}}%
807360\\960\\-960
\end{tabular}\endgroup%
\kern3pt%
\begingroup \smaller\smaller\smaller\begin{tabular}{@{}c@{}}%
960\\-8\\0
\end{tabular}\endgroup%
\kern3pt%
\begingroup \smaller\smaller\smaller\begin{tabular}{@{}c@{}}%
-960\\0\\1
\end{tabular}\endgroup%
{$\left.\llap{\phantom{%
\begingroup \smaller\smaller\smaller\begin{tabular}{@{}c@{}}%
0\\0\\0
\end{tabular}\endgroup%
}}\!\right]$}%
\EasyButWeakLineBreak%
{$\left[\!\llap{\phantom{%
\begingroup \smaller\smaller\smaller\begin{tabular}{@{}c@{}}%
0\\0\\0
\end{tabular}\endgroup%
}}\right.$}%
\begingroup \smaller\smaller\smaller\begin{tabular}{@{}c@{}}%
-29\\-3120\\-27360
\end{tabular}\endgroup%
\HardButStrongLineBreak\kern3pt%
\begingroup \smaller\smaller\smaller\begin{tabular}{@{}c@{}}%
-1\\-107\\-940
\end{tabular}\endgroup%
\HardButStrongLineBreak\kern3pt%
\begingroup \smaller\smaller\smaller\begin{tabular}{@{}c@{}}%
-1\\-105\\-930
\end{tabular}\endgroup%
\HardButStrongLineBreak\kern3pt%
\begingroup \smaller\smaller\smaller\begin{tabular}{@{}c@{}}%
0\\1\\4
\end{tabular}\endgroup%
\HardButStrongLineBreak\kern3pt%
\begingroup \smaller\smaller\smaller\begin{tabular}{@{}c@{}}%
1\\120\\960
\end{tabular}\endgroup%
\HardButStrongLineBreak\kern3pt%
\begingroup \smaller\smaller\smaller\begin{tabular}{@{}c@{}}%
0\\0\\-1
\end{tabular}\endgroup%
\HardButStrongLineBreak\kern3pt%
\begingroup \smaller\smaller\smaller\begin{tabular}{@{}c@{}}%
-1\\-110\\-960
\end{tabular}\endgroup%
\HardButStrongLineBreak\kern3pt%
\begingroup \smaller\smaller\smaller\begin{tabular}{@{}c@{}}%
-1\\-108\\-946
\end{tabular}\endgroup%
{$\left.\llap{\phantom{%
\begingroup \smaller\smaller\smaller\begin{tabular}{@{}c@{}}%
0\\0\\0
\end{tabular}\endgroup%
}}\!\right]$}%

\medskip%
%
\leavevmode\llap{}%
$W_{370}$%
\qquad\llap{4} lattices, $\chi=24$%
\hfill%
$2|22|22|22|2\rtimes D_{4}$%
\nopagebreak\smallskip\hrule\nopagebreak\medskip%
%
%
\leavevmode%
${L_{370.1}}$%
{} : {$1\above{1pt}{1pt}{1}{7}8\above{1pt}{1pt}{1}{1}64\above{1pt}{1pt}{1}{7}{\cdot}1\above{1pt}{1pt}{2}{}3\above{1pt}{1pt}{-}{}{\cdot}1\above{1pt}{1pt}{-2}{}5\above{1pt}{1pt}{-}{}$}\EasyButWeakLineBreak%
{${960}\above{1pt}{1pt}{l}{2}{8}\above{1pt}{1pt}{}{2}{15}\above{1pt}{1pt}{r}{2}{64}\above{1pt}{1pt}{*}{2}{60}\above{1pt}{1pt}{l}{2}{8}\above{1pt}{1pt}{}{2}{960}\above{1pt}{1pt}{r}{2}{4}\above{1pt}{1pt}{b}{2}$}%
\nopagebreak\par%
\nopagebreak\par\leavevmode%
{$\left[\!\llap{\phantom{%
\begingroup \smaller\smaller\smaller\begin{tabular}{@{}c@{}}%
0\\0\\0
\end{tabular}\endgroup%
}}\right.$}%
\begingroup \smaller\smaller\smaller\begin{tabular}{@{}c@{}}%
960\\0\\0
\end{tabular}\endgroup%
\kern3pt%
\begingroup \smaller\smaller\smaller\begin{tabular}{@{}c@{}}%
0\\8\\0
\end{tabular}\endgroup%
\kern3pt%
\begingroup \smaller\smaller\smaller\begin{tabular}{@{}c@{}}%
0\\0\\-1
\end{tabular}\endgroup%
{$\left.\llap{\phantom{%
\begingroup \smaller\smaller\smaller\begin{tabular}{@{}c@{}}%
0\\0\\0
\end{tabular}\endgroup%
}}\!\right]$}%
\EasyButWeakLineBreak%
{$\left[\!\llap{\phantom{%
\begingroup \smaller\smaller\smaller\begin{tabular}{@{}c@{}}%
0\\0\\0
\end{tabular}\endgroup%
}}\right.$}%
\begingroup \smaller\smaller\smaller\begin{tabular}{@{}c@{}}%
-11\\-120\\-480
\end{tabular}\endgroup%
\HardButStrongLineBreak\kern3pt%
\begingroup \smaller\smaller\smaller\begin{tabular}{@{}c@{}}%
-1\\-9\\-40
\end{tabular}\endgroup%
\HardButStrongLineBreak\kern3pt%
\begingroup \smaller\smaller\smaller\begin{tabular}{@{}c@{}}%
-2\\-15\\-75
\end{tabular}\endgroup%
\HardButStrongLineBreak\kern3pt%
\begingroup \smaller\smaller\smaller\begin{tabular}{@{}c@{}}%
-1\\-4\\-32
\end{tabular}\endgroup%
\HardButStrongLineBreak\kern3pt%
\begingroup \smaller\smaller\smaller\begin{tabular}{@{}c@{}}%
-1\\0\\-30
\end{tabular}\endgroup%
\HardButStrongLineBreak\kern3pt%
\begingroup \smaller\smaller\smaller\begin{tabular}{@{}c@{}}%
0\\1\\0
\end{tabular}\endgroup%
\HardButStrongLineBreak\kern3pt%
\begingroup \smaller\smaller\smaller\begin{tabular}{@{}c@{}}%
1\\0\\0
\end{tabular}\endgroup%
\HardButStrongLineBreak\kern3pt%
\begingroup \smaller\smaller\smaller\begin{tabular}{@{}c@{}}%
0\\-1\\-2
\end{tabular}\endgroup%
{$\left.\llap{\phantom{%
\begingroup \smaller\smaller\smaller\begin{tabular}{@{}c@{}}%
0\\0\\0
\end{tabular}\endgroup%
}}\!\right]$}%
%
%
%
%
%
%
%
%
%
%
%
%
%
%

\medskip%
%
\leavevmode\llap{}%
$W_{371}$%
\qquad\llap{4} lattices, $\chi=144$%
\hfill%
$\infty|\infty222\infty|\infty222\infty|\infty222\infty|\infty222\rtimes D_{4}$%
\nopagebreak\smallskip\hrule\nopagebreak\medskip%
%
%
\leavevmode%
${L_{371.1}}$%
{} : {$1\above{1pt}{1pt}{-2}{{\rm II}}8\above{1pt}{1pt}{-}{5}{\cdot}1\above{1pt}{1pt}{1}{}17\above{1pt}{1pt}{1}{}289\above{1pt}{1pt}{1}{}$}\spacer%
\instructions{2}%
\EasyButWeakLineBreak%
{${136}\above{1pt}{1pt}{34,1}{\infty z}{34}\above{1pt}{1pt}{68,13}{\infty a}{136}\above{1pt}{1pt}{b}{2}{2}\above{1pt}{1pt}{b}{2}{2312}\above{1pt}{1pt}{b}{2}{34}\above{1pt}{1pt}{68,1}{\infty b}{136}\above{1pt}{1pt}{34,21}{\infty z}{34}\above{1pt}{1pt}{b}{2}{8}\above{1pt}{1pt}{b}{2}{578}\above{1pt}{1pt}{b}{2}$}\relax$\,(\times2)$%
\nopagebreak\par%
\nopagebreak\par\leavevmode%
{$\left[\!\llap{\phantom{%
\begingroup \smaller\smaller\smaller
\endgroup%
}}\!\right]$}%

\medskip%
%
\leavevmode\llap{}%
$W_{372}$%
\qquad\llap{8} lattices, $\chi=72$%
\hfill%
$22|2222|2222|2222|22\rtimes D_{4}$%
\nopagebreak\smallskip\hrule\nopagebreak\medskip%
%
%
\leavevmode%
${L_{372.1}}$%
{} : {$1\above{1pt}{1pt}{1}{1}8\above{1pt}{1pt}{-}{5}64\above{1pt}{1pt}{1}{7}{\cdot}1\above{1pt}{1pt}{2}{}3\above{1pt}{1pt}{1}{}{\cdot}1\above{1pt}{1pt}{2}{}7\above{1pt}{1pt}{1}{}$}\EasyButWeakLineBreak%
{${448}\above{1pt}{1pt}{r}{2}{12}\above{1pt}{1pt}{b}{2}{56}\above{1pt}{1pt}{s}{2}{12}\above{1pt}{1pt}{b}{2}{448}\above{1pt}{1pt}{s}{2}{4}\above{1pt}{1pt}{*}{2}{64}\above{1pt}{1pt}{l}{2}{1}\above{1pt}{1pt}{}{2}$}\relax$\,(\times2)$%
\nopagebreak\par%
\nopagebreak\par\leavevmode%
{$\left[\!\llap{\phantom{%
\begingroup \smaller\smaller\smaller
\endgroup%
}}\!\right]$}%

\medskip%
%
\leavevmode\llap{}%
$W_{373}$%
\qquad\llap{8} lattices, $\chi=72$%
\hfill%
$2222|2222|2222|2222|\rtimes D_{4}$%
\nopagebreak\smallskip\hrule\nopagebreak\medskip%
%
%
\leavevmode%
${L_{373.1}}$%
{} : {$1\above{1pt}{1pt}{-2}{2}16\above{1pt}{1pt}{-}{3}{\cdot}1\above{1pt}{1pt}{-}{}3\above{1pt}{1pt}{1}{}9\above{1pt}{1pt}{-}{}{\cdot}1\above{1pt}{1pt}{-2}{}5\above{1pt}{1pt}{-}{}{\cdot}1\above{1pt}{1pt}{2}{}7\above{1pt}{1pt}{-}{}$}\EasyButWeakLineBreak%
{${18}\above{1pt}{1pt}{b}{2}{560}\above{1pt}{1pt}{*}{2}{12}\above{1pt}{1pt}{l}{2}{315}\above{1pt}{1pt}{}{2}{48}\above{1pt}{1pt}{}{2}{35}\above{1pt}{1pt}{r}{2}{12}\above{1pt}{1pt}{*}{2}{5040}\above{1pt}{1pt}{b}{2}{2}\above{1pt}{1pt}{l}{2}{5040}\above{1pt}{1pt}{}{2}{3}\above{1pt}{1pt}{r}{2}{140}\above{1pt}{1pt}{*}{2}{48}\above{1pt}{1pt}{*}{2}{1260}\above{1pt}{1pt}{l}{2}{3}\above{1pt}{1pt}{}{2}{560}\above{1pt}{1pt}{r}{2}$}%
\nopagebreak\par%
\nopagebreak\par\leavevmode%
{$\left[\!\llap{\phantom{%
\begingroup \smaller\smaller\smaller
\endgroup%
}}\!\right]$}%

\medskip%
%
\leavevmode\llap{}%
$W_{374}$%
\qquad\llap{4} lattices, $\chi=72$%
\hfill%
$222222|2222222\slashinfty2\rtimes D_{2}$%
\nopagebreak\smallskip\hrule\nopagebreak\medskip%
%
%
\leavevmode%
${L_{374.1}}$%
{} : {$1\above{1pt}{1pt}{1}{1}8\above{1pt}{1pt}{1}{7}256\above{1pt}{1pt}{1}{1}{\cdot}1\above{1pt}{1pt}{2}{}9\above{1pt}{1pt}{1}{}$}\EasyButWeakLineBreak%
{${2304}\above{1pt}{1pt}{}{2}{1}\above{1pt}{1pt}{r}{2}{256}\above{1pt}{1pt}{*}{2}{36}\above{1pt}{1pt}{s}{2}{256}\above{1pt}{1pt}{s}{2}{4}\above{1pt}{1pt}{*}{2}{2304}\above{1pt}{1pt}{l}{2}{1}\above{1pt}{1pt}{}{2}{256}\above{1pt}{1pt}{}{2}{9}\above{1pt}{1pt}{r}{2}{256}\above{1pt}{1pt}{*}{2}{4}\above{1pt}{1pt}{s}{2}{2304}\above{1pt}{1pt}{b}{2}{8}\above{1pt}{1pt}{48,1}{\infty b}{8}\above{1pt}{1pt}{l}{2}$}%
\nopagebreak\par%
\nopagebreak\par\leavevmode%
{$\left[\!\llap{\phantom{%
\begingroup \smaller\smaller\smaller
\endgroup%
}}\!\right]$}%
\endgroup

\end{document}